\documentclass{cmim-l}
\usepackage{amsmath, amsbsy, amsthm, amsfonts, array,rlepsf}
 \newtheorem{thm}{Theorem}[chapter]
 \newtheorem{cor}[thm]{Corollary}
 \newtheorem{lem}[thm]{Lemma}
 \newtheorem{prop}[thm]{Proposition}
 \newtheorem{claim}[thm]{Claim}
 \theoremstyle{definition}
 \newtheorem{defn}[thm]{Definition}
 \newtheorem{exam}[thm]{Example}
 
 \newtheorem{ex}[thm]{Exercise}
 \newtheorem{rem}[thm]{Remark}
 \newtheorem{assumption}[thm]{Assumption}
 
 \numberwithin{equation}{chapter}

 \newcommand{\A}{\mathcal{A}}

 \newcommand{\abs}[1]{\left\vert#1\right\vert}

 \newcommand{\ii}[1]{{\it #1}}

\def\RR{{\mathbb R}}

\def\bbeta{{\mathbb \beta}}

\def\Zee{\mathbb{Z}}

\def\Ar{\mathbb{R}}
\def\Cee{\mathbb{C}}

\def\RR{{\mathbb R}}

\def\be{\begin{equation}}
\def\ee{\end{equation}}
\def\ba{\begin{eqnarray}}
\def\ea{\end{eqnarray}}
\def\tilde{\widetilde}

\def\e1{\epsilon}
\def\AAl{\mathcal{A}_{\lambda}}
\def\A0{\stackrel{\circ}{\AAl}}

\def\o1{\omega}
\def\01{\Omega}
\def\c1{\gamma}

\def\g1{\Sigma}

\def\bigcup{\cup}

\def\l1{\Lambda}

\def\v1{\varphi}

\def\d1{\delta}
\def\part{\partial}

\def\f1{\frac}
\def\t1{\theta}
\def\b1{\beta}
\def\bar{\overline}
\def\bs{\begin{eqnarray*}}
\def\es{\end{eqnarray*}}
\def\m1{\Theta}
\def\w1{\wedge}

\makeindex

\begin{document}

\title{\Huge{\textbf{Ricci Flow and the Poincar\'e Conjecture}}}

\author
{\huge{John W. Morgan and Gang Tian}\\}
 \maketitle
 \thispagestyle{empty}

\titlepage

\tableofcontents \pagestyle{headings}
\chapter*{Introduction} 
 In this book we present a complete
and detailed proof of the

\begin{quote}{\bf Poincar\'e Conjecture:\index{Poincar\'e Conjecture|ii}
a closed, smooth, simply connected $3$-manifold is
diffeomorphic}\footnote{Every topological $3$-manifold admits a
differentiable structure and every homeomorphism between smooth
$3$-manifolds can be approximated by a diffeomorphism. Thus,
classification results about topological $3$-manifolds up to
homeomorphism and about smooth $3$-manifolds up to diffeomorphism
are equivalent. In this book `manifold' means `smooth manifold.'}
{\bf to $S^3$.}
\end{quote}

\noindent 
This conjecture was formulated by Henri Poincar\'e \cite{Poincare}
in 1904 and has remained open until the recent work of Perelman. The
arguments we give here are a detailed version of those that appear
in Perelman's three preprints~\cite{Perelman1,Perelman2,Perelman3}.
Perelman's arguments rest on a foundation built by Richard Hamilton
with his study of the Ricci flow equation for Riemannian metrics.
Indeed, Hamilton believed that Ricci flow could be used to establish
the Poincar\'e Conjecture and more general topological
classification results in dimension $3$, and laid out a program to
accomplish this. The difficulty was to deal with singularities in
the Ricci flow. Perelman's breakthrough was to understand the
qualitative nature of the singularities sufficiently to allow him to
prove the Poincar\'e Conjecture (and Theorem~\ref{Theorem1} below
which implies the Poincar\'e Conjecture).   For
a detailed history of the Poincar\'e Conjecture, see Milnor's survey
article \cite{Milnorsurvey}.

\bigskip
A  class of examples closely related to the $3$-sphere are {\em the
$3$-dimensional spherical space-forms}\index{space-form!$3$-dimensional}, i.e.,
the quotients of $S^3$ by free, linear actions of finite subgroups of the
orthogonal group $O(4)$. There is a generalization of the Poincar\'e
Conjecture, called the {\bf $3$-dimensional spherical space-form
conjecture}\index{space-form!$3$-dimensional!conjecture}, which conjectures
that any closed $3$-manifold with finite fundamental group is diffeomorphic to
a $3$-dimensional spherical space-form. Clearly, a special case of the
$3$-dimensional spherical space-form conjecture is the Poincar\'e Conjecture.

 As indicated in Remark 1.4 of \cite{Perelman3},  the arguments we present
 here not only prove the Poincar\'e Conjecture, they prove the $3$-dimensional space-form conjecture.
 In fact,  the purpose of this book is to prove the following
 more general theorem.

\begin{thm}\label{Theorem1} Let $M$ be a closed, connected $3$-manifold and
suppose that the fundamental group of $M$ is a free product of
finite groups and infinite cyclic groups. Then $M$ is diffeomorphic
to a connected sum of spherical space-forms, copies of $S^2\times
S^1$, and copies of the unique (up to diffeomorphism) non-orientable
$2$-sphere bundle over $S^1$.
\end{thm}

This immediately implies an affirmative resolution of the Poincar\'e
Conjecture and of the $3$-dimensional spherical space-form
conjecture.

\begin{cor}
(a) A closed, simply connected $3$-manifold is diffeomorphic to
$S^3$. (b) A closed  $3$-manifold with finite fundamental
group is diffeomorphic to a $3$-dimensional spherical space-form.
\end{cor}

Before launching into a more detailed description of the contents of
this book, one remark on the style of the exposition is in order.
Because of the importance and visibility of the results discussed
here, and because of the number of incorrect claims of proofs of
these results in the past, we felt that it behooved us to work out
and present the arguments in great detail. Our goal was to make the
arguments clear and convincing and also to make them more easily
accessible to a wider audience. As a result, experts may find some
of the points are overly elaborated.

\section{Overview of Perelman's argument}

In dimensions less than or equal to three, any Riemannian metric of constant
Ricci curvature has constant sectional curvature. Classical results in
Riemannian geometry show that the universal cover of a closed manifold of
constant positive curvature is diffeomorphic to the sphere and that the
fundamental group is identified with a finite subgroup of the orthogonal group
acting linearly and freely on the universal cover. Thus, one can approach the
Poincar\'e Conjecture and the more general $3$-dimensional spherical space-form
problem by asking the following question. Making the appropriate fundamental
group assumptions on $3$-manifold $M$, how does one establish the existence of
a metric of constant Ricci curvature on $M$? The essential ingredient in
producing such a metric is the Ricci flow equation\index{Ricci flow!equation}
introduced by Richard Hamilton in \cite{Hamilton3MPRC}:
$$\frac{\partial g(t)}{\partial t}=-2{\rm Ric}(g(t)),$$
where ${\rm Ric}(g(t))$ is the Ricci curvature of the metric $g(t)$.
The fixed points (up to rescaling) of this equation are the
Riemannian metrics of constant Ricci curvature. For a general
introduction to the subject of the Ricci flow see Hamilton's survey
paper \cite{Hamiltonsurvey}, the book by Chow-Knopf
\cite{ChowKnopf}, or the book by Chow, Lu, and Ni \cite{ChowLuNi}.
The Ricci flow equation is a (weakly) parabolic partial differential
flow equation for Riemannian metrics on a smooth manifold. Following
Hamilton, one defines a Ricci flow \index{Ricci flow}to be a family
of Riemannian metrics $g(t)$ on a fixed smooth manifold,
parameterized by $t$ in some interval, satisfying this equation. One
considers $t$ as time and studies the equation as an initial value
problem: Beginning with any Riemannian manifold $(M,g_0)$ find a
Ricci flow with $(M,g_0)$ as initial metric; that is to say find a
one-parameter family $(M,g(t))$ of Riemannian manifolds with
$g(0)=g_0$ satisfying the Ricci flow equation. This equation is
valid in all dimensions but we concentrate here on dimension three.
In a sentence, the method of proof is to begin with any Riemannian
metric on the given smooth $3$-manifold and flow it using the Ricci
flow equation to obtain the constant curvature metric for which one
is searching. There are two examples where things work in exactly
this way, both due to Hamilton. (i) If the initial metric has
positive Ricci curvature, Hamilton proved over twenty years ago,
\cite{Hamilton3MPRC}, that under the Ricci flow the manifold shrinks
to a point in finite time, that is to say, there is a finite-time
singularity, and, as we approach the singular time, the diameter of
the manifold tends to zero and the curvature blows up at every
point. Hamilton went on to show that, in this case, rescaling by a
time-dependent function so that the diameter is constant produces a
one-parameter family of metrics converging smoothly to a metric of
constant positive curvature. (ii) At the other extreme, in
\cite{HamiltonNSRF3M} Hamilton showed that if the Ricci flow exists
for all time and if there is an appropriate curvature bound together
with another geometric bound, then as $t\rightarrow\infty$, after
rescaling to have a fixed diameter, the metric converges to a metric
of constant negative curvature.

The results in the general case are much more complicated to
formulate and much more difficult to establish. While Hamilton
established that the Ricci flow equation  has short-term existence
properties, i.e., one can  define $g(t)$ for $t$ in some interval
$[0,T)$ where $T$ depends on the initial metric, it turns out that
if the topology of the manifold is sufficiently complicated, say it
is a non-trivial connected sum, then no matter what the initial
metric is one must encounter finite-time singularities, forced by
the topology. More seriously, even if the manifold has simple
topology, beginning with an arbitrary metric one expects to (and
cannot rule out the possibility that one will) encounter finite-time
singularities in the Ricci flow. These singularities, unlike in the
case of positive Ricci curvature, occur along proper subsets of the
manifold, not the entire manifold.
 Thus, to derive the topological
consequences stated above, it is not sufficient in general to stop
the analysis the first time a singularity arises in the Ricci flow.
One is led to study a more general evolution process
 called {\em Ricci flow with surgery},\index{Ricci flow!with surgery}
 first introduced by Hamilton in the context
 of four-manifolds, \cite{Hamilton4MPIC}.
 This evolution process
 is still parameterized by an interval in time, so that for each $t$
 in the interval of definition there is a compact Riemannian $3$-manifold
 $M_t$. But there is a discrete set of
times at which the manifolds and metrics undergo topological and
metric discontinuities (surgeries). In each of the complementary
intervals to the singular times, the evolution is the usual Ricci
flow, though, because of the surgeries, the topological type of the
manifold $M_t$ changes as $t$ moves from one complementary interval
to the next. From an analytic point of view, the surgeries at the
discontinuity times are introduced in order to `cut away' a
neighborhood of the singularities as they develop and insert by
hand, in place of the `cut away' regions,  geometrically nice
regions. This allows one to continue  the Ricci flow (or more
precisely, restart the Ricci flow with the new metric constructed at
the discontinuity time). Of course, the surgery process also changes
the topology. To be able to say anything useful topologically about
such a process, one needs results about Ricci flow, and one also
needs to control both the topology and the geometry of the surgery
process at the singular times. For example, it is crucial for the
topological applications that we do surgery along $2$-spheres rather
than surfaces of higher genus. Surgery along $2$-spheres produces
the connected sum decomposition, which is well-understood
topologically, while, for example, Dehn surgeries along tori can
completely destroy the topology, changing any $3$-manifold into any
other.

The change in topology turns out to be completely understandable and
amazingly, the surgery processes produce exactly the topological
operations needed to cut the manifold into pieces on which the Ricci
flow can produce the metrics sufficiently controlled so that the
topology can be recognized.

The bulk of this book (Chapters 1-17 and the Appendix) concerns the
establishment of the following long-time existence result for Ricci
flow with surgery\index{Ricci flow!with surgery!long-time
existence}.

\begin{thm}\label{surgery} Let $(M,g_0)$ be a closed Riemannian $3$-manifold.
Suppose that there is no embedded, locally separating $\Ar P^2$
contained\footnote{I.e., no embedded $\Ar P^2$ in $M$ with trivial
normal bundle. Clearly, all orientable manifolds satisfy this
condition.} in $M$. Then there is a Ricci flow with surgery defined
for all $t\in [0,\infty)$ with initial metric $(M,g_0)$. The set of
discontinuity times for this Ricci flow with surgery is a discrete
subset of $[0,\infty)$. The topological change in the $3$-manifold
as one crosses a surgery time is a connected sum decomposition
together with removal of connected components, each of which is
diffeomorphic to one of $S^2\times S^1$, $\Ar P^3\#\Ar P^3$, the
non-orientable $2$-sphere bundle over $S^1$, or a manifold admitting
a metric of constant positive curvature.
\end{thm}

While Theorem~\ref{surgery} is central for all applications of Ricci
flow to the topology of three-dimensional manifolds, the argument
for the $3$-manifolds described in Theorem~\ref{Theorem1} is
simplified, and avoids all references to the nature of the flow  as
time goes to infinity, because of the following finite-time
extinction result\index{Ricci flow!with surgery!finite-time
extinction}.

\begin{thm}\label{finiteext}
Let $M$ be a closed $3$-manifold whose fundamental group is a free
product of finite groups and infinite cyclic groups\footnote{In
\cite{Perelman3} Perelman states the result for $3$-manifolds
without prime factors that are acyclic. It is a standard exercise in
$3$-manifold topology to show that Perelman's condition is
equivalent to the group theory hypothesis stated here; see
Corollary~\ref{corequiv}.}. Let $g_0$ be any Riemannian metric on
$M$. Then $M$ admits no locally separating $\Ar P^2$, so that there
is a Ricci flow with surgery defined for all positive time with
$(M,g_0)$ as initial metric as described in Theorem~\ref{surgery}.
This Ricci flow with surgery becomes extinct after some time
$T<\infty$, in the sense that the manifolds $M_t$ are empty for all
$t\ge T$.
\end{thm}

This result is established in Chapter 18 following the argument
given by Perelman in \cite{Perelman3}, see also
\cite{ColdingMinicozzi}.

We immediately deduce Theorem~\ref{Theorem1} from
Theorems~\ref{surgery} and~\ref{finiteext} as  follows: Let $M$ be a
$3$-manifold satisfying the hypothesis of Theorem~\ref{Theorem1}.
Then there is a finite sequence $M=M_0,M_1,\ldots,M_k=\emptyset$
such that for each $i,\ 1\le i\le k$, $M_i$ is obtained from
$M_{i-1}$ by  a connected sum decomposition or $M_i$ is obtained
from $M_{i-1}$ by removing  a component diffeomorphic to one of
$S^2\times S^1$, $\Ar P^3\#\Ar P^3$, a non-orientable $2$-sphere
bundle over $S^1$, or a $3$-dimensional spherical space-form.
Clearly, it follows by downward induction on $i$ that  each
connected component of $M_i$ is diffeomorphic to a connected sum of
$3$-dimensional spherical space-forms, copies of $S^2\times S^1$,
and copies of the non-orientable $2$-sphere bundle over $S^1$. In
particular, $M=M_0$ has this form. Since $M$ is connected by
hypothesis, this proves the theorem. In fact, this argument proves
the following:

\begin{cor}\label{corequiv}
Let $(M_0,g_0)$ a connected Riemannian manifold with no locally
separating $\Ar P^2$. Let $({\mathcal M},G)$ be a Ricci flow with
surgery defined for $0\le t<\infty$ with $(M_0,g_0)$ as initial
manifold. Then the following four conditions are equivalent:
\begin{enumerate} \item[(1)] $({\mathcal M},G)$ becomes extinct after a
finite time, i.e., $M_T=\emptyset$ for all $T$ sufficiently large,
\item[(2)] $M_0$ is diffeomorphic to a connected sum of three-dimensional
spherical space-forms and $S^2$-bundles over $S^1$,\item[(3)] the
fundamental group of $M_0$ is a free product of finite groups and
infinite cyclic groups,
\item[(4)] no prime\footnote{A three-manifold $P$ is prime
if  every separating two-sphere in $P$ bounds a three-ball in $P$.
Equivalently, $P$ is prime if it admits no non-trivial connected sum
decomposition. Every closed three-manifold decomposes as a connected
sum of prime factors with the decomposition being unique up to
diffeomorphism of the factors and the order of the factors.} factor
of $M_0$ is acyclic, i.e., every prime factor of $M_0$ has either
non-trivial $\pi_2$ or non-trivial $\pi_3$.
\end{enumerate}
\end{cor}

\begin{proof}
Repeated application of Theorem~\ref{surgery} shows that (1) implies
(2). The implication (2) implies (3) is immediate from van Kampen's
theorem. The fact that (3) implies (1) is Theorem~\ref{finiteext}.
This shows that (1), (2) and (3) are all equivalent. Since
three-dimensional spherical space-forms and $S^2$-bundles over $S^1$
are easily seen to be prime,  (2) implies  (4). Thus, it remains
only to see that (4) implies (3). We consider a manifold $M$
satisfying (4), a prime factor $P$ of $M$, and  universal covering
$\widetilde P$ of $P$. First suppose that $\pi_2(P)=\pi_2(\widetilde
P)$ is trivial. Then, by hypothesis $\pi_3(P)=\pi_3(\widetilde P)$
is non-trivial. By the Hurewicz theorem this means that
$H_3(\widetilde P)$ is non-trivial, and hence that $\widetilde P$ is
a compact, simply connected three-manifold. It follows that
$\pi_1(P)$ is finite. Now suppose that $\pi_2(P)$ is non-trivial.
Then $P$ is not diffeomorphic to $\Ar P^3$. Since $P$ is prime and
contains no locally separating $\Ar P^2$, it follows that $P$
contains no embedded $\Ar P^2$. Then by the sphere theorem there is
an embedded two-sphere in $P$ that is homotopically non-trivial.
Since $P$ is prime, this sphere cannot separate, so cutting $P$ open
along it produces a connected manifold $P_0$ with two boundary
two-spheres. Since $P_0$ is prime, it follows that $P_0$ is
diffeomorphic to $S^2\times I$ and hence $P$ is diffeomorphic to a
two-sphere bundle over the circle.
\end{proof}

\begin{rem}
(i) The use of the sphere theorem is unnecessary in the above
argument for what we actually prove is that if every prime factor of
$M$ has non-trivial $\pi_2$ or non-trivial $\pi_3$, then the Ricci
flow with surgery with $(M,g_0)$ as initial metric becomes extinct
after a finite time. In fact, the sphere theorem for closed
three-manifolds follows from the results here.

\noindent (ii)  If the initial manifold is simpler then all the
time-slices are simpler: If $({\mathcal M},G)$ is a Ricci flow with
surgery whose initial manifold is prime, then every time-slice is a
disjoint union of connected components, all but at most one being
diffeomorphic to a three-sphere and if there is one not
diffeomorphic to a three-sphere, then it is diffeomorphic to the
initial manifold. If the initial manifold is a simply connected
manifold $M_0$, then every component of every time-slice $M_T$ must
be simply connected, and thus {\em a posteriori} every time-slice is
a disjoint union of manifolds diffeomorphic to the three-sphere.
Similarly, if the initial manifold has finite fundamental group,
then every connected component of every time-slice is either simply
connected or has the same fundamental group as the initial manifold.

\noindent (iii) The conclusion of this result is a natural
generalization of Hamilton's conclusion in analyzing the Ricci flow
on manifolds of positive Ricci curvature in \cite{Hamilton3MPRC}.
Namely, under appropriate hypotheses, during the evolution process
of Ricci flow with surgery the manifold breaks into components each
of which disappears in finite time. As a component disappears at
some finite time, the metric on that component is well enough
controlled to show that the disappearing component admits a
non-flat, homogeneous Riemannian metric of non-negative sectional
curvature, i.e., a metric locally isometric to either a round $S^3$
or to a product of a round $S^2$ with the usual metric on $ \Ar$.
The existence of such a metric on a component immediately gives the
topological conclusion of Theorem~\ref{Theorem1} for that component,
i.e., that it is diffeomorphic to a $3$-dimensional spherical
space-form, to $S^2\times S^1$ to a non-orientable $2$-sphere bundle
over $S^1$, or to $\Ar P^3\#\Ar P^3$. The biggest difference between
these two results is that Hamilton's hypothesis is geometric
(positive Ricci curvature) whereas Perelman's is homotopy theoretic
(information about the fundamental group).

\noindent (iv) It is also worth pointing out that it follows from
Corollary~\ref{corequiv} that the manifolds that satisfy the four equivalent
conditions in that corollary are exactly the closed, connected, three-manifolds
that admit a Riemannian metric of positive scalar curvature, cf,
\cite{SchoenYau2} and \cite {GromovLawson}.
\end{rem}

One can use Ricci flow in a more general study of three-manifolds
than the one we carry out here. There is a conjecture due to
Thurston, see \cite{Thurstongeom},  known as Thurston's
Geometrization Conjecture\index{Thurston's Geometrization
Conjecture|ii} or simply as the Geometrization Conjecture for
three-manifolds. It conjectures that every $3$-manifold without
locally separating $\Ar P^2$'s (in particular every orientable
$3$-manifold) is a connected sum of prime $3$-manifolds each of
which admits a decomposition along incompressible\footnote{I.e.,
embedded by a map that is injective on $\pi_1$.} tori into pieces
that admit locally homogeneous geometries of finite volume.  Modulo
questions about cofinite-volume lattices in $SL_2(\Cee)$, proving
this conjecture leads to a complete classification
 of  $3$-manifolds without locally separating $\Ar P^2$'s,
and in particular to a complete classification of all orientable
$3$-manifolds. (See Peter Scott's survey article \cite{Scott}.) By
passing to the orientation double cover and working equivariantly,
these results can be extended to all $3$-manifolds.

 Perelman in
\cite{Perelman2} has stated results which imply a positive
resolution of Thurston's Geometrization conjecture. Perelman's
proposed proof of Thurston's Geometrization Conjecture relies in an
essential way on Theorem~\ref{surgery}, namely the existence of
Ricci flow with surgery for all positive time. But it also involves
a further analysis of the limits of these Ricci flows as time goes
to infinity. This further analysis involves analytic arguments which
are exposed in Sections 6 and 7 of Perelman's second paper
(\cite{Perelman2}), following earlier work of Hamilton
(\cite{HamiltonNSRF3M}) in a simpler case of bounded curvature. They
also involve a result (Theorem 7.4 from \cite{Perelman2}) from the
theory of manifolds with curvature locally bounded below that are
collapsed, related to results of Shioya-Yamaguchi \cite{SY}. The
Shioya-Yamaguchi results in turn rely on an earlier, unpublished
work of Perelman proving the so-called `Stability Theorem.'
Recently, Kapovich, \cite{Kap} has put a preprint on the archive
giving a proof of the stability result. We have been examining
another approach, one suggested by Perelman in \cite{Perelman2},
avoiding the stability theorem, cf, \cite{KleinerLott2} and
\cite{MorganTian2}. It is our view that the collapsing results
needed for the Geometrization Conjecture are in place, but that
before a definitive statement that the Geometrization Conjecture has
been resolved can be made these arguments must be subjected to the
same close scrutiny that the arguments proving the Poincar\'e
Conjecture have received. This process is underway.

In this book we do not attempt to explicate any of the results
beyond Theorem~\ref{surgery} described in the previous paragraph
that are needed for the Geometrization Conjecture. Rather, we
content ourselves with presenting a proof of Theorem~\ref{Theorem1}
above which, as we have indicated, concerns initial Riemannian
manifolds for which the Ricci flow with surgery becomes extinct
after finite time. We are currently preparing a detailed proof,
along the lines suggested by Perelman, of the further results that
will complete the proof of the Geometrization Conjecture.

As should be clear from the above overview, Perelman's argument did
not arise in a vacuum. Firstly, it resides in a context provided by
the general theory of Riemannian manifolds. In particular,  various
notions of convergence of sequences of manifolds play a crucial
role. The most important is geometric convergence (smooth
convergence on compact subsets). Even more importantly, Perelman's
argument resides in the context of the theory of the Ricci flow
equation, introduced by Richard Hamilton and extensively studied by
him and others. Perelman makes use of almost every previously
established result for $3$-dimensional Ricci flows. One exception is
Hamilton's  proposed classification results for three-dimensional
singularities. These are replaced by Perelman's strong qualitative
description of singularity development for Ricci flows on compact
three-manifolds.

The first five chapters of the book review the necessary background
material from these two subjects. Chapters 6 through 11 then explain
Perelman's advances. In Chapter 12 we introduce the standard
solution\index{Ricci flow!standard solution}, which is the manifold
constructed by hand that one `glues in' in doing surgery. Chapters
13 through 17 describe in great detail the surgery process and prove
the main analytic and topological estimates that are needed to show
that one can continue the process for all positive time. At the end
of Chapter 17 we have established Theorem~\ref{surgery}. Chapter 18
discusses the finite-time extinction result. Chapter 19 is an
appendix on some topological results that were needed in the surgery
analysis in Chapters 13-17.

\section{Background material from Riemannian geometry}

\subsection{Volume and injectivity radius}

One important general concept that is used throughout is the notion
of a manifold being non-collapsed at a point. Suppose that we have a
point $x$ in a complete Riemannian $n$-manifold. Then we say that
the manifold is {\em
$\kappa$-non-collapsed}\index{$\kappa$-non-collapsed} at $x$
provided that the following holds: For any $r$ such that the norm of
the Riemannian curvature tensor, $|{\rm Rm}|$, is $\le r^{-2}$ at
all points of the metric ball, $B(x,r)$, of radius $r$ centered at
$x$,  we have ${\rm Vol}\, B(x,r)\ge \kappa r^n$. There is a
relationship between this notion and the injectivity
radius\index{injectivity radius} of $M$ at $x$. Namely, if $|{\rm
Rm}|\le r^{-2}$ on $B(x,r)$ and if $B(x,r)$ is
$\kappa$-non-collapsed then the injectivity radius of $M$ at $x$ is
greater than or equal to  a positive constant that depends only on
$r$ and $\kappa$. The advantage of working with the volume
non-collapsing condition is that, unlike for the injectivity radius,
there is a simple equation for the evolution of volume under Ricci
flow.

Another important general result is the Bishop-Gromov volume
comparison\index{Bishop-Gromov Theorem}\index{volume
comparison|see{Bishop-Gromov Theorem}} result that says that if the
Ricci curvature of a complete Riemannian $n$-manifold $M$ is bounded
below by a constant $(n-1)K$  then for any $x\in M$ the ratio of the
volume of $B(x,r)$ to the volume of the ball of radius $r$ in the
space of constant curvature $K$ is a non-increasing function whose
limit as $r\rightarrow 0$ is $1$.

All of these basic facts from Riemannian geometry are reviewed in
the first chapter.

\subsection{Manifolds of non-negative curvature}

For reasons that should be clear from the above description and in
any event will become much clearer shortly, manifolds of
non-negative curvature play an extremely important role in the
analysis of Ricci flows with surgery. We need several general
results about them. The first is the soul theorem for manifolds of
non-negative sectional curvature. A {\em soul} is a compact, totally
geodesic submanifold. The entire manifold is diffeomorphic to the
total space of a vector bundle over any of its souls. If a
non-compact $n$-manifold has positive sectional curvature, then any
soul for it is a point, and in particular, the manifold is
diffeomorphic to Euclidean space. In addition, the distance function
$f$ from a soul has the property that for every $t>0$ the pre-image
$f^{-1}(t)$ is homeomorphic to an $(n-1)$-sphere and the pre-image
under this distance function of any non-degenerate interval
$I\subset \Ar ^+$ is homeomorphic to $S^{n-1}\times I$.

Another important result is the splitting theorem\index{splitting
theorem}, which says that, if a complete manifold of non-negative
sectional curvature has a geodesic line (an isometric copy of $\Ar$)
that is distance minimizing between every pair of its points, then
that manifold is a metric product of a manifold of one lower
dimension and $\Ar$. In particular, if a complete $n$-manifold of
non-negative sectional curvature has two ends then it is a metric
product $N^{n-1}\times \Ar$ where $N^{n-1}$ is a compact manifold.

Also, we need some of the elementary comparison results from Toponogov
theory\index{Toponogov theory}. These compare ordinary triangles in the
Euclidean plane with triangles in a manifold of non-negative sectional
curvature whose sides are minimizing geodesics in that manifold.

\subsection{Canonical neighborhoods}

\begin{figure}[ht]
                {\centerline{\epsfbox{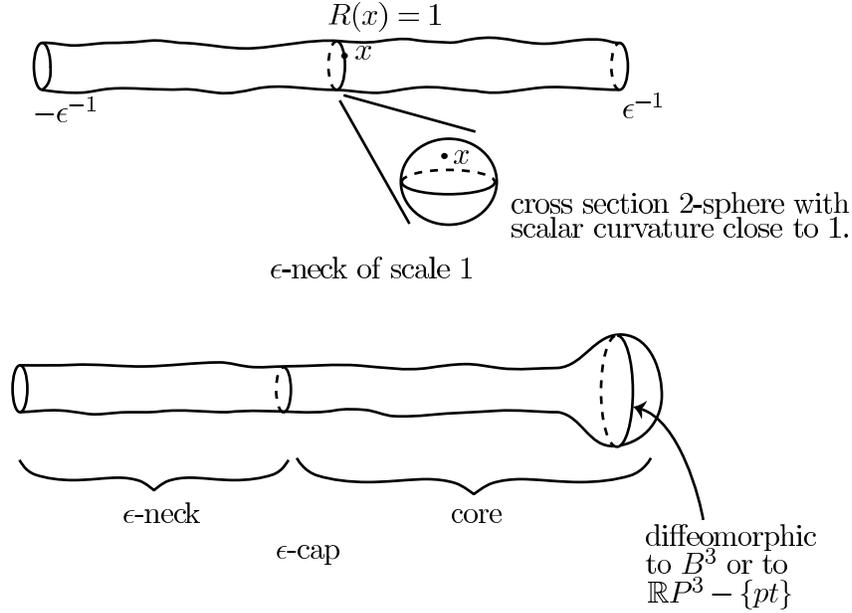}}}
                \caption{Canonical neighborhoods.}
                \label{fig:epsneck}
        \end{figure}

Much of the analysis of the geometry of Ricci flows revolves around
the notion of canonical neighborhoods\index{canonical neighborhood}.
Fix some $\epsilon>0$ sufficiently small. There are two types of
non-compact canonical neighborhoods: $\epsilon$-necks
\index{$\epsilon$-neck}and $\epsilon$-caps\index{$\epsilon$-cap}. An
$\epsilon$-neck in a Riemannian $3$-manifold $(M,g)$ centered at a
point $x\in M$ is a submanifold $N\subset M$ and a diffeomorphism
$\psi\colon S^2\times (-\epsilon^{-1},\epsilon^{-1})\to N$ such that
$x\in \psi(S^2\times \{0\})$ and such that the pullback of the
rescaled metric, $\psi^*(R(x)g)$, is within $\epsilon$ in the
$C^{[1/\epsilon]}$-topology of the product of the round metric of
scalar curvature $1$ on $S^2$ with the usual metric on the interval
$(-\epsilon^{-1},\epsilon^{-1})$. (Throughout, $R(x)$ denotes the
scalar curvature of $(M,g)$ at the point $x$.) An $\epsilon$-cap is
a non-compact submanifold ${\mathcal C}\subset M$ with the property
that a neighborhood $N$ of infinity in ${\mathcal C}$ is an
$\epsilon$-neck, such that every point of $N$ is the center of an
$\epsilon$-neck in $M$, and such that the {\em
core}\index{$\epsilon$-cap!core of}, ${\mathcal C}\setminus
\overline N$, of the $\epsilon$-cap is diffeomorphic to either a
$3$-ball or a punctured $\Ar P^3$. It will also be important to
consider $\epsilon$-caps that, after rescaling to make $R(x)=1$ for
some point $x$ in the cap, have bounded geometry (bounded diameter,
bounded ratio of the curvatures at any two points, and bounded
volume). If $C$ represents the bound for these quantities, then we
call the cap an $(C,\epsilon)$-cap\index{$(C,\epsilon)$-cap}. See
{\sc Fig.}~\ref{fig:epsneck}.
 An
$\epsilon$-tube\index{$\epsilon$-tube} in $M$ is a submanifold of
$M$ diffeomorphic to $S^2\times (0,1)$  which is a union of
$\epsilon$-necks and with the property that each point of the
$\epsilon$-tube is the center of an $\epsilon$-neck in $M$.

There are two other types of canonical neighborhoods in
$3$-manifolds -- (i) a $C$-component\index{$C$-component} and (ii)
an $\epsilon$-round component\index{$\epsilon$-round component}. The
$C$-component is a compact, connected Riemannian manifold of
positive sectional curvature diffeomorphic to either $S^3$ or $\Ar
P^3$ with the property that rescaling  the metric by $R(x)$ for any
$x$ in the component produces a Riemannian manifold whose  diameter
 is at most $C$, whose sectional curvature at any point and in any
$2$-plane direction is between $C^{-1}$ and $C$, and whose volume is
between $C^{-1}$ and $C$. An $\epsilon$-round component is a
component on which the metric rescaled by $R(x)$ for any $x$ in the
component is within $\epsilon$ in the $C^{[1/\epsilon]}$-topology of
a round metric of scalar curvature one.

As we shall see, the singularities at time $T$ of a $3$-dimensional
Ricci flow are contained in subsets that are unions of canonical
neighborhoods with respect to the metrics at nearby, earlier times
$t'<T$. Thus, we need to understand the topology of manifolds that
are unions of $\epsilon$-tubes and $\epsilon$-caps. The fundamental
observation is that, provided that $\epsilon$ is sufficiently small,
when two $\epsilon$-necks intersect (in more than a small
neighborhood of the boundaries) their product structures almost line
up, so that the two $\epsilon$-necks can be glued together to form a
manifold fibered by $S^2$'s. Using this idea we show that, for
$\epsilon>0$ sufficiently small, if a connected manifold is a union
of $\epsilon$-tubes and $\epsilon$-caps  then it is diffeomorphic to
$\Ar^3$,  $S^2\times \Ar$, $S^3$, $S^2\times S^1$, $\Ar P^3\#\Ar
P^3$,  the total space of a line bundle over $\Ar P^2$, or  the
non-orientable $2$-sphere bundle over $S^1$. This topological result
is proved in the appendix at the end of the book. {\bf We shall fix
$\epsilon>0$ sufficiently small so that these results hold.}

There is one result relating canonical neighborhoods and manifolds
of positive curvature of which we make repeated use: Any complete
$3$-manifold of positive curvature does not admit $\epsilon$-necks
of arbitrarily high curvature. In particular, if $M$ is a complete
Riemannian $3$-manifold with the property that every point of scalar
curvature greater than $r_0^{-2}$ has a canonical neighborhood, then
$M$ has bounded curvature. This turns out to be of central
importance and is used repeatedly.

All of these basic facts about Riemannian manifolds of non-negative
curvature are recalled in the second chapter.

\section{Background material from Ricci flow}

Hamilton \cite{Hamilton3MPRC} introduced the Ricci flow
equation\index{Ricci flow!equation},
$$\frac{\partial g(t)}{\partial t}=-2{\rm Ric}(g(t)).$$
This is an evolution equation for a one-parameter family of Riemannian metrics
$g(t)$ on a smooth manifold $M$. The Ricci flow equation is weakly parabolic
and is strictly parabolic modulo the `gauge group', which is the group of
diffeomorphisms of the underlying smooth manifold. One should view this
equation as a non-linear, tensor version of the heat equation. From it, one can
derive the evolution equation for the Riemannian metric tensor, the Ricci
tensor, and the scalar curvature function. These are all parabolic equations.
For example, the evolution equation for scalar curvature $R(x,t)$ is
\begin{equation}\label{scalarflow} \frac{\partial R}{\partial t}(x,t)=\triangle
R(x,t)+2|{\rm Ric}(x,t)|^2, \end{equation} illustrating the
similarity with the heat equation. (Here $\triangle$ is the
Laplacian with non-positive spectrum.)

\subsection{First results}

Of course, the first results we need are uniqueness and short-time existence
for solutions to the Ricci flow equation for compact manifolds. These results
were proved by Hamilton (\cite{Hamilton3MPRC}) using the Nash-Moser inverse
function theorem,
 (\cite{HamiltonIFTNM}). These results are standard for
strictly parabolic equations. By now there is a fairly standard method for
working `modulo' the gauge group (the group of diffeomorphisms) and hence
arriving at a strictly parabolic situation where the classical existence,
uniqueness and smoothness results apply. The method for the Ricci flow equation
goes under the name of `DeTurck's trick.'

There is also a result that allows us to patch together local solutions
$(U,g(t)),\ a\le t\le b$, and $(U,h(t)),\ b\le t\le c$, to form a smooth
solution defined on the interval $a\le t\le c$ provided that $g(b)=h(b)$.

Given a Ricci flow $(M,g(t))$ we can always translate, replacing $t$
by $t+t_0$ for some fixed $t_0$, to produce a new Ricci flow. We can
also rescale by any positive constant $Q$ by setting
$h(t)=Qg(Q^{-1}t)$ to produce a new Ricci flow.

\subsection{Gradient shrinking solitons}

Suppose that $(M,g)$ is a complete Riemannian manifold, and suppose
that there is a constant $\lambda>0$ with the property that
$${\rm Ric}(g)=\lambda g.$$
In this case, it is easy to see that there is a Ricci flow given by
$$g(t)=(1-2\lambda t)g.$$ In particular, all the metrics in this flow
differ by a constant factor depending on time and the metric is a
decreasing function of time. These are called {\em shrinking
solitons}\index{gradient shrinking soliton}. Examples are compact
manifolds of constant positive Ricci curvature.

There is a closely related, but more general, class of examples: the
{\em gradient shrinking solitons}. Suppose that $(M,g)$ is
 a complete Riemannian manifold, and suppose that  there is a constant
$\lambda>0$ and a function $f\colon M\to \Ar$ satisfying
$${\rm Ric}(g)=\lambda g-{\rm Hess}^{g}f.$$
 In this case, there is a Ricci flow which is a shrinking
 family after we pull back by the one-parameter family of diffeomorphisms generated by
 the time-dependent vector field $\frac{1}{1-2\lambda t}\nabla_{g} f$.
 An example of a gradient shrinking soliton is the manifold
 $S^2\times \Ar$ with the family of metrics being the product of the
 shrinking family of round metrics on $S^2$ and the constant
 family of standard metrics on $\Ar$. The function $f$ is $s^2/4$ where $s$ is the
 Euclidean parameter on $\Ar$.

\subsection{Controlling higher derivatives of curvature}

Now let us discuss the smoothness results for geometric limits. The
general result along these lines is Shi's theorem\index{Shi's
derivative estimates}\index{Shi's Theorem|see{Shi's derivative
estimates}}, see \cite{Shi1,Shi2}. Again, this is a standard type of
result for parabolic equations. Of course, the situation here is
complicated somewhat by the existence of the gauge group. Roughly,
Shi's theorem says the following. Let us denote by $B(x,t_0,r)$ the
metric ball in $(M,g(t_0))$ centered at $x$ and of radius $r$. If we
can control the norm of the Riemannian curvature tensor on a
backward neighborhood of the form $B(x,t_0,r)\times [0,t_0]$, then
for each $k>0$ we can control the $k^{th}$ covariant derivative of
the curvature on $B(x,t_0,r/2^k)\times [0,t_0]$ by a constant over
$t^{k/2}$. This result has many important consequences in our study
because it tells us that geometric limits are smooth limits. Maybe
the first result to highlight is the fact (established earlier by
Hamilton) that if $(M,g(t))$ is a Ricci flow defined on $0\le
t<T<\infty$, and if the Riemannian curvature is uniformly bounded
for the entire flow, then the Ricci flow extends past time $T$.

In the third chapter this material is reviewed and, where necessary,
slight variants of results and arguments in the literature are
presented.

\subsection{Generalized Ricci flows}

Because we cannot restrict our attention to Ricci flows, but rather
must consider more general objects, Ricci flows with surgery, it is
important to establish the basic analytic results and estimates in a
context more general than that of Ricci flow. We choose to do this
in the context of generalized Ricci flows\index{Ricci
flow!generalized}.

A generalized three-dimensional Ricci flow consists of a smooth
four-dimensional manifold ${\mathcal M}$ (space-time)\index{Ricci
flow!space-time of} with a time function ${\bf t}\colon {\mathcal M}\to \Ar$
and a smooth vector field $\chi$. These are required to satisfy:
\begin{enumerate}
\item Each $x\in {\mathcal M}$ has a neighborhood of the form
$U\times J$, where $U$ is an open subset in $\Ar^3$ and $J\subset
\Ar$ is an interval, in which ${\bf t}$ is the projection onto $J$
and $\chi$ is the unit vector field tangent to the one-dimensional
foliation $\{u\}\times J$ pointing in the direction of increasing
${\bf t}$. We call ${\bf t}^{-1}(t)$ the $t$ time-slice. It is a
smooth $3$-manifold.
\item The image ${\bf t}({\mathcal M})$ is a connected interval $I$ in
$\Ar$, possibly infinite. The boundary of ${\mathcal M}$ is the
pre-image under ${\bf t}$ of the boundary of $I$.
\item The level sets ${\bf t}^{-1}(t)$ form a codimension-one foliation of
${\mathcal M}$, called the horizontal foliation, with the boundary
components of ${\mathcal M}$ being leaves.
\item There is a metric $G$ on the horizontal distribution, i.e., the distribution
tangent to the level sets of ${\bf t}$. This metric induces a
Riemannian metric on each $t$ time-slice varying smoothly as we vary
the time-slice. We define the curvature of $G$ at a point
$x\in{\mathcal M}$ to be the curvature of the Riemannian metric
induced by $G$ on the time-slice $M_t$ at $x$.
\item Because of the first property the integral curves of $\chi$
preserve the horizontal foliation and hence the horizontal
distribution. Thus, we can take the Lie derivative of $G$ along
$\chi$. The Ricci flow equation is then
$${\mathcal L}_\chi(G)=-2{\rm Ric}(G).$$
\end{enumerate}

Locally in space-time the horizontal metric is simply a smoothly varying family
of Riemannian metrics on a fixed smooth manifold and the evolution equation is
the ordinary Ricci flow equation. This means that the usual formulas for the
evolution of the curvatures as well as much of the analytic analysis of Ricci
flows still hold in this generalized context. In the end, a Ricci flow with
surgery is a more singular type of space-time, but it will have an open dense
subset which is a generalized Ricci flow, and all the analytic estimates take
place in this open subset.

The notion of canonical neighborhoods make sense in the context of
generalized Ricci flows. There is also the notion of a strong
$\epsilon$-neck. Consider an embedding $\psi\colon \left(S^2\times
(-\epsilon^{-1},\epsilon^{-1})\right)\times (-1,0]$ into space-time
such that the time function pulls back to the projection onto
$(-1,0]$ and the vector field $\chi$ pulls back to
$\partial/\partial t$. If there is such an embedding into an
appropriately shifted and rescaled version of the original
generalized Ricci flow so that the pull-back of the rescaled
horizontal metric is within $\epsilon$ in the
$C^{[1/\epsilon]}$-topology of the product of the shrinking family
of round $S^2$'s with the Euclidean metric on
$(-\epsilon^{-1},\epsilon^{-1})$, then we say that $\psi$ is a
strong $\epsilon$-neck\index{$\epsilon$-neck!strong} in the
generalized Ricci flow.

\subsection{The maximum principle}\label{sectmax}

The Ricci flow equation satisfies various forms of the maximum
principle. The fourth chapter explains this principle, which is due
to Hamilton (see Section 4 of \cite{Hamiltonsurvey}), and derives
many of its consequences,  which are also due to Hamilton (cf.
\cite{HamiltonNSRF3M}). This principle and its consequences are at
the core of all the detailed results about the nature of the flow.
We illustrate the idea by considering the case of the scalar
curvature. A standard scalar maximum principle argument applied to
Equation~(\ref{scalarflow}) proves that the minimum of the scalar
curvature is a non-decreasing function of time. In addition, it
shows that if the minimum of scalar curvature at time $0$ is
positive then we have
$$R_{\rm min}(t)\ge R_{\rm min}(0)\left(\frac{1}{1-\frac{2t}{n}R_{\rm
min}(0)}\right),$$ and thus the equation must develop a singularity at or
before time $n/\left(2R_{\rm min}(0)\right)$.

While the above result about the scalar curvature is important and
is used repeatedly, the most significant uses of the maximum
principle involve the tensor version, established by Hamilton, which
applies for example to the Ricci tensor and the full curvature
tensor.  These have given the most significant understanding of the
Ricci flows, and they form the core of the arguments that Perelman
uses in his application of Ricci flow to $3$-dimensional topology.
Here are the main results established by Hamilton:

\begin{enumerate}
\item[(1)] For $3$-dimensional flows, if the Ricci curvature is positive,
then the family of metrics becomes singular at finite time and as
the family becomes singular, the metric becomes closer and closer to
round; see \cite{Hamilton3MPRC}.
\item[(2)] For $3$-dimensional flows, as the scalar curvature goes to $+\infty$
the ratio of the absolute value of any negative eigenvalue of the
Riemannian curvature to the largest positive eigenvalue goes to
zero; see \cite{HamiltonNSRF3M}. This condition is called {\em
pinched toward positive curvature}\index{curvature!pinched toward
positive}.
\item[(3)] Motivated by a Harnack inequality for the heat equation established
by Li-Yau \cite{LiYau}, Hamilton established a Harnack
inequality\index{Ricci flow!Harnack inequality for} for the
curvature tensor under the Ricci flow for complete manifolds
$(M,g(t))$  with bounded, non-negative curvature operator; see
\cite{Hamiltonharnack}. In the applications to three dimensions, we
shall need the following consequence for the scalar curvature:
Suppose that $(M,g(t))$ is  a Ricci flow defined for all $t\in
[T_0,T_1]$ of complete manifolds of non-negative curvature operator
with bounded curvature. Then
$$\frac{\partial R}{\partial t}(x,t)+\frac{R(x,t)}{t-T_0}\ge 0.$$
In particular,
 if $(M,g(t))$ is an ancient solution (i.e.,
defined for all $t\le 0$) of bounded, non-negative curvature then
$\partial R(x,t)/\partial t\ge 0$.
\item[(4)] If a complete $3$-dimensional Ricci flow
$(M,g(t)),\ 0\le t\le T$, has non-negative curvature, if $g(0)$ is
not flat, and if there is at least one point $(x,T)$ such that the
Riemannian curvature tensor of $g(T)$ has a flat direction in
$\wedge^2TM_x$, then $M$ has a cover $\widetilde M$ so that for each
$t>0$ the Riemannian manifold $(\widetilde M,g(t))$ splits as a
Riemannian product of a surface of positive curvature and a
Euclidean line. Furthermore, the flow on the cover $\widetilde M$ is
the product of a $2$-dimensional flow and the trivial
one-dimensional Ricci flow on the line; see Sections 8 and 9 of
\cite{Hamilton4MPCO}.
\item[(5)] In particular, there is no Ricci flow of non-negative curvature tensor
$(U,g(t)),$ defined for $\ 0\le t\le T$ with $T>0$, such that
$(U,g(T))$ is isometric to an open subset in a non-flat,
$3$-dimensional metric cone.
\end{enumerate}

\subsection{Geometric limits}

In\index{geometric limit} the fifth chapter we discuss geometric
limits of Riemannian manifolds and of Ricci flows. Let us review the
history of these ideas. The first results about geometric limits of
Riemannian manifolds go back to Cheeger in his thesis in 1967; see
\cite{Cheeger}. Here Cheeger obtained topological results. In
\cite{Gromov} Gromov proposed that geometric limits should exist in
the Lipschitz topology and suggested a result along these lines,
which also was known to Cheeger.
 In
\cite{GreeneWu}, Greene-Wu gave a rigorous proof of the compactness
theorem suggested by Gromov and also enhanced the convergence to be
$C^{1,\alpha}$-convergence by using harmonic coordinates; see also
\cite{SPeters}. Assuming that all the derivatives of curvature are
bounded, one can apply elliptic theory to the expression of
curvature in harmonic coordinates and deduce $C^\infty$-convergence.
These ideas lead to various types of compactness results  that go
under the name Cheeger-Gromov compactness for Riemannian manifolds.
Hamilton in \cite{Hamiltonlimits} extended these results to Ricci
flows. We shall use the compactness results for both Riemannian
manifolds and for Ricci flows. In a different direction, geometric
limits were extended to the non-smooth context by Gromov in
\cite{Gromov} where he introduced a weaker topology, called the
Gromov-Hausdorff topology and proved a compactness theorem.

 Recall that a
sequence of based Riemannian manifolds $(M_n,g_n,x_n)$ is said to
{\em converge geometrically}\index{converge!geometrically} to a
based, complete Riemannian manifold $(M_\infty,g_\infty,x_\infty)$
if there is a sequence of open subsets $U_n\subset M_\infty$ with
compact closures, with $x_\infty\in U_1\subset \overline U_1\subset
U_2\subset \overline U_2\subset U_3\subset \cdots $ with
$\cup_nU_n=M_\infty$, and embeddings $\varphi_n\colon U_n\to M_n$
sending $x_\infty$ to $x_n$ so that the pull back metrics,
$\varphi_n^*g_n$, converge uniformly on compact subsets of
$M_\infty$ in the $C^\infty$-topology to $g_\infty$. Notice that the
topological type of the limit can be different from the topological
type of the manifolds in the sequence. There is a similar notion of
geometric convergence for a sequence of based Ricci flows.

Certainly, one of the most important consequences of Shi's results,
cited above, is that, in concert with  Cheeger-Gromov  compactness,
it allows us to form smooth geometric limits of sequences of based
Ricci flows. We have the following result of Hamilton's; see
\cite{Hamiltonlimits}:

\begin{thm}
Suppose we have a sequence of based Ricci flows
$(M_n,g_n(t),(x_n,0))$ defined for $t\in (-T,0]$ with the
$(M_n,g_n(t))$ being complete. Suppose that:
\begin{enumerate}
\item[(1)] There is $r>0$ and $\kappa>0$ such that for every $n$ the metric ball $B(x_n,0,r)\subset
(M_n,g_n(0))$ is $\kappa$-non-collapsed.
\item[(2)] For each $A<\infty$ there is $C=C(A)<\infty$ such that the Riemannian
curvature on $B(x_n,0,A)\times (-T,0]$ is bounded by $C$.
\end{enumerate}
Then after passing to a subsequence there is a geometric limit which
is a based Ricci flow $(M_\infty,g_\infty(t),(x_\infty,0))$ defined
for $t\in (-T,0]$.
\end{thm}

To emphasize, the two conditions that we must check in order to
extract a geometric limit of a subsequence based at points at time
zero are: (i) uniform non-collapsing at the base point in the time
zero metric, and (ii) for each $A<\infty$ uniformly bounded
curvature for the restriction of the flow to the metric balls of
radius $A$ centered at the base points.

Most steps in Perelman's argument require invoking this result in order to form
limits of appropriate sequences of Ricci flows, often rescaled to make the
scalar curvatures at the base point equal to $1$. If, before rescaling, the
scalar curvature at the base points goes to infinity as we move through the
sequence, then the resulting limit of the rescaled flows has non-negative
sectional curvature. This is a consequence of the fact that the sectional
curvatures of the manifolds in the sequence are uniformly pinched toward
positive. It is for exactly this reason that non-negative curvature plays such
an important role in the study of singularity development in three-dimensional
Ricci flows.

\section{Perelman's advances}

So far we have been discussing the results that were known before
Perelman's work. They concern almost exclusively Ricci flow (though
Hamilton in \cite{Hamilton4MPIC} had introduced the notion of
surgery and proved that surgery can be performed preserving the
condition that the curvature is pinched toward positive, as in (2)
above). Perelman extended in two essential ways the analysis of
Ricci flow -- one involves the introduction of a new analytic
functional, {\em the reduced length}, which is the tool by which he
establishes the needed non-collapsing results, and the other is a
delicate combination of geometric limit ideas and consequences of
the maximum principle together with the non-collapsing results in
order to establish bounded curvature at bounded distance results.
These are used to prove in an inductive way the existence of
canonical neighborhoods, which is a crucial ingredient in proving
that is possible to do surgery iteratively, creating a flow defined
for all positive time.

While it is easiest to formulate and consider these techniques in
the case of Ricci flow, in the end one needs them in the more
general context of Ricci flow with surgery since we inductively
repeat the surgery process, and in order to know at each step that
we can perform surgery we need to apply these results to the
previously constructed Ricci flow with surgery. We have chosen to
present these new ideas only once -- in the context of generalized
Ricci flows -- so that we can derive the needed consequences in all
the relevant contexts from this one source.

\subsection{The reduced length function}

\begin{figure}[ht]
  \relabelbox{
  \centerline{\epsfbox{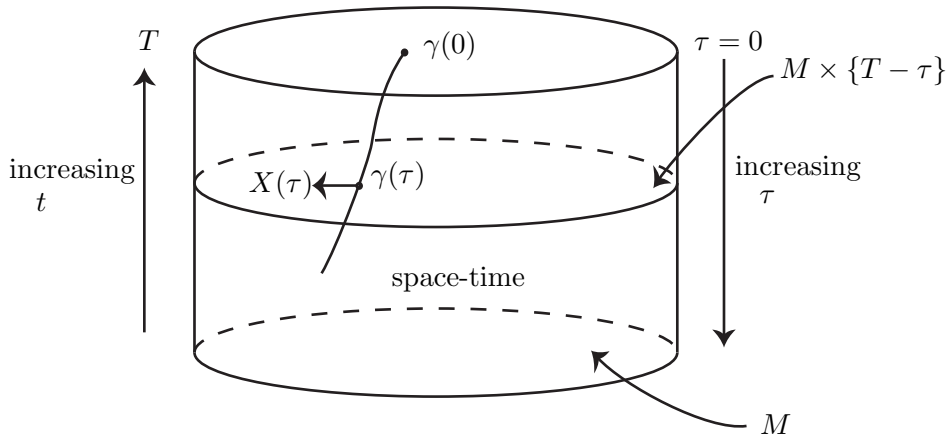}}}
  \relabel{1}{$T$}
  \relabel{2}{increasing}
  \relabel{3}{$t$}
  \relabel{4}{$\gamma(0)$}
  \relabel{5}{$\gamma(\tau)$}
  \relabel{6}{$X(\tau)$}
  \relabel{7}{space-time}
  \relabel{8}{$\tau=0$}
  \relabel{9}{increasing}
  \relabel{a}{$\tau$}
  \relabel{b}{$M$}
  \relabel{c}{$M\times\{T-\tau\}$}
  \endrelabelbox
  \caption{Curves in space-time parameterized by
  $\tau$.}\label{fig:path}
\end{figure}

In\index{reduced length function} Chapter~\ref{lengthfn} we come to
the first of Perelman's major contributions. Let us first describe
it in the context of an ordinary three-dimensional Ricci flow, but
viewing the Ricci flow as a horizontal metric on a space-time which
is the manifold $M\times I$, where $I$ is the interval of definition
of the flow. Suppose that $I=[0,T)$ and fix $(x,t)\in M\times
(0,T)$. We consider paths $\gamma(\tau),\ 0\le \tau\le
\overline\tau$, in space-time with the property that for every
$\tau\le \overline\tau$ we have $\gamma(\tau)\in M\times \{t-\tau\}$
and $\gamma(0)=x$. These paths are said to be {\em parameterized by
backward time}. See {\sc Fig.}~\ref{fig:path}. The ${\mathcal
L}$-{\em length} of such a path is given by
$${\mathcal L}(\gamma)=\int_0^{\overline
\tau}\sqrt{\tau}\left(R(\gamma(\tau))+|\gamma'(\tau)|^2\right)d\tau,$$
where the derivative on $\gamma$ refers to the spatial derivative.
There is also the closely related {\em reduced length}
$$\ell(\gamma)=\frac{{\mathcal L}(\gamma)}{2\sqrt{\overline\tau}}.$$
There is a theory for the functional ${\mathcal L}$ analogous to the
theory for the usual energy function\footnote{Even though this
functional is called a length, the presence of the
$|\gamma'(\tau)|^2$ in the integrand means that it behaves more like
the usual energy functional for paths in a Riemannian manifold.}.
 In
particular, there is the notion of an ${\mathcal
L}$-geodesic\index{${\mathcal L}$-geodesic}, and the reduced length
as a function on space-time $\ell_{(x,t)}\colon M\times [0,t)\to
\Ar$. One establishes a crucial monotonicity\index{reduced length
function!monotonicity of} for this reduced length along ${\mathcal
L}$-geodesics. Then one defines the reduced volume
$$\widetilde
V_{(x,t)}(U\times\{\bar t\})=\int_{U\times\{\bar
t\}}\bar\tau^{-3/2}e^{-\ell_{(x,t)}(q,\bar\tau)}d{\rm
vol}_{g(\bar\tau}(q),$$ where, as before $\bar\tau=t-\bar t$.
Because of the monotonicity of $\ell_{(x,t)}$ along ${\mathcal
L}$-geodesics, the reduced volume\index{reduced volume} is also
non-increasing under the flow (forward in $\bar\tau$ and hence
backward in time) of open subsets along ${\mathcal L}$-geodesics.
This is the fundamental tool which is used to establish
non-collapsing results which in turn are essential in proving the
existence of geometric limits.

The definitions and the analysis of the reduced length function and
the reduced volume as well as the monotonicity results are valid in
the context of the generalized Ricci flow. The only twist to be
aware of is that in the more general context one cannot always
extend ${\mathcal L}$-geodesics; they may run `off the edge' of
space-time. Thus, the reduced length function and reduced volume
cannot be defined globally, but only on appropriate open subsets of
a time-slice (those reachable by minimizing ${\mathcal
L}$-geodesics). But as long as one can flow an open set $U$ of a
time-slice along minimizing ${\mathcal L}$-geodesics in the
direction of decreasing $\bar\tau$, the reduced volumes of the
resulting family of open sets form a monotone non-increasing
function of $\bar\tau$. This turns out to be sufficient to extend
the non-collapsing results to Ricci flow with surgery provided that
we are careful in how we choose the parameters that go into the
definition of the surgery process.

\subsection{Application to non-collapsing results}

As we indicated in the previous paragraph, one of the main
applications of the reduced length function is to prove
non-collapsing results for three-dimensional Ricci flows with
surgery. In order to make this  argument work, one takes a weaker
notion of $\kappa$-non-collapsed\index{$\kappa$-non-collapsed} by
making a stronger curvature bound assumption: one considers
 points $(x,t)$ and constants $r$ with the property that $|{\rm Rm}|\le r^{-2}$ on
$P(x,t,r,-r^2)=B(x,t,r)\times (t-r^2,t]$. The
$\kappa$-non-collapsing condition applies to these balls and says
that ${\rm Vol}(B(x,t,r))\ge \kappa r^3.$ The basic idea in proving
non-collapsing is to use the fact that as we flow forward in time
via minimizing ${\mathcal L}$-geodesics the reduced volume is a
non-decreasing function. Hence, a lower bound of the reduced volume
of an open set at an earlier time implies the same lower bound for
the corresponding open subset at a later time. This is contrasted
with direct computations (related to the heat kernel in $\Ar^3$)
that say if the manifold is highly collapsed near $(x,t)$ (i.e.,
satisfies the curvature bound above but is not
$\kappa$-non-collapsed for some small $\kappa$) then the reduced
volume $\widetilde V_{(x,t)}$ is small at times close to $t$. Thus,
to show that the manifold is non-collapsed at $(x,t)$ we need only
find an open subset at an earlier time that is reachable by
minimizing ${\mathcal L}$-geodesics and that has a reduced volume
bounded away from zero.

One case where it is easy to do this is when we have a Ricci flow of
compact manifolds or of complete manifolds of non-negative
curvature. Hence, these manifolds are non-collapsed at all points
with a  non-collapsing constant that depends only on the geometry of
the initial metric of the Ricci flow. Non-collapsing results are
crucial and are used repeatedly in dealing with Ricci flows with
surgery in Chapters 10 -- 17, for these give one of the two
conditions required in order to take geometric limits.

\subsection{Application to ancient $\kappa$-non-collapsed solutions}

There is another important application of the length function, which
is to the study of non-collapsed, ancient solutions\index{ancient
solution} in dimension three. In the case that the generalized Ricci
flow is an ordinary Ricci flow either on a compact manifold or on a
complete manifold (with bounded curvatures) one can say much more
about the reduced length function and the reduced volume. Fix a
point $(x_0,t_0)$ in space-time. First of all, one shows that every
point $(x,t)$ with $t<t_0$ is reachable by a minimizing ${\mathcal
L}$-geodesic and thus that the reduced length is defined as a
function on all points of space at all times $t<t_0$. It turns out
to be a locally Lipschitz function in both space and time and hence
its gradient and its time derivative exist as $L^2$-functions and
satisfy important differential inequalities in the weak sense.

\begin{figure}[ht]
 {\relabelbox
  \centerline{\epsfbox{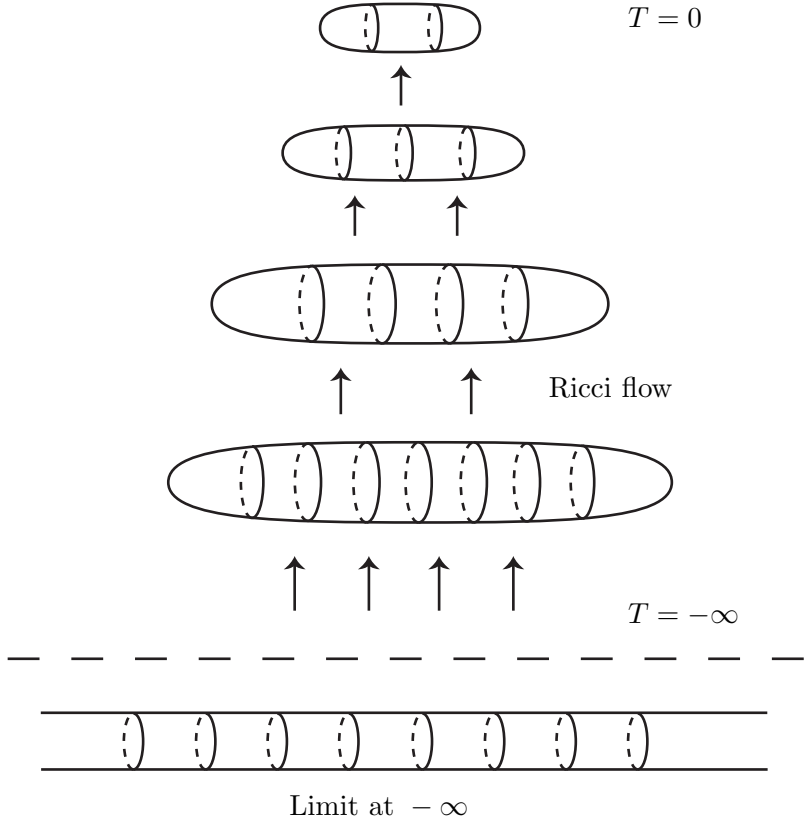}}
   \relabel{1}{Ricci flow}
  \relabel{2}{$T=0$}
  \relabel{3}{$T=-\infty$}
    \relabel{6}{$\text{Limit at }-\infty$}
  \endrelabelbox}
\caption{The asymptotic Soliton.}\label{fig:asymp}
\end{figure}

These results apply to a class of Ricci flows called
$\kappa$-solutions\index{$\kappa$-solution}, where $\kappa$ is a
positive constant. By definition a $\kappa$-solution is a Ricci flow
defined for all $t\in (-\infty,0]$, each time-slice is a non-flat,
complete $3$-manifold of non-negative, bounded curvature and each
time-slice is $\kappa$-non-collapsed. The differential inequalities
for the reduced length from any point $(x,0)$ imply that, for any
$t<0$, the minimum value of $\ell_{(x,0)}(y,t)$ for all $y\in M$ is
at most $3/2$. Furthermore, again using the differential
inequalities for the reduced length function, one shows that for any
sequence $t_n\rightarrow -\infty$, and any points $(y_n,t_n)$ at
which the reduced length function is bounded above by $3/2$, there
is a subsequence of based Riemannian manifolds,
$(M,\frac{1}{|t_n|}g(t_n),y_n)$, with a geometric limit, and this
limit is a gradient shrinking soliton. This gradient shrinking
soliton is called an {\em asymptotic
soliton}\index{$\kappa$-solution!asymptotic
soliton}\index{asymptotic soliton|see{$\kappa$-solution}} for the
original $\kappa$-solution, see {\sc Fig.}~\ref{fig:asymp}.

The point is that there are only two types of gradient shrinking
solitons in dimension three -- (i) those finitely covered by a
family of shrinking round $3$-spheres and (ii) those finitely
covered by a family of shrinking round cylinders $S^2\times \Ar$. If
a $\kappa$-solution has a gradient shrinking soliton of the first
type then it is in fact isomorphic to its gradient shrinking
soliton. More interesting is the case when the $\kappa$-solution has
a gradient shrinking soliton which is of the second type. If the
$\kappa$-solution does not have strictly positive curvature, then it
is isomorphic to its gradient shrinking soliton. Furthermore, there
is a constant $C_1<\infty$ depending on $\epsilon$ (which remember
is taken sufficiently small) such that a $\kappa$-solution of
strictly positive curvature is either a $C_1$-component, or is a
union of cores of $(C_1,\epsilon)$-caps and points that are the
center points of $\epsilon$-necks.

In order to prove the above results (for example the uniformity of $C_1$ as
above over all $\kappa$-solutions)  one needs the following result:
\begin{thm}
The space of based $\kappa$-solutions,  based at points $(x,0)$ with
$R(x,0)=1$, is compact.
\end{thm}

This result does not generalize to ancient solutions that are not
non-collapsed because, in order to prove compactness, one has to
take limits of subsequences, and in doing this the non-collapsing
hypothesis is essential. See Hamilton's work \cite{Hamiltonsurvey}
for more on general ancient solutions (i.e., those that are not
necessarily non-collapsed).

Since $\epsilon>0$ is sufficiently small so that all the results from the
appendix about manifolds covered by $\epsilon$-necks and $\epsilon$-caps hold,
the above results about gradient shrinking solitons lead to a rough qualitative
description of all $\kappa$-solutions. There are those which do not have
strictly positive curvature. These are gradient shrinking solitons, either an
evolving family of round $2$-spheres times $\Ar$ or the quotient of this family
by an involution. Non-compact $\kappa$-solutions of strictly positive curvature
are diffeomorphic to $\Ar^3$ and are the union of an $\epsilon$-tube and a core
of a $(C_1,\epsilon)$-cap. The compact ones of strictly positive curvature are
of two types. The first type are positive, constant curvature shrinking
solitons. Solutions of the second type are diffeomorphic to either $S^3$ or
$\Ar P^3$. Each time-slice of a $\kappa$-solution of the second type is either
of uniformly bounded geometry (curvature, diameter, and volume) when rescaled
so that the scalar curvature at a point is one, or admits an $\epsilon$-tube
whose complement is either a disjoint union of the cores of two
$(C_1,\epsilon)$-caps.

This gives a rough qualitative understanding of  $\kappa$-solutions. Either
they are round, or they are finitely covered by the product of a round surface
and a line, or they are a union of $\epsilon$-tubes and cores of
$(C_1,\epsilon)$-caps , or they are diffeomorphic to $S^3$ or $\Ar P^3$ and
have bounded geometry (again after rescaling so that there is a point of scalar
curvature $1$). This is the source of canonical neighborhoods for Ricci flows:
 the point is that this qualitative result remains true for any point $x$ in a
Ricci flow that has an appropriate size neighborhood within
$\epsilon$ in the $C^{[1/\epsilon]}$-topology of a neighborhood in a
$\kappa$-solution. For example, if we have a sequence of based
generalized flows $({\mathcal M}_n,G_n,x_n)$ converging to a based
$\kappa$-solution, then for all $n$ sufficiently large $x$ will have
a canonical neighborhood, one that is either an $\epsilon$-neck
centered at that point, a $(C_1,\epsilon)$-cap whose core contains
the point, a $C_1$-component, or an $\epsilon$-round component.

\subsection{Bounded curvature at bounded distance}

Perelman's other major breakthrough is his result establishing
bounded curvature at bounded distance for blow-up limits of
generalized Ricci flows. As we have alluded to several times, many
steps in the argument require taking (smooth) geometric limits of a
sequence of based generalized flows about points of curvature
tending to infinity.  To study  such a sequence we rescale each term
in the sequence so that its curvature at the base point becomes one.
Nevertheless, in taking such limits we face the problem that even
though the curvature at the point we are focusing on (the points we
take as base points) was originally large and has been  rescaled to
be one, there may be other points in the same time-slice of much
larger curvature, which, even after the rescalings, can tend to
infinity.
 If
these points are  at uniformly bounded (rescaled) distance from the
base points, then they would preclude the existence of a smooth
geometric limit of the based, rescaled flows. In his arguments,
Hamilton avoided this problem by always focusing on points of
maximal curvature (or almost maximal curvature). That method will
not work in this case. The way to deal with this possible problem is
to show that a generalized Ricci flow satisfying appropriate
conditions satisfies the following. For each $A<\infty$ there are
constants $Q_0=Q_0(A)<\infty$ and $Q(A)<\infty$ such that any point
$x$ in such a generalized flow for which the scalar curvature
$R(x)\ge Q_0$ and for any $y$ in the same time-slice as $x$ with
$d(x,y)<AR(x)^{-1/2}$ satisfies $R(y)/R(x)<Q(A)$. As we shall see,
this and the non-collapsing result are the fundamental tools that
allow Perelman to study neighborhoods of points of sufficiently
large curvature by taking smooth limits of rescaled flows, so
essential in studying the prolongation of Ricci flows with surgery.

The basic idea in proving this result is to assume the contrary and
take an incomplete geometric limit of the rescaled flows based at
the counterexample points. The existence of points at bounded
distance with unbounded, rescaled curvature means that there is a
point at infinity at finite distance from the base point where the
curvature blows up. A neighborhood of this point at infinity is
cone-like in a manifold of non-negative curvature. This contradicts
Hamilton's maximum principle result (5) in Chapter~\ref{sectmax})
that the result of a Ricci flow of manifolds of non-negative
curvature is never an open subset of a cone. (We know that any
`blow-up limit' like this has non-negative curvature because of the
curvature pinching result.) This contradiction establishes the
result.

\section{The standard solution and the surgery process}

Now we are ready to discuss three-dimensional Ricci flows with surgery.

\subsection{The standard solution}

\begin{figure}[ht]
  \centerline{\epsfbox{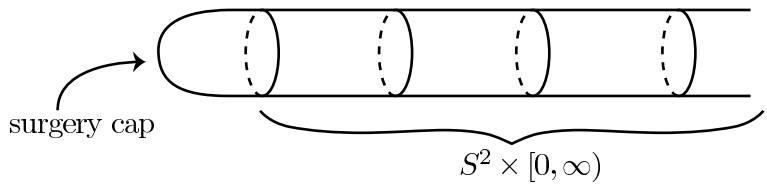}}
  \caption{The standard metric.}\label{fig:stsoln}
\end{figure}

In preparing the way for defining the surgery process, we must
construct a metric on the $3$-ball that we shall glue in when we
perform surgery. This we do in Chapter~\ref{stdsolnsect}. We fix a
non-negatively curved, rotationally symmetric metric on $\Ar^3$ that
is isometric near infinity to $S^2\times [0,\infty)$ where the
metric on $S^2$ is the round metric of scalar curvature $1$, and
outside this region has positive sectional curvature, see {\sc
Fig.}~\ref{fig:stsoln}. Any such metric will suffice for the gluing
process, and we fix one and call it the {\em standard
metric}\index{Ricci flow!standard metric}\index{standard
metric|see{Ricci flow}}. It is important to understand Ricci flow
with the standard metric as initial metric.  Because of the special
nature of this metric (the rotational symmetry and the asymptotic
nature at infinity), it is fairly elementary to show that there is a
unique solution of bounded curvature on each time-slice to the Ricci
flow equation with the standard metric as the initial metric; this
flow is defined for $0\le t<1$; and for any $T<1$ outside of a
compact subset $X(T)$ the restriction of the flow to $[0,T]$ is
close to the evolving round cylinder. Using the length function, one
shows that the Ricci flow is non-collapsed, and that the bounded
curvature and bounded distance result applies to it.
 This allows one to prove that every point $(x,t)$ in this flow has one of the
following types of neighborhoods:
\begin{enumerate}
\item $(x,t)$ is contained in the core of a $(C_2,\epsilon)$-cap,
where $C_2<\infty$ is a given universal constant depending only on $\epsilon$.
\item $(x,t)$ is the center of a strong $\epsilon$-neck.
\item $(x,t)$ is the center of an evolving $\epsilon$-neck whose
initial slice is at time zero.
\end{enumerate}

These form the second source of models for canonical neighborhoods in a Ricci
flow with surgery. Thus, we shall set $C=C(\epsilon)={\rm
max}(C_1(\epsilon),C_2(\epsilon))$ and we shall find $(C,\epsilon)$-canonical
neighborhoods in Ricci flows with surgery.

\subsection{Ricci flows with surgery}

\begin{figure}[ht]
  \centerline{\epsfbox{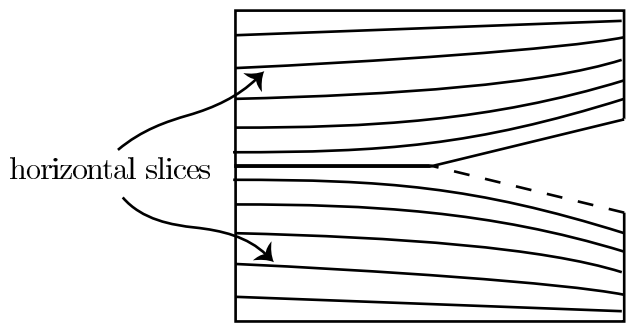}}
  \caption{Model for singularities in
  space-time.}\label{fig:spacetime}
\end{figure}

Now it is time to introduce the notion of a Ricci flow with
surgery\index{Ricci flow!with surgery}. To do this we formulate an
appropriate notion of $4$-dimensional space-time that allows for the
surgery operations. We define {\em space-time} to be a
$4$-dimensional Hausdorff singular space with a time function ${\bf
t}$ with the property that each time-slice is a compact, smooth
$3$-manifold, but level sets at different times are not necessarily
diffeomorphic. Generically space-time is a smooth $4$-manifold, but
there are {\em exposed regions}\index{exposed region} at a discrete
set of times. Near a point in the exposed region  space-time is a
$4$-manifold with boundary. The singular points of space-time are
the boundaries of the exposed regions. Near these, space-time is
modeled on the product of $\Ar^2$ with the square $(-1,1)\times
(-1,1)$, the latter having a topology in which the open sets are, in
addition to the usual open sets, open subsets of $(0,1)\times
[0,1)$, see {\sc Fig.}~\ref{fig:spacetime}. There is a natural
notion of smooth functions on space-time. These are smooth in the
usual sense on the open subset of non-singular points. Near the
singular points, and in the local coordinates described above, they
are required to be pull-backs from smooth functions on $\Ar^2\times
(-1,1)\times (-1,1)$ under the natural map. Space-time\index{Ricci
flow!space-time of} is equipped with a smooth vector field $\chi$
with $\chi({\bf t})=1$.

A Ricci flow with surgery is a smooth horizontal metric $G$ on a
space-time  with the property that the restriction of $G$, ${\bf t}$
and $\chi$ to the open subset of smooth points forms a generalized
Ricci flow.  We call this the {\em associated generalized Ricci
flow} for the Ricci flow with surgery.

\subsection{The inductive conditions necessary for doing surgery}

With all this preliminary work out of the way, we are ready to show
that one can construct Ricci flow with surgery which is precisely
controlled both topologically and metrically. This result is proved
inductively, one interval of time after another, and it is important
to keep track of various properties as we go along to ensure that we
can continue to do surgery. Here we discuss the conditions we verify
at each step.

Fix $\epsilon>0$ sufficiently small and let $C={\rm
max}(C_1,C_2)<\infty$, where $C_1$ is the constant associated to
$\epsilon$ for $\kappa$-solutions and $C_2$ is the constant
associated to $\epsilon$ for the standard solution. We say that a
point $x$ in a generalized Ricci flow has a $(C,\epsilon)$-canonical
neighborhood if one of the following holds:
\begin{enumerate}
\item $x$ is contained in a connected component of a
time-slice that is a $C$-component.
\item $x$ is contained in a connected component of its time-slice
that is within $\epsilon$ of round in the
$C^{[1/\epsilon]}$-topology.
\item $x$ is contained in the core of a $(C,\epsilon)$-cap.
\item $x$ is the center of a strong $\epsilon$-neck.
\end{enumerate}

We shall study  Ricci flows with surgery defined for $0\le
t<T<\infty$ whose associated generalized Ricci flows satisfy the
following properties:
\begin{enumerate}
\item The initial metric is normalized, meaning that for the
metric at time zero the norm the Riemann curvature is bounded above
by one and the volume of any ball of radius one is at least half the
volume of the unit ball in Euclidean space.
\item The curvature of the flow is pinched toward positive.
\item There is $\kappa>0$ so that the associated generalized Ricci flow is
$\kappa$-non-collapsed on scales at most $\epsilon$, in the sense
that we require only that balls of radius $r\le \epsilon$ be
$\kappa$-non-collapsed.
\item There is $r_0>0$ such that any point of space-time at which
the scalar curvature is $\ge r_0^{-2}$ has an
$(C,\epsilon)$-canonical neighborhood.
\end{enumerate}

The main result is that, having a Ricci flow with surgery defined on
some time interval satisfying these conditions, it is possible to
extend it to a longer time interval in such a way that it still
satisfies the same conditions, possibly allowing the constants
$\kappa$ and $r_0$ defining these conditions to get closer to zero,
but keeping them bounded away from $0$ on each compact time
interval. We repeat this construction inductively, and, since it is
easy to see that on any compact time interval there can only be a
bounded number of surgeries. In the end, we create a Ricci flow with
surgery defined for all positive time. As far as we know, it may be
the case that in the entire flow defined all the way to infinity
there are infinitely many surgeries.

\subsection{Surgery}

Let us describe how we extend a Ricci flow with surgery satisfying
all the conditions listed above and becoming singular at time
$T<\infty$. Fix $T^-<T$ so that there are no surgery  times in the
interval $[T^-,T)$. Then we can use the Ricci flow to identify all
the time-slices $M_t$ for $t\in [T^-,T)$, and hence view this part
of the Ricci flow with surgery as an ordinary Ricci flow. Because of
the canonical neighborhood assumption, there is an open subset
$\Omega\subset M_{T^-}$ on which the curvature stays bounded as
$t\rightarrow T$. Hence, by Shi's results, there is a limiting
metric at time $T$ on $\Omega$. Furthermore, the scalar curvature is
a proper function, bounded below, from $\Omega$ to $\Ar$, and each
end of $\Omega$ is an $\epsilon$-tube where the cross-sectional area
of the $2$-spheres goes to zero as we go to the end of tube. We call
such tubes {\em $\epsilon$-horns}. We are interested in
$\epsilon$-horns whose boundary is contained in the part of $\Omega$
where the scalar curvature is bounded above by some fixed finite
constant $\rho^{-2}$. We call this region $\Omega_\rho$. Using the
bounded curvature at bounded distance result and using the
non-collapsing hypothesis, one shows that given any $\delta>0$ there
is $h=h(\delta,\rho,r_0)$ such that for any $\epsilon$-horn
${\mathcal H}$ whose boundary lies in $\Omega_\rho$ and for any
$x\in{\mathcal H}$ with $R(x)\ge h^{-2}$, the point $x$ is the
center of a strong $\delta$-neck.

Now we are ready to describe the surgery procedure. It depends on
our choice of standard solution on $\Ar ^3$ and on a choice of
$\delta>0$ sufficiently small. For each $\epsilon$-horn in $\Omega$
whose boundary is contained in $\Omega_\rho$ fix a point of
curvature $(h(\delta, \rho,r_0))^{-2}$ and fix a strong
$\delta$-neck centered at this point. Then we cut the
$\epsilon$-horn open along the central $2$-sphere $S$ of this neck
and remove the end of the $\epsilon$-horn that is cut off by $S$.
Then we glue in a ball of a fixed radius around the tip from the
standard solution, after scaling the metric on this ball by
$(h(\delta,\rho,r_0))^2$. To glue these two metrics together we must
use a partition of unity near the $2$-spheres that are matched.
There is also a delicate point that we first bend in the metrics
slightly so as to achieve positive curvature near where we are
gluing. This is an idea due to Hamilton, and it is needed in order
to show that the condition of curvature pinching toward positive  is
preserved. In addition, we remove all components of $\Omega$ that do
not contain any points of $\Omega_\rho$.

\begin{figure}[ht]
  \centerline{\epsfbox{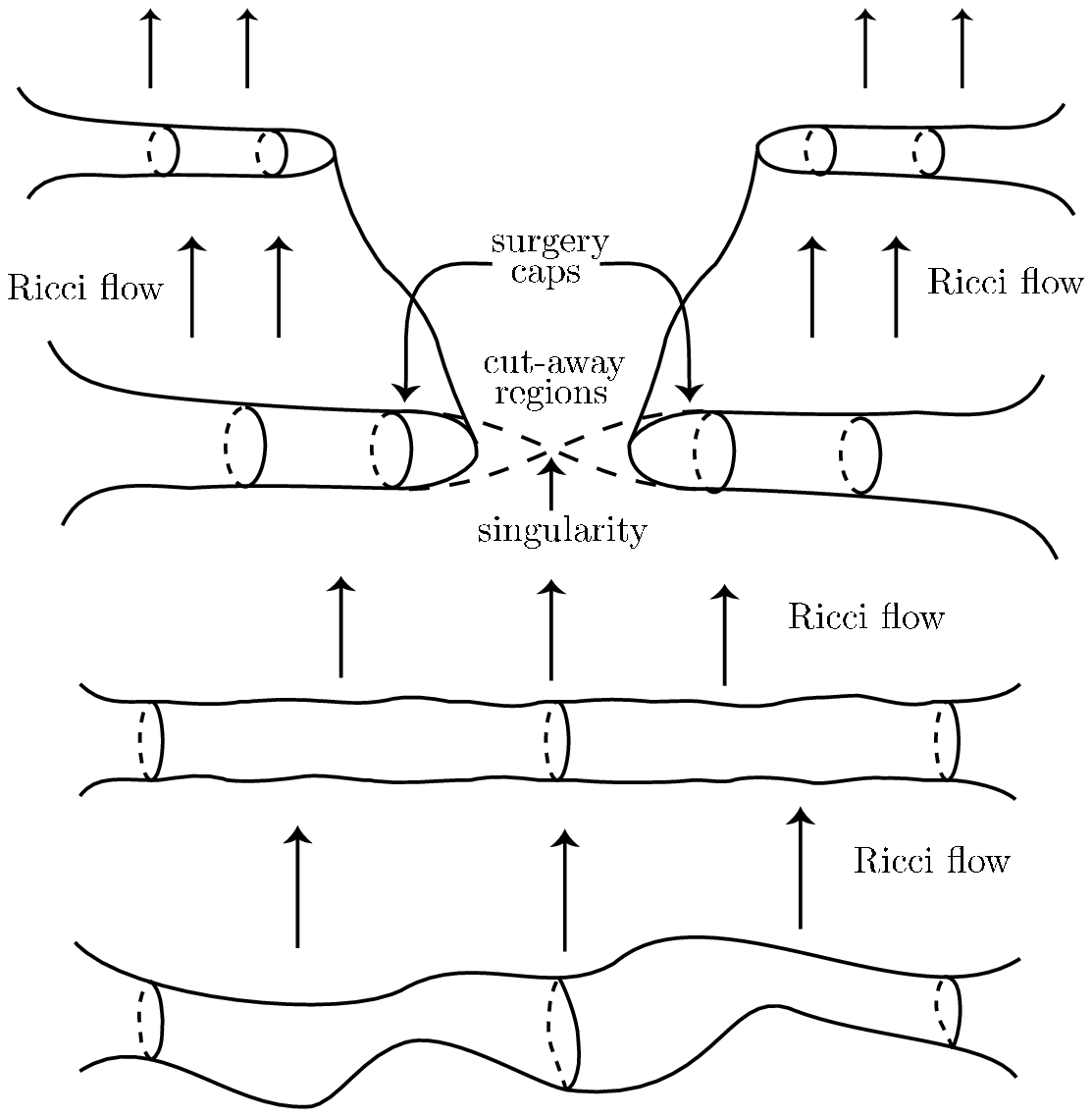}}
  \caption{Surgery.}\label{fig:rfwsurgery}
\end{figure}

This operation produces a new compact $3$-manifold.  One continues
the Ricci flow with surgery by letting this Riemannian manifold at
time $T$ evolve under the Ricci flow. See {\sc
Fig.}~\ref{fig:rfwsurgery}.

\subsection{Topological effect of surgery}

Looking\index{Ricci flow!with surgery!topological effect} at the
situation just before the surgery time, we see a finite number of
disjoint submanifolds, each diffeomorphic to either $S^2\times I$ or
the $3$-ball, where the curvature is large. In addition there may be
entire components of where the scalar curvature is large. The effect
of $2$-sphere surgery is to do a finite number of ordinary
topological surgeries along $2$-spheres in the $S^2\times I$. This
simply effects a partial connected-sum decomposition and may
introduce new components diffeomorphic to $S^3$. We also remove
entire components, but these are covered by $\epsilon$-necks and
$\epsilon$-caps  so that they have standard topology (each one is
diffeomorphic to $S^3$, $\Ar P^3$, $\Ar P^3\# \Ar P^3$, $S^2\times
S^1$, or the non-orientable $2$-sphere bundle over $S^1$). Also, we
remove $C$-components and $\epsilon$-round components (each of these
is either diffeomorphic to $S^3$ or $\Ar P^3$ or admits a metric of
constant positive curvature). Thus, the topological effect of
surgery is to do a finite number of ordinary $2$-sphere topological
surgeries and to remove a finite number of topologically standard
components.

\section{Extending Ricci flows with surgery}

Now we come to the crux of the argument. We must show that if we
have a Ricci flow with surgery defined for some time $0\le
t<T<\infty$ satisfying the four conditions: normalized initial
metric, curvature pinched toward positive, all points of scalar
curvature $\ge r^{-2}$ have canonical neighborhoods, and the flow is
$\kappa$-non-collapsed on scales $\le \epsilon$; then it is possible
to extend to a Ricci flow with surgery defined past $T$ to a larger
time keeping all these conditions satisfied (possibly with different
constants $r_0'<r_0$ and $\kappa'<\kappa$). In order to do this we
need to choose the surgery parameter $\delta>0$ sufficiently small.
There is also the issue of whether the surgery times can accumulate.

Of course, the initial metric does not change as we extend surgery
so that the condition that the normalized initial metric is clearly
preserved as we extend surgery. As we have already remarked,
Hamilton had proved earlier that one can do surgery in such a way as
to preserve the condition that the curvature is pinched toward
positive. The other two conditions require more work, and, as we
indicated above, the constants may decay to zero as we extend the
Ricci flow with surgery.

If we have all the conditions for the Ricci flow with surgery up to time $T$,
then the analysis of the open subset on which the curvature remains bounded
holds, and given $\delta>0$ sufficiently small, we do surgery on the central
$S^2$ of a strong $\delta$-neck in each $\epsilon$-horn meeting $\Omega_\rho$.
In addition we remove entirely all components that do not contain points of
$\Omega_\rho$. We then glue in the cap from the standard solution. This gives
us a new compact $3$-manifold and we restart the flow from this manifold.

The $\kappa$-non-collapsed result is extended to the new part of the
Ricci flow with surgery using the fact that it holds at times
previous to $T$. To establish this extension one uses ${\mathcal
L}$-geodesics in the associated generalized Ricci flow and reduced
volume as indicated before. In order to get this argument to work,
one must require $\delta>0$ to be sufficiently small; how small is
determined by $r_0$.

The other thing that we must establish is the existence of canonical
neighborhoods for all points of sufficiently large scalar curvature.
Here the argument is by contradiction. We consider all Ricci flows
with surgery that satisfy all four conditions on $[0,T)$ and we
suppose that we can find a sequence of such containing points
(automatically at times $T'>T$) of arbitrarily large curvature where
there are not canonical neighborhoods. In fact, we take the points
at the first such time violating this condition. We base our flows
at these points. Now we consider rescaled versions of the
generalized flows so that the curvature at these base points is
rescaled to one. We are in a position to apply the bounded curvature
and bounded distance results to this sequence, and of course the
$\kappa$-non-collapsing results which have already been established.
There are two possibilities. The first is that the rescaled sequence
converges to an ancient solution. This ancient solution has
non-negative curvature by the pinching hypothesis. General results
about three-manifolds  of non-negative curvature imply that it also
has bounded curvature. It is $\kappa$-non-collapsed. Thus, in this
case the limit is a $\kappa$-solution. This produces the required
canonical neighborhoods for the base points of the tail of the
sequence modeled on the canonical neighborhoods of points in a
$\kappa$-solution. This contradicts the assumption that none of
these points has a canonical neighborhood.

The other possibility is that one can take a partial smooth limit
but that this limit does not extend all the way back to $-\infty$.
The only way this can happen is if there are surgery caps that
prevent extending the limit back to $-\infty$. This means that the
base points in our sequence are all within a fixed distance and time
(after the rescaling) of a surgery region. But in this case results
from the nature of the standard solution show that if we have taken
$\delta>0$ sufficiently small, then the base points have canonical
neighborhoods modeled on the canonical neighborhoods in the standard
solution, again contradicting our assumption that none of the base
points has a canonical neighborhood. In order to show that our base
points have neighborhoods near those of the standard solution, one
appeals to a geometric limit argument as $\delta\rightarrow 0$. This
argument uses the uniqueness of the Ricci flow for the standard
solution. (Actually, Bruce Kleiner pointed out to us that one only
needs a compactness result for the space of all Ricci flows with the
standard metric as initial metric, not uniqueness, and the
compactness result can be proved by the same arguments that prove
the compactness of the space of $\kappa$-solutions.)

Interestingly enough, in order to establish the uniqueness of the
Ricci flow for the standard solution, as well as to prove that this
flow is defined for time $[0,1)$ and to prove that at infinity it is
asymptotic to an evolving cylinder requires the same results --
non-collapsing and the bounded curvature at bounded distance that we
invoked above. For this reason, we order the material described here
as follows. First, we introduce  generalized Ricci flows, and then
introduce the length function in this context and establish the
basic monotonicity results. Then we have a chapter on stronger
results for the length function in the case of complete manifolds
with bounded curvature. At this point we are in a position to prove
the needed results about the Ricci flow from the standard solution.
Then we are ready to define the surgery process and prove the
inductive non-collapsing results and the existence of canonical
neighborhoods.

In this way, one establishes the existence of canonical neighborhoods. Hence,
one can continue to do surgery, producing a Ricci flow with surgery defined for
all positive time.  Since these arguments are inductive, it turns out that the
constants in the non-collapsing and in the canonical neighborhood statements
decay in a predetermined rate as time goes to infinity.

Lastly, there is the issue of ruling out the possibility that the surgery times
accumulate. The idea here is very simple: Under Ricci flow during an elapsed
time $T$, volume increases at most by a  multiplicative factor which is a fixed
exponential of the time $T$. Under each surgery there is a removal of at least
a fixed positive amount of volume depending on the surgery scale $h$, which in
turns depends on $\delta$ and $r_0$. Since both $\delta$ and $r_0$ are bounded
away from zero on each finite interval, there can be at most finitely many
surgeries in each finite interval. Notice that this argument allows for the
possibility that in the entire flow all the way to infinity there are
infinitely many surgeries. It is still unknown whether that possibility ever
happens.

This completes our outline of the proof of Theorem~\ref{surgery}.

\section{Finite-time extinction}

The last topic we discuss is the  proof of the finite-time
extinction\index{finite-time extinction} for Ricci flows with
initial metrics satisfying the hypothesis  of
Theorem~\ref{finiteext}.

As we present it, the finite extinction result has two steps. The
first step is to show that there is $T<\infty$ (depending on the
initial metric) such that for all $t\ge T$, all connected components
of the $t$-time-slice $M_t$ have trivial $\pi_2$.  First, an easy
topological argument shows that only finitely many of the $2$-sphere
surgeries in a Ricci flow with surgery can be along homotopically
non-trivial $2$-spheres. Thus, after some time $T_0$ all $2$-sphere
surgeries are along homotopically trivial $2$-spheres. Such a
surgery does not affect $\pi_2$. Thus, after time $T_0$, the only
way that $\pi_2$ can change is by removal of components with
non-trivial $\pi_2$. (An examination of the topological types of
components that are removed shows that there are only two types of
such components with non-trivial $\pi_2$: $2$-sphere bundles over
$S^1$ and $\Ar P^3\#\Ar P^3$.) We suppose that at every $t\ge T_0$
there is a component of $M_t$ with non-trivial $\pi_2$. Then we can
find a connected open subset ${\mathcal X}$ of ${\bf
t}^{-1}([T_0,\infty))$ with the property that for each $t\ge  T_0$
the intersection ${\mathcal X}(t)= {\mathcal X}\cap M_t$ is a
component of $M_t$ with non-trivial $\pi_2$. We define a function
$W_2\colon [T_0,\infty)\to \Ar$ associated with such an ${\mathcal
X}$. The value $W_2(t)$ is the minimal area of all homotopically
non-trivial $2$-spheres mapping into ${\mathcal X}(t)$. This minimal
area $W_2(t)$ is realized by a harmonic map of $S^2$ into ${\mathcal
X}(t)$. The function $W_2$ varies continuously under Ricci flow and
at a surgery is lower semi-continuous. Furthermore, using an idea
that goes back to Hamilton (who applied it to minimal disks) one
shows that the forward difference quotient of the minimal area
satisfies
$$\frac{dW_2(t)}{dt}\le -4\pi +\frac{3}{(4t+1)}W_2(t).$$
(Here, the explicit form of the bound for the forward difference
quotient depends on the way we have chosen to normalize initial
metric and also on Hamilton's curvature pinching result.)

But any function $W_2(t)$ with these properties and defined for all
$t>T_0$,  becomes negative at some finite $T_1$ (depending on the
initial value). This is absurd since $W_2(t)$ is the minimum of
positive quantities. This contradiction shows that such a path of
components with non-trivial $\pi_2$ cannot exist for all $t\ge T_0$.
In fact, it even gives a computable upper bound on how long such a
component ${\mathcal X}$, with every time-slice having non-trivial
$\pi_2$, can exist in terms of the minimal area of a homotopically
non-trivial $2$-sphere mapping into ${\mathcal X}(T_0)$. It follows
that there is $T<\infty$ with the property that every component of
$M_T$ has trivial $\pi_2$. This condition then persists for all
$t\ge T$.

Three remarks are in order. This argument showing that eventually
every component of the time-slice $t$ has trivial $\pi_2$ is not
necessary for the topological application (Theorem~\ref{finiteext}),
or indeed, for any other topological application. The reason is the
sphere theorem (see \cite{Hempel}), which says that if $\pi_2(M)$ is
non-trivial then either $M$ is diffeomorphic to an $S^2$ bundle over
$S^1$ or $M$ has a non-trivial connected sum decomposition. Thus, we
can establish results for all $3$-manifolds if we can establish them
for $3$-manifolds with $\pi_2=0$. Secondly, the reason for giving
this argument is that it is pleasing to see Ricci flow with surgery
implementing the connected sum decomposition required for
geometrization of $3$-manifolds. Also, this argument is a simpler
version of the one that we use to deal with components with
non-trivial $\pi_3$. Lastly, these results on Ricci flow do not use
the sphere theorem so that establishing the cutting into pieces with
trivial $\pi_2$ allows us to give a different proof of this result
(though admittedly one using much deeper ideas).

Let us now fix $T<\infty$ such that for all  $t\ge T$ all the
time-slices $M_t$ have trivial $\pi_2$. There is a simple
topological consequence of this and our assumption on the initial
manifold. If $M$ is a compact $3$-manifold whose fundamental group
is either a non-trivial free product or an infinite cyclic group,
then $M$ admits a homotopically non-trivial embedded $2$-sphere.
Since we began with a manifold $M_0$ whose fundamental group is a
free product of finite groups and infinite cyclic groups, it follows
that for $t\ge T$ every component of $M_t$ has finite fundamental
group.  Fix $t\ge T$. Then  each component of $M_t$ has a finite
cover that is simply connected, and thus, by an elementary argument
in algebraic topology, each component of $M_t$ has non-trivial
$\pi_3$. The second step in the finite-time extinction argument is
to use a non-trivial element in this group analogously to the way we
used homotopically non-trivial $2$-spheres to show that eventually
the manifolds have trivial $\pi_2$.

 There are two
approaches to this second step: the first is due to Perelman in
\cite{Perelman3} and the other due to Colding-Minicozzi in
\cite{ColdingMinicozzi}. In their approach Colding-Minicozzi
associate to a non-trivial element in $\pi_3(M)$ a non-trivial
element in $\pi_1({\rm Maps}(S^2,M))$. This element is represented
by  a one-parameter family of $2$-spheres (starting and ending at
the constant map) representing a non-trivial element $\xi\in
\pi_3(M_0)$. They define the {\em width} of this homotopy class by
$W(\xi,t)$ by associating to each representative the maximal energy
of the $2$-spheres in the family and then minimizing over all
representatives of the homotopy class. Using results of Jost
\cite{Jost}, they show that this function satisfies the same forward
difference inequality that $W_2$ satisfies (and has the same
continuity property under Ricci flow and the same semi-continuity
under surgery). Since $W(\xi,t)$ is always $\ge 0$ if it is defined,
this forward difference quotient inequality implies that the
manifolds $M_t$ must eventually become empty.

While this approach seemed completely natural to us, and while we
believe that it works, we found the technical details
daunting\footnote{Colding and Minicozzi tell us they plan to give an
expanded version of their argument with a more detailed proof.}
(because one is forced to consider index-one critical points of the
energy functional rather than minima). For this reason we chose to
follow Perelman's approach. He represents a non-trivial element in
$\pi_3(M)$ as a non-trivial element in $\xi\in \pi_2(\Lambda M,*)$
where $\Lambda M$ is the free loop space of $M$. He then associates
to a family $\Gamma\colon S^2\to \Lambda M$ of homotopically trivial
loops an invariant $W(\Gamma)$ which is the maximum of the areas of
minimal spanning disks for the loops $\Gamma(c)$ as $c$ ranges over
$S^2$. The invariant of a non-trivial homotopy class $\xi$ is then
the infimum over all representatives $\Gamma$ for $\xi$ of
$W(\Gamma)$. As before, this function is continuous under Ricci flow
and is lower semi-continuous under surgery (unless the surgery
removes the component in question). It also satisfies a forward
difference quotient
$$\frac{dW(\xi)}{dt}\le -2\pi+\frac{3}{4t+1}W(\xi).$$
The reason for the term $-2\pi$ instead of $-4\pi$ which occurs in
the other cases is that we are working with minimal $2$-disks
instead of minimal $2$-spheres. Once this forward difference
quotient estimate (and the continuity) have been established the
argument finishes in the same way as the other cases: a function $W$
with the properties we have just established cannot be non-negative
for all positive time. This means the component in question, and
indeed all components at later time derived from it, must disappear
in finite time. Hence, under the hypothesis on the fundamental group
in Theorem~\ref{finiteext} the entire manifold must disappear at
finite time.

Because this approach uses only minima for the energy or area
functional, one does not have to deal with higher index critical
points. But one is forced to face other difficulties though --
namely boundary issues. Here, one must prescribe the deformation of
the family of boundary curves before computing the forward
difference quotient of the energy. The obvious choice is the
curve-shrinking flow (see \cite{AG}). Unfortunately, this flow can
only be defined when the curve in question is immersed and even in
this case the curve-shrinking flow can develop singularities even if
the Ricci flow does not. Following Perelman, or indeed \cite{AG},
one uses the device of taking the product with a small circle and
using loops, called {\em ramps}\index{ramp}, that go around that
circle once. In this context the curve-shrinking flow remains
regular as long as the Ricci flow does. One then projects this flow
to a flow of families of $2$-spheres in the free loop space of the
time-slices of the original Ricci flow. Taking the length of the
circle sufficiently small yields the boundary deformation needed to
establish the forward difference quotient result. This requires a
compactness result which holds under local total curvature bounds.
This compactness result holds off of a set of time of small total
measure, which is sufficient for the argument. At the very end of
the argument we need an elementary but complicated result on annuli,
which we could not find in the literature. For more details on these
points see Chapter~\ref{energy}.

\section{Acknowledgements}

In sorting out Perelman's arguments we have had the aid of many
people. First of all, during his tour of the United States in the
Spring of 2003, Perelman spent enormous amounts of time explaining
his ideas and techniques both in public lectures and in private
sessions to us and to many others. Since his return to Russia, he
has also freely responded to our questions by e-mail. Richard
Hamilton has also given unstintingly of his time, explaining to us
his results and his ideas about Perelman's results. We have
benefitted tremendously from the work of Bruce Kleiner and John
Lott. They produced a lengthy set of notes filling in the details in
Perelman's arguments \cite{KleinerLott}. We have referred to those
notes countless times as we came to grips with Perelman's ideas. In
late August and early September of 2004, Kleiner, Lott and the two
of us participated in  a workshop at Princeton University, supported
by the Clay Math Institute, going through Perelman's second paper
(the one on Ricci flow with surgery) in detail. This workshop played
a significant role in convincing us that Perelman's arguments were
complete and correct and also in convincing us to write this book.
We thank all the participants of this workshop and especially
Guo-Feng Wei, Peng Lu, Yu Ding, and X.-C. Rong  who, together with
Kleiner and Lott, made significant contributions to the workshop.
Before, during, and after this workshop, we have benefitted from
private conversations too numerous to count with Bruce Kleiner and
John Lott.

Several of the analytic points were worked out in detail by others,
and we have freely adapted their work. Rugang Ye wrote a set of
notes proving the Lipschitz properties of the length function and
proving the fact that the cut locus for the reduced length function
is of measure zero. Closely related to this, Ye gave detailed
arguments proving the requisite limiting results for the length
function required to establish the existence of a gradient shrinking
soliton. None of these points was directly addressed by Perelman,
and it seemed to us that they needed careful explanation. Also, the
proof of the uniqueness of the standard solution was established
jointly by the second author and Peng Lu, and the proof of the
refined version of Shi's theorem where one has control on a certain
number of derivatives at time zero was also shown to us by Peng Lu.
We also benefitted from Ben Chow's expertise and vast knowledge of
the theory of Ricci flow especially during the ClayMath/MSRI Summer
School on Ricci Flow in 2005. We had several very helpful
conversations with Tom Mrowka and also with Robert Bryant,
especially about annuli of small area.

The second author gave courses at Princeton University and ran
seminars on this material at MIT and Princeton. Natasa Sesum and
Xiao-Dong Wang, see \cite{STW}, wrote notes for the seminars at MIT,
and Edward Fan and Alex Subotic took notes for the seminars and
courses at Princeton, and they produced preliminary manuscripts. We
have borrowed freely from these manuscripts, and it is a pleasure to
thank each of them for their efforts in the early stages of this
project.

The authors thank all the referees for their time and effort which
helped us to improve the presentation. In particular, we wish to
thank Cliff Taubes who carefully read the entire manuscript and gave
us enumerable comments and suggestions for clarifications and
improvements. We also thank Terry Tao for much helpful feedback
which improved the exposition. We also thank Colin Rourke for his
comments on the finite-time extinction argument.

We thank Lori Lejeune for drawing the figures and John Etnyre for technical
help in inserting these figures into the manuscript.

Lastly, during this work we were both generously supported by the Clay
Mathematical Institute, and it is a pleasure to thank the Clay Mathematics
Institute for its support and to thank its staff, especially Vida Salahi, for
their help in preparing this manuscript. We also thank the National Science Foundation
for their continued support and the second author thanks the Jim Simons
for his support during the period that he was a faculty member at MIT.

\section{List of related papers}

For the readers' convenience we gather here references to all the
closely related articles.

First and foremost are Perelman's three preprints, \cite{Perelman1},
\cite{Perelman2}, and \cite{Perelman3}. The first of these
introduces the main techniques in the case of Ricci flow, the second
discusses the extension of these techniques to Ricci flow with
surgery, and the last gives the short-cut to the Poincar\'e
Conjecture and the $3$-dimensional spherical space-form conjecture,
avoiding the study of the limits as time goes to infinity and
collapsing space arguments. There are the detailed notes by Bruce
Kleiner and John Lott, \cite{KleinerLott}, which greatly expand and
clarify Perelman's arguments from the first two preprints. There is
also a note on Perelman's second paper by Yu Ding  \cite{yuding}.
There is the article by Colding-Minicozzi \cite{ColdingMinicozzi},
which gives their alternate approach to the material in Perelman's
third preprint. Collapsing space arguments needed for the full
geometrization conjecture are discussed in Shioya-Yamaguchi
\cite{SY}. Lastly, after we had submitted a preliminary version of
this manuscript for refereeing, H.-D. Cao and X.-P. Zhu published an
article on the Poincar\'e Conjecture and Thurston's Geometrization
Conjecture; see \cite{CaoZhu}.

\chapter{Preliminaries from Riemannian
geometry}\label{sectprelim}

In this chapter we will recall some basic facts in Riemannian
geometry. For more details we refer the reader to \cite{doCarmo} and
\cite{Petersen}. Throughout, we always adopt Einstein's summation
convention on repeated indices and `manifold' means a paracompact,
Hausdorff, smooth manifold.


\section{Riemannian metrics and the Levi-Civita connection}

Let $M$ be a manifold and let $p$ be a point of $M$. Then $TM$
denotes the tangent bundle of $M$ and $T_pM$ is the tangent space at
$p$. Similarly, $T^*M$ denotes the cotangent bundle of $M$ and
$T^*_pM$ is the cotangent space at $p$. For any vector bundle
${\mathcal V}$ over $M$ we denote by $\Gamma({\mathcal V})$ the
vector space of smooth sections of ${\mathcal V}$.

\begin{defn}
Let $M$ be an $n$-dimensional manifold. A \it Riemannian metric \rm
$g$ on $M$ is a smooth section of $T^*M\otimes T^*M$ defining a
positive definite symmetric bilinear form on $T_pM$ for each $p\in
M$. In local coordinates $(x^1,\cdots, x^n)$, one has a natural
local basis $\{\partial_1,\cdots,\partial_n\}$ for $TM$, where
$\partial_i=\frac{\partial}{\partial x^i}$. The metric tensor
$g=g_{ij} dx^i\otimes dx^j$ is represented by a smooth matrix-valued
function
\[ g_{ij}=g(\partial_i,\partial_j). \]
 The
pair $(M,g)$ is a {\sl Riemannian manifold}. We denote by $(g^{ij})$
the inverse of the matrix $(g_{ij})$.
\end{defn}

Using a partition of unity one can easily see that any manifold
admits a Riemannian metric. A Riemannian metric on $M$ allows us to
measure lengths of smooth paths in $M$ and hence to define a
distance function by setting $d(p,q)$ equal to the infimum of the
lengths of smooth paths from $p$ to $q$. This makes $M$ a metric
space. For a point $p$ in a Riemannian manifold $(M,g)$ and for $r>
0$ we denote the metric ball of radius $r$ centered at $p$ in $M$ by
$B(p,r)$ or by $B_g(p,r)$ if the metric needs specifying or
emphasizing. It is defined by
$$B(p,r)=\{q\in M\left|\right. d(p,q)<r\}.$$

\begin{thm}
Given a Riemannian metric $g$ on $M$, there uniquely exists a
torsion-free connection on $TM$ making $g$ parallel, i.e., there is
a unique $\mathbb{R}$-linear mapping $\nabla\colon
\Gamma(TM)\rightarrow \Gamma(T^*M\otimes TM)$ satisfying the Leibniz
formula
\[ \nabla(fX)=df\otimes X+f\nabla X, \]
and the following two additional conditions for all vector fields
$X$ and $Y$:
\begin{itemize}
\item[$\bullet$] ($g$ orthogonal) $d(g(X,Y))=g(\nabla X,Y)+g(X,\nabla
Y)$.
\item[$\bullet$] (Torsion-free) $\nabla_XY-\nabla_YX-[X,Y]=0$
(where, as is customary, we denote $\nabla Y(X)$ by $\nabla_XY$);
\end{itemize}
\end{thm}

 We call the
above connection the {\sl Levi-Civita connection}\index{Levi-Civita
connection|ii} of the metric and $\nabla X$ the {\sl covariant
derivative} of $X$. On a Riemannian manifold we always use the
Levi-Civita connection.

In local coordinates $(x^1,\ldots,x^n)$ the Levi-Civita connection
$\nabla$ is given by the $\nabla_{\partial_i}(\partial_j) =
\Gamma^{k}_{ij}\partial_k$, where the Christoffel symbols
$\Gamma^k_{ij}$ are the smooth functions
\begin{equation}\label{Gamma}
\Gamma^{k}_{ij}= \frac{1}{2}g^{kl}(\partial_i g_{l j} +\partial_j
g_{il}-\partial_l g_{ij}). \end{equation}
 Note that the above two
additional conditions for the Levi-Civita connection $\nabla$
correspond respectively to

\smallskip

$\bullet$ $\Gamma^k_{ij}=\Gamma^k_{ji}$,

$\bullet$ $\partial_k g_{ij}=g_{l j}\Gamma^l_{ki}+g_{il
}\Gamma^l_{kj}$.

The covariant derivative extends to all tensors. In the special case
of a function $f$ we have $\nabla( f)=df$. Note that there is a
possible confusion between this and  the notation in the literature
since one often sees $\nabla f$ written for the gradient of $f$,
which is the vector field dual to $df$. We always use $\nabla f$ to
mean $df$, and we will denote the gradient of $f$ by $(\nabla f)^*$,

The covariant derivative allows us to define the Hessian of a smooth
function at any point, not just a critical point. Let $f$ be a
smooth real-valued function on $M$.  We define the
Hessian\index{Hessian|ii}  of $f$, denoted ${\rm Hess}(f)$, as
follows:
\begin{equation}\label{Hessian}
{\rm Hess}(f)(X,Y)=X(Y(f))-\nabla_XY(f).\end{equation}
\begin{lem}\label{Hessformula}
 The Hessian is a contravariant, symmetric two-tensor, i.e., for
 vector fields $X$ and $Y$ we have
$${\rm Hess}(f)(X,Y)={\rm Hess}(f)(Y,X)$$
and
$${\rm Hess}(f)(\phi X,\psi Y)=\phi\psi {\rm Hess}(f)(X,Y)$$
for all smooth functions $\phi,\psi$. Other formulas for the Hessian
are
$${\rm Hess}(f)(X,Y)=\langle \nabla_X(\nabla f),Y\rangle=\nabla_X(\nabla_Y(f))=\nabla^2f(X,Y).$$
Also, in local coordinates  we have
$${\rm Hess}(f)_{ij}=\partial_i\partial_jf-(\partial_kf)\Gamma^k_{ij}.$$
\end{lem}

\begin{proof}
The proof of symmetry is direct from the torsion-free assumption:
$${\rm Hess}(f)(X,Y)-{\rm Hess}(f)(Y,X)=[X,Y](f)-(\nabla_XY-\nabla_YX)(f)=0.$$
 The fact that
${\rm Hess}(f)$ is a tensor is also established by direct
computation. The equivalence of the various formulas is also
immediate:
\begin{eqnarray}\label{hessform}
\langle\nabla_X(\nabla f),Y\rangle & = & X(\langle \nabla
f,Y\rangle)-\langle\nabla f,\nabla_XY\rangle
\\ & = & X(Y(f))-\nabla_XY(f)={\rm Hess}(f)(X,Y).\nonumber
\end{eqnarray}
Since $df=(\partial_rf) dx^r$ and
$\nabla(dx^k)=-\Gamma_{ij}^kdx^i\otimes dx^j$, it follows that
$$\nabla(df)=\left(\partial_i\partial_jf-(\partial_kf)\Gamma^k_{ij}\right)dx^i\otimes
dx^j.$$ It is direct from the definition that
$${\rm Hess}(f)_{ij}={\rm Hess}(f)(\partial_i,\partial_j)=\partial_i\partial_jf-(\partial_kf)\Gamma_{ij}^k.$$
\end{proof}

When the metric that we are using to define the Hessian is not clear
from the context, we introduce it into the notation and write ${\rm
Hess}_g(f)$ to denote the Hessian of $f$ with respect to the metric
$g$.

The Laplacian\index{Laplacian|ii} $\triangle f$ is defined as the
trace of the Hessian: That is to say, in local coordinates near $p$
we have
$$\triangle f(p)=\sum_{ij}g^{ij}{\rm Hess}(f)(\partial_i,\partial_j).$$
Thus, if $\{X_i\}$ is an orthonormal basis for $T_pM$ then
\begin{equation}\label{laplacformula}
\triangle f(p)=\sum_i{\rm Hess}(f)(X_i,X_i).\end{equation} Notice that this is
the form of the Laplacian that is non-negative at a local minimum, and
consequently has a non-positive spectrum.

\section{Curvature of a Riemannian manifold}\index{curvature}

For the rest of this chapter $(M,g)$ is a Riemannian manifold.

\begin{defn}
The \textit{Riemann curvature tensor} of $M$\index{Riemann
curvature!tensor|ii}\index{curvature!Riemann|see{Riemann curvature
tensor}} is the  $(1,3)$-tensor on $M$
\[ {\mathcal R}(X,Y)Z = \nabla^{2}_{X,Y}Z
- \nabla^{2}_{Y,X}Z = \nabla_{X}\nabla_{Y}Z-\nabla_{Y}\nabla_{X}Z
-\nabla_{[X,Y]}Z, \] where
$\nabla^2_{X,Y}Z=\nabla_{X}\nabla_{Y}Z-\nabla_{\nabla_XY}Z$.
\end{defn}

In local coordinates the curvature tensor can be represented as \[
{\mathcal R}(\partial_i,\partial_j)\partial_k=
{{{R}_{ij}}^l}_k\partial_l, \] where
$${{R_{ij}}^l}_k=
\partial_i \Gamma^{l}_{jk}- \partial_j
\Gamma^{l}_{ik}+
\Gamma^{s}_{jk}\Gamma^{l}_{is}-\Gamma^{s}_{ik}\Gamma^{l}_{js}.$$
Using the metric tensor $g$, we can change ${\mathcal R}$ to a
$(0,4)$-tensor as follows:
\[ {\mathcal R}(X,Y,Z,W)=g({\mathcal R}(X,Y)W,Z). \]
(Notice the change of order in the last two variables.) Notice that
we use the same symbol and the same name for both the $(1,3)$ tensor
and the $(0,4)$ tensor; which one we are dealing with in a given
context is indicated by the index structure or the variables to
which the tensor is applied.  In local coordinates, the Riemann
curvature tensor can be represented as
\begin{align*} {\mathcal R}(\partial_i,\partial_j,\partial_k,\partial_l) & = R_{ijkl}\\
&=g_{ks}{{R_{ij}}^s}_l\\
&=g_{ks}(\partial_i \Gamma^{s}_{jl}- \partial_j \Gamma^{s}_{il}+
\Gamma^{t}_{jl}\Gamma^{s}_{it}-\Gamma^{t}_{il}\Gamma^{s}_{jt}).
\end{align*}

One can easily verify the following:

\begin{claim}\label{Bianchi}
The Riemann curvature tensor ${\mathcal R}$ satisfies the following
properties:

$\bullet$ (Symmetry) $R_{ijkl}=-R_{jikl}$, $R_{ijkl}=-R_{ijlk}$,
$R_{ijkl}=R_{klij}$,

$\bullet$ (1st Bianchi identity)\index{Bianchi identity} The sum of
$R_{ijkl}$ over the cyclic permutation of

\quad any three indices vanishes,

$\bullet$ (2nd Bianchi identity)
$R_{ijkl,h}+R_{ijlh,k}+R_{ijhk,l}=0$, where
$$R_{ijkl,h}=(\nabla_{\partial_h}{\mathcal R})_{ijkl}.$$
\end{claim}

There are many important related curvatures.

\begin{defn}
The \textit{sectional curvature}\index{curvature!sectional|ii} of a
2-plane $P\subset T_pM$ is defined as
\[ K(P) = {\mathcal R}(X,Y,X,Y), \]
where $\{X,Y\}$ is an orthonormal basis of $P$. We say that $(M,g)$
has \textit{positive sectional curvature} (resp., \textit{negative
sectional curvature}) if $K(P)>0$ (resp., $K(P)<0$) for every
2-plane $P$.  There are analogous notions of non-negative and
non-positive sectional curvature.
\end{defn}

In local coordinates, suppose that $X = X^{i}\partial_i$ and $Y =
Y^{i}\partial_i$. Then we have
\[ K(P) = R_{ijkl}X^{i}Y^{j}X^{k}Y^{l}.\]
A Riemannian manifold is said to have {\sl constant sectional
curvature} if $K(P)$ is the same for all $p\in M$ and all $2$-planes
$P\subset T_pM$. One can show that a manifold $(M,g)$ has constant
sectional curvature $\lambda$ if and only if \[ R_{ijkl} =
\lambda(g_{ik}g_{jl}-g_{il}g_{jk}). \] Of course, the sphere of
radius $r$ in $\Ar^n$ has constant sectional curvature $1/r^2$,
$\Ar^n$ with the Euclidean metric has constant sectional curvature
$0$, and the hyperbolic space ${\mathbb H}^n$, which, in the
Poincar\'e model, is given by the unit disk with the metric
$$\frac{4(dx^{2}_{1}+\cdots+
dx^{2}_{n})}{(1-\abs{x}^{2})^{2}},$$ or in the upper half-space
model with coordinates $(x^1,\ldots,x^n)$ is given by
$$\frac{ds^2}{(x^n)^2}$$ has constant sectional curvature $-1$. In all
three cases we denote the constant curvature metric by $g_{\rm st}$.

\begin{defn}
Using the metric, one can replace the Riemann curvature tensor
${\mathcal R}$ by  a symmetric bilinear form ${\rm Rm}$ on
$\wedge^2TM$. In local coordinates let
$\varphi=\varphi^{ij}\partial_i\wedge\partial _j$ and
$\psi=\psi^{kl}\partial _k\wedge\partial_l$ be local sections of
$\wedge^2TM$. The formula for ${\rm Rm}$ is
\[ {\rm Rm}(\varphi, \psi)=R_{ijkl}\varphi^{ij}\psi^{kl}. \]
We call ${\rm Rm}$ the \textit{curvature
operator}\index{curvature!operator|ii}. We say $(M,g)$ has
\textit{positive curvature operator} if ${\rm
Rm}(\varphi,\varphi)>0$ for any nonzero 2-form $\varphi =
\varphi^{ij}\partial_{i}\wedge
\partial_{j}$ and has \textit{nonnegative curvature operator} if
${\rm Rm}(\varphi,\varphi)\geq0$ for any  $\varphi\in \wedge^2TM$.
\end{defn}

Clearly, if the curvature operator is a positive (resp.,
non-negative) operator then the manifold is positively (resp.,
non-negatively) curved.

\begin{defn}
The \textit{Ricci curvature tensor}\index{curvature!Ricci|ii},
denoted ${\rm Ric}$ or ${\rm Ric}_g$ when it is necessary to specify
the metric, is a symmetric contravariant two-tensor. In local
coordinates it is defined by
\[ {\rm Ric}(X,Y)=g^{kl}R(X,\partial_k,Y,\partial_l). \]
The value of this tensor at a point $p\in M$ is  given by
$\sum_{i=1}^nR(X(p),e_i,Y(p),e_i)$ where $\{e_{1},\cdots,e_{n}\}$ is
an orthonormal basis of $T_pM$. Clearly ${\rm Ric}$ is a symmetric
bilinear form on $TM$, given in local coordinates by
$${\rm Ric}={\rm Ric}_{ij}dx^i\otimes dx^j,$$ where
${\rm Ric}_{ij}={\rm Ric}(\partial_i,\partial _j)$. The
\textit{scalar curvature}\index{curvature!scalar|ii} is defined by:
$$R=R_g={\rm tr}_g {\rm Ric}=g^{ij}{\rm Ric}_{ij}.$$ We will say that ${\rm Ric} \geq
k$ (or $\leq k$) if all the eigenvalues of ${\rm Ric}$ are $\geq k$
(or $\leq k$).
\end{defn}

Clearly, the curvatures are natural in the sense that if $F\colon N
\rightarrow M$ is a diffeomorphism and if $g$ is a Riemannian metric
on $M$, then $F^*g$ is a Riemannian metric on $N$ and we have ${\rm
Rm}(F^*g)=F^*({\rm Rm}(g))$, ${\rm Ric}(F^{*}g) = F^{*}({\rm
Ric}(g))$, and $R(F^*g)=F^*(R(g))$.

\subsection{Consequences of the Bianchi identities}

There is one consequence of the second Bianchi identity that will be
important later. For any contravariant two-tensor $\omega$ on $M$
(such as ${\rm Ric}$ or ${\rm Hess}(f)$) we define the contravariant
one-tensor ${\rm div}(\omega)$ as follows: For any vector field $X$
we set
$${\rm div}(\omega)(X)=\nabla^*\omega(X)=g^{rs}\nabla_r(\omega)(X,\partial_s).$$

\begin{lem}\label{divRic}
$$dR=2{\rm div}({\rm Ric})=2\nabla^*{\rm Ric}.$$
\end{lem}

For a proof see Proposition 6 of Chapter 2 on page 40 of \cite{Petersen}.

We shall also need a formula relating the connection Laplacian on
contravariant one-tensors with the Ricci curvature. Recall that for
a smooth function $f$, we defined the symmetric two-tensor
$\nabla^2f$ by
$$\nabla^2f(X,Y)=\nabla_X\nabla_Y(f)-\nabla_{\nabla_X(Y)}(f)={\rm Hess}(f)(X,Y),$$
and then defined the Laplacian
$$\triangle f={\rm
tr}\nabla^2f=g^{ij}(\nabla^2f)_{ij}.$$ These operators extend to tensors of any
rank. Suppose that $\omega$ is a contravariant tensor of rank $k$. Then we
define $\nabla^2\omega$ to be a contravariant tensor of rank $k+2$ given by
$$\nabla^2\omega(\cdot,X,Y)=(\nabla_{X}\nabla_{Y}\omega)(\cdot)
-\nabla_{\nabla_X(Y)}\omega(\cdot).$$ This expression is not
symmetric in the vector fields $X,Y$ but the commutator is given by
evaluating the curvature operator ${\mathcal R}(X,Y)$ on $\omega$.
We define the connection Laplacian on the tensor $\omega$ to be
$$\triangle \omega=g^{ij}\nabla^2(\omega)(\partial_i,\partial_j).$$
Direct computation gives the standard Bochner formula relating these
Laplacians with the Ricci curvature; see for example Proposition
4.36 on page 168 of \cite{Gallot}.

\begin{lem}\label{lapformula}
 Let $f$ be a smooth function on a Riemannian manifold.
Then we have the following formula for contravariant one-tensors:
$$ \triangle df=d\triangle f+{\rm Ric}((\nabla f)^*,\cdot).$$
\end{lem}

\subsection{First examples}

The most homogeneous Riemannian manifolds are those of constant
sectional curvature. These are easy to classify; see Corollary 10 of
Chapter 5 on page 147 of \cite{Petersen}.

\begin{thm}
{\bf (Uniformization Theorem)}\index{Uniformization Theorem} If
$(M^n,g)$ is a complete, simply-connected Riemannian manifold of
constant sectional curvature $\lambda$, then:
\begin{enumerate}
\item If $\lambda=0$, then $M^n$ is isometric to Euclidean $n$-space.
\item If $\lambda>0$ there is a diffeomorphism $\phi \colon M \rightarrow S^{n}$
 such that g = $\lambda^{-1}\phi^{*}(g_{\rm
st})$ where $g_{st}$ is the usual metric on the unit sphere in $\Ar^{n+1}$.
\item If $\lambda<0$ there is a diffeomorphism $\phi \colon M \rightarrow {\mathbb H}^{n}$
 such that g = $\abs{\lambda}^{-1}\phi^{*}(g_{\rm
st})$ where $g_{st}$ is the Poincar\'e metric of constant curvature
$-1$ on ${\mathbb H}^n$.
\end{enumerate}
\end{thm}

Of course, if $(M^n,g)$ is a complete manifold of constant sectional
curvature then its universal covering satisfies the hypothesis of
the theorem and hence is one of $S^n, \Ar^n$, or ${\mathbb H}^n$, up
to a constant scale factor. This implies that $(M,g)$ is isometric
to a quotient of one of these simply connected spaces of constant
curvature by the free action of a discrete group of isometries. Such
a Riemannian manifold is called a
\textit{space-form}\index{space-form}.

\begin{defn}
 The Riemannian manifold
$(M,g)$ is said to be an \textit{Einstein manifold with Einstein
constant}\index{Einstein manifold} $\lambda$ if ${\rm
Ric}(g)=\lambda g$.
\end{defn}

\begin{exam} Let $M$ be an $n$-dimensional manifold with $n$ being either $2$ or $3$.
If $(M,g)$ is Einstein with Einstein constant $\lambda$, one can
easily show that $M$ has constant sectional curvature
$\frac{\lambda}{n-1}$, so that in fact $M$ is  a space-form.
\end{exam}

\subsection{Cones}

Another class of examples that will play an important role in our
study of the Ricci flow is that of cones.

\begin{defn}\label{conedefn}
Let $(N,g)$ be a Riemannian manifold. We define the {\em open cone}
over $(N,g)$ to be the manifold $N\times (0,\infty)$ with the metric
$\tilde g$ defined as follows: For any $(x,s)\in N\times (0,\infty)$
we have
$$\tilde g(x,s)=s^2g(x)+ds^2.$$
\end{defn}

Fix local coordinates $(x^1,\ldots,x^n)$ on $N$. Let
$\Gamma_{ij}^k;\ 1\le i,j,k\le n$, be the Christoffel symbols for
the Levi-Civita connection on $N$. Set $x^0=s$. In the local
coordinates $(x^0,x^1,\ldots,x^n)$ for the cone we have the
Christoffel symbols $\widetilde \Gamma_{ij}^k;\ 0\le i,j,k\le n$,
for $\tilde g$. The relation between the metrics gives the following
relations between the two sets of Christoffel symbols:
\begin{eqnarray*}\widetilde
\Gamma_{ij}^k & = & \Gamma_{ij}^k; \ \ \ \ 1\le i,j,k\le n \\
\widetilde \Gamma_{ij}^0 & = & -sg_{ij}; \ \ \ \ 1\le i,j\le n \\
\widetilde\Gamma_{i0}^j & = & \widetilde\Gamma_{0i}^j =
s^{-1}\delta^j_i; \ \ \ \ 1\le i,j\le n \\
\widetilde \Gamma_{i0}^0 & = & 0;\ \ \ \ 0\le i\le n \\
\widetilde \Gamma_{00}^i & = & 0;\ \  \ \ 0\le i\le n.
\end{eqnarray*}

Denote by ${\mathcal R}_g$ the curvature tensor for $g$ and by
${\mathcal R}_{\tilde g}$ the curvature tensor for $\tilde g$. Then
the above formulas lead directly to:
\begin{eqnarray*}
{\mathcal R}_{\tilde g}(\partial_i,\partial_j)(\partial_0) & = & 0;\ \ \ \  0\le i,j\le n \\
{\mathcal R}_{\tilde g}(\partial_i,\partial_j)(\partial_i) & = &
{\mathcal
R}_g(\partial_i,\partial_j)(\partial_i)+g_{ii}\partial_j-g_{ji}\partial_i;\
\ \ \  1\le i,j\le n \\
\end{eqnarray*}

This allows us to compute the Riemann curvatures of the cone in
terms of those of $N$.

\begin{prop}\label{conecurv}
Let $N$ be a Riemannian manifold of dimension $n-1$. Fix $(x,s)\in
c(N)=N\times (0,\infty)$. With respect to the coordinates
$(x^0,\ldots,x^n)$ the curvature operator ${{\rm Rm}}_{\tilde
g}(p,s)$ of the cone decomposes as
$$\begin{pmatrix} 0 & 0 \\ s^2({\rm Rm}_g(p)-\wedge^2g(p)) & 0
\end{pmatrix},$$
where $\wedge^2g(p)$ is the symmetric form on $\wedge^2T_pN$ induced
by $g$.
\end{prop}

\begin{cor}
For any $p\in N$ let $\lambda_1,\ldots,\lambda_{(n-1)(n-2)/2}$ be
the eigenvalues of ${\rm Rm}_g(p)$. Then for any $s>0$  there are
$(n-1)$ zero eigenvalues of ${\rm Rm}_{\tilde g}(p,s)$. The other
$(n-1)(n-2)/2$ eigenvalues of ${\rm Rm}_{\tilde g}(p,s)$ are
$s^{-2}(\lambda_i-1)$.
\end{cor}

\begin{proof}
Clearly from Proposition~\ref{conecurv}, we see that under the
orthogonal decomposition $\wedge^2T_{(p,s)}c(N)=\wedge^2T_pN\oplus
T_pN$ the second subspace is contained in the null space of $ {\rm
Rm}_{\tilde g}(p,s)$, and hence contributes $(n-1)$ zero
eigenvalues. Likewise, from this proposition we see that the
eigenvalues of the restriction of ${\rm Rm}_{\tilde g}(p,s)$ to the
subspace $\wedge^2T_pN$ are  given by
$s^{-4}(s^2(\lambda_i-1))=s^{-2}(\lambda_i-1)$.
\end{proof}

\section{Geodesics and the exponential map}

Here we review standard material about geodesics, Jacobi fields, and
the exponential map.

\subsection{Geodesics and the energy functional}

\begin{defn} Let $I$ be an open interval.
A smooth curve $\gamma\colon I \rightarrow M$ is called a
\textit{geodesic}\index{geodesic} if $\nabla_{\dot \gamma}\dot
\gamma =0$.
\end{defn}

In local coordinates, we write $\gamma(t)=(x^1(t),\ldots,x^n(t))$
and this equation becomes
$$0=\nabla_{\dot \gamma}\dot \gamma(t) =
\left(\sum\limits_k\left(\ddot x^k(t) +\dot x^i(t)\dot
x^j(t)\Gamma^{k}_{ij}(\gamma(t))\right)\partial_k\right).$$ This is
a system of $2^{nd}$ order ODE's. The local existence, uniqueness
and smoothness of a geodesic through any point $p\in M$ with initial
velocity vector $v\in T_pM$ follow from the classical ODE theory.
Given any two points in a complete manifold, a standard limiting
argument shows that there is a rectifiable curve of minimal length
between these points. Any such curve is a geodesic. We call
geodesics that minimize the length between their endpoints {\em
minimizing geodesics}\index{geodesic!minimizing}.

We have the classical theorem showing that on a complete manifold
all geodesics are defined for all time (see Theorem 16 of Chapter 5
on p. 137 of \cite{Petersen}).

\begin{thm}
(Hopf-Rinow) If $(M,g)$ is complete as a metric space, then every
geodesic extends to a geodesic defined for all time.
\end{thm}

Geodesics are critical points of the energy functional. Let $(M,g)$
be a complete Riemannian manifold. Consider the space of $C^1$-paths
in $M$ parameterized by the unit interval. On this space we have the
energy functional
$$E(\gamma)=\frac{1}{2}\int_{0}^{1}\langle \gamma'(t),\gamma'(t)\rangle dt.$$
 Suppose that we have a one-parameter family of  paths parameterized by
 $[0,1]$,
 all having the same initial point $p$ and the same final point $q$.
 By this we mean that we have a surface $\tilde\gamma(t,u)$ with the property that for each
 $u$ the path $\gamma_u=\tilde\gamma(\cdot,u)$ is a path from $p$ to $q$ parameterized by $[0,1]$.
 Let $\widetilde X=\partial\tilde \gamma/\partial t$ and
 $\widetilde Y=\partial\tilde \gamma/\partial u$ be the corresponding vector
 fields along the surface swept out by $\tilde \gamma$, and denote by $X$ and $Y$ the restriction
 of these vector fields along $\gamma_0$.
 We compute
 \begin{eqnarray*}
 \frac{dE(\gamma_u)}{du}\Bigl|_{u=0}\Bigr.& = & \left( \int_0^1\langle \nabla_{\widetilde Y}\widetilde X,\widetilde
 X\rangle dt\right)\left|_{u=0}\right. \\
 & = & \left(\int_0^1\langle\nabla_{\widetilde X}\widetilde Y,\widetilde
 X\rangle dt \right)\left|_{u=0}\right.\\
 & = & -\left( \int_0^1\langle\nabla_{\widetilde X}\widetilde X,\widetilde
 Y\rangle dt\right)\left|_{u=0}\right.=-\int_0^1\langle\nabla_XX,Y\rangle,
\end{eqnarray*}
where the first equality in the last line comes from integration by
parts and the fact that $\widetilde Y$ vanishes at the endpoints.
Given any vector field $Y$ along $\gamma_0$ there is a one-parameter
family $\tilde \gamma(t,u)$ of paths from $p$ to $q$ with $\tilde
\gamma(t,0)=\gamma_0$ and with $\widetilde Y(t,0)=Y$. Thus, from the
above expression we see that $\gamma_0$ is a critical point for the
energy functional on the space of paths from $p$ to $q$
parameterized by the interval $[0,1]$ if and only if $\gamma_0$ is a
geodesic.

Notice that it follows immediately from the geodesic equation that
the length of a tangent vector along a geodesic is constant. Thus,
if a geodesic is parameterized by $[0,1]$ we have
$$E(\gamma)=\frac{1}{2}L(\gamma)^2.$$
It is immediate from the Cauchy-Schwarz inequality that for any
curve $\mu$ parameterized by $[0,1]$ we have
$$E(\mu)\ge \frac{1}{2}L(\mu)^2$$
with equality if and only if $|\mu'|$ is constant. In particular, a
curve parameterized by $[0,1]$ minimizes distance between its
endpoints if it is a minimum for the energy functional on all paths
parameterized by $[0,1]$ with the given endpoints.

\subsection{Families of geodesics and Jacobi fields}

Consider a family of geodesics $\tilde\gamma(u,t)=\gamma_u(t)$
parameterized by the interval  $[0,1]$ with $\gamma_u(0)=p$ for all
$u$. Here, unlike the discussion above, we allow $\gamma_u(1)$ to
vary with $u$. As before define vector fields along the surface
swept out by $\tilde \gamma$:  $\widetilde X=\partial
\tilde\gamma/\partial t$ and let $\tilde
Y=\partial\tilde\gamma/\partial u$.  We denote by $X$ and $Y$ the
restriction of these vector fields to the geodesic
$\gamma_0=\gamma$.
 Since each
$\gamma_u$ is a geodesic, we have $\nabla_{\widetilde X}{\widetilde
X}=0$. Differentiating this equation in the $\widetilde Y$-direction
yields $\nabla_{\widetilde Y}\nabla_{\widetilde X}\widetilde X=0$.
Interchanging the order of differentiation, using
$\nabla_{\widetilde X}\widetilde Y=\nabla_{\widetilde Y}{\widetilde
X}$, and then restricting to $\gamma$, we get the {\em Jacobi
equation}:
$$\nabla_X\nabla_XY+{\mathcal R}(Y,X)X=0.$$
Notice that the left-hand side of the equation depends only on the
value of $Y$ along $\gamma$, not on the entire family. We denote the
left-hand side of this equation by ${rm Jac}(Y)$, so that the Jacobi
equation now reads
$${rm Jac}(Y)=0.$$
 The fact that all the geodesics begin at the same point at time
$0$ means that $Y(0)=0$. A vector field $Y$ along a geodesic
$\gamma$ is said to be a {\em Jacobi field}\index{Jacobi field} if
it satisfies this equation and vanishes at the initial point $p$. A
Jacobi field is determined by its first derivative at $p$, i.e., by
$\nabla_XY(0)$. We have just seen that this is the equation
describing, to first order, variations of $\gamma$ by a family of
geodesics with the same starting point.

Jacobi fields are also determined by the energy functional. Consider
the space of paths parameterized by $[0,1]$ starting at a given
point $p$ but free to end anywhere in the manifold. Let $\gamma$ be
a geodesic (parameterized by $[0,1]$) from $p$ to $q$. Associated to
any one-parameter family $\tilde \gamma(t,u)$ of paths parameterized
by $[0,1]$ starting at $p$ we associate the second derivative of the
energy at $u=0$. Straightforward computation gives
$$\frac{d^2E(\gamma_u)}{du^2}\Bigl|_{u=0}\Bigr.=\langle \nabla_XY(1),Y(1)\rangle+\langle
X(1),\nabla_Y\widetilde Y(1,0)\rangle-\int_0^1\langle {rm
Jac}(Y),Y\rangle dt.$$

Notice that the first term is a boundary term from the integration
by parts, and it depends not just on the value of $Y$ (i.e., on
$\widetilde Y$ restricted to $\gamma$) but also on the first-order
variation of $\widetilde Y$ in the $Y$ direction. There is the
associated bilinear form that comes from two-parameter families
$\tilde \gamma(t,u_1,u_2)$ whose value at $u_1=u_1=0$ is $\gamma$.
It is
$$\frac{d^2E}{du_1du_2}\Bigl|_{u_1=u_2=0}\Bigr.=\langle \nabla_XY_1(1),Y_2(1)\rangle+\langle
X(1),\nabla_{Y_1}\widetilde Y_2(1,0)\rangle-\int_0^1\langle {rm
Jac}(Y_1),Y_2\rangle dt.$$ Notice that restricting to the space of
vector fields that vanish at both endpoints, the second derivatives
depend only on $Y_1$ and $Y_2$ and the formula is
$$\frac{d^2E}{du_1du_2}\Bigl|_{u_1=u_2=0}\Bigr.=-\int_0^1\langle {rm Jac}(Y_1),Y_2\rangle dt,$$
so that this expression is symmetric in $Y_1$ and $Y_2$. The
associated quadratic form on the space of vector fields along
$\gamma$ vanishing at both endpoints
$$-\int_0^1\langle {rm Jac}(Y),Y\rangle dt$$
is the second derivative of the energy function at $\gamma$ for any
one-parameter family whose value at $0$ is $\gamma$ and whose first
variation is given by $Y$.

\subsection{Minimal geodesics}

\begin{defn}
Let $\gamma$ be a geodesic beginning at $p\in M$. For any $t>0$ we
say that $q=\gamma(t)$ is a {\sl conjugate point along
$\gamma$}\index{conjugate point} if there is a non-zero Jacobi field
along $\gamma$ vanishing at $\gamma(t)$.
\end{defn}

\begin{prop}\label{jacmin}
Suppose that $\gamma\colon [0,1]\to M$ is a minimal geodesic. Then
for any $t<1$ the restriction of $\gamma$ to $[0,t]$ is the unique
minimal geodesic between its endpoints and there are no conjugate
points on $\gamma([0,1))$, i.e., there is no non-zero Jacobi field
along $\gamma$ vanishing at any $t\in [0,1)$.
\end{prop}

We shall sketch the proof. For more details see Proposition 19 and Lemma 14 of
Chapter 5 on pp. 139 and 140 of \cite{Petersen}.

\begin{proof} (Sketch)
Fix $0<t_0<1$. Suppose that there were a different geodesic
$\mu\colon [0,t_0]\to M$ from $\gamma(0)$ to $\gamma(t_0)$, whose
length was at most that of $\gamma|_{[0,t_0]}$. The fact that $\mu$
and $\gamma|_{[0,t_0]}$ are distinct means that
$\mu'(t_0)\not=\gamma'(t_0)$. Then the curve formed by concatenating
$\mu$ with $\gamma|_{[t_0,1]}$ is a curve from $\gamma(0)$ to
$\gamma(1)$ whose length is at most that of $\gamma$. But this
concatenated curve is not smooth at $\mu(t_0)$, and hence it is not
a geodesic, and in particular there is a curve with shorter length
(a minimal geodesic) between these points. This is contrary to our
assumption that $\gamma$ was minimal.

To establish that there are no conjugate points at $\gamma(t_0)$ for
$t_0<1$ we need the following claim.

\begin{claim}
Suppose that $\gamma$ is a minimal geodesic and $Y$ is a field
vanishing at both endpoints. Let $\tilde\gamma(t,u)$ be any
one-parameter family of curves parameterized by $[0,1]$, with
$\gamma_0=\gamma$ and with $\gamma_u(0)=\gamma_0(0)$ for all $u$.
Suppose that the first-order variation of $\tilde \gamma$ at $u=0$
is given by $Y$. Then
$$\frac{d^2E(\gamma_u)}{du^2}\Bigl|_{u=0}\Bigr.=0$$
if and only if $Y$ is a Jacobi field.
\end{claim}

\begin{proof}
Suppose that $\tilde\gamma(u,t)$ is a one-parameter family of curves
from $\gamma(0)$ to $\gamma(1)$ with $\gamma_0=\gamma$ and $Y$ is
the first-order variation of this family along $\gamma$.  Since
$\gamma$ is a minimal geodesic we have
$$-\int_0^1\langle {rm Jac}(Y),Y\rangle dt=\frac{d^2E(\gamma_u)}{du^2}\Bigl|_{u=0}\Bigr.\ge 0.$$
 The associated symmetric bilinear form is
$$B_\gamma(Y_1,Y_2)=-\int_\gamma\langle {rm Jac}(Y_1),Y_2\rangle dt$$
is symmetric when $Y_1$ and $Y_2$ are constrained to vanish at both
endpoints. Since the associated quadratic form is non-negative, we
see by the usual argument for symmetric bilinear forms that
$B_\gamma(Y,Y)=0$ if and only if $B_\gamma(Y,\cdot)=0$ as a linear
functional on the space of vector fields along $\gamma$ vanishing at
point endpoints. This of course occurs if and only if ${rm
Jac}(Y)=0$.
\end{proof}

Now let us use this claim to show that there are no conjugate points
on $\gamma|_{(0,1)}$. If for some $t_0<1$, $\gamma(t_0)$ is a
conjugate point along $\gamma$, then there is a non-zero Jacobi
field $Y(t)$ along $\gamma$ with $Y(t_0)=0$. Notice that since $Y$
is non-trivial $\nabla_XY(t_0)\not= 0$. Extend $Y(t)$ to a vector
field $\hat Y$ along all of $\gamma$ by setting it equal to $0$ on
$\gamma|_{[t_0,1]}$. Since the restriction of $Y$ to
$\gamma([0,t_0])$ is a Jacobi field vanishing at both ends and since
$\gamma|_{[0,t_0]}$ is a minimal geodesic, the second-order
variation of length of $\gamma|_{[0,t_0]}$ in the $Y$-direction is
zero. It follows that the second-order variation of length along
$\hat Y$ vanishes. But $\hat Y$ is not smooth (at $\gamma(t_0)$) and
hence it is not a Jacobi field along $\gamma$. This contradicts the
fact discussed in the previous paragraph that for minimal geodesics
the null space of the quadratic form is exactly the space of Jacobi
fields.
\end{proof}

\subsection{The exponential mapping}

\begin{defn}
For any $p\in M$, we can define the \textit{exponential
map}\index{exponential mapping} at $p$, ${\rm exp}_{p}$. It is
defined on an open neighborhood $O_{p}$ of the origin in $T_{p}M$
and is defined by ${\rm exp}_{p}(v)=\gamma_v(1)$, the endpoint of
the unique geodesic $\gamma_v\colon [0,1]\rightarrow M$ starting
from $p$ with initial velocity vector $v$. We always take
$O_{p}\subset T_{p}M$ to be the maximal domain on which ${\rm
exp}_p$ is defined, so that $O_{p}$ is a star-shaped open
neighborhood of $0\in T_pM$. By the Hopf-Rinow Theorem, if $M$ is
complete, then the exponential map is define on all of $T_pM$.
\end{defn}

By the inverse function theorem there exists $r_0=r_{0}(p,M)
> 0$, such that the restriction of ${\rm exp}_p$ to the ball
$B_{g|T_pM}(0,r_0)$ in $T_pM$  is a diffeomorphism  onto
$B_g(p,r_0)$. Fix $g$-orthonormal linear coordinates on $T_pM$.
Transferring these coordinates via ${\rm exp}_p$ to coordinates on
$B(p,r_0)$ gives us {\sl Gaussian normal coordinates}\index{Gaussian
normal coordinates} on $B(p,r_0)\subset M$.

Suppose now that  $M$ is complete, and fix a point $p\in M$. For
every $q\in M$, there is a length-minimizing path from $p$ to $q$.
When parameterized at constant speed equal to its length, this path
is a geodesic with domain  interval $[0,1]$. Consequently, ${\rm
exp}_p\colon T_pM\to M$ is onto. The differential of the exponential
mapping is given by Jacobi fields:
 Let
$\gamma\colon [0,1]\to M$ be a geodesic from $p$ to $q$, and let
 $X\in T_pM$ be $\gamma'(0)$. Then the exponential
mapping at $p$ is a smooth map from $T_p(M)\to M$ sending $X$ to
$q$. Fix $Z\in T_pM$. Then there is a unique Jacobi field $Y_Z$
along $\gamma$ with $\nabla_XY_Z(0)=Z$. The association $Z\mapsto
Y_Z(1)\in T_qM$ is a linear map from $T_p(M)\to T_qM$. Under the
natural identification of $T_pM$ with the tangent plane to $T_pM$ at
the point $Z$, this linear mapping is the differential of ${\rm
exp}_p\colon T_pM\to M$ at the point $X\in T_pM$.

\begin{cor}\label{star}
Suppose that $\gamma$ is a minimal geodesic parameterized by $[0,1]$
starting at $p$. Let $X(0)=\gamma'(0)\in T_pM$. Then for each
$t_0<1$ the restriction $\gamma|_{[0,t_0]}$ is a minimal geodesic
and ${\rm exp}_p\colon T_pM\to M$ is a local diffeomorphism near
$t_0X(0)$.
\end{cor}

\begin{proof}
Of course, ${\rm exp}_p(t_0X(0))=\gamma(t_0)$. According to the
previous discussion, the kernel of the differential of the
exponential mapping at $t_0X(0)$ is identified with the space of
Jacobi fields along $\gamma$ vanishing at $\gamma(t_0)$. According
to  Proposition~\ref{jacmin} the only such Jacobi field is the
trivial one. Hence, the differential of ${\rm exp}_p$ at $t_0X(0)$
is an isomorphism, completing the proof.
\end{proof}

\begin{defn}
There is an open neighborhood  $U_p\subset T_pM$ of $0$ consisting
of all $v\in T_pM$ for which: (i) $\gamma_v$ is the unique minimal
geodesic from $p$ to $\gamma_v(1)$, and (ii) $ {\rm exp}_p$ is a
local diffeomorphism at $v$.  We set ${\mathcal C}_p\subset M$ equal
to $M\setminus {\rm exp}_p(U_p)$. Then ${\mathcal C}_p$ is called
the {\sl cut locus from $p$}\index{cut locus}. It is a closed subset
of measure $0$.
\end{defn}

It follows from Corollary~\ref{star} that  $U\subset T_pM$ is a
star-shaped open neighborhood  of $0\in T_pM$.

\begin{prop}
The map
$${\rm exp}_p\colon U_p\to M\setminus {\mathcal C}_p$$
is a diffeomorphism. \end{prop} For a proof see p. 139 of \cite{Petersen}.

\begin{defn}
The {\sl injectivity radius ${\rm inj}_M(p)$ of $M$ at $p$} is the
supremum of the $r> 0$ for which the restriction of ${\rm
exp}_p\colon T_pM\to M$ to the ball $B(0,r)$ of radius $r$ in $T_pM$
is a diffeomorphism into $M$. Clearly, ${\rm inj}_M(p)$ is the
distance in $T_pM$ from $0$ to the frontier of $U_p$. It is also the
distance in $M$ from $p$ to the cut locus ${\mathcal C}_p$.
\end{defn}

Suppose that  ${\rm inj}_M(p)=r$. There are two possibilities:
Either there is a broken, closed geodesic through $p$, broken only
at $p$, of length $2r$, or there is a geodesic $\gamma$ of length
$r$ emanating from $p$ whose endpoint is a conjugate point along
$\gamma$. The first case happens when the exponential mapping is not
one-to-one of the closed ball of radius $r$ in $T_pM$, and the
second happens when there is a tangent vector in $T_pM$ of length
$r$ at which ${\rm exp}_p$ is not a local diffeomorphism.

\section{Computations in Gaussian normal coordinates}

In this section we compute the metric and the Laplacian (on
functions) in local Gaussian coordinates. A direct computation shows
that in Gaussian normal coordinates on a metric ball about $p\in M$
the metric takes the form
\begin{eqnarray}\label{metricexp}
g_{ij}(x)& = & \delta_{ij} + \frac{1}{3}R_{iklj}x^kx^l +
\frac{1}{6}R_{iklj,s}x^kx^lx^s \\
& & +(\frac{1}{20}R_{iklj,st} +\frac{2}{45}\sum_m
R_{iklm}R_{jstm})x^kx^lx^sx^t + O(r^{5}), \nonumber
\end{eqnarray}
where $r$ is the distance from $p$. (See, for example Proposition
3.1 on page 41 of~\cite{Sakai}, with the understanding that, with
the conventions there, the quantity $R_{ijkl}$ there differs by sign
from ours.)

Let $\gamma$ be a geodesic in $M$ emanating from $p$ in the direction $v$.
Choose local coordinates $\theta^1,\ldots,\theta^{n-1}$ on the unit sphere in
$T_pM$ in a neighborhood of $v/|v|$. Then $(r, \theta^{1}, ... ,\theta^{n-1})$
are local coordinates at any point of the ray emanating from the origin in the
$v$ direction (except at $p$). Transferring these via ${\rm exp}_p$ produces
local coordinates $(r,\theta^{1}, ... ,\theta^{n-1})$ along $\gamma$. Using
Gauss's lemma (Lemma 12 of Chapter 5 on p. 133 of \cite{Petersen}), we can
write the metric locally as
\[ g=dr^{2}+
r^{2}h_{ij}(r,\theta)d\theta^{i}\otimes d\theta^{j}. \] Then the
volume form
\begin{align*}
dV &= \sqrt{{\rm det}(g_{ij})}dr\wedge d\theta^{1}
\wedge\cdots\wedge
d\theta^{n-1}\\
& = r^{n-1}\sqrt{{\rm det}(h_{ij})}dr\wedge d\theta^{1}
\wedge\cdots\wedge d\theta^{n-1}.
\end{align*}

\begin{lem}
The \textit{Laplacian operator}\index{Laplacian} acting on scalar
functions on $M$ is given in local coordinates by
\[ \triangle  = \frac{1}{\sqrt{{\rm det}(g)}}\partial_i
\left(g^{ij}\sqrt{{\rm det}(g)}\partial_j\right). \]
\end{lem}

\begin{proof}
Let us compute the derivative at a point $p$. We have
$$\frac{1}{\sqrt{{\rm det}(g)}}\partial_i
\left(g^{ij}\sqrt{{\rm
det}(g)}\partial_j\right)f=g^{ij}\partial_i\partial_jf+\partial_ig^{ij}\partial_jf+
\frac{1}{2}g^{ij}\partial_iTr(\tilde g)\partial_jf,$$ where $\tilde
g=g(p)^{-1}g$. On the other hand from the definition of the
Laplacian, Equation~(\ref{laplacformula}), and
Equation~(\ref{Hessformula}) we have
$$\triangle
f=g^{ij}{\rm
Hess}(f)(\partial_i,\partial_j)=g^{ij}\left(\partial_i\partial_j(f)-
\nabla_{\partial_i}\partial_jf\right)
=g^{ij}\partial_i\partial_jf-g^{ij}\Gamma_{ij}^k\partial_kf.$$ Thus,
to prove the claim it suffices to show that
$$g^{ij}\Gamma_{ij}^k=-(\partial_ig^{ik}+\frac{1}{2}g^{ik}{\rm Tr}(\partial_i\tilde g)).$$
From the definition of the Christoffel symbols we have
$$g^{ij}\Gamma_{ij}^k=\frac{1}{2}g^{ij}g^{kl}(\partial_ig_{jl}+\partial_jg_{il}-\partial_lg_{ij}).$$
Of course, $g^{ij}\partial_ig_{jl}=-\partial_ig^{ij}g_{jl}$, so that
 $g^{ij}g^{kl}\partial_ig_{jl}=-\partial_ig^{ik}$. It follows by symmetry that $g^{ij}g^{jl}\partial_jg_{il}
 =-\partial_ig^{ik}$. The last term is clearly $-\frac{1}{2}g^{ik}{\rm Tr}(\partial_i\tilde g).$
\end{proof}

Using Gaussian local coordinates near $p$, we have
\begin{align*}
\triangle r &=
\frac{1}{r^{n-1}\sqrt{{\rm det}(h)}}\partial_r\left(r^{n-1}\sqrt{{\rm det}(h)}\right)\\
&= \frac{n-1}{r}
 + \partial_r\log\left(\sqrt{{\rm det}(h)}\right).
\end{align*}
From this one computes directly that
 \[ \triangle r
 = \frac{n-1}{r} - \frac{r}{3}{\rm Ric}(v,v)+ O(r^{2}), \]
 where $v = \dot r(0)$, cf,  p.265-268 of \cite{Petersen}. So
 $$\triangle r \leq
 \frac{n-1}{r}\ \ {\rm when\ \ } r \ll 1\ \ {\rm and\ } {\rm Ric} > 0.$$
This local computation has the following global analogue.
\smallskip

\begin{ex}(E.Calabi, 1958)
Let $f(x)=d(p,x)$ be the distance function from $p$. If $(M, g)$ has
${\rm Ric} \geq 0$, then
\[ \triangle f \leq \frac{n-1}{f} \]
 in the sense of distributions.
\end{ex}
\noindent [Compare \cite{Petersen}, p. 284 Lemma 42].

\begin{rem}\label{calabi}
The statement that $\triangle f\leq \frac{n-1}{f}$ in {\sl the sense
of distributions} (or equivalently {\sl in the weak sense}) means
that for any non-negative test function $\phi$, that is to say for
any compactly supported $C^\infty$-function $\phi$, we have
$$\int_Mf\triangle \phi d{\rm vol}\le \int_M\left(\frac{n-1}{f}\right)\phi d{\rm vol}.$$
Since  the triangle inequality implies that $|f(x)-f(y)|\le d(x,y)$, it follows
that $f$ is Lipschitz, and hence that the restriction of $\nabla f$ to any
compact subset of $M$ is an $L^2$ one-form. Integration by parts then shows
that
$$\int_Mf\triangle
\phi d{\rm vol}=-\int_M\langle \nabla f,\nabla\phi\rangle d{\rm
vol}.$$
\end{rem}

\medskip

Since $\abs{\nabla f} = 1$ and $\triangle f$ is the mean curvature
of the geodesic sphere $\partial B(x,r)$, ${\rm Ric}(v,v)$ measures
the difference of the mean curvature between the standard Euclidean
sphere and the geodesic sphere in the direction $v$. Another
important geometric object  is the shape operator associated to $f$,
denoted $S$. By definition it is the Hessian of $f$; i.e.,
$S=\nabla^2f={\rm Hess}(f)$.

\section{Basic curvature comparison results}
In this section we will recall some of the basic curvature
comparison results in Riemannian geometry. The reader can refer to
\cite{Petersen}, Section 1 of Chapter 9 for details.

\smallskip

We fix a point $p\in M$. For any real number $k\ge 0$ let $H^n_k$ denote the
simply connected, complete Riemannian $n$-manifold of constant sectional
curvature $-k$. Fix a point $q_k\in H^n_k$, and consider the exponential map
${\rm exp}_{q_k}\colon T_{q_k}(H^n_k)\to H^n_k$. This map is a global
diffeomorphism. Let us consider the pullback, $\tilde h_k$, of the Riemannian
metric on $H^n_k$ to $T_{q_k}H^n_k$. A formula for this tensor is easily given
in polar coordinates on $T_{q_k}(H^n_k)$ in terms of the following function.

\begin{defn}
We define a function ${\rm sn}_k$ as follows:$${\rm
sn}_k(r)=\begin{cases}
r \ \  & {\rm if\ \ }k=0 \\
\frac{1}{\sqrt{k}}{\rm sinh}(\sqrt{k}r) \ \ & {\rm if \ \ } k>0.
\end{cases}$$
\end{defn}

The function ${\rm sn}_{k}(r)$ is the solution to the equation
\begin{align*}
 \varphi'' - k\varphi &= 0,\\
\varphi(0) &= 0,\\
\varphi'(0) &= 1.
\end{align*}
We define
 ${\rm ct}_{k}(r) = \frac{{\rm sn}_{k}^{'}(r)}{\sqrt{k}{\rm sn}_{k}(r)}$.

Now we can compare manifolds of varying sectional curvature with
those of constant curvature.

\begin{thm}\label{SCC}
(Sectional Curvature Comparison) Fix $k\ge 0$. Let $(M, g)$ be a
Riemannian manifold with the property that $-k \leq K(P)$ for every
$2$-plane $P$ in $TM$. Fix a minimizing geodesic $\gamma\colon
[0,r_0)\to M$ parameterized at unit speed with $\gamma(0)=p$. Impose
Gaussian polar coordinates $(r,\theta^1,\ldots,\theta^{n-1})$ on a
neighborhood of $\gamma$ so that $g=dr^2+g_{ij}\theta^i\otimes
\theta^j$. Then for all $0<r<r_0$ we have
$$ (g_{ij}(r,\theta))_{1 \leq
i,j \leq n-1} \leq {\rm sn}_{k}^{2}(r),$$ and the shape operator
associated to  the distance function from $p$, $f$, satisfies
$$(S_{ij}(r,\theta))_{1 \leq i,j \leq n-1} \leq \sqrt{k}{\rm ct}_{k}(r).$$
\end{thm}

There is also an analogous result for a positive upper bound to the
sectional curvature, but in fact all we shall need is the local
diffeomorphism property of the exponential mapping.

\begin{lem}\label{localdiffeo} Fix $K\ge 0$.
If $|{\rm Rm}(x)|\le K$ for all  $x\in B(p,\pi/\sqrt{K})$, then
${\rm exp}_p$ is a local diffeomorphism from the ball
$B(0,\pi/\sqrt{K})$ in $T_pM$  to the  ball $B(p,\pi/\sqrt{K})$ in
$M$.
\end{lem}

There is a crucial comparison result for volume which involves the
Ricci curvature.

\begin{thm}\label{riccurvcomp}
{\bf (Ricci curvature comparison)} Fix $k\ge 0$. Assume that $(M,
g)$ satisfies ${\rm Ric} \geq -(n-1)k$. Let $\gamma\colon [0,r_0)\to
M$ be a minimal geodesic of unit speed. Then for any $r<r_0$ at
$\gamma(r)$ we have
$$\sqrt{{\rm det}\, g(r,\theta)} \leq {\rm sn}^{n-1}_{k}(r)$$ and
$${\rm Tr}(S)(r,\theta) \leq (n-1)\frac{{\rm sn}_{k}^{'}(r)}{{\rm sn}_{k}(r)}.$$
\end{thm}

Note that the inequality in Remark~\ref{calabi} follows from this
theorem.

The comparison result in Theorem~\ref{riccurvcomp} holds out to
every radius, a fact that will be used repeatedly in our arguments.
This result evolved over the period 1964-1980 and now is referred to
as the Bishop-Gromov inequality\index{Bishop-Gromov Theorem|ii}; see
Proposition 4.1 of \cite{CheegerGromovTaylor}

\begin{thm}\label{BishopGromov}
(Relative Volume Comparison, Bishop-Gromov 1964-1980) Suppose $(M,
g)$ is a Riemannian manifold. Fix a point $p\in M$, and suppose that
$B(p,R)$ has compact closure in $M$. Suppose that for some $k\ge 0$
we have ${\rm Ric} \geq -(n-1)k$ on $B(p,R)$. Recall that $H^n_k$ is
the simply connected, complete manifold of constant curvature $-k$
and $q_k\in H^n_k$ is a point. Then
\[ \frac{{\rm Vol}\,B(p,r)}{{\rm Vol}\, B_{H^n_k}B(q_k,r)} \]
is a non-increasing function of $r$ for $r<R$, whose limit as $r\rightarrow 0$
is 1. In particular, if the Ricci curvature of $(M,g)$ is $\ge 0$ on $B(p,R)$,
then ${\rm Vol}\,B(p,r)/r^n$ is a non-increasing function of $r$ for $r<R$.
\end{thm}

\section{Local volume and the injectivity radius}

As the following results show, in the presence of bounded curvature the volume
of a ball $B(p,r)$ in $M$ is bounded away from zero if and only if the
injectivity radius of $M$ at $p$ is bounded away from zero.

\begin{prop}\label{injvol}
Fix an integer $n>0$. For every $\epsilon>0$ there is $\delta>0$
depending on $n$ and $\epsilon$ such that the following holds.
Suppose that $(M^n,g)$ is a complete Riemannian manifold of
dimension $n$ and that $p\in M$. Suppose that $|{\rm Rm}(x)|\le
r^{-2}$ for all $x\in B(p,r)$. If the injectivity radius of $M$ at
$p$ is at least $\epsilon r$, then ${\rm Vol}(B(p,r))\ge \delta
r^n$.
\end{prop}

\begin{proof}
Suppose that $|{\rm Rm}(x)|\le r^{-2}$ for all $x\in B(p,r)$.
Replacing $g$ by $r^2 g$ allows us to assume that $r=1$. Without
loss of generality we can assume that $\epsilon\le 1$. The map ${\rm
exp}_p$ is a diffeomorphism on the ball $B(0,\epsilon)$ in the
tangent space, and by Theorem~\ref{SCC} the volume of
$B(p,\epsilon)$ is at least that of the ball of radius $\epsilon$ in
the $n$-sphere of radius $1$. This gives a lower bound to the volume
of $B(p,\epsilon)$, and a fortiori to $B(p,1)$, in terms of $n$ and
$\epsilon$.
\end{proof}

 We shall normally work with volume, which behaves nicely under Ricci
flow, but in order to take limits we need to bound the injectivity
radius away from zero. Thus, the more important, indeed crucial,
result for our purposes is the converse to the previous proposition;
see Theorem 4.3, especially Inequality (4.22), on page 46 of
\cite{CheegerGromovTaylor}, or see Theorem 5.8 on page 96 of
\cite{CheegerEbin}.

\begin{thm}\label{volinj}
Fix an integer $n>0$. For every $\epsilon>0$ there is $\delta>0$
depending on $n$ and $\epsilon$ such that the following holds.
Suppose that $(M^n,g)$ is a complete Riemannian manifold of
dimension $n$ and that $p\in M$. Suppose that $|{\rm Rm}(x)|\le
r^{-2}$ for all $x\in B(p,r)$. If $Vol(B(p,r))\ge \epsilon r^n$ then
the injectivity radius of $M$ at $p$ is at least $\delta r$.
\end{thm}

\chapter{Manifolds of non-negative curvature}\label{nonnegcurv}

In studying singularity development in $3$-dimensional Ricci flows
one forms blow-up limits. By this we mean the following. One
considers a sequence of points $x_k$ in the flow converging to the
singularity. It will be the case that $R(x_k)$ tends to $\infty$ as
$k$ tends to $\infty$. We form a sequence of based Riemannian
manifolds labeled by $k$, where the $k^{th}$ Riemannian manifold is
obtained by taking the time-slice of $x_k$, rescaling its metric by
$R(x_k)$, and then taking $x_k$ as the base point. This creates a
sequence with the property that for each member of the sequence the
scalar curvature at the base point is one. Because of a pinching
result of Hamilton's (see Chapter~\ref{secmaxprin}), if there is a
geometric limit of this sequence, or of any subsequence of it, then
that limit is non-negatively curved. Hence, it is important to
understand the basic properties of Riemannian manifolds of
non-negative curvature in order to study singularity development. In
this chapter we review the properties that we shall need. We suppose
that $M$ is non-compact and of positive (resp., non-negative)
curvature. The key to understanding these manifolds is the Busemann
function associated to a minimizing geodesic ray.

\section{Busemann functions}

A geodesic ray $\lambda\colon [0,\infty)\to M$  is said to be {\em
minimizing} if the restriction of $\lambda$ to every compact
subinterval of $[0,\infty)$ is a length-minimizing geodesic arc,
i.e., a geodesic arc whose length is equal to the distance between
its endpoints. Likewise, a geodesic line $\lambda\colon
(-\infty,\infty)\to M$ is said to be {\em minimizing} if its
restriction to every compact sub-interval of $\Ar$ is a length
minimizing geodesic arc.

Clearly, if a sequence of minimizing geodesic arcs $\lambda_k$
converges to a geodesic arc, then the limiting geodesic arc is also
minimizing. More generally, if $\lambda_k$ is a sequence of length
minimizing geodesic arcs whose initial points converge and whose
lengths go to infinity, then, after passing to a subsequence, there
is a limit which is a minimizing geodesic ray. (The existence of a
limit of a subsequence is a consequence of the fact that a geodesic
ray is determined by its initial point and its initial tangent
direction.) Similarly, if $I_k$ is an sequence of compact intervals
with the property that every compact subset of $\Ar$ is contained in
$I_k$ for all sufficiently large $k$, if for each $k$ the map
$\lambda_k\colon I_k\to M$ is a minimizing geodesic arc, and if
${\rm lim}_{k\rightarrow \infty}\lambda_k(0)$ exists, then, after
passing to a subsequence there is a limit which is a minimizing
geodesic line. Using these facts one establishes the following
elementary lemma.

\begin{lem}\label{ends}
Suppose that $M$ is a complete, connected, non-compact Riemannian
manifold and let $p$ be a point of $M$. Then $M$ has a minimizing
geodesic ray emanating from $p$. If $M$ has more than one end, then
it has a minimizing line.
\end{lem}

\begin{defn}
Suppose that $\lambda\colon [0,\infty)\to M$ is a minimizing
geodesic ray with initial point $p$. For each $t\ge 0$ we consider
$B_{\lambda,t}(x)=d(\lambda(t),x)-t$. This is a family of functions
satisfying $|B_{\lambda,t}(x)-B_{\lambda,t}(y)|\le d(x,y)$. Since
$\lambda$ is a minimizing geodesic, $B_{\lambda,t}(p)=0$ for all
$t$. It follows that $B_{\lambda,t}(x)\ge -d(x,p)$ for all $x\in M$.
Thus,  the family of functions $B_{\lambda,t}$ is pointwise bounded
below. The triangle inequality shows that for each $x\in M$ the
function $B_{\lambda,t}(x)$ is a non-increasing function of $t$.
 It follows that, for
each $x\in M$, ${\rm lim}_{t\rightarrow \infty} B_{\lambda,t}(x)$
exists. We denote this limit by $B_\lambda(x)$. This is the {\sl
Busemann function}\index{Busemann function|ii} for $\lambda$.
\end{defn}

Clearly, $B_\lambda(x)\ge -d(x,\lambda(0))$. By equicontinuity
$B_\lambda(x)$ is a continuous function of $x$ and in fact a
Lipschitz function satisfying $|B_\lambda(x)-B_\lambda(y)|\le
d(x,y)$ for all $x,y\in X$. Clearly $B_\lambda(\lambda(s))=-s$ for
all $s\ge 0$. Since $B_\lambda$ is Lipschitz, $\nabla B_\lambda$ is
well-defined as an $L^2$-vector field.

\begin{prop}\label{Blambda}
Suppose that $M$ is complete and of non-negative Ricci curvature.
Then, for any minimizing geodesic ray $\lambda$, the Busemann
function $B_\lambda$ satisfies $\Delta B_\lambda\le 0$ in the weak
sense.
\end{prop}

\begin{proof}
First notice that since $B_\lambda$ is Lipschitz, $\nabla B_\lambda$
is an $L^2$-vector field on $M$. That is to say, $B_\lambda\in
W_{\rm loc}^{1,2}$, i.e., $B_\lambda$ locally has one derivative in
$L^2$. Hence, there is a sequence of $C^\infty$-functions $f_n$
converging to $B_\lambda$ in $W^{1,2}_{loc}$. Let $\varphi$ be a
test function (i.e., a compactly supported $C^\infty$-function).
Integrating by parts yields
$$-\int_M\langle \nabla
f_n,\nabla\varphi\rangle d{\rm vol}=\int_Mf_n\triangle \varphi d{\rm
vol}.$$ Using the fact that $f_n$ converges to $B_\lambda$ in
$W^{1,2}_{loc}$ and taking limits yields
$$-\int_M\langle \nabla B_\lambda,\nabla \varphi\rangle d{\rm vol}= \int_M
B_\lambda\triangle \varphi d{\rm vol}.$$

Thus, to prove the proposition we need only show that if $\varphi$
is a non-negative test function, then
$$-\int_M \langle\nabla B_\lambda,\nabla\varphi\rangle d{\rm vol}\le 0.$$
For a proof of this see
 Proposition 1.1  and its proof on pp. 7 and 8 in \cite{SchoenYau}.
\end{proof}

\section{Comparison results in non-negative curvature}

Let us review some elementary comparison results for manifolds of
non-negative curvature. These form the basis for Toponogov
theory\index{Toponogov theory|(ii},  \cite{Toponogov}. For any pair
of points $x,y$ in a complete Riemannian manifold $s_{xy}$ denotes a
minimizing geodesic from $x$ to $y$. We set $|s_{xy}|=d(x,y)$ and
call it the {\sl length} of the side. A {\sl triangle} in a
Riemannian manifold consists of three vertices $a,b,c$ and three
sides $s_{ab}$,$s_{ac}$,$s_{bc}$. We denote by $\angle_a$ the angle
of the triangle at $a$, i.e., the angle at $a$ between the geodesic
rays $s_{ab}$ and $s_{ac}$.

\begin{thm}{\bf (Length comparison)}\label{lengthcompar}
Let $(M,g)$ be a manifold of non-negative curvature. Suppose that
$\triangle(a,b,c)$ is a triangle in $M$  and let
$\triangle(a',b',c')$ be a Euclidean triangle.
\begin{enumerate}
\item  Suppose that the corresponding
sides of $\triangle(a,b,c)$ and $\triangle(a',b',c')$ have the same
lengths. Then the angle at each vertex of the Euclidean triangle is
no larger than the corresponding angle of $\triangle(a,b,c)$.
Furthermore, for any $\alpha$ and $\beta$ less than $|s_{ab}|$ and
$|s_{ac}|$ respectively, let $x$, resp. $x'$, be the point on
$s_{ab}$, resp. $s_{a'b'}$, at distance $\alpha$ from $a$, resp.
$a'$, and let $y$, resp. $y'$, be the point on $s_{ac}$, resp.
$s_{a'c'}$, at distance $\beta$ from $a$, resp. $a'$. Then
$d(x,y)\ge d(x',y')$.
\item Suppose that $|s_{ab}|=|s_{a'b'}|$, that
$|s_{ac}|=|s_{a'c'}|$ and that $\angle_a=\angle_{a'}$. Then
$|s_{b'c'}|\ge |s_{bc}|$.
\end{enumerate}
\end{thm}

\begin{figure}[ht]
  \relabelbox{
  \centerline{\epsfbox{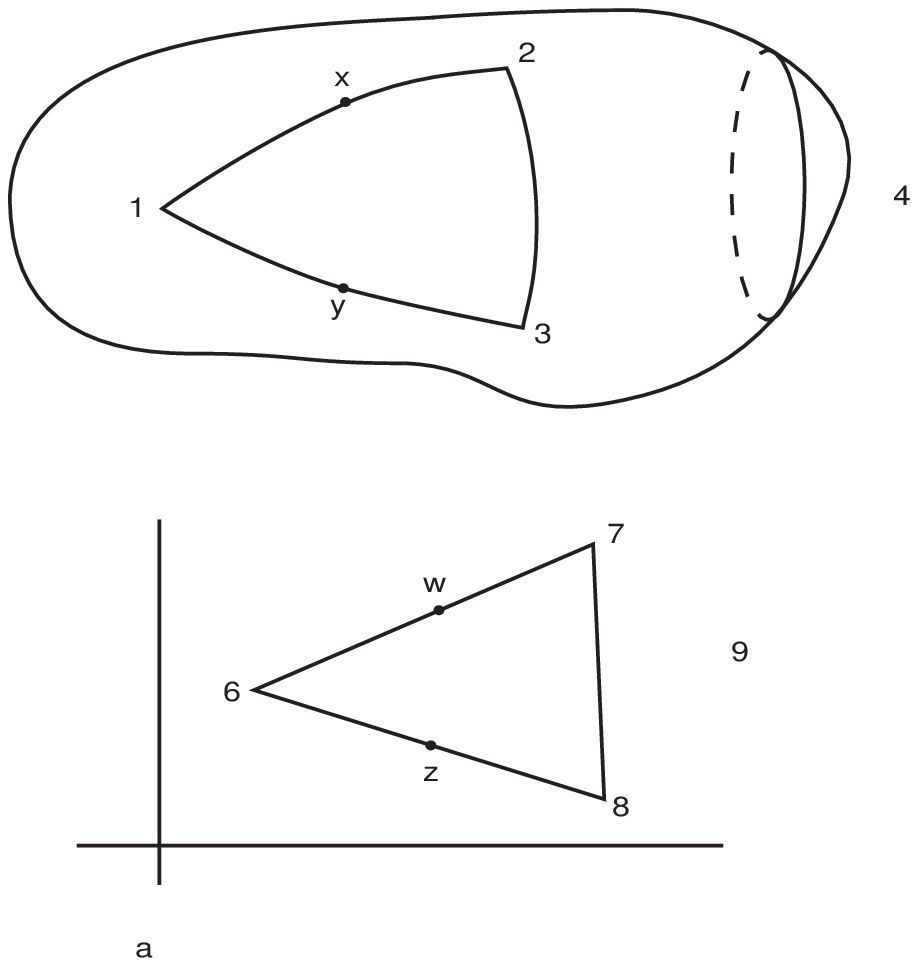}}}
  \relabel{1}{$a$}
  \relabel{2}{$b$}
  \relabel{3}{$c$}
  \relabel{4}{$M$}
  \relabel{6}{$a'$}
  \relabel{7}{$b'$}
  \relabel{8}{$c'$}
  \relabel{9}{$\mathbb{R}^2$}
  \relabel{w}{$x'$}
  \relabel{z}{$y'$}
  \relabel{a}{$d(x,y)\geq d(x',y')  \quad \text{and}\quad \angle_{bac}\geq \angle_{b'a'c'}$}
  \relabel{x}{$x$}
  \relabel{y}{$y$}
  \endrelabelbox
  \caption{Toponogov comparison.}\label{fig:Topog}
\end{figure}

See \textsc{Fig.}~\ref{fig:Topog}. For a proof of this result see
Theorem 4.2 on page 161 of \cite{Sakai}, or Theorem 2.2 on page 42
of \cite{CheegerEbin}.

One corollary  is a monotonicity result. Suppose that
$\triangle(a,b,c)$ is a triangle in a complete manifold of
non-negative curvature. Define a function $EA(u,v)$ defined for
$0\le u\le |s_{ab}|$ and $0\le v\le |s_{ac}|$ as follows. For $u$
and $v$ in the indicated ranges, let $x(u)$ be the point on $s_{ab}$
at distance $u$ from $a$ and let $y(v)$ be the point of $s_{ac}$ at
distance $v$ from $a$. Let  $EA(u,v)$ be the angle at $a'$ of the
Euclidean triangle with side lengths $|s_{a'b'}|=u$, $|s_{a'c'}|=v$
and $|s_{b'c'}|=d(x(u),y(v))$.

\begin{cor}\label{anglemono}
Under the assumptions of the previous theorem, $EA(u,v)$ is a
monotone non-increasing function of each variable $u$ and $v$ when
the other variable is held fixed.
\end{cor}

 Suppose that $\alpha, \beta,\gamma$ are three geodesics
emanating from a point $p$ in a Riemannian manifold. Let $\angle_p
(\alpha,\beta)$, $\angle_p(\beta,\gamma)$ and
$\angle_p(\alpha,\gamma)$ be the angles of these geodesics at $p$ as
measured by the Riemannian metric. Then of course
$$\angle_p(\alpha,\beta)+\angle_p(\beta,\gamma)+\angle_p(\alpha,\gamma)\le
2\pi$$ since this inequality holds for the angles between straight
lines in Euclidean $n$-space. There is a  second corollary of
Theorem~\ref{lengthcompar} which gives an analogous result for the
associated Euclidean angles.

\begin{cor}\label{anglecompar}
Let $(M,g)$ be a complete Riemannian manifold of non-negative
curvature. Let $p,a,b,c$ be four points in $M$ and let
$\alpha,\beta,\gamma$ be minimizing geodesic arcs  from the point
$p$ to $a,b,c$ respectively. Let $T(a,p,b)$, $T(b,p,c)$ and
$T(c,p,a)$ be the triangles in $M$ made out of these minimizing
geodesics and minimizing geodesics between $a,b,c$. Let
$T(a',p',b')$, $T(b',p',c')$ and $T(c',p',a')$ be planar triangles
with the same side lengths. Then
$$\angle_{p'}T(a',p',b')+\angle_{p'}T(b',p',c')+\angle_{p'}T(c',p',a')\le 2\pi.$$
\end{cor}

\begin{proof}
Consider the sum of these angles as the geodesic arcs in $M$ are
shortened without changing their direction. By the first property of
Theorem~\ref{lengthcompar} the sum of the angles of these triangles
is a monotone decreasing function of the lengths. Of course, the
limit as the lengths all go to zero is the corresponding Euclidean
angle. The result is now clear\index{Toponogov theory|)}.
\end{proof}

\section{The soul theorem}

A subset $X$  of a Riemannian manifold $(M,g)$ is said to be {\em
totally convex} if every geodesic segment with endpoints in $X$ is
contained in $X$. Thus, a point $p$ in $M$ is totally convex if and
only if there is no broken geodesic arc in $M$ broken exactly at
$x$.

\begin{thm}\label{soul} (Cheeger-Gromoll, see
\cite{CheegerGromollbul} and \cite{CheegerGromoll}) Suppose that
$(M,g)$ is a connected, complete, non-compact Riemannian manifold of
non-negative sectional curvature. Then $M$ contains a soul $S\subset
M$. By definition a {\em soul}\index{soul|ii} is a compact, totally
geodesic, totally convex submanifold (automatically of positive
codimension). Furthermore, $M$ is diffeomorphic to the total space
of the normal bundle of the $S$ in $M$. If $(M,g)$ has positive
curvature, then any soul for it is a point, and consequently $M$ is
diffeomorphic to $\Ar^n$.
\end{thm}

\begin{rem}
We only use the soul theorem for manifolds with positive curvature
and the fact that any soul of such a manifold is a point. A proof of
this result  first appears in \cite{GromollMeyer}.
\end{rem}

The rest of this section is devoted to a sketch of the proof of this
result. Our discussion follows closely that in \cite{Petersen}
starting on  p.~349. We shall need more information about complete,
non-compact manifolds of non-negative curvature, so we review a
little of their theory as we sketch the proof of the soul theorem.

\begin{lem}
Let $(M,g)$ be a complete, non-compact Riemannian manifold of
non-negative sectional curvature and let $p\in M$. For every
$\epsilon>0$ there is a compact subset $K=K(p,\epsilon)\subset M$
such that for all points $q\notin K$, if $\gamma$ and $\mu$ are
minimizing geodesics from $p$ to $q$, then the angle that $\gamma$
and $\mu$ make at $q$ is less than $\epsilon$.
\end{lem}

\begin{figure}[ht]
  \relabelbox{
  \centerline{\epsfbox{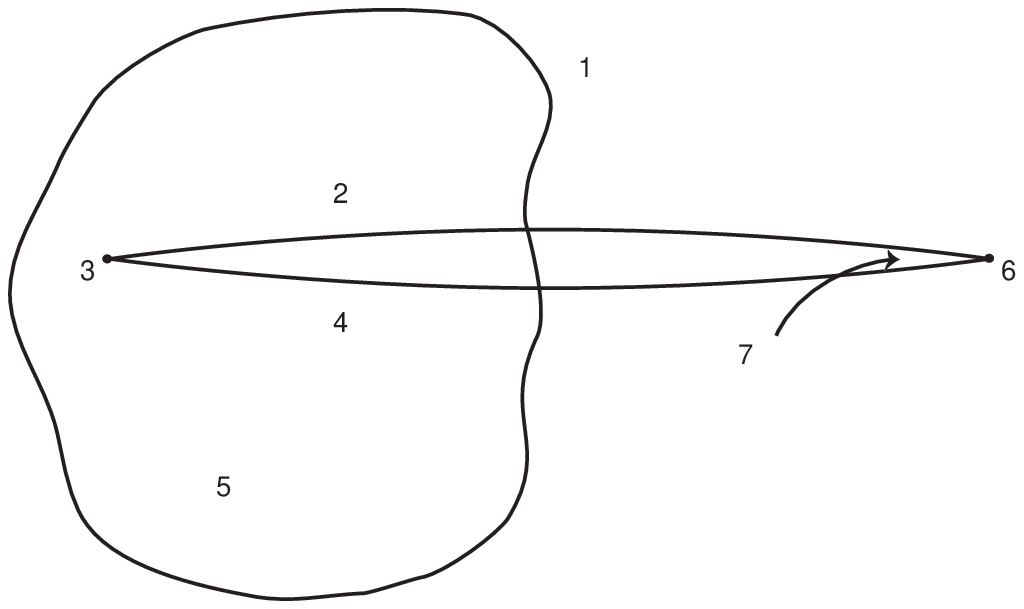}}}
  \relabel{1}{$\gamma, \gamma' $ minimal geodesics}
  \relabel{2}{$\gamma$}
  \relabel{3}{$p$}
  \relabel{4}{$\gamma'$}
  \relabel{5}{$K(p,\epsilon)$}
  \relabel{6}{$q$}
  \relabel{7}{angle $<\epsilon$}
  \endrelabelbox
  \caption{Shallow angles of minimal geodesics.}\label{fig:angles}
\end{figure}

See \textsc{Fig.}~\ref{fig:angles}.

\begin{proof}
The proof is by contradiction. Fix $0<\epsilon<1$ sufficiently small
so that ${\rm cos}(\epsilon/2)<1-\epsilon^2/12$. Suppose that there
is a sequence of points $q_n$ tending to infinity such that for each
$n$ there are minimizing geodesics $\gamma_n$ and $\mu_n$ from $p$
to $q_n$ making angle at least $\epsilon$ at $q_n$. For each $n$ let
$d_n=d(p,q_n)$. By passing to a subsequence we can suppose that for
all $n$ and $m$ the cosine of the angle at $p$ between $\gamma_n$
and $\gamma_m$ at least $1-\epsilon^2/24$, and the cosine of the
angle at $p$ between $\mu_n$ and $\mu_m$ is at least
$1-\epsilon^2/24$. We can also assume that for all $n\ge 1$ we have
$d_{n+1}\ge (100/\epsilon^2) d_n$. Let $\delta_n=d(q_n,q_{n+1})$.
Applying the first Toponogov property at $p$, we see that
$\delta_n^2\le d_n^2+d_{n+1}^2-2d_nd_{n+1}(1-\epsilon^2/24)$.
Applying the same property at $q_n$ we have
$$d_{n+1}^2\le d_n^2+\delta_n^2-2d_n\delta_n{\rm cos}(\theta),$$ where
$\theta\le \pi$ is the angle at $q_n$ between $\gamma_n$ and a
minimal geodesic joining $q_n$ to $q_{n+1}$. Thus,
$${\rm cos}(\theta)\le
\frac{d_n-d_{n+1}(1-\epsilon^2/24)}{\delta_n}.$$ By the triangle
inequality (and the fact that $\epsilon<1$) we have $\delta_n\ge
(99/\epsilon)d_n$ and $\delta_n\ge d_{n+1}(1-(\epsilon^2/100))$.
Thus,
$${\rm cos}(\theta)\le
\epsilon^2/99-(1-\epsilon^2/24)/(1-(\epsilon^2/100))<-(1-\epsilon^2/12).$$
This implies that ${\rm cos}(\pi-\theta)>(1-\epsilon^2/12)$, which
implies that $\pi-\theta<\epsilon/2$. That is to say, the angle at
$q_n$ between $\gamma_n$ and a shortest geodesic from $q_n$ to
$q_{n+1}$ is between $\pi-\epsilon/2$ and $\pi$. By symmetry, the
same is true for the angle between $\mu_n$ and the same shortest
geodesic from $q_n$ to $q_{n+1}$. Thus, the angle between $\gamma_n$
and $\mu_n$ at $q_n$ is less than $\epsilon$, contradicting our
assumption.
\end{proof}

\begin{cor}\label{homeo}
Let $(M,g)$ be a complete, non-compact manifold of non-negative
sectional curvature. Let $p\in M$ and define a function $f\colon
M\to \Ar$ by $f(q)=d(p,q)$. Then there is $R<\infty$ such that for
$R\le s<s'$ we have:
\begin{enumerate}
\item $f^{-1}([s,s'])$ is homeomorphic to $f^{-1}(s)\times[s,s']$ and
in particular the level sets $f^{-1}(s)$ and $f^{-1}(s')$ are
homeomorphic; \item  $f^{-1}([s,\infty)$ is homeomorphic to
$f^{-1}(s)\times[s,\infty)$.
\end{enumerate}
\end{cor}

\begin{proof}
Given $(M,g)$ and $p\in M$ as in the statement of the corollary, choose  a
constant $R<\infty$ such that any two minimal geodesics from $p$ to a point $q$
with $d(p,q)\ge R/2$ make an angle at most $\pi/6$ at $q$. Now following
\cite{Petersen} p. 335, it is possible to find a smooth unit vector field $X$
on $U=M-\overline{B(p,R/2)}$ with the property that $f(\cdot)=d(p,\cdot)$ is
increasing along any integral curve for $X$ at a rate bounded below by ${\rm
cos}(\pi/3)$. In particular, for any $s\ge R$ each integral curve of $X$
crosses the level set $f^{-1}(s)$  in a single point. Using this vector field
we see that for any $s,s'>R$, the pre-image $f^{-1}([s,s'])$ is homeomorphic to
$f^{-1}(s)\times[s,s']$ and that the end $f^{-1}\left([s,\infty)\right)$ is
homeomorphic to $f^{-1}(s)\times [s,\infty)$.
\end{proof}

In a complete, non-compact $n$-manifold of positive curvature any soul is a
point.  While the proof of this result uses the same ideas as discussed above,
we shall not give a proof. Rather we refer the reader to Theorem 84 of
\cite{Petersen} on  p. 349. A soul has the property that if two minimal
geodesics emanate from $p$ and end at the same point $q\not= p$, then the angle
that they make at $q$ is less than $\pi/2$. Also, of course, the exponential
mapping is a diffeomorphism sufficiently close to the soul. Applying the above
lemma and a standard compactness argument, we see that in fact there is
$\epsilon>0$ such that all such pairs of minimal geodesics from $p$ ending at
the same point make angle less than $\pi/2-\epsilon$ at that point. Hence, in
this case there is a vector field $X$ on all of $M$ vanishing only at the soul,
and agreeing with the gradient of the distance function near the soul, so that
the distance function from $p$ is strictly increasing to infinity along each
flow line of $X$ (except the fixed point). Using $X$ one establishes that $M$
is diffeomorphic to $\Ar^n$. It also follows that all the level surfaces
$f^{-1}(s)$ for $s>0$ are homeomorphic to $S^{n-1}$ and for $0<s<s'$ the
preimage $f^{-1}([s,s'])$ is homeomorphic to $S^{n-1}\times[s,s']$.

There is an analogue of this result for the distance function from
any point, not just a soul.

\begin{cor}\label{htpytype}
Let $(M,g)$ be a complete, non-compact Riemannian $n$-manifold of positive
curvature. Then for any point $p\in M$ there is a constant $R=R(p)$ such that
for any $s<s'$ with $R\le s$ both $f^{-1}(s,s')$ and $f^{-1}(s,\infty)$ are
homotopy equivalent to $S^{n-1}$.
\end{cor}

\begin{proof}
Given $(M,g)$ and $p$ fix $R<\infty$ sufficiently large so that
Corollary~\ref{homeo} holds. Since $M$ is diffeomorphic to $\Ar^n$
it has only one end and hence the level sets $f^{-1}(s)$ for $s\ge
R$ are connected. Given any compact subset $K\subset M$ there is a
larger compact set $B$ (a ball) such that $M\setminus B$ has trivial
fundamental group and trivial homology groups $H_i$ for $i<n-1$.
Hence for any subset $Z\subset M\setminus B$, the inclusion of $Z\to
M\setminus K$ induces the trivial map on $\pi_1$ and on $H_i$ for
$i<n-1$. Clearly, for any $R\le s<b$ the inclusion
$f^{-1}(b,\infty)\to f^{-1}(s,\infty)$ is a homotopy equivalence.
Thus, it must be the case that $f^{-1}(b,\infty)$ has trivial
fundamental group and $H_i$ for $i<n-1$. Hence, the same is true for
$f^{-1}(s,\infty)$ for any $s\ge R$. Lastly, since
$f^{-1}(s,\infty)$ is connected and simply connected ( hence
orientable) and has two ends, it follows by the non-compact form of
Poincar\'e duality that $H_{n-1}(f^{-1}(s,\infty))\cong \Zee$.
Hence, by the Hurewicz theorem $f^{-1}(s,\infty)$ is homotopy
equivalent to $S^{n-1}$ for any $s\ge R$. Of course, it is also true
for $R\le s\le s'$ that $f^{-1}(s,s')$ is homotopy equivalent to
$S^{n-1}$.
\end{proof}

\section{Ends of a manifold}\label{ENDS}

Let us review the basic notions about ends of a manifold.

\begin{defn}
Let $M$ be a connected manifold.  Consider the inverse system of
spaces indexed by the compact, codimension-$0$ submanifolds
$K\subset M$, where the space associated to $K$ is the finite set
$\pi_0(M\setminus K)$ with the discrete topology. The inverse limit
of this inverse system is the {\em space of ends} of $M$. It is a
compact space. An {\em end}\index{end, of a manifold} of $M$ is a
point of the space of ends. An end ${\mathcal E}$ determines a
complementary component of each compact, codimension-$0$ submanifold
$K\subset M$, called a {\em neighborhood} of the end. Conversely, by
definition these neighborhoods are cofinal in the set of all
neighborhoods of the end. A sequence $\{x_n\}$ in $M$ {\em converges
to the end} ${\mathcal E}$ if it is eventually in every neighborhood
of the end. In fact, what we are doing is defining a topology on the
union of $M$ and its space of ends that makes this union a compact,
connected Hausdorff space which is a compactification of $M$.
\end{defn}

A proper map between topological manifolds induces a map on the
space of ends, and in fact induces a map on the compactifications
sending the subspace of ends of the compactification of the domain
to the subspace of ends of the compactification of the range.

We say that a path $\gamma\colon [a,b)\to M$ is a path {\em to the
end ${\mathcal E}$} if it is a proper map and it sends the end
$\{b\}$ of $[a,b)$ to the end ${\mathcal E}$ of $M$. This condition
is equivalent to saying that given a neighborhood $U$ of ${\mathcal
E}$ there is a neighborhood of the end $\{b\}$ of $[a,b)$ that maps
to $U$.

Now suppose that $M$ has a Riemannian metric $g$. Then we can
distinguish between ends at finite and infinite distance. An end is
at {\em finite distance} if there is a rectifiable path of finite
length to the end. Otherwise, the end is at infinite distance. If an
end is at finite distance we have the notion of the distance from a
point $x\in M$ to the end. It is the infimum of the lengths of
rectifiable paths from $x$ to the end. This distance is always
positive. Also, notice that the Riemannian manifold is complete if
and only if has no end at finite distance.

\section{The splitting theorem}

In this section we give a proof of the following theorem which is
originally due to Cheeger-Gromoll \cite{CheegerGromoll2}. The weaker
version giving the same conclusion under the stronger hypothesis of
non-negative sectional curvature (which is in fact all we need in
this work) was proved earlier by Toponogov, see \cite{Toponogov}.

\begin{thm}\label{splitting}
Suppose\index{splitting theorem|ii} that $M$ is complete, of
non-negative Ricci curvature and suppose that $M$ has at least two
ends. Then $M$ is isometric to a product $N\times \Ar$ where $N$ is
a compact manifold.
\end{thm}

\begin{proof}
We begin the proof by establishing a result of independent interest,
which was formulated as the main theorem in \cite{CheegerGromoll2}..

\begin{lem}\label{line}
Any complete Riemannian manifold $X$ of non-negative Ricci curvature
containing a minimizing line is isometric to a product $N\times \Ar$
for some Riemannian manifold $N$.
\end{lem}

\begin{proof}
Given a minimizing line $\lambda\colon \Ar\to X$, define
$\lambda_\pm\colon [0,\infty)\to X$ by $\lambda_+(t)=\lambda(t)$ and
$\lambda_-(t)=\lambda(-t)$. Then we have the Busemann functions
$B_+=B_{\lambda_+}$ and $B_-=B_{\lambda_-}$.
Proposition~\ref{Blambda} applies to both $B_+$ and $B_-$ and shows
that $\Delta (B_+ +B_-)\le 0$. On the other hand, using the fact
that $\lambda$ is distance minimizing, we see that for any $s,t>0$
and for any $x\in M$ we have $d(x,\lambda(t))+d(x,\lambda(-s))\ge
s+t$, and hence $B_+(x)+B_-(x)\ge 0$. Clearly, $B_+(x)+B_-(x)=0$ for
any $x$ in the image of $\lambda$. Thus, the function $B_++B_-$  is
everywhere $\ge 0$, vanishes at at a least one point and satisfies
$\Delta (B_++B_-)\le 0$ in the weak sense. This is exactly the
set-up for the maximum principle\index{maximum principle!scalar
functions}, cf. \cite{Petersen}, p. 279.

\begin{thm}{\bf (The Maximum Principle)}
Let $f$ be a real-valued continuous function on a connected
Riemannian manifold with $\Delta f\ge 0$ in the weak sense. Then $f$
is locally constant near any local maximum. In particular, if $f$
achieves its maximum then it is a constant.
\end{thm}

Applying this result to $-(B_++B_-)$, we see that $B_++B_-=0$, so that
$B_-=-B_+$. It now follows that $\Delta B_+=0$ in the weak sense. By standard
elliptic regularity results this implies that $B_+$ is a smooth harmonic
function.

Next, we show that  for all $x\in M$ we have $|\nabla B_+(x)|=1$.
Fix $x\in M$. Take a sequence $t_n$ tending to infinity and consider
minimizing geodesics $\mu_{+,n}$ from $x$ to $\lambda_+(t_n)$. By
passing to a subsequence we can assume that there is a limit as
$n\rightarrow \infty$. This limit  is a minimizing geodesic ray
$\mu_+$ from $x$, which we think of as being `asymptotic at
infinity' to $\lambda_+$. Similarly, we construct a minimizing
geodesic ray $\mu_-$ from $x$ asymptotic at infinity to $\lambda_+$.
 Since
$\mu_+$ is a minimal geodesic ray, it follows that for any $t$ the
restriction $\mu_+|_{[0,t]}$ is the unique length minimizing
geodesic from $x$ to $\mu_+(t)$ and that $\mu_+(t)$ is not a
conjugate point along $\mu_+$. It follows  by symmetry that $x$ is
not a conjugate point along the reversed geodesic $-\mu_+|_{[0,t]}$
and hence that $x\in U_{\mu_+(t)}$. This means that the function
$d(\mu_+(t),\cdot)$ is smooth at $x$ with gradient equal to the unit
tangent vector in the negative direction at $x$ to $\mu_+$, and
consequently that $B_{\mu_+,t}$ is smooth at $x$. Symmetrically, for
any $t>0$ the function $B_{\mu_-,t}$ is smooth at $x$ with the
opposite gradient. Notice that these gradients have norm one. We
have
$$B_{\mu_+,t}+B_+(x)\ge B_+=-B_-\ge -(B_{\mu_-,t}+B_-(x)).$$
Of course, $B_{\mu_+,t}(x)=0$ and $B_{\mu_-,t}(x)=0$, so that
$$B_{\mu_+,t}(x)+B_+(x)=-(B_{\mu_-,t}(x)+B_-(x)).$$
This squeezes $B_+$ between two smooth functions with the same value
and same gradient at $x$ and hence shows that  $B_+$ is $C^1$ at $x$
and $|\nabla B_+(x)|$ is of norm one.

Thus, $B$ defines a smooth Riemannian submersion from $M\to \Ar$
which implies that $M$ is isometric to a product of the fiber over
the origin with $\Ar$.
\end{proof}

This result together with Lemma~\ref{ends} shows  that if $M$
satisfies the hypothesis of  the theorem, then it can be written as
a Riemannian product $M=N\times \Ar$. Since $M$ has at least two
ends, it follows immediately that $N$ is compact. This completes the
proof of the theorem.
\end{proof}

\section{$\epsilon$-necks}

Certain types of (incomplete) Riemannian manifolds play an
especially important role in our analysis. The purpose of this
section is to introduce these manifolds and use them to prove one
essential result in Riemannian geometry.

For all of the following definitions we fix $0<\epsilon<1/2$. Set
$k$ equal to the greatest integer less than or equal to
$\epsilon^{-1}$. In particular, $k\ge 2$.

\begin{defn}\label{epsclose}
Suppose that we have a fixed metric $g_0$ on a manifold $M$ and an
open submanifold $X\subset M$. We say that another metric $g$ on $X$
is {\em within $\epsilon$ of $g_0|_X$ in the
$C^{[1/\epsilon]}$-topology} if, setting $k=[1/\epsilon]$ we have
\begin{equation}\label{eqnepsclose}{\rm sup}_{x\in X}\left(|g(x)-g_0(x)|_{g_0}^2+
\sum_{\ell=1}^k|\nabla_{g_0}^\ell g(x)|_{g_0}^2\right)<
\epsilon^2,\end{equation} where the covariant derivative
$\nabla^\ell_{g_0}$ is the Levi-Civita connection of  $g_0$ and
norms are  the pointwise $g_0$-norms on
$${\rm
Sym}^2T^*M\otimes \underbrace{T^*M\otimes\cdots\otimes
T^*M}_{\ell-{\rm times}}.$$

More generally, given two smooth families of metrics $g(t)$ and
$g_0(t)$ on  $M$ defined for $t$ in some interval $I$ we say that
the family $g(t)|_X$ is within $\epsilon$ of the family $g_0(t)|_X$
in the $C^{[1/\epsilon]}$-topology if  we have
$${\rm sup}_{(x,t)\in X\times I}\left(\left| g(x,t)
-g_0(x,t)\right|_{g_0(t)}^2 +\sum_{\ell=1}^k\left|\nabla_{g_0}^\ell
g(x,t)\right|_{g_0}^2\right)<\epsilon^2.$$
\end{defn}

\begin{rem}
Notice that if we view a one-parameter family of metrics $g(t)$ as a
curve in the space of metrics on $X$ with the
$C^{[1/\epsilon]}$-topology then this is the statement that the two
paths are pointwise within $\epsilon$ of each other. It says nothing
about the derivatives of the paths, or equivalently about the time
derivatives of the metrics and of their covariant derivatives. We
will always be considering paths of metrics satisfying the Ricci
flow equation. In this context two one-parameter families of metrics
that are close in the $C^{2k}$-topology exactly when  the $r^{th}$
time derivatives of the $s^{th}$-covariant derivatives are close for
all $r,s$ with $s+2r\le 2k$.
\end{rem}

The first object of interest is one that, up to scale, is close to a
long, round cylinder.

\begin{defn}\label{epsneck}
Let $(N,g)$ be a Riemannian manifold and $x\in N$ a point. Then {\em
an $\epsilon$-neck\index{$\epsilon$-neck|ii} structure on $(N,g)$
centered at $x$} consists of a diffeomorphism
$$\varphi\colon  S^2\times (-\epsilon^{-1},\epsilon^{-1})\to N,$$
with $x\in \varphi(S^2\times\{0\})$, such that the metric
$R(x)\varphi^*g$  is within $\epsilon$ in the
$C^{[1/\epsilon]}$-topology of the product of the usual Euclidean
metric on the open interval with the metric of constant Gaussian
curvature $1/2$ on $S^2$. We also use the terminology {\em $N$ is an
$\epsilon$-neck centered at $x$}.  The image under $\varphi$ of the
family of submanifolds $S^2\times \{t\}$ is called the {\em family
of $2$-spheres of the $\epsilon$-neck}. The submanifold
$\varphi(S^2\times \{0\})$ is called {\em the central $2$-sphere} of
the $\epsilon$-neck structure.  We denote by $s_N\colon N\to \Ar$
the composition $p_2\circ \varphi^{-1}$, where $p_2$ is the
projection of $S^2\times (-\epsilon^{-1},\epsilon^{-1})$ to the
second factor. There is also the vector field $\partial/\partial
s_N$ on $N$ which is $\varphi_*$ of the standard vector field in the
interval-direction of the product. We also use the terminology of
the {\em plus} and {\em minus} end of the $\epsilon$-neck in the
obvious sense. The opposite (or reversed) $\epsilon$-neck structure
is the one obtained by composing the structure map with ${\rm
Id}_{S^2}\times -1$. We define the
 {\em positive half of the neck} to be the region
$s_N^{-1}(0,\epsilon^{-1})$ and the {\em negative half} to be the
region $s_N^{-1} (-\epsilon^{-1},0)$. For any other fraction, e.g.,
the left-hand three-quarters, the right-hand one-quarter, there are
analogous notions, all measured with respect to $s_N\colon N\to
(-\epsilon^{-1},\epsilon^{-1})$. We also use the terminology the
middle one-half, or middle one-third of the $\epsilon$-neck; again
these regions have their obvious meaning when measured via $s_N$.

{\em An $\epsilon$-neck} in a Riemannian manifold $X$ is a
codimension-zero submanifold $N$ and an $\epsilon$-structure on $N$
centered at some point $x\in N$.

The {\em scale}\index{$\epsilon$-neck!scale} of an $\epsilon$-neck $N$ centered
at $x$ is $R(x)^{-1/2}$. The scale of $N$ is denoted $r_N$. Intuitively, this
is a measure of the radius of the cross-sectional $S^2$ in the neck. In fact,
the extrinsic diameter of any $S^2$ factor in the neck is close to $\sqrt{2}\pi
r_N$. See \textsc{Fig.}~\ref{fig:epsneck} in the introduction.
\end{defn}

Here is the result that will be so important in our later arguments.

\begin{prop}\label{narrows} The following
holds for any $\epsilon>0$ sufficiently small.
 Let $(M,g)$ be a complete, positively curved Riemannian
$3$-manifold. Then $(M,g)$ does not contain $\epsilon$-necks of
arbitrarily small scale.
\end{prop}

\begin{proof}
The result is obvious if $M$ is compact, so we assume that $M$ is
non-compact. Let $p\in M$ be a soul for $M$ (Theorem~\ref{soul}),
and let $f$ be the distance function from $p$. Then $f^{-1}(s)$ is
connected for all $s>0$.

\begin{lem}\label{neckseparate}
Suppose that $\epsilon>0$ is sufficiently small that Lemma~\ref{topiso}  from
the appendix holds. Let $(M,g)$ be a non-compact $3$-manifold of positive
curvature and let $p\in M$ be a soul for it. Then for any $\epsilon$-neck $N$
disjoint from $p$
 the central $2$-sphere of $N$  separates the soul from the end of the manifold. In particular, if
 two $\epsilon$-necks $N_1$ and $N_2$ in $M$ are disjoint from each other and from $p$, then
the central $2$-spheres of $N_1$ and $N_2$ are the boundary
components of a region in $M$ diffeomorphic to $S^2\times I$.
\end{lem}

\begin{proof}
Let $N$ be an $\epsilon$-neck disjoint from $p$. By
Lemma~\ref{topiso} for any point $z$ in the middle third of $N$, the
boundary of the metric ball $B(p,d(z,p))$ is a topological
$2$-sphere in $N$ isotopic in $N$ to the central $2$-sphere of $N$.
Hence, the central $2$-sphere separates the soul from the end of
$M$. The second statement follows immediately by applying this to
$N_1$ and $N_2$.
\end{proof}

Let $N_1$ and $N_2$ be disjoint $\epsilon$-necks, each disjoint from
the soul. By the previous lemma, the central $2$-spheres $S_1$ and
$S_2$ of these necks are smoothly isotopic to each other and they
are the boundary components of a region  diffeomorphic to $S^2\times
I$. Reversing the indices if necessary we can assume that $N_2$ is
closer to $\infty$ than $N_1$, i.e., further from the soul.
Reversing the directions of the necks if necessary, we can arrange
that for $i=1,2$ the function $s_{N_i}$ is increasing as we go away
from the soul. We define $C^\infty$- functions $\psi_i$ on $N_i$,
functions depending only on $s_{N_i}$, as follows. The function
$\psi_1$ is zero on the negative side of the middle third of $N_1$
and increases to be identically one on the positive side of the
middle third. The function $\psi_2$ is one on the negative side of
the middle third of $N_2$ and decreases to be zero on the positive
side. We extend $\psi_1,\psi_2$ to a function $\psi$ defined on all
of $M$ by requiring that it be identically one on the region $X$
between $N_1$ and $N_2$ and to be identically zero on $M\setminus
(N_1\cup X\cup N_2)$.

Let $\lambda$ be a geodesic ray from the soul of $M$ to infinity,
and $B_\lambda$ its Busemann function. Let $N$ be any
$\epsilon$-neck disjoint from the soul, with $ s_N$ direction chosen
so that it points away from the soul. At any point of the middle
third of $N$ where $B_\lambda$ is smooth,  $\nabla B_\lambda$ is a
unit vector in the direction of the unique minimal geodesic ray from
the end of $\lambda$ to this point. Invoking Lemma~\ref{directions}
from the appendix we see that at such points $\nabla B_\lambda$ is
close to $-R(x)^{1/2}\partial/\partial s_N$, where $x\in N$ is the
center of the $\epsilon$-neck. Since $\nabla B_\lambda$ is $L^2$ its
non-smooth points have measure zero and hence, the restriction of
$\nabla B_\lambda$ to the middle third of $N$ is close in the
$L^2$-sense to $-R(x)^{1/2}\partial/\partial s_N$.

Applying this to $N_1$ and $N_2$ we see that
\begin{equation}\label{locallabel}
\int_M\langle \nabla B_\lambda,\nabla\psi\rangle d{\rm
vol}=\left(\alpha_2R(x_2)^{-1}-\alpha_1R(x_1)^{-1}\right){\rm
Vol}_{h_0}(S^2)),\end{equation} where $h(0)$ is the round metric of
scalar curvature $1$ and where each of $\alpha_1$ and $\alpha_2$
limits to $1$ as $\epsilon$ goes to zero. Since $\psi\ge 0$,
Proposition~\ref{Blambda} tells us that the left-hand side of
Equation~(\ref{locallabel}) must be $\ge 0$. This shows that,
provided that $\epsilon$ is sufficiently small, $R(x_2)$ is bounded
above by $2R(x_1)$. This completes the proof of the proposition.
\end{proof}

\begin{cor}\label{strnarrows}
Fix $\epsilon>0$ sufficiently small so that Lemma~\ref{topiso} holds. Then
there is a constant $C<\infty$ depending on $\epsilon$ such that the following
holds. Suppose that $M$ is a non-compact $3$-manifold of positive sectional
curvature. Suppose that $N$ is an $\epsilon$-neck in $M$  centered at a point
$x$ and disjoint from a soul $p$ of $M$. Then for any $\epsilon$-neck $N'$ that
is separated from $p$ by $N$ with center $x'$ we have $R(x')\le CR(x)$.
\end{cor}

\section{Forward difference quotients}

Let us review quickly some standard material on forward difference quotients.

Let $f\colon [a,b]\to \Ar$ be a continuous function on an interval.
We say that the {\em forward difference quotient of $f$ at a point
$t\in [a,b)$, denoted $\frac{df}{dt}(t)$, is less than
$c$}\index{forward difference quotient|ii} provided that
$$\overline{\rm lim}_{\triangle t\rightarrow
0^+}\frac{f(t+\triangle t)-f(t)}{\triangle t}\le c.$$ We say that it
is greater than or equal to $c'$ if $$c'\le  \underline{\rm
lim}_{\triangle t\rightarrow 0^+}\frac{f(t+\triangle
t)-f(t)}{\triangle t}.$$ Standard comparison arguments show:

\begin{lem}\label{fordiffquot}
Suppose that $f\colon [a,b]\to \Ar$ is a continuous function. Suppose that
$\psi$ is a $C^1$-function on $[a,b]\times \Ar$ and suppose that
$\frac{df}{dt}(t)\le \psi(t,f(t))$ for  every $t\in [a,b)$ in the sense of
forward difference quotients. Suppose also that there is a function  $G(t)$
defined on $[a,b]$ that satisfies the differential equation
$G'(t)=\psi(t,G(t))$ and has $f(a)\le G(a)$. Then $f(t)\le G(t)$ for all $t\in
[a,b]$.
\end{lem}

The application we shall make of these results is the following.

\begin{prop}\label{fordiffmax}
Let $M$ be a smooth manifold with a smooth vector field $\chi$ and a
smooth function ${\bf t}\colon M\to [a,b]$ with $\chi({\bf t})=1$.
Suppose also that $F\colon M\to \Ar$ is a smooth function with the
properties:
\begin{enumerate}
\item  for each $t_0\in [a,b]$ the restriction of $F$ to the level set ${\bf
t}^{-1}(t_0)$ achieves its maximum, and \item the subset ${\mathcal Z}$ of $M$
consisting of all $x$ for which $F(x)\ge F(y)$ for all $y\in {\bf t}^{-1}({\bf
t}(x))$ is a compact set.
\end{enumerate}
 Suppose also that at each $x\in {\mathcal Z}$ we have $\chi(F(x))\le \psi({\bf
t}(x),F(x))$. Set $F_{\rm max}(t)={\rm max}_{x\in {\bf
t}^{-1}(t)}F(x)$. Then $F_{\rm max}(t)$ is a continuous function and
$$\frac{dF_{\rm max}}{dt}(t)\le \psi(t,F_{\rm max}(t))$$
in the sense of forward difference quotients.
 Suppose that $G(t)$ satisfies the
differential equation
$$G'(t)=\psi(t,G(t))$$ and has initial condition $F_{\rm
max}(a)\le G(a)$.  Then for all $t\in [a,b]$ we have
$$F_{\rm max}(t)\le G(t).$$
\end{prop}

\begin{proof}
Under the given hypothesis it is a standard and easy exercise to
establish the statement about the forward difference quotient of
$F_{\rm max}$. The second statement then is an immediate corollary
of the previous result.
\end{proof}

\chapter{Basics of Ricci flow}\label{flowbasics}
 In this chater we introduce the Ricci flow equation due to R.
Hamilton \cite{Hamilton3MPRC}. For the basic material on the Ricci flow
equation see \cite{ChowKnopf}.

\section{The definition of the Ricci flow}
\quad
\begin{defn}
The {\em Ricci flow equation}\index{Ricci flow!equation|ii} is the
following evolution equation for a Riemannian metric:
\begin{equation}\label{Ricciflow}
\frac{\partial g}{\partial t}= -2{\rm Ric}(g). \end{equation} A solution to
this equation (or a {\em Ricci flow})\index{Ricci flow|ii} is a one-parameter
family of metrics $g(t)$, parameterized by $t$ in a non-degenerate interval
$I$, on a smooth manifold $M$ satisfying Equation~(\ref{Ricciflow}). If $I$ has
an initial point $t_0$ then $(M,g(t_0))$ is called  {\sl the initial condition
of} or {\em the initial metric for}\index{Ricci flow!initial metric for}  the
Ricci flow (or of the solution).
\end{defn}

Let us give a quick indication of what the Ricci flow equation
means. In harmonic coordinates $(x^1,\ldots,x^n)$ about $p$,
 that is to say coordinates where $\triangle x^i=0$ for all $i$,
we have
\[ {\rm Ric}_{ij}={\rm Ric}(\frac{\partial}{\partial x^i},\frac{\partial}{\partial x^j})=
-\frac{1}{2}\triangle g_{ij} + Q_{ij}(g^{-1},\partial g)\]
where $Q$ is a quadratic form in $g^{-1}$ and $\partial g$, and so in particular is a lower order term
in the derivatives of $g$. See Lemma 3.32 on page 92 of \cite{ChowKnopf}.
So, in these coordinates, the Ricci flow equation is actually a heat
equation  for the Riemannian metric
\[ \frac{\partial}{\partial t}g =
\triangle g +2Q(g^{-1},\partial g). \]

\begin{defn}\label{defnRicci}
We introduce some notation that will be used throughout. Given a
Ricci flow $(M^n,g(t))$ defined for $t$ contained in an interval
$I$, then the {\sl space-time} for this flow is  $M\times I$.  The
{\sl $t$ time-slice} of space-time is the Riemannian manifold
$M\times\{t\}$ with the Riemannian metric $g(t)$. Let ${\mathcal
HT}(M\times I)$ be the {\em horizontal tangent
bundle}\index{horizontal tangent space} of space-time, i.e., the
bundle of tangent vectors to the time-slices. It is a smooth,
rank-$n$ subbundle the tangent bundle of space-time. The evolving
metric $g(t)$ is then a smooth section of ${\rm Sym}^2{\mathcal
HT}^*(M\times I)$. We denote points of space-time as pairs $(p,t)$.
Given $(p,t)$ and any $r>0$ we denote by $B(p,t,r)$ the metric ball
of radius $r$ centered at $(p,t)$ in the $t$ time-slice. For any
$\Delta t>0$ for which $[t-\Delta t,t]\subset I$, we define the {\sl
backwards parabolic neighborhood}\index{parabolic neighborhood}
$P(x,t,r,-\Delta t)$ to be the product $B(x,t,r)\times
[t-\Delta t,t]$ in space-time. Notice that the intersection of
$P(x,t,r,-\Delta t)$ with a time-slice other that the $t$ time-slice
need not be a metric ball in that time-slice. There is the
corresponding notion of a forward parabolic neighborhood
$P(x,t,r,\Delta t)$ provided that $[t,t+\Delta t]\subset I$.
\end{defn}

\section{Some exact solutions to the Ricci flow}

\subsection{Einstein manifolds}
Let $g_0$ be an Einstein\index{Einstein manifold|ii} metric: ${\rm
Ric}(g_0) = \lambda g_0$, where $\lambda$ is a constant. Then for
any positive constant $c$, setting $g = cg_0$ we have ${\rm Ric}(g)
= {\rm Ric}(g_0) = \lambda g_0= \frac{\lambda}{c} g.$ Using this we
can construct solutions to the Ricci flow equation as follows.
Consider $g(t) = u(t)g_{0}$. If this one-parameter family of
metrics is a solution of the Ricci flow, then
\begin{align*}
\frac{\partial g}{\partial t}& =u'(t)g_{0}\\& = -2{\rm Ric}(u(t)g_{0})\\&
= -2{\rm Ric}(g_{0})\\& = -2\lambda g_{0}.
\end{align*}
So $u'(t)= -2\lambda$,  and hence $u(t) = 1-2 \lambda t$. Thus
$g(t)= (1-2 \lambda t)g_{0}$ is a solution of the Ricci flow. The
cases $\lambda
>0, \lambda=0,$ and $\lambda< 0$ correspond to \textit{shrinking}, \textit{steady}
and \textit{expanding} solutions. Notice that in the shrinking case the
solution exists for $t\in [0,\frac{1}{2\lambda})$ and goes singular at
$t=\frac{1}{2\lambda}$.

\begin{exam}
The standard metric on each of $S^{n}, \Ar^{n},$ and ${\mathbb
H}^{n}$ is Einstein. Ricci flow is contracting on $S^n$, constant on
$\Ar^n$, and expanding on ${\mathbb H}^n$. The Ricci flow on $S^n$
has a finite-time singularity where the diameter of the manifold
goes to zero and the curvature goes uniformly to $+\infty$. The
Ricci flow on ${\mathbb H}^n$ exists for all $t\ge 0$ and as $t$
goes to infinity the distance between any pair of points grows
without bound and the curvature goes uniformly to zero.
\end{exam}

\begin{exam}
$\mathbb{C}P^{n}$ equipped with the Fubini-Study metric, which is induced
from the standard metric of $S^{2n+1}$ under the Hopf fibration with the fibers
of great circles, is Einstein.
\end{exam}

\begin{exam}
Let $h_0$ be the round metric on $S^2$ with constant Gausssian
curvature $1/2$. Set $h(t)=(1-t)h_0$. Then the flow $$(S^2,h(t)),\
-\infty<t<1,$$ is a Ricci flow. We also have the product of this flow
with the trivial flow on the line: $(S^2\times \Ar,h(t)\times
ds^2),\ -\infty<t<1$. This is called the {\em standard shrinking
round cylinder}\index{shrinking round cylinder}.
\end{exam}

The standard shrinking round cylinder is a model for evolving
$\epsilon$-necks. In Chapter~\ref{sectprelim} we introduced the
notion of an $\epsilon$-neck. In the case of flows in order to take
smooth geometric limits, it is important to have a stronger version
of this notion. In this stronger notion, the neck not only exists in
one time-slice but it exists backwards in the flow for an
appropriate amount of time and is close to the standard shrinking
round cylinder on the entire time interval. The existence of
evolving necks
 is exploited when we study limits of Ricci flows.

\begin{figure}[ht]
  \relabelbox{
  \centerline{\epsfbox{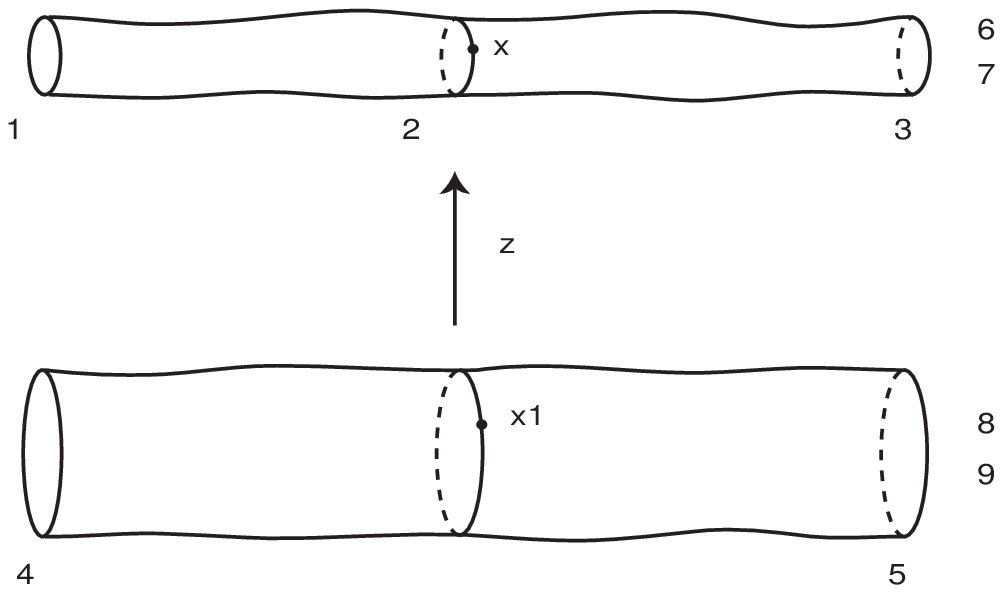}}}
  \relabel{1}{$-\epsilon^{-1}$}
  \relabel{2}{$S^2\times\{0\}$}
  \relabel{3}{$\epsilon^{-1}$}
  \relabel{4}{$-\epsilon^{-1}$}
  \relabel{5}{$\epsilon^{-1}$}
  \relabel{6}{$t=0$}
  \relabel{7}{$R\sim 1$}
  \relabel{8}{$t=-1$}
  \relabel{9}{$R\sim \frac 12$}
  \relabel{x}{$x$}
  \relabel{x1}{$x$}
  \relabel{z}{Ricci flow}
  \endrelabelbox
  \caption{Strong $\epsilon$-neck of scale 1.}\label{fig:strneck}
\end{figure}

\begin{defn}\label{strongepsneck}
Let $(M,g(t))$ be a Ricci flow. An {\em evolving $\epsilon$-neck
centered at $(x,t_0)$ and defined for rescaled time
$t_1$}\index{$\epsilon$-neck!evolving|ii} is an $\epsilon$-neck
$$\varphi\colon S^2\times
(-\epsilon^{-1},\epsilon^{-1})\buildrel\cong\over\longrightarrow
N\subset (M,g(t))$$ centered at $(x,t_0)$ with the property that
pull-back via $\varphi$ of the family of metrics $R(x,t_0)g(t')|_N,\
-t_1< t'\le 0$, where $t_1=R(x,t_0)^{-1}(t-t_0)$, is within
$\epsilon$ in the $C^{[1/\epsilon]}$-topology of the product of the
standard metric on the interval with evolving round metric on $S^2$
with scalar curvature $1/(1-t')$ at time $t'$. A {\em strong
$\epsilon$-neck centered at
$(x,t_0)$}\index{$\epsilon$-neck!strong|ii} in a Ricci flow is an
evolving $\epsilon$-neck centered at $(x,t_0)$ and defined for
rescaled time $1$, see {\sc Fig.}~\ref{fig:strneck}.
\end{defn}

\subsection{Solitons}

A {\sl Ricci soliton}\index{Ricci flow!soliton}\index{soliton|see{Ricci flow}}
is a Ricci flow $(M,g(t)),\ 0\le t<T\le \infty$, with the property that for
each $t\in [0,T)$ there is a diffeomorphism $\varphi_t\colon M\to M$ and a
constant $\sigma(t)$ such that $\sigma(t)\varphi_t^*g(0)=g(t)$. That is to say,
in a Ricci soliton all the Riemannian manifolds $(M,g(t))$ are isometric up to
a scale factor that is allowed to vary with $t$. The soliton is said to be {\em
shrinking} if $\sigma'(t)<0$ for all $t$. One way to generate Ricci solitons is
the following: Suppose that we have a vector field $X$ on $M$ and a constant
$\lambda$ and a metric $g(0)$ such that
\begin{equation}\label{soliton}
-{\rm Ric}(g(0))=\frac{1}{2}{\mathcal L}_Xg(0)-\lambda g(0).
\end{equation} We set $T=\infty$ if $\lambda\le 0$ and equal to $(2\lambda)^{-1}$
if $\lambda>0$. Then, for all $t\in [0,T)$ we define a function
$$\sigma(t)=1-2\lambda t,$$
and a vector field
$$Y_t(x)=\frac{X(x)}{\sigma(t)}.$$
 Then we define $\varphi_t$ as the one-parameter
family of diffeomorphisms generated by the time-dependent vector
fields $Y_t$.
\begin{claim}\label{solitonclaim}
 The flow
$(M,g(t)),\ 0\le t<T$, where $g(t)=\sigma(t)\varphi_t^*g(0)$, is a soliton. It
is a shrinking soliton if $\lambda>0$.
\end{claim}

\begin{proof}
We check that this flow satisfies the Ricci flow equation; from
that, the result follows immediately. We have
\begin{eqnarray*} \frac{\partial g(t)}{\partial t} & = &
\sigma'(t)\varphi_t^*g(0)+\sigma(t)\varphi_t^*{\mathcal L}_{Y(t)}g(0)\\
& = & \varphi_t^*(-2\lambda+{\mathcal L}_X)g(0)\\
& = & \varphi_t^*(-2{\rm Ric}(g(0))) =-2{\rm Ric}(\varphi_t^*(g(0))).
\end{eqnarray*}
Since ${\rm Ric}(\alpha g)={\rm Ric}(g)$ for any $\alpha>0$, it follows that
$$\frac{\partial g(t)}{\partial t} = -2{\rm Ric}(g(t)).$$
\end{proof}

There is one class of shrinking solitons which are of special
importance to us. These are the gradient shrinking solitons.

\begin{defn}
A shrinking soliton $(M,g(t)),\ 0\le t<T$, is said to be a {\em
gradient shrinking soliton}\index{gradient shrinking soliton|ii} if
the vector field $X$ in Equation~(\ref{soliton}) is the gradient of
a smooth function $f$ on $M$.
\end{defn}

\begin{prop}\label{GSSgeneration}
 Suppose we have a complete Riemannian manifold $(M,g(0))$, a
smooth function $f\colon M\to \Ar$, and a constant $\lambda>0$ such
that
\begin{equation}\label{GSS}
-{\rm Ric}(g(0))={\rm Hess}(f)-\lambda g(0). \end{equation}
 Then there is $T>0$ and a gradient shrinking soliton
$(M,g(t))$ defined for  $0\le t<T.$ \end{prop}

\begin{proof}
Since
$${\mathcal L}_{\nabla f}g(0)=2Hess(f),$$
 Equation~(\ref{GSS}) is the soliton equation,
Equation~(\ref{soliton}), with the vector field $X$ being the
gradient vector field $\nabla f$. It is a shrinking soliton by
assumption since $\lambda>0$.
\end{proof}

\begin{defn}
In this case we say that $(M,g(0))$ and $f\colon M\to \Ar$ {\em
generate} a gradient shrinking soliton.
\end{defn}

\section{Local existence and uniqueness}\label{3.3} The following is the
first basic result in the theory -- local existence and uniqueness for Ricci
flow in the case of compact manifolds.

\begin{thm}{\bf (Hamilton, cf.
\cite{Hamilton3MPRC}.)}\label{compuniq} Let $(M,g_0)$ be a compact
Riemannian manifold of dimension $n$.
\begin{enumerate}
\item There is a $T>0$ depending on $(M,g_0)$ and  a Ricci flow
$(M,g(t)),\  0\le t<T$, with $g(0)=g_0$.
\item Suppose that we have Ricci flows with initial conditions $(M,g_0)$ at time
$0$ defined respectively on time intervals $I$ and $I'$. Then these
flows agree on $I\cap I'$. \end{enumerate}
\end{thm}

\smallskip

We remark that the Ricci flow is a weakly parabolic system where
degeneracy comes from the gauge invariance of the equation under
diffeomorphisms. Therefore the short-time existence does not come
from general theory. R. Hamilton's original proof of the short-time
existence was involved and used the Nash-Moser inverse function
theorem, \cite{HamiltonIFTNM}. Soon after, DeTurck \cite{DeTurck}
substantially simplified the short-time existence proof by breaking
the diffeomorphism invariance of the equation. For the reader's
convenience, and also because in establishing the uniqueness for
Ricci flows from the standard solution in Section~\ref{12.4}
we use a version of this idea
in the non-compact case, we sketch DeTurck's argument.

\begin{proof}
Let's sketch the proof due to DeTurck \cite{DeTurck}, cf, Section
3 of Chapter 3 starting on page 78 of \cite{ChowKnopf} for more details. First, we compute the first variation
at a Riemannian metric $g$ of minus twice the Ricci curvature tensor in the direction $h$:
\begin{align*}
\delta_g(-2{\rm Ric})(h)=\triangle h-{\rm Sym}(\nabla V)+S
\end{align*}
where:
\begin{enumerate}
\item $V$ is the one-form given by
$$V_k=\frac{1}{2}g^{pq}(\nabla_ph_{qk}+\nabla_qh_{pk}-\nabla_kh_{pq}),$$
\item ${\rm Sym}(\nabla V)$ is the symmetric two-tensor obtained by
symmetrizing the covariant derivative of
$V$, and
\item  $S$ is a symmetric two-tensor constructed from the inverse of the metric,
the Riemann curvature tensor and $h$, but involves no derivatives of $h$.
\end{enumerate}

Now let $g_0$ be the initial metric. For any metric $g$ we define a one-form
$\hat W$ by taking the trace, with respect to $g$, of the matrix-valued
one-form that is the difference of the connections of $g$ and $g_0$. Now we
form a second-order operator of $g$ by setting
$$P(g)={\mathcal L}_Wg,$$
the Lie derivative of $g$ with respect to the vector field $W$ dual to $\hat
W$. Thus, in local coordinates we have
$P(g)_{ij}=\nabla_i\hat W_j+\nabla_j\hat W_i$.
The linearization at $g$ of the second-order operator $P$ in the direction
$h$ is symmetric and is given by
$$\delta_gP(h)={\rm Sym}(\nabla V)+T$$
where $T$ is a first-order operator in $h$.
Thus, defining $Q=-2{\rm Ric}+P$ we have
$$\delta_g(Q)(h)=\triangle h+U$$
where $U$ is a first-order operator in $h$.
Now we introduce the Ricci-DeTurck  flow

\begin{align}\label{brokeneqn}
\frac{\partial g}{\partial t}&=-2{\rm Ric}(g)+P.
\end{align}
The computations above show that the Ricci-DeTurck flow is  strictly parabolic.
 Thus, Equation~(\ref{brokeneqn})  has a
short-time solution $\bar g(t)$ with $\bar g(0)=g_0$ by the standard PDE theory.
Given this solution $\bar g(t)$ we
define the time-dependent vector field $W(t)=W(\bar g(t),g_0)$ as above.
 Let $\phi_t$ be a one-parameter family of diffeomorphisms, with $\phi_0={\rm Id}$,
generated by this time-dependent vector field, i.e.,
$$\frac{\partial\phi_t}{\partial t}=W(t).$$
Then, direct computation shows that $g(t)=\phi_t^*\bar g(t)$ solves the Ricci flow equation.
\end{proof}

In performing surgery at time $T$, we will have an open submanifold $\Omega$ of
the compact manifold with the following property.  As $t$ approaches $T$ from
below, the metrics $g(t)|_\Omega$ converge smoothly to a limiting metric $g(T)$
on $\Omega$. We will `cut away' the rest of the manifold $M\setminus \Omega$
where the metrics are not converging and glue in a piece $E$ coming from the
standard solution to form a new compact manifold $M'$. Then we extend the
Riemannian metric $g(T)$ on $\Omega$ to one defined on $M'=\Omega\cup E$. The
resulting Riemannian manifold forms the initial manifold at time $T$ for
continuing the Ricci flow $\tilde g(t)$ on an interval $T\le t<T'$. It is
important to know that the two Ricci flows $(\Omega,g(t)),\ t\le T$ and
$(\Omega, \tilde g(t)),\ T\le t<T'$ glue together to make a smooth solution
spanning across the surgery time $T$. That this is true is a consequence of the
following elementary result.

\begin{prop}\label{patch}
 Suppose that $(U,g(t)),\ a\le t<b$, is a Ricci flow and suppose that
 there is a Riemannian metric $g(b)$ on $U$ such that
 as $t\rightarrow b$ the metrics $g(t)$ converge in the $C^\infty$-topology, uniformly
 on compact subsets, to $g(b)$. Suppose also that $(U,g(t)),\ b\le t<c$, is a
 Ricci flow. Then the one-parameter family of metrics $g(t), \ a\le
 t<c$, is a $C^\infty$-family and is a solution to the Ricci flow
 equation on the entire interval $[a,c)$.
 \end{prop}


\section{Evolution of curvatures}\label{evolcurv}

Let\index{curvature!evolution under Ricci flow} us fix a set
$(x^1,\ldots,x^n)$ of local coordinates. The Ricci flow equation,
written in local coordinates
\[ \frac{\partial g_{ij}}{\partial t} =
-2{\rm Ric}_{ij}
\]  implies a heat equation for the Riemann
curvature tensor $R_{ijkl}$ which we now derive. Various second-order
derivatives of the curvature tensor are likely to differ by
terms quadratic in the curvature tensors. To this end we introduce
the tensor $$B_{ijkl} = g^{pr}g^{qs}R_{ipjq}R_{krls}.$$ Note that we
have the obvious symmetries
$$B_{ijkl} = B_{jilk} = B_{klij},$$ but the other symmetries of
the curvature tensor $R_{ijkl}$ may fail to hold for $B_{ijkl}$.

\smallskip

\begin{thm}\label{thmRmevol}
The curvature tensor $R_{ijkl}$, the Ricci curvature ${\rm Ric}_{ij}$, the scalar
curvature $R$, and the volume form $d{\rm vol}(x,t)$ satisfy the following evolution
equations under Ricci flow:
\begin{eqnarray}
\frac{\partial R_{ijkl}}{\partial t} &= & \triangle R_{ijkl} +2(B_{ijkl} -
B_{ijlk}-B_{iljk} + B_{ikjl})\nonumber \\&  &  - g^{pq}(R_{pjkl}{\rm Ric}_{qi} +
R_{ipkl}{\rm Ric}_{qj} + R_{ijpl}{\rm Ric}_{qk} +
R_{ijkp}{\rm Ric}_{ql}) \label{Rmevol} \\
\label{Ricevol} \frac{\partial}{\partial t}{\rm Ric}_{jk} & = & \triangle
{\rm Ric}_{jk}+2g^{pq}g^{rs}R_{pjkr}{\rm Ric}_{qs}-2g^{pq}{\rm Ric}_{jp}{\rm Ric}_{qk} \\
\label{Revol}\frac{\partial}{\partial t}R & = & \Delta R+2|{\rm
Ric}|^2 \\
\label{dvolevol} \frac{\partial}{\partial t}d{\rm vol}(x,t) & = &
-R(x,t)d{\rm vol}(x,t).
\end{eqnarray}
\end{thm}

These equations are contained in Lemma 6.15  on page 179, Lemma 6.9
 on page 176,  Lemma 6.7 on page 176, and Equation (6.5) on
page 175 of \cite{ChowKnopf}, respectively.

Let us derive some consequences of these evolution equations.
The first result is obvious from the Ricci flow equation and will be used implicitly throughout
the paper.

\begin{lem}\label{posdistdec}
Suppose that $(M,g(t)),\ a<t<b$ is a Ricci flow of non-negative Ricci curvature
with $M$ a connected manifold.
Then for any points $x,y\in M$ the function $d_{g(t)}(x,y)$ is a non-increasing function of $t$.
\end{lem}

\begin{proof}
The Ricci flow equation tells us that non-negative Ricci curvature implies that
$\partial g/\partial t\le 0$. Hence, the length of any tangent vector in $M$, and consequently the length
of any path in $M$, is a non-increasing function of $t$. Since the distance between points is the infimum
over all rectifiable paths from $x$ to $y$ of the length of the path, this function is also a non-increasing function
of $t$.
\end{proof}

\begin{lem}
Suppose that $(M,g(t)),\ 0\le t\le T$, is a Ricci flow and
$|Rm(x,t)|\le K$ for all $x\in M$ and all $t\in [0,T]$. Then there
are constants $A,A'$  depending on $K,T$ and the dimension such
that:
\begin{enumerate}
\item For any non-zero tangent vector $v\in T_xM$ and any $t\le T$ we have
$$A^{-1}\langle v,v\rangle_{g(0)}\le \langle v,v\rangle_{g(t)}\le A\langle
v,v\rangle_{g(0)}.$$
\item For any open subset $U\subset M$ and any $t\le T$ we have
$$(A')^{-1}{\rm Vol}_0(U)\le {\rm Vol}_t(U)\le A'{\rm Vol}_0(U).$$
\end{enumerate}
\end{lem}

\begin{proof}
The Ricci flow equation yields
$$\frac{d}{dt}\left(\langle v,v\rangle_{g(t)}\right)=-2{\rm Ric}(v,v).$$
The bound on the Riemann curvature gives a bound on ${\rm Ric}$.
Integrating yields the result. The second statement is proved
analogously using Equation~(\ref{dvolevol}).
\end{proof}

\section{Curvature evolution in an evolving orthonormal
frame}\label{evolframe}

It is often best to study the evolution of the representative of the
tensor in an orthonormal frame $F$. Let $(M,g(t)),\ 0\le t<T$, be a
Ricci flow, and suppose that ${\mathcal F}$ is a frame on an open
subset $U\subset M$ consisting of vector fields $\{F_{1},
F_{2},\cdots,F_{n}\}$ on $U$ that are $g(0)$-orthonormal at every
point. Since the metric evolves by the Ricci flow, to keep the frame
orthonormal we must evolve it by an equation involving Ricci
curvature. We evolve this local frame  according to the formula
\begin{equation}\label{frameevol} \frac{\partial F_{a}}{\partial t}
= {\rm Ric}(F_{a},\cdot)^*,\end{equation} i.e.,  assuming that in local
coordinates $(x^1,\ldots,x^n)$, we have
\[ F_{a}=F^{i}_{a}\frac{\partial}{\partial x_{i}}, \]
then the evolution equation is
\[\frac{\partial F^{i}_{a}}{\partial t} = g^{ij}{\rm Ric}_{jk}F^{k}_{a}.\]
 Since this is a linear ODE, there are unique
solutions for all times $t\in [0,T)$.

The next remark to make is that this frame remains orthonormal:

\begin{claim}\label{ortho}
 Suppose that ${\mathcal F}(0)=\{F_a\}_a$ is a local  $g(0)$-orthonormal frame,
  and suppose that ${\mathcal F}(t)$ evolves according to Equation~(\ref{frameevol}). Then
 for all $t\in [0,T)$ the frame ${\mathcal F}(t)$ is a local $g(t)$-orthonormal frame.
\end{claim}

\begin{proof}
\begin{eqnarray*}
\frac{\partial}{\partial t}\langle F_a(t),F_b(t)\rangle_{g(t)} & = & \langle
\frac{\partial F_a}{\partial t},F_b\rangle + \langle F_b,\frac{\partial
F_b}{\partial t}\rangle +
\frac {\partial g}{\partial t}(F_a,F_b) \\
 & = & {\rm Ric}(F_a,F_b)+{\rm Ric}(F_b,F_a)-2{\rm Ric}(F_a,F_b)=0.
\end{eqnarray*}
\end{proof}

Notice that if ${\mathcal F}'(0)=\{F'_a\}_a$ is another
frame related to ${\mathcal F}(0)$ by, say,
$$F'_a=A^b_aF_b$$
then $$F_a(t)=A^b_aF_b(t).$$ This means that the evolution of frames
actually defines a bundle automorphism
$$\Phi\colon TM|_U\times [0,T)\to TM|_U\times [0,T)$$
covering the identity map of $U\times[0,T)$ which is independent of
the choice of initial frame and is the identity at time $t=0$. Of
course, since the resulting bundle automorphism is independent of
the initial frame it globalizes to produce a bundle isomorphism
$$\Phi\colon TM\times[0,T)\to TM\times [0,T)$$
covering the identity on $M\times [0,T)$. We view this as an
evolving identification $\Phi_t$ of $TM$ with itself which is the
identity at $t=0$. The content of Claim~\ref{ortho} is:

\begin{cor}
$$\Phi_t^*(g(t))=g(0).$$
\end{cor}

Returning to the local situation of the orthonormal frame ${\mathcal
F}$, we set ${\mathcal F}^*=\{F^1,\ldots, F^n\}$ equal the dual
coframe to $\{F_1,\ldots,F_n\}$. In this coframe the Riemann
curvature tensor is given by $R_{abcd}F^aF^bF^cF^d$ where
\begin{equation}\label{framecoord} R_{abcd} =
R_{ijkl}F^{i}_{a}F^{j}_{b}F^{k}_{c}F^{l}_{d}.\end{equation} One
advantage of working in the evolving frame is that the evolution
equation for the Riemann curvature tensor simplifies:

\begin{lem}\label{M}
Suppose that the orthonormal frame ${\mathcal F}(t)$  evolves by
Formula~(\ref{frameevol}). Then we have the evolution equation
$$\frac{\partial R_{abcd}}{\partial t} = \triangle R_{abcd} +
2(B_{abcd} +B_{acbd}-B_{abdc}-B_{adbc}),$$ where $B_{abcd} =
\sum_{e,f}R_{aebf}R_{cedf}$.
\end{lem}

\begin{proof}
For a proof see Theorem 2.1 in \cite{Hamiltonharnack}.
\end{proof}

Of course, the other way to describe all of this  is to consider the
four-tensor $\Phi_t^*({\mathcal R}_{g(t)})=R_{abcd}F^aF^bF^cF^d$ on
$M$. Since $\Phi_t$ is a bundle map but not a bundle map induced by
a diffeomorphism, even though the pullback of the metric
$\Phi_t^*g(t)$ is constant, it is not the case that the pullback of
the curvature $\Phi_t^*{\mathcal R}_{g(t)}$ is constant.
 The next proposition gives the evolution equation for the pullback
 of the Riemann curvature tensor.

It simplifies the notation somewhat to work directly with a basis of
$\wedge^2TM$. We chose an orthonormal basis
$$\{\varphi^{1},\ldots,\varphi^{\frac{n(n-1)}{2}}\},$$
of $\wedge^2T^*_pM$ where  we have \[ \varphi^{\alpha}(F_{a},F_{b})
= \varphi^{\alpha}_{ab} \] and write the curvature tensor in this
basis as ${\mathcal T} = ({\mathcal T}_{\alpha\beta})$ so that
\begin{equation}\label{Mdefn}
R_{abcd}={\mathcal
T}_{\alpha\beta}\varphi^{\alpha}_{ab}\varphi^{\beta}_{cd}.\end{equation}

\begin{prop}\label{Mevol}
The evolution of the curvature operator ${\mathcal
T}(t)=\Phi_t^*{\rm Rm}(g(t))$ is given by
\[ \frac{\partial {\mathcal T}_{\alpha\beta}}{\partial t} = \triangle
{\mathcal T}_{\alpha\beta} + {\mathcal T}^{2}_{\alpha\beta} +
{\mathcal T}^{\sharp}_{\alpha\beta},
\] where ${\mathcal T}^{2}_{\alpha\beta}={\mathcal T}_{\alpha\gamma}{\mathcal
T}_{\gamma\beta}$ is the operator square; ${\mathcal T}^{\sharp}_{\alpha\beta}
 =
c_{\alpha\gamma\zeta}c_{\beta\delta\eta}{\mathcal T}_{\gamma\delta}{\mathcal
T}_{\zeta\eta}$ is the Lie algebra square; and $c_{\alpha\beta\gamma} =
\langle[\varphi^{\alpha},\varphi^{\beta}],\varphi^{\gamma}\rangle $ are the
structure constants of the Lie algebra ${\rm so}(n)$ relative to the basis
$\{\varphi^{\alpha}\}$. The structure constants $c_{\alpha\beta\gamma}$  are
fully antisymmetric in the three indices.
\end{prop}

\begin{proof} We work in local coordinates that are orthonormal at
the point. By the first Bianchi identity $$R_{abcd}+ R_{acdb}+
R_{adbc} = 0,$$ we get
\begin{align*}
\sum_{e,f}R_{abef}R_{cdef} &=
\sum_{e,f}(-R_{aefb}-R_{afbe})(-R_{cefd}-R_{cfde})\\& =
\sum_{e,f}2R_{aebf}R_{cedf}-2R_{aebf}R_{cfde}\\& = 2(B_{abcd} -
B_{adbc}).
\end{align*}
Note that
\begin{align*}
\sum_{e,f}R_{abef}R_{cdef}& =
\sum_{e,f}{\mathcal T}_{\alpha\beta}\varphi^{\alpha}_{ab}\varphi^{\beta}_{ef}
{\mathcal T}_{\gamma\lambda}\varphi^{\gamma}_{cd}\varphi^{\lambda}_{ef}\\
&={\mathcal T}_{\alpha\beta}\varphi^{\alpha}_{ab}{\mathcal
T}_{\gamma\lambda}\varphi^{\gamma}_{cd}\delta^{\beta\lambda}\\& =
{\mathcal
T}^{2}_{\alpha\beta}\varphi^{\alpha}_{ab}\varphi^{\beta}_{cd}.
\end{align*}
Also,
\begin{align*}
2(B_{acbd}-B_{adbc}) &= 2\sum_{e,f}(R_{aecf}R_{bedf} -
R_{aedf}R_{becf})\\& = 2\sum_{e,f}({\mathcal
T}_{\alpha\beta}\varphi^{\alpha}_{ae}\varphi^{\beta}_{cf}{\mathcal
T}_{\gamma\lambda}\varphi^{\gamma}_{be}\varphi^{\lambda}_{df} -
{\mathcal
T}_{\alpha\beta}\varphi^{\alpha}_{ae}\varphi^{\beta}_{df}{\mathcal
T}_{\gamma\lambda}\varphi^{\gamma}_{be}\varphi^{\lambda}_{cf})
\\&=
2\sum_{e,f}{\mathcal T}_{\alpha\beta}{\mathcal
T}_{\gamma\lambda}\varphi^{\alpha}_{ae}\varphi^{\gamma}_{be}(\varphi^{\beta}_{cf}\varphi^{\lambda}_{df}-\varphi^{\beta}_{df}\varphi^{\lambda}_{cf})
\\&=2\sum_e{\mathcal T}_{\alpha\beta}{\mathcal T}_{\gamma\lambda}\varphi^{\alpha}_{ae}
\varphi^{\gamma}_{be}[\varphi^{\beta},\varphi^{\lambda}]_{cd}
\\&=
\sum_e{\mathcal T}_{\alpha\beta}{\mathcal
T}_{\gamma\lambda}[\varphi^{\beta},\varphi^{\lambda}]_{cd}
(\varphi^{\alpha}_{ae}\varphi^{\gamma}_{be}
- \varphi^{\alpha}_{be}\varphi^{\gamma}_{ae})
\\&={\mathcal T}_{\alpha\beta}{\mathcal T}_{\gamma\delta}[\varphi^{\beta},\varphi^{\lambda}]_{cd}[\varphi^{\alpha},\varphi^{\gamma}]_{ab}
\\&=
{\mathcal
T}^{\sharp}_{\alpha\beta}\varphi^{\alpha}_{ab}\varphi^{\beta}_{cd}.
\end{align*}
So we can rewrite the equation for the evolution of the curvature
tensor given in Lemma~\ref{M} as \[ \frac{\partial R_{abcd}
}{\partial t} = \triangle R_{abcd} + {\mathcal
T}^{2}_{\alpha\beta}\varphi^{\alpha}_{ab}\varphi^{\beta}_{cd} +
{\mathcal
T}^{\sharp}_{\alpha\beta}\varphi^{\alpha}_{ab}\varphi^{\beta}_{cd},
\] or equivalently as \[ \frac{\partial {\mathcal T}_{\alpha\beta}}{\partial t}
= \triangle {\mathcal T}_{\alpha\beta} + {\mathcal
T}^{2}_{\alpha\beta} + {\mathcal T}^{\sharp}_{\alpha\beta}. \] We
abbreviate the last equation as
\[ \frac{\partial {\mathcal T}}{\partial t}= \triangle {\mathcal T} + {\mathcal T}^{2} +
{\mathcal T}^{\sharp}.\]
\end{proof}

\begin{rem}
Notice that neither ${\mathcal T}^{2}$ nor ${\mathcal T}^{\sharp}$
satisfies the Bianchi identity, but their sum does.
\end{rem}

\section{Variation of distance under Ricci flow}

There is one result that we will use several times in the arguments
to follow. Since it is an elementary result (though the proof is
somewhat involved), we have chosen to include it here.

\begin{prop}\label{I.8.3}
 Let $t_0\in \Ar$ and let $(M,g(t))$ be a Ricci flow defined for $t$ in an interval containing $t_0$
 with $(M,g(t))$ complete for every $t$ in this interval. Fix a constant $K<\infty$.  Let $x_0,x_1$ be two points
of $M$ and let $r_0>0$ such that $d_{t_0}(x_0,x_1)\ge 2r_0$. Suppose
that ${\rm Ric}(x,t_0)\le (n-1)K$ for all $x\in B(x_0,r_0,t_0)\cup
B(x_1,r_0,t_0)$. Then
$$\frac{d(d_t(x_0,x_1))}{dt}\Bigl|_{t=t_0}\Bigr.\ge -2(n-1)\left(\frac{2}{3}Kr_0+r_0^{-1}\right).$$
 If the distance function $d_{t}(x_0,x_1)$ is not a differentiable
function of $t$ at $t=t_0$ then this inequality is understood as an
inequality for the forward difference quotient.
\end{prop}

\begin{rem}
Of course, if the distance function is differentiable at $t=t_0$
then the derivative statement is equivalent to the forward
difference quotient statement. Thus, in the proof of this result we
shall always work with the forward difference quotients.
\end{rem}

\begin{proof}
The first step in the proof is to replace the distance function by
the length of minimal geodesics. The following is standard.

\begin{claim}\label{geosuff}
Suppose that for every minimal $g(t_0)$-geodesic $\gamma$ from $x_0$
to $x_1$ the function $\ell_t(\gamma)$ which is the $g(t)$-length of
$\gamma$ satisfies
$$\frac{d(\ell_t(\gamma))}{dt}\Bigl|_{t=t_0}\Bigr.\ge C.$$
Then
$$\frac{d(d_t(x_0,x_1))}{dt}\Bigl|_{t=t_0}\Bigr.\ge C,$$
where, as in the statement of the proposition, if the distance
function is not differentiable at $t_0$ then the inequality in the
conclusion is interpreted by replacing the derivative on the
left-hand side with the ${\rm liminf}$ of the forward difference
quotients of $d_t(x_0,x_1)$ at  $t_0$.
\end{claim}

The second step in the proof is to estimate the time derivative of a
minimal geodesic under the hypothesis of the proposition.

\begin{claim}\label{3.23} Assuming the hypothesis of the proposition, for any
 minimal $g(t_0)$-geodesic $\gamma$ from $x_0$ to $x_1$, we have
$$\frac{d(\ell_t(\gamma))}{dt}\Bigl|_{t=t_0}\Bigr.\ge -2(n-1)
\left(\frac{2}{3}Kr_0+r_0^{-1}\right).$$
\end{claim}

\begin{proof}
Fix a minimal $g(t_0)$-geodesic $\gamma(u)$ from $x_0$ to $x_1$,
parameterized by arc length. We set $d=d_{t_0}(x_0,x_1)$, we set
$X(u)=\gamma'(u)$, and we take tangent vectors $Y_1,\ldots,Y_{n-1}$
in $T_{x_0}M$ which together with $X(0)=\gamma'(0)$ form an
orthonormal basis. We let $Y_i(u)$ be the parallel translation of
$Y_i$ along $\gamma$. Define $f\colon [0,d]\to [0,1]$ by:
$$f(u)=\begin{cases}u/r_0  & 0\le u\le r_0 \\ 1 & r_0\le u\le d-r_0 \\
(d-u)/r_0 & d-u\le r_0\le d, \end{cases}$$ and define $$\widetilde
Y_i(u)=f(u)Y_i(u).$$ See \textsc{Fig.}~\ref{fig:jacobi}.
 For $1\le i\le n-1$, let $s''_{\widetilde
Y_i}(\gamma)$ be the second variation of the $g(t_0)$-length of
$\gamma$ along $\widetilde Y_i$. Since $\gamma$ is a minimal
$g(t_0)$-geodesic, for all $i$ we have
\begin{equation}\label{geoineq}
s''_{\widetilde Y_i}(\gamma)\ge 0.\end{equation}

\begin{figure}[ht]
  \relabelbox
  \centerline{\epsfbox{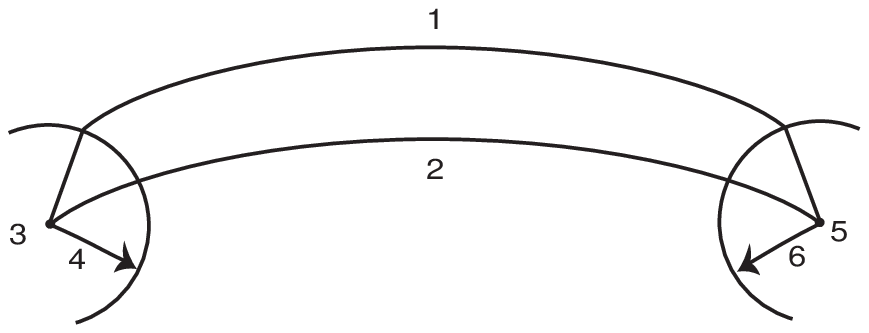}}
  \relabel{1}{$\widetilde{Y}_i$}
  \relabel{2}{$\gamma$}
  \relabel{3}{$x_0$}
  \relabel{4}{$r_0$}
  \relabel{5}{$x_1$}
  \relabel{6}{$r_0$}
  \endrelabelbox
  \caption{$\widetilde{Y}_i$ along $\gamma$.}\label{fig:jacobi}
\end{figure}

Let us now compute $s''_{\widetilde Y_i}(\gamma)$ by taking a
two-parameter family $\gamma(u,s)$ such that the curve $\gamma(u,0)$
is the original minimal geodesic and $\frac{\partial}{\partial
s}(\gamma(u,s))|_{s=0}=\widetilde Y_i(u)$.  We denote by $X(u,s)$
the image $D\gamma_{(u,s)}(\partial/\partial u)$ and by $\widetilde
Y_i(u,s)$ the image $D\gamma_{(u,s)}(\partial/\partial s)$. We wish
to compute
\begin{eqnarray}
\lefteqn{s''_{\widetilde Y_i}(\gamma) =
\frac{d^2}{ds^2}\left(\int_0^d\sqrt{X(u,s),X(u,s)}du\right)\Bigl|_{s=0}\Bigr.}
& &
 \nonumber \\
& = & \frac{d}{ds}\left(\int_0^d\langle X(u,s),X(u,s)\rangle^{-1/2}
\langle
X(u,s),\nabla_{\widetilde Y_i}X(u,s)\rangle du\right)\Bigl|_{s=0}\Bigr. \nonumber \\
& = &  \int_0^d-\langle X(u,0),X(u,0)\rangle^{-3/2}\langle
X(u,0),\nabla_{\widetilde
Y_i}X(u,0)\rangle^2du \label{eq2ndvar} \\
& &  +\int_0^d\frac{\langle \nabla_{\widetilde
Y_i}X(u,0),\nabla_{\widetilde Y_i}X(u,0)\rangle +\langle
X(u,0),\nabla_{\widetilde Y_i}\nabla_{\widetilde
Y_i}X(u,0)\rangle}{\langle X(u,0),X(u,0)\rangle^{1/2}}du \nonumber.
\end{eqnarray}
Using the fact that $X$ and $\widetilde Y_i$ commute (since they are
the coordinate partial derivatives of a map of a surface into $M$)
and using the fact that $Y_i(u)$ is parallel along $\gamma$, meaning
that $\nabla_X(Y_i)(u)=0$, we see that $\nabla_{\widetilde
Y_i}X(u,0)=\nabla_X\widetilde Y_i(u,0)=f'(u)Y_i(u)$. By construction
$\langle Y_i(u),X(u,0)\rangle=0$. It follows that
$$\langle \nabla_{\widetilde Y_i}X(u,0),X(u,0)\rangle=\langle \nabla_X(\widetilde
Y_i)(u,0),X(u,0)\rangle = \langle f'(u)Y_i(u),X(u,0)\rangle =0.$$
Also, $\langle X(u,0),X(u,0)\rangle=1$, and by construction $\langle
Y_i(u,0),Y_i(u,0)\rangle=1$. Thus, Equation~(\ref{eq2ndvar})
simplifies to
\begin{eqnarray}
\lefteqn{s''_{\widetilde
Y_i}(\gamma)=\frac{d^2}{ds^2}\left(\int_0^d\sqrt{X(u,s),X(u,s)}du\right)\Bigl|_{s=0}\Bigr.
 \label{2ndvaria}}
\\& = & \int_0^d\left((f'(u))^2\langle
Y_i(u),Y_i(u)\rangle+\langle \nabla_{\widetilde
Y_i}\nabla_X(\widetilde
Y_i(u,0)),X(u,0)\rangle\right)du \nonumber\\
& = & \int_0^d\left(\langle R(\widetilde Y_i,X)\widetilde
Y_i(u,0),X(u,0)\rangle-\langle \nabla_X\nabla_{\widetilde
Y_i}\widetilde Y_i(u,0),X(u,0)\rangle+(f'(u))^2\right)
du.\nonumber\end{eqnarray} Now we restrict to $s=0$ and for simplicity of notation
we leave the variable $u$ implicit. We have
$$\langle
\nabla_X\nabla_{\widetilde Y_i}\widetilde
Y_i,X\rangle=\frac{d}{du}\langle \nabla_{\widetilde Y_i}\widetilde
Y_i,X\rangle-\langle \nabla_{\widetilde Y_i}\widetilde
Y_i,\nabla_XX\rangle=\frac{d}{du}\langle \nabla_{\widetilde
Y_i}\widetilde Y_i,X\rangle,$$ where the last equality is a
consequence of the geodesic equation, $\nabla_XX=0$. It follows that
$$\int_0^d\langle
\nabla_X\nabla_{\widetilde Y_i}\widetilde Y_i,X\rangle
du=\int_0^d\frac{d}{du}\langle \nabla_{\widetilde Y_i}\widetilde
Y_i,X\rangle=0,$$ where the last equality is a consequence of the
fact that $\widetilde Y_i$ vanishes at the end points.

Consequently, plugging these into Equation~(\ref{2ndvaria}) we have
\begin{equation} \label{sformula} s''_{\widetilde Y_i}(\gamma) =
  \int_0^d\left(\langle R(\widetilde Y_i,X)\widetilde
Y_i(u,0),X(u,0)\rangle+(f'(u))^2\right)du. \end{equation} Of course,
it is immediate from the definition that $f'(u)^2=1/r_0^2$ for $u\in
[0,r_0]$ and for $u\in [d-r_0,d]$ and is zero otherwise. Also, from
the definition of the vector fields $Y_i$ we have
$$\sum_{i=1}^{n-1}\langle
R(Y_i,X)Y_i(u),X(u)\rangle=-{\rm Ric}_{g(t_0)}(X(u),X(u)),$$ so that
$$\sum_{i=1}^{n-1}\langle
R(\widetilde Y_i,X)\widetilde
Y_i(u),X(u)\rangle=-f^2(u){\rm Ric}_{g(t_0)}(X(u),X(u)).$$

Hence, summing Equalities~(\ref{sformula}) for $i=1,\ldots, n-1$ and
using Equation~(\ref{geoineq}) gives
\begin{eqnarray*}
0\le\sum_{i=1}^{n-1}s''_{\widetilde Y_i}(\gamma) & = &
\int_0^{r_0}\left[\frac{u^2}{r_0^2}\left(-{\rm Ric}_{g(t_0)}\left(X(u),X(u)\right)\right)
+\frac{n-1}{r_0^2}\right]du
\\
& & + \int_{r_0}^{d-r_0}-{\rm Ric}_{g(t_0)}(X(u),X(u))du
\\
 & & +\int_{d-r_0}^d\left[\frac{(d-u)^2}{r_0^2}
\left(-{\rm Ric}_{g(t_0)}(X(u),X(u))\right)+\frac{n-1}{r_0^2}\right]du.\end{eqnarray*}
Rearranging the terms yields
\begin{eqnarray*} 0 & \le & -\int_0^d {\rm Ric}_{g(t_0)}(X(u),X(u))du \\
& &+\int_0^{r_0}\left[\left(1-\frac{u^2}{r_0^2}\right)
\left({\rm Ric}_{g(t_0)}(X(u),X(u))\right) +\frac{n-1}{r_0^2}\right]du\\
& & +\int_{d-r_0}^d
\left[\left(1-\frac{(d-u)^2}{r_0^2}\right)\left({\rm Ric}_{g(t_0)}(X(u),X(u))\right)
+\frac{n-1}{r_0^2}\right]du.
\end{eqnarray*}
Since
\begin{eqnarray*}
\frac{d(\ell_t(\gamma))}{dt}\Bigl|_{t=t_0}\Bigr. & = &
\frac{d}{dt}\left[\left(\int_0^d\sqrt{\langle X(u),X(u)\rangle}dt\right)^{1/2}\right]|_{t=t_0} \\
& = & -\int_0^d {\rm Ric}_{g(t_0)}(X(u),X(u))du,
\end{eqnarray*}
we have
\begin{eqnarray*}
\lefteqn{\frac{d(\ell_t(\gamma))}{dt}\Bigl|_{t=t_0}\Bigr.\ge
-\left\{\int_0^{r_0}\left[\left(1-\frac{u^2}{r_0^2}\right)\left({\rm Ric}_{g(t_0)}(X(u),X(u))\right)
+\frac{n-1}{r_0^2}\right]du\right.}
\\
& & \left. +\int_{d-r_0}^d
\left[\left(1-\frac{(d-u)^2}{r_0^2}\right)\left({\rm Ric}_{g(t_0)}(X(u),X(u))\right)
+\frac{n-1}{r_0^2}\right]du\right\}
\end{eqnarray*}

Now, since $|X(u)|=1$, by the hypothesis of the proposition we have
the estimate ${\rm Ric}_{g(t_0)}(X(u),X(u))\le (n-1)K$ on the regions of
integration on the right-hand side of the above inequality. Thus,
\begin{eqnarray*}
\frac{d(\ell_t(\gamma))}{dt}\Bigl|_{t=t_0}\Bigr.\ge
-2(n-1)\left(\frac{2}{3}r_0K+r_0^{-1}\right).
\end{eqnarray*}
This completes the proof of Claim~\ref{3.23}.
\end{proof}

Claims~\ref{geosuff} and~\ref{3.23} together prove the proposition.
\end{proof}

\begin{cor}\label{2ndI.8.3}
Let $t_0\in \Ar$ and let $(M,g(t))$ be a Ricci flow defined for $t$ in an
interval
containing $t_0$ and with $(M,g(t))$ complete for every $t$ in this
interval. Fix a constant $K<\infty$. Suppose that ${\rm Ric}(x,t_0)\le (n-1)K$
for all $x\in M$.
 Then for any points $x_0,x_1\in M$ we have
$$\frac{d(d_t(x_0,x_1))}{dt}\Bigl|_{t=t_0}\Bigr.\ge -4(n-1)\sqrt{\frac{2K}{3}}$$
in the sense of forward difference quotients.
\end{cor}

\begin{proof}
There are two cases: Case (i): $d_{t_0}(x_0,x_1)\ge
\sqrt{\frac{6}{K}}$ and Case (ii) $d_{t_0}(x_0,x_1)<
\sqrt{\frac{6}{K}}$. In Case (i) we  take $r_0=\sqrt{3/2K}$ in
Proposition~\ref{I.8.3}, and we conclude that the ${\rm liminf}$ at
$t_0$ of the difference quotients for $d_{t}(x_0,x_1)$ is at most
$-4(n-1)\sqrt{\frac{2K}{3}}$. In Case (ii) w let $\gamma(u)$ be any
minimal $g(t_0)$-geodesic from $x_0$ to $x_1$ parameterized by arc
length. Since
$$\frac{d}{dt}(\ell_t(\gamma))|_{t=t_0}=-\int_\gamma
{\rm Ric}_{g(t_0)}(\gamma'(u),\gamma'(u))du,$$ we see that
$$\frac{d}{dt}(\ell_t(\gamma))|_{t=t_0}\ge
-(n-1)K\sqrt{6/K}=-(n-1)\sqrt{6K}.$$ By Claim~\ref{geosuff}, this
implies that the ${\rm liminf}$ of the forward difference quotient of
$d_t(x_0,x_1)$ at $t=t_0$ is at least $-(n-1)\sqrt{6K}\ge
-4(n-1)\sqrt{2K/3}$.
\end{proof}

\begin{cor}\label{corI.8.3}
 Let $(M,g(t)),\ a\le t\le b$, be a Ricci flow with $(M,g(t))$ complete for every $t
\in [0,T)$. Fix a positive function $K(t)$, and suppose that ${\rm
Ric}_{g(t)}(x,t)\le (n-1)K(t)$ for all $x\in M$ and all $t\in [a,b]$. Let
$x_0,x_1$ be two points of $M$. Then
$$d_a(x_0,x_1)\le d_b(x_0,x_1)+4(n-1)\int_a^b\sqrt{\frac{2K(t)}{3}}dt.$$
\end{cor}

\begin{proof} By Corollary~\ref{2ndI.8.3} we have
\begin{equation}\label{liminf}
\frac{d}{dt}d_t(x_0,x_1)|_{t=t'}\ge -4(n-1)\sqrt{\frac{2K(t')}{3}}
\end{equation}
in the sense  of forward difference quotients. Thus, this result is
an immediate consequence of Lemma~\ref{fordiffquot}.
\end{proof}

\section{Shi's derivative estimates}

The last `elementary' result we discuss is Shi's result\index{Shi's
derivative estimates|(ii} controlling all derivatives in terms of a
bound on curvature. This is a consequence of the parabolic nature of
the Ricci flow equation. More precisely, we can control all
derivatives of the curvature tensor at a point $p\in M$ and at a
time $t$ provided that we have an upper bound for the curvature on
an entire backward parabolic neighborhood of $(p,t)$ in space-time.
The estimates become weaker as the parabolic neighborhood shrinks,
either in the space direction or the time direction.

Recall that for any $K<\infty$  if $(M,g)$ is a Riemannian manifold
with $|{\rm Rm}|\le K$ and if for some $r\le \pi/\sqrt{K}$ the metric ball
$B(p,r)$ has compact closure in $M$, then the exponential mapping
${\rm exp}_p$ is defined on the ball $B(0,r)$ of radius $r$ centered
at the origin of $T_pM$ and ${\rm exp}_p\colon B(0,r)\to M$ is a
local diffeomorphism onto $B(p,r)$.

The first of Shi's derivative estimates controls the first
derivative of ${\rm Rm}$.

\begin{thm}\label{firstshi}
There is a constant $C=C(n)$, depending only on the dimension $n$,
such that the following holds for every $K<\infty$, for every $T>0$,
and for every $r>0$. Suppose that $(U,g(t)),\ 0\le t\le T$, is an $n$-dimensional
Ricci flow with $|{\rm Rm}(x,t)|\le K$ for all $x\in U$ and $t\in [0,T]$.
Suppose that $p\in U$ has the property that $B(p,0,r)$ has compact
closure in $U$. Then
$$|\nabla {\rm Rm}(p,t)|\le
CK\left(\frac{1}{r^2}+\frac{1}{t}+K\right)^{1/2}.$$
\end{thm}

For a proof of this result see Chapter 6.2, starting on page 212, of
\cite{ChowLuNi}.

We also need higher derivative estimates. These are also due to Shi,
but they take a slightly different form. (See Theorem 6.9 on page
210 of \cite{ChowLuNi}.)

\begin{thm}{\bf(Shi's Derivative Estimates)}\label{shi}
Fix the dimension $n$ of the Ricci flows under consideration. Let
$K<\infty$ and $\alpha>0$ be positive constants.  Then for each
non-negative integer $k$ and each $r>0$ there is a constant
$C_k=C_k(K,\alpha,r,n)$ such that the following holds. Let
$(U,g(t)),\ 0\le t\le T$, be a Ricci flow with $T\le \alpha/K$.
Fix $p\in U$ and suppose that the metric ball $B(p,0,r)$ has compact
closure in $U$. If
$$|{\rm Rm}(x,t)|\le K\ \ \text{for all }\ (x,t)\in P(x,0,r,T),$$
then
$$|\nabla^k({\rm Rm}(y,t))|\le \frac{C_k}{t^{k/2}}$$
for all $y\in B(p,0,r/2)$ and all $t\in(0,T]$.
\end{thm}

For a proof of this result see Chapter 6.2 of \cite{ChowLuNi} where
these estimates are proved for the first and second derivatives of
${\rm Rm}$. The proofs of the higher derivatives follow similarly.  Below, we
shall prove a stronger form of this result below including the proof
for all derivatives.

We shall need a stronger version of this result, a version which is
well-known but for which there seems to be no good reference. The
stronger version takes as hypothesis $C^k$-bounds on the
initial conditions and produces a better bound on the derivatives of
the curvature at later times. The argument is basically the same as
that of the result cited above, but since there is no good reference
for it we include the proof, which was shown to us by Lu Peng.

\begin{thm}\label{shiw/deriv}
Fix the dimension $n$ of the Ricci flows under consideration. Let
$K<\infty$ and $\alpha>0$ be given positive constants. Fix an
integer $l\ge 0$. Then for each integer $k\ge 0$ and for each $r>0$
there is a constant $C'_{k,l}=C'_{k,l}(K,\alpha,r,n)$ such that the
following holds. Let $(U,g(t)),\ 0\le t\le T$, be a Ricci flow with
$T\le \alpha/K$. Fix $p\in U$ and suppose that the metric ball
$B(p,0,r)$ has compact closure in $U$. Suppose that
\begin{eqnarray*}
\left\vert \operatorname*{Rm}\left(  x,t\right)  \right\vert  &
\leq K\ \ \text{ for all }x\in U\text{ and all}\  t\in[0,T],\\
\left\vert \nabla^{\beta}\operatorname*{Rm}\left(  x,0\right)
\right\vert  & \leq K\ \ \text{ for all }x\in U\ \ \text{ and all }\
\ \beta\leq l.
\end{eqnarray*}
Then
\[
\left\vert \nabla^{k}\operatorname*{Rm}\left(  y,t\right)
\right\vert
\leq\frac{C'_{k,l}}{t^{\max\left\{  k-l,0\right\}  /2}}%
\]
for all $y\in$ $B(p,0,r/2)  $ and all $t\in(0,T].$ In particular if
$k\leq l$, then for $y\in B(p,0,r/2)$ and $t\in (0,T]$ we have
\[
\left\vert \nabla^{k}\operatorname*{Rm}\left(  y,t\right)
\right\vert \leq C_{k,l}'.
\]
\end{thm}

\begin{rem}
Clearly, the case $l=0$ of Theorem~\ref{shiw/deriv} is Shi's theorem
(Theorem~\ref{firstshi}).
\end{rem}

Theorem~\ref{shiw/deriv} leads immediately to the following:

\begin{cor}\label{Cinftybd} Suppose that $(M,g(t)),\ 0\le t\le T$, is a Ricci flow
with $(M,g(t))$ being complete and with $T<\infty$. Suppose that
$\operatorname*{Rm}(x,0)$ is bounded in the $C^\infty$-topology
independent of $x\in M$ and suppose that
$\left\vert\operatorname*{Rm}(x,t)\right\vert$ is bounded
independent of $x\in M$ and $t\in[0,T]$. Then  the operator
$\operatorname*{Rm}(x,t)$  is bounded in the $C^\infty$-topology
independent of $(x,t)\in M\times [0,T]$.
\end{cor}

 For a proof of Theorem~\ref{shi}  see \cite{Shi1, Shi2}. We give the proof of
a stronger result, Theorem~\ref{shiw/deriv}.\index{Shi's derivative
estimates|)}

\begin{proof}
The first remark is that establishing Theorem~\ref{shiw/deriv} for
one value of $r$ immediately gives it for all $r'\ge 2r$. The reason
is that for such $r'$ any point  $y\in B(p,0,r'/2)$ has the property
that $B(y,0,r)\subset B(p,0,r')$ so that a curvature bound on
$B(p,0,r')$ will imply one on $B(y,0,r)$ and hence by the result for
$r$ will imply the higher derivative bounds at $y$.

 Thus, without loss of
generality we can suppose that $r\le \pi/2\sqrt{K}$. We shall assume this
from now on in the proof. Since
$B(p,0,r)$ has compact closure in $M$, for some $r<r'<\pi/\sqrt{K}$
the ball $B(p,0,r')$ also has compact closure in $M$. This means
that the exponential mapping from the ball of radius $r'$ in $T_pM$
is a local diffeomorphism onto $B(p,0,r')$.

The proof is by induction: We assume that we have established the
result for $k=0,\ldots,m$, and then we shall establish it for
$k=m+1$. The inductive hypothesis tells us that there are constants
$A_j,\ 0\le j\le m$, depending on $(l,K,\alpha,r,n)$ such that for
all $(x,t)\in B(p,0,r/2)\times (0,T]$ we have
\begin{equation}\label{strongineq}
\left\vert \nabla^{j}\operatorname*{Rm} \left(  x,t\right)
\right\vert \leq A_{j}t^{-\max\left\{  j-l,0\right\}  /2}.
\end{equation}
 Even better, applying the inductive result to $B(y,0,r/2)$ with
$y\in B(p,0,r/2)$ we see, after replacing the $A_j$ by the
larger constants associated with $(l,K,\alpha,r/2,n)$, that we have
the same inequality for all $y\in B(x,0,3r/4)$.

We fix a constant $C\geq \max(4A^2_m,1)$ and consider
\[
F_{m}(x,t)=\left(  C+t^{\max\left\{ m-l,0\right\} }  \left\vert \nabla
^{m}\operatorname*{Rm}(x,t)\right\vert ^{2}\right)  t^{\max\left\{
m+1-l,0\right\}  }\left\vert
\nabla^{m+1}\operatorname*{Rm}(x,t)\right\vert ^{2}.
\]
Notice that bounding $F_m$ above by a constant $(C'_{m+1,l})^2$ will
yield
$$|\nabla^{m+1}{\rm Rm}(x,t)|^2\le \frac{(C'_{m+1,l})^2}{t^{{\rm max}\{m+1-l,0\}}},$$
and hence will complete the proof of the result.

 Bounding $F_m$
above (assuming the inductive hypothesis) is what is accomplished in
the rest of this proof. The main calculation is the proof of the
following claim under the inductive hypothesis.

\begin{claim}\label{Fm}
 With $F_m$ as defined above and with $C\ge \max(4A^2_m,1)$,
there are constants $c_1$ and $C_0,C_1$ depending on $C$ as well as
$K,\alpha,A_1,\ldots,A_m$ for which the following holds on
$B(p,0,3r/4)\times (0,T]$:
$$
\left(  \frac{\partial}{\partial t}-\Delta\right)
F_{m}(x,t)\leq-\frac {c_1}{t^{{\rm s}\left\{  \max\left\{
m-l+1,0\right\}  \right\} }}\left( F_{m}(x,t)-C_0\right)
^{2}+\frac{C_1}{t^{{\rm s}\{\max\left\{ m-l+1,0\right\} \} }},
$$
where
$${\rm s}(n)=\begin{cases} +1 & \text{if } n>0 \\ 0 & \text{if }n=0 \\ -1&
\text{if }n<0.\end{cases}$$
\end{claim}

Let us assume this claim and use it to prove Theorem~\ref{shiw/deriv}. We fix
$C=\max\{4A_m^2,1\}$, and consider the resulting function $F_m$. The
constants $c_1,C_0,C_1$ from Claim~\ref{Fm} depend only on $K,\alpha$, and
$A_1,\ldots,A_m$. Since $r\le \pi/2\sqrt{K}$, and $B(p,0,r)$ has
compact closure in $U$, there is some $r'>r$ so that the exponential
mapping ${\rm exp}_p\colon B(0,r')\to U$ is a local diffeomorphism
onto $B(p,0,r')$. Pulling back by the exponential map, we  replace
the Ricci flow on $U$ by a Ricci flow on $B(0,r')$ in $T_pM$.
Clearly, it suffices to establish the estimates in the statement of
the proposition for $B(0,r/2)$. This remark allows us to assume that
the exponential mapping is a diffeomorphism onto $B(p,0,r)$. Bounded
curvature then comes into play in the following crucial proposition,
which goes back to Shi. The function given in the next proposition allows us to localize
the computation in the ball $B(p,0,r)$.

\begin{prop}\label{etaprop}
Fix constants $0<\alpha$ and the dimension $n$. Then there is a constant
$C_2'=C_2'(\alpha,n)$ and for each $r>0$ and $K<\infty$ there is a constant
$C_2=C_2(K,\alpha,r,n)$ such that the following holds. Suppose that $(U,g(t)),\
0\le t\le T$,  is an $n$-dimensional Ricci flow  with $T\le \alpha/K$. Suppose
that $p\in U$ and that $B(p,0,r)$ has compact closure in $U$ and that the
exponential mapping from the ball of radius $r$ in $T_pU$ to $B(p,0,r)$ is a
diffeomorphism. Suppose that $|{\rm Rm}(x,0)|\le K$ for all $x\in B(p,0,r)$. There is
a smooth function $\eta\colon B(p,0,r)\to [0,1]$ satisfying the following for
all $t\in [0,T]$:
\begin{enumerate}
\item $\eta$ has compact support in $B(p,0,r/2)$
\item The restriction of $\eta$ to $B(p,0,r/4)$ is identically $1$.
\item $|\Delta_{g(t)} \eta|\le C_2(K,\alpha,r,n)$.
\item $\frac{|\nabla \eta|_{g(t)}^2}{\eta}\le \frac{C'_2(\alpha,n)}{r^2}$
\end{enumerate}
\end{prop}

For a proof of this result see Lemma 6.62 on page 225 of
\cite{ChowLuNi}.

We can apply this proposition to our situation, because we are
assuming that $r\le \pi/2\sqrt{K}$ so that the exponential mapping
is a local diffeomorphism onto $B(p,0,r)$ and we have pulled the
Ricci flow back to the ball in the tangent space.

Fix any $y\in B(p,0,r/2)$ and choose $\eta$ as in the previous
proposition for the constants $C_2(\alpha,n)$ and
$C_2'(K,\alpha,r/4,n)$. Notice that $B(y,0,r/4)\subset B(p,0,3r/4)$
so that the conclusion of Claim~\ref{Fm} holds for every $(z,t)$
with $z\in B(y,0,r/4)$ and $t\in [0,T]$.  We shall show that the restriction of $\eta
F_m$ to $P(y,0,r/4,T)$ is bounded by a constant that depends only on
$K,\alpha,r,n,A_1,\ldots,A_m$. It will then follow immediately that
the restriction of $F_m$ to $P(y,0,r/8,T)$ is bounded by the same
constant. In particular, the values of $F_m(y,t)$ are bounded by the
same constant for all $y\in B(p,0,r/2)$ and $t\in [0,T]$.

 Consider a point $(x,t)\in B(y,0,r/2)\times
[0,T]$ where $\eta F_m$ achieves its maximum; such a point exists
since the ball $B(y,0,r/2)\subset B(p,0,r)$, and hence $B(y,0,r/2)$
has compact closure in $U$. If $t=0$, then $\eta F_m$ is bounded by
$(C+K^2)K^2$ which is a constant depending only on $K$ and $A_m$.
This, of course, immediately implies the result. Thus we can assume
that the maximum is achieved at some $t>0$. When $s\left\{
\max\left\{ m+1-l,0\right\} \right\}  -0,$ according to the
Claim~\ref{Fm}, we have

\[
\left(  \frac{\partial}{\partial t}-\Delta\right)
F_{m}\leq-c_1\left( F_{m}-C_0\right)  ^{2}+C_1.
\]

We compute
\[
\left(  \frac{\partial}{\partial t}-\Delta\right)  \left(  \eta
F_{m}\right) \leq\eta\left(  -c_1\left(
F_{m-C_0}\right)^{2}+C_1\right) -\Delta\eta\cdot
F_m-2\nabla\eta\cdot\nabla F_m.
\]

 Since $(x,t)$ is a maximum point for $\eta F_m$ and since
$t>0$,  a simple maximum principle argument shows that
$$\left(  \frac{\partial}{\partial t}-\Delta\right)\eta  F_{m}(x,t)\ge 0.$$
Hence, in this case we conclude that
\begin{eqnarray*} 0& \le & \left(
\frac{\partial}{\partial t}-\Delta\right)  \left( \eta(x)
F_{m}(x,t)\right) \leq\eta(x)\left(  -c_1\left(
F_{m}(x,t)-C_0\right)^{2}+C_1\right) \\
& &  -\Delta\eta(x)\cdot F_m(x,t)-2\nabla\eta(x)\cdot\nabla
F_m(x,t). \end{eqnarray*}
 Hence, $$c_1\eta (F_m(x,t)-C_0)^2\le \eta(x)
C_1-\Delta\eta(x)\cdot F_m(x,t)-2\nabla\eta(x)\cdot \nabla F_m(x,t).$$ Since we
are proving that $F_m$ is bounded, we are free to argue by contradiction and
assume that $F_m(x,t)\ge 2C_0$, implying that $F_m(x,t)-C_0\ge F_m(x,t)/2$.
Using this inequality yields \begin{eqnarray*} \eta(x) (F_m(x,t)-C_0) & \le &
\frac{2\eta C_1}{c_1F_m(x,t)}-\frac{2\Delta\eta(x)}{c_1} -\frac{4}{c_1F_m(x,t)}
\nabla\eta(x)\cdot\nabla F_m(x,t) \\
&\le &  \frac{\eta C_1}{c_1C_0}-\frac{2\Delta\eta(x)}{c_1}
-\frac{4}{c_1F_m(x,t)} \nabla\eta(x)\cdot\nabla F_m(x,t)\end{eqnarray*}
Since
$(x,t)$ is a maximum for $\eta F_m$ we have
$$0=\nabla(\eta(x) F_m(x,t))=\nabla\eta(x) F_m(x,t)+\eta(x)\nabla F_m(x,t),$$
so that
$$\frac{\nabla\eta(x)}{\eta(x)}=-\frac{\nabla F_m(x,t)}{F_m(x,t)}.$$
Plugging this in gives
$$\eta(x) F_m(x,t)\le
\frac{C_1}{c_1C_0}-\frac{2\Delta\eta(x)}{c_1}+4\frac{|\nabla\eta(x)|^2}{c_1\eta(x)}+\eta
C_0.$$

Of course, the gradient and Laplacian of $\eta$ are taken at the
point $(x,t)$.  Thus, because of the properties of $\eta$ given in
Proposition~\ref{etaprop}, it immediately follows that $\eta F_m(x,t)$ is
bounded by a constant depending
 only on $K,n,\alpha,r,c_1,C_0,C_1$, and as we have already seen, $c_1,C_0,C_1$ depend only on
 $K,\alpha,A_1,\ldots,A_m$.

Now suppose that $s\left\{  \max\left\{  m-l+1,0\right\}
\right\}=1.$ Again we compute the evolution inequality for $\eta
F_{m}$. The result is $$\left(  \frac{\partial}{\partial
t}-\Delta\right) \left(  \eta F_{m}\right) \leq\eta\left(
-\frac{c_1}{t} (F_{m}-C_0)^{2}+\frac{C_1}{t}\right) -\Delta\eta\cdot
F_m-2\nabla\eta\cdot\nabla F_m.$$ Thus,  using the maximum principle
as before, we have
$$\left(  \frac{\partial}{\partial t}-\Delta\right)\eta  F_{m}(x,t)\ge 0.$$
Hence,
$$\frac{\eta(x) c_1(F_m(x,t)-C_0)^2}{t}\le
\frac{\eta(x)C_1}{t}-\Delta\eta(x)F_m(x,t)-2\nabla\eta(x)\cdot\nabla F_m(x,t).$$
Using the assumption that $F_m(x,t)\ge 2C_0$ as before, and rewriting the last
term as before, we have
$$\eta F_m(x,t)\le
\frac{\eta(x)C_1}{c_1C_0}-\frac{2t\Delta\eta(x)}
{c_1}+\frac{4t|\nabla\eta(x)|^2}{c_1\eta(x)}+\eta C_0.$$ The right-hand side is
bounded by a constant depending only on $K,n,\alpha,r,c_1,C_0,C_1$. We conclude
that in all cases $\eta F_{m}$ is bounded by a constant depending only on
$K,n,\alpha,r,c_1,C_0,C_1$, and hence on $K,n,\alpha,r,A_1,\ldots,A_m$.

This proves that for any $y\in B(p,0,r/2)$, the value
 $\eta F_m(x,t)$ is bounded by a
constant $A_{m+1}$ depending only on $(m+1,l,K,n,\alpha,r)$ for all
$(x,t)\in B(y,0,r/2)\times [0,T]$. Since $\eta(y)=1$, for all $0\le t\le T$ we have
$$t^{{\rm max}\{m+1-l,0\}}|\nabla^{m+1}{\rm Rm}(y,t)|^2\le F_m(y,t)=\eta(y)F_m(y,t)\le A_{m+1}.$$

This completes the inductive proof that the result holds for
$k=m+1$ and hence establishes Theorem~\ref{shiw/deriv}, modulo the proof of Claim~\ref{Fm}.
\end{proof}

Now we turn to the proof of Claim~\ref{Fm}.

\begin{proof} In this argument we fix
$(x,t)\in B(p,0,3r/4)\times (0,T]$ and we drop $(x,t)$ from the
notation. Recall that by Equations (7.4a) and (7.4b) on p. 229 of
\cite{ChowKnopf}  we have
\begin{align}
\frac{\partial}{\partial t}\left\vert \nabla^{\ell
}\operatorname*{Rm}\right\vert ^{2} & \leq\Delta\left\vert
\nabla^{\ell }\operatorname*{Rm}\right\vert ^{2} -2\left\vert
\nabla^{\ell +1}\operatorname*{Rm}\right\vert ^{2}+\sum_{i=0} ^{\ell
}c_{\ell,j}\left\vert \nabla^{i}\operatorname*{Rm}\right\vert
\left\vert \nabla^{\ell -i}\operatorname*{Rm}\right\vert \left\vert
\nabla^{\ell } \operatorname*{Rm}\right\vert, \label{curderfor}
\end{align}
where the constants $c_{\ell,j}$ depend only on $\ell$ and $j$.

Hence, setting $m_l= \max\left\{ m+1-l,0\right\}$ and denoting
$c_{m+1,i}$ by $\tilde c_i$, we have

\begin{eqnarray}
\lefteqn{\frac{\partial}{\partial t}\left(  t^{m_l}\left\vert
\nabla^{m+1}\operatorname*{Rm}\right\vert ^{2}\right)
 \leq\Delta\left(  t^{m_l }\left\vert \nabla
^{m+1}\operatorname*{Rm}\right\vert ^{2}\right)  -2t^{m_l }
\left\vert \nabla^{m+2}\operatorname*{Rm}\right\vert ^{2}}\label{1steqn}\\
& & +t^{m_l  }\sum_{i=0}^{m+1}\tilde c_{i}\left\vert \nabla
^{i}\operatorname*{Rm}\right\vert \left\vert
\nabla^{m+1-i}\operatorname*{Rm} \right\vert \left\vert
\nabla^{m+1}\operatorname*{Rm}\right\vert
 +m_l  t^{m_l -1}\left\vert \nabla^{m+1}\operatorname*{Rm}\right\vert ^{2}\nonumber\\
&  \leq & \Delta\left(  t^{m_l }\left\vert \nabla
^{m+1}\operatorname*{Rm}\right\vert ^{2}\right)  -2t^{m_l }
\left\vert \nabla^{m+2}\operatorname*{Rm}\right\vert ^{2} +(\tilde
c_{0}+\tilde
c_{m+1})t^{m_l}\left\vert\operatorname*{Rm}\right\vert\left\vert\nabla^{m+1}\operatorname*{Rm}
\right\vert^2\nonumber\\ & &+ t^{m_l }\sum_{i=1}^{m}\tilde
c_{i}\left\vert \nabla
^{i}\operatorname*{Rm}\right\vert \left\vert \nabla^{m+1-i}\operatorname*{Rm}%
\right\vert \left\vert \nabla^{m+1}\operatorname*{Rm}\right\vert
+m_l  t^{m_l -1}\left\vert \nabla^{m+1}\operatorname*{Rm}\right\vert
^{2}.\nonumber
\end{eqnarray}

Using the inductive hypothesis, Inequality~(\ref{strongineq}), there is a constant
$A<\infty$ depending only on $K,\alpha,A_1,\ldots,A_m$ such that
$$
\sum_{i=1}^{m}\tilde c_{i}\left\vert \nabla
^{i}\operatorname*{Rm}\right\vert \left\vert
\nabla^{m+1-i}\operatorname*{Rm} \right\vert \leq At^{-m_l/2}.
$$
Also, let $c=\tilde c_{0}+\tilde c_{m+1}$ and define a new  constant
$B$
by
$$B=c(\alpha+K)+m_l.$$ Then, since $t\le T\le \alpha/K$ and $m_l\ge 0$, we have
$$((\tilde c_{0}+\tilde c_{m+1})t\left\vert\operatorname*{Rm}\right\vert +m_l)t^{m_l-1}\le
\frac{Bt^{m_l}}{t^{s(m_l)}}.$$ Putting this together allows us to
rewrite Inequality~(\ref{1steqn}) as
\begin{align*}
\frac{\partial}{\partial t}\left(  t^{m_l}\left\vert
\nabla^{m+1}\operatorname*{Rm}\right\vert ^{2}\right)
 \leq &\Delta\left(  t^{m_l }\left\vert \nabla
^{m+1}\operatorname*{Rm}\right\vert ^{2}\right)  -2t^{m_l }
\left\vert \nabla^{m+2}\operatorname*{Rm}\right\vert ^{2}\\
 &
+At^{m_l /2}\left\vert \nabla^{m+1} \operatorname*{Rm}\right\vert
+\left(ct\left\vert \operatorname*{Rm}\right\vert +m_l \right)
t^{m_l -1}\left\vert \nabla^{m+1}\operatorname*{Rm}\right\vert
^{2}\\
 \leq &\Delta\left(  t^{m_l }\left\vert \nabla
^{m+1}\operatorname*{Rm}\right\vert ^{2}\right)  -2t^{m_l }
\left\vert \nabla^{m+2}\operatorname*{Rm}\right\vert ^{2}\\
& +
\frac{B}{t^{s(m_l)}}t^{m_l}\left\vert\nabla^{m+1}\operatorname*{Rm}\right\vert^2
+ At^{m_l/2}\left\vert\nabla^{m+1}\operatorname*{Rm}\right\vert.
\end{align*}
Completing the square gives
\begin{align*}
\frac{\partial}{\partial t}\left(  t^{m_l}\left\vert
\nabla^{m+1}\operatorname*{Rm}\right\vert ^{2}\right)
 &
\leq\Delta\left( t^{m_l }\left\vert \nabla
^{m+1}\operatorname*{Rm}\right\vert ^{2}\right) -2t^{m_l  }
\left\vert \nabla^{m+2}\operatorname*{Rm}\right\vert ^{2}\\
&  +(B+1)t^{m_l  -s(m_l)}\left\vert \nabla^{m+1}
\operatorname*{Rm}\right\vert ^{2}+\frac{A^2}{4}t^{s(m_l)}.
\end{align*}

Let $\hat{m}_l = \max\left\{  m-l,0\right\} $. From
(\ref{curderfor}) and the induction hypothesis, there is a constant
$D$, depending on $K,\alpha,A_1,\ldots,A_m$ such that

\begin{eqnarray*}
\frac{\partial}{\partial t}\left(  t^{\hat{m}_l }\left\vert
\nabla^{m}\operatorname*{Rm}\right\vert ^{2}\right) & \leq &
\Delta\left( t^{\hat{m}_l} \left\vert
\nabla^{m}\operatorname*{Rm}\right\vert ^{2}\right) -2t^{\hat{m}_l
}\left\vert
\nabla^{m+1}\operatorname*{Rm}\right\vert ^{2}\\
& & +\hat{m}_l  t^{\hat{m}_l  -1}\left\vert
\nabla^{m}\operatorname*{Rm}\right\vert ^{2}+D.
\end{eqnarray*}

Now, defining new constants $\widetilde B=B+1$ and $\widetilde
A=A^2/4$ we have

\begin{eqnarray*}
\lefteqn{ \left(  \frac{\partial}{\partial t}-\Delta\right)F_m=
\left( \frac{\partial}{\partial t}-\Delta\right)  \left[ \left(
C+t^{\hat{m}_l  }\left\vert \nabla^{m}\operatorname*{Rm} \right\vert
^{2}\right)  t^{m_l }\left\vert
\nabla^{m+1}\operatorname*{Rm}\right\vert ^{2}\right]  \leq}   \\
& & \left(  C+t^{\hat{m}_l  }\left\vert \nabla^{m}
\operatorname*{Rm}\right\vert ^{2}\right)  \left( -2t^{m_l
}\left\vert \nabla^{m+2} \operatorname*{Rm} \right\vert ^{2}
+\frac{\widetilde B}{t^{s\{m_l\}}} t^{m_l}\left\vert \nabla^{m+1}
\operatorname*{Rm} \right\vert ^{2}+\widetilde At^{s(m_l)} \right) \\
& &  +\left( -2t^{\hat{m}_l }\left\vert
\nabla^{m+1}\operatorname*{Rm} \right\vert ^{2} +\hat{m}_l
t^{\hat{m}_l  -1}\left\vert \nabla^{m}\operatorname*{Rm}\right\vert
^{2}+D \right) t^{m_
l} \left\vert \nabla^{m+1} \operatorname*{Rm}\right\vert ^{2}\\
& &  -2t^{\hat{m}_l  +m_l } \nabla\left(\left\vert
\nabla^{m}\operatorname*{Rm}\right\vert ^{2}\right)\cdot
\nabla\left(\left\vert \nabla^{m+1}\operatorname*{Rm}\right\vert
^{2}\right). \end{eqnarray*}

Since $C\ge 4t^{\hat
m_l}\left\vert\nabla^m\operatorname*{Rm}\right\vert^2$, this implies

\begin{eqnarray}
\lefteqn{ \left(  \frac{\partial}{\partial t}-\Delta\right) F_m \leq
-10t^{\hat{m}_l  +m_l} \left\vert
\nabla^{m}\operatorname*{Rm}\right\vert ^{2}\left\vert
\nabla^{m+2}\operatorname*{Rm}\right\vert ^{2}}\label{2ndeqn}\\
 & & -8t^{\hat{m}_l  +m_l  }\left\vert
\nabla^{m}\operatorname*{Rm}\right\vert \left\vert \nabla^{m+1}%
\operatorname*{Rm}\right\vert ^{2}\left\vert \nabla^{m+2}\operatorname*{Rm}%
\right\vert  -2t^{\hat m_l  +m_l  }\left\vert
\nabla^{m+1}\operatorname*{Rm}\right\vert ^{4}\nonumber\\
 & & +\left(  C+t^{\hat m_l}\left\vert \nabla^{m}
\operatorname*{Rm}\right\vert ^{2}\right)  \left( \widetilde
Bt^{m_l-s(m_l) }\left\vert \nabla^{m+1}\operatorname*{Rm}\right\vert
^{2}+\widetilde At^{s(m_l)}\right) \nonumber\\
 & & +\left(\hat m_l  t^{\hat{m}_l -1}\left\vert
\nabla^{m}\operatorname*{Rm}\right\vert ^{2}+D\right) t^{m_l
}\left\vert \nabla^{m+1}\operatorname*{Rm} \right\vert
^{2}.\nonumber
\end{eqnarray}

Now we can write the first three terms on the right-hand side of
Inequality~(\ref{2ndeqn}) as \begin{equation}\label{3rdeqn} -t^{\hat
m_l+m_l}\left(\sqrt{10}\left\vert\nabla^{m+2}\operatorname*{Rm}\right\vert
\left\vert\nabla^m\operatorname*{Rm}\right\vert+\frac{4}{\sqrt{10}}\left\vert\nabla^{m+1}
\operatorname*{Rm}\right\vert^2\right)^2-\frac{2}{5}t^{\hat
m_l+m_l}\left\vert\nabla^{m+1}\operatorname*{Rm}\right\vert^4.\end{equation}
In addition we have
\begin{equation}\label{4theqn}
C+t^{\hat m_l}\left\vert \nabla^{m} \operatorname*{Rm}\right\vert
^{2}  \leq  C+A_m^2.
\end{equation}

Let us set $\widetilde D=\max(\alpha/K,1)D$. If $\hat m_l=0$, then
\begin{equation}\label{5theqn}
\hat m_lt^{\hat
m_l-1}\left\vert\nabla^m\operatorname*{Rm}\right\vert^2+D =D\le
\frac{\widetilde D}{t^{s(m_l)}}=\hat m_lA_m^2+D\le \frac{\hat
m_lA_m^2+\widetilde D}{t^{s(m_l)}}.\end{equation}  On the other
hand, if $\hat m_l>0$, then $s(\hat m_l)=s(m_l)=1$ and hence
$$\hat m_lt^{\hat m_l-1}\left\vert\nabla^m\operatorname*{Rm}\right\vert^2+D\le
\frac{1}{t^{s(m_l)}}\hat m_lA_m^2+D\le \frac{\hat
m_lA_m^2+\widetilde D}{t^{s(m_l)}}.$$

Since $\hat m_l =m_l-s(m_l)$,
Inequalities~(\ref{3rdeqn}),~(\ref{4theqn}), and~(\ref{5theqn}) then
allow us  rewrite Inequality~(\ref{2ndeqn}) as
\begin{eqnarray*}
\lefteqn{ \left(  \frac{\partial}{\partial t}-\Delta\right)F_m
 \leq   -\frac{2}{5t^{s( m_l)  }
 }t^{2m_l}\left\vert\nabla^{m+1}\operatorname*{Rm}\right\vert^4}\\
& &  +(C+A_m^2)\left(\frac{\widetilde B}{t^{s(m_l)}}
t^{m_l}\left\vert\nabla^{m+1}\operatorname*{Rm}\right\vert^2+\widetilde
At^{s(m_l)}\right) +\frac{\hat m_lA_m^2+\widetilde
D}{t^{s(m_l)}}t^{m_l}\left\vert\nabla^{m+1}\operatorname*{Rm}
\right\vert^2.
\end{eqnarray*}

Setting
$$B'=(C+A_m^2)\widetilde B+(\hat m_lA_m^2+\widetilde D),$$
and $A'=\widetilde A(C+A_m^2)$ we have
\begin{eqnarray*}
\left(  \frac{\partial}{\partial t}-\Delta\right) F_m  & \leq &
-\frac{2}{5t^{s( m_l)  }
 }\left(t^{m_l}\left\vert\nabla^{m+1}\operatorname*{Rm}\right\vert^2\right)^2 \\
 & & +
\frac{B'}{t^{s(m_l)}}t^{m_l}\left\vert\nabla^{m+1}\operatorname*{Rm}\right\vert^2
+A't^{s(m_l)}.
\end{eqnarray*}
We rewrite this as
\begin{eqnarray*}
\left(  \frac{\partial}{\partial t}-\Delta\right) F_m & \leq &
-\frac{2}{5t^{s(m_l)}}\left(t^{m_l}\left\vert\nabla^{m+1}\operatorname*{Rm}
\right\vert^2   -\frac{5B'}{4}\right)^2
\\ & & +\frac{5(B')^2}{8t^{s(m_l)}} +A't^{s(m_l)},
\end{eqnarray*}
and hence
\begin{equation*}
\left(  \frac{\partial}{\partial t}-\Delta\right) F_m \leq
-\frac{2}{5t^{s(m_l)}}\left(t^{m_l}\left\vert\nabla^{m+1}\operatorname*{Rm}
\right\vert^2-B''\right)^2  +\frac{A''}{t^{s(m_l)}}
\end{equation*}
where the constants $B''$ and $A''$ are defined by $ B''=5 B'/4$ and
$$A''=(\max\{\alpha/K,1\})^2+5( B')^2/8.$$ (Recall that $t\le T\le
\alpha/K$.) Let
$$Y=(C+t^{\hat
m_l}\left\vert\nabla^m\operatorname*{Rm}\right\vert^2).$$ (Notice
that $Y$ is not a constant.) Of course, by definition
$$F_m=Yt^{m_l}|\nabla^{m+1}{\rm Rm}|^2.$$ Then the previous inequality
becomes
\begin{equation*}
\left(\frac{\partial}{\partial t}-\Delta\right)  F_m \leq
-\frac{2}{5t^{s(m_l)}Y^2}\left(Yt^{m_l}\left\vert\nabla^{m+1}\operatorname*{Rm}
\right\vert^2-B'' Y\right)^2
 +\frac{A''}{t^{s(m_l)}}
\end{equation*}
Since $C\le Y\le 5C/4$ we have
\begin{eqnarray*}
\left(  \frac{\partial}{\partial t}-\Delta\right)  F_m \leq
-\frac{32}{125t^{s(m_l)}C^2}\left(F_m-B''
Y\right)^2+\frac{A''}{t^{s(m_l)}}
\end{eqnarray*}
At any point where $F_m\ge 5CB''/4$,  the last inequality gives
$$\left(  \frac{\partial}{\partial t}-\Delta\right)  F_m \leq
-\frac{32}{125t^{s(m_l)}C^2}\left(F_m-5CB''/4\right)^2+\frac{A''}{t^{s(m_l)}}.$$
At any point where $F_m\le 5CB''/4$,  since $F_m\ge 0$ and $0\le
B''Y\le 5CB''/4$, we have $(F_m-B''Y)^2\le 25C^2(B'')^2/16$, so that
$$-\frac{32}{125t^{s(m_l)}C^2}\left(F_m-5CB''/4\right)^2\ge
-2(B'')^2/5t^{s(m_l)}.$$ Thus, in this case we have
$$\left(\frac{\partial}{\partial t}-\Delta\right)  F_m \leq
\frac{A''}{t^{s(m_l)}}\le
-\frac{32}{125t^{s(m_l)}C^2}\left(F_m-5CB''/4\right)^2+\frac{A''+2(B'')^2/5}{t^{s(m_l)}}.$$

These two cases together prove Claim~\ref{Fm}.
\end{proof}

\section{Generalized Ricci flows}

In this section we introduce a generalization of the Ricci flow
equation. The generalization does not involve changing the PDE that
gives the flow. Rather it allows for the global topology of
space-time to be different from a product.
\subsection{Space-time}

There are two basic ways to view an $n$-dimensional Ricci flow: (i)
as a one-parameter family of metrics $g(t)$ on a fixed smooth
$n$-dimensional manifold $M$, and (ii) as a partial metric (in the
horizontal directions) on the $(n+1)$-dimensional manifold $M\times
I$. We call the latter $(n+1)$-dimensional manifold {\sl
space-time}\index{Ricci flow!generalized!space-time|ii} and the
horizontal slices are the {\em time-slices}. In defining the
generalized Ricci flow, it is the second approach that we
generalize.

\begin{defn}
By {\em space-time} we mean a smooth $(n+1)$-dimensional manifold
${\mathcal M}$ (possibly with boundary), equipped with a smooth
function ${\bf t}\colon {\mathcal M}\to \Ar$, called {\em time} and
a smooth vector field $\chi$ subject to the following axioms:
\begin{enumerate}
\item The image of ${\bf t}$ is an interval $I$ (possibly infinite) and
the boundary of ${\mathcal M}$ is the preimage under ${\bf t}$ of
$\partial I$.
\item For each $x\in {\mathcal M}$ there is an open neighborhood $U\subset
{\mathcal M}$ of $x$ and a diffeomorphism $f\colon  V\times J\to U$,
where $V$ is an open subset in $\Ar ^n$ and $J$ is an interval with
the property that (i) ${\bf t}$ is the composition of $f^{-1}$
followed by the projection onto the interval $J$ and (ii) $\chi$ is
the image under $f$ of the unit vector field in the positive
direction tangent to the foliation by the lines $\{v\}\times J$ of
$V\times J$.
\end{enumerate}
\end{defn}

 Notice that it follows
that $\chi({\bf t})=1$.

\begin{defn}
The {\em time-slices}\index{Ricci flow!generalized!time-slices|ii}
of space-time are the level sets ${\bf t}$. These form a
codimension-one foliation of ${\mathcal M}$. For each $t\in I$ we
denote by $M_t\subset {\mathcal M}$ the $t$ time-slice, that is to
say ${\bf t}^{-1}(t)$. Notice that each boundary component of
${\mathcal M}$ is contained in a single time-slice. The {\em
horizontal distribution}, ${\mathcal H}T{\mathcal M}$ is the
distribution tangent to this foliation. A {\em horizontal
metric}\index{Ricci flow!generalized!horizontal metric|ii} on
space-time is a smoothly varying positive definite inner product on
${\mathcal H}T{\mathcal M}$.
\end{defn}

Notice that a horizontal metric on space-time induces an ordinary
Riemannian metric on each time-slice. Conversely, given a Riemannian
metric on each time-slice $M_t$, the condition that they fit
together to form a horizontal metric on space-time is that they vary
smoothly on space-time. We define the curvature of a horizontal
metric $G$ to be the section of the dual of the symmetric square of
$\wedge^2{\mathcal H}T{\mathcal M}$ whose value at each point $x$ with ${\bf
t}(x)=t$ is the usual Riemannian curvature tensor of the induced
metric on $M_t$  at the point $x$. This is a smooth section of ${\rm
Sym}^2(\wedge^2{\mathcal H}T^*{\mathcal M})$.
 The Ricci curvature and the scalar curvature of a horizontal metric
 are given in the usual way from its Riemannian curvature.
 The Ricci curvature is a smooth section of ${\rm Sym}^2({\mathcal
 H}T^*{\mathcal M})$ while the scalar curvature is a smooth function on ${\mathcal M}$.

\subsection{The generalized Ricci flow equation}

Because of the second condition in the definition of space-time, the
vector field $\chi$ preserves the horizontal foliation and hence the
horizontal distribution. Thus, we can form the Lie derivative of a
horizontal metric with respect to $\chi$.

\begin{defn}
An {\em $n$-dimensional generalized Ricci flow}\index{Ricci
flow!generalized|ii} consists of a space-time ${\mathcal M}$ that is
$(n+1)$-dimensional and a horizontal metric $G$ satisfying the
generalized Ricci flow equation:
$${\mathcal L}_\chi(G)=-2{\rm Ric}(G).$$
\end{defn}

\begin{rem}
Let $({\mathcal M},G)$ be a generalized Ricci flow and let $x\in
{\mathcal M}$. Pulling $G$ back to  the local coordinates $V\times
J$ defined near any point gives a one-parameter family of metrics
$(V,g(t)),\  t\in J$, satisfying the usual Ricci flow equation. It
follows that all the usual evolution formulas for Riemannian
curvature, Ricci curvature, and scalar curvature hold in this more
general context.
\end{rem}

Of course, any ordinary Ricci flow is a generalized Ricci flow where
space-time is a product $M\times I$ with time being the projection
to $I$ and $\chi$ being the unit vector field in the positive
$I$-direction.

\subsection{More definitions for generalized Ricci flows}

\begin{defn}
Let ${\mathcal M}$ be a space-time. Given a space $C$ and an
interval $I\subset \Ar$  we say that an embedding $C\times I\to
{\mathcal M}$ is {\em compatible with}\index{Ricci
flow!generalized!compatible embedding} the time and the vector field
if: (i) the restriction of ${\bf t}$ to the image agrees with the
projection onto the second factor and (ii) for each $c\in C$ the
image of $\{c\}\times I$ is the integral curve for the vector field
$\chi$.
 If in addition $C$ is a subset of $M_{t}$
we require that $t\in I$ and that the map $C\times\{t\}\to M_{t}$ be
the identity. Clearly, by the uniqueness of integral curves for vector fields,
 two such embeddings agree on their common interval of definition, so that,
given $C\subset M_{t}$ there is a maximal interval $I_C$ containing
$t$ such that such an embedding, compatible with time and the vector field,
 is defined on $C\times I$. In the
special case when $C=\{x\}$ for a point $x\in M_{t}$ we say that
such an embedding is {\em the flow line} through $x$. The embedding
of the maximal interval through $x$ compatible with time and the vector
field $\chi$ is called {\em the domain of definition} of the flow
line through $x$. For a more general subset $C\subset M_{t}$ there
is an embedding $C\times I$ compatible with time and the vector
field $\chi$ if an only if for every $x\in C$, $I$ is contained in
the domain of definition of the flow line through $x$.
\end{defn}

\begin{defn}
We say that $t$ is a {\em regular time} if there is $\epsilon>0$ and
a diffeomorphism $M_{t}\times (t-\epsilon,t+\epsilon)\to {\bf
t}^{-1}((t-\epsilon,t+\epsilon))$ compatible with time
and the vector field. A time is {\em singular} if it is not regular.
Notice that if all times are regular, then space-time is a product
$M_{t}\times I$ with ${\bf t}$ and $\chi$ coming from the second
factor. If the image ${\bf t}({\mathcal M})$ is an interval $I$
bounded below, then the {\em initial time} for the flow is the
greatest lower bound for $I$. If $I$ includes $(-\infty,A]$ for some
$A$, then the initial time for the generalized Ricci flow is
$-\infty$.
\end{defn}

\begin{defn}
Suppose that $({\mathcal M},G)$ is a generalized Ricci flow and that
$Q>0$ is a positive constant. Then we can define a new generalized
Ricci flow by setting $G'=QG$, ${\bf t'}=Q{\bf t}$ and
$\chi'=Q^{-1}\chi$. It is easy to see that the result still
satisfies the generalized Ricci flow equation. We denote this new
generalized Ricci flow by $(Q{\mathcal M},QG)$ where the changes in
${\bf t}$ and $\chi$ are denoted by the factor of $Q$ in front of
${\mathcal M}$.

It is also possible to translate a generalized solution $({\mathcal M},G)$ by
replacing the time function ${\bf t}$ by ${\bf t'}={\bf t}+a$ for any constant
$a$, leaving $G$ and $\chi$ unchanged.
\end{defn}

\begin{figure}[ht]
  \relabelbox{
  \centerline{\epsfbox{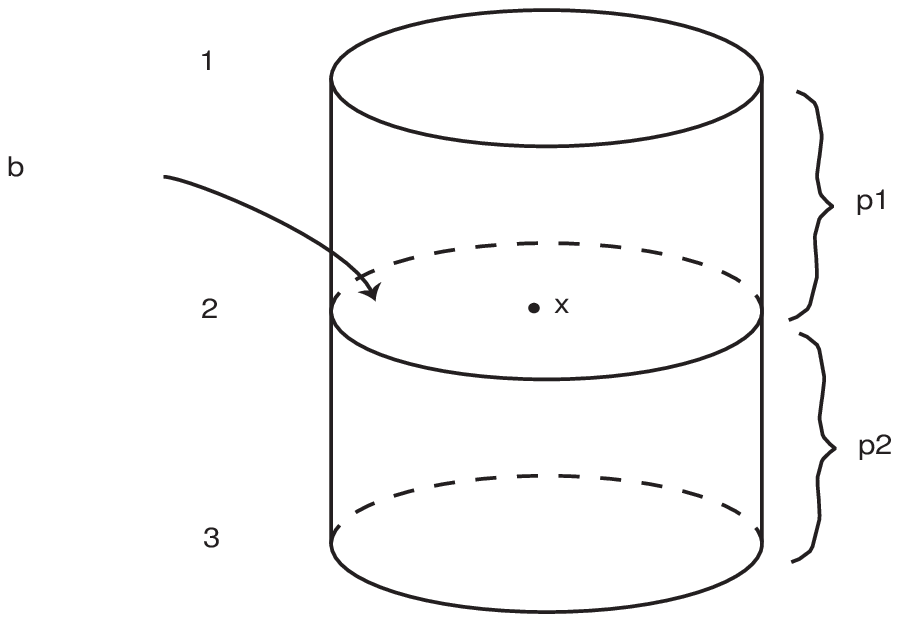}}}
  \relabel{1}{$t+\Delta t$}
  \relabel{2}{$t$}
  \relabel{3}{$t-\Delta t$}
  \relabel{b}{$B(x,t,r)$}
  \relabel{x}{$x$}
  \relabel{p1}{$P(x,t,r,\Delta t)$}
  \relabel{p2}{$P(x,t,r,-\Delta t)$}
  \endrelabelbox
  \caption{Parabolic neighborhoods}\label{fig:paranbhd}
\end{figure}

\begin{defn}
Let $({\mathcal M},G)$ be a generalized Ricci flow and let $x$ be a
point of space-time. Set $t={\bf t}(x)$. For any  $r>0$ we define
$B(x,t,r)\subset M_t$ to be the metric ball of radius $r$ centered
at $x$ in the Riemannian manifold $(M_t,g(t))$. For any $\Delta t>0$
we say that $P(x,t,r,\Delta t)$,
respectively, $P(x,r,t,-\Delta
t)$, {\em exists} in ${\mathcal M}$ if there is an embedding
$B(x,t,r)\times [t,t+\Delta t]$, respectively, $B(x,t,r)\times
[t-\Delta t,t]$, into ${\mathcal M}$ compatible with time and the
vector field. When this embedding exists, its image is defined to be the {\em
forward parabolic neighborhood}\index{parabolic neighborhood|ii}
$P(x,t,r,\Delta t)$, respectively the {\em backward parabolic
neighborhood} $P(x,t,r,-\Delta t)$.
See {\sc Fig.}~\ref{fig:paranbhd}.
\end{defn}

\chapter{The maximum principle}\label{secmaxprin}
 Recall that the maximum principle\index{maximum principle!heat
equation} for the heat equation says that if $h$ is a  solution to
the heat equation $$\frac{\partial h}{\partial t}=\Delta h$$ on a
compact manifold and if $h(x,0)\ge 0$ for all $x\in M$, then
$h(x,t)\ge 0$ for all $(x,t)$. In this chapter we discuss analogues
of this result for the scalar curvature, the Ricci curvature, and
the sectional curvature under Ricci flow. Of course, in all three
cases we are working with quasi-linear versions of the heat equation
so it is important to control the lower order (non-linear) terms and
in particular show that at zero curvature they have the appropriate
sign. Also, in the latter two cases we are working with tensors
rather than with scalars and hence we require a tensor version of
the maximum principle, which was established by Hamilton in
\cite{Hamilton4MPIC}.

As further applications of these results beyond just establishing
non-negativity, we indicate Hamilton's result that if the initial
conditions have positive Ricci curvature then the solution becomes
singular at finite time and as it does it becomes round (pinching to
round). We also give Hamilton's result showing that at points where
the scalar curvature is sufficiently large the curvature is pinched
toward positive. This result is crucial for understanding
singularity development. As a last application, we give Hamilton's
Harnack inequality for Ricci flows of non-negative curvature.

The maximum principle is used here in two different ways. The first
assumes non-negativity of something (e.g., a curvature) at time zero
and uses the maximum principle to establish non-negativity of this
quantity at all future times. The second assumes non-negativity of
something at all times and positivity at one point, and then uses
the maximum principle to establish positivity at all points and all
later times. In the latter application one compares the solution
with a solution to the linear heat equation where such a property is
known classically to hold.

\section{Maximum principle for scalar curvature}

Let\index{maximum principle!scalar curvature|(ii} us begin with the
easiest evolution equation, that for the scalar curvature, where the
argument uses only the (non-linear) version of the maximum
principle. This result is valid in all dimensions:

\begin{prop}\label{scalarevol}
Let $(M,g(t)),\ 0\le t<T$, be a Ricci flow with $M$ a compact
$n$-dimensional manifold. Denote by $R_{\rm min}(t)$ the minimum
value of the scalar curvature of $(M,g(t))$. Then:
\begin{itemize} \item $R_{\rm min}(t)$ is a non-decreasing function
of $t$.
\item If $R_{\rm
min}(0)\ge 0$, then
$$R_{\rm min}(t)\ge R_{\rm min}(0)\left(\frac{1}{1-\frac{2t}{n}R_{\rm
min}(0)}\right),$$ in particular,
$$T\le\frac{n}{2R_{\rm min}(0)}.$$
\item If $R_{\rm min}(0)<0$, then
$$R_{\rm min}(t)\ge
-\frac{n\bigl| R_{\rm min}(0)\bigr|}{2t\bigl|R_{\rm
min}(0)\bigr|+n}.$$
\end{itemize}
\end{prop}

\begin{proof}
According to Equation~(\ref{Revol}), the evolution equation for $R$
is
$$\frac{\partial}{\partial t}R(x,t)=\Delta
R(x,t)+2|{\rm Ric}(x,t)|^2.$$ Since $M$ is compact, the function
$R_{\rm min}(t)$ is continuous but may not be $C^1$ at points where
the minimum of the scalar curvature is achieved at more than one
point.

The first thing to notice is the following:

\begin{claim}\label{firstclaim}  If $R(x,t)=R_{\rm
min}(t)$ then $(\partial R/\partial t)(x,t)\ge \frac{2}{n}R^2(x,t)$.
\end{claim}

\begin{proof}
This is immediate from  the evolution equation for $R$, the fact
that if $R(x,t)=R_{\rm min}(t)$, then $\Delta R(x,t)\ge 0$, and the
fact that $R$ is the trace of ${\rm Ric}$ which implies by the
Cauchy-Schwarz inequality that $|R|^2\le n|{\rm Ric}|^2$.
\end{proof}

Now it follows that:

\begin{claim}\label{Rminprime}
$$\frac{d}{d t}(R_{\rm min}(t))\ge \frac{2}{n}R_{\rm min}^2(t),$$
where, at times $t$ where $R_{\rm min}(t)$ is not smooth, this
inequality is interpreted as an inequality for the forward
difference quotients.
\end{claim}

\begin{proof}
This is immediate from the first statement in
Proposition~\ref{fordiffmax}.
\end{proof}

If follows immediately from Claim~\ref{Rminprime} and
Lemma~\ref{fordiffquot} that $R_{\rm min}(t)$ is a non-decreasing
function of $t$. This establishes the first item and also the second
item in the case when $R_{\rm min}(0)=0$.

Suppose that $R_{\rm min}(0)\not= 0$.
 Consider the function
$$S(t)=\frac{-1}{R_{\rm min}(t)}-\frac{2t}{n}+\frac{1}{R_{\rm min}(0)}.$$
 Clearly, $S(0)=0$ and $S'(t)\ge 0$ (in the sense of forward difference quotients),
  so that by  Lemma~\ref{fordiffquot} we have $S(t)\ge 0$ for all $t$. This means that
\begin{equation}\label{Rminineq} \frac{1}{R_{\rm min}(t)}\le
\frac{1}{R_{\rm min}(0)}-\frac{2t}{n}\end{equation} provided that
$R_{\rm min}$ is not ever zero on the interval $[0,t]$. If $R_{\rm
min}(0)>0$, then by the first item, $R_{\rm min}(t)>0$ for all $t$
for which the flow is defined, and the inequality in the second item
of the proposition is immediate from Equation~(\ref{Rminineq}). The
third inequality in the proposition also follows easily from
Equation~(\ref{Rminineq}) when $R_{\rm min}(t)<0$. But if $R_{\rm
min}(t)\ge 0$, then the third item is obvious.
\end{proof}\index{maximum principle!scalar curvature|)}

\section{The maximum principle for tensors}

For\index{maximum principle!tensors|(ii} the applications to the
Ricci curvature and the curvature tensor we need a version of the
maximum principle for tensors that is due to Hamilton; see
\cite{Hamilton4MPCO}.

Suppose that $V$ is a finite-dimensional real vector space and $Z\subset V$ is
a closed convex set. For each $z$ in the frontier of $ Z$ we define the {\em
tangent cone to} $Z$ at $z$, denoted $T_zZ$, to be the intersection of all
closed half-spaces $H$ of $V$ such that $z\in
\partial H$ and $Z\subset H$. For $z\in {\rm int}\, Z$ we define
$T_zZ=V$. Notice that $v\notin T_zZ$ if and only if there is a
affine linear function $\ell$ vanishing at $z$ non-positive on $Z$
and positive on $v$.

\begin{defn}
Let $Z$ be a closed convex subset of a finite-dimensional real vector space
$V$. We say that a smooth vector field $\psi$ defined on an open neighborhood
$U$ of $Z$ in $V$ {\em preserves} $Z$ if for every $z\in Z$ we have $\psi(z)\in
T_zZ$.
\end{defn}

It is an easy exercise to show the following; see Lemma 4.1 on page
183 of \cite{Hamilton4MPCO}:

\begin{lem}
Let $Z$ be a closed convex subset in a finite dimensional real
vector space $V$. Let $\psi$ be a smooth vector field defined on an
open neighborhood of $Z$ in $V$. Then $\psi$ preserves $Z$ if and
only if  every integral curve $\gamma\colon [0,a)\to V$ for $\psi$
with $\gamma(0)\in Z$ has $\gamma(t)\in Z$ for all $t\in [0,a)$.
Said more informally, $\psi$ preserves $Z$ if and only if every
integral curve for $\psi$ that starts in $Z$ remains in $Z$.
\end{lem}

\subsection{The global
version}\label{4.3.3}

The maximum principle for tensors generalizes this to tensor flows
evolving by parabolic equations. First we introduce a generalization
of the notion of a vector field preserving a closed convex set to
the context of vector bundles.

\begin{defn}
 Let $\pi\colon {\mathcal V}\to M$ be a vector bundle and let
${\mathcal Z}\subset {\mathcal V}$ be  a closed subset. We say that
${\mathcal Z}$ is {\em convex} if  for every $x\in M$ the fiber
$Z_x$ of ${\mathcal Z}$ over $x$ is a convex subset of the vector
space fiber $V_x$ of ${\mathcal V}$ over $x$. Let $\psi$ be a
fiberwise vector field on an open neighborhood ${\mathcal U}$ of
${\mathcal Z}$ in ${\mathcal V}$. We say that $\psi$ {\em preserves}
${\mathcal Z}$ if for each $x\in M$ the restriction of $\psi$ to the
fiber $U_x$ of ${\mathcal U}$ over $x$ preserves $Z_x$.
\end{defn}

The following global version of the maximum principle for tensors is
Theorem 4.2 of  \cite{Hamilton4MPCO}.

\begin{thm}\label{MP}{\bf (The maximum principle for tensors)}
Let $(M,g)$ be a compact Riemannian manifold. Let ${\mathcal V}\to M$ be a
tensor bundle and let ${\mathcal Z}\subset {\mathcal V}$ be a closed, convex
subset invariant under the parallel translation induced by the Levi-Civita
connection. Suppose that $\psi$ is a fiberwise vector field defined on an open
neighborhood of ${\mathcal Z}$ in ${\mathcal V}$ that preserves ${\mathcal Z}$.
Suppose that ${\mathcal T}(x,t),\ 0\le t\le T$, is a one-parameter family of
sections of ${\mathcal V}$ that evolves according to the parabolic equation
$$\frac{\partial {\mathcal T}}{\partial t}=\Delta {\mathcal T}
+\psi({\mathcal T}).$$ If ${\mathcal T}(x,0)$ is contained in
${\mathcal Z}$ for all $x\in M$, then ${\mathcal T}(x,t)$ is
contained in ${\mathcal Z}$ for all $x\in M$ and for all $0\le t\le
T$.
\end{thm}

For a proof we refer the reader to  Theorem 4.3 and its proof (and
the related Theorem 4.2 and its proof) in  \cite{Hamilton4MPCO}.

There is a slight improvement of this result where the convex set ${\mathcal
Z}$ is allowed to vary with $t$. It is proved by the same argument; see Theorem
4.8 on page 101 of \cite{ChowKnopf}.

\begin{thm}\label{varyingZ}
Let $(M,g)$ be a compact Riemannian manifold. Let ${\mathcal V}\to M$ be a
tensor bundle and let ${\mathcal Z}\subset {\mathcal V}\times [0,T]$ be a
closed subset with the property that for each $t\in [0,T]$ the time-slice
${\mathcal Z}(t)$ is a  convex subset of ${\mathcal V}\times\{t\}$ invariant
under the parallel translation induced by the Levi-Civita connection. Suppose
that $\psi$ is a fiberwise vector field defined on an open neighborhood of
${\mathcal Z}$ in ${\mathcal V}\times [0,T]$ that preserves the family
${\mathcal Z}(t)$ in the sense that any integral curve $\gamma(t),\ t_0\le t\le
t_1$, for $\psi$ with the property that $\gamma(t_0)\in {\mathcal Z}(t_0)$ has
$\gamma(t)\in {\mathcal Z}(t)$ for every $t\in [t_0,t_1]$. Suppose that
${\mathcal T}(x,t),\ 0\le t\le T$, is a one-parameter family of sections of
${\mathcal V}$ that evolves according to the parabolic equation
$$\frac{\partial {\mathcal T}}{\partial t}=\Delta {\mathcal T}+\psi({\mathcal T}).$$
If ${\mathcal T}(x,0)$ is contained in
${\mathcal Z}(0)$ for all $x\in M$, then ${\mathcal T}(x,t)$ is
contained in ${\mathcal Z}(t)$ for all $x\in M$ and for all $0\le
t\le T$.
\end{thm}

\subsection{The local version}

Here is the local result. It is proved by the same argument as given
in the proof of Theorem 4.3 in \cite{Hamilton4MPCO}.

\begin{thm}\label{localMP}
Let $(M,g)$ be a Riemannian manifold. Let $\overline U\subset M$ be a compact,
smooth, connected, codimension-$0$ submanifold. Let ${\mathcal V}\to M$ be a
tensor bundle and let ${\mathcal Z}\subset {\mathcal V}$ be a closed, convex
subset. Suppose that $\psi$ is a fiberwise vector field defined on an open
neighborhood of ${\mathcal Z}$ in ${\mathcal V}$ preserving ${\mathcal Z}$.
Suppose that ${\mathcal Z}$ is invariant under the parallel translation induced
by the Levi-Civita connection. Suppose that ${\mathcal T}(x,t),\ 0\le t\le T$,
is a one-parameter family of sections of ${\mathcal V}$ that evolves according
to the parabolic equation
$$\frac{\partial {\mathcal T}}{\partial t}=\Delta {\mathcal T}
+\psi({\mathcal T}).$$ If ${\mathcal T}(x,0)$ is contained in
${\mathcal Z}$ for all $x\in \overline U$ and if ${\mathcal
T}(x,t)\in {\mathcal Z}$ for all $x\in \partial \overline U$ and all
$0\le t\le T$, then ${\mathcal T}(x,t)$ is contained in ${\mathcal
Z}$ for all $x\in \overline U$ and all $0\le t\le T$.
\end{thm}\index{maximum principle!tensors|)}

\section{Applications of the maximum principle}

Now let us give some applications of these results to Riemann and
Ricci curvature. In order to do this we first need to specialize the
above general maximum principles for tensors to the situation of the
curvature.

\subsection{Ricci flows with normalized initial conditions}

As we have already seen, the Ricci flow equation is invariant under multiplying
space and time by the same scale. This means that there can be no absolute
constants in the results about Ricci surgery. To break this gauge symmetry and
make the constants absolute we impose scale fixing (or rather scale bounding)
conditions on the initial metrics of the flows that we shall consider. The
following definition makes precise the exact conditions that we shall use.

\begin{defn}\label{norminitcond}
We say that a that a Ricci flow $(M,g(t))$ has {\em normalized
initial conditions}\index{Ricci flow!normalized initial conditions}
if $0$ is the initial time for the flow and if the compact
Riemannian manifold $(M^n,g(0))$ satisfies:
\begin{enumerate}
\item $|{\rm Rm}(x,0)|\le 1$ for all $x\in M$.
\item Let $\omega_n$ be the volume of the ball of radius $1$ in $n$ -dimensional Euclidean
space. Then ${\rm Vol}(B(x,0,r))\ge (\omega_n/2)r^n$ for any $p\in
M$ and any $r\le 1$.
\end{enumerate}
We also use the terminology $(M,g(0))$ {\em is normalized} to
indicate that it satisfies these two conditions.
\end{defn}

The evolution equation for the Riemann curvature and a standard maximum
principle argument show that if $(M,g(0))$ has an upper bound on the Riemann
curvature and a lower bound on the volume of balls of a fixed radius, then the
flow has Riemann curvature bounded above and volumes of balls bounded below on
a fixed time interval. Here is the result in the context of normalized initial
condition.

\begin{prop}\label{kappa0r0t0}
 There is $\kappa_0>0$  depending only on the dimension $n$
 such that the following holds.
Let $(M^n,g(t)),\ 0\le t\le T$, be a Ricci flow with bounded curvature, with
each $(M,g(t))$ being complete, and with normalized initial conditions. Then
$|{\rm Rm}(x,t)|\le 2$ for all $x\in M$ and all $t\in [0,{\rm min}(T,2^{-4})]$.
Furthermore, for any $t\in [0,{\rm min}(T,2^{-4})]$ and any $x\in M$ and any
$r\le 1$ we have ${\rm Vol}\,B(x,t,r)\ge \kappa_0r^n$.
\end{prop}

\begin{proof}
The bound on the Riemann curvature follows directly from Lemma 6.1
on page 207 of \cite{ChowLuNi} and the definition of normalized
initial conditions. Once we know that the Riemann curvature is
bounded by $2$ on $[0,2^{-4}]$, there is an $0<r_0$ depending on $n$
such that for every $x\in M$ and every $r\le r_0$ we have
$B(x,0,r_0r)\subset B(x,t,r)\subset B(x,0,1)$. Also, from the bound
on the Riemann curvature and the evolution equation for volume given
in Equation~(\ref{Revol}), we see that there is $A<\infty$ such that
${\rm Vol}_t\,(B(x,0,s))\ge A^{-1}{\rm Vol}_0\,(B(x,0,s))$. Putting
this together we see that
$${\rm Vol}_t\,(B(x,t,r)\ge A^{-1}(\omega_n/2)r_0^nr^n.$$
This proves the result.
\end{proof}

\subsection{Extending flows}

There is one other consequence that will be important for us. For a
reference see \cite{ChowLuNi} Theorem 6.3 on page 208.

\begin{prop}\label{flowextend}
Let $(M,g(t)),\ 0\le t<T<\infty$, be a Ricci flow with $M$ a compact manifold.
Then either the flow extends to an interval $[0,T')$ for some $T'>T$ or $|{\rm
Rm}|$ is unbounded on $M\times[0,T)$.
\end{prop}

\subsection{Non-negative curvature is preserved}

We need to consider the tensor versions of the maximum principle
when the tensor in question is the Riemann or Ricci curvature and
the evolution equation is that induced by the Ricci flow. This part
of the discussion is valid in dimension three only. We begin by
evaluating the expressions in Equation~(\ref{Mevol}) in the
$3$-dimensional case. Fix a symmetric bilinear form ${\mathcal S}$
on a $3$-dimensional real vector space $V$ with a positive definite
inner product. The inner product determines an identification of
$\wedge^2V$ with $V^*$. Hence, $\wedge^2{\mathcal S}^*$ is
identified with a symmetric automorphism of $V$,  denoted by
${\mathcal S}^\sharp$.

\begin{lem}\label{quad}
Let $(M,g)$ be a Riemannian $3$-manifold. Let ${\mathcal T}\in {\rm
Sym}^2(\wedge^2T_x^*M)$ be the curvature operator written with
respect to the evolving frame as in Proposition~\ref{Mevol}. Then
the evolution equation given in Proposition~\ref{Mevol} is:
$$\frac{\partial {\mathcal T}}{\partial t}=\triangle {\mathcal T}+\psi({\mathcal T})$$
where
$$\psi({\mathcal T})={\mathcal T}^2+{\mathcal T}^\sharp.$$
In particular, in an orthonormal basis in which
$${\mathcal T}=\begin{pmatrix} \lambda & 0 & 0 \\ 0 & \mu & 0 \\ 0 &
0  & \nu\end{pmatrix}$$ with $\lambda\ge \mu\ge \nu$, the vector
field is given by
$$\psi({\mathcal T})={\mathcal T}^2+{\mathcal T}^\sharp=\begin{pmatrix}
\lambda^2+\mu\nu & 0 & 0 \\ 0 & \mu^2+\lambda\nu & 0 \\ 0 & 0 &
\nu^2+\lambda\mu\end{pmatrix}.$$
\end{lem}

\begin{cor}\label{rm>0}
Let $(M,g(t)),\ 0\le t\le T$, be a Ricci flow with $M$ a compact,
connected $3$-manifold. Suppose that  ${\rm Rm}(x,0)\ge 0$ for all
$x\in M$. Then ${\rm Rm}(x,t)\ge 0$ for all $x\in M$ and all $t\in
[0,T]$.
\end{cor}

\begin{proof}
Let $\nu_x\colon {\rm Sym}^2(\wedge^2T^*_xM)\to \Ar$ associate to
each endomorphism its smallest eigenvalue. Then $\nu_x({\mathcal
T})$ is the minimum over all lines in $\wedge^2T_xM$ of the trace of
the restriction of ${\mathcal T}$ to that line. As a minimum of
linear functions, $\nu_x$ is a convex function. In particular,
$Z_x=\nu_x^{-1}([0,\infty))$ is a convex subset. We let ${\mathcal
Z}$ be the union over all $x$ of $Z_x$. Clearly, ${\mathcal Z}$ is a
closed convex subset of the tensor bundle. Since parallel
translation is orthogonal, ${\mathcal Z}$ is invariant under
parallel translation. The expressions in Lemma~\ref{quad} show that
if ${\mathcal T}$ is an endomorphism of $\wedge^2T^*_xM$ with
$\nu({\mathcal T})\ge 0$, then the symmetric matrix $\psi({\mathcal
T})$ is non-negative. This implies that $\nu_x$ is non-decreasing in
the direction $\psi({\mathcal T})$ at the point ${\mathcal T}$. That
is to say, for each $x\in M$, the vector field $\psi({\mathcal T})$
preserves the set $\{\nu_x^{-1}([c,\infty))\}$ for any $c\ge 0$. The
hypothesis that ${\rm Rm}(x,0)\ge 0$ means that ${\rm Rm}(x,0)\in
{\mathcal Z}$ for all $x\in M$. Applying Theorem~\ref{MP} proves the
result.
\end{proof}

\begin{cor}\label{Ric>0}
Suppose that $(M,g(t)),\ 0\le t\le T$, is a Ricci flow with $M$ a
compact, connected $3$-manifold with ${\rm Ric}(x,0)\ge 0$ for all
$x\in M$. Then ${\rm Ric}(x,t)\ge 0$ for all $t>0$.
\end{cor}

\begin{proof}
The statement that ${\rm Ric}(x,t)\ge 0$ is equivalent to the statement that
for every two-plane in $\wedge^2 T_xM$ the trace of the Riemann curvature
operator on this plane is $\ge 0$. For ${\mathcal T}\in {\rm
Sym}^2(\wedge^2T_x^*M)$, we define $s({\mathcal T})$ as the minimum over all
two-planes $P$ in $\wedge^2TM$ of the trace of ${\mathcal T}$ on $P$. The
restriction $s_x$ of $s$ to the fiber over $x$ is the minimum of a collection
of linear functions and hence is convex. Thus, the subset ${\mathcal
S}=s^{-1}([0,\infty))$ is convex. Clearly, $s$ is preserved by orthogonal
isomorphisms, so ${\mathcal S}$ is invariant under parallel translation. Let
$\lambda\ge\mu\ge\nu$ be the eigenvalues of ${\mathcal T}$. According to
Lemma~\ref{quad} the derivative of $s_x$ at ${\mathcal T}$ in the
$\psi({\mathcal T})$-direction is
$(\mu^2+\lambda\nu)+(\nu^2+\lambda\mu)=(\mu^2+\nu^2)+\lambda(\mu+\nu)$. The
condition that $s({\mathcal T})\ge 0$, is the condition that $\nu+\mu\ge 0$,
and hence $\mu\ge 0$, implying that $\lambda\ge 0$. Thus, if $s({\mathcal
T})\ge 0$, it is also the case that the derivatifve of $s_x$ in the
$\psi({\mathcal T})$-direction is non-negative. This implies that $\psi$
preserves ${\mathcal S}$. Applying Theorem~\ref{MP} gives the result.
\end{proof}

\section{The strong maximum principle for curvature}

First let us state the strong maximum principle for the heat
equation\index{maximum principle!heat equation!strong}.

\begin{thm}\label{SMPheat}
Let $\overline U$ be a compact, connected manifold, possibly with boundary. Let
$h(x,t),\ 0\le t\le T$, be a solution to the heat equation
$$\frac{\partial h(x,t)}{\partial t}=\Delta h(x,t).$$
 Suppose that  $h$ has Dirichlet boundary conditions in the sense
that $h(x,t)=0$ for all $(x,t)\in\partial \overline U\times [0,T]$.
If $h(x,0)\ge 0$ for all $x\in \overline U$, then $h(x,t)\ge 0$ for
all $(x,t)\in \overline U\times [0,T]$. If, in addition, there is
$y\in \overline U$ with $h(y,0)>0$, then $h(x,t)>0$ for all
$(x,t)\in {\rm int}(\overline U)\times(0,T]$.
\end{thm}

We shall use this strong maximum principle to establish an analogous
result for the curvature tensors. The hypotheses are in some ways
more restrictive -- they are set up to apply to the Riemann and
Ricci curvature.

\begin{prop}\label{SMPcurv}
Let $(M,g)$ be a Riemannian manifold and let ${\mathcal V}$ be a
tensor bundle. Suppose that $\overline U$ is a compact, connected,
smooth codimension-$0$ submanifold of $M$. Consider a one-parameter
family of sections ${\mathcal T}(x,t),\ 0\le t\le T$, of ${\mathcal
V}$. Suppose that ${\mathcal T}$ evolves according to the equation
$$\frac{\partial {\mathcal T}}{\partial t}=\Delta {\mathcal T}+\psi({\mathcal T})$$
for some smooth, fiberwise vector field $\psi({\mathcal T})$
defined on ${\mathcal V}$. Suppose that $s\colon {\mathcal V}\to
\Ar$ is a function satisfying the following properties:
\begin{enumerate}
\item For each $x\in M$ the restriction $s_x$ to the fiber $V_x$ of
${\mathcal V}$ over $x$ is a  convex function.
\item For any $A$ satisfying $s_x(A)\ge 0$ the vector $\psi(A)$ is contained
in the tangent cone of the convex set $\{y|s_x(y)\ge s_x(A)$ at the
point $A$.
\item $s$ is invariant under parallel translation.
\end{enumerate}
Suppose that $s({\mathcal T}(x,0))\ge 0$ for all $x\in \overline U$
and that $s({\mathcal T}(x,t))\ge 0$ for all $x\in \partial\overline
U$ and all $t\in [0,T]$. Suppose also that there is $x_0\in {\rm
int}(\overline U)$ with $s({\mathcal T}(x_0,0))>0$. Then
$s({\mathcal T}(x,t))> 0$ for all $(x,t)\in {\rm int}(\overline
U)\times(0,T]$.
\end{prop}

\begin{proof}
Let $h\colon \overline U\times\{0\}\to \Ar$ be a smooth function
with $h(x,0)=0$ for all $x\in \partial \overline U$ and with
$s({\mathcal T}(x,0))\ge h(x,0)\ge 0$ for all $x\in \overline U$. We
choose $h$ so that $h(x_0,0)>0$. Let $h(x,t),\ 0\le t<\infty$,  be
the solution to the heat equation on $\overline U$
$$\frac{\partial h}{\partial t}=\Delta h$$
with Dirichlet boundary conditions $h(x,t)=0$ for all $x\in \partial
\overline U$ and all $t\ge 0$ and with the given initial conditions.

Consider the tensor bundle ${\mathcal V}\oplus \Ar$ over $M$. We
define
$$Z_x=\left\{({\mathcal T},h)\in V_x\oplus
\Ar\bigl|\bigr. s_x({\mathcal T})\ge h\ge 0\right\}.$$ The union
over all $x\in M$ of the $Z_x$ defines a closed convex subset
${\mathcal Z}\subset {\mathcal V}\oplus \Ar$ which is invariant
under parallel translation since $s$ is. We consider the family of
sections $({\mathcal T}(x,t),h(x,t)),\ 0\le t\le T$, of ${\mathcal
V}\oplus \Ar$. These evolve by
$$\frac{d\left({\mathcal T}(x,t),h(x,t)\right)}{dt}=\left(\Delta {\mathcal T}(x,t),\Delta
h(x,t)\right)+\widetilde \psi\left({\mathcal T}(x,t),h(x,t)\right)$$
where $\widetilde \psi({\mathcal T},h)=\left(\psi({\mathcal
T}),0\right)$. Clearly, by our hypotheses, the vector field
$\widetilde \psi$ preserves the convex set ${\mathcal Z}$. Applying
the local version of the maximum principle (Theorem~\ref{localMP}),
we conclude that ${\mathcal T}(x,t)\ge h(x,t)$ for all $(x,t)\in
\overline U\times [0,T]$.

The result then follows immediately from Theorem~\ref{SMPheat}.
\end{proof}

\subsection{Applications of the strong maximum principle}

We have the following applications of the strong maximum principle.

\begin{thm}\label{firststrongmax}
Let $(U,g(t)),\ 0\le t\le T$, be a $3$-dimensional Ricci flow with
non-negative sectional curvature with $U$ connected but not
necessarily complete and with $T>0$. If $R(p,T)=0$ for some $p\in
U$, then $(U,g(t))$ is flat for every $t\in [0,T]$.
\end{thm}

\begin{proof}
We suppose that there is $p\in U$ with $R(p,T)=0$. Since all the
metrics in the flow are of non-negative sectional curvature, if the
flow does not consist entirely of flat manifolds then there is
$(q,t)\in U\times [0,T]$ with $R(q,t)>0$. Clearly, by continuity, we
can assume $t<T$. By restricting to the time interval $[t,T]$ and
shifting by $-t$ we can arrange that $t=0$. Let $V$ be a compact,
connected smooth submanifold with boundary whose interior contains
$q$ and $p$. Let $h(y,0)$ be a smooth non-negative function with
support in $V$, positive at $q$, such that $R(y,0)\ge h(y,0)$ for
all $y\in V$. Let $h(y,t)$ be the solution to the heat equation on
$V\times [0,T]$ that vanishes on $\partial V$. Of course, $h(y,T)>0$
for all $y\in {\rm int}(V)$. Also, from Equation~(\ref{Revol}) we
have
$$\frac{\partial}{\partial t}(R-h)=\triangle (R-h)+2|{\rm Ric}|^2,$$
so that $(R-h)(y,0)\ge 0$ on $(V\times \{0\})\cup (\partial V\times [0,T])$. It
follows from the  maximum principle that $(R-h)\ge 0$ on all of $V\times
[0,T]$. In particular, $R(p,T)\ge h(p,T)>0$. This is a contradiction,
establishing the theorem.
\end{proof}

\begin{cor}\label{localprod}
Fix $T>0$. Suppose that $(U,g(t)),\ 0\le t\le T$, is a  Ricci flow
such that for each $t$, the Riemannian manifold  $(U,g(t))$ is a
(not necessarily complete) connected, $3$-manifold of non-negative
sectional curvature. Suppose that $(U,g(0))$ is not flat and that
for some $p\in M$ the Ricci curvature at $(p,T)$ has a zero
eigenvalue. Then for each $t\in (0,T]$ the Riemannian manifold
$(U,g(t))$ splits locally as a product of a surface of positive
curvature and a line, and under this local splitting the flow is
locally the product of a Ricci flow on the surface and the trivial
flow on the line.
\end{cor}

\begin{proof}
First notice that it follows from Theorem~\ref{firststrongmax} that
because $(U,g(0))$ is not flat, we have $R(y,t)>0$ for every
$(y,t)\in U\times (0,T]$.

 We consider the function $s$ on
${\rm Sym}^2(\wedge^2T^*_yU)$ that associates to each endomorphism
the sum of the smallest $2$ eigenvalues. Then $s_y$ is the minimum
of the traces on $2$-dimensional subsets  in $\wedge^2T_yU$. Thus,
$s$ is a convex function, and the subset ${\mathcal
S}=s^{-1}([0,\infty))$ is a convex subset. Clearly, this subset is
invariant under parallel translation. By the computations in the
proof of Corollary~\ref{Ric>0} it is invariant under the vector
field $\psi({\mathcal T})$. The hypothesis of the corollary tells us
that $s(p,T)=0$. Suppose that $s(q,t)>0$ for some $(q,t)\in U\times
[0,T]$. Of course, by continuity we can take $t< T$. Shift the time
parameter so that $t=0$, and fix a compact connected,
codimension-$0$ submanifold $V$ containing $p,q$ in its interior.
Then by Theorem~\ref{SMPcurv} $s(y,T)>0$ for all $y\in {\rm int}(V)$
and in particular $s(p,T)>0$. This is a contradiction, and we
conclude that $s(q,t)=0$ for all $(q,t)\in U\times [0,T]$.

Since we have already established that each $R(y,t)>0$ for all $(y,t)\in
U\times (0,T]$, so that ${\rm Rm}(y,t)$ is not identically zero, this means
that for all $y\in U$ and all $t\in (0,T]$ that the null space of the  operator
${\rm Rm}(y,t)$ is a $2$-dimensional subspace of $\wedge^2T_yU$. This
$2$-dimensional subspace is dual to a line in $T_xM$. Thus, we have a
one-dimensional distribution (a line bundle in the tangent bundle) ${\mathcal
D}$ in $U\times (0,T]$ with the property that the sectional curvature ${\rm
Rm}(y,t)$ vanishes on any $2$-plane containing the line ${\mathcal D}(y,t)$.
The fact that the sectional curvature of $g(t)$ vanishes on all two-planes in
$T_yM$ containing ${\mathcal D}(y,t)$ means that its eigenvalues are
$\{\lambda,0,0\}$ where $\lambda>0$ is the sectional curvature of the
$g(t)$-orthogonal $2$-plane to ${\mathcal D}(y,t)$. Hence ${\mathcal
R}(V(y,t),\cdot,\cdot,\cdot)=0$.

Locally in space and time, there is a unique (up to sign) vector field $V(y,t)$
that generates ${\mathcal D}$ and satisfies $|V(y,t)|^2_{g(t)}=1$. We wish to
show that this local vector field is invariant under parallel translation and
time translation; cf. Lemma 8.2 in \cite{Hamilton4MPCO}. Fix a point $x\in M$,
a direction $X$ at $x$, and a time $t$. Let $\tilde V(y,t)$ be a parallel
extension of $V(x,t)$ along a curve $C$ passing through $x$ in the
$X$-direction, and let $\tilde W(y,t)$ be an arbitrary parallel vector field
along $C$. Since the sectional curvature is non-negative, we have ${\mathcal
R}(\tilde V,\tilde W,\tilde V,\tilde W)(y)\ge 0$ for all $y\in C$; furthermore,
this expression vanishes at $x$. Hence, its first variation vanishes at $x$.
That is to say
$$\nabla \left({\mathcal R}(\tilde V,\tilde W,\tilde V,\tilde W)\right)(x,t)=(\nabla{\mathcal R})(\tilde
V,\tilde W,\tilde V,\tilde W)$$ vanishes at $(x,t)$. Since this is true for all
$\tilde W$, it follows that the null space of the quadratic form
$\nabla{\mathcal R}(x,t)$ contains the null space of ${\mathcal R}(x,t)$, and
thus $$(\nabla {\mathcal R})(V(x,t),\cdot,\cdot,\cdot)=0.$$
 Now let us consider three parallel
vector fields $\tilde W_1,\tilde W_2,$ and $\tilde W_3$ along $C$.
We compute $0=\nabla_X\left({\mathcal R}( V(y,t),\tilde
W_1(y,t),\tilde W_2(y,t),\tilde W_3(y,t))\right)$. (Notice that
while the $\tilde W_i$ are parallel along $C$, $V(y,t)$ is defined
to be the vector field spanning ${\mathcal D}(y,t)$ rather than a
parallel extension of $V(x,t)$.) Given the above result we find that
$$0=2{\mathcal R}(\nabla_X V(x,t),\tilde W_1(x,t),\tilde W_2(x,t),
\tilde W_3(x,t)).$$ Since this is true for all triples of vector
fields $\tilde W_i(x,t)$, it follows that $\nabla_XV(x,t)$ is a real
multiple of $V(x,t)$. But since $|V(y,t)|^2_{g(t)}=1$, we see that
$\nabla_XV(x,t)$ is orthogonal to $V(x,t)$. We conclude that
$\nabla_XV(x,t)=0$. Since $x$ and $X$ are general, this shows that
the local vector field $V(x,t)$ is invariant under the parallel
translation associated to the metric $g(t)$.

It follows that locally $(M,g(t))$ is a Riemannian product of a surface of
positive curvature with a line. Under this product decomposition, the curvature
is the pullback of the curvature of the surface. Hence, by
Equation~(\ref{Rmevol}), under Ricci flow on the $3$-manifold, the time
derivative of the curvature at time $t$ also decomposes as the pullback of the
time derivative of the curvature of the surface under Ricci flow on the
surface. In particular, $(\partial {\mathcal R}/\partial
t)(V,\cdot,\cdot,\cdot)=0$. It now follows easily that $\partial
V(x,t)/\partial t=0$.

This completes the proof that the unit vector field in the direction ${\mathcal
D}(x,t)$ is invariant under parallel translation and under time translation.
Thus, there is a local Riemannian splitting of the $3$-manifold into a surface
and a line, and this splitting is invariant under the Ricci flow. This
completes the proof of the corollary.
\end{proof}

In the complete case, this local product decomposition globalizes in
some cover; see Lemma 9.1 in \cite{Hamilton4MPCO}.

\begin{cor}\label{nullspace}
Suppose that $(M,g(t)),\ 0\le t\le T$, is a Ricci flow of complete, connected
Riemannian $3$-manifolds with ${\rm Rm}(x,t)\ge 0$ for all $(x,t)$ and with
$T>0$. Suppose that $(M,g(0))$ is not flat and that for some $x\in M$  the
endomorphism ${\rm Rm}(x,T)$ has a zero eigenvalue. Then $M$ has a cover
$\widetilde M$ such that, denoting the induced family of metrics on this cover
by $\tilde g(t)$, we have that $(\widetilde M,\tilde g(t))$ splits as a product
$$(N,h(t))\times (\Ar,ds^2)$$
where $(N,h(t))$ is a surface of positive curvature for all $0<t\le
T$. The Ricci flow is a product of the Ricci flow $(N,h(t)),\ 0\le
t\le T$, with the trivial flow on $\Ar$.
\end{cor}

\begin{rem}
Notice that there are only four possibilities for the cover required
by the corollary. It can be trivial, or a normal $\Zee$-cover or it
can be a two-sheeted cover or a normal infinite dihedral group
cover. In the first two cases, there is a unit vector field on $M$
parallel under $g(t)$ for all $t$ spanning the null direction of
${\rm Ric}$. In the last two cases, there is no such vector field,
only a non-orientable line field.
\end{rem}

 Let $(N,g)$ be a Riemannian manifold. Recall from Definition~\ref{conedefn}
  that the open cone
on $(N,g)$ is the space $N\times (0,\infty)$ with the Riemannian
metric $\tilde g(x,s)=s^2g(x)+ds^2$. An extremely important result
for us is that open pieces in non-flat cones cannot arise as the
result of Ricci flow with non-negative curvature.

\begin{prop}\label{nocones}
Suppose that $(U,g(t)),\ 0\le t\le T$, is a $3$-dimensional Ricci
flow with non-negative sectional curvature, with $U$ being connected
but not necessarily complete and $T>0$. Suppose that $(U,g(T))$ is
isometric to a non-empty open subset of a cone over a Riemannian
manifold. Then $(U,g(t))$ is flat for every $t\in [0,T]$.
\end{prop}

\begin{proof}
If $(U,g(T))$ is flat, then by Theorem~\ref{firststrongmax}
 for every $t\in [0,T]$ the
Riemannian manifold $(U,g(t))$ is flat.

We must rule out the possibility that $(U,g(T))$ is  non-flat. Suppose that
$(U,g(T))$ is an open subset in a non-flat cone. According to
Proposition~\ref{conecurv},  for each $x\in U$ the Riemann curvature tensor of
$(U,g(T))$ at $x$ has a $2$-dimensional null space in $\wedge^2T_xU$. Since we
are assuming that $(U,g(T))$ is not flat, the third eigenvalue of the Riemann
curvature tensor is not identically zero. Restricting to a smaller open subset
if necessary, we can assume that the third eigenvalue is never zero. By the
computations in Proposition~\ref{conecurv} the non-zero eigenvalue is not
constant, and in fact it scales by $s^{-2}$ in the terminology of that
proposition, as we move along the cone lines. Of course, the $2$-dimensional
null-space for the Riemann curvature tensor at each point is equivalent to a
line field in the tangent bundle of the manifold. Clearly, that line field is
the line field along the cone lines. Corollary~\ref{localprod} says that since
the Riemann curvature of $(U,g(T))$ has a $2$-dimensional null-space in
$\wedge^2T_xU$ at every point $x\in U$, the Riemannian manifold $(U,g(T))$
locally splits as a Riemannian product of a line with a surface of positive
curvature, and the $2$-dimensional null-space for the Riemannian curvature
tensor is equivalent to the line field in the direction of the second factor.
Along these lines the non-zero eigenvalue of the curvature is constant. This is
a contradiction and establishes the result.
\end{proof}

Lastly, we have Hamilton's result (Theorem 15.1 in
\cite{Hamilton3MPRC}) that compact $3$-manifolds of non-negative
Ricci curvature become round under Ricci flow:

\begin{thm}\label{flowtoround}
Suppose that $(M,g(t)),\ 0\le t<T$, is a Ricci flow with $M$ being a
compact $3$-dimensional manifold. If ${\rm Ric}(x,0)\ge 0$ for all
$x\in M$, then either ${\rm Ric}(x,t)>0$ for all $(x,t)\in M\times
(0,T)$ or ${\rm Ric}(x,t)=0$ for all $(x,t)\in M\times[0,T)$.
Suppose that ${\rm Ric}(x,t)>0$ for some $(x,t)$ and that the flow
is maximal in the sense that there is no $T'>T$ and an extension of
the given flow to a  flow defined on the time interval $[0,T')$. For
each $(x,t)$, let $\lambda(x,t)$, resp. $\nu(x,t)$, denote the
largest, resp. smallest, eigenvalue of ${\rm Rm}(x,t)$ on
$\wedge^2T_xM$. Then as $t$ tends to $T$ the Riemannian manifolds
$(M,g(t))$ are becoming round in the sense that
$${\rm lim}_{t\rightarrow T}\frac{{\rm max}_{x\in M}\lambda(x,t)}{{\rm
min}_{x\in M}\nu(x,t)}=1.$$ Furthermore, for any $x\in M$ the
largest eigenvalue $\lambda(x,t)$ tends to $\infty$ as $t$ tends to
$T$, and rescaling $(M,g(t))$ by $\lambda(x,t)$ produces a family of
Riemannian manifolds converging smoothly as $t$ goes to $T$ to a
compact round manifold. In particular, the underlying smooth
manifold supports a Riemannian metric of constant positive curvature
so that the manifold is diffeomorphic to a $3$-dimensional spherical
space-form.
\end{thm}

Hamilton's proof in \cite{Hamilton3MPRC} uses the maximum principle
and Shi's derivative estimates.

\subsection{Solitons of positive curvature}

One nice application of this pinching result is the following
theorem.

\begin{thm}\label{sphsf}
Let $(M,g)$ be a compact $3$-dimensional soliton of positive Ricci
curvature. Then $(M,g)$ is round. In particular, $(M,g)$ is the
quotient of $S^3$ with a round metric by a finite subgroup of $O(4)$
acting freely; that is to say, $M$ is a $3$-dimensional spherical
space-form.
\end{thm}

\begin{proof}
Let $(M,g(t)),\ 0\le t<T$, be the maximal Ricci flow with initial
manifold $(M,g)$. Since ${\rm Ric}(x,0)>0$ for all $x\in M$, it
follows from Theorem~\ref{flowtoround} that $T<\infty$ and that as
$t$ tends to $T$ the metrics $g(t)$ converge smoothly to a round
metric. Since all the manifolds $(M,g(t))$ are isometric up to
diffeomorphism and a constant conformal factor, this implies that
all the $g(t)$ are of constant positive curvature.

The last statement is a standard consequence of the fact that the
manifold has constant positive curvature.
\end{proof}

\begin{rem}
After we give a stronger pinching result in the next section, we
shall improve this result, replacing the positive Ricci curvature
assumption by the {\em a priori} weaker assumption that the soliton
is a shrinking soliton.
\end{rem}

\section{Pinching toward positive
curvature}\index{curvature!pinched toward positive|(}

As the last application of the maximum principle for tensors we give
a theorem due to R. Hamilton (Theorem 4.1 in \cite{HamiltonNSRF3M})
and T. Ivey \cite{Ivey} which shows that, in dimension three, as the
scalar curvature gets large, the sectional curvatures pinch toward
the positive. Of course, if the sectional curvatures are
non-negative, then the results in the previous section apply. Here,
we are considering the case when the sectional curvature is not
everywhere positive. The pinching result says roughly the following:
At points where the Riemann curvature tensor has a negative
eigenvalue, the smallest (thus negative) eigenvalue of the Riemann
curvature tensor divided by the largest eigenvalue limits to zero as
the scalar curvature grows. This result is central in the analysis
of singularity development in finite time for a $3$-dimensional
Ricci flow.

\begin{thm}\label{pinch}{\bf (Pinching toward positive curvature)}
Let $(M,g(t)),\ 0\le t<T$, be a Ricci flow with $M$ a compact
$3$-manifold. Assume at for every $x\in M$, the eigenvalues,
$\lambda(x,0)\geq\mu(x,0)\geq\nu(x,0)$, of ${\rm Rm}(x,t)$ are all
at least $-1$. Set $X(x,t)={\rm max}(-\nu(x,t),0)$. Then we have:
\begin{enumerate}
\item $R(x,t)\ge \frac{-6}{4t+1}$, and
\item for all $(x,t)$ for which $0<X(x,t)$
$$R(x,t)\geq 2X(x,t)\left({\rm log}X(x,t)+{\rm
log}(1+t)-3\right).$$
\end{enumerate}
\end{thm}

For any fixed $t$, the limit as $X$ goes to $0$ from above of
$X({\rm log}(X)+{\rm log}(1+t)-3)$ is zero, so that it is natural to
interpret this expression to be zero when $X=0$. Of course, when
$X(x,t)=0$ all the eigenvalues of ${\rm Rm}(x,t)$ are non-negative
so that $R(x,t)\ge 0$ as well. Thus, with this interpretation of the
expression in Part 2 of the theorem, it remains valid even when
$X(x,t)=0$.

\begin{rem}
This theorem tells us, among other things, that as the scalar
curvature goes to infinity  then  absolute values of all the
negative eigenvalues (if any) of ${\rm Rm}$ are  arbitrarily small
with respect to the scalar curvature.
\end{rem}\index{curvature!pinched toward positive|)}

We proof we give below follows Hamilton's original proof in
\cite{HamiltonNSRF3M} very closely.

\begin{proof}
First note that by Proposition~\ref{scalarevol}, if $R_{\rm
min}(0)\ge 0$, then the same is true for $R_{\rm min}(t)$ for every
$t>0$ and thus the first inequality stated in the theorem is clearly
true. If $R_{\rm min}(0)<0$, the first inequality stated in the
theorem follows easily from the last inequality in
Proposition~\ref{scalarevol}.

We turn now to the second inequality in the statement of the
theorem. Consider the tensor bundle ${\mathcal V}={\rm
Sym}^2(\wedge^2T^*M)$. Then the curvature operator written in the
evolving frame, ${\mathcal T}(x,t)$, is a one-parameter family of
smooth sections of this bundle, evolving by
$$\frac{\partial{\mathcal T}}{\partial t}=
\Delta{\mathcal T}+\psi({\mathcal T}).$$ We consider two subsets of ${\mathcal
V}$. There are two solutions to $x({\rm log}(x)+({\rm log}(1+t)-3)=-3/(1+t)$.
One is $x=1/(1+t)$; let $\xi(t)>1/(1+t)$ be the other. We set $S({\mathcal
T})={\rm tr}({\mathcal T})$, so that $R=2S$, and we set $X({\mathcal T})={\rm
max}(-\nu({\mathcal T}),0)$. Define
\begin{eqnarray*}
 {\mathcal Z}_1(t) & = & \{{\mathcal T}\in {\mathcal V}\bigl|\bigr.
 S({\mathcal T})\geq-\frac{3}{(1+t)}\} \\
{\mathcal Z}_2(t) & = & \{ {\mathcal T}\in {\mathcal
V}\bigl|\bigr.S({\mathcal T})\geq f_t(X({\mathcal T})),\ \ \
\text{if }\ \ X( {\mathcal T})\geq\xi(t)\},
\end{eqnarray*}
where $f_t(x)=x({\rm log}x+{\rm log}(1+t)-3)$. Then we define
$${\mathcal Z}(t)={\mathcal Z}_1(t)\cap {\mathcal Z}_2(t).$$

\begin{claim}
For each $x\in M$ and each $t\ge 0$, the fiber $Z(x,t)$ of
${\mathcal Z}(t)$ over $x$  is a convex subset of ${\rm
Sym}^2(\wedge^2T^*M)$.
\end{claim}

\begin{proof}
First consider the function $f_t(x)=x({\rm log}(x)+{\rm
log}(1+t)-3)$ on the interval $[\xi(t),\infty)$. Direct computation
shows that $f'(x)>0$ and $f''(x)>0$ on this interval. Hence, for
every $t\ge 0$ the region ${\mathcal C}(t)$ in the $S$-$X$ plane
 defined by $S\ge -3/(1+t)$
and $S\ge f_t(X)$ when $X\ge \xi(t)$ is convex and has the property
that if $(S,X)\in {\mathcal C}(t)$ then so is $(S,X')$ for all
$X'\le X$. (See {\sc Fig.}~\ref{fig:SX}). By definition an element
${\mathcal T}\in {\mathcal V}$ is contained in ${\mathcal Z}(t)$ if
and only if $(S({\mathcal T}),X({\mathcal T})\in {\mathcal C}(t)$.
Now fix $t\ge 0$ and suppose that ${\mathcal T}_1$ and ${\mathcal
T}_2$ are elements of ${\rm Sym}^2(\wedge^2T^*M_x)$ such that
setting $S_i={\rm tr}({\mathcal T}_i)$ and $X_i=X({\mathcal T}_i)$
we have $(S_i,X_i)\in {\mathcal C}(t)$ for $i=1,2$. Then we consider
${\mathcal T}=s{\mathcal T}_1+(1-s){\mathcal T}_2$ for some $s\in
[0,1]$. Let $S={\rm tr}({\mathcal T})$ and $X=X({\mathcal T})$.
Since ${\mathcal C}(t)$ is convex, we know that
$(sS_1+(1-s)S_2,sX_1+(1-s)X_2)\in{\mathcal C}(t)$, so that
${\mathcal T}\in {\mathcal Z}(t)$.
 Clearly,
$S=sS_1+(1-s)S_2$, so that we conclude that $(S,(sX_1+(1-s)X_2))\in{\mathcal
C}(t)$.  But since $\nu$ is a convex function, $X$ is a concave function, i.e.,
$X\le sX_1+(1-s)X_2$. Hence $(S,X)\in {\mathcal C}(t)$.
\end{proof}

\begin{figure}[ht]
  \centerline{\epsfbox{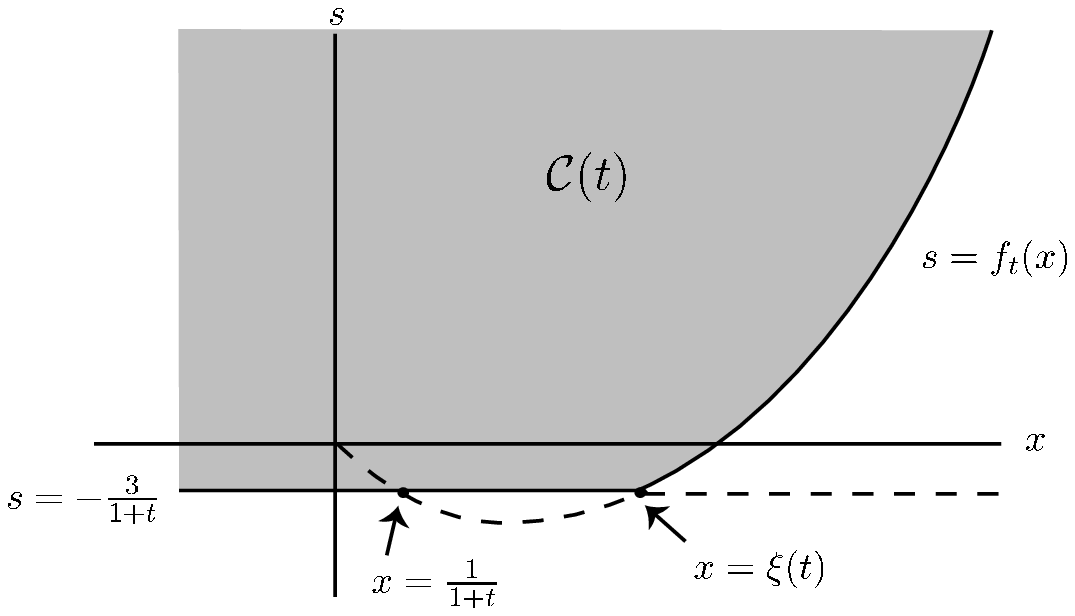}}
  \caption{Curvature convex set}\label{fig:SX}
\end{figure}

\begin{claim}
${\mathcal T}(x,0)\in Z(x,0)$ for all $x\in M$.
\end{claim}

\begin{proof}
 Note that by the
hypothesis of the theorem we have
$$\nu(x,0)+\mu(x,0)+\lambda(x,0)\geq -3$$ so $(S(x,0),X(x,0))\in
{\mathcal C}(0)$ for all $x\in M$. On the other hand, if $0<
X(x,0)$, then since $X(x,0)\le 1$ we have $S(x,0)\geq-3X(x,0)\geq
X({\rm log}X-3)$. This completes the proof that ${\mathcal
T}(x,0)\in {\mathcal C}(0)$ for all $x\in M$.
\end{proof}

\begin{claim}
The vector field $\psi({\mathcal T})={\mathcal T}^2+{\mathcal
T}^\sharp$ preserves the family ${\mathcal Z}(t)$ of convex sets.
\end{claim}

\begin{proof}
Fix $x\in M$ and suppose that we have an integral curve $\gamma(t),\
t_0\le t\le T$, for $\psi$ with $\gamma(t_0)\in Z(x,t_0)$. We wish
to show that $\gamma(t)\in Z(x,t)$ for all $t\in [t_0,T]$. The
function $S(t)=S(\gamma(t))$ satisfies
$$\frac{dS}{dt}=\lambda^2+\mu^2+\nu^2+\lambda\mu+\lambda\nu+\mu\nu=\frac{1}{2}\left((\lambda+\mu)^2+
(\lambda+\nu)^2+(\mu+\nu)^2\right).$$ By Cauchy-Schwarz we have
$$(\lambda+\mu)^2+
(\lambda+\nu)^2+(\mu+\nu)^2)\ge \frac{4S^2}{3}\ge \frac{2S^2}{3}.$$ Since
$\gamma(t_0)\in Z(x,t_0)$ we have $S(t_0)\ge -3/(1+t_0)$. It then follows that
\begin{equation}\label{Sineq}
S(t)\ge -3/(1+t)\ \ \ {\rm for\ all}\ \ \  t\ge t_0.
\end{equation}

Now let us consider the evolution of $X(t)=X(\gamma(t))$. Assume
that we are at a point $t$ for which $X(t)>0$. For this computation
we set $Y=-\mu$.
\begin{align*}
\frac{dX}{dt}&=-\frac{d\nu}{dt}=-\nu^2-\mu\lambda =-X^2+Y\lambda,\\
\frac{dS}{dt}&=\frac{d(\nu+\mu+\lambda)}{dt}
=\nu^2+\mu^2+\lambda^2+\mu\lambda+\nu\lambda+\nu\mu\\
&=X^2+Y^2+\lambda^2+XY-\lambda(X+Y).
\end{align*}
Putting this together yields
\begin{equation}\label{dX}
X\frac{dS}{dt}-(S+X)\frac{dX}{dt} =X^3+I, \end{equation}
 where
$I=XY^2+\lambda Y(Y-X)+\lambda^2(X-Y)$.
\begin{claim}  $I\geq 0$.
\end{claim}

\begin{proof}

 First we consider the case when
 $Y\leq 0$. This means that  $\mu\geq 0$ and
 hence that $\lambda\geq0$. Since by definition $X\ge 0$, we have $X\ge  Y$.
This immediately gives  $I\geq 0$. Now let us consider the case when
$Y>0$ which means that $\nu\le \mu<0$. In this case, we have
\[ I=Y^3+(X-Y)(\lambda^2-\lambda Y+Y^2)>0 \]
since $X\geq Y$ and $\lambda^2-\lambda
Y+Y^2=(\lambda-\frac{Y}{2})^2+\frac{3Y^2}{4}>0$.
\end{proof}

The above claim and Equation~(\ref{dX}) immediately imply that
\begin{equation}\label{dX1}
 X\frac{dS}{dt}-(S+X)\frac{dX}{dt}\geq
X^3. \end{equation}
 Set $W=\frac{S}{X}-{\rm log}\,X$, then rewriting
 Equation~(\ref{dX1})
 in terms of $W$ gives
\begin{equation}\label{dW}
\frac{dW}{dt}\geq X.
\end{equation}

Now suppose that $\gamma(t)\not\in Z(x,t)$ for some $t\in [t_0,T]$.
Let $t_1<T$ be maximal subject to the condition that $\gamma(t)\in
Z(x,t)$ for all $t_0\le t\le t_1$. Of course, $\gamma(t_1)\in
\partial Z(x,t_1)$ which implies that $(S(t_1),X(t_1))\in \partial
{\mathcal C}(t_1)$. There are two possibilities: either
$S(t_1)=-3/(1+t_1)$ and $X(t_1)<\xi(t_1)$ or $X(t_1)\ge
\xi(t_1)>1/(1+t_1)$  and $S(t_1)=f_{t_1}(X(t_1))$. But
Equation~(\ref{Sineq}) implies that $S(t)\ge -3/(1+t)$ for all $t$.
Hence, if the first case holds then $\gamma(t)\in Z(x,t)$ for $t$ in
some interval $[t_0,t_1']$ with $t_1'>t_1$. This contradicts the
maximality of $t_1$. Thus, it must be the case that $X(t_1)\ge
\xi(t_1)$. But then $X(t)> \frac{1}{1+t}$ for all $t$ sufficiently
close to $t_1$. Hence, by Equation~(\ref{dW}) we have
\[ \frac{dW}{dt}(t)\geq X(t)>\frac{1}{1+t}, \]
for all $t$ sufficiently close to $t_1$. Also, since
$S(t_1)=f_{t_1}(X(t_1))$, we have $W(t_1)=({\rm log}(1+t_1)-3)$. It
follows immediately that $W(t)\ge ({\rm log}(1+t)-3)$ for all
$t>t_1$ sufficiently close to $t_1$. This proves that $S(t)\ge
f_t(X(t))$ for all $t\ge t_1$ sufficiently close to $t_1$, again
contradicting the maximality of $t_1$.

This contradiction proves that $\psi$ preserves the family
${\mathcal Z}(t)$.
\end{proof}

By Theorem~\ref{varyingZ}, the previous three claims imply that ${\mathcal
T}(x,t)\in {\mathcal Z}(t)$ for all $x\in M$ and all $t\in [0,T)$. That is to
say, $S(x,t)\ge -3/(1+t)$ and $S(x,t)\ge f_t(X(x,t))$ whenever $X(x,t)\ge
\xi(t)$. For $X\in [1/(1+t),\xi(t)]$ we have $f_t(X)\le -3/(1+t)$, and thus in
fact $S(x,t)\ge f_t(X(x,t))$ as long as $X(x,t)\ge 1/(1+t)$. On the other hand,
if $0< X(x,t)\le 1/(1+t)$ then $f_t(X(x,t))< -3X(x,t)\le S(x,t)$. On the other
hand, since $X(x,t)$ is the negative of the smallest eigenvalue of ${\mathcal
T}(x,t)$ and $S(x,t)$ is the trace of this matrix, we have $S(x,t)\ge
-3X(x,t)$. Thus, $S(x,t)\ge f_t(X(x,t))$ in this case as well. This completes
the proof of Theorem~\ref{pinch}.
\end{proof}

Actually, the proof establishes a stronger result which we shall
need.

\begin{thm}\label{pincha}
Fix $a\ge 0$. Let $(M,g(t)),\ a\le t<T$, be a Ricci flow with $M$ a
compact $3$-manifold. Suppose the eigenvalues of ${\rm Rm}(x,t)$ are
$\lambda(x,t)\ge \mu(x,t)\ge \nu(x,t)$ and set $X(x,t)={\rm
max}(-\nu(x,t),0)$. Assume that for every $x\in M$ we have
$$R(x,a) \ge  \frac{-6}{4a+1}$$
and $$ R(x,a) \geq 2X(x,a)\left({\rm log}X(x,a)+{\rm
log}(1+a)-3\right),$$ where the second inequality holds whenever
$X(x,a)>0$. Then for all $a\le t<T$ we have:
\begin{eqnarray}\label{Rlower}
R(x,t) & \ge & \frac{-6}{4t+1} \\\label{RXineq} R(x,t) & \geq &
2X(x,t)\left({\rm log}X(x,t)+{\rm log}(1+t)-3\right),
\end{eqnarray}
whenever  $X(x,t)> 0$.
\end{thm}

Once again it is natural to interpret the right-hand side of the
inequalities relating $R$ and $X$  to be zero when $X(x,t)=0$. With
this convention the result remains true even when $X(x,t)=0$.

\begin{cor}
Fix $a\ge 0$. Suppose that $(M,g(t)),\ a\le t<T$, is a Ricci flow with $M$ a
compact $3$-manifold, and suppose that the two hypotheses of the previous
theorem hold. Then there is a continuous function $\phi$ such that for all
$R_0<\infty$, if $R(x,t)\le R_0$ then $|{\rm Rm}(x,t)|\le \phi(R_0)$.
\end{cor}

\begin{proof}
Fix $R_0\ge e^4$ sufficiently large, and suppose that $R(x,t)\le
R_0$. If $X(x,t)=0$, then $|{\rm Rm}(x,t)|\le R(x,t)/2$. If
$X(x,t)>0$, then by Theorem~\ref{pincha} it is bounded by $R_0$.
Thus, $\lambda(x,t)\le 3R_0$. Thus, we have an upper bound on
$\lambda(x,t)$ and a lower bound on $\nu(x,t)$ in terms of $R_0$.
\end{proof}

This theorem leads to a definition.

\begin{defn}\label{pinchd}
Let $({\mathcal M},G)$ be a generalized Ricci flow whose domain of definition
is contained in $[0,\infty)$. Then we say that $({\mathcal M},G)$ has {\em
curvature pinched toward positive}\index{curvature!pinched toward positive} if
for every $x\in {\mathcal M}$ the following two conditions hold:
\begin{enumerate}
\item $$R(x)\ge \frac{-6}{4{\bf t}(x)+1}$$
\item $$R(x)\geq 2X(x)\left({\rm log}X(x)+{\rm
log}(1+{\bf t}(x))-3\right),
$$ whenever
$0<X(x)$
\end{enumerate}
where, as in the statement of Theorem~\ref{pinch}, $X(x)$ is the
maximum of zero and the negative of the smallest eigenvalue of ${\rm
Rm}(x)$.
\end{defn}

The content of Theorem~\ref{pincha} is that if $(M,g(t)),\ 0\le a\le
t<T$, is a Ricci flow with $M$ a compact $3$-manifold and if the
curvature of $(M,g(a))$ is pinched toward positive, then the same is
true for the entire flow.

\subsection{Application of the pinching result}

As an application of this pinching toward positive curvature result
we establish a strengthening of Theorem~\ref{sphsf}.

\begin{thm}
Let $(M,g)$ be a compact $3$-dimensional shrinking soliton, i.e.,
there is a Ricci flow $(M,g(t)),\ 0\le t<T$, so that for each $t\in
[0,T)$ there is a constant $c(t)$ with  ${\rm lim}_{t\rightarrow
T}c(t)=0$ and with the property that there is an isometry from
$(M,g(t))$ to $(M,c(t)g)$. Then $(M,g)$ is round.
\end{thm}

\begin{proof}
By rescaling we can assume that for all $x\in M$ all the eigenvalues of ${\rm
Rm}(x,0)$ have absolute value $\le 1$.  This implies that $(M,g(0))$ satisfies
the hypothesis of Theorem~\ref{pinch}. Our first goal is to show that ${\rm
Rm}(x,0)\ge 0$ for all $x\in M$. Suppose that this is not true; then there is a
point $x$ with $X(x,0)>0$. Consider $A=X(x,0)/R(x,0)$. For each $t<T$ let
$x_t\in M$ be the image of $x$ under the isometry from $(M,g(0))$ to
$(M,c(t)g(t))$.
 Then $X(x_t,t)=c^{-1}(t)X(x,0)$ and  $X(x_t,t)/R(x_t,t)=A$. Since
$c(t)$ tends to $0$ as $t$ approaches $T$, this contradicts
Theorem~\ref{pinch}. Now, according to Theorem~\ref{rm>0} either all
$(M,g(t))$ are flat or ${\rm Rm}(x,t)>0$ for all $(x,t)\in M\times
(0,T)$. But if the $(M,g(t))$ are all flat, then the flow is trivial
and hence the diameters of the $(M,g(t))$ do not go to zero as $t$
approaches $T$, contradicting the hypothesis. Hence, ${\rm
Rm}(x,t)>0$ for all $(x,t)\in M\times (0,T)$. According to
Theorem~\ref{flowtoround} this means that as the singularity
develops the metrics are converging to round. By the shrinking
soliton hypothesis, this implies that all the metrics $(M,g(t)),\
0<t<T$, are in fact round. Of course, it then follows that $(M,g)$
is round.
\end{proof}

The following more general result was first given by T. Ivey
\cite{Ivey}.

\begin{thm}\label{ivey}
Any 3-dimensional compact Ricci soliton\index{Ricci flow!soliton}
$g_0$ is Einstein\index{Einstein manifold}.
\end{thm}

Since we do not need this result, we do not include a proof.

\subsection{The Harnack inequality}

The last consequence of the maximum principle that we need is
Hamilton's version of the  Harnack inequality for Ricci
flows\index{Ricci flow!Harnack inequality for|ii}, see Theorem 1.1
and Corollary 1.2 of \cite{Hamiltonharnack}.

\begin{thm}\label{Harnack}
Suppose that $(M,g(t))$ is a Ricci flow defined for $(T_0,T_1)$ with
$(M,g(t))$ a complete manifold of non-negative curvature operator
with bounded curvature  for each $t\in (T_0,T_1)$. Then for any
time-dependent vector field $\chi(x,t)$ on $M$ we have:
$$\frac{\partial R(x,t)}{\partial t}+\frac{R(x,t)}{t-T_0}+2\langle
\chi(x,t),\nabla R(x,t)\rangle +2{\rm
Ric}(x,t)(\chi(x,t),\chi(x,t))\ge 0.$$ In particular, we have
$$\frac{\partial R(x,t)}{\partial t}+\frac{R(x,t)}{t-T_0}\ge 0.$$
\end{thm}

\begin{rem}
Notice that the second result follows from the first by taking
$\chi=0$.
\end{rem}

\begin{cor}\label{posderiv}
If $(M,g(t))$ is a Ricci flow defined for $-\infty<t\le 0$ with
$(M,g(t))$ a complete manifold of bounded, non-negative curvature
operator for each $t$, then
$$\frac{\partial R(x,t)}{\partial t}\ge 0.$$
\end{cor}

\begin{proof}
Apply the above theorem with $\chi(x,t)=0$ for all $(x,t)$ and for a
sequence of $T_0\rightarrow -\infty$.
\end{proof}

The above is the differential form of Hamilton's Harnack inequality.
There is also the integrated version, also due to Hamilton; see
Corollary 1.3 of \cite{Hamiltonharnack}.

\begin{thm}\label{Harnack2}
Suppose that $(M,g(t))$ is a Ricci flow defined for $t_1\le t\le
t_2$ with $(M,g(t))$ a complete manifold of non-negative, bounded
curvature operator for all $t\in [t_1,t_2]$. Let $x_1$ and $x_2$ be
two points of $M$.  Then $${\rm
log}\left(\frac{R(x_2,t_2)}{R(x_1,t_1)}\right)\ge
-\frac{1}{2}\frac{d^2_{t_1}(x_2,x_1)}{(t_2-t_1)}.$$
\end{thm}

\begin{proof}
Apply the differential form of the Harnack inequality to
$\chi=-\nabla({\rm log} R)/2=-\nabla R/2R$, and divide by $R$. The
result is
$$R^{-1}(\partial R/\partial t)-|\nabla ({\rm log} R)|^2+\frac{{\rm Ric}(\nabla({\rm log} R),\nabla({\rm
log} R))}{2R}\ge 0.$$ Since ${\rm Ric}(A,A)/R\le |A|^2$, it follows
that
$$\frac{\partial}{\partial t}({\rm log} R)-\frac{|\nabla({\rm log} R)|^2}{2}\ge 0.$$

Let $d$ be the $g(t_1)$-distance from $x_1$ to $x_2$ and let
 $\gamma\colon [t_1,t_2]\to M$ be a $g(t_1)$-geodesic from $x_1$
 to $x_2$, parameterized at speed $d/(t_2-t_1)$. Then
 let $\mu(t)=(\gamma(t),t)$ be a path in space-time.
 We compute
 \begin{eqnarray*}
 {\rm log}\left(\frac{R(x_2,t_2)}{R(x_1,t_1)}\right) & = & \int_{t_1}^{t_2}
 \frac{d}{dt}{\rm
 log}(R(\mu(t))dt \\
 & = & \int_{t_1}^{t_2}\frac{\frac{\partial R}{\partial t}(\mu(t))}{R(\mu(t))}+\langle\nabla({\rm
 log} R)(\mu(t)),\frac{d\mu}{dt}(\mu(t))dt \\
 & \ge & \int_{t_1}^{t_2}\frac{1}{2}|\nabla ({\rm log} R)(\mu(t))|^2-|\nabla ({\rm
 log} R)(\mu(t))|\cdot \left|\frac{d\gamma}{dt}\right|dt \\
& \ge  & -\frac{1}{2}\int_{t_1}^{t_2}\left|\frac
{d\gamma}{dt}\right|^2dt,\end{eqnarray*} where the last inequality
comes form completing the square.
 Since ${\rm Ric}(x,t)\ge 0$,
$|d\gamma/dt|_{g(t)}\le |d\gamma/dt|_{g(t_1)}$, thus
$${\rm log}\left(\frac{R(x_2,t_2)}{R(x_1,t_1)}\right)\ge
-\frac{1}{2}\int_{t_1}^{t_2}\left|\frac{d\gamma}{dt}\right|^2_{g(t_1)}dt.$$
Since $\gamma$ is a $g(t_1)$-geodesic, this latter integral is
$$-\frac {1}{2}\frac{d_{g(t_1)}^2(x_1,x_2)}{t_2-t_1}.$$
\end{proof}

\chapter{Convergence results for Ricci flow}

The most obvious notion of smooth convergence of Riemannian
manifolds is the $C^\infty$-version of Cheeger-Gromov compactness:
We have a sequence of Riemannian metrics $g_n$ on a fixed smooth
manifold $M$ converging uniformly on compact subsets of $M$ in the
$C^\infty$-topology to a limit metric $g_\infty$. There is also a
version of this compactness for based, complete Riemannian
manifolds. The most common starts with a sequence of based complete
Riemannian manifolds $(M_n,g_n,x_n)$, typically of unbounded
diameter. Then a geometric limit is a based complete
$(M_\infty,g_\infty,x_\infty)$ so that for every $R<\infty$ the
metric balls $B(x_n,R)\subset M_n$ converge uniformly in the
$C^\infty$-topology to the metric ball $B(x_\infty,R)\subset
M_\infty$. This allows the topology to change -- even if all the
$M_n$ are diffeomorphic, $M_\infty$ can have a different topological
type; for example the $M_n$ could all be compact and $M_\infty$
could be non-compact.

But we also need to be able to deal with incomplete limits. In the
case of incomplete limits, the basic idea remains the same, but it
requires some care to give a definition of a geometric limit that
makes it unique up to canonical isometry. One must somehow impose
conditions that imply that the limit eventually fills up most of
each of the manifolds in the sequence.

\section{Geometric convergence of Riemannian manifolds}

Above we referred to filling up `most' of the manifold. The measure of most of
the manifold is in terms of the $\delta$-regular points as defined below.

\begin{defn}\label{deltareg}
Let $(U,g)$ be a Riemannian manifold. Let $\delta>0$ be given. We
say that $p\in U$ is a {\em $\delta$-regular}
point\index{$\delta$-regular point|ii} if for every $r'<\delta$ the
metric ball $B(p,r')$ has compact closure in $U$. Equivalently, $p$
is $\delta$-regular if the exponential mapping at $p$ is defined on
the open ball of radius $\delta$ centered at the origin in $T_pU$,
i.e., if each geodesic ray emanating from $p$ extends to a geodesic
defined on $[0,\delta)$. We denote by ${\rm Reg}_\delta(U,g)$ the
subset of $\delta$-regular points in $(U,g)$. For any $x\in {\rm
Reg}_\delta(U,g)$ we denote by ${\rm Reg}_\delta(U,g,x)$ the
connected component of ${\rm Reg}_\delta(U,g)$ containing $x$.
\end{defn}

Intuitively, the $\delta$-regular points of $(U,g)$ are at distance
at least $\delta$ from the boundary on $U$.

\begin{lem}
 ${\rm Reg}_\delta(U,g)$  is a closed subset of $U$.
\end{lem}

\begin{proof}
Suppose that $p_n$ converges to $p$ as $n$ tends to $\infty$ and
suppose that $p_n\in {\rm Reg}_\delta(U,g)$ for all $n$. Fix
$r'<\delta$ and consider the
 ball $B(p,r')$. For all $n$ sufficiently large, this
ball is contained in $B(p_n,(\delta+r')/2)$, and hence has compact
closure.
\end{proof}

Now we are ready for the basic definition of geometric convergence
of Riemannian manifolds.

\begin{defn}\label{smoothconv} For each $k$
let $(U_k,g_k,x_k)$ be a based, connected Riemannian manifold. A
{\sl geometric limit} of the sequence $\{U_k,g_k,x_k\}_{k=0}^\infty$
is a based, connected Riemannian manifold
$(U_\infty,g_\infty,x_\infty)$
 with the extra data:
\begin{enumerate}
\item[(1)] An increasing sequence $V_k\subset U_\infty$ of connected open subsets
of $U_\infty$ whose union is $U_\infty$ and which satisfy the
following for all $k$:
\begin{enumerate}
\item the closure $\overline V_k$ is compact, \item
$\overline V_k\subset V_{k+1}$,  \item $V_k$ contains $x_\infty$.
\end{enumerate}
\item For each $k\ge 0$ a smooth embedding $\varphi_k\colon (V_k,x_\infty)\to
(U_k,x_k)$ with the properties that: \begin{enumerate} \item ${\rm
lim}_{k\rightarrow \infty}\varphi_k^*g_k=g_\infty$, where the limit
is in the uniform $C^\infty$-topology on compact subsets of
$U_\infty$.
\item For any $\delta>0$ and any $R<\infty$ for all $k$ sufficiently large,
$x_k\in {\rm Reg}_\delta(U_k,g_k)$ and for any $\ell\ge k$ the image
$\varphi_\ell(V_k)$ contains $B(x_\ell,R)\cap{\rm
Reg}_\delta(U_\ell,g_\ell,x_\ell)$.
\end{enumerate}
\end{enumerate}

We also say that the sequence {\em converges geometrically
to}\index{converge!geometrically} $(U_\infty,g_\infty,x_\infty)$ if
there exist $(V_k,\varphi_k)$ as required in the above definition.
We also say that $(U_\infty,g_\infty,x_\infty)$ is the {\em
geometric limit}\index{geometric limit} of the sequence.

More generally, given $(U_\infty,g_\infty,x_\infty)$, a sequence of open
subsets and $\{V_k\}_{k=1}^\infty$ satisfying (1) above, and smooth maps
$\varphi_k\colon V_k\to U_k$ satisfying (2a) above, we say that
$(U_\infty,g_\infty,x_\infty)$ is a {\em partial geometric
limit}\index{geometric limit!partial|ii} of the sequence.
\end{defn}

\begin{rem}
Conditions (1) and (2a) in the definition above also appear in the definition
in the case of complete limits. It is Condition (2b) that is extra in this
incomplete case. It says that once $k$ is sufficiently large then the image
$\varphi_\ell(V_k)$ contains all points satisfying two conditions: they are at
most a given bounded distance from $x_\ell$, and also they are at least a fixed
distance from the boundary of $U_\ell$.
\end{rem}

Notice that if the $(U_k,g_k)$ have uniformly bounded volume by, say, $V$, then
any geometric limit has volume $\le V$.

\begin{lem}
The geometric limit of a sequence $(U_k,g_k,x_k)$ is unique up to
based isometry.
\end{lem}

\begin{proof}
Suppose that we have two geometric limits
$(U_\infty,g_\infty,x_\infty)$ and
$(U'_\infty,g'_\infty,x'_\infty)$. Let $\{V_k,\varphi_k\}$ and
$\{V'_k,\varphi_k'\}$ be the sequences of open subsets and maps as
required by the definition of the limit.

Fix $k$. Since $V_k$ is connected and has compact closure, there are
$R<\infty$ and $\delta>0$ such that $V_k\subset B(x_\infty,R)\cap
{\rm Reg}_\delta(U_\infty,g_\infty,x_\infty)$. Let $x$ be contained
in the closure of $V_k$. Then by the triangle inequality the closed
ball $\overline{B(x,\delta/3)}$ is contained in
$B(x_\infty,R+\delta)\cap{\rm
Reg}_{\delta/2}(U_\infty,g_\infty,x_\infty)$. Since the union of
these closed balls as $x$ ranges over $\overline V_k$ is a compact
set, for all $\ell$ sufficiently large, the restriction of
$\varphi_\ell^*g_\ell$ to the union of these balls is close to the
restriction of $g_\infty$ to the same subset. In particular, for all
$\ell$ sufficiently large and any $x\in V_k$ we see that
$\varphi_\ell\left(B(x,\delta/3)\right)$ contains
$B(\varphi_\ell(x),\delta/4)$. Thus, for all $\ell$ sufficiently
large $\varphi_\ell(V_k)\subset B(x_\ell,R+2\delta)\cap {\rm
Reg}_{\delta/4}(U_\ell,g_\ell,x_\ell)$. This implies that, for given
$k$, for all $\ell$ sufficiently large $\varphi_\ell(V_k)\subset
\varphi'_\ell(V'_\ell)$. Of course, $(\varphi'_\ell)^{-1}\circ
\varphi_\ell(x_\infty)=x_\infty'$. Fix $k$ and pass to a subsequence
of $\ell$, such that as $\ell\rightarrow\infty$, the compositions
$(\varphi'_\ell)^{-1}\circ\left(\varphi_\ell|_{V_k}\right)\colon
V_k\to U'_\infty$ converge to a base-point preserving isometric
embedding of $V_k$ into $U'_\infty$. Clearly, as we pass from $k$ to
$k'>k$ and take a further subsequence of $\ell$ these limiting
isometric embeddings are compatible. Their union is then a
base-point preserving isometric embedding of $U_\infty$ into
$U'_\infty$.

The last thing we need to see is that the embedding of $U_\infty$ into
$U'_\infty$ constructed in the previous paragraph is onto. For each $n$ we have
$\overline V'_n\subset V'_{n+1}$. Since $\overline V'_n$ is compact and
connected, it follows that there are $R<\infty$ and $\delta>0$ (depending on
$n$) such that $\overline V'_n\subset B(x'_\infty,R)\cap {\rm
Reg}_\delta(V_{n+1},g'_\infty,x_\infty')$. Since $V'_{n+1}$ has compact closure
in $U'_\infty$, as $\ell$ tends to $\infty$ the metrics
$(\varphi'_\ell)^*g_\ell$ converge uniformly on $V_{n+1}$ to
$g'_\infty|_{V_{n+1}}$. This means that there are $R'<\infty$ and $\delta'>0$
(depending on $n)$ such that for all $\ell$ sufficiently large,
$\varphi'_\ell(V_n)\subset B(x_k,R')\cap {\rm
Reg}_{\delta'}(U_\ell,g_\ell,x_\ell)$. This implies that for all $k$
sufficiently large and any $\ell\ge k$ the image $\varphi'_\ell(V'_n)$ is
contained in the image of $\varphi_\ell(V_k)$. Hence, for all $k$  sufficiently
large and any $\ell\ge k$ we have
$V'_n\subset(\varphi'_\ell)^{-1}(\varphi_\ell(V_k))$. Hence, the isometric
embedding $U_\infty\to U'_\infty$ constructed above contains $V'_n$. Since this
is true for every $n$, it follows that this isometric embedding is in fact an
isometry $U_\infty\to U'_\infty$.
\end{proof}

Here is the basic existence result.

\begin{thm}\label{basicconv}
Suppose that $\{(U_k,g_k,x_k)\}_{k=1}^\infty$ is a sequence of
based, connected, $n$-dimensional Riemannian manifolds. In addition,
suppose the following:

\begin{enumerate}
\item[(1)] There is $\delta>0$ such that $x_k\in {\rm
Reg}_\delta(U_k,g_k)$ for all $k$.
\item[(2)] For each $R<\infty$ and $\delta>0$ there is a constant $V(R,\delta)<\infty$ such
that ${\rm Vol}(B(x_k,R)\cap {\rm Reg}_\delta(U_k,x_k))\le
V(R,\delta)$ for all $k$ sufficiently large.
\item[(3)] For each non-negative integer $\ell$, each $\delta>0$, and each $R<\infty$,
 there is a constant
$C(\ell,\delta,R)$ such that for every $k$ sufficiently large we
have
$$|\nabla^\ell {\rm Rm}(g_k)|\le C(\ell,\delta,R)$$
on all of $B(x_k,R)\cap{\rm Reg}_\delta(U_k,g_k)$.
\item[(4)] For every $R<\infty$ there are $r_0>0$ and $\kappa>0$ such that
for every $k$ sufficiently large, for every $\delta\le r_0$ and
every $x\in B(x_k,R)\cap {\rm Reg}_\delta(U_k,g_k,x_k)$ the volume
of the metric ball $B(x,\delta)\subset U_k$ is at least
$\kappa\delta^n$.
\end{enumerate}
Then, after passing to a subsequence, there exists a based
Riemannian manifold $(U_\infty,g_\infty,x_\infty)$ that is a
geometric limit of the sequence $\{(U_k,g_k,x_k)\}_{k=1}^\infty$.
\end{thm}

Before giving the proof of this result, we begin with a standard
lemma.

\begin{lem}
Suppose that we have a sequence of $n$-dimensional balls $(B_k,h_k)$
of radius $r$ in Riemannian $n$-manifolds. Suppose that for each
$\ell$ there is a constant $C(\ell)$ such that for every $k$, we
have $|\nabla ^\ell {\rm Rm}(h_k)|\le C(\ell)$ throughout $B_k$.
Suppose also that for each $n$ the exponential mapping from the
tangent space at the center of $B_k$ induces a diffeomorphism from a
ball in the tangent space onto $B_k$. Then choosing an isometric
identification of the tangent spaces at the central points of the
$B_k$ with $\Ar^n$ and pulling back the metrics $h_k$ via the
exponential mapping to metrics $\tilde h_k$ on the ball $B$ of
radius $r$ in $\Ar^n$ gives us a family of metrics on $B$ that,
after passing to a subsequence, converge in the $C^\infty$-topology,
uniformly on compact subsets of $B$,  to a limit.
\end{lem}

The basic point in proving this lemma is to `find the right gauge,' which in
this case means find local coordinates so that the metric tensor is controlled
by the curvature. The correct local coordinates are the Gaussian coordinates
centered at the center of the ball.

\begin{proof}(of the theorem).
Fix $R<\infty$ and $\delta>0$. Let
$$X(\delta,R)=B(x_k,R)\cap
{\rm Reg}_{2\delta}(U_k,g_k,x_k).$$
 From the non-collapsing assumption
and the curvature bound assumption if follows from
Theorem~\ref{volinj} that there is a uniform positive lower bound
(independent of $k$) to the injectivity radius of every point in
$X(\delta,R)$. Fix $0<\delta'\le {\rm min}(r_0,\delta/2)$ much less
than this injectivity radius. We also choose $\delta'>0$
sufficiently small so that any ball of radius $2\delta'$ in
$B(x_k,R+\delta)\cap {\rm Reg}_\delta(U_k,g_k,x_k)$ is geodesically
convex. (This is possible because of the curvature bound.)  We cover
$X(\delta,R)$ by balls $B_1',\ldots,B_N'$
 of radii $\delta'/2$ centered at points of $X(\delta,R)$ with the
property that the sub-balls of radius $\delta'/4$ are disjoint. We
denote by $B_i'\subset B_i\subset \widetilde B_i$ the metric balls
with the same center and radii $\delta'/2$, $\delta'$, and
$2\delta'$ respectively.
 Notice that each of the balls $\widetilde B_i$  is contained in $B(x_k,R+\delta)\cap {\rm
Reg}_{\delta}(U_k,g_k,x_k)$. Because $\delta'\le r_0$, because ${\rm
Vol}\,B(x_k,R+\delta)$ is bounded independent of $k$, and because the
concentric balls of radius $\delta'/4$ are disjoint, there is a uniform bound
(independent of $k$) to the number of such balls. Passing to a subsequence we
can assume that the number of balls in these coverings is the same for all $k$.
We number them $\widetilde B_1,\ldots, \widetilde B_N$. Next, using the
exponential mapping at the central point, identify each of these balls with the
ball of radius $2\delta'$ in $\Ar ^n$. By passing to a further subsequence we
can arrange that the metrics on each $\widetilde B_i$ converge uniformly. (This
uses the fact that the concentric balls of radius $2\delta\ge 4\delta'$ are
embedded in the $U_k$ by the exponential mapping.) Now we pass to a further
subsequence so that the distance between the centers of the balls converges,
and so that for any pair $\widetilde B_i$ and $\widetilde B_j$ for which the
limiting distance between their centers is less than $4\delta'$, the overlap
functions in the $U_k$ also converge. The limits of the overlap functions
defines a limiting equivalence relation on $\coprod_i\widetilde B_i$.

This allows us to form a limit manifold $\widehat U_\infty$. It is
the quotient of the disjoint union of the $\widetilde B_i$ with the
limit metrics under the limit equivalence relation. We set
$(U_\infty(\delta,R),g_\infty(\delta,R),x_\infty(\delta,R))$ equal
to the submanifold of $\widehat U_\infty$ that is the union of the
sub-balls $B_i\subset \widetilde B_i$ of radii $\delta'$. A standard
argument using partitions of unity and the geodesic convexity of the
balls $\widetilde B_i$ shows that, for all $k$ sufficiently large,
there are smooth embeddings $\varphi_k(\delta,R)\colon
U_\infty(\delta,R)\to B(x_k,R+\delta)\cap {\rm
Reg}_\delta(U_k,g_k,x_k)$ sending $x_\infty(\delta,R)$ to $x_k$ and
converging as $k\rightarrow \infty$, uniformly in the
$C^\infty$-topology on each $B_i$, to the identity. Furthermore, the
images of each of these maps contains $B(x_k,R)\cap {\rm
Reg}_{2\delta}(U_k,g_k,x_k)$; compare \cite{Cheeger}. Also, the pull
backs under these embeddings of the metrics $g_k$ converge uniformly
to $g_\infty(\delta,R)$.

Repeat the process with $R$ replaced by $2R$ and $\delta=\delta_1$
replaced by $\delta_2\le \delta_1/2$. This produces
$$\left(U_\infty(\delta_2,2R),g_\infty(\delta_2,2R),x_\infty(\delta_2,2R)\right)$$
and, for all $k$ sufficiently large, embeddings
$\varphi_k(\delta_2,2R)$ of this manifold into
$$B(x_k,2R+\delta_2)\cap {\rm Reg}_{\delta_2}(U_k,g_k,x_k).$$ Hence,
the image of these embeddings contains the images of the original
embeddings. The compositions $(\varphi_k(\delta_2,2R))^{-1}\circ
\varphi_k(\delta,R)$ converge to an isometric embedding
$$\left(U_\infty(\delta,R),g_\infty(\delta,R),x_\infty(\delta,R)\right)\to
\left(U_\infty(\delta_2,2R),g_\infty(\delta_2,2R),x_\infty(\delta_2,2R)\right).$$
Repeating  this construction infinitely often produces a manifold
$(U_\infty,g_\infty,x_\infty)$ which is written as an increasing union of open
subsets $V_k=U_\infty(\delta_k,2^kR)$, where the $\delta_k$ tend to zero as $k$
tends to $\infty$. For each $k$ the open subset $V_k$ has compact closure
contained in $V_{k+1}$. By taking a subsequence of the original sequence we
have maps $\varphi_k\colon V_k\to U_k$ so that (2a) in the definition of
geometric limits holds. Condition (2b) clearly holds by construction.
\end{proof}

Now let us turn to complete Riemannian manifolds, where the result
is the $C^\infty$-version of the classical Cheeger-Gromov
compactness.

\begin{lem}
Suppose that $(U_k,g_k,x_k)$ is a sequence of based Riemannian
manifolds and that there is a partial geometric limit
$(U_\infty,g_\infty,x_\infty)$ that is a complete Riemannian
manifold. Then this partial geometric limit is a geometric limit.
\end{lem}

\begin{proof}
Since the balls $B(x_\infty,R)$ have compact closure in $U_\infty$
and since
$${\rm Reg}_\delta(U_\infty,g_\infty,x_\infty)=U_\infty$$ for
every $\delta>0$, it is easy to see that the extra condition, (2b), in
Definition~\ref{smoothconv} is automatic in this case.
\end{proof}

Now as an immediate corollary of Theorem~\ref{basicconv} we have the
following.

\begin{thm}\label{mfdconv}
Let $\{(M_k,g_k,x_k)\}_{k=1}^\infty$ be a sequence  of connected,
based Riemannian manifolds. Suppose that:
\begin{enumerate}
\item[(1)] For every $A<\infty$ the ball
$B(x_k,A)$ has compact closure in $M_k$ for all $k$ sufficiently
large.
\item[(2)] For each integer $\ell\ge 0$ and each $A<\infty$ there is a constant
$C=C(\ell,A)$ such that for each $y_k\in B(x_k,A)$ we have
$$\left|\nabla^\ell{\rm Rm}(g_k)(y_k)\right|\le C$$
for all $k$ sufficiently large.
\item[(3)] Suppose also that there is a constant $\delta>0$ such that ${\rm
inj}_{(M_k,g_k)}(x_k)\ge \delta$ for all $k$ sufficiently large.
\end{enumerate}
 Then after passing to a subsequence there is a geometric limit which
 is a complete Riemannian manifold.
\end{thm}

\begin{proof}
By the curvature bounds, it follows from the Bishop-Gromov theorem
(Theorem~\ref{BishopGromov}) that for each $A<\infty$ there is a uniform bound
to the volumes of the balls $B(x_k,A)$ for all $k$ sufficiently large. It also
follows from the same result that  the uniform lower bound on the injectivity
radius at the central point implies that for each $A<\infty$ there is a uniform
lower bound for the injectivity radius on the entire ball $B(x_k,A)$, again for
$k$ sufficiently large. Given these two facts, it follows immediately from
Theorem~\ref{basicconv} that there is a geometric limit.

Since, for every $A<\infty$, the $B(x_k,A)$ have compact closure in
$M_k$ for all $k$ sufficiently large, it follows that for every
$A<\infty$ the ball $B(x_\infty,A)$ has compact closure in
$M_\infty$. This means that $(M_\infty,g_\infty)$ is complete.
\end{proof}

\begin{cor}\label{2ndmfdconv}
Suppose that $\{(M_k,g_k,x_k)\}_{k=1}^\infty$ is a sequence of based, connected
Riemannian manifolds. Suppose that the first two conditions in
Theorem~\ref{mfdconv} hold and suppose also that there are constants $\kappa>0$
and $\delta>0$ such that ${\rm Vol}_{g_k}B(x_k,\delta)\ge \kappa\delta^n$ for
all $k$. Then after passing to a subsequence there is a geometric limit which
is a complete Riemannian manifold.
\end{cor}

\begin{proof}
Let $A={\rm max}(\delta^{-2},C(0,\delta))$, where $C(0,\delta)$ is the constant
given in the second condition in Theorem~\ref{mfdconv}. Rescale, replacing the
Riemannian metric $g_k$ by $Ag_k$. Of course, the first condition of
Theorem~\ref{mfdconv} still holds as does the second with different constants,
and we have $\left|{\rm Rm}_{Ag_k}(y_k)\right|\le 1$ for all $y_k\in
B_{Ag_k}(x_k,\sqrt{A}\delta)$. Also, ${\rm
Vol}\,B_{Ag_k}(x_k,\sqrt{A}\delta)\ge \kappa (\sqrt{A}\delta)^n$. Thus, by the
Bishop-Gromov inequality (Theorem~\ref{BishopGromov}), we have ${\rm
Vol}_{Ag_k}B(x_k,1)\ge \kappa/\Omega$ where
$$\Omega=\frac{V(\sqrt{A}\delta)}{(\sqrt{A}\delta)^nV(1)},$$
where $V(a)$ is the volume of the ball of radius $a$ in hyperbolic
$n$-space (the simply connected $n$-manifold of constant curvature
$-1$).
 Since
$\sqrt{A}\delta\ge 1$, this proves that for the rescaled manifolds the absolute
values of the sectional curvatures on $B_{Ag_k}(x_k,1)$ are bounded by $1$ and
the ${\rm Vol}_{Ag_k}B(x_k,1)$ are bounded below by a positive constant
independent of $k$.
 According to Theorem~\ref{volinj} the lower bound on the volume of the ball
of radius $1$ and the curvature bound on the ball of radius $1$ yield a uniform
positive lower bound $r>0$ for the injectivity radius of the rescaled manifolds
at $x_k$. Hence, the injectivity radii at the base points of the original
sequence are bounded below by $\delta/\sqrt{A}$. This means that the original
sequence of manifolds satisfies the third condition in Theorem~\ref{mfdconv}.
Invoking this theorem gives the result.
\end{proof}

\subsection{Geometric convergence of manifolds in the case of Ricci
flow}

As the next theorem shows, because of Shi's theorem, it is much
easier to establish the geometric convergence manifolds in the
context of Ricci flows than in general.

\begin{thm}\label{sliceconv}
Suppose that $({\mathcal M}_k,G_k,x_k)$ is a sequence of based
generalized $n$-dimensional Ricci flows with ${\bf t}(x_k)=0$. Let
$(M_k,g_k)$ be the $0$ time-slice of $({\mathcal M}_k,G_k)$. Suppose
that for each $A<\infty$ there are constants $C(A)<\infty$ and
$\delta(A)>0$ such that for all $k$ sufficiently large the following
hold:
\begin{enumerate}
\item[(1)] the ball $B(x_k,0,A)$ has compact closure in $M_k$,
\item[(2)] there is an embedding $B(x_k,0,A)\times (-\delta(A),0]\to {\mathcal M}_k$ compatible with
the time function and with the vector field, \item[(3)]   $|{\rm Rm}|\le C(A)$
on the image of the embedding in the Item (2), and \item[(4)] there is $r_0>0$
and $\kappa>0$ such that ${\rm Vol}\,B(x_k,0,r_0)\ge \kappa r_0^n$ for all $k$
sufficiently large.
\end{enumerate} Then, after passing to a
subsequence, there is a geometric limit
$(M_\infty,g_\infty,x_\infty)$ of the $0$ time-slices
$(M_k,g_k,x_k)$. This limit is a complete Riemannian manifold.
\end{thm}

\begin{proof}
The first condition in Theorem~\ref{mfdconv} holds by our first
assumption. It is immediate from Shi's theorem (Theorem~\ref{shi})
that the second condition of Theorem~\ref{mfdconv} holds. The result
is then immediate from Corollary~\ref{2ndmfdconv}.
\end{proof}

\section{Geometric convergence of Ricci flows}

 In this section we extend this notion of geometric convergence for based
 Riemannian manifolds in the obvious way to geometric convergence of based
 Ricci flows. Then we give Hamilton's theorem about the existence of
 such geometric limits.

\begin{defn}
Let $\{({\mathcal M}_k,G_k,x_k)\}_{k=1}^\infty$ be a sequence of based
generalized Ricci flows. We suppose that ${\bf t}(x_k)=0$ for all $k$ and we
denote by $(M_k,g_k)$ the time-slice of $({\mathcal M}_k,G_k)$. For some
$0<T\le \infty$, we say that a based Ricci flow
$(M_\infty,g_\infty(t),(x_\infty,0))$ defined for $t\in (-T,0]$ is a {\sl
partial geometric limit Ricci flow}\index{geometric limit!for Ricci
flows!partial} if:
\begin{enumerate}
\item
There are open subsets $x_\infty\in V_1\subset V_2\subset \cdots \subset
M_\infty$ satisfying (1) of Definition~\ref{smoothconv} with $M_\infty$ in
place of $U_\infty$,
\item there is a sequence $0<t_1<t_2<\cdots$ with  ${\rm
lim}_{k\rightarrow\infty}t_k=T$, \item and maps
$$\widetilde\varphi_k\colon (V_k\times [-t_k,0])\to {\mathcal M}_k$$
compatible with time and the vector field
\end{enumerate}
such that the sequence of horizontal families of metrics
$\widetilde\varphi_k^*G_k$ converges uniformly on compact subsets of
$M_\infty\times (-T,0]$ in the $C^\infty$-topology to the horizontal
family of metrics $g_\infty(t)$ on $M_\infty\times (-T,0]$.
\end{defn}

Notice that the restriction to the $0$ time-slices of a partial
geometric limit of generalized Ricci flows is a partial geometric
limit of the $0$ time-slices.

\begin{defn}
For $0<T\le \infty$, if $(M_\infty,g_\infty(t),x_\infty),\ -T<t\le 0$, is a
partial geometric limit Ricci flow  of the based generalized Ricci flows
$({\mathcal M}_k,G_k,x_k)$ and if $(M_\infty,g_\infty(0),x_\infty)$ is a
geometric limit of the $0$ time-slices, then we say that the partial geometric
limit is a {\em geometric limit  Ricci flow defined on the time interval
$(-T,0]$}\index{geometric limit!for Ricci flows}.
\end{defn}

Again Shi's theorem, together with a computation of Hamilton, allows us to form
geometric limits of generalized Ricci flows. We have the following result due
originally to Hamilton \cite{Hamiltonlimits}.

\begin{prop}\label{partialflowlimit}
Fix constants $-\infty\le T'\le 0\le T\le \infty$ and suppose that $T'<T$. Let
$\{({\mathcal M}_k,G_k,x_k)\}_{k=1}^\infty$ be a sequence of based generalized
Ricci flows. Suppose that ${\bf t}(x_k)=0$ for all $k$, and denote by
$(M_k,g_k)$ the $0$ time-slice of $({\mathcal M}_k,G_k)$. Suppose that there is
a partial geometric limit $(M_\infty,g_\infty,x_\infty)$ for the
$(M_k,g_k,x_k)$ with open subsets $\{V_k\subset M_\infty\}$ and maps
$\varphi_k\colon V_k\to M_k$ as in Definition~\ref{smoothconv}. Suppose that
for every compact subset $K\subset M_\infty$ and every compact interval
$I\subset (T',T)$ containing $0$, for all $k$ sufficiently large, there is an
embedding $\widetilde \varphi_k(K,I)\colon K\times I\to {\mathcal M}_k$
compatible with time and the vector field and extending the map $\varphi_k$ on
the $0$ time-slice. Suppose in addition that for every $k$ sufficiently large
there is a uniform bound (independent of $k$) to the norm of Riemann curvature
on the image of $\widetilde\varphi_k(K,I)$. Then after passing to a subsequence
the flows $\widetilde\varphi_k^*G_k$ converge  to a partial geometric limit
Ricci flow $g_\infty(t)$ defined for $t\in (T',T)$.
\end{prop}

\begin{proof}
 Suppose that we have a partial geometric limit of the
time-zero slices as stated in the proposition. Fix a compact subset $K\subset
M_\infty$ and a compact sub-interval $I\subset (T',T)$. For all $k$
sufficiently large we have embeddings $\widetilde \varphi_k(K,I)$ as stated. We
consider the flows $g_k(K,I)(t)$ on $K\times I$ defined by pulling back the
horizontal metrics $G_k$ under the maps $\widetilde \varphi_k(K,I)$. These of
course satisfy the Ricci flow equation on $K\times I$. Furthermore, by
assumption the flows $g_k(K,I)(t)$ have uniformly bounded curvature. Then under
these hypothesis, Shi's theorem can be used to show that the curvatures of the
$g_k(K,I)$ are uniformly bounded $C^\infty$-topology. The basic computation
 done by Hamilton in \cite{Hamiltonlimits} shows that after passing
to a further subsequence, the Ricci flows $g_k(K,I)$ converge
uniformly in the $C^\infty$-topology to a limit flow on $K\times I$.
A standard diagonalization argument allows us to pass to a further
subsequence so that the pullbacks $\widetilde\varphi_k^*G_k$
converge uniformly in the $C^\infty$-topology on every compact
subset of $M_\infty\times (T',T)$. Of course, the limit satisfies
the Ricci flow equation.
\end{proof}

This `local' result leads immediately to the following result for
complete limits.

\begin{thm}\label{flowlimit}
Fix $-\infty\le T'\le 0\le T\le \infty$ with $T'<T$. Let $\{({\mathcal
M}_k,G_k,x_k)\}_{k=1}^\infty$ be a sequence of based generalized Ricci flows.
Suppose that ${\bf t}(x_k)=0$ for all $k$, and denote by $(M_k,g_k)$ the $0$
time-slice of $({\mathcal M}_k,G_k)$. Suppose that for each $A<\infty$ and each
compact interval $I\subset (T',T)$ containing $0$ there is a constant $C(A,I)$
such that the following hold for all $k$ sufficiently large:
\begin{enumerate}
\item[(1)]  the ball
$B_{g_k}(x_k,0,A)$ has compact closure in $M_k$, \item[(2)]  there is an
embedding $B_{g_k}(x_k,0,A)\times I\to {\mathcal M}_k$ compatible with time and
with the vector field,
\item[(3)] the norms of the Riemann curvature of $G_k$ on the image of the embedding in the
previous item are  bounded by $C(A,I)$, and
\item[(4)] there is $r_0>0$ and $\kappa>0$ with ${\rm Vol}\,B(x_k,0,r_0)\ge \kappa r_0^n$
for all $k$ sufficiently large.
\end{enumerate}
 Then after passing to a subsequence there is a
flow  $(M_\infty,g_\infty(t),(x_\infty,0))$ which is the geometric
limit. It is a solution to the Ricci flow equation defined for $t\in
(T',T)$. For every $t\in (T',T)$ the Riemannian manifold
$(M_\infty,g_\infty(t))$ is complete.
\end{thm}

\begin{proof}
By Theorem~\ref{sliceconv} there is a geometric limit $(M_\infty,g_\infty(0))$
of the $0$ time-slices, and the limit  is a complete Riemannian manifold. Then
by Proposition~\ref{partialflowlimit} there is a geometric limit flow defined
on the time interval $(T',T)$. Since for every $t\in (T',T)$ there is a compact
interval $I$ containing $0$ and $t$, it follows that the Riemann curvature of
the limit is bounded on $M_\infty\times I$. This means that the metrics
$g_\infty(0)$ and $g_\infty(t)$ are commensurable with each other. Since
$g_\infty(0)$ is complete so is $g_\infty(t)$.
\end{proof}

\begin{cor}
Suppose that $(U,g(t)),\ 0\le t<T<\infty$,\ is a Ricci flow. Suppose
that $|{\rm Rm}(x,t)|$ is bounded independent of $(x,t)\in U\times
[0,T)$. Then for any open subset $V\subset U$ with compact closure
in $U$, there is an extension of the Ricci flow $(V,g(t)|_V)$ past
time $T$.
\end{cor}

\begin{proof}
 Take a sequence $t_n\rightarrow T$ and consider the
sequence of Riemannian manifolds $(V,g(t_n))$. By Shi's theorem and the fact
that $V$ has compact closure in $U$, the restriction of this sequence of
metrics to $V$ has uniformly bounded curvature derivatives. Hence, this
sequence has a convergent subsequence with limit $(V,g_\infty)$, where the
convergence is uniform in the $C^\infty$-topology. Now by Hamilton's result
\cite{Hamiltonlimits} it follows that, passing to a further subsequence, the
flows $(V,g(T+t-t_n),(p,0))$ converge to a flow $(V,g_\infty(t),(p,0))$ defined
on $(0,T]$. Clearly, for any $0<t<T$ we have $g_\infty(t)=g(t)$. That is to
say, we have extended the original Ricci flow smoothly to time $T$. Once we
have done this, we extend it to a Ricci flow on $[T,T_1)$ for some $T_1> T$
using the local existence results. The extension to $[T,T_1)$ fits together
smoothly with the flow on $[0,T]$ by Proposition~\ref{patch}.
\end{proof}

\section{Gromov-Hausdorff convergence}

Let us begin with the notion of the Gromov-Hausdorff distance between based
metric spaces of finite diameter. Let $Z$ be a metric space. We define the
Hausdorff distance between subsets of $Z$ as follows: $d^Z_H(X,Y)$ is the
infimum of all $\delta\ge 0$ such that $X$ is contained in the
$\delta$-neighborhood of $Y$ and $Y$ is contained in the $\delta$-neighborhood
of $X$. For metric spaces $X$ and $Y$ we define the Gromov-Hausdorff distance
between them, denoted $D_{GH}(X,Y)$, to be the infimum over all metric spaces
$Z$ and isometric embeddings $f\colon X\to Z$ and $g\colon Y\to Z$ of the
Hausdorff distance between $f(X)$ and $g(Y)$. For pointed metric spaces $(X,x)$
and $(Y,y)$ of finite diameter, we define the Gromov-Hausdorff
distance\index{Gromov-Hausdorff distance|ii} between them, denoted
$D_{GH}((X,x),(Y,y))$, to be the infimum of $D_H^Z(f(X),g(Y))$ over all triples
$((Z,z),f,g)$ where $(Z,z)$ is a pointed metric space and $f\colon (X,x)\to
(Z,z)$ and $g\colon (Y,y)\to (Z,z)$ are base-point preserving isometries.

To see that $D_{GH}$ is a distance function we must establish the
triangle inequality. For this it is convenient to introduce
$\delta$-nets in metric spaces.

\begin{defn}\label{epsnet}
A $\delta$-net in $(X,x)$ is a subset $L$ of $X$ containing $x$
whose $\delta$-neighborhood covers $X$ and for which there is some
$\delta'>0$ with $d(\ell_1,\ell_2)\ge \delta'$ for all $\ell_1\not=
\ell_2$ in $L$.
\end{defn}

 Clearly, the Gromov-Hausdorff distance from a based
metric space $(X,x)$ to a $\delta$-net $(L,x)$ contained in it is at
most $\delta$. Furthermore, for every $\delta>0$ the based space
$(X,x)$ has a $\delta$-net: Consider subsets $L\subset X$ containing
$x$ with the property that the $\delta/2$-balls centered at the
points of $L$ are disjoint. Any maximal such subset (with respect to
the inclusion relation) is a $\delta$-net in $X$.

\begin{lem}
The Gromov-Hausdorff distance satisfies the triangle inequality.
\end{lem}

\begin{proof}
Suppose that $D_{GH}((X,x),(Y,y))=a$ and $D_{GH}((Y,y),(Z,z))=b$.
Fix any $\delta>0$. Then there is a metric $d_1$ on $X\vee Y$ such
that $d_1$ extends the metrics on $X,Y$ and the
$(a+\delta)$-neighborhood of $X$ is all of $X\vee Y$ as is the
$(a+\delta)$-neighborhood of $Y$. Similarly, there is a metric $d_2$
on $Y\vee Z$ with the analogous properties (with $b$ replacing $a$).
Take a $\delta$-net $(L,y)\subset (Y,y)$, and define
$$d(x',z')={\rm inf}_{\ell\in L}d(x',\ell)+d(\ell,z').$$
We claim that $d(x',z')>0$ unless $x'=z'$ is the common base point. The reason
is that if ${\rm inf}_{\ell\in L}d(x',\ell)=0$, then by the triangle
inequality, any sequence of $\ell_n\in L$ with $d(x',\ell_n)$ converging to
zero is a Cauchy sequence, and hence is eventually constant. This means that
for all $n$ sufficiently large, $x'=\ell_n\in L\cap X$ and hence $x'$ is the
common base point. Similarly for $z'$.

A straightforward computation shows that the function $d$ above,
together with the given metrics on $X$ and $Z$, define a metric on
$X\vee Z$ with the property that the $(a+b+3\delta)$-neighborhood of
$X$ is all of $X\vee Z$ and likewise for $Z$. Since we can do this
for any $\delta>0$, we conclude that $D_{GH}((X,x),(Z,z))\le a+b$.
\end{proof}

Thus, the Gromov-Hausdorff distance is a pseudo-metric. In fact, the
restriction of the Gromov-Hausdorff distance to complete metric
spaces of bounded diameter is a metric. We shall not establish this
result, though we prove below closely related results about the
uniqueness of Gromov-Hausdorff limits.

\begin{defn}\label{GHconverge}
We say that a sequence of  based metric spaces $(X_k,x_k)$ of
uniformly bounded diameter {\sl converges in the Gromov-Hausdorff
sense}\index{converge!Gromov-Hausdorff sense}\index{Gromov-Hausdorff
convergence|see{converge}} to a based metric space $(Y,y)$ of finite
diameter if
$${\rm lim}_{k\rightarrow\infty}D_{GH}((X_k,x_k),(Y,y))=0.$$
Thus, a based metric space $(X,x)$ of bounded diameter is the limit
of a sequence of $\delta_n$-nets $L_n\subset X$ provided that
$\delta_n\rightarrow 0$ as $n\rightarrow\infty$.
\end{defn}

\begin{exam}
A sequence of compact $n$-manifolds of diameter tending to zero has
a point as  Gromov-Hausdorff limit.
\end{exam}

\begin{defn}\label{seqconv}
Suppose that $\{(X_k,x_k)\}_k$ converges in the Gromov-Hausdorff sense to
$(Y,y)$. Then a {\sl realization sequence} for this convergence is a sequence
of triples $((Z_k,z_k),f_k,g_k)$ where, for each $k$, the pair $(Z_k,z_k)$ is a
based metric space,
$$f_k\colon (X_k,x_k)\to
(Z_k,z_k) \ \ \ \ {\rm and} \ \ \ g_k\colon (Y,y)\to (Z_k,z_k)$$ are
isometric embeddings  and $D_{GH}(f_k(X_k),g_k(Y))\rightarrow 0$ as
$k\rightarrow\infty$. Given a realization sequence for the
convergence, we say that a sequence $\ell_k\in X_k$ converges to
$\ell\in Y$ (relative to the given realization sequence) if
$d(f_k(\ell_k),g_k(\ell))\rightarrow 0$ as $i\rightarrow\infty$.
\end{defn}

Notice that, with a different realization sequence for the
convergence, a sequence $\ell_k\in X_k$ can converge to a different
point of $Y$. Also notice that, given a realization sequence for the
convergence, every $y\in Y$ is the limit of some sequence $x_k\in
X_k$, a sequence $x_k\in X_k$ has at most one limit in $Y$, and if
$Y$ is compact then every sequence $x_k\in X_k$ has a subsequence
converging to a point of $Y$. Lastly, notice that under any
realization sequence for the convergence, the base points $x_k$
converge to the base point $y$.

\begin{lem}
Let $(X_k,x_k)$ be a sequence of metric spaces whose diameters are uniformly
bounded. Then the $(X_k,x_k)$ converge in the Gromov-Hausdorff
sense\index{converge!Gromov-Hausdorff sense} to $(X,x)$ if and only if the
following holds for every $\delta>0$. For every $\delta$-net $L\subset X$, for
every $\eta>0$, and for every $k$ sufficiently large, there is a
$(\delta+\eta)$-net $L_k\subset X_k$ and a bijection $L_k\to L$ sending $x_k$
to $x$ so that the push forward of the metric on $L_k$ induced from that of
$X_k$ is $(1+\eta)$-bi-Lipschitz equivalent to the metric on $L$ induced from
$X$.
\end{lem}

For a proof see Proposition 3.5 on page 36 of \cite{Gromov}.

\begin{lem}
Let $(X_k,x_k)$ be a sequence of  based metric spaces whose diameters are
uniformly bounded. Suppose that $(Y,y)$ and $(Y',y')$ are limits in the
Gromov-Hausdorff sense of this sequence and each of $Y$ and $Y'$ are compact.
Then $(Y,y)$ is isometric to $(Y',y')$.
\end{lem}

\begin{proof}
By the triangle inequality for Gromov-Hausdorff distance, it follows
from the hypothesis of the lemma that $D_{GH}((Y,y),(Y',y'))=0$. Fix
$\delta>0$. Since $D_{GH}((Y,y),(Y',y'))=0$, for any $n>0$ and
finite $1/n$-net $L_n\subset Y$ containing $y$ there is an embedding
$\varphi_n\colon L_n\to Y'$ sending $y$ to $y'$ such that the image
is a $2/n$-net in $Y'$ and such that the map from $L_n$ to its image
is a $(1+\delta)$-bi-Lipschitz homeomorphism. Clearly, we can
suppose that in addition the $L_n$ are nested: $L_n\subset
L_{n+1}\subset \cdots$. Since $Y'$ is compact and $L_n$ is finite,
and we can pass to a subsequence so that ${\rm
lim}_{k\rightarrow\infty}\varphi_k|_{L_n}$ converges to a map
$\psi_n\colon L_n\to Y'$ which is a $(1+\delta)$-bi-Lipschitz map
onto its image which is a $2/n$ net in $Y'$. By a standard
diagonalization argument, we can arrange that
$\psi_{n+1}|_{L_{n}}=\psi_{n}$ for all $n$. The $\{\psi_n\}$ then
define an embedding $\cup_n L_n\to Y'$ that is a
$(1+\delta)$-bi-Lipschitz map onto its image which is a dense subset
of $Y'$. Clearly, using the compactness of $Y'$ this map extends to
a $(1+\delta)$-bi-Lipschitz embedding $\psi_\delta\colon (Y,y)\to
(Y',y')$ onto a dense  subset of $Y'$. Since $Y$ is also compact,
this image is in fact all of $Y'$. That is to say, $\psi_\delta$ is
a $(1+\delta)$-bi-Lipschitz homeomorphism $(Y,y)\to (Y',y')$. Now
perform this construction for  a sequence of $\delta_n\rightarrow 0$
and $(1+\delta_n)$-bi-Lipschitz homeomorphisms
$\psi_{\delta_n}\colon (Y,y)\to (Y',y')$. These form an
equicontinuous family so that by passing to a subsequence we can
extract a limit $\psi\colon (Y,y)\to (Y',y')$. Clearly, this limit
is an isometry.
\end{proof}

Now let us consider the more general case of spaces of not
necessarily bounded diameter. It turns out that the above definition
is too restrictive when applied to such spaces. Rather one takes:

\begin{defn}
For  based metric spaces $(X_k,x_k)$ (not necessarily of finite
diameter) to {\sl converge in the Gromov-Hausdorff sense to a based
metric space}\index{converge!Gromov-Hausdorff sense} $(Y,y)$ means
that for each $r>0$ there is a sequence $\delta_k\rightarrow 0$ such
that the sequence of balls $B(x_k,r+\delta_k)$ in $(X_k,x_k)$
converges in the Gromov-Hausdorff sense to the ball $B(y,r)$ in $Y$.
\end{defn}

Thus, a sequence of cylinders $S^{n-1}\times \Ar$ with any base points and with
the radii of the cylinders going to zero has the real line as Gromov-Hausdorff
limit.

\begin{lem}
Let $(X_k,x_k)$ be a sequence of locally compact metric spaces. Suppose that
$(Y,y)$ and $(Y',y')$ are complete, locally compact, based metric spaces that
are limits of the sequence in the Gromov-Hausdorff sense. Then there is an
isometry $(Y,y)\to (Y',y')$.
\end{lem}

\begin{proof}
We show that for each $r<\infty$ there is an isometry between the
closed balls $\overline{B(y,r)}$ and $\overline{B(y',r)}$. By the
local compactness and completeness,  these closed balls are compact.
Each is the limit in the Gromov-Hausdorff sense of a sequence
$B(x_k,r+\delta_k)$ for some $\delta_k\rightarrow 0$ as
$k\rightarrow\infty$. Thus, invoking the previous lemma we see that
these closed balls are isometric. We take a sequence
$r_n\rightarrow\infty$ and isometries $\varphi_n\colon
(B(y,r_n),y)\to (B(y',r_n),y')$. By a standard diagonalization
argument, we pass to a subsequence such that for each $r<\infty$ the
sequence $\varphi_n|_{B(y,r)}$ of isometry converges to an isometry
$\varphi_r\colon B(y,r)\to B(y',r)$. These then fit together to
define a global isometry $\varphi\colon (Y,y)\to (Y',y')$.
\end{proof}

If follows from this that if a sequence of points $\ell_k\in X_k$ converges to
$\ell\in Y$ under one realization sequence for the convergence and to $\ell'\in
Y$ under another, then there is an isometry of $(Y,y)$ to itself carrying
$\ell$ to $\ell'$.

\begin{exam}
Let $(M_n,g_n,x_n)$ be a sequence of based Riemannian manifolds
converging geometrically to $(M_\infty,g_\infty,x_\infty)$. Then the
sequence also converges in the Gromov-Hausdorff sense to the same
limit.
\end{exam}

\subsection{Precompactness}

There is a fundamental compactness result due to Gromov. We begin
with a definition.

\begin{defn}
A {\sl length space}\index{length space} is a connected metric space
$(X,d)$ such that for any two points $x,y$ there is a rectifiable
arc $\gamma$ with endpoints $x$ and $y$ and with the length of
$\gamma$ equal to $d(x,y)$.
\end{defn}

For any based metric space $(X,x)$ and constants $\delta>0$ and
$R<\infty$ let $N(\delta,R,X)$ be the maximal number of disjoint
$\delta$-balls in $X$ that can be contained in $B(x,R)$.

\begin{thm}\label{GHcompact}
Suppose that $(X_k,x_k)$ is a sequence of based length spaces. Then
there is a based length space $(X,x)$ that is the limit in the
Gromov-Hausdorff sense of a subsequence of the $(X_k,x_k)$ if for
every $\delta>0$ and $R<\infty$ there is an $N<\infty$ such that
$N(\delta,R,X_k)\le N$ for all $k$. On the other hand, if the
sequence $(X_k,x_k)$ has a Gromov-Hausdorff limit, then for every
$\delta>0$ and $R<\infty$ the $N(\delta,R,X_k)$ are bounded
independent of $k$.
\end{thm}

For a proof of this result see Proposition 5.2 on page 63 of
\cite{Gromov}.

\subsection{The Tits cone}

\begin{figure}[ht]
  \relabelbox{
  \centerline{\epsfbox{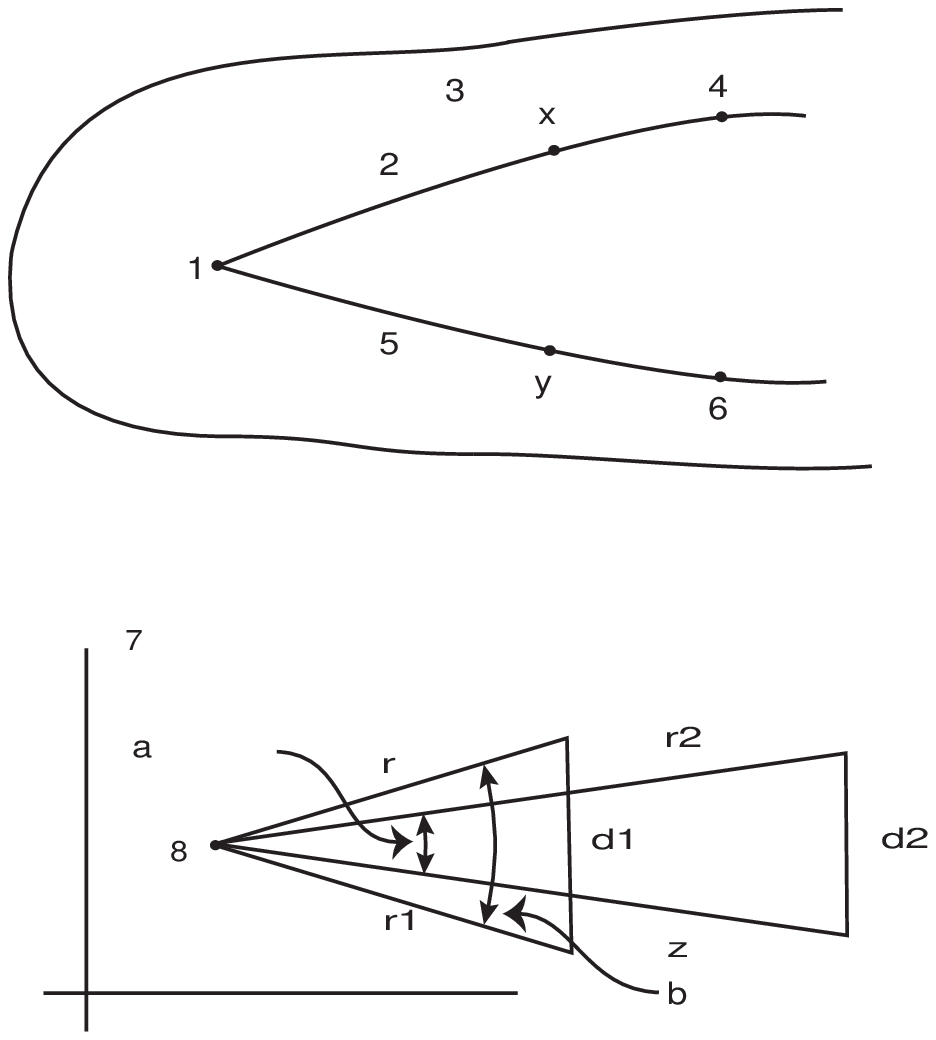}}}
  \relabel{1}{$p$}
  \relabel{2}{$\gamma$}
  \relabel{3}{$M$}
  \relabel{4}{$\gamma(r')$}
  \relabel{5}{$\mu$}
  \relabel{6}{$\mu(r')$}
  \relabel{7}{$\mathbb{R}^2$}
  \relabel{8}{$p'$}
  \relabel{x}{$\gamma(r)$}
  \relabel{y}{$\mu(r)$}
  \relabel{a}{$\Theta_{r'}(\gamma,\mu)$}
  \relabel{b}{$\Theta_r(\gamma,\mu)$}
  \relabel{r}{$r$}
  \relabel{r1}{$r$}
  \relabel{r2}{$r'$}
  \relabel{z}{$r'$}
    \relabel{d1}{$d(\gamma(r),\mu(r))$}
  \relabel{d2}{$d(\gamma'(r'),\mu(r'))$}
  \endrelabelbox
  \caption{Angles at infinity}\label{fig:anglesinf}
\end{figure}

Let $(M,g)$ be a complete, non-compact Riemannian manifold of
non-negative sectional curvature. Fix a point $p\in M$, and let
$\gamma$ and $\mu$ be minimal geodesic rays emanating from $p$. For
each $r>0$ let $\gamma(r)$ and $\mu(r)$ be the points along these
geodesic rays at distance $r$ from $p$. Then by  Part 1 of
Theorem~\ref{lengthcompar} we see that
$$\ell(\gamma,\mu,r)=\frac{d(\gamma(r),\mu(r))}{r}$$
is a non-increasing function of $r$. Hence, there is a limit
$\ell(\gamma,\mu)\ge 0$ of $\ell(\gamma,\mu,r)$ as
$r\rightarrow\infty$. We define the {\sl angle at infinity} between
$\gamma$ and $\mu$, $0\le\theta_\infty(\gamma,\mu)\le \pi$, to be
the angle at $b$ of the Euclidean triangle $a,b,c$ with side lengths
$|ab|=|bc|=1$ and $|bc|=\ell(\gamma,\mu)$, see {\sc
Fig.}~\ref{fig:anglesinf}.  If $\nu$ is a third geodesic ray
emanating from $p$, then clearly,
$\theta_\infty(\gamma,\mu)+\theta_\infty(\mu,\nu)\ge
\theta_\infty(\gamma,\nu)$.

\begin{defn}
Now we define a metric space whose underlying space is the quotient
space of the equivalence classes of minimal geodesic rays emanating
from $p$, with two rays equivalent if and only if the angle at
infinity between them is zero. The pseudo-distance function
$\theta_\infty$ descends to a metric on this space. This space is a
 length space \cite{BGP}. Notice that the distance between
any two points in this metric space is at most $\pi$. We denote this
space by $S_\infty(M,p)$.
\end{defn}

\begin{claim}
$S_\infty(M,p)$ is a compact space.
\end{claim}

\begin{proof}
Let $\{[\gamma_n]\}_n$ be a sequence of points in $S_\infty(M,p)$. We show that
there is a subsequence with a limit point. By passing to a subsequence we can
arrange that the unit tangent vectors to the $\gamma_n$ at $p$ converge to a
unit tangent vector $\tau$, say. Fix $d<\infty$, and let $x_n$ be the point of
$\gamma_n$ at distance $d$ from $p$. Then by passing to a subsequence we can
arrange that the $x_n$ converge to a point $x$. The minimizing geodesic
segments $[p,x_n]$ on $\gamma_n$ then converge to a minimizing geodesic segment
connecting $p$ to $x$. Performing this construction for a sequence of $d$
tending to infinity and then taking a diagonal subsequence produces a
minimizing geodesic ray $\gamma$ from $p$ whose class is the limit of a
subsequence of the $\{[\gamma_n]\}$
\end{proof}

 We define the {\em Tits cone}\index{Tits cone} of $M$ at $p$, denoted ${\mathcal T}(M,p)$,
 to be the cone over
$S_\infty(M,p)$, i.e.,  the quotient of the space $S_\infty(M,p)\times
[0,\infty)$ where all points $(x,0)$ are identified together (to become the
cone point). The cone metric on this space is given as follows: Let $(x_1,a_1)$
and $(x_2,a_2)$ be points of $S_\infty(M,p)\times [0,\infty)$. Then the
distance between their images in the cone is determined by
$$d^2([x_1,a_1],[x_2,a_2])=a_1^2+a_2^2-2a_1a_2{\rm cos}(\theta_\infty(x_1,x_2)).$$
 It is an easy exercise to show that the Tits cone of $M$ at $p$ is in fact
independent of the choice of $p$. From the previous claim, it follows that the
Tits cone of $M$ is locally compact and complete.

\begin{prop}\label{conelimit}
Let $(M,g)$ be a non-negatively curved, complete, non-compact
Riemannian manifold of dimension $k$. Fix a point $p\in M$ and let
$\{x_n\}_{n=1}^\infty$ be a sequence tending to infinity in $M$. Let
$\lambda_n=d^2(p,x_n)$ and consider the sequence of based Riemannian
manifolds $(M,g_n,p)$, where $g_n=\lambda_n^{-1}g$. Then there is a
subsequence converging in the Gromov-Hausdorff sense. Any
Gromov-Hausdorff limit
 $(X,g_\infty,x_\infty)$  of a
subsequence $(X,g_\infty)$ is isometric to the Tits cone ${\mathcal
T}(M,p)$ with base point the cone point.
\end{prop}

\begin{proof}
Let $c$ be the cone point of ${\mathcal T}(M,p)$, and denote by $d$ the
distance function on ${\mathcal T}(M,p)$. Consider the ball $B(c,R)\subset
{\mathcal T}(M,p)$. Since $S_\infty(M,p)$ is the metric completion of the
quotient space of minimal geodesic rays emanating from $p$, for any $\delta>0$
there is a $\delta$-net $L\subset B(c,R)$ consisting of the cone point together
with points of the form $([\gamma],t)$ where $\gamma$ is a minimal geodesic ray
emanating from $p$ and $t>0$. We define a map from $\psi_n\colon L \to (M,g_n)$
by sending the cone point to $p$ and sending $([\gamma],t)$ to the point at
$g_n$-distance $t$ from $p$ along $\gamma$. Clearly, $\psi_n(L)$ is contained
in $B_{g_n}(p,R)$. From the second item of Theorem~\ref{lengthcompar} and the
monotonicity of angles it follows that the map $\psi_n\colon L\to (M,g_n)$ is a
distance non-decreasing map; i.e., $\psi_n^*(g_n|\psi_n(L))\ge d|_L$. On the
other hand, by the monotonicity, $\psi_{n+1}^*(g_{n+1}|\psi_{n+1}(L))\le
\psi_n^*(g_n|\psi_n(L))$ and this non-increasing sequence of metrics converges
to $d|_L$. This proves that for any $\delta>0$ for all $n$ sufficiently large,
the embedding $\psi_n$ is a $(1+\delta)$-bi-Lipschitz homeomorphism.

 It remains to show that for any $\eta>0$ the images $\psi_n(L)$ are eventually
$\delta+\eta$-nets in $B_{g_n}(p,R)$. Suppose not. Then after
passing to a subsequence, for each $n$ we have a point $x_n\in
B_{g_n}(p,R)$ whose distance from $\psi_n(L)$ is at least
$\delta+\eta$. In particular, $d_{g_n}(x_n,p)\ge \delta$. Consider a
sequence of minimal geodesic rays $\mu_n$ connecting $p$ to the
$x_n$. Since the $g$-length of $\mu_n$ is at least $n\delta$,  by
passing to a further subsequence, we can arrange that the $\mu_n$
converge to a minimal geodesic ray $\gamma$ emanating from $p$. By
passing to a further subsequence if necessary, we arrange that
$d_{g_n}(x_n,p)$ converges to $r>0$. Now consider the points $\tilde
x_n$ on $\gamma$ at $g$-distance $\sqrt{\lambda_n}r$ from $p$.
Clearly, from the second item of Theorem~\ref{lengthcompar} and the
fact that the angle at $p$ between the $\mu_n$ and $\mu$ tends to
zero as $n\rightarrow\infty$ we have $d_{g_n}(x_n,\tilde
x_n)\rightarrow 0$ as $n\rightarrow\infty$. Hence, it suffices to
show that for all $n$ sufficiently large, $\tilde x_n$ is within
$\delta$ of $\psi_n(L)$ to obtain a contradiction. Consider the
point $z=([\mu],r)\in {\mathcal T}(M,p)$. There is a point
$\ell=([\gamma],t')\in L$ within distance $\delta$ of $z$ in the
metric $d$. Let $\tilde y_n\in M$ be the point in $M$ at
$g$-distance $\sqrt{\lambda_n}t'$ along $\gamma$. Of course, $\tilde
y_n=\psi_n(\ell)$. Then $d_{g_n}(\tilde x_n,\tilde y_n)\rightarrow
d(\ell,z)<\delta$. Hence, for all $n$ sufficiently large,
$d_{g_n}(\tilde x_n,\tilde y_n)<\delta$. This proves that for all
$n$ sufficiently large $\tilde x_n$ is within $\delta$ of
$\psi_n(L)$ and hence for all $n$ sufficiently large $x_n$ is within
$\delta+\eta$ of $\psi_n(L)$.

We have established that for every $\delta,\eta>0$ and every $R<\infty$ there
is a finite $\delta$-net $L$ in $({\mathcal T}(M,p),c)$ and for all $n$
sufficiently large an $(1+\delta)$-bi-Lipschitz embedding $\psi_n$ of $L$ into
$(M,g_n,p)$ with image a $\delta+\eta$-net for $(M,g_n,p)$. This proves that
the sequence $(M,g_n,p)$ converges in the Gromov-Hausdorff sense to ${\mathcal
T}(M,p),c))$.
\end{proof}

\section{Blow-up limits}

Here we introduce a type of geometric limit. These were  originally
introduced and studied by Hamilton in \cite{Hamiltonsurvey}, where,
among other things, he showed that  $3$-dimensional blow-up limits
have non-negative sectional curvature. We shall use repeatedly
blow-up limits and the positive curvature result in the arguments in
the later sections.

\begin{defn} Let $({\mathcal M}_k,G_k,x_k)$ be a sequence of based
generalized Ricci flows. We suppose that ${\bf t}(x_k)=0$ for all
$n$. We set $Q_k$ equal to $R(x_k)$. We denote by $(Q_k{\mathcal
M}_k,Q_kG_k,x_k)$ the family of generalized flows that have been
rescaled so that $R_{Q_kG_k}(x_k)=1$.
 Suppose that ${\rm
lim}_{k\rightarrow\infty}Q_k=\infty$ and that after passing to a
subsequence there is a geometric limit of the sequence
$(Q_k{\mathcal M}_k,Q_kG_k,x_k)$ which is a Ricci flow defined for
$-T<t\le 0$. Then we call this limit a {\em blow-up
limit}\index{blow-up limit|ii} of the original based sequence. In
the same fashion, if there is a geometric limit for a subsequence of
the zero time-slices of the $(Q_k{\mathcal M}_k,Q_kG_k,x_k)$, then
we call this limit
 the  blow-up limit of the $0$ time-slices.
\end{defn}

The significance of the condition that the generalized Ricci flows
have curvature pinched toward positive is that, as Hamilton
originally established in \cite{Hamiltonsurvey}, the latter
condition implies that any blow-up limit has non-negative curvature.

\begin{thm}\label{blowupposcurv}
Let $({\mathcal M}_k,G_k,x_k)$ be a sequence of generalized $3$-dimensional
Ricci flows, each of which has time interval of definition contained in
$[0,\infty)$ and each of which has curvature pinched toward positive. Suppose
that $Q_k=R(x_k)$ tends to infinity as $k$ tends to infinity. Let $t_k={\bf
t}(x_k)$ and let $({\mathcal M}_k',G_k',x_k)$ be the result of shifting time by
$-t_k$ so that ${\bf t}'(x_k)=0$. Then any blow-up limit of the sequence
$({\mathcal M}_k,G_k',x_k)$ has non-negative Riemann curvature. Similarly, any
blow-up limit of the zero time-slices of this sequence has non-negative
curvature.
\end{thm}

\begin{proof}
Let us consider the case of the geometric limit of the zero time-slice first.
Let $( M_\infty,g_\infty(0),x_\infty)$ be a blow-up limit of the zero
time-slices in the sequence. Let $V_k\subset M_\infty$ and $\varphi_k\colon
V_k\to (M_k)_0$ be as in the definition of the geometric limit. Let $y\in
{\mathcal M}_\infty$ be a point and let $\lambda(y)\ge \mu(y)\ge \nu(y)$ be the
eigenvalues of the Riemann curvature operator for $g_\infty$ at $y$. Let
$\{y_k\}$ be a sequence in $Q_k{\mathcal M}_k'$ converging to $y$, in the sense
that $y_k=\varphi_k(y)$ for all $k$ sufficiently large. Then
\begin{eqnarray*}
\lambda(y) & = & {\rm lim}_{n\rightarrow\infty}Q_k^{-1}\lambda(y_k) \\
\mu(y) & = & {\rm lim}_{n\rightarrow\infty}Q_k^{-1}\mu(y_k) \\
\nu(y) & = & {\rm lim}_{n\rightarrow\infty}Q_k^{-1}\nu(y_k)
\end{eqnarray*}
Since by Equation~(\ref{Rlower}) we have $R(y_k)\ge -6$ for all $k$ and since
by hypothesis $Q_k$ tends to infinity as $n$ does, it follows that $R(y)\ge 0$.
Thus if $\lambda(y)=0$, then ${\rm Rm}(y)=0$ and the result is established at
$y$. Hence, we may assume that $\lambda(y)>0$, which means that $\lambda(y_k)$
tends to infinity as $k$ does. If $\nu(y_k)$ remains bounded below as $k$ tends
to infinity, then $Q_k^{-1}\nu(y_k)$ converges to a limit which is $\ge 0$, and
consequently $Q_k^{-1}\mu(y_k)\ge Q_k^{-1}\nu(y_k)$ has a non-negative limit.
Thus, in this case the Riemann curvature of $g_\infty$ at $y$ is non-negative.
On the other hand, if $\nu(y_k)$ goes to $-\infty$ as $k$ does, then according
to Equation~(\ref{RXineq}) the ratio of $X(y_k)/R(y_k)$ goes to zero. Since
$Q_k^{-1}R(y_k)$ converges to the finite limit $R(y)$, the product
$Q_k^{-1}X(y_k)$ converges to zero as $k$ goes to infinity. This means that
$\nu(y)=0$ and consequently that $\mu(y)\ge 0$. Thus, once again we have
non-negative curvature for $g_\infty$ at $y$.

The argument in the case of a geometric limit flow is identical.
\end{proof}

\begin{cor}
Suppose that $(M_k,g_k(t))$ is a sequence of Ricci flows each of which has time
interval of definition contained in $[0,\infty)$ with each $M_k$ being a
compact $3$-manifold. Suppose further that, for each $k$, we have $|{\rm
Rm}(p_k,0)|\le 1$ for all $p_k\in M_k$. Then any blow-up limit of this sequence
of Ricci flows has non-negative curvature.
\end{cor}

\begin{proof}
According to Theorem~\ref{pinch} the hypotheses imply that for every $k$ the
Ricci flow $(M_k,g_k(t))$ has curvature pinched toward positive. From this, the
corollary follows immediately from the previous theorem.
\end{proof}

\section{Splitting limits at infinity}

In our later arguments we shall need a splitting result at infinity
in the non-negative curvature case. Assuming that a geometric limit
exists, the splitting result is quite elementary. For this reason we
present it here, though it will not be used until
Chapter~\ref{kappasect}.

\begin{figure}[ht]
  \relabelbox{
  \centerline{\epsfbox{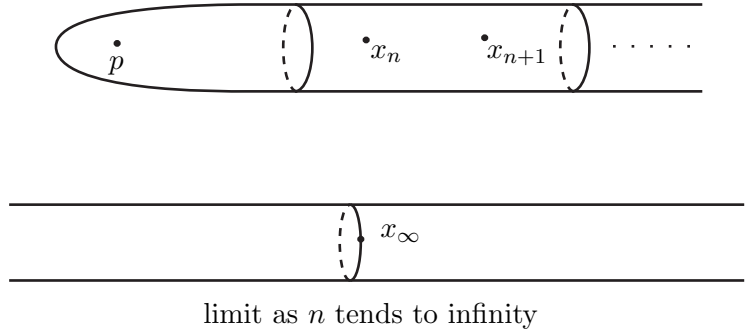}}}
  \relabel{1}{$p$}
  \relabel{2}{$x_n$}
  \relabel{3}{$x_{n+1}$}
  \relabel{4}{$x_\infty$}
  \relabel{5}{limit as $n$ tends to infinity}
  \endrelabelbox
  \caption{Splitting at infinity.}\label{fig:splinfty}
\end{figure}

The main result of this section gives a condition under which a geometric limit
automatically splits\index{splitting at infinity} off a line; see {\sc
Fig.}~\ref{fig:splinfty}.

\begin{thm}\label{topsplit}
Let $(M,g)$ be a complete, connected manifold of non-negative
sectional curvature. Let $\{x_n\}$ be a sequence of points going off
to infinity, and suppose that we can find scaling factors
$\lambda_n>0$ such that the based Riemannian manifolds
$(M,\lambda_ng,x_n)$ have a geometric limit
$(M_\infty,g_\infty,x_\infty)$. Suppose that there is a point $p\in
M$ such that $\lambda_nd^2(p,x_n)\rightarrow \infty$ as
$n\rightarrow\infty$.  Then, after passing to a subsequence,
minimizing geodesic arcs $\gamma_n$ from $x_n$ to $p$ converge to a
minimizing geodesic ray in $M_\infty$. This minimizing geodesic ray
is part of a minimizing geodesic line $\ell$ in $M_\infty$. In
particular, there is a Riemannian product decomposition
$M_\infty=N\times \Ar$ with the property that $\ell$ is $\{x\}\times
\Ar$ for some $x\in N$.
\end{thm}

\begin{proof}
Let $d_n$ be the distance from $p$ to $x_n$. Consider minimizing
geodesic arcs $\gamma_n$ from $p$ to $x_n$. By passing to a
subsequence we can assume that tangent directions at $p$ of these
arcs converge. Hence, for every $0<\delta<1$ there is $N$ such that
for all $n,m\ge N$ the angle between $\gamma_n$ and $\gamma_m$ at
$p$ is less than $\delta$. For any $n$ we can choose $m(n)$ such
that $d_{m(n))}\ge d_n(1+1/\delta)$. Let $\mu_n$ be a minimizing
geodesic from $x_n$ to $x_{m(n)}$.
 Now applying the Toponogov comparison
(first part of Theorem~\ref{lengthcompar}) and the usual law of
cosines in Euclidean space, we see that the distance $d$ from $x_n$
to $x_{m(n)}$ satisfies
$$d_{m(n)}-d_n\le d\le \sqrt{d_n^2+d_{m(n)}^2-2d_nd_{m(n)}{\rm
cos}(\delta)}.$$ Let $\theta_n=\angle_{x_n'}$ of the Euclidean triangle
$\triangle(x_n',p',x_{m(n)}')$ with $|s_{x_n'p'}|=d_n,|s_{x_n'x_{m(n)}'}|=d$
and $|s_{p'x_{m(n)}'}|=d_{m(n)}$. Then for any $\alpha<d_n$ and $\beta<d$ let
$x$ and $y$ be the points on $s_{x_np}$ and on $s_{x_nx_{m(n)}}$ at distances
$\alpha$ and $\beta$ respectively from $x_n$. Given this, according to the
Toponogov comparison result (first part of Theorem~\ref{lengthcompar}), we have
$$d(x,y)\ge \sqrt{\alpha^2+\beta^2-2\alpha\beta{\rm cos}(\theta_n)}.$$
The angle $\theta_n$ satisfies:
$$d_n^2+d^2-2d_nd{\rm cos}(\theta_n)=d_{m(n)}^2.$$
Thus, \begin{eqnarray*} {\rm cos}(\theta_n) & = &
\frac{d_n^2+d^2-d_{m(n)}^2}{2d_nd} \\
& \le & \frac{2d_n^2-2d_nd_{m(n)}{\rm cos}(\delta)} {2d_nd} \\
& = & \frac{d_n}{d}-\frac{d_{m(n)}}{d}{\rm cos}(\delta)\\
& \le & \delta-(1-\delta){\rm cos}(\delta).
\end{eqnarray*}
Since $\delta\rightarrow 0$ as $n\rightarrow\infty$, it follows that
given any $\delta>0$, for all  $n$ sufficiently large, $1+{\rm
cos}(\theta_n)<\delta$.


We are assuming that the based Riemannian manifolds
$\{(M,\lambda_ng,x_n)\}_{n=1}^\infty$ converge to a geometric limit
$(M_\infty,g_\infty,x_\infty)$. Also, by assumption,
$d_{\lambda_ng_n}(p,x_n)\rightarrow \infty$ as $n\rightarrow
\infty$, so that the lengths of the $\gamma_n$ tend to infinity in
the metrics $\lambda_ng_n$. This also means that the lengths of
$\mu_n$, measured in the metrics $\lambda_ng_n$, tend to infinity.
Thus, by passing to a subsequence we can assume that each of these
families, $\{\gamma_n\}$ and $\{\mu_n\}$, of minimizing geodesic
arcs converges to a minimizing geodesic arc, which we denote $\tilde
\gamma$ and $\tilde \mu$, respectively, in $M_\infty$ emanating from
$x_\infty$. The above computation shows that the angle between these
arcs is $\pi$ and hence that their union is a geodesic, say $\ell$.
The same computation shows that  $\ell $ is minimizing.

The existence of the minimizing geodesic line $\ell$ together with
the fact that the sectional curvatures of the limit are $\ge 0$
implies by Lemma~\ref{line} that the limit manifold is a Riemannian
product $N\times \Ar$ in such a way that $\ell$ is of the form
$\{x\}\times \Ar$ for some $x\in N$.
\end{proof}

\chapter{A comparison geometry approach to the Ricci
flow}\label{lengthfn}

In this section we discuss Perelman's notions, introduced in
\cite{Perelman1}, of the ${\mathcal L}$-length\index{${\mathcal
L}$-length} in the context of generalized Ricci flows. This is a
functional defined on paths in space-time parameterized by backward
time, denoted $\tau$. The ${\mathcal L}$-length is the analogue in
this context of the energy for paths in a Riemannian manifold. We
derive the associated Euler-Lagrange equation for the ${\mathcal
L}$-length; the critical paths are then ${\mathcal
L}$-geodesics\index{${\mathcal L}$-geodesic}. Using ${\mathcal
L}$-geodesics we define the ${\mathcal L}$-exponential
mapping\index{${\mathcal L}$-exponential map}. We derive the
${\mathcal L}$-Jacobi equation\index{${\mathcal L}$-Jacobi equation}
and relate ${\mathcal L}$-Jacobi fields\index{${\mathcal L}$-Jacobi
field} to the differential of the ${\mathcal L}$-exponential
mapping. There is the analogue of the interior of the cut locus. It
is the open subset, depending on the parameter $\tau$, of the
tangent space of initial vectors for ${\mathcal L}$-geodesics which
are minimizing out to time $\tau$ and at which the ${\mathcal
L}$-geodesic map is a local diffeomorphism at time $\tau$. The
difference between this situation and that of geodesics in a
Riemannian manifold is that there is such an open set in the tangent
space for each positive $\tau$. The analogue of the fact that, for
ordinary geodesics, the interior of the cut locus in the tangent
space is star-shaped from the origin is that the open set of `good'
initial conditions at $\tau$ is contained the open subset of `good'
initial conditions at time $\tau'$ for any $\tau'<\tau$. All of
these results are local and are established in the context of
generalized Ricci flows. In the next section we consider the case of
ordinary Ricci flows, where we are able to extend our results over
the entire manifold.

There are two applications of this theory in our study. In
Section~\ref{noncoll} we use the theory of ${\mathcal L}$-geodesics and the
associated notion of reduced volume to establish non-collapsing results. These
are crucial when we wish to take blow-up limits in studying singularities in
Ricci flows and Ricci flows with surgery. The second application will be in
Section~\ref{kappasect} to $\kappa$-solutions (ancient, $\kappa$-non-collapsed
solutions of bounded non-negative curvature). Here the second-order
inequalities on the length function that we establish in this section are used
to prove the existence of an asymptotic soliton for any $\kappa$-solution. This
asymptotic soliton is important for giving qualitative results on
$\kappa$-solutions.

\section{${\mathcal L}$-length and ${\mathcal L}$-geodesics}

The running assumption throughout this section is that we have an
$n$-dimensional generalized Ricci flow $({\mathcal M},G)$. In
particular, the space-time ${\mathcal M}$ is a smooth manifold of
dimension $n+1$ whose boundary lies at the initial and final times
(if they exist). Recall that its tangent bundle naturally decomposes
as the direct sum of the sub-line bundle spanned by the vector field
$\chi$ and the horizontal tangent bundle, denoted ${\mathcal
HT}{\mathcal M}$. We also fix a time $T$ in the time interval of
definition of the flow distinct from the initial time.

\begin{defn}
Let $0\le \tau_1<\tau_2$ be given and let $\gamma\colon
[\tau_1,\tau_2]\to {\mathcal M}$ be a continuous map. We say that
$\gamma$ is {\em parameterized by backward time}\index{parameterized
by backward time} provided that $\gamma(\tau)\in M_{T-\tau}$ for all
$\tau\in [\tau_1,\tau_2]$
\end{defn}

Throughout this section the paths $\gamma$ that we consider shall be
parameterized by backward time. We begin with the definition of
${\mathcal L}$-length of such a path.

\begin{defn}\label{Llength}
Let  $\gamma\colon [\tau_{1},\tau_{2}]\rightarrow {\mathcal M},\ 0
\le \tau_{1} < \tau_{2}$, be a  $C^{1}$-path parameterized by
backward time. We define $X_\gamma(\tau)$ to be the horizontal
projection of the tangent vector $d\gamma(\tau)/d\tau$, so that
$d\gamma/d\tau=-\chi+X_\gamma(\tau)$ with $X_\gamma(\tau)\in
{\mathcal HT}{\mathcal M}$. We define the $\mathcal{L}$-{\em
length}\index{${\mathcal L}$-length|ii} of $\gamma$ to be:
$$\mathcal{L}(\gamma) = \int_{\tau_{1}}^{\tau_{2}}
\sqrt{\tau}\left(R(\gamma(\tau)) +
\abs{X_\gamma(\tau)}^{2}\right)d\tau,\index{${\mathcal L}$|ii}$$
where the norm of $X_\gamma(\tau)$ is measured using the metric
$G_{T-\tau}$ on ${\mathcal HT}{\mathcal M}$. When $\gamma$ is clear
from the context, we write $X$ for $X_\gamma$; see {\sc
Fig.}~\ref{fig:path} from the Introduction.
\end{defn}

With a view toward better understanding the properties of the paths
that are critical points of this functional, the so-called
$\mathcal{L}$-geodesics, especially near $\tau=0$, it is helpful to
introduce a convenient reparameterization. We set $s=\sqrt{\tau}$.
We use the notation $A(s)$ to denote the horizontal component of the
derivative of $\gamma$ with respect to the variable $s$. One sees
immediately by the chain rule that
\begin{equation}\label{changeofvar}
A(s^2) = 2s X(s^2)\ \ {\rm or} \ \
A(\tau)=2\sqrt{\tau}X(\tau).\end{equation} With respect to the
variable $s$, the ${\mathcal L}$-functional is
\begin{equation}\label{seqn}
\mathcal{L}(\gamma)=
\int_{\sqrt{\tau_1}}^{\sqrt{\tau_2}}\left(\frac{1}{2}|A(s)|^{2} +
2R(\gamma(s))s^{2}\right)ds.
\end{equation}

Let's consider the simplest example.

\begin{exam}
Suppose that our generalized Ricci flow is a constant family of
Euclidean metrics on $\Ar^n\times [0,T]$. That is to say, $g(t) =
g_{0}$ is the usual Euclidean metric. Then we have
$R(\gamma(\tau))\equiv 0$. Using the change of variables
$s=\sqrt{\tau}$, we have
$$\mathcal{L}(\gamma) = \frac{1}{2}\int_{\sqrt{\tau_1}}^{\sqrt{\tau_2}} \abs{A(s)}^{2}ds,
$$ which is the standard energy functional in Riemannian geometry
for the path $\gamma(s)$. The minimizers for this functional are the
maps $s\mapsto (\alpha(s),T-s^2)$ where $\alpha(s)$ is a straight
line in $\Ar^n$ parameterized at constant speed. Written in the
$\tau$ variables the minimizers are
$$\gamma(\tau)=(x+\sqrt{\tau}v,T-\tau),$$
straight lines parameterized at speed varying linearly with
$\sqrt{\tau}$.

\end{exam}

\subsection{${\mathcal L}$-geodesics}

\begin{lem}\label{EL}  The
Euler-Lagrange equation\index{${\mathcal L}$-length!Euler Lagrange
equation} for critical paths for the $\mathcal{L}$-length is
\begin{equation}\label{EulerLag}
 \nabla_{X}X - \frac{1}{2}\nabla R
+ \frac{1}{2\tau}X + 2{\rm Ric}(X,\cdot)^* = 0.
\end{equation}
\end{lem}

\begin{rem}
${\rm Ric}(X,\cdot)$ is a horizontal one-form along $\gamma$ and its
dual ${\rm Ric}(X,\cdot)^*$ is a horizontal tangent vector field
along $\gamma$.
\end{rem}

\begin{proof} First, let us suppose that the
generalized Ricci flow is an ordinary Ricci flow $(M,g(t))$. Let
$\gamma_u(\tau)=\gamma(\tau,u)$ be a family of curves parameterized
by backward time. Let
$$\widetilde Y(\tau,u)=\frac{\partial\gamma}{\partial u}.$$ Then $\widetilde X(\tau,u)=X_{\gamma_u}(\tau,u)$
 and $\widetilde Y(\tau,u)$
are  the coordinate vector fields along the surface obtained by taking the
projection of $\gamma(\tau,u)$ into $M$. Thus, $[\widetilde X,\widetilde Y]
=0$. We denote by $X$ and $Y$ the restrictions of $\widetilde X$ and
$\widetilde Y$, respectively to $\gamma_0$. We have
\begin{align*}
\frac{d}{du}\mathcal{L}(\gamma_{u})\bigl|_{u=0}\bigr.  & =
\frac{d}{du}\left(\int_{\tau_{1}}^{\tau_{2}}
\sqrt{\tau}(R(\gamma_u(\tau)) + \abs{\widetilde
X(\tau,u)}^{2})d\tau\right)\Bigl|_{u=0}\Bigr.
\\&=\int_{\tau_{1}}^{\tau_{2}} \sqrt{\tau}(\langle \nabla R, Y\rangle  +
2\langle (\nabla_{Y}\widetilde X)|_{u=0},X\rangle )d\tau
\end{align*}
On the other hand, since $\partial g/\partial \tau =2{\rm Ric}$ and
since $[\widetilde X,\widetilde Y]=0$, we have
\begin{eqnarray*} 2\frac{d}{d\tau}(\sqrt{\tau}\langle Y,X\rangle_{g(T-\tau)} )
& = & \frac{1}{\sqrt{\tau}}\langle Y,X\rangle +
2\sqrt{\tau}\langle\nabla_XY,X\rangle
+2\sqrt{\tau}\langle Y,\nabla_XX\rangle \\ & & +4\sqrt{\tau}{\rm Ric}(Y,X) \\
& = & \frac{1}{\sqrt{\tau}}\langle Y,X\rangle  +
2\sqrt{\tau}\langle(\nabla_Y\widetilde X)|_{u=0},X\rangle
+2\sqrt{\tau}\langle Y,\nabla_XX\rangle \\ & &  +4\sqrt{\tau}{\rm
Ric}(Y,X)
\end{eqnarray*}
Using this we obtain \begin{eqnarray}\nonumber \frac{d}{du}
\mathcal{L}(\gamma_u)\bigl|_{u=0}\bigr. & = & \int_{\tau_1}^{\tau_2}
\Bigl(2\frac{d}{d\tau}\left[(\sqrt{\tau})\langle
Y,X\rangle\right]-\frac{1}{\sqrt{\tau}}\langle
Y,X\rangle\Bigr. \\
& & +\sqrt{\tau}\bigl(\langle\nabla R,Y\rangle
 - 2\langle Y,\nabla_{X}X\rangle  - 4{\rm Ric}(X,Y)\bigr)\Bigr)d\tau \nonumber \\
 & = & 2\sqrt{\tau}\langle Y,X\rangle |_{\tau_{1}}^{\tau_{2}} \nonumber \\
& & + \int_{\tau_{1}}^{\tau_{2}} \sqrt{\tau}\langle Y,\bigl(\nabla R
- \frac{1}{\tau}X - 2\nabla_{X}X - 4{\rm
Ric}(X,\cdot)^*\bigr)\rangle d\tau.\label{variation}
\end{eqnarray}

Now we drop the assumption that the generalized Ricci flow is an
ordinary Ricci flow. Still we can partition the interval
$[\tau_1,\tau_2]$ into finitely many sub-intervals with the property
that the restriction of $\gamma_0$ to each of the sub-intervals is
contained in a patch of space-time on which the generalized Ricci
flow is isomorphic to an ordinary Ricci flow. The above argument
then applies to each of the sub-intervals. Adding up
Equation~(\ref{variation}) over these sub-intervals shows that the
same equation for the first variation of length for the entire
family $\gamma_u$ holds.

We consider a variation $\gamma(\tau,u)$ with fixed endpoints, so
that  $Y(\tau_{1})=Y(\tau_{2})= 0.$ Thus, the condition that
$\gamma$ be a critical path for the ${\mathcal L}$-length is that
the integral expression vanish for all variations $Y$ satisfying
$Y(\tau_1)=Y(\tau_2)=0$. Equation~(\ref{variation}) holds for all
such $Y$ if and only if $\gamma$ satisfies
Equation~(\ref{EulerLag}).
\end{proof}

\begin{rem}
In the Euler-Lagrange equation, $\nabla R$ is the horizontal
gradient, and the equation is an equation of horizontal vector
fields along $\gamma$.
\end{rem}

\begin{defn}
A curve $\gamma$, parameterized by backward time, that is a critical
point of the ${\mathcal L}$-length is called an $\mathcal{L}$-{\em
geodesic}\index{${\mathcal L}$-geodesic|ii}.
Equation~(\ref{EulerLag}) is the {\em ${\mathcal L}$-geodesic
equation}\index{${\mathcal L}$-geodesic equation|ii}.
\end{defn}

Written with respect to the variable $s=\sqrt{\tau}$ the
$\mathcal{L}$-geodesic equation becomes
\begin{equation}\label{reparam}
\nabla_{A(s)}A(s) -2s^{2}\nabla R + 4s{\rm Ric}(A(s),\cdot)^*=0.\end{equation}
 Notice that in this form the ODE is
regular even at $s=0$.

\begin{lem}
Let $\gamma\colon [0,\tau_2]\to {\mathcal M}$ be an ${\mathcal L}$-geodesic.
Then ${\rm lim}_{\tau\rightarrow 0}\sqrt{\tau}X_\gamma(\tau)$ exists. The
${\mathcal L}$-geodesic $\gamma$ is completely determined by this limit (and by
$\tau_2$).
\end{lem}

\begin{proof}
Since the ODE in Equation~(\ref{reparam})  is non-singular even at
zero, it follows that $A(s)$ is a smooth function of $s$ in a
neighborhood of $s=0$. , The lemma follows easily by the change of
variables formula, $A(\tau)=2\sqrt{\tau}X_\gamma(\tau)$.
\end{proof}

\begin{defn}
An ${\mathcal L}$-geodesic is said to be {\em
minimizing}\index{${\mathcal L}$-geodesic!minimizing} if there is no
curve parameterized by backward time with the same endpoints and
with smaller ${\mathcal L}$-length.
\end{defn}

\subsection{The ${\mathcal L}$-Jacobi equation}

Consider a family $\gamma(\tau,u)$ of ${\mathcal L}$-geodesics parameterized by
$u$ and defined on $[\tau_1,\tau_2]$ with $0\le \tau_1<\tau_2$. Let $Y(\tau)$
be the horizontal vector field along $\gamma$ defined by
$$Y(\tau)=\frac{\partial}{\partial u}\gamma(\tau,u)|_{u=0}.$$

\begin{lem}
$Y(\tau)$ satisfies the {\em ${\mathcal L}$-Jacobi
equation}\index{${\mathcal L}$-Jacobi equation|ii}:
\begin{equation}\label{jaceqn}\nabla_X\nabla_XY+{\mathcal
R}(Y,X)X-\frac{1}{2}\nabla_Y(\nabla
R)+\frac{1}{2\tau}\nabla_XY+2(\nabla_Y{\rm Ric})(X,\cdot)^*+2{\rm
Ric}(\nabla_XY,\cdot)^*=0.\end{equation} This is a second-order
linear equation for $Y$. Supposing that $\tau_1>0$, there is a
unique horizontal vector field $Y$ along $\gamma$ solving this
equation vanishing at $\tau_1$ with a given first-order derivative
along $\gamma$ at $\tau_1$. Similarly, there is a unique solution
$Y$ to this equation vanishing at $\tau_2$ and with a given
first-order derivative at $\tau_2$.
\end{lem}

\begin{proof}
Given a family $\gamma(\tau,u)$ of ${\mathcal L}$-geodesics, then
from Lemma~\ref{EL} we have
$$ \nabla_{\widetilde X}\widetilde X = \frac{1}{2}\nabla R(\gamma)
-\frac{1}{2\tau}\widetilde X - 2{\rm Ric}(\widetilde X,\cdot)^*.$$
Differentiating this equation in the $u$-direction along the curve
$u=0$ yields
$$\nabla_Y
\nabla_{\widetilde X}\widetilde X|_{u=0}=\frac{1}{2}\nabla_Y(\nabla
R)-\frac{1}{2\tau}\nabla_Y(\widetilde X)|_{u=0} -2\nabla_Y({\rm
Ric}(\widetilde X,\cdot))^*|_{u=0}.$$ Of course, we have
$$\nabla_Y({\rm Ric}(\widetilde X,\cdot)^*)|_{u=0}=(\nabla_Y{\rm Ric})(X,\cdot)^*
+{\rm Ric}(\nabla_Y\widetilde X|_{u=0},\cdot)^*.$$ Plugging this in,
interchanging the orders of differentiation on the left-hand side, using
$\nabla_{\widetilde Y}\widetilde X=\nabla_{\widetilde X}\widetilde Y$, and
restricting to $u=0$ yields the equation given in the statement of the lemma.
This equation is a regular, second-order linear equation  for all $\tau>0$, and
hence is determined by specifying the value and first derivative at any
$\tau>0$.
\end{proof}

Equation~(\ref{jaceqn}) is obtained by applying $\nabla_Y$ to
Equation~(\ref{EulerLag}) and exchanging orders of differentiation.
The result Equation~(\ref{jaceqn}) is a second-order differential
equation for $Y$ that makes no reference to an extension of
$\gamma(\tau)$ to an entire family of curves.

\begin{defn} A field $Y(\tau)$ along an ${\mathcal L}$-geodesic is
called an {\em ${\mathcal L}$-Jacobi field}\index{${\mathcal
L}$-Jacobi field|ii} if it satisfies the ${\mathcal L}$-Jacobi
equation, Equation~(\ref{jaceqn}), and if it vanishes at $\tau_1$.
For any horizontal vector field $Y$ along $\gamma$ we denote by
${\rm Jac}(Y)$\index{${\rm Jac}$|ii} the expression on the left-hand
side of Equation~(\ref{jaceqn}).
\end{defn}

In fact, there is a similar result even for $\tau_1=0$.

\begin{lem}
Let $\gamma$ be an ${\mathcal L}$-geodesic defined on $[0,\tau_2]$
and let $Y(\tau)$ be an ${\mathcal L}$-Jacobi field along $\gamma$.
Then
$${\rm lim}_{\tau\rightarrow 0}\sqrt{\tau}\nabla_XY$$
exists. Furthermore, $Y(\tau)$ is completely determined by this
limit.
\end{lem}

\begin{proof}
We use the variable $s=\sqrt{\tau}$, and let $A(s)$ be the
horizontal component of $d\gamma/ds$. Then differentiating the
${\mathcal L}$-geodesic equation written with respect to this
variable we see
$$\nabla_A\nabla_AY=-{\mathcal R}(Y,A)A+2s^2\nabla_Y(\nabla R)-4s(\nabla_Y{\rm Ric})(A,\cdot)
-4s{\rm Ric}(\nabla_AY,\cdot).$$ Hence, for each tangent vector $Z$,
there is a unique solution to this equation with the two initial
conditions $Y(0)=0$ and $\nabla_AY(0)=Z$.

On the other hand, from Equation~(\ref{changeofvar}) we have
$\nabla_X(Y)=\frac{1}{2\sqrt{\tau}}\nabla_A(Y)$, so that
$$\sqrt{\tau}\nabla_X(Y)=\frac{1}{2}\nabla_A(Y).$$
\end{proof}

\subsection{Second order variation of ${\mathcal
L}$} We shall need the relationship of the ${\mathcal L}$-Jacobi
equation to the second-order variation of ${\mathcal L}$. This is
given in the next proposition.

\begin{prop}\label{2ndvarj}
Suppose that $\gamma$ is a minimizing ${\mathcal L}$-geodesic. Then,
for any vector field $Y$ along $\gamma$, vanishing at both
endpoints, and any family $\gamma_u$ of curves parameterized by
backward time with $\gamma_0=\gamma$ and with the $u$-derivative of
the family at $u=0$ being the vector field $Y$ along $\gamma$, we
have
$$\frac{d^2}{du^2}{\mathcal L}(\gamma_u)|_{u=0}=-\int_{\tau_1}^{\tau_2}2\sqrt{\tau}
\langle {\rm Jac}(Y),Y\rangle d\tau.$$ This quantity vanishes if and
only if $Y$ is an ${\mathcal L}$-Jacobi field.
\end{prop}

Let us begin the proof of this proposition with the essential computation.

\begin{lem}\label{2ndvari}
Let $\gamma$ be an ${\mathcal L}$-geodesic defined on $
[\tau_1,\tau_2]$, and let $Y_1$ and $Y_2$ be horizontal vector
fields along $\gamma$ vanishing at $\tau_1$. Suppose that
$\gamma_{u_1,u_2}$ is any family of curves parameterized by backward
time with the property that $\gamma_{0,0}=\gamma$ and the derivative
of the family in the $u_i$-direction at $u_1=u_2=0$ is $Y_i$. Let
$\widetilde Y_i$ be the image of $\partial/\partial u_i$ under
$\gamma_{u_1,u_2}$ and let $\widetilde X$ be the image of the
horizontal projection of $\partial/\partial \tau$ under this same
map, so that the restrictions of these three vector fields to the
curve $\gamma_{0,0}=\gamma$ are $Y_1,Y_2$ and $X$ respectively. Then
we have \begin{eqnarray*} \frac{\partial}{\partial
u_1}\frac{\partial}{\partial u_2}{\mathcal
L}(\gamma_{u_1,u_2})|_{u_1=u_2=0} & = & 2\sqrt{\tau_2}
Y_1(\tau_2)\langle \widetilde Y_2(\tau_2,u_1,0),\widetilde
X(\tau_2,u_1,0)\rangle|_{u_1=0}\\
& & -\int_{\tau_1}^{\tau_2}2\sqrt{\tau}\langle {\rm
Jac}(Y_1),Y_2\rangle d\tau.\end{eqnarray*}
\end{lem}

\begin{proof}
According to Equation~(\ref{variation}) we have
\begin{eqnarray*}
\frac{\partial}{\partial u_2}{\mathcal L}(\gamma)(u_1,u_2)
 & = & 2\sqrt{\tau_2}\langle \widetilde Y_2(\tau_2,u_1,u_2),\widetilde
X(\tau_2,u_1,u_2)\rangle
\\ & & -\int_{\tau_1}^{\tau_2}2\sqrt{\tau}\langle EL(\widetilde
X(\tau,u_1,u_2),\widetilde Y_2(\tau,u_1,u_2)\rangle d\tau,\end{eqnarray*} where
$EL(\widetilde X(\tau,u_1,u_2))$ is the Euler-Lagrange expression for
geodesics, i.e., the left-hand side of Equation~(\ref{EulerLag}).
Differentiating again yields:
\begin{eqnarray}
\lefteqn{\frac{\partial }{\partial u_1}\frac{\partial }{\partial u_2}{\mathcal
L}(\gamma_{u_1,u_2}\bigl)|_{u_1=u_2=0}\big.
 =  2\sqrt{\tau_2}Y_1(\tau_2)\langle \widetilde
Y_2(\tau_2,u_1,0),\widetilde X(\tau_2,u_1,0)\rangle\bigl|_{u_1=0}\bigr.  } \nonumber\\
& & -
\int_{\tau_1}^{\tau_2}2\sqrt{\tau}\left(\langle\nabla_{Y_1}EL(\widetilde
X), Y_2\rangle +\langle EL(X),\nabla_{Y_1}\widetilde
Y_2\rangle\right)(\tau,0,0)d\tau. \label{crosspart} \end{eqnarray}
Since $\gamma_{0,0}=\gamma$ is a geodesic, the second term in the
integrand vanishes, and since $[\widetilde X,\widetilde Y_1]=0$, we
have $\nabla_{Y_1}EL(\widetilde X(\tau,0,0))={\rm Jac}(Y_1)(\tau)$.
This proves the lemma.
\end{proof}

\begin{rem}\label{variremark}
Let $\gamma(\tau,u)$ be a family of curves as above with
$\gamma(\tau,0),\ \tau_1\le \tau\le \bar \tau$, being an ${\mathcal
L}$-geodesic. It follows from Lemma~\ref{2ndvari} and the remark
after the introduction of the ${\mathcal L}$-Jacobi equation that
the second-order variation of length at $u=0$ of this family is
determined by the vector field $Y(\tau)=\partial \gamma/\partial u$
along $\gamma(\cdot,0)$ and by the second-order information about
the curve $\gamma(\bar\tau,u)$ at $u=0$.
\end{rem}

\begin{cor}\label{symmY}
Let $\gamma$ be an ${\mathcal L}$-geodesic and let $Y_1,Y_2$ be vector fields
along $\gamma$ vanishing at $\tau_1$. Suppose $Y_1(\tau_2)=Y_2(\tau_2)=0$.
 Then the
bilinear pairing
$$-\int_{\tau_1}^{\tau_2}2\sqrt{\tau}\langle
{\rm Jac}(Y_1),Y_2\rangle d\tau$$ is a symmetric function of $Y_1$
and $Y_2$.
\end{cor}

\begin{proof}
Given $Y_1$ and $Y_2$ along $\gamma$ we construct a two-parameter
family of curves parameterized by backward time as follows. Let
$\gamma(\tau,u_1)$ be the value at $u_1$ of the geodesic through
$\gamma(\tau)$ with tangent vector $Y_1(\tau)$. This defines a
family of curves parameterized by backward time, the family being
parameterized by  $u_1$ sufficiently close to $0$. We extend
 $Y_1$ and $X$ to vector fields on this entire family by defining them to be
 $\partial/\partial u_1$ and the
horizontal projection of $\partial/\partial \tau$, respectively.
 Now we
extend the vector field $Y_2$ along $\gamma$ to a vector field on
this entire one-parameter family of curves. We do this so that
$Y_2(\tau_2,u_1)=Y_1(\tau_2,u_1)$. Now given this extension
$Y_2(\tau,u_1)$ we define a two-parameter family of curves
parameterized by backward time
 by setting $\gamma(\tau,u_1,u_2)$ equal to the value at $u_2$ of the geodesic
 through $\gamma(\tau,u_1)$ in the direction $Y_2(\tau,u_1)$.
 We then extend  $Y_1$, $Y_2$, and $X$ over this entire family by letting them be
 $\partial/\partial u_1$,
$\partial/\partial u_2$, and the horizontal projection of $\partial/ \partial
\tau$, respectively. Applying Lemma~\ref{2ndvari} and using the fact that
$Y_i(\bar\tau)=0$ we conclude that
$$\frac{\partial}{\partial u_1}\frac{\partial}{\partial u_2}{\mathcal L}(\gamma)|_{u_1=u_2=0}= -\int_{\tau_1}^{\tau_2}2\sqrt{\tau}\langle
{\rm Jac}(Y_1),Y_2\rangle d\tau$$ and symmetrically that
$$\frac{\partial}{\partial u_2}\frac{\partial}{\partial u_1}{\mathcal L}(\gamma)|_{u_1=u_2=0}= -\int_{\tau_1}^{\tau_2}2\sqrt{\tau}\langle
{\rm Jac}(Y_2),Y_1\rangle d\tau.$$ Since the second cross partials
are equal, the corollary follows.
\end{proof}

Now we are in a position to establish Proposition~\ref{2ndvarj}.

\begin{proof} (Of Proposition~\ref{2ndvarj})
From the equation in Lemma~\ref{2ndvari}, the equality of the second variation
of ${\mathcal L}$-length at $u=0$ and the integral is immediate from the fact
that $Y(\tau_2)=0$. It follows immediately that, if  $Y$ is an ${\mathcal
L}$-Jacobi field vanishing at $\tau_2$, then the second variation of the length
vanishes at $u=0$. Conversely, suppose given a family $\gamma_u$ with
$\gamma_0=\gamma$ with the property that the second variation of length
vanishes at $u=0$, and that the vector field $Y=(\partial\gamma/\partial
u)|_{u=0}$ along $\gamma$ vanishes at the end points. It follows that the
integral also vanishes. Since $\gamma$ is a minimizing ${\mathcal L}$-geodesic,
for any variation $W$, vanishing at the endpoints, the first variation of the
length vanishes and the second variation of length is non-negative. That is to
say,
$$-\int_{\tau_1}^{\tau_2}2\sqrt{\tau}\langle {\rm Jac}(W),W\rangle d\tau\ge 0$$
for all vector fields $W$ along $\gamma$ vanishing at the endpoints.
Hence, the restriction to the space of vector fields along $\gamma$
vanishing at the endpoints of the bilinear form
$$B(Y_1,Y_2)=-\int_{\tau_1}^{\tau_2}2\sqrt{\tau}\langle
{\rm Jac}(Y_1),(Y_2)d\tau,$$ which is symmetric by Corollary~\ref{symmY}, is
positive semi-definite. Since $B(Y,Y)=0$, it follows that $B(Y,\cdot)=0$; that
is to say, ${\rm Jac}(Y)=0$\index{${\rm Jac}$}.
\end{proof}

\section{The ${\mathcal L}$-exponential map and its first-order
properties}\index{${\mathcal L}$-exponential
map|ii}\label{sect:Lexp}

We use  ${\mathcal L}$-geodesics in order to define the ${\mathcal
L}$-exponential map.

For Section~\ref{sect:Lexp} we fix $\tau_1\ge 0$ and a point $x\in
{\mathcal M}$ with ${\bf t}(x)=T-\tau_1$. We suppose that $T-\tau_1$
is greater than the initial time of the generalized Ricci flow.
Then, for every $Z \in T_{x}M_{T-\tau_1},$ there is a maximal
$\mathcal{L}$-geodesic, denoted $\gamma_{Z}$, defined on some
positive $\tau$-interval, with $\gamma_{Z}(\tau_1)=x$ and with
$\sqrt{\tau_1}X(\tau_1)=Z$. (In the case $\tau_1=0$, this equation
is interpreted to mean $\lim_{\tau\rightarrow 0} \sqrt{\tau}X(\tau)
= Z$.)

\begin{figure}[ht]
  \relabelbox{
  \centerline{\epsfbox{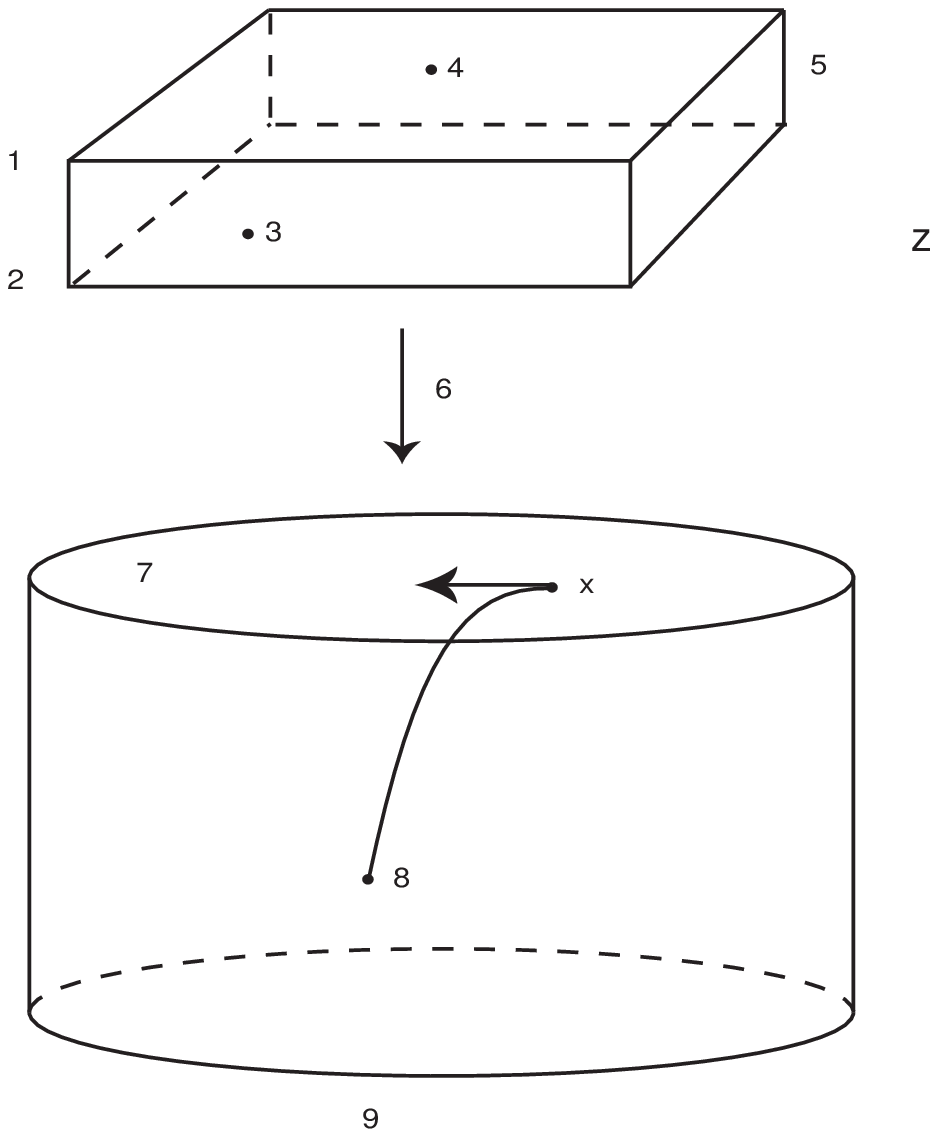}}}
  \relabel{1}{$0$}
  \relabel{2}{$\tau$}
  \relabel{3}{$(Z,\tau)$}
  \relabel{4}{$0$}
  \relabel{5}{$\mathcal{H}T_xM\times [0,\infty)$}
  \relabel{6}{$\mathcal{L}\text{exp}_x$}
  \relabel{7}{$\displaystyle{\lim_{\tau\to 0} \sqrt{\tau}X_\tau=Z}$}
  \relabel{8}{$\gamma_Z(\tau)=\mathcal{L}\text{exp}_x(Z,\tau)$}
  \relabel{9}{Space-time}
  \relabel{x}{$x$}
  \endrelabelbox
  \caption{The  map $\mathcal{L}\text{exp}_x$.}\label{fig:Lexp}
\end{figure}

\begin{defn}\label{defnD}
We define {\em the domain of definition of} ${\mathcal L}{\rm
exp}_x$\index{${\mathcal L}{\rm exp}$|ii},  denoted ${\mathcal D}_x$, to be the
subset of $T_xM_{T-\tau_1}\times (\tau_1,\infty)$ consisting of all pairs
$(Z,\tau)$ for which $\tau> \tau_1$ is in the maximal domain of definition of
$\gamma_Z$.  Then we define $\mathcal{L}{\rm exp}_{x}\colon {\mathcal D}_x\to
{\mathcal M}$ by setting $\mathcal{L}{\rm exp}_{x}(Z,\tau) = \gamma_{Z}(\tau)$
for all $(Z,\tau)\in {\mathcal D}_x$. (See {\sc Fig.}~\ref{fig:Lexp}.) We
define the map\index{$\widetilde L$} $\widetilde L\colon {\mathcal D}_x\to \Ar$
by $\widetilde L(Z,\tau)={\mathcal L}\left(\gamma_Z|_{[\tau_1,\tau]}\right)$.
Lastly, for any $\tau> \tau_1$ we denote by ${\mathcal L}{\rm
exp}_x^\tau$\index{${\mathcal L}{\rm exp}_x^\tau$|ii} the restriction of
${\mathcal L}{\rm exp}_x$ to the slice
$${\mathcal
D}^\tau_x={\mathcal D}_x\cap \left(T_xM_{T-\tau_1}\times \{\tau\}\right),$$
which is {\em the domain of definition of} ${\mathcal L}{\rm exp}_x^\tau$. We
also denote by  $\widetilde L^\tau$\index{$\widetilde L^\tau$|ii} the
restriction of $\widetilde L$ to this slice. We will implicitly identify
${\mathcal D}_x^\tau$ with a subset of $T_xM_{T-\tau_1}$.
\end{defn}

\begin{lem}\label{lcdiffeo}
${\mathcal D}_x$ is an open subset of $T_xM_{T-\tau_1}\times (\tau_1,\infty)$;
its intersection with each line $\{Z\}\times (\tau_1,\infty)$ is a non-empty
interval whose closure contains $\tau_1$. Furthermore, ${\mathcal L}{\rm
exp}_x\colon {\mathcal D}_x\to {\mathcal M}$ is a smooth map, and $\widetilde
L$ is a smooth function.
\end{lem}

\begin{proof}
The tangent vector in space-time of  the ${\mathcal L}$-geodesic
$\gamma$ is the vector field $-\chi+X_\gamma(\tau)$ along $\gamma$,
where $X_\gamma(\tau)$ satisfies Equation~(\ref{EulerLag}). As
above, in the case $\tau_1=0$, it is convenient to replace the
independent variable $\tau$ by $s=\sqrt{\tau}$, so that the ODE
becomes Equation~(\ref{reparam}) which is regular at $0$. With this
change, the lemma then follows immediately by the usual results on
existence, uniqueness and $C^\infty$-variation with parameters of
ODE's.
\end{proof}

\subsection{The differential of ${\mathcal L}{\rm exp}$}
Now we compute the differential of ${\mathcal L}{\rm exp}$.

\begin{lem}\label{DLJacobi}
Let $Z\in {\mathcal D}_x^{\bar\tau}\subset T_xM_{T-\tau_1}$. The
differential of ${\mathcal L}{\rm
exp}^{\bar\tau}_x$\index{${\mathcal L}{\rm exp}$!differential of|ii}
at the point $Z$ is given as follows: For each $W\in
T_x(M_{T-\tau_1})$ there is a unique ${\mathcal L}$-Jacobi
field\index{${\mathcal L}$-Jacobi field} $Y_W(\tau)$ along
$\gamma_Z$ with the property that $\sqrt{\tau_1}Y_W(\tau_1)=0$ and
$\sqrt{\tau_1}\nabla_X(Y_W)(\tau_1)=W$. We have
$$d_Z{\mathcal L}{\rm exp}^{\bar\tau}_x(W)=Y_W(\bar\tau).$$
Again, in case $\tau_1=0$, both of the conditions on $Y_W$ are interpreted as
the limits as $\tau\rightarrow 0$.
\end{lem}

\begin{proof}
Let $Z(u)$ be a curve in ${\mathcal D}_x^{\bar\tau}$ with $Z(0)=Z$.
Let $\gamma_u$ be the ${\mathcal L}$-geodesic starting at $x$ with
$\sqrt{\tau_1}X_{\gamma_u}(\tau_1)=Z(u)$. Then, clearly,
$$d_Z{\mathcal L}{\rm exp}^{\bar\tau}_x\left(\frac{dZ}{du}(0)\right)
=\frac{\partial}{\partial u}\left(\gamma_u(\bar\tau)\right)|_{u=0}.$$ On the
other hand, the vector field $Y(\tau)=\left(\partial\gamma_u(\tau)/\partial
u\right)|_{u=0} $ is an ${\mathcal L}$-Jacobi field along $\gamma_Z$.  Thus, to
complete the proof in the case when $\tau_1>0$ we need only see that
$\nabla_X\widetilde Y(\tau_1)=\nabla_Y\widetilde X(\tau_1)$. This is clear
since, as we have already seen, $[\widetilde X,\widetilde Y]=0$.

  When $\tau_1=0$, we complete the
argument using the following claim.

\begin{claim} If $\tau_1=0$, then
$$\frac{\partial}{\partial u}\left({\rm lim}_{\tau\rightarrow 0}\sqrt{\tau}X(\tau,u)\right)|_{u=0}=
{\rm lim}_{\tau\rightarrow 0}\sqrt{\tau}\frac{d}{d\tau}Y(\tau).$$
\end{claim}

\begin{proof} This follows immediately by changing variables,
replacing $\tau$ by $s=\sqrt{\tau}$.
\end{proof}
This completes the proof of Lemma~\ref{DLJacobi}.
\end{proof}

\subsection{Positivity of the second variation at a minimizing
${\mathcal L}$-geodesic}

If $\gamma$ is a minimizing ${\mathcal L}$-geodesic, then variations
of $\gamma$ fixing the endpoints  give curves whose ${\mathcal
L}$-length is no less than that of $\gamma$. In fact, there is a
second-order version of this inequality which we shall need later.

\begin{cor}\label{posform}
Let $Z\in T_xM_{T-\tau_1}$. Suppose that the associated ${\mathcal
L}$-geodesic $\gamma_Z$ minimizes ${\mathcal L}$-length between its
endpoints, $x$ and $\gamma_Z(\bar\tau)$, and that $d_Z{\mathcal
L}{\rm exp}^{\bar\tau}_x$ is an isomorphism. Then for any family
$\gamma_u$ of curves parameterized by backward time with
$Y=(\partial \gamma/\partial u)|_{u=0}$ vanishing at both endpoints,
we have
$$\frac{d^2}{du^2}{\mathcal L}(\gamma_u)|_{u=0}\ge 0,$$
with equality if and only if $Y=0$.
\end{cor}

\begin{proof}
According to Proposition~\ref{2ndvarj} the second variation  in the
$Y$-direction is non-negative and vanishes if and only if $Y$ is an ${\mathcal
L}$-Jacobi field. But since $d_Z{\mathcal L}{\rm exp}_x^{\bar\tau}$ is a
diffeomorphism, by Lemma~\ref{DLJacobi} there are no non-zero ${\mathcal
L}$-Jacobi fields vanishing at both endpoints of $\gamma_Z$.
\end{proof}

\subsection{The gradient of $\widetilde L^\tau$}

Recall\index{$\widetilde L^\tau$!gradient} that $\widetilde
L^{\tau}$ is the map from ${\mathcal D}^{\tau}_x$ to $\Ar$ that
assigns to each $Z$ the ${\mathcal L}$-length of
$\gamma_Z|_{[\tau_1,\tau]}$. We compute its gradient.

\begin{lem}\label{DtildeL}
Suppose that $Z\in {\mathcal D}^{\tau}_x$. Then for any $\widetilde
Y\in T_xM_{T-\tau_1}=T_Z({\mathcal D}^{\tau}_x)$ we have
$$ \langle \nabla\widetilde L^{\tau},\widetilde Y\rangle =2\sqrt{\tau}\langle
X(\tau),d_Z\left({\mathcal L}{\rm exp}_x^{\tau}\right)(\widetilde
Y)\rangle.$$
\end{lem}

\begin{proof}
Since ${\mathcal D}^{\tau}_x$ is an open subset of $T_x(M_{T-\tau_1})$, it
follows that for any $\widetilde Y\in T_x(M_{T-\tau_1})$ there is a
one-parameter family $\gamma_u(\tau')=\gamma(\tau',u)$ of
 ${\mathcal L}$-geodesics, defined for $\tau_1\le \tau'\le \tau$,
starting at $x$ with $\gamma(\cdot,0)=\gamma_Z$ and with
$\frac{\partial}{\partial u}\left(\sqrt{\tau_1}X(\tau_1)\right)=\widetilde Y$.
(Again, when $\tau_1=0$, this equation is interpreted to mean
$\frac{\partial}{\partial u}{\rm lim}_{\tau'\rightarrow
0}(\sqrt{\tau'}X(\tau',u))=\widetilde Y$.) Let
$Y(\tau')=\frac{\partial}{\partial u}(\gamma(\tau',u))|_{u=0}$ be the
corresponding ${\mathcal L}$-Jacobi field along $\gamma_Z$. Since
$\gamma(\tau_1,u)=x$ for all $u$, we have $Y(\tau_1)=0$. Since
$\gamma(\cdot,u)$ is an ${\mathcal L}$-geodesic for all $u$, according to
Equation~(\ref{variation}), and in the case $\tau_1=0$, using the fact that
$\sqrt{\tau}X(\tau')$ approaches a finite limit as $\tau\rightarrow 0$, we have
$$\frac{d}{du}{\mathcal
L}(\gamma_u)|_{u=0}=2\sqrt{\tau}\langle X(\tau),Y(\tau)\rangle.$$ By
Lemma~\ref{DLJacobi} we have $Y(\tau)=d_Z{\mathcal L}{\rm exp}_x^{\tau}(\tilde
Y)$. Thus,
$$\langle\nabla\widetilde L^{\tau},\widetilde Y\rangle=\frac{d}{du}{\mathcal
L}(\gamma_u)|_{u=0}=2\sqrt{\tau}\langle
X(\tau),Y(\tau)\rangle=2\sqrt{\tau}\langle X(\tau),d_Z({\mathcal
L}{\rm exp}^{\tau}_x)(\widetilde Y)\rangle.$$
\end{proof}

\subsection{Local diffeomorphism near the initial $\tau$}

Now let us use the nature of the ${\mathcal L}$-Jacobi equation to
study ${\mathcal L}{\rm exp}_x$ for $\tau>\tau_1$ but $\tau$
sufficiently close to $\tau_1$.

\begin{lem}\label{locstr}
For any $x$ in ${\mathcal M}$ with ${\bf
t}(x)=T-\tau_1$ and any $Z\in T_xM_{t-\tau_1}$, there is $\delta>0$
such that for any $\tau$ with $\tau_1<\tau<\tau_1+\delta$ the map
${\mathcal L}{\rm exp}_x^{\tau}$ is a local diffeomorphism from a
neighborhood of $Z$ in $T_xM_{T-\tau_1}$ to $M_{T-\tau}$.
\end{lem}

\begin{proof}
Fix $x$ and $Z$ as in the statement of the lemma. To establish the result it
suffices to prove that there is $\delta>0$ such that $d_Z{\mathcal L}{\rm
exp}_x^{\tau}$ is an isomorphism for all $\tau_1<\tau<\tau_1+\delta$. By
Lemma~\ref{DLJacobi} it is enough to find a $\delta>0$ such that any ${\mathcal
L}$-Jacobi field $Y$ along
 $\gamma_Z$ with $\sqrt{\tau_1}\nabla_XY(\tau_1)\not= 0$ does
 not vanish on the interval $(\tau_1,\tau_1+\delta)$.
 Because the ${\mathcal L}$-Jacobi equation is linear, it suffices to consider the
 case of ${\mathcal L}$-Jacobi fields with
 $|\nabla_XY(\tau_1)|=1$. The space of such fields is identified with the
 unit sphere in $T_xM_{T-\tau_1}$.
 Let us consider first the case when $\tau_1\not= 0$. Then for any
 such tangent vector $\nabla_XY(\tau_1)\not=0$. Since $Y(\tau_1)=0$,
 it follows that $Y(\tau)\not=0$ in some interval
 $(\tau_1,\tau_1+\delta)$, where $\delta$ can depend on $Y$.  Using
 compactness of the unit sphere in the tangent space, we see that there is $\delta>0$
 independent of $Y$ so that the above holds.

 In case when $\tau_1=0$, it is convenient to shift to the
 $s=\sqrt{\tau}$ parameterization. Then the geodesic equation and
 the ${\mathcal L}$-Jacobi equation are non-singular at the origin. Also, letting
 $A=d\gamma_Z/ds$ we have $\nabla_AY=2{\rm lim}_{\tau\rightarrow
 0}\sqrt{\tau}\nabla_XY$. In these variables, the argument for $\tau_1=0$
 is the same as the one above for $\tau_1>0$.
\end{proof}

\begin{rem}\label{backwards}
When $\tau_1>0$ it is possible to consider the ${\mathcal L}{\rm
exp}^{\tau}_x$ defined for $0<\tau<\tau_1$. In this case, the curves
are moving backward in $\tau$ and hence are moving forward with
respect to the time parameter ${\bf t}$. Two comments are in order.
First of all, for $\tau<\tau_1$, the gradient of ${\widetilde
L}^{\tau}_x$ is $-2\sqrt{\tau}X(\tau)$. The reason for the sign
reversal is that the length is given by the integral from $\tau$ to
$\tau_1$ and hence its derivative with respect to $\tau$ is the
negative of the integrand. The second thing to remark is that
Lemma~\ref{locstr} is true for $\tau<\tau_1$ with $\tau$
sufficiently close to $\tau_1$.
\end{rem}

\section{Minimizing ${\mathcal L}$-geodesics and the injectivity
domain}\label{sect:InjDom}

Now we discuss the analogue of the interior of the cut locus for the
usual exponential map of a Riemannian manifold. For
Section~\ref{sect:InjDom} we keep the assumption that $x\in
{\mathcal M}$ with ${\bf t}(x)=T-\tau_1$ for some $\tau_1\ge 0$.

\begin{defn}\label{injdefn}
 The {\sl injectivity set}\index{injectivity set|ii}
$\widetilde{\mathcal U}_x\subset {\mathcal D}_x\subset
\left(T_xM_{T-\tau_1}\times (\tau_1,\infty)\right)$\index{$\widetilde{\mathcal
U}_x$|ii} is the subset of all $(Z,\tau)\in {\mathcal D}_x$ with the following
properties:
\begin{enumerate}
\item The map ${\mathcal L}{\rm exp}_x^{\tau}$ is a local
diffeomorphism near $Z$ from $T_x(M_{T-\tau_1})$ to $M_{T-\tau}$.
\item There is a neighborhood ${\mathcal Z}$ of $Z$ in ${\mathcal D}^{\tau}_x$
 such that for every $Z'\in {\mathcal Z}$ the ${\mathcal L}$-geodesic
$\gamma_{Z'}|_{[\tau_1,\tau]}$ is the unique minimizing path
parameterized by backward time for the ${\mathcal L}$-length. That
is to say, the ${\mathcal L}$-length of
$\gamma_{Z'}|_{[\tau_1,\tau]}$ is less than the ${\mathcal
L}$-length of any other path parameterized by backward time between
the same endpoints.
\end{enumerate}
For any $\tau>\tau_1$, we set $\widetilde {\mathcal
U}_x(\tau)\subset T_xM_{T-\tau_1}$\index{$\widetilde{\mathcal
U}_x(\tau)$|ii} equal to the slice of $\widetilde{\mathcal U}_x$ at
$\tau$, i.e., $\widetilde {\mathcal U}_x(\tau)$ is determined by the
equation
$$\widetilde{\mathcal
U}_x(\tau)\times\{\tau\}=\widetilde{\mathcal
U}_x\cap\left(T_xM_{T-\tau_1}\times\{\tau\}\right).$$
\end{defn}

It is clear from the definition that $\widetilde{\mathcal
U}_x\subset {\mathcal D}_x$ is an open subset and hence
$\widetilde{\mathcal U}_x$ is an open subset of
$T_xM_{T-\tau_1}\times (\tau_1,\infty)$. Of course, this implies
that $\widetilde {\mathcal U}_x(\tau)$ is an open subset of
${\mathcal D}_x^\tau$ for every $\tau>\tau_1$.

\begin{defn}\label{L_x}
We set ${\mathcal U}_x\subset {\mathcal M}$\index{${\mathcal
U}_x$|ii} equal to $ {\mathcal L}{\rm exp}_x(\widetilde{\mathcal
U}_x)$. We call this subset of ${\mathcal M}$ the {\em injectivity
domain (of $x$)}\index{injectivity domain|ii}.
 For any
$\tau>\tau_1$ we set ${\mathcal U}_x(\tau)={\mathcal U}_x\cap
M_{T-\tau}$\index{${\mathcal U}_x(\tau)$|ii}.
\end{defn}

By definition, for every point $q\in {\mathcal U}_x$ for any
$(Z,\tau)\in \widetilde{\mathcal U}_x$ with ${\mathcal L}{\rm
exp}_x(Z,\tau)=q$, the ${\mathcal L}$-geodesic
$\gamma_{Z}|_{[\tau_1,\tau]}$  is a minimizing
 ${\mathcal L}$-geodesic to $q$. In particular, there is a minimizing
 ${\mathcal L}$-geodesic from $x$ to $q$.

\begin{defn}
The function $L_x\colon {\mathcal U}_x\to \Ar$\index{$L_x$|ii}
assigns to each $q$ in ${\mathcal U}_x$ the length of any minimizing
${\mathcal L}$-geodesic from $x$ to $q$. For any $\tau>\tau_1$, we
denote by $L_x^\tau$\index{$L_x^\tau$|ii} the restriction of $L_x$
to the $T-\tau$ time-slice of ${\mathcal U}_x$, i.e., the
restriction of $L_x$ to ${\mathcal U}_x(\tau)$.
\end{defn}

 This brings us to the
analogue of the fact that in Riemannian geometry the restriction to
the interior of the cut locus of the exponential mapping is a
diffeomorphism onto an open subset of the manifold.

\begin{prop}\label{diffeo}
The map
$${\mathcal L}{\rm exp}_x\colon\widetilde {\mathcal
U}_x\to {\mathcal M}$$ is a diffeomorphism onto the open subset
${\mathcal U}_x$ of ${\mathcal M}$. The function $L\colon {\mathcal
U}_x\to \Ar$ that associates to each $q\in {\mathcal U}_x$ the
length of the unique minimizing ${\mathcal L}$-geodesic from $x$ to
$q$ is a smooth function and
$$L_x\circ {\mathcal L}{\rm exp}_x|_{\widetilde{\mathcal U}_x}
=\widetilde L|_{\widetilde{\mathcal U}_x}.$$
\end{prop}

\begin{proof}
We consider the differential of ${\mathcal L}{\rm exp}_x$ at any
 $(Z,\tau)\in\widetilde {\mathcal U}_x$. By construction
 the restriction of this differential to $T_xM_{T-\tau_1}$ is a local
 isomorphism onto ${\mathcal H}T{\mathcal M}$ at the image point. On the other
 hand, the differential of ${\mathcal L}{\rm exp}_x$ in the $\tau$
 direction is $\gamma_Z'(\tau)$, whose `vertical' component is
 $-\chi$. By the inverse function theorem this
 shows that ${\mathcal L}{\rm exp}_x$ is a local diffeomorphism at
 $(Z,\tau)$, and
 its image is an open subset of ${\mathcal M}$.
The uniqueness in Condition 2, of the definition immediately implies
that the restriction of ${\mathcal L}{\rm exp}_x$ to
$\widetilde{\mathcal U}_x$ is one-to-one, and hence that it is a
global diffeomorphism onto its image ${\mathcal U}_x$.

Since for every $(Z,\tau)\in\widetilde{\mathcal U}_x$ the ${\mathcal
L}$-geodesic  $\gamma_Z|_{[\tau_1,\tau]}$ is ${\mathcal
L}$-minimizing, we see that $L_x\circ {\mathcal L}{\rm
exp}_x|_{\widetilde{\mathcal U}_x}=\widetilde L|_{\widetilde
{\mathcal U}_x}$ and that $L_x\colon {\mathcal U}_x\to \Ar$ is a
smooth function.
\end{proof}

According to Lemma~\ref{DtildeL} we have:

\begin{cor}\label{DL}
At any $q\in {\mathcal U}_x(\tau)$ we have
$$\nabla L_x^{\tau}(q)=2\sqrt{\tau}X(\tau)$$
where  $X(\tau)$ is the horizontal component of $\gamma'(\tau)$,
where $\gamma$ is the unique minimizing ${\mathcal L}$-geodesic
connecting $x$  to $q$. (See {\sc Fig.}~\ref{fig:path} in the
Introduction.)
\end{cor}

At the level of generality that we are working (arbitrary
generalized Ricci flows) there is no result analogous to the fact in
Riemannian geometry that the image under the exponential mapping of
the interior of the cut locus is an open dense subset of the
manifold. There is an analogue in the special case of Ricci flows on
compact manifolds or on complete manifolds of bounded curvature.
These will be discussed in Section~\ref{newcomp2}.

\subsection{Monotonicity of the $\widetilde{\mathcal U}_x(\tau)$ with respect to
$\tau$}\index{$\widetilde{\mathcal U}_x(\tau)$!monotonicty}

Next, we have the analogue of the fact in Riemannian geometry that
the cut locus is star-shaped.

\begin{prop}\label{inclusion}
Let $\bar\tau'>\bar\tau$. Then $\widetilde{\mathcal
U}_x(\bar\tau')\subset \widetilde{\mathcal U}_x(\bar\tau)\subset
T_xM_{T-\tau_1}$.
\end{prop}

\begin{proof}
For $Z\in \widetilde{\mathcal U}_x(\bar\tau')$,
 we shall show that:
(i) the ${\mathcal L}$-geodesic $\gamma_{Z'}|_{[\tau_1,\bar\tau]}$
is the unique minimizing ${\mathcal L}$-geodesic from $x$ to
$\gamma_{Z}(\bar\tau)$, and (ii) the differential $d_Z{\mathcal
L}{\rm exp}_x^{\bar\tau}$ is an isomorphism. Given these two
conditions, it follows from the definition that $\widetilde{\mathcal
U}_x(\bar\tau')$ is contained in $\widetilde{\mathcal
U}_x(\bar\tau)$.

 We  show that the ${\mathcal L}$-geodesic
$\gamma_{Z}|_{[\tau_1,\bar\tau]}$ is  the unique minimizing
${\mathcal L}$- geodesic to its endpoint.  If there is an ${\mathcal
L}$-geodesic $\gamma_1$, distinct from
$\gamma_{Z}|_{[\tau_1,\bar\tau]}$, from $x$ to $\gamma_Z(\bar\tau)$
whose ${\mathcal L}$-length is at most that of
$\gamma_Z|_{[\tau_1,\bar\tau]}$, then we can create a broken path
$\gamma_1*\gamma_Z|_{[\bar\tau,\bar\tau']}$ parameterized by
backward time whose ${\mathcal L}$-length is at most that of
$\gamma_{Z}$. Since this latter path is not smooth, its ${\mathcal
L}$-length cannot be the minimum, which is a contradiction.

 Now suppose  that $d_Z{\mathcal L}{\rm exp}_x^{\bar\tau}$ is not an isomorphism.
The argument is similar to the one above, using a non-zero
${\mathcal L}$-Jacobi field vanishing at both endpoints rather than
another geodesic. Let $\tau'_2$ be the first $\tau$ for which
$d_Z{\mathcal L}{\rm
 exp}_x^\tau$ is not an isomorphism.
According to Lemma~\ref{locstr}, $\tau_1<\tau'_2\le \bar\tau$. Since
${\mathcal L}{\rm exp}_x^{\tau'_2}$ is not a local diffeomorphism at
$(Z,\tau'_2)$, by Lemma~\ref{DLJacobi}  there is a non-zero
${\mathcal L}$-Jacobi field $Y$ along $\gamma_Z|_{[\tau_1,\tau'_2]}$
vanishing at both ends. Since $\gamma_Z|_{[\tau_1,\tau_2']}$ is
${\mathcal L}$-minimizing, according to Proposition~\ref{2ndvarj},
the second variation of the length of $\gamma_Z|_{\tau_1,\tau'_2]}$
in the $Y$-direction vanishes, in the sense that if $\gamma(u,\tau)$
is any one-parameter family of paths parameterized by backward time
from $x$ to $\gamma_Z(\tau'_2)$ with $(\partial \gamma/\partial
u)|_{u=0}=Y$ then
$$\frac{\partial^2{\mathcal L}(\gamma_u)}{\partial
u^2}\bigl|_{u=0}\bigr.=0.$$ Extend $Y$ to a horizontal vector field
$\widehat Y$ along $\gamma_Z$ by setting $\widehat Y(\tau)=0$ for
all $\tau\in [\tau'_2,\bar\tau]$. Of course, the extended horizontal
vector field $\widehat Y$ is not $C^2$ at $\tau'_2$ since $Y$, being
a non-zero ${\mathcal L}$-Jacobi field, does not vanish to second
order there. This is the first-order variation of the family
$\hat\gamma(u,\tau)$ that agrees with $\gamma(u,\tau)$ for all
$\tau\le \tau'_2$ and has $\hat\gamma(u,\tau)=\gamma_Z(\tau)$ for
all $\tau\in [\tau'_2,\bar\tau]$. Of course, the second-order
variation of this extended family at $u=0$ agrees with the
second-order variation of the original family at $u=0$, and hence
vanishes. But according to Proposition~\ref{2ndvarj} this means that
$\widehat Y$ is an ${\mathcal L}$-Jacobi field, which is absurd
since it is not a $C^2$-vector field.
\end{proof}

We shall also need a closely related result.

\begin{prop}\label{tausmooth}
Let $\gamma$ be a minimizing ${\mathcal L}$-geodesic defined for
$[\tau_1,\bar\tau]$. Fix $0\le\tau_1<\tau_2<\bar\tau$, and set
$q_2=\gamma(\tau_2)$, and $Z_2=\sqrt{\tau_2}X_{\gamma}(\tau_2)$.
Then, the map ${\mathcal L}{\rm exp}_{q_2}$ is  diffeomorphism from
a neighborhood of $\{Z_2\}\times (\tau_2,\bar\tau]$ in
$T_{q}M_{T-\tau_2}\times (\tau_2,\infty)$ onto a neighborhood of the
image of  $\gamma|_{(\tau_2,\bar\tau]}$.
\end{prop}

\begin{proof}
It suffices to show that the differential of ${\mathcal L}{\rm
exp}_{q_2}^\tau$ is an isomorphism for all $\tau\in
(\tau_2,\bar\tau]$. If this is not the case, then there is a
$\tau'\in (\tau_2,\bar\tau]$ and a non-zero ${\mathcal L}$-Jacobi
field $Y$ along $\gamma_Z|_{[\tau_2,\tau']}$ vanishing at both ends.
We extend $Y$ to a horizontal vector field $\widehat Y$ along all of
$\gamma_Z|_{[\tau_1,\tau']}$  by setting it equal to zero on
$[\tau_1,\tau_2]$. Since $Y$ is an ${\mathcal L}$-Jacobi field, the
second-order variation of ${\mathcal L}$-length in the direction of
$Y$ is zero, and consequently the second-order variation of the
length of $\gamma_Z|_{[\tau_1,\tau']}$ vanishes. Hence by
Proposition~\ref{2ndvarj} it must be the case that $\widehat Y$ is a
${\mathcal L}$-Jacobi field. This is impossible since $\widehat Y$
is not smooth at $\tau'$.
\end{proof}

We finish this section with a computation of the $\tau$-derivative
of $L_x$.

\begin{lem}\label{Ltau}
 Suppose that
$q\in  {\mathcal U}_x$ with ${\bf t}(q)=T-\bar\tau$ for some $\bar\tau>\tau_1$.
Let $\gamma\colon [\tau_1,\bar\tau]\to {\mathcal M}$ be the unique minimizing
${\mathcal L}$-geodesic from $x$ to $q$. Then we have
\begin{equation}\label{Ltaueqn}
\frac{\partial L_x}{\partial \tau}(q)=2\sqrt{\bar\tau}R(
q)-\sqrt{\bar\tau}\left(R( q)+|X(\bar\tau)|^2\right).
\end{equation}
\end{lem}

\begin{proof}
By definition and the Fundamental Theorem of Calculus, we have
$$\frac{d}{d\tau}L_x(\gamma(\tau))= \sqrt{\tau}\left(R(\gamma(\tau))+|X(\tau)|^2\right).$$
On the other hand since $\gamma'(\tau)=-\partial/\partial t+X(\tau)$
the chain rule implies
$$\frac{d}{d\tau}L_x(\gamma(\tau))=\langle \nabla
L_x,X(\tau)\rangle+\frac{\partial L_x}{\partial
\tau}(\gamma(\tau)),$$ so that
$$\frac{\partial L_x}{\partial \tau}(\gamma(\tau))=\sqrt{\tau}\left(R(\gamma(\tau))+|X(\tau)|^2\right)-\langle\nabla
L_x,X(\tau)\rangle.$$ Now using Corollary~\ref{DL}, and rearranging
the terms gives the result.
\end{proof}

\section{Second-order differential inequalities for $\widetilde L^{\bar\tau}$
and $L_x^{\bar\tau}$}\label{secondorder}

Throughout Section~\ref{secondorder} we fix $x\in {\mathcal M}$ with $x\in
M_{T-\tau_1}$.

\subsection{The second variation formula for $\widetilde L^{\bar\tau}$}

Our goal here is to compute the second variation of $\widetilde
L^{\bar\tau}$ in the direction of a horizontal vector field
$Y(\tau)$ along an ${\mathcal L}$-geodesic $\gamma$. Here is the
main result of this subsection.

\begin{prop}\label{2ndvariation} Fix $0\le \tau_1<\bar \tau$.
Let $\gamma$ be an ${\mathcal L}$-geodesic defined on
$[\tau_1,\bar\tau]$ and let $\gamma_u=\widetilde\gamma(\tau,u)$ be a
smooth family of curves parameterized by backward time with
$\gamma_0=\gamma$. Let $\widetilde Y(\tau,u)$ be $\partial
\widetilde\gamma/\partial u$ and let $\widetilde X$ be the
horizontal component of $\partial \widetilde \gamma/\partial \tau$.
These are horizontal vector fields along the image of $\widetilde
\gamma$. We set $Y$ and $X$ equal to the restrictions of $\widetilde
Y$ and $\widetilde X$, respectively, to $\gamma$. We assume that
$Y(\tau_1)=0$. Then
\begin{eqnarray*}
\lefteqn{\frac{d^2}{du^2}\left({\mathcal L}(\gamma_u)\right)|_{u=0}
 =  2\sqrt{\bar\tau}\langle \nabla_{Y(\tau)}\widetilde
Y(\bar\tau,u)|_{u=0},X(\bar\tau)\rangle } \\
& & + \int_{\tau_1}^{\bar\tau}\sqrt{\tau}({\rm Hess}(R)(Y,Y)
        +2\langle{\mathcal R}(Y,X)Y,X\rangle
        - 4(\nabla_{Y}{\rm Ric})(X,Y) \\
        & &  + 2(\nabla_{X}{\rm Ric})(Y,Y)+ 2\abs{\nabla_{X} Y}^{2})d\tau.
\end{eqnarray*}
\end{prop}

As we shall see, this is simply a rewriting of the equation in
Lemma~\ref{2ndvari} in the special case when $u_1=u_2$.

We begin the proof of this result with the following computation.

\begin{claim}  Let
$\gamma(\tau)$ be a curve parameterized by backward time. Let $Y$ be
a horizontal vector field along $\gamma$ and let $X$ be the
horizontal component of $\partial\widetilde\gamma/\partial \tau$.
Then
\begin{align*}\frac{\partial}{\partial\tau}\langle \nabla_{ X}Y, Y\rangle
 &= \langle \nabla_{ X} Y, \nabla_{ X} Y\rangle  + \langle
\nabla_{ X}\nabla_{ X} Y, Y\rangle
\\&\hspace{0.5cm} +2{\rm Ric}(\nabla_XY,Y)) +(\nabla_X{\rm Ric})(Y,Y))\end{align*}
\end{claim}

\begin{proof}
We can break $\frac{\partial}{\partial\tau}\langle \nabla_{X}Y,Y\rangle $ into
two parts: the first  assumes that the metric is constant and the second deals
with the variation with $\tau$ of the metric. The first contribution is the
usual formula
$$\frac{\partial}{\partial\tau}\langle \nabla_{X}Y,Y\rangle _{G(T-\tau_0)}=
\langle \nabla_{X}Y,\nabla_{X}Y\rangle _{G(T-\tau_0)}+ \langle
\nabla_{X}\nabla_{X}Y,Y\rangle _{G(T-\tau_0)}.$$ This gives us the first two
terms of the right-hand side of the equation in the claim.

We show that the last two terms in that equation come from differentiating the
metric with respect to $\tau$. To do this recall that in local coordinates,
writing the metric $G(T-\tau)$ as $g_{ij}$, we have
$$\langle \nabla_{X}Y,Y \rangle=g_{ij}\bigl(X^k\partial_k
Y^i+\Gamma_{kl}^iX^kY^l\bigr)Y^j.$$ There are two contributions coming from
differentiating the metric with respect to $\tau$. The first is when we
differentiate $g_{ij}$. This leads to
$$2{\rm Ric}_{ij}\bigl(X^k\partial_k
Y^i+\Gamma_{kl}^iX^kY^l\bigr)Y^j=2{\rm Ric}(\nabla_XY,Y\rangle.$$
The other contribution is from differentiating the Christoffel
symbols. This yields
$$g_{ij}\frac{\partial\Gamma^{i}_{kl}}{\partial \tau}X^{k}Y^{l}Y^{j}.$$
Differentiating the formula
$\Gamma_{kl}^i=\frac{1}{2}g^{si}(\partial_kg_{sl}+\partial_lg_{sk}-\partial_s
g_{kl})$ leads to \begin{eqnarray*} g_{ij}\frac{\partial
\Gamma_{kl}^i}{\partial \tau} & = & -2{\rm
Ric}_{ij}\Gamma_{kl}^i+g_{ij}g^{si}(\partial_k{\rm Ric}_{sl}+\partial_l{\rm
Ric}_{sk}
-\partial_s{\rm Ric}_{kl})\\
& = & -2{\rm Ric}_{ij}\Gamma_{kl}^i+
\partial_k{\rm Ric}_{jl}+\partial_l{\rm Ric}_{jk} -\partial_j{\rm Ric}_{kl}.
\end{eqnarray*}
Thus, we have
\begin{eqnarray*}
g_{ij}\frac{\partial\Gamma^{i}_{kl}}{\partial \tau}X^{k}Y^{l}Y^{j}
 & = & \bigl(-2{\rm Ric}_{ij}\Gamma^i_{kl}+\partial_k{\rm Ric}_{jl})\bigr)X^kY^lY^j \\
& = & (\nabla_X{\rm Ric})(Y,Y)
\end{eqnarray*}
This completes the proof of the claim.
\end{proof}

Now we return to the proof of the second variational formula in
Proposition~\ref{2ndvariation}.

\begin{proof}
According to Lemma~\ref{2ndvari} we have
$$\frac{d^2}{du^2}{\mathcal L}_{u=0}=2\sqrt{\bar
\tau}Y(\bar\tau)(\langle \widetilde Y(\bar\tau,u),\widetilde
X(\bar\tau,u)\rangle)|_{u=0}-\int_{\tau_1}^{\tau_2}2\sqrt{\tau}\langle
{\rm Jac}(Y),Y\rangle d\tau.$$\index{${\rm Jac}$} We plug in
Equation~\ref{jaceqn} for ${\rm Jac}(Y)$ and this makes the
integrand
\begin{eqnarray*}\sqrt{\tau}\langle\nabla_Y(\nabla
R),Y\rangle+2\sqrt{\tau}\langle{\mathcal R}(Y,X)Y,X\rangle
-\bigl(2\sqrt{\tau}\langle\nabla_X\nabla_XY,Y\rangle+\frac{1}{\sqrt{\tau}}\langle\nabla_XY,Y\rangle\bigr)\\
-4\sqrt{\tau}(\nabla_Y{\rm Ric})(X,Y)-4\sqrt{\tau}{\rm
Ric}(\nabla_XY,Y)
\end{eqnarray*}
The first term is $\sqrt{\tau}{\rm Hess}(R)(Y,Y)$. Let us deal with
the third and fourth terms, which are grouped together within
parentheses. According to the previous claim, we have
\begin{eqnarray*}
\frac{\partial}{\partial \tau}\bigl(2\sqrt{\tau}\langle
\nabla_XY,Y\rangle\bigr) & = &
\frac{1}{\sqrt{\tau}}\langle\nabla_XY,Y\rangle+2\sqrt{\tau}\langle\nabla_X\nabla_XY,Y\rangle
 +2\sqrt{\tau}\langle\nabla_XY,\nabla_XY\rangle \\ & & +4\sqrt{\tau}{\rm Ric}(\nabla_XY,Y)
 +2\sqrt{\tau} (\nabla_X{\rm Ric})(Y,Y).
\end{eqnarray*}
This allows us to replace the two terms under consideration by
$$-\frac{\partial}{\partial t}\bigl(2\sqrt{\tau}\langle
\nabla_XY,Y\rangle\bigr)+2\sqrt{\tau}\langle\nabla_XY,\nabla_XY\rangle+4\sqrt{\tau}{\rm
Ric}(\nabla_XY,Y) +2\sqrt{\tau} (\nabla_X{\rm Ric})(Y,Y).$$
Integrating the total derivative out of the integrand  and canceling
terms leaves the integrand as
\begin{eqnarray*}\sqrt{\tau}{\rm Hess}(R)(Y,Y)+2\sqrt{\tau}\langle{\mathcal R}(Y,X)Y,X\rangle
+2\sqrt{\tau}|\nabla_XY|^2 \\
-4\sqrt{\tau}(\nabla_Y{\rm Ric})(X,Y)+2\sqrt{\tau} (\nabla_X{\rm
Ric})(Y,Y),
\end{eqnarray*}
and makes the boundary term (the one in front of the integral) equal
to
$$2\sqrt{\bar\tau}\bigl(Y(\bar\tau)\langle
\widetilde Y(\bar\tau,u),\widetilde
X(\bar\tau,u)\rangle|_{u=0}-\langle
\nabla_XY(\bar\tau),Y(\bar\tau)\rangle\bigr)=2\sqrt{\bar\tau}\langle
X(\bar\tau),\nabla_Y\widetilde Y(\bar\tau,u)|_{u=0}\rangle.$$ This
completes the proof of the proposition.
\end{proof}

\subsection{Inequalities for the Hessian of $L_x^{\bar\tau}$}

Now we shall specialize the type of vector fields along $\gamma$.
This will allow us to give an inequality for the Hessian of
${\mathcal L}$ involving the integral of the vector field along
$\gamma$. These lead to inequalities for the Hessian of
$L_x^{\bar\tau}$\index{Hessian!of $L_x^{\bar\tau}$}. The main result
of this section is Proposition~\ref{Hessineq} below. In the end we
are interested in the case when the $\tau_1=0$. In this case the
formulas simplify. The reason for working here in the full
generality of all $\tau_1$ is in order to establish differential
inequalities at points not in the injectivity domain. As in the case
of the theory of geodesics, the method is to establish weak
inequalities at these points by working with smooth barrier
functions. In the geodesic case the barriers are constructed by
moving the initial point out the geodesic a small amount. Here the
analogue is to move the initial point of an ${\mathcal L}$-geodesic
from $\tau_1=0$ to a small positive $\tau_1$. Thus, the case of
general $\tau_1$ is needed so that we can establish the differential
inequalities for the barrier functions that yield the weak
inequalities at non-smooth points.

\begin{defn}
Let $q\in {\mathcal U}_x(\bar\tau)$ and let $\gamma\colon
[\tau_1,\bar\tau]\to {\mathcal M}$ be the unique minimizing
${\mathcal L}$-geodesic from $x$ to $q$. We say that a horizontal
vector field $Y(\tau)$ along $\gamma$ is {\em adapted} if it solves
the following ODE on $[\tau_1,\bar \tau]$:
\begin{equation}\label{firstYeqn}
\nabla_{X}Y(\tau) = -{\rm Ric}(Y(\tau),\cdot)^* +
\frac{1}{2\sqrt{\tau}(\sqrt{\tau}-\sqrt{\tau_1})}Y(\tau).
\end{equation}
 \end{defn}

Direct computation shows the following:

 \begin{lem}\label{adaptform}
 Suppose that $Y(\tau)$ is an adapted vector field along $\gamma$.
 Then
\begin{eqnarray}\label{secondYeqn}
\frac{d}{d\tau}\langle Y(\tau),Y(\tau)\rangle & = &
2{\rm Ric}(Y(\tau),Y(\tau)) + 2\langle \nabla_{X}Y(\tau),Y(\tau)\rangle \\
& = &  \frac{1}{\sqrt{\tau}(\sqrt{\tau}-\sqrt{\tau_1})}\langle
Y(\tau),Y(\tau)\rangle\nonumber.
\end{eqnarray} It follows that
$$|Y(\tau)|^2 = C\frac{(\sqrt{\tau}-\sqrt{\tau_1})^2}{(\sqrt{\bar\tau}-\sqrt{\tau_1})^2},$$
where $C=|Y(\bar\tau)|^2$.
\end{lem}

Now we come to the main result of this subsection, which is an
extremely important inequality for the Hessian of $L_x^{\bar\tau}$.

\begin{prop}\label{Hessineq}
Suppose that $q\in {\mathcal U}_x(\bar\tau)$, that $Z\in
\widetilde{\mathcal U}_x(\bar\tau)$ is the pre-image of $q$,  and
that $\gamma_Z$ is the  ${\mathcal L}$-geodesic to $q$ determined by
$Z$. Suppose  that $Y(\tau)$ is an adapted vector field along
$\gamma_Z$. Then
\begin{equation}\label{Hessin}
{\rm Hess}(L_x^{\bar\tau})(Y(\bar\tau),Y(\bar\tau)) \leq
\left(\frac{|Y(\bar\tau)|^2}{\sqrt{\bar\tau}-\sqrt{\tau_1}}\right)-
2\sqrt{\bar\tau}{\rm Ric}(Y(\bar\tau),Y(\bar\tau)) -
\int_{\tau_1}^{\bar\tau}\sqrt{\tau}H(X,Y)d\tau,
\end{equation} where
\begin{eqnarray} \nonumber
H(X,Y) & =&-{\rm Hess}(R)(Y,Y) - 2\langle  {\mathcal R}(Y,X)Y,X\rangle
\\& &-4(\nabla_{X}{\rm Ric})(Y,Y) +4(\nabla_{Y}{\rm Ric})(Y,X)\label{Heqn}
\\& & - 2\frac{\partial{\rm Ric}}{\partial\tau}(Y,Y) + 2\abs{{\rm Ric}(Y,\cdot)}^{2} -
\frac{1}{\tau}{\rm Ric}(Y,Y),\nonumber\index{$H(X,Y)$|ii}.
\end{eqnarray}
 We
have equality in Equation~(\ref{Hessin}) if and only if the adapted vector
field $Y$ is also a ${\mathcal L}$-Jacobi field.
\end{prop}

\begin{rem}
In spite of the notation, $H(X,Y)$ is a purely quadratic function of
the vector field $Y$ along $\gamma_Z$.
\end{rem}

We begin the proof of this proposition with three elementary lemmas. The first
is an immediate consequence of the definition of $\widetilde {\mathcal
U}_x(\bar\tau)$.

\begin{lem}
Suppose that $q\in {\mathcal U}_x(\bar\tau)$  and that $\gamma\colon
[\tau_1,\bar\tau]\to{\mathcal M}$ is the minimizing ${\mathcal
L}$-geodesic from $x$ to $q$. Then for every tangent vector $
Y(\bar\tau)\in T_qM_{T-\bar\tau}$ there is a one-parameter family of
${\mathcal L}$-geodesics   $\tilde\gamma(\tau,u)$ defined on
$[\tau_1,\bar\tau]$ with $\tilde\gamma(0,u)=x$ for all $u$, with
$\tilde\gamma(\tau,0)=\gamma(\tau)$ and $\partial\tilde
\gamma(\bar\tau,0)/\partial u=Y(\bar \tau)$. Also, for every $Z\in
T_xM_{T-\tau_1}$ there is a family of ${\mathcal L}$-geodesics
$\tilde\gamma(\tau,u)$ such that $\gamma(0,u)=x$ for all $u$,
$\tilde\gamma(\tau,0)=\gamma(\tau)$ and such that, setting
$Y(\tau)=\frac{\partial}{\partial u} \tilde\gamma_u(\tau)|_{u=0}$,
we have
$$\nabla_{\sqrt{\tau_1}X(\tau_1)}Y(\tau_1)=Z.$$
\end{lem}

\begin{lem}\label{Jacobi}
Let  $\gamma$ be a minimizing ${\mathcal L}$-geodesic from $x$, and
let $Y(\tau)$ be an ${\mathcal L}$-Jacobi field along $\gamma$. Then
$$2\sqrt{\bar\tau}\langle\nabla_XY(\bar\tau),Y(\bar\tau)\rangle = {\rm Hess}(L_x^{\bar\tau})(Y(\bar\tau),
Y(\bar\tau)).$$
\end{lem}

\begin{proof}
Let $\gamma(\tau,u)$ be a one-parameter family of ${\mathcal
L}$-geodesics emanating from $x$ with $\gamma(u,0)$ being the
${\mathcal L}$-geodesic in the statement of the lemma  and with
$\frac{\partial}{\partial u}\gamma(\tau,0)=Y(\tau)$. We have the
extensions of $X(\tau)$ and $Y(\tau)$ to vector fields $\widetilde
X(\tau,u)$ and $\widetilde Y(\tau,u)$ defined at $\gamma(\tau,u)$
for all $\tau$ and $u$. Of course,
\begin{eqnarray*}
\lefteqn{2\sqrt{\bar\tau}\langle \nabla_Y\widetilde
X(\bar\tau,u)|_{u=0},Y(\bar\tau)\rangle  } \\
& = & Y(\langle 2\sqrt{\bar\tau}\widetilde X(\bar\tau,u),\widetilde
Y(\bar\tau,u)\rangle)|_{u=0}-\langle
2\sqrt{\bar\tau}X(\bar\tau),\nabla_Y\widetilde
Y(\bar\tau,u)|_{u=0}\rangle.\end{eqnarray*} Then by
Corollary~\ref{DL} we have
\begin{align*}
2\sqrt{\bar\tau}\langle \nabla_Y\widetilde
X(\bar\tau,u)|_{u=0},Y(\bar\tau)\rangle & = Y(\bar\tau)(\langle
\nabla L_x^{\bar\tau},\widetilde
Y(\bar\tau,u)\rangle)|_{u=0}-\langle \nabla
L_x^{\bar\tau},\nabla_{Y(\bar\tau)}\widetilde
Y(\bar\tau,u)|_{u=0}\rangle
\\
& = Y(\bar\tau)(\widetilde Y(\bar\tau,u)
L_x^{\bar\tau})|_{u=0}-\nabla_{Y(\bar\tau)}\widetilde
Y(\bar\tau,u)|_{u=0}(L_x^{\bar\tau}) \\ & =  {\rm
Hess}(L_x^{\bar\tau})(Y(\bar\tau),Y(\bar\tau)).
\end{align*}
\end{proof}

Now suppose that we have a horizontal vector field that is both
adapted and ${\mathcal L}$-Jacobi. We get:

\begin{lem}\label{adapjab}
Suppose that $q\in {\mathcal U}_x(\bar\tau)$, that $Z\in \widetilde
{\mathcal U}_x(\bar\tau)$ is the pre-image of $q$, and that
$\gamma_Z$ is the  ${\mathcal L}$-geodesic to $q$ determined by $Z$.
Suppose further that $Y(\tau)$ is a horizontal vector field along
$\gamma$ that is both adapted and an ${\mathcal L}$-Jacobi field.
Then,  we have
$$\frac{1}{2\sqrt{\bar\tau}(\sqrt{\bar\tau}-\sqrt{\tau_1})}|Y(\bar\tau)|^2
=\frac{1}{2\sqrt{\bar\tau}}{\rm
Hess}(L_x^{\bar\tau})(Y(\bar\tau),Y(\bar\tau))+{\rm
Ric}(Y(\bar\tau),Y(\bar\tau)).$$
\end{lem}

\begin{proof}
From the definition of an adapted vector field  $Y(\tau)$  we have
\begin{equation*}
{\rm Ric}(Y(\tau),Y(\tau)) + \langle
\nabla_{X}Y(\tau),Y(\tau)\rangle =
\frac{1}{2\sqrt{\tau}(\sqrt{\tau}-\sqrt{\tau_1})}\langle
Y(\tau),Y(\tau)\rangle.
\end{equation*}
Since $Y(\tau)$ is an ${\mathcal L}$-Jacobi field, according to
Lemma~\ref{Jacobi} we have
$$\langle\nabla_XY(\bar\tau),Y(\bar\tau)\rangle =
\frac{1}{2\sqrt{\bar\tau}}{\rm
Hess}(L_x^{\bar\tau})(Y(\bar\tau),Y(\bar\tau)).$$ Putting these
together gives the result.
\end{proof}

Now we are ready to begin the proof of Proposition~\ref{Hessineq}.

\begin{proof}
Let $\tilde\gamma(\tau,u)$ be a family  of curves with
$\gamma(\tau,0)=\gamma_Z$ and with $\frac{\partial}{\partial
u}\gamma(\tau,u)=\widetilde Y(\tau,u)$. We denote by $Y$ the
horizontal vector field which is the restriction of $\widetilde Y$
to $\gamma_0=\gamma_Z$. We denote by
$q(u)=\tilde\gamma(\bar\tau,u)$. By restricting to a smaller
neighborhood of $0$ in the $u$-direction, we can assume that
$q(u)\in {\mathcal U}_x(\bar\tau)$ for all $u$. Then ${\mathcal
L}(\tilde\gamma_u)\ge L_x^{\bar\tau}(q(u))$. Of course,
$L_x^{\bar\tau}(q(0))={\mathcal L}(\gamma_Z)$. This implies that
$$\frac{d}{d u}L_x^{\bar\tau}(q(u))\bigl|\bigr._{u=0}=
\frac{d}{du}{\mathcal L}(\gamma_u)\bigl|\bigr._{u=0},$$
 and
$$
Y(\bar\tau)(\widetilde
Y(\bar\tau,u)(L^{\bar\tau}_x))|_{u=0}=\frac{d^2}{du^2}L_x^{\bar\tau}(q(u))\bigl|\bigr._{u=0}\le
\frac{d^2}{du^2}{\mathcal L}(\gamma_u)\bigl|\bigr._{u=0}.$$
 Recall that
$\nabla L_x^{\bar\tau}(q) = 2\sqrt{\bar\tau}X(\bar\tau)$, so that
$$\nabla_{Y(\bar\tau)}\widetilde Y(\bar\tau,u)|_{u=0}(L_x^{\bar\tau})=
\langle \nabla_{Y(\bar\tau)}\widetilde Y(\bar\tau,u)|_{u=0},\nabla
L^{\bar\tau}\rangle = 2\sqrt{\bar\tau}\langle
\nabla_{Y(\bar\tau)}\widetilde
Y(\bar\tau,u)|_{u=0},X(\bar\tau)\rangle .$$ Thus, by
Proposition~\ref{2ndvariation}, and using the fact that
$Y(\tau_1)=0$, we have
\begin{eqnarray*}
\lefteqn{{\rm Hess}(L^{\bar\tau})(Y(\bar\tau),Y(\bar\tau))  =
Y(\bar\tau)\left(\widetilde
Y(\bar\tau,u)(L_x^{\bar\tau})\right)|_{u=0}-
\nabla_{Y(\bar\tau)}\widetilde Y(\bar\tau,u)|_{u=0}(L_x^{\bar\tau})} \\
& \le & \frac{d^2}{du^2}{\mathcal
L}(\gamma_u)-2\sqrt{\bar\tau}\langle
\nabla_{Y(\bar\tau)}\widetilde Y(\bar\tau,u)|_{u=0},X(\bar\tau)\rangle \\
 & = & \int_{\tau_1}^{\bar\tau}
\sqrt{\tau}\bigl({\rm Hess}(R)(Y,Y) +
2\langle {\mathcal R}(Y,X)Y,X\rangle - 4(\nabla_{Y}{\rm Ric})(X,Y)\bigr. \\
& & \ \ \ \ \bigl.+ 2(\nabla_{X}{\rm Ric})(Y,Y) +
2\abs{\nabla_XY}^2\bigr)d\tau.
\end{eqnarray*}
 Plugging in Equation~(\ref{firstYeqn}), and using the fact that
 $|Y(\tau)|^2=
|Y(\bar\tau)|^2\frac{(\sqrt{\tau}-\sqrt{\tau_1})^2}{(\sqrt{\bar\tau}-\sqrt{\tau_1})^2}$,
gives
\begin{eqnarray*}
\lefteqn{{\rm Hess}(L^{\bar\tau})(Y(\bar\tau),Y(\bar\tau)) } \\ &
\le & \int_{\tau_1}^{\bar\tau} \sqrt{\tau}\bigl({\rm Hess}(R)(Y,Y) +
2\langle {\mathcal R}(Y,X)Y,X\rangle  - 4(\nabla_{Y}{\rm Ric})(X,Y)\bigr. \\
& &  \ \ \ \ \bigl. +2(\nabla_X{\rm Ric})(Y,Y) +
 2\abs{{\rm Ric}(Y,\cdot)}^{2}\bigr)d\tau \\ & &
 +\int_{\tau_1}^{\bar\tau}\left[
\frac{|Y(\bar\tau)|^2}{2\sqrt{\tau}(\sqrt{\bar\tau}-\sqrt{\tau_1})^2}
- \frac{2}{(\sqrt{\tau}-\sqrt{\tau_1})}{\rm Ric}(Y,Y)\right] d\tau
\end{eqnarray*}
 Using the definition of $H(X,Y)$\index{$H(X,Y)$} given in the statement, Equation~(\ref{Heqn}),
 allows us to write
\begin{eqnarray*}
\lefteqn{{\rm Hess}(L^{\bar\tau})(Y(\bar\tau),Y(\bar\tau))} \\
 & \le &
-\int_{\tau_1}^{\bar\tau}\sqrt{\tau}H(X,Y)d\tau \\ & &
+\int_{\tau_1}^{\bar\tau}\Bigl[\sqrt{\tau} \bigl(-2(\nabla_{X}{\rm
Ric})(Y,Y)
-2\frac{\partial{\rm Ric}}{\partial\tau}(Y,Y) +4|{\rm Ric}(Y,\cdot)|^2 \bigr) \Bigr.\\
& & \Bigl.+
\frac{|Y(\bar\tau)|^2}{2\sqrt{\tau}(\sqrt{\bar\tau}-\sqrt{\tau_1})^2}
-\left(\frac{2}{(\sqrt{\tau}-
\sqrt{\tau_1})}+\frac{1}{\sqrt{\tau}}\right) {\rm
Ric}(Y,Y)\Bigr]d\tau,
\end{eqnarray*}\index{$H(X,Y)$}

 To simplify
further, we compute, using Equation~(\ref{firstYeqn})
\begin{eqnarray*}
    \frac{d}{d\tau}\bigl({\rm Ric}(Y(\tau),Y(\tau))\bigr) &= & \frac{\partial{\rm Ric}}{\partial{\tau}}(Y,Y) +
        2{\rm Ric}(\nabla_{X}Y,Y) + (\nabla_{X}{\rm Ric})(Y,Y)\\
        & = & \frac{\partial{\rm Ric}}{\partial{\tau}}(Y,Y)+ (\nabla_{X}{\rm Ric})(Y,Y) \\ & & +
\frac{1}{\sqrt{\tau}(\sqrt{\tau}-\sqrt{\tau_1})}{\rm Ric}(Y,Y) -
2|{\rm Ric}(Y,\cdot)|^{2}.
\end{eqnarray*}
 Consequently, we have

\begin{eqnarray*}
\frac{d\left(2\sqrt{\tau}{\rm Ric}(Y(\tau),Y(\tau))\right)}{d\tau} &
= & 2\sqrt{\tau}\left(\frac{\partial{\rm Ric}}{\partial{\tau}}(Y,Y)
+ (\nabla_{X}{\rm
Ric})(Y,Y)- 2|{\rm Ric}(Y,\cdot)|^{2}\right)\\
 & & + \left(\frac{2}{(\sqrt{\tau}-\sqrt{\tau_1})} +\frac{1}{\sqrt{\tau}}\right){\rm Ric}(Y,Y)
\end{eqnarray*}
Using this, and the fact that $Y(\tau_1)=0$, gives
\begin{eqnarray}\label{Hessineq1}
\lefteqn{{\rm Hess}(L_x^{\bar\tau})(Y(\bar\tau),Y(\bar\tau)) \le}  \\
 &  &
-\int_{\tau_1}^{\bar\tau}\left(\sqrt{\tau}H(X,Y)
-\frac{d}{d\tau}\left(2\sqrt{\tau}{\rm Ric}(Y,Y)\right)
-\frac{|Y(\bar\tau)|^2}{2\sqrt{\tau}(\sqrt{\bar\tau}-\sqrt{\tau_1})^2}\right)d\tau \nonumber\\
 & = &
\frac{|Y(\bar\tau)|^2}{\sqrt{\bar\tau}-\sqrt{\tau_1}}
-2\sqrt{\bar\tau}{\rm
Ric}(Y(\bar\tau),Y(\bar\tau))-\int_{\tau_1}^{\bar\tau}\sqrt{\tau}H(X,Y)d\tau.
\nonumber
\end{eqnarray}\index{$H(X,Y)$}

This proves Inequality~(\ref{Hessin}). Now we examine when equality
holds in this expression. Given an adapted vector field $Y(\tau)$
along $\gamma$, let $\mu(v)$ be a geodesic through
$\gamma(\bar\tau,0)$ with tangent vector $Y(\bar\tau)$. Then there
is a one-parameter family $\mu(\tau,v)$ of minimizing ${\mathcal
L}$-geodesics with the property that $\mu(\bar\tau,v)=\mu(v)$. Let
$\widetilde Y'(\tau,v)$ be $\partial\mu(\tau,v)/\partial v$. It is
an ${\mathcal L}$-Jacobi field with $\widetilde
Y'(\bar\tau,0)=Y(\bar\tau)$. Since $L_x\circ {\mathcal L}{\rm
exp}_x=\widetilde L$, we see that
$$\frac{d^2}{d v^2}{\mathcal
L}(\mu_v)|_{v=0}=\frac{d^2}{du^2}L^{\bar \tau}_x(\mu(u))|_{u=0}.$$
Hence, the assumption that we have equality in~(\ref{Hessin})
implies that
$$\frac{d^2}{d v^2}{\mathcal
L}(\mu_v)|_{v=0}=\frac{d^2}{d u^2}{\mathcal L}(\tilde\gamma_u)|_{u=0}.$$

 Now we extend this
one-parameter family to a two-parameter family $\mu(\tau,u,v)$ so
that $\partial\mu(\tau,0,0)/\partial v=\widetilde Y'$ and
$\partial\mu(\tau,0,0)/\partial u=Y(\tau)$. Let $w$ be the variable
$u-v$, and let $\widetilde W$ be the tangent vector in this
coordinate direction, so that $\widetilde W=\widetilde Y-\widetilde
Y'$. We denote by $W$ the restriction of $\widetilde W$ to
$\gamma_{0,0}=\gamma_Z$. By Remark~\ref{variremark} the second
partial derivative of the length of this family in the $u$-direction
at $u=v=0$ agrees with the second derivative of the length of the
original family $\widetilde \gamma$ in the $u$-direction.

\begin{claim}
$$\frac{\partial}{\partial v}\frac{\partial}{\partial w}{\mathcal L}(\mu)|_{u=v=0}=
\frac{\partial}{\partial w}\frac{\partial}{\partial v}{\mathcal
L}(\mu)|_{v=w=0}=0.$$
\end{claim}

\begin{proof}
Of course, the second partial derivatives are equal. According to
Lemma~\ref{2ndvari} we have
$$\frac{\partial}{\partial v}\frac{\partial}{\partial w}{\mathcal L}(\mu)|_{v=w=0}
=2\sqrt{\bar\tau} \widetilde Y'(\bar\tau)\langle \widetilde
W(\bar\tau),X(\bar\tau)\rangle-\int_{\tau_1}^{\bar\tau}2\sqrt{\tau}\langle
{\rm Jac}(\widetilde Y'),W\rangle d\tau.$$ Since $W(\bar\tau)=0$ and
since $\nabla_{\widetilde Y'}(\widetilde W)=\nabla_W(\widetilde
Y')$, we see that the boundary term in the above expression
vanishes. The integral vanishes since $\widetilde Y'$ is an
${\mathcal L}$-Jacobi field.
\end{proof}

If Inequality~(\ref{Hessin}) is an equality, then
$$\frac{\partial^2}{\partial v^2}{\mathcal L}(\mu)|_{u=v=0}=\frac{\partial^2}{\partial u^2}{\mathcal
L}(\mu)|_{u=v=0}.$$ We write $\partial/\partial u=\partial/\partial
v+\partial/\partial w$. Expanding out the right-hand side and
canceling the common terms gives
$$0=\left(\frac{\partial}{\partial v}\frac{\partial}{\partial w}+
\frac{\partial}{\partial w}\frac{\partial}{\partial
v}+\frac{\partial^2}{\partial w^2}\right){\mathcal
L}(\mu)|_{u=v=0}.$$ The previous claim tells us that the first two
terms on the right-hand side of this equation vanish, and hence we
conclude
$$\frac{\partial^2}{\partial w^2}{\mathcal L}(\mu)|_{u=v=0}=0$$
Since $W$ vanishes at both endpoints this implies, according to
Proposition~\ref{2ndvarj}, that $\widetilde W(\tau,0,0)=0$ for all
$\tau$, or in other words $Y(\tau)=\widetilde Y'(\tau,0,0)$ for all
$\tau$. Of course by construction $\widetilde Y'(\tau,0,0)$ is an
${\mathcal L}$-Jacobi field. This shows that equality holds only if
the adapted vector field $Y(\tau)$ is also an ${\mathcal L}$-Jacobi
field.

Conversely, if the adapted vector field $Y(\tau)$ is also an
${\mathcal L}$-Jacobi field, then inequality between the second
variations at the beginning of the proof is an equality. In the rest
of the argument we dealt only with equalities. Hence, in this case
Inequality~(\ref{Hessin}) is an equality.

This shows that we have equality in~(\ref{Hessin}) if and only if
the adapted vector field $Y(\tau)$ is also an ${\mathcal L}$-Jacobi
field.
\end{proof}

\subsection{Inequalities for $\triangle L_x^{\bar\tau}$}

The inequalities for the Hessian of $L_x^{\bar\tau}$ lead to
inequalities for $\triangle L_x^{\bar\tau}$\index{$\triangle
L_x^{\bar\tau}$} which we establish in this section. Here is the
main result.

\begin{prop}\label{triL} Suppose that $q\in {\mathcal U}_x(\bar\tau)$, that $Z\in \widetilde
{\mathcal U}_x(\bar\tau)$ is the pre-image of $q$ and that
$\gamma_Z$ is the ${\mathcal L}$-geodesic determined by $Z$. Then
\begin{equation}\label{lineq}
 \triangle L_x^{\bar\tau}(q) \leq
\frac{n}{\sqrt{\bar\tau}-\sqrt{\tau_1}} -2\sqrt{\bar\tau}R(q) -
\frac{1}{(\sqrt{\bar\tau}-\sqrt{\tau_1})^2}{\mathcal
K}_{\tau_1}^{\bar\tau}(\gamma_Z),
\end{equation}
where, for any path $\gamma$ parameterized by backward time on the
interval
 $[\tau_1,\bar\tau]$ taking
 value $x$ at $\tau=\tau_1$ we define
$${\mathcal K}_{\tau_1}^{\bar\tau}(\gamma)= \int_{\tau_1}^{\bar\tau}
\sqrt{\tau}(\sqrt{\tau}-\sqrt{\tau_1})^2H(X)d\tau,$$\index{${\mathcal
K}_{\tau_1}^{\bar\tau}$} with
\begin{equation}\label{Hdefnn}\index{$H(X)$|ii}
H(X)= -\frac{\partial R}{\partial \tau} - \frac{1}{\tau}R - 2\langle
\nabla R,X\rangle + 2{\rm Ric}(X,X), \end{equation} where $X$ is the
horizontal projection of $\gamma'(\tau)$. Furthermore,
Inequality~(\ref{lineq}) is an equality if and only if for every
$Y\in T_q(M_{T-\bar\tau})$ the unique adapted vector field $Y(\tau)$
along $\gamma$ satisfying $Y(\bar\tau)=Y$ is an ${\mathcal
L}$-Jacobi field. In this case
$${\rm Ric}+\frac{1}{2\sqrt{\bar\tau}}{\rm Hess}(L_x^{\bar\tau})
=\frac{1}{2\sqrt{\bar\tau}(\sqrt{\bar\tau}-\sqrt{\tau_1})}G(T-\bar\tau).$$
\end{prop}

\begin{proof}
Choose an orthonormal basis $\{Y_{\alpha}\}$ for
$T_q(M_{T-\bar\tau})$. For each $\alpha$, extend $\{Y_{\alpha}\}$ to
an adapted vector field along the $ \mathcal{L}$-geodesic $\gamma_Z$
by solving
$$\nabla_{X}Y_{\alpha} = \frac{1}{2\sqrt{\tau}(\sqrt{\tau}-\sqrt{\tau_1})}Y_{\alpha}
-{\rm Ric}(Y_{\alpha},\cdot)^*.$$ As in Equation~(\ref{secondYeqn}),
we have
\begin{align*}
\frac{d}{d\tau}\langle Y_{\alpha},Y_{\beta}\rangle  &= \langle
\nabla_{X}Y_{\alpha},Y_{\beta}\rangle  + \langle
\nabla_{X}Y_{\beta}, Y_{\alpha}\rangle + 2{\rm
Ric}(Y_{\alpha},Y_{\beta})
\\&= \frac{1}{\sqrt{\tau}(\sqrt{\tau}-\sqrt{\tau_1})}\langle Y_{\alpha},Y_{\beta}\rangle .
\end{align*}
 By integrating we get $$\langle Y_{\alpha},Y_{\beta}\rangle (\tau) =
\frac{(\sqrt{\tau}-\sqrt{\tau_1})^2}{(\sqrt{\bar\tau}-\sqrt{\tau_1})^2}\delta_{\alpha\beta}.$$
To simplify the notation we set
$$I(\tau)=\frac{\sqrt{\bar\tau}-\sqrt{\tau_1}}{\sqrt{\tau}-\sqrt{\tau_1}}$$
and $W_\alpha(\tau)=I(\tau)Y_\alpha(\tau)$.
 Then $\{W_{\alpha}(\tau)\}_\alpha$
form an orthonormal basis at $\tau$. Consequently, summing
Inequality~(\ref{Hessineq1}) over $\alpha$ gives
\begin{equation}\label{triLinequ}
\triangle L_x^{\bar\tau}(q)\le
\frac{n}{\sqrt{\bar\tau}-\sqrt{\tau_1}}
-2\sqrt{\bar\tau}R(q)-\sum_\alpha\int_{\tau_1}^{\bar\tau}\sqrt{\tau}H(X,Y_\alpha)d\tau.\end{equation}\index{$H(X,Y)$}
To establish Inequality~(\ref{lineq}) it remains to prove the
following claim.

\begin{claim}\label{Hclaim}
$$\sum_\alpha
H(X,Y_\alpha)=\frac{(\sqrt{\tau}-\sqrt{\tau_1})^2}{(\sqrt{\bar\tau}-\sqrt{\tau_1})^2}H(X).$$\index{$H(X,Y)$}\index{$H(X)$}
\end{claim}

\begin{proof}
To prove the claim we sum Equation~(\ref{Heqn}) giving
\begin{eqnarray*}
I^2(\tau) \sum_{\alpha}H(X,Y_{\alpha}) &= & \sum_\alpha
H(X,W_\alpha)
\\
& = & -\triangle R + 2{\rm Ric}(X,X) - 4\langle \nabla R,X\rangle
+4\sum_{\alpha}(\nabla_{W_\alpha} {\rm Ric})(W_{\alpha},X) \\ & &
-2\sum_{\alpha}{\rm Ric}_{\tau}( W_{\alpha},W_{\alpha})  +2|{\rm
Ric}|^{2} -\frac{1}{\tau}R.
\end{eqnarray*}\index{$H(X,Y)$}

Taking the trace of the second Bianchi identity, we get
$$\sum_{\alpha}(\nabla_{W_\alpha}
{\rm Ric})(W_{\alpha},X) =\frac{1}{2}\langle \nabla R,X\rangle .$$ In addition
by~(\ref{Revol}), recalling that $\partial R/\partial\tau=-\partial R/\partial
t$, we have
$$\frac{\partial R}{\partial \tau} = -\triangle R - 2|{\rm Ric}|^{2}.$$
On the other hand, $$\frac{\partial R}{\partial\tau} =
\partial(g^{ij}R_{ij})/\partial{\tau}= -2|{\rm Ric}|^{2} + \sum_{\alpha}\frac{\partial{\rm
Ric}}{\partial\tau}(W_{\alpha}, W_{\alpha}),$$ and so
$\sum_{\alpha}\frac{\partial{\rm
Ric}}{\partial\tau}(W_{\alpha},W_{\alpha}) = -\triangle R$. Putting
all this together gives\index{$H(X,Y)$}\index{$H(X)$}
$$I^2(\tau)\sum_{\alpha}H(X,Y_{\alpha})= H(X).$$
\end{proof}

Clearly, Inequality~(\ref{lineq}) follows immediately from
Inequality~(\ref{triLinequ}) and the claim. The last statement of
Proposition~\ref{triL} follows directly from the last statement of
Proposition~\ref{Hessineq} and Lemma~\ref{adapjab}. This completes
the proof of Proposition~\ref{triL}.
\end{proof}

\section{Reduced length}\label{sect:redln}

We introduce the reduced length function both on the tangent space
and on space-time. The reason that the reduced length $l_x$ is
easier to work with is that it is scale invariant when $\tau_1=0$.
Throughout Section~\ref{sect:redln} we fix $x\in{\mathcal M}$ with
${\bf t}(x)=T-\tau_1$. We shall always suppose that $T-\tau_1$ is
greater than the initial time of the generalized Ricci flow.

\subsection{The reduced length function $l_x$ on space-time}

\begin{defn} We define
 the $\mathcal{L}$-{\sl reduced length}\index{${\mathcal L}$-length!reduced} (from $x$)
$$l_x\colon {\mathcal U}_x\to \Ar$$\index{$l_x$|ii}
by setting $$l_x(q)= \frac{L_x(q)}{2\sqrt{\tau}},$$ where $\tau=T-{\bf t}(q)$.
We denote by $l_x^\tau$ the restriction of $l_x$ to the slice ${\mathcal
U}_x(\tau)$.
\end{defn}

In order to understand the differential inequalities that $l_x$
satisfies, we first need to introduce a quantity closely related to
the function ${\mathcal K}^{\bar\tau}_{\tau_1}$ defined in
Proposition~\ref{triL}.

\begin{defn}
For any ${\mathcal L}$-geodesic $\gamma$ parameterized by $[\tau_1,\bar\tau]$
we define
$$K_{\tau_1}^{\bar\tau}(\gamma)=
\int_{\tau_1}^{\bar\tau}\tau^{3/2}H(X)d\tau.$$\index{$K_{\tau_1}^{\bar\tau}$|ii}\index{$H(X)$}
In the special case when $\tau_1=0$ we denote this integral by
$K^{\bar\tau}(\gamma)$\index{$K^{\bar\tau}$|ii}.
\end{defn}

The following is immediate from the definitions.

\begin{lem}\label{Kequal}
For any ${\mathcal L}$-geodesic $\gamma$ defined on $[0,\bar\tau]$ both
$K^{\bar\tau}_{\tau_1}(\gamma)$ and ${\mathcal K}^{\bar\tau}_{\tau_1}(\gamma)$
are continuous in $\tau_1$ and at $\tau_1=0$ they take the same value. Also,
$$\left(\frac{\tau_1}{\bar\tau}\right)^{3/2}\left(R(\gamma(\tau_1))+|X(\tau_1)|^2\right)$$
is continuous for all $\tau_1> 0$ and has limit $0$ as
$\tau_1\rightarrow 0$. Here, as always, $X(\tau_1)$ is the
horizontal component of $\gamma'$ at $\tau=\tau_1$.
\end{lem}

\begin{lem}\label{Keqn} Let $q\in {\mathcal
U}_x(\bar\tau)$, let $Z\in \widetilde {\mathcal U}_x$ be the
pre-image of $q$ and let $\gamma_Z$ be the ${\mathcal L}$-geodesic
determined by $Z$. Then we have
\begin{equation}\label{Kequat}
\bar\tau^{-\frac{3}{2}}K^{\bar\tau}_{\tau_1}(\gamma_Z)=\frac{
l_x(q)}{\bar\tau}-(R(q)
+|X(\bar\tau)|^2)+\left(\frac{\tau_1}{\bar\tau}\right)^{3/2}\left(R(x)+|X(\tau_1)|^2\right).
\end{equation}\index{$K_{\tau_1}^{\bar\tau}$|(}

In the case when $\tau_1=0$, the last term on the right-hand side of
Equation~(\ref{Kequat}) vanishes.
\end{lem}

\begin{proof}
Using the ${\mathcal L}$-geodesic equation and the definition of $H$ we have
\begin{eqnarray*}
\lefteqn{\frac{d}{d\tau}(R(\gamma_Z(\tau))+ \abs{X(\tau)}^{2})} \\ &
= & \frac{\partial R}{\partial\tau}(\gamma_Z(\tau)) + \langle \nabla
R(\gamma_Z(\tau)),
X(\tau)\rangle + 2\langle \nabla_{X}X(\tau),X(\tau)\rangle \\
& & + 2{\rm Ric}(X(\tau),X(\tau))
\\&= & \frac{\partial R}{\partial\tau}(\gamma_Z(\tau)) + 2X(\tau)(R) - \frac{1}{\tau}\abs{X(\tau)}^{2} -
2{\rm Ric}(X(\tau),X(\tau)) \\&= &-H(X(\tau))
-\frac{1}{\tau}(R(\gamma_Z(\tau)+ \abs{X(\tau)}^{2}).
\end{eqnarray*}
Thus
$$\frac{d}{d\tau}(\tau^{\frac{3}{2}}(R(\gamma_Z(\tau)+\abs{X(\tau)}^{2}) )=
\frac{1}{2}\sqrt{\tau}(R(\gamma_Z(\tau)+\abs{X(\tau)}^{2}) -
\tau^{\frac{3}{2}}H(X(\tau)).$$ Integration from $\tau_1$ to $\bar\tau$ gives
$$\bar\tau^{3/2}\left(R(q))+|X(\bar\tau)|^2\right)-\tau_1^{3/2}(R(x)+|X(\tau_1)|^2)
=\frac{L_x^{\bar\tau}(q)}{2}-K_{\tau_1}^{\bar\tau}(\gamma_Z),$$
which is equivalent to Equation~(\ref{Kequat}). In the case when
$\tau_1=0$, the last term on the right-hand side vanishes since
$${\rm lim}_{\tau\rightarrow 0}\tau^{3/2}|X(\tau)|^2=0.$$
\end{proof}

Now we come to the most general of the differential inequalities for
$l_x$ that will be so important in what follows. Whenever the
expression
$\left(\frac{\tau_1}{\bar\tau}\right)^{3/2}\left(R(x)+|X(\tau_1)|^2\right)$
appears in a formula, it is interpreted to be zero in the case when
$\tau_1=0$.

\begin{lem}\label{linequal}
For any  $q\in{\mathcal U}_x(\bar\tau)$, let $Z\in
\widetilde{\mathcal U}_x(\bar\tau)$ be the pre-image of $q$ and let
$\gamma_Z$ be the ${\mathcal L}$-geodesic determined by $Z$. Then we
have
\begin{eqnarray*} \frac{\partial l_x}{\partial \tau}(q) & = &
R(q)-\frac{l_x(q)}{\bar\tau}+\frac{K_{\tau_1}^{\bar\tau}(\gamma_Z)}{2\bar\tau^{3/2}}-\frac{1}{2}\left(\frac{\tau_1}{\bar\tau}
\right)^{3/2}\left(R(x)+|X(\tau_1)|^2\right) \\
|\nabla l^{\bar\tau}_x(q)|^2 & = & |X(\bar\tau)|^2
=\frac{l^{\bar\tau}_x(q)}{\bar\tau}-\frac{K_{\tau_1}^{\bar\tau}
(\gamma_Z)}{\bar\tau^{3/2}}-R(q) +\left(\frac{\tau_1}{\bar\tau}
\right)^{3/2}\left(R(x)+|X(\tau_1)|^2\right) \\
\triangle l^{\bar\tau}_x(q) & = & \frac{1}{2\sqrt{\bar\tau}}\triangle
L^{\bar\tau}_x(q)\le \frac{n}{2\sqrt{\bar\tau}(\sqrt{\bar\tau}-\sqrt{\tau_1})}
-R(q)-\frac{{\mathcal
K}_{\tau_1}^{\bar\tau}(\gamma_Z)}{2\sqrt{\bar\tau}(\sqrt{\bar\tau}-\sqrt{\tau_1})^{2}}.
\end{eqnarray*}\index{${\mathcal K}_{\tau_1}^{\bar\tau}$}
\end{lem}

\begin{proof}
 It follows immediately from Equation~(\ref{Ltaueqn}) that
$$\frac{\partial l_x}{\partial \tau}=R-\frac{1}{2}(R+|X|^2)-\frac{l_x}{2\tau}.$$
Using Equation~(\ref{Kequat}) this gives the first equality stated
in the lemma.
 It follows immediately from Corollary~\ref{DL} that
$\nabla l_x^\tau=X(\tau)$ and hence $|\nabla
l_x^\tau|^2=|X(\tau)|^2$. From this and Equation~(\ref{Kequat}) the
second equation follows. The last inequality is immediate from
Proposition~\ref{triL}.

When $\tau_1=0$, the last terms on the right-hand sides of the first
two equations vanish, since the last term on the right-hand side of
Equation~(\ref{Kequat}) vanishes in this case.
\end{proof}

When $\tau_1=0$, which is the case of main interest, all these
formulas simplify and we get:

\begin{thm}\label{Ueqn}
Suppose that $x\in M_T$ so that $\tau_1=0$. For any $q\in{\mathcal
U}_x(\bar\tau)$, let $Z\in \widetilde{\mathcal U}_x(\bar\tau)$ be
the pre-image of $q$ and let $\gamma_Z$ be the ${\mathcal
L}$-geodesic determined by $Z$. As usual, let $X(\tau)$ be the
horizontal projection of $\gamma_Z'(\tau)$. Then we have
\begin{eqnarray*} \frac{\partial l_x}{\partial \tau}(q) & = &
R(q)-\frac{l_x(q)}{\bar\tau}+\frac{K^{\bar\tau}(\gamma_Z)}{2\bar\tau^{3/2}}
\\
|\nabla l_x^{\bar\tau}(q)|^2 & = & |X(\bar\tau)|^2  =
\frac{l_x^{\bar\tau}(q)}{\bar\tau}-\frac{K^{\bar\tau}
(\gamma_Z)}{\bar\tau^{3/2}}-R(q)  \\
\triangle l_x^{\bar\tau}(q) & = &
\frac{1}{2\sqrt{\bar\tau}}\triangle L_x^{\bar\tau}(q)\le \frac{n}{2
\bar\tau}-R(q)-\frac{K^{\bar\tau}(\gamma_Z)}{2\bar\tau^{3/2}}.
\end{eqnarray*}
\end{thm}

\begin{proof}
This is immediate from the formulas in the previous lemma.
\end{proof}

Now let us reformulate the differential inequalities in
Theorem~\ref{Ueqn} in a way that will be useful later.

\begin{cor}\label{lformula}
Suppose that $x\in M_T$ so that $\tau_1=0$. Then for $q\in{\mathcal
U}^{\bar\tau}_x$ we have
$$\frac{\partial l_x}{\partial \tau}(q)+\triangle l_x^{\bar\tau}(q)\le \frac{(n/2)-l_x^{\bar\tau}(q)}{\bar\tau}.$$
$$\frac{\partial l_x}{\partial \tau}(q)-\triangle l_x^{\bar\tau}(q)+|\nabla l_x^{\bar\tau}(q)|^2-R(q)+\frac{n}{2\bar\tau}\ge 0 .$$
$$2\triangle l_x^{\bar\tau}(q)-|\nabla l_x^{\bar\tau}(q)|^2+R(q)+\frac{l_x^{\bar\tau}(q)-n}{\bar\tau}
\le 0.$$

In fact, setting $$\delta=\frac{n}{2
\bar\tau}-R(q)-\frac{K^{\bar\tau}(\gamma_Z)}{2\bar\tau^{3/2}}-\triangle
l_x^{\bar\tau}(q),$$ then $\delta\ge 0$ and
$$\frac{\partial l_x}{\partial \tau}(q)-\triangle l_x^{\bar\tau}(q)
+|\nabla l_x^{\bar\tau}(q)|^2-R(q)+\frac{n}{2\bar\tau}=\delta$$
$$2\triangle l_x^{\bar\tau}(q)-|\nabla l_x^{\bar\tau}(q)|^2+R(q)
+\frac{l_x^{\bar\tau}(q)-n}{\bar\tau}=-2\delta.$$
\end{cor}\index{$K_{\tau_1}^{\bar\tau}$|)}

\subsection{The tangential version $\widetilde l$ of the reduced length function}

For any path $\gamma\colon [\tau_1,\bar\tau]\to ({\mathcal M},G)$
parameterized by backward time we define
$$l(\gamma)=\frac{1}{2\sqrt{\bar\tau}}{\mathcal L}(\gamma).$$
This leads immediately to a reduced length on $\widetilde {\mathcal
U}_x$.

\begin{defn}
We define $\widetilde l\colon \widetilde {\mathcal U}_x\to
\Ar$\index{$\widetilde l$}  by
$$\widetilde l(Z,\tau)=\frac{\widetilde L(Z,\tau)}{2\sqrt{\tau}}=l(\gamma_Z|_{[\tau_1,\bar\tau]}).$$
\end{defn}

At first glance it may appear that the computations of the gradient
and $\tau$-derivatives for $l_x$ pass immediately to those for
$\widetilde l$. For the spatial derivative this is correct, but for
the $\tau$-derivative it is not true. As the computation below
shows, the $\tau$-derivatives of $\widetilde l$ and $l_x$ do not
agree under the identification ${\mathcal L}{\rm exp}_x$. The reason
is that this identification does not line up the $\tau$-vector field
in the domain with $-\partial/\partial t$ in the range. So it is an
entirely different computation with a different answer.

\begin{lem}\label{tildeltau}
$$\frac{\partial \widetilde
l(Z,\tau)}{\partial
\tau}=\frac{1}{2}\left(R(\gamma_Z(\tau))+X(\tau)|^2\right)-\frac{\widetilde
l(Z,\tau)}{2\tau}.$$
\end{lem}

\begin{proof}
By the Fundamental Theorem of Calculus
$$\frac{\partial}{\partial\tau}\widetilde L(Z,\tau)=\sqrt{\tau}\left(R(\gamma_z(\tau))+|X(\tau)|^2\right).$$
Thus,
$$\frac{\partial}{\partial\tau}\widetilde l(Z,\tau)=\frac{1}{2}\left(R(\gamma_z(\tau))+|X(\tau)|^2\right)
-\frac{\widetilde l(Z,\tau)}{2\tau}.$$
\end{proof}

\begin{cor}\label{tildelK} Suppose that $x\in M_T$ so that
$\tau_1=0$. Then
$$\frac{\partial}{\partial\tau}\widetilde
l(Z,\tau)=-\frac{K^\tau(\gamma_Z)}{2\tau^{\frac{3}{2}}}.$$
\end{cor}

\begin{proof}
This is immediate from Lemma~\ref{tildeltau} and Lemma~\ref{Keqn} (after the
latter is rewritten using $\widetilde L$ instead of $L_x$).
\end{proof}

\section{Local Lipschitz estimates for $l_x$}\label{sect:lips}

It is important for the applications to have results on the
Lipschitz properties of $l_x$, or equivalently $L_x$. Of course,
these are the analogues of the fact that in Riemannian geometry the
distance function from a point is Lipschitz. The proof of the
Lipschitz property given here is based on the exposition in
\cite{Ye}. In Section~\ref{sect:lips}, we fix $x\in
M_{T-\tau_1}\subset {\mathcal M}$.

\subsection{Statement and corollaries}

\begin{defn}
Let $({\mathcal M},G)$ be a generalized Ricci flow and let $x\in
M_{T-\tau_1}\subset {\mathcal M}$. The reduced length function $l_x$ is defined
on the subset of ${\mathcal M}$ consisting of all points $y\in {\mathcal M}$
for which there is a minimizing ${\mathcal L}$-geodesic from $x$ to $y$. The
value $l_x(y)$ is the quotient of ${\mathcal L}$-length of any such minimizing
${\mathcal L}$-geodesic divided by $2\sqrt{\tau}$.
\end{defn}

Here is the main result of this subsection.

\begin{prop}\label{lips}
Let $({\mathcal M},G)$ be a generalized Ricci flow and let $x\in
M_{T-\tau_1}\subset {\mathcal M}$. Let $\epsilon>0$ be given and let $A\subset
{\mathcal M}\cap {\bf t}^{-1}(-\infty,T-\tau_1+\epsilon)$. Suppose that there
is a subset $F\subset {\mathcal M}$ on which $|{\rm Ric}|$ and $|\nabla R|$ are
bounded and a neighborhood $\nu(A)$ of $A$ contained in $F$ with the property
that for every point $z\in \nu(A)$ there is a minimizing ${\mathcal
L}$-geodesic from $x$ to $z$ contained in $F$. Then $l_x$ is defined on  all of
$\nu(A)$. Furthermore, there is a smaller neighborhood $\nu_0(A)\subset \nu(A)$
of $A$ on which $l_x$ is a locally Lipschitz function with respect to the
Riemannian metric, denoted $\widehat G$, on ${\mathcal M}$ which is defined as
the orthogonal sum of the Riemannian metric $G$ on ${\mathcal H}T{\mathcal M}$
and the metric $dt^2$ on the tangent line spanned by $\chi$.
\end{prop}

\begin{cor}
With $A$ and $\nu_0(A)$ as in Proposition~\ref{lips}, the
restriction of $l_x$ to  $\nu_0(A)\cap M_{T-\bar\tau}$ is a locally
Lipschitz function with respect to the metric $G_{T-\bar\tau}$.
\end{cor}

\subsection{The proof of Proposition~\protect{\ref{lips}}}

Proposition~\ref{lips} follows from a much more precise, though more
complicated to state, result. In order to state this more technical
result we introduce the following definition.

\begin{defn}
Let $y\in {\mathcal M}$ with ${\bf t}(y)=t$ and suppose that for
some $\epsilon>0$ there is an embedding $\iota\colon B(y,t,r)\times
(t-\epsilon,t+\epsilon)\to {\mathcal M}$ that is compatible with
time and the vector field. Then we denote by $\widetilde
P(y,r,\epsilon)\subset {\mathcal M}$ the image of $\iota$. Whenever
we introduce $\widetilde P(y,r,\epsilon)\subset {\mathcal M}$
implicitly we are asserting that such an embedding exists.

For  $A\subset {\mathcal M}$, if $\widetilde
P(a,\epsilon,\epsilon)\subset {\mathcal M}$ exists for every $a\in
A$, then we denote by $\nu_\epsilon(A)$ the union over all $a\in A$
of $\widetilde P(a,\epsilon,\epsilon)$.
\end{defn}

Now we are ready for the more precise, technical result.

\begin{prop}\label{lipaty}
Given constants  $\epsilon>0$, $\bar\tau_0<\infty$, $l_0<\infty$,
and $C_0<\infty$, there are constants $ C<\infty$ and
$0<\delta<\epsilon$ depending only the  given constants  such that
the following holds. Let $({\mathcal M},G)$ be a generalized Ricci
flow and let $x\in{\mathcal M}$ be a point with ${\bf
t}(x)=T-\tau_1$. Let $y\in {\mathcal M}$ be a point with ${\bf
t}(y)=t=T-\bar\tau$ where $\tau_1+\epsilon\le \bar\tau\le
\bar\tau_0$. Suppose that there is a minimizing ${\mathcal
L}$-geodesic $\gamma$ from $x$ to $y$ with $l(\gamma)\le l_0$.
Suppose that the ball $B(y,t,\epsilon)$ has compact closure in
$M_{t}$ and that $\widetilde P(y,\epsilon,\epsilon)\subset {\mathcal
M}$ exists and that the sectional curvatures of the restriction of
$G$ to this submanifold are bounded by $C_0$. Lastly, suppose that
for every point of the form $z\in \widetilde P(y,\delta,\delta)$
there is a minimizing ${\mathcal L}$-geodesic from $x$ to $z$ with
$|{\rm Ric}|$ and $|\nabla R|$ bounded by $C_0$ along this geodesic.
Then for all $(b,t')\in B(y,t,\delta)\times(t-\delta,t+\delta)$ we
have
$$|l_x(y)-l_x(\iota(b,t'))|\le C\sqrt{d_{t}(y,b)^2+|t-t'|^2}.$$
\end{prop}

Before proving Proposition~\ref{lipaty}, let us show how it implies
Proposition~\ref{lips}.

\begin{proof} (that Proposition~\ref{lipaty} implies
Proposition~\ref{lips}) Suppose given $\epsilon>0$, $A$, $\nu(A)$ and $F$ as in
the statement of Proposition~\ref{lips}. For each $y\in A$ there is
$0<\epsilon'<\epsilon$ and a neighborhood $\nu'(y)$ with (i) the closure
$\overline{\nu}'(y)$ of $\nu'(y)$ being a compact subset of $\nu(A)$ and (ii)
for each $z\in \overline{\nu}'(y)$ the parabolic neighborhood $\widetilde
P(z,\epsilon',\epsilon')$ exists and has compact closure in $\nu(A)$. It
follows that for every $z\in \overline{\nu}'(A)$, ${\rm Rm}_G$ is bounded on
$\widetilde P(z,\epsilon',\epsilon')$ and every point of $\widetilde
P(z,\epsilon',\epsilon')$ is connected to $x$ by a minimizing ${\mathcal
L}$-geodesic in $F$. Thus, Proposition~\ref{lipaty}, with $\epsilon$ replaced
by $\epsilon'$, applies to $z$. In particular, $l_x$ is continuous at $z$, and
hence is continuous on all of $\overline{\nu}'(y)$. Thus, $l_x$ is bounded on
$\overline{\nu}'(y)$. Since we have uniform bounds for the curvature on
$\widetilde P(z,\epsilon',\epsilon')$ according to Proposition~\ref{lipaty}
there are constants $C<\infty$ and $0<\delta<\epsilon'$ such that for any $z\in
\overline{\nu}'(y)$ and any $z'\in \widetilde P(z,\delta,\delta)$, we have
$$|l_x(z)-l_x(z')|\le C|z-z'|_{G({\bf t}(z)+dt^2}.$$
Since we have a uniform bound for the curvature on $\widetilde
P(z,\epsilon',\epsilon')$ independent of $z\in \nu'(y)$, the metrics $\widehat
G=G+dt^2$ and $G({\bf t}(z))+dt^2$ are uniformly comparable on all of
$\widetilde P(z,\delta,\delta)$. It follows that there is a constant
$C'<\infty$ such that for all $z\in \nu'(y)$ and all $z'\in \widetilde
P(z,\delta,\delta)$ we have
$$|l_x(z)-l_x(z')|\le C'|Z-z'|_{\widehat G}.$$

We set $\nu_0(A)=\cup_{y\in A}\nu'(y)$. This is an open neighborhood of $A$
contained in $\nu(A)$ on which $l_x$ is locally Lipschitz with respect to the
metric $\widehat G$.
\end{proof}

Now we turn to the proof of Proposition~\ref{lipaty}. We begin with
several preliminary results.

\begin{lem}\label{min-max}
 Suppose  that
$\gamma$ is an ${\mathcal L}$-geodesic defined on
$[\tau_1,\bar\tau]$, and suppose that for all $\tau\in
[\tau_1,\bar\tau]$ we have $|\nabla R(\gamma(\tau))|\le C_0$ and
$|{\rm Ric}(\gamma(\tau))|\le C_0$. Then
$${\rm max}_{\tau}\left(\sqrt{\tau}|X_\gamma(\tau)|\right)\le C_1{\rm
min}_{\tau}\left(\sqrt{\tau}|X_\gamma(\tau)|\right)+\frac{(C_1-1)}{2}\sqrt{\bar\tau},$$
where $C_1=e^{2C_0\bar\tau}$.
\end{lem}

\begin{proof}
The geodesic equation in terms of the variable $s$,
Equation~(\ref{reparam}), gives
\begin{eqnarray}
\frac{d|\gamma'(s)|^2}{ds} & = & 2\langle
\nabla_{\gamma'(s)}\gamma'(s),\gamma'(s)\rangle+4s{\rm Ric}(\gamma'(s),\gamma'(s))\nonumber \\
& = & 4s^2\langle\nabla R,\gamma'(s)\rangle -4s{\rm
Ric}(\gamma'(s),\gamma'(s))\label{nablas}.
\end{eqnarray}

Thus,, by our assumption on $|\nabla R|$ and $|{\rm Ric}|$ along
$\gamma$, we have
$$\left|\frac{d|\gamma'(s)|^2}{ds}\right| \le
4C_0s^2|\gamma'(s)|+4C_0s|\gamma'(s)|^2.$$ It follows that
$$\left|\frac{d|\gamma'(s)|}{ds}\right|\le 2C_0s^2+2C_0s|\gamma'(s)|\le 2C_0\bar\tau+2
C_0\sqrt{\bar\tau}|\gamma'(s)|,$$ and hence that
$$-2C_0\sqrt{\bar\tau}ds\le \frac{d|\gamma'(s)|}{\sqrt{\bar\tau}+|\gamma'(s)|}
\le 2C_0\sqrt{\bar\tau}ds.$$ Suppose that $s_0<s_1$. Integrating from $s_0$ to
$s_1$ gives
\begin{eqnarray*}
|\gamma'(s_1)| & \le &  C|\gamma'(s_0)|+(C-1)\sqrt{\bar\tau} \\
|\gamma'(s_0)| & \le & C|\gamma'(s_1)|+(C-1)\sqrt{\bar\tau}
\end{eqnarray*}
where
\begin{equation*}
C=e^{2C_0\sqrt{\bar\tau}(s_1-s_0)}.
\end{equation*}
 Since
$\sqrt{\tau}X_\gamma(\tau)=\frac{1}{2}\gamma'(s)$, this completes
the proof of the lemma.
\end{proof}

\begin{cor}\label{minmax}
Given $\bar\tau_0<\infty$, $C_0<\infty$, $\epsilon>0$, and $l_0<\infty$, there
is a constant $C_2$ depending only on $C_0$, $l_0$, $\epsilon$ and $\bar\tau_0$
such that the following holds. Let $\gamma$ be an ${\mathcal L}$-geodesic
defined on $[\tau_1,\bar\tau]$ with $\tau_1+\epsilon\le\bar\tau\le \bar\tau_0$
and with $|\nabla R(\gamma(\tau))|\le C_0$ and $|{\rm Ric}(\gamma(\tau))|\le
C_0$ for all $\tau\in [\tau_1,\bar\tau]$. Suppose also that $ l(\gamma)\le
l_0$. Then, we have
$${\rm max}_{\tau}\left(\sqrt{\tau}|X_\gamma(\tau)|\right)\le C_2.$$
\end{cor}

\begin{proof} From the definition
${\mathcal L}(\gamma)=\int_{\sqrt{\tau_1}}^{\sqrt{\bar\tau}}
(2s^2R+\frac{1}{2}|\gamma'(s)|^2)ds$. Because of the bound on $|{\rm
Ric}|$ (which implies that $|R|\le 3C_0$) we have
$$\frac{1}{2}\int_{\sqrt{\tau_1}}^{\sqrt{\bar\tau}}|\gamma'(s)|^2ds\le {\mathcal L}(\gamma)+
2C_0\bar\tau^{3/2}.$$ Thus,
$$(\sqrt{\bar \tau}-\sqrt{\tau_1}){\rm min}(|\gamma'(s)|^2)\le 2{\mathcal L}(\gamma)+
4C_0\bar\tau^{3/2}.$$ The bounds $\tau_1+\epsilon\le \bar\tau\le\bar\tau_0$,
then imply that ${\rm min}|\gamma'(s)|^2\le C''$ for some $C''$ depending on
$C_0,l_0,\epsilon,$ and $\bar\tau_0$. Since
$\sqrt{\tau}X_\gamma(\tau)=\frac{1}{2}\gamma'(s)$, we have
$${\rm min}_{\tau}\left(\sqrt{\tau}|X_\gamma(\tau)|\right)\le C'$$ for some constant $C'$ depending
only on $C_0$, $l_0$, $\epsilon$ and $\bar\tau_0$. The result is now immediate
from Lemma~\ref{min-max}.
\end{proof}

Now we are ready to show that, for $z$ sufficiently close to  $y$,
the reduced length $l_x(z)$ is bounded above by a constant depending
on the curvature bounds, on $l_x(y)$, and on the distance in
space-time from $z$ to $y$.

\begin{lem}\label{upperbd}
Given constants  $\epsilon>0$, $\bar\tau_0<\infty$, $C_0<\infty$,
and $l_0<\infty$, there are $C_3<\infty$ and $0<\delta_2\le
\epsilon/4$ depending only on the  given constants  such that the
following holds. Let $y\in {\mathcal M}$ be a point with ${\bf
t}(y)=t_0=T-\bar\tau$ where $\tau_1+\epsilon\le \bar\tau\le
\bar\tau_0$. Suppose that there is a minimizing ${\mathcal
L}$-geodesic $\gamma$ from $x$ to $y$ with $l_x(\gamma)\le l_0$.
Suppose that $|\nabla R|$ and $|{\rm Ric}|$ are bounded by $C_0$
along $\gamma$. Suppose also that the ball $B(y,t_0,\epsilon)$ has
compact closure in $M_{t_0}$ and that there is an embedding
$$\iota\colon B(y,t_0,\epsilon)\times (t_0-\epsilon,t_0+\epsilon)\stackrel
{\cong}{\longrightarrow} \tilde P(y,\epsilon,\epsilon)\subset
{\mathcal M}$$ compatible with  time  and the vector field so that
the sectional curvatures of the restriction of $G$ to the image of
this embedding are bounded by $C_0$. Then for any point $b\in
B(y,t_0,\delta_2)$ and for any $t'\in (t_0-\delta_2,t_0+\delta_2)$
there is a curve $\gamma_1$ from $x$ to the point $z=\iota(b,t')$,
parameterized by backward time, such that
$$l(\gamma_1)\le  l(\gamma)+C_3\sqrt{d_{t_0}(y,b)^2+|t_0-t'|^2}.$$
\end{lem}

\begin{proof}
Let $C_2$ be the constant depending on $C_0$, $l_0$, $\epsilon$, and
$\bar\tau_0$ from Corollary~\ref{minmax}, and set
$$C'=\frac{\sqrt{2}}{\sqrt{\epsilon}}C_2.$$
Since $\bar\tau\ge \epsilon$, it follows that
$\bar\tau-\epsilon/2\ge \epsilon/2$, so that by
Corollary~\ref{minmax} we have $|X_\gamma(\tau)|\le C'$ for all
$\tau\in [\bar\tau-\epsilon/2,\bar\tau]$. Set $0<\delta_0$
sufficiently small (how small depends only on $C_0$) such that for
all $(z,t)\in \tilde P(y,\epsilon,\delta_0)$ we have
$$\frac{1}{2}g(z,t)\le g(z,t_0)\le 2g(z,t),$$ and define
$$\delta_2={\rm
min}\left(\frac{\epsilon}{8},\frac{\epsilon}{8C'},\frac{\delta_0}{4}\right).$$
 Let $b\in B(y,t_0,\delta_2)$ and $t'\in
(t_0-\delta_2,t_0+\delta_2)$ be given. Set
$\alpha=\sqrt{d_{t_0}(y,b)^2+|t_0-t'|^2}$, set $t_1=t_0-2\alpha$,
and let $\tau_1=T-t_1$. Notice that
$\alpha<\sqrt{2}\delta_2<\epsilon/4$, so that the norm of the Ricci
curvature is bounded by $C'$ on
$\iota(B(y,t_0,\epsilon)\times(t_1,t_0+2\alpha))$.

\begin{claim}\label{6.66}
$\gamma(\tau_1)\in \widetilde P(y,\epsilon,\epsilon)$ and writing
$\gamma(\tau_1)=\iota(c,t_1)$ we have $d_{t_0}(c,b)\le
(4C'+1)\alpha$.
\end{claim}

\begin{proof}
Since $|X_\gamma(\tau)|\le C'$ for all $\tau\in
[\bar\tau-2\alpha,\bar\tau]$, and  $\delta_2\le \delta_0/4$, it
follows that $2\alpha\le \delta_0$ and hence that
$|X_\gamma(\tau)|_{g(t_0)}\le 2C'$ for all  $\tau\in
[\bar\tau-2\alpha,\bar\tau]$. Since $\gamma(\bar\tau)=y$, this
implies that
$$d_{t_0}(y,c)\le 4C'\alpha.$$
The claim then follows from the triangle inequality.
\end{proof}

Now let $\overline\mu\colon [\bar\tau-2\alpha,T-t']\to
B(y,t_0,\epsilon)$ be a shortest $g(t_0)$-geodesic from $c$ to $b$,
parameterized at constant $g(t_0)$-speed, and let $\mu$ be the path
parameterized by backward time defined by
$$\mu(\tau)=\iota(\overline\mu(\tau),T-\tau)$$
for all $\tau\in [\bar\tau-2\alpha,T-t']$. Then the concatenation
$\gamma_1=\gamma|_{[\tau_1,\bar\tau-2\alpha]}*\mu$ is a path
parameterized by backward time from $x$ to $\iota(b,t')$.

\begin{claim} There is a constant $C'_1$ depending only on
$C_0$, $C'$, and $\bar\tau_0$ such that
$$l(\gamma_1)\le l(\gamma|_{[\tau_1,\bar\tau-2\alpha]})+C'_1\alpha$$
\end{claim}

\begin{proof}
First notice that since $\bar\tau=T-t_0$ and $|t'-t_0|\le \alpha$ we have
$(T-t')-(\bar\tau-2\alpha)=2\alpha+(t'-t_0)\ge \alpha$. According to
Claim~\ref{6.66} this implies that the $g(t_0)$-speed of $\mu$ is at most
$(4C'+1)$, and hence  that $|X_{\mu}(\tau)|_{g(T-\tau)}\le 8C'+2$ for all
$\tau\in [\bar\tau-2\alpha,T-t']$. Consequently, $R+|X_\mu|^2$ is bounded above
along $\mu$ by a constant $\widetilde C$ depending only on $C'$ and $C_0$. This
implies that ${\mathcal L}(\mu)\le \widetilde C\alpha\sqrt{T-t'}$. Of course,
$T-t'\le \bar\tau+\epsilon<2\bar\tau\le 2\bar\tau_0$. This completes the proof
of the claim.
\end{proof}

On the other hand, since $R\ge -3C_0$ in $P(y,\epsilon,\epsilon)$
and $|X|^2\ge 0$, we see that
$${\mathcal L}(\gamma|_{[\tau_1,\bar\tau-2\alpha]})\le {\mathcal L}(\gamma)
+ 6C_0\alpha\sqrt{\bar\tau_0}.$$  Together with the previous claim
this establishes Lemma~\ref{upperbd}.
\end{proof}

This is a one-sided inequality which says that the nearby values of
$l_x$ are bounded above in terms of  $l_x(y)$, the curvature bounds,
and the distance in space-time from $y$. In order to complete the
proof of Proposition~\ref{lipaty} we must establish inequalities in
the opposite direction. This requires reversing the roles of the
points.

\begin{proof}(of Proposition~\ref{lipaty})
Let $\delta_2$ and $C_3$ be the constants given by
Lemma~\ref{upperbd} associated to $\epsilon/2$, $\bar\tau_0$, $C_0$,
and $l_0$. We shall choose $C\ge C_3$ and $\delta\le \delta_2$ so
that by Lemma~\ref{upperbd} we will automatically have
$$l_x(\iota(b,t'))\le l_x(y)+C_3\sqrt{d_{t_0}(y,b)^2+|t_0-t'|^2}
\le l_x(y)+C\sqrt{d_{t_0}(y,b)^2+|t_0-t'|^2}$$ for all
$\iota(b,t')\in P(y,\delta,\delta)$. It remains to choose $C$ and
$\delta$ so that
$$l_x(y)\le l_x(\iota(b,t'))+C\sqrt{d_{t_0}(y,b)^2+|t_0-t'|^2}.$$
 Let
$\delta'_2$ and $ C'_3$ be the constants given by
Lemma~\ref{upperbd} for the following set of input constants:
$C_0'=C_0$, $\bar\tau_0$ replaced by
$\bar\tau_0'=\bar\tau_0+\epsilon/2$, and $l_0$ replaced by
$l_0'=l_0+\sqrt{2}C_3\delta_2$, and $\epsilon$ replaced by
$\epsilon'=\epsilon/4$. Then set $C={\rm max}(2C_3',C_3)$.

Let $z=\iota(b,t')\in \widetilde P(y,\delta,\delta)$.

\begin{claim}
For $\delta$ sufficiently small (how small depending on $\delta_2$
and $\delta'_2$) we have $B(z,t',\epsilon/4)\subset
B(y,t_0,\epsilon)$.
\end{claim}

\begin{proof}
Since $|t_0-t'|<\delta\le \delta_2$, and by construction
$\delta_2<\delta_0$, it follows  that for any $c\in
B(y,t_0,\epsilon)$ we have $d_{t'}(b,c)\le 2d_{t_0}(b,c)$. Since
$d_{t_0}(y,b)<\delta\le \epsilon/4$, the result is immediate from
the triangle inequality.
\end{proof}

 By
the above and the fact that $ \delta\le\epsilon/4$, the sectional
curvatures on $\widetilde P(z,\epsilon/4,\epsilon/4)$ are bounded by
$C_0$. By Lemma~\ref{upperbd} there is a curve parameterized by
backward time from $x$ to $z$ whose $l$-length is at most $l_0'$.
Thus the $l$-length any minimizing ${\mathcal L}$-geodesic from $x$
to $z$ is at most $l_0'$. By assumption we have a minimizing
${\mathcal L}$-geodesic with the property that  $|{\rm Ric}|$
and$|\nabla R|$ are bounded by $C_0$ along the ${\mathcal
L}$-geodesic.

Of course, $t_0-\delta<t'<t_0+\delta$, so that
$\tau_1+\epsilon/2<T-t'\le \bar\tau_0+\epsilon/4$. This means that
Lemma~\ref{upperbd} applies to show that for very $w=\iota(c,t)\in
\widetilde P(z,\delta_2',\delta'_2)$,  we have
$$l_x(w)\le l_x(z)+C'_3\sqrt{d_{t'}(b,c)^2+|t-t'|^2}.$$

The proof is then completed by showing the following:

\begin{claim}
$y\in \widetilde P(z,\delta_2',\delta'_2)$.
\end{claim}

\begin{proof}
By construction $|t'-t_0|<\delta\le \delta_2'$. Also,
$d_{t_0}(y,b)<\delta\le \delta'_2/2$. Since $d_{t_0}\le 2d_{t'}$, we
have $d_{t'}(y,b)< \delta'_2$ the claim is then immediate.
\end{proof}

It follows immediately that \begin{eqnarray*} l_x(y) & \le &
l_x(z)+C'_3\sqrt{d_{t'}(b,y)^2+|t_0-t'|^2} \\ & \le & l_x(z)+2
C'_3\sqrt{d_{t_0}(b,y)^2+|t_0-t'|^2}\le
l_x(z)+C\sqrt{d_{t_0}(b,y)^2+|t_0-t'|^2}.
\end{eqnarray*}

This completes the proof of Proposition~\ref{lipaty}.
\end{proof}

\begin{cor}\label{fullmeasure}
Let $({\mathcal M},G)$ be a generalized Ricci flow and let $x\in {\mathcal M}$
with ${\bf t}(x)=T-\tau_1$. Let
 $A\subset {\mathcal
M}\cap {\bf t}^{-1}(-\infty,T-\tau_1)$ be a subset whose intersection with each
time-slice $M_t$ is measurable. Suppose that there is a  subset $F\subset
{\mathcal M}$ such that $|\nabla R|$ and $|{\rm Ric}|$ are bounded on $F$ and
such that every minimizing ${\mathcal L}$ geodesic from $x$ to any point in a
neighborhood, $\nu(A)$, of $A$ is contained in $F$. Then for each $\tau\in
(\tau_1,\bar\tau]$ the intersection of $A$ with ${\mathcal U}_x(\tau)$ is an
open subset of full measure in $A\cap M_{T-\tau}$.
\end{cor}

\begin{proof}
Since ${\mathcal U}_x(\tau)$ is an open subset of $M_{T-\tau}$, the complement
of $\nu(A)\cap{\mathcal U}_x(\tau)$ in $\nu(A)\cap M_{T-\tau}$ is a closed
subset of $\nu(A)\cap M_{T-\tau}$. Since there is a minimizing ${\mathcal
L}$-geodesic to every point of $\nu(A)\cap M_{T-\tau}$, the ${\mathcal
L}$-exponential map ${\mathcal L}{\rm exp}^\tau_x$ is onto $\nu(A)\cap
M_{T-\tau}$.

\begin{claim}\label{twosets} The complement of $\nu(A)\cap {\mathcal U}_{x}(\tau)$
in $\nu(A)$ is contained in the union of two sets: The first is the set of
points $z$ where there is more than one minimizing ${\mathcal L}$-geodesic from
$x$ ending at $z$ and if $Z$ is the initial condition for any minimizing
${\mathcal L}$-geodesic to $z$ then the differential of ${\mathcal L}{\rm
exp}_x^\tau$ at any $Z$ is an isomorphism. The second is the intersection of
the set of critical values of ${\mathcal L}{\rm exp}^\tau$ with $\nu(A)\cap
M_{T-\tau}$.
\end{claim}

\begin{proof}
Suppose that $q\in \nu(A)\cap M_{T-\tau}$ is not contained in ${\mathcal U}_x$.
Let $\gamma_Z$ be a minimal ${\mathcal L}$-geodesic  from $x$ to $q$. If the
differential of ${\mathcal L}{\rm exp}_x$ is not an isomorphism at $Z$, then
$q$ is contained in the second set given in the claim. Thus, we can assume that
the differential of ${\mathcal L}{\rm exp}_x$ at $Z$, and hence ${\mathcal
L}{\rm exp}_x$ identifies a neighborhood $\widetilde V$ of $Z$ in ${\mathcal
H}T_z{\mathcal M}$ with a neighborhood $V\subset \nu(A)$  of $q$ in
$M_{T-\tau}$. Suppose that there is no neighborhood $\widetilde V'\subset
\widetilde V$ of $Z$ so that the ${\mathcal L}$-geodesics are unique minimal
${\mathcal L}$-geodesics to their endpoints in $M_{T_\tau}$. Then there is a
sequence of minimizing ${\mathcal L}$-geodesics $\gamma_n$ whose endpoints
converge to $q$, but so that no $\gamma_n$ has initial condition contained in
$\widetilde V'$. By hypothesis all of these geodesics are contained in $F$ and
hence $|{\rm Ric}|$ and $|\nabla R|$ are uniformly bounded on these geodesics.
Also,  by the continuity of ${\mathcal L}$, the ${\mathcal L}$-lengths of
$\gamma_n$ are uniformly bounded as $n$ tends to infinity. By
Corollary~\ref{minmax} we see that the initial conditions
$Z_n=\sqrt{\tau_1}X_{\gamma_n}(\tau_1)$ (meaning the limit as $\tau\rightarrow
0$ of these quantities in the case when $\tau_1=0$) are of uniformly bounded
norm. Hence, passing to a subsequence we can arrange that the $Z_n$ converge to
some $Z_\infty$. The tangent vector $Z_\infty$ is the initial condition of an
${\mathcal L}$-geodesic $\gamma_\infty$. Since the $\gamma_n$ are minimizing
${\mathcal L}$-geodesics to a sequence of points $q_n$ converging to $q$, by
continuity it follows that $\gamma_\infty$ is a minimizing ${\mathcal
L}$-geodesic to $q$. Since none of the $Z_n$ is contained in $\widetilde V'$,
it follows that $Z_\infty\not= Z$. This is a contradiction, showing that
throughout some neighborhood $\widetilde V'$ of $Z$ the ${\mathcal
L}$-geodesics are unique minimizing ${\mathcal L}$-geodesics and completing the
proof of the claim.
\end{proof}

 According to the next claim, the first subset given in Claim~\ref{twosets} is contained in the
set of points of $\nu(A)\cap M_{T-\tau}$ where $L_x^\tau$ is
non-differentiable. Since $L_x^\tau$ is a locally Lipschitz function on
$\nu(A)$, this subset is of measure zero in $\nu(A)$; see Rademacher's Theorem
on p. 81 of \cite{EvansGariepy}. The second set is of measure zero by Sard's
theorem. This proves, modulo the next claim, that ${\mathcal U}_x(\tau)\cap A$
is full measure in $A\cap M_{T-\tau}$.

\begin{claim}
Let $z\in M_{T-\tau}$. Suppose that there is a neighborhood of $z$ in
$M_{T-\tau}$ with the property that every point of the neighborhood is the
endpoint of a minimizing ${\mathcal L}$-geodesic from $x$, so that $L^\tau_x$
is defined on this neighborhood of $z$. Suppose that there are two distinct,
minimizing ${\mathcal L}$-geodesics $\gamma_{Z_1}$ and $\gamma_{Z_2}$ from $x$
ending at $z$ with the property that the differential of ${\mathcal L}{\rm
exp}^{\tau}$ is an isomorphism at both $Z_1$ and $Z_2$. Then the function
$L^\tau_x$ is non-differentiable at $z$.
\end{claim}

\begin{proof}
 Suppose that
$\gamma_{Z_0}|_{[0,\tau]}$ is an ${\mathcal L}$-minimal ${\mathcal L}$-geodesic
and that $d_{Z_0}{\mathcal L}{\rm exp}_x^\tau$ is an isomorphism. Then use
${\mathcal L}{\rm exp}_x^{\tau}$ to identify a neighborhood of $Z_0\in T_xM$
with a neighborhood of $z$ in $M_{T-\tau}$, and push  the function
$\widetilde{\mathcal L}^\tau_x$ on this neighborhood of $Z_0$ down to a
function $L_{Z_0}$ on a neighborhood in $M_{T-\tau}$ of $z$. According to
Lemma~\ref{DtildeL} the resulting function $L_{Z_0}$ is smooth and its gradient
at $z$ is $2\sqrt{\tau}X(\tau)$. Now suppose that there is a second ${\mathcal
L}$-minimizing ${\mathcal L}$-geodesic to $z$ with initial condition
$Z_1\not=Z_0$ and with $d_{Z_1}{\mathcal L}{\rm exp}_x^\tau$ being an
isomorphism. Then near $z$ the function $L_x^{\tau}$ is less than or equal to
the minimum of two smooth functions $L_{Z_0}$ and $L_{Z_1}$. We have
$L_{Z_0}(z)=L_{Z_1}(z)=L_x^\tau(z)$,  and furthermore, $L_{Z_0}$ and $L_{Z_1}$
have distinct gradients at $z$. It follows that $L_x^{\tau}$ is not smooth at
$z$.
\end{proof}

This completes the proof of the Corollary~\ref{fullmeasure}.
\end{proof}

\section{Reduced volume}

Here, we  assume that $x\in M_T\subset {\mathcal M}$, so that
$\tau_1=0$ in this subsection.

\begin{defn}\label{redvol}
Let $A\subset {\mathcal U}_x(\tau)$ be a  measurable subset of
$M_{T-\tau}$. The $\mathcal{L}$-{\em reduced volume of $A$ from $x$}
(or the {\em reduced volume}\index{reduced volume} for short) is
defined to be
$$\widetilde V_x(A) = \int_{A} \tau^{-\frac{n}{2}}{\rm
exp}(-l_x(q))dq$$\index{$\widetilde V_x$} where $dq$ is the volume
element of the metric $G(T-\tau)$.
\end{defn}

\begin{lem}\label{redvoltilde}
Let $A\subset {\mathcal U}_x(\tau)$ be a  measurable subset. Define
$\widetilde A\subset \widetilde{\mathcal U}_x(\tau)$ to be the
pre-image under ${\mathcal L}{\rm exp}_x^{\tau}$ of $A$. Then
$$\widetilde V_x(A)=\int_{\widetilde A}\tau^{-\frac{n}{2}}{\rm
exp}(-\tilde l(Z, \tau)){\mathcal J}(Z,\tau)dZ,$$ where $dZ$ is the usual
Euclidean volume element and ${\mathcal J}(Z,\tau)$ is the Jacobian determinant
of ${\mathcal L}{\rm exp}_x^{\tau}$ at $Z\in T_xM_T$.
\end{lem}

\begin{proof}
This is simply the change of variables formula for integration.
\end{proof}

Before we can study the reduced volume we must study the function
that appears as the integrand in its definition. To understand the
limit as $\tau\rightarrow 0$ requires a rescaling argument.

\subsection{Rescaling}

Fix $Q>0$. We rescale to form $(Q{\mathcal M},QG)$ and then we shift
the time by $T-QT$ so that the time-slice $M_T$ in the original flow
is the $T$ time-slice of the new flow. We call the result
$({\mathcal M}',G')$. Recall that $\tau=T-{\bf t}$ is the parameter
for ${\mathcal L}$-geodesics in $({\mathcal M},G)$ The corresponding
parameter in the rescaled flow $({\mathcal M}',G')$ is $\tau'=T-{\bf
t'}=Q\tau$. We denote by ${\mathcal L}'{\rm exp}_x$ the ${\mathcal
L}$-exponential map from $x$ in $({\mathcal M}',G')$, and by $l_x'$
the reduced length function for this Ricci flow. The associated
function on the tangent space is denoted $\tilde l'$.

\begin{lem}\label{Lrescale}
Let $({\mathcal M},G)$ be a generalized Ricci flow and let $x\in
M_T\subset {\mathcal M}$. Fix $Q>0$ and let $({\mathcal M}',G')$ be
the $Q$ scaling and shifting of $({\mathcal M},G)$ as described in
the previous paragraph. Let $\iota\colon{\mathcal M}\to {\mathcal
M}'$ be the identity map. Suppose that $\gamma\colon [0,\bar\tau]\to
{\mathcal M}$ is a path parameterized by backward time with
$\gamma(0)=x$. Let $\beta\colon [0,Q\bar\tau]\to Q{\mathcal M}$ be
defined by
$$\beta(\tau')=\iota(\gamma(\tau'/Q)).$$
Then $\beta(0)=x$ and $\beta$ is parameterized by backward time in
$({\mathcal M}',G')$, and ${\mathcal L}(\beta)=\sqrt{Q}{\mathcal
L}(\gamma)$. Furthermore, $\beta$ is an ${\mathcal L}$-geodesic if
and only if $\gamma$ is. In this case, if $Z={\rm
lim}_{\tau\rightarrow 0}\sqrt{\tau}X_\gamma(\tau)$ then
$\sqrt{Q^{-1}}Z={\rm lim}_{\tau'\rightarrow
0}\sqrt{\tau'}X_\beta(\tau')$
\end{lem}

\begin{rem}
Notice that $|Z|_G^2=|\sqrt{Q^{-1}}Z|^2_{G'}$.
\end{rem}

\begin{proof}
It is clear that $\beta(0)=x$ and that $\beta$ is parameterized by
backward time in $({\mathcal M}',G')$. Because of the scaling of
space and time by $Q$, we have $R_{G'}=R_G/Q$ and
$X_\beta(\tau')=d\iota(X_\gamma(\tau))/Q$, and hence
$|X_\beta(\tau')|^2_{G'}=\frac{1}{Q}|X_\gamma(\tau)|^2_{G}$. A
direct change of variables in the integral then shows that
$${\mathcal L}(\beta)=\sqrt{Q}{\mathcal L}(\gamma).$$
It follows that $\beta$ is an ${\mathcal L}$-geodesic if and only if
$\gamma$ is. The last statement follows directly.
\end{proof}

Immediately from the definitions we see the following:

\begin{cor}\label{COR}
With notation as above, and with the substitution $\tau'=Q\tau$, for
any $Z\in {\mathcal H}T_x{\mathcal M}$ and any $\tau>0$ we have
$${\mathcal
L}'{\rm exp}_x(\sqrt{Q^{-1}}Z,\tau')=\iota({\mathcal L}{\rm
exp}_x(Z,\tau))$$ and
$$ \tilde l'(\sqrt{Q^{-1}}Z,\tau')=\tilde l(Z,\tau),$$
whenever these are defined.
\end{cor}

\subsection{The integrand in the reduced volume integral}

Now we turn our attention to the integrand\index{reduced
volume!integrand} (over $\widetilde {\mathcal U}_x(\tau)$) in the
reduced volume integral. Namely, set
$$f(\tau)=\tau^{-n/2}e^{-\tilde l(Z,\tau)}{\mathcal J}(Z,\tau),$$
where ${\mathcal J}(Z,\tau)$\index{${\mathcal J}(Z,\tau)$\ii} is the Jacobian
determinant of ${\mathcal L}{\rm exp}_x^\tau$ at the point $Z\in \widetilde
U_x(\tau)\subset T_xM_T$. We wish to see that this quantity is invariant under
the rescaling.

\begin{lem}\label{Qformula}
With the notation as above let ${\mathcal J}'(Z,\tau')$ denote the
Jacobian determinant of  ${\mathcal L}'{\rm exp}_x$. Then, with the
substitution $\tau'=Q\tau$, we have
$$(\tau')^{-n/2}e^{-\tilde l'(\sqrt{Q^{-1}}Z,\tau')}{\mathcal J}'(\sqrt{Q^{-1}}Z,\tau')=
\tau^{-n/2}e^{-\tilde l(Z,\tau)}{\mathcal J}(Z,\tau).$$
\end{lem}

\begin{proof}
It follows from the first equation in Corollary~\ref{COR} that
$$J(\iota){\mathcal J}(Z,\tau)=J(\sqrt{Q^{-1}}){\mathcal
J}'(\sqrt{Q^{-1}}Z,\tau'),$$ where $J(\iota)$ is the Jacobian
determinant of $\iota$ at ${\mathcal L}{\rm exp}_x(Z,\tau)$ and
$J(\sqrt{Q^{-1}})$ is the Jacobian determinant of multiplication by
$\sqrt{Q^{-1}}$ as a map from $T_xM_T$ to itself, where the domain
has the metric $G$ and the range has metric $G'=QG$. Clearly, with
these conventions, we have $J(\iota)=Q^{n/2}$ and
$J(\sqrt{Q^{-1}})=1$. Hence, we conclude
$$Q^{n/2}{\mathcal J}(Z,\tau)={\mathcal J}'(\sqrt{Q^{-1}}Z,\tau').$$

Letting $\gamma$ be the ${\mathcal L}$-geodesic in $({\mathcal
M},G)$ with initial condition $Z$ and $\beta$ the ${\mathcal
L}$-geodesic in $({\mathcal M}',G')$ with initial condition
$\sqrt{Q^{-1}}Z$,  by Lemma~\ref{Lrescale} we have
$\gamma(\tau)=\beta(\tau')$. From Corollary~\ref{COR} and the
definition of the reduced length, we get $$\tilde
l'(\sqrt{Q^{-1}}Z,\tau')=\tilde l(\gamma,\tau).$$ Plugging these in
gives the result.
\end{proof}

Let us evaluate $f(\tau)$ in the case of $\Ar^n$ with the Ricci flow being the
constant family of Euclidean metrics.

\begin{exam}
Let the Ricci flow be the constant family of standard metrics on
$\Ar^n$. Fix $x=(p,T)\in \Ar ^n\times (-\infty,\infty)$. Then
$${\mathcal L}{\rm exp}_x(Z,\tau)=(p+2\sqrt{\tau}Z,T-\tau).$$
In particular, the Jacobian determinant of ${\mathcal L}{\rm
exp}_{(x,T)}^\tau$ is constant and equal to $2^n\tau^{n/2}$. The
$\widetilde l$-length of the ${\mathcal L}$-geodesic
$\gamma_Z(\tau)=(p+2\sqrt{\tau}Z, T-\tau),\ 0\le \tau\le \bar\tau$,
is $|Z|^2$.
\end{exam}

Putting these computations together gives the following.

\begin{claim}\label{Arnquantity}
In the case of the constant flow on Euclidean space we have
$$f(\tau)=\tau^{-n/2}e^{-\tilde l(Z,\tau)}{\mathcal J}(Z,\tau)=2^ne^{-|Z|^2}.$$
\end{claim}

This computation has consequences for all Ricci flows.

\begin{prop}\label{Jlimit}
Let $({\mathcal M},G)$ be a generalized Ricci flow and let $x\in M_T\subset
{\mathcal M}$. Then, for any $A<\infty$, there is $\delta>0$ such that the map
${\mathcal L}{\rm exp}_x$ is defined on $B(0,A)\times (0,\delta)$, where
$B(0,A)$ is the ball of radius $A$ centered at the origin in $T_xM_T$.
Moveover, ${\mathcal L}{\rm exp}_x$ defines a diffeomorphism of $B(0,A)\times
(0,\delta)$ onto an open subset of ${\mathcal M}$. Furthermore,
$${\rm lim}_{\tau\rightarrow 0}\tau^{-n/2}e^{-\tilde l(Z,\tau)}{\mathcal J}(Z,\tau)=2^ne^{-|Z|^2},$$
where the convergence is uniform on each compact subset of $T_xM_T$.
\end{prop}

\begin{proof}
First notice that since $T$ is greater than the initial time of ${\mathcal M}$,
there is $\epsilon>0$, and an embedding $\rho\colon B(x,T,\epsilon)\times
[T-\epsilon,T]\to {\mathcal M}$ compatible with time  and the vector field. By
taking $\epsilon>0$ smaller if necessary, we can assume that the image of
$\rho$ has compact closure in ${\mathcal M}$. By compactness every higher
partial derivative (both spatial and temporal) of the metric is bounded on the
image of $\rho$.

Now take a sequence of positive constants $\tau_k$ tending to $0$ as
$k\rightarrow\infty$, and set $Q_k=\tau_k^{-1}$. We let $({\mathcal
M}_k,G_k)$ be the $Q_k$-rescaling and shifting of $({\mathcal M},G)$
as described at the beginning of this section. The rescaled version
of $\rho$ is an embedding
$$\rho_k\colon B_{G_k}(x,T,\sqrt{Q_k}\epsilon)\times [T-Q_k\epsilon,T]\to
{\mathcal M}_k$$ compatible with the time function ${\bf t}_k$ and the vector
field. Furthermore, uniformly on the image of $\rho_k$, every higher partial
derivative of the metric is bounded by a constant that goes to zero with $k$.
Thus, the generalized Ricci flows $({\mathcal M}_k,G_k)$ based at $x$ converge
geometrically to the constant family of Euclidean metrics on $\Ar^n$. Since the
ODE given in Equation~(\ref{reparam}) is regular even at $0$, this implies that
the ${\mathcal L}$-exponential maps for these flows converge uniformly on the
balls of finite radius centered at the origin of the tangent spaces at $x$ to
the ${\mathcal L}$-exponential map of $\Ar^n$ at the origin. Of course, if
$Z\in T_xM_T$ is an initial condition for an ${\mathcal L}$-geodesic in
$({\mathcal M},G)$, then $\sqrt{Q_k^{-1}}Z$ is the initial condition for the
corresponding ${\mathcal L}$-geodesic in $({\mathcal M}_k,G_k)$. But
$|Z|_G=|\sqrt{Q_k^{-1}}Z|_{G_k}$, so that if $Z\in B_G(0,A)$ then
$\sqrt{Q_k^{-1}}Z\in B_{G_k}(0,A)$.
 In particular, we see
that for any $A<\infty$, for all $k$ sufficiently large, the
${\mathcal L}$-geodesics are defined on $B_{G_k}(0,A)\times (0,1]$
and the image is contained in the image of $\rho_k$. Rescaling shows
that for any $A<\infty$  there is $k$ for which the ${\mathcal
L}$-exponential map is defined on $B_G(0,A)\times (0,\tau_k]$ and
has image contained in $\rho$.

Let $Z\in B_Q(0,A)\subset T_xM_T$, and let $\gamma$ be the
${\mathcal L}$-geodesic with ${\rm lim}_{\tau\rightarrow
0}\sqrt{\tau}X_\gamma(\tau)=Z$. Let $\gamma_k$ be the corresponding
${\mathcal L}$-geodesic in $({\mathcal M}_k,G_k)$. Then ${\rm
lim}_{\tau\rightarrow
0}\sqrt{\tau}X_{\gamma_k}(\tau)=\sqrt{\tau_k}Z=Z_k$. Of course,
$|Z_k|^2_{G_k}=|Z|^2_G$, meaning that $Z_k$ is contained in the ball
$B_{G_k}(0,A)\subset T_xM_T$ for all $k$. Hence, by passing to a
subsequence we can assume that, in the geometric limit, the
$\sqrt{\tau_k}Z$ converge to a tangent vector $Z'$ in the ball of
radius $A$ centered at the origin in the tangent space to Euclidean
space. Of course $|Z'|^2=|Z|_G^2$. By Claim~\ref{Arnquantity}, this
means that we have
$${\rm lim}_{k\rightarrow \infty}1^{-n/2}e^{-\tilde l_k(\sqrt{Q_k^{-1}}Z,1)}{\mathcal
J}_k(\sqrt{Q_k^{-1}}Z,1)=2^ne^{-|Z|^2},$$ where ${\mathcal J}_k$ is
the Jacobian determinant of the ${\mathcal L}$-exponential map for
$({\mathcal M}_k,G_k)$.
 Of course, since
$\tau_k=Q_k^{-1}$, by Lemma~\ref{Qformula} we have
$$1^{-n/2}e^{-\tilde l_k(\sqrt{Q_k^{-1}}Z,1)}{\mathcal
J}_k(\sqrt{Q_k^{-1}}Z,1)=\tau_k^{-n/2}e^{-\tilde
l(Z,\tau_k)}{\mathcal J}(Z,\tau_k).$$ This establishes the limiting
result.

Since the geometric limits are uniform on balls of finite radius
centered at the origin in the tangent space, the above limit also is
uniform over each of these balls.
\end{proof}

\begin{cor}\label{Jlimit1}
Let $({\mathcal M},G)$ be a generalized Ricci flow whose sectional
curvatures are bounded. For any $x\in M_T$ and any $R<\infty$ for
all $\tau>0$ sufficiently small, the ball of radius $R$ centered at
the origin in $T_xM_T$ is contained in $\widetilde {\mathcal
U}_x(\tau)$.
\end{cor}

\begin{proof}
According to the last result, given $R<\infty$, for all $\delta>0$ sufficiently
small the ball of radius $R$ centered at the origin in $T_xM_T$ is contained in
${\mathcal D}_x^\delta$, in the domain of definition of ${\mathcal L}{\rm
exp}_x^\delta$ as given in Definition~\ref{defnD}, and ${\mathcal L}{\rm
exp}_x$ is a diffeomorphism on this subset. We shall show that if $\delta>0$ is
sufficiently small, then the resulting ${\mathcal L}$-geodesic $\gamma$ is the
unique minimizing ${\mathcal L}$-geodesic. If not then there must be another,
distinct ${\mathcal L}$-geodesic to this point whose ${\mathcal L}$-length is
no greater than that of $\gamma$. According to Lemma~\ref{min-max} there is a
constant $C_1$ depending on the curvature bound and on $\delta$ such that if
$Z$ is an initial condition for an ${\mathcal L}$-geodesic then for all
$\tau\in (0,\delta)$ we have
$$C_1^{-1}\left(|Z|-\frac{(C_1-1)}{2}\sqrt{\delta}\right)\le \sqrt{\tau}|X(\tau)|\le
C_1|Z|+\frac{(C_1-1)}{2}\sqrt{\delta}.$$ From the formula given in
Lemma~\ref{min-max} for $C_1$, it follows that, fixing the bound of
the curvature and its derivatives,
 $C_1\rightarrow 1$ as $\delta\rightarrow 0$. Thus, with
a given curvature bound, for $\delta$ sufficiently small,
$\sqrt{\tau}|X(\tau)|$ is almost a constant along ${\mathcal
L}$-geodesics. Hence, the integral of $\sqrt{\tau}|X(\tau)|^2$ is
approximately $2\sqrt{\delta}|Z|^2$. On the other hand, the absolute
value of the integral of $\sqrt{\tau}R(\gamma(\tau))$ is at most
$2C_0\delta^{3/2}/3$ where $C_0$ is an upper bound for the absolute
value of the scalar curvature.

Given $R<\infty$, choose  $\delta>0$ sufficiently small such that ${\mathcal
L}{\rm exp}_x$ is a diffeomorphism on the ball of radius $9R$ centered at the
origin and such that the following estimate holds: The ${\mathcal L}$-length of
an ${\mathcal L}$-geodesic defined on $[0,\delta]$ with initial condition $Z$
is between $\sqrt{\delta}|Z|^2$ and $3\sqrt{\delta}|Z|^2$. To ensure the latter
estimate we need only take $\delta$ sufficiently small given the curvature
bounds and the dimension. Hence, for these $\delta$ no ${\mathcal L}$-geodesic
with initial condition outside the ball of radius $9R$ centered at the origin
in $T_xM_T$ can be as short as any ${\mathcal L}$-geodesic with initial
condition in the ball of radius $R$ centered at the same point. This means that
the ${\mathcal L}$-geodesics defined on $[0,\delta]$ with initial condition
$|Z|$ with $|Z|<R$ are unique minimizing ${\mathcal L}$-geodesics.
\end{proof}

\subsection{Monotonicity of reduced volume}\index{reduced
volume!monotonicity}

Now we are ready to state and prove our main result concerning the
reduced volume.

\begin{thm}\label{Amono}
Fix $x\in M_T\subset {\mathcal M}$. Let $A\subset{\mathcal
U}_x\subset {\mathcal M}$ be an open  subset. We suppose that for
any $0<\tau\le \bar\tau$ and any $y\in A_\tau=A\cap M_{T-\tau}$ the
minimizing ${\mathcal L}$-geodesic from $x$ to $y$ contained in
$A\cup \{x\}$.  Then $\widetilde V_x(A_\tau)$ is a non-increasing
function of $\tau$ for all $0<\tau\le \bar\tau$.
\end{thm}

\begin{proof}
Fix $\tau_0\in (0,\bar\tau]$. To prove the theorem we shall show
that for any $0<\tau<\tau_0$ we have $\widetilde V_x(A_\tau)\ge
\widetilde V_x(A_{\tau_0})$. Let $\widetilde A_{\tau_0}\subset
\widetilde {\mathcal U}_x(\tau_0)$ be the pre-image under ${\mathcal
L}{\rm exp}_x^{\tau_0}$ of $A_{\tau_0}$. For each $0<\tau\le \tau_0$
we set
$$A_{\tau,\tau_0}={\mathcal L}{\rm exp}_x^\tau(\widetilde
A_{\tau_0})\subset M_{T-\tau}.$$ It follows from the assumption on
$A$ that $A_{\tau,\tau_0}\subset A_\tau$, so that $\widetilde
V_x(A_{\tau,\tau_0})\le \widetilde V_x(A_\tau)$. Thus, it suffices
to show that for all $0<\tau\le \tau_0$ we have
$$\widetilde V_x(A_{\tau,\tau_0})\ge \widetilde V_x(\tau_0).$$

Since
$$ \widetilde V_x(A_{\tau,\tau_0})=\int_{\widetilde A_{\tau_0}} \tau^{-\frac{n}{2}}
{\rm exp}(-\tilde l(Z,\tau)){\mathcal J}(Z,\tau)dZ,$$ the theorem
follows from:

\begin{prop}\label{fmonotone} For each $Z\in \widetilde
{\mathcal U}_x(\bar\tau)\subset T_xM_T$ the function
$$f(Z,\tau)=\tau^{-\frac{n}{2}}e^{-\tilde l(Z,\tau)}{\mathcal J}(Z,\tau)$$
is a non-increasing function of $\tau$ on the interval
$(0,\bar\tau]$ with ${\rm lim}_{\tau\rightarrow
0}f(Z,\tau)=2^ne^{-|Z|^2}$, the limit being uniform on any compact
subset of $T_xM_T$.
\end{prop}

\begin{proof}
 First, we analyze the Jacobian
${\mathcal J}(Z,\tau)$. We know that $\mathcal{L}{\rm exp}_x^{\tau}$
is smooth in a neighborhood of $Z$. Choose a basis
$\{\partial_\alpha\}$ for $T_xM_T$ such that $\partial_\alpha$
pushes forward under the differential at $Z$ of ${\mathcal L}{\rm
exp}_x^\tau$  to an orthonormal basis $\{Y_\alpha\}$ for
$M_{T-\tau}$ at $\gamma_Z(\tau)$. Notice that, letting $\tau'$ range
from $0$ to $\tau$ and taking the push-forward of the
$\partial_\alpha$ under the differential at $Z$ of ${\mathcal L}{\rm
exp}_x^{\tau'}$ produces a basis of ${\mathcal L}$-Jacobi fields
$\{Y_\alpha(\tau')\}$ along $\gamma_Z$. With this understood, we
have:
\begin{eqnarray*}
\frac{\partial}{\partial \tau}{\rm ln}\mathcal{J}|_\tau & = &
\frac{d}{d\tau}{\rm ln}(\sqrt{{\rm det}(\langle
Y_{\alpha},Y_{\beta}\rangle
)})\\
 & = &
\frac{1}{2}\Bigl(\frac{d}{d\tau}\sum_{\alpha}\abs{Y_{\alpha}}^{2}\Bigr)\Bigl|_{\tau}\Bigr.
.
\end{eqnarray*}
By Lemma~\ref{Jacobi} and by Proposition~\ref{Hessineq} (recall that
$\tau_1=0$) we have
\begin{eqnarray}\nonumber
\frac{1}{2}\frac{d}{d\tau}|Y_\alpha(\tau)|^2 & = &
\frac{1}{2\sqrt{\tau}}{\rm Hess}(L)(Y_\alpha,Y_\alpha) +{\rm
Ric}(Y_\alpha,Y_\alpha) \\ & \le &
\frac{1}{2\tau}-\frac{1}{2\sqrt{\tau}}\int_0^{\tau}\sqrt{\tau'}H(X,\tilde
Y_\alpha(\tau'))d\tau',\label{basic}\end{eqnarray} where $\tilde
Y_\alpha(\tau')$ is the adapted vector field along $\gamma$ with
$\tilde Y(\tau)=Y_\alpha(\tau)$. Summing over $\alpha$ as in the
proof of Proposition~\ref{triL} and Claim~\ref{Hclaim} yields
\begin{eqnarray}
\frac{\partial}{\partial\tau}{\rm ln}\mathcal{J}(Z,\tau)|_{\tau} &
\le & \frac{n}{2\tau} -
\frac{1}{2\sqrt{\tau}}\sum_{\alpha}\int_{0}^{\tau}\sqrt{\tau'}H(X,\tilde
Y_{\alpha(\tau')})d\tau' \nonumber
\\&=  & \frac{n}{2\tau} - \frac{1}{2}\tau^{-\frac{3}{2}}K^\tau(\gamma_Z).\label{Jinequ}
\end{eqnarray}

On $\widetilde {\mathcal U}_x(\tau)$ the expression
$\tau^{-\frac{n}{2}}e^{-\widetilde l(Z,\tau)}{\mathcal J}(Z,\tau)$
is positive, and so we have
$$\frac{\partial}{\partial\tau}{\rm ln}\left(\tau^{-\frac{n}{2}}e^{-\widetilde
l(Z,\tau)}{\mathcal J}(Z,\tau)\right)\le
\left(-\frac{n}{2\tau}-\frac{d\widetilde l}{d
\tau}+\frac{n}{2\tau}-\frac{1}{2}\tau^{-\frac{3}{2}}K^\tau(\gamma_Z)\right).$$

Corollary~\ref{tildelK} says that the right-hand side of the
previous inequality is zero. Hence, we conclude
\begin{equation}\label{Jineq}
\frac{d}{d\tau}\left(\tau^{-\frac{n}{2}}e^{-\tilde
l(X,\tau)}{\mathcal J}(X,\tau)\right)\le 0. \end{equation}

This proves the inequality given in the statement of the
proposition. The limit statement as $\tau\rightarrow 0$ is contained
in Proposition~\ref{Jlimit}.
\end{proof}

As we have already seen, this proposition implies
 Theorem~\ref{Amono}, and hence the proof of this theorem is complete.
\end{proof}

Notice that we have established the following:

\begin{cor}\label{finiteredvol}
For any measurable subset $A\subset {\mathcal U}_x(\tau)$ the
reduced volume  $\widetilde V_x(A)$ is at most $(4\pi)^{n/2}$.
\end{cor}

\begin{proof} Let $\widetilde A\subset \widetilde {\mathcal U}_x(\tau)$ be the pre-image
of $A$. We have seen that
$$\widetilde V_x(A)=\int_{A}\tau^{-n/2}e^{-l(q,\tau)}dq=\int_{\widetilde
A}\tau^{-n/2}e^{-\tilde l(Z,\tau)}{\mathcal J}(Z,\tau)dz.$$ By
Theorem~\ref{Amono} we see that $\tau^{-n/2}e^{-\tilde
l(Z,\tau)}{\mathcal J}(Z,\tau)$ is a non-increasing function of
$\tau$ whose limit as $\tau\rightarrow 0$ is the restriction of
$2^ne^{-|Z|^2}$ to $\widetilde A$. The result is immediate from
Lebesgue dominated convergence.
\end{proof}

\chapter{Complete Ricci flows of bounded curvature}\label{newcomp2}

 In this chapter we establish strong results for ${\mathcal L}{\rm
exp}_x$ in the case of ordinary Ricci flow on complete $n$-manifolds
with appropriate curvature bounds. In particular, for these flows we
show that there is a minimizing ${\mathcal L}$-geodesic to every
point. This means that $l_x$\index{$l_x$|(} is everywhere defined.
We extend the differential inequalities for $l_x$ established in
Section~\ref{secondorder} at the `smooth points' to weak
inequalities (i.e., inequalities in the distributional sense) valid
on the whole manifold. Using this we prove an upper bound for the
minimum of $l_x^\tau$.

Let us begin with a definition that captures the necessary curvature
bound for these results.

\begin{defn}
Let $(M,g(t)),\ a\le t\le b$, be a  Ricci flow. We say that the flow
is {\em complete of bounded curvature} if for each $t\in [a,b]$ the
Riemannian manifold $(M,g(t))$ is complete and if there is
$C<\infty$ such that $|{\rm Rm}|(p,t)\le C$ for all $p\in M$ and all
$t\in [a,b]$. Let $I$ be an interval and let $(M,g(t)),\ t\in I$, be
a Ricci flow. Then we say that the flow is {\em complete with
curvature locally bounded in time} if for each compact subinterval
$J\subset I$ the restriction of the flow to $(M,g(t)),\ t\in J$, is
complete of bounded curvature.
\end{defn}

\section{The functions $L_{x}$ and $l_{x}$}

Throughout Chapter~\ref{newcomp2} we have a Ricci flow $(M,g(t)),\
0\le t\le T<\infty$, and we set $\tau=T-t$. All the results of the
last chapter apply in this context,
 but in fact in this context there are much stronger results, which
we develop here.

\subsection{Existence of ${\mathcal L}$-geodesics}

We assume here that $(M,g(t)),\ 0\le t\le T<\infty$, is a Ricci flow
which is complete of bounded curvature. In Shi's Theorem
(Theorem~\ref{shi}) we take $K$ equal to the bound of the norm of
the Riemannian curvature on $M\times [0,T]$, we take $\alpha=1$, and
we take $t_0=T$. It follows from Theorem~\ref{shi} that there is a
constant $C(K,T)$ such that $|\nabla R(x,t)|\le C/t^{1/2}$. Thus,
for any $\epsilon>0$ we have a uniform bound for $|\nabla R|$ on
$M\times [\epsilon,T]$. Also, because of the uniform bound for the
Riemann curvature and the fact that $T<\infty$, there is a constant
$C$, depending on the curvature bound and $T$ such that
\begin{equation}\label{compmetric}
C^{-1}g(x,t)\le g(x,0)\le Cg(x,t)
\end{equation}
for all $(x,t)\in M\times [0,T]$.

\begin{lem}\label{geoexist} Assume that $M$ is connected.
Given $p_1,p_2\in M$ and $0\le\tau_1<\tau_2\le T$, there is
 a minimizing ${\mathcal L}$-geodesic:
$\gamma\colon[\tau_{1},\tau_{2}]\to M\times [0,T]$ connecting
$(p_1,\tau_1)$ to $(p_2,\tau_2)$.
\end{lem}

\begin{proof} For any curve $\gamma$ parameterized by backward time, we set
$\bar \gamma$ equal to the path in $M$ that is the image under
projection of $\gamma$. We set $A(s)=\bar\gamma'(s)$. Define
$$c((p_1,\tau_1),(p_2,\tau_{2})) = {\rm inf}
\{\mathcal{L}(\gamma)|\gamma\colon [\tau_{1},\tau_{2}] \rightarrow
M\times [0,T], \bar\gamma(\tau_{1}) = p_1, \bar\gamma(\tau_{2})=
p_2\}.$$ From Equation~(\ref{seqn}) we see that the infimum exists
since, by assumption, the curvature is uniformly bounded (below).
Furthermore, for a minimizing sequence $\gamma_{i}$, we have
$\int_{s_1}^{s_{2}}\abs{A_i(s)}^{2}ds \leq C_0$, for some constant
$C_0$, where $s_i=\sqrt{\tau_i}$ for $i=1,2$. It follows from this
and the inequality in Equation~(\ref{compmetric}) that there is a
constant $C_1<\infty$ such that for all $i$ we have
$$\int_{s_{1}}^{s_{2}} \abs{
A_{i}}^{2}_{g(0)}d\tau \leq C_1.$$ Therefore the sequence $\{\gamma_{i}\}$ is
uniformly continuous with respect to the metric $g(0)$; by Cauchy-Schwarz we
have
$$\abs{\bar\gamma_{i}(s) - \bar\gamma_{i}(s')}_{g(0)} \leq
\int_{s'}^{s} \abs{A_{i}}_{g(0)}ds \leq \sqrt{C_1}\sqrt{s-s'}.$$ By
the uniform continuity, we see that a subsequence of the $\gamma_i$
converges uniformly pointwise to a continuous curve $\gamma$
parameterized by $s$, the square root backward time. By passing to a
subsequence we can arrange that the $\gamma_i$ converge weakly in
$H^{2,1}$. Of course, the limit in $H^{2,1}$ is represented by the
continuous limit $\gamma$. That is to say, after passing to a
subsequence, the $\gamma_i$ converge uniformly and weakly in
$H^{2,1}$ to a continuous curve $\gamma$. Let $A(s)$ be the
$L^2$-derivative of $\gamma$. Weak convergence in $H^{2,1}$ implies
that $\int_{s'}^s|A(s)|^2ds\le {\rm
lim}_{i\rightarrow\infty}\int_{s'}^s|A_i(s)|^2ds$, so that
${\mathcal L}(\gamma)\le {\rm lim}_{i\rightarrow\infty}{\mathcal
L}(\gamma_i)$. This means that $\gamma$ minimizes the ${\mathcal
L}$-length. Being a minimizer of ${\mathcal L}$-length, $\gamma$
satisfies the Euler-Lagrange equation and is smooth by the
regularity theorem of differential equations. This then is the
required minimizing ${\mathcal L}$-geodesic from $(p_1,\tau_1)$ to
$(p_2,\tau_2)$.
\end{proof}

Let us now show that it is always possible to uniquely extend
${\mathcal L}$-geodesics up to time $T$.

\begin{lem} For any $0\le \tau_1<\tau_2<T$ suppose that
 $\gamma\colon [\tau_1,\tau_2]\to M\times [0,T]$ is an ${\mathcal L}$-geodesic.
 Then $\gamma$ extends uniquely to an ${\mathcal L}$-geodesic
$\gamma\colon [0,T)\to M\times (0,T]$.
\end{lem}

\begin{proof}
We work with the parameter $s=\sqrt{\tau}$. According to
Equation~(\ref{reparam}), we have
$$\nabla_{\gamma'(s)}\gamma'(s)=2s^2\nabla
R-4s{\rm Ric}(\gamma'(s),\cdot).$$ This is an everywhere
non-singular ODE. Since the manifolds $(M,g(t))$ are complete and
their metrics are uniformly related as in
Inequality~(\ref{compmetric}), to show that the solution is defined
on the entire interval $s\in [0,\sqrt{T})$ we need only show that
there is a uniform bound to the length, or equivalently the energy
of $\gamma$ of any compact subinterval of $[0,T)$ on which it is
defined. Fix $\epsilon>0$. It follows immediately from
Lemma~\ref{min-max}, and the fact that the quantities $R$, $|\nabla
R|$ and $|{\rm Rm}|$ are bounded on $M\times [\epsilon,T]$, that
there is a bound on ${\rm max}|\gamma'(s)|$ in terms of
$|\gamma'(\tau_1)|$, for all $s\in [0,\sqrt{T-\epsilon}]$ for which
$\gamma$ is defined. Since $(M,g(0))$ is complete, this, together
with a standard extension result for second-order ODE's, implies
that $\gamma$ extends uniquely to the entire interval
$[0,\sqrt{T-\epsilon}]$. Changing the variable from $s$ to
$\tau=s^2$ shows that the ${\mathcal L}$-geodesic extends uniquely
to the entire interval $[0,T-\epsilon]$. Since this is true for
every $\epsilon>0$, this completes the proof.
\end{proof}

Let $p\in M$ and set $x=(p,T)\in M\times [0,T]$. Recall that from
Definition~\ref{injdefn} for every $\tau>0$, the injectivity set
$\widetilde{\mathcal U}_{x}(\tau)\subset T_pM$ consists of all $Z\in
T_pM$ for which (i) the ${\mathcal L}$-geodesic
$\gamma_Z|_{[0,\tau]}$ is the unique minimizing ${\mathcal
L}$-geodesic from $x$ to its endpoint, (ii) the differential of
${\mathcal L}{\rm exp}_{x}^\tau$ is an isomorphism at $Z$, and (iii)
for all  $Z'$ sufficiently close to
 $Z$ the ${\mathcal L}$-geodesic $\gamma_{Z'}|_{[0,\tau]}$ is the unique minimizing
 ${\mathcal L}$-geodesic to its endpoint. The image
of $\widetilde{\mathcal U}_{x}(\tau)$ is denoted ${\mathcal U}_{x}(\tau)\subset
M$

The existence of minimizing ${\mathcal L}$-geodesics from $x$ to
every point of $M\times (0,T)$ means that the functions $L_{x}$ and
$l_{x}$ are defined on all of $M\times (0,T)$. This leads to:

\begin{defn} Suppose that $(M,g(t)),\ 0\le t\le T<\infty$, is a Ricci
flow, complete of bounded curvature.  We define the function
$L_{x}\colon M\times [0,T)\to \Ar$ by assigning to each $(q,t)$ the
length of any ${\mathcal L}$-minimizing ${\mathcal L}$-geodesic from
$x$ to $y=(q,t)\in M\times [0,T)$. Clearly, the restriction of this
function to ${\mathcal U}_{x}$ agrees with the smooth function
$L_{x}$ given in Definition~\ref{L_x}. We define
 $L_{x}^\tau\colon M\to \Ar$ to be the restriction of $L_x$ to
 $M\times \{T-\tau\}$.
Of course, the restriction of $L_{x}^\tau$ to ${\mathcal
U}_{x}(\tau)$ agrees with the smooth function $L_{x}^\tau$ defined
in the last chapter. We define $l_{x}\colon M\times [0,T)\to \Ar$ by
$l_{x}(y)=L_{x}(y)/2\sqrt{\tau}$, where, as always $\tau=T-t$, and
we define $l_{x}^\tau(q)=l_{x}(q,T-\tau)$.
\end{defn}

\subsection{Results about $l_{x}$ and ${\mathcal
U}_{x}(\tau)$}

Now we come to our main result about the nature of ${\mathcal
U}_{x}(\tau)$\index{${\mathcal U}_x(\tau)$} and the function $l_{x}$
in the context of Ricci flows which are complete and of bounded
curvature.

\begin{prop}\label{lipcomplete}
Let $(M,g(t)),\ 0\le t\le T<\infty$, be a Ricci flow that is
complete and of bounded curvature. Let $p\in M$, let $x=(p,T)\in
M\times [0,T]$, and let $\tau\in (0,T)$.
\begin{enumerate}
\item[(1)] The functions $L_{x}$ and $l_{x}$ are locally Lipschitz functions
on $M\times (0,T)$.
\item[(2)] ${\mathcal L}{\rm exp}_{x}^\tau$
 is a diffeomorphism from $\widetilde{\mathcal U}_{x}(\tau)$ onto
an open subset ${\mathcal U}_{x}(\tau)$\index{${\mathcal
U}_x(\tau)$} of $M$.
\item[(3)] The complement of ${\mathcal U}_{x}(\tau)$ in $M$ is a closed
subset of $M$ of zero Lebesgue measure. \item[(4)]
 For every
$\tau<\tau'<T$ we have
$${\mathcal U}_{x}(\tau')\subset {\mathcal U}_{x}(\tau).$$
\end{enumerate}
\end{prop}

\begin{proof}
By Shi's Theorem (Theorem~\ref{shi}) the curvature bound on $M\times
[0,T]$ implies that for each $\epsilon>0$ there is a bound for
$|\nabla R|$ on $M\times (\epsilon,T]$. Thus,
Proposition~\ref{lipaty} shows that $L_{x}$ is a
 locally Lipschitz function on $M\times (\epsilon,T)$. Since this is
true for every $\epsilon>0$, $L_x$ is a locally Lipschitz function
on $M\times (0,T)$. Of course, the same is true for $l_{x}$. The
second statement is contained in Proposition~\ref{diffeo}, and the
last one is contained in Proposition~\ref{inclusion}. It remains to
prove the third statement, namely that the complement of ${\mathcal
U}_{x}(\tau)$ is closed nowhere dense. This follows immediately from
Corollary~\ref{fullmeasure} since $|{\rm Ric}|$ and $|\nabla R|$ are
bounded on $F=M\times [T-\tau,T]$.
\end{proof}

\begin{cor}
The function $l_x$ is a continuous function on $M\times (0,T)$ and
is smooth on the complement of a closed subset ${\mathcal C}$ that
has the property that its intersection with each $M\times \{t\}$ is
of zero Lebesgue measure in $M\times \{t\}$. For each $\tau\in
(0,T)$ the gradient $\nabla l^\tau_x$ is then a smooth vector field
on the complement of ${\mathcal C}\cap M_{T_\tau}$, and it is a
locally essentially bounded vector field in the following sense. For
each $q\in M$ there is a neighborhood $V\subset M$ of $q$ such that
the restriction of $|\nabla l_x^\tau|$ to $V\setminus\left(V\cap
{\mathcal C}\right)$ is a bounded smooth function. Similarly,
$\partial l_x/\partial t$ is an essentially bounded smooth vector
field on $M\times (0,T)$.
\end{cor}

\section{A bound for ${\rm min}\, l^\tau_x$}

We continue to assume that we have a Ricci flow $(M,g(t)),\ 0\le
t\le T<\infty$, complete and of bounded curvature and a point
$x=(p,T)\in M\times [0,T]$. Our purpose here is to extend the first
differential inequality given in
 Corollary~\ref{lformula} to a differential
inequality in the weak or distributional sense for $l_{x}$ valid on
all of $M\times (0,T)$. We then use this to establish that ${\rm
min}_{q\in M} l^\tau_{x}(q)\le n/2$ for all $0<\tau<T$.

In establishing inequalities in the non-smooth case the notion of a
support function or a barrier function is often convenient.

\begin{figure}[ht]
  \centerline{\epsfbox{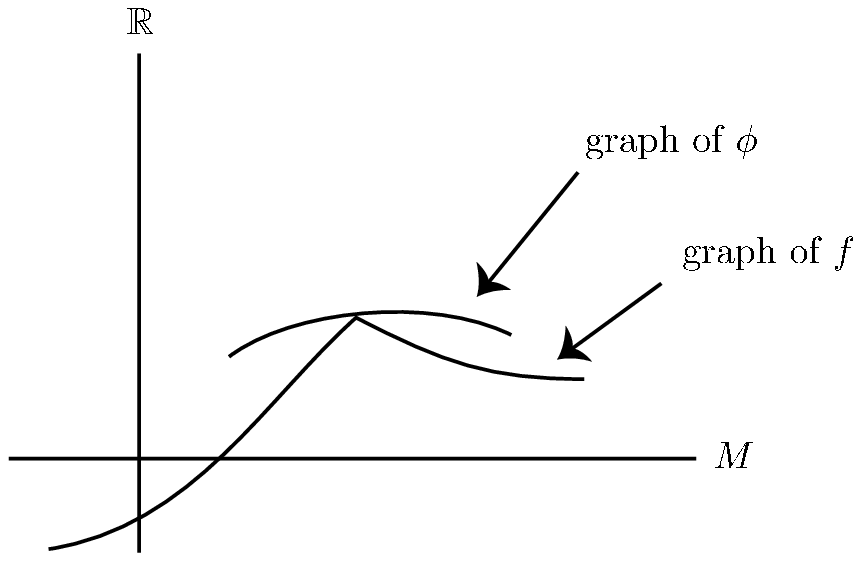}}
  \caption{Upper barrier.}\label{fig:upperbarrier}
\end{figure}

\begin{defn}
Let $P$ be a smooth manifold and let $f\colon P\to \Ar$ be a continuous
function. An {\em upper barrier}\index{upper barrier} for $f$ at $p\in P$ is a
smooth function $\varphi$ defined on a neighborhood of $p$ in $P$, say $U$,
satisfying $\varphi(p)=f(p)$ and $\varphi(u)\ge f(u)$ for all $u\in U$, see
{\sc Fig.}~\ref{fig:upperbarrier}.
\end{defn}

\begin{prop}\label{weaksense}
Let $(M,g(t)),\ 0\le t\le T<\infty$, be an $n$-dimensional Ricci
flow,  complete  of bounded curvature. Fix a point $x=(p,T)\in
M\times [0,T]$, and for any $(q,t)\in M\times [0,T]$, set
$\tau=T-t$. Then for any $(q,t)$, with $0<t<T$, we have
$$\frac{\partial l_x}{\partial \tau}(q,\tau) +\triangle l_x(q,\tau)\le
\frac{(n/2)-l_x(q,\tau)}{\tau}$$ in the barrier sense. This means
that
 for each $\epsilon>0$ there is a neighborhood $U$ of $(q,t)$ in
 $M\times [0,T]$
 and an upper barrier $\varphi$ for $l_x$ at this point defined on $U$
satisfying
$$\frac{\partial \varphi}{\partial \tau}(q,\tau)+\triangle \varphi(q,\tau)
\le \frac{(n/2)-l_x(q,\tau)}{\tau}+\epsilon.$$
\end{prop}

\begin{rem}
The operator $\triangle$ in the above statement is the horizontal
Laplacian, i.e., the Laplacian of the restriction of the indicated
function to the slice $M\times\{t=T-\tau\}$ as defined using the
metric $g(T-\tau)$ on this slice.
\end{rem}

\begin{proof}
If $(q,T-\tau)\in {\mathcal U}_{x}$, then $l_x$ is smooth near $(q,T-\tau)$,
and the result is immediate from the first inequality in
Corollary~\ref{lformula}.

Now consider a general point $(q,t=T-\tau)$ with $0<t<T$. According to
Lemma~\ref{geoexist} there is  a minimizing ${\mathcal L}$-geodesic $\gamma$
from $x=(p,T)$ to $(q,t=T-\tau)$. Let $\gamma$ be any minimizing ${\mathcal
L}$-geodesic between these points. Fix $0<\tau_1<\tau$ let $q_1=\gamma(\tau_1)$
and set $t_1=T-\tau_1$. Even though $q_1$ is contained in the $t_1$ time-slice,
we keep $\tau=T-t$ so that paths beginning at $q_1$ are parameterized by
intervals in the $\tau$-line of the form $[\tau_1,\tau']$ for some $\tau'<T$.
Consider ${\mathcal L}{\rm exp}_{(q_1,t_1)}\colon T_{q_1}M\times (\tau_1,T)\to
M\times (0,t_1)$. According to Proposition~\ref{tausmooth} there is a
neighborhood $\widetilde V$ of $\{\sqrt{\tau_1}X_\gamma(\tau_1)\}\times
(\tau_1,\tau]$ which is mapped diffeomorphically  by ${\mathcal L}{\rm
exp}_{(q_1,t_1)}\colon T_{q_1}M\times (\tau_1,\tau)\to M\times (0,t_1)$ onto a
neighborhood $V$ of $\gamma((\tau_1,\tau))$. (Of course, the neighborhood $V$
depends on $\tau_1$.) Let $L_{(q_1,t_1)}$ be the length function on $V$
obtained by taking the ${\mathcal L}$-lengths of geodesics parameterized by
points of $\widetilde V$. Let $\varphi_{\tau_1}\colon V\to \Ar$ be defined by
$$\varphi_{\tau_1}(q',\tau')=\frac{1}{2\sqrt{\tau'}}\left({\mathcal
L}(\gamma|_{[0,\tau_1]})+L_{(q_1,t_1)}(q',T-\tau')\right).$$
Clearly, $\varphi_{\tau_1}$ is an upper barrier for $l_x$ at
$(q,\tau)$. According to Lemma~\ref{linequal} we have
\begin{eqnarray*}
\frac{\partial \varphi_{\tau_1}}{d\tau}(q,\tau)+\triangle
\varphi_{\tau_1}(q,\tau) & \le &
\frac{n}{2\sqrt{\tau}(\sqrt{\tau}-\sqrt{\tau_1})}-\frac{\varphi_{\tau_1}(q,\tau)}{\tau}
+\frac{1}{2\tau^{3/2}}{\mathcal L}(\gamma|_{[0,\tau_1]})\\
& & +\frac{K^{\tau}_{\tau_1}(\gamma)}{2\tau^{3/2}}-\frac{{\mathcal
K}_{\tau_1}^{\tau}(\gamma)}{2\sqrt{\tau}(\sqrt{\tau}-\sqrt{\tau_1})^{2}}
\\ & & -\frac{1}{2}\left(\frac{\tau_1}{\tau}\right)^{3/2}\left(R(q_1,t_1)+|X(\tau_1)|^2\right).
\end{eqnarray*}
By Lemma~\ref{Kequal},  it follows easily that
$${\rm lim}_{\tau_1\rightarrow 0^+}
\frac{\partial \varphi_{\tau_1}}{\partial \tau}(q,t)+\triangle
\varphi_{\tau_1}(q,t)\le \frac{(n/2)-l_x(q,t)}{\tau}.$$ This
establishes the result.
\end{proof}

\begin{thm}\label{n/2} Suppose that $(M,g(t)),\ 0\le t\le T<\infty$,
is an $n$-dimensional Ricci flow, complete of bounded curvature.
Then for any $x=(p,T)\in M\times [0,T]$ and for every $0<\tau<T$
there is a point $q_\tau\in M$ such that $l_{x}(q_\tau,\tau)\le
\frac{n}{2}$.
\end{thm}

\begin{proof}
We set $l_{\rm min}(\tau)={\rm inf}_{q\in
M}l_x(q,\tau)$\index{$l_{\rm min}$|ii}. (We are not excluding the
possibility that this infimum is $-\infty$.) To prove this corollary
we first need to establish the following claim.

\begin{claim}\label{8.1}
For every
 $\tau\in (0,T)$ the function
$l_x(\cdot,\tau)$\index{$l_x$|)} achieves its minimum. Furthermore,
for every compact interval $I\subset (0,T)$ the subset of
$(q,\tau)\in M\times I$
 for which $l_x(q,\tau)=l_{\rm min}(\tau)$ is a compact set.
\end{claim}

 First, let us assume this claim and use it to prove the
theorem. We set $l_{\rm min}(\tau)={\rm min}_{q\in M}l_x(q,\tau)$.
(This minimum exists by the first statement in the claim.) From the
compactness result in the claim, it follows (see for example
Proposition~\ref{fordiffmax}) that $l_{\rm min}(\tau)$ is a
continuous function of $\tau$.

Suppose that $l_x(\cdot,\tau)$ achieves its minimum at $q$. Then by
the previous result, for any $\epsilon>0$ there is an upper barrier
$\varphi$ for $l_x$ at $(q,\tau)$ defined on an open subset $U$ of
$(q,\tau)\in M\times (0,T)$ and satisfying
$$\frac{d\varphi}{d\tau}(q,\tau)+\triangle \varphi(q,\tau)\le \frac{(n/2)-l_x(q,t)}{\tau}+\epsilon.$$
Since $l_x(q,\tau)=l_{\rm min}(\tau)$, it follows that
$\varphi(q',\tau)\ge \varphi(q,\tau)$ for all $(q',\tau)\in U\cap
M_{T-\tau}$. This means that $\triangle \varphi(q,\tau)\ge 0$, and
we conclude that
$$\frac{d\varphi}{d\tau}(q,\tau)\le \frac{(n/2)-l_{\rm min}(\tau)}{\tau}+\epsilon.$$
Since $\varphi$ is an upper barrier for $l_x$ at $(q,\tau)$ it
follows immediately that
$${\rm limsup}_{\tau'\rightarrow
\tau^+}\frac{l_x(q,\tau')-l_x(q,\tau)}{\tau'-\tau}\le
\frac{(n/2)-l_x(q,\tau)}{\tau}+\epsilon.$$ Since this is true for every
$\epsilon>0$, we see that
$${\rm limsup}_{\tau'\rightarrow
\tau^+}\frac{l_x(q,\tau')-l_x(q,\tau)}{\tau'-\tau}\le
\frac{(n/2)-l_x(q,\tau)}{\tau}.$$ Since $l_{\rm
min}(\tau)=l_x(q,\tau)$, the same inequality holds for the forward
difference quotient of $l_{\rm min}$ at $\tau$. That is to say, we
have
$${\rm limsup}_{\tau'\rightarrow \tau^+}\frac{l_{\rm min}(\tau')-l_{\rm
min}(\tau)}{\tau'-\tau}\le \frac{(n/2)-l_{\rm min}(\tau)}{\tau}.$$
The preceding equation implies that if $l_{\rm min}(\tau)\le n/2$
then $l_{\rm min}(\tau')\le n/2$ for every $\tau'\ge \tau$. On the
other hand ${\rm lim}_{\tau\rightarrow 0}l_{\rm min}(\tau)=0$. Then
reason for this is that the path $\tau'\mapsto (P,T-\tau')$ for
$\tau'\in [0,\tau]$ has ${\mathcal L}$-length $ O(\tau^{3/2})$ as
$\tau\rightarrow 0$. It follows that $l_{\rm min}(\tau)<n/2$ when
$\tau$ is small.

To complete the proof of Theorem~\ref{n/2}, it remains to prove
Claim~\ref{8.1}.

\begin{proof}
In the case when $M$ is compact, the claim is obvious. We consider
the case when $M$ is complete and the flow has bounded curvature.
Since the curvature on $M\times [0,T]$ is bounded, according to
Inequality~(\ref{compmetric}) there is a constant $C$ such that for
all $t,t'\in [0,T]$ we have
$$C^{-1}g(t')\le g(t)\le Cg(t').$$
For any compact interval $I\subset (0,T)$, there is $l_0<\infty$
such that $l_{\rm min}(\tau)\le l_0$ for all $\tau\in I$. According
to Corollary~\ref{minmax}, for every $\bar\tau\in I$ and all
${\mathcal L}$-geodesics from $x$ to points $(q,T-\bar\tau)$ of
lengths at most $2|l_0|$ there is an upper bound, say $C_2$, to
$|\sqrt{\tau}X_\gamma(\tau)|$. Thus, $|X_\gamma(\tau|\le
\frac{C_2}{\sqrt{\tau}}$, and hence
$$|X_\gamma(\tau)|_{g(T)}\le C_2\sqrt{C}/\sqrt{\tau}$$
 for these geodesics.
Thus, $$\int_0^{\bar\tau}|X_\gamma(\tau)|_{g(T)}d\tau\le
2\sqrt{\bar\tau}C_2\sqrt{C}.$$ This shows that there is $A<\infty$
such that for each $\bar\tau\in I$ and for any ${\mathcal
L}$-geodesic $\gamma$ defined on $[0,\bar\tau]$ of length at most
$2|l_0|$ the following holds. Letting $q\in M$ be such that
$\gamma(\bar\tau)=(q,T-\tau)$, the point $q$ lies in $B_T(p,A)$.
This implies  that  the endpoints of all such ${\mathcal
L}$-geodesics lie in a fixed compact subset of $M$ independent of
$\bar\tau\in I$ and the geodesic. Since the set of $(q,\tau)$ where
$l_x(q,\tau)=l_{\rm min}(\tau)$ is clearly a closed set, it follows
that the subset of $M\times I$ of all $(q,\tau)\in M\times I$ for
which $l_x(q,\tau)=l_{\rm min}(\tau)$ is compact. The last thing to
show is that for every $\tau\in I$ the function $l_x(\cdot,\tau)$
achieves it minimum. Fix $\tau\in I$ and let $q_n$ be a minimizing
sequence for $l_x(\cdot,\tau)$. We have already established that the
$q_n$ are contained in a compact subset of $M$, and hence we can
assume that they converge to a limit $q\in M$. Clearly, by the
continuity of $l_x$ we have $l_x(q,\tau)={\rm lim}_{n\rightarrow
\infty}l_x(q_n,\tau)={\rm inf}_{q'\in M}l_x(q',\tau)$, so that
$l_x(\cdot,\tau)$ achieves its minimum at $q$.
\end{proof}

Having established the claim, we have now completed the proof of
Theorem~\ref{n/2}.
\end{proof}

Actually, the proof given here also shows the following, which will
be useful later.

\begin{cor}\label{compactcoro}
Suppose that $({\mathcal M},G)$ is a generalized $n$-dimensional Ricci flow and
that $x\in {\mathcal M}$ is given and set $t_0={\bf t}(x)$. We suppose that
there is an open subset $U\subset {\bf t}^{-1}(-\infty,t_0)$ with the following
properties:
\begin{enumerate}
\item For every $y\in U$ there is a minimizing ${\mathcal L}$-geodesic from $x$
to $y$.
\item There are $r>0$ and $\Delta t>0$ such that the backward parabolic neighborhood
$P(x,t_0,r,-\Delta t)$ of $x$ exists in ${\mathcal M}$ and has  the property
that $P\cap {\bf t}^{-1}(-\infty,t_0)$ is contained in $U$.
\item For each compact interval (including the case of degenerate intervals
consisting of a single point) $I\subset (-\infty,t_0)$ the subset of points
$y\in {\bf t}^{-1}(I)\cap U$ for which ${\mathcal L}(y)={\rm inf}_{z\in {\bf
t}^{-1}({\bf t}(y))\cap U}{\mathcal L}(z)$ is compact and non-empty.
\end{enumerate}
 Then
for every $t<t_0$ the minimum of the restriction of $l_x$ to the time-slice
${\bf t}^{-1}(t)\cap U$ is at most $n/2$.
\end{cor}

\subsection{Extension of the other inequalities in
 Corollary~\protect{\ref{lformula}}}

The material in this subsection is adapted from \cite{Ye}. It
captures (in a weaker way) the fact that, in the case of geodesics
on a Riemannian manifold, the interior of the cut locus in $T_xM$ is
star-shaped from the origin.

\begin{thm}\label{weak}
Let $(M,g(t)),\ 0\le t\le T<\infty$, be a Ricci flow, complete and
of bounded curvature, and let $x=(p,T)\in M\times[0,T]$. The last
two inequalities in Corollary~\ref{lformula}, namely
$$\frac{\partial l_x}{\partial \tau}+|\nabla l^\tau_x|^2-R+\frac{n}{2\tau}
- \triangle l^\tau_x\ge 0$$
$$-|\nabla l^\tau_x|^2+R+\frac{l^\tau_x-n}{\tau}+2\triangle l^\tau_x\le 0$$
hold in the weak or distributional sense on all of $M\times\{\tau\}$
for all $\tau>0$. This means that for any $\tau>0$ and for any
non-negative, compactly supported, smooth function $\phi(q)$ on $M$
we have the following two inequalities:
$$\int_{M\times \{\tau\}}\Bigl[ \phi\cdot\left(\frac{\partial l_x}{\partial \tau}+|\nabla
l_x^\tau|^2-R+\frac{n}{2\tau}\right) - l_x^\tau\triangle \phi\Bigr]\
d{\rm vol}(g(t))\ge 0$$
$$\int_{M\times \{\tau\}}\Big[ \phi\cdot\left(-|\nabla
l_x^\tau|^2+R+\frac{l_x^\tau-n}{\tau}\right)+2l^\tau_x\triangle
\phi\Bigr]\ d{\rm vol}(g(t))\le 0.$$ Furthermore, equality holds in
either of these weak inequalities for all functions $\phi$ as above
and all $\tau$ if and only if it holds in both. In that case $l_x$
is a smooth function on space-time and the equalities hold in the
usual smooth sense.
\end{thm}

\begin{rem}
The terms in these inequalities are interpreted in the following
way: First of all, $\nabla l^\tau_x$ and $\triangle l^\tau_x$ are
computed using only the spatial derivatives (i.e., they are
horizontal differential operators). Secondly, since $l_x$ is a
locally Lipschitz functions defined on all of $M\times (0,T)$ we
have seen that $\partial l_x/\partial t$ and $|\nabla l^\tau_x|^2$
are continuous functions on the open subset ${\mathcal U}_x(\tau)$
of full measure in $M\times\{\tau\}$ and furthermore, that they are
locally bounded on all of $M\times\{\tau\}$ in the sense that for
any $q\in M$ there is a neighborhood $V$ of $q$ such that the
restriction of $|\nabla l_x^\tau|^2$ to $V\cap {\mathcal U}_x(\tau)$
is bounded. This means that $\partial l_x/\partial t$ and $|\nabla
l_x^\tau|^2$ are elements of $L_{\rm loc}^\infty(M)$ and hence can
be integrated against any smooth function with compact support. In
particular, they are distributions.

 Since
$\nabla l^\tau_x$ is a smooth, locally bounded vector field on an
open subset of full measure, for any compactly supported test
function $\phi$, integration by parts yields
$$\int\triangle \phi\cdot l^\tau_x\ d{\rm vol}(g(t))=-\int \langle \nabla \phi,\nabla
l_x^\tau\rangle\ d{\rm vol}(g(t)).$$ Thus, formulas in
Theorem~\ref{weak} can also be taken to mean:
$$\int_{M\times \{\tau\}} \Bigl[\phi\cdot(\frac{\partial l_x}{\partial \tau}
+|\nabla l^\tau_x|^2-R+\frac{n}{2\tau})+ \langle \nabla
l_x^\tau,\nabla \phi \rangle\Bigr] d{\rm vol}(g(t))\ge 0$$
$$\int_{M\times\{\tau\}} \Bigl[\phi\cdot(-|\nabla
l^\tau_x|^2+R+\frac{l^\tau_x-n}{\tau})-2\langle\nabla
l^\tau_x,\nabla\phi \rangle\Bigr] d{\rm vol}(g(t))\le 0.$$
\end{rem}

The rest of this subsection is devoted to the proof of these
inequalities. We fix $(M,g(t)),\ 0\le t\le T<\infty$, as in the
statement of the theorem.
 We
fix $x$ and denote by $L$ and $l$ the functions $L_{x}$ and $l_{x}$.
We also fix $\tau$, and we denote by $L^\tau$ and $l^\tau$ the
restrictions of $L$ and $l$ to the slice $M\times \{T-\tau\}$. We
begin with a lemma.

\begin{lem}\label{phisupport}
There is a continuous function $C\colon M\times (0,T)\to \Ar$ such
that for each point $(q,t)\in M\times (0,T)$, setting $\tau=T-t$,
the following holds. There is an upper barrier $\varphi_{(q,t)}$ for
$L_x^\tau$ at the point $q$ defined on a neighborhood $U_{(q,t)}$ of
$q$ in $M$ satisfying $|\nabla \varphi_{(q,t)}(q)|\le C(q,t)$ and
$${\rm Hess}(\varphi)(v,v)\le C(q,t)|v|^2$$
for all tangent vectors $v\in T_qM$.
\end{lem}

\begin{proof}
By Proposition~\ref{lipcomplete}, $L$ is a locally Lipschitz
function on $M\times (0,T)$, and in particular is continuous. The
bound $C(q,t)$ will depend only on the bounds on curvature and its
first two derivatives and on the function $L(q,t)$. Fix $(q,t)$ and
let $\gamma$ be a minimizing ${\mathcal L}$-geodesic from $x$ to
$(q,t)$. (The existence of such a minimizing geodesic is established
in Lemma~\ref{geoexist}.) Fix $\tau_1>0$, with $\tau_1<(T-t)/2$, let
$t_1=T-\tau_1$, and let $q_1=\gamma(\tau_1)$. Consider
$\varphi_{(q,t)}={\mathcal L}(\gamma_{[0,\tau_1]})+
L^\tau_{(q_1,t_1)}$. This is an upper barrier  for $L_x^\tau$ at $q$
defined in some neighborhood $V\subset M$ of $q$.  Clearly, $\nabla
\varphi_{(q,t)}=\nabla L^\tau_{(q_1,t_1)}$ and ${\rm
Hess}(\varphi_{(q,t)}={\rm Hess}( L^\tau_{(q_1,t_1)})$.

According to Corollary~\ref{DL} we have $\nabla
L^\tau_{(q_1,t_1)}(q)=2\sqrt{\tau}X_{\gamma}(\tau)$. On the other
hand, by Corollary~\ref{minmax} there is a bound on
$\sqrt{\tau}|X_{\gamma}(\tau)|$ depending only on the bounds on
curvature and its first derivatives, on $\tau$ and $\tau_1$ and on
$l_x(q,\tau)$. Of course, by Shi's theorem (Theorem~\ref{firstshi})
for every $\epsilon>0$ the norms of the first derivatives of
curvature on $M\times [\epsilon,T]$ are bounded in terms of
$\epsilon$ and the bounds on curvature.
 This
proves that $|\nabla \varphi_{(q,t)}(q)|$ is bounded by a continuous
function $C(q,t)$ defined on all of $M\times (0,T)$.

 Now consider Inequality~(\ref{Hessin}) for $\gamma$ at $\bar\tau=\tau$.
 It is clear that the first
 two terms on the right-hand side are bounded by $C|Y(\tau)|^2$,
 where $C$ depends on the curvature bound and on $T-t$. We consider
 the last term, $\int_{\tau_1}^{\tau}\sqrt{\tau'}H(X,Y)d\tau'$.
 We claim that this integral is also bounded by $C'|Y(\tau)|^2$
 where $C'$ depends on the bounds on curvature and its first and
 second derivatives along $\gamma_1$ and on $T-t$.
We consider $\tau'\in [\tau_1,\tau]$. Of course,
 $\sqrt{\tau'}|X(\tau')|$ is bounded on this interval. Also,
 $$|Y(\tau')|=\left(\frac{\sqrt{\tau'}-\sqrt{\tau_1}}{\sqrt{\tau}
 -\sqrt{\tau_1}}\right)|Y(\tau)|\le \frac{\sqrt{\tau'}}{\sqrt{\tau}}|Y(\tau)|.$$
Hence $|Y(\tau')|/\sqrt{\tau'}$ and $|Y(\tau')||X(\tau')|$ are
bounded in terms of $T-t$, $|Y(\tau)|$, and the bound on
$\sqrt{\tau'}|X(\tau')|$ along the ${\mathcal L}$-geodesic. From
this it follows immediately from Equation~(\ref{Heqn}) that $H(X,Y)$
is bounded along the ${\mathcal L}$-geodesic  by $C|Y(\bar\tau)|^2$
 where the constant $C$ depends on $T-t$ and the
 bounds on curvature and its first two derivatives.
\end{proof}

Of course, if $(q,t)\in {\mathcal U}_{x}(\tau)$, then this argument
shows that the Hessian of $L^\tau_{x}$ is bounded near $(q,t)$.

At this point in the proof of Theorem~\ref{weak} we wish to employ
arguments using convexity. To carry these out we find it convenient
to work with a Euclidean metric and usual convexity rather than the
given metric $g(t)$ and convexity measured using $g(t)$-geodesics.
In order to switch to a Euclidean metric we must find one that well
approximates $g(t)$. The following is straightforward to prove.

\begin{claim}
For each point $(q,t)\in M\times(0,T)$ there is an open metric ball
$B_{(q,t)}$ centered at $q$ in $(M,g(t))$ which is the diffeomorphic
image of a ball $\widetilde B\subset T_qM$ under the exponential map
for $g(t)$ centered at $q$ such that the following hold:
\begin{enumerate}
\item $B_{(q,t)}\subset U_{(q,t)}$ so that the upper barrier $\varphi_{(q,t)}$
from Lemma~\ref{phisupport} is defined on all of $B_{(q,t)}$.
\item The constants $C(z,t)$ of Lemma~\ref{phisupport} satisfy
$C(z,t)\le 2C(q,t)$ for all $z\in B_{(q,t)}$.
\item The push-forward, $h$, under the exponential mapping of the Euclidean
metric on $T_qM$ satisfies $$h/2\le g\le 2h.$$
\item The Christoffel symbols $\Gamma_{ij}^k$ for the
metric $g(t)$ written using the Gaussian normal coordinates (the
image under the exponential mapping of orthonormal linear
coordinates on $T_qM$) are bounded in absolute value by
$1/(8n^3C(q,t))$ where $n$ is the dimension of $M$.
\end{enumerate}
\end{claim}

Instead of working in the given metric $g(t)$ on $B_{(q,t)}$ we
shall use the Euclidean metric $h$ as in the above claim. For any
function $f$ on $B_{(q,t)}$ we denote by ${\rm Hess}(f)$ the Hessian
of $f$ with respect to the metric $g(t)$ and by ${\rm Hess}^h(f)$
the Hessian of $f$ with respect to the metric $h$. By
Formula~(\ref{Hessian}), for any $z\in B_{(q,t)}$ and any $v\in
T_zM$, we have
$${\rm Hess}(\varphi_{(z,t)})(v,v)={\rm Hess}^h(\varphi_{(z,t)})(v,v)
-\sum_{i,j,k}v^iv^j\Gamma_{ij}^k\frac{\partial
\varphi_{(z,t)}}{\partial x^k}.$$ Thus, it follows from the above
assumptions on the $\Gamma_{ij}^k$ and the bound on $|\nabla
\varphi_{(z,t)}|$ that for all $z\in B_{(q,t)}$ we have
\begin{equation}\label{Hessdiff}
\left|{\rm Hess}(\varphi_{(z,t)})(v,v)-{\rm
Hess}^h(\varphi_{(z,t)})(v,v)\right|\le
\frac{1}{4}|v|^2_h,\end{equation} and hence for every $z\in
B_{(q,t)}$ we have
$${\rm Hess}^h(\varphi_{(z,t)})(v,v)\le 2C(q,t)|v|^2_g+\frac{|v|_h^2}{4}\le
\left(4C(q,t)+\frac{1}{4}\right)|v|_h^2.$$

This means:

\begin{claim}
For each $(q,t)\in M\times(0,T)$ there is a smooth function
$$\psi_{(q,t)}\colon B_{(q,t)}\to \Ar$$
 with the property that at each
$z\in B_{(q,t)}$ there is an upper barrier $b_{(z,t)}$ for
$L^\tau+\psi_{(q,t)}$ at $z$ with
$${\rm Hess}^h(b_{(z,t)})(v,v)\le -3|v|_h^2/2$$
for all $v\in T_zM$.
\end{claim}

\begin{proof}
 Set
$$\psi_{(q,t)}=-(2C(q,t)+1)d^2_h(q,\cdot).$$
Then  for any $z\in B_{(q,t)}$ the function
$b_{(z,t)}=\varphi_{(z,t)}+\psi_{(q,t)}$ is an upper barrier  for
$L^\tau+\psi_{(q,t)}$ at $z$. Clearly, for all $v\in T_zM$ we have
$${\rm Hess}^h(b_{(z,t)})(v,v)={\rm Hess}^h(\varphi_{(z,t)})(v,v)+{\rm Hess}^h(\psi_{(q,t)})(v,v)\le -3|v|_h^2/2.$$
\end{proof}

This implies that if $\alpha\colon [a,b]\to B_{(q,t)}$ is any
Euclidean straight-line segment in $B_{(q,t)}$ parameterized by
Euclidean arc length and if $z= \alpha(s)$ for some $s\in (a,b)$,
then
$$(b_{(z,t)}\circ\alpha)''(s)\le -3/2.$$

\begin{claim}
Suppose that $\beta\colon [-a,a]\to \Ar$ is a continuous function
and that at each $s\in (-a,a)$ there is an upper barrier $\hat b_s$
for $\beta$  at $s$ with $\hat b_s''\le -3/2$. Then
$$\frac{\beta(a)+\beta(-a)}{2}\le \beta(0)-\frac{3}{4}a^2.$$
\end{claim}

\begin{proof}
Fix $c<3/4$ and define a continuous function
$$A(s)=\frac{(\beta(-s)+\beta(s))}{2}+cs^2-\beta(0)$$ for $s\in [0,a]$.
Clearly, $A(0)=0$. Also, using the upper barrier at $0$ we see that
for $s>0$ sufficiently small $A(s)<0$. For any $s\in (0,a)$ there is
an upper barrier $c_s=(\hat b_{s}+\hat b_{-s})/2+cs^2-\beta(0)$ for
$A(s)$ at $s$, and $c_s''(t)\le 2c-3/2<0$. By the maximum principle
this implies that $A$ has no local minimum in $(0,a)$, and
consequently that it is a non-increasing function of $s$ on this
interval. That is to say, $A(s)<0$ for all $s\in (0,a)$ and hence
$A(a)\le 0$, i.e., $(\beta(a)+\beta(-a))/2+ca^2\le \beta(0)$. Since
this is true for every $c<3/4$, the result follows.
\end{proof}

Now applying this to Euclidean intervals in $B_{(q,t)}$ we conclude:

\begin{cor}\label{lineineq}
For any $(q,t)\in M\times (0,T)$, the function
$$\beta_{(q,t)}=L^\tau+\psi_{(q,t)}\colon B_{(q,t)}\to\Ar$$ is
uniformly strictly convex with respect to $h$. In fact, let
$\alpha\colon [a,b]\to B_{(q,t)}$ be a Euclidean geodesic arc. Let
$y,z$ be the endpoints of $\alpha$, let $w$ be its midpoint, and let
$|\alpha|$ denote the length of this arc (all defined using the
Euclidean metric). We have
$$\beta_{(q,t)}(w)\ge \frac{\left(\beta_{(q,t)}(y)+\beta_{(q,t)}(z)\right)}{2}+\frac{3}{16}|\alpha|^2.$$
\end{cor}

What follows is a simple interpolation result (see \cite{GreeneWu}).
For each $q\in M$ we let $B_{(q,t)}'\subset B_{(q,t)}$ be a smaller
ball centered at $q$, so that $B'_{(q,t)}$ has compact closure in
$B_{(q,t)}$.

\begin{claim}
Fix $(q,t)\in M\times(0,T)$, and let $\beta_{(q,t)}\colon B_{(q,t)}\to \Ar$ be
as above. Let $S\subset M$ be the singular locus of $L^\tau$, i.e.,
$S=M\setminus {\mathcal U}_{x}(\tau)$. Set $S_{(q,t)}=B_{(q,t)}\cap S$. Of
course, $\beta_{(q,t)}$ is smooth on $B_{(q,t)}\setminus S_{(q,t)}$. Then there
is a sequence of smooth functions $\{f_k\colon B'_{(q,t)}\to
\Ar\}_{k=1}^\infty$ with the following properties:
\begin{enumerate}
\item As $k\rightarrow \infty$ the functions $f_k$ converge uniformly to $\beta_{(q,t)}$ on
$B_{(q,t)}'$.
\item For any $\epsilon>0$ sufficiently small, let $\nu_\epsilon(S_{(q,t)})$ be the
$\epsilon$-neighborhood (with respect to the Euclidean metric) in
$B_{(q,t)}$ of $S_{(q,t)}\cap B_{(q,t)}$. Then, as $k\rightarrow
\infty$ the restrictions of $f_k$ to $B_{(q,t)}'\setminus
\left(B_{(q,t)}'\cap \nu_\epsilon(S_{(q,t)})\right)$ converge
uniformly in the $C^\infty$-topology to the restriction of
$\beta_{(q,t)}$ to this subset.
\item For each $k$, and for any $z\in B_{(q,t)}'$ and any $v\in T_zM$ we have
$${\rm Hess}(f_k)(v,v)\le -|v|_{g(t)}^2/2.$$
That is to say, $f_k$ is strictly convex with respect to the metric
$g(t)$.
\end{enumerate}
\end{claim}

\begin{proof}
Fix $\epsilon>0$ sufficiently small so that for any $z\in
B_{(q,t)}'$ the Euclidean $\epsilon$-ball centered at $z$ is
contained in $B_{(q,t)}$. Let $B_0$ be the ball of radius $\epsilon$
centered at the origin in $\Ar^n$ and let $\xi\colon B_0\to \Ar$ be
a non-negative $C^\infty$-function with compact support and with
$\int_{B_0}\xi d{\rm vol}_h=1$. We define
$$\beta_{(q,t)}^\epsilon(z)=\int_{B_0}\xi(y)\beta_{(q,t)}(z+y)dy,$$
for all $z\in B_{(q,t)}'$. It is clear that for each $\epsilon>0$
sufficiently small, the function $\beta_{(q,t)}^\epsilon\colon
B_{(q,t)}'\to \Ar$ is $C^\infty$ and that as $\epsilon\rightarrow 0$
the $\beta_{(q,t)}^\epsilon$ converge uniformly on $B_{(q,t)}'$ to
$\beta_{(q,t)}$. It is also clear that for every $\epsilon>0$
sufficiently small, the conclusion of Corollary~\ref{lineineq} holds
for $\beta_{(q,t)}^\epsilon$ and for each Euclidean straight-line
segment $\alpha$ in $B_{(q,t)}'$. This implies that ${\rm
Hess}^h(\beta_{(q,t)}^\epsilon)(v,v)\le -3|v|_h^2/2$, and hence that
by Inequality~(\ref{Hessdiff}) that
$${\rm Hess}(\beta_{(q,t)}^\epsilon)(v,v)\le -|v|^2_h=-|v|^2_{g(t)}/2.$$
This means that $\beta_{(q,t)}^\epsilon$ is convex with respect to
$g(t)$. Now take a sequence $\epsilon_k\rightarrow 0$ and let
$f_k=\beta_{(q,t)}^{\epsilon_k}$.

Lastly, it is a standard fact that $f_k$ converge uniformly in the
$C^\infty$-topology to $\beta_{(q,t)}$ on any subset of $B_{(q,t)}'$
whose closure is disjoint from $S_{(q,t)}$.
\end{proof}

\begin{defn}
For any continuous function $\psi$ defined on $B_{(q,t)}'\setminus
\left(S_{(q,t)}\cap B_{(q,t)}'\right)$ we define
$$\int_{(B_{(q,t)}')^*}\psi d{\rm vol}(g(t))={\rm lim}_{\epsilon\rightarrow
0}\int_{B_{(q,t)}'\setminus \nu_\epsilon(S_{(q,t)})\cap
B_{(q,t)}'}\psi d{\rm vol}(g(t)).$$
\end{defn}

We now have:

\begin{claim}\label{wineq}
Let $\phi\colon B_{(q,t)}'\to \Ar$ be a non-negative, smooth
function with compact support. Then
$$\int_{B_{(q,t)}'}\beta_{(q,t)}\triangle \phi d{\rm vol}(g(t))\le \int_{(B_{(q,t)}')^*}\phi
\triangle \beta_{(q,t)} d{\rm vol}(g(t)).$$
\end{claim}

\begin{rem}
Here $\triangle$ denotes the Laplacian with respect to the metric
$g(t)$.
\end{rem}

\begin{proof}
Since $f_k\rightarrow \beta_{(q,t)}$ uniformly on $B_{(q,t)}'$ we
have
$$\int_{B_{(q,t)}'}\beta_{(q,t)}\triangle \phi d{\rm vol}(g(t))={\rm lim}_{k\rightarrow
\infty}\int_{B_{(q,t)}'}f_k\triangle \phi  d{\rm vol}(g(t)).$$ Since
$f_k$ is strictly convex with respect to the metric $g(t)$,
$\triangle f_k\le 0$ on all of $B_{(q,t)}'$. Since $\phi\ge 0$, for
every $\epsilon$ and $k$ we have
$$\int_{\nu_\epsilon(S_{(q,t)})\cap B_{(q,t)}'}\phi\triangle f_k  d{\rm vol}(g(t))\le
0.$$ Hence, for every $k$ and for every $\epsilon$ we have
\begin{eqnarray*}
\int_{B_{(q,t)}'}f_k\triangle \phi d{\rm vol}(g(t)) & = &
\int_{B_{(q,t)}'}\phi\triangle f_k  d{\rm vol}(g(t)) \\
& \le & \int_{B_{(q,t)}'\setminus \left(B_{(q,t)}'\cap
\nu_\epsilon(S_{(q,t)})\right)}\phi\triangle f_k d{\rm
vol}(g(t)).\end{eqnarray*}
 Taking the limit as $k\rightarrow\infty$,
using the fact that $f_k\rightarrow \beta_{(q,t)}$ uniformly on
$B_{(q,t)}'$ and that restricted to
$B_{(q,t)}'\setminus(B_{(q,t)}'\cap \nu_\epsilon (S_{(q,t)}))$ the
$f_k$ converge uniformly in the $C^\infty$-topology to
$\beta_{(q,t)}$ yields
$$\int_{B_{(q,t)}'}\beta_{(q,t)}\triangle \phi d{\rm vol}(g(t))
\le \int_{B_{(q,t)}'\setminus \left(B_{(q,t)}'\cap
\nu_\epsilon(S_{(q,t)})\right)}\phi\triangle \beta_q d{\rm
vol}(g(t)).$$ Now taking the limit as $\epsilon\rightarrow 0$
establishes the claim.
\end{proof}

\begin{cor}\label{ltauineq}
Let $\phi\colon B_{(q,t)}'\to \Ar$ be a non-negative, smooth
function with compact support. Then
$$\int_{B_{(q,t)}'}l^\tau\triangle \phi d{\rm vol}(g(t))\le \int_{(B_{(q,t)}')^*}\phi
\triangle l^\tau d{\rm vol}(g(t)).$$
\end{cor}

\begin{proof}
Recall that $\beta_{(q,t)}=L^\tau+\psi_{(q,t)}$ and that
$\psi_{(q,t)}$ is a $C^\infty$-function. Hence,
$$\int_{B_{(q,t)}'}\psi_{(q,t)}\triangle \phi d{\rm vol}(g(t))= \int_{(B_{(q,t)}')^*}\phi
\triangle \psi_{(q,t)} d{\rm vol}(g(t)).$$ Subtracting this equality
from the inequality in the previous claim and dividing by
$2\sqrt{\tau}$ gives the result.
\end{proof}

Now we turn to the proof proper of Theorem~\ref{weak}.

\begin{proof}
Let $\phi\colon M\to \Ar$ be a non-negative, smooth function of
compact support. Cover $M$ by open subsets of the form $B_{(q,t)}'$
as above. Using a partition of unity we can write
$\phi=\sum_i\phi_i$ where each $\phi_i$ is a non-negative smooth
function supported in some $B_{(q_i,t)}'$. Since the inequalities we
are trying to establish are linear in $\phi$, it suffices to prove
the result for each $\phi_i$. This allows us to assume (and we shall
assume) that $\phi$ is supported in $B_{(q,t)}'$ for some $q\in M$.

Since $l_x^\tau$ is a locally Lipschitz function, the restriction of
$|\nabla l_x^\tau|^2$ to $B_{(q,t)}'$ is an $L_{\rm
loc}^\infty$-function. Similarly, $\partial l_x/\partial \tau$ is an
$L_{\rm loc}^\infty$-function. Hence
\begin{eqnarray*} \lefteqn{\int_{B_{(q,t)}'}
\phi\cdot\left(\frac{\partial l_x}{\partial \tau}+|\nabla
l_x^\tau|^2-R+\frac{n}{2\tau}\right) d{\rm vol}(g(t))} \\ & = &
\int_{(B_{(q,t)}')^*}\phi\left(\frac{\partial l_x}{\partial
\tau}+|\nabla l_x^\tau|^2-R+\frac{n}{2\tau}\right) d{\rm
vol}(g(t)).\end{eqnarray*} On the other hand, by
Corollary~\ref{ltauineq} we have
$$\int_{B_{(q,t)}'}l_x^\tau\triangle \phi  d{\rm vol}(g(t))\le \int_{(B_{(q,t)}')^*}\phi\triangle l_x^\tau
d{\rm vol}(g(t)).$$ Putting these together we see
\begin{eqnarray*}
\lefteqn{\int_{B_{(q,t)}'}\phi\left(\frac{\partial l_x}{\partial
\tau}+|\nabla
l_x^\tau|^2-R+\frac{n}{2\tau}\right) -l_x^\tau\triangle \phi  d{\rm vol}(g(t))} \\
& \ge & \int_{(B_{(q,t)}')^*}\phi\left(\frac{\partial l_x}{\partial
\tau}+|\nabla l_x^\tau|^2-R+\frac{n}{2\tau} -\triangle
l_x^\tau\right) d{\rm vol}(g(t)).\end{eqnarray*}

It follows immediately from the second inequality in
Corollary~\ref{lformula} that, since $\phi\ge 0$ and
$\left(B'_{(q,t)}\right)^*\subset {\mathcal U}_x(\tau)$, we have
$$\int_{(B_{(q,t)}')^*}\phi\left(\frac{\partial l_x}{\partial \tau}+|\nabla l_x^\tau|^2-R+
\frac{n}{2\tau}-\triangle l_x^\tau\right) d{\rm vol}(g(t))\ge 0.$$
This proves the first inequality in the statement of the theorem.

The second inequality in the statement of the theorem is proved in
exactly the same way using the third inequality in
Corollary~\ref{lformula}.

Now let us consider the distributions
$$D_1=\frac{\partial l_x}{\partial \tau}+|\nabla l_x^\tau|^2-R+\frac{n}{2\tau}- \triangle
l_x^\tau $$ and
$$D_2=-|\nabla
l_x^\tau|^2+R+\frac{l_x^\tau-n}{\tau}+2\triangle l_x^\tau
$$ on $M\times\{\tau\}$. According to Corollary~\ref{lformula} the following equality
holdes on ${\mathcal U}_x(\tau)$:
$$2\frac{\partial l_x}{\partial \tau}+|\nabla l_x^\tau|^2-R+\frac{l_x^\tau}{\tau}=0.$$
By Proposition~\ref{lipcomplete} the open set ${\mathcal U}_x(\tau)$
has full measure in $M$ and $|\nabla l_x^\tau|^2$ and $\partial
l_x/\partial\tau$ are locally essentially bounded. Thus, this
equality is an equality of locally essentially bounded, measurable
functions, i.e., elements of $L_{\rm loc}^\infty(M)$, and hence is
an equality of distributions on $M$. Subtracting $2D_1$ from this
equality yields $D_2$. Thus,
$$D_2=-2D_1,$$
as distributions on $M$. This shows that $D_2$ vanishes as a
distribution if and only if $D_1$ does. But if $D_2=0$ as a
distribution for some $\tau$, then by elliptic regularity $l_x^\tau$
is smooth on $M\times\{\tau\}$ and the equality is the na\"ive one
for smooth functions. Thus, if $D_2=0$ for all $\tau$, then
$l_x^\tau$ and $\partial l/\partial \tau$ are $C^\infty$ functions
on each slice $M\times\{\tau\}$ and both $D_1$ and $D_2$ hold in the
na\"ive sense on each slice $M\times \{\tau\}$. It follows from a
standard bootstrap argument that in this case $l_x^\tau$ is smooth
on all of space-time.
\end{proof}

\section{Reduced volume}

We have established that for a Ricci flow $(M,g(t)),\ 0\le t\le T$,
and a point  $x=(p,T)\in M\times [0,T]$ the reduced length function
$l_x$ is defined on all of $M\times (0,T)$. This allows us to
defined the reduced volume of $M\times \{\tau\}$ for any $\tau\in
(0,T)$ Recall that the {\sl reduced volume}\index{reduced volume} of
$M$ is defined to be
$$\widetilde V_x(M,\tau) = \int_{M} \tau^{-\frac{n}{2}}exp(-l_x(q,\tau))dq.$$
This function is defined for $0<\tau<T$.

There is one simple case where we can make an explicit computation.
\begin{lem}\label{flat}
If $(M,g(t))$ is flat Euclidean $n$-space (independent of $t$), then
for any $x\in \Ar^n\times (-\infty,\infty)$ we have
$$\widetilde V_x(M,\tau)=(4\pi)^{n/2}$$ for all $\tau>0$.
\end{lem}

\begin{proof}
By symmetry we can assume that $x=(0,T)\in \Ar^n\times [0,T]$, where
$0\in \Ar^n$ is the origin. We have already seen that the ${\mathcal
L}$-geodesics in flat space are the usual geodesics when
parameterized by $s=\sqrt{\tau}$. Thus, for any $X\in
\Ar^n=T_0\Ar^n$ $\gamma_X(\tau)=2\sqrt{\tau}X$, and hence ${\mathcal
L}{\rm exp}(X,\bar\tau)=2\sqrt{\bar\tau}X$. This means that for any
$\tau>0$ and any $X\in T_0\Ar^n$ we have ${\mathcal U}(\tau)=T_pM$,
and ${\mathcal J}(X,\tau)=2^n\tau^{n/2}$. Also,
 $L_x(X,\tau)=2\sqrt{\tau}|X|^2$, so that $l_x(X,\tau)=|X|^2$.
 Thus, for any $\tau>0$
 $$\widetilde
 V_x(\Ar^n,\tau)=\int_{\Ar^n}\tau^{-n/2}e^{-|X|^2}2^n\tau^{n/2}dX=(4\pi)^{n/2}.$$
\end{proof}

In the case when $M$ is non-compact, it is not clear {\em a priori}
that the integral defining the reduced volume is finite in general.
In fact, as the next proposition shows, it is always finite and
indeed, it is bounded above by the integral for  $\Ar^n$.

\begin{thm}\label{4pi} Let $(M,g(t)),\ 0\le t\le T$, be a Ricci flow
of bounded curvature with the property that for each $t\in [0,T]$
the Riemannian manifold $(M,g(t))$ is complete. Fix a point
$x=(p,T)\in M\times [0,T]$. For every $0<\tau<T$ the reduced volume
$$\widetilde V_x(M,\tau)=\int_{M} \tau^{-\frac{n}{2}}{\rm exp}(-l_x(q,\tau))dq$$
is absolutely convergent and  $\widetilde V_x(M,\tau)\le
(4\pi)^{\frac{n}{2}}$. The function $\widetilde V_x(M,\tau)$ is a
non-increasing function of $\tau$ with
$${\rm lim}_{\tau\rightarrow 0}\widetilde V_x(M,\tau)=(4\pi)^{\frac{n}{2}}.$$
\end{thm}

\begin{proof}
By Proposition~\ref{lipcomplete} ${\mathcal U}_x(\tau)$ is an open
subset of full measure in $M$. Hence, $$\widetilde
V_x(M,\tau)=\int_{{\mathcal U}_x(\tau)}\tau^{-\frac{n}{2}}{\rm
exp}(-l_x(q,\tau))dq.$$
 Take
linear orthonormal coordinates $(z^1,\ldots,z^n)$ on $T_pM$. It
follows from  the previous equality and Lemma~\ref{redvoltilde} that
$$\widetilde V_x(M,\tau)=\int_{\widetilde {\mathcal U}_{x}(\tau)}f(Z,\tau)dz^1\cdots dz^n,$$
where $f(Z,\tau)=\tau^{-\frac{n}{2}}e^{-\widetilde
l(Z,\tau)}{\mathcal J}(Z,\tau)$. By Proposition~\ref{fmonotone} for
each $Z$ the integrand, $f(Z,\tau)$, is a non-increasing function of
$\tau$ and the function converges uniformly on compact sets as
$\tau\rightarrow 0$ to $2^ne^{-|Z|^2}$. This implies that
$f(Z,\tau)\le 2^ne^{-|Z|^2}$ for all $\tau>0$, and hence that
$$\int_{\widetilde{\mathcal U}_x(\tau)}f(Z,\tau)dz^1\ldots dz^n$$
converges absolutely for each $\tau>0$, and the integral has value
at most $(4\pi)^{n/2}$.

Fix $0<\tau_0<T$. According to Theorem~\ref{Amono} (with $A=M\times
(T-\tau_0,T)$), the reduced volume $\widetilde V_x(M,\tau)$ is a
non-increasing function of $\tau$ on $(0,\tau_0]$. Since this is
true for any $0<\tau_0<T$, it follows that $\widetilde V_x(M,\tau)$
is a non-increasing function of $\tau$ for all $\tau\in (0,T)$.
 (This of course is a consequence of the
monotonicity of $f(Z,\tau)$ in $\tau$ and the fact that $\widetilde {\mathcal
U}_x(\tau)\subset \widetilde {\mathcal U}_x(\tau')$ for $\tau'<\tau$.)

To show that ${\rm lim}_{\tau\rightarrow 0}\widetilde
 V_x(M,\tau)= (4\pi)^{n/2}$ we need only see that for
  each $A<\infty$ for all $\tau>0$ sufficiently small
 $\widetilde {\mathcal U}_x(\tau)$ contains the ball of radius $A$ centered at the origin
 in $T_pM$. Since the  curvature is bounded, this is exactly
 the content of Corollary~\ref{Jlimit1}.
\end{proof}

\subsection{Converse to Lemma~\protect{\ref{flat}}}

In Lemma~\ref{flat} we showed that for the trivial flow on flat
Euclidean $n$-space and for any point $x\in \Ar^n\times \{T\}$ the
reduced volume $\widetilde V_x(\Ar^n,\tau)$ is independent of
$\tau>0$  and is equal to $(4\pi)^{n/2}$. In this subsection we use
the monotonicity results of the last subsection to establish the
converse to Lemma~\ref{flat}, namely to show that if $(M,g(t)),0\le
t\le T$, is a Ricci flow complete with bounded curvature and if
$\widetilde V_x(M,\bar\tau)=(4\pi)^{n/2}$ for some $\bar\tau>0$ and
some $x\in M\times\{T\}$, then the flow on the interval
$[T-\bar\tau,T]$ is the trivial flow on flat Euclidean $n$-space.

\begin{prop}\label{4picase} Suppose that $(M,g(\tau)),\ 0\le \tau\le T$,
 is a solution to the backward Ricci flow
equation, complete of bounded curvature. Let $x=(p,T)\in M\times
\{T\}$, and suppose that $0<\bar\tau<T$. If $\widetilde
V_x(M,\bar\tau)=(4\pi)^{n/2}$, then the backward Ricci flow on the
interval $[0,\bar\tau]$ is the trivial flow on flat Euclidean space.
\end{prop}

\begin{proof}
If $\widetilde V_x(M,\bar\tau)=(4\pi)^{n/2}$, then by
Lemma~\ref{flat}, $\widetilde V_x(M,\tau)$ is constant on the
interval $(0,\bar\tau]$. Hence, it follows from the proof of
Theorem~\ref{4pi} that the closure of $\widetilde{\mathcal U}(\tau)$
is all of $T_pM$ for all $\tau\in (0,\bar\tau]$ and that
$f(Z,\tau)=e^{-|Z|^2}2^n$ for all $Z\in T_pM$ and all $\tau\le
\bar\tau$. In particular,
$$\frac{\partial {\rm ln}(f(Z,\tau))}{\partial \tau}=0.$$
From the proof of
 Proposition~\ref{fmonotone} this means that
Inequality~(\ref{Jinequ}) is an equality and consequently, so is
Inequality~(\ref{basic}). Thus, by Proposition~\ref{Hessineq}
 (with $\tau_1=0$) each of the vector fields $Y_\alpha(\tau)=\tilde
Y_\alpha(\tau)$ is both a Jacobi field and adapted. By
Proposition~\ref{triL} we then have
$${\rm Ric}+{\rm Hess}(l^\tau_x)=\frac{g}{2\tau}.$$
In particular, $l_x$ is smooth. Let $\varphi_\tau\colon M\to M,\
0<\tau\le \bar\tau$,\ be the one-parameter family of diffeomorphisms
obtained by solving
$$\frac{d\varphi_\tau}{d\tau}=\nabla
l_x(\cdot,\tau)\ \ \ \ {\rm and} \ \ \ \varphi_{\bar\tau}={\rm
Id}.$$ We now consider
$$h(\tau)=\frac{\bar\tau}{\tau}\varphi_\tau^*g(\tau).$$
We compute \begin{eqnarray*}\frac{\partial h}{\partial \tau} & = &
-\frac{\bar\tau}{\tau^2}\varphi_\tau^*g(\tau)
+\frac{\bar\tau}{\tau}\varphi_\tau^*{\mathcal
L}_{\frac{d\varphi_\tau}{d\tau}}(g(\tau))
+\frac{\bar\tau}{\tau}\varphi_t^*2{\rm Ric}(g(\tau)) \\
& = & -\frac{\bar\tau}{\tau^2}\varphi_\tau^*g(\tau) +
\frac{\bar\tau}{\tau}\varphi_\tau^*2{\rm
Hess}(l_x^\tau)+\frac{\bar\tau}{\tau}
\varphi_\tau^*\left(\frac{1}{\tau}g(\tau)-2{\rm
Hess}(l_x^\tau)\right)=0.\end{eqnarray*} That is to say the family
of metrics $h(\tau)$ is  constant in $\tau$: for all $\tau\in
(0,\bar\tau]$ we have $h(\tau)=h(\bar\tau)=g(\bar\tau)$. It then
follows that
$$g(\tau)=\frac{\tau}{\bar\tau}(\varphi_\tau^{-1})^*g(\bar\tau),$$
which means that the entire flow in the interval $(0,\bar\tau]$
differs by diffeomorphism and scaling from $g(\bar\tau)$. Suppose
that $g(\bar\tau)$ is not flat, i.e., suppose that there is some
$(x,\bar\tau)$ with $|{\rm Rm}(x,\bar\tau)|=K>0$. Then from the flow
equation we see that $|{\rm
Rm}(\varphi^{-1}_\tau(x),\tau)|=K\bar\tau^2/\tau^2$, and these
curvatures are not bounded as $\tau\rightarrow 0$. This is a
contradiction. We conclude that $g(\bar\tau)$ is flat, and hence,
again by the flow equation so are all the $g(\tau)$ for $0<\tau\le
\bar\tau$, and by continuity, so is $g(0)$. Thus, $(M,g(\tau))$ is
isometric to a quotient of $\Ar^n$ by a free, properly discontinuous
group action. Lastly, since $\widetilde V_x(M,\tau)=(4\pi)^{n/2}$,
it follows that $(M,g(\tau))$ is isometric to $\Ar^n$ for every
$\tau\in [0,\bar\tau]$. Of course, it then follows that the flow is
the constant flow.
\end{proof}

\chapter{Non-collapsed results}\label{noncoll}

In this chapter we apply the results for the reduced length function
and reduced volume established in the last two sections to prove
non-collapsing results. In the first section we give a general
result that applies to generalized Ricci flows and will eventually
be applied to Ricci flows with surgery to prove the requisite
non-collapsing. In the second section we give a non-collapsing
result for Ricci flows on compact $3$-manifolds with normalized
initial metrics.

\section{A non-collapsing result for generalized Ricci flows}

The main result of this chapter is a
$\kappa$-non-collapsed\index{$\kappa$-non-collapsed} result.

 \begin{thm}\label{THM}
 Fix  positive constants  $\bar\tau_0<\infty$, $l_0<\infty$, and  $V>0$.
 Then there is $\kappa>0$ depending on $\bar\tau_0$, $V$, and $l_0$
 and the dimension $n$ such
that the following holds. Let $({\mathcal M},G)$ be a generalized
$n$-dimensional Ricci flow, and let $0<\tau_0\le \bar\tau_0$. Let
$x\in {\mathcal M}$ be fixed. Set $T={\bf t}(x)$. Suppose that
$0<r\le \sqrt{\tau_0}$ is given. These data are required to satisfy:
\begin{enumerate} \item[(1)] The ball $B(x,T,r)\subset M_T$ has compact
closure. \item[(2)] There is an embedding $B(x,T,r)\times
[T-r^2,T]\subset {\mathcal M}$ compatible with ${\bf t}$ and with
the vector field. \item[(3)]  $|{\rm Rm}|\le r^{-2}$ on the image of
the embedding in (2). \item[(4)] There is an open subset $\widetilde
W\subset \widetilde {\mathcal U}_x(\tau_0)\subset T_x M_T$ with the
property that for every ${\mathcal L}$-geodesic $\gamma\colon
[0,\tau_0]\to {\mathcal M}$ with initial condition contained in
$\widetilde W$, the $l$-length of $\gamma$ is at most $l_0$.
\item[(5)] For each $\tau\in [0,\tau_0]$, let $W(\tau)={\mathcal
L}{\rm exp}_x^\tau(\widetilde W)$. The volume of the image
$W(\tau_0)\subset M_{T-\tau_0}$ is at least $V$.
\end{enumerate}
Then
$${\rm Vol}(B(x,T,r))\ge \kappa r^n.$$
 \end{thm}

\begin{figure}[ht]
  \centerline{\epsfbox{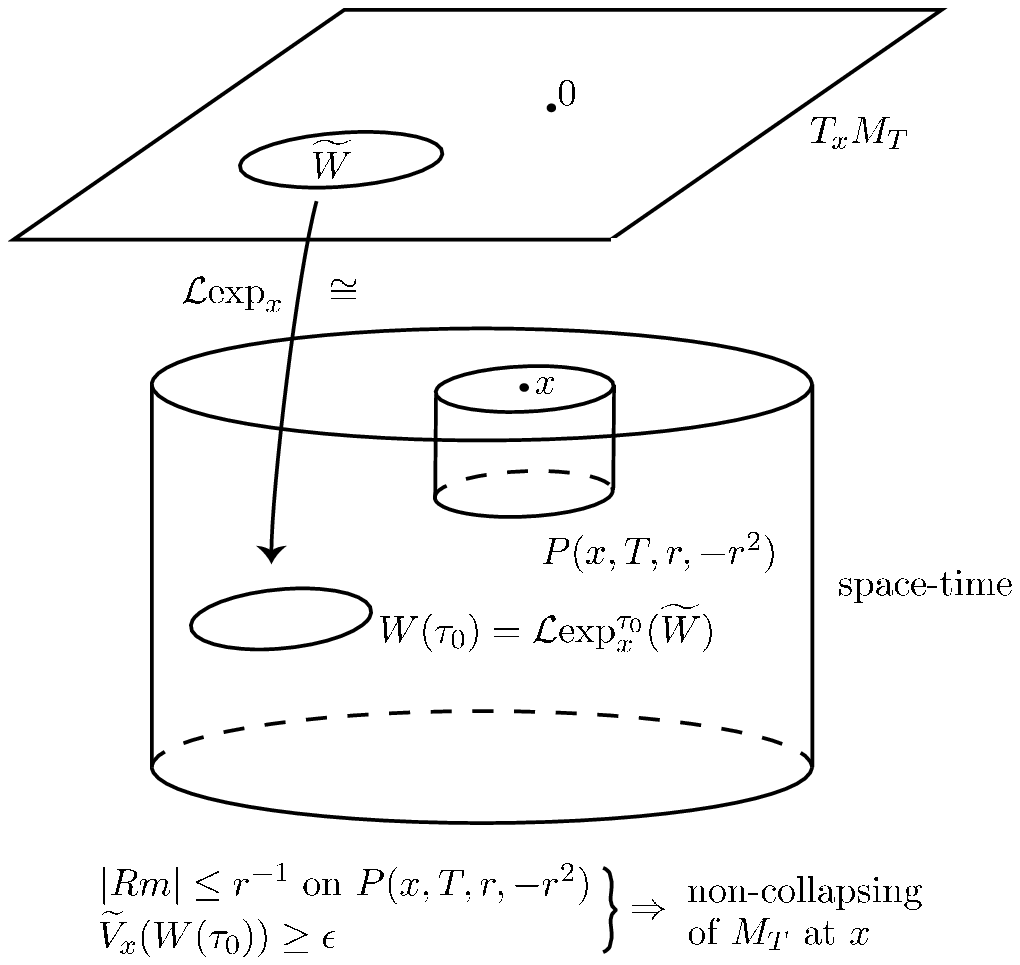}}
  \caption{Non-collapsing.}\label{fig:8.1}
\end{figure}

See {\sc Fig.}~\ref{fig:8.1}.

In this section we denote by $g(\tau), \ 0\le \tau\le r^2$, the
family of metrics on $B(x,T,r)$ induced from pulling back $G$ under
the embedding $B(x,T,r)\times [T-r^2,T]\to {\mathcal M}$. Of course,
this family of metrics satisfies the backward Ricci flow equation.

\begin{proof}
Clearly from the definition of the reduced volume\index{reduced
volume}, we have
$$\widetilde V_x(W(\tau_0))\ge \tau_0^{-n/2}Ve^{-l_0}\ge
\bar\tau_0^{-n/2}Ve^{-l_0}.$$\index{$\widetilde V_x$|(} By the
monotonicity result (Theorem~\ref{Amono}) it follows that for any
$\tau\le \tau_0$, and in particular for any $\tau\le r^2$, we have
\begin{equation}\label{redvoleqn} \widetilde V_x(W(\tau))\ge
\bar\tau_0^{-n/2}Ve^{-l_0}.\end{equation}

Let $\varepsilon=\sqrt[n]{{\rm Vol}(B(x,T,r))}/r$, so that ${\rm
Vol}\,B(x,T,r)=\varepsilon^nr^n$.  The basic result we need to
establish in order to prove this theorem is the following:

\begin{prop}\label{step} There is a positive constant
$\varepsilon_0\le 1/4n(n-1)$ depending on $\bar\tau_0$ and $l_0$
such that if $\varepsilon\le \varepsilon_0$ then, setting
$\tau_1=\varepsilon r^2$, we have
 $\widetilde V_x(W(\tau_1)) < 3\varepsilon^{\frac{n}{2}}$.
\end{prop}

Given this proposition, it follows immediately that either
$\varepsilon>\varepsilon_0$ or
$$\varepsilon\ge \left(\frac{\widetilde V_x(W(\tau_1))}{3}\right)^{2/n}\ge
\frac{1}{3^{2/n}\bar\tau_0}V^{2/n}e^{-2l_0/n}.$$ Since
$\kappa=\varepsilon^n$, this proves the theorem.

\begin{proof}  We divide $\widetilde W$
into
$$\widetilde W_{\rm sm}=\widetilde W\cap \left\{Z\in
T_xM_T\bigl|\bigr. |Z|\le \frac{1}{8}\varepsilon^{-1/2}\right\}$$
and
$$\widetilde W_{\rm lg}=\widetilde W\setminus \widetilde W_{\rm sm},$$
(see {\sc Fig. 8.2}).

\begin{figure}[ht]
  \centerline{\epsfbox{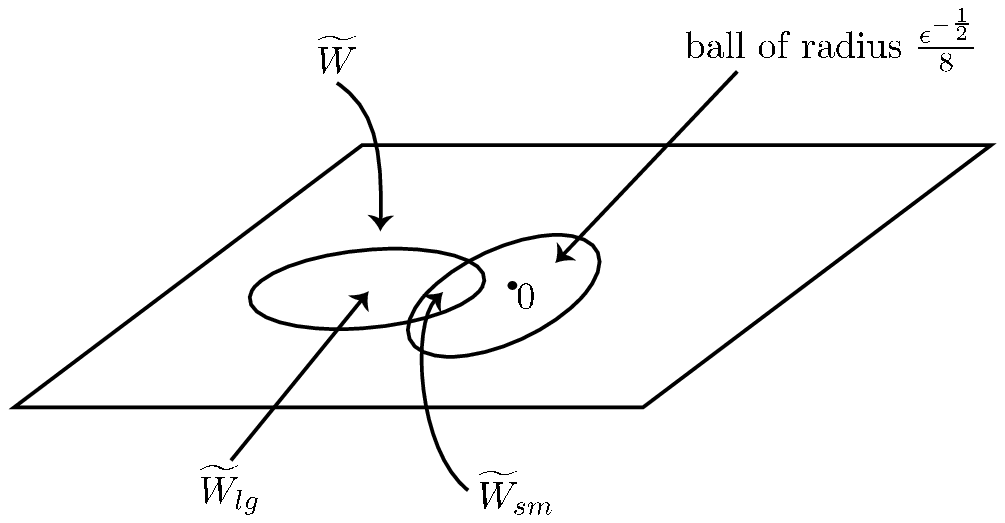}}
  \caption{$\widetilde{W}_{lg}$ and $\widetilde{W}_{sm}$.}
\end{figure}

We set $W_{\rm sm}(\tau_1)={\mathcal L}{\rm
exp}_x^{\tau_1}(\widetilde W_{\rm sm})$ and $W_{\rm
lg}(\tau_1)={\mathcal L}{\rm exp}_x^{\tau_1}(\widetilde W_{\rm
lg}).$ Clearly, since $W(\tau_1)$ is the union of $W_{\rm
sm}(\tau_1)$ and $W_{\rm lg}(\tau_1)$ and since these subsets are
disjoint measurable subsets, we have
$$\widetilde V_x(W(\tau_1))=\widetilde V_x(W_{\rm sm}(\tau_1))+\widetilde V_x(W_{\rm
lg}(\tau_1)).$$  We shall show that there is $\varepsilon_0$ such that either
$\varepsilon>\varepsilon_0$ or  $\widetilde V_x(W_{\rm sm}(\tau_1))\le
2\varepsilon^{n/2}$ and $\widetilde V_x(W_{\rm lg}(\tau_1))\le
\varepsilon^{n/2}$. This will establish Proposition~\ref{step} and hence
Theorem~\ref{THM}.

\subsection{Upper bound for  $\widetilde V_x(W_{\rm sm}(\tau_1))$}

The idea here is that ${\mathcal L}$-geodesics with initial vector
in $\widetilde W_{\rm sm}$ remain in the parabolic neighborhood
$P=B(x,T,r)\times [T-r^2,T]$ for $\tau\in [0,r^2]$. Once we know
this it is easy to see that their ${\mathcal L}$-lengths are bounded
from below. Then if the volume of $B(x,T,r)$ was arbitrarily small,
the reduced volume of  $W_{\rm sm}(\tau_1)$ would be arbitrarily
small.

\begin{lem}\label{small}
Setting $\tau_1=\varepsilon r^2$, there is a constant
$\varepsilon_0>0$ depending on $\bar\tau_0$ such that, if
$\varepsilon\le \varepsilon_0$, we have
$$\int_{\widetilde
W_{\rm sm}}\tau_1^{-n/2}e^{-\tilde l(Z,\tau_1)}{\mathcal
J}(Z,\tau_1)dZ\le 2\varepsilon^{\frac{n}{2}}.$$
\end{lem}

Of course,  we have
$$\widetilde V_x(W_{\rm sm}(\tau_1))=\int_{W_{\rm
sm}(\tau_1)}\tau_1^{-n/2}e^{-l(q,\tau_1)}d{\rm
vol}_{g(\tau_1)}=\int_{\widetilde W_{\rm sm}}\tau_1^{-n/2}e^{-\tilde
l(Z,\tau_1)}{\mathcal J}(Z,\tau_1)dZ,$$ so that it will follow
immediately from the lemma that:

\begin{cor}
There is a  constant $\varepsilon_0>0$ depending on $\bar\tau_0$
such that, if $\varepsilon\le \varepsilon_0$, we have
$$\widetilde V_x(W_{\rm sm}(\tau_1))\le 2\varepsilon^{\frac{n}{2}}.$$
\end{cor}

\begin{proof}{\em (Of Lemma~\ref{small})}
In order to establish Lemma~\ref{small} we need two preliminary
estimates:

\begin{claim}\label{nabR}
There is a universal positive constant $\varepsilon'_0$ such that, if
$\varepsilon\le \varepsilon'_0$, then there is a constant $C_1<\infty$
depending only on the dimension $n$ such that the following hold for all $y\in
B(x,T,r/2)$, and for all $t\in[T-\tau_1,T]$:
\begin{enumerate}
\item[(1)] $$|\nabla R(y,t)|\le \frac{C_1}{r^3}$$
\item[(2)]  $$(1-C_1\varepsilon)\le \frac{g(y,t)}{g(y,T)}\le (1+C_1\varepsilon).$$
\end{enumerate}
\end{claim}

\begin{proof}
Recall that by hypothesis $|{\rm Rm}(y,t)|\le 1/r^2$ on
$B(x,T,r)\times [T-r^2,T]$. Rescale the flow by multiplying the
metric and time by $r^{-2}$ resulting in a ball $\widetilde B$ of
radius one and a flow defined for a time interval of length one with
$|{\rm Rm}|\le 1$ on the entire parabolic neighborhood
$B(x,T,1)\times [T-1,T]$. Then according to Theorem~\ref{shi} there
is a universal constant $C_1$ such that $|\nabla R(y,t)|\le C_1$ for
all $(y,t)\in B(x,T,1/2)\times [T-1/2,T]$. Rescaling back by $r^2$
to the original flow, we see that on this flow $|\nabla R(y,t)|\le
C_1/r^3$ for all $(y,t)\in B(x,T,r/2)\times [T-r^2/2,T]$. Taking
$\varepsilon'_0\le 1/2$ gives the first item in the claim.

 Since  $|{\rm Ric}|\le (n-1)/r^2$ for all
$(y,t)\in B\times[T-r^2,T]$ it follows by integrating that
$$e^{-2(n-1)(T-t)/r^2}\le \frac{g(x,t)}{g(x,T)}\le e^{2(n-1)(T-t)/r^2}.$$
Thus, for $t\in [T-\tau_1,T]$ we have
$$e^{-2(n-1)\varepsilon}\le \frac{g(x,t)}{g(x,T)}\le e^{2(n-1)\varepsilon}.$$
From this the second item in the claim is immediate.
\end{proof}

At this point we view the ${\mathcal L}$-geodesics as paths
$\gamma\colon [0,\tau_1]\to B(x,T,r)$ (with the understanding that
the path in space-time is given by the composition of the path $
(\gamma(\tau),T-\tau)$ in $B(x,T,r)\times [T-r^2,T]$ followed by the
given inclusion of this product into ${\mathcal M}$.

The next step in the  proof is to show that for any $Z\in \widetilde
W_{\rm sm}$ the ${\mathcal L}$-geodesic $\gamma_Z$ (the one having
${\rm lim}_{\tau\rightarrow 0}\sqrt\tau X_{\gamma_Z}(\tau)=Z$)
remains in $B(x,T,r/2)$ up to time $\tau_1$. Because of this, as we
shall see, these paths contribute a small amount to the reduced
volume since $B(x,T,r/2)$ has small volume. We set
$X(\tau)=X_{\gamma_Z}(\tau)$

\begin{claim} There is  a
positive constant $\varepsilon_0\le 1/4n(n-1)$ depending on
$\bar\tau_0$, such that the following holds. Suppose that
 $\varepsilon\le \varepsilon_0$  and $\tau_1' \leq
\tau_1=\varepsilon r^2$. Let $Z\in T_xM_T$ and let $\gamma_Z$ be the
associated ${\mathcal L}$-geodesic from $x$. Suppose  that
$\gamma_Z(\tau) \in B(x,T,r/2)$ for all $\tau < \tau_1'$. Then for
all $\tau<\tau_1'$ we have
$$\abs{|\sqrt{\tau}X(\tau)|_{g(T)} - |Z|} \leq 2\varepsilon (1+|Z|).$$
\end{claim}

\begin{proof}
First we make sure that $\varepsilon_0$ is less than or equal to the
universal constant $\varepsilon_0'$ of the last claim. For all
$(y,t)\in B(x,T,r)\times [T-r^2,T]$ we have $|{\rm Rm}(y,t)|\le
r^{-2}$ and $|\nabla R(y,t)|\le C_1/r^3$ for some universal constant
$C_1$. Of course, $r^2\le \bar\tau$. Thus, at the expense of
replacing $C_1$ by a larger constant, we can (and shall) assume that
$C_1/r^3>(n-1)r^{-2}\ge |{\rm Ric}(y,t)|$ for all $(y,t)\in
B(x,T,r)\times [T-r^2,T]$. Thus, we can take the constant $C_0$ in
the hypothesis of Lemma~\ref{min-max} to be $C_1/r^3$. We take the
constant $\bar\tau$ in the hypothesis of that lemma to be
$\varepsilon r^2$. Then, we have that
$${\rm max}_{0\le \tau\le \tau_1'}\sqrt{\tau}|X(\tau)|\le
e^{2C_1\varepsilon^2}|Z|+\frac{e^{2C_1\varepsilon^2}-1}{2}\sqrt{\varepsilon}
r,$$ and
$$|Z|\le e^{2C_1\varepsilon^2}{\rm min}_{0\le \tau\le
\tau_1'}\sqrt{\tau}|X(\tau)|+\frac{e^{2C_1\varepsilon^2}-1}{2}\sqrt{\varepsilon}r.$$
By choosing $\varepsilon_0>0$ sufficiently small (as determined by
the universal constant $C_1$ and by $\bar\tau_0$), we have
$${\rm max}_{0\le \tau\le \tau_1'}\sqrt{\tau}|X(\tau)|_{g(T-\tau)}\le
(1+\frac{\varepsilon}{2})|Z|+\frac{\varepsilon}{2},$$ and
$$|Z|\le (1+\frac{\varepsilon}{2}){\rm min}_{0\le \tau\le
\tau_1'}\sqrt{\tau}|X(\tau)|_{g(T-\tau)}+\frac{\varepsilon}{2}.$$ It
is now immediate that
$$\abs{|\sqrt{\tau}X(\tau)|_{g(T-\tau)} - |Z|} \leq \varepsilon (1+|Z|).$$

Again choosing $\varepsilon_0$ sufficiently small the result now
follows from the second inequality in Claim~\ref{nabR}
\end{proof}

Now we are ready to establish that the ${\mathcal L}$-geodesics
whose initial conditions are elements of $\widetilde W_{\rm sm}$ do
not leave $B(x,T,r/2)\times [T-r^2,T]$ for any  $\tau\le\tau_1$.

\begin{claim}\label{inside}
Suppose $\varepsilon_0\le 1/4n(n-1)$ is the constant from the last
claim. Set $\tau_1=\varepsilon r^2$, and suppose that
$\varepsilon\le \varepsilon_0$. Lastly, assume that $|Z|\le
\frac{1}{8\sqrt{\varepsilon}}$. Then $\gamma_Z(\tau) \in B(x,T,r/2)$
for all $\tau\le \tau_1$.
\end{claim}

\begin{proof}
Since $\varepsilon\le \varepsilon_0\le 1/4n(n-1)\le 1/8$, by the
last claim we have
$$|\sqrt{\tau}X(\tau)|_{g(T)}\le (1+2\varepsilon )|Z|+2\varepsilon\le
\frac{5}{4}|Z|+\frac{3}{32\sqrt\varepsilon},$$ provided that
$\gamma|_{[0,\tau)}$ is contained in $B(x,T,r/2)\times [T-\tau,T]$.
Since $|Z|\le (8\sqrt\varepsilon)^{-1}$ we conclude that
$$|\sqrt{\tau}X(\tau)|_{g(T)}\le \frac{1}{4\sqrt\varepsilon},$$
as long as $\gamma([0,\tau))$ is contained in $B(x,T,r/2)\times
[T-\tau,T]$.

 Suppose that there is
$\tau'<\tau_1=\varepsilon r^2$ for which $\gamma_Z$ exits
$B(x,T,r/2)\times [T-r^2,T]$. We take $\tau'$ to be the first such
time. Then we have
$$
|\gamma_Z(\tau')- x|_{g(T)} \leq
\int_{0}^{\tau'}|X(\tau)|_{g(T)}d\tau \le
\frac{1}{4\varepsilon^{\frac{1}{2}}}\int_{0}^{\tau'}\frac{d\tau}{\sqrt{\tau}}
=\frac{1}{2\varepsilon^{\frac{1}{2}}}\sqrt{\tau'} < r/2.
$$
This contradiction implies that  $\gamma_Z(\tau) \in B(x,T,r/2)$ for
all $\tau< \tau_1=\varepsilon r^2$.
\end{proof}

Now we assume that $\varepsilon_0>0$ depending on $\bar\tau_0$ is as
above and that $\varepsilon\le \varepsilon_0$, and we shall estimate
$$
\widetilde V_x(W_{\rm sm}(\tau_1)) = \int_{W_{\rm sm}(\tau_1)}
(\tau_1)^{-\frac{n}{2}}e^{- l(q,\tau_1)}d{\rm vol}_{g(\tau_1)}.$$ In
order to do this we estimate $ l_x(q,\tau_1)$ on $ W_{\rm
sm}(\tau_1)$. By hypothesis $|{\rm Rm}|\le 1/r^2$ on
$B(x,T,r/2)\times [0,\tau_1]$ and by Lemma~\ref{inside} every
${\mathcal L}$-geodesic $\gamma_Z$, defined on $[0,\tau_1]$, with
initial conditions $Z$ satisfying $|Z|\leq
\frac{1}{8}\varepsilon^{-\frac{1}{2}}$ remains in $B(x,T,r/2)$.
Thus, for such $\gamma_Z$ we have $R(\gamma_Z(\tau))\ge
-n(n-1)/r^2$. Thus, for any $q\in W_{\rm sm}(\tau_1)$ we have
$$ L_x(q,\tau_1) =
\int_{0}^{\tau_1}\sqrt{\tau}(R+|X(\tau)|^{2})d \tau \geq
-\frac{2n(n-1)}{3r^{2}}(\tau_1)^{\frac{3}{2}}=
-\frac{2n(n-1)}{3}\varepsilon^{\frac{3}{2}}r,$$ and hence
$$l+x(q,\tau_1)=\frac{L_x(q,\tau_1)}{2\sqrt{\tau_1}}\ge -\frac{n(n-1)}{3}\varepsilon.$$ Since $W_{\rm
sm}(\tau)\subset B(x,T,r/2)\subset B(x,T,r)$, we have:
\begin{eqnarray}\label{volineq}
 \widetilde V_x(W_{\rm sm}(\tau_1)) & \leq &
\varepsilon^{-\frac{n}{2}}r^{-n}e^{n(n-1)\varepsilon/3}{\rm
Vol}_{g(T-\tau_1)}\, W_{\rm sm}(\tau) \\
& \le & \varepsilon^{-\frac{n}{2}}r^{-n}e^{n(n-1)\varepsilon/3}{\rm
Vol}_{g(T-\tau_1)}\,B(x,T,r).\nonumber\end{eqnarray}

\begin{claim}\label{volvomp}
There is a universal constant $\varepsilon_0>0$ such that if
$\varepsilon\le \varepsilon_0$, for any open subset $U$ of
$B(x,T,r)$, and for any $0\le \tau_1\le \tau_0$, we have
$$0.9\le {\rm Vol}_{g(T)}U/{\rm Vol}_{g(T-\tau_1)}U\le 1.1.$$
\end{claim}

\begin{proof}
This is immediate from the second item in Claim~\ref{nabR}.
\end{proof}

Now assume that $\varepsilon_0$ also satisfies this claim. Plugging
this into Equation~(\ref{volineq}), and using the fact that
$\varepsilon\le \varepsilon_0\le 1/4n(n-1)$, , so that
$n(n-1)\varepsilon/3\le 1/12$, and the fact that from the definition
we have ${\rm Vol}_{g(T)}\,B(x,T,r)=\varepsilon^nr^n$, gives
$$\widetilde V_x(W_{\rm sm}(\tau_1))\leq \varepsilon^{-\frac{n}{2}}r^{-n}e^{n(n-1)\varepsilon/3}
(1.1){\rm Vol}_{g(T)}B(x,T,r)\le(1.1)
\varepsilon^{\frac{n}{2}}e^{\frac{1}{12}}.$$ Thus,
$$\widetilde V_x(W_{\rm sm}(\tau_1))\le
2\varepsilon^{\frac{n}{2}}.$$

This completes the proof of Lemma~\ref{small}.
\end{proof}

\subsection{Upper bound for $\widetilde V_x( W_{\rm lg}(\tau_1))$}

Here the basic point is to approximate the reduced volume integrand
by the heat kernel, which drops off exponentially fast as we go away
from the origin.

Recall that ${\rm Vol}\,B(x,T,r)=\varepsilon^nr^n$ and
$\tau_1=\varepsilon r^2$.

\begin{lem}\label{large}
There is a universal positive constant $\varepsilon_0>0$ such that
if $\varepsilon\le \varepsilon_0$, we have
$$
\widetilde V_x(W_{\rm lg}(\tau_1)) \le  \int_{{\widetilde
U}(\tau_1)\cap \{Z\bigl|\bigr. |Z|\geq
\frac{1}{8}\varepsilon^{-\frac{1}{2}}\}}(\tau_1)^{-\frac{n}{2}}
e^{-\tilde l(q,\tau_1)}\mathcal{J}(Z,\tau_1)dZ\le
\varepsilon^{\frac{n}{2}}.$$
\end{lem}

\begin{proof}
By the monotonicity result (Proposition~\ref{fmonotone}), we see
that the restriction of the function
$\tau_1^{-\frac{n}{2}}e^{-\tilde l(Z,\tau_1)}{\mathcal J}(Z,\tau_1)$
to ${\widetilde U}(\tau_1)$ is less than or equal to the restriction
of the function $2^ne^{-|Z|^2}$ to the same subset. This means that
$$\widetilde V_x(W_{\rm lg}(\tau_1))\le \int_{{\widetilde U}(\tau_1)\setminus
\widetilde U(\tau_1)\cap B(0,\frac{1}{8}\varepsilon^{-1/2})}2^n
e^{-|Z|^2}dZ\le \int_{T_p{\mathcal M}_T\setminus
B(0,\frac{1}{8}\varepsilon^{-1/2})}2^ne^{-|Z|^2}dZ.$$ So it suffices
to estimate this latter integral.

Fix some $a>0$ and  let $I(a)=\int_{B(0,a)}2^ne^{-|Z|^2}dZ$. Let
$R(a/\sqrt{n})$ be the $n$-cube centered at the origin with side
lengths $2a/\sqrt{n}$. Then $R(a/\sqrt{n})\subset B(0,a)$, so that
\begin{eqnarray*}
I(a) & \ge  & \int_{R(a/\sqrt{n})}2^ne^{-|Z|^2}dZ \\
& = & \prod_{i=1}^n
\left(\int_{-a/\sqrt{n}}^{a/\sqrt{n}}2e^{-z_i^2}dz_i\right)
\\
& = &
\left(\int_0^{2\pi}\int_0^{a/\sqrt{n}}4e^{-r^2}rdrd\theta\right)^{n/2}.
\end{eqnarray*}

Now
$$\int_0^{2\pi}\int_0^{a/\sqrt{n}}4e^{-r^2}rdrd\theta=4\pi(1-e^{-\frac{a^2}{n}}).$$
Applying this with $a=(8\sqrt\varepsilon)^{-1}$ we have
$$\widetilde V_x(W_{\rm lg}(\tau_1))\le \int_{\Ar^n}2^ne^{-|Z|^2}dZ-I(1/8\sqrt\varepsilon)
\le(4\pi)^{n/2}\left(1-\left(1-e^{-1/(64n\varepsilon)}\right)^{n/2}\right).$$
Thus,
$$\widetilde V_x(W_{\rm lg}(\tau_1))\le
(4\pi)^{n/2}\frac{n}{2}e^{-1/(64n\varepsilon)}.$$\index{$\widetilde
V_x$} There is $\varepsilon_0>0$ such that the expression on the
right-hand side is less than $\varepsilon^{n/2}$ if $\varepsilon\le
\varepsilon_0$. This completes the proof of Lemma~\ref{large}.
\end{proof}

Putting Lemmas~\ref{small} and~\ref{large} together, establishes
Proposition~\ref{step}. \end{proof}

As we have already remarked, Proposition~\ref{step} immediately
implies Theorem~\ref{THM}. This completes the proof of
Theorem~\ref{THM}.
\end{proof}

\section{Application to compact Ricci flows}

Now let us apply this result to Ricci flows with normalized initial
metrics to show that they are universally $\kappa$-non-collapsed on
any fixed, finite time interval. In this section we specialize to
$3$-dimensional Ricci flows. We do not need this result in what
follows for we shall prove a more delicate result in the context of
Ricci flows with surgery. Still, this result is much simpler and
serves as a paradigm of what will come.

\begin{thm}
Fix positive constants $\omega>0$ and $T_0<\infty$. Then there is
$\kappa>0$ depending only on these constants such that the following
holds. Let $(M,g(t)),\ 0\le t<T\le T_0$, be a $3$-dimensional Ricci
flow with $M$ compact. Suppose that $|{\rm Rm}(p,0)|\le 1$ and also
that ${\rm Vol}\,B(p,0,1)\ge \omega$ for all $p\in M$. Then for any
$t_0\le T$, any  $r>0$ with $r^2\le t_0$ and any $(p,t_0)\in M\times
\{t_0\}$, if $|{\rm Rm}(q,t)|\le r^{-2}$ on $B(p,t_0,r)\times
[t_0-r^2,t_0]$ then ${\rm Vol}\,B(p,t_0,r)\ge \kappa r^3$.
\end{thm}

\begin{figure}[ht]
  \relabelbox{
  \centerline{\epsfbox{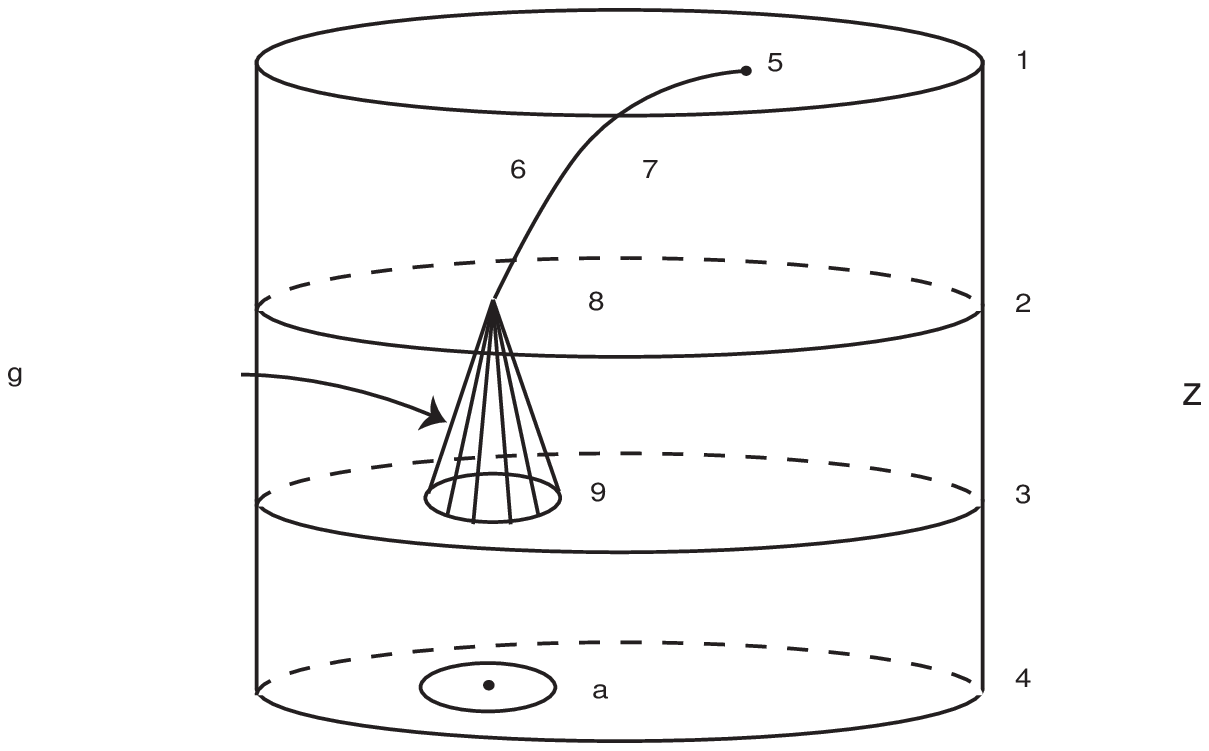}}}
  \relabel{1}{$M\times \{t_0\}$}
  \relabel{2}{$M\times \{2^{-4}\}$}
  \relabel{3}{$M\times \{2^{-5}\}$}
  \relabel{4}{$M\times \{0\}$}
  \relabel{5}{$(q,t_0)$}
  \relabel{6}{$\gamma$}
  \relabel{7}{$l(\gamma)\leq \frac 32$}
  \relabel{8}{$(q,2^{-4})$}
  \relabel{9}{$B(q,0,1)\times \{2^{-5}\}$}
  \relabel{a}{$B(q,0,1)$}
  \relabel{g}{$g(0)$-geodesics}
  \endrelabelbox
  \caption{Non-collapsing of Ricci flows.}\label{fig:ncrf}
\end{figure}

\begin{proof}
Fix any $x=(p,t_0)\in M\times [0,T]$. First, we claim that we can
suppose that $t_0\ge 1$. For if not, then rescale the flow by
$Q=1/t_0$. This does not affect the curvature inequality at time
zero. Furthermore, there is $\omega'>0$ depending only on $\omega$
such that for any ball $B$ at time zero and of radius one in the
rescaled flow we have ${\rm Vol}\,B\ge \omega'$. The reason for the
latter fact is the following: By the Bishop-Gromov inequality
(Theorem~\ref{BishopGromov}) there is $\omega'>0$ depending only on
$\omega$ such that for any $q\in M$ and any $r\le 1$ we have ${\rm
Vol}\,B(q,0,r)\ge \omega'r^3$. Of course, the rescaling increases
$T$, but simply restrict to the rescaled flow on $[0,1]$.

Next, we claim that we can assume that $r\le \sqrt{t_0}/2$. If $r$
does not satisfy this inequality, then we replace $r$ with
$r'=\sqrt{t_0}/2$. Of course, the curvature inequalities hold for
$r'$ if they hold for $r$. Suppose that we have established the
result for $r'$. Then
$${\rm Vol}\,B(p,T,r)\ge {\rm Vol}\,B(p,T,r')\ge \kappa(r')^3\ge
\kappa\left(\frac{r}{2}\right)^3=\frac{\kappa}{8}r^3.$$

From now on we assume that $t_0\ge 1$ and $r\le\sqrt{t_0}/2$.
According to Proposition~\ref{kappa0r0t0} for any $(p,t)\in M\times
[0,2^{-4}]$ we have $|{\rm Rm}(p,t)|\le 2$ and ${\rm
Vol}\,B(p,t,r)\ge \kappa_0 r^3$ for all $r\le 1$.

Once we know that $|{\rm Rm}|$ is universally bounded on $M\times
[0,2^{-4}]$ it follows that there is a universal constant $C_1$ such
that $C_1^{-1}g(q,0)\le g(q,t) \le C_1g(q,0)$ for all $q\in M$ and
all $t\in [0,2^{-4}]$. This means that there is a universal constant
$C<\infty $ such that the following holds. For any points $q_0,q\in
M$ with $d_0(q_0,q)\le 1$ let $\gamma_{q_0,q}$ be the path in
$M\times [2^{-5},2^{-4}]$ given by
$$\gamma_{q_0,q}(\tau)=(A_{q_0,q}(\tau),2^{-4}-\tau),\ 0\le \tau\le 2^{-5},$$
where $A_{q_0,q}$ is a shortest $g(0)$-geodesic from $q_0$ to $q$.
Then ${\mathcal L}(\gamma_{q_0,q})\le C$.

By Theorem~\ref{n/2} there is a point $q_0\in M$ and an ${\mathcal
L}$-geodesic $\gamma_0$ from $x=(p,t_0)$ to $(q_0,2^{-4})$ with
$l(\gamma_0)\le 3/2$. Since $t_0\ge 1$, this means that there is a
universal constant $C'<\infty$ such that for each point $q\in
B(q_0,0,1)$ the path which is the composite of $\gamma_0$ followed
by $\gamma_{q_0,q}$ has $\tilde l$-length at most $C'$. Setting
$\tau_0=t_0-2^{-5}$, this implies that $l_x(q,\tau_0)\le C'$ for
every $q\in B(q_0,0,1)$. This ball has volume at least $\kappa_0$.
By Proposition~\ref{lipcomplete}, the open subset ${\mathcal
U}_x(\tau_0)$ is of full measure in $M\times\{2^{-5}\}$. Hence,
$W(\tau_0)=\left(B(q_0,0,1)\times \{2^{-5}\}\right)\cap {\mathcal
U}_x(\tau_0)$ also has volume at least $\kappa_0$. Since $r^2\le
t_0/4<\tau_0$, Theorem~\ref{THM} now gives the result. (See {\sc
Fig.}~\ref{fig:ncrf}.)
\end{proof}

\chapter{$\kappa$-non-collapsed ancient solutions}\label{kappasect}

 In this chapter we discuss the qualitative properties of
$\kappa$-non-collapsed, ancient solutions. One of the most important is the
existence of a gradient shrinking soliton that is asymptotic at $-\infty$ to
the solution. The other main qualitative result is the compactness result (up
to scaling) for these solutions. Also extremely important for us is
classification of $3$-dimensional gradient shrinking solitons -- up to finite
covers there are only two: a shrinking family of round $S^3$'s and a shrinking
family of products of round $S^2$'s with $\Ar$. This leads to a rough
classification of all $3$-dimensional $\kappa$-non-collapsed, ancient
solutions. The $\kappa$-solutions are in turn the models for singularity
development in $3$-dimensional Ricci flows on compact manifolds, and eventually
for singularity development in $3$-dimensional Ricci flows with surgery.

\section{Preliminaries}

Our objects of study are Ricci flows $(M,g(t)),\ -\infty<t\le 0$,
with each $(M,g(t))$ being a complete manifold of bounded
non-negative curvature. The first remark to make is that the
appropriate notion of non-negative curvature is that the Riemann
curvature operator\index{Riemann curvature operator}
$${\rm Rm}\colon \wedge^2TM\to \wedge^2TM,$$
which is a symmetric operator, is non-negative. In general, this
implies, but is stronger than, the condition that the sectional
curvatures are all non-negative. In case the dimension of $M$ is at
most three, every element of $\wedge^2TM$ is represented by a
$2$-plane (with area form) and hence the Riemann curvature operator
is non-negative if and only if all the sectional curvatures are
non-negative.  In the case of non-negative curvature operator,
bounded curvature operator is equivalent to bounded scalar
curvature.

It follows immediately from the Ricci flow equation that since the $(M,g(t))$
have non-negative Ricci curvature,  the metric is non-increasing in time in the
sense that for any point $p\in M$ and any $v\in T_pM$ the function
$|v|^2_{g(t)}$ is a non-increasing function of $t$.

There are  stronger results under the assumption of bounded, non-negative
curvature operator. These are consequences of the Harnack inequality (see
\cite{Hamiltonharnack}). As was established in Corollary~\ref{posderiv}, since
the flow exists for $t\in (-\infty,0]$ and since the curvature operator is
non-negative and bounded for each $(q,t)\in M\times (-\infty,0]$, it follows
that $\partial R(q,t)/\partial t\ge 0$ for all $q$ and $t$. That is to say, for
each $q\in M$ the scalar curvature $R(q,t)$  is a non-decreasing function of
$t$.

\subsection{Definition}

Now we turn to the definition of what it means for a Ricci flow to
be $\kappa$-non-collapsed.

\begin{defn}\label{kappadefn}
Fix $\kappa>0$. Let $(M,g(t)),\ a<t\le b$, be a Ricci flow of
complete $n$-manifolds. Fix $r_0>0$. We say that $(M,g(t))$ is {\em
$\kappa$-non-collapsed\index{$\kappa$-non-collapsed|ii} on scales at
most $r_0$} if the following holds for any $(p,t)\in M\times (a,b]$
and any $0<r\le r_0$  with the property that $a\le t-r^2$. Whenever
$|{\rm Rm}(q,t')|\le r^{-2}$ for all $q\in B(p,t,r)$ and all $t'\in
(t-r^2,t]$, then ${\rm Vol}\,B(p,t,r)\ge \kappa r^n$. We say that
$(M,g(t))$ is {\em $\kappa$-non-collapsed}, or equivalently {\em
$\kappa$-non-collapsed on all scales} if it is
$\kappa$-non-collapsed on scales at most $r_0$ for every
$r_0<\infty$.
\end{defn}

\begin{defn}
An {\em ancient solution}\index{ancient solution|ii} is a Ricci flow
 $(M,g(t))$ defined for $-\infty <t \leq 0$ such
that for each $t$, $(M,g(t))$ is a connected, complete, non-flat
Riemannian manifold whose curvature operator is bounded and
nonnegative. For any $\kappa>0$, an ancient solution is {\sl
$\kappa$-non-collapsed} if it is $\kappa$-non-collapsed on all
scales. We also use the terminology $\kappa$-{\sl solution} for a
$\kappa$-non-collapsed, ancient solution.
\end{defn}

Notice that a $\kappa$-solution is a $\kappa'$-solution for any
$0<\kappa'\le \kappa$.

\subsection{Examples}

Here are some examples of $\kappa$-solutions:
\begin{exam}
Let $(S^{2},g_{0})$ be the standard round $2$-sphere of scalar
curvature $1$ (and hence Ricci tensor $g_0/2$). Set
$g(t)=(1-t)g_{0}$. Then $\partial g(t)/\partial t=-2{\rm
Ric}(g(t))$, $-\infty<t\le 0$. This Ricci flow is an ancient
solution which is $\kappa$-non-collapsed on all scales for any
$\kappa$ at most the volume of the ball of radius one in the unit
$2$-sphere.
\end{exam}

According to a result of Hamilton which we shall prove below
(Corollary~\ref{2DGSS}):

\begin{thm}
Every orientable, $2$-dimensional $\kappa$-solution is a rescaling
of the previous example, i.e., is a family of shrinking round
$2$-spheres.
\end{thm}

\begin{exam}
Let $(S^n,g_0)$ be the standard round $n$-sphere of scalar curvature
$n/2$. Set $g(t)=(1-t)g_0$. This is a $\kappa$-solution for any
$\kappa$ which is at most the volume of the ball of radius one in
the unit $n$-sphere. If $\Gamma$ is a finite subgroup of the
isometries of $S^n$ acting freely on $S^n$, then the quotient
$S^n/\Gamma$ inherits an induced family of metrics $\bar g(t)$
satisfying the Ricci flow equation. The result is a
$\kappa$-solution for any $\kappa$ at most $1/|\Gamma|$ times the
volume of the ball of radius one in the unit sphere.
\end{exam}

\smallskip

\begin{exam}
Consider the product $S^{2}\times \mathbb{R} $, with the metric $
g(t) =(1-t)g_{0}+ds^{2}.$ This is a $\kappa$-solution for any
$\kappa$ at most the volume of a ball of radius one in the product
of the unit $2$-sphere with $\Ar$.
\end{exam}

\smallskip

\begin{exam}
The quotient $S^{2}\times \mathbb{R}/\langle\iota\rangle$, where the
involution $\iota$ is the product of the antipodal map on $S^2$ with
$s\mapsto -s$ on the $\Ar$ factor,  is an orientable
$\kappa$-solution for some $\kappa>0$.
\end{exam}

\begin{exam}
Consider the metric product $(S^{2},g_0) \times (S^{1}_{R},ds^2)$
where $(S^1_R,ds^2)$ is the circle of radius $R$. We define
$g(t)=(1-t)g_0+ds^2$. This is an ancient solution to the Ricci flow.
But it is not $\kappa$-non-collapsed for any $\kappa>0$. The reason
is that
$$|{\rm Rm}(p,t)|=\frac{1}{1-t},$$ and
$$\frac{{\rm Vol}_{g(t)}\,B(p,\sqrt{1-t})}{(1-t)^{3/2}}
\le \frac{{\rm Vol}_{g(t)}(S^{2}\times S^{1}_{R})}{(1-t)^{3/2}}
=\frac{2\pi R(1-t)4\pi}{(1-t)^{3/2}}=\frac{8\pi^2R}{\sqrt{1-t}}.$$
Thus, as $t\rightarrow-\infty$ this ratio goes to zero.
\end{exam}

\subsection{A consequence of Hamilton's Harnack inequality}

In order to prove the existence of an asymptotic gradient shrinking
soliton associated to every $\kappa$-solution, we need the following
inequality which is a consequence of Hamilton's Harnack
inequality\index{Ricci flow!Harnack inequality for} for Ricci flows
with non-negative curvature operator.

\begin{prop}\label{uniformlip}
Let $(M,g(t)),\ -\tau_0\le t\le 0$, be an $n$-dimensional Ricci flow
such that for each $t\in [-\tau_0,0]$ the Riemannian manifold
$(M,g(t))$ is complete with non-negative, bounded curvature
operator. Let $\tau=-t$. Fix a point $p\in M$ and let $x=(p,0)\in
M\times [-\tau_0,0]$. Then for any $0<c<1$ and any $\tau\le
(1-c)\tau_0$ we have
$$|\bigtriangledown l_x(q,\tau))|^{2} +R(q,\tau) \leq \frac{(1+2c^{-1})l_x(q,\tau)}{\tau},
\ \ \ {\rm and}$$
$$ R(q,\tau)-\frac{(1+c^{-1})l_x(q,\tau)}{\tau}\le
\frac{\partial l_x}{\partial \tau} $$\index{$l_x$|(} where these
inequalities hold on the open subset of full measure of $M\times
[-(1-c)\tau_0,0)$ on which $l_x$ is a smooth function.
\end{prop}

\begin{proof}
Recall that from Definition~\ref{Hdefnn} we have
$$H(X)=-\frac{\partial R}{\partial
\tau}-\frac{R}{\tau}-2\langle \nabla R,X\rangle+2{\rm Ric}(X,X).$$
Using Hamilton's Harnack's inequality (Theorem~\ref{Harnack}) with
$\chi=-X$, we have
\begin{equation*}
-\frac{\partial R}{\partial \tau}-\frac{R}{\tau_0-\tau}-2\langle
\nabla R,X\rangle+2{\rm Ric}(X,X)\ge 0.\end{equation*}
 Together these imply
 $$H(X)\ge \left(\frac{1}{\tau-\tau_0}-\frac{1}{\tau}\right)R
 =\frac{\tau_0}{\tau(\tau-\tau_0)}R.$$
Restricting to $\tau\le (1-c)\tau_0$ gives
$$H(X)\ge -\frac{c^{-1}}{\tau}R.$$\index{$H(X)$}
Take a minimal ${\mathcal L}$-geodesic from $x$ to $(q,\bar\tau)$,
we have
\begin{equation}\label{Kbartaueqn}
K^{\bar\tau}(\gamma)=\int_0^{\bar\tau}\tau^{3/2}H(X)d\tau\ge
-c^{-1}\int_0^{\bar\tau}\sqrt{\tau}Rd\tau\ge
-2c^{-1}\sqrt{\bar\tau}l_x(q,\bar\tau). \end{equation} Together with
the second equality in Theorem~\ref{Ueqn}, this gives
\begin{align*}
4\bar\tau |\bigtriangledown l_x(q,\bar\tau)|^{2} & = -4\bar\tau
R(q,\bar\tau) + 4l_x(q,\bar\tau) -\frac{4}{\sqrt{\bar\tau}}
K^{\bar\tau}(\gamma)\\ & \le -4\bar\tau R(q,\bar\tau)
+4l_x(q,\bar\tau) + 8c^{-1}l_x(q,\bar\tau)
\end{align*}

Dividing through by $4\bar\tau$, and replacing $\bar\tau$ with
$\tau$ yields the first inequality in the statement of the
proposition:
$$|\nabla l_x(q,\tau)|^2+R(q,\tau)\le \frac{(1+2c^{-1})l_x(q,\tau)}{\tau}$$
for all $0<\tau\le (1-c)\tau_0.$ This is an equation of smooth
functions on the open dense subset ${\mathcal U}(\bar \tau)$ but it
extends as an equation of $L^\infty_{\rm loc}$-functions on all of
$M$.

As to the second inequality in the statement, by the first equation
in Theorem~\ref{Ueqn} we have
$$\frac{\partial l_x(q,\tau)}{\partial \tau}
=R(q,\tau)-\frac{l_x(q,\tau)}{\tau}+\frac{1}{2\tau^{3/2}}K^{\tau}(\gamma).$$
The estimate on $K^\tau$ in Equation~(\ref{Kbartaueqn}) then gives
$$R(q,\tau)-\frac{(1+c^{-1})l_x(q,\tau)}{\tau}\le
\frac{\partial l_x(q,\tau)}{\partial \tau}.$$ This establishes the
second inequality.
\end{proof}

\begin{cor}\label{unianc}
Let $(M,g(t)),\ -\infty< t\le 0$, be a Ricci flow on a complete,
$n$-dimensional manifold with bounded, non-negative curvature
operator. Fix a point $p\in M$ and let $x=(p,0)\in M\times
(-\infty,0]$. Then for any $\tau>0$ we have
$$|\bigtriangledown l_x(q,\tau))|^{2} +R(q,\tau) \leq \frac{3l_x(q,\tau)}{\tau},$$
$$-\frac{2l_x(q,\tau)}{\tau}\le \frac{\partial l_x(q,\tau)}{\partial \tau}
\le \frac{l_x(q,\tau)}{\tau}.$$ where these inequalities are valid
in the sense of smooth functions on the open subset of full measure
of $M\times \{\tau\}$ on which $l_x$ is a smooth function, and are
valid as inequalities of $L^\infty_{\rm loc}$-functions on all of
$M\times \{\tau\}$.
\end{cor}

\begin{proof}
Fix $\tau$ and take a sequence of $\tau_0\rightarrow \infty$,
allowing us to take $c\rightarrow 1$, and apply the previous
proposition. This gives the first inequality and gives the lower
bound for $\partial l_x/\partial \tau$ in the second inequality.

 To establish the upper bound in the second inequality we
consider the path that is the concatenation of a minimal ${\mathcal
L}$-geodesic $\gamma$ from $x$ to $(q,\tau)$ followed by the path
$\mu(\tau')=(q,\tau')$ for $\tau'\ge \tau$. Then
$$l_x(\gamma*\mu|_{[\tau,\tau_1]})=\frac{1}{2\sqrt{\tau_1}}\left({\mathcal L}(\gamma)
+\int_{\tau}^{\tau_1} \sqrt{\tau'}R(q,\tau')d\tau'\right).$$
Differentiating at $\tau_1=\tau$ gives
\begin{eqnarray*}
\frac{\partial
l_x(\gamma*\mu)}{\partial\tau}\Bigl|_{\tau_1=\tau}\Bigr. & = &
-\frac{1}{4\tau^{3/2}}{\mathcal L}(\gamma)+
\frac{1}{2\sqrt{\tau}}\sqrt{\tau}R(q,\tau) \\
& = & -\frac{l_x(q,\tau)}{2\tau}+\frac{R(q,\tau)}{2}.\end{eqnarray*}
By the first inequality in this statement, we have
$$-\frac{l_x(q,\tau)}{2\tau}+\frac{R(q,\tau)}{2}\le
 \frac{l_x(q,\tau)}{\tau}.$$ Since $l_x(q,\tau')\le \tilde l(\gamma*\mu|_{[\tau,\tau']})$
 for all $\tau'\ge \tau$, this establishes the claimed upper
bound for $\partial l_x/\partial \tau$.
\end{proof}

\section{The asymptotic gradient shrinking soliton for
$\kappa$-solutions}\label{sect9.2}

We fix $\kappa>0$ and we consider an $n$-dimensional
$\kappa$-solution $(M,g(t)),\ -\infty <t\leq 0.$ Our goal in this
section is to establish the existence of an asymptotic gradient
shrinking soliton associated to this $\kappa$-solution. Fix a
reference point $p \in M$ and set $x=(p,0)\in M\times (-\infty,0]$.
By Theorem~\ref{n/2} for every $\tau>0$ there is a point $q(\tau)\in
M$ at which the function $l_x(\cdot,^\tau)$ achieves its minimum,
and furthermore, we have
$$l_x(q(\tau),\tau) \leq \frac{n}{2}.$$

 For $\bar\tau>0,$ define
$$g_{\bar\tau}(t) = \frac{1}{\bar\tau}g(\bar\tau t),\ \  -\infty<t\le 0.$$

Now we come to one of the main theorems about $\kappa$-solutions, a
result that will eventually provide a qualitative description of all
$\kappa$-solutions.

\begin{thm}\label{asympGSS} Let $(M,g(t)),\ -\infty<t\le 0$,
be a $\kappa$-solution of dimension $n$. Fix $x=(p,0)\in
M\times(-\infty,0]$. Suppose that $\{\bar\tau_{k}\}_{k=1}^\infty $
is a sequence tending to $ \infty$ as $k\rightarrow\infty$. Then,
after replacing $\{\bar\tau_k\}$ by a subsequence, the following
holds. For each $k$ denote by $M_k$ the manifold $M$, by $g_k(t)$
the family of metrics $g_{\bar\tau_k}(t)$ on $M_k$, and by $q_k\in
M_k$ the point $q(\bar\tau_k)$. The sequence of pointed flows
$(M_{k}, g_k(t), (q_{k},-1))$ defined for $t \in (-\infty,0)$
converges smoothly to a non-flat based Ricci flow $(M_{\infty},
g_{\infty}(t), (q_{\infty},-1))$ defined for $t\in (-\infty,0)$.
This limiting Ricci flow  satisfies the gradient shrinking
soliton\index{gradient shrinking soliton|(} equation in the sense
that there is a smooth function $f\colon
M_\infty\times(-\infty,0)\to \Ar$ such that for every
$t\in(-\infty,0)$ we have
\begin{equation}\label{GSSequation}
 {\rm Ric}_{g_{\infty}(t)} + {\rm Hess}^{g_{\infty}(t)}(f(t)) +
\frac{1}{2t}g_{\infty}(t) = 0.
\end{equation}
Furthermore, $(M_\infty,g_\infty(t))$ has non-negative curvature operator, is
$\kappa$-non-collapsed, and satisfies $\partial R_{g_\infty}(x,t)/\partial t\ge
0$ for all $x\in M_\infty$ and all $t<0$.
\end{thm}

\begin{figure}[ht]
  \centerline{\epsfbox{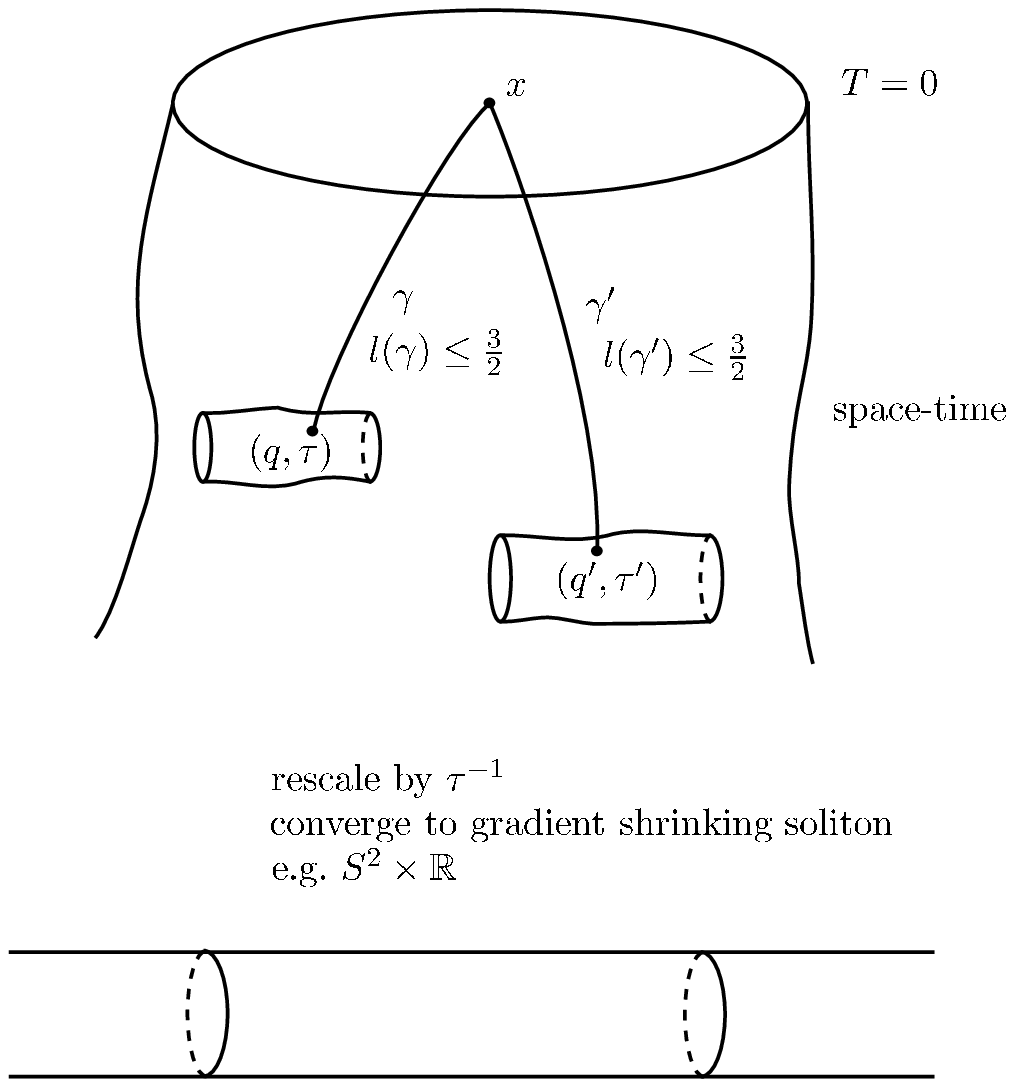}}
  \caption{Gradient shrinking soliton.}\label{fig:GSS}
\end{figure}

See {\sc Fig.}~\ref{fig:GSS}.

\begin{rem}
We are not claiming that the gradient shrinking soliton is a
$\kappa$-solution (or more precisely an extension forward in time of
a time-shifted version of a $\kappa$-solution) because we are not
claiming that the time-slices have bounded curvature operator.
Indeed, we do not know if this is true in general. We shall
establish below (see Corollary~\ref{2DGSS} and
Corollary~\ref{3DGSSkappa}) that in the case $n=2,3$, the gradient
shrinking soliton does indeed have time-slices of bounded curvature,
and hence is an extension of a $\kappa$-solution. We are also not
claiming at this point that the limiting flow is a gradient
shrinking soliton in the sense that there is a one-parameter family
of diffeomorphisms $\varphi_t\colon M_\infty\to M_\infty,\ t<0$,
with the property that $|t|\varphi_t^*g_\infty(-1)=g_\infty(t)$ and
with the property that the $\varphi_t$ are generated by the gradient
vector field of a function. We shall also
 establish this result in dimensions $2$ and $3$ later in this chapter.
\end{rem}

We will divide the proof of Theorem~\ref{asympGSS} into steps.
First, we will show that the reduced length and norm of the
curvature $|{\rm Rm}|$ are bounded throughout the sequence in some
way. Then using the $\kappa$-non-collapsed assumption, by the
compactness theorem (Theorem~\ref{flowlimit}), we conclude that a
subsequence of the sequence of flows converges geometrically to a
limiting flow. Then, using the fact that the limit of the reduced
volumes\index{reduced volume}, denoted $\widetilde
V_{\infty}(M_\infty\times\{t\})$\index{$\widetilde V_x$}, is
constant we show that the limit flow is a gradient shrinking
soliton\index{gradient shrinking soliton|)}. Finally we argue that
the limit is non-flat. The proof occupies the rest of
Section~\ref{sect9.2}.

\subsection{Bounding the reduced length and the curvature}
Now let's carry this procedure out in detail. The first remark is
that since rescaling does not affect the $\kappa$-non-collapsed
hypothesis, all the Ricci flows $(M_k,g_k(t))$ are
$\kappa$-non-collapsed on all scales. Next, we have the effect on
reduced volume.

\begin{claim}\label{scaleinv2} For each $k\ge 1$ denote by $x_k\in M_k$ the point $(p,0)\in
M_k$. Let $\widetilde V_{x_k}(\tau)=\widetilde
V_{x_k}(M_k\times\{\tau\})$ denote the reduced volume function for
the Ricci flow $(M_k,g_k(t))$ from the point $x_k$, and let
$\widetilde V_x(\tau)$ denote the reduced volume of
$M\times\{\tau\}$ for the Ricci flow $(M,g(t))$ from the point $x$.
Then
$$\widetilde V_{x_k}(\tau) = \widetilde V_x (\bar \tau_{k}\tau).$$
\end{claim}

\begin{proof}
This is a special case of the reparameterization equation for
reduced volume (Lemma~\ref{Qformula}).
\end{proof}

By Theorem~\ref{4pi} the reduced volume function $\widetilde
V_x(\tau)$ is a non-increasing function of $\tau$ with ${\rm
lim}_{\tau\rightarrow 0}\widetilde V_x(\tau)=(4\pi)^{\frac{n}{2}}$.
Since the integrand for $\widetilde V_x(\tau)$ is everywhere
positive, it is clear that $\widetilde V_x(\tau)>0$ for all $\tau$.
Hence, ${\rm lim}_{\tau\rightarrow\infty}\widetilde V_x(\tau)$
exists. By Corollary~\ref{4picase} either this limit as $\tau$ goes
to infinity is less than $(4\pi)^{n/2}$ or the flow is the constant
flow on flat Euclidean space. The latter is ruled out by our
assumption that the manifolds are non-flat. It follows immediately
from this and Claim~\ref{scaleinv2} that:

\begin{cor}\label{redvolconst}
There is a non-negative constant $V_\infty<(4\pi)^{n/2}$ such that
for all $\tau\in(0,\infty)$, we have
\begin{equation}\label{limVeqn}
{\rm lim}_{k\rightarrow\infty}\widetilde V_{x_k}(\tau)=V_\infty.
\end{equation}
\end{cor}

Now let us turn to the length functions $l_{x_k}$.

\begin{claim}\label{tausqtau}
For any $\tau>0$ we have
$$l_{x_k}(q_k,\tau)\le \frac{n}{2\tau^2}+\frac{n\tau}{2}.$$
\end{claim}

\begin{proof}
 By the choice of $q_k$ we have $l_{x_k}(q_k,\tau_k)\le \frac{n}{2}$.
 By the scale invariance of $l$ (Corollary~\ref{COR}) we have
 $l_{x_k}(q_k,-1)\le n/2$ for all $k$. Fix $0<\tau<1$. Integrating the inequality
$$\frac{-2l_x(q_k,\tau)}{\tau}\le \frac{\partial l_{x_k}(q_k,\tau)}{\partial \tau}$$
 from $\tau$ to $1$  yields
$$l_{x_k}(q_k,\tau)\le \frac{n}{2\tau^2}.$$
If $\tau>1$, then
 integrating the second inequality in the second displayed line of Corollary~\ref{unianc} gives
 $l_{x_k}(q_k,\tau)\le \frac{n\tau}{2}$.
\end{proof}

\begin{cor}\label{lxkest}
There is a  positive continuous function $C_1(\tau)$ defined for
$\tau>0$ such that for any $q\in M_k$ we have:
$$l_{x_k}(q,\tau)\le
\left(\sqrt{\frac{3}{\tau}}d_{g_k(-\tau)}(q_k,q)+C_1(\tau)\right)^2,$$
$$|\nabla l_{x_k}(q,\tau)|\le
\frac{3}{\tau}d_{g_k(-\tau)}(q_k,q)+\sqrt{\frac{3}{\tau}}C_1(\tau).$$
\end{cor}

\begin{proof}
 By Corollary~\ref{unianc}, for any $q\in M_k$ we have $|\nabla
l_{x_k}(q,\tau)|^2\le 3l_{x_k}(q,\tau)/\tau$. Since
$l_{x_k}(q_k,\tau)\le \frac {n}{2\tau_0^2}+\frac{n\tau}{2}$,
integrating yields
$$l_{x_k}(q,\tau)\le
\left(\sqrt{\frac{3}{\tau}}d_{g_k(-\tau)}(q_k,q)+C_1(\tau)\right)^2,$$
with $C_1(\tau)$ being $\sqrt{(n/2\tau^2)+(n\tau/2)}$. The second
statement follows from this and Proposition~\ref{uniformlip}.
\end{proof}

It follows immediately from Corollary~\ref{lxkest} that for each
$A<\infty$ and $\tau_0>0$, the functions $l_{x_k}$ are uniformly
bounded (by a bound that is independent of $k$ but depends on
$\tau_0$ and $A$) on the balls $B(q_k,-\tau_0,A)$. Once we know that
the $l_{x_k}$ are uniformly bounded on $B(q_k,-\tau_0,A)$, it
follows from Corollary~\ref{unianc} that $R_{g_k}$ are also
uniformly bounded on the $B(q_k,-\tau_0,A)$. Invoking
Corollary~\ref{posderiv},
 we see
that for any $A<\infty$ the scalar curvatures of the metrics $g_k$
are uniformly bounded on $B_{g_k}(q_k,-\tau_0,A)\times (-\infty,
-\tau_0]$. Since the metrics have non-negative curvature operator,
this implies that the eigenvalues of this operator are  uniformly
bounded on these regions.
 Since we are assuming
that the original Ricci flows are $\kappa$-non-collapsed on all
scales, it follows from Theorem~\ref{flowlimit} that after passing
to a subsequence there is a geometric limit
$(M_\infty,g_\infty(t),(q_\infty,-1)),\ -\infty<t\le -\tau_0$, which
is a Ricci flow which is $\kappa$-non-collapsed on all scales.

Since this is true for every $\tau_0>0$, by a standard
diagonalization argument passing to a further subsequence we get a
geometric limit flow $(M_\infty,g_\infty(t),(q_\infty,-1)),\
-\infty<t<0$.

Let us summarize our progress to this point.

\begin{cor}
After passing to a subsequence of the $\tau_k$ there is a smooth
limiting flow of the $(M_k,g_k(t),(q_k,-1)),-\infty<t\le 0$,
$$(M_\infty,g_\infty(t),(q_\infty,-1)),$$ defined for  $-\infty<t<0$. For every $t<0$
the Riemannian manifold $(M_\infty,g_\infty(t))$ is complete of
non-negative curvature. The flow is $\kappa$-non-collapsed on all
scales and satisfies $\partial R/\partial t\ge 0$.
\end{cor}

\begin{proof}
 Since the  flows in the sequence are all
$\kappa$-non-collapsed on all scales and have non-negative curvature
operator, the limiting flow is $\kappa$-non-collapsed on all scales
and has non-negative curvature operator. By the consequence of
Hamilton's Harnack inequality (Corollary~\ref{posderiv}), we have
$\partial R/\partial t\ge 0$ for the original $\kappa$-solution.
This condition also passes to the limit.
\end{proof}

\subsection{The limit function}

The next step in the proof is to construct the limiting function
$l_\infty$ of the $l_{x_k}$ and show that it satisfies the gradient
shrinking soliton equation.

By definition of the geometric limit, for any compact connected set
$K\subset M_\infty$ containing $q_\infty$ and any compact
subinterval $J$ of $(-\infty,0)$ containing $-1$, for all $k$
sufficiently large we have smooth embeddings $\psi_k\colon K\to M_k$
sending $q_\infty$ to $q_k$ so that the pullbacks of the
restrictions of the family of metrics $g_k(t)$ for $t\in J$ to $K$
converge uniformly in the $C^\infty$-topology to the restriction of
$g_\infty(t)$ on $K\times J$. Take an exhausting sequence $K_k\times
J_k$ of such products of compact sets with closed intervals, and
pass to a subsequence so that for all $k$ the diffeomorphism
$\psi_k$ is defined on $K_k\times J_k$. We denote by $l_k$ the
pullback of $l_{x_k}$ under these embeddings and by $h_k(t)$ the
pullback of the family of metrics $g_k(t)$. We denote by $
\nabla^{h_k}$ the gradient with respect to $h_k(t)$, and similarly
$\triangle^{h_k}$ denotes the Laplacian for the metric $h_k(t)$. By
construction, for any compact subset of $M_\infty\times (-\infty,0)$
for all $k$ sufficiently large the function $ l_k$ is defined on the
compact set. We use $\nabla$ and $\triangle$ to refer to the
covariant derivative and the Laplacian in the limiting metric
$g_\infty$.

Now let us consider the functions $l_{x_k}$. According to
Corollary~\ref{lxkest}, for any $A<\infty$ and any $0<\tau_0<T$,
both $l_{x_k}$ and $|\nabla l_{x_k}|$ are uniformly bounded on
$B(q_k,-1,A)\times [-T,-\tau_0]$ independent of $k$. Hence, the
$l_{x_k}$ are uniformly Lipschitz on these subspaces. Doing this for
each $A$, $\tau_0$, and $T$ and using a standard diagonalization
argument then shows that, after transferring to the limit, the
functions $l_k$ are uniformly locally bounded and uniformly locally
Lipschitz on $M_\infty\times (-\infty,0)$ with respect to the
limiting metric $g_\infty$.

Fix $0<\alpha<1$. Passing to a further subsequence if necessary, we
can arrange that the $l_k$ converge strongly 
in $C^{0,\alpha}_{\rm loc}$ to a function $l_\infty$ defined on
$M_\infty\times(-\infty,0)$. Furthermore, it follows that the
restriction of $l_\infty$ is locally Lipschitz, and hence the
function $l_\infty$ is an element of $W^{1,2}_{\rm
loc}\left(M_\infty\times(-\infty,0)\right)$. Also, by passing to a
further subsequence if necessary, we can assume that the $l_k$
converge weakly in $W^{1,2}_{\rm loc}$ to $l_\infty$.

\begin{cor}\label{linfest}
For any $\tau>0$ and any $q$ we have
$$|\nabla l_\infty(q,\tau)|\le
\frac{3}{\tau}d_{g_\infty}(-\tau)(q_\infty,q)+\sqrt{\frac{3}{\tau}}C_1(\tau),$$
where $C_1(\tau)$ is the continuous function from
Corollary~\ref{lxkest}.
\end{cor}

\begin{proof}
This is immediate from Corollary~\ref{lxkest} and Fatou's lemma.
\end{proof}

\begin{rem}
{\bf N.B.} We are not claiming that $l_\infty$ is the reduced length
function from a point of $M_\infty\times (-\infty,0)$.
\end{rem}

\subsection{Differential inequalities for $l_\infty$}

The next step is to establish differential equalities for $l_\infty$
related to, but stronger than, those that we established in
Chapter~\ref{newcomp2} for $l_x$. Here is a crucial result.

\begin{prop}\label{weakeqn} The function $l_\infty$ is a smooth function
on $M\times (-\infty,0)$ and satisfies the following two
differential equalities: \begin{equation}\label{weak1}
\frac{\partial l_\infty}{\partial \tau}+|\nabla
l_\infty|^2-R+\frac{n}{2\tau}-\triangle l_\infty=0\end{equation}
 and
\begin{equation}\label{weak2}2\triangle l_\infty-|\nabla
l_\infty|^2+R+\frac{l_\infty-n}{\tau}=0.\end{equation}
\end{prop}

The proof of this result is contained in Sections~\ref{sect9.2.3}
through~\ref{sect9.2.5}

\subsection{Preliminary results toward the proof of
Proposition~\protect{\ref{weakeqn}}}\label{sect9.2.3}

In this subsection we shall prove  that the left-hand side of
Equation~(\ref{weak1}) is a distribution and is $\ge 0$ in the
distributional sense. We shall also show that this distribution
extends to a continuous linear functional on compactly supported
functions in $W^{1,2}$.

 The first step in the proof of
this result is the following, somewhat delicate lemma.

\begin{lem}\label{nablaconv}
For any $t\in (-\infty,0)$ we have
$${\rm lim}_{k\rightarrow\infty}|\nabla^{h_k} l_k|_{h_k}^2d{\rm vol}(h_k)
=|\nabla l_\infty|_{g_\infty}^2 d{\rm vol}(g_\infty)$$ in the sense
of distributions on $M_\infty\times \{t\}$.
\end{lem}

\begin{proof}
It suffices to fix $0<\tau_0<|t|$.
 The inequality in one direction ($\ge$) is a general
result. Here is the argument.  Since the $|\nabla^{g_k} l_{x_k}|_{g_k}$ are
uniformly essentially bounded on every $B(x_k,-\tau_0,A)\times [-T,-\tau_0]$,
the $|\nabla^{h_k}l_k|_{h_k}$ are uniformly essentially bounded on
$B(x_\infty,-\tau_0,A)\times [-T,-\tau_0]$. (Of course,
$\nabla^{h_k}l_k=dl_k=\nabla l_k$.) Since the $h_k$ converge uniformly on
compact sets to $g_\infty$, it is clear that
\begin{equation}\label{firstdist} {\rm
lim}_{k\rightarrow\infty}\left(|\nabla^{h_k}l_k|^2_{h_k}d{\rm
vol}(h_k)-|\nabla l_k|_{g_\infty}^2d{\rm
vol}(g_\infty)\right)=0\end{equation} in the sense of distributions
on $M\times \{t\}$.
 Since the $
l_k$ converge uniformly on compact subsets to $l_\infty$, it follows
immediately from Fatou's lemma that
$${\rm lim}_{k\rightarrow \infty}|\nabla
l_k|^2_{g_\infty}d{\rm vol}(g_\infty)\ge|\nabla
l_\infty|^2_{g_\infty}d{\rm vol} (g_\infty)$$ in the sense of
distributions on $M_\infty\times \{t\}$. Thus, we have the following
inequality of distributions:
$${\rm
lim}_{k\rightarrow\infty}|\nabla^{h_k}l_k|^2_{h_k}d{\rm vol}(h_k)\ge
|\nabla l_\infty|^2_{g_\infty}d{\rm vol} (g_\infty).$$

 We need to establish the opposite
inequality which is not a general result, but rather relies on  the
bounds on $\triangle^{g_k}l_{x_k}$ (or equivalently on
$\triangle^{h_k}l_k$) given in second inequality in
Theorem~\ref{weak}. We must show that for each $t\le -\tau_0$ and
for any $\varphi$, a non-negative, smooth function with compact
support in $M_\infty\times \{t\}$, we have
$${\rm lim}_{k\rightarrow\infty}\int_{M\times\{t\}}\varphi
\left(|\nabla^{h_k} l_k|_{h_k}^2d{\rm vol}(h_k)-|\nabla
l_\infty|_{g_\infty}^2d{\rm vol}(g_\infty)\right)\le 0.$$ First,
notice that since, on the support of $\varphi$, the metrics $h_k$
converge uniformly in the $C^\infty$-topology to $g_\infty$ and
since $|\nabla^{h_k}l_k|_{h_k}^2$ and $|\nabla
l_\infty|_{g_\infty}^2$ are essentially bounded on the support of
$\varphi$, we have
\begin{eqnarray}\nonumber
\lefteqn{{\rm lim}_{k\rightarrow\infty}\int_{M\times\{t\}}\varphi
\left(|\nabla^{h_k}l_k|_{h_k}^2d{\rm vol}(h_k)-|\nabla
l_\infty|_{g_\infty}^2d{\rm vol}(g_\infty) \right)} & & \\
&  = & {\rm lim}_{k\rightarrow\infty}\int_{M\times\{t\}}\varphi (|
\nabla^{h_k} l_k|_{h_k}^2-|\nabla
l_\infty|_{h_k}^2)d{\rm vol}(h_k) \nonumber \\
\nonumber& = & {\rm
lim}_{k\rightarrow\infty}\int_{M\times\{t\}}\langle\nabla^{h_k}
l_k-\nabla l_\infty),\varphi\nabla^{h_k} l_k\rangle_{h_k}d{\rm
vol}(h_k)
\\ & & \hskip.55in+\int_{M\times\{t\}}\langle\nabla^{h_k}l_k-\nabla l_\infty),
\varphi\nabla l_\infty\rangle_{h_k}d{\rm vol}(h_k).
\label{11}\end{eqnarray} We claim that, in the limit, the last term
in this expression vanishes. Using the fact that the $h_k$ converge
uniformly in the $C^\infty$-topology to $g_\infty$ on the support of
$\varphi$, and $|\nabla l_\infty|$ is bounded on this support we can
rewrite the last term as
\begin{equation}\label{12}
{\rm lim}_{k\rightarrow\infty}\int_{M\times\{t\}}\langle\nabla(
l_k-l_\infty),\varphi\nabla l_\infty\rangle_{g_\infty}d{\rm vol}(
g_\infty).\end{equation} Since $ l_k-l_\infty$ goes to zero weakly
in $W^{1,2}$ on the support of $\varphi$  whereas $ l_\infty$ is an
element of $W^{1,2}$ of this compact set, we see that the expression
given in~(\ref{12}) vanishes and hence that
$${\rm
lim}_{k\rightarrow\infty}\int_{M\times\{t\}}\langle\nabla^{h_k}(
l_k-l_\infty),\varphi\nabla l_\infty\rangle_{h_k}d{\rm
vol}(h_k)=0.$$

It remains to consider the first term in the last expression in
Equation~(\ref{11}). (This is where we shall need the differential
inequality for the $\triangle^{g_k}l_{x_k}$.) Since the $l_k$
converge uniformly to $l_\infty$ on the support of $\varphi$, we can
choose positive constants $\epsilon_k$ tending to  $0$ as $k$ tends
to $\infty$ so that $l_\infty- l_k+\epsilon_k>0$ on the support of
$\varphi$. We can rewrite
\begin{eqnarray*}\lefteqn{{\rm
lim}_{k\rightarrow\infty}\int_{M\times\{t\}}\langle\left(\nabla^{h_k}
l_k-\nabla l_\infty\right),\varphi\nabla^{h_k} l_k\rangle_{h_k}
d{\rm vol}(h_k) = } \\
& & {\rm
lim}_{k\rightarrow\infty}\int_{M\times\{t\}}\langle\nabla^{h_k}(
l_k-l_\infty-\epsilon_k),\varphi\nabla^{h_k}l_k\rangle_{h_k} d{\rm
vol}(h_k).\end{eqnarray*}
\begin{claim}
$${\rm lim}_{k\rightarrow\infty}\int_{M\times\{t\}}\langle\nabla^{h_k}(l_k-l_\infty-\epsilon_k),
\varphi\nabla^{h_k}  l_k\rangle_{h_k} d{\rm vol}(h_k)\le 0.$$
\end{claim}

\begin{proof}
Since $\varphi$ is a compactly supported, non-negative smooth
function, it follows from Theorem~\ref{weak} that we have the
following inequality of distributions:
$$\varphi\triangle^{h_k} l_k\le
\frac{\varphi}{2}\left(|\nabla^{h_k} l_k|_{h_k}^2-R_{h_k}-\frac{
l_k-n}{\tau}\right).$$ (Here $R_{h_k}$ is the scalar curvature of
$h_k$.)
 That is to say, for any non-negative
$C^\infty$-function $f$ we have
\begin{eqnarray*}
\lefteqn{\int_{M\times\{t\}}-\langle\nabla^{h_k}
l_k,\nabla^{h_k}(\varphi \cdot f)\rangle_{h_k} d{\rm vol}(h_k) \le}
\\
& & \int_{M\times\{t\}}\frac{\varphi f}{2}\left(|\nabla^{h_k}
l_k|_{h_k}^2-R_{h_k}-\frac{l_k-n}{\tau}\right)d{\rm vol}(h_k).
\end{eqnarray*}

We claim that the same inequality holds as long as $f$ is a
non-negative, locally  Lipschitz function. The point is that given
such a function $f$, we can find a sequence of non-negative
$C^\infty$-functions $f_k$ on the support of $\varphi$ (by say
mollifying $f$)  that converge to $f$ strongly in the $W^{1,2}$-norm
on the support of $\varphi$. The sought-after inequality holds for
every $f_k$. Since both sides of the inequality are continuous in
the $W^{1,2}$-norm of the function, the result holds for the limit
function $f$ as well.

Now we apply this with $f$ being the non-negative locally Lipschitz
function $l_\infty-l_k+\epsilon_k$. We conclude that
\begin{eqnarray*}
\lefteqn{\int_{M\times\{t\}}\langle \nabla^{h_k}(\varphi(
l_k-l_\infty-\epsilon_k)),\nabla^{h_k} l_k\rangle_{h_k}d{\rm vol}(h_k) \le}  \\
& & \int_{M\times\{t\}}\frac{\varphi(l_\infty-
l_k+\epsilon_k)}{2}\left(|\nabla^{h_k}
l_k|_{h_k}^2-R_{h_k}-\frac{l_k-n}{\tau}\right)d{\rm vol}(h_k).
\end{eqnarray*}
Now taking the limit as $k\rightarrow\infty$, we see that the
right-hand side of this inequality tends to zero since $(l_\infty-
l_k+\epsilon_k)$ tends uniformly to zero on the support of $\varphi$
and $|\nabla^{h_k}l_k|_{h_k}^2$, $R_k$ and $l_k$ are all uniformly
essentially bounded on the support of $\varphi$. Thus, the term
$$\int_{M\times\{t\}}\langle\nabla^{h_k}(\varphi(
l_k-l_\infty-\epsilon_k)),\nabla^{h_k} l_k\rangle_{h_k}d{\rm
vol}(h_k)$$ has a limsup $\le 0$ as $k$ tends to $\infty$. Now we
expand
$$\nabla^{h_k}(\varphi(l_k-l_\infty-\epsilon_k))=\nabla^{h_k}(\varphi)(l_k-l_\infty-\epsilon_k)+
\varphi\nabla^{h_k}(l_k-l_\infty-\epsilon_k).$$ The first term on
the right-hand side converges to zero as $k\rightarrow\infty$ since
$l_k-l_\infty-\epsilon_k$ tends uniformly to zero on the support of
$\varphi$. This completes the proof of the claim.
\end{proof}
We have now established the  inequalities in both directions and
hence completed the proof of Lemma~\ref{nablaconv}.
\end{proof}

\begin{lem}
Consider the distribution
$${\mathcal D}=\frac{\partial
l_\infty}{\partial \tau}+|\nabla
l_\infty|^2-R+\frac{n}{2\tau}-\triangle l_\infty$$ on
$M_\infty\times (-\infty,0)$. Then ${\mathcal D}$ extends to a
continuous linear functional on the space of compactly supported
$W^{1,2}$-functions on $M_\infty\times (-\infty,0)$. Furthermore, if
$\psi$ is a non-negative Lipschitz function on $M_\infty\times
(-\infty,0)$ with compact support, then ${\mathcal D}(\psi)\le 0$.
\end{lem}

\begin{proof}
Clearly, since the $l_k$ converge uniformly on compact subsets of
$M_{\infty}\times (-\infty,0)$ to $l_\infty$ and the metrics $h_k$
converge smoothly to $g_\infty$, uniformly on compact sets, it
follows that the $\triangle ^{h_k}l_k $ converge in the weak sense
to $\triangle l_\infty$ and similarly, the $\partial l_k/\partial
\tau$ converge in the weak sense to $\partial l_\infty/\partial
\tau$. Hence, by taking limits from Theorem~\ref{weak}, using
Lemma~\ref{nablaconv}, we see that
\begin{equation}\label{1stlinftyeqn}
{\mathcal D}=\frac{\partial l_\infty}{\partial \tau}+|\nabla
l_\infty|^2-R+\frac{n}{2\tau}-\triangle l_\infty\ge 0\end{equation}
 in the weak sense on $M\times (-\infty,0)$.

Since $R$ and $\frac{n}{2\tau}$ are $C^\infty$-functions, it is
clear that the distributions given by these terms extend to
continuous linear functionals on the space of compactly supported
$W^{1,2}$-functions. Similarly, since $|\nabla l_\infty|^2$ is an
element of $L^\infty_{\rm loc}$, it also extends to a continuous
linear functional on compactly supported $W^{1,2}$-functions. Since
$|\partial l_\infty/\partial \tau|$ is an locally essentially
bounded function, $\partial l_\infty/\partial \tau$ extends to a
continuous functional on the space of compactly supported $W^{1,2}$
functions. Lastly, we consider $\triangle l_\infty$. As we have
seen, the value of the associated distribution on $\varphi$ is given
by
$$\int_{M\times (-\infty,,0)}-\langle \nabla \varphi,\nabla
l_\infty\rangle_{g_\infty} d{\rm vol}(g_\infty)d\tau.$$ Since
$|\nabla l_\infty|$ is a locally essentially bounded function, this
expression also extends to a continuous linear functional on
compactly supported $W^{1,2}$-functions.

Lastly, if $\psi$ is an element of $W^{1,2}$ with compact support
and hence can be approximated in the $W^{1,2}$-norm by non-negative
smooth functions. The last statement  is now immediate from
Equation~(\ref{1stlinftyeqn}).
\end{proof}

This leads immediately to:

\begin{cor}\label{etimesdis}
The functional
$$\varphi\mapsto{\mathcal D}(e^{-l_\infty}\varphi)$$ is a distribution and
its value on any non-negative, compactly supported $C^\infty$-function
$\varphi$ is $\ge 0$.
\end{cor}

\begin{proof}
If $\varphi$ is a compactly supported non-negative
$C^\infty$-function, then $e^{-l_\infty}\varphi$ is a compactly
supported non-negative Lipschitz function. Hence, this result is an
immediate consequence of the previous corollary.
\end{proof}

\subsection{Extension to non-compactly supported functions}

The next step in this proof is to estimate the $l_{x_k}$ uniformly from below
in order to show that  the integrals involved in the distributions in
Proposition~\ref{weakeqn} are absolutely convergent so that they extend to
continuous functionals on a certain  space of functions that includes
non-compactly supported functions.

\begin{lem}\label{lowerbdd}
There is a constant $c_1>0$ depending only on the dimension $n$ such
that for any $p,q\in M_k$ we have
$$l_{x_k}(p,\tau)\ge -l_{x_k}(q,\tau)-1+c_1\frac{d^2_{g(-\tau)}(p,q)}{\tau}.$$
\end{lem}

\begin{proof}
Since both sides of this inequality and also Ricci flow are
invariant if the metric and time are simultaneously rescaled, it
suffices to consider the case when $\tau=1$. Also, since ${\mathcal
U}_x(1)$ is a dense subset, it suffices to assume that $p,q\in
{\mathcal U}_x(1)$. Also, by symmetry, we can suppose that
$l_{x_k}(q,1)\le l_{x_k}(p,1)$.

 Let
$\gamma_1$ and $\gamma_2$ be the minimizing ${\mathcal L}$-geodesics
from $x$ to $(p,1)$ and $(q,1)$ respectively. We define a function
$f\colon M_k\times M_k\times [0,\infty)\to \Ar$ by
$$f(a,b,\tau)=d_{g_k(-\tau)}(a,b).$$
Since $\gamma_1(0)=\gamma_2(0)$ we have
\begin{eqnarray}\nonumber
d_{g_k(-1)}(p,q) & = & f(p,q,1) \\ \nonumber & =&
\int_0^1\frac{d}{d\tau}f(\gamma_1(\tau),\gamma_2(\tau),\tau)d\tau \\
\nonumber & = & \int_0^1\bigl(\frac{\partial f}{\partial
\tau}(\gamma_1(\tau),\gamma_2(\tau),\tau)+\langle \nabla
f_a,\gamma_1'(\tau)\rangle\bigr. \\
\label{fform} & & \hskip.2in+\bigl.\langle \nabla
f_b,\gamma_2'(\tau)\rangle\bigr)d\tau,
\end{eqnarray}
where $\nabla_af$ and $\nabla_bf$ refer respectively to  the
gradient of $f$ with respect to the first copy of $M_k$ in the
domain and the second copy of $M_k$ in the domain. Of course,
$|\nabla f_a|=1$ and $|\nabla f_b|=1$.

By Corollary~\ref{DL}, we have $\gamma_1'(\tau)=\nabla
l_{x_k}(\gamma_1(\tau),\tau)$ and $\gamma_2'(\tau)=\nabla
l_{x_k}(\gamma_2(\tau),\tau)$. Since $R\ge 0$ we have
$$l_{x_k}(\gamma_1(\tau),\tau)=\frac{1}{2\sqrt{\tau}}{\mathcal
L}_{x_k}(\gamma_1|_{[0,\tau]})\le \frac{1}{2\sqrt{\tau}}{\mathcal
L}_{x_k}(\gamma_1)=\frac{1}{\sqrt{\tau}}l_{x_k}(p,1).$$
Symmetrically, we have
$$l_{x_k}(\gamma_2(\tau),\tau)\le \frac{1}{\sqrt{\tau}}l_{x_k}(q,1).$$
From this inequality, Corollary~\ref{unianc}, and the fact that
$R\ge 0$, we have
\begin{eqnarray}\nonumber \bigl|\langle \nabla
f_a(\gamma_1(\tau),\gamma_2(\tau),\tau),\gamma_1'(\tau)\rangle
\bigr|& \le &
|\gamma_1'(\tau)| =|\nabla l_{x_k}(\gamma_1(\tau),\tau) | \nonumber \\
& \le &
\frac{\sqrt{3}}{\tau^{3/4}}\sqrt{l_{x_k}(p,1)}  \nonumber \\
& \le &
\frac{\sqrt{3}}{\tau^{3/4}}\sqrt{l_{x_k}(p,1)+1}.\label{1eqn}\end{eqnarray}
Symmetrically, we have \begin{equation}\label{2eqn} \left|\langle
\nabla
f_b(\gamma_1(\tau),\gamma_2(\tau),\tau),\gamma'_2(\tau)\rangle\right|\le
\frac{\sqrt{3}}{\tau^{3/4}}\sqrt{l_{x_k}(q,1)}\le
\frac{\sqrt{3}}{\tau^{3/4}}\sqrt{l_{x_k}(q,1)+1}.\end{equation}

 It follows from
Corollary~\ref{unianc} that for any $p$
$$|\nabla (\sqrt{l_{x_k}(p,\tau)})|\le
\frac{\sqrt{3}}{2\sqrt{\tau}}.$$ Set
$r_0(\tau)=\tau^{3/4}(l_{x_k}(q,1)+1)^{-1/2}$. For any $p'\in
B_{g_k}(\gamma_1(\tau),\tau,r_0(\tau))$ integrating gives
$$l_{x_k}^{1/2}(p',\tau)\le
l_{x_k}^{1/2}(\gamma_1(\tau),\tau)+\frac{\sqrt{3}}{2\sqrt{\tau}}r_0(\tau)\le
\left(\tau^{-1/4}+\frac{\sqrt{3}}{2}\tau^{1/4}\right)\sqrt{l_{x_k}(p,1)+1},$$
where in the last inequality we have used the fact that $1\le
l_{x_k}(q,1)+1\le l_{x_k}(p,1)+1$. Again using
Corollary~\ref{unianc} we have
$$R(p',\tau)\le
\frac{3}{\tau}\left(\tau^{-1/4}+\frac{\sqrt{3}}{2}\tau^{1/4}\right)^2(l_{x_k}(p,1)+1).$$
Now consider $q'\in B_{g_k(\tau)}(\gamma_2(\tau),\tau,r_0(\tau))$.
Similarly to the above computations, we have
$$l^{1/2}_{x_k}(q',\tau)\le
l_{x_k}^{1/2}(q,1)+\frac{\sqrt{3}}{2\sqrt{\tau}}r_0(\tau),$$ so that
$$l_{x_k}^{1/2}(q',\tau)\le \left(\tau^{-1/4}
+\frac{\sqrt{3}}{2}\tau^{1/4}\right)\sqrt{l_{x_k}(q,1)+1},$$ and
$$|{\rm Ric}(q',\tau)\le R(q',\tau)\le
\frac{3}{\tau}\left(\tau^{-1/4}+\frac{\sqrt{3}}{2}\tau^{1/4}\right)^2(l_{x_k}(q,1)+1).$$
We set
$$K=\frac{3}{\tau}\left(\tau^{-1/4}+\frac{\sqrt{3}}{2}\tau^{1/4}\right)^2(l_{x_k}(q,1)+1).$$

Now, noting that $\partial/\partial \tau$ here is
$-\partial/\partial t$ of Proposition~\ref{I.8.3}, we apply
Proposition~\ref{I.8.3} to see that
\begin{eqnarray*}
\left|\frac{\partial}{\partial
\tau}f(\gamma_1(\tau),\gamma_2(\tau),\tau)\right|
& \le & 2(n-1)\left(\frac{2}{3(n-1)}Kr_0(\tau)+r_0(\tau)^{-1}\right) \\
& \le &
\left(C_1\tau^{-3/4}+C_2\tau^{-1/4}+C_3\tau^{1/4}\right)\sqrt{l_{x_k}(q,1)+1},
\end{eqnarray*}
where $C_1,C_2,C_3$ are constants depending only on the dimension
$n$.

Now plugging Equation~(\ref{1eqn}) and~(\ref{2eqn}) and  the above
inequality into Equation~(\ref{fform}) we see that
\begin{eqnarray*}
d_{g(-1)}(p,q) & \le &
\int_0^1\left(\left(C_1\tau^{-3/4}+C_2\tau^{-1/4}+C_3\tau^{1/4}\right)\sqrt{l_{x_k}(q,1)+1}\right.
\\ & & +
\left.\sqrt{3}\tau^{-3/4}\sqrt{l_{x_k}(q,1)+1}+\sqrt{3}\tau^{-3/4}\sqrt{l_{x_k}(p,1)+1}\right)d\tau.
\end{eqnarray*}
This implies that
$$d_{g(-1)}(p,q)\le
C\left(\sqrt{l_{x_k}(q,1)+1}+\sqrt{l_{x_k}(p,1)+1}\right),$$ for
some constant depending only on the dimension. Thus, since we are
assuming that $l_{x_k}(p,1)\ge l_{x_k}(q,1)$ we have
$$d^2_{g(-1)}(p,q)\le
C^2\left(3(l_{x_k}(p,1)+1)+(l_{x_k}(q,1)+1)\right)\le
4C^2(l_{x_k}(p,1)+1+l_{x_k}(q,1)),$$ for some constant $C<\infty$
depending only on the dimension. The result now follows immediately.
\end{proof}

\begin{cor}
For  any $q'\in M$ and any $0<\tau_0\le\tau'$ we have
$$l_{x_k}(q',\tau') \ge -\frac{n}{2(\tau')^2}-\frac{\tau'}{2}-1+c_1\frac{d_{g^2_k(-\tau_0)}(q_k,q')}{\tau'},$$
where  $c_1$ is the constant from Lemma~\ref{lowerbdd}.
\end{cor}

\begin{proof}
By Claim~\ref{tausqtau}
$$l_{x_k}(q_k,\tau')\le \frac{n}{2(\tau')^2}+\frac{n\tau'}{2}.$$
Now applying Lemma~\ref{lowerbdd} we see that for any $0<\tau'$ and any $q'\in
M_k$ we have \begin{eqnarray*} l_{x_k}(q',\tau') & \ge &
-\frac{n}{2(\tau')^2}-\frac{n\tau'}{2}-1+c_1\frac{d^2_{g_k(-\tau')}(q_k,q')}{\tau'}
\\
& \ge &
-\frac{n}{2(\tau')^2}-\frac{n\tau'}{2}-1+c_1\frac{d^2_{g_k(-\tau_0)}(q_k,q')}{\tau'}.
\end{eqnarray*}\index{$l_x$|)} In the last inequality, we use the fact that the
Ricci curvature is positive so that the metric is decreasing under
the Ricci flow.
\end{proof}

Since the time slices of all the flows in question have non-negative
curvature, by Theorem~\ref{BishopGromov} the volume of the ball of
radius $s$ is at most $\omega s^{n}$ where $\omega$ is the volume of
the ball of radius one in $\Ar^n$. Since the $l_k$ converge
uniformly to $l_\infty$ on compact sets and since the metrics $h_k$
converge uniformly in the $C^\infty$-topology on compact sets to
$g_\infty$, it follows that for any $\epsilon>0$, for any
$0<\tau_0\le \tau'<\infty$ there is a radius $r$ such that for every
$k$ and any $\tau\in [\tau_0,\tau']$ the integral
$$\int_{M_\infty\setminus B_{h_k(-\tau_0)}(q_k,r)}e^{-l_k(q,\tau)}dq<\epsilon.$$
It follows by Lebesgue dominated convergence that
$$\int_{M_\infty\setminus
B_{g_\infty(-\tau_0)}(q_\infty,r)}e^{-l_\infty(q,\tau)}dq\le
\epsilon.$$

\begin{claim}\label{absconv1}
 Fix a compact interval $[-\tau,-\tau_0]\subset
(-\infty,0)$. Let $f$ be a locally Lipschitz function that is
defined on $M_\infty\times [-\tau,-\tau_0]$ and such that there is a
constant $C$ with the property that $f(q,\tau')$ by $C$ times ${\rm
max}(l_\infty(q,\tau'),1)$ . Then the distribution ${\mathcal
D}_1=fe^{-l_\infty}$ is absolutely convergent in the following
sense. For any bounded smooth function $\varphi$ defined on all of
$M_\infty\times [-\tau,-\tau_0]$ and any sequence of compactly
supported, non-negative smooth functions $\psi_k$, bounded above by
$1$ everywhere that are eventually $1$ on every compact subset, the
following limit exists and is finite:
$${\rm lim}_{k\rightarrow\infty}{\mathcal D}_1(\varphi\psi_k).$$
 Furthermore, the limit is independent of the choice of the
$\psi_k$ with the given properties.
\end{claim}

\begin{proof}
It follows from the above discussion that there are constants $c>0$
 and a ball $B\subset M_\infty$ centered at $q_\infty$ such that on
$M_\infty\times [-\tau,-\tau_0]\setminus B\times [\tau,-\tau_0]$ the
function $l_\infty$ is greater than
$cd^2_{g_\infty(-\tau_0)}(q_\infty,\cdot)-C'$. Thus,
$fe^{-l_\infty}$ has fixed exponential decay at infinity. Since the
Riemann curvature of $M_\infty\times\{\tau'\}$ is non-negative for
every $\tau'$, the flow is distance decreasing, and there is a fixed
polynomial upper bound to the growth rate of volume at infinity.
This leads to the claimed convergence property.
\end{proof}

\begin{cor}\label{absconv}
The distributions $|\nabla l_\infty|^2e^{-l_\infty}$,
$Re^{-l_\infty}$, $|(\partial l_\infty/\partial \tau)|e^{-l_\infty}$
are absolutely convergent in the sense of the above claim. \end{cor}

\begin{proof}
 By
Corollary~\ref{unianc}, each of the Lipschitz functions $|\nabla
l_\infty|^2$, $|\partial l_\infty/\partial \tau|$ and $R$
  is at most a constant multiple of $l_\infty$.
Hence, the corollary follows from the previous claim.
\end{proof}

There is a slightly weaker statement that is true for $\triangle
e^{-l_\infty}$.

\begin{claim}\label{absconv2}
Suppose that $\varphi$ and $\psi_k$ are as in Claim~\ref{absconv1},
but in addition $\varphi$ and all the $\psi_k$ are uniformly
Lipschitz. Then
$${\rm lim}_{k\rightarrow\infty}\int_{M_\infty}\varphi\psi_k\triangle
e^{-l_\infty}d{\rm vol}_{g_\infty}$$ converges absolutely.
\end{claim}

\begin{proof}
This time the value of the distribution on a compactly supported
smooth function $\rho$ is given by the integral of $$-\langle
\nabla\rho,\nabla e^{-l_\infty}\rangle=\langle \nabla \rho,\nabla
l_\infty\rangle e^{-l_\infty}.$$ Since $|\nabla l_\infty|$ is less
than or equal to the maximum of $1$ and $|\nabla l_\infty|^2$, it
follows immediately, that if $|\nabla \rho|$ is bounded, then the
integral is absolutely convergent. From this the claim follows
easily.
\end{proof}

\begin{cor}\label{weakneg}
Fix $0<\tau_0<\tau_1<\infty$. Let $f$ be a non-negative, smooth
bounded function on $M_\infty\times [\tau_0,\tau_1]$ with (spatial)
gradient of bounded norm. Then
$$\int_{\tau_0}^{\tau_1}\int_{M_\infty\times\{-\tau\}}\left(\frac{\partial l_\infty}{\partial \tau}+|\nabla
l_\infty|^2-R+\frac{n}{2\tau}-\triangle
l_\infty\right)f\tau^{-n/2}e^{-l_\infty}d{\rm
vol}_{g_\infty}d\tau\ge 0.$$
\end{cor}

\begin{proof}
For the interval $[\tau_0,\tau']$  we construct a sequence of
uniformly Lipschitz functions $\psi_k$ on $M_\infty\times
[\tau_0,\tau']$ that are non-negative, bounded above by one and
eventually one on every compact set. Let $\rho(x)$ be a smooth bump
function which is one for $x$ less than $1/4$ and is zero from $x\ge
3/4$ and is everywhere between $0$ and $1$. For any $k$ sufficiently
large let $\psi_k$ be the composition of
$\rho(d_{g_\infty(-\tau_0)}(q_\infty,\cdot)-k)$. Being compositions
of $\rho$ with Lipschitz functions with Lipschitz constant $1$, the
$\psi_k$ are a uniformly Lipschitz family of functions on
$M_\infty\times\{-\tau_0\}$. Clearly then they form a uniformly
Lipschitz family on $M_\infty\times [\tau_0,\tau']$ as required.
This allows us to define any of the above distributions on Lipschitz
functions on $M_\infty\times [\tau_0,\tau']$.

Take a family $\psi_k$ of uniformly Lipschitz functions, each
bounded between $0$ and $1$ and eventually one of every compact
subset of $M_\infty\times [\tau_0,\tau_1]$. Then the family
$f\psi_k$ is a uniformly Lipschitz family of compactly supported
functions. Hence, we can apply Claims~\ref{absconv1}
and~\ref{absconv2} to establish that the integral in question is the
limit of an absolutely convergence sequence. By
Corollary~\ref{etimesdis} each term in the sequence is non-positive.
\end{proof}

 \subsection{Completion of the proof of
Proposition~\protect{\ref{weakeqn}}}\label{sect9.2.5}

Lebesgue dominated convergence implies that the following limit
exists
$${\rm lim}_{k\rightarrow\infty}\widetilde V_k(\tau)\equiv\widetilde V_\infty(\tau)=\int_{M_\infty\times\{-\tau\}}\tau^{-n/2}e^{-
l_\infty(q,\tau)}d{\rm vol}_{g_\infty(\tau)}.$$\index{$\widetilde
V_x$} By Corollary~\ref{redvolconst}, the function $\tau\rightarrow
\widetilde V_\infty(\tau)$ is constant. On the other hand, note that
for any $0< \tau_0<\tau_1<\infty$, we have
\begin{eqnarray*}
\widetilde V_\infty(\tau_1)-\widetilde V_\infty(\tau_0) & = &
\int_{\tau_0}^{\tau_1}\frac{d\widetilde
V_\infty}{d\tau}d\tau \\
& = & \int_{\tau_0}^{\tau_1}\int_{M_\infty}\left(\frac{\partial
l_\infty}{\partial \tau}-R+\frac{n}{2\tau}\right)
\left(\tau^{-n/2}e^{-l_\infty(q,\tau)}d{\rm
vol}_{g_\infty(\tau)}\right).\end{eqnarray*} According to
Corollary~\ref{absconv} this is an absolutely convergent integral,
and so  this integral is zero.

\begin{claim}
\begin{eqnarray*}
\int_{\tau_0}^{\tau_1}\int_{M_\infty\times\{-\tau\}}\triangle
e^{-l_\infty}d{\rm vol}_{g_\infty}d\tau & =
&\int_{\tau_0}^{\tau_1}\int_{M_\infty\times\{-\tau\}}\left(|\nabla
l_\infty|^2-\triangle l_\infty\right)e^{-l_\infty}d{\rm vol}_{g_\infty}d\tau \\
& = & 0.
\end{eqnarray*}
\end{claim}

\begin{proof}
Since we are integrating against the constant function $1$, this
result is clear, given the convergence result,
Corollary~\ref{absconv}, necessary to show that this integral is
well defined.
\end{proof}

Adding these two results together gives us the following
\begin{equation}\label{tildeD}
\int_{\tau_0}^{\tau_1}\int_{M_\infty\times\{-\tau\}}\left(\frac{\partial
l_\infty}{\partial \tau}+|\nabla l_\infty|^2
-R+\frac{n}{2\tau}-\triangle
l_\infty\right)\tau^{-n/2}e^{-l_\infty}d{\rm
vol}_{g_\infty}=0.\end{equation}

Now let $\varphi$ be any compactly supported, non-negative smooth
function. By scaling by a positive constant, we can assume that
$\varphi\le 1$ everywhere. Let $\widetilde{\mathcal D}$ denote the
distribution given by
$$\widetilde{\mathcal D}(\varphi)=\int_{\tau_0}^{\tau_1}\int_{M_\infty\times\{-\tau\}}
\varphi\left(\frac{\partial l_\infty}{\partial \tau}+|\nabla
l_\infty|^2 -R+\frac{n}{2\tau}-\triangle
l_\infty\right)\tau^{-n/2}e^{-l_\infty}d{\rm vol}_{g_\infty}.$$ Then
we have seen that $\widetilde{\mathcal D}$ extends to a functional
on bounded smooth functions of bounded gradient. Furthermore,
according to Equation~(\ref{tildeD}), we have $\widetilde {\mathcal
D}(1)=0$. Thus,
$$0={\mathcal D}(1)={\mathcal D}(\varphi)+{\mathcal D}(1-\varphi).$$
Since both $\varphi$ and $1-\varphi$ are non-negative, it follows from
Corollary~\ref{weakneg}, that ${\mathcal D}(\varphi)$ and ${\mathcal
D}(1-\varphi)$ are each $\ge 0$. Since their sum is zero, it must be the case
that each is individually zero.

This proves that the Inequality~(\ref{1stlinftyeqn}) is actually an
equality in the weak sense, i.e., an equality of distributions on
$M_\infty\times [\tau_0,\tau')$. Taking limits we see:
\begin{equation}\label{constcon}
\widetilde{\mathcal D}=\left(\frac{\partial l_\infty}{\partial
\tau}+|\nabla l_\infty|^2-R+\frac{n}{2\tau}-\triangle
l_\infty\right)\tau^{-n/2}e^{-l_\infty}=0,\end{equation} in the weak
sense on all of $M\times(-\infty,0)$. Of course, this implies that
$$\frac{\partial l_\infty}{\partial \tau}+|\nabla
l_\infty|^2-R+\frac{n}{2\tau}-\triangle l_\infty=0$$ in the weak
sense.

It now follows by parabolic regularity that $l_\infty$ is a smooth
function on $M_\infty\times(-\infty,0)$ and that
Equation~(\ref{constcon}) holds in the usual sense.

Now from the last two equations in Corollary~\ref{lformula} and the
convergence of the $l_{x_k}$ to $l_\infty$, we conclude that the
following equation also holds:
\begin{equation}\label{wk=}2\triangle l_\infty-|\nabla
l_\infty|^2+R+\frac{l_\infty-n}{\tau}=0.\end{equation}

This completes the proof of Proposition~\ref{weakeqn}.

\subsection{The gradient shrinking soliton equation}

Now we return to the proof of Theorem~\ref{asympGSS}. We have shown that the
limiting Ricci flow   referred to in that result exists, and we have
established that the limit $l_\infty$ of the length functions $l_{x_k}$ is a
smooth function and satisfies the differential equalities given in
Proposition~\ref{weakeqn}. We shall use these to establish the gradient
shrinking soliton\index{gradient shrinking soliton} equation,
Equation~(\ref{GSSequation}), for the limit for $f=l_\infty$.

\begin{prop}\label{gradshrink}
The following equation holds on $M_\infty\times (-\infty,0)$:
$$ {\rm Ric}_{g_{\infty}(t)} + {\rm Hess}^{g_\infty(t)}(l_\infty(\cdot,\tau))
-\frac{1}{2\tau}g_{\infty}(t) = 0,$$ where $\tau=-t$,
\end{prop}

\begin{proof}  This result will follow immediately
from:

\begin{lem}\label{9.1}
Let $(M,g(t)),\ 0\le t\le T$, be an $n$-dimensional Ricci flow, and
let $f\colon M\times [0,T]\to \Ar$ be a smooth function. As usual
set $\tau=T-t$. Then the function
$$u=(4\pi\tau)^{-\frac{n}{2}}e^{-f}$$ satisfies the conjugate heat
equation
$$-\frac{\partial u}{\partial t}-\triangle u+Ru=0,$$ if and only if we have
$$\frac{\partial f}{\partial t}+\triangle f-|\nabla
f|^2+R-\frac{n}{2\tau}=0.$$ Assuming that $u$ satisfies the
conjugate heat equation, then setting
$$v=\left[\tau\left(2\triangle f-|\nabla
f|^2+R\right)+f-n\right]u,$$ we have
$$-\frac{\partial v}{\partial t}-\triangle v+Rv=-2\tau
\bigl|{\rm Ric}_g+{\rm Hess}^g(f)-\frac{1}{2\tau}g\bigr|^2u.$$
\end{lem}

Let us assume the lemma for a moment and use it to complete the
proof of the proposition.

We apply the lemma to the limiting Ricci flow
$(M_\infty,g_\infty(t))$ with the function $f=l_\infty$. According
to Proposition~\ref{weakeqn} and the first statement in
Lemma~\ref{9.1}, the function $u$ satisfies the conjugate heat
equation. Thus, according to the second statement in
Lemma~\ref{9.1}, setting
$$v=\left[\tau\left(2\triangle f-|\nabla
f|^2+R\right)+f-n\right]u,$$ we have
$$\frac{\partial v}{\partial\tau}-\triangle v+Rv=-2\tau
\bigl|{\rm Ric}_g+{\rm Hess}(f)-\frac{1}{2\tau}g\bigr|^2u.$$ On the
other hand, the second equality in Proposition~\ref{weakeqn} shows
that $v=0$. Since $u$ is nowhere zero, this implies that
$${\rm Ric}_{g_\infty}+{\rm Hess}^{g_\infty}(f)-\frac{1}{2\tau}g_\infty=0.$$
This completes the proof of the proposition assuming the lemma.
\end{proof}

Now we turn to the proof of the lemma.

\begin{proof}{\em (of Lemma~\ref{9.1})}
Direct computation shows that $$-\frac{\partial u}{\partial
t}-\triangle u+Ru= \left(-\frac{n}{2\tau}+\frac{\partial f}{\partial
t}+\triangle f-|\nabla f|^2+R\right)u.$$ From this, the first
statement of the lemma is clear.
 Let
$$H=\bigl[\tau(2\triangle f-|\nabla f|^2+R)+f-n\bigr]$$ so that $v=Hu$.
Then, of course,
$$\frac{\partial v}{\partial t}=\frac{\partial H}{\partial t}u+H\frac{\partial u}{\partial t}\ \ \ {\rm and}$$
$$\triangle v=\triangle H\cdot u+2\langle\nabla H,\nabla u\rangle+H\triangle u.$$
Since $u$ satisfies the conjugate heat equation, we have
$$-\frac{\partial v}{\partial t}-\triangle v+Rv=\left(-\frac{\partial H}{\partial t}-\triangle
H\right)u-2\langle\nabla H,\nabla u\rangle.$$ Differentiating the
definition of $H$ yields \begin{equation}\label{dHeqn}
\frac{\partial H}{\partial t} = -(2\triangle f-|\nabla
f|^2+R)+\frac{\partial f}{\partial t}+ \tau \left(2\frac{\partial
}{\partial t}\triangle f -\frac{\partial }{\partial t}(|\nabla
f|^2)+\frac{\partial R}{\partial t}\right)
\end{equation}

\begin{claim}\label{11.20}
$$\frac{\partial }{\partial t}\triangle
f=\triangle(\frac{\partial f}{\partial t})+2\langle {\rm Ric},{\rm
Hess}(f)\rangle.$$
\end{claim}

\begin{proof}
We work in local coordinates. We have $$\triangle
f=g^{ij}\nabla_i\nabla_j
f=g^{ij}(\partial_i\partial_jf-\Gamma_{ij}^k\partial_kf),$$ so that
from the Ricci flow equation we have
\begin{eqnarray*}
\frac{\partial }{\partial t}\triangle f & = & 2{\rm Ric}^{ij}{\rm
Hess}(f)_{ij}+g^{ij}\frac{\partial}{\partial t}\left({\rm
Hess}(f)_{ij}\right)
\\
& = & 2{\rm Ric}^{ij}{\rm Hess}(f)_{ij}+g^{ij}{\rm
Hess}\left({\frac{\partial f}{\partial
t}}\right)_{ij}-g^{ij}\frac{\partial \Gamma_{ij}^k}{\partial
t}\partial_kf.
\end{eqnarray*}
Since the first term is $2\langle {\rm Ric},{\rm Hess}(f)\rangle$
and the second is $\triangle(\frac{\partial f}{\partial t})$, to
complete the proof of the claim, we must show that the last term of
this equation vanishes. In order to simplify the computations, we
assume that the metric is standard to second order at the point and
time under consideration. Then, using the Ricci flow equation, the
definition of the Christoffel symbols in terms of the metric, and
the fact that $g_{ij}$ is the identity matrix at the given point and
time and that its covariant derivatives in all spatial directions
vanish at this point and time, we get
$$g^{ij}\frac{\partial \Gamma_{ij}^k}{\partial t}=g^{kl}g^{ij}\left(-(\nabla_j {\rm Ric})_{li}-(\nabla_i
 {\rm Ric})_{lj}+(\nabla_l{\rm Ric})_{ij}\right).$$
 This expression vanishes by the second Bianchi identity (Claim~\ref{Bianchi}).
 This completes the proof of the claim.
\end{proof}

We also have
$$\frac{\partial }{\partial t}(|\nabla f|^2)=2{\rm Ric}(\nabla f,\nabla f)+2\langle
\nabla\frac{\partial f}{\partial t},\nabla f\rangle.$$ (Here $\nabla
f$ is a one-form, which explains the positive sign in the Ricci
term.)

Plugging this and Claim~\ref{11.20} into Equation~(\ref{dHeqn})
yields
\begin{eqnarray*}
\frac{\partial H}{\partial t} & = & -2\triangle f+|\nabla f|^2-R+\frac{\partial f}{\partial t} \\
& & + \tau \left(4\langle {\rm Ric},{\rm Hess}(f)\rangle +2\triangle
\frac{\partial f}{\partial t} -2{\rm Ric}(\nabla f,\nabla f)
-2\langle \nabla\frac{\partial f}{\partial t},\nabla
f\rangle+\frac{\partial R}{\partial t}\right).
\end{eqnarray*}
Also,
\begin{equation*} \triangle H  =  \triangle
f+\tau\left(2\triangle^2f-\triangle( |\nabla f|^2)+\triangle
R\right).
\end{equation*}

Since $u$ satisfies the conjugate heat equation, from the first part
of the lemma we have
\begin{equation}\label{feqn}\frac{\partial f}{\partial t}=-\triangle f+|\nabla
f|^2-R+\frac{n}{2\tau}.\end{equation}
 Putting all this together and
using the Equation~(\ref{Revol}) for $\partial R/\partial t$ yields
\begin{eqnarray*}
\frac{\partial H}{\partial t}+\triangle H & = & -\triangle f+|\nabla
f|^2+\frac{\partial f}{\partial t}-R \\ & & +\tau\bigl(4\langle {\rm
Ric},{\rm Hess}(f)\rangle+2\triangle \frac{\partial f}{\partial t}
+2\triangle^2f-2{\rm Ric}(\nabla f,\nabla f)\bigr. \\
& & \bigl.-\triangle(|\nabla f|^2) -2\langle \nabla\frac{\partial
f}{\partial t},\nabla f\rangle+2\triangle
R +2|{\rm Ric}|^2\bigr) \\
& = & -\triangle f+|\nabla f|^2+\frac{\partial f}{\partial t}-R
+\tau\bigl(4\langle {\rm Ric},{\rm
Hess}(f)\rangle+2\triangle(|\nabla f|^2-R)\bigr.
\\
& & \bigl.-2{\rm Ric}(\nabla f,\nabla f)-\triangle(|\nabla
f|^2)-2\langle\nabla\frac{\partial f}{\partial t},\nabla
f\rangle+2\triangle R+2|{\rm Ric}|^2\bigr)
\\
& = & -\triangle f+|\nabla f|^2+\frac{\partial f}{\partial
t}-R+\tau\Bigl[4\langle {\rm Ric},{\rm
Hess}(f)\rangle+\triangle(|\nabla f|^2)\Bigr.
\\
& & -2{\rm Ric}(\nabla f,\nabla f)  +2\langle\nabla(\triangle
f),\nabla f\rangle -2\langle \nabla(|\nabla f|^2),\nabla
f\rangle \\
& & \Bigl. +2\langle\nabla R,\nabla f\rangle+2|{\rm Ric}|^2\Bigr].
\end{eqnarray*}
Similarly, we have
\begin{eqnarray*}
\frac{2\langle\nabla u,\nabla H\rangle}{u} & = & -2\langle\nabla f,\nabla H\rangle \\
& = & -2|\nabla f|^2-2\tau\langle\nabla f,\left(\nabla(2\triangle
f)-|\nabla f|^2+R\right)\rangle \\
& = & -2|\nabla f|^2-\tau\left(4\langle\nabla f,\nabla(\triangle
f)\rangle-2\langle\nabla f,\nabla(|\nabla f|^2)\rangle+\langle\nabla
f,\nabla R\rangle\right).
\end{eqnarray*}

Thus,
\begin{eqnarray*}
\frac{\partial H}{\partial t} +\triangle H+\frac{2\langle\nabla
u,\nabla H\rangle}{u} & = & -\triangle f-|\nabla f|^2+\frac{\partial
f}{\partial t}-R + \tau\Bigl[4\langle {\rm Ric},{\rm Hess}(f)\rangle
\Bigr.
\\
& & +\triangle(|\nabla f|^2)-2{\rm Ric}(\nabla f,\nabla f)+2|{\rm
Ric}|^2
\\
& & \Bigl.-2\langle\nabla f,\nabla(\triangle f)\rangle\Bigr].
\end{eqnarray*}

\begin{claim} The following equality holds:
\begin{eqnarray*} \triangle(|\nabla f|^2) & = &
2\langle\nabla(\triangle f),\nabla f\rangle+2{\rm Ric}(\nabla
f,\nabla f)+2|{\rm Hess}(f)|^2,
\end{eqnarray*}
\end{claim}

\begin{proof} We have
$$\triangle(|\nabla f|^2)=\triangle \langle \nabla f,\nabla f\rangle =\triangle \langle df,df\rangle=2\langle \triangle
df,df\rangle+2\langle \nabla df,\nabla df\rangle.$$ The last term is
$|{\rm Hess}(f)|^2$. According to Lemma~\ref{lapformula} we have
$\triangle df=d(\triangle f)+{\rm Ric}(\nabla f,\cdot)$. Plugging
this in gives
$$\triangle(|\nabla f|^2)=2\langle d(\triangle f),df\rangle+2\langle
{\rm Ric}(\nabla f,\cdot),df\rangle+2|{\rm Hess}(f)|^2,$$ which is
clearly another way of writing the claimed result.
\end{proof}

Using this we can simplify the above to
\begin{eqnarray*}
\frac{\partial H}{\partial t} +\triangle H+\frac{2\langle\nabla
u,\nabla H\rangle}{u} & = &-\triangle f-|\nabla f|^2+\frac{\partial
f}{\partial t}-R
\\
& & +\tau\left(4\langle {\rm Ric},{\rm Hess}(f)\rangle +2|{\rm
Hess}(f)|^2+2|{\rm Ric}|^2\right).
\end{eqnarray*}
Now using Equation~(\ref{feqn}) we have
\begin{eqnarray*}
\frac{\partial H}{\partial t} +\triangle H+\frac{2\langle\nabla
u,\nabla H\rangle}{u} & = &
-2\triangle f-2R+\frac{n}{2\tau} \\
& & +\tau\left(4\langle {\rm Ric},{\rm Hess}(f)\rangle+2|{\rm Ric}|^2+2|{\rm Hess}(f)|^2\right) \\
& = & 2\tau\bigl(2\langle {\rm Ric},{\rm Hess}(f)\rangle+|{\rm Ric}|^2+|{\rm Hess}(f)|^2\bigr. \\
& & -\frac{\triangle f}{\tau}-\frac{R}{\tau}+\frac{n}{4\tau^2}\bigr) \\
& = & 2\tau\bigl|{\rm Ric}+{\rm
Hess}(f)-\frac{1}{2\tau}g_\infty\bigr|^2
\end{eqnarray*}
Since
$$-\frac{\partial v}{\partial t}-\triangle v+Rv=-u\left(\frac{\partial H}{\partial t}
 +\triangle H+\frac{2\langle\nabla u,\nabla H\rangle}{u}\right),$$
 this proves the lemma.
\end{proof}

At this point, setting $f=l_\infty$, we have established all the
results claimed in Theorem~\ref{asympGSS} except for the fact that
the limit is not flat. This we establish in the next chapter.

\subsection{Completion of the proof of
Theorem~\protect{\ref{asympGSS}}}

 To complete the
proof of Theorem~\ref{asympGSS} it remains to show that for no $t\in
(-\infty,0)$ is $(M_\infty,g_\infty(t))$ flat.

\begin{claim}
If, for some $t\in (-\infty,0)$, the Riemannian manifold
$(M_\infty,g_\infty(t))$ is flat, then there is an isometry from
$\Ar^n$ to $(M_\infty,g_\infty(t))$  and the pullback under this
isometry of the function $l_\infty(x,\tau)$ is the function
$|x|^2/4\tau+\langle x,a\rangle+b\cdot \tau$ for some $a\in \Ar^n$
and $b\in \Ar$.
\end{claim}

\begin{proof}
We know that $f=l_\infty(\cdot,\tau)$ solves the equation given in
Lemma~\ref{9.1} and hence by the above argument, $f$ also satisfies
the equation given in  Proposition~\ref{gradshrink}. If the limit is
flat, then the equation becomes
$${\rm Hess}(f)=\frac{1}{2\tau}g.$$ The universal covering of $(M_\infty,g_\infty(t))$
is isometric to $\Ar^n$. Choose an identification with $\Ar^n$, and lift $f$ to
the universal cover. Call the result $\tilde f$. Then $\tilde f$ satisfies
${\rm Hess}(\tilde f)=\frac{1}{2\tau}\tilde g$, where $\tilde g$ is the usual
Euclidean metric on $\Ar^n$. This means that $\tilde f-|x|^2/4\tau$ is an
affine linear function. Clearly, then $\tilde f$ is not invariant under any
free action of a non-trivial group, so that the universal covering in question
is the trivial cover. This completes the proof of the claim.
\end{proof}

If $(M_\infty,g_\infty(t))$ is flat for some $t\in (-\infty,0)$, then by the
above $(M_\infty,g_\infty(t))$ is isometric to $\Ar^n$. According to
Theorem~\ref{4picase} this implies that $\widetilde
V_\infty(\tau)=(4\pi)^{n/2}$. This contradicts Corollary~\ref{redvolconst}, and
the contradiction establishes that $(M_\infty,g_\infty(t))$ is not flat for any
$t<0$.  Together with Proposition~\ref{gradshrink}, this completes the proof of
Theorem~\ref{asympGSS}, namely of the fact that $(M_\infty,g_\infty(t)),\
-\infty<t<0$, is a non-flat, $\kappa$-non-collapsed Ricci flow with
non-negative curvature operator that satisfies the gradient shrinking soliton
equation, Equation~(\ref{GSSequation}).

To emphasize once again, we do not claim that $(M_\infty,g_\infty(t))$ is a
$\kappa$-solution, since we do not claim that each time-slice has bounded
curvature operator.

\section{Splitting results at infinity}

\subsection{Point-picking}\index{point-picking}

There is a very simple, general result about Riemannian manifolds
that we shall use in various contexts to prove that certain types of
Ricci flows split at infinity as a product with $\Ar$.

\begin{lem}\label{4times}
Let $(M,g)$ be a Riemannian manifold  and let $p\in M$ and $r>0$ be
given. Suppose that $B(p,2r)$ has compact closure in $M$ and suppose
that $f\colon B(p,2r)\times (-2r,0]\to \Ar$ is a continuous, bounded
function with $f(p,0)>0$. Then there is a point $(q,t)\in
B(p,2r)\times (-2r,0]$ with the following properties:
\begin{enumerate}
\item[(1)] $f(q,t)\ge f(p,0)$.
\item[(2)] Setting $\alpha=f(p,0)/f(q,t)$ we have $d(p,q)\le 2r(1-\alpha)$
and $t\ge -2r(1-\alpha)$.
\item[(3)] $f(q',t')<2f(q,t)$ for all $(q',t')\in B(q,\alpha r)\times (t-\alpha r,t]$.
\end{enumerate}
\end{lem}

\begin{proof}
Consider sequences of points
$x_0=(p,0),x_1=(p_1,t_1),\ldots,x_j=(p_j,t_j)$ in $B(p,2r)\times
(-2r,0]$ with the following properties:
\begin{enumerate}
\item[(1)] $f(x_i)\ge 2f(x_{i-1})$;
\item[(2)] Setting $r_i=rf(x_0)/f(x_{i-1})$, then $r_i\le 2^{i-1}r$, and
 we have that $$x_i\in
B(p_{i-1},r_i)\times (t_{i-1}-r_i,t_{i-1}].$$
\end{enumerate}

Of course, there is exactly one such sequence with  $j=0$: it has
$x_0=(p,0)$. Suppose we have such a sequence defined for some $j\ge
0$. If follows immediately from the properties of the sequence that
$f(p_j,t_j)\ge 2^jf(p,0)$, that
$$t_j\ge -r(1+2^{-1}+\cdots+2^{1-j}),$$
and that $r_{j+1}\le 2^{-j}r$. It also follows immediately from the
triangle inequality that $d(p,p_j)\le r(1+2^{-1}+\cdots +2^{1-j})$.
 This means that
 $$B(p_j,r_{j+1})\times(t_j-r_{j+1},t_j]\subset B(p,2r)\times (-2r,0].$$
  Either the point
$x_j$ satisfies the conclusion of the lemma, or we can find
$x_{j+1}\in B(p_j,r_{j+1})\times (t_j-r_{j+1},t_j]$ with
$f(x_{j+1})\ge 2f(x_j)$. In the latter case we extend our sequence
by one term. This shows that either the process terminates at some
$j$, in which case $x_j$ satisfies the conclusion of the lemma, or
it continues indefinitely. But it cannot continue indefinitely since
$f$ is bounded on $B(p,2r)\times (-2r,0]$.
\end{proof}

One special case worth stating separately is when $f$ is independent
of $t$.

\begin{cor}\label{pp}
Let $(M,g)$ be a Riemannian manifold  and let $p\in M$ and $r>0$ be given.
Suppose that $B(p,2r)$ has compact closure in $M$ and suppose that $f\colon
B(p,2r)\to \Ar$ is a continuous, bounded function with $f(p)>0$. Then there is
a point $q\in B(x,2r)$ with the following properties:
\begin{enumerate}
\item[(1)] $f(q)\ge f(p)$.
\item[(2)] Setting $\alpha=f(p)/f(q)$ we have $d(p,q)\le 2r(1-\alpha)$
and $f(q')<2f(q)$ for all $q'\in B(q,\alpha r)$.
\end{enumerate}
\end{cor}

\begin{proof}
Apply the previous lemma to $\widehat f\colon B(p,2r)\times
(-2r,0]\to \Ar$ defined by $\widehat f(p,t)=f(p)$.
\end{proof}

\subsection{Splitting results}

Here we prove a splitting result\index{splitting result} for ancient solutions
of non-negative curvature. They are both based on Theorem~\ref{topsplit}.

\begin{prop}\label{prodatinf}
Suppose that $(M,g(t)),\ -\infty <t< 0$, is a $\kappa$-non-collapsed
 Ricci flow of dimension\footnote{This result in fact
 holds in all dimensions.} $n\le 3$. Suppose that $(M,g(t))$ is a
 complete, non-compact, non-flat Riemannian manifold
with non-negative curvature operator for each $t$. Suppose  that
$\partial R(q,t)/\partial t\ge 0$ for all $q\in M$ and all $t<0$.
Fix $p\in M$. Suppose that there is a sequence of points $p_i\in M$
going to infinity with the property that
$${\rm lim}_{i\rightarrow \infty}R(p_i,-1)d^2_{g(-1)}(p,p_i)=\infty.$$
  Then
there is a sequence of points $q_i\in M$ tending to infinity such
that, setting $Q_i=R(q_i,-1)$, we have ${\rm
lim}_{i\rightarrow\infty}d^2(p,q_i)Q_i=\infty$. Furthermore, setting
$g_i(t)=Q_ig(Q_i^{-1}(t+1)-1)$,
 the
sequence of  based flows $(M,g_i(t),(q_i,-1)),\ -\infty<t\le -1$,
converges smoothly to $(N^{n-1},h(t))\times (\Ar,ds^2)$, a product
Ricci flow defined for $-\infty<t\le -1$ with $(N^{n-1},h(-1))$
being non-flat and of bounded, non-negative curvature.
\end{prop}

\begin{cor}\label{compact2D}
There is no two-dimensional flow satisfying the hypotheses of
Proposition~\ref{prodatinf}.
\end{cor}

\begin{proof} (of Proposition~\ref{prodatinf})
Take a sequence $p_i\in M$ such that
$$d_{g(-1)}^2(p,p_i)R(p_i,-1)\rightarrow\infty$$ as $i\rightarrow \infty$.
  We set $d_i=d_{g(-1)}(p,p_i)$ and we set
$B_i=B(p_i,-1,d_i/2)$, and we let $f\colon B_i\to \Ar$ be the square
root of the scalar curvature. Since $(M,g(-1))$ is complete, $B_i$
has compact closure in $M$, and consequently $f$ is a bounded
continuous function on $B_i$. Applying Corollary~\ref{pp} to
$(B_i,g(-1))$ and $f$, we conclude that there is a point $q_i\in
B_i$ with the following properties:
\begin{enumerate}
\item[(1)] $R(q_i,-1)\ge R(p_i,-1)$
\item[(2)] $B'_i=B(q_i,-1,(d_iR(p_i,-1)^{1/2})/(4R(q_i,t_i)^{1/2})\subset B(p_i,-1,d_i/2)$.
\item[(3)] $R(q',-1)\le 4R(q_i,-1)$ for
all $(q',-1)\in B'_i$.
\end{enumerate}
Since $d_{g(-1)}(p,q_i)\ge d_i/2$, it is also the case that
$d^2_{g(-1)}(p,q_i)R(q_i,-1)$ tends to infinity as $i$ tends to infinity.
Because of our assumption on the time derivative of $R$, it follows that
$R(q',t)\le 4R(q_i,-1)$ for all $q'\in B'_i$ and for all $t\le -1$.

 Set
$Q_i=R(q_i,-1)$. Let $M_i=M$, and set $x_i=(q_i,-1)$. Lastly, set
$g_i(t)=Q_ig(Q_i^{-1}(t+1)-1)$. We consider the based Ricci flows
$(M_i,g_i(t),x_i),\ -\infty<t\le -1$. We see that $R_{g_i}(q',t)\le
4$ for all $(q',t)\in B_{g_i}(q_i,-1,d_iR(p_i,-1)^{1/2}/4)\times
(-\infty,-1]$. Since the original Ricci flows are
$\kappa$-non-collapsed, the same is true for the rescaled flows.
Since $d_iR(p_i,-1)^{1/2}/4\rightarrow \infty$, by
Theorem~\ref{flowlimit} there is a geometric limit flow
$(M_\infty,g_\infty(t),(q_\infty,-1))$ defined for $t\in
(-\infty,-1]$. Of course, by taking limits we see that
$(M_\infty,g_\infty(t))$ is $\kappa$-non-collapsed, its scalar
curvature is bounded above by $4$, and its curvature operator is
non-negative. It follows that $(M_\infty,g_\infty(t))$ has bounded
curvature.

To complete the proof we show that the Ricci flow $(M_\infty,g_\infty(t))$
splits as a product of a line with a Ricci flow of one lower dimension. By
construction $(M_\infty,g_\infty(-1))$ is the geometric limit constructed from
$(M,g(-1))$ in the following manner. We have a sequence of points $q_i$ tending
to infinity in $M$ and constants $\lambda_i=R(q_i,-1)$ with the property that
$\lambda_id^2_{g(-1)}(p,q_i)$ tending to infinity such that
$(M_\infty,g_\infty(-1))$ is the geometric limit of $(M,\lambda_ig(-1),q_i)$.
 Thus, according
to Theorem~\ref{topsplit}, the limit $(M_\infty,g_\infty(-1))$ splits as a
Riemannian product with a line. If the dimension of $M_\infty$ is two, then
this is a contradiction: We have that $(M_\infty,g_\infty(-1))$ splits as the
Riemannian product of a line and a one-manifold and hence is flat, but
$R(q_\infty,-1)=1$. Suppose that the dimension of $M_\infty$ is three. Once we
know that $(M_\infty,g_\infty(-1))$ splits as a product with a line, it follows
from the maximum principle (Corollary~\ref{localprod}) that the entire flow
splits as a product with a line, and the Ricci flow on the surface has strictly
positive curvature.
\end{proof}

\section{Classification of gradient shrinking solitons in dimensions $2$ and $ 3$}

In this section we fix $\kappa>0$ and we classify all
$\kappa$-solutions $(M,g_\infty(t)),\ -\infty<t<0$, that satisfy the
gradient shrinking soliton equation at the time-slice $t=-1$ in the
sense that there is a function $f\colon M \to \Ar$ such that
\begin{eqnarray}\label{GSS-1}{\rm
Ric}_{g_\infty(-1)}+{\rm Hess}^{g_\infty(-1)}(f)-\frac{1}{2}
g_\infty(-1)=0.\end{eqnarray} This will give a classification of the
two- and three-dimensional asymptotic gradient shrinking solitons
constructed in Theorem~\ref{asympGSS}.

 Let us give some examples in
dimensions two and three of ancient solutions that have such
functions. It turns out, as we shall see below, that in dimensions
two and three the only such are compact manifolds of constant
positive curvature -- i.e., Riemannian manifolds finitely covered by
the round sphere. We can create another, non-flat gradient shrinking
soliton in dimension three by taking $(M,g_{-1})$ equal to the
product of $(S^2,h_{-1})$, the round sphere of Gaussian curvature
$1/2$, with the real line (with the metric on the real line denoted
$ds^{ 2}$) and setting $g(t)=|t|h_{-1}+ ds^2$ for all $t<0$. We
define $f\colon M\times (-\infty,0)\to\Ar$ by $f(p,t)=s^2/4|t|$
where $s\colon M\to \Ar$ is the projection onto the second factor.
Then it is easy to see that
$${\rm Ric}_{g(t)}+{\rm Hess}^{g(t)}(f)-\frac{1}{2|t|}g(t)=0,$$
so that this example is a gradient shrinking soliton. There is a
free, orientation-preserving involution on this Ricci flow: the
product of the sign change on $\Ar$ with the antipodal map on $S^2$.
This preserves the family of metrics and hence there is an induced
Ricci flow on the quotient. Since this involution also preserves the
function $f$, the quotient  is also a gradient shrinking soliton.
These are the basic $3$-dimensional examples. As the following
theorem shows, they are all the $\kappa$-non-collapsed gradient
shrinking solitons in dimension three.

First we need a definition for a single Riemannian manifold analogous to a
definition we have already made for Ricci flows.

\begin{defn}
Let $(M,g)$ be an $n$-dimensional complete Riemannian manifold and
fix $\kappa>0$. We say that $(M,g)$ is {\em $\kappa$-non-collapsed}
if for every $p\in M$ and any $r>0$, if $|{\rm Rm}_g|\le r^{-2}$ on
$B(p,r)$ then ${\rm Vol}\, B(p,r)\ge \kappa r^n$.
\end{defn}

Here is the theorem that we shall prove:

\begin{thm}\label{GSSclass}
Let $(M,g)$ be a  complete, non-flat Riemannian manifold of bounded
non-negative curvature of dimension $2$ or $3$. Suppose that the
Riemannian manifold $(M,g)$ is $\kappa$-non-collapsed. Lastly,
suppose that there is a $C^2$-function $f\colon M \to \Ar$ such that
$${\rm Ric}_{g}+{\rm Hess}^{g}(f)=\frac{1}{2}g.$$
 Then there is a Ricci flow $(M,G(t)),\ -\infty<t<0$, with $G(-1)=g$ and with $(M,G(t))$ isometric
 to $(M,|t|g)$ for every $t<0$. In addition, $(M,G(t))$ is of one
 of the following three types:
\begin{enumerate}
\item[(1)] The flow $(M,G(t)),\ -\infty<t<0$, is a shrinking family of
compact, round (constant positive curvature) manifolds.
\item[(2)] The flow $(M,G(t)),\ -\infty<t<0$, is a product of a shrinking family
of round $2$-spheres with the real line. \item[(3)] $(M,G(t))$ is
isomorphic to the quotient family of metrics of the product of a
shrinking family of round $2$-spheres and the real line under the
action of an isometric involution.
\end{enumerate}
\end{thm}

Now let us begin the proof of Theorem~\ref{GSSclass}

\subsection{Integrating $\nabla f$}

Since the curvature of $(M,g)$ is bounded, it follows immediately from the
gradient shrinking soliton equation that ${\rm Hess}^g(f)$ is bounded. Fix a
point $p\in M$. For any $q\in M$ let $\gamma(s)$ be a minimal geodesic from $p$
to $q$ parameterized at unit length. Since
$$\frac{d}{ds}\left(|\nabla f(\gamma(s))|\right)^2=2\langle {\rm
Hess}(f)(\gamma'(s),\nabla f(\gamma(s))\rangle,$$ it follows that
$$\frac{d}{ds}\left(|\nabla f(\gamma(s))|\right)\le C,$$
where $C$ is an upper bound for $|{\rm Hess}(f)|$.
 By integrating, it
follows  that
$$|\nabla f(q)|\le Cd_g(p,q)+|\nabla f(p)|.$$
This means that any flow line $\lambda(t)$ for $\nabla f$ satisfies
$$\frac{d}{dt}d_g(p,\lambda(t))\le Cd_g(p,\lambda (t))+|\nabla f(p)|,$$
and hence these flow lines do not escape to infinity in finite time.
It follows that there is a flow $\Phi_t\colon M\to M$ defined for
all time with $\Phi_0={\rm Id}$ and $\partial\Phi_t/\partial
t=\nabla f$. We consider the one-parameter family of diffeomorphisms
$\Phi_{-{\rm log}(|t|)}\colon M\to M$ and define
\begin{equation}\label{Gdefn} G(t)=|t|\Phi_{-{\rm log}(|t|)}^*g.\ \
\ -\infty<t<0.\end{equation} We compute
$$\frac{\partial G}{\partial t}=-\Phi_{h(t)}^*g+2\Phi_{h(t)}^*{\rm
Hess}^g(f)=-2\Phi_{h(t)}^*{\rm Ric}(g)=-2{\rm Ric}(G(t)),$$ so that
$G(t)$ is a Ricci flow. Clearly, every time-slice is a complete,
non-flat manifold of non-negative bounded curvature. It is clear
from the construction that $G(-1)=g$ and that $(M,G(t))$ is
isometric to $(M,|t|g)$. This shows that $(M,g)$ is the $-1$
time-slice of a Ricci flow $(M,G(t))$ defined for all $t<0$, and
that, furthermore, all the manifolds $(M,G(t))$ are equivalent up to
diffeomorphism and scaling by $|t|$.

\subsection{Case 1: $M$ is compact and the curvature is strictly
positive}\label{compactcs}

\begin{claim} Suppose that $(M,g)$ and $f\colon M\to \Ar$ satisfies the hypotheses of
Theorem~\ref{GSSclass} and that $M$ is compact and of positive curvature. Then
the Ricci flow $(M,G(t))$ with $G(-1)=g$ given in Equation~(\ref{Gdefn}) is a
shrinking family of compact round manifolds.
\end{claim}

\begin{proof} The
manifold $(M,G(t))$ given in Equation~(\ref{Gdefn}) is equivalent up to
diffeomorphism and scaling by $|t|$ to $(M,g)$. If the dimension of $M$ is
three, then according to Hamilton's pinching toward positive curvature result
(Theorem~\ref{flowtoround}), the Ricci flow becomes singular in finite time and
as it becomes singular the metric approaches constant curvature in the sense
that the ratio of the largest sectional curvature to the smallest goes to one.
But this ratio is invariant under scaling and diffeomorphism, so that it must
be the case that for each $t$, all the sectional curvatures of the metric
$G(t)$ are equal; i.e., for each $t$ the metric $G(t)$ is round. If the
dimension of $M$ is two, then the results go back to Hamilton in
\cite{Hamiltonsurface}. According to Proposition 5.21 on p.. 118 of
\cite{ChowKnopf}, $M$ is a shrinking family of constant positive curvature
surfaces, which must be either $S^2$ or $\Ar P^2$. This completes the analysis
in the compact case.
\end{proof}

From this result, we can easily deduce a complete classification of
$\kappa$-solutions with compact asymptotic gradient shrinking
soliton.

\begin{cor}\label{compactcase}
Suppose that $(M,g(t))$ is a $\kappa$-solution of dimension $3$ with a compact
asymptotic gradient shrinking soliton. Then the Ricci flow $(M,g(t))$ is
isomorphic to a time-shifted version of its asymptotic gradient shrinking
soliton.
\end{cor}

\begin{proof}
We suppose that the compact asymptotic gradient shrinking soliton is
the limit of the $(M,g_{\tau_n}(t),(q_n,-1))$ for some sequence of
$\tau_n\rightarrow\infty$. Since by the discussion in the compact
case, this limit is of constant positive curvature. It follows that
for all $n$ sufficiently large, $M$ is diffeomorphic to the limit
manifold and the metric $g_{\tau_n}(-1)$ is close to a metric of
constant positive curvature.  In particular, for all $n$
sufficiently large, $(M,g_{\tau_n}(-1))$ is compact and of strictly
positive curvature. Furthermore, as $n\rightarrow\infty$
$\tau_n\rightarrow\infty$ and Riemannian manifolds
$(M,g_{\tau_n}(-1))$ become closer and closer to round in the sense
that the ratio of its largest sectional curvature to its smallest
sectional curvature goes to one. Since this is a scale invariant
ratio, the same is true for the sequence of Riemannian manifolds
$(M,g(-\tau_n))$. In the case when the dimension of $M$ is three, by
Hamilton's pinching toward round result or Ivey's theorem (see
Theorem~\ref{flowtoround}), this implies that the $(M,g(t))$ are all
exactly round.

This proves that $(M,g(t))$ is a shrinking family of round metrics.
The only invariants of such a family are the diffeomorphism type of
$M$ and the time $\Omega$ at which the flow becomes singular. Of
course, $M$ is diffeomorphic to its asymptotic soliton. Hence, the
only remaining invariant is the singular time, and hence $(M,g(t))$
is equivalent to a time-shifted version of its asymptotic soliton.
\end{proof}

\subsection{Case 2: Non-strictly positively curved.}

\begin{claim}\label{flatssplit}
Suppose that $(M,g)$ and $f\colon M\to \Ar$ are as in the statement
of Theorem~\ref{GSSclass} and that $(M,g)$ does not have strictly
positive curvature. Then $n=3$ and the Ricci flow $(M,G(t))$ with
$G(-1)=g$ given in Equation~(\ref{Gdefn}) has a one- or two-sheeted
covering that is a product of a two-dimensional
$\kappa$-non-collapsed Ricci flow of positive curvature and a
constant flat copy of $\Ar$. The curvature is bounded on each
time-slice.
\end{claim}

\begin{proof}
According to Hamilton's strong maximum principle
(Corollary~\ref{nullspace}), the Ricci flow $(M,G(t))$ has a one- or
two-sheeted covering that splits as a product of an evolving family
of manifolds of one dimension less of positive curvature and a
constant one-manifold.  It follows immediately that $n=3$.
 Let
$\tilde f$ be the lifting of $f$ to this one- or two-sheeted
covering. Let $Y$ be a unit tangent vector in the direction of the
one-manifold. Then it follows from Equation~(\ref{GSS-1}) that the
value of the Hessian of $\tilde f$ of $(Y,Y)$ is one. If the flat
one-manifold factor is a circle then there can be no such function
$\tilde f$. Hence, it follows that the one- or two-sheeted covering
is a product of an evolving surface with a constant copy of $\Ar$.
 Since
$(M,g)$ is $\kappa$-non-collapsed and of bounded curvature,
$(M,G(t))$ is $\kappa$-non-collapsed and each time-slice has
positive bounded curvature. These statements are also true for the
flow of surfaces.
\end{proof}

\subsection{Case 3: $M$ is non-compact and
strictly positively curved}\label{noncomp}

Here the main result is that this case does not occur.

\begin{prop}\label{noncompsol}
There is no  two- or three-dimensional  Ricci flow satisfying the hypotheses of
Theorem~\ref{GSSclass} with $(M,g)$ non-compact and of positive curvature.
\end{prop}

We suppose that we have $(M,g)$ as in Theorem~\ref{GSSclass} with
$(M,g)$ being non-compact and of positive curvature. Let $n$ be the
dimension of $M$, so that $n$ is either $2$ or $3$. Taking the trace
of the gradient shrinking soliton equation yields
$$R+\triangle f-\frac{n}{2}=0,$$
and consequently that $$ dR+d(\triangle f)=0.$$ Using
Lemma~\ref{lapformula} we rewrite this equation as
\begin{equation}\label{Rfeqn}
dR+\triangle(df)-{\rm Ric}(\nabla f,\cdot)=0.\end{equation}
 On the other
hand, taking the divergence of the gradient shrinking soliton
equation and using the fact that $\nabla^* g=0$ gives
$$\nabla^*{\rm Ric}+\nabla^*{\rm Hess}(f)=0.$$
Of course, $$\nabla^*{\rm Hess}(f)=\nabla^*(\nabla\nabla
f)=(\nabla^*\nabla)\nabla f=\triangle (df),$$ so that
$$\triangle(df)=-\nabla^*{\rm Ric}.$$ Plugging this into Equation~\ref{Rfeqn}
 gives
$$dR-\nabla^*{\rm Ric}-{\rm Ric}(\nabla f,\cdot)=0.$$
Now invoking Lemma~\ref{divRic} we have  \begin{equation}\label{dR}
dR=2{\rm Ric}(\nabla f,\cdot).\end{equation}

 Fix a point $p\in M$. Let
$\gamma(s);\ 0\le s\le \bar s$, be a shortest geodesic (with respect to the
metric $g$), parameterized at unit speed, emanating from $p$, and set
$X(s)=\gamma'(s)$.

\begin{claim} There is a constant $C$ independent of the choice of
$\gamma$ and of $\bar s$ such that
$$\int_0^{\bar s}{\rm Ric}(X,X)ds\le C.$$
\end{claim}

\begin{proof}
Since the curvature is bounded, clearly it suffices to assume that
$\bar s>>1$. Since $\gamma$ is length-minimizing and parameterized
at unit speed, it follows that it is a local minimum for the energy
functional $E(\gamma)=\frac{1}{2}\int_0^{\bar s}|\gamma'(s)|^2ds$
among all paths with the same end points. Thus, letting
$\gamma_u(s)=\gamma(s,u)$ be a one-parameter family of variations
(fixed at the endpoints) with $\gamma_0=\gamma$ and with
$d\gamma/du|_{u=0}=Y$, we see
$$0\le \delta_{Y}^2E(\gamma_u)=\int_0^{\bar
s}|\nabla_XY|^2+\langle {\mathcal R}(Y,X)Y,X\rangle ds.$$ We
conclude that
\begin{equation}\label{Yeqn}
\int_0^{\bar s}\langle -{\mathcal R}(Y,X)Y,X\rangle ds\le
\int_0^{\bar s}|\nabla_XY|^2ds.
\end{equation} Fix an orthonormal basis $\{E_i\}_{i=1}^n$  at $p$ with $E_n=X$, and let
$\widetilde E_i$ denote the parallel translation of $E_i$ along $\gamma$. (Of
course, $\widetilde E_n=X$.) Then, for $i\le n-1$, we define
$$Y_i=\left\{\begin{array}{ll} s\widetilde E_i & \mbox{if $0\le s\le
1$} \\ \widetilde E_i & \mbox{if $1\le s\le \bar s-1$} \\
(\bar s-s)E_i & \mbox{if $\bar s-1\le s\le \bar s$.} \end{array}
\right. $$ Adding up Equation~(\ref{Yeqn}) for each $i$ gives
$$-\sum_{i=1}^{n-1}\int_0^{\bar s}\langle{\mathcal R}(Y_i,X)Y_i,X\rangle ds\le
\sum_{i=1}^{n-1}\int_0^{\bar s}|\nabla_XY_i|^2ds.$$ Of course, since the
$\widetilde E_i$ are parallel along $\gamma$, we have
$$|\nabla_XY_i|^2=\left\{\begin{array}{ll} 1 & \mbox{if $0\le s\le 1$}
\\ 0 & \mbox{if $1\le s\le \bar s-1$} \\ 1 & \mbox{if $\bar s-1\le s\le
\bar s$}\end{array}\right. ,$$ so that
$$\sum_{i=1}^{n-1}\int_0^{\bar s}|\nabla_XY_i|^2=2(n-1).$$
On the other hand,
$$-\sum_{i=1}^{n-1}\langle{\mathcal R}(Y_i,X)(Y_i),X\rangle =\left\{
\begin{array}{ll} s^2{\rm Ric}(X,X) & \mbox{if $0\le s\le 1$} \\
{\rm Ric}(X,X) & \mbox{if $1\le s\le \bar s-1$} \\ (\bar s-s)^2{\rm
Ric}(X,X) & \mbox{if $\bar s-1\le s\le \bar s$.} \end{array} \right.
$$ Since the curvature is bounded and $|X|=1$, we see that
$\int_0^{\bar s}(1-s^2){\rm Ric}(X,X)ds+\int_{\bar s-1}^{\bar
s}(\bar s-s)^2{\rm Ric}(X,X)$ is bounded independent of $\gamma$ and
of $\bar s$. This concludes the proof of the claim.
\end{proof}

\begin{claim}
$|{\rm Ric}(X,\cdot)|^2\le R \cdot {\rm Ric}(X,X)$.
\end{claim}

\begin{proof}
This is obvious if $n=2$, so we may as well assume that $n=3$. We
diagonalize ${\rm Ric}$ in an orthonormal basis $\{e_i\}$. Let
$\lambda_i\ge 0$
 be the eigenvalues. Write $X=X^ie_i$ with $\sum_i(X^i)^2=1$. Then
$${\rm Ric}(X,\cdot)=X^i\lambda_i(e_i)^*,$$ so that
$|{\rm Ric}(X,\cdot)|^2=\sum_i(X^i)^2\lambda_i^2$. Of course, since
the $\lambda_i\ge 0$, this gives
$$R\cdot {\rm Ric}(X,X)=(\sum_i\lambda_i)\sum_i\lambda_i(X^i)^2\ge
\sum_i\lambda_i^2(X^i)^2,$$ establishing the claim.
\end{proof}

Now we compute, using Cauchy-Schwarz,
\begin{eqnarray*}
\left(\int_0^{\bar s}|{\rm Ric}(X,\widetilde E_i)|ds\right)^2 & \le
&\bar s\int_0^{\bar s}|{\rm Ric}(X,\widetilde E_i)|^2ds\le \bar
s\int_0^{\bar s}|{\rm Ric}(x,\cdot)|^2ds\\ & \le &  \bar
s\int_0^{\bar s}R\cdot {\rm Ric}(X,X)ds.\end{eqnarray*} Since $R$ is
bounded, it follows from the first claim that there is a constant
$C'$ independent of $\gamma$ and $\bar s$ with
\begin{equation}\label{sqrts}
\int_0^{\bar s}|{\rm Ric}(X,\widetilde E_i)|ds\le C'\sqrt{\bar
s}.\end{equation} Since $\gamma$ is a geodesic in the metric $g$, we have
$\nabla_XX=0$. Hence,
$$\frac{d^2f(\gamma(s))}{ds^2}=X(X(f))={\rm Hess}(f)(X,X).$$
Applying  the gradient shrinking soliton equation to the pair
$(X,X)$ gives
$$\frac{d^2f(\gamma(s))}{ds^2}=\frac{1}{2}-{\rm Ric}_{g}(X,X).$$
Integrating we see
$$\frac{df(\gamma(s))}{ds}|_{s=\bar
s}=\frac{df(\gamma(s))}{ds}|_{s=0}+\frac{\bar s}{2}-\int_0^{\bar
s}{\rm Ric}(X,X)ds.$$ It follows that
\begin{equation}\label{X}
X(f)(\gamma(\bar s))\ge \frac{\bar s}{2}-C'',\end{equation} for some
constant $C''$ depending only on $(M,g)$ and $f$. Similarly,
applying the gradient shrinking soliton equation to the pair
$(X,\widetilde E_i)$, using Equation~(\ref{sqrts}) and the fact that
$\nabla_X\widetilde E_i=0$ gives
\begin{equation}\label{sqrts2}
|\widetilde E_i(f)(\gamma(\bar s))|\le C''(\sqrt{\bar
s}+1).\end{equation}

These two inequalities imply that for $\bar s$ sufficiently large,
$f$ has no critical points  and that $\nabla f$ makes a small angle
with the gradient of the distance function from $p$, and $|\nabla
f|$ goes to infinity as the distance from $p$ increases. In
particular, $f$ is a proper function going off to $+\infty$ as we
approach infinity in $M$.

Now apply Equation~(\ref{dR}) to see that $R$ is increasing along the gradient
curves of $f$. Hence, there is a sequence $p_k$ tending to infinity in $M$ with
${\rm lim}_kR(p_k)={\rm limsup}_{q\in M} R_g(q)>0$.

The Ricci flow $(M,G(t)),\ -\infty<t<0$, given in
Equation~(\ref{Gdefn}) has the property that $G(-1)=g$ and that
$(M,G(t))$ is isometric to $(M,|t|g)$. Since the original Riemannian
manifold $(M,g)$ given in the statement of Theorem~\ref{GSSclass} is
$\kappa$-non-collapsed, it follows that, for every $t<0$, the
Riemannian manifold $(M,G(t))$ is $\kappa$-non-collapsed.
Consequently, the Ricci flow $(M,G(t))$ is $\kappa$-non-collapsed.
It clearly has bounded non-negative curvature on each time-slice and
is non-flat.  Fix a point $p\in M$. There is a sequence of points
$p_i$ tending to infinity with $R(p_i,-1)$ bounded away from zero.
It follows  that ${\rm lim}_{i\rightarrow
\infty}R(p_i,-1)d_{g(-1)}^2(p,p_i)=\infty$. Thus, this flow
satisfies all the hypotheses of Proposition~\ref{prodatinf}. Hence,
by Corollary~\ref{compact2D} we see that $n$ cannot be equal to two.
Furthermore, by Proposition~\ref{prodatinf}, when $n=3$ there is
another subsequence $q_i$ tending to infinity in $M$ such that there
is a geometric limit $(M_\infty,g_\infty(t),(q_\infty,-1)), \
-\infty<t\le -1$, of the flows $(M,G(t),(q_i,-1))$ defined for all
$t<0$ and this limit splits as a product of a surface flow
$(\Sigma^2,h(t))$ times the real line where the surfaces
$(\Sigma^2,h(t))$ are all of positive, bounded curvature and the
surface flow is $\kappa$-non-collapsed. Since there is a constant
$C<\infty$ such that the curvature of $(M,G(t)), -\infty<t\le
t_0<0$, is bounded by $C/|t_0|$, this limit actually exists for
$-\infty<t<0$ with the same properties.

Let us summarize our progress to date.

\begin{cor}\label{compact2Dsol}
There is no non-compact, two-dimensional Riemannian manifold $(M,g)$
satisfying the hypotheses of Theorem~\ref{GSSclass}. For any
non-compact three-manifold $(M,g)$ of positive curvature satisfying
the hypotheses of Theorem~\ref{GSSclass}, there is a sequence of
points $q_i\in M$ tending to infinity such that ${\rm
lim}_{i\rightarrow \infty}R_g(q_i)={\rm sup}_{p\in M}$ such that the
based Ricci flows $(M,G(t),(q_i,-1))$ converge to a Ricci flow
$(M_\infty,G_\infty(t),(q_\infty,-1))$ defined for $-\infty<t<0$
that splits as a product of a line and a family of surfaces, each of
positive, bounded curvature $(\Sigma^2,h(t))$. Furthermore, the flow
of surfaces is $\kappa$-non-collapsed.
\end{cor}

\begin{proof}
In Claim~\ref{flatssplit} we saw that every two-dimensional $(M,g)$ satisfying
the hypotheses of Theorem~\ref{GSSclass} has strictly positive curvature. The
argument that we just completed shows that there is no non-compact
two-dimensional example of strictly positive curvature.

The final statement  is exactly what we just established.
\end{proof}

\begin{cor}\label{2DGSS}
\begin{enumerate}
\item[(1)] Let $(M,g(t))$ be a two-dimensional Ricci flow satisfying all
the hypotheses of Proposition~\ref{2DGSS} except possible the
non-compactness hypothesis. Then $M$ is compact and for any $a>0$
the restriction of the flow to any interval of the form
$(-\infty,-a]$ followed by a shift of time by $+a$ is a
$\kappa$-solution.
\item[(2)] Any asymptotic gradient
shrinking soliton for a two-dimensional $\kappa$-solution is a shrinking family
of round surfaces.
\item[(3)] Let $(M,g(t)), -\infty<t\le 0$, be a
two-dimensional $\kappa$-solution. Then $(M,g(t))$ is a shrinking family of
compact, round surfaces.
\end{enumerate}
\end{cor}

\begin{proof}
Let $(M,g(t))$ be a two-dimensional Ricci flow satisfying all the hypotheses of
Proposition~\ref{prodatinf} except possibly non-compactness. It then follows
from Corollary~\ref{compact2D} that $M$ is compact. This proves the first item.

Now suppose that $(M,g(t))$ is an asymptotic soliton for a
$\kappa$-solution of dimension two. If $(M,g(-1))$ does not have
bounded curvature, then there is a sequence $p_i\rightarrow \infty$
so that ${\rm lim}_{i\rightarrow \infty}R(p_i,-1)=\infty$. By this
and
 Theorem~\ref{asympGSS} the Ricci flow $(M,g(t))$ satisfies all the hypotheses
 of Proposition ~\ref{prodatinf}. But this contradicts Corollary~\ref{compact2D}.
 We conclude that $(M,g(-1))$ has bounded curvature. According to
 Corollary~\ref{compact2Dsol} this means that
 $(M,g(t))$ is compact.
Results going back to Hamilton in \cite{Hamiltonsurface} imply that
this compact asymptotic shrinking soliton is a shrinking family of
compact, round surfaces. For example, this result is contained in
Proposition 5.21 on p. 118 of \cite{ChowKnopf}.
  This proves the second item.

 Now suppose that $(M,g(t))$ is a two-dimensional $\kappa$-solution.
 By the second item any asymptotic gradient shrinking soliton for
 this $\kappa$-solution is compact.
It follows that $M$ is compact. We know that as $t$ goes to $-\infty$ the
Riemannian surfaces $(M,g(t))$ are converging to compact, round surfaces.
Extend the flow forward from $0$ to a maximal time $\Omega<\infty$. By Theorem
5.64 on p. 149 of \cite{ChowKnopf} the surfaces $(M,g(t))$ are also becoming
round as $t$ approaches $\Omega$ from below. Also, according to Proposition
5.39 on p. 134 of \cite{ChowKnopf} the entropy of the flow is weakly monotone
decreasing and is strictly decreasing unless the flow is a gradient shrinking
soliton. But we have seen that the limits at both $-\infty$ and $\Omega$ are
round manifolds, and hence of the same entropy. It follows that the
$\kappa$-solution is a shrinking family of compact, round surfaces.
\end{proof}

Now that we have shown that every two-dimensional $\kappa$-solution
is a shrinking family of round surfaces, we can complete the proof
of Proposition~\ref{noncompsol}. Let $(M,g)$ be a non-compact
manifold of positive curvature satisfying the hypotheses of
Theorem~\ref{GSSclass}. According to Corollary~\ref{2DGSS} the
limiting Ricci flow $(M_\infty,G_\infty(t))$ referred to in
Corollary~\ref{compact2Dsol} is the product of a line and a
shrinking family of round surfaces. Since $(M,g)$ is non-compact and
has positive curvature, it is diffeomorphic to $\Ar ^3$ and hence
does not contain an embedded copy of a projective plane. It follows
that the round surfaces are in fact  round two-spheres. Thus,
$(M_\infty,G_\infty(t)),\ -\infty<t<0$,  splits as the product of a
shrinking family $(S^2,h(t)),\ -\infty<t<0$, of round two-spheres
and the real line.

\begin{claim}
The scalar curvature of $(S^2,h(-1))$ is equal to $1$.
\end{claim}

\begin{proof}
Since the shrinking family of round two-spheres $(S^2,h(t))$ exists for all
$-\infty<t<0$, it follows that the scalar curvature of $(S^2,h(-1))$ is at most
$1$. On the other hand, since the scalar curvature is increasing along the
gradient flow lines of $f$, the infimum of the scalar curvature of $(M,g)$,
$R_{\rm inf}$, is positive. Thus, the infimum of the scalar curvature of
$(M,G(t))$ is $R_{\rm inf}/|t|$ and goes to infinity as $|t|$ approaches $0$.
Thus, the infimum of the scalar curvature of  $(S^2,h(t))$ goes to infinity as
$t$ approaches zero. This means that the shrinking family of two-spheres
becomes singular as $t$ approaches zero, and consequently the scalar curvature
of $(S^2,h(-1))$ is equal to $1$.
\end{proof}

 It follows that for any $p$ in a
neighborhood of infinity of $(M,g)$, we have
$$R_g(p)< 1.$$
For any unit vector $Y$ at any point of $M\setminus K$ we have
$${\rm Hess}(f)(Y,Y)=\frac{1}{2}-{\rm Ric}(Y,Y)\ge \frac{1}{2}-\frac{R}{2}>0.$$
(On a manifold with non-negative curvature ${\rm Ric}(Y,Y)\le R/2$ for any unit
tangent vector $Y$.) This means that for $u$ sufficiently large the level
surfaces of $N_u=f^{-1}(u)$ are convex and hence have increasing area as $u$
increases.

According to Equations~(\ref{X}) and~(\ref{sqrts2}) the angle between $\nabla
f$ and the gradient of the distance function from $p$ goes to zero as we go to
infinity. According to Theorem~\ref{topsplit} the gradient of the distance
function from $p$ converges to the unit vector field in the $\Ar$-direction of
the product structure. It follows that the unit vector in the $\nabla
f$-direction converges to the unit vector in the $\Ar$-direction. Hence, as $u$
tends to $\infty$ the level surfaces $f^{-1}(u)$ converge in the $C^1$-sense to
$\Sigma\times \{0\}$. Thus, the areas of these level surfaces converge to the
area of $(\Sigma,h(-1))$ which is $8\pi$ since the scalar curvature of this
limiting surface is ${\rm limsup}_{p\in M} R(p,-1)= 1$. It follows that the
area of $f^{-1}(u)$ is less than $8\pi$ for all $u$ sufficiently large.

Now let us estimate the intrinsic curvature of $N=N_u=f^{-1}(u)$.
Let $K_N$ denote the sectional curvature of the induced metric on
$N$, whereas $K_M$ is the sectional curvature of $M$. We also denote
by $R_N$ the scalar curvature of the induced metric on $N$. Fix an
orthonormal basis $\{e_1,e_2,e_3\}$ at a point of $N$, where
$e_3=\nabla f/|\nabla f|$. Then by the Gauss-Codazzi formula we have
$$R_N=2K_N(e_1,e_2)=2(K_M(e_1,e_2)+{\rm det}\, S)$$
where $S$ is the shape operator
$$S=\frac{{\rm Hess}(f|TN)}{|\nabla f|}.$$
Clearly, we have $R-2{\rm Ric}(e_3,e_3)=2K_M(e_1,e_2)$, so that
$$R_N=R-2{\rm Ric}(e_3,e_3)+2{\rm det}\, S.$$
We can assume that the basis is chosen so that ${\rm Ric}|_{TN}$ is
diagonal; i.e., in the given basis we have
$${\rm Ric}=\begin{pmatrix} r_1 & 0  & c_1 \\ 0 & r_2 & c_2 \\ c_1 & c_2 &
r_3\end{pmatrix}.$$ From the gradient shrinking soliton equation we
have ${\rm Hess}(f)=(1/2)g-{\rm Ric}$ so that \begin{eqnarray*} {\rm
det}({\rm
Hess}(f|TN)) & = & \left(\frac{1}{2}-r_1\right)\left(\frac{1}{2}- r_2\right) \\
& \le &
\frac{1}{4}(1-r_1-r_2)^2 \\
& = & \frac{1}{4} (1-R+{\rm Ric}(e_3,e_3))^2.
\end{eqnarray*}
Thus, it follows that \begin{equation}\label{R_N} R_N\le R-2{\rm
Ric}(e_3,e_3)+\frac{(1-R+{\rm Ric}(e_3,e_3))^2}{2|\nabla f|^2}.
\end{equation}
 It follows from Equation~(\ref{X}) that $|\nabla
f(x)|\rightarrow\infty$ as $x$ goes to infinity in $M$. Thus, since
the curvature of $(M,g(-1))$ is bounded, provided that $u$ is
sufficiently large, we have $1-R+{\rm Ric}(e_3,e_3)<2|\nabla f|^2$.
Since the left-hand side of this inequality is positive (since $R<
1$), it follows that
\begin{equation*}
(1-R+{\rm Ric}(e_3,e_3))^2  <  2(1-R+{\rm Ric}(e_3,e_3))|\nabla
f|^2.
\end{equation*}
Plugging this into Equation~(\ref{R_N}) gives that
$$R_N<1-{\rm Ric}(e_3,e_3)\le 1,$$
assuming that  $u$ is sufficiently large.

This contradicts the Gauss-Bonnet theorem for the surface $N$: Its
area is less than $8\pi$, and the  scalar curvature of the induced
metric is less than $1$, meaning that its Gaussian curvature is less
than $1/2$; yet $N$ is diffeomorphic to a $2$-sphere. This completes
the proof of Proposition~\ref{noncompsol}, that is to say this shows
that there are no non-compact positive curved examples satisfying
the hypotheses of Theorem~\ref{GSSclass}.

\subsection{Case of non-positive curvature revisited}

We return now to the second case of Theorem~\ref{GSSclass}.  We
extend $(M,g)$ to a Ricci flow $(M,G(t))$ defined for $-\infty<t<0$
as given in Equation~(\ref{Gdefn}). By Claim~\ref{flatssplit} $M$
has either a one- or $2$-sheeted covering $\widetilde M$ such that
$(\widetilde M,\tilde G(t))$ is a metric product of a surface and a
one-manifold for all $t<0$. The evolving metric on the surface is
itself a $\kappa$-solution and hence by Corollary~\ref{2DGSS} the
surfaces are compact and the metrics are all round. Thus, in this
case, for any $t<0$, the manifold $(\widetilde M,\tilde G(t))$ is a
metric product of a round $S^2$ or $\Ar P^2$ and a flat copy of
$\Ar$.  The conclusion in this case is that the one- or two-sheeted
covering $(\widetilde M,\tilde G(t))$ is a product of a round $S^2$
or $\Ar P^2$ and the line for all $t<0$.

\subsection{Completion of the proof of Theorem~\protect{\ref{GSSclass}}}

\begin{cor}\label{3Dinfty}
Let $(M,g(t))$ be a three-dimensional Ricci flow satisfying the hypotheses of
Proposition~\ref{prodatinf}. Then the limit constructed in that proposition
splits as a product of a shrinking family of compact round surfaces with a
line. In particular, for any non-compact gradient shrinking soliton of a
three-dimensional $\kappa$-solution the limit constructed in
Proposition~\ref{prodatinf} is the product of a shrinking family of round
surfaces and the real line.
\end{cor}

\begin{proof}
Let $(M,g(t))$ be a three-dimensional Ricci flow satisfying the
hypotheses of Proposition~\ref{prodatinf} and let $(N^2,h(t))\times
(\Ar,ds^2)$ be the limit constructed in that proposition. Since this
limit is $\kappa$-non-collapsed, $(N,h(t))$ is
$\kappa'$-non-collapsed for some $\kappa'>0$ depending only on
$\kappa$. Since the limit is not flat and has non-negative
curvature, the same is true for $(N,h(t))$. Since $\partial
R/\partial t\ge 0$ for the limit, the same is true for $(N,h(t))$.
That is to say $(N,h(t))$ satisfies all the hypotheses of
Proposition~\ref{prodatinf} except possibly non-compactness. It now
follows from Corollary~\ref{2DGSS} that $(N,h(t))$ is a shrinking
family of compact, round surfaces.
\end{proof}

\begin{cor}\label{3DGSSkappa}
Let $(M,g(t)),\ -\infty<t<0$, be an asymptotic gradient shrinking soliton for a
three-dimensional $\kappa$-solution. Then for each $t<0$, the Riemannian
manifold $(M,g(t))$ has bounded curvature. In particular, for any $a>0$ the
flow $(M,g(t)),\ -\infty<t\le -a$, followed by a shift of time by $+a$ is a
$\kappa$-solution.
\end{cor}

\begin{proof}
If an asymptotic gradient shrinking soliton $(M,g(t))$ of a three-dimensional
$\kappa$-solution does not have strictly positive curvature, then according to
Corollary~\ref{nullspace}, $(M,g(t))$ has a covering that splits as a product
of a a two-dimensional Ricci flow and a line. The two-dimensional Ricci flow
satisfies all the hypotheses of Proposition~\ref{prodatinf} except possibly
compactness, and hence by Corollary~\ref{2DGSS} it is a shrinking family of
round surfaces. In this case, it is clear that each time-slice of $(M,g(t))$
has bounded curvature.

Now we consider the remaining case when $(M,g(t))$ has strictly positive
curvature.  Assume that $(M,g(t))$ has unbounded curvature. Then there is a
sequence of points $p_i$ tending to infinity in $M$ such that $R(p_i,t)$ tends
to infinity. By Corollary~\ref{3Dinfty} we can replace the points $p_i$ by
points $q_i$ with $Q_i=R(q_i,t)\ge R(p_i,t)$ so that the based Riemannian
manifolds $(M,Q_ig(t),q_i)$ converge to a product of a round surface $(N,h(t))$
with $\Ar$. The surface $N$ is either diffeomorphic to $S^2$ or $\Ar P^2$.
Since $(M,g(t))$ has positive curvature, by Theorem~\ref{soul}, it is
diffeomorphic to $\Ar ^3$, and hence it contains no embedded $\Ar P^2$. It
follows that $(N,h(t))$ is a round two-sphere.

Fix $\epsilon>0$ sufficiently small as in Proposition~\ref{narrows}. Then the
limiting statement means that, for every $i$ sufficiently large, there is an
$\epsilon$-neck in $(M,g(t))$ centered at $q_i$ with scale $Q_i^{-1/2}$. This
contradicts Proposition~\ref{narrows}, establishing that for each $t<0$ the
curvature of $(M,g(t))$ is bounded.
\end{proof}

\begin{cor}
Let $(M,g(t)),\ -\infty<t\le 0$, be a $\kappa$-solution of dimension three.
Then any asymptotic gradient shrinking soliton $(M_\infty,g_\infty(t))$ for
this $\kappa$-solution, as constructed in Theorem~\ref{asympGSS}, is of one of
the three types listed in Theorem~\ref{GSSclass}.
\end{cor}

\begin{proof}
Let $(M_\infty,g_\infty(t)),\ -\infty<t<0$, be an asymptotic
gradient shrinking soliton for $(M,g(t))$. According to
Corollary~\ref{3DGSSkappa}, this soliton is a $\kappa$-solution,
implying that $(M_\infty,g_\infty(-1))$ is a complete Riemannian
manifold of bounded, non-negative curvature. Suppose that
$B(p,-1,r)\subset M_\infty$ is a metric ball and $|{\rm
Rm}_{g_\infty}|(x,-1)\le r^{-2}$ for all $x\in B(p,-1,r)$. Since
$\partial R_{g_\infty}(x,t)/\partial t\ge 0$, it follows that
$R(x,t)\le 3r^{-2}$ on $B(p,-1,r)\times (-1-r^2,-1]$, and hence that
$|{\rm Rm}_{g_\infty}|\le 3r^{-2}$ on this same region. Since the
Ricci flow $(M_\infty,g_\infty(t))$ is $\kappa$-non-collapsed, it
follows that ${\rm Vol}\, B(p,-1,r/\sqrt{3})\ge
\kappa(r/\sqrt{3})^3$. Hence, ${\rm Vol}, B(p,-1,r)\ge
(\kappa/3\sqrt{3})r^3$. This proves that the manifold
$(M_\infty,g_\infty(-1))$ is $\kappa'$-non-collapsed for some
$\kappa'>0$ depending only on $\kappa$. On the other hand, according
to Theorem~\ref{asympGSS} there is a function $f(\cdot, -1)$ from
$M_\infty$ to $\Ar$ satisfying the gradient shrinking soliton
equation at the time-slice $-1$.
 Thus,
Theorem~\ref{GSSclass} applies to $(M_\infty,g_\infty(-1))$ to
produce a Ricci flow $G(t),\ -\infty<t<0$, of one of the three types
listed in that theorem and with $G(-1)=g_\infty(-1)$.

Now we must show that $G(t)=g_\infty(t)$ for all $t<0$. In the first
case when $M_\infty$ is compact, this is clear by uniqueness of the
Ricci flow in the compact case. Suppose that $(M_\infty,G(t))$ is of
the second type listed in Theorem~\ref{GSSclass}. Then
$(M_\infty,g_\infty(-1))$ is a product of a round two-sphere and the
real line. By Corollary~\ref{nullspace} this implies that the entire
flow $(M_\infty,g_\infty(t))$ splits as the product of a flow of
compact two-spheres and the real line. Again by uniqueness in the
compact case, this family of two-spheres must be a shrinking family
of round two-spheres. In the third case, one passes to a finite
sheeted covering of the second type, and applies the second case.
\end{proof}

\subsection{Asymptotic curvature}

There is one elementary result that will be needed in what follows.

\begin{defn}
Let $(M,g)$ be a complete, connected, non-compact Riemannian manifold of
non-negative curvature. Fix a point $p\in M$. We define the {\em asymptotic
scalar curvature}\index{asymptotic scalar curvature}
$${\mathcal R}(M,g)={\rm
limsup}_{x\rightarrow\infty}R(x)d^2(x,p).$$ Clearly, this limit is
independent of $p$.
\end{defn}

\begin{prop}\label{asympscal}
Suppose that $(M,g(t)),\ -\infty<t<0$, is a connected, non-compact
$\kappa$-solution of dimension at most three\footnote{This result, in fact,
holds in all dimensions.}.  Then ${\mathcal R}(M,g(t))=+\infty$ for every
$t<0$.
\end{prop}

\begin{proof}
By Corollary~\ref{2DGSS} the only two-dimensional $\kappa$-solutions are
compact, so that the result is vacuously true in this case. Suppose that
$(M,g(t))$ is three-dimensional
 If $(M,g(t))$ does
not have strictly positive curvature, then, since it is not flat, by
Corollary~\ref{nullspace} it must be three-dimensional and it has a
finite-sheeted covering space that splits as a product $(Q,h(0))\times
(\Ar,ds^2)$ with $(Q,h(0))$ being a surface of strictly positive curvature and
$T$ being a flat one-manifold. Clearly, in this case the asymptotic curvature
is infinite.

Thus, without loss of generality we can assume that $(M,g(t))$ has
strictly positive curvature. Let us first consider the case when
${\mathcal R}(M,g(t))$ has a finite, nonzero value. Fix a point
$p\in M$. Take a sequence of points $x_n$ tending to infinity and
set $\lambda_n=d^2_0(x_n,p)$ and $Q_n=R(x_n,t)$. We choose this
sequence such that
$${\rm lim}_{n\rightarrow\infty}Q_n\lambda_n={\mathcal R}(M,g(t)).$$
We consider the sequence of Ricci flows $(M,h_n(t),(x_n,0))$, where
$$h_n(t)=Q_ng(Q_n^{-1}t).$$

 Fix $0<a<\sqrt{{\mathcal R}(M,g(t))}<b<\infty$. Consider  the annuli
$$A_n=\{y\in M\left|\right. a<d_{h_n(0)}(y,p)<b\}.$$
Because of the choice of sequence, for all $n$ sufficiently large,
the scalar curvature of the restriction of $h_n(0)$ to $A_n$ is
bounded independent of $n$. Furthermore, since $d_{h_n}(p,x_n)$
converges to $\sqrt{{\mathcal R}(M,g(t))}$, there is $\alpha>0$ such
that for all $n$ sufficiently large, the annulus $A_n$ contains
$B_{h_n}(x_n,0,\alpha)$. Consequently, we have a bound, independent
of $n$, for the scalar curvature of $h_n(0)$ on these balls. By the
hypothesis that $\partial R/\partial t\ge 0$, there is a bound,
independent of $n$, for the scalar curvature of $h_n$ on
$B_{h_n}(x_n,0,\alpha)\times (-\infty,0]$. Using the fact that the
flows have non-negative curvature, this means that there is a bound,
independent of $n$, for $|{\rm Rm}_{h_n}(y,0)|$ on
$B_{h_n}(x_n,0,\alpha)\times (-\infty,0]$.
 This means
that by Shi's theorem (Theorem~\ref{shi}), there are bounds,
independent of $n$, for every covariant derivative of the curvature
on $B_{h_n}(x_n,0,\alpha/2)\times (-\infty,0]$.

Since the original flow is $\kappa$-non-collapsed on all scales, it
follows that the rescaled flows are also $\kappa$ non-collapsed on
all scales. Since the curvature is bounded, independent of $n$, on
$B_{h_n}(x_n,0,\alpha)$, this implies that there is $\delta>0$,
independent of $n$, such that for all $n$ sufficiently large, every
ball of radius $\delta$ centered at any point of
$B_{h_n}(x_n,0,\alpha/2)$ has volume at least $\kappa \delta^3$, Now
applying Theorem~\ref{basicconv} we see that  a subsequence
converges geometrically to a limit which will automatically be a
metric ball $B_{g_\infty}(x_\infty,0,\alpha/2)$. In fact, by
Hamilton's result (Proposition~\ref{partialflowlimit}) there is a
limiting flow on $B_{g_\infty}(x_\infty,0,\alpha/4)\times
(-\infty,0]$. Notice that the limiting flow is not flat since
$R(x_\infty,0)=1$.

On the other hand, according to Lemma~\ref{conelimit} the
Gromov-Hausdorff limit\index{Gromov-Hausdorff limit} of a
subsequence $(M,\lambda_n^{-1}g_n(0),x_n)$ is the Tits
cone\index{Tits cone}, i.e., the cone over $S_\infty(M,p)$. Since
$Q_n={\mathcal R}(M,g(t))\lambda_n^{-1}$, the rescalings
$(M,Q_ng_n(0),x_n)$  also converge to  a cone, say $(C,h,y_\infty)$,
which is in fact simply a rescaling of the Tits cone by a factor
${\mathcal R}(M,g(t))$.  Pass to a subsequence so that both the
geometric limit on the ball of radius $\alpha/2$ and the
Gromov-Hausdorff limit exist. Then the geometric limit
$B_{g_\infty}(x_\infty,0,\alpha/2)$ is isometric to an open ball in
the cone. Since we have a limiting Ricci flow
$$(B_{g_\infty}(x_\infty,0,\alpha/2),g_\infty(t)),\ -\infty<t\le 0,$$
 this contradicts Proposition~\ref{nocones}. This completes the proof that
 it is not possible for the asymptotic curvature to be finite
and nonzero.

Lastly, we consider the possibility that the asymptotic curvature is zero.
Again we fix $p\in M$. Take any sequence of points $x_n$ tending to infinity
and let $\lambda_n=d^2_0(p,x_n)$. Form the sequence of based Ricci flows
$(M,h_n(t),(x_n,0))$ where $h_n(t)=\lambda_n^{-1}g(\lambda_nt)$. On the one
hand, the Gromov-Hausdorff limit (of a subsequence) is the Tits cone. On the
other hand, the curvature condition tells us the following: For any $0<a<1<b$
on the regions
$$\{y\in M\left|\right. a< d_{h_n(0)}(y,p)<b\},$$
the curvature tends uniformly to zero as $n$ tends to infinity. Arguing as in
the previous case, Shi's theorem, Hamilton's result,
Theorem~\ref{partialflowlimit}, and the fact that the original flow is $\kappa$
non-collapsed on all scales tells us that we can pass to a subsequence so that
these annuli centered at $x_n$ converge geometrically to a limit. Of course,
the limit is flat. Since this holds for all $0<a<1<b$, this implies that the
Tits cone is smooth and flat except possibly at its cone point. In particular,
the sphere at infinity, $S_\infty(M,p)$, is a smooth surface of constant
curvature $+1$.

\begin{claim}
 $S_\infty(M,p)$ is isometric to a round $2$-sphere.
 \end{claim}

\begin{proof}
Since $M$ is orientable the complement of the cone point in the Tits
cone is an orientable manifold and hence $S_\infty(M,p)$ is an
orientable surface. Since we have already established that it has a
metric of constant positive curvature, it must be diffeomorphic to
$S^2$, and hence isometric to a round sphere. (In higher dimensions
one can prove that $S_\infty(M,p)$ is simply connected, and hence
isometric to a round sphere.)
\end{proof}

It follows that the Tits cone is a smooth flat manifold even at the origin, and
hence is isometric to Euclidean $3$-space. This means that in the limit, for
any $r>0$ the volume of the ball of radius $r$ centered at the cone point
  is exactly $\omega_3 r^3$, where
$\omega_3$ is the volume of the unit ball in $\Ar^3$. Consequently,
$${\rm lim}_{n\rightarrow\infty}{\rm
Vol}\left(B_{g}(p,0,\sqrt{\lambda_n}r)\setminus
B_g(p,0,1)\right)\rightarrow \omega_3\lambda_n^{3/2}r^3.$$ By
Theorem~\ref{BishopGromov} and the fact that the Ricci curvature is
non-negative, this implies that
$${\rm Vol}\, B_g(p,0,R)=\omega_3R^3$$
for all $R<\infty$. Since the Ricci curvature is non-negative, this
means that $(M,g(t))$ is Ricci-flat, and hence flat. But this
contradicts the fact that $(M,g(t))$ is a $\kappa$-solution and
hence is not flat.

Having ruled out the other two cases, we are left with only one possibility:
${\mathcal R}(M,g(t))=\infty$.
\end{proof}

\section{Universal $\kappa$}

The first consequence of the existence of an asymptotic gradient
shrinking soliton is that there is a universal $\kappa$ for all
$3$-dimensional $\kappa$-solutions, except those of constant
positive curvature.

\begin{prop}\label{kappa0}
There is a $\kappa_0>0$ such that any non-round $3$-dimensional
$\kappa$-solution is a $\kappa_0$-solution.
\end{prop}

\begin{proof}
Let $(M,g(t))$ be a non-round $3$-dimensional $\kappa$-solution. By
Corollary~\ref{compactcase} since $(M,g(t))$ is not a family of
round manifolds, the asymptotic soliton for the $\kappa$-solution
cannot be compact. Thus, according to Corollary~\ref{GSSclass} there
are only two possibilities for the asymptotic soliton
$(M_\infty,g_\infty(t))$ -- either $(M_\infty,g_\infty(t))$ is the
product of a round $2$-sphere of Gaussian curvature $1/2|t|$ with a
line or has a two-sheeted covering by such a product. In fact, there
are three possibilities: $S^2\times \Ar$, $\Ar P^2\times \Ar$ or the
twisted $\Ar$-bundle over $\Ar P^2$ whose total space is
diffeomorphic to the complement of a point in $\Ar P^3$.

Fix a point $x=(p,0)\in M\times\{0\}$. Let $\bar\tau_k$ be a
sequence converging to $\infty$, and $q_k\in M$ a point with
$l_{x}(q_k,\bar\tau_k)\le 3/2$. The existence of an asymptotic
soliton means that, possibly after passing to a subsequence, there
is a gradient shrinking soliton $(M_\infty,g_\infty(t))$ and a ball
$B$ of radius $1$ in $(M_\infty,g_\infty(-1))$ centered at a point
$q_\infty\in M_\infty$ and a sequence of embeddings $\psi_k\colon
B\to M$ such that $\psi_k(q_\infty)=q_k$ and such that the map
$$B\times [-2,-1]\to M\times [-2\bar\tau_k,-\bar\tau_k]$$
given by $(b,t)\mapsto (\psi_k(b),\bar\tau_kt)$ has the property
that the pullback of $\bar\tau^{-1}_kg(\bar\tau_k t)$ converges
smoothly and uniformly as $k\rightarrow\infty$ to the restriction of
$g_\infty(t),\ -2\le t\le -1$, to $B$. Let $(M_k,g_k(t))$ be this
rescaling of the the $\kappa$-solution by $\tau_k$. Then the
embeddings $\psi_k\times {\rm id}\colon B\times (-2,-1]\to
(M_k\times [-2,-1]$ converge as $k\rightarrow \infty$ to a
one-parameter family of isometries. That is to say, the image
$\psi_k(B\times [-2,-1])\subset M_k\times [-2,-1]$ is an almost
isometric embedding. Since the reduced length function from $x$ to
$(\psi_k(a),-1)$ is at most $3/2$ (from the invariance of reduced
length under rescalings, see Corollary~\ref{COR}), it follows easily
that the reduced length function on $\psi_k(B\times \{-2\})$ is
bounded independent of $k$. Similarly, the volume of
$\psi_k(B\times\{-2\})$ is bounded independent of $k$. This means
the reduced volume of $\psi_k(B\times\{-2\})$ in $(M_k,g_k(t))$ is
bounded independent of $k$. Now by Theorem~\ref{THM} this implies
that $(M_k,g_k(t))$ is $\kappa_0$-non-collapsed at $(p,0)$ on scales
$\le \sqrt{2}$ for some $\kappa_0$ depending only on the geometry of
the three possibilities for $(M_\infty,g_\infty(t)),\ -2\le t\le
-1$. Being $\kappa_0$-non-collapsed is invariant under rescalings,
so that it follows immediately that $(M,g(t))$ is
$\kappa_0$-non-collapsed on scales $\le \sqrt{2\bar\tau_k}$. Since
this is true for all $k$, it follows that $(M,g(t))$ is
$\kappa_0$-non-collapsed on all scales at $(p,0)$.

This result holds of course for every $p\in M$, showing that at
$t=0$ the flow is $\kappa_0$-non-collapsed. To prove this result at
points of the form $(p,t)\in M\times(-\infty,0]$ we simply shift the
original $\kappa$-solution upward by $|t|$ and remove the part of
the flow at positive time. This produces a new $\kappa$-solution and
the point in question has been shifted to the time-zero slice, so
that we can apply the previous results.
\end{proof}

\section{Asymptotic volume}

Let $(M,g(t))$ be an $n$-dimensional $\kappa$-solution. For any
$t\le 0$ and any point $p\in M$ we consider $({\rm
Vol}\,B_{g(t)}(p,r))/r^n$. According to the Bishop-Gromov Theorem
(Theorem~\ref{BishopGromov}), this is a non-increasing function of
$r$. We define the {\em asymptotic volume}\index{asymptotic
volume|ii} ${\mathcal V}_\infty(M,g(t))$, or ${\mathcal
V}_\infty(t)$ if the flow is clear from the context, to be the limit
as $r\rightarrow\infty$ of this function. Clearly, this limit is
independent of $p\in M$.

\begin{thm}\label{asympvol}
For\footnote{This theorem and all the other results of this section are valid
in all dimensions. Our proofs use Theorem~\ref{asympscal} and
Proposition~\ref{prodatinf} which are also valid in all dimensions but which we
proved only in dimensions $2$ and $3$. Thus, while we state the results of this
section for all dimensions, strictly speaking we give proofs only for
dimensions $2$ and $3$. These are the only cases we need in what follows.} any
$\kappa>0$ and any $\kappa$-solution $(M,g(t))$ the asymptotic volume
${\mathcal V}_\infty(M,g(t))$ is identically zero.
\end{thm}

\begin{proof}
The proof is by induction on the dimension $n$ of the solution. For $n=2$ by
Corollary~\ref{2DGSS} there are only compact $\kappa$-solutions, which clearly
have zero asymptotic volume. Suppose that we have established the result for
$n-1\ge 2$ and let us prove it for $n$.

According to Proposition~\ref{prodatinf} there is a sequence of points $p_n\in
M$ tending to infinity such that setting $Q_n=R(p_n,0)$ the sequence of Ricci
flows
$$(M,Q_ng(Q_n^{-1}t),(q_n,0))$$ converges geometrically to a limit
$(M_\infty,g_\infty(t),(q_\infty,0))$, and this limit splits off a
line: $(M_\infty,g_\infty(t))= (N,h(t))\times \Ar$. Since the ball
of radius $R$ about a point $(x,t)\in N\times \Ar$ is contained in
the product of the ball of radius $R$ about $x$ in $N$ and an
interval of length $2R$, it follows that $(N,h(t))$ is a
$\kappa/2$-ancient solution. Hence, by induction, for every $t$, the
asymptotic volume of $(N,h(t))$ is zero, and hence so is that of
$(M,g(t))$.
\end{proof}

\subsection{Volume comparison}

One important consequence of the asymptotic volume result is a
volume comparison result.

\begin{prop}\label{volcomp} Fix the dimension $n$.
For every $\nu>0$ there is $A<\infty$ such that the following holds.
Suppose that $(M_k,g_k(t)),\ -t_k\le t\le 0$, is a sequence of (not
necessarily complete) $n$-dimensional Ricci flows of non-negative
curvature operator. Suppose in addition we have points $p_k\in M_k$
and radii $r_k>0$ with the property that for each $k$ the ball
$B(p_k,0,r_k)$ has compact closure in $M_k$. Let $Q_k=R(p_k,0)$ and
suppose that $R(q,t)\le 4Q_k$ for all $q\in B(p_k,0,r_k)$ and for
all $t\in [-t_k,0]$, and suppose that $t_kQ_k\rightarrow\infty$ and
$r_k^2Q_k\rightarrow\infty$ as $k\rightarrow\infty$. Then ${\rm
Vol}\,B(p_k,0,A/\sqrt{Q_k})< \nu(A/\sqrt{Q_k})^n$ for all $k$
sufficiently large.
\end{prop}

\begin{proof}
Suppose that the result fails for some $\nu>0$. Then there is a
sequence $(M_k,g_k(t)),\ -t_k\le t\le 0$, of $n$-dimensional Ricci
flows, points $p_k\in M_k$, and radii $r_k$ as in the statement of
the lemma such that for every $A<\infty$ there is an arbitrarily
large $k$ with ${\rm Vol}\,B(p_k,0,A/\sqrt{Q_k})\ge
\nu(A/\sqrt{Q_k})^n$. Pass to a subsequence so that for each
$A<\infty$ we have
$${\rm
Vol}\,B(p_k,0,A/\sqrt{Q_k})\ge\nu(A/\sqrt{Q_k})^n$$ for all $k$
sufficiently large. Consider now the flows
$h_k(t)=Q_kg_k(Q_k^{-1}t)$, defined for $-Q_kt_k\le t\le 0$. Then
for every $A<\infty$ for all $k$ sufficiently large we have
$R_{h_k}(q,t)\le 4$ for all $q\in B_{h_k}(p_k,0,A)$ and all $t\in
(-t_kQ_k,0]$. Also, for every $A<\infty$ for all $k$ sufficiently
large we have ${\rm Vol}\,B(p_k,0,A)\ge\nu A^n$. According to
Theorem~\ref{flowlimit} we can then pass to a subsequence that has a
geometric limit which is an ancient flow of complete Riemannian
manifolds. Clearly, the time-slices of the limit have non-negative
curvature operator, and the scalar curvature is bounded (by $4$) and
is equal to $1$ at the base point of the limit. Also, the asymptotic
volume ${\mathcal V}(0)\ge \nu$.

\begin{claim}
 Suppose that $(M,g(t))$ is an ancient Ricci flow such that for each
 $t\le 0$ the Riemannian manifold $(M,g(t))$ is complete and has
bounded, non-negative  curvature operator. Let ${\mathcal V}(t)$ be
the asymptotic volume of the manifold $(M,g(t))$.
\begin{enumerate}
\item[(1)] The asymptotic volume ${\mathcal V}(t)$ is a non-increasing function
of $t$. \item[(2)] If ${\mathcal V}(t)=V>0$ then every metric ball
$B(x,t,r)$ has volume at least $Vr^n$.
\end{enumerate}
\end{claim}

\begin{proof}
We begin with the proof of the first item. Fix $a<b\le 0$. By
hypothesis there is a constant $K<\infty$ such that the scalar
curvature of $(M,g(0))$ is bounded by $(n-1)K$. By the Harnack
inequality (Corollary~\ref{posderiv}) the scalar curvature of
$(M,g(t))$ is bounded by $(n-1)K$ for all $t\le 0$. Hence, since the
$(M,g(t))$ have non-negative curvature, we have ${\rm Ric}(p,t)\le
(n-1)K$ for all $p$ and $t$. Set $A=4(n-1)\sqrt{\frac{2K}{3}}$. Then
by Corollary~\ref{corI.8.3} we have
$$d_a(p_0,p_1)\le d_b(p_0,p_1)+A(b-a).$$
 This means that for any $r>0$ we have
$$B(p_0,b,r)\subset B(p_0,a,r+A(b-a)).$$

On the other hand, since $d{\rm Vol}/dt=-Rd{\rm Vol}$, it follows
that in the case of non-negative curvature that the volume of any
open set is non-increasing in time. Consequently,
$${\rm Vol}_{g(b)}B(p_0,b,r)\le {\rm
Vol}_{g(a)}B(p_0,a,r+A(b-a)),$$ and hence
$$\frac{{\rm Vol}_{g(b)}B(p_0,b,r)}{r^n}\le\frac{{\rm
Vol}_{g(a)}B(p_0,a,r+A(b-a))}{(r+A(b-a))^n}\frac{(r+A(b-a))^n}{r^n}.$$
Taking the limit as $r\rightarrow \infty$ gives
$${\mathcal V}(b)\le {\mathcal V}(a).$$

The second item of the claim is immediate from the Bishop-Gromov
inequality (Theorem~\ref{BishopGromov}).
\end{proof}

Now we return to the proof of the proposition. Under the assumption
that there is a counterexample to the proposition for some $\nu>0$,
we have constructed a limit that is an ancient Ricci flow with
bounded, non-negative curvature with ${\mathcal V}(0)\ge \nu$. Since
${\mathcal V}(0)\ge \nu$, it follows from the claim that ${\mathcal
V}(t)\ge \nu$ for all $t\le 0$ and hence, also by the claim,  we see
that $(M,g(t))$ is $\nu$-non-collapsed for all $t$. This completes
the proof that the limit is a $\nu$-solution. This contradicts
Theorem~\ref{asympvol} applied with $\kappa=\nu$, and proves the
proposition.
\end{proof}

This proposition has two useful corollaries about balls in
$\kappa$-solutions with volumes bounded away from zero. The first
says that the normalized curvature is bounded on such balls.

\begin{cor}\label{bdvolbdcurv}
For any $\nu>0$ there is a $C=C(\nu)<\infty$ depending only on the
dimension $n$ such that the following holds. Suppose that
$(M,g(t)),\ -\infty<t\le 0$, is an $n$-dimensional Ricci flow with
each $(M,g(t))$ being complete and with bounded, non-negative
curvature operator. Suppose $p\in M$, and $r>0$ are such that ${\rm
Vol}\,B(p,0,r)\ge \nu r^n$. Then $r^2R(q,0)\le C$ for all $q\in
B(p,0,r)$.
\end{cor}

\begin{proof}
Suppose that the result fails for some $\nu>0$. Then there is a
sequence $(M_k,g_k(t))$ of $n$-dimensional Ricci flows, complete,
with bounded non-negative curvature operator and  points $p_k\in
M_k$, constants $r_k>0$, and points $q_k\in B(p_k,0,r_k)$ such that:
\begin{enumerate}
\item[(1)] ${\rm Vol}\,B(p_k,0,r_k)\ge \nu r_k^n$, and
\item[(2)] setting $Q_k=R(q_k,0)$ we have $r_k^2Q_k\rightarrow\infty$ as $k\rightarrow\infty$.
\end{enumerate}
Using Lemma~\ref{4times} we can find points $q_k'\in B(p_k,0,2r_k)$
and constants $s_k\le r_k$, such that setting $Q_k'=R(q_k',0)$ we
have $Q_k's_k^2=Q_kr_k^2$ and  $R(q,0)<4Q_k'$ for all $q\in
B(q_k',0,s_k)$.  Of course, $Q'_ks_k^2\rightarrow\infty$ as
$k\rightarrow\infty$. Since $d_0(p_k,q'_k)<2r_k$, we have
$B(p_k,0,r_k)\subset B(q_k',0,3r_k)$ so that
$${\rm Vol}\,B(q_k',0,3r_k)\ge {\rm Vol}\,B(p_k,0,r_k)\ge \nu
r_k^n=(\nu/3^n)(3r_k)^n.$$ Since the sectional curvatures of
$(M,g_k(0))$ are non-negative, it follows from the Bishop-Gromov
inequality (Theorem~\ref{BishopGromov}) that ${\rm
Vol}\,B(q_k',0,s)\ge (\nu/3^n)s^n$ for any $s\le s_k$.

Of course, by Corollary~\ref{posderiv}, we have $R(q,t)<4Q_k'$ for
all $t\le 0$ and all $q\in B(q'_k,0,s_k)$. Now consider the sequence
of based, rescaled flows $$(M_k,Q'_kg(Q_k'^{-1}t),(q_k',0)).$$ In
these manifolds all balls centered at $(q_k',0)$ of radii at most
$\sqrt{Q_k}s_k$ are
 $(\nu/3^n)$ non-collapsed. Also, the curvatures of these manifolds
 are non-negative and the scalar curvature is bounded by $4$.
 It follows that by passing to a subsequence we can extract a
geometric limit.
 Since $Q_k's_k^2\rightarrow\infty$
as $k\rightarrow\infty$ the asymptotic volume of this limit is at
least $\nu/3^n$. But this geometric limit is a
$\nu/3^n$-non-collapsed ancient solution with non-negative curvature
operator with scalar curvature bounded by $4$.  This contradicts
Theorem~\ref{asympvol}.
\end{proof}

The second corollary gives curvature bounds at all points in terms
of the distance to the center of the ball.

\begin{cor}\label{global}
Fix the dimension $n$. Given $\nu>0$,  there is a function
$K(A)<\infty$, defined for $A\in (0,\infty)$, such that if
$(M,g(t)),\ -\infty<t\le 0$, is an $n$-dimensional Ricci flow,
complete of bounded, non-negative curvature operator, $p\in M$ is a
point and $0<r<\infty$ is such that ${\rm Vol}\,B(p,0,r)\ge \nu r^n$
then for all $q\in M$ we have
$$(r+d_0(p,q))^2R(q,0)\le K(d_0(p,q)/r).$$
\end{cor}

\begin{proof}
Fix $q\in M$ and let $d=d_0(p,q)$. We have
$${\rm Vol}\,B(q,0,r+d)\ge {\rm
Vol}\,B(p,0,r)\ge \nu r^n=\frac{\nu}{(1+(d/r))^n}(r+d)^n.$$ Let
$K(A)=C(\nu/^n)$, where $C$ is the constant provided by the previous
corollary. The result is immediate from the previous corollary.
\end{proof}

\section{Compactness of the space of $3$-dimensional
$\kappa$-solutions}\index{$\kappa$-solution!compactness}

This section is devoted to proving the following result.

\begin{thm}\label{compactkappa}
Let $(M_k,g_k(t),(p_k,0))$ be a sequence of based $3$-dimensional
$\kappa$-solutions satisfying $R(p_k,0)=1$. Then there is a
subsequence converging smoothly to a based $\kappa$-solution.
\end{thm}

The main point in proving this theorem is to establish the uniform
curvature bounds given in the next lemma.

\begin{lem}\label{proposition}
For each $r<\infty$ there is a constant $C(r)<\infty$, such that the
following holds. Let $(M,g(t),(p,0))$ be a  based $3$-dimensional
$\kappa$-solution satisfying $R(p,0)=1$. Then $R(q,0)\le C(r)$ for
all $q\in B(p,0,r)$.
\end{lem}

\begin{proof}
Fix a based $3$-dimensional $\kappa$-solution $(M,g(t),(p,0))$. By
Theorem~\ref{asympscal} we have
$${\rm sup}_{q\in M}d_{0}(p,q)^2R(q,0)=\infty.$$
Let $q$ be a closest point to $p$ satisfying
$$d_0(p,q)^2R(q,0)=1.$$

We set $d=d_{0}(p,q)$, and we set $Q=R(q,0)$. Of course, $d^2Q=1$.
We carry this notation and these assumptions through the next five
claims. The goal of these claims is to show that $R(q',0)$ is
uniformly bounded for $q'$ near $(p,0)$ so that in fact the distance
$d$ from the point $q$ to $p$ is uniformly bounded from below by a
positive constant (see Claim~\ref{bdnearxk} for a more precise
statement). Once we have this the lemma will follow easily. To
establish this uniform bound requires a sequence of claims.

\begin{claim}\label{2Rbound}
There is a universal (i.e., independent of the $3$-dimensional
$\kappa$-solution) upper bound $C$ for $R(q',0)/R(q,0)$ for all
$q'\in B(q,0,2d)$.
\end{claim}

\begin{proof}
Suppose not. Then there is a sequence $(M_k,g_k(t),(p_k,0))$ of
$3$-dimensional $\kappa$-solutions with $R(p_k,0)=1$, points $q_k$
in $(M_k,g_k(0))$ closest to $p_k$ satisfying $d_k^2R(q_k,0)=1$,
where $d_k=d_0(p_k,q_k)$, and points $q_k'\in B(q_k,0,2d_k)$ with
$${\rm lim}_{k\rightarrow\infty}(2d_k)^2R(q_k',0)=\infty.$$ Then according
to Corollary~\ref{bdvolbdcurv} for every $\nu>0$ for all $k$
sufficiently large, we have \begin{equation}\label{V3}{\rm
Vol}\,B(q_k,0,2d_k)<\nu (2d_k)^3.\end{equation} Therefore, by
passing to a subsequence, we can assume that for each $\nu>0$
\begin{equation}\label{lowerbd}
{\rm Vol}\,B(q_k,0,2d_k)< \nu (2d_k)^3 \end{equation} for all $k$
sufficiently large. Let $\omega_3$ be the volume of the unit ball in
$\Ar^3$. Then for all $k$ sufficiently large, ${\rm
Vol}\,B(q_k,0,2d_k)<[\omega_3/2](2d_k)^3$. Since the sectional
curvatures of $(M_k,g_k(0))$ are non-negative, by the Bishop-Gromov
inequality (Theorem~\ref{BishopGromov}),  it follows that for every
$k$ sufficiently large there is $r_k<2d_k$ such that
\begin{equation}\label{V/2}
{\rm Vol}\,B(q_k,0,r_k)=[\omega_3/2]r_k^3.\end{equation} Of course,
because of Equation~(\ref{lowerbd}) we see that ${\rm
lim}_{k\rightarrow\infty}r_k/d_k=0$. Then, according to
Corollary~\ref{global}, we have for all $q\in M_k$

$$(r_k+d_{g_k(0)}(q_k,q))^2R(q,0)\le K(d_{g_k(0)}(q_k,q)/r_k),$$
where $K$ is as given in Corollary~\ref{global}. Form the sequence
$(M_k,g'_k(t),(q_k,0))$, where $g'_k(t)=r_k^{-2}g_k(r_k^2t)$. This
is a sequence  of based Ricci flows. For each $A<\infty$ we have
$$(1+A)^2R_{g'_k}(q,0)\le K(A)$$ for all
$q\in B_{g_k'(0)}(q_k,0,A)$. Hence, by the consequence of Hamilton's
Harnack inequality (Corollary~\ref{posderiv})
$$R_{g'_k}(q,t)\le K(A),$$
for all  $(q,t)\in B_{g_k'(0)}(q_k,0,A)\times(-\infty,0]$. Using
this and the fact that all the flows are $\kappa$-non-collapsed,
Theorem~\ref{flowlimit} implies that, after passing to a
subsequence, the sequence $(M_k,g_k'(t),(q_k,0))$ converges
geometrically to a limiting Ricci flow
$(M_\infty,g_\infty(t),(q_\infty,0))$ consisting of non-negatively
curved, complete manifolds $\kappa$-non-collapsed on all scales
(though possibly with unbounded curvature).

Furthermore, Equation~(\ref{V/2}) passes to the limit to give
\begin{equation}\label{V3/2}
{\rm Vol}\,B_{g_\infty}(q_\infty,0,1)=\omega_3/2. \end{equation}
 Since $r_k/d_k\rightarrow 0$ as $k\rightarrow\infty$ and since $R_{g_k}(q_k,0)=d_k^{-2}$, we see that
$R_{g_\infty}(q_\infty,0)=0$. By the strong maximum principle for
scalar curvature (Theorem~\ref{firststrongmax}), this implies that
the limit $(M_\infty,g_\infty(0))$ is flat. But
Equation~(\ref{V3/2}) tells us that this limit is not $\Ar^3$. Since
it is complete and flat, it must be a quotient of $\Ar^3$ by an
action of a non-trivial group of isometries acting freely and
properly discontinuously. But the quotient of $\Ar^3$ by any
non-trivial group of isometries acting freely and properly
discontinuously has zero asymptotic volume. [Proof: It suffices to
prove the claim in the special case when the group is infinite
cyclic. The generator of this group has an axis $\alpha$ on which it
acts by translation and on the orthogonal subspace its acts by an
isometry. Consider the circle in the quotient that is the image of
$\alpha$, and let $L_\alpha$ be its length. The volume of the ball
of radius $r$ about $L_\alpha$ is $\pi r^{2}L_\alpha$. Clearly then,
for any $p\in \alpha$, the volume of the ball of radius $r$ about
$p$ is at most $\pi L_\alpha r^2$. This proves that the asymptotic
volume of the quotient is zero.]

We have now shown that $(M_\infty,g_\infty(0))$ has zero curvature and zero
asymptotic volume. But this implies that it is not $\kappa$-non-collapsed on
all scales, which is a contradiction. This contradiction completes the proof of
Claim~\ref{2Rbound}.
\end{proof}

This claim establishes the existence of a universal constant $C<\infty$
(universal in the sense that it is independent of the $3$-dimensional
$\kappa$-solution) such that $R(q',0)\le CQ$ for all $q'\in B(q,0,2d)$. Since
the curvature of $(M,g(t))$ is non-negative and bounded, we know from the
Harnack inequality (Corollary~\ref{posderiv}) that  $R(q',t)\le CQ$ for all
$q'\in B(q,0,2d)$ and all $t\le 0$. Hence, the Ricci curvature ${\rm
Ric}(q',t)\le CQ$ for all $q'\in B(q,0,2d)$ and all $t\le 0$.

\begin{claim}
Given any constant $c>0$ there is a constant $\widetilde
C=\widetilde C(c)$, depending only on $c$ and not on the
$3$-dimensional $\kappa$-solution, so that
$$d_{g(-cQ^{-1})}(p,q)\le \tilde CQ^{-1/2}.$$
\end{claim}

\begin{proof}
Let $\gamma\colon[0,d]\to M$ be a $g(0)$-geodesic from $p$ to $q$,
parameterized at unit speed. Denote by $\ell_t(\gamma)$ the length
of $\gamma$ under the metric $g(t)$. We have $d_t(p,q)\le
\ell_t(\gamma)$. We estimate $\ell_t(\gamma)$ using the fact that
$|{\rm Ric}|\le CQ$ on the image of $\gamma$ at all times.
\begin{eqnarray*}
\frac{d}{dt}\ell_{t_0}(\gamma) & = &
\frac{d}{dt}\left(\int_0^d\sqrt{\langle\gamma'(s),\gamma'(s)\rangle_{g(t)}}ds\right)
\Bigl|_{t=t_0}\Bigr.
\\
& = & \int_0^{d}\frac{-{\rm Ric}_{g(t_0)}(\gamma'(s),\gamma'(s))}
{\sqrt{\langle\gamma'(s),\gamma'(s)\rangle_{g(t_0)}}}ds \\
& \ge & -CQ\int_0^d|\gamma'(s)|_{g(t_0)}ds=-CQ\ell_{t_0}(\gamma).
\end{eqnarray*}
Integrating yields
$$\ell_{-t}(\gamma)\le e^{CQt}\ell_0(\gamma)=e^{CQt}Q^{-1/2}.$$
(Recall $d^2Q=1$.) Plugging in $t=cQ^{-1}$ gives us
$$d_{-cQ^{-1}}(p,q)\le \ell_{-cQ^{-1}}(\gamma)\le e^{cC}Q^{-1/2}.$$
Setting $\widetilde C=e^{cC}$ completes the proof of the claim.
\end{proof}

The integrated form of  Hamilton's Harnack inequality,
Theorem~\ref{Harnack2}, tells us that
$${\rm log}\left(\frac{R(p,0)}{R(q,-cQ^{-1})}\right)\ge
-\frac{d_{-cQ^{-1}}^2(p,q)}{2cQ^{-1}}.$$ According to the above
claim, this in turn tells us
$${\rm log}\left(\frac{R(p,0)}{R(q,-cQ^{-1})}\right)\ge -\widetilde
C^2/2c.$$ Since $R(p,0)=1$, it immediately follows that
$R(q,-cQ^{-1})\le {\rm exp}(\widetilde C^2/(2c))$.

\begin{claim}\label{zkbound}
There is a universal (i.e., independent of the $3$-dimensional
$\kappa$-solution) upper bound  for $Q=R(q,0)$.
\end{claim}

\begin{proof}
Let $G'=QG$ and $t'=Qt$. Then $R_{G'}(q',0)\le C$ for all $q'\in
B_{G'}(q,0,2)$. Consequently, $R_{G'}(q',t')\le C$ for all $q'\in
B_{G'}(q,0,2)$ and all $t'\le 0$.  Thus, by Shi's derivative
estimates (Theorem~\ref{shi})applied  with $T=2$ and $r=2$, there is
a universal constant $C_1$ such that for all $-1\le t'\le 0$
$$|\triangle {\rm Rm}_{G'}(q,t')|_{G'}\le C_1,$$
(where the Laplacian is taken with respect to the metric $G'$).
Rescaling by $Q^{-1}$ we see that for all $-Q^{-1}\le t\le 0$ we
have
$$|\triangle {\rm Rm}_{G}(q,t)|\le C_1Q^2,$$
where the Laplacian is taken with respect to the metric $G$. Since
the metric is non-negatively curved, by Corollary~\ref{posderiv} we
have  $2|{\rm Ric}(q,t)|^2\le 2Q^2$ for all $t\le 0$.
 From these two facts we conclude from
the flow equation~(\ref{Revol}) that there is a constant $1<C''<\infty$ with
the property that $\partial R/\partial t(q,t)\le C''Q^2$ for all $-Q^{-1}<t\le
0$.
 Thus for any $0<c<1$, we have $Q=R(q,0)\le
cC''Q+R(q,-cQ^{-1})\le cC''Q+e^{(\widetilde C^2(c)/2c)}$. Now we
take $c=(2C'')^{-1}$ and $\widetilde C=\widetilde C(c)$. Plugging
these values into the previous inequality yields
$$Q\le 2e^{(\tilde C^2C'')}.$$
\end{proof}

This leads immediately to:

\begin{claim}\label{bdnearxk}
There are universal constants $\delta>0$ and  $C_1<\infty$
(independent of the based $3$-dimensional $\kappa$-solution
$(M,g(t),(p,0))$ with $R(p,0)=1$) such that $d(p,q)\ge \delta$. In
addition,  $R(q',t)\le C_1$ for all $q'\in B(p,0,d)$ and all $t\le
0$.
\end{claim}

\begin{proof}
Since, according to the previous claim, $Q$ is universally bounded
above and $d^2Q=1$, the existence of $\delta>0$ as required is
clear. Since $B(p,0,d)\subset B(q,0,2d)$, since $R(q',0)/R(q,0)$ is
universally bounded on $B(q,0,2d)$ by Claim~\ref{2Rbound}, and since
$R(q,0)$ is universally bounded by Lemma~\ref{zkbound}, the second
statement is clear for all $(q',0)\in B(p,0,d)\subset B(q,0,2d)$.
Given this, the fact that the second statement holds for all
$(q',t)\in B(p,0,d)\times(-\infty,0]$ then follows immediately from
the derivative inequality for $\partial R(q,t)/\partial t$,
Corollary~\ref{posderiv}.
\end{proof}

This, in turn, leads immediately to:

\begin{cor}\label{rho} Fix $\delta>0$ the universal constant of the last
claim. Then $R(q',t)\le \delta^{-2}$ for all $q'\in B(p,0,\delta)$
and all $t\le 0$.
\end{cor}

Now we return to the proof of Lemma~\ref{proposition}. Since
$(M,g(t))$ is $\kappa$-non-collapsed, it follows from the previous
corollary that ${\rm Vol}\,B(p,0,\delta)\ge \kappa \delta^3$. Hence,
according to Corollary~\ref{global} for each $A<\infty$ there is a
constant $K(A)$ such that $R(q',0)\le K(A/\delta)/(\delta+A)^2$ for
all $q'\in B(p,0,A)$. Since $\delta$ is a universal positive
constant, this completes the proof of Lemma~\ref{proposition}.
\end{proof}

Now let us turn to the proof of Theorem~\ref{compactkappa}, the
compactness result for $\kappa$-solutions.

\begin{proof}
 Let
$(M_k,g_k(t),(p_k,0))$ be a sequence of based $3$-dimensional
$\kappa$-solutions with $R(p_k,0)=1$ for all $k$.
 The immediate
consequence of Lemma~\ref{proposition} and Corollary~\ref{posderiv}
is the following. For every $r<\infty$ there is a constant
$C(r)<\infty$ such that $R(q,t)\le C(r)$ for all $q\in B(p_k,0,r)$
and for all $t\le 0$. Of course, since, in addition, the elements in
the sequence are $\kappa$-non-collapsed, by Theorem~\ref{flowlimit}
this implies that there is a subsequence of the
$(M_k,g_k(t),(p_k,0))$ that converges geometrically to an ancient
flow $(M_\infty,g_\infty(t),(p_\infty,0))$. Being a geometric limit
of $\kappa$-solutions, this limit is complete and
$\kappa$-non-collapsed, and each time-slice is of non-negative
curvature. Also, it is not flat since, by construction,
$R(p_\infty,0)=1$. Of course, it also follows from the limiting
procedure that $\partial R(q,t)/\partial t\ge 0$ for every $(q,t)\in
M_\infty\times(-\infty,0]$. Thus, according to
Corollary~\ref{3DGSSkappa} the limit
 $(M_\infty,g_\infty(t))$ has bounded curvature for each $t\le 0$.
Hence, the limit is a $\kappa$-solution. This completes the proof of
Theorem~\ref{compactkappa}.
\end{proof}

\begin{cor}\label{derivcor}
Given $\kappa>0$, there is $C<\infty$ such that for any
$3$-dimensional $\kappa$-solution $(M,g(t)),\ -\infty<t\le 0$, we
have
\begin{eqnarray}\label{gradest} {\rm sup}_{(x,t)}\frac{|\nabla R(x,t)|}{R(x,t)^{3/2}}< C\\
\label{timeest} {\rm
sup}_{(x,t)}\frac{\left|\frac{d}{dt}R(x,t)\right|}{R(x,t)^2}<C.
\end{eqnarray}

\end{cor}

\begin{proof}
Notice that the two inequalities are scale invariant. Thus, this
result is immediate from the compactness theorem,
Theorem~\ref{compactkappa}.
\end{proof}

Because of Proposition~\ref{kappa0}, and the fact that the previous
corollary obviously holds for any shrinking family of round metrics,
we can take the constant $C$ in the above corollary to be
independent of $\kappa>0$.

Notice that, using Equation~(\ref{Revol}), we can rewrite the second
inequality in the above corollary as
\begin{eqnarray*}
{\rm sup}_{(x,t)}\frac{|\triangle R+2|{\rm Ric}|^2|}{R(x,t)^2}<C.
\end{eqnarray*}

\section{Qualitative description of
$\kappa$-solutions}\index{$\kappa$-solution!qualitative description}

In Chapter~\ref{nonnegcurv} we defined the notion of an
$\epsilon$-neck. In this section we define a stronger version of
these, called strong $\epsilon$-necks. We also introduce other types
of canonical neighborhoods -- $\epsilon$-caps, $\epsilon$-round
components and $C$-components. These definitions pave the way for a
qualitative description of $\kappa$-solutions.

\subsection{Strong canonical neighborhoods}

The next manifold we introduce is one with controlled topology
(diffeomorphic either to the $3$-disk or a punctured $\Ar P^3$) with
the property that the complement of a compact submanifold is an
$\epsilon$-neck.

\begin{defn}\label{epscap}
Fix constants $0<\epsilon<1/2$ and $C<\infty$. Let $(M,g)$ be a
Riemannian $3$-manifold. A $(C,\epsilon)$-cap in $(M,g)$ is an open
submanifold $({\mathcal C},g|_{\mathcal C})$ together with an open
submanifold $N\subset {\mathcal C}$ with the following properties:
\begin{enumerate}
\item[(1)] ${\mathcal C}$ is diffeomorphic either to an open $3$-ball or to a
punctured $\Ar P^3$.
\item[(2)] $N$ is an $\epsilon$-neck with compact complement in
${\mathcal C}$.
\item[(3)] $\overline Y={\mathcal C}\setminus N$ is a compact submanifold with boundary.
Its  interior, $Y$, is called
 the {\em core} of ${\mathcal C}$. The frontier of $Y$, which is $\partial \overline Y$,
  is a central
 $2$-sphere of an $\epsilon$-neck contained in ${\mathcal C}$.
\item[(4)] The scalar curvature $R(y)>0$ for every $y\in {\mathcal C}$ and
$$ {\rm diam}({\mathcal C},g|_{\mathcal C})< C \left({\rm sup}_{y\in
{\mathcal C}}R(y)\right)^{-1/2}.$$
\item[(5)] ${\rm sup}_{x,y\in {\mathcal C}}\left[ R(y)/R(x)\right]< C$.
\item[(6)]  ${\rm Vol}\,{\mathcal C} < C({\rm sup}_{y\in {\mathcal C}}R(y))^{-3/2}$.
\item[(7)] For any $y\in Y$  let $r_y$ be defined so
that ${\rm sup}_{y'\in B(y,r_y)}R(y')=r_y^{-2}$. Then for each $y\in
Y$, the ball $B(y,r_y)$ lies in ${\mathcal C}$ and indeed has
compact closure in ${\mathcal C}$. Furthermore,
$$C^{-1}<{\rm inf}_{y\in
Y}\frac{{\rm Vol}\,B(y,r_y)}{r_y^3}.$$
\item[(8)] Lastly,
$${\rm sup}_{y\in {\mathcal C}}\frac{|\nabla R(y)|}{R(y)^{3/2}}<C$$
and
$${\rm sup}_{y\in {\mathcal C}}\frac{\left|\triangle R(y)+2|{\rm Ric}|^2\right|}{R(y)^2}<C.$$
\end{enumerate}
\end{defn}

\begin{rem}
If the ball $B(y,r_y)$ meets the complement of the core of
${\mathcal C}$ then it contains a point whose scalar curvature is
close to $R(x)$, and hence $r_y$ is bounded above by, say
$2R(x)^{-1}$. Since $\epsilon<1/2$, using the fact that $y$ is
contained in the core of ${\mathcal C}$ it follows that $B(y,r_y)$
is contained in ${\mathcal C}$ and has compact closure in ${\mathcal
C}$.

Implicitly, we always orient the $\epsilon$-neck structure on $N$ so
that the closure of its negative end meets the core of ${\mathcal
C}$. See {\sc Fig.}~\ref{fig:epsneck} in the Introduction.
\end{rem}

Condition (8) in the above definition may seem unnatural, but here
is the reason for it.

\begin{claim}\label{partialtR}
Suppose that $(M,g(t))$ is a Ricci flow and that $\left({\mathcal
C},g(t)|_{\mathcal C}\right)$ is a subset of a $t$ time-slice. Then
Condition (8) above is equivalent to
$${\rm sup}_{(x,t)\in {\mathcal C}}\frac{\left|\frac{\partial R(x,t)}{\partial
t}\right|}{R^2(x,t)}<C.$$
\end{claim}

\begin{proof}
This is immediate from Equation~(\ref{Revol}).
\end{proof}

Notice that the definition of a $(C,\epsilon)$-cap is a scale
invariant notion.

\begin{defn}\label{Ccomponent}
Fix a positive constant $C$. A compact connected Riemannian manifold
$(M,g)$ is called {\em a $C$-component} if
\begin{enumerate}
\item[(1)] $M$ is diffeomorphic to either $S^3$ or $\Ar P^3$.
\item[(2)] $(M,g)$ has positive sectional curvature.
\item[(3)] $$C^{-1}<\frac{{\rm inf}_{P}K(P)}{{\rm sup}_{y\in M}R(y)}$$
where $P$ varies over all $2$-planes in $TX$ (and $K(P)$ denotes the
sectional curvature in the $P$-direction).
\item[(4)] $$C^{-1}{\rm sup}_{y\in M}\left(R(y)^{-1/2}\right)<{\rm diam}(M)<C{\rm
inf}_{y\in M}\left(R(y)^{-1/2}\right).$$
\end{enumerate}
\end{defn}

\begin{defn}\label{defnepsround}
Fix $\epsilon>0$. Let $(M,g)$ be a compact, connected $3$-manifold.
Then $(M,g)$ is {\em within $\epsilon$ of round in the
$C^{[1/\epsilon]}$-topology} if there exist a constant $R>0$, a
compact manifold $(Z,g_0)$ of constant curvature $+1$, and a
diffeomorphism $\varphi\colon Z\to M$ with the property that the
pull back under $\varphi$ of $Rg$ is within $\epsilon$ in the
$C^{[1/\epsilon]}$-topology of $g_0$.
\end{defn}

Notice that both of these notions are scale invariant notions.

\begin{defn}
Fix $C<\infty$ and $\epsilon>0$. For any Riemannian manifold
$(M,g)$,  an open neighborhood $U$ of a point $x\in M$ is a {\em
$(C,\epsilon)$-canonical neighborhood}\index{canonical
neighborhood|ii} if one of the following holds:
\begin{enumerate}
\item[(1)] $U$ is an $\epsilon$-neck\index{$\epsilon$-neck|ii} in $(M,g)$ centered at $x$.
\item[(2)] $U$ a $(C,\epsilon)$-cap\index{$(C,\epsilon)$-cap|ii} in $(M,g)$ whose core contains $x$.
\item[(3)] $U$ is a $C$-component\index{$C$-component|ii} of $(M,g)$.
\item[(4)] $U$ is  an $\epsilon$-round component\index{$\epsilon$-round component|ii} of $(M,g)$.
\end{enumerate}
\end{defn}

Whether or not a point $x\in M$ has a $(C,\epsilon)$-canonical
neighborhood in $M$ is a scale invariant notion.

The notion of $(C,\epsilon)$-canonical neighborhoods is sufficient
for some purposes, but often we need a stronger notion.

\begin{defn}\label{defncanonnbhd}
Fix constants $C<\infty$ and $\epsilon>0$. Let $({\mathcal M},G)$ be
a generalized Ricci flow. An {\em evolving  $\epsilon$-neck defined
for an interval of normalized time of length
$t'>0$}\index{$\epsilon$-neck!evolving|ii} centered at a point $x\in
{\mathcal M}$ with ${\bf t}(x)=t$ is an embedding $\psi\colon
S^2\times (-\epsilon^{-1},\epsilon^{-1})\stackrel{\cong}{\to}
N\subset M_t$ with $x\in \psi(S^2\times\{0\})$ satisfying the
following properties:
\begin{enumerate}
\item[(1)] There is an embedding $N\times (t-R(x)^{-1}t',t]\to {\mathcal
M}$ compatible with time and the vector field.
\item[(2)] The pullback under $\psi$ of the one-parameter family of
metrics on $N$ determined by restricting $R(x)G$ to the image of
this embedding is within $\epsilon$ in the
$C^{[1/\epsilon]}$-topology of the standard family $(h(t),ds^2),
-t'<t\le 0$, where $h(t)$ is the round metric of scalar curvature
$1/(1-t)$ on $S^2$ and $ds^2$ is the usual Euclidean metric on the
interval (see Definition~\ref{epsclose} for the notion of two
families of metrics being close).
\end{enumerate}
A {\em strong $\epsilon$-neck}\index{$\epsilon$-neck!strong|ii} is
the image of an evolving $\epsilon$-neck which is defined for an
interval of normalized time of length  $1$.

Both of these notions are scale invariant notions.

Let $({\mathcal M},G)$ be a generalized Ricci flow. Let $x\in
{\mathcal M}$ be a point with ${\bf t}(x)=t$. We say that an open
neighborhood $U$ of $x$ in $M_t$  is a {\em strong
$(C,\epsilon)$-canonical neighborhood}\index{canonical
neighborhood!strong|ii} of $x$ if one of the following holds
\begin{enumerate}
\item[(1)]   $U$ is a strong $\epsilon$-neck in $({\mathcal M},G)$ centered at $x$.
\item[(2)] $U$ is a $(C,\epsilon)$-cap in $M_t$ whose  core contains $x$.
\item[(3)] $U$ is a $C$-component of $M_t$.
\item[(4)] $U$ is an $\epsilon$-round component of $M_t$.
\end{enumerate}
\end{defn}

Whether or not a point $x$ in a generalized Ricci flow has a strong
$(C,\epsilon)$-canonical neighborhood is a scale invariant notion.

\begin{prop}\label{canonvary}
The following holds for any  $\epsilon<1/4$ and any $C<\infty$. Let
$({\mathcal M},G)$ be a generalized Ricci flow and let $x\in
{\mathcal M}$ be a point with ${\bf t}(x)=t$.

(1) Suppose that $U\subset M_t$ is a  $(C,\epsilon)$-canonical
neighborhood for $x$. Then for any horizontal metric $G'$
sufficiently close to $G|_U$ in the $C^{[1/\epsilon]}$-topology,
$(U,G'|_U)$ is a  $(C,\epsilon)$ neighborhood for any $x'\in U$
sufficiently close to $x$.

(2) Suppose that in $({\mathcal M},G)$ there is an evolving
$\epsilon$-neck $U$ centered at $(x,t)$ defined for an interval of
normalized time of length  $a>1$. Then any Ricci flow on $U\times
(t-aR(x,t)^{-1},t]$ sufficiently close in the
$C^{[1/\epsilon]}$-topology to the pullback of $G$ contains a strong
$\epsilon$-neck centered at $(x,t)$.

(3) Given $(C,\epsilon)$ and $(C',\epsilon')$ with $C'>C$ and
$\epsilon'>\epsilon$ there is $\delta>0$ such that the following
holds. Suppose that $R(x)\le 2$. If $(U,g)$ is a
$(C,\epsilon)$-canonical neighborhood of $x$ then for any metric
$g'$ within $\delta$ of $g$ in the $C^{[1/\epsilon]}$-topology
$(U,g')$ contains a $(C',\epsilon')$-neighborhood of $x$.

(4) Suppose that $g(t),\ -1<t\le 0$, is a one-parameter family of
metrics on $(U,g)$ that is a strong $\epsilon$-neck centered at
$(x,0)$ and $R_g(x,0)=1$. Then any one-parameter family $g'(t)$
within $\delta$ in the $C^{[1/\epsilon]}$-topology of $g$ with
$R_{g'}(x,0)=1$ is a strong $\epsilon'$-neck.
\end{prop}

\begin{proof}
Since $\epsilon<1/4$,  the diameter of $(U,g)$, the volume of
$(U,g)$, the supremum over $x\in U$ of $R(x)$, the supremum over $x$
and $y$ in $U$ of $R(y)/R(x)$, and the infimum over all $2$-planes
$P$ in ${\mathcal H}TU$ of $K(P)$ are all continuous functions of
the horizontal metric $G$ in the $C^{[1/\epsilon]}$-topology.

Let us consider the first statement. Suppose $(U,G|_U)$ is a
$C$-component or an $\epsilon$-round component. Since the defining
inequalities are strict, and, as we just remarked, the quantities in
these inequalities vary continuously with the metric in the
$C^{[1/\epsilon]}$-topology the result is clear in this case.

Let us consider the case when $U\subset M_t$ is an $\epsilon$-neck
centered at $x$. Let $\psi\colon S^2\times
(-\epsilon^{-1},\epsilon^{-1})\to U$ be the map giving the
$\epsilon$-neck structure. Then for all horizontal metrics $G'$
sufficiently close to $G$ in the $C^{[1/\epsilon]}$-topology the
same map $\psi$ is determines an $\epsilon$-neck centered at $x$ for
the structure $(U,G'|_U)$. Now let us consider moving $x$ to a
nearby point $x'$, say $x'$ is the image of $(a,s)\in S^2\times
(-\epsilon^{-1},\epsilon^{-1})$. We pre-compose $\psi$ by a map
which is the product of the identity in the $S^2$-factor with a
diffeomorphism $\alpha$ on $(-\epsilon^{-1},\epsilon^{-1})$ that is
the identity near the ends and moves $0$ to $s$. As $x'$ approaches
$x$, $s$ tends to zero, and hence we can choose $\alpha$ so that it
tends to the identity in the $C^\infty$-topology. Thus, for $x'$
sufficiently close to $x$, this composition will determine an
$\epsilon$-neck structure centered at $x'$. Lastly, let us consider
the case when $(U,G|_U)$ is a $(C,\epsilon)$-cap whose core $Y$
contains $x$. Let $G'$ be a horizontal metric sufficiently close to
$G|_U$ in the $C^{[1/\epsilon]}$-topology. Let $N\subset U$ be the
$\epsilon$-neck $U\setminus \bar Y$. We have just seen that
$(N,G'|_N)$ is an $\epsilon$-neck. Similarly, if $N'\subset U$ is an
$\epsilon$-neck with central $2$-sphere $\partial \overline Y$, then
$(N',G'|_{N'})$ is an $\epsilon$-neck if $G'$ is sufficiently close
to $G$ in the $C^{[1/\epsilon]}$-topology.

Thus, Conditions (1),(2), and (3) in the definition of a
$(C,\epsilon)$-cap hold for $(U,G'_U)$.
 Since the
curvature, volume and diameter inequalities in Conditions (4), (5),
and (6) are strict, they also hold for $g'$. To verify that
Condition (7) holds for $G'$, we need only remark that $r_y$ is a
continuous function of the metric. Lastly, since the derivative
inequalities for the curvature in Condition (8) are strict
inequalities and $\epsilon^{-1}>4$, if these inequalities hold for
all horizontal metrics $G'$ sufficiently close to $G$ in the
$C^{[1/\epsilon]}$-topology. This completes the examination of all
cases and proves the first statement.

The second statement is proved in the same way using the fact that
if $g'(t)$ is sufficiently close to $g$ in the
$C^{[1/\epsilon]}$-topology and if $x'$ is sufficiently close to $x$
then $R_{g'(x')}^{-1}<aR_g(x)^{-1}$.

Now let us turn to the third statement. The result is clear for
$\epsilon$-necks. Also, since $R(x)\le 2$ the result is clear for
$\epsilon$-round components and $C$-components as well. Lastly, we
consider a $(C,\epsilon)$-cap $U$ whose core $Y$ contains $x$.
 Clearly,
since $R(x)$ is bounded above by  $2$, for $\delta>0$ sufficiently
small, any metric $g'$ within $\delta$ of $g$ will satisfy the
diameter, volume and curvature and the derivative of the curvature
inequalities with $C'$ replacing $C$. Let $N$ be the $\epsilon$-neck
in $(U,g)$ containing the end of $U$. Assuming that $\delta$ is
sufficiently small, let $N'$ be the image of $S^2\times
\left(-\epsilon^{-1},2(\epsilon')^{-1}-\epsilon^{-1}\right)$. Then
$(N',g')$ becomes an $\epsilon'$-neck structure once we shift the
parameter in the $s$-direction by $\epsilon^{-1}-(\epsilon')^{-1}$.
We let $U'=\overline Y\cup N'$. Clearly, the $\epsilon$-neck with
central $2$-sphere $\partial \overline Y$ will also determine an
$\epsilon'$-neck with the same central $2$-sphere provided that
$\delta>0$ is sufficiently small. Thus, for $\delta>0$ sufficiently
small, for any $(C,\epsilon)$ the result of this operation is a
$(C',\epsilon')$-cap with the same core.

The fourth statement is immediate.
\end{proof}

\begin{cor}\label{strepsopen}
In an ancient solution $(M,g(t))$ the set of points that are centers
of strong $\epsilon$-necks is an open subset
\end{cor}

\begin{proof}
Let $T$ be the final time of the flow. Suppose that $(x,t)$ is the
center of a strong $\epsilon$-neck $U\times (t-R(x,t)^{-1},t]\subset
M\times (-\infty,0]$. This neck extends backwards for all time and
forwards until the final time $T$ giving an embedding $U\times
(-\infty,T]\to M\times (-\infty,T]$. There is $a>1$ such that for
all $t'$ sufficiently close to $t$ the restriction of this embedding
determines an evolving $\epsilon$-neck centered at $(x,t')$ defined
for an interval of normalized time of length  $a$. Composing this
neck structure with a self-diffeomorphism of $U$ moving $x'$ to $x$,
as described above, shows that all $(x',t')$ sufficiently close to
$(x,t)$ are centers of strong $\epsilon$-necks.
\end{proof}

\begin{defn}\label{epstube}
{\em An $\epsilon$-tube ${\mathcal T}$}\index{$\epsilon$-tube|ii} in
a Riemannian $3$-manifold $M$ is a submanifold diffeomorphic to the
product of $S^2$ with a non-degenerate interval with the following
properties:
\begin{enumerate}
\item[(1)] Each boundary component $S$ of ${\mathcal T}$ is the central
$2$-sphere of an $\epsilon$-neck $N(S)$ in $M$.
\item[(2)] ${\mathcal T}$
is a union of $\epsilon$-necks and the closed half $\epsilon$-necks
whose boundary sphere is a component of $\partial{\mathcal T}$.
Furthermore, the central $2$-sphere of each of the $\epsilon$-necks
is isotopic in ${\mathcal T}$ to the $S^2$-factors of the product
structure.
\end{enumerate}

An {\em open $\epsilon$-tube} is one without boundary. It is a union
of $\epsilon$-necks with the central spheres that are isotopic to
the $2$-spheres of the product structure.

A {\em $C$-capped $\epsilon$-tube}\index{$\epsilon$-tube!capped|ii}
in $M$ is a connected submanifold that is the union of a
$(C,\epsilon)$-cap ${\mathcal C}$
 and an open $\epsilon$-tube where the intersection of ${\mathcal C}$ with the
$\epsilon$-tube   is diffeomorphic to $S^2\times (0,1)$ and contains
an end of the $\epsilon$-tube and an end of the cap.
 A {\em doubly $C$-capped
$\epsilon$-tube}\index{$\epsilon$-tube!doubly capped|ii} in $M$ is a
closed, connected submanifold of $M$ that is the union of two
$(C,\epsilon)$-caps ${\mathcal C}_1$ and ${\mathcal C}_2$  and  an
open $\epsilon$-tube. Furthermore, we require (i) that the cores
$Y_1$ and $Y_2$ of ${\mathcal C}_1$ and ${\mathcal C}_2$ have
disjoint closures, (ii) that the union of either ${\mathcal C}_i$
with the $\epsilon$-tube is  a capped $\epsilon$-tube and ${\mathcal
C}_1$ and ${\mathcal C}_2$ contain the opposite ends of the
$\epsilon$-tube. There is one further closely related notion, that
of an {\em $\epsilon$-fibration}\index{$\epsilon$-fibration|ii}. By
definition an $\epsilon$-fibration is a closed, connected manifold
that fibers over the circle with fibers $S^2$ that is also a union
of $\epsilon$-necks with the property that the central $2$-sphere of
each neck is isotopic to a fiber of the fibration structure. We
shall not see this notion again until the appendix, but because it
is clearly closely related to the notion of an $\epsilon$-tube, we
introduce it here.
\end{defn}

\begin{figure}[ht]
  \centerline{\epsfbox{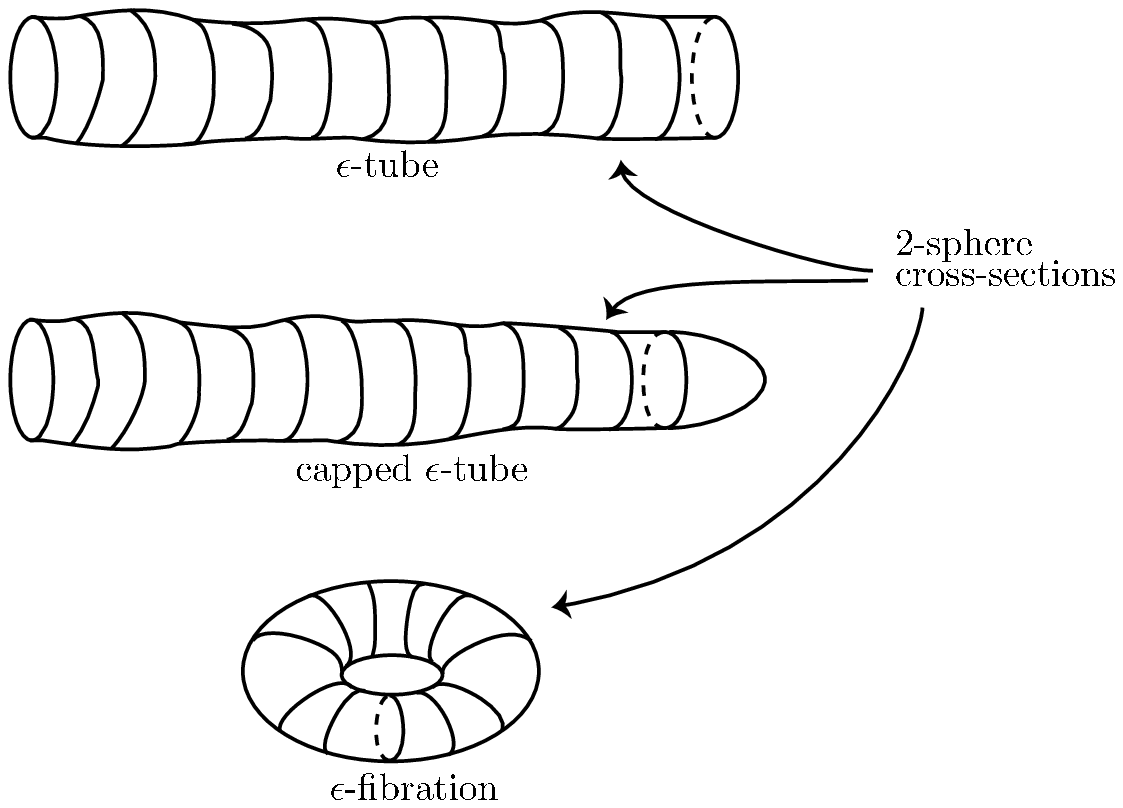}}
  \caption{$\epsilon$-canonical neighborhoods}\label{fig:canonnbhd}
\end{figure}

See {\sc Fig.}~\ref{fig:canonnbhd}.

\begin{defn}
A {\em strong $\epsilon$-tube}\index{$\epsilon$-tube!strong|ii} in a
generalized Ricci flow is an $\epsilon$-tube with the property that
each point of the tube is the center of a strong $\epsilon$-neck in
the generalized flow.
\end{defn}

\subsection{Canonical neighborhoods for $\kappa$-solutions}

\begin{prop}\label{infprod2}
Let $(M,g(t))$ be a  $3$-dimensional $\kappa$-solution.  Then one of
the following hold:
\begin{enumerate}
\item[(1)] For every $t\le 0$ the manifold $(M,g(t))$ has positive
curvature.
\item[(2)] $(M,g(t))$ is the product of an evolving family of
round $S^2$'s with a line.
\item[(3)] $M$ is diffeomorphic to a line  bundle over $\Ar P^2$, and
there is a finite covering of $(M,g(t))$ that is  a flow as in (2).
\end{enumerate}
\end{prop}

\begin{proof}
Suppose that $(M,g(t))$ does not have positive curvature for some
$t$. Then, by the application of the strong maximum principle given
in Corollary~\ref{nullspace}, there is a covering $\widetilde M$ of
$M$, with either one or $2$-sheets, such that $(\widetilde M,\tilde
g(t))$ is the product of an evolving family of round surfaces with a
flat one-manifold (either a circle or the real line).  Of course,
the covering must be $\kappa$-solution. In the case in which
$(\widetilde M,\tilde g(t))$ is isometric to the product of an
evolving family of round surfaces  and a circle, that circle has a
fixed length, say $L<\infty$. Since the curvature of the surface in
the $t$ time-slice goes to zero as $t\rightarrow -\infty$, we see
that the flow is not $\kappa$-non-collapsed on all scales for any
$\kappa>0$. Thus, $(M,g(t))$ has either a trivial cover or a double
cover isometric to the product of a shrinking family of round
surfaces with $\Ar$. If the round surface is $S^2$, then we have
established the result. If the round surface is $\Ar P^2$ a further
double covering is a product of round two-spheres with $\Ar$.
 This proves the proposition.
\end{proof}

\begin{lem}
Let $(M,g(t))$ be a non-compact $3$-dimensional $\kappa$-solution of
positive curvature and let $p\in M$. Then there is $D'<\infty$,
possibly depending on $(M,g(0))$ and $p$, such that every point of
$M\times\{0\}\setminus B(p,0,D'R(p,0)^{-1/2})$ is the center of an
evolving $\epsilon$-neck in $(M,g(t))$ defined for an interval of
normalized time of length  $2$. Furthermore, there is $D'_1<\infty$
such that for any point $x\in B(p,0,D'R(p,0)^{-1/2})$ and any
$2$-plane $P_x$ in $T_xM$
 we have $(D'_1)^{-1}< K(P_x)/R(p,0)< D'_1$ where $K(P_x)$ denotes the sectional
 curvature in the direction of the $2$-plane $P_x$.
\end{lem}

\begin{proof}
 Given $(M,g(t))$ and $p$, suppose that no such $D'<\infty$ exists.
Because the statement is scale invariant, we can arrange that
$R(p,0)=1$.
 Then we
can find a sequence of points $p_k\in M$ with $d_0(p,p_k)\rightarrow
\infty$ as $k\rightarrow \infty$ such that no $p_k$ is the center of
an evolving $\epsilon$-neck in $(M,g(0))$ defined for an interval of
normalized time of length  $2$. By passing to a subsequence we can
assume that one of two possibilities holds: either
$d_0^2(p,p_k)R(p_k,0)\rightarrow\infty$ as $k\rightarrow\infty$ or
${\rm lim}_{k\rightarrow\infty}d_0^2(p,p_k)R(p_k,0)=\ell<\infty$. In
the first case,  set $\lambda_k=R(p_k,0)$ and consider the based
flows $(M,\lambda_kg(\lambda_k^{-1}t),(p_k,0))$. According to
Theorem~\ref{compactkappa}, after passing to a subsequence there is
a geometric limit. Thus, by Theorem~\ref{topsplit} and
Corollary~\ref{localprod} the limit splits as a product of a
$2$-dimensional $\kappa$-solution and $\Ar$. By
Corollary~\ref{2DGSS} it follows that the limit is the standard
round evolving cylinder. This implies that for all $k$ sufficiently
large $(p_k,0)$ is the center of an evolving $\epsilon$-neck in
$(M,g(t))$ defined for an interval of normalized time of length
$2$. This contradiction establishes the existence of $D'$ as
required in this case.

Now suppose that ${\rm
lim}_{k\rightarrow\infty}d_0^2(p,p_k)R(p_k,0)=\ell<\infty$. Of
course, since $d_0(p,p_k)\rightarrow\infty$, it must be the case
that $R(p_k,0)\rightarrow 0$ as $k\rightarrow \infty$. Set
$Q_k=R(p_k,0)$. By passing to a subsequence we can arrange that
$d_0^2(p,p_k)Q_k< \ell+1$ for all $k$. Consider the
$\kappa$-solutions $(M_k,g_k(t))=(M,Q_kg(Q_k^{-1}t))$. For each $k$
we have $p\in B_{g_k}(p_k,0,\ell+1)$, and
$R_{g_k}(p,0)=Q_k^{-1}\rightarrow\infty$ as $k\rightarrow\infty$.
This contradicts Lemma~\ref{proposition}, and completes the proof of
the existence of $D$ as required in this case as well.

 The existence of $D'_1$ is immediate since the
closure of the ball is compact and the manifold has positive
curvature.
\end{proof}

In fact a much stronger result is true. The constants $D'$ and
$D_1'$ in the above lemma can be chosen independent of the
non-compact $\kappa$-solutions.

\begin{prop}\label{uniformD} For any $0<\epsilon$ sufficiently small there are constants
 $D=D(\epsilon)<\infty$ and $D_1=D_1(\epsilon)<\infty$
 such that the following holds for  any non-compact $3$-dimensional
 $\kappa$-solution $(M,g(t))$ of positive
curvature. Let $p\in M$ be a soul of $(M,g(0))$. Then:
\begin{enumerate}
\item[(1)] Every point in $M\setminus B(p,0,DR(p,0)^{-1/2})$ is the center
of a strong $\epsilon$-neck in $(M,g(t))$. Furthermore, for any
$x\in B(p,0,DR(p,0)^{-1/2})$ and any $2$-plane $P_x$ in $T_xM$ we
have
$$D_1^{-1}<K(P_x)/R(p,0)<D_1.$$
Also,
$$D_1^{-3/2}R(p,0)^{-3/2}< {\rm
Vol}(B(p,0,DR(p,0)^{-1/2})<D_1^{3/2}R(p,0)^{-3/2}.$$
 \item[(2)] Let $f$
denote the distance function from $p$. For any $\epsilon$-neck
$N\subset (M,g(0))$, the middle two-thirds of $N$ is disjoint from
$p$, and the central $2$-sphere $S_N$ of $N$ is (topologically)
isotopic in $M\setminus \{p\}$ to $f^{-1}(a)$ for any $a>0$. In
particular, given two disjoint central $2$-spheres of
 $\epsilon$-necks in $(M,g(0))$ the region of $M$ bounded by these $2$-spheres is
 diffeomorphic
 to $S^2\times [0,1]$.
 \end{enumerate}
\end{prop}

\begin{rem}
In Part 1 of this theorem one can  replace $p$ by any point $p'\in
M$ that is not the center of a strong $\epsilon$-neck.
\end{rem}

\begin{proof}
First suppose that no $D$ exists so that the first statement holds.
Then there is a sequence of such solutions $(M_k,g_k(t))$, with
$p_k\in M_k$ being a soul of $(M_k,g_k(0))$ and points $q_k\in M_k$
with $d^2_0(p_k,q_k)R(p_k,0)\rightarrow\infty$ as
$k\rightarrow\infty$ such that $q_k$ is not the center of a strong
$\epsilon$-neck in $(M_k,g_k(0))$. By rescaling we can assume that
$R(p_k,0)=1$ for all $k$, and hence that
$d_0(p_k,q_k)\rightarrow\infty$. Then, according to
Theorem~\ref{compactkappa}, by passing to a subsequence we can
assume that there is a geometric limit
$(M_\infty,g_\infty(t),(p_\infty,0))$ with $R(p_\infty,0)=1$.
 By
Lemma~\ref{topiso}, provided that $\epsilon$ is sufficiently small
for all $k$ the soul $(p_k,0)$ is not the center of a strong
$2\epsilon$-neck in $(M_k,g_k(t))$.  Hence, invoking Part 4 of
Proposition~\ref{canonvary} and using the fact that $R(p_k,0)=1$ for
all $k$ and hence $R(p_\infty,0)=1$, we see that $(p_\infty,0)$ is
not the center of a strong $\epsilon$-neck in
$(M_\infty,g_\infty(t))$. Since the manifolds $M_k$ are non-compact
and have metrics of positive curvature they are diffeomorphic to
$\Ar^3$ and in particular, do not contain embedded copies of $\Ar
P^2$. Thus, the limit $(M_\infty,g_\infty(t))$ is a non-compact
$\kappa$-solution containing no embedded copy of $\Ar P^2$. Thus, by
Proposition~\ref{infprod2} either it is positively curved or it is a
Riemannian product $S^2$ times $\Ar$. In the second case every point
is the center of a strong $\epsilon$-neck. Since we have seen that
the point $(p_\infty,0)$ is not the center of a strong
$\epsilon$-neck, it follows that the limit is a positively curved
$\kappa$-solution.

 Then according to the previous lemma there is $D'$, depending on
$(M_\infty,g_\infty(0))$ and $p_\infty$, such that every point
outside $B(p_\infty,0,D')$ is the center of an evolving
$\epsilon/2$-neck defined for an interval of normalized time of
length $2$.

Now since $(M_k,g_k(t),(p_k,0))$ converge geometrically to
$(M_\infty,g_\infty(t),(p_\infty,0))$, by Part 2 of
Proposition~\ref{canonvary} for any $L<\infty$, for all $k$
sufficiently large, all points of $B(p_k,0,L)\setminus B(p_k,0,2D')$
are centers of strong $\epsilon$-necks in $(M_k,g_k(t))$. In
particular, for all $k$ sufficiently large, $d_0(p_k,q_k)>L$. Let
$L_k$ be a sequence tending to infinity as $k\rightarrow \infty$.
Passing to a subsequence, we can suppose that every point of
$\left(B(p_k,0,L_k)\setminus B(p_k,0,2D')\right)\subset M_k$ is the
center of a strong $\epsilon$-neck in $(M_k,g_k(0))$.  Of course,
for all $k$ sufficiently large, $q_k\in M_k\setminus B(p_k,0,2D'))$.
By Corollary~\ref{strepsopen} the subset of points in $M_k\times
\{0\}$ that are centers of strong $\epsilon$-necks is an open set.
Thus, replacing $q_k$ with another point if necessary we can suppose
that it $q_k$ is a closest point to $p_k$ contained in $M_k\setminus
B(p_k,0,2D')$ with the property $q_k$ is not the center of a strong
$\epsilon$-neck. Then $q_k\in M_k\setminus B(p_k,0,L_k)$ and
$(q_k,0)$ is in the closure of the set of points in $M_k$ that are
centers of strong $\epsilon$-necks in $(M_k,g_k(t))$, and hence by
Part 3 of Proposition~\ref{canonvary} each $(q_k,0)$ is the center
of a $2\epsilon$-neck in $(M_k,g_k(t))$.

Let $\gamma_k$ be a minimizing geodesic connecting $(p_k,0)$ to
$(q_k,0)$, and let $\mu_k$ be a minimizing geodesic ray from
$(q_k,0)$ to infinity. Set $Q_k=R(q_k,0)$. Since $(q_k,0)$ is the
center of a $2\epsilon$-neck, from Lemma~\ref{neckseparate} we see
that, provided that $\epsilon$ is sufficiently small,  the
$2\epsilon$-neck centered at $q_k$ separates $p$ from $\infty$, so
that $\gamma_k$ and $\mu_k$ exit this $2\epsilon$-neck at opposite
ends. According to Theorem~\ref{compactkappa}, after passing to a
subsequence, the based, rescaled flows
$$(M_k,Q_kg(Q_k^{-1}t),(q_k,0))$$
converge geometrically to a limit. Let $(q_\infty,0)$ be the base
point of the resulting limit. By Part 3 of
Proposition~\ref{canonvary}, it is the center of a $4\epsilon$-neck
in the limit.

\begin{claim}
$d_0^2(p_k,q_k)Q_k\rightarrow \infty$ as $k\rightarrow \infty$.
\end{claim}

\begin{proof}
Suppose not. Then by passing to a subsequence we can suppose that
these products are bounded independent of $k$. Then since
$d_0(p_k,q_k)\rightarrow\infty$, we see that $Q_k\rightarrow 0$.
Thus, in the rescaled flows $(M_k,Q_kg_k(Q_k^{-1}t))$ the curvature
at $(p_k,0)$ goes to infinity. But this is impossible since the
$Q_kg_k$-distance from $(p_k,0)$ to $(q_k,0)$ is
$\sqrt{Q_k}d_0(p_k,q_k)$ which is bounded independent of $k$ and the
scalar curvature of $(p,0)$ in the metric $Q_kg_k(0)$ is
$R(p_k,0)Q_k^{-1}=Q_k^{-1}$ tends to $\infty$. Unbounded curvature
at bounded distance contradicts Lemma~\ref{proposition}, and this
establishes the claim.
\end{proof}

A subsequence of the based flows $(M_k,Q_kg_k(Q_k^{-1}t),(q_k,0))$
converge geometrically to a $\kappa$-solution. According
Theorem~\ref{topsplit} and Corollary~\ref{localprod}, this limiting
flow is the product of a $2$-dimensional $\kappa$-solution with a
line. Since $M$ is orientable, this $2$-dimensional
$\kappa$-solution is an evolving family of round $2$-spheres. This
implies that for all $k$ sufficiently large, $(q_k,0)$ is the center
of a strong $\epsilon$-neck in $(M_k,g_k(t))$. This is a
contradiction and proves the existence of $D<\infty$ as stated  in
the proposition.

Let $(M_k,g_k(t),(p_k,0))$ be a  sequence of non-compact Ricci flows
based at a soul  $p_k$ of $(M_k,g_k(0))$. We rescale so that
$R(p_k,0)=1$. By Lemma~\ref{topiso}, if $\epsilon$ is sufficiently
small, then $p_k$ cannot be the center of an $\epsilon$-neck. It
follows from Proposition~\ref{canonvary} that for any limit of a
subsequence the point $p_\infty$, which is the limit of the $p_k$,
is not the center of an $2\epsilon$-neck in the limit. Since the
limit manifold is orientable, it is either contractible with
strictly positive curvature or is a metric product of a round
$2$-sphere and the line. It follows that the limit manifold has
strictly positive curvature at $(p_\infty,0)$, and hence positive
curvature everywhere. The existence of $D_1<\infty$ as required is
now immediate from Theorem~\ref{compactkappa}.

The fact that any soul is disjoint from the middle two-thirds of any
$\epsilon$-neck and the fact that the central $2$-spheres of all
$\epsilon$-necks are isotopic in $M\setminus\{p\}$ are contained in
Lemma~\ref{topiso} and Corollary~\ref{neckseparate}.
\end{proof}

\begin{cor}\label{noncompkappa}
There is $\bar\epsilon_2>0$ such that for any $0<\epsilon\le \bar
\epsilon_2$ the following holds.
 There
is $C_0=C_0(\epsilon)$ such that for any $\kappa>0$ and any non-compact
$3$-dimensional $\kappa$-solution not containing an embedded $\Ar P^2$ with
trivial normal bundle, the zero time-slice is either a strong $\epsilon$-tube
or a $C_0$-capped strong $\epsilon$-tube.
\end{cor}

\begin{proof}
For $\epsilon>0$ sufficiently small let $D(\epsilon)$ and
$D_1(\epsilon)$ be as in the previous corollary. At the expenses of
increasing these, we can assume that they are at least the constant
$C$ in Corollary~\ref{derivcor}. We set
$$C_0(\epsilon)={\rm max}(D(\epsilon),D_1(\epsilon)).$$
 If the non-compact $\kappa$-solution has
positive curvature, then the corollary follows immediately from
Proposition~\ref{uniformD} and Corollary~\ref{derivcor}. If the
$\kappa$-solution is the product of an evolving round $S^2$ with the
line, then every point of the zero time-slice is the center of a
strong $\epsilon$-neck for every $\epsilon>0$ so that the zero
time-slice of the solution is a strong $\epsilon$-tube. Suppose the
solution is double covered by the product of an evolving round
$2$-sphere and the line. Let $\iota$ be the involution and take the
product coordinates so that $S^2\times\{0\}$ is the invariant
$2$-sphere of $\iota$ in the zero time-slice. Then any point in the
zero time-slice at distance at least $3\epsilon^{-1}$ from
$P=(S^2\times \{0\})/\iota$ is the center of a strong
$\epsilon$-neck. Furthermore, an appropriate neighborhood of $P$ in
the time zero slice is a $(C,\epsilon)$-cap whose core contains the
$3\epsilon^{-1}$ neighborhood of $P$. The derivative bounds in this
case come from the fact that the metric is close in the
$C^{[1/\epsilon]}$-topology to the standard evolving flow. This
proves the corollary in this case and hence completes the proof.
\end{proof}

Now let us consider compact $\kappa$-solutions.

\begin{thm}\label{thmep3}
There is $\bar \epsilon_3>0$ such that for every $0<\epsilon\le \bar
\epsilon_3$  there is $C_1=C_1(\epsilon)<\infty$ such that one of
the following holds for any $\kappa>0$ and any compact
$3$-dimensional $\kappa$-solution $(M,g(t))$.
\begin{enumerate}
\item[(1)] The manifold $M$ is compact and of constant positive sectional curvature.
\item[(2)] The diameter of $(M,g(0))$ is  less than $C_1\cdot ({\rm max}_{x\in
M}R(x,0))^{-1/2}$, and $M$ is diffeomorphic to either $S^3$ or $\Ar
P^3$.
\item[(3)] $(M,g(0))$
 is  a double  $C_1$-capped  strong $\epsilon$-tube.
\end{enumerate}
\end{thm}

\begin{proof}
First notice that if $(M,g(t))$ is not of strictly positive
curvature, then the universal covering of $(M,g(0))$ is a Riemannian
product $S^2\times \Ar$, and hence $(M,g(0))$ is either non-compact
or finitely covered by the product flow on $S^2\times S^1$. The
former case is ruled out since we are assuming that $M$ is compact
and the latter is ruled out because such flows are not
$\kappa$-non-collapsed for any $\kappa>0$. We conclude that
$(M,g(t))$ is of positive curvature. This implies that the
fundamental group of $M$ is finite. If there were an embedded $\Ar
P^2$ in $M$ with trivial normal bundle, that $\Ar P^2$ cannot
separate (since the Euler characteristic of $\Ar P^2$ is one, it is
not the boundary of a compact $3$-manifold). But a non-separating
surface in $M$ induces a surjective homomorphism of $H_1(M)$ onto
$\Zee$. We conclude from this that $M$ does not contain an embedded
$\Ar P^2$ with trivial normal bundle.

 We assume that $(M,g(0))$ is not round so that by Proposition~\ref{kappa0}
 there is a universal $\kappa_0>0$ such that $(M,g(0))$ is a
$\kappa_0$-solution. Let $C_0(\epsilon)$ be the constant from
 Corollary~\ref{noncompkappa}.

 \begin{claim}
Assuming that $(M,g(0))$ is compact but not of constant positive
sectional curvature, for each $\epsilon>0$ there is $C_1$ such that
if the diameter of $(M,g(0))$ is greater than $C_1({\rm max}_{x\in
M}R(x,0))^{-1/2}$ then every point of $(M,g(0))$ is either contained
in the core of $(C_0(\epsilon),\epsilon)$-cap or is the center of a
strong $\epsilon$-neck in $(M,g(t))$.
\end{claim}

\begin{proof}
Suppose that for some $\epsilon>0$ there is no such $C_1$. Then we
take a sequence of constants $C'_k$ that diverges to $+\infty$ as
$k\rightarrow \infty$ and a sequence $(M_k,g_k(t), (p_k,0))$ of
based $\kappa_0$-solutions such that the diameter of $(M_k,0)$ is
greater than $C'_kR^{-1/2}(p_k,0)$ and yet $(p_k,0)$ is not
contained in the core of a $(C_0(\epsilon),\epsilon)$-cap nor is the
center of a strong $\epsilon$-neck.  We scale $(M_k,g_k(t))$ by
$R(p_k,0)$. This allows us to assume that $R(p_k,0)=1$ for all $k$.
According to Theorem~\ref{compactkappa}, after passing to a
subsequence we can assume these based $\kappa$-solutions converge to
a based $\kappa$-solution $(M_\infty,g_\infty(t),(p_\infty,0))$.
Since the diameters of the $(M_k,g_k(0))$ go to infinity, $M_\infty$
is non-compact. According to Corollary~\ref{noncompkappa} the point
$p_\infty$ is either the center of a strong $\epsilon$-neck, or is
contained in the core of  a $(C_0(\epsilon),\epsilon)$-cap. Since
$R(p_k,0)=1$ for all $k$, it follows from Parts 1 and 4 of
Proposition~\ref{canonvary} that for all $k$ sufficiently large,
$(p_k,0)$ is either the center of a strong $\epsilon$-neck in
$(M_k,g_k(t))$ or is contained in the core of a
$(C_0(\epsilon),\epsilon)$-cap. This is a contradiction, proving the
claim.
\end{proof}

Now it follows from Proposition~\ref{epstopology} that if the
diameter of $(M,g(0))$ is greater than $C_1({\rm max}_{x\in
M}R(x,0))^{-1/2}$ and if it is not of constant positive curvature,
then $M$ is diffeomorphic to either $S^3$, $\Ar P^3$,  $\Ar P^3\#\Ar
P^3$ or is a $S^2$-fibration over $S^1$. On the other hand, since
$M$ is compact of positive curvature its fundamental group is
finite, see Theorem 4.1 on p. 154 of \cite{Petersen}. This rules out
the last two cases. This implies that when $(M,g(0))$ has diameter
greater than $C_1({\rm max}_{x\in M}R(x,0))^{-1/2}$ and is not of
constant positive curvature, it is a double $C_0$-capped
$\epsilon$-tube.

We must consider the case when $(M,g(0))$ is not of constant
positive curvature and its diameter is less than or equal to
$C_1({\rm max}_{x\in M}R(x,0))^{-1/2}$. Since $(M,g(0))$ is not
round, by Corollary~\ref{compactcase} its asymptotic soliton is not
compact. Thus, by Theorem~\ref{GSSclass} its asymptotic soliton is
either $S^2\times \Ar$ or is double covered by this product. This
means that for $t$ sufficiently negative the diameter of $(M,g(t))$
is greater than $C_1({\rm max}_{x\in M}R(x,0))^{-1/2}$. Invoking the
previous result for this negative time tells us that $M$ is
diffeomorphic to $S^3$ or $\Ar P^3$.
\end{proof}

\begin{prop}
Let $\bar\epsilon_2$ and $\bar \epsilon_3$ be as in
Corollary~\ref{noncompkappa} and Theorem~\ref{thmep3}, respectively.
For each $0<\epsilon\le {\rm min}(\bar \epsilon_2,\bar\epsilon_3)$
let $C_1=C_1(\epsilon)$ be as in Theorem~\ref{thmep3}. There is
$C_2=C_2(\epsilon)<\infty$ such that for any $\kappa>0$ and any
compact $\kappa$-solution $(M,g(t))$ the following holds. If
$(M,g(0))$ is not of constant positive curvature and if $(M,g(0))$
is of diameter less than $C_1({\rm max}_{x\in M}R(x,0))^{-1/2}$ then
for any $x\in M$ we have
$$C_2^{-1}R(x,0)^{-3/2}< {\rm Vol}(M,g(0))< C_2R(x,0)^{-3/2}.$$
In addition, for any $y\in M$ and any $2$-plane $P_y$ in $T_yM$  we
have
$$C_2^{-1}< \frac{K(P_y)}{R(x,0)}< C_2,$$
where $K(P_y)$ is the sectional curvature in the $P_y$-direction.
\end{prop}

\begin{proof}
The result is immediate from Corollary~\ref{kappa0} and
Theorem~\ref{compactkappa}.
\end{proof}

\begin{rem}
For a round $\kappa$-solution $(M,g(t))$ we have $R(x,0)=R(y,0)$ for
all $x,y\in M$, and the volume of $(M,g(0))$ is bounded above by a
constant times $R(x,0)^{-3/2}$. There is no universal lower bound to
the volume in terms of the curvature. The lower bound takes the form
$C_2|\pi_1(M)|^{-1}R(x,0)^{-3/2}$, where $|\pi_1(M)|$ is the order
of the fundamental group $\pi_1(M)$.
\end{rem}

Let us summarize our results.

\begin{thm}\label{kappasummary}
There is  $\bar \epsilon>0$ such that the following is true for any
$0<\epsilon<\bar \epsilon$.  There is $C=C(\epsilon)$ such that for
any $\kappa>0$ and any $\kappa$-solution $(M,g(t))$ one of the
following holds.
\begin{enumerate}
\item[(1)] $(M,g(t))$ is round for all $t\le 0$. In this case $M$ is
diffeomorphic to the quotient of $S^3$ by a finite subgroup of
$SO(4)$ acting freely.
\item[(2)] $(M,g(0))$ is compact and of positive curvature. For any $x,y\in M$
and any $2$-plane $P_y$ in $T_yM$
 we have \begin{eqnarray*}
 C^{-1/2}R(x,0)^{-1} &  <{\rm diam}(M,g(0))
 &  <CR(x,0)^{-1/2} \\
 C^{-1}R(x,0)^{-3/2} & < {\rm Vol}(M,g(0)) & < CR(x,0)^{-3/2} \\
C^{-1}R(x,0) & < K(P_y) & < CR(x,0).
\end{eqnarray*}
  In this case $M$ is diffeomorphic either to
 $S^3$ or to $\Ar P^3$.
 \item[(3)] $(M,g(0))$ is of positive curvature and is a double  $C$-capped strong
 $\epsilon$-tube, and in
 particular $M$ is diffeomorphic to $S^3$ or to $\Ar P^3$.
 \item[(4)] $(M,g(0))$ is of positive curvature and is a
 $C$-capped strong $\epsilon$-tube and $M$ is diffeomorphic to $\Ar ^3$.
 \item[(5)] $(M,g(0))$ is isometric to the quotient of the product of
 a round $S^2$ and $\Ar$ by a
 free, orientation-preserving involution. It is a  $C$-capped
 strong $\epsilon$-tube and is diffeomorphic to a punctured $\Ar P^3$.
 \item[(6)] $(M,g(0))$ is isometric to the product of a round $S^2$ and $ \Ar$
  and is a strong $\epsilon$-tube.
 \item[(7)] $(M,g(0))$ is isometric to a product $\Ar P^2\times \Ar$, where the metric on $\Ar P^2$
 is of constant Gaussian curvature.
\end{enumerate}

In particular, in all cases except the first two and the last one,
all points of $(M,g(0))$ are either contained in the core of a
$(C,\epsilon)$-cap or are the centers of a strong $\epsilon$-neck in
$(M,g(0))$.

Lastly, in all cases we have \begin{eqnarray}{\rm sup}_{p\in M, t\le
0}\frac{|\nabla
R(p,t)|}{R(p,t)^{3/2}} & <C \label{nablaR} \\
{\rm sup}_{p\in M,t\le 0}\frac{\left|\partial R(p,t)/\partial
t\right|}{R(p,t)^2} & <C.\label{dRdt}
\end{eqnarray}
\end{thm}

An immediate consequence of this result is:

\begin{cor}\label{kappacannbhd}
For every $0<\epsilon\le\bar\epsilon'$ there is
$C=C(\epsilon)<\infty$ such that every point in a $\kappa$-solution
has a strong $(C,\epsilon)$-canonical neighborhood unless the
$\kappa$-solution is a product $\Ar P^2\times \Ar$.
\end{cor}

\begin{cor}\label{limitcannbhd}
Fix $0<\epsilon\le \bar\epsilon'$, and let $C(\epsilon)$ be as in
the last corollary. Suppose that $({\mathcal M}_n,G_n,x_n)$ is a
sequence of based, generalized Ricci flows with ${\bf t}(x_n)=0$ for
all $n$. Suppose that none of the time-slices of the ${\mathcal
M}_n$ contain embedded $\Ar P^2$'s with trivial normal bundle.
Suppose also that there is a smooth limiting flow
$(M_\infty,g_\infty(t),(x_\infty,0))$ defined for $-\infty<t\le 0$
that is a $\kappa$-solution. Then for all $n$ sufficiently large
$x_n$ has a strong $(C,\epsilon)$-canonical neighborhood in
$({\mathcal M}_n,G_n,x_n)$.
\end{cor}

\begin{proof}
The limiting manifold $M_\infty$ cannot contain an embedded $\Ar
P^2$ with trivial normal bundle. Hence, by the previous corollary,
the point $(x_\infty,0)$ has a strong $(C,\epsilon)$-canonical
neighborhood in the limiting flow. If the limiting $\kappa$-solution
is round, then for all $n$ sufficiently large $x_n$ is contained in
a component of the zero time-slice that is $\epsilon$-round. If
$(x_\infty,0)$ is contained in a $C$-component of the zero
time-slice of the limiting $\kappa$-solution, then for all $n$
sufficiently large $x_n$ is contained in a $C$-component of the zero
time-slice of ${\mathcal M}_n$. Suppose that $(x_\infty,0)$ is the
center of a strong $\epsilon$-neck in the limiting flow. This neck
extends backwards in the limiting solution some amount past an
interval of normalized time of length  $1$, where by continuity it
is an evolving $\epsilon$-neck defined backwards for an interval of
normalized time of length  greater than $1$. Then by Part 2 of
Proposition~\ref{canonvary}, any family of metrics on this neck
sufficiently close to the limiting metric will determine an strong
$\epsilon$-neck. This implies that for all $n$ sufficiently large
$x_n$ is the center of a strong $\epsilon$-neck in $({\mathcal
M}_n,G_n)$. Lastly, if $(x_\infty,0)$ is contained in the core of a
$(C,\epsilon)$-cap in the limiting flow, then by Part 1 of
Proposition~\ref{canonvary} for all $n$ sufficiently large $x_n$ is
contained in the core of a $(C,\epsilon)$-cap in $({\mathcal
M}_n,G_n)$.
\end{proof}

\chapter{Bounded curvature at bounded distance} This chapter is
devoted to Perelman's result about bounded curvature at bounded
distance for blow-up limits. Crucial to the argument is that each
member of the sequence of generalized Ricci flows has curvature
pinched toward positive and also has strong canonical neighborhoods.

\section{Pinching toward positive: the definitions}\label{10.1}

In this section we give the definition of what it means for a
generalized Ricci flow to have curvature pinched toward
positive\index{curvature!pinched toward positive|ii}. This is the
obvious generalization of the corresponding notion for Ricci flows.

\begin{defn}\label{pinchdef}
Let $({\mathcal M},G)$ be a generalized three-dimensional Ricci flow
whose domain of definition is contained in $[0,\infty)$. For each
$x\in {\mathcal M}$, let $\nu(x)$ be the smallest eigenvalue of
${\rm Rm}(x)$ on $\wedge^2T_xM_{{\bf t}(x)}$, as measured with
respect to a $G(x)$-orthonormal basis for the horizontal space at
$x$, and set $X(x)={\rm max}(0,-\nu(x))$. We say that $({\mathcal
M},G)$ has curvature {\em pinched toward positive} if, for all $x\in
{\mathcal M}$, if the following two inequalities hold:
\begin{enumerate}
\item[(1)] $$R(x)\ge \frac{-6}{1+4{\bf t}(x)},$$
\item[(2)] $$R(x)\geq 2X(x)\left({\rm log}X(x)+{\rm
log}(1+{\bf t}(x))-3\right),
$$ whenever
$0<X(x)$.
\end{enumerate}
\end{defn}

According to Theorem~\ref{pincha} if $(M,g(t)),\ 0\le a\le t<T$, is
Ricci flow with $M$  a compact three-manifold, and if the two
conditions given in the definition hold at the initial time $a$,
then they hold for all $t\in [a,T)$. In particular, if $a=0$ and if
$|{\rm Rm}(p,0)|\le 1$ for all $p\in M$, then the curvature of the
flow is pinched toward positive.

Next we fix $\epsilon_0>0$ sufficiently small such that for any
$0<\epsilon\le \epsilon_0$ all the results of the Appendix hold for
$2\epsilon$ and $\alpha=10^{-2}$, and Proposition~\ref{narrows}
holds for $2\epsilon$.

\section{The statement of the theorem}

Here is the statement of the main theorem of this chapter, the
theorem that establishes bounded curvature at bounded distance for
blow-up limits.

\begin{thm}\label{bcbd}
Fix $0<\epsilon\le \epsilon_0$ and  $C<\infty$. Then for each
$A<\infty$ there are $D_0<\infty$ and $D<\infty$ depending on $A$,
$\epsilon$ and $C$ such that the following holds. Suppose that
$({\mathcal M},G)$ is a generalized three-dimensional Ricci flow
whose interval of definition is contained in $[0,\infty)$, and
suppose that $x\in {\mathcal M}$. Set $t={\bf t}(x)$. We suppose
that these data satisfy the following:
\begin{enumerate}
\item[(1)] $({\mathcal M},G)$ has curvature pinched toward positive.
\item[(2)] Every point $y\in{\mathcal  M}$ with
$R(y)\ge 4R(x)$  and ${\bf t}(y)\le t$ has a strong
$(C,\epsilon)$-canonical neighborhood.
\end{enumerate}
If $R(x)\ge D_0$, then $R(y)\le DR(x)$ for all $y\in
B(x,t,AR(x)^{-1/2})$.
\end{thm}

This chapter is devoted to the proof of this theorem. The proof is
by contradiction. Suppose that there is some $A_0<\infty$ for which
the result fails. Then there are a sequence of generalized
three-dimensional Ricci flows $({\mathcal M}_n,G_n)$  whose
intervals of definition are contained in $[0,\infty)$ and whose
curvatures are pinched toward positive. Also, there are points
$x_n\in {\mathcal M}_n$ satisfying the second condition given in the
theorem and points $y_n\in {\mathcal M}_n$ such that for all $n$ we
have:
\begin{enumerate}
\item[(1)] ${\rm lim}_{n\rightarrow \infty}R(x_n)=\infty$.
\item[(2)] ${\bf t}(y_n)={\bf t}(x_n)$,
\item[(3)] $d(x_n,y_n)< A_0R(x_n)^{-1/2}$,
\item[(4)] $${\rm lim}_{n\rightarrow\infty}\frac{R(y_n)}{R(x_n)}=\infty.$$
\end{enumerate}
For the rest of this chapter we assume that such a sequence of
generalized Ricci flows exists. We shall eventually derive a
contradiction.

Let us sketch how the argument goes. We show that there is a
(partial) geometric blow-up limit of the sequence $({\mathcal
M}_n,G_n)$ based at the $x_n$. We shall see that the following hold
for this limit. It is an incomplete manifold $U_\infty$
diffeomorphic to $S^2\times (0,1)$ with the property that the
diameter of $U_\infty$ is finite and the curvature goes to infinity
at one end of $U_\infty$, an end denoted ${\mathcal E}$, while
remaining bounded at the other end. (The non-compact manifold in
question is diffeomorphic to $S^2\times (0,1)$ and, consequently, it
has two ends.)
 Every point of $U_\infty$
sufficiently close to ${\mathcal E}$ is the center of a
$2\epsilon$-neck in $U_\infty$. In fact, there is a partial
geometric limiting flow on $U_\infty$ so that these points are
centers of evolving $2\epsilon$-necks. Having constructed this
incomplete blow-up limit of the original sequence we then consider
further blow-up limits about the end ${\mathcal E}$, the end where
the scalar curvature goes to infinity. On the one hand, a direct
argument shows that a sequence of rescalings of $U_\infty$ around
points converging to the end ${\mathcal E}$ converge in the
Gromov-Hausdorff sense to a cone. On the other hand, a slightly
different sequence of rescalings at the same points converges
geometrically to a limiting non-flat Ricci flow. Since both limits
are non-degenerate three-dimensional spaces, we show that the ratio
of the rescaling factors used to construct them converges to a
finite, non-zero limit. This means that the two limits differ only
by an overall constant factor. That is to say the geometric blow-up
limit is isometric to an open subset of a non-flat cone. This
contradicts Hamilton's result (Theorem~\ref{nocones}) which says
that it is not possible to flow under the Ricci flow to an open
subset of a non-flat cone. Now we carry out all the steps in this
argument.

\section{The incomplete geometric limit}\label{sectincomp}

We fix a sequence $({\mathcal M}_n,G_n,x_n)$ of generalized Ricci
flows as above. The first step  is to shift and rescale this
sequence of generalized Ricci flows so that we can form an
(incomplete) geometric limit which will be a tube of finite length
with scalar curvature going to infinity at one end.

We shift the time parameter of $({\mathcal M}_n,G_n)$ by $-{\bf
t}(x_n)$. We change notation and denote these shifted flows by
$({\mathcal M}_n,G_n)$. This allows us to arrange that ${\bf
t}(x_n)=0$ for all $n$. Since shifting leaves the curvature
unchanged, the shifted flows  satisfy a weaker version of curvature
pinched toward positive. Namely, for the shifted flows we have
\begin{eqnarray}\label{weakerpos}
R(x) & \ge & -6 \nonumber \\
R(x) & \ge & 2X(x)\left({\rm log}(X(x))-3\right).
\end{eqnarray}

 We set $Q_n=R(x_n)$, and we denote by $M_n$ the $0$-time-slice of
${\mathcal M}_n$. We rescale $({\mathcal M}_n,G_n)$ by $Q_n$. Denote
by $({\mathcal M}'_n,G_n')$ the rescaled (and shifted) generalized
flows. For the rest of this argument we implicity use the metrics
$G'_n$. If we are referring to $G_n$ we mention it explicitly.

\subsection{The sequence of tubes}

 Let
$\gamma_n$ be a smooth path in $B_{G_n}(x_n,0,A_0Q_n^{-1/2})$ from
$x_n$ to $y_n$. For all $n$ sufficiently large we have
$R_{G'_n}(y_n) \gg 1$. Thus, there is a point $z_n\in \gamma_n$ such
that $R_{G'_n}(z_n)=4$ and such that on the sub-path
$\gamma_n|_{[z_n,y_n]}$ we have $R_{G'_n}\ge 4$. We replace
$\gamma_n$ by this sub-path. Now, with this replacement, according
to the second condition in the statement of the theorem, every point
of $\gamma_n$ has a strong $(C,\epsilon)$ canonical
neighborhood\index{canonical neighborhood}. As $n$ tends to infinity
the ratio of $R(y_n)/R(z_n)$ tends to infinity. This means that for
all $n$ sufficiently large, no point of  $\gamma_n$ can be contained
in an $\epsilon$-round component\index{$\epsilon$-round component}
or a $C$-component\index{$C$-component}, because if it were then all
of $\gamma_n$ would be contained in that component, contradicting
the fact that the curvature ratio is arbitrarily large for large
$n$. Hence, for $n$ sufficiently large, every point of $\gamma_n$ is
either contained in the core of a
$(C,\epsilon)$-cap\index{$(C,\epsilon)$-cap} or is the center of a
strong $\epsilon$-neck\index{$\epsilon$-neck!strong}. According to
Proposition~\ref{Xcontainedin}, for all $n$ sufficiently large
$\gamma_n$ is contained an open submanifold $X_n$ of the zero
time-slice of ${\mathcal M}'_n$ that is one of the following:
\begin{enumerate}
\item[(1)] an $\epsilon$-tube  and both endpoints of $\gamma_n$ are centers
of $\epsilon$-necks
contained in $X_n$,
\item[(2)] a $C$-capped $\epsilon$-tube\index{$\epsilon$-tube!capped} with cap ${\mathcal C}$,
and each endpoint of $\gamma_n$ either is contained in the core $Y$
of ${\mathcal C}$ or is the center of an $\epsilon$-neck contained
in $X_n$,
\item[(3)] a double $C$-capped $\epsilon$-tube\index{$\epsilon$-tube!doubly capped}, or finally
\item[(4)] the union of two $(C,\epsilon)$-caps.
\end{enumerate}

The fourth possibility is incompatible with the fact that the ratio
of the curvatures at the endpoints of $\gamma_n$ grows arbitrarily
large as $n$ tends to infinity. Hence, this fourth possibility
cannot occur for $n$ sufficiently large. Thus, for all $n$
sufficiently large $X_n$ is one of the first three types listed
above.

\begin{claim}
There is a geodesic $\hat\gamma_n$ in $X_n$ with endpoints $z_n$ and
$y_n$. This geodesic is minimizing among all paths in $X_n$ from
$z_n$ to $y_n$.
\end{claim}

\begin{proof}
This is clear in the third case since $X_n$ is a closed manifold.

Let us consider the first case. There are $\epsilon$-necks $N(z_n)$
and $N(y_n)$ centered at $z_n$ and $y_n$ and contained in $X_n$.
Suppose first that the central $2$-spheres $S(z_n)$ and $S(y_n)$ of
these necks are disjoint. Then they are the boundary of a compact
sumanifold $X_n'$ of $X_n$. It follows easily from
Lemma~\ref{curveshort} that any sequence of minimizing paths from
$z_n$ to $y_n$ is contained in the union of $X_n'$ with the  middle
halves of $N(z_n)$ and $N(y_n)$. Since this manifold has compact
closure in $X_n$, the usual arguments show that one can extract a
limit of a subsequence which is a minimizing geodesic in $X_n$ from
$z_n$ to $y_n$. If $S(z_n)\cap S(y_n)\not=\emptyset$, then $y_n$ is
contained in the middle half of $N(z_n)$, and again it follows
immediately from Lemma~\ref{curveshort} that there is a minimizing
geodesic in $N(z_n)$ between these points.

Now let us consider the second case. If each of $z_n$ and $y_n$ is
the center of an $\epsilon$-neck in $X_n$, the argument as in the
first case applies. If both points are contained in the core of
${\mathcal C}$ then, since that core has compact closure in $X_n$,
the result is again immediate. Lastly, suppose that one of the
points, we can assume by the symmetry of the roles of the points
that it is $z_n$, is the center of an $\epsilon$-neck $N(z_n)$ in
$X_n$ and the other is contained in the core of ${\mathcal C}$.
Suppose that the central $2$-sphere $S(z_n)$ of $N(z_n)$  meets the
core $Y$ of ${\mathcal C}$. Then $z_n$ lies in the  half of the neck
$N={\mathcal C}\setminus \overline Y$ whose closure contains the
frontier of $Y$. Orient $s_N$ so that this half is the positive
half. Thus, by Lemma~\ref{curveshort} any minimizing sequence of
paths from $z_n$ to $y_n$ is eventually contained in the union of
the core of ${\mathcal C}$ and the the positive three-quarters of
this neck. Hence, as before we can pass to a limit and construct  a
minimizing geodesic in $X_n$ connecting $z_n$ to $y_n$. On the other
hand, if $S(z_n)$ is disjoint from $Y$, then $S(z_n)$ separates
$X_n$ into a compact complementary component and a non-compact
complementary component and the compact complementary component
contains $Y$. Orient the $s_N$-direction so that the compact
complementary component lies on the positive side of $S(z_n)$. Then
any minimizing sequence of paths in $X_n$ from $z_n$ to $y_n$ is
eventually contained in the union of the compact complementary
component of $N(z_n)$ and the positive $3/4$'s of $N(z_n)$. As
before, this allows us to pass to a limit to obtain a minimizing
geodesic in $X_n$.
\end{proof}

This claim allows us to assume (as we now shall) that $\gamma_n$ is
a minimizing geodesic in $X_n$ from $z_n$ to $y_n$.

\begin{claim}
For every $n$ sufficiently large, there is a sub-geodesic
$\gamma'_n$ of $\gamma_n$ with end points $z_n'$ and $y_n'$ such
that the following hold:
\begin{enumerate}
\item[(1)] The length of $\gamma'_n$ is bounded independent of $n$.
\item[(2)] $R(z_n')$ is bounded independent of $n$.
\item[(3)] $R(y_n')$ tends to infinity as $n$ tends to infinity.
\item[(4)] $\gamma_n'$ is contained in a strong $\epsilon$-tube
$T_n$ that is the union of a balanced
chain\index{$\epsilon$-neck!strong!balanced chain of} of strong
$\epsilon$-necks centered at points of $\gamma'_n$. The first
element in this chain is a strong $\epsilon$-neck $N(z'_n)$ centered
at $z'_n$. The last element is a strong $\epsilon$-neck containing
$y'_n$.
\item[(5)] For every  $x\in T_n$, we have $R(x)>3$ and $x$ is the center of a strong
$\epsilon$-neck in the flow $({\mathcal M}'_n,G'_n)$.
\end{enumerate}
\end{claim}

\begin{proof}
The first item is clear since, for all $n$, the geodesic $\gamma_n$
has $G_n$-length at most $A_0Q_n^{-1/2}$ and hence $G'_n$-length at
most $A_0$. Suppose that we have a $(C,\epsilon)$-cap ${\mathcal C}$
whose core $Y$ contains a point of $\gamma_n$. Let $N$ be the
$\epsilon$-neck that is the complement of the closure of $Y$ in
${\mathcal C}$, and let $\widehat Y$ be the union of $Y$ and the
closed negative half of $N$. We claim that $\widehat Y$ contains
either $z_n$ or $y_n$. By Corollary~\ref{coreint}, since $Y$
contains a point of $\gamma_n$, the intersection of $\widehat Y$
with $\gamma_n$ is a subinterval containing one of the end points of
$\gamma_n$, i.e., either $z_n$ or $y_n$. This means that any point
$w$  which is contained in a $(C,\epsilon)$-cap whose core contains
a point of $\gamma_n$ must satisfy one of the following:
$$R(w)< CR(z'_n)\ \ \ {\rm or}\ \ \ R(w)> C^{-1}R(y'_n).$$

We pass to a subsequence so that $R(y_n)/R(z_n)>4C^2$ for all $n$,
and we pass to a subinterval $\gamma'_n$ of $\gamma_n$ with
endpoints $z_n'$ and $y_n'$ such that:
\begin{enumerate}
\item[(1)] $R(z'_n)=2CR(z_n)$
\item[(2)] $R(y'_n)=(2C)^{-1}R(y_n)$
\item[(3)] $R(z'_n)\le R(w)\le R(y'_n)$ for all $w\in \gamma'_n$.
\end{enumerate}
Clearly, with these choices $R(z'_n)$ is bounded independent of $n$
and $R(y'_n)$ tends to infinity as $n$ tends to infinity. Also, no
point of $\gamma'_n$ is contained in the core of a
$(C,\epsilon)$-cap. Since every point of $\gamma'_n$ has a strong
$(C,\epsilon)$-canonical neighborhood, it follows that every point
of $\gamma'_n$ is the center of a strong $\epsilon$-neck. It now
follows from Proposition~\ref{epschains} that there is a balanced
$\epsilon$-chain consisting of strong $\epsilon$-necks centered at
points of $\gamma'_n$ whose union contains $\gamma'_n$. (Even if the
$2$-spheres of these necks do not separate the zero time-slice of
${\mathcal M}'_n$, as we build the balanced $\epsilon$-chain as
described in Proposition~\ref{epschains} the new necks we add can
not meet the negative end of $N(z'_n)$ since the geodesic
$\gamma'_n$ is minimal.) We can take the first element in the
balanced chain to be a strong $\epsilon$-neck $N(z'_n)$ centered at
$z'_n$, and the last element to be a strong $\epsilon$-neck $N^+_n$
containing $y'_n$. The union of this chain is $T_n$. (See {\sc
Fig.}~\ref{fig:N}.)

\begin{figure}[ht]
  \relabelbox{
  \centerline{\epsfbox{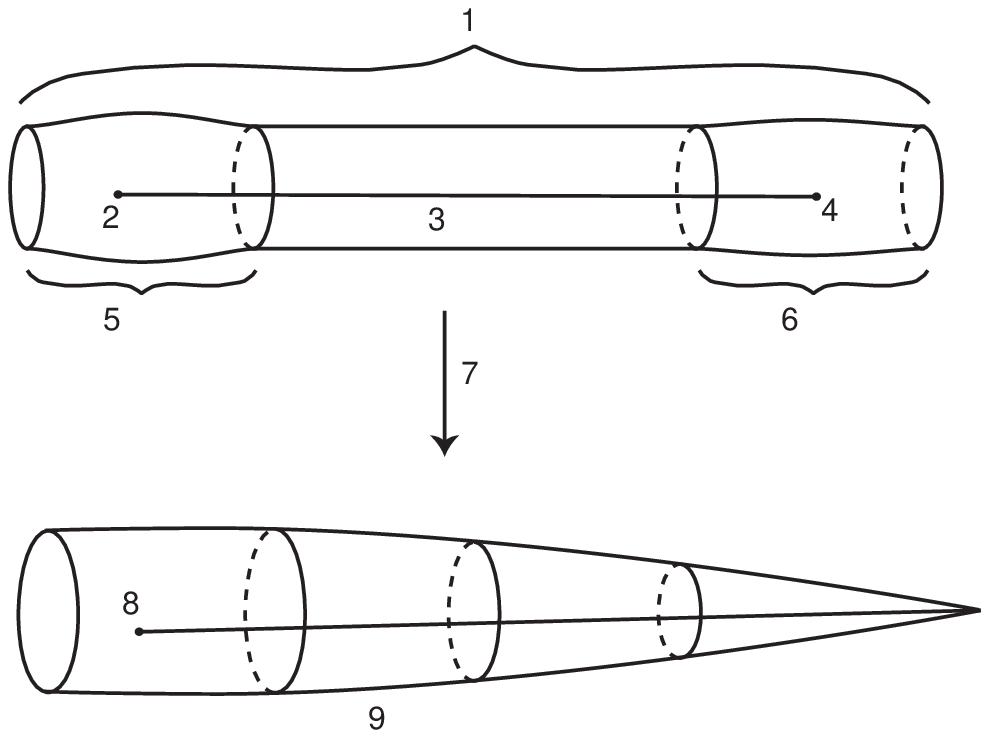}}}
  \relabel{1}{$T_n$}
  \relabel{2}{$z_n'$}
  \relabel{3}{$\gamma_n'$}
  \relabel{4}{$y_n'$}
  \relabel{5}{$N(z_n')$}
  \relabel{6}{$N(y_n')$}
  \relabel{7}{converge to}
  \relabel{8}{$z_\infty'$}
  \relabel{9}{$(U_\infty,g_\infty, z_\infty)$}
  \endrelabelbox
  \caption{Limiting tube}\label{fig:N}
\end{figure}

Next, we show that every point of $T_n$ is the center of a strong
$\epsilon$-neck in $({\mathcal M}_n,G_n)$. We must rule out the
possibility that there is a
 point of $T_n$ that is contained in the core of a $(C,\epsilon)$-cap\index{$(C,\epsilon)$-cap!core of}.
Since $T_n$ is a union of $\epsilon$-necks centered at points of
$\gamma_n'$ we see that every point  $w\in T_n$ has
$$(3C/2)R(z_n)<R(w)<(2/3C)R(y_n).$$
This tells us that no point of $T_n$ is contained in a
$(C,\epsilon)$-cap whose core contains a point of $\gamma_n$. Thus,
to complete the argument we need only see that if there is a point
of $T_n$ contained in the core of a $(C,\epsilon)$-cap then the core
of that $(C,\epsilon)$-cap also contains a point of $\gamma_n$.
 The
scalar curvature inequality implies that both $z_n$ and $y_n$ are
outside $T_n$. This means that $\gamma_n$ traverses $T_n$ from one
end to the other. Let $w_-$, resp. $w_+$, be the point of $\gamma_n$
that lies in the frontier of $T_n$ contained in the closure of the
$N(z_n')$, resp. $N^+_n$. Since the scalar curvatures at these two
points of $\gamma$ satisfy the weak version of the above
inequalities, we see that there are strong $\epsilon$-necks $N(w_-)$
and $N(w_+)$ centered at them. Let $\widehat T_n$ be the union of
$T_n$, $N(w_-)$ and $N(w_+)$. It is also a strong $\epsilon$-tube,
and every point $\hat w$ of $\widehat T_n$ satisfies
$$(1.1)CR(z_n)<R(\hat w)<(0.9)C^{-1}R(y_n).$$
Thus, $z_n$ and $y_n$ are disjoint from $\widehat T_n$ and hence
$\gamma$ crosses $\widehat T_n$ from one end to the other.

 Now
suppose that $T_n$ meets the core $Y$ of a $(C,\epsilon)$-cap
${\mathcal C}$.
 Consider the boundary $S$ of the closure of $Y$.
If it is disjoint from $T_n$ then $T_n$ is contained in the core
$Y$. For large $n$ this is inconsistent with the fact that the ratio
of the scalar curvature at the endpoints of $\gamma'_n$ goes to
infinity. Thus, we are left to consider the case when
 $S$ contains a point of the tube $T_n$. In this case $S$ is
completely contained in $\widehat T_n$ and by
Corollary~\ref{s2isotopic} $S$ is isotopic to the $2$-spheres of the
product decomposition of $\widehat T_n$. Hence, $S$ meets a point of
$\gamma_n$ and consequently the core $Y$ contains a point of
$\gamma_n$. But we have already seen that this is not possible.

Lastly, we must show that $R(x)>3$ for every $x\in T_n$. We have
just seen that every $x\in T_n$ is the center of an $\epsilon$-neck.
If $x$ is contained in the $\epsilon$-neck centered at $z'_n$ or
$y'_n$, then since $R(z'_n)\ge 4$ and $R(y'_n)\ge 4$, clearly
$R(x)>3$. We must consider the case when $x$ is not contained in
either of these $\epsilon$-necks. In this case the central
$2$-sphere $S_x$of the $\epsilon$-neck centered at $x$ is contained
in the compact submanifold of $T_n'$ bounded by the central
$2$-spheres of the necks centered at $z'_n$ and $y'_n$. These
$2$-spheres are disjoint and by Condition 4 in
Proposition~\ref{S2intersection} each is a homotopically non-trivial
$2$-sphere in $T_n'$. Hence, the compact manifold with their
disjoint union as boundary is diffeomorphic to $S^2\times [0,1]$
and, again according to Condition 4 of
Proposition~\ref{S2intersection}, $S_x$ is isotopic to the
$2$-sphere factor in this product decomposition. Since the
intersection of $\gamma'_n$ with this submanifold is an arc spanning
from one boundary component to the other, $S_x$ must meet
$\gamma'_n$, in say $w$.
 By construction, since $w\in \gamma'_n$ we have $R(w)\ge 4$.
  This implies that $R(x)>3$.
 This completes the proof of the claim.
\end{proof}

\subsection{Extracting a limit of a subsequence of the tubes}

Passing to a subsequence we arrange that the $R(z_n')$ converge. Now
consider the subset ${\mathcal A}\subset \Ar$ consisting of all
$A>0$ such that there is a uniform bound, independent of $n$, for
the curvature on $B(z'_n,A)\cap T_n$. The set ${\mathcal A}$ is
non-empty since $R(z'_n)$ is bounded independent of $n$ and for
every $n$ there is a strong $\epsilon$-neck $N(z'_n)$ centered at
$z'_n$ contained in $T_n$. On the other hand, since
$d_{G_n'}(z_n',y_n')$ is uniformly bounded and $R(y_n')\rightarrow
\infty$,  there is a finite upper bound for ${\mathcal A}$. Let
$A_1$ be the least upper bound of ${\mathcal A}$. We set
$U_n=T_n\cap B(z'_n,A_1)$. This is an open subset of $T_n$
containing $z'_n$. We let $g_n'=G_n'|U_n$.

\begin{claim}
For all $n$ sufficiently large,  $3R(z'_n)^{-1/2}\epsilon^{-1}/2$ is
less than $A_1$, and hence $U_n$ contains the strong $\epsilon$-neck
$N(z_n')$ centered at $z'_n$.
\end{claim}

\begin{proof}
The curvature on  $N(z'_n)$ is bounded independent of $n$. Consider
a point $w$ near the end of $N(z'_n)$ that separates $y_n'$ from
$z_n'$. It is also the center of a strong $\epsilon$-neck $N(w)$. By
Proposition~\ref{S2intersection} and our assumption that
$\epsilon\le \bar\epsilon(10^{-2})$, the scalar curvature on
$N(z'_n)\cup N(w)$ is between $(0.9)R(z'_n)$ and $(1.1)R(z'_n)$.
Since, by construction, the negative end of $N(z'_n)$ contains an
end of $T_n$, this implies that
$$N(z'_n)\cup N(w)\supset B(z'_n,7R(z'_n)^{-1/2}\epsilon^{-1}/4)\cap T_n,$$
 so that we see that
$A_1\ge 7\epsilon^{-1}{\rm lim}_{n\rightarrow\infty}
R(z'_n)^{-1/2}/4$. Thus, $A_1>3R(z'_n)^{-1}\epsilon^{-1}/2$ for all
$n$ sufficiently large. Obviously then $U_n$ contains $N(z'_n)$.
\end{proof}

The next claim uses terminology from Definition~\ref{deltareg}.

\begin{claim}
For any $\delta>0$ there is a uniform bound, independent of $n$, for
the curvature on ${\rm Reg}_\delta(U_n,g'_n)$.
\end{claim}

\begin{proof}
To prove this it suffices to show that given $\delta>0$ there is
$A<A_1$ such that ${\rm Reg}_\delta(U_n,g'_n)\subset B(z'_n,A)$ for
all $n$ sufficiently large. Of course, if we establish this for
every $\delta>0$ sufficiently small, then it follows for all
$\delta>0$. First of all, by Corollary~\ref{neckgeo} and
Lemma~\ref{neckcurv}, the fact that $\epsilon\le
\bar\epsilon(10^{-2})$ implies that any point $w$ with the property
that the strong $\epsilon$-neighborhood centered at $w$ contains
$z_n'$ is contained in the ball of radius
$(1.1)R(z'_n)^{-1/2}\epsilon^{-1}<A_1$ centered at $z_n'$. Thus, it
suffices to consider points $w_n$ in ${\rm Reg}_\delta(U_n,g'_n)$
with the property that the strong $\epsilon$-neck centered at $w_n$
does not contain $z_n'$. Fix such a $w_n$. Take a path $\mu_n(s)$
starting at $w_n$ moving in the $s$-direction at unit speed measured
in the $s$-coordinate of the $\epsilon$-neck centered at $w_n$ away
from $z'_n$ and ending at the frontier of this neck. Let $u_1$ be
the final point of this path. The rescaled version of
Lemma~\ref{neckdistcomp} implies that the forward difference
quotient for the distance from $z'_n$ satisfies
$$(0.99)R(w_n)^{-1/2}\le \frac{d}{ds}d(z'_n,\mu_n(s))
\le (1.01)R(w_n)^{-1/2}.$$ Of course, since we are working in an
$\epsilon$-neck we also have
$$(1-\epsilon)R(w_n)^{-1/2}\le \frac{d(d(w_n,\mu_n(s)))}{ds}\le
(1+\epsilon)R(w_n)^{-1/2}.$$  We continue the path $\mu_n$ moving in
the $s$-direction of a neck centered at $u_1$. Applying
Lemma~\ref{neckdistcomp} again both to the distance from $w_n$ and
the distance from $z'_n$ yields:
$$(0.99)R(u_1)^{-1/2}\le \frac{d(d(z'_n,\mu_n(s)))}
{ds}\le (1.01)R(u_1)^{-1/2}$$
$$(0.99)R(u_1)^{-1/2}\le \frac{d(d(w_n,\mu_n(s)))}
{ds}\le (1.01)R(u_1)^{-1/2}$$ on this part of the path $\mu_n$. We
repeat this process as many times as necessary until we reach a
point $w'_n\in U_n$ at distance $\delta/2$ from $w_n$. This is
possible since the ball of radius $\delta$ centered at $w_n$ is
contained in $U_n$. By the difference quotient inequalities, it
follows that $d(z'_n,w'_n)-d(z'_n,w_n)>\delta/4$.
 Since $w_n'\in U_n$ and
consequently that $d(z_n',w_n')<A_1$. It follows that
$d(z_n',w_n)\le A_1-\delta/4$. This proves that, for all $n$
sufficiently large, ${\rm Reg}_\delta(U_n,g'_n)\subset
B(z_n',A_1-\delta/4)$, and consequently that the curvature on ${\rm
Reg}_\delta(U_n,g'_n)$ is bounded independent of $n$.
\end{proof}

By Shi's theorem (Theorem~\ref{shi}), the fact that each point of
$U_n$ is the center of a strong $\epsilon$-neck means that there is
a bound, independent of $n$, on all covariant derivatives of the
curvature at any point of $U_n$ in terms of the bound on the
curvature at the center point. In particular, because of the
previous result, we see that for any $\epsilon>0$ and any $\ell\ge
0$ there is a uniform bound for $|\nabla^\ell {\rm Rm}|$ on ${\rm
Reg}_\delta(U_n,g'_n)$. Clearly, since the base point $z_n'$ has
bounded curvature it lies in ${\rm Reg}_\delta(U_n,g'_n)$ for
sufficiently small $\delta$ (how small being independent of $n$).
Lastly, the fact that every point in $U_n$ is the center of an
$\epsilon$-neighborhood implies that $(U_n,g'_n)$ is $\kappa$
non-collapsed on scales $\le r_0$ where both $\kappa$ and $r_0$ are
universal. Since the $\gamma'_n$ have uniformly bounded lengths, the
$\epsilon$-tubes $T_n'$ have uniformly bounded diameter. Also, we
have seen that their have curvatures are bounded from below by $3$.
It follows that their volumes are uniformly bounded. Now invoking
Theorem~\ref{basicconv} we see that after passing to a subsequence
we have a geometric limit $(U_\infty,g_\infty,z_\infty)$ of a
subsequence of $(U_n,g'_n,z'_n)$.

\subsection{Properties of the limiting tube}

Now we come to a result establishing all the properties we need for
the limiting manifold.

\begin{prop}\label{Uinfty}
The geometric limit $(U_\infty,g_\infty,z_\infty)$ is an incomplete
Riemannian $3$-manifold of finite diameter. There is a
diffeomorphism $\psi\colon U_\infty\to S^2\times (0,1)$. There is a
$2\epsilon$-neck centered at $z_\infty$ whose central $2$-sphere
$S^2(z_\infty)$ maps under $\psi$ to a $2$-sphere isotopic to a
$2$-sphere  factor in the product decomposition. The scalar
curvature is bounded at one end of $U_\infty$ but tends to infinity
at the other end, the latter end  which is denoted ${\mathcal E}$.
Let ${\mathcal U}_\infty\subset U_\infty\times(-\infty,0]$ be the
open subset consisting of all $(x,t)$ for which $-R(x)^{-1}<t\le 0$.
We have a generalized Ricci flow on ${\mathcal U}_\infty$ which is a
partial geometric limit of a subsequence of the generalized Ricci
flows $({\mathcal M}'_n,G'_n,z'_n)$. In particular, the zero-time
slice of the limit flow is $(U_\infty,g_\infty)$. The Riemannian
curvature is non-negative at all points of the limiting smooth flow
on ${\mathcal U}_\infty$. Every point $x\in U_\infty\times\{0\}$
which is not separated from ${\mathcal E}$ by $S^2(z_\infty)$ is the
center of an evolving $2\epsilon$-neck $N(x)$ defined for an
interval of normalized time of length $1/2$. Furthermore, the
central $2$-sphere of $N(x)$ is isotopic to the $2$-sphere factor of
$U_\infty$ under the diffeomorphism $\psi$ (see {\sc
Fig.}~\ref{fig:N}).
\end{prop}

The proof of this proposition occupies the rest of
Chapter~\ref{sectincomp}.

\begin{proof}
Let $V_1\subset V_2 \subset \cdots \subset U_\infty$ be the open
subsets and $\varphi_n\colon V_n\to U_n$ be the maps having all the
properties stated in Definition~\ref{smoothconv} so as to exhibit
$(U_\infty,g_\infty,z_\infty)$ as the  geometric limit of the
$(U_n,g_n',z_n')$.

Since the $U_n$ are all contained in $B(z_n',A_1)$, it follows that any point
of $U_\infty$ is within $A_1$ of the limiting base point $z_\infty$. This
proves that the diameter of $U_\infty$ is bounded.

For each $n$ there is the  $\epsilon$-neck $N(z_n')$ centered at
$z_n'$ contained in $U_n$. The middle two-thirds, $N_n'$, of this
neck has closure contained in ${\rm Reg}_\delta(U_n,g_n)$ for some
$\delta>0$ independent of $n$ (in fact, restricting to $n$
sufficiently large, $\delta$ can be taken to be approximately equal
to $R(z_\infty)^{-1/2}\epsilon^{-1}/3$). This means that for some
$n$ sufficiently large and for all $m\ge n$ the image
$\varphi_m(V_n)\subset U_m$ contains $N'_m$. For any fixed $n$ as
$m$ tends to infinity  the metrics $\varphi_m^*g_m|_{V_n}$ converge
uniformly in the $C^\infty$-topology to $g_\infty|_{V_n}$. Thus, it
follows from Proposition~\ref{canonvary} that for all $m$
sufficiently large, $\varphi^{-1}_m(N'_m)$ is a $3\epsilon/2$-neck
centered at $z_\infty$. We fix such a neck $N'(z_\infty)\subset
U_\infty$. Let $S(z_\infty)$ be the central $2$-sphere of
$N'(z_\infty)$. For each $n$ sufficiently large,
$\varphi_n(S(z_\infty))$ separates $U_n$ into two components, one,
say $W^-_n$ contained in $N(z'_n)$ and the other, $W^+_n$ containing
all of $U_n\setminus N(z'_n)$. It follows that $S(z_\infty)$
separates $U_\infty$ into two components, one, denoted $W^-_\infty$,
where the curvature is bounded (and where, in fact, the curvature is
close to $R(z_\infty)$) and the other, denoted $W^+_\infty$, where
it is unbounded.

\begin{claim}
Any point $q\in W^+_\infty$ is the center of a $2\epsilon$-neck in
$U_\infty$.
\end{claim}

\begin{proof}
Fix a point $q\in W^+_\infty$.
 For all $n$ sufficiently large denote by $q_n=\varphi_n(q)$.
 Then for all $n$ sufficiently large, $q_n\in W^+_n$ and
${\rm lim}_{n\rightarrow\infty}R(q_n)=R(q)$. This means that for all
$n$ sufficiently large $R(y_n')>>R(q_n))$, and hence
 the $3\epsilon/2$-neck centered at
 $q_n\in U_n$ is disjoint from $N(y'_n)$. Thus, by the rescaled version of
 Corollary~\ref{neckgeo}, we see that  the distance from
 the $3\epsilon/2$-neck centered at $q_n$ to $N(y'_n)$ is bounded below by
 $(0.99)\epsilon^{-1}R(q_n)^{-1/2}/4\ge \epsilon^{-1}R(q_\infty)^{-1/2}/12$.
 Also, since $q_n\in W_n$, this $3\epsilon/2$-neck $N'(q_n)$ centered at
 $q_n$ does not extend past the $2$-sphere at $s^{-1}(-3\epsilon^{-1}/4)$ in the
 $ \epsilon$-neck $N(z'_n)$. It follows that for all $n$ sufficiently large that
 this $3\epsilon/2$-neck has
 compact closure contained in ${\rm Reg}_\delta(U_n,g_n)$ for some $\delta$
 independent of $n$, and hence there is $m$ such that for all $n$ sufficiently large $N'(q_n)$ is
 contained in the image $\varphi_n(V_m)$. Again using the fact that
 $\varphi_n^*(g_n|_{V_m})$ converges in the $C^\infty$-topology to
 $g_\infty|_{V_m}$ as $n$ tends to infinity, we see, by Proposition~\ref{canonvary}
 that for all $n$
 sufficiently large $\varphi_n^{-1}(N_m)$ contains a $2\epsilon$-neck in $U_\infty$ centered
 at $q$.
\end{proof}

It now follows from Proposition~\ref{Xcontainedin} that $W^+_\infty$
is contained in an $2\epsilon$-tube $T_\infty$ that is contained in
$U_\infty$. Furthermore, the frontier of $W^+_\infty$ in $T_\infty$
is the $2$-sphere $S(z_\infty)$ which is isotopic to the central
$2$-spheres of the $2\epsilon$-necks making up $T_\infty$. Hence,
the closure $\overline W^+_\infty$ of $W^+_\infty$ is a
$2\epsilon$-tube with boundary $S(z_\infty)$. In particular,
$\overline W^+_\infty$ is diffeomorphic to $S^2\times [0,1)$.

Now we consider the closure $\overline W^-_\infty$ of  $W^-_\infty$.
Since the closure of each $W^-_n$ is the closed negative half of the
$\epsilon$-neck $N(z_n')$ and the curvatures of the $z_n'$ have a
finite, positive limit, the limit $\overline W^-_\infty$ is
diffeomorphic to a product $S^2\times (-1,0]$. Hence, $U_\infty$ is
the union of  $\overline W^+_\infty$ and  $\overline W^-_\infty$
along their common boundary. It follows immediately that $U_\infty$
is diffeomorphic to $S^2\times (0,1)$.

\begin{claim}
The curvature is bounded in a neighborhood of one end of $U_\infty$
and goes to infinity at the other end. \end{claim}

\begin{proof}
 A neighborhood of one end of $U_\infty$, the end $\overline W^-_\infty$,
 is the limit of the
 negative halves of $\epsilon$-necks centered at $z_n'$.
 Thus, the curvature is bounded on this neighborhood, and in fact is
approximately equal to $R(z_\infty)$. Let $x_k$ be any sequence of
points in $U_\infty$ tending to the other end. We show that $R(x_k)$
tends to $\infty$ as $k$ does. The point is that since the sequence
is tending to the end, the distance from $x_k$ to the end of
$U_\infty$ is going to zero. Yet, each $x_k$ is the center of an
$\epsilon$-neck in $U_\infty$. The only way this is possible is if
the scales of these $\epsilon$-necks are converging to zero as $k$
goes to infinity. This is equivalent to the statement that $R(x_k)$
tends to $\infty$ as $k$ goes to infinity.
\end{proof}

The next step in the proof of Proposition~\ref{Uinfty} is to extend
the flow backwards a certain amount. As stated in the proposition,
the amount of backward time that we can extend the flow is not
uniform over all of $U_\infty$, but rather depends on the curvature
of the point at time zero.

\begin{claim}
For each $x\in U_n\subset M_n$ there is a flowline $\{x\}\times
(-R(x)^{-1},0]$ in ${\mathcal M}_n$. Furthermore, the scalar
curvature at any point of this flow line is less than or equal to
the scalar curvature at $x$.
\end{claim}

\begin{proof}
Since $x\in U_n\subset T_n$, there is a strong $\epsilon$-neck in
${\mathcal M}_n$ centered at $x$. Both statements follow immediately
from that.
\end{proof}

Let $X\subset U_\infty$ be an open submanifold with compact closure
and set
$$t_0(X)= {\rm sup}_{x\in X}(-R_{g_\infty}(x)^{-1}).$$ Then for
all $n$ sufficiently large $\varphi_n$ is defined on $X$ and the
scalar curvature of the flow $g_n(t)$ on $\varphi_n(X)\times
(t_0,0]$ is uniformly bounded independent of $n$. Thus, according to
Proposition~\ref{partialflowlimit} by passing to a subsequence we
can arrange that there is a limiting flow defined on $X\times
(t_0,0]$. Let ${\mathcal U}_\infty\subset U_\infty\times
(-\infty,0]$ consist of all pairs $(x,t)$ with the property that
$-R_{g_\infty}(x,0)^{-1}<t\le 0$.
 Cover ${\mathcal U}_\infty$ by countably many such boxes
of the type $X\times (-t_0(X),0]$ as described above, and take a
diagonal subsequence. This allows us to pass to a subsequence so
that the limiting flow exists (as a generalized Ricci flow) on
${\mathcal U}_\infty$.

\begin{claim}
 The curvature of the generalized Ricci flow on ${\mathcal
U}_\infty$ is non-negative.
\end{claim}

\begin{proof}
 This claim follows from the fact that the
original sequence $({\mathcal M}_n,G_n)$ consists of generalized
flows whose curvatures are pinched toward positive in the weak sense
given in Equation~\ref{weakerpos} and the fact that $Q_n\rightarrow
\infty$ as $n\rightarrow \infty$. (See Theorem~\ref{blowupposcurv}.)
\end{proof}

This completes the proof that all the properties claimed in
Proposition~\ref{Uinfty} hold for the geometric limit
$(U_\infty,g_\infty,z_\infty)$. This  completes the proof of that
proposition. \end{proof}

\section{Cone limits near the end ${\mathcal E}$ for rescalings of $U_\infty$}\label{GromHaus}

The next step is to study the nature  of the limit $U_\infty$ given
in Proposition~\ref{Uinfty}. We shall show that an appropriate
blow-up limit (limit in the Gromov-Hausdorff sense) around the end
is a cone.

Let $(X,d_X)$ be a metric space. Recall that the  cone on $X$,
denoted $C(X)$, is the quotient space $X\times [0,\infty)$ under the
identification $(x,0)\cong (y,0)$ for all $x,y\in X$. The image of
$X\times\{0\}$ is the {\sl cone point} of the cone. The metric on
$C(X)$ is given by
\begin{equation}\label{conemetric}
d((x,s_1),(y,s_2))=s_1^2+s_2^2-2s_1s_2{\rm cos}({\rm
min}(d_X(x,y),\pi)). \end{equation} The open cone $C'(X)$ is the
complement of the cone point in $C(X)$ with the induced metric.

The purpose of this section is to prove the following result.

\begin{figure}[ht]
  \relabelbox{
  \centerline{\epsfbox{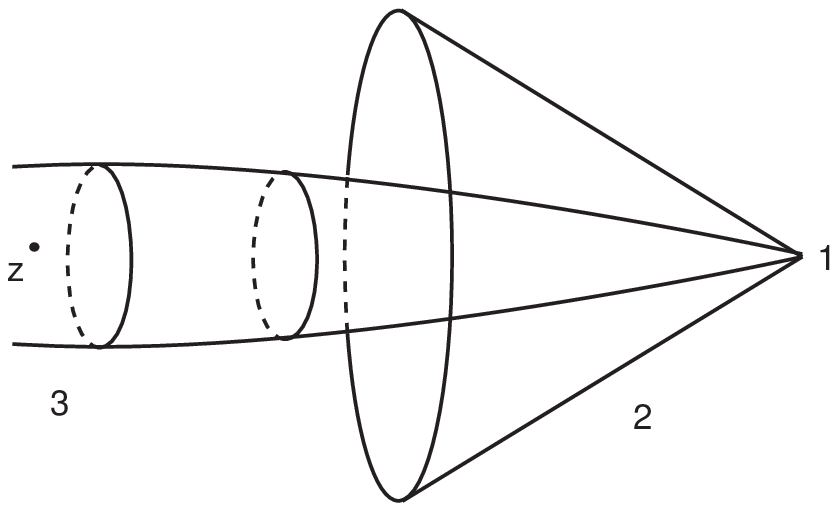}}}
  \relabel{1}{\Large{$\mathcal{E}$}}
  \relabel{z}{$z_\infty$}
  \relabel{2}{$(C_\mathcal{E},g_\mathcal{E})$}
  \relabel{3}{$(U_\infty, g_\infty, z_\infty)$}
  \endrelabelbox
  \caption{Limiting cone.}\label{fig:conelimit}
\end{figure}

\begin{prop}\label{statement}
Let $(U_\infty,g_\infty,z_\infty)$  be as in the conclusion of
Proposition~\ref{Uinfty}. Let $Q_\infty=R_{g_\infty}(z_\infty)$ and
let ${\mathcal E}$ be the end of $U_\infty$ where the scalar
curvature is unbounded. Let $\lambda_n$ be any sequence of positive
numbers with ${\rm lim}_{n\rightarrow\infty}\lambda_n=+\infty$. Then
there is a sequence $x_n$ in $U_\infty$ such that for each $n$ the
distance from $x_n$ to ${\mathcal E}$ is $\lambda_n^{-1/2}$, and
such that the pointed Riemannian manifolds
$(U_\infty,\lambda_ng_\infty,x_n)$ converge in the Gromov-Hausdorff
sense to an open cone, an open cone\index{cone} not homeomorphic to
an open ray (i.e., not homeomorphic to the open cone on a point).
 (see {\sc Fig.}~\ref{fig:conelimit}).
\end{prop}

The rest of this section is devoted to the proof of this result.

\subsection{Directions at ${\mathcal E}$}

 We orient the direction down the tube $U_\infty$ so that
${\mathcal E}$ is at the positive end. This gives an $s_N$-direction
for each $2\epsilon$-neck $N$ contained in $U_\infty$.

Fix a point $x\in U_\infty$. We say a ray $\gamma$ with endpoint $x$
limiting to ${\mathcal E}$ is {\sl a minimizing geodesic ray} if for
every $y\in \gamma$ the segment on $\gamma$ from $x$ to $y$ is a
minimizing geodesic segment; i.e., the length of this geodesic
segment is equal to $d(x,y)$.

\begin{claim}
There is a minimizing geodesic ray to ${\mathcal E}$ from each $x\in
U_\infty$ with $R(x)\ge 2Q_\infty$.
\end{claim}

\begin{proof}
Fix $x$ with $R(x)\ge 2Q_\infty$ and fix a $2\epsilon$-neck $N_x$
centered at $x$. Let $S_x^2$ be the central $2$-sphere of this neck.
Take a sequence of points $q_n$ tending to the end ${\mathcal E}$,
each being closer to the end than $x$ in the sense that $S^2_x$ does
not separate any $q_n$ from the end ${\mathcal E}$. We claim that
there is a minimizing geodesic from $x$ to each $q_n$. The reason is
that by Lemma~\ref{genint} any minimizing sequence of arcs from $x$
to $q_n$ cannot exit from the minus end of $N_x$ nor the plus end of
a $2\epsilon$-neck centered at $q_n$. Consider a sequence of paths
from $x$ to $q_n$ minimizing the distance.  Hence these paths all
lie in a fixed compact subset of $U_n$. After replacing the sequence
by a subsequence, we can pass to a limit, which is clearly a
minimizing geodesic from $x$ to $q_n$. Consider minimizing geodesics
$\mu_n$ from $x$ to $q_n$. The same argument shows that, after
passing to a subsequence, the $\mu_n$ converge to a minimizing
geodesic ray from $x$ to ${\mathcal E}$.
\end{proof}

\begin{claim}
(1) Any minimizing geodesic ray from $x$ to the end ${\mathcal E}$
is a shortest ray from $x$ to the end ${\mathcal E}$, and conversely
any shortest ray from $x$ to the end ${\mathcal E}$ is a minimizing
geodesic ray.

(2) The length of a shortest ray from $x$ to ${\mathcal E}$ is the
distance (see Section~\ref{ENDS}) from $x$ to ${\mathcal E}$.
\end{claim}

\begin{proof}
The implication in (1) in one direction is clear: If $\gamma$ is a
ray from $x$ to the end ${\mathcal E}$, and for some $y\in \gamma$
the segment on $\gamma$ from $x$ to $y$ is not minimizing, then
there is a shorter geodesic segment $\mu$ from $x$ to $y$. The union
of this together with the ray on $\gamma$ from  $y$ to the end is a
shorter ray from $x$ to the end.

Let us establish the opposite implication. Suppose that $\gamma$ is
a minimizing geodesic ray from $x$ to the end ${\mathcal E}$ and
that there is a $\delta>0$ and a shortest geodesic ray $\gamma'$
from $x$ to the end ${\mathcal E}$ with $|\gamma'|=|\gamma|-\delta$.
As we have just seen, $\gamma'$ is a minimizing geodesic ray. Take a
sequence of points $q_i$ tending to the end ${\mathcal E}$ and let
$S^2_i$ be the central $2$-sphere in the $2\epsilon$-neck centered
at $q_i$. Of course, for all $i$ sufficiently large, both $\gamma'$
and $\gamma$ must cross $S^2_i$. Since the scalar curvature tends to
infinity at the end ${\mathcal E}$, it follows from
Lemma~\ref{directions} for all $i$ sufficiently large, the extrinsic
diameter of $S^2_i$ is less than $\delta/3$. Let $p_i$ be a point of
intersection of $\gamma$ with $S^2_i$. For all $i$ sufficiently
large the length $d_i$ of the sub-ray in $\gamma$ from $p_i$ to the
end ${\mathcal E}$ is at most $\delta/3$. Let $p'_i$ be a point of
intersection of $\gamma'$ with $S^2_i$ and let $d'_i$ be the length
of the ray in $\gamma'$ from $p'_i$ to the end ${\mathcal E}$. Let
$\lambda$ be the sub-geodesic of $\gamma$ from $x$ to $p_i$ and
$\lambda'$ the sub-geodesic of $\gamma'$ from $x$ to $p'_i$. Let
$\beta$ be a minimizing geodesic from $p'_i$ to $p_i$. Of course,
$|\beta|<\delta/3$ so that by the minimality of $\lambda$ and
$\lambda'$ we have
$$-\delta/3<|\lambda|-|\lambda'|<\delta/3.$$ Since
$|\lambda'|+d'_i=|\lambda|+d_i-\delta$, we have
$$2\delta/3\le d_i-d'_i.$$
This is absurd since $d'_i>0$ and $d_i<\delta/3$.

(2) follows immediately from (1) and the definition.
\end{proof}

Given this result, the usual arguments show:

\begin{cor}\label{distlim}
If $\gamma$ is a minimizing geodesic ray from $x$ to the end
${\mathcal E}$, then for any $y\in \gamma\setminus\{x\}$ the sub-ray
of $\gamma$ from $y$ to the end, is the unique shortest geodesic
from $y$ to the end.
\end{cor}

Also, we have a version of the triangle inequality for distances to
${\mathcal E}$.

\begin{lem}
Let $x$ and $y$ be points of $M$. Then the three distances $d(x,y)$,
$d(x,{\mathcal E})$ and $d(y,{\mathcal E})$ satisfy the triangle
inequality.
\end{lem}

\begin{proof}
From the definitions it is clear that $d(x,y)+d(y,{\mathcal E})\ge
d(x,{\mathcal E})$, and symmetrically, reversing the roles of $x$
and $y$. The remaining inequality that we must establish is the
following: $d(x,{\mathcal E})+d(y,{\mathcal E})\ge d(x,y)$. Let
$q_n$ be any sequence of points converging to ${\mathcal E}$. Since
the end is at finite distance, it is clear that $d(x,{\mathcal E})=
{\rm lim}_{n\rightarrow\infty}d(x,q_n)$. The remaining inequality
follows from this and  the usual triangle inequality applied to
$d(x,q_n)$, $d(y,q_n)$ and $d(x,y)$.
\end{proof}

\begin{defn}
We say that two minimizing geodesic rays limiting to ${\mathcal E}$
are {\em equivalent} if one is contained in the other. From the
unique continuation of geodesics it is easy to see that this
generates an equivalence relation. An equivalence class is a {\sl
direction at ${\mathcal E}$}, and the set of equivalence classes is
the {\sl set of directions at ${\mathcal E}$}.
\end{defn}

\begin{lem}\label{morethanone}
There is more than one direction at ${\mathcal E}$.
\end{lem}

\begin{proof}
Take a minimal geodesic ray $\gamma$ from a point $x$ limiting to
the end and let $y$ be a point closer to ${\mathcal E}$ than $x$ and
not lying on $\gamma$. Then a minimal geodesic ray from $y$ to
${\mathcal E}$ gives a direction at ${\mathcal E}$ distinct from the
direction determined by $\gamma$.
\end{proof}

\begin{rem}
In fact, the general theory of positively curved spaces implies that
the space of directions is homeomorphic to $S^2$. Since we do not
need this stronger result we do not prove it.
\end{rem}

\subsection{The Metric on the space of directions at ${\mathcal
E}$}

\begin{defn}\label{defntheta}
Let $\gamma$ and $\mu$ be minimizing geodesic rays limiting to
${\mathcal E}$, of lengths $a$ and $b$, parameterized by the
distance from the end. For $0<s\le a$ and $0<s'\le b$ construct a
triangle $\alpha_se\beta_{s'}$ in the Euclidean plane with
$|\alpha_se|=s, |e\beta_{s'}|=s'$ and
$|\alpha_s\beta_{s'}|=d(\gamma(s),\mu(s'))$. We define
$\theta(\gamma,s,\mu,s')$ to be the angle at $e$ of the triangle
$\alpha_se\beta_{s'}$.
\end{defn}

\begin{lem}\label{thetaconv}
For all $\gamma,s,\mu,s'$ as in the previous definition we have
$$0\le \theta(\gamma,s,\mu,s')\le \pi.$$ Furthermore,  $\theta(\gamma,s,\mu,s')$
is a non-increasing function of $s$  when $\gamma,\mu,s'$ are held
fixed, and symmetrically it is a non-increasing function of $s'$
when $\gamma,s,\mu$ are held fixed. In particular, fixing $\gamma$
and $\mu$, the function $\theta(\gamma,s,\mu,s')$ is non-decreasing
as $s$ and $s'$ tend to zero. Thus, there is a well-defined limit as
$s$ and $s'$ go to zero, denoted $ \theta(\gamma,\mu)$. This limit
is greater than or equal to $\theta(\gamma,s,\mu,s')$ for all $s$
and $s'$ for which the latter is defined. We have
$0\le\theta(\gamma,\mu)\le \pi$. The angle $\theta(\gamma,\mu)=0$ if
and only if $\gamma$ and $\mu$ are equivalent. Furthermore, if
$\gamma$ is equivalent to $\gamma'$ and $\mu$ is equivalent to
$\mu'$, then $\theta(\gamma,\mu)=\theta(\gamma',\mu')$.
\end{lem}

\begin{proof}
By restricting $\gamma$ and $\mu$ to slightly smaller rays, we can
assume that each is the unique shortest ray from its endpoint to the
end ${\mathcal E}$. Let $x$, resp., $y$ be the endpoint of $\gamma$,
resp., $\mu$. Now let $q_n$ be any sequence of points in $U_\infty$
limiting to the end ${\mathcal E}$, and consider minimizing geodesic
rays $\gamma_n$ from $q_n$ to $x$ and $\mu_n$ from $q_n$ to $y$,
each parameterized by the distance from $q_n$. By passing to a
subsequence we can assume that each of the sequences $\{\gamma_n\}$
and $\{\mu_n\}$ converge to a minimizing geodesic ray, which  by
uniqueness, implies that the first sequence limits to $\gamma$ and
the second to $\mu$. For $s,s'$ sufficiently small, let
$\theta_n(s,s')$ be the angle at $\tilde q_n$ of the Euclidean
triangle $\alpha_n\tilde q_n\beta_n$, where $|\alpha_n\tilde
q_n|=d(\gamma_n(s),q_n)$, $|\beta_n\tilde q_n|=d(\mu_n(s'),q_n)$ and
$|\alpha_n\beta_n|=d(\gamma_n(s),\mu_n(s'))$. Clearly, for fixed $s$
and $s'$ sufficiently small, $\theta_n(s,s')$ converges as
$n\rightarrow\infty$ to $\theta(\gamma,s,\mu,s')$. By the Toponogov
property (Theorem~\ref{lengthcompar}) for manifolds with
non-negative curvature, for each $n$ the function $\theta_n(s,s')$
is a non-increasing function of each variable, when the other is
held fixed. This property then passes to the limit, giving the first
statement in the lemma.

By the monotonicity, $\theta(\gamma,\mu)=0$ if and only if for all
$s,s'$ sufficiently small we have $\theta(\gamma,s,\mu,s')=0$, which
means one of $\gamma$ and $\mu$ is contained in the other.

It is obvious that the last statement holds.
\end{proof}

It follows that $\theta(\gamma,\mu)$ yields a well-defined function
on the set of pairs of directions at ${\mathcal E}$. It is clearly a
symmetric, non-negative function which is positive off of the
diagonal. The next lemma shows that it is a metric by establishing
the triangle inequality for $\theta$.

\begin{lem}
If $\gamma,\mu,\nu$ are minimizing geodesic rays limiting to
${\mathcal E}$, then
$$\theta(\gamma,\mu)+\theta(\mu,\nu)\ge \theta(\gamma,\nu).$$
\end{lem}

\begin{proof}
By Corollary~\ref{distlim}, after replacing $\gamma, \mu,\nu$ by
equivalent, shorter geodesic arcs, we can assume that they are the
unique minimizing geodesics from their end points, say $x,y,z$
respectively, to ${\mathcal E}$. Let $q_n$ be a sequence of points
limiting to ${\mathcal E}$, and let $\gamma_n,\mu_n,\nu_n$ be
minimizing geodesics from $x,y,z$ to $q_n$. Denote by
$\theta_n(x,y), \theta_n(y,z)$, and $\theta_n(x,z)$, respectively,
the angles  at $\tilde q_n$ of the Euclidean triangles with the
following edge lengths: $\{d(x,y), d(x,q_n), d(y,q_n)\}$,
$\{d(y,z),d(y,q_n),d(z,q_n)\}$, and $\{d(z,x),d(z,q_n),d(x,q_n)\}$.
According to Corollary~\ref{anglecompar} we have
$\theta_n(x,y)+\theta_n(y,z)\ge \theta_n(x,z)$. Passing to the limit
as $n$ goes to $\infty$ and then the limit as $x$, $y$ and $z$ tend
to ${\mathcal E}$, gives the result.
\end{proof}

\begin{defn}
Let $X({\mathcal E})$ denote the set of directions at ${\mathcal
E}$. We define the metric on $X({\mathcal E})$ by setting
$d([\gamma],[\mu])=\theta(\gamma,\mu)$. We call this the {\sl
(metric) space of realized directions at} ${\mathcal E}$. The {\sl
metric space of directions at} ${\mathcal E}$ is the completion
$\bar X({\mathcal E})$ of $X({\mathcal E})$ with respect to the
given metric. We denote by $(C_{\mathcal E},g_{\mathcal E})$ the
cone on $\bar X({\mathcal E})$ with the cone metric as given in
Equation~(\ref{conemetric}). (See {\sc Fig.}~\ref{fig:conelimit}.)
\end{defn}

\begin{prop}\label{notray}
$(C_{\mathcal E},g_{\mathcal E})$ is a metric cone that is not
homeomorphic to a ray.
\end{prop}

\begin{proof}
By construction $(C_{\mathcal E},g_{\mathcal E})$ is a metric cone.
That it is not homeomorphic to a ray follows immediately from
Lemma~\ref{morethanone}.
\end{proof}

\subsection{Comparison results for distances}

\begin{lem}\label{distcomp}
Suppose that $\gamma$ and $\mu$ are unique shortest geodesic rays
from points $x$ and $y$ to the end ${\mathcal E}$. Let $[\gamma]$
and $[\mu]$ be the points of $X({\mathcal E})$ represented by these
two geodesics rays. Let $a$, resp. $b$, be the distance from $x$,
resp. $y$, to ${\mathcal E}$. Denote by  $x'$, resp. $y'$, the image
in $C_{\mathcal E}$ of the point $([\gamma],a)$, resp. $([\mu],b)$,
of $X({\mathcal E})\times [0,\infty)$. Then
$$d_{g_\infty}(x,y)\le d_{g_{\mathcal E}}(x',y').$$
\end{lem}

\begin{proof}
By the definition of the cone metric we have
$$d_{g_{\mathcal E}}(x',y')=a^2+b^2-2ab\,{\rm cos}(\theta(\gamma,\mu)).$$
On the other hand by Definition~\ref{defntheta} and the law of
cosines for Euclidean triangles, we have
$$d_{g_\infty}(x,y)=a^2+b^2-2ab\,{\rm cos}(\theta(\gamma,a,\mu,b)).$$
The result is now immediate from the fact, proved in
Lemma~\ref{defntheta} that
$$0\le \theta(\gamma,a,\mu,b)\le \theta(\gamma,\mu)\le \pi,$$
and the fact that the cosine is a monotone decreasing function on
the interval $[0,\pi]$.
\end{proof}

\begin{cor}\label{refdistcomp}
Let $\gamma,\mu,x,y$ be as in the previous lemma. Fix $\lambda>0$.
Let $a=d_{\lambda g_\infty}(x,{\mathcal E})$ and $b=d_{\lambda
g_\infty}(y,{\mathcal E})$. Set $x'_\lambda$ and $y'_\lambda$ equal
to the points in the cone $([\gamma],a)$ and $([\mu],b)$. Then we
have
$$d_{\lambda g_\infty}(x,y)\le d_{g_{\mathcal E}}(x'_\lambda,y'_\lambda).$$
\end{cor}

\begin{proof}
This is immediate by applying the previous lemma to the rescaled
manifold $(U_\infty,\lambda g_\infty)$, and noticing that rescaling
does not affect the cone $C_{\mathcal E}$ nor its metric.
\end{proof}

\begin{lem}
For any $\delta>0$ there is $K=K(\delta)<\infty$ so that for any set
of realized directions at ${\mathcal E}$ of cardinality $K$,
$\ell_1,\ldots,\ell_{K}$, it must be the case that there are $j$ and
$j'$ with $j\not=j'$ such that $\theta(\ell_j,\ell_{j'})<\delta$.
\end{lem}

\begin{proof}
Let $K$ be such that, given $K$ points in the central $2$-sphere of
any $2\epsilon$-tube of scale $1$, at least two are within distance
$\delta/2$ of each other. Now suppose that we have $K$ directions
$\ell_1,\ldots,\ell_K$ at ${\mathcal E}$. Let
$\gamma_1,\ldots,\gamma_K$ be minimizing geodesic rays limiting to
${\mathcal E}$ that represent these directions. Choose a point $x$
sufficiently close to the end ${\mathcal E}$ so that all the
$\gamma_j$ cross the central $2$-sphere $S^2$ of the
$2\epsilon$-neck centered at $x$. By replacing the $\gamma_j$ with
sub-rays we can assume that for each $j$ the endpoint $x_j$ of
$\gamma_j$ lies in $S^2$. Let $d_j$ be the length of $\gamma_j$. By
taking $x$ sufficiently close to ${\mathcal E}$ we can also assume
the following. For each $j$ and $j'$, the angle at $e$ of the
Euclidean triangle $\alpha_je\alpha_{j'}$, where $|\alpha_je|=d_j;
|\alpha_{j'}e|=d_{j'}$ and $|\alpha_j\alpha_{j'}|=d(x_j,x_{j'})$ is
within $\delta/2$ of $\theta(\ell_j,\ell_{j'})$. Now there must be
$j\not= j'$ with $d(x_i,x_j)<(\delta/2)r_i$ where $r_i$ is the scale
of $N_i$. Since $d_j,d_{j'}>\epsilon^{-1}r_i/2$, it follows that the
angle at $e$ of $\alpha_je\alpha_{j'}$ is less than $\delta/2$.
Consequently, $\theta(\ell_j,\ell_{j'})<\delta$.
\end{proof}

Recall that a $\delta$-net in a metric space $X$ is a finite set of
points such that $X$ is contained in the union of the
$\delta$-neighborhoods of these points. The above lemma immediately
yields:

\begin{cor}\label{deltanet}
The metric completion $\bar X({\mathcal E})$ of the space of
directions at ${\mathcal E}$ is a compact space. For every
$\delta>0$ this space has a  $\delta$-net consisting of realized
directions. For every $0<r<R<\infty$ the annular region $A_{\mathcal
E}(r,R)=\bar{X}({\mathcal E})\times [r,R]$ in $C_{\mathcal E}$ has a
$\delta$-net consisting of points $(\ell_i,s_i)$ where for each $i$
we have $\ell_i$ is a realizable direction and $r<s_i<R$.
\end{cor}

\subsection{Completion of the proof of a cone limit at ${\mathcal E}$}

Now we are ready to prove Proposition~\ref{statement}. In fact, we
prove a version of the proposition that identifies the sequence of
points $x_n$ and also identifies the cone\index{cone} to which the
rescaled manifolds converge.

\begin{prop}\label{secondstatment}
Let $(U_\infty,g_\infty)$  be an incomplete Riemannian $3$-manifold
of non-negative curvature with an end ${\mathcal E}$ as in the
hypothesis of Proposition~\ref{statement}. Fix a minimizing geodesic
ray $\gamma$ limiting to ${\mathcal E}$. Let $\lambda_n$ be any
sequence of positive numbers tending to infinity. For each $n$
sufficiently large let $x_n\in \gamma$ be the point at distance
$\lambda_n^{-1/2}$ from the end ${\mathcal E}$. Then the based
metric spaces $(U_\infty,\lambda_ng_\infty,x_n)$ converge in the
Gromov-Hausdorff sense\index{converge!Gromov-Hausdorff sense} to
$\left( C'_{\mathcal E},g_{\mathcal E},([\gamma],1)\right)$. Under
this convergence the distance function from the end ${\mathcal E}$
in $(U_\infty,\lambda_ng_\infty)$ converges to the distance function
from the cone point in the open cone.
\end{prop}

\begin{proof}
It suffices to prove that given any subsequence of the original
sequence, the result holds for a further subsequence. So let us
replace the given sequence by a subsequence. Recall that for each
$0<r<R<\infty$ we have $A_{\mathcal E}(r,R)\subset C'_{\mathcal E}$,
the compact annulus which is the image of $\bar X({\mathcal
E})\times [r,R]$. The statement about the non-compact spaces
converging in the Gromov-Hausdorff topology, means that for each
compact subspace $K$ of $ C'_{\mathcal E}$ containing the base
point, for all $n$ sufficiently large, there are compact subspaces
$K_n\subset (U_\infty,\lambda_ng_\infty)$ containing $x_n$ with the
property that the $(K_n,x_n)$ converge in the Gromov-Hausdorff
topology to $(K,x)$ (see Section D of Chapter 3, p. 39, of
\cite{Gromov}).

Because of this, it suffices to  fix $0<r<1<R<\infty$ arbitrarily
and prove the convergence result for $A_{\mathcal E}(r,R)$. Since
the Gromov-Hausdorff distance from a compact pointed metric space to
a $\delta$-net in it containing the base point is at most $\delta$,
it suffices to prove that for $\delta>0$ there is a $\delta$-net
$({\mathcal N},h)$ in $A_{\mathcal E}(r,R)$, with $([\gamma],1)\in
{\mathcal N}$ such that for all $n$ sufficiently large there are
embeddings $\varphi_n$ of ${\mathcal N}$ into $A_n(r,R)=\bar
B_{\lambda_ng_\infty}({\mathcal E},R)\setminus
B_{\lambda_ng_\infty}({\mathcal E},r)$  with the following four
properties:
\begin{enumerate}
\item[(1)]  $\varphi_n^*(\lambda_ng_\infty)$ converge to $h$ as $n\rightarrow\infty$,
\item[(2)] $\varphi_n([\gamma],1)=x_n$,
and
\item[(3)]  $\varphi_n({\mathcal N})$ is a $\delta$-net in $A_n(r,R)$, and
\item[(4)] denoting the cone point by $c\in C_{\mathcal E}$,
if $d(p,c)=r$ then $d(\varphi_n(p),{\mathcal E})=r$.
\end{enumerate}

According to Corollary~\ref{deltanet} there is a $\delta$-net
${\mathcal N}\subset A_{\mathcal E}(r,R)$ consisting of points
$(\ell_i,s_i)$ where the $\ell_i$ are realizable directions and
$r<s_i< R$. Add $([\gamma],1)$ to ${\mathcal N}$ if necessary so
that we can assume that $([\gamma],1)\in {\mathcal N}$. Let
$\gamma_i$ be a minimizing geodesic realizing $\ell_i$ and let $d_i$
be its length.

Fix $n$ sufficiently large so that $\lambda_n^{-1/2}R\le d_i$ for all $i$. We
define $\varphi_n\colon {\mathcal N}\to A_n(r,R)$ as follows. For any
$a_i=([\gamma_i],s_i)\in {\mathcal N}$ we let
$\varphi_n(a_i)=\gamma_i(\lambda_n^{-1/2}s_i)$. (Since $\lambda_n^{-1/2}s\le
\lambda_n^{-1/2}R\le d_i$, the geodesic $\gamma_i$ is defined at
$\lambda_n^{-1/2}s_i$.) This defines the embeddings $\varphi_n$ for all $n$
sufficiently large. Notice that
$$d_{g_\infty}\left(\varphi_n(\ell_i,s_i),\varphi_n(\ell_j,s_j)\right)=\lambda_n^{-1}s_i^2+\lambda_n^{-1}
s_j^2-2\lambda_n^{-1}s_is_j\theta(\gamma_i,\lambda_n^{-1/2}s_i,\gamma_j,\lambda_n^{-1/2}s_j),$$
or equivalently
$$d_{\lambda_ng_\infty}\left(\varphi_n(\ell_i,s_i),\varphi_n(\ell_j,s_j)\right)=s_i^2+s_j^2-2s_is_j
\theta(\gamma_i,\lambda_n^{-1/2}s_i,\gamma_j,\lambda_n^{-1/2}s_j).$$
Because of the convergence result on angles (Lemma~\ref{thetaconv}),
for all $i$ and $j$ we have \begin{eqnarray*} {\rm
lim}_{n\rightarrow\infty}d_{\lambda_ng_\infty}\left(\varphi_n(\ell_i,s_i),\varphi_n(\ell_j,s_j)\right)
 & = & s_i^2+s_j^2-2s_is_j{\rm
cos}(\theta(\gamma_i,\gamma_j)) \\ & = & d_{g_{\mathcal
E}}\left((\ell_i,s_i),(\ell_j,s_j)\right).\end{eqnarray*}
 This
establishes the existence of the $\varphi_n$ for all $n$
sufficiently large satisfying the first condition. Clearly, from the
definition $\varphi_n([\gamma],1)=x_n$, and for all $p\in {\mathcal
N}$ we have $d(\varphi_n(p),{\mathcal E})=d(p,c)$.

It remains to check that for all $n$ sufficiently large
$\varphi_n({\mathcal N})$ is a $\delta$-net in $A_n(r,R)$. For $n$
sufficiently large let $z\in A_n(r,R)$ and let $\gamma_z$ be a
minimizing geodesic ray from $z$ to ${\mathcal E}$ parameterized by
the distance from the end. Set
$d_n=d_{\lambda_ng_\infty}(z,{\mathcal E})$, so that $r\le d_n\le
R$. Fix $n$ sufficiently large so that $\lambda_n^{-1/2}R<d_i$ for
all $i$. The point $([\gamma_z],d_n)\in C_{\mathcal E}$ is contained
in $A_{\mathcal E}(r,R)$ and hence there is an element
$a=([\gamma_i],s_i)\in {\mathcal N}$ within distance $\delta$ of
$([\gamma_z],d_n)$ in $C_{\mathcal E}$. Since $s_i\le R$,
$\lambda_n^{-1/2}s_i\le d_i$ and hence
$x=\gamma_i(\lambda_n^{-1/2}s_i)$ is defined. By
Corollary~\ref{refdistcomp} we have
$$d_{\lambda_ng_\infty}(x,z)\le
d_{g_{\mathcal E}}\left(([\gamma],d_n),([\gamma_i],s_i)\right)\le
\delta.$$ This completes the proof that for $n$ sufficiently large
the image $\varphi_n({\mathcal N})$ is a $\delta$-net in $A_n(r,R)$.

This shows that the $(U_\infty,\lambda_ng_\infty,x_n)$ converge in
the Gromov-Hausdorff topology to $(C'_{\mathcal E},g_{\mathcal
E},([\gamma],1))$.
\end{proof}

\begin{rem}
Notice that since the manifolds $(U_\infty,\lambda_ng_\infty,x_n)$
are not complete, there can be more than one Gromov-Hausdorff limit.
For example we could take the full cone as a limit. Indeed, the cone
is the only Gromov-Hausdorff limit that is complete as a metric
space.
\end{rem}

\section{Comparison of the Gromov-Hausdorff limit and the smooth
limit}

Let us recap the progress to date. We constructed an incomplete
geometric blow-up limit $({\mathcal U}_\infty,G_\infty,z_\infty)$
for our original sequence. It has non-negative Riemann curvature. We
showed that the zero time-slice $U_\infty$ of the limit is
diffeomorphic to a tube $S^2\times (0,1)$ and that at one end of the
tube the scalar curvature goes to infinity. Also, any point
sufficiently near this end is the center of an evolving
$2\epsilon$-neck defined for an interval of normalized  time of
length $1/2$ in the limiting flow. Then we took a further blow-up
limit. We chose a sequence of points $x_n\in U_\infty$ tending to
the end ${\mathcal E}$ where the scalar curvature goes to infinity.
Then we formed $(U_\infty,\lambda_ng_\infty,x_n)$ where the distance
from $x_n$ to the end ${\mathcal E}$ is $\lambda_n^{-1/2}$. By
fairly general principles (in fact it is a general theorem about
manifolds of non-negative curvature) we showed that this sequence
converges in the Gromov-Hausdorff
sense\index{converge!Gromov-Hausdorff sense} to a cone.

The next step is to show that this second blow-up limit also exists
as a geometric limit away from the cone point. Take a sequence of
points $x_n\in U_\infty$ tending to ${\mathcal E}$. We let
$\lambda'_n=R(x_n)$, and we consider the based Riemannian manifolds
$(U_\infty,\lambda_n'g_\infty(0),x_n)$. Let $B_n\subset U_n$ be the
metric ball of radius $\epsilon^{-1}/3$ centered at $x_n$ in
$(U_\infty,\lambda'_ng_\infty(0))$. Since this ball is contained in
a $2\epsilon$-neck centered at $x_n$, the curvature on this ball is
bounded, and this ball has compact closure in $U_\infty$. Also, for
each $y\in B_n$, there is a rescaled flow $\lambda'g(t)$  defined on
$\{y\}\times (-1/2,0]$ whose  curvature on $B_n\times (-1/2,0]$ is
bounded. Hence, by Theorem~\ref{sliceconv} we can pass to a
subsequence and extract a geometric limit. In fact, by
Proposition~\ref{partialflowlimit} there is even a geometric
limiting flow defined on the time interval $(-1/2,0]$.

We must compare the zero time-slice of this geometric limiting flow
with the corresponding open subset of the Gromov-Hausdorff
limit\index{Gromov-Hausdorff limit} constructed in the previous
section. Of course, one obvious difference is that we have used
different blow-up factors: $d(x_n,{\mathcal E})^{-2}$ in the first
case and $R(x_n)$ in the second case. So one important ingredient in
comparing the limits will be to compare these factors, at least in
the limit.

\subsection{Comparison of the blow-up factors}

Now let us compare the two limits: (i) the Gromov-Hausdorff limit of
the sequence $(U_\infty,\lambda_ng_\infty,x_n)$ and (ii) the
geometric limit of the sequence $(U_\infty,\lambda_n'g_\infty,x_n)$
 constructed above.

\begin{claim}
The ratio $\rho_n=\lambda_n'/\lambda_n$ is bounded above and below
by positive constants.
\end{claim}

\begin{proof}
First of all, since there is a $2\epsilon$-neck centered at $x_n$,
by Proposition~\ref{S2intersection} we see that the distance
$\lambda_n^{-1/2}$ from $x_n$ to ${\mathcal E}$ is at least
$R(x_n)^{-1/2}\epsilon^{-1}/2=(\lambda'_n)^{-1/2}\epsilon^{-1}/2$.
Thus,
$$\rho_n^{-1}=\lambda_n/\lambda'_n\le 4\epsilon^2.$$

On the other hand, suppose that
$\rho_n=\lambda_n'/\lambda_n\rightarrow \infty$ as $n\rightarrow
\infty$. Rescale by $\lambda'_n$ so that $R(x_n)=1$. The distance
from $x_n$ to ${\mathcal E}$ is $\sqrt{\rho_n}$. Then by
Lemma~\ref{directions} with respect to this metric there is a sphere
of diameter at most $2\pi$  through  $x_n$ that separates all points
at distance at most $\sqrt{\rho_n}-\epsilon^{-1}$ from ${\mathcal
E}$ from all points at distance at least
$\sqrt{\rho_n}+\epsilon^{-1}$ from ${\mathcal E}$. Now rescale the
metric by $\rho_n$. In the rescaled metric there is a $2$-sphere of
diameter at most $2\pi/\sqrt{\rho_n}$ through $x_n$ that separates
all points at distance at most $1-\epsilon^{-1}/\sqrt{\rho_n}$ from
${\mathcal E}$ from all points at distance at least
$1+\epsilon^{-1}/\sqrt{\rho_n}$ from ${\mathcal E}$. Taking the
Gromov-Hausdorff limit of these spaces, we see that the base point
$x_\infty$ separates all points of distance less than one from
${\mathcal E}$ from all points of distance greater than one from
${\mathcal E}$. This is impossible since the Gromov-Hausdorff limit
is a cone that is not the cone on a single point.
\end{proof}

\subsection{Completion of the comparison of the blow-up limits}

Once we know that the $\lambda_n/\lambda_n'$ are bounded above and
below by positive constants, we can pass to a subsequence so that
these ratios converge to a finite positive limit. This means that
the Gromov-Hausdorff limit of the sequence of based metric spaces
$(U_\infty,\lambda_n'g_\infty,x_n)$ is a  cone, namely the
Gromov-Hausdorff limiting cone constructed is Section~\ref{GromHaus}
rescaled by ${\rm lim}_{n\rightarrow\infty}\rho_n$. In particular,
the balls of radius $\epsilon^{-1}/2$ around the base points in this
sequence converge in the Gromov-Hausdorff sense to the ball of
radius $\epsilon^{-1}/2$ about the base point of a cone.

But we have already seen that the balls of radius $\epsilon^{-1}/2$
centered at the base points converge geometrically to a limiting
manifold. That is to say, on every   ball of radius less than
$\epsilon^{-1}/2$ centered at the base point the metrics converge
uniformly in the $C^\infty$-topology to a limiting smooth metric.
Thus, on every ball of radius less than $\epsilon^{-1}/2$ centered
at the base point the limiting smooth metric is isometric to the
metric of the Gromov-Hausdorff limit. This means that the limiting
smooth metric on the ball $B_\infty$ of radius $\epsilon^{-1}/2$
centered at the base point is isometric to an open subset of a cone.
Notice that the scalar curvature of the limiting smooth metric at
the base point is $1$, so that this cone is a non-flat cone.

\section{The final contradiction}

We have now shown that the smooth limit of the balls of radius
$\epsilon^{-1}/2$ centered at the base points of
$(U_\infty,\lambda_n'g_\infty,x_n)$ is isometric to an open subset
of a non-flat cone\index{cone}, and is also the zero time-slice of a
Ricci flow defined for the time interval $(-1/2,0]$. This
contradicts Proposition~\ref{nocones},  one of the consequences of
the maximum principle established by Hamilton. The contradiction
shows that the limit $(U_\infty,g_\infty,x_\infty)$ cannot exist.
The only assumption that we made in order to construct this limit
was that Theorem~\ref{bcbd} did not hold for some $A_0<\infty$.
Thus, we have established Theorem~\ref{bcbd}  by contradiction.

\chapter{Geometric limits of generalized Ricci flows}

In this chapter we apply the main result of the last section,
bounded curvature at bounded distance, to blow-up limits in order to
establish the existence of a smooth limit for sequences of
generalized Ricci flows. In the first section we establish a blow-up
limit that is defined for some interval of time of positive length,
where the length of the interval of time is allowed to depend on the
limit. In the second section we give conditions under which this
blow-up limit can be extended backwards to make an ancient Ricci
flow. In the third section we construct limits at the singular time
of a generalized Ricci flow satisfying appropriate conditions. We
characterize the ends of the components of these limits. We show
that they are $\epsilon$-horns -- the ends are diffeomorphic to
$S^2\times [0,1)$ and the scalar curvature goes to infinity at the
end. In the fourth section we prove  for any $\delta>0$  that there
are $\delta$-necks sufficiently deep in any $\epsilon$-horn,
provided that the curvature at the other end of the horn is not too
large. Throughout this chapter we fix $\epsilon>0$ sufficiently
small such that all the results of the Appendix hold for $2\epsilon$
and $\alpha=10^{-2}$, and Proposition~\ref{narrows} holds for
$2\epsilon$.

\section{A smooth blow-up limit defined for a small time}

We begin with a theorem that produces a blow-up limit\index{blow-up
limit} flow defined on some small time interval.

\begin{thm}\label{smlmtflow}
Fix  canonical neighborhood constants $(C,\epsilon)$, and
non-collapsing constants $r>0,\kappa>0$. Let $({\mathcal
M}_n,G_n,x_n)$ be a sequence of based generalized $3$-dimensional
Ricci flows. We set $t_n={\bf t}(x_n)$ and $Q_n=R(x_n)$. We denote
by $M_n$ the $t_n$ time-slice of ${\mathcal M}_n$. We suppose that:
\begin{enumerate}
\item[(1)] Each $({\mathcal M}_n,G_n)$ either has a time interval of definition contained in $[0,\infty)$
and has  curvature pinched toward positive\index{curvature!pinched
toward positive}, or has  non-negative curvature.
\item[(2)] Every point $y_n\in ({\mathcal M}_n,G_n)$ with ${\bf t}(y_n)\le
t_n$ and with $R(y_n)\ge 4R(x_n)$ has a strong
$(C,\epsilon)$-canonical neighborhood\index{canonical
neighborhood!strong}.
\item[(3)] ${\rm lim}_{n\rightarrow\infty}Q_n=\infty$.
\item[(4)] For each $A<\infty$ the following holds for all $n$ sufficiently large. The ball
$B(x_n,t_n,AQ_n^{-1/2})$ has compact closure in $M_n$ and the flow
is $\kappa$-non-collapsed on scales $\le r$ at each point of
$B(x_n,t_n,AQ_n^{-1/2})$.
\item[(5)] There is $\mu>0$ such that for every $A<\infty$ the following holds
 for all $n$ sufficiently large.
For every $y_n\in B(x_n,t_n,AQ_n^{-1/2})$ the maximal flow line
through $y_n$ extends backwards for a time at least $\mu\left({\rm
max}(Q_n,R(y_n))\right)^{-1}$.
\end{enumerate}

Then, after passing to a subsequence and shifting the times of each
of the generalized flows so that $t_n=0$ for every $n$, there is a
geometric limit $(M_\infty,g_\infty,x_\infty)$ of the sequence of
based Riemannian manifolds $(M_n,Q_nG_n(0),x_n)$. This limit is a
complete $3$-dimensional Riemannian manifold of bounded,
non-negative curvature. Furthermore, for some $t_0>0$ which depends
on the curvature bound for $(M_\infty,g_\infty)$ and on $\mu$, there
is a geometric limit Ricci flow defined on $(M_\infty,g_\infty(t)),
-t_0\le t\le 0$, with $g_\infty(0)=g_\infty$.
\end{thm}

Before beginning the proof of this theorem we establish a lemma that
we shall need both in its proof and also for later applications.

\begin{lem}\label{controlnbhd}
Let $({\mathcal M},G)$ be a generalized $3$-dimensional Ricci flow.
Suppose that $r_0>0$ and that any $z\in {\mathcal M}$ with $R(z)\ge
r_0^{-2}$ has a strong $(C,\epsilon)$-canonical
neighborhood\index{canonical neighborhood!strong}. Suppose $z\in
{\mathcal M}$ and ${\bf t}(z)=t_0$. Set
$$r=\frac{1}{2C\sqrt{{\rm max}(R(z),r_0^{-2})}}$$ and
$$\Delta t=\frac{1}{16C\left(R(z)+r_0^{-2}\right)}.$$
Suppose that $r'\le r$ and that $|t'-t_0|\le \Delta t$ and let $I$
be the interval with endpoints $t_0$ and $t'$. Suppose that there is
an embedding of $j\colon B(z,t_0,r')\times I$ into ${\mathcal M}$
compatible with time and with the vector field. Then $R(y)\le
2\left(R(z)+r_0^{-2}\right)$ for all $y$ in the image of $j$.
\end{lem}

\begin{proof}
We first  prove that for any $y\in B(z,t_0,r)$ we have
\begin{equation}\label{Rineq}
R(y)\le \frac{16}{9}(R(z)+r_0^{-2}).\end{equation} Let $\gamma\colon
[0,s_0]\to B(z,t_0,r)$ be a path of length $s_0<r$ connecting
$z=\gamma(0)$ to $y=\gamma(s_0)$. We take $\gamma$ parameterized by
arc length. For any $s\in [0,s_0]$ let $R(s)=R(\gamma(s))$.
According to the strong $(C,\epsilon)$-canonical neighborhood
assumption at any point where $R(s)\ge r_0^{-2}$ we have $|R'(s)|\le
CR^{3/2}(s)$. Let $J\subset [0,s_0]$ be the closed subset consisting
of $s\in [0,s_0]$ for which $R(s)\ge r_0^{-2}$. There are three
possibilities. If $s_0\not\in J$ then $R(y)\le r_0^{-2}$ and we have
established Inequality~(\ref{Rineq}). If $J=[0,s_0]$, then we have
$|R'(s)|\le CR^{3/2}(s)$ for all $s$ in $J$. Using this differential
inequality and the fact that the interval has length at most
$\frac{1}{2C\sqrt{R(z)}}$, we see that $R(y)\le 16R(z)/9$, again
establishing Inequality~(\ref{Rineq}). The last possibility is that
$J\not= [0,s_0]$ but $s_0\in J$. We restrict attention to the
maximal interval of $J$ containing $s_0$. This interval has length
at most $\frac{r_0}{2C}$ and at its initial point $R$ takes the
value $r_0^{-2}$. For every $s$ in this interval by our assumptions
we again have the inequality $|R'(s)|\le CR^{3/2}(s)$, it follows
immediately that $R(y)\le 16r_0^{-2}/9$. This establishes
Inequality~(\ref{Rineq}) in all cases.

Now consider the vertical path $j(\{y\}\times I)$. Let
$R(t)=R(j(y,t))$. Again by the strong canonical neighborhood
assumption $|R'(t)|\le CR^2(t)$ at all points where $R(t)\ge
r_0^{-2}$. Consider the closed subset $K$ of $I$ where $R(t)\ge
r_0^{-2}$. There are three cases to consider: $t'\not\in K$, $t'\in
K\not=I$, or $K=I$.  In the first case, $R(y,t')\le r_0^{-2}$ and we
have established the result. In the second case, let $K'$ be the
maximal subinterval of $K$ containing $t'$. On the interval $K'$ we
have $|R'(t)|\le CR^2(t)$ and at one endpoint $R(t)=r_0^{-2}$. Since
this interval has length at most $r_0^2/16C$, it follows easily that
$R(t')\le 16r_0^{-2}/15$, establishing the result. In the last case
where $K=I$, then, by what we established above, the initial
condition is $R(t_0)=R(y)\le 16(R(z)+r_0^{-2})/9$, and the
differential inequality $|R'(t)|\le CR^2(t)$ holds for all $t\in I$.
Since the length of $I$ is at most $\frac{1}{16C(R(y)+r_0^{-2})}$ we
see directly that $R(t')\le 2(R(z)+r_0^{-2})$, completing the proof
in this case as well.
\end{proof}

Now we begin the proof of Theorem~\ref{smlmtflow}.

\begin{proof}(of Theorem~\ref{smlmtflow})
We shift the times for the flows so that $t_n=0$ for all $n$. Since
$Q_n$ tends to $\infty$ as $n$ tends to $\infty$, according to
Theorem~\ref{bcbd} for any $A<\infty$, there is a bound
$Q(A)<\infty$ on the scalar curvature of $Q_nG_n(0)$ on
$B_{Q_nG_n}(x_n,0,A)$  for all $n$ sufficiently large. According to
the hypothesis of Theorem~\ref{smlmtflow}, this means that there is
$t_0(A)>0$ and, for each $n$ sufficiently large, an embedding of
$B_{Q_nG_n}(x_n,0,A)\times [-t_0(A),0]$ into ${\mathcal M}_n$
compatible with time and with the vector field. In fact, we can
choose $t_0(A)$ so that more is true.

\begin{cor}\label{cnnbhd2}
For each $A<\infty$, let $Q(A)$ be a bound  on the  scalar curvature
of the restriction of $Q_nG_n$ to  $B_{Q_nG_n}(x_n,0,A)$ for all $n$
sufficiently large. Then there exist a constant  $t'_0(A)>0$
depending on $t_0(A)$ and $Q(A)$, and a constant $Q'(A)<\infty$
depending only on $Q(A)$, and, for all $n$ sufficiently large, an
embedding
$$B_{Q_nG_n}(x_n,0,A)\times (-t'_0(A),0]\to{\mathcal
M}_n$$ compatible with time and with the vector field with the
property that  the scalar curvature of the restriction of $Q_nG_n$
to the image of this subset is bounded by $Q'(A)$.
\end{cor}

\begin{proof}
This is immediate from  Lemma~\ref{controlnbhd} and Assumption (5)
in the hypothesis of the theorem.
\end{proof}

Now since the curvatures of the $Q_nG_n$ are pinched toward positive
or are non-negative, bounding the scalar curvature above gives a
bound on $|{\rm Rm}_{Q_nG_n}|$ on the product
$B_{Q_nG_n}(x_n,0,A)\times (-t'_0(A),0]$. Now we invoke Shi's
theorem (Theorem~\ref{shi}):

\begin{cor}
For each $A<\infty$ and for each integer $\ell\ge 0$, there is a
constant $C_2$ such that for all $n$ sufficiently large we have
$$|\nabla^\ell{\rm {\rm Rm}}_{Q_nG_n}(x)|\le C_2$$ for all $x\in
B_{Q_nG_n}(x_n,0,A) $.
\end{cor}

Also, by the curvature bound and  the $\kappa$-non-collapsed
hypothesis we have the following:

\begin{claim}
There is $\eta>0$ such that for all $n$ sufficiently large
$${\rm Vol}(B_{Q_nG_n}(x_n,0,\eta))\ge \kappa\eta^3.$$
\end{claim}

Now we are in a position to apply Corollary~\ref{2ndmfdconv}. This
implies that, after passing to a subsequence, there is a geometric
limit $(M_\infty,g_\infty,x_\infty)$ of the sequence of based
Riemannian manifolds $(M_n,Q_nG_n(0),x_n)$. The geometric limit is a
complete Riemannian manifold. If the $({\mathcal M}_n,G_n)$ satisfy
the curvature pinched toward positive hypothesis, by
Theorem~\ref{blowupposcurv}, the limit Riemannian manifold
$(M_\infty,g_\infty)$ has non-negative curvature. If the $({\mathcal
M}_n,G_n)$ have non-negative curvature, then it is obvious that the
limit has non-negative curvature. By construction $R(x_\infty)=1$.

In fact, by Proposition~\ref{partialflowlimit} for each $A<\infty$,
there is $t(A)>0$ and, after passing to a subsequence, geometrically
limit flow defined on $B(x_\infty,0,A)\times (-t(A),0]$.

\begin{claim}
Any point in $(M_\infty,g_\infty)$ of curvature greater than $4$ has
a $(2C,2\epsilon)$-canonical neighborhood\index{canonical
neighborhood}.
\end{claim}

\begin{proof}
The fact that $(M_\infty,g_\infty,x_\infty)$ is the geometric limit
of the $(M_n,Q_nG_n(0),x_n)$ means that we have the following. There
is an exhausting sequence  $V_1\subset V_2\subset\cdots \subset
M_\infty$ of open subsets of $M_\infty$, with compact closure, each
containing $x_\infty$, and for each $n$  an embedding $\varphi_n$ of
$V_n$ into the zero time-slice of ${\mathcal M}_n$ such that
$\varphi_n(x_\infty)=x_n$ and such that the Riemannian metrics
$\varphi_n^*G_n$ converge uniformly on compact sets to $g_\infty$.
Let $q\in M_\infty$ be a point with $R_{g_\infty}(q)>4$. Then for
all $n$ sufficiently large, $q\in V_n$, so that $q_n=\varphi_n(q)$
is defined, and $R_{Q_nG_n}(q_n)>4$. Thus, $q_n$ has an
$(C,\epsilon)$-canonical neighborhood, $U_n$, in ${\mathcal M}_n$;
and, since $R(q_n)>4$ for all $n$,  there is a uniform bound to the
distance from any point of $U_n$ to $q_n$. Thus, there exists $m$
such that for all $n$ sufficiently large $\varphi_n(V_m)$ contains
$U_n$. Clearly as $n$ goes to infinity the Riemannian metrics
$\varphi_n^*(G_n)|_{\varphi_n^{-1}(U_m)}$ converge smoothly to
$g_\infty|_{\varphi_n^{-1}(U_n)}$. Thus, by
Proposition~\ref{canonvary} for all $n$ sufficiently large the
restriction of $g_\infty$ to $\varphi_n^{-1}(U_n)$ contains a
$(2C,2\epsilon)$-canonical neighborhood of $q$.
\end{proof}

\begin{claim}
The limit Riemannian manifold $(M_\infty,g_\infty)$ has bounded
curvature.
\end{claim}

\begin{proof}
First, suppose that $(M_\infty,g_\infty)$ does not have strictly
positive curvature. Suppose that $y\in M_\infty$ has the property
that ${\rm Rm}(y)$ has a zero eigenvalue. Fix $A<\infty$ greater
than $d_{g_\infty}(x_\infty,y)$. Then applying
Corollary~\ref{localprod} to the limit flow on
$B(x_\infty,0,A)\times (-t(A),0]$,  we see that the Riemannian
manifold $(B(x_\infty,0,A),g_\infty)$  is locally a Riemannian
product of a compact surface of positive curvature with a
one-manifold. Since this is true for every $A<\infty$ sufficiently
large, the same is true for $(M_\infty,g_\infty)$. Hence
$(M_\infty,g_\infty)$ has a one- or two-sheeted covering that is a
global Riemannian product of a compact surface and one-manifold.
Clearly, in this case the curvature of $(M_\infty,g_\infty)$ is
bounded.

If $M_\infty$ is compact, then it is clear that the curvature is
bounded.

It remains to consider the  case where $(M_\infty,g_\infty)$ is
non-compact and of strictly positive curvature. Since any point of
curvature greater than $4$ has a $(2C,2\epsilon)$-canonical
neighborhood, and since $M_\infty$ is non-compact, it follows that
the only possible canonical neighborhoods for $x\in M_\infty$ are  a
$2\epsilon$-neck centered at $x$  or $(2C,2\epsilon)$-cap whose core
contains $x$.  Each of these canonical neighborhoods contains a
$2\epsilon$-neck. Thus, if $(M_\infty,g_\infty)$ has unbounded and
positive Riemann curvature or equivalently, it has  unbounded scalar
curvature, then it has $(2C,2\epsilon)$-canonical neighborhoods of
arbitrarily small scale, and hence $2\epsilon$-necks of arbitrarily
small scale. But this contradicts Proposition~\ref{narrows}. It
follows from this contradiction that the curvature of
$(M_\infty,g_\infty)$ is bounded.
\end{proof}

To complete the proof of Theorem~\ref{smlmtflow} it remains to
extend the limit for the $0$ time-slices of the $({\mathcal
M}_n,G_n)$ that we have just constructed to a limit flow defined for
some positive amount of time backward. Since the curvature of
$(M_\infty,g_\infty)$ is bounded, this implies that there is a
bound, $Q$, such that for any $A<\infty$
 the curvature of the restriction of
$Q_nG_n$ to $B_{Q_nG_n}(x_n,0,A)$ is bounded by $Q$ for all $n$
sufficiently large. Thus, we can take the constant $Q(A)$ in
Corollary~\ref{cnnbhd2} to be independent of $A$. According to that
corollary this implies that there is a $t'_0>0$ and $Q'<\infty$ such
that for every $A$ there is an embedding $B_{Q_nG_n}(x_N,0,A)\times
(-t'_0,0]\to {\mathcal M}_n$ compatible with time and with the
vector field so that the scalar curvature of the restriction of
$Q_nG_n$ to the image is bounded by $Q'$ for all $n$ sufficiently
large. This uniform bound on the scalar curvature yields a uniform
bound, uniform in the sense of being independent of $n$, on $|{\rm
Rm}_{Q_nG_n}|$ on the image of the embedding
$B_{Q_nG_n}(x_N,0,A)\times (-t'_0,0]$.

 Then by Hamilton's result,
Proposition~\ref{partialflowlimit}, we see that, after passing to a
further subsequence, there is a limit flow defined on $(-t_0',0]$.
Of course, the zero time-slice of this limit flow is the limit
$(M_\infty,g_\infty)$. This completes the proof of
Theorem~\ref{smlmtflow}.
\end{proof}

\section{Long-time blow-up limits}

Now we wish to establish conditions under which we can, after
passing to a further subsequence, establish the existence of a
geometric limit flow defined on $-\infty<t\le 0$. Here is the main
result.

\begin{thm}\label{kaplimit}
Suppose that $\{({\mathcal M}_n,G_n,x_n)\}_{n=1}^\infty$ is a
sequence of generalized $3$-dimensional Ricci flows satisfying all
the hypothesis of Theorem~\ref{smlmtflow}. Suppose in addition that
there is  $T_0$ with $0<T_0\le \infty$ such that the following
holds. For any $T<T_0$, for each $A<\infty$, and all $n$
sufficiently large, there is an embedding
$B(x_n,t_n,AQ_n^{-1/2})\times (t_n-TQ_n^{-1},t_n]$ into ${\mathcal
M}_n$ compatible with time and with the vector field and at every
point of the image the generalized flow is
$\kappa$-non-collapsed\index{$\kappa$-non-collapsed} on scales $\le
r$. Then, after shifting the times of the generalized flows so that
$t_n=0$ for all $n$ and passing to a subsequence there is a
geometric limit Ricci flow
$$(M_\infty,g_\infty(t),x_\infty),\ -T_0<t\le 0,$$
for the rescaled generalized flows $(Q_n{\mathcal M}_n,Q_nG_n,x_n)$.
This limit flow is complete and of  non-negative curvature.
Furthermore, the curvature is locally bounded in time. If in
addition $T_0=\infty$, then it is a $\kappa$-solution.
\end{thm}

\begin{rem}
Let us point out the differences between this result and
Theorem~\ref{smlmtflow}. The hypotheses of this theorem include all
the hypotheses of Theorem~\ref{smlmtflow}. The main difference
between the conclusions is that in Theorem~\ref{smlmtflow} the
amount of backward time for which the limit flow is defined depends
on the curvature bound for the final time-slice of the limit (as
well as how far back the flows in the sequence are defined). This
amount of backward time tends to zero as the curvature of the final
time-slice limit tends to infinity. Here, the amount of backward
time for which the limit flow is defined depends only on how far
backwards the flows in the sequence are defined.
\end{rem}

\begin{proof}
In Theorem~\ref{smlmtflow} we proved that, after passing to a
subsequence, there is a geometric limit Ricci flow, complete of
bounded non-negative curvature,
$$(M_\infty,g_\infty(t),x_\infty),\ -t_0\le t\le 0,$$
defined for some $t_0>0$. Our next step is to
extend the limit flow all the way back to time $-T_0$.

 \begin{prop}\label{maxT}
With the notation of, and under the hypotheses of
Theorem~\ref{kaplimit}, suppose that there is a geometric limit flow
$(M_\infty,g_\infty(t))$ defined for $-T<t\le 0$ which has
non-negative curvature locally bounded in time. Suppose that
$T<T_0$. Then the curvature of the limit flow is bounded and the
geometric limit flow can be extended to a flow with bounded
curvature defined on $(-(T+\delta),0]$ for some $\delta>0$.
 \end{prop}

\begin{proof}
 The
argument is by contradiction, so we suppose that there is a $T<T_0$
as in the statement of the proposition. Then the geometric limit
flow on $(-T,0]$ is complete of non-negative curvature and with the
curvature locally bounded in time. First suppose that the scalar
curvature is bounded by, say $Q<\infty$. Fix $T'<T$.
 The Riemannian manifold
$(M_\infty,g_\infty(T'))$ is complete of non-negative curvature with
the scalar curvature, and hence the norm of the Riemann curvature,
bounded by $Q$. Thus, for any $A<\infty$ for all $n$ sufficiently
large, the norm of the Riemann curvature of $Q_nG_n(-T')$ on
$B_{Q_nG_n}(x_n,-T,A)$ is bounded above by $2Q$. Also, arguing as in
the proof of Theorem~\ref{smlmtflow} we see that any point $y\in
M_\infty$ with $R(y,-T')> 4$ has a $(2C,2\epsilon)$-canonical
neighborhood. Hence, applying Lemma~\ref{controlnbhd} as in the
argument in the proof of
 Corollary~\ref{cnnbhd2} shows that for all $n$ sufficiently large,
every point in $B_{Q_nG_n}(x_n,-T',A)$ has a uniform size parabolic
neighborhood on which the Riemann curvature is uniformly bounded,
where both the time interval in the parabolic neighborhood and the
curvature bound on this neighborhood depend only on $C$ and the
curvature bound on $Q$ for the limit flow. According to Hamilton's
result (Proposition~\ref{partialflowlimit}) this implies that, by
passing to a further subsequence, we can extend the limit flow
backward
 beyond $-T'$ a uniform amount of time, say $2\delta$. Taking $T'>T-\delta$
 then gives the desired extension under the condition that the
scalar curvature is bounded on $(-T,0]$.

It remains to show that, provided that $T<T_0$, the scalar curvature
of the limit flow $(M_\infty,g_\infty(t)),\ -T<t\le 0$, is bounded.
 To establish this we need a couple of preliminary results.

\begin{lem}\label{compactbd}
Suppose that there is a geometric limit flow defined on $(-T,0]$ for
some $0<T\le T_0$ with $T<\infty$. We suppose that this limit is
complete with non-negative curvature, and with curvature locally
bounded in time. Suppose that $X\subset M_\infty$ is a compact,
connected subset. If ${\rm min}_{x\in X}(R_{g_\infty} (x,t))$ is
bounded, independent of $t$, for all $t\in (-T,0]$, then there is a
finite upper bound on $R_{g_\infty}(x,t)$ for all $x\in X$ and all
$t\in(-T,0]$.
\end{lem}

\begin{proof}
Let us begin with:

\begin{claim}\label{distch}
 Let $Q$ be an upper bound on $R(x,0)$ for all $x\in
M_\infty$. Then for any points $x,y\in M_\infty$ and any $t\in
(-T,0]$ we have
$$d_{t}(x,y)\le d_0(x,y)+16\sqrt{\frac{Q}{3}}T.$$
\end{claim}

\begin{proof}
Fix $-t_0\in (-T,0]$. Then for any $\epsilon>0$ sufficiently small,
by the Harnack inequality\index{Ricci flow!Harnack inequality for}
(the second result in Theorem~\ref{Harnack})
 we have
$$\frac{\partial R}{\partial t}(x,t)\ge -\frac{R(x,t)}{t+T-\epsilon}.$$
Taking the limit as $\epsilon\rightarrow 0$ gives
$$\frac{\partial R}{\partial t}(x,t)\ge -\frac{R(x,t)}{t+T},$$
and hence, fixing $x$,
$$\frac{dR(x,t)}{R(x,t)}\ge \frac{-dt}{(t+T)}.$$
Integrating from $-t_0$ to $0$ shows that
$${\rm log}(R(x,0))-{\rm log}(R(x,-t_0))\ge {\rm log}(T-t_0)-{\rm log}(T),$$
and since $R(x,0)\le Q$, this implies
$$R(x,-t_0)\le Q\frac{T}{T-t_0}.$$
Recalling that $n=3$ and that the curvature is non-negative we see that
$${\rm Ric}(x,-t_0)\le (n-1)\frac{QT}{2}\frac{1}{T-t_0}.$$
Hence by Corollary~\ref{corI.8.3}, for all $-t_0\in (-T,0]$ we have
that
$${\rm dist}_{-t_0}(x,y)\le {\rm dist}_0(x,y)+8\int_{-t_0}^0\sqrt{\frac{QT}{3(T+t)}}
\le {\rm dist}_0(x,y)+16\sqrt{\frac{Q}{3}}T.$$
\end{proof}

It follows immediately from this claim that any compact subset
$X\subset M_\infty$ has uniformly bounded diameter under all the
metrics $g_\infty(t);\ -T<t\le 0$.

By the hypothesis of the lemma there is a constant $C'<\infty$ such
that for each $t\in (-T,0]$ there is $y_{t}\in X$ with
$R_{g_\infty}(y_{t},t)\le C'$. Suppose that the conclusion of the
lemma does not hold. Then there is a sequence $t_m\rightarrow -T$ as
$m\rightarrow\infty$ and points $z_m\in X$ such that
$R_{g_\infty}(z_m,t_m)\rightarrow\infty$ as $m\rightarrow\infty$. In
this case, possibly after redefining the constant $C'$, we can also
assume that there is a point $y_m$ such that $2\le R(y_m,t_m)\le
C'$. Since the sequence $({\mathcal M}_n,Q_nG_n,x_n)$ converges
smoothly to $(M_\infty,g_\infty(t),x_\infty)$ for $t\in (-T,0]$, it
follows that for each $m$ there are sequences $\{y_{m,n}\in
{\mathcal M}_n\}_{n=1}^\infty$ and $\{z_{m,n}\in {\mathcal
M}_n\}_{n=1}^\infty$ with ${\bf t}(y_{m,n})={\bf t}(z_{m,n})=t_m$
converging to $(y_m,t_m)$ and $(z_m,t_m)$ respectively. Thus, for
all $m$ there is $n_0=n_0(m)$ such that for all $n\ge n_0$ we have:
\begin{enumerate}
\item[(1)] $1\le R_{Q_nG_n}(y_{m,n})\le 2C'$,
\item[(2)] $R_{Q_nG_n}(z_{m,n})\ge R_{g_\infty}(z_m,t_m)/2$,
\item[(3)] $d_{Q_nG_n}((y_{m,n}),(z_{m,n}))\le 2\,{\rm diam}_{g_\infty(t_m)}(X).$
\end{enumerate}
Because of the third condition and the fact that $X$ has uniformly
bounded diameter under all the metrics $g_\infty(t)$ for $t\in
(-T,0]$, the distance $d_{Q_{m,n}G_{m,n}}(z_{m,n},y_{m,n})$ is
bounded independent of $m$ and $n$ as long as $n\ge n_0$. Because of
the fact that $R_{Q_nG_n}(y_{m,n})\ge 1=R(x_n)$, it follows that any
point $z\in {\mathcal M}_n$ with ${\bf t}(z_m)\le t_{m}$ and with
$R(z)\ge 4R(y_{m,n})$ has a strong $(C,\epsilon)$-canonical
neighborhood\index{canonical neighborhood!strong}. This then
contradicts Theorem~\ref{bcbd} and completes the proof of the lemma.
\end{proof}

Clearly, this argument will be enough to handle the case when
$M_\infty$ is compact. The case when $M_\infty$ is non-compact uses
additional results.

\begin{lem}
Let $(M,g)$ be a complete, connected, non-compact manifold of
non-negative sectional curvature and let $x_0\in M$ be a point. Then
there is $D>0$, such that for any $y\in M$ with $d(x_0,y)=d\ge D$,
there is $x\in M$ with $d(y,x)=d$ and with $d(x_0,x)>3d/2$.
\end{lem}

\begin{proof}
Suppose that the result is false for $(M,g)$ and $x\in M$. Then
there is a sequence $y_n\in M$ such that setting $d_n=d(x,y_n)$ we
have ${\rm lim}_{n\rightarrow\infty}d_n=\infty$ and yet
$B(y_n,d_n)\subset B(x,3d_n/2)$ for every $n$. Let $\gamma_n$ be a
minimal geodesic from $x$ to $y_n$. By passing to a subsequence we
arrange that the $\gamma_n$ converge to a minimal geodesic ray
$\gamma$ from $x$ to infinity in $M$. In particular, the angle at
$x$ between $\gamma_n$ and $\gamma$ tends to zero as
$n\rightarrow\infty$. Let $w_n$ be the point on $\gamma$ at distance
$d_n$ from $x$, and let $\alpha_n=d(y_n,w_n)$. Because $(M,g)$ has
non-negative curvature, by Corollary~\ref{anglemono}, ${\rm
lim}_{n\rightarrow\infty}\alpha_n/d_n=0$. In particular, for all $n$
sufficiently large, $\alpha_n<d_n$. This implies that there is a
point $z_n$ on the sub-ray of $\gamma$ with endpoint $w_n$ at
distance $d_n$ from $y_n$. By the triangle inequality,
$d(w_n,z_n)\ge d_n-\alpha_n$. Since $\gamma$ is a minimal geodesic
ray,  $d(z,z_n)=d(z,w_n)+d(w_n,z_n)\ge 2d_n-\alpha_n$. Since
$\alpha_n/d_n\rightarrow 0$ as $n\rightarrow \infty$, it follows
that for all $n$ sufficiently large $d(z,z_n)>3d_n/2$. This
contradiction proves the lemma.
\end{proof}

\begin{claim}\label{outerbd}
Fix $D<\infty$ greater than or equal to the constant given in the
previous lemma
 for the Riemannian manifold
$(M_\infty,g_\infty(0))$ and the point $x_\infty$. We also choose
$D\ge 32\sqrt{\frac{Q}{3}}T$. Then for any $y\in M_\infty\setminus
B(x_\infty,0,D)$ the scalar curvature $R_{g_\infty}(y,t)$ is
uniformly bounded for all $t\in (-T,0]$.
\end{claim}

\begin{proof}
Suppose this does not hold for some $y\in M_\infty\setminus
B(x_\infty,0,D)$. Let $d=d_0(x_\infty,y)$. Of course, $d\ge D$.
Thus, by the lemma there is $z\in M_\infty$ with $d_0(y,z)=d$ and
$d_0(x_\infty,z)>3d/2$. Since the scalar curvature $R(y,t)$ is not
uniformly bounded for all $t\in (-T,0]$, there is $t$ for which
$R(y,t)$ is arbitrarily large and hence $(y,t)$ has an
$(2C,2\epsilon)$-canonical neighborhood of arbitrarily small scale.
By Claim~\ref{distch}  we have $d_{t}(x_\infty,y)\le
d+8\sqrt{\frac{Q}{3}}T$ and $d_{t}(y,z)\le d+8\sqrt{\frac{Q}{3}}T$.
Of course, since ${\rm Ric}\ge 0$ the metric is non-increasing in
time and hence $d\le {\rm min}(d_{t}(y,z),d_{t}(x_\infty,y))$ and
$3d/2\le d_{t}(x_0,z)$. Since $y$ has a $(2C,2\epsilon)$-canonical
neighborhood in $(M_\infty,g_\infty(t))$, either $y$ is the center
of an $2\epsilon$-neck in $(M_\infty,g_\infty(t))$ or $y$ is
contained in the core of a $(2C,2\epsilon)$-cap in
$(M_\infty,g_\infty(t))$. (The other two possibilities for canonical
neighborhoods require that $M_\infty$ be compact.)

\begin{claim}
$y$ cannot lie in the core of a $(2C,2\epsilon)$-cap in
$(M_\infty,g_\infty(t))$, and hence it is the center of a
$2\epsilon$-neck\index{$\epsilon$-neck} $N$ in
$(M_\infty,g_\infty(t))$. Furthermore, minimal $g(t)$-geodesics from
$y$ to $x_\infty$ and $z$ exit out of opposite ends of $N$ (see {\sc
Fig.}~\ref{fig:mingeo}).
\end{claim}

\begin{figure}[ht]
  \centerline{\epsfbox{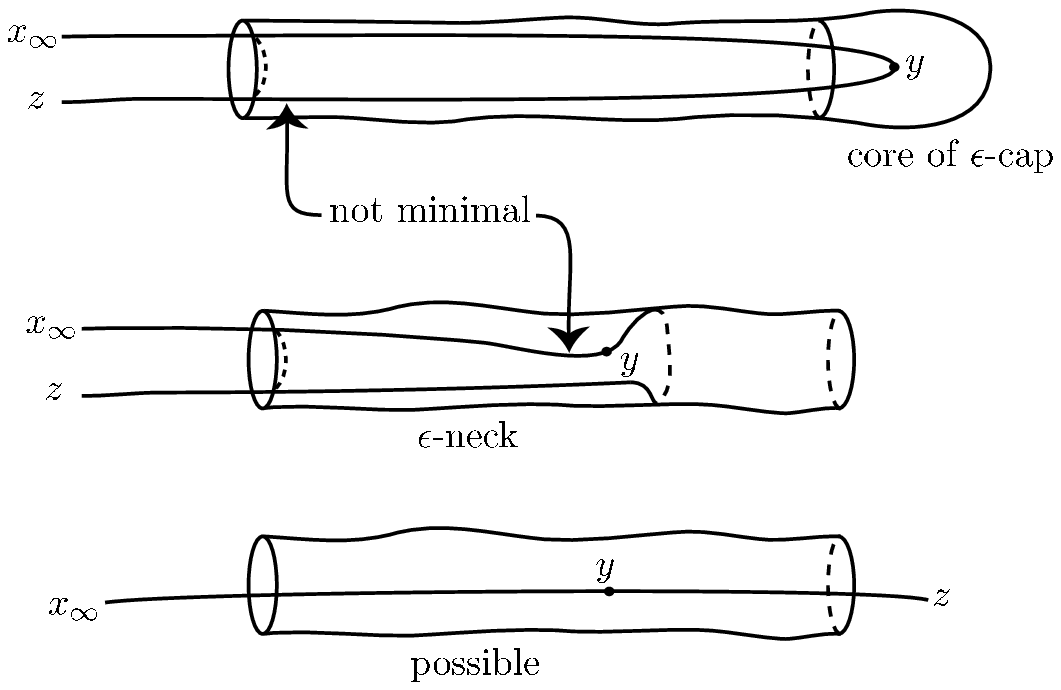}}
  \caption{Minimal geodesics in necks and caps}\label{fig:mingeo}
\end{figure}

\begin{proof}
Let ${\mathcal C}$ be a $(2C,2\epsilon)$-canonical
neighborhood\index{canonical neighborhood} of $y$ in
$(M_\infty,g_\infty(t))$. Since $R(y,t)$ can be arbitrarily large,
we can assume that $d \gg 2CR(y)^{-1/2}$, which is a bound on the
diameter of ${\mathcal C}$. This means that minimal $g(t)$-geodesics
$\gamma_{x_\infty}$ and $\gamma_z$ connecting $y$ to $x_\infty$ and
to $z$, respectively, must exit from ${\mathcal C}$. Let $a$ be a
point on $\gamma_{x_\infty}\cap{\mathcal C}$ close to the complement
of ${\mathcal C}$. Let $b$ be a point at the same $g(t)$-distance
from $y$ on $\gamma_z$. In the case that  ${\mathcal C}$ is a cap
or that it is a $2\epsilon$-neck and $\gamma_{x_\infty}$ and
$\gamma_z$ exit from the same end, then
$d_t(b,y)/d_t(a,y)<4\pi\epsilon$. This means that the angle $\theta$
of the Euclidean triangle with these side lengths at the point
corresponding to $y$ satisfies
$${\rm cos}(\theta)\ge 1-\frac{(4\pi\epsilon)^2}{2}.$$

Recall that $Q$ is the maximum value of $R(x,0)$,  and that by
Claim~\ref{distch} we have
 $$d\le d_t(x_\infty,y)\le d+16\sqrt{\frac{Q}{3}}T,$$
 with the same inequalities holding with $d_t(z,y)$ replacing $d_t(x_\infty,y)$.
 Also, by construction $d\ge 32\sqrt{\frac{Q}{3}}T$.
We set $a_0=d_t(x_\infty,y)$ and $a_1=d_t(z,y)$. Then by the
Toponogov property we have
$$d_t(x,z)^2\le a_0^2+a_1^2-2a_0a_1\left(
1-\frac{(4\pi\epsilon)^2}{2}\right)=(a_0-a_1)^2+(4\pi\epsilon)^2a_0a_1.$$
Since $|a_0-a_1|\le d/2$ and $a_0,a_1\le 3d/2$ and
$\epsilon<1/8\pi$, it follows that $d_t(x,z)<d$. Since distances do
not increase under the flow, it follows that $d_0(x,z)<d$. This
contradicts the fact that $d_0(x,z)=d$.
\end{proof}

It follows that  the point $y$ is the center of a
$(2C,2\epsilon)$-neck $N$ in $(M_\infty,g_\infty(t))$ and minimal
$g(t)$-geodesics from $y$ to $z$ and to $x_\infty$ exit out of
opposite ends of $N$. This implies that $B(y,t,4\pi R(y,t)^{-1/2})$
separates $x_\infty$ and $z$. Since the curvature of the time-slices
is non-negative, the Ricci flow does not increase distances. Hence,
$B(y,0,4\pi R(y,t)^{-1/2})$ separates $z$ from $x_\infty$. (Notice
that since $d>4\pi R(y,t)^{-1/2}$, neither $z$ nor $x_\infty$ lies
in this ball.) Thus, if  $R(y,t)$ is unbounded as $t\rightarrow -T$
then arbitrarily small $g(0)$-balls centered at $y$ separate $z$ and
$x_\infty$. Since $y$ is distinct from $x_\infty$ and $z$, this is
clearly impossible.
\end{proof}

Next we establish that the curvature near the base point $x_\infty$
is bounded for all $t\in (-T,0]$.

\begin{cor}\label{ancientlimit}
Suppose that there is a geometric limit flow
$(M_\infty,g_\infty(t))$ of a subsequence defined on $(-T,0]$ for
some $T<\infty$. Suppose that the limit flow is complete of
non-negative curvature with the curvature locally bounded in time.
Then for every $A<\infty$ the scalar curvature $R_{g_\infty}(y,t)$
is uniformly bounded for all $(y,t)\in B(x_\infty,0,A)\times
(-T,0]$.
\end{cor}

\begin{proof}
First we pass to a subsequence so that a geometric limit flow
$$(M_\infty,g_\infty(t),(x_\infty,0))$$ exists on  $(-T,0]$. We let
$Q$ be the upper bound for $R(x,0)$ for all $x\in M_\infty$. We now
divide the argument into two cases: (i) $M_\infty$ is compact, and
(ii) $M_\infty$ is non-compact.

Suppose that $M_\infty$ is compact. By Proposition~\ref{scalarevol}
we know that
$${\rm min}_{x\in M_\infty}(R_{g_\infty}(x,t))$$
is a non-decreasing function of $t$. Since
$R_{g_\infty}(x_\infty,0)=1$, it follows that for each $t\in
(-T,0]$, we have ${\rm min}_{x\in M_\infty}R(x,t)\le 1$, and hence
there is a point $x_{t}\in M_\infty$ with $R(x_{t},t)\le 1$. Now we
can apply Lemma~\ref{compactbd} to see that the scalar curvature of
$g_\infty$ is bounded on all of $M_\infty\times (-T,0]$.

If $M_\infty$ is non-compact, choose $D$ as in Lemma~\ref{outerbd}.
According to that lemma every point in the boundary of
$B(x_\infty,0,D)$ has bounded curvature under $g_\infty(t)$ for all
$t\in (-T,0]$. In particular, for each $t\in (-T,0]$ the minimum of
$R(x,t)$ over $\overline B(x_\infty,0,D)$ is bounded independent of
$t$. Now apply Lemma~\ref{compactbd} to the closure of
$B(x_\infty,0,D)$. We conclude that the curvature of
$B(x_\infty,0,D)$ is uniformly bounded for all $g_\infty(t)$ for all
$t\in (-T,0]$. In particular, $R(x_\infty,t)$ is uniformly bounded
for all $t\in (-T,0]$.

Now for any $A<\infty$ we apply Lemma~\ref{compactbd} to the compact
subset $\overline B(x_\infty,0,A)$ to conclude that the curvature is
uniformly bounded on $B(x_\infty,0,A)\times (-T,0]$. This completes
the proof of the corollary.
\end{proof}

Now let us return to the proof of Proposition~\ref{maxT}.

\begin{claim}
For each $A<\infty$ and for all $n$ sufficiently large, there are
$\delta>0$ with $\delta\le T_0-T$ and a bound, independent of $n$,
on the scalar curvature of the restriction of $Q_nG_n$ to
$B_{Q_nG_n}(x_n,0,A)\times [-(T+\delta),0]$.
\end{claim}

\begin{proof}
Fix $A<\infty$ and let $K$ be the bound for the scalar curvature of
$g_\infty$ on $B(x_\infty,0,2A)\times (-T,0]$ from
Corollary~\ref{ancientlimit}.  Lemma~\ref{controlnbhd} shows that
there are $\delta>0$ and a bound in terms of $K$ and $C$ on the
scalar curvature of the restriction of $Q_nG_n$ to
$B_{Q_nG_n}(x_n,0,A)\times [-(T+\delta),0]$.
\end{proof}

Since the scalar curvature is bounded, by the assumption that either
the curvature is pinched toward positive or the Riemann curvature is
non-negative, this implies that the sectional curvatures of $Q_nG_n$
are also uniformly bounded on the products
$B_{Q_nG_n}(x_n,0,A)\times [-(T+\delta),0]$ for all $n$ sufficiently
large. Consequently, it follows that by passing to a further
subsequence we can arrange that the $-T$ time-slices of the
$({\mathcal M}_n,G_n,x_n)$ converge to a limit
$(M_\infty,g_\infty(-T))$.
 This limit manifold
satisfies the hypothesis of Proposition~\ref{narrows} and hence, by
that proposition, it has bounded sectional curvature. This means
that there is a uniform $\delta>0$ such that for all $n$
sufficiently large and for any $A<\infty$ the scalar curvatures (and
hence the Riemann curvatures) of the restriction of $Q_nG_n$ to
$B_{Q_nG_n}(x_n,0,A)\times [-(T+\delta),0]$ are uniformly bounded.
This allows us to pass to a further subsequence for which there is a
geometric limit defined on $(-(T+\delta/2),-T]$. This geometric
limit is complete of  bounded, non-negative curvature. Hence, we
have now constructed a limit flow on $(-(T+\delta/2),0]$ with the
property that for each $t\in (-(T+\delta/2),0]$ the Riemannian
manifold $(M,g(t))$ is complete and of bounded non-negative
curvature. (We still don't know whether the entire flow is of
bounded curvature.) But now invoking Hamilton's Harnack inequality
(Theorem~\ref{Harnack}), we see that the curvature is bounded on
$[-T,0]$. Since we already know it is bounded in $(-T+\delta/2,-T]$,
this completes the proof of the proposition.
\end{proof}

 It follows immediately from Proposition~\ref{maxT} that there is a
geometric  limit flow defined on  $(-T_0,0]$. The geometric limit
flow on $(-T_0,0]$ is complete of non-negative curvature, locally
bounded in time.

It remains to prove the last statement in the theorem. So let us
suppose that $T_0=\infty$. We have just established the existence of
a geometric limit flow defined for $t\in (-\infty,0]$.  Since the
$({\mathcal M}_n,G_n)$ either have curvature pinched toward positive
or are of non-negative curvature, it follows from
Theorem~\ref{blowupposcurv} that all time-slices of the limit flow
are complete manifolds of non-negative curvature. Since points of
scalar curvature greater than $4$ have $(2C,2\epsilon)$-canonical
neighborhoods, it follows from Proposition~\ref{narrows} that the
curvature is bounded on each time-slice, and hence universally
bounded by the Harnack inequality (Theorem~\ref{Harnack}). Since for
any $A<\infty$ and every $T<\infty$ the parabolic neighborhoods
$B_{Q_nG_n}(x_n,0,A)\times [-T,0]$ are $\kappa$-non-collapsed on
scales $Q_nr$ for every $n$ sufficiently large,  the limit  is
$\kappa$-non-collapsed on scales $\le {\rm lim}_{n\rightarrow
\infty}Q_nr$. Since  $r>0$ and ${\rm lim}_{n\rightarrow
\infty}Q_n=\infty$, it follows that the limit flow is
$\kappa$-non-collapsed on all scales. Since $R_{Q_n G_n}(x_n)=1$,
$R_{g_\infty}(x_\infty,0)=1$ and the limit flow is non-flat. This
establishes all the properties need to show that the limit is a
$\kappa$-solution. This completes the proof of
Theorem~\ref{kaplimit}.
\end{proof}

\section{Incomplete smooth limits at singular times}

Now we wish to consider smooth limits where we do not blow up, i.e.,
do not rescale the metric. In this case the limits that occur can be
incomplete, but we have strong control over their ends.

\subsection{Assumptions}

We shall assume the following about the generalized Ricci flow
$({\mathcal M},G)$:
\begin{assumption}\label{assumptions}

\begin{enumerate}
\item[(a)] The singular times form a discrete subset of $\Ar$, and
each time slice of the flow at a non-singular time is a compact
$3$-manifold. \item[(b)] The time interval of definition of the
generalized Ricci flow $({\mathcal M},G)$ is contained in
$[0,\infty)$ and its curvature is pinched toward positive.
\item[(c)] There are $r_0>0$ and $C<\infty$, such that any point
$x\in {\mathcal M}$ with $R(x)\ge r_0^{-2}$ has a strong
$(C,\epsilon)$-canonical neighborhood. In particular, for every
$x\in {\mathcal M}$  with $R(x)\ge r_0^{-2}$ the following two
inequalities hold:
$$\left|\frac{\partial R(x)}{\partial t}\right|< CR^2(x),$$
$$|\nabla R(x)|< CR^{3/2}(x).$$
\end{enumerate}
\end{assumption}
With these assumptions we can say quite a bit about the limit metric
at time $T$.

\begin{thm}\label{omega}
Suppose that $({\mathcal M},G)$ is a generalized Ricci flow defined
for $0\le t<T<\infty$ satisfying the three assumptions given
in~\ref{assumptions}. Let $T^-<T$ be such that there is a
diffeomorphism $\rho\colon M_{T^-}\times [T^-,T)\to {\bf
t}^{-1}([T^-,T))$ compatible with time and with the vector field.
Set $M=M_{T^-}$ and let $g(t),\ T^-\le t<T$, be the family of
metrics $\rho^*G(t)$ on $M$. Let $\Omega\subset
M$\index{$\Omega$|(ii} be the subset of defined by
$$\Omega=\left\{x\in M\bigl|\bigr. {\rm liminf}_{t\rightarrow T}R_g(x,t)<\infty\right\}.$$
Then $\Omega\subset M$ is an open subset and there is a Riemannian
metric $g(T)$ with the following properties:
\begin{enumerate}
\item[(1)] As $t\rightarrow T$ the metrics $g(t)|_\Omega$ limit to $g(T)$
 uniformly in the $C^\infty$-topology
on every compact subset of $\Omega$.
\item[(2)] The scalar curvature $R(g(T))$ is a proper function from
$\Omega\to \Ar$ and is bounded below.
\item[(3)] Let
$$\widehat{\mathcal M} = {\mathcal M}\cup_{\Omega\times [T^-,T)}\left(\Omega\times [T^-,T]\right).$$
Then the generalized Ricci flow $({\mathcal M},G)$ extends to a
generalized Ricci flow  $(\widehat {\mathcal M},\widehat G)$.
\item[(4)] Every end of a connected component of
$\Omega$ is contained in a strong
$2\epsilon$-tube\index{$\epsilon$-tube!strong}.
\item[(5)] Any point $x\in \Omega\times\{T\}$ with $R(x)>r_0^{-2}$
has a strong $(2C,2\epsilon)$-canonical neighborhood\index{canonical
neighborhood!strong} in $\widehat{\mathcal M}$.
\end{enumerate}
\end{thm}

\begin{rem}
Recall that by definition a function $f$ is proper if the pre-image
under $f$ of every compact set is compact.
\end{rem}

In order to prove this result we establish a sequence of lemmas. The
first in the series establishes that $\Omega$ is an open subset and
also establishes the first two of the above five conclusions.

\begin{lem}\label{Rproper}
Suppose that $({\mathcal M},G)$ is a generalized Ricci flow defined
for $0\le t<T<\infty$ satisfying the three assumptions given
in~\ref{assumptions}. Let $T'<T$ be as in the previous theorem, set
$M=M_{T^-}$, and let $g(t)$ be the family of metrics on $M$ and let
$\Omega\subset M$, each being  as defined in the previous theorem.
Then $\Omega\subset M$ is an open subset of $M$.
 Furthermore, the restriction of the family $g(t)$ to $\Omega$
 converges in the  $C^\infty$-topology, uniformly on compact sets of
 $\Omega$,
 to a Riemannian metric $g(T)$.
Lastly,  $R(g(T))$ is a proper function, bounded below, from
$\Omega$ to $\Ar$.
\end{lem}

\begin{proof}
 We pull back $G$ to $M\times [T^-,T)$ to define a Ricci flow
$(M,g(t)),\ T^{-}\le t<T$. Suppose that $x\in \Omega$. Then there is
a sequence $t_n\rightarrow T$ as $n\rightarrow \infty$ such that
$R(x,t_n)$ is bounded above, independent of $n$, by say $Q$. For all
$n$ sufficiently large we have $T-t_n\le
\frac{1}{16C(Q^2+r_0^{-2})}$. Fix such an $n$. Then, according to
the Lemma~\ref{controlnbhd}, there is $r>0$ such that $R(y,t)$ is
uniformly bounded for $y\in B(x,t_n,r)\times [t_n,T)$. This means
that $B(x,t_n,r)\subset \Omega$, proving that $\Omega$ is open in
$M$.

Furthermore, since $R(y,t)$ is bounded on $B(x,t_n,r)\times
[t_n,T)$, it follows from the curvature pinching toward positive
hypothesis that $|{\rm Rm}(y,t)|$ is bounded on
$B(x,t_n,r)\times[t_n,T)$. Now applying Theorem~\ref{shi} we see
that in fact ${\rm Rm}$ is bounded in the $C^\infty$-topology on
$B(x,t_n,r)\times[(t_n+T)/2 ,T)$. The same is of course also true
for ${\rm Ric}$ and hence for $\frac{\partial g}{\partial t}$ in the
$C^\infty$-topology. It then follows that there is a continuous
extension of $g$ to $B(x,t_n,r)\times [t_n,T]$.  Since this is true
for every $x\in \Omega$ we see that $g(t)$ converges in the
$C^\infty$-topology, uniformly on compact subsets of $\Omega$, to
$g(T)$.

Lastly, let us consider the function $R(g(T))$ on $\Omega$. Since
the metric $g(T)$ is a smooth metric on $\Omega(T)$, this is a
smooth function. Clearly, by the curvature pinching toward positive
hypothesis, this function is bounded below. We must show that it is
proper. Since $M$ is compact, it suffices to show that if $x_n$ is a
sequence in $\Omega\subset M$ converging to a point $x\in M\setminus
\Omega$ then $R(x_n,T)$ is unbounded. Suppose that $R(x_n,T)$ is
bounded independent of $n$. It follows from Lemma~\ref{controlnbhd}
that there is a positive constant $\Delta t$ such that $R(x_n,t)$ is
uniformly bounded for all $n$ and all $t\in [T-\Delta t,T)$, and
hence, by the same result, there is $r>0$ such that $R(y_n, t)$ is
bounded for all $n$,  all $y_n\in B(x_n,T-\Delta t,r)$, and all
$t\in [T-\Delta t,T)$. Since the $x_n\rightarrow x\in M$, it follows
that for all $n$ sufficiently large that $x\in B(x_n,T-\Delta t,r)$,
and hence $R(x,t)$ is uniformly bounded as $t\rightarrow T$. This
contradicts the fact that $x\not\in \Omega$.
\end{proof}

\begin{defn}\label{hatM}
Let
$$\widehat {\mathcal M}= {\mathcal
M}\cup_{\Omega\times[T^-,T)} \left(\Omega\times
[T^-,T]\right).$$\index{$\widehat{\mathcal M}$} Since both
${\mathcal M}$ and $\Omega\times [T^-,T]$ have the structure of
space-times and the time functions and vector fields agree on the
overlap, $\widehat {\mathcal M}$ inherits the structure of a
space-time. Let $G'(t),\ T^-\le t\le T$, be the smooth family of
metrics on $\Omega$. The horizontal metrics, $G$, on ${\mathcal M}$
and this family of metrics on $\Omega$ agree on the overlap and
hence define a horizontal metric $\widehat G$ on $\widehat {\mathcal
M}$. Clearly, this metric satisfies the Ricci flow equation, so that
$(\widehat {\mathcal M},\widehat G)$ is a generalized Ricci flow
extending $({\mathcal M},G)$. We call this the {\em maximal
extension of } $({\mathcal M},G)$ to time $T$. Notice that even
though the time-slices $M_t$ of ${\mathcal M}$ are compact, it will
not necessarily be the case that the time-slice $\Omega$ is
complete.
\end{defn}

At this point we have established the first three of the five
conclusions stated in Theorem~\ref{omega}. Let us turn to the last
two.

\subsection{Canonical neighborhoods for $(\widehat{\mathcal M},\widehat G)$}

We continue with the notation and assumptions of the previous
subsection. Here we establish the fifth conclusion in
Theorem~\ref{omega}, namely the existence of strong canonical
neighborhoods\index{canonical neighborhood!strong} for
$(\widehat{\mathcal M},\widehat G)$

\begin{lem}\label{omegacanon}
For any $x\in \Omega\times\{T\}$ with $R(x,T)>r_0^{-2}$ one of the
following holds:
\begin{enumerate}
\item[(1)] $(x,T)$ is the center of a strong $2\epsilon$-neck in
$(\widehat{\mathcal M},\widehat G)$. \item[(2)] There is a
$(2C,2\epsilon)$-cap in $(\Omega(T),\widehat G(T))$ whose core
contains $(x,T)$. \item[(3)] There is a $2C$-component of
$\Omega(T)$ that contains $(x,T)$.
\item[(4)] There  is a $2\epsilon$-round component of $\Omega(T)$ that
contains $(x,T)$.
\end{enumerate}
\end{lem}

\begin{proof}
We fix $x\in \Omega(T)$ with $R(x,T)>r_0^{-2}$. First notice that
 for all $t<T$ sufficiently close to $T$ we
have $R(x,t)>r_0^{-2}$. Thus, for all such $t$ the point $(x,t)$ has
a strong $(C,\epsilon)$-canonical neighborhood\index{canonical
neighborhood!strong} in $({\mathcal M},G)\subset (\widehat{\mathcal
M},\widehat G)$. Furthermore, since ${\rm lim}_{t\rightarrow
T}R(x,t)=R(x,T)<\infty$, for all $t<T$ sufficiently close to $T$,
there is a constant $D<\infty$ such that for any point $y$ contained
in  a strong $(C,\epsilon)$-canonical neighborhood containing
$(x,t)$, we have $D^{-1}R(x,T)\le R(y,t)\le DR(x,T)$. Again assuming
that $t<T$ is sufficiently close to $T$, by Lemma~\ref{controlnbhd}
there is $D'<\infty$ depending only on $D$, $t$, and $r_0$ such that
the curvature $R(y,T)$ satisfies $(D')^{-1}R(x,T)\le R(y,T)\le
D'R(x,T)$. By Lemma~\ref{Rproper} this implies that there is a
compact subset $K\subset \Omega(T)$ containing all the
$(C,\epsilon)$-canonical neighborhoods for $(x,t)$.
 By the same lemma, the metrics $G(t)|_K$
converge uniformly in the $C^\infty$-topology to $G(T)|_K$. If there
is a sequence of $t$ converging to $T$ for which the canonical
neighborhood of $(y,t)$ is an $\epsilon$-round component, resp. a
$C$-component, then $(y,T)$ is contained in a $2\epsilon$-round,
resp. a $2C$-component of $\widehat \Omega$. If there is a sequence
of $t_n$ converging to $T$ so that each $(y,t_n)$ has a canonical
neighborhood ${\mathcal C}_n$ that is a $(C,\epsilon)$-cap whose
core contains $(y,t_n)$, then by Proposition~\ref{canonvary} since
these caps are all contained in a fixed compact subset $K$ and since
the $G(t_n)|_K$ converge uniformly in the $C^\infty$-topology to
$G(T)|_K$, it follows that for any $n$  sufficiently large, the
metric $G(T)$ restricted to ${\mathcal C}_n$ contains a
$(2C,2\epsilon)$-cap ${\mathcal C}$ whose core contains $(y,T)$.

Now we examine the case of strong $\epsilon$-necks.

\begin{claim}\label{epslimit}
Fix a point $x\in\Omega$. Suppose that there is a sequence
$t_n\rightarrow T$ such that for every $n$, the point $(x,t_n)$ is
the center of a strong $\epsilon$-neck in $\widehat{\mathcal M}$.
Then $(x,T)$ is the center of a strong $2\epsilon$-neck in $\widehat
{\mathcal M}$.
\end{claim}

\begin{proof}
By an overall rescaling we can assume that $R(x,T)=1$. For each $n$
let $N_n\subset \Omega$ and let $\psi_n\colon S^2\times
(-\epsilon^{-1},\epsilon^{-1})\to N_n\times\{t\}$ be a strong
$\epsilon$-neck centered at $(x,t_n)$. Let
$B=B(x,T,2\epsilon^{-1}/3)$. Clearly, for all $n$ sufficiently large
$B\subset N_n$. Thus, for each point $y\in B$ and each $n$ there is
a flow line through $y$ defined on the interval
$(t_n-R(x,t_n)^{-1},t_n]$. Since the $t_n\rightarrow T$ and since
$R(x,t_n)\rightarrow R(x,T)=1$ as $n\rightarrow \infty$, it follows
that there is a flow line through $y$ defined on $(T-1,T]$.

Consider the maps
$$\alpha_n\colon B\times (-1,0]\to \widehat{\mathcal M}$$
that send $(y,t)$ to the value at time $t_n-tR(x,t_n)^{-1}$ of the
flow line through $y$. Pulling back the metric $R(x,t_n)\widehat G$
by $\alpha_n$ produces the restriction of a strong $\epsilon$-neck
structure to $B$. The maps $\alpha_n$ converge uniformly in the
$C^\infty$-topology to the map $\alpha\colon B\times (-1,0]\to
\widehat M$ defined by sending $(y,t)$ to the value of the flowline
through $(y,T)$ at the time $T-t$. Hence, the sequence of metrics
$\alpha^*_n(R(x,t_n))\widehat G$ on $B\times (-1,0]$ converges
uniformly on compact subsets of $B\times (-1,0]$ in the
$C^\infty$-topology to the family $\alpha^*(\widehat G)$.  Then, for
all $n$ sufficiently large, the image $\psi_n(S^2\times
(-\epsilon^{-1}/2,\epsilon^{-1}/2))$ is contained in $B$ and has
compact closure in $B$. Since the family of  metrics
$\psi_n^*\widehat G$ on $B$ converge smoothly to $\psi^*\widehat G$,
it follows that
 for every $n$ sufficiently
large, the restriction of $\psi_n$ to $S^2\times
(-\epsilon^{-1}/2,\epsilon^{-1}/2)$ gives the coordinates showing
that the restriction of the family of metrics $\psi^*(\widehat G)$
to the image $\psi_n(S^2\times (-\epsilon^{-1}/2,\epsilon^{-1}/2))$
is a strong $2\epsilon$-neck at 
time $T$.
\end{proof}

This completes the proof of the lemma.
\end{proof}

The lemma tells us that every point  $x\in\Omega\times\{T\}$ with
$R(x)>r_0^{-2}$ has a strong $(2C,2\epsilon)$-canonical
neighborhood. Since, by assumption, points at time before $T$ with
scalar curvature at least $r_0^{-2}$ have strong
$(C,\epsilon)$-canonical neighborhoods, this completes the proof of
the fifth conclusion of Theorem~\ref{omega}. It remains to establish
the fourth conclusion of that theorem.

\subsection{The ends of $(\Omega,g(T))$}

\begin{defn}
 A {\em strong $2\epsilon$-horn}\index{horn|ii} in $(\Omega,g(T))$ is a submanifold of $\Omega$ diffeomorphic
to $S^2\times [0,1)$ with the following properties:
\begin{enumerate}
\item[(1)] The embedding $\psi$ of $S^2\times [0,1)$ into $\Omega$ is a proper map.
\item[(2)] Every point of  the image of this map is the center of a strong
$2\epsilon$-neck in $(\widehat{\mathcal M},\widehat G)$.
\item[(3)] The image of the boundary $S^2\times \{0\}$
is the central sphere of a strong $2\epsilon$-neck.
\end{enumerate}
\end{defn}

\begin{defn}
 A {\em strong double $2\epsilon$-horn}\index{horn!double|ii} in
$(\Omega,g(T))$ is a component of $\Omega$ diffeomorphic to
$S^2\times (0,1)$ with the property that every point of this
component is the center of a strong $2\epsilon$-neck in $\widehat
{\mathcal M}$. This means that a strong double $2\epsilon$-horn is a
$2\epsilon$-tube and hence is a component of $\Omega$ diffeomorphic
to $S^2\times(-1,1)$. Notice that each end of a strong double
$2\epsilon$-horn contains a strong $2\epsilon$-horn.

 For any $C'<\infty$, a $C'$-{\em capped
$2\epsilon$-horn}\index{horn!capped|ii} in $(\Omega,g(T))$ is a
component of $\Omega$ that is a the union of a the core of a
$(C',2\epsilon)$-cap and a strong $2\epsilon$-horn. Such a component
is diffeomorphic to an open $3$-ball or to a punctured $\Ar P^3$.
\end{defn}

\begin{figure}[ht]
  \centerline{\epsfbox{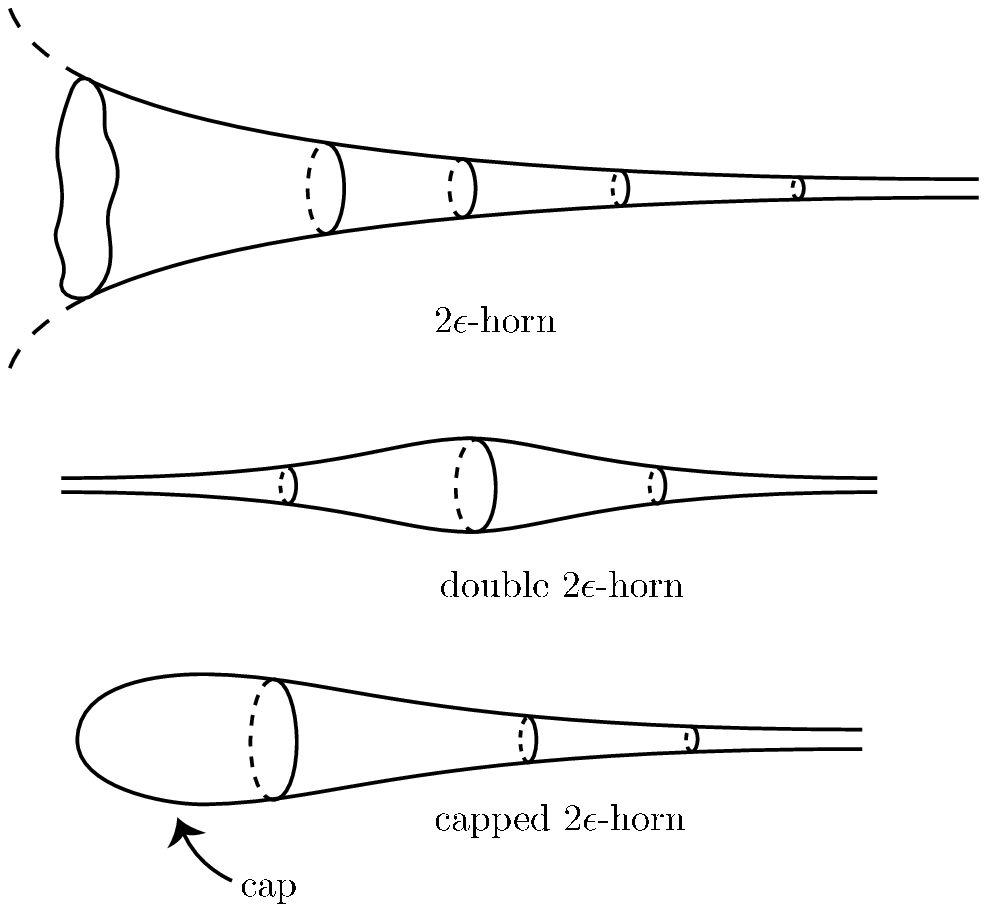}}
  \caption{Horns.}\label{fig:horns}
\end{figure}

See {\sc Fig.}~\ref{fig:horns}.

\begin{defn}
Fix any $\rho,\ 0<\rho<r_0$. We define $\Omega_\rho\subset \Omega$
to be the closed subset of all $x\in \Omega$ for which $R(x,T)\le
\rho^{-2}$. We say that a strong $2\epsilon$-horn $\psi\colon
S^2\times[0,1)\to \Omega$ has {\em boundary contained in}
$\Omega_\rho$\index{$\Omega_\rho$} if its boundary,
$\psi(S^2\times\{0\})$, is contained in $\Omega_\rho$.
\end{defn}

\begin{figure}[ht]
  \centerline{\epsfbox{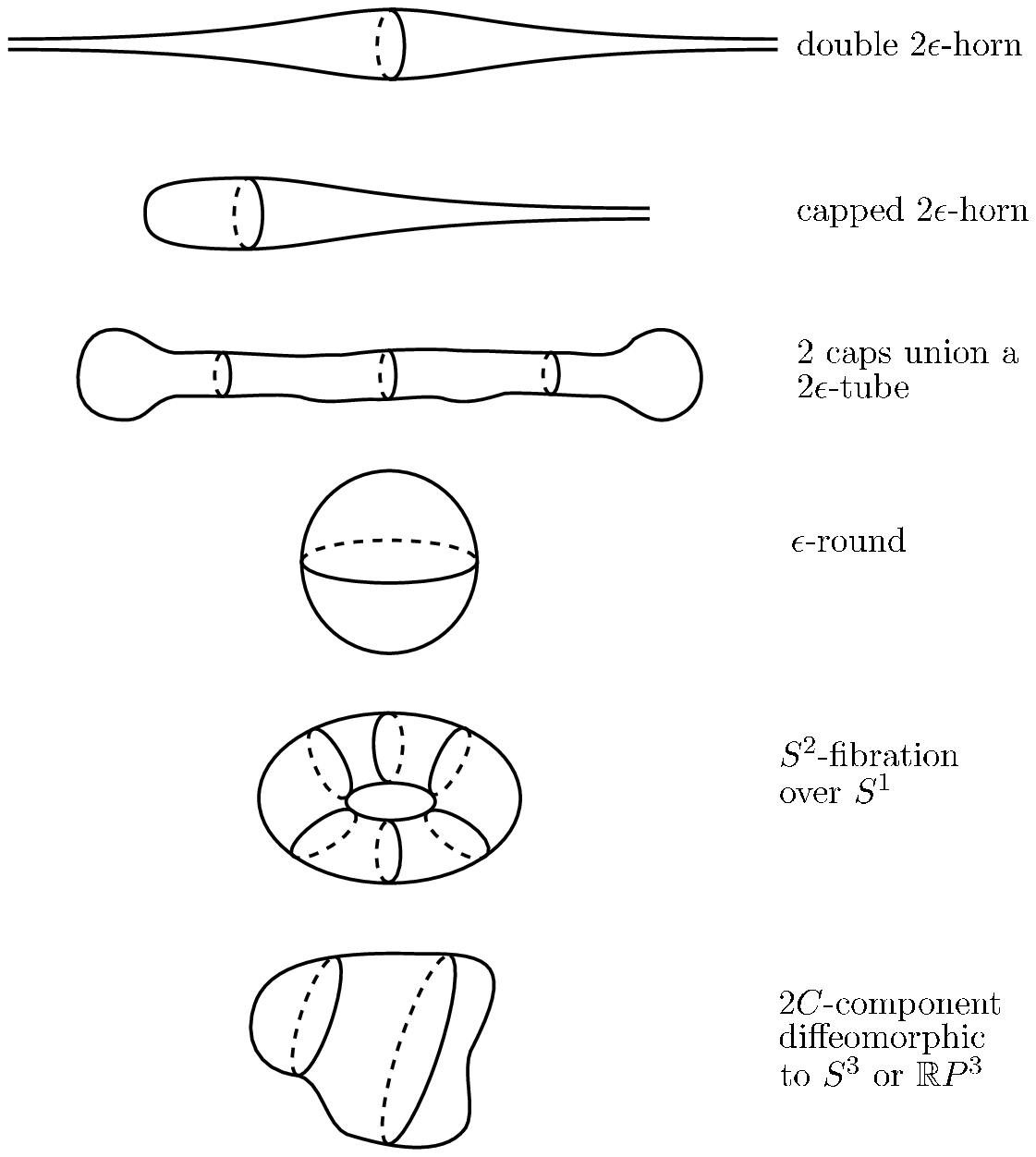}}
  \caption{Components of $\Omega$ disjoint from
  $\Omega_\rho$.}\label{fig:Omega}
\end{figure}

\begin{lem}
Suppose that $0<\rho<r_0$ and that $\Omega^0$ is a component of
$\Omega$ which contains no point of $\Omega_\rho$. Then one of the
following holds:
\begin{enumerate}
\item[(1)] $\Omega^0$ is a strong double $2\epsilon$-horn\index{horn!double} and is diffeomorphic to $S^2\times \Ar$.
\item[(2)] $\Omega^0$ is a  $2C$-capped $2\epsilon$-horn\index{horn!capped} and is diffeomorphic to $\Ar^3$
or to a punctured $\Ar P^3$.
\item[(3)] $\Omega^0$ is a compact component and is the union of the cores of
two $(2C,2\epsilon)$-caps and a strong $2\epsilon$-tube. It is
diffeomorphic to $S^3$, $\Ar P^3$ or $\Ar P^3\#\Ar P^3$.
\item[(4)] $\Omega^0$ is a compact $2\epsilon$-round component and is diffeomorphic to a compact
manifold of constant positive curvature.
\item[(5)] $\Omega^0$ is a compact component that fibers over $S^1$
with fibers $S^2$.
\item[(6)] $\Omega^0$ is a compact
$2C$-component and is diffeomorphic to $S^3$ or to $\Ar P^3$.
\end{enumerate}
See {\sc Fig.}\ref{fig:Omega}.
\end{lem}

\begin{proof}
Let $\Omega^0$ be a component of $\Omega$ containing no point of
$\Omega_\rho$. Then for every  $x\in \Omega^0$, we have
$R(x,T)>r_0^{-2}$.   Therefore, by Lemma~\ref{omegacanon} $(x,T)$
has a $(2C,2\epsilon)$-canonical neighborhood. Of course, this
entire neighborhood is contained in $\widehat{\mathcal M}$ and hence
is contained in $\Omega^0$ (or, more precisely, in the case of
strong $2\epsilon$-necks in the union of maximum backward flow lines
ending at points of $\Omega^0$). If the canonical neighborhood of
$(x,T)\in \Omega^0$ is a $2C$-component or is an $2\epsilon$-round
component, then of course $\Omega^0$ is that $2C$-component or
$2\epsilon$-round component. Otherwise, each point of $\Omega^0$ is
either the center of a strong $2\epsilon$-neck or is contained in
the core of a $(2C,2\epsilon)$-cap. We have chosen $2\epsilon$
sufficiently  small so that the result follows from
Proposition~\ref{epstopology}.
\end{proof}

\begin{rem}
We do not claim that there are only finitely many such components;
in particular, as far as we know there may be infinitely double
$2\epsilon$-horns.
\end{rem}

It follows immediately from this lemma that if $X$ is a component of $\Omega$
not containing any point of $\Omega_\rho$, then every end of $X$ is contained
in a strong $2\epsilon$-tube. To complete the proof of Theorem~\ref{omega}, it
remains only to establish the same result for the components of $\Omega$ that
meet $\Omega_\rho$. That is part of the content of the next lemma.

\begin{lem}\label{2epshorns}
Let $({\mathcal M},G)$ be a generalized $3$-dimensional Ricci flow
defined for $0\le t<T<\infty$ satisfying
Assumptions~\ref{assumptions}. Fix $0<\rho<r_0$.  Let
$\Omega^0(\rho)$ be the union of all components of $\Omega$
containing points of $\Omega_\rho$. Then $\Omega^0(\rho)$ has
finitely many components and is a union of a compact set and
finitely many strong $2\epsilon$-horns each of which is disjoint
from $\Omega_\rho$ and has its boundary contained in
$\Omega_{\rho/2C}$.
\end{lem}

\begin{proof}
Since $R\colon \Omega\times\{T\}\to \Ar$ is a proper function
bounded below, $\Omega_\rho$ is compact. Hence, there are only
finitely many components of $\Omega$ containing points of
$\Omega_\rho$. Let $\Omega^0$ be a non-compact component of $\Omega$
containing a point of $\Omega_\rho$, and let ${\mathcal E}$ be an
end of $\Omega^0$. Let
$$X=\{x\in \Omega^0\bigl|\bigr.
R(x)\ge 2C^2\rho^{-2}\}.$$ Then $X$ is a closed set and contains a
neighborhood of the end ${\mathcal E}$. Since $\Omega^0$ contains a
point of $\Omega_\rho$, $\Omega^0\setminus X$ is non-empty. Let
$X_0$ be the connected component of $X$ that contains a neighborhood
of ${\mathcal E}$. This is a closed, connected set every point of
which has  a $(2C,2\epsilon)$-canonical neighborhood. Since $X_0$
includes an end of $\Omega^0$, no point of $X_0$ can be contained in
an $\epsilon$-round component nor in a $C$-component. Hence, every
point of $X_0$ is either the center of a strong $2\epsilon$-neck or
is contained in the core of a $(2C,2\epsilon)$-cap. Since
$2\epsilon$ is sufficiently small to invoke
Proposition~\ref{Xcontainedin}, the latter implies that $X_0$ is
contained either in a $2\epsilon$-tube which is a union of strong
$2\epsilon$-necks centered at points of $X_0$ or $X_0$ is contained
in a $2C$-capped $2\epsilon$-tube where the core of the cap contains
a point of $X_0$. ($X_0$ cannot be contained in a double capped
$2\epsilon$-tube since the latter is compact.) In the second case,
since this capped tube contains an end of $\Omega^0$, it is in fact
equal to $\Omega^0$. Since a point of $X_0$ is contained in the core
of the $(2C,2\epsilon)$-cap, the curvature of this point is at most
$2C^2\rho^{-2}$ and hence the curvature at any point of the cap is
at least $2C\rho^{-2}>\rho^{-2}$. This implies that the cap is
disjoint from $\Omega_{\rho}$. Of course, any $2\epsilon$-neck
centered at a point of $X_0$ has curvature at least $C^2\rho^{-2}$
and hence is also disjoint from $\Omega_\rho$. Hence, if $\Omega^0$
is a $2C$-capped $2\epsilon$-tube and there is a point of $X_0$ in
the core of the cap, then this component is disjoint from
$\Omega_\rho$, which is a contradiction. Thus,  $X_0$ is contained
in a $2\epsilon$-tube made up of strong $2\epsilon$-necks centered
at points of $X_0$.

This proves that $X_0$ is contained in a strong $2\epsilon$-tube,
$Y$, every point of which has curvature $\ge C^2\rho^{-2}$. Since
$X_0$ is closed but not the entire component $\Omega^0$, it follows
that $X_0$ has a frontier point $y$. Of course,
$R(y)=2C^2\rho^{-2}$. Let $N$ be the strong $2\epsilon$-neck
centered at $y$ and let $S_N^2$ be its central $2$-sphere. Clearly,
every $y'\in S^2_N$ satisfies $R(y')\le 4C^2\rho^{-2}$, so that
$S^2_N$ is contained in $\Omega_{\rho/2C}$. Let $Y'\subset Y$ be the
complementary component of $S^2_N$ in $Y$ that contains a
neighborhood of the end ${\mathcal E}$. Then the closure of $Y'$ is
the required strong $2\epsilon$-horn containing a neighborhood of
${\mathcal E}$, disjoint from $\Omega_\rho$ and with boundary
contained in $\Omega_{\rho/2C}$.

 The last thing to see is
that there are only finitely many such ends in a given component
$\Omega^0$. First suppose that the boundary $2$-sphere of one of the
$2\epsilon$-horns is homotopically trivial in $\Omega^0$. Then this
$2$-sphere separates $\Omega^0$ into two components one of which is
compact and hence $\Omega^0$ has only one boundary component. Thus,
we can assume that all the boundary $2$-spheres of the
$2\epsilon$-horns are homotopically non-trivial. Suppose that two of
these $2\epsilon$-horns containing different ends of $\Omega^0$ have
non-empty intersection. Let $N$ be the $2\epsilon$-neck whose
central $2$-sphere is the boundary of one of the $2\epsilon$-horns.
Then the boundary of the other $2\epsilon$-horn is also contained in
$N$. This means that
 the
union of the two $2\epsilon$-horns and $N$ is a component of
$\Omega$\index{$\Omega$|)}. Clearly, this component has exactly two
ends. Thus, we can assume that all the $2\epsilon$-horns with
boundary in $\Omega_{\rho/2C}$ are disjoint. If two of the
$2\epsilon$-horns have boundary components that are topologically
parallel in $\Omega^0\cap \Omega_{\rho/2C}$ (meaning that they are
the boundary components of a compact submanifold diffeomorphic to
$S^2\times I$), then $\Omega^0$ is diffeomorphic to $S^2\times
(0,1)$ and has only two ends.
 By compactness of
$\Omega_{\rho/2C}$, there can only be finitely many disjoint
$2\epsilon$-horns with non-parallel, homotopically noon-trivial
boundaries in $\Omega^0\cap \Omega_{\rho/2C}$. This completes the
proof of the fact that each component of $\Omega_{\rho/2C}$ has only
finitely many ends.
\end{proof}

This completes the proof of Theorem~\ref{omega}.

\section{Existence of strong $\delta$-necks sufficiently deep in
a $2\epsilon$-horn}

We keep the notation and assumptions of the previous section.

\begin{figure}[ht]
  \relabelbox{
  \centerline{\epsfbox{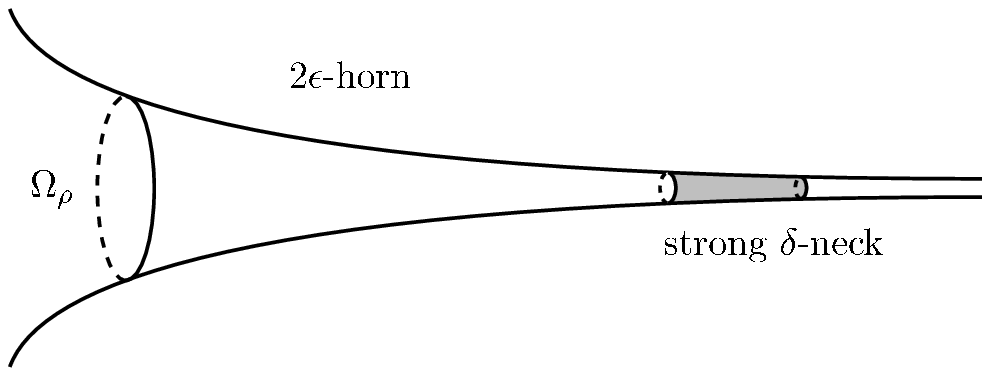}}}
  \relabel{1}{$2\epsilon$-horn}
  \relabel{2}{$\Omega_\rho$}
  \relabel{3}{strong $\delta$-neck}
  \endrelabelbox
  \caption{$\delta$-necks deep in a $2\epsilon$-horn.}\label{fig:deep}
\end{figure}

\begin{thm}\label{hexist}
Fix $\rho>0$. Then for any $\delta>0$ there is an
$0<h=h(\delta,\rho)\le {\rm min}(\rho\cdot\delta,\rho/2C)$,
implicitly depending on $r$ and $(C,\epsilon)$ which are fixed, such
that for any generalized Ricci flow $({\mathcal M},G)$ defined for
$0\le t<T<\infty$ satisfying Assumptions~\ref{assumptions} and for
any $2\epsilon$-horn ${\mathcal H}$ of $(\Omega,g(T))$ with boundary
contained in $\Omega_{\rho/2C}$, every point $x\in {\mathcal H}$
with $R(x,T)\ge h^{-2}$ is the center of a strong
$\delta$-neck\index{strong $\delta$-neck} in $(\widehat{\mathcal
M},\widehat G)$ contained in ${\mathcal H}$. Furthermore, there is a
point $y\in {\mathcal H}$ with $R(y)=h^{-2}$ with the property that
the central $2$-sphere of the $\delta$-neck centered at $y$ cuts off
an end of the ${\mathcal H}$ disjoint from $\Omega_\rho$. See {\sc
Fig.}~\ref{fig:deep}.
\end{thm}

\begin{proof}
The proof of the first statement is by contradiction. Fix $\rho>0$
and $\delta>0$ and suppose that there is no $0<h\le {\rm
min}(\rho\cdot\delta,\rho/2C)$ as required. Then there is a sequence
of generalized Ricci flows $({\mathcal M}_n,G_n)$ defined for $0\le
t<T_n<\infty$ satisfying Assumptions~\ref{assumptions} and points
$x_n\in {\mathcal M}_n$ with ${\bf t}(x_n)=T_n$ contained in
$2\epsilon$-horns ${\mathcal H}_n$ in $\Omega_n$ with boundary
contained in $(\Omega_n)_{\rho/2C}$ with $Q_n=R(x_n)\rightarrow
\infty$ as $n\rightarrow \infty$ but such that no $x_n$ is the
center of a strong $\delta$-neck in $({\mathcal M}_n, G_n)$. Form
the maximal extensions, $(\widehat{\mathcal M}_n,\widehat G_n)$, to
time $T$ of the $({\mathcal M}_n,G_n)$.

\begin{claim}\label{hyp}
The sequence $(\widehat{\mathcal M}_n,\widehat G_n,x_n)$ satisfies
the five hypothesis of Theorem~\ref{smlmtflow}.
\end{claim}

\begin{proof}
By our assumptions, Hypotheses (1) and (3) of
Theorem~\ref{smlmtflow} hold for this sequence.
 Also, we are assuming that any point $y\in {\mathcal M}_n$ with
 $R(y)\ge r_0^{-2}$ has a strong $(C,\epsilon)$-canonical neighborhood. Since
$R(x_n)=Q_n\rightarrow \infty$ as $n\rightarrow\infty$ this means
that for all $n$ sufficiently large, any point $y\in \widehat
{\mathcal M}_n$ with $R(y)\ge R(x_n)$ has a strong
$(C,\epsilon)$-canonical neighborhood. This establishes Hypothesis
(2) in the statement of Theorem~\ref{smlmtflow}.

Next, we have:

\begin{claim}
For any $A<\infty$ for all $n$ sufficiently large,
$B(x_n,0,AQ_n^{-1/2})$ is contained in the $2\epsilon$-horn
${\mathcal H}_n$ and has compact closure in ${\mathcal M}_n$.
\end{claim}

\begin{proof}
Any point  $z\in\partial {\mathcal H}_n$ has scalar curvature at
most $16C^2\rho^{-2}$ and there is a $2\epsilon$-neck centered at
$z$. This means that for all $y$ with
$d_{G_n}(z,y)<\epsilon^{-1}\rho/2C$ we have $R(y)\le
32C^2\rho^{-2}$. Hence, for all $n$ sufficiently large,
$d_{G_n}(x_n,z)>\epsilon^{-1}\rho/2C$, and thus
$d_{Q_nG_n}(x_n,z)>Q^{1/2}_n\epsilon^{-1}\rho/2C$. This implies
that, given $A<\infty$, for all $n$ sufficiently large, $z\not\in
B_{Q_nG_n}(x_n,0,A)$. Since this is true for all $z\in \partial
{\mathcal H}_n$, it follows that for all $n$ sufficiently large
$B_{Q_n,G_n}(x_n,0,A)\subset {\mathcal H}_n$. Next, we must show
that, for all $n$ sufficiently large, this ball has compact closure.
That is to say, we must show that for every $A$ for all $n$
sufficiently large the distance from $x_n$ to the end of the horn
${\mathcal H}_n$ is greater than $AQ_n^{-1/2}$. If not, then since
the curvature at the end of ${\mathcal H}_n$ goes to infinity for
each $n$, this sequence would violate Theorem~\ref{bcbd}.
\end{proof}

Because $B(x_n,0,AQ_n^{-1/2})$ is contained in a $2\epsilon$-horn,
it is $\kappa$-non-collapsed on scales $\le r$ for a universal
$\kappa>0$ and $r>0$. Also, because every point in the horn is the
center of a strong $2\epsilon$-neck, for every $n$ sufficiently
large and every $y\in B(x_n,0,AQ_n^{-1/2})$ the flow is defined on
an interval through $y$ defined for backward time $R(y)^{-1}$.

This completes the proof that all the hypotheses of
Theorem~\ref{smlmtflow} hold and establishes Claim~\ref{hyp}.
\end{proof}

We form a new sequence of generalized Ricci flows from the
$({\mathcal M}_n, G_n)$ by translating by $-T_n$, so that the final
time-slice is at ${\bf t}_n=0$, where ${\bf t}_n$ is the time
function for ${\mathcal M}_n$.

Theorem~\ref{smlmtflow} implies that, after passing to a
subsequence, there is a limit flow
$(M_\infty,g_\infty(t),(x_\infty,0)),\ t\in[-t_0,0])$ defined for
some $t_0>0$ for the  sequence $(Q_n\widehat{\mathcal
M}_n,Q_n\widehat G_n,x_n)$. Because of the curvature pinching toward
positive assumption, by Theorem~\ref{blowupposcurv}, the limit Ricci
flow has non-negative sectional curvature. Of course,
$R(x_\infty)=1$ so that the limit $(M_\infty,g_\infty(0))$ is
non-flat.

\begin{claim}
$(M_\infty,g_\infty(0))$ is isometric to the product $(S^2,h)\times
(\Ar,ds^2)$, where $h$ is a metric of non-negative curvature on
$S^2$ and $ds^2$ is the usual Euclidean metric on the real line.
\end{claim}

\begin{proof} Because of the fact that the $({\mathcal M}_n,G_n)$
have curvature pinched toward positive, and since $Q_n$ tend to
$\infty$ as $n$ tends to infinity, it follows that the geometric
limit $(M_\infty,g_\infty)$ has non-negative curvature. In
${\mathcal H}_n$ take a minimizing geodesic ray $\alpha_n$ from
$x_n$ to the end of ${\mathcal H}_n$ and a minimizing geodesic
$\beta_n$ from $x_n$ to $\partial {\mathcal H}_n$. As we have seen,
the lengths of both $\alpha_n$ and $\beta_n$ tend to $\infty$ as
$n\rightarrow \infty$. By passing to a subsequence, we can assume
that the $\alpha_n$ converge to a minimizing geodesic ray $\alpha$
in $(M_\infty,g_\infty)$ and that the $\beta_n$ converge to a
minimizing geodesic ray $\beta$ in $(M_\infty,g_\infty)$. Since, for
all $n$, the union of $\alpha_n$ and $\beta_n$ forms a piecewise
smooth ray in ${\mathcal H_n}$ meeting the central $2$-sphere of a
$2\epsilon$-neck centered at $x_n$ in a single point and at this
point crossing from one side of this $2$-sphere to the other, the
union of $\alpha$ and $\beta$ forms a proper, piecewise smooth map
of $\Ar$ to $M_\infty$ that meets  the central $2$-sphere of a
$2\epsilon$-neck centered at $x_\infty$ in a single point and
crosses from one side to the other at the point. This means that
$M_\infty$ has at least two ends. Since $(M_\infty,g_\infty)$ has
non-negative curvature, according to Theorem~\ref{splitting}, this
implies that $M_\infty$ is a product of a surface with $\Ar$. Since
$M$ has non-negative curvature, the surface has non-negative
curvature. Since $M$ has positive curvature at at least one point,
the surface is diffeomorphic to the $2$-sphere.
\end{proof}

According to Theorem~\ref{smlmtflow}, after passing to a subsequence
there is a limit flow defined on some interval of the form
$[-t_0,0]$ for $t_0>0$. Suppose that, after passing to a subsequence
there is a limit flow defined on $[-T,0]$ for some $0<T<\infty$. It
follows that for any $t\in [-T,0]$, the Riemannian manifold
$(M_\infty,g_\infty(t))$ is of non-negative curvature and has two
ends. Again by Theorem~\ref{splitting}, this implies that for every
$t\in [-T,0]$ the Riemannian manifold $(M_\infty,g_\infty(t))$ is a
Riemannian product of a metric of non-negative curvature on $S^2$
with $\Ar$. Thus, by Corollary~\ref{localprod} the Ricci flow is a
product of a Ricci flow  $(S^2,h(t))$ with the trivial flow on
$(\Ar,ds^2)$. It now follows from Corollary~\ref{rm>0} that for
every $t\in (-T,0]$ the curvature of $g_\infty(t)$ on $S^2$ is
positive.

Let $M_n$ be the zero time-slice of ${\mathcal M}_n$. Since
$(M_\infty,g_\infty,x_\infty)$ is the geometric limit of the
$(M_n,Q_nG_n(0),x_n)$, there is an exhausting sequence $x_\infty\in
V_1\subset V_2\subset\cdots$ of open subsets of $M_\infty$ with
compact closure and embeddings $\varphi_n\colon V_n\to M_n$ sending
$x_\infty$ to $x_n$ such that $\varphi_n^*(Q_nG_n(0))$ converges in
the $C^\infty$-topology, uniformly on compact sets, to $g_\infty$.

\begin{claim}
For any $z\in M_\infty$ for all $n$ sufficiently large, $z\in V_n$,
so that $\varphi_n(z)$ is defined. Furthermore, for all $n$
sufficiently large, there is a backward flow line through
$\varphi_n(z)$ in the generalized Ricci flow $(Q_n{\mathcal
M}_n,Q_nG_n)$ defined on the interval
$(-T-(R_{Q_nG_n}^{-1}(\varphi_n(z),0)/2),0]$. The scalar curvature
is bounded above on this entire flow line by $R(\varphi_n(z),0)$.
\end{claim}

\begin{proof}
Of course, for any compact subset $K\subset M_\infty$ and any $t'<T$
for all $n$ sufficiently large, $K\subset V_n$, and there is an
embedding $\varphi_n(K)\times [-t',0]\subset Q_n{\mathcal M}_n$
compatible with time and the vector field. The map $\varphi_n$
defines a map $Q_n^{-1}\varphi_n\colon K\times [-Q_n^{-1}t',0]\to
{\mathcal M}_n$. Since the scalar curvature of the limit is
positive, and hence bounded away from zero on the compact set
$K\times [-t',0]$ and since $Q_n\rightarrow \infty$ as $n$ tends to
infinity the following is true: For any compact subset $K\subset M$
and any $t'<T$, for all $n$ sufficiently large, the scalar curvature
of $G_n$ on the image $Q_n^{-1}\varphi_n(K)\times [-Q_n^{-1}t',0]$
is greater than $r_0^{-2}$, and hence for all $n$ sufficiently
large, every point in $Q_n^{-1}\varphi_n(K)\times [-Q_n^{-1}t',0]$
has a strong $(C,\epsilon)$-canonical neighborhood in ${\mathcal
M}_n$. Since having a strong $(C,\epsilon)$-canonical neighborhood
is invariant under rescaling, it follows that for all $n$
sufficiently large, every point of $\varphi_n(K)\times [-t',0]$ has
a strong $(C,\epsilon)$-canonical neighborhood.

Next we claim that, for all $n$ sufficiently large and for any $t\in
[-t',0]$, the point $(\varphi_n(z),t)$ is the center of a strong
$\epsilon$-neck. We have already seen that for all $n$ sufficiently
large $(\varphi_n(z),t)$ has a strong $(C,\epsilon)$-canonical
neighborhood. Of course, since $M_\infty$ is non-compact, for $n$
sufficiently large, the canonical neighborhood of $(\varphi_n(z),t)$
must either be a $(C,\epsilon)$-cap or a strong $\epsilon$-neck. We
shall rule out the possibility of a $(C,\epsilon)$-cap, at least for
all $n$ sufficiently large.

To do this,  take $K$ to be a neighborhood of  $(z,0)$ in the limit
$(M_\infty,g_\infty(0))$ with the topology of $S^2\times I$ and with
the metric being the product of a positively curved metric on $S^2$
with the Euclidean metric on $I$. We take $K$ to be sufficiently
large to contain the $2C$-ball centered at $(z,0)$. Because the
limit flow is the product of a positively curved flow on $S^2$ with
the trivial flow on $\Ar$, the flow is distance decreasing. Thus,
for every $t\in [-t',0]$ the submanifold $K\times\{t\}$ contains the
ball in $(M_\infty,g_\infty(t))$ centered at $(z,t)$ of radius $2C$.
 For every $n$
sufficiently large, consider the submanifolds $\varphi_n(K)\times
\{t\}$ of $(M_n,Q_nG_n(t))$.  Since the metrics
$\varphi_n^*Q_nG_n(t)$ are converging uniformly for all $t\in
[-t',0]$ to the product flow on $K$, for all $n$ sufficiently large
and any $t\in [-t',0]$, this submanifold contains the $C$-ball
centered at $(\varphi_n(z),t)$ in $(M_n,Q_nG_n(t))$. Furthermore,
the maximal curvature two-plane at any point of
$\varphi_n(K)\times\{t\}$ is almost tangent to the $S^2$-direction
of $K$. Hence, by Lemma~\ref{neckcurv} the central $2$-sphere of any
$\epsilon$-neck contained $\varphi_n(K)\times\{t\}$ is almost
parallel to the $S^2$-factors in the product structure on $K$ at
every point. This implies that the central $2$-sphere of any such
$\epsilon$-neck is isotopic to the $S^2$-factor of
$\varphi_n(K)\times\{t\}$. Suppose that $(\varphi_n(z),t)$ is
contained in the core of a $(C,\epsilon)$-cap ${\mathcal C}$. Then
${\mathcal C}$ is contained in $\varphi_n(K)\times\{t\}$. Consider
the $\epsilon$-neck $N\subset {\mathcal C}$ that is the complement
of the core of ${\mathcal C}$. Its central $2$-sphere, $\Sigma$, is
isotopic in $K$ to the $2$-sphere factor of $K$, but this is absurd
since $\Sigma$ bounds a $3$-ball in the ${\mathcal C}$. This
contradiction shows that for all $n$ sufficiently large and all
$t\in [-t',0]$, it is not possible for $(\varphi_n(z),t)$ to be
contained in the core of a $(C,\epsilon)$-cap. The only other
possibility is that for all $n$ sufficiently large and all $t\in
[-t',0]$ the point $(\varphi_n(z),t)$ is the center of a strong
$\epsilon$-neck in $(M_n,Q_nG_n(t))$.

Fix $n$ sufficiently large. Since, for all $t\in [-t',0]$, the point
$(\varphi_n(z),t)$ is the center of a strong $\epsilon$-neck, it
follows from Definition~\ref{defncanonnbhd} that for all $t\in
[-t',0]$ we have $R(\varphi_n(z),t)\le R(\varphi_n(z),0)$ (this
follows from the fact that the partial derivative in the
time-direction of the scalar curvature of a strong $\epsilon$-neck
of scale one is positive and bounded away from $0$). It also follows
from Definition~\ref{defncanonnbhd} that the flow near
$(\varphi_n(z),-t')$ extends backwards to time
$$-t'-R^{-1}_{Q_nG_n}(\varphi_n(z,t'))<-t'-R^{-1}_{Q_nG_n}(\varphi_n(z,0)),$$
with the same inequality for scalar curvature holding for all $t$ in
this extended interval. Applying this for $t'<T$ but sufficiently
close to $T$ establishes the last statement in the claim, and
completes the proof of the claim.
\end{proof}

Let $Q_0$ be the upper bound of the scalar curvature of
$(M_\infty,g_\infty(0))$.  By the previous claim, $Q_0$ is also an
upper bound for the curvature of $(M_\infty,g_\infty(-t'))$ for any
$t'<T$. Applying Theorem~\ref{smlmtflow} to the flows
$(M_n,Q_nG_n(t)),\, -t'-Q_0^{-1}/2< t\le -t'$, we conclude that
there is $t_0$ depending only on the bound of the scalar curvature
of $(M_\infty,g_\infty(-t'))$, and hence depending only on $Q_0$,
such that, after passing to a further subsequence the limit flow
exists for $t\in [-t'-t_0,-t']$. Since the limit flow already exists
on $[-t',0]$, we conclude that, for this further subsequence, the
limit flow exists on $[-t'-t_0,0]$. Now apply this with
$t'=T-t_0/2$. This proves that if, after passing to a subsequence,
there is a limit flow defined on $[-T,0]$, then, after passing to a
further subsequence there is a limit flow defined on $[-T-t_0/2,0]$
where $t_0$ depends only on $Q_0$, and in particular, is independent
of $T$. Repeating this argument with $T+(t_0/2)$ replacing $T$, we
pass to a further subsequence so that the limit flow is defined on
$[-T-t_0,0]$. Repeating this inductively, we can find a sequence of
subsequences so that for the $n$ subsequence the limit flow is
defined on $[-T-nt_0,0]$. Taking a diagonal subsequence produces a
subsequence for which the limit is defined on $(-\infty,0]$.

The limit flow is the product of a flow on $S^2$ of positive
curvature defined for $t\in (-\infty,0]$ and the trivial flow on
$\Ar$. Now, invoking Hamilton's result (Corollary~\ref{2DGSS}), we
see that the ancient solution of positive curvature on $S^2$ must be
a shrinking round $S^2$. This means that the limit flow is the
product of the shrinking round $S^2$ with $\Ar$, and implies that
for all $n$ sufficiently large there is a strong $\delta$-neck
centered at $x_n$. This contradiction proves the existence of $h$ as
required.

Now let us establish the last statement in Theorem~\ref{hexist}. The
subset of ${\mathcal H}$ consisting of all $z\in {\mathcal H}$ with
$R(z)\le \rho^{-2}$ is compact (since $R$ is a proper function), and
disjoint from any $\delta$-neck of scale $h$ since $h<\rho/2C$. On
the other hand, for any point $z\in {\mathcal H}$ with $R(z)\le
\rho^{-2}$ take a minimal geodesic from $z$ to the end of ${\mathcal
H}$. There must be a point $y$ on this geodesic with $R(y)=h^{-2}$.
The $\delta$-neck centered at $y$ is disjoint from $z$ (since
$h<\rho/2C$) and hence this neck separates $z$ from the end of
${\mathcal H}$. It now follows easily that there is a point $y\in
{\mathcal H}$ with $R(y)=h^{-2}$ and such that the central
$2$-sphere of the $\delta$-neck centered at $y$ divides ${\mathcal
H}$ into two pieces with the non-compact piece disjoint from
$\Omega_\rho$.
\end{proof}

\begin{cor}\label{hmonotone}
We can take the function $h(\rho,\delta)$ in the last lemma to be
$\le \delta\rho$, to be a weakly monotone non-decreasing function of
$\delta$ when $\rho$ is fixed, and to be a weakly monotone
non-decreasing function of $\rho$ when $\delta$ is held fixed.
\end{cor}

\begin{proof}
If $h$ satisfies the conclusion of Theorem~\ref{hexist} for $\rho$
and $\delta$ and if $\rho'\ge \rho$ and $\delta'\ge \delta$ then $h$
also satisfies the conclusion of Theorem~\ref{hexist} for $\rho'$
and $\delta'$. Also, any $h'\le h$ also satisfies the conclusion of
Theorem~\ref{hexist} for $\delta$ and $\rho$. Take a sequence
$(\delta_n,\rho_n)$ where each of the sequences $\{\delta_n\}$ and
$\{\rho_n\}$ is a monotone decreasing sequence with limit $0$. Then
we choose $h_n=h(\rho_n,\delta_n)\le \rho_n\delta_n$ as in the
statement of Theorem~\ref{hexist}. We of course can assume that
$\{h_n\}_n$ is a non-increasing sequence of positive numbers with
limit $0$. Then for any $(\rho,\delta)$ we take the largest $n$ such
that $\rho\ge \rho_n$ and $\delta\ge \delta_n$, and we define
$h(\rho,\delta)$ to be $h_n$ for this value of $n$. This constructs
the function $h(\delta,\rho)$ as claimed in the corollary.
\end{proof}

\chapter{The standard solution}\label{stdsolnsect}

The process of surgery involves making a choice of the metric on a
three-ball to `glue in'. In order to match approximatively with the
metric coming from the flow, the metric we glue in must be
asymptotic to the product of a round two-sphere and an interval near
the boundary. There is no natural choice for this metric; yet it is
crucial to the argument that we choose an initial metric so that the
Ricci flow with these initial conditions has several properties.
Conditions on the initial metric that  ensure the required
properties for the subsequence flow are contained in the following
definition.

\begin{figure}[ht]
  \relabelbox{
  \centerline{\epsfbox{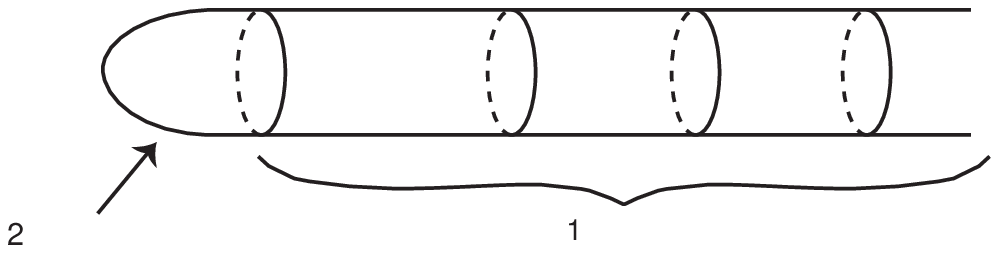}}}
  \relabel{1}{$S^2\times [0,\infty)$}
  \relabel{2}{positive curvature}
  \endrelabelbox
  \caption{A standard initial metric.}\label{fig:stdsoln}
\end{figure}

\begin{defn}\label{defn:stinmetr}
A {\em standard initial metric}\index{initial metric, standard|ii}
is a metric $g_0$ on $\Ar^3$ with the following properties:
\begin{itemize}
\item $g_0$ is a complete metric.
\item $g_0$ has non-negative sectional curvature at every point.
\item $g_0$ is invariant under the usual $SO(3)$-action on $\Ar^3$.
\item there is a compact ball $B\subset \Ar^3$ so
that the restriction of the metric $g_0$ to the complement of this
ball is isometric to the product $(S^2,h)\times (\Ar^+,ds^2)$ where
$h$ is the round metric of scalar curvature $1$ on $S^2$.
\item $g_0$ has constant sectional curvature $1/4$ near the origin. (This point will be
denoted $p$ and is called the {\em tip} of the initial metric.)
\end{itemize} See {\sc
Fig.}~\ref{fig:stdsoln}.
\end{defn}

Actually, one can work with an alternative weaker version of the
fourth condition, namely:

\noindent (iv) $g_0$ is asymptotic at infinity in the
$C^\infty$-topology to the product of the round metric $h_0$ on
$S^2$ of  scalar curvature $1$ with the usual metric $ds^2$ on the
real line. By this we mean  that if $x_n\in \Ar^3$ is any sequence
converging to infinity, then the based Riemannian manifolds
$(\Ar^3,g_0,x_n)$ converge smoothly to $(S^2,h_0)\times (\Ar,ds^2)$.
But we shall only use standard initial metrics as given in
Definition~\ref{defn:stinmetr}.

\begin{lem}\label{stdinit}
There is a standard initial metric.
\end{lem}

\begin{proof}
We construct our Riemannian manifold as follows. Let
$(x_0,x_1,x_2,x_3)$ be Euclidean coordinates on $\Ar^4$. Let
$y=f(s)$ be a function defined for $s\ge 0$ and satisfying:
\begin{enumerate}
\item[(1)] $f$ is $C^\infty$ on $(0,\infty)$
\item[(2)] $f(s)>0$ for all $s>0$.
\item[(3)] $f''(s)\le 0$ for all $s>0$.
\item[(4)] There is $s_1>0$ such that $f(s)=\sqrt{2}$ for all $s\ge
s_1$.
\item[(5)] There is $s_0>0$ such that
$f(s)=\sqrt{4s-s^2}$ for all $s\in [0,s_0]$.
\end{enumerate}
Given such a function $f$, consider the graph
$$\Gamma=\{(x_0,x_1)\bigl|\bigr. x_0\ge 0 \text{ and }
x_1=f(x_0)\}$$ in the $(x_0,x_1)$-plane. We define $\Sigma(f)$ by
rotating $\Gamma$ about the $x_0$-axis in four-space:
$$\Sigma(f)=\{(x_0,x_1,x_2,x_3)\bigl|\bigr. x_0\ge 0 \text{ and }
x_1^2+x_2^2+x_3^2=f(x_0)^2\}.$$ Because of the last condition on
$f$, there is a neighborhood of $0\in \Sigma(f)$ that is isometric
to a neighborhood of the north pole in the three-sphere of radius
$2$. Because of this and the first item, we see that $\Sigma(f)$ is
a smooth submanifold of $\Ar^4$. Hence, it inherits a Riemannian
metric $g_0$. Because of the fourth item, a neighborhood of infinity
of $(\Sigma(f),g_0)$ is isometric to $(S^2,h)\times (0,\infty)$, and
in particular, $(\Sigma(f),g_0)$ is complete. Clearly, the rotation
action of $S0(3)$ on $\Sigma(f)$, induced by the orthogonal action
on the last three coordinates in $\Ar^4$, is an isometric action
with the origin as the only fixed point. It is also clear that
$\Sigma(f)$ is diffeomorphic to $\Ar^3$ by a diffeomorphism that
send the $SO(3)$ action to the standard one on $\Ar ^3$.

It remains to compute the sectional curvatures of $g_0$. Let $q\in
\Sigma(f)$ be a point distinct from the fixed point of the
$SO(3)$-action. Direct computation shows that the tangent plane to
the two-dimensional $SO(3)$-orbit through $q$ is a principal
direction for the curvature, and the sectional curvature on this
 tangent two-plane  is given by
$$\frac{1}{f(q)^2(1+f'(q)^2)}.$$
On the subspace in $\wedge^2T_q\Sigma(f)$ perpendicular to the line
given by this two-plane, the curvature is constant  with eigenvalue
$$\frac{-f''(q)}{f(q)(1+f'(q)^2)^2}.$$
Under our assumptions about $f$, it is clear that $\Sigma(f)$ has
non-negative curvature and has constant sectional curvature $1/4$
near the origin. It remains to choose the function $f$ satisfying
Items (1) -- (5) above.

Consider the function $h(s)=(2-s)/\sqrt{4s-s^2}$. This function is
integrable from $0$ and the definite integral from  zero to $s$ is
equal to $\sqrt{4s-s^2}$. Let $\lambda(s)$ be a non-increasing
$C^\infty$-function defined on $[0,1/2]$, with $\lambda$ identically
one near $0$ and identically equal to $0$ near $1/2$. We extend
$\lambda$ to be identically $1$ for $s<0$ and identically $0$ for
$s>1/2$. Clearly,
$$\int_0^2h(s)\lambda(s-3/2)ds>\int_0^{3/2}h(s)ds>\sqrt{2}$$
and $$\int_0^2h(s)\lambda(s)ds<\int_0^{1/2}h(s)<\sqrt{2}.$$ Hence,
for some $s_0\in (1/2,3/2)$ we have
$$\int_0^2h(s)\lambda(s-s_0)ds=\sqrt{2}.$$
We define $$f(s)=\int_0^sh(\sigma)\lambda(\sigma-s_0)d\sigma.$$ It
is easy to see that $f$ satisfies all the above conditions.
\end{proof}

The following lemma is clear from the construction.

\begin{lem}\label{A_0}
There is $A_0<\infty$ such that
$$(\Ar^3\setminus B(0,A_0),g(0))$$
is isometric to the product of a round metric on $S^2$ of scalar
curvature $1$ with the Euclidean metric on $[0,\infty)$. There is a
constant $K<\infty$ such that the volume of $B_{g(0)}(0,A_0)$ is at
most $K$. Furthermore, there is a constant $D<\infty$ so that the
scalar curvature of standard initial metric $(\Ar^3,g(0))$ is
bounded above by $D$ and below by $D^{-1}$.
\end{lem}

\subsection{Uniqueness and properties: The statement}

{\bf Fix once and for all a standard initial metric $g_0$ on
$\Ar^3$.}

\begin{defn}
A {\em partial standard Ricci flow}\index{standard Ricci
flow!partial|ii} is a Ricci flow $(\Ar^3,g(t)),\ 0\le t<T$, such
that $g(0)=g_0$ and such that the curvature is locally bounded in
time. We say that a partial standard Ricci flow is a {\em standard
Ricci flow}\index{standard Ricci flow|ii} if it has the property
that $T$ is maximal in the sense that there is no extension of the
flow to a flow on a larger time interval $[0,T')$ with $T'>T$ with
the property that the extension has curvature locally bounded in
time.
\end{defn}

Here is the main result of this chapter.

\begin{thm}\label{stdsoln} There is a standard Ricci flow defined for some
positive amount of time. Let $(\Ar^3,g(t)),\ 0\le t<T$, be a
standard Ricci flow. Then the following hold.
\begin{enumerate}
\item[(1)] {\bf (Uniqueness):}
If $(\Ar^3,g'(t)),\ 0\le t<T'$, is a standard Ricci flow, then
$T'=T$ and $g'(t)=g(t)$.
\item[(2)] {\bf (Time Interval):} $T=1$.
\item[(3)] {\bf (Positive curvature):} For each $t\in (0,1)$ the metric $g(t)$ on $\Ar^3$ is complete
of strictly positive curvature.
\item[(4)] {\bf ($SO(3)$-invariance):} For each $t\in [0,T)$ the
Riemannian manifold $(\Ar^3,g(t))$  is invariant under the
$SO(3)$-action on $\Ar^3$.
\item[(5)] {\bf (Asymptotics at $\infty$):}
For any $t_0<1$ and any $\epsilon>0$ there is a compact subset $X$
of $\Ar^3$ such that for any $x\in \Ar ^3\setminus X$ the
restriction of the standard flow to an appropriate neighborhood of
$x$ for time $t\in [0,t_0]$ is within $\epsilon$ in the
$C^{[1/\epsilon]}$-topology of the product Ricci flow $(S^2\times
(-\epsilon^{-1},\epsilon^{-1})),h(t)\times ds^2,\ 0\le t\le t_0$,
where $h(t)$ is the round metric with scalar curvature $1/(1-t)$ on
$S^2$.
\item[(6)] {\bf (Non-collapsing):} There are $r>0$ and $\kappa>0$ such that $(\Ar^3,g(t)),\ 0\le
t<1$, is $\kappa$-non-collapsed on scales less than $r$.
\end{enumerate}
\end{thm}

The proof of this result occupies the next few subsections. All the
properties except the uniqueness are fairly straightforward to
prove. We establish uniqueness by reducing the Ricci flow to the
Ricci-DeTurck by establishing the existence of a solution to the
harmonic map flow in this case. This technique can be made to work
more generally in the case of complete manifolds of bounded
curvature, see \cite{ChenZhu}, but we preferred to give the more
elementary argument that covers this case, where the symmetries
allow us to reduce the existence problem for the harmonic map flow
to a problem that is the essentially one-dimensional. Also, in the
rest of the argument one does not need uniqueness, only a
compactness result for the space of all Ricci flows of bounded
curvature on each time-slice with the given initial conditions.
Kleiner and Lott pointed out to us that this uniqueness can be
easily derived from the other properties by arguments similar to
those used to establish the compactness of the space of
$\kappa$-solutions.

\section{Existence of a standard flow}

For any $R<\infty$, denote by $B_R\subset \Ar^3$, the ball of radius
$R$ about the origin in the metric $g_0$. For $R\ge A_0+1$,  a
neighborhood of the boundary of this ball is isometric to
$(S^2,h)\times ([0,1],ds^2)$. Thus, in this case, we can double the
ball, gluing the boundary to itself by the identity, forming a
manifold we denote by $S^3_R$. The doubled metric will be a smooth
Riemannian metric $g_R$ on $S^3_R$. Let $p\in S^3_R$ be the image of
the origin in the first copy of $B_R$. Now take a sequence, $R_n$,
tending to infinity to construct based Riemannian manifolds
$(S^3_{R_n},g_{R_n},p)$ that converge geometrically to
$(\Ar^3,g_0,p)$. For each $n$, let $(S^3_{R_n},g_{R_n}(t)),\ 0\le
t<T_n$ be maximal Ricci flow with $(S^3_{R_n},g_{R_n})$ as initial
metric.
 The maximum principle
applied to Equation~(\ref{Revol}), $\partial R/\partial t=\triangle
R+|{\rm Ric}|^2$, then implies by Proposition~\ref{fordiffmax} that
the maximum of $R$ at time $t$, $R_{\rm max}(t)$ obeys the
inequality $\partial R_{\rm max}/\partial t\le R_{\rm max}(t)^2$,
and integrating this inequality (i.e., invoking
Lemma~\ref{fordiffquot}) one finds a positive constants $t_0$ and
$Q_0$ such that for each $n$, the norm of the scalar curvature of
$g_{R_n}(t)$ are bounded by $Q_0$ on the interval $[0,{\rm
max}(t_0,T_n))$. By Corollary~\ref{rm>0}, for each $n$ the sectional
curvature of the flow $(S^3_{R_n},g_{R_n}(t)),\ 0\le t<T_n$ is
non-negative, and hence the sectional curvature of  this flow is
also bounded by $Q_0$ on $[0,{\rm max}(t_0,T_n))$. It now follows
from Proposition~\ref{flowextend} and the fact that the $T_n$ are
maximal that $T_n>t_0$ for all $n$.
 Since the Riemann
curvatures of the $(S^3_{R_n},g_{R_n}(t)),\ 0\le t<t_0$, are bounded
independent of $n$, and since the $(S^3_{R_n},g_{R_n},p)$ converge
geometrically to $(\Ar^3,g_0,p)$, it follows from
Theorem~\ref{flowlimit} that there is a geometric limiting flow
defined on $[0,t_0)$. Since this flow is the geometric limit of
flows of uniformly bounded curvature, it has uniformly bounded
curvature. Taking a maximal extension of this flow to one of locally
bounded curvature gives a standard flow.

\section{Completeness, positive curvature, and asymptotic behavior}

 Let $\left( \mathbb{R}^{3},g(t)
\right)$, $t\in [0,T)$, be a partial standard solution. Let
$y_{i}\rightarrow\infty$ be the sequence of points in
$\mathbb{R}^{3}$ converging to infinity. From the  definition we see
that the based Riemannian manifolds $\left(
\mathbb{R}^{3},g_{0},y_{i}\right)$ converge smoothly to
$(S^2\times\mathbb{R},h(0)\times ds^2)$ where $h(0)$ is the round
metric of scalar curvature $1$ on $S^2$.

Let us begin by proving the third item in the statement of
Theorem~\ref{stdsoln}:

\begin{lem}\label{third}
For each $t_0\in [0,T)$ the Riemannian manifold $(\Ar^3,g(t_0))$ is
complete and of positive curvature.
\end{lem}

\begin{proof}
Fix $t_0\in [0,T)$. By hypothesis $(\Ar^3,g(t)),\ 0\le t\le t_0$ has bounded
curvature. Hence, there is a constant $C<\infty$ such that $ g(0)\le Cg(t_0)$,
so that for any points $x,y\in \Ar^3$, we have $d_0(x,y)\le
\sqrt{C}d_{t_0}(x,y)$. Since $g(0)$ is complete, this implies that $g(t_0)$ is
also complete.

Now let us show that $(M,g(t_0))$ has non-negative curvature. Here,
the argument is the analogue of the proof of Corollary~\ref{quad}
with one additional step, the use of a function $\varphi$ to
localize the argument. Suppose this is false, i.e., suppose that
there is $x\in M$ with ${\rm Rm}(x,t_0)$ having an eigenvalue less
than zero. Since the restriction of the flow to $[0,t_0]$ is
complete and of bounded curvature, according to
\cite{Hamiltonlimits} for any constants $C<\infty$ and $\eta>0$ and
any compact subset $K\subset M\times [0,t_0]$ there is $\epsilon>0$
and a function $\varphi\colon M\times [0,t_0]\to \Ar$ with the
following properties:
\begin{enumerate}
\item[(1)] $\varphi|_K\le \eta$.
\item[(2)] $\varphi\ge \epsilon$ everywhere.
\item[(3)] For each $t\in [0,t_0]$ the restriction of $\varphi$  to $M\times \{t\}$
goes to infinity at infinity in the sense that  for any $A<\infty$
the pre-image $\varphi^{-1}([0,A]\cap (M\times\{t\})$ is compact.
\item[(4)] On all of $M\times [0,t_0]$ we have
$
  \left(\frac{\partial}{\partial t}-\triangle\right)\varphi\ge C\varphi.
$
\end{enumerate}
Recall from Section~\ref{4.3.3} that ${\mathcal T}$ is the curvature
tensor written with respect to an evolving orthonormal frame
$\{F_\alpha\}$ for the tangent bundle. Consider the symmetric,
horizontal two-tensor $\widehat{\mathcal T}={\mathcal T}+\varphi g$.
Let $\hat\mu(x,t)$ denote the smallest eigenvalue of this symmetric
two-tensor at $(x,t)$. Clearly, since the curvature is bounded, it
follows from the third property of $\varphi$ that for each $t\in
[0,t_0]$ the restriction of $\hat\mu$ to $M\times\{t\}$ goes to
infinity at infinity in $M$. In particular, the subset of $(x,t)\in
M\times [0,t_0]$ with the property that $\hat\mu(x,t)\le
\hat\mu(y,t)$ for all $y\in M$ is a compact subset of $M\times
[0,t_0]$. It follows from Proposition~\ref{fordiffmax} that
$f(t)={\rm min}_{x\in M}\hat\mu(x,t)$ is a continuous function of
$t$. Choosing $\eta>0$ sufficiently small and $K$ to include
$(x,t_0)$, then $\widehat{\mathcal T}$ will have a negative
eigenvalue at $(x,t_0)$. Clearly, it has only positive eigenvalues
on $M\times \{0\}$. Thus, there is $0<t_1<t_0$ so that
$\widehat{\mathcal T}$ has only positive eigenvalues on $M\times
[0,t_1)$ but has a zero eigenvalue at $(y,t_1)$ for some $y\in M$.
That is to say, ${\mathcal T}\ge -\varphi g$ on $M\times [0,t_1]$.
Diagonalizing ${\mathcal T}$ at any point $(x,t)$ with $t\le t_1$,
all its eigenvalues are at least $-\varphi(x,t_1)$. It follows
immediately that on $M\times [0,t_1]$ the smallest eigenvalue of the
symmetric form ${\mathcal T}^2+{\mathcal T}^\#$ is bounded below by
$2\varphi$. Thus, choosing $C\ge 4$  we see that for $t\le t_1$
every eigenvalue of ${\mathcal T}^2+{\mathcal T}^\#$ is at least
$-C\varphi/2$.

We compute the evolution equation using the formula in
Lemma~\ref{quad} for the evolution of ${\mathcal T}$ in an evolving
orthonormal frame:
\begin{eqnarray*}\frac{\partial \widehat{\mathcal T}}{\partial t} & = &
\frac{\partial{\mathcal T}}{\partial t}+\frac{\partial
\varphi}{\partial t}g-2\varphi {\rm Ric}(g) \\
& = & \triangle {\mathcal T}+{\mathcal T}^2+{\mathcal T}^\#+\frac{\partial \varphi}{\partial t}g-2\varphi {\rm Ric}(g) \\
& = & \triangle \widehat{\mathcal T}+{\mathcal T}^2+{\mathcal
T}^\#+\left(\frac{\partial \varphi}{\partial t}-\triangle
\varphi\right)g-2\varphi {\rm Ric}(g)\\
&\ge & \triangle \widehat{\mathcal T}+{\mathcal T}^2+{\mathcal
T}^\#+\left(Cg-2{\rm Ric}(g)\right)\varphi.
\end{eqnarray*}
Since every eigenvalue of ${\mathcal T}^2+{\mathcal T}^\#$  on
$M\times [0,t_1]$ is at least $-C\varphi/2$, it follows that on
$M\times [0,t_1]$
$$\frac{\partial \widehat{\mathcal T}}{\partial t}\ge
 \triangle \widehat{\mathcal T}+(Cg/2-2{\rm Ric}_g)\varphi.$$ Once again assuming that $C$ is
sufficiently large, we see that for any $t\le t_1$
$$\frac{\partial \widehat{\mathcal T}}{\partial t}\ge \triangle
\widehat {\mathcal T}.$$ Thus, at any local minimum $x\in M$ for
$\hat\mu(\cdot,t)$, we have
$$\frac{\partial \hat\mu}{\partial t}\ge 0.$$
This immediately implies by Proposition~\ref{fordiffmax} that
$\psi(t)={\rm min}_{x\in M}\hat\mu(x,t)$ is a non-decreasing
function of $t$. Since its value at $t=0$ is at least $\epsilon>0$
and its value at $t_1$ is zero, this is a contradiction. This
establishes that the solution has non-negative curvature everywhere.
Indeed, by Corollary~\ref{nullspace} it has strictly positive
curvature for every $t>0$.
\end{proof}

Now let us turn to the asymptotic behaviour of the flow.

Fix $T'<T$. Let $y_k$ be a sequence tending to infinity in
$(\Ar^3,g_0)$. Fix $R<\infty$. Then there is $k_0(R)$ such that for
all $k\ge k_0(R)$ there is an isometric embedding $\psi_k\colon
(S^2,h)\times (-R,R)\to (\Ar^3,g_0)$ sending $(x,0)$ to $y_k$. These
maps realize the product $(S^2,h)\times (\Ar,ds^2)$ as the geometric
limit of the $(\Ar^3,g_0,y_i)$. Furthermore, for each $R<\infty$
there is a uniform $C^\infty$ point-wise bound to the curvatures of
$g_0$ restricted to the images of the $\psi_k$ for $k\ge k_0(R)$.
 Since the flow $g(t)$ has bounded
curvature on $\Ar^3\times [0,T']$, it follows from
Theorem~\ref{shiw/deriv} that there are uniform $C^\infty$
point-wise bounds for the curvatures of $g(t)$ restricted to
$\psi_k(S^2\times (-R,R))$. Thus,
 by Theorem~\ref{partialflowlimit}, after passing to a subsequence, the flows
$\psi_k^*g(t)$ converge to a  limiting flow on $S^2\times \Ar$. Of
course, since the curvature of $g(t)$ is everywhere $\ge 0$, the
same is true of this limiting flow. Since the time-slices of this
flow have two ends, it follows from Theorem~\ref{splitting} that
every manifold in the flow is a product of a compact surface with
$\Ar$. According to Corollary~\ref{nullspace} this implies that the
flow is the product $(S^2,h(t))\times (\Ar,ds^2)$. This means that
given $\epsilon>0$, for all $k$ sufficiently large, the restriction
of the flow to the cylinder of length $2R$ centered at $y_k$  is
within $\epsilon$ in the $C^{[1/\epsilon]}$-topology of the
shrinking cylindrical flow on time $[0,T']$. Given $\epsilon>0$ and
$R<\infty$ this statement is true for all $y$ outside a compact ball
$B$ centered at the origin.

We have now established the following

\begin{prop}\label{item5}
Given $T'<T$ and $\epsilon>0$ there is a compact ball $B$ centered
at the origin of $\Ar^3$ such that the restriction of the flow
$(\Ar^3\setminus B,g(t)),\ 0\le t\le T'$, is within $\epsilon$ in
the $C^{[1/\epsilon]}$-topology of the standard evolving cylinder
$(S^2,h(t))\times (\Ar^+,ds^2)$.
\end{prop}

\begin{cor}\label{1}
The maximal time $T$ is $\le 1$.
\end{cor}

\begin{proof}
If $T>1$, then we can apply the above result to $T'$ with
$1<T'<T$,and see that the solution at infinity is asymptotic to the
evolving cylinder $(S^2,h(t))\times (\Ar,ds^2)$ on the time interval
$[0,T']$. But this is absurd since this evolving cylindrical flow
becomes completely singular at time $T=1$.
\end{proof}

\section{Standard solutions are rotationally symmetric}

Next, we consider the fourth item in the statement of the theorem.
Of course, rotational symmetry would follow immediately from
uniqueness. But here we shall use the rotational symmetry to reduce
the uniqueness problem to a one-dimensional problem which we then
solve. One can also use the general uniqueness theorem for complete,
non-compact manifolds due to Chen and Zhu (\cite{ChenZhu}), but we
have chosen to present a more elementary, self-contained argument in
this special case which we hope will be more accessible.

Let ${\rm Ric}_{ij}$ be the Ricci tensor and ${\rm Ric}^i_k= g^{ij}
{\rm Ric}_{jk}$ be the dual tensor. Let $X$ be a vector field
evolving by
\begin{equation}\label{eq1}\frac{\partial}{\partial t}X=\triangle X+{\rm Ric}(X,\cdot)^*.
\end{equation}
 In local coordinates $(x^1,\ldots,x^n)$, if
 $X=X^i\partial_i$, then the equation becomes
\begin{equation}
\frac{\partial}{\partial t}X^{i}=(\Delta X)^{i}+{\rm
Ric}_{k}^{i}X^{k} \label{eq.3}
\end{equation}
Let $X^*$ denote the dual one-form to $X$. In local coordinates we
have $X^*=X^*_idx^i$ with $X^*_i=g_{ij}X^j$. Since the evolution
equation for the metric is the Ricci flow, the evolution equation
for $X^*$ is
$$\frac{\partial X^*}{\partial t}=\triangle X^*-{\rm Ric}(X,\cdot),$$
or in local coordinates
$$\frac{\partial X^*_i}{\partial t}=(\triangle X^*)_i-{\rm Ric}_{ij}X^j.$$

\begin{lem}\label{12.9}
With $X$ and its dual $X^*$ evolving by the above equations, set
$V=\nabla X^*$, so that $V$ is a contravariant two-tensor. In local
coordinates we have $V=V_{ij}dx^i\otimes dx^j$ with
$$V_{ij}=(\nabla_iX)_j=g_{jk}(\nabla_iX)^k.$$
This symmetric two-tensor satisfies
\begin{equation}
\frac{\partial}{\partial t}V=\triangle
V-\left(2{{R_k}^{rl}}_jV_{rl}+{\rm Ric}_k^lV_{lj} +{\rm
Ric}_j^lV_{kl}\right)dx^k\otimes dx^j. \label{eq.4}
\end{equation}
\end{lem}

\begin{rem}
The covariant derivative acts on one-forms $\omega$ in such a way
that the following equation holds:
$$\langle\nabla(\omega),\xi\rangle=\langle
\omega,\nabla(\xi)\rangle$$ for every vector field $\xi$. This means
that in local coordinates we have
$$\nabla_{\partial_r}(dx^k)=-\Gamma_{rl}^kdx^l.$$ Similarly, the Riemann
curvature acts on one-forms $\omega$ satisfying
$${\rm Rm}(\xi_1,\xi_2)(\omega)(\xi)=-\omega\left({\rm Rm}(\xi_1,\xi_2)(\xi)\right).$$
Recall that in local coordinates  $$R_{ijkl}=\langle {\rm
Rm}(\partial_i,\partial_j)(\partial_l),\partial_k\rangle.$$ Thus, we
have
$${\rm Rm}(\partial_i,\partial_j)(dx^k)=-g^{ka}R_{ijal}dx^l=-{{R_{ij}}^k}_ldx^l,$$
where as usual we use the inverse metric tensor to raise the index.

Also, notice that  $\Delta X_{i}-{\rm
Ric}_{ik}X^{k}=-\Delta_{d}X_{i\text{ }}$, where by $\Delta_d$ we
mean the Laplacian associated to the operator $d$ from vector fields
to one-forms with values in the vector field. Since
\begin{align*}
 -\left(  d\delta+\delta d\right)  X_{i}
&  =-\nabla_{i}\left(  -\nabla^{k}X_{k}\right)  -\left(
-\nabla^{k}\right)
\left(  \nabla_{k}X_{i}-\nabla_{i}X_{k}\right) \\
&
=\nabla_{i}\nabla^{k}X_{k}+\nabla^{k}\nabla_{k}X_{i}-\nabla^{k}\nabla
_{i}X_{k}\\
&  ={{{R_{i}}^k}_k}^jX_{j}+\nabla^{k}\nabla_{k}X_{i}=\Delta
X_{i}-{\rm Ric}_i^jX_{j}.
\end{align*}
\end{rem}

\begin{proof}(of Lemma~\ref{12.9})
The computation is routine, if complicated. We make the computation
at a point $(p,t)$ of space-time. We fix local $g(t)$-Gaussian
coordinates $(x^1,\ldots,x^n)$ centered at $p$ for space, so that
the Christoffel symbols vanish at $(p,t)$.

 We
compute
\begin{eqnarray}\nonumber
\frac{\partial}{\partial t}V  &  = &\frac{\partial}{\partial
t}\left( \nabla X^* \right) =-\left( \frac{\partial}{\partial
t}\Gamma_{kj}^{l}\right)  X^*_{l}dx^k\otimes dx^j +\nabla\left(
\frac{\partial}{\partial t}X^*\right) \\ \nonumber
 &= &
\left(-\nabla^{l}{\rm Ric}_{kj}+\nabla_{k}{\rm
Ric}_j^l+\nabla_{j}{\rm Ric}_k^l\right)X^*_{l}dx^k\otimes
dx^j \\
& & +\nabla\left(  \Delta X^*-{\rm Ric}(X,\cdot)\right)
.\label{eqn1}
\end{eqnarray}

We have \begin{eqnarray*} \nabla(\triangle X^*) & = &
\nabla\left((g^{rs}\left(\nabla_r\nabla_s(X^*)-\Gamma_{rs}^l\nabla_lX^*\right)\right)
\\
&= &
g^{rs}\left(\nabla\left(\nabla_r\nabla_s(X^*)-\Gamma_{rs}^l\nabla_lX^*\right)\right).
\end{eqnarray*}

Let us recall the formula for commuting $\nabla$ and $\nabla_r$. The
following is immediate from the definitions.

\begin{claim}
For any tensor $\phi$ we have
$$\nabla(\nabla_r\phi)=\nabla_r(\nabla\phi)+dx^k\otimes
{\rm Rm}(\partial_k,\partial_r)(\phi)-\nabla_r(dx^l)\otimes
\nabla_l(\phi).$$
\end{claim}

Applying this to our formula gives
\begin{eqnarray*}\nabla(\triangle
X^*) & = &  g^{rs}\left(\nabla_r\nabla\nabla_sX^*+dx^k\otimes {\rm
Rm}(\partial_k,\partial_r)(\nabla_sX^*)-\nabla_r(dx^l)\otimes
\nabla_l\nabla_sX^* \right. \\
& & \left. -\nabla(\Gamma_{rs}^l\nabla_lX^*)\right) .
\end{eqnarray*}
 Now we
apply the same formula to commute $\nabla$ and $\nabla_s$. The
result is
\begin{eqnarray*}
\nabla(\triangle X^*) & = & g^{rs}\left(\nabla_r\nabla_s\nabla
X^*+\nabla_r\left(dx^k\otimes {\rm
Rm}(\partial_k,\partial_s)X^*-\nabla_sdx^l\otimes \nabla_lX^*\right)\right. \\
& & \left. +dx^k\otimes {\rm
Rm}(\partial_k,\partial_r)(\nabla_sX^*)-\nabla_r(dx^l)\otimes
\nabla_l\nabla_sX^*-\nabla(\Gamma_{rs}^l\nabla_lX^*)\right).
\end{eqnarray*}

Now we expand
\begin{eqnarray*}
\lefteqn{\nabla_r\left(dx^k\otimes {\rm
Rm}(\partial_k,\partial_s)X^*-\nabla_sdx^l\otimes
\nabla_lX^*\right)}
\\ & = & \nabla_r(dx^k)\otimes
{\rm Rm}(\partial_k,\partial_s)X^*+dx^k\otimes\nabla_r({\rm
Rm}(\partial_k,\partial_s))X^*
\\
& & + dx^k\otimes {\rm
Rm}(\partial_k,\partial_s)\nabla_rX^*-\nabla_r\nabla_sdx^l\otimes\nabla_lX^*-\nabla_sdx^l\otimes
\nabla_r\nabla_lX^*
\end{eqnarray*}
Invoking the fact that the Christoffel symbols vanish at the point
of space-time where we are making the computation, this above
expression simplifies to
\begin{eqnarray*}
\lefteqn{\nabla_r\left(dx^k\otimes {\rm
Rm}(\partial_k,\partial_s)X^*-\nabla_sdx^l\otimes
\nabla_lX^*\right)}
\\ & = & dx^k\otimes\nabla_r({\rm Rm}(\partial_k,\partial_s))X^* + dx^k\otimes
{\rm
Rm}(\partial_k,\partial_s)\nabla_rX^*-\nabla_r\nabla_sdx^l\otimes\nabla_lX^*.
\end{eqnarray*}

Also, expanding and using the vanishing of the Christoffel symbols
we have
\begin{eqnarray*}
-\nabla(\Gamma_{rs}^l\nabla_lX^*)& = &
-d\Gamma_{rs}^l\otimes\nabla_lX^*-\Gamma_{rs}^l\nabla\nabla_lX^* \\
& = & -d\Gamma_{rs}^l\otimes\nabla_lX^*.
\end{eqnarray*}

Plugging these computations into equation above and using once more the
vanishing of the Christoffel symbols gives
\begin{eqnarray*}
\nabla(\triangle X^*) & = & \triangle(\nabla X^*)+g^{rs}\left(
dx^k\otimes\nabla_r({\rm Rm}(\partial_k,\partial_s))X^* +
dx^k\otimes
{\rm Rm}(\partial_k,\partial_s)\nabla_rX^* \right.\\
& & \left.
  -\nabla_r\nabla_sdx^l\otimes\nabla_lX^* +dx^k\otimes
{\rm Rm}(\partial_k,\partial_r)(\nabla_sX^*)
 -d\Gamma_{rs}^l\otimes\nabla_lX^*\right).
\end{eqnarray*}

 Now by the symmetry of $g^{rs}$ we can amalgamate the second and
 fourth terms on the right-hand side to give
\begin{eqnarray*}
\nabla(\triangle X^*) & = & \triangle(\nabla X^*)+g^{rs}\left(
dx^k\otimes\nabla_r({\rm Rm}(\partial_k,\partial_s))X^* \right. \\
 & & \left. + 2dx^k\otimes
{\rm Rm}(\partial_k,\partial_s)\nabla_rX^*
-\nabla_r\nabla_sdx^l\otimes\nabla_lX^*
 -d\Gamma_{rs}^l\otimes\nabla_lX^*\right).
\end{eqnarray*}

We expand
$${\rm Rm}(\partial_k,\partial_s)\nabla_rX^* =-{{R_{ks}}^{l}}_jV_{rl}dx^j.$$
Also we have (again using the vanishing of the Christoffel symbols)
\begin{eqnarray*}
-\nabla_r\nabla_sdx^l-d\Gamma_{rs}^l & = &
\nabla_r\Gamma_{ks}^ldx^k-\partial_k\Gamma_{rs}^ldx^k \\
& = & {{R_{rk}}^l}_s dx^k.
\end{eqnarray*}

Lastly,
\begin{eqnarray*}
\nabla_r({\rm
Rm}(\partial_k,\partial_s))X^*=-{{(\nabla_rR)_{ks}}^l}_jX^*_ldx^j.
\end{eqnarray*}

Plugging all this in and raising indices yields
\begin{eqnarray*}
\nabla(\triangle X^*) & = & \triangle(\nabla X^*)-
g^{rs}{{(\nabla_rR)_{ks}}^l}_jX^*_ldx^k\otimes
dx^j-2{{R_k}^{rl}}_jV_{rl}dx^k\otimes dx^j \\
& &  +
g^{rs}{{R_{rk}}^l}_sV_{lj}dx^k\otimes dx^j \\
& = & \triangle(\nabla X^*)- g^{rs}{{(\nabla_rR)_{ks}}^j}_lX^*_jdx^k\otimes
dx^l-2{{R_k}^{rl}}_jV_{rl}dx^k\otimes dx^j  \\
& & -{\rm Ric}_k^lV_{lj}dx^k\otimes dx^j.
\end{eqnarray*}

Thus, we have
\begin{eqnarray*}
\lefteqn{\nabla(\triangle X^*) -\nabla({\rm Ric}(X,\cdot)^*)=} \\ &  &
\triangle(\nabla X^*)- g^{rs}{{(\nabla_rR)_{ks}}^l}_jX^*_ldx^k\otimes
dx^j-2{{R_k}^{rl}}_jV_{rl}dx^k\otimes dx^j
\\
& & -\left({\rm Ric}_k^lV_{lj}+\nabla_k({\rm Ric})_j^lX^*_l+{\rm
Ric}_j^lV_{kl}\right)dx^k\otimes dx^j,
\end{eqnarray*}
and consequently, plugging back into Equation~(\ref{eqn1}), and
canceling the two like terms appearing with opposite sign, we have
\begin{eqnarray*}
\frac{\partial}{\partial t}V  & = & \left(-\nabla^{l}{\rm
Ric}_{kj}+\nabla_{j}{\rm Ric}_k^l\right)X^*_{l}dx^k\otimes
dx^j +\triangle(\nabla X^*) \\
& & -g^{rs}{{(\nabla_rR)_{ks}}^l}_jX^*_ldx^k\otimes
dx^j-2{{R_k}^{rl}}_jV_{rl}dx^k\otimes dx^j
\\
& & -\left({\rm Ric}_k^lV_{lj}+{\rm Ric}_j^lV_{kl}\right)dx^k\otimes
dx^j.
\end{eqnarray*}

The last thing we need to see in order to complete the proof is that
$$-g^{rs}{{(\nabla_rR)_{ks}}^l}_j-\nabla^{l}{\rm Ric}_{kj}+\nabla_{j}{\rm Ric}_k^l=0.$$
This is obtained by  contracting $g^{rs}$ against the Bianchi
identity
$$\nabla_r{{R_{ks}}^l}_j+\nabla^lR_{ksjr}+\nabla_j{R_{ksr}}^l=0.$$
\end{proof}

Let $h_{ij}$ be defined by $h_{ij}= V_{ij}+V_{ji}$. It follows from
(\ref{eq.4}) that
\begin{equation}
\frac{\partial}{\partial t}h_{ij}=\Delta_{L}h_{ij}, \label{eq.5}
\end{equation}
where by definition $\Delta_{L}h_{ij}=\Delta
h_{ij}+2{{R_i}^{kl}j}h_{kl}-{\rm Ric}_i^k h_{kj}-{\rm
Ric}_j^kh_{ki}$ is the \textbf{Lichnerowicz
Laplacian}\index{Lichnerowicz Laplacian}. A simple calculation shows
that there is a constant $C>0$ such that
\begin{align}
\left(  \frac{\partial}{\partial t}-\Delta\right)  \left\vert h_{ij}
\right\vert ^{2}  &  =-2\left\vert \nabla_{k}h_{ij}\right\vert ^{2}
+4R^{ijkl}h_{jk}h_{il} \label{eq.6}\\
\frac{\partial}{\partial t}\left\vert h_{ij}\right\vert ^{2}  & \leq
\Delta\left\vert h_{ij}\right\vert ^{2}-2\left\vert \nabla_{k}h_{ij}
\right\vert ^{2}+C\left\vert h_{ij}\right\vert ^{2}. \label{eq.7}
\end{align}
Note that $X(t)$ is a Killing vector field for $g(t)$ if and only if
$h_{ij}(t)=0$. Since Equation~(\ref{eq1}) is linear and since the
curvature is bounded on each time-slice, for any given bounded
Killing vector field $X(0)$ for metric $g(0)$, there is a bounded
solution $X^{i}\left( t\right) $ of Equation~(\ref{eq1}) for
$t\in\lbrack0,T]$. Then $\left \vert h_{ij}\left( t\right) \right
\vert^2 $ is a bounded function satisfying (\ref{eq.7}) and $\left
\vert h_{ij} \right \vert^2 (0)=0$. One can apply the maximum
principle to (\ref{eq.7}) to conclude that $h_{ij}(t)=0$ for all
$t\geq 0$. This is done as follows: Let $h(t)$ denote the maximum of
$|h_{ij}(x,t)|^2$ on the $t$ time-slice. Note that, for any fixed
$t$ the function $|h_{ij}(x,t)|^2$ approaches $0$ as $x$ tends to
infinity  since the metric is asymptotic at infinity to the product
of a round metric on $S^2$ and the standard metric on the line. By
virtue of (\ref{eq.7}) and Proposition~\ref{fordiffmax}, the
function $h(t)$ satisfies $dh/dt\le Ch$ in the sense of forward
difference quotients, so that $d(e^{-Ct}h)/dt\le 0$, also in the
sense of forward difference quotients. Thus, by
Corollary~\ref{fordiffquot}, since $h(0)=0$ and $h\ge 0$, it follows
that $e^{-Ct}h(t)=0$ for all $t\ge 0$, and consequently, $h(t)=0$
for all $t\ge 0$.

Thus, the evolving vector field $X( t) $ is a Killing vector field
for $g(t)$ for all $t\in [0,T)$. The following is a very nice
observation of Bennett Chow; we thank him for allowing us to use it
here. From $h_{ij}=0$ we have $\nabla_{j} X^i +\nabla_{i} X^{j}=0$.
Taking the $\nabla_{j}$ derivative and summing over $j$ we get
$\Delta X^i +R^i_kX^k =0$ for all $t$. Hence (\ref{eq.3}) gives
$\frac{\partial}{\partial t}X^{i}=0$ and $X(t)=X(0)$, i.e., the
Killing vector fields are stationary and remain Killing vector
fields for the entire flow $g(t)$.
 Since at $t=0$ the Lie algebra $so(3)$ of the standard rotation action
 consists of Killing vector fields, the same is true for all the metrics $g(t)$
 in the standard solution.
 Thus, the
rotation group $SO(3)$ of $\mathbb{R}^3$ is  contained in the
isometry group of $g(t)$ for every $t\in [0,T)$. We have shown:

\begin{cor} The standard solution $g(t),\ t\in [0,T)$, consists of a family of metrics
all of 
which are
rotationally symmetric by the standard action of $SO(3)$ on
$\mathbb{R}^3$.
\end{cor}

\subsection{Non-collapsing}

\begin{prop}\label{stdkappa}
For any $r>0$ sufficiently small, there is a $\kappa>0$ such that
the standard flow is $\kappa$-non-collapsed on all scales $\le r$.
\end{prop}

\begin{proof}
Since the curvature of the standard solution is non-negative, it
follows directly that $2|{\rm Ric}|^2\le R^2$. By
Equation~(\ref{Revol}) this gives
$$\frac{\partial R}{\partial t}=\triangle R+2|{\rm Ric}|^2\le \triangle
R+R^2.$$ Let $C={\rm max}(2,{\rm max}_{x\in \Ar^3}R(x,0))$. Suppose
that $t_0< T$ and $t_0<1/C$.

\begin{claim}
For all $x\in \Ar^3$ and $t\in [0,t_0]$ we have
$$R(x,t)\le \frac{C}{1-Ct}.$$
\end{claim}

\begin{proof}
By the asymptotic condition, there is a compact subset $X\subset
\Ar^3$ such that for any point $p\in \Ar^3\setminus X$ and for any
$t\le t_0$ we have $R(p,t)<2/(1-t)$. Since $C\ge 2$, for all $t$ for
which ${\rm sup}_{x\in \Ar^3}R(x,t)\le 2/(1-t)$,  we also have
$$R(x,t)\le \frac{C}{1-Ct}.$$

Consider the complementary subset of $t$, that is to say the subset
of $[0,t_0]$ for which there is $x\in \Ar ^3$ with
$R(x,t)>C/(1-Ct)$. This is an open subset of $[0,t_0]$, and hence is
a disjoint union of relatively open intervals. Let $\{t_1<t_2\}$ be
the endpoints of one such interval. If $t_1\not=0$, then clearly
$R_{\rm max}(t_1)=C/(1-Ct_1)$. Since $C\ge {\rm sup}_{x\in
\Ar^3}R(x,0)$, this is also true if $t_1=0$. For every  $t\in
[t_1,t_2]$ the maximum of $R$ on the $t$ time-slice is achieved, and
the subset of $\Ar^3\times [t_1,t_2]$ of all points where maxima are
achieved is compact. Furthermore, at any maximum point we have
$\partial R/\partial t\le R^2$. Hence, according to
Proposition~\ref{fordiffmax} for all $t\in [t_1,t_2]$ we have
$$R_{\rm max}(t)\le G(t)$$
where $G'(t)=G^2(t)$ and $G(t_1)=C/(1-Ct_1)$. It is easy to see that
$$G(t)= \frac{C}{1-Ct}.$$
This shows that for all $t\in [t_1,t_2]$ we have $R(x,t)\le
\frac{C}{1-Ct}$, completing the proof of the claim.
\end{proof}

This shows that for $t_0<T$ and $t_0<1/C$ the scalar curvature is
bounded on $M\times [0,t_0]$ by a constant depending only on $C$ and
$t_0$. Since we are assuming that our flow is maximal, it follows
that $T\ge 1/C$.

Since  $(\Ar^3,g_0)$ is asymptotic to $(S^2\times \Ar,h(0)\times
ds^2)$, by compactness there is $V>0$ such that for any metric ball
$B(x,0,r)$ on which $|{\rm Rm}|\le r^{-2}$ we have ${\rm
Vol}\,B(x,r)\ge V r^3$. Since there is a uniform bound on the
curvature on $[0,1/2C]$, it follows that there is $V'>0$ so that any
ball $B(q,t,r)$ with $t\le 1/2C$ on which $|{\rm Rm}|\le r^{-2}$
satisfies ${\rm Vol}\,B(q,t,r)\ge V'r^3$. Set $t_0=1/4C$. For any
point $x=(p,t)$ with $t\ge 1/2C$ there is a point $(q,t_0)$ such
that $l_x(q,t_0)\le 3/2$; this by Theorem~\ref{n/2}. Since
$B(q,0,1/\sqrt{R_{\rm max}(0)})\subset \Ar^3$ has volume at least
$V/R_{\rm max}(0)^{3/2}$, and clearly $l_x$ is bounded above on
$B(q,0,1/\sqrt{R_{\rm max}(0)})$ by a uniform constant, we see that
the reduced volume of $B(q,0,1/\sqrt{R_{\rm max}(0)})$ is uniformly
bounded from below. It now follows from Theorem~\ref{THM} that there
is $\kappa_0>0$ such that if $|{\rm Rm}|$ is bounded by $r^{-2}$ on
the parabolic neighborhood $P(p,t,r,-r^2)$ and $r\le \sqrt{1/4C}$,
then the volume of this neighborhood is at least $\kappa_0r^3$.
Putting all this together we see that there is a universal
$\kappa>0$ such that the standard solution is $\kappa$-non-collapsed
on all scales at most $\sqrt{1/4C}$.
\end{proof}

\section{Uniqueness}\label{12.4}

Now we turn to the proof of uniqueness. The idea is to mimic the
proof of uniqueness in the compact case, by replacing the Ricci flow
by a strictly parabolic flow. The material we present here is
closely related to and derived from the presentation given in
\cite{LuTian}. The presentation here is the analogy in the context
of the standard solution of DeTurck's argument presented in
Section~\ref{3.3}.

\subsection{From Ricci flow to Ricci-DeTurck flow}

In this subsection we discuss the Ricci-DeTurck flow and the
harmonic map flow. Let $(M^n,g(t)),t \in [t_0,T]$ be a solution of
the Ricci flow and let $\psi_{t}\colon M \rightarrow M,\ t \in
\left[ t_0,T_{1}\right]$ be a solution of the harmonic map
flow\index{harmonic map flow|ii}
\begin{align}
\frac{\partial\psi_{t}}{\partial t}  &  =\Delta_{g\left(  t\right)
,g(t_0) }\psi_{t},\label{eq2.1} \quad \psi_{t_0}  ={\rm Id}.
\end{align}

Here, $\Delta_{g\left(  t\right) ,g(t_0) }$ is the Laplacian for
maps from the Riemannian manifold $(M,g(t))$ to the Riemannian
manifold $(M,g(t_0))$. In local coordinates $(x^i)$ on the domain
$M$ and $(y^\alpha)$ on the target $M$, the harmonic map flow
(\ref{eq2.1}) can be written as
\begin{equation}
\left( \frac{\partial}{\partial t}-\Delta_{g\left(  t\right)
}\right) \psi^{\alpha}\left(  x,t\right)  =g^{ij}\left( x,t\right)
\Gamma _{\beta\gamma}^{\alpha}\left(  \psi\left(  x,t\right) \right)
\frac {\partial\psi^{\beta}\left(  x,t\right)  }{\partial
x^{i}}\frac{\partial \psi^{\gamma}\left(  x,t\right)  }{\partial
x^{j}} \label{eq2.2}
\end{equation}
where $\Gamma_{\beta\gamma}^{\alpha}$ are the Christoffel symbols of
$g(t_0)$. Suppose $\psi\left(x, t\right)$ is a bounded smooth
solution of~\ref{eq2.2} with $\psi_{t_0}={\rm Id}$. Then $\psi\left(
t\right),\ t \in[t_0, T_1]$ are diffeomorphisms when $T_{1}>t_0$ is
sufficiently close to $t_0$. For any such $T_1$ and for $t_0 \leq t
\leq T_1$, define $\hat{g}\left( t\right) =\left(
\psi_{t}^{-1}\right) ^{\ast}g\left( t\right)$. Then $\hat g(t)$
satisfies the following equation:
\begin{equation}\label{RDfloweqn}
\frac{\partial}{\partial t}\hat{g}_{ij} =-2\widehat{{\rm
Ric}}_{ij}+\hat{\nabla} _{i}W_{j}(t)+\hat{\nabla}_{j}W_{i}(t) \quad
\quad \hat{h}\left( 0\right) =h(0),
\end{equation}
where $\widehat{{\rm Ric}}_{ij}$ and $\hat{\nabla}_{i}$ are the
Ricci curvature and Levi-Civita connection of $\hat{g}(t)$
respectively and $W(t)$ is the time-dependent $1$-form defined by
\[
W(t) _{j}=\hat{g}_{jk}(t)\hat{g}\,^{pq}(t)\left(
\hat{\Gamma}_{pq}^{k}(t)-\Gamma_{pq}^{k}(t_0)   \right).
\]
Here, $\hat{\Gamma}_{pq}^{k}(t)$ denotes the Christoffel symbols of
the metric $\hat g(t)$ and $\Gamma_{pq}^k(t_0)$ denotes the
Christoffel symbols of the metric $g(t_0)$. (See, for example,
(\cite{Shi1} Lemma 2.1).)
 We call a solution to this flow equation a
 \textbf{Ricci-DeTurck flow}\index{Ricci-DeTurck flow|ii} (see \cite{DeTurck}, or \cite{ChowKnopf}
Chapter 3 for details). In local coordinates we have

\begin{align}
& \frac{\partial\hat{g}_{ij}}{\partial t}  =\hat{g}^{kl}\nabla
_{k}\nabla_{l}\hat{g}_{ij}-\hat{g}^{kl}g(t_0)_{ ip}\hat{g}^{pq}R_{jk
ql}\left(  g(t_0)\right)  -\hat{g} ^{kl}g(t_0)_{
jp}\hat{g}^{pq}R_{ik ql}\left(  g(t_0)\right) \nonumber
\\
&  +\frac{1}{2}\hat{g}^{kl}\hat{g}^{pq}\left[
\begin{array}
[c]{c} \nabla_{i}\hat{g}_{pk}\nabla_{ j}\hat{g}_{ql}+2\nabla
_{k}\hat{g}_{jp}\nabla_{ q}\hat{g}_{il}\\
-2\nabla_{k}\hat{g}_{jp}\nabla_{l}\hat{g}_{iq}-2\nabla_{
j}\hat{g}_{pk}\nabla_{l}\hat{g}_{iq}-2\nabla_{ i}\hat
{g}_{pk}\nabla_{l}\hat{g}_{jq}
\end{array}
\right]  . \label{eqRcDeTurck}
\end{align}
where $\nabla$ is the Levi-Civita connection of $g(t_0)$. This is a
strictly parabolic equation.

\begin{lem}\label{12.7}
Suppose that $g(t)$ solves the Ricci flow equation and suppose that
$\psi_t$ solves the harmonic map flow equation,
Equation~(\ref{eq2.1}); then $\hat g(t)=(\psi_t^{-1})^*g(t)$ solves
the Ricci-DeTurck flow, Equation~(\ref{RDfloweqn}) and $\psi_t$
satisfies the following ODE:
$$\frac{\partial \psi_t}{\partial t}=-\hat g^{ij}(t)W(t).$$
\end{lem}

\begin{proof}
The first statement follows from the second statement and a standard
Lie derivative computation. For the second statement, we need to
show
$$\triangle_{g(t),g(0)} \psi^\alpha = -\hat
g^{pq}\left(\hat\Gamma_{pq}^\alpha(t)-\Gamma_{pq}^\alpha(t_0)\right).$$
Notice that this equation is a tensor equation, so that we can
choose coordinates in the domain and range so that $\Gamma(t)$
vanishes at the point $p$ in question and $\Gamma(t_0)$ vanishes at
$\psi_t(p)$.  With these assumptions we need to show
$$g^{pq}(t)\frac{\partial^2 \psi^\alpha}{\partial x^p\partial x^q}=-\hat g^{pq}(t)\hat\Gamma_{pq}^\alpha(t).$$
This is a direct computation using the change of variables formula
relating $\hat \Gamma$ and $\Gamma$.
\end{proof}

\begin{cor}\label{dtuniqtouniq}
Suppose that $(M,g_1(t)),\ t_0\le t\le T$, and $(M,g_2(t)),\ t_0\le
t\le T$, are solutions to the Ricci flow equation for which there
are solutions
$$\psi_{1,t}\colon (M,g_1(t))\to (M,g_1(0))$$ and
$$\psi_{2,t}\colon
(M,g_2(t))\to (M,g_2(0))$$ to the harmonic map equation with
$\psi_{1,t_0}=\psi_{2,t_0}={\rm Id}$. Let $\hat
g_1(t)=(\psi_{1,t}^{-1})^*g_1(t)$ and $\hat
g_2(t)=(\psi_{2,t}^{-1})^*g_2(t)$ be the corresponding solutions to
the Ricci-DeTurck flow. Suppose that $\hat g_1(t)=\hat g_2(t)$ for
all $t\in [t_0,T]$. Then $g_1(t)=g_2(t)$ for all $t\in [t_0,T]$.
\end{cor}

\begin{proof}
Since $\psi_{a,t}$ satisfies the equation
$$\frac{\partial \psi_{a,t}}{\partial t}=-\hat g_a^{ij}W(t)_j$$
where the time-dependent vector field $W(t)$ depends only on $\hat
g_a$, we see that $\psi_{1,t}$ and $\psi_{2,t}$ both solve the same
time-dependent ODE and since $\psi_{1,t_0}=\psi_{2,t_0}={\rm Id}$,
it follows that $\psi_{1,t}=\psi_{2,t}$ for all $t\in [t_0,T]$. On
the other hand, $g_a(t)=\psi_{a,t}^*\hat g_a(t)$, so that it follows
that $g_1(t)=g_2(t)$ for $t\in [t_0,T]$.
\end{proof}

Our strategy of proof is to begin with a standard solution $g(t)$
and show that there is a solution to the harmonic map equation for
this Ricci flow with appropriate decay conditions at infinity. It
follows that the solution to the Ricci-DeTurck flow constructed is
well-controlled at infinity. Suppose that we have two standard
solutions $g_1(t)$ and $g_2(t)$ (with the same initial conditions
$g_0$) that agree on the interval $[0,t_0]$ which is a proper
subinterval of the intersection of the intervals of definition of
$g_1(t)$ and $g_2(t)$. We construct solutions to the harmonic map
flow equation from $g_a(t)$ to $g_a(t_0)$ for $a=1,2$. We show that
solutions always exist for some amount of time past $t_0$. The
corresonding Ricci-DeTurck flows $\hat g_a(t)$ starting at $g_{t_0}$
are well-controlled at infinity. Since the Ricci-DeTurck flow
equation is a purely parabolic equation, it has a unique solution
with appropriate control at infinity and given initial condition
$g_1(t_0)=g_2(t_0)$. This implies that the two Ricci-DeTurck flows
we have constructed are in fact equal. Invoking the above corollary,
we conclude that $g_1(t)$ and $g_2(t)$ agree on a longer interval
extending  past $t_0$. From this it follows easily that $g_1(t)$ and
$g_2(t)$ agree on their common domain of definition. Hence, if they
are both maximal flows, they must be equal.

\section{Solution of the harmonic map flow}

In order to pass from a solution to the Ricci flow equation to a
solution of the Ricci-DeTurck flow we must prove the existence of a
solution of the harmonic map flow\index{harmonic map flow}
associated with the Ricci flow. In this section we study the
existence of the harmonic flow (\ref{eq2.1}) and its asymptotic
behavior at the space infinity when $h(t)=g(t)$ is a standard
solution. Here we use in an essential way the rotationally symmetric
property and asymptotic property at infinity of $g(t)$. In this
argument there is no reason, and no advantage, to restricting to
dimension three, so we shall consider rotationally symmetric
complete metrics on $\Ar^n$, i.e., complete metrics on $\Ar^n$
invariant under the standard action of $SO(n)$. Let
$\theta=(\theta^1,\cdots,\theta^{n-1})$ be local coordinates on the
round $(n-1)$-sphere of radius $1$, and let $d\sigma$ be the metric
on the sphere. We denote by $\hat r$ the standard radial coordinate
in $\Ar^n$. Since $g(t)$ is rotationally symmetric and $n \geq 3$,
we can write

\begin{align}
g(t)=dr^2+f(r,t)^2 d \sigma \label{eq2.4}
\end{align}
Here $r=r(\hat r,t)$ is the (time-dependent) radial coordinate on
$\mathbb{R}^n$ for the metric $g(t)$.

\begin{claim}\label{cl12.20}
For any fixed $t$ the function $r\colon \Ar^n \to[0,\infty)$ is a
function only of $\hat r$. Considered as a function of two
variables,
 $r(\hat r,t)$ is a smooth function defined for $\hat r\ge 0$. It is
 an odd function of $\hat r$. For fixed $t$ it is an increasing
 function of $\hat r$.
\end{claim}

\begin{proof}
Write the metric $g(t)=g_{ij}dx^idx^j$ and let
\begin{equation}
x^1  =  \hat{r} \cos \theta^1, \ \ x^2 =  \hat{r} \sin \theta^1 \cos
\theta^2, \cdots,
 x^n  =  \sin \theta^1 \cdots \sin \theta^{n-1}.
 \end{equation}
 We compute $f(r,t)$
by restricting attention to the ray $\hat{r}=x^1$ and $\theta^1
=\cdots=\theta^{n-1}=0$, i.e., $x^2 =\cdots =x^n=0$. Then
\[
g(t)=g_{11}(\hat{r},0,\cdots,0,t)d\hat{r}^2+g_{22}(\hat{r},0,\cdots,0,t)
\hat{r}^2 d \sigma.
\]
Both $g_{11}$ and $g_{22}$ are positive smooth and even in $\hat r$.
Clearly,$\sqrt{g_{11}(\hat r,0,\cdots,0,t)}$ is a positive smooth
function defined for all $(\hat r,t)$ and is invariant under the
involution $\hat r\mapsto -\hat r$. Hence its restriction to $\hat
r\ge0$ is an even function. Since
\[
r= \int_0^{\hat{r}} \sqrt{g_{11}(\hat{s},0,\cdots,0,t)}d\hat{s} =
\hat{r} \int_0^{1} \sqrt{g_{11}(\hat{r}s,0,\cdots,0,t)}ds,
\]
we see that $r(\hat r,t)$ is of the form $\hat r\cdot \phi(\hat
r,t)$ where $\phi(\hat r,t)$ is an even smooth function. This shows
that $r(\hat r,t)$ is an odd function. It is also clear from this
formula that $\partial r/\partial \hat r>0$.
\end{proof}

 Since, for each $t_0$, the function $r(\hat
r,t_0)$ is an increasing function of $\hat r$,  it can be inverted
to give a function $\hat r(r,t_0)$. In Equation~(\ref{eq2.4}), we
have chosen to write $f$ as a function of $r$ and $t$, rather than a
function of $\hat r$ and $t$. We look for rotationally symmetric
solutions to Equation~(\ref{eq2.1}), i.e., solutions of the form:

\begin{align}\label{ansatz}
\psi(t)\colon \mathbb{R}^n \rightarrow  \mathbb{R}^n \qquad
\psi(t)(r,\theta) & =(\rho(r,t),\theta) \ \ \ {\rm for\ \ } t\ge t_0 \\
\psi(r,t_0) & = {\rm Id} \nonumber
\end{align}

We shall adopt the following conventions: we shall consider
functions $f(w,t)$ defined in the closed half-plane $w\ge 0$. When
we say that such a function is {\em smooth} we mean that for each
$n,m\ge 0$ we have a continuous function $f_{nm}(w,t)$ defined for
all $w\ge 0$ with the following properties:
\begin{enumerate}
\item[(1)] $$f_{00}=f$$
\item[(2)] $$\frac{\partial f_{nm}}{\partial t}=f_{n(m+1)}$$
\item[(3)] $$\frac{\partial f_{nm}}{\partial w}=f_{(n+1)m},$$
\end{enumerate}
where in Item (3) the partial derivative along the boundary $w=0$ is
a right-handed derivative only. We say such a function is {\em even}
if $f_{(2k+1)m}(0,t)=0$ for all $k\ge 0$.

We have the following elementary lemma:

\begin{lem} (a) Suppose that $f(w,t)$ is a smooth function defined for $w\ge 0$.
Define $\phi(r,t)=f(r^2,t)$. Then $\phi(r,t)$ is a smooth function
defined for all $r\in \Ar$. Now fix $k$ and let $\hat r\colon
\Ar^k\to [0,\infty)$ be the usual radial coordinate.  Then we have a
smooth family of smooth functions on $\Ar^k$ defined by
$$\hat\phi(x^1,\ldots,x^k,t)=\phi(\hat r(x^1,\ldots,x^k),t)=f(\sum_{i=1}^k(x^i)^2,t).$$

\noindent (b) If $\psi(r,t)$ is a smooth function defined for $r\ge
0$ and if it is even in the sense that its Taylor expansion to all
orders along the line $r=0$ involves only even powers of $r$, then
there is a smooth function $f(w,t)$ defined for $w\ge 0$ such that
$\psi(r,t)=f(r^2,t)$. In particular, for any $k\ge 2$ the function
$\hat\psi((x^1,\ldots,x^k),t)=\psi(r(x^1,\ldots,x^k),t)$ is a smooth
family of smooth functions on $\Ar^k$.
\end{lem}

\begin{proof}
Item (a) is obvious, and Item (b) is obvious away from $r=0$. We
establish Item (b) along the line $r=0$. Consider the Taylor theorem
with remainder to order $2N$ in the $r$-direction for $\psi(r,t)$ at
a point $(0,t)$. By hypothesis it takes the form
$$\sum_{i=0}c_i(t)w^{2i}+w^{2N+1}R(w,t).$$
Now we replace $w$ by $\sqrt{r}$ to obtain
$$f(r,t)=\sum_{i=0}c_ir^i+\sqrt{r}^{2N+1}R(\sqrt{r},t).$$
Applying the usual chain rule and taking limits as $r\rightarrow
0^+$ we see that $f(r,t)$ is $N$ times differentiable along the line
$r=0$. Since this is true for every $N<\infty$, the result follows.
\end{proof}

Notice that an even function $f(r,t)$ defined for $r\ge 0$ extends
to a smooth function on the entire plane invariant under $r\mapsto
-r$. When we say a function $f(r,t)$ {\em defines a smooth family of
smooth functions on $\Ar^n$} we mean that, under the substitution
$\hat f((x^1,\ldots,x^n),t)=f(r(x^1,\ldots,x^n),t)$, the function
$\hat f$ is a smooth function on $\Ar^n$ for each $t$.

We shall also consider odd functions $f(r,t)$, i.e., smooth
functions defined for $r\ge 0$ whose Taylor expansion in the
$r$-direction along the line $r=0$ involves only odd powers of $r$.
These do not define smooth functions on $\Ar^n$. On the other hand,
by the same argument as above with the Taylor expansion one sees
that they can be written as $rg(r,t)$ where $g$ is even, and hence
define  smoothly varying families of smooth functions on $\Ar^n$.
Notice also that the product of two odd functions $f_1(r,t)f_2(r,t)$
is an even function and hence this product defines a smoothly
varying family of smooth function on $\Ar^n$.

\subsection{The properties of $r$ as a function of $\hat r$ and $t$}

We shall make a change of variables and write the harmonic map flow
equation in terms of $r$ and $\theta$. For this we need some basic
properties of $r$ as a function of $\hat r$ and $t$. Recall that we
are working on $\Ar^n$ with its usual Euclidean coordinates
$(x^1,\ldots,x^n)$. We shall also employ spherical coordinates $\hat
r,\theta^1,\ldots,\theta^{n-1}$. (We denote the fixed radial
coordinate on $\Ar^n$ by $\hat r$ to distinguish it from the varying
radial function $r=r(t)$ that measures the distance from the tip in
the metric $g(t)$.)

As a corollary of Claim~\ref{cl12.20} we have:

\begin{cor}\label{rcor}
$r^2(\hat r,t)$ is a smoothly varying family of smooth functions on
$\Ar^n$. Also, $\hat r$ is a smooth  function of $(r,t)$ defined for
$r\ge 0$ and odd in $r$. In particular, any smooth even function of
$r$ is a smooth even function of $\hat r$ and thus defines a smooth
function on $\Ar^n$. Moreover, there is a smooth function $\xi(w,t)$
such that $d({\rm log}\, r)/dt=r^{-1}(dr/dt)=\xi(r^2,t)$.
\end{cor}

For future reference we define
\begin{equation}\label{Bdefn}
B(w,t)=\frac{1}{2}\int_0^w\xi(u,t)du. \end{equation} Then $B(r^2,t)$
is a smooth function even in $r$ and hence, as $t$ varies, defines a
smoothly varying family of smooth functions on $\Ar^n$. Notice that
$$\frac{\partial B(r^2,t)}{\partial r}=2r\frac{\partial B}{\partial
w}(w,t)\left|_{w=r^2}\right.=2r\left(\frac{1}{2}\xi(r^2,t)\right)=\frac{dr}{dt}.$$

Now let us consider $f(r,t)$.

\begin{claim}
$f(r,t)$ is a smooth function defined for $r\ge 0$. It is an odd
function of $r$.
\end{claim}

\begin{proof}
We have
$$f(r,t)=\hat r(r,t)\sqrt{g_{22}(\hat r(r,t),0,\ldots,0,t)}.$$
Since $\sqrt{g_{22}(\hat r,0,\ldots,0,t)}$ is a smooth function of
$(\hat r,t)$ defined for $\hat r\ge 0$ and since it is an  even
function of $\hat r$, it follows immediately from the fact that
$\hat r$ is a smooth odd function of $r$, that $f(r,t)$ is a smooth
odd function of $r$.
\end{proof}

\begin{cor}
There is a smooth function $h(w,t)$ defined for $w\ge 0$ so that
$f(r,t)=rh(r^2,t)$. In particular, $h(r^2,t)$ defines a smooth
function on all of $\Ar^n$. Clearly, $h(w,t)>0$ for all $w\ge 0$ and
all $t$.
\end{cor}

We set $\tilde h(w,t)={\rm log}(h(w,t))$, so that $f(r,t)=re^{\tilde
h(r^2,t)}$. Notice that $\tilde h(r^2,t)$ defines a smooth function
of $\hat r^2$ and $t$ and hence is a smoothly varying family of
smooth functions on $\Ar^n$.

\subsection{The harmonic map flow equation}.

Let $\psi(t)\colon \Ar^n\to \Ar^n$ be a smoothly varying family of
smooth functions as given in Equation~(\ref{ansatz}). Using
(\ref{eq2.4}) and (\ref{ansatz}) it is easy to calculate the energy
functional using spherical coordinates with $r$ as the radial
coordinate.

\begin{align*}
E(\psi(t)) & = \frac{1}{2} \int_{\mathbb{R}^n} \left \vert \nabla
\psi(t)
\right \vert^2_{g(t),g_{t_0}} dV_{g(t)}  \\
&= \frac{1}{2} \int_{\mathbb{R}^n} \left[ \left( \frac{\partial
\rho}{\partial r} \right)^2 + (n-1) f^2(\rho,t_0) f^{-2}(r,t)
\right] dV_{g(t)}.
\end{align*}

If we have a compactly supported variation $\delta \rho =w$, then
letting $d{\rm vol}_\sigma$ denote the standard volume element on
$S^{n-1}$, we have
\begin{align*}
& \delta E(\psi(t))(w) = \frac{1}{2} \int_{\mathbb{R}^n} \left[2
\frac{\partial \rho}{\partial r} \frac{\partial w}{\partial r}
+2(n-1) f(\rho,t_0)
\frac{\partial f(\rho,t_0)}{\partial \rho}f^{-2}(r,t)w \right] dV_g(t) \\
&= \int_0^{+\infty} \left[  f^{n-1}(r,t)\frac{\partial
\rho}{\partial r} \frac{\partial w}{\partial r} +(n-1) f(\rho,t_0)
\frac{\partial f(\rho,t_0)}{\partial
\rho}f^{n-3}(r,t)w \right] dr \cdot \int_{S^{n-1}} d{\rm vol}_\sigma \\
& = \int_0^{+\infty} \left[-\frac{\partial}{\partial
r}\left(f^{n-1}\frac{\partial\rho}{\partial r}\right)w+(n-1)
f(\rho,t_0) \frac{\partial f(\rho,t_0)}{\partial \rho}f^{n-3}(r,t)w
\right] dr \cdot \int_{S^{n-1}} d{\rm vol}_\sigma\\ &=
\int_{\mathbb{R}^n}\left[ - f^{1-n} \frac{\partial }{\partial r}
\left( \frac{\partial \rho}{\partial r} f^{n-1}\right) +(n-1)
f(\rho,t_0) \frac{\partial f(\rho,t_0)}{\partial \rho}f^{-2}(r,t)
\right] w dV_{g(t)}.
\end{align*}

The usual argument shows that for a compactly supported variation
$w$ we have
$$\delta_w\left(\frac{1}{2}\int_{\Ar^n}|\nabla_{g(t),g(t_0)}\psi|^2d{\rm vol}=\int_{\Ar^n}\langle
w,-\triangle_{g(t),g(t_0)}\psi\rangle d{\rm vol}\right).$$ Thus,
$$\triangle_{g(t),g(t_0)}\psi= \left[
f^{1-n}\frac{\partial }{\partial r} \left( \frac{\partial
\rho}{\partial r} f^{n-1}\right) -(n-1) f(\rho,t_0) \frac{\partial
f(\rho,t_0)}{\partial \rho}f^{-2}(r,t) \right]
\frac{\partial}{\partial r}$$ where we have written this expression
using the coordinates $(r,\theta)$ on the range $\Ar ^n$ (rather
than the fixed coordinates $(\hat r,\theta)$).

Now let us compute $\partial \psi/\partial t(\hat r,t)$ in these
same coordinates. (We use $\hat r$ for the coordnates for $\psi$ in
the domain to emphasize that this must be the time derivative at a
fixed point in the underlying space.) Of course, by the chain  rule,
\begin{eqnarray*}
\frac{\partial \psi(\hat r,t)}{\partial t} & = & \frac{\partial
\psi(r,t)}{\partial r}\frac{\partial r}{\partial t}+\frac{\partial
\psi(r,t)}{\partial t} \\
& = & \frac{\partial \rho(r,t)}{\partial r}\frac{\partial r(\hat
r,t)}{\partial t}+\frac{\partial \rho(r,t)}{\partial
t}.\end{eqnarray*}

Consequently,   for rotationally symmetric maps as in
Equation~(\ref{ansatz}) the harmonic map flow equation (\ref{eq2.1})
has the following form:

\begin{equation*}
 \frac{\partial \rho}{\partial t}+\frac{\partial \rho}{\partial r}\frac{\partial r}{\partial t}
 = \frac{1}{f^{n-1}(r,t)}
\frac{\partial }{\partial r} \left(  f^{n-1}(r,t) \frac{\partial
\rho}{\partial r} \right)-(n-1)f^{-2}(r,t) f(\rho,t_0)
\frac{\partial f(\rho,t_0)}{\partial \rho}
\end{equation*}
or equivalently

\begin{equation}\label{eqrothar}
\frac{\partial \rho}{\partial t}=\frac{1}{f^{n-1}(r,t)}
\frac{\partial }{\partial r} \left(  f^{n-1}(r,t) \frac{\partial
\rho}{\partial r} \right)-(n-1)f^{-2}(r,t) f(\rho,t_0)
\frac{\partial f(\rho,t_0)}{\partial \rho} -\frac{\partial
\rho}{\partial r}\frac{\partial r}{\partial t}
\end{equation}

The point of  rewriting the harmonic map equation in this way is to
find an equation for the functions $\rho(r,t), f(r,t)$ defined on
$r\ge 0$. Even though the terms in this rewritten equation involve
odd functions of $r$, as we shall see, solutions to these equations
will be even in $r$ and hence will produce a smooth solution to the
harmonic map flow equation on $\Ar ^n$.

\subsection{An equation equivalent to the harmonic map flow
equation}

We will solve (\ref{eqrothar})  for solutions of the form

\[ \rho(r,t)= re^{\tilde{\rho}(r,t)}, \ \ \ t\ge t_0;\ \ \ \tilde \rho(r,t_0)=0.
\]

For $\psi$ as in Equation~(\ref{ansatz}) to define a diffeomorphism,
it must be the case that $\rho(r,t)$ is a smooth function for $r\ge
0$ which is odd in $r$. It follows from the above expression that
$\tilde\rho(r,t)$ is a smooth function of $r$ and $t$ defined for
$r\ge 0$ and even in $r$, so that it defines a smoothly varying
family of smooth functions on $\Ar^n$. Then some straightforward
calculation shows that (\ref{eqrothar}) becomes

\begin{eqnarray}\label{newform}
\frac{\partial\tilde{\rho} }{\partial t} & = &
\frac{\partial^{2}\tilde{\rho}}{\partial r^{2}}+\frac{n+1}{r}\frac
{\partial\tilde{\rho}}{\partial r}
 +\left(  n-1\right) \frac{\partial\tilde{h}(
r^{2},t)}{\partial r}  \frac{\partial\tilde{\rho} }{\partial r}
+\left( \frac{\partial\tilde{\rho}
}{\partial r}\right)  ^{2}    \\
& & +\frac{n-1}{r^{2}}\left[  1-e^{2\tilde{h}\left(
\rho^{2},t_0\right) -2\tilde{h}\left(  r^{2},t\right)  }\right]
 +2\left(  n-1\right)  \frac{\partial\tilde{h}}{\partial w}\left(  r^{2},t\right) \nonumber \\
&  & -2\left(  n-1\right)  e^{2\tilde{h}\left(  \rho ^{2},t_0\right)
+2\tilde{\rho}  -2\tilde{h}\left( r^{2},t\right)
}\frac{\partial\tilde{h}}{\partial w}(\rho^2,t_0)
-\frac{2}{r}\frac{\partial r}{\partial t}-\frac{\partial r}{\partial
t} \frac{\partial\tilde{\rho} }{\partial r} .\nonumber
\end{eqnarray}

Note that from the definition, $\tilde{h}\left( 0,t \right)=0$, we
can write $\tilde{h}\left( w,t \right)=w \tilde{h}^*(w,t)$ where
$\tilde{h}^*(w,t)$ is a smooth function of $w\ge 0$ and $t$. So
$$
\frac{n-1}{r^{2}}\left[  1-e^{2\tilde{h}\left( \rho^{2},t_0\right)
-2\tilde{h}\left(  r^{2},t\right)  }\right] =\frac{n-1}{r^{2}}\left[
1-e^{2r^2\left[ e^{2\tilde{\rho}} \tilde{h}^*\left(
\rho^{2},t_0\right)  -\tilde{h}^*\left( r^{2},t\right) \right]
}\right]
$$
which is a smooth function of $\tilde{\rho},r^2,t$.

Let
\begin{align}
G(\tilde{\rho},w,t) = & \frac{n-1}{w}\left[ 1-e^{2\tilde{h}\left(
\rho^{2},t_0\right)  -2\tilde{h}\left( w,t\right)  }\right]
 +2\left(  n-1\right)  \frac{\partial\tilde{h}}{\partial w}\left(  w,
 t\right) \nonumber \\
& -2\left(  n-1\right)  e^{2\tilde{h}\left(  \rho ^{2},t_0\right)
+2\tilde{\rho}  -2\tilde{h}\left( w,t\right) }\frac{\partial\tilde
h}{\partial w}(\rho^2,t_0) -2\xi(w,t), \label{defofG}
\end{align}
where $\xi$ is the function from Corollary~\ref{rcor}. Then
$G(\tilde{\rho},w,t)$ is a smooth function defined for $w\ge 0$.
Notice that when $r$ and $\tilde \rho$ are the functions associated
with the varying family of metrics $g(t)$ and the solutions to the
harmonic map flow, then $G(\tilde\rho,r^2,t)$ defines a smoothly
varying family of smooth functions on $\Ar^n$.

We have the following form of equation (\ref{newform}):

\begin{align*}
\frac{\partial\tilde{\rho} }{\partial t} = &
\frac{\partial^{2}\tilde{\rho}}{\partial r^{2}}+\frac{n+1}{r}\frac
{\partial\tilde{\rho}}{\partial r} +\left[ ( n-1)
\frac{\partial\tilde{h}}{\partial r} -\frac{\partial B}{\partial r}
\right] \left( r^{2},t\right) \frac{\partial\tilde{\rho} }{\partial
r} +\left( \frac{\partial\tilde{\rho} }{\partial r}\right)  ^{2}
+G(\tilde\rho,r^2,t).
\end{align*}

Now we think of $\tilde{\rho}$ as a rotationally symmetric function
defined on $\mathbb{R}^{n+2}$ and let $\widehat
G(\tilde\rho,(x^1,\ldots,x^{n+2}),t)=
G(\tilde\rho,\sum_{i=1}^{n+2}(x^i)^2,t)$ and then the above equation
can be written as

\begin{align}
\frac{\partial\tilde{\rho} }{\partial t}
 = \Delta \tilde{\rho} +  \nabla [ (n-1)\tilde{h} -B] \cdot \nabla\tilde{\rho}
+ \left | \nabla \tilde{\rho} \right |^2 +G(\tilde{\rho},x,t)
\label{betterharMapeq}
\end{align}
where $\nabla$ and $\Delta$ are the Levi-Civita connection and
Laplacian defined by the Euclidean metric on  $\mathbb{R}^{n+2}$
respectively and where $B$ is the function defined in
Equation~(\ref{Bdefn}).

\begin{rem}
The whole purpose of this rewriting of the PDE for $\tilde \rho$ is
to present this equation in such a form that all its coefficients
represent smooth functions of $\hat r$ and $t$ that are even in
$\hat r$ and hence define smooth functions on Euclidean space of any
dimension. We have chosen to work on $\Ar^{n+2}$ because the
expression for the Laplacian in this dimension has the term
$((n+1)/r)\partial \tilde \rho/\partial r$.
\end{rem}

It is important to understand the asymptotic behavior of our
functions at spatial infinity.

\begin{claim}\label{asympt}
For any fixed $t$ we have the following asymptotic expansions at
spatial infinity.
\begin{enumerate}
\item[(1)]
 $e^{\tilde{h}(r^2,t)}$ is asymptotic to $\frac{1}{(1-t)r}$.
\item[(2)] $\tilde{h}(r^2,t)$ is asymptotic to $-\log r$.
\item[(3)] $\frac{\partial \tilde{h}}{\partial w}(r^2,t)$ is asymptotic to
 $-\frac{1}{2r^2}$.
\item[(4)]  $r^{-1} \frac{\partial r}{\partial t}$ is asymptotic to
$\frac{C}{r}$.
\item[(5)] $\frac{\partial B(r^2,t)}{\partial r}$ is asymptotic to  $C$.
\item[(6)] $|G( \tilde{\rho},r^2,t)| \leq C_*< \infty$ where
$C_* =C_*\left(\sup\{|\tilde{\rho}|,\tilde{h}\}\right)$ is a
constant depending only on $\sup\{|\tilde{\rho}|,\tilde{h}\}$.
\end{enumerate}
\end{claim}

\begin{proof}
The first item is immediate from  Proposition~\ref{item5}. The
second and third follow immediately from the first. The fourth is a
consequence of the fact that by Proposition~\ref{item5} $dr/dt$ is
asymptotic to a constant at infinity on each time-slice. The fifth
follows immediately from the fourth and the definition of
$B(r^2,t)$. Given all these asymptotic expressions, the last is
clear from the expression for $G$ in terms of $\tilde \rho$, $r^2$,
and $t$.
\end{proof}

\subsection{The short time existence}

The purpose of this subsection is to prove the following short-time
existence theorem for the harmonic map flow equation.

\begin{prop}
For any $t_0\ge 0$ for which there is a standard solution $g(t)$
defined on $[0,T_1]$ with $t_0<T_1$ there is $T>t_0$ and a solution
to Equation~(\ref{betterharMapeq}) with initial condition
$\tilde{\rho}(r,t_0)=0$  defined on the time-interval $[t_0,T]$.
\end{prop}

 At this point to simplify the notation we shift time by $-t_0$ so
that our initial time is $0$, though our initial metric is not $g_0$
but rather is the time $t_0$-slice of the standard solution we are
considering, so that now $t_0=0$ and our initial condition is
$\tilde \rho(r,0)=0$.

Let $x=(x^1,\cdots,x^{n+2})$ and $y=(y^1,\cdots,y^{n+2})$ be two
points in $\mathbb{R}^{n+1}$ and

\[
H(x,y,t)=\frac{1}{(4 \pi t)^{(n+2)/2}} e^{-\frac{|x-y|^2}{4t}}
\]
be the heat kernel. We solve (\ref{betterharMapeq}) by successive approximation
\cite{LiTam}.

Define
$$F(x, \tilde{\rho},\nabla\tilde{\rho},t)
=  \nabla \left[ (n-1)\tilde{h}-B \right] \cdot \nabla\tilde{\rho} +
\left | \nabla \tilde{\rho} \right |^2 +G( \tilde{\rho},x,t)
$$

Let $\tilde{\rho}_0(x,t)=0$ and for $i \geq 1$ we define
$\tilde{\rho}_i$ by

\begin{align}
\tilde{\rho}_i=\int_0^t \int_{\mathbb{R}^{n+2}}
H(x,y,t-s)F(y,\tilde{\rho}_{i-1}, \nabla\tilde{\rho}_{i-1} ,t)dyds
\label{rhoIntg}
\end{align}
which solves

\begin{align}
\frac{\partial\tilde{\rho}_i}{\partial t} = \Delta \tilde{\rho}_i +
F(x, \tilde{\rho}_{i-1},\nabla\tilde{\rho}_{i-1},t) \qquad
\tilde{\rho}_i(x,0)=0. \label{rhodiffeq}
\end{align}

To show the existence of $\tilde{\rho}_i$ by induction, it suffices
to prove the following statement: For any $i \geq 1$, if
$|\tilde{\rho}_{i-1}| , |\nabla \tilde{\rho}_{i-1}| $ are bounded,
then $\tilde{\rho}_i$ exists and $|\tilde{\rho}_{i}|, |\nabla
\tilde{\rho}_{i}|$ are bounded. Assume  $|\tilde{\rho}_{i-1}| \leq
C_1, |\nabla \tilde{\rho}_{i-1}| \leq C_2$ are bounded on
$\mathbb{R}^{n+2} \times [0,T]$; then it follows from
Claim~\ref{asympt} that $G(\tilde{\rho}_{i-1},{\bf x},t)$ is bounded
on $\mathbb{R}^{n+2} \times [0,T]$

\[
|G(\tilde{\rho}_{i-1},x,t)| \leq C_*(C_1,\tilde{h} ),
\]
and also because of Claim~\ref{asympt} both $|\nabla B|$ and
$|\nabla \tilde h|$ are bounded on all of $\Ar^{n+2}\times [0,T]$,
it follows that $F(x,
\tilde{\rho}_{i-1},\nabla\tilde{\rho}_{i-1},t)$ is bounded:
\begin{align*}
& |F(x, \tilde{\rho}_{i-1},\nabla\tilde{\rho}_{i-1},t)|  \\
\leq & \left[(n-1)\sup |\nabla \tilde{h}| + \sup |\nabla B|
\right]C_2+C_2^2+ C_*(C_1,\tilde{h} ) = C_3
\end{align*}
Hence $\tilde{\rho}_i$ exists.

The bounds on $|\tilde{\rho}_{i}|$ and $|\nabla \tilde{\rho}_{i}|$
follow from the following estimates

\[
| \tilde{\rho}_i | \leq \int_0^t \int_{\mathbb{R}^{n+2}}
H(x,y,t-s)C_3 dyds \leq C_3t,
\]
and

\begin{align*}
& | \nabla \tilde{\rho}_i | = | \int_0^t \int_{\mathbb{R}^{n+2}}
[\nabla_x H(x,y,t-s)]
F(y,\tilde{\rho}_{i-1}, \nabla\tilde{\rho}_{i-1} ,t)dyds| \\
& \leq \int_0^t \int_{\mathbb{R}^{n+2}} | \nabla_x H(x,y,t-s)| C_3 dyds \\
& = \int_0^t \int_{\mathbb{R}^{n+2}} \frac{1}{(4 \pi
{(t-s)})^{(n+2)/2}}
e^{-\frac{|x-y|^2}{4{(t-s)}} } \frac{|x-y|}{2{(t-s)}}  C_3 dyds \\
& \leq \frac{(n+2)C_3}{\sqrt{\pi}} \int_0^t  \frac{1}{\sqrt{t-s}} ds
= \frac{2(n+2)C_3}{\sqrt{\pi}} \sqrt{t}.
\end{align*}

Assuming, as we shall, that $T\le   \min \{\frac{C_3}{C_1},\frac{\pi
C_2^2}{ 4(n+2)^2C_3^2} \}$, then for $0 \leq t \leq T$ we have for
all $i$,
\begin{align}
|\tilde{\rho}_{i}| \leq C_1 \qquad {\rm and} \qquad  |\nabla
\tilde{\rho}_{i}| \leq C_2. \label{rho_i est}
\end{align}

We prove the convergence of $\tilde{\rho}_i$ to a solution of
(\ref{betterharMapeq}) via proving that it is a Cauchy sequence in
$C^1$-norm.  Note that $\tilde{\rho}_i- \tilde{\rho}_{i-1}$
satisfies

\begin{align}
& \frac{\partial (\tilde{\rho}_i- \tilde{\rho}_{i-1})}{\partial t} =
\Delta (\tilde{\rho}_i- \tilde{\rho}_{i-1}) + F(x,
\tilde{\rho}_{i-1},\nabla\tilde{\rho}_{i-1},t)-F(x,\tilde{\rho}_{i-2},
\nabla\tilde{\rho}_{i-2},t) \nonumber \\
& (\tilde{\rho}_i- \tilde{\rho}_{i-1})(x,0)=   0. \label{rhoiminus}
\end{align}
where
\begin{align*}
& F(x,
\tilde{\rho}_{i-1},\nabla\tilde{\rho}_{i-1},t)-F(x,\tilde{\rho
}_{i-2},\nabla\tilde{\rho}_{i-2},t) \\
= & [(n-1) \nabla \tilde{h}- \nabla B +\nabla (\tilde{\rho}_{i-1}+
\tilde{\rho}_{i-2})]
 \cdot \nabla(\tilde{\rho}_{i-1} - \tilde{\rho}_{i-2}) \\
& + G( \tilde{\rho}_{i-1},{\bf x},t)-G(\tilde{\rho}_{i-2},{\bf x},t)
\end{align*}

By lengthy but straightforward calculations one can verify the
Lipschitz property of $G(\tilde{\rho},{\bf x},t)$

\[
|G(\tilde{\rho}_{i-1},{\bf x},t)-G(\tilde{\rho}_{i-2},{\bf x},t)|
\leq C_\&(C_1,C _2, \tilde{f},\tilde{f}_0) \cdot |
\tilde{\rho}_{i-1}-\tilde{\rho}_{i-2}|.
\]
This and (\ref{rho_i est}) implies

\begin{align}
& |F(x,
\tilde{\rho}_{i-1},\nabla\tilde{\rho}_{i-1},t)-F(x,\tilde{\rho}_{i-2},
\nabla\tilde{\rho}_{i-2},t)|
\nonumber \\
\leq &  C_4 \cdot | \tilde{\rho}_{i-1}-\tilde{\rho}_{i-2}| +C_5
\cdot |\nabla \tilde{\rho}_{i-1}- \nabla \tilde{\rho}_{i-2} |
\label{Flip}
\end{align}
where $C_4 = C_\&(C_1,C_2, \tilde{f},\tilde{f}_0)$ and $C_5 = [(n-1)
\sup |\nabla \tilde{f}| +\sup |\nabla B|+2C_2]$.

Let

\begin{align*}
& A_i(t)= \sup_{0 \leq s \leq t, x \in \mathbb{R}^{n+2}}
|\tilde{\rho}_i-
\tilde{\rho}_{i-1}|(x,s) \\
& B_i(t)=\sup_{0 \leq s \leq t,x \in \mathbb{R}^{n+2}}
|\nabla(\tilde{\rho}_i-
 \tilde{\rho}_{i-1})|(x,s).
\end{align*}

From Equations~(\ref{rhoiminus}) and~(\ref{Flip}) we can estimate
$|\tilde{\rho}_i- \tilde{\rho}_{i-1}|$ and $ |\nabla(\tilde{\rho}_i-
\tilde{\rho}_{i-1})|$ in the same way as we estimate
$|\tilde{\rho}_{i}|$ and $|\nabla\tilde{\rho}_i|$ above; we conclude

\begin{align*}
& A_i(t)\leq [ C_4  A_{i-1}(t) +C_5 B_{i-1}(t)] \cdot t  \\
& B_i(t) \leq \frac{2(n+2)[C_4  A_{i-1}(t) +C_5
B_{i-1}(t)]}{\sqrt{\pi}} \cdot \sqrt{t}.
\end{align*}

Let $C_6 = \max \{C_4,C_5 \}$; then we get

\[
A_i(t)+B_i(t) \leq \left( C_6t+ \frac{2(n+2)C_6
\sqrt{t}}{\sqrt{\pi}} \right) \cdot \left(  A_{i-1}(t)+B_{i-1}(t)
\right).
\]

Now suppose that $T\le T_2$ where $T_2$ satisfies  $C_6T_2+
\frac{2(n+2)C_6 \sqrt{T_2}}{\sqrt{\pi}} =\frac{1}{2}$; then for all
$t\le T$ we have
\[
A_i(t)+B_i(t) \leq \frac{1}{2} \left(  A_{i-1}(t)+B_{i-1}(t)
\right).
\]
This proves that $\tilde{\rho}_{i}$ is a Cauchy sequence in
$C^1(\mathbb{R}^{n+2})$. Let $\lim_{i \rightarrow +\infty}
\tilde{\rho}_{i} =\tilde{\rho}_{\infty}$. Then $\nabla
\tilde{\rho}_{i} \rightarrow \nabla \tilde{\rho}_{\infty}$ and
$F(x,\tilde{\rho}_{i-1},\nabla \tilde{\rho}_{i-1},t) \rightarrow
F(x, \tilde{\rho}_{\infty}, \nabla \tilde{\rho}_{\infty},t)$
uniformly. Hence we get from (\ref{rhoIntg}),

\begin{align} \tilde{\rho}_\infty=\int_0^t \int_{\mathbb{R}^{n+2}} H(x,y ,t-s)F(y,\tilde{ \rho}_{\infty},
\nabla\tilde{\rho}_{\infty} ,t)dyds \label{rhoInfIntg}
\end{align}

The next argument is similar to the argument in \cite{LiTam}, p.21.
The function $\tilde{\rho}_{i}$ is a smooth solution of
(\ref{rhodiffeq}) with $\tilde{\rho}_i(x,0)=0$. Also, both
 $\tilde{\rho}_{i}$
and $F(x,\tilde{\rho}_{i_1},\nabla \tilde{\rho}_{i-1},t)$ are
uniformly bounded on $\mathbb{R}^{n+2} \times[0,T]$. Thus, by
Theorem 1.11 \cite{LSU}, p.211 and Theorem 12.1 \cite{LSU}, p.223,
for any compact $K \subset \mathbb{R}^{n+2}$ and any $0<t_*<T$,
there is $C_7$ and $\alpha \in(0,1)$ independent of $i$ such that

\[
| \nabla \tilde{\rho}_{i}(x,t) -\nabla \tilde{\rho}_{i}(y,s)| \leq
C_7 \cdot \left( |x-y|^\alpha+ |t-s|^{\alpha/2} \right)
\]
where $x,y \in K$ and $0 \leq t< s \leq t_*$.

Letting $i \rightarrow \infty$ we get

\begin{align}
| \nabla \tilde{\rho}_{\infty}(x,t) -\nabla
\tilde{\rho}_{\infty}(y,s)| \leq C_7 \cdot \left( |x-y|^\alpha+
|t-s|^{\alpha/2} \right). \label{rhoHolder}
\end{align}

Hence $\nabla \tilde{\rho}_{\infty}\in C^{\alpha,\alpha/2}$, i.e.,
it  is $\alpha$-H\"{o}lder continuous in space and
$\alpha/2$-H\"older continuous.

From (\ref{rhoInfIntg}) we conclude that $\tilde{\rho}_{\infty}$ is
a solution of (\ref{betterharMapeq}) on $\mathbb{R}^{n+2} \times
[0,T]$ with $\tilde{\rho}_{\infty}(x,0)=0$.

\subsection{The asymptotic behavior of the
solutions}\label{harmsec}

In the rest of this subsection we study the asymptotic behavior of
solution $\tilde{\rho}(x,t)$  as $x \rightarrow \infty$. First we
prove inductively that there is a constant $\lambda$ and $T_3$ such
that, provided that $T\le T_3$, for $x \in \mathbb{R}^{n+2}, t \in
[0,T]$, we have

\begin{align}
| \tilde{\rho}_i (x,t) | \leq  \frac{\lambda }{(1+|x|)^2} \qquad
{\rm and} \qquad |\nabla \tilde{\rho}_i (x,t) | \leq
\frac{\lambda}{(1+|x|)^2} \label{rho_iDcay}
\end{align}

Clearly, since $\tilde \rho_0=0$, these estimates hold for $i=0$. It
follows from (\ref{rho_i est}) and Claim~\ref{asympt} that there is
a constant $C_8$ independent of $i$ such that

\begin{align*}
& | G(\tilde{\rho}_i,{\bf x},t) | \leq \frac{C_8}{(1+|x|)^2} \\
& \left[ (n-1) |\nabla \tilde{h}|+ |\nabla B| \right](x,t) \leq C_8.
\end{align*}

Now we assume these estimates hold for $i$. Then for $0 \leq t \leq
T$ we have

\begin{align*}
| \tilde{\rho}_i (x,t) | & \leq  \int_0^t \int_{\mathbb{R}^{n+2}}
H(x,y,t-s) \left[ \frac{C_8 \lambda }{(1+|y|)^2} + \frac{\lambda^2
}{(1+|y|)^2}
+\frac{C_8}{(1+|y|)^2} \right ] dyds \\
& = \int_0^t \int_{\mathbb{R}^{n+2}}
\frac{1}{(4\pi(t-s))^{(n+2)/2}}e^{ -\frac{|x-y|^2}{4(t-s)}} \left[
\frac{C_8 \lambda +\lambda^2+C_8}{(1+|y|)^2}  \right ] dyds \\
& \leq (C_8 \lambda +\lambda^2+C_8) \cdot \frac{C(n)t}{(1+|x|)^2}.
\end{align*}

Also, we have

\begin{align*}
& |\nabla \tilde{\rho}_i (x,t) | \leq  \int_0^t
\int_{\mathbb{R}^{n+2}} |\nabla_x H(x,y,t-s)| \left[
\frac{C_8 \lambda + \lambda^2 + C_8}{(1+|y|)^2} \right ]  dyds \\
& = \int_0^t \int_{\mathbb{R}^{n+2}}
\frac{|x-y|}{2(t-s)}\frac{1}{(4\pi(t-s))^{(n+2)/2}}e^{-\frac{|x-y|^2}{4(t-s)}}
\left[
\frac{C_8 \lambda +\lambda^2+C_8}{(1+|y|)^2}  \right ] dyds \\
& \leq (C_8 \lambda +\lambda^2+C_8) \cdot
\frac{C(n)\sqrt{t}}{(1+|x|)^2}.
\end{align*}

If we choose $T_3$ such that

\[
(C_8 \lambda +\lambda^2+C_8) \cdot C(n)T_3 \leq \lambda \quad
\text{and} \quad (C_8 \lambda +\lambda^2+C_8) \cdot C(n) \sqrt{T_3}
\leq \lambda,
\]
then (\ref{rho_iDcay}) hold for all $i$. From the definition of
$\tilde{ \rho}_\infty$ we conclude

\begin{align}
| \tilde{\rho}_\infty (x,t) | \leq  \frac{\lambda }{(1+|x|)^2}
\qquad \qquad |\nabla \tilde{\rho}_\infty (x,t) | \leq
\frac{\lambda}{(1+|x|)^2} \label{rhodecay}
\end{align}

Recall that $\tilde\rho_\infty$ is a solution of the following
linear equation (in $\upsilon$):

\begin{align*}
& \frac{\partial \upsilon }{\partial t} = \Delta \upsilon + \nabla [
(n-1)\tilde{h} -B] \cdot \nabla \upsilon
+G( \tilde{\rho}_\infty,{\bf x},t) \\
& \upsilon(x,0)=0.
\end{align*}

From (\ref{rhoHolder}) and Claim~\ref{asympt} we know that $\nabla [
(n-1)\tilde{h} -B+ \tilde{\rho}_\infty ]$ has
$C^{\alpha,\alpha/2}$-H\"older-norm bounded (this means
$\alpha$-H\"older norm in space and the $\alpha/2$-H\"older norm in
time). By some lengthy calculation we
 get

\[
| G(\tilde{\rho}_\infty,{\bf x},t)|_{C^{\alpha,\alpha/2}} \leq
\frac{C_9}{(1+|x|)^2}.
\]

By local Schauder estimates for parabolic equations we conclude

\[
|\tilde{\rho}_\infty |_{C^{2+\alpha,1+\alpha/2}} \leq
\frac{C_{10}}{(1+|x|)^2}.
\]

Using this estimate one can further show by calculation that

\begin{align*}
& |\nabla \nabla [ (n-1)\tilde{f} -B+
\tilde{\rho}_\infty ] |_{C^{\alpha,\alpha/2}} \leq C_{11} \\
& |\nabla G(\tilde{\rho}_\infty,{\bf x},t)|_{C^{\alpha,\alpha/2}}
\leq \frac{C_{12}}{(1+|x|)^2}.
\end{align*}

By local high order Schauder estimates for parabolic equations we
conclude

\[
|\nabla \tilde{\rho}_\infty |_{C^{2+\alpha,1+\alpha/2}} \leq
\frac{C_{13}}{(1+|x|)^2}.
 \]

We have proved the following:

\begin{prop}\label{harmapflow}
For a standard solution $(\mathbb{R}^n,g(t)),\ 0\le t<T$, and for
any $t_0\in [0,T)$ there is a rotationally symmetric solution
$\psi_t({\bf x)} = x e^{\tilde{\rho}({\bf x},t)}$ to the harmonic
map flow
\[
\frac{\partial \psi_t}{\partial t}=\Delta_{g(t),g(t_0)}\psi(t)
\qquad \psi(t_0)({\bf x})={\bf x},
\]
and $ |\nabla^i \tilde{\rho}|({\bf x},t)\leq \frac{C_{14}}{(1+|{\bf
x}|)^2}$ for $0 \leq i \leq 3$ defined on some non-degenerate
interval $[t_0,T']$.
\end{prop}

\subsection{The uniqueness for the solutions of Ricci-DeTurck flow}

We prove the following general uniqueness result for Ricci-DeTurck
flow on open manifolds.

\begin{prop}\label{RDuniq}
Let $\hat g_1(  t)  $ and $\hat g_2(  t),\ 0 \leq t\leq T $, be two
bounded solutions of the Ricci-DeTurck flow on complete and
noncompact manifold $M^{n}$ with initial metric $g_1( t_0) =g_2(t_
0)  =g$. Suppose that for some $1<C<\infty$ we have
\begin{align*}
& C^{-1}  g    \leq \hat g_1( t) \leq C  g \\
& C^{-1}  g   \leq \hat g_2(  t) \leq C g .\\
\end{align*}
Suppose that in addition we have

\begin{align*}
& \left\Vert \hat g_1(t) \right\Vert _{C^{2}\left(M\right),g} \le C \\
& \left\Vert \hat g_2(t) \right\Vert _{C^{2}\left(M\right),g} \le
C.
\end{align*}

Lastly, suppose there is an exhausting sequence of compact, smooth
submanifolds of $\Omega_{k}\subset M$, i.e., $\Omega _{k}\subset{\rm
int}\Omega_{k+1}$ and $\cup\Omega_{k}=M$ such that $\hat g_1\left(
t\right)$ and $\hat g_2(t)$ have the same sequential asymptotic
behavior at $\infty$ in the sense that for any $\epsilon>0,$ there
is a $k_0$ arbitrarily large with

\[
\left\vert \hat g_1(  t)  -\hat g_2(  t) \right\vert
_{C^{1}\left(\partial\Omega_{k_0}\right) ,g}\leq\epsilon,
\]
Then $ \hat g_1(  t)  =\hat g_2(  t)$.
\end{prop}

\begin{proof}
Letting $\tilde \nabla$ be the covariant derivative determined by
$g$, then, using the Ricci-DeTurck flow (\ref{eqRcDeTurck}) for
$\hat g_1$ and $\hat g_2$, we can  make the following estimate for
an appropriate constant $C_{14}$ depending on $g$.
\begin{align*}
& \frac{\partial}{\partial t}\left\vert \hat g_1(  t) -\hat g_2( t)
\right\vert _{g}^{2} =2\left\langle \frac{\partial}{\partial
t}\left(  \hat g_1(  t) -\hat g_2( t)  \right) ,\hat g_1( t)
-\hat g_2(  t)  \right\rangle _{g}   \\
& \leq 2\left\langle\hat
g_1^{\alpha\beta}\tilde{\nabla}_{\alpha}\tilde{\nabla}_{\beta}\left(
\hat g_1(t)-\hat g_2(t)\right),  \left(  \hat g_1( t) -\hat g_2(t)
\right)\right\rangle_g
\\ &  +C_{14}\left\vert \hat g_1(  t)  -\hat g_2(
t) \right\vert _{g}^{2}  +C_{14}\left\vert \tilde\nabla\left( \hat
g_1( t) - \hat g_2( t)\right)  \right\vert _{g} \left\vert \hat g_1(
t)
-\hat g_2(  t)  \right\vert _{g}\\
& \leq \hat
g_1^{\alpha\beta}\tilde{\nabla}_{\alpha}\tilde{\nabla}_{\beta}\left(\left\vert
\hat g_1(  t)  -\hat g_2(  t)  \right\vert _{g}^{2}\right) -2\hat
g_1^{\alpha\beta}\left\langle\tilde{\nabla}_{\beta}\left( \hat
g_1(t)-\hat g_2(t)\right),\tilde{\nabla}_{\alpha}\left(  \hat g_1(
t)  -\hat g_2 (t)  \right)\right\rangle_g \\
&  +C_{14}\left\vert \hat g_1(t)-\hat g_2(t) \right\vert _{g}^{2}
+C_{14}\left\vert \tilde\nabla\left(  \hat g_1(t)-\hat g_2(t)
\right) \right\vert
_{g}\left\vert \hat g_1(t) - \hat g_2(t)\right\vert _{g}\\
&  \leq   \hat
g_1^{\alpha\beta}\tilde{\nabla}_{\alpha}\tilde{\nabla}_{\beta}\left(\left\vert
\hat g_1\left(  t\right)  -\hat g_2\left(  t\right)  \right\vert
_{g}^{2}\right)-2C^{-1}  \left\vert \tilde\nabla\left( \hat
g_1(t)-\hat g_2(t)\right)
\right\vert _{g}^{2}\\
&  +C_{14}\left\vert \hat g_1(t)-\hat g_2(t)\right\vert
_{g}^{2}+C^{-1} \left\vert \tilde\nabla\left( \hat g_1(t)-\hat
g_2(t) \right) \right\vert _{g}^{2}+\frac{C_{14}^{2}}{4C^{-1}
}\left\vert \hat g_1(t)-\hat g_2(t) \right\vert _{g}^{2},
\end{align*}
where the last inequality comes from completing the square to
replace the last term in the previous expression. Thus, we have
proved

\begin{equation}
\frac{\partial}{\partial t}\left\vert \hat g_1\left(  t\right) -\hat
g_2\left( t\right) \right\vert _{g}^{2}\leq 2 \hat
g_1^{\alpha\beta}\tilde{\nabla}_{\alpha}
\tilde{\nabla}_{\beta}\left\vert \hat g_1\left(  t\right)  -\hat
g_2\left( t\right) \right\vert _{g}^{2}+C_{15}\left\vert \hat
g_1\left( t\right)  -\hat g_2\left( t\right)  \right\vert _{g}^{2}
\label{eq est for unique}
\end{equation}
pointwise on $\Omega_{k}$ with $C_{15}$ a constant that depends only
on $n$, $C$ and $g$.

Suppose that $\hat g_1\left(  t\right)  \neq  \hat g_2 \left(
t\right) $ for some $t$. Then there is a point $x_{0}$ such that
$\left\vert \hat g_1(x_0, t) - \hat g_2(x_0, t) \right\vert _{g}^{2}
>\epsilon_{0}$ for  some $\epsilon_{0}>0$.

We choose a $k_0$ sufficiently large that $x_{0}\in \Omega_{k_0}$
and for all $t'\in [t_0,T]$ we have

\begin{equation}
\sup_{x\in\partial\Omega_{b}}\left\vert \hat g_1(x,  t') -\hat
g_2(x, t') \right\vert _{g}^{2} \leq \epsilon \label{eq est for uniq
bddy}
\end{equation}
where $\epsilon >0$ is a constant to be chosen later.

Recall we have the initial condition $\left\vert \hat g_1( 0) -\hat
g_2(0)  \right\vert _{g}^{2}=0$. Using Equation~(\ref{eq est for
unique}) and applying the maximum principle to $\left\vert \hat g_1(
t) -\hat g_2( t) \right\vert _{g}^{2}$ ) on the domain
$\Omega_{k_0}$, we get
\[
e^{-C_{15}t}\left\vert \hat g_1 (  t)  -\hat g_2 ( t) \right\vert
_{g}^{2}(  x) \leq\epsilon.\text{ for all }x\in\Omega _{k_0}.
\]

This is a contradiction if we choose $\epsilon \leq
\epsilon_{0}e^{-C_{15}T}$. This contradiction establishes the
proposition.
\end{proof}

 Let $g_1(t),\ 0\le t<T_1$, and $g_2(t),\ 0\le t<T_2$, be standard solutions that
 agree on the interval $[0,t_0]$ for some $t_0\ge 0$.
By Proposition~\ref{harmapflow} there are $\psi_1(t)$  and
$\psi_2(t)$  which are solutions of the harmonic map flow defined
for $t_0 \leq t \leq T$ for some $T>t_0$ for the Ricci flows
$g_1(t)$ and $g_2(t)$. Let $\hat{g}_1(t)= (\psi^{-1}(t))^*g_1(t)$
and $\hat{g}_2(t) = (\psi^{-1}(t))^*g_2(t)$. Then $\hat{g}_1(t)$ and
$\hat{g}_2(t)$ are two solutions of the Ricci-DeTurck flow with
$\hat{g}_1(t_0)=\hat{g}_2(t_0)$. Choose $T' \in (t_0,T]$ such that
$\hat{g}_1(t)$ and $\hat{g}_2(t)$ are $\delta$-close to
$\hat{g}_1(t_0)$ as required in Proposition~\ref{RDuniq}. It follows
from Lemma~\ref{item5} and the decay estimate in
Proposition~\ref{harmapflow} that $\hat{g}_1(t)$ and $\hat{g}_2(t)$
are bounded solutions and that they have same sequential asymptotic
behavior at infinity. We can apply Proposition~\ref{RDuniq} to
conclude $\hat{g}_1(t)=\hat{g}_2(t)$ on $t_0 \leq t \leq T'$.  We
have proved:

\begin{cor}\label{duniq}
Let $g_1(t)$ and $g_2(t)$  be standard solutions. Suppose that
$g_1(t)=g_2(t)$ for all $t\in [0,t_0]$ for some $t_0\ge 0$. The
Ricci-DeTurck solutions $\hat{g}_1(t)$ and $\hat{g}_2(t)$
constructed from standard solutions $g_1(t)$ and $g_2(t)$ with
$g_1(t_0)=g_2(t_0)$ exist and satisfy $\hat{g}_1(t)=\hat{g}_1(t)$
for $t \in [t_0,T']$ for some $T'>t_0$.
\end{cor}

\section{Completion of the proof of uniqueness}

Now we are ready to prove the uniqueness of the standard solution.
Let $g_1(t),\ 0\le t<T_1$, and $g_2(t),\ 0\le t<T_2$, be a standard
solutions. Consider the maximal interval $I$ (closed or half-open)
containing $0$ on which $g_1$ and $g_2$ agree.

{\bf Case 1: $T_1< T_2$ and $I=[0,T_1)$}

In this case since $g_1(t)=g_2(t)$ for all $t<T_1$ and $g_2(t)$
extends smoothly past time $T_1$, we see that the curvature of
$g_1(t)$ is bounded as $t$ tends to $T_1$. Hence, $g_1(t)$ extends
past time $T_1$, contradicting the fact that it is a maximal flow.

{\bf Case 2: $T_2<T_1$ and $I=[0,T_2)$}

The argument in this case is the same as the previous one with the
roles of $g_1(t)$ and $g_2(t)$ reversed.

There is one more case to rule out.

{\bf Case 3: $I$ is a closed interval $I=[0,t_0]$.}

In this case, of course, $t_0<{\rm min}(T_1,T_2)$. Hence we apply
Proposition~\ref{harmapflow} to construction solutions $\psi_1$ and
$\psi_2$ to the harmonic map flow for $g_1(t)$ and $g_2(t)$ with
$\psi_1$ and $\psi_2$ being the identity at time $t_0$. These
solutions will be defined on an interval of the form $[t_0,T]$ for
some $T>t_0$. Using these harmonic map flows we construct solutions
$\hat g_1(t)$ and $\hat g_2(t)$ to the Ricci-DeTurck flow defined on
the interval $[t_0,T]$. According to Corollary~\ref{duniq}, there is
a uniqueness theorem for these Ricci-DeTurck flows, which implies
that $\hat g_1(t)=\hat g_2(t)$ for all $t\in [t_0,T']$ for some
$T'>t_0$. Invoking Corollary~\ref{dtuniqtouniq} we conclude that
$g_1(t)=g_2(t)$ for all $t\in [0,T']$, contradicting the maximality
of the interval $I$.

If none of these three cases can occur, then the only remaining
possibility is that $T_1=T_2$ and $I=[0,T_1)$, i.e., the flows are
the same. This then completes the proof of the uniqueness of the
standard flow.

\subsection{$T=1$ and existence of canonical neighborhoods}

At this point we have established all the properties claimed in
Theorem~\ref{stdsoln} for the standard flow except for the fact that
$T$, the endpoint of the time-interval of definition, is equal to
$1$. We have shown that $T\le 1$. In order to establish the opposite
inequality, we must show the existence of canonical neighborhoods
for the standard solution.

Here is the result about the existence of canonical neighborhoods
for the standard solution.

\begin{thm}\label{R_0cannbhd}
Fix $0<\epsilon<1$.   Then there is $r>0$ such that for any point
$(x_0,t_0)$ in the standard flow with $R(x_0,t_0)\ge r^{-2}$ the
following hold.
\begin{enumerate}
\item[(1)] $t_0>r^2$.
\item[(2)] $(x_0,t_0)$ has a strong canonical
$(C(\epsilon),\epsilon)$-neighborhood\index{canonical
neighborhood!strong}. If this canonical neighborhood is a strong
$\epsilon$-neck centered at $(x_0,t_0)$, then the strong neck
extends to an evolving neck defined for backward rescaled time
$(1+\epsilon)$.
\end{enumerate}
\end{thm}

\begin{proof}
 Take an increasing sequence of times $t'_n$ converging to $T$.
Since the curvature of $(\Ar^3,g(t))$ is locally bounded in time,
for each $n$, there is a bound on the scalar curvature on
$\Ar^3\times [0,t'_n]$. Hence, there is a finite upper bound $R_n$
on $R(x,t)$ for all points $(x,t)$ with $t\le t'_n$ for which the
conclusion of the theorem does not hold. (There clearly are such
points since the conclusion of the theorem fails for all $(x,0)$.)
Pick $(x_n,t_n)$ with $t_n\le t'_n$, with $R(x_n,t_n)\ge R_n/2$ and
such that the conclusion of the theorem does not hold for
$(x_n,t_n)$. To prove the theorem we must show that $\overline{\rm
lim}_{n\rightarrow\infty}R(x_n,t_n)<\infty$. Suppose the contrary.
By passing to a subsequence we can suppose that ${\rm
lim}_{n\rightarrow \infty}R(x_n,t_n)=\infty$. We set
$Q_n=R(x_n,t_n)$. We claim that all the hypotheses of
Theorem~\ref{kaplimit} apply to the sequence
$(\Ar^3,g(t),(x_n,t_n))$. First, we show that all the hypotheses of
Theorem~\ref{smlmtflow} (except the last) hold. Since $(\Ar^3,g(t))$
has non-negative curvature all these flows have curvature pinched
toward positive. By Theorem~\ref{stdkappa} there are $r>0$ and
$\kappa>0$ so that all these flows are $\kappa$-non-collapsed on
scales $\le r$. By construction if $t\le t_n$ and $R(y,t)>2Q_n\ge
R_n$ then  the point $(y,t)$ has a strong canonical
$(C(\epsilon),\epsilon)$-neighborhood. We are assuming that
$Q_n\rightarrow \infty$ as $n\rightarrow \infty$ in order to achieve
the contradiction. Since all time-slices are complete, all balls of
finite radius have compact closure.

Lastly, we need to show that the extra hypothesis of
Theorem~\ref{kaplimit} (which includes the last hypothesis of
Theorem~\ref{smlmtflow}) is satisfied. This is clear since
$t_n\rightarrow T$ as $n\rightarrow \infty$ and $Q_n\rightarrow
\infty$ as $n\rightarrow \infty$. Applying Theorem~\ref{kaplimit} we
conclude that after passing to a subsequence there is a limiting
flow which is a $\kappa$-solution. Clearly, this and
Corollary~\ref{limitcannbhd} imply that for all sufficiently large
$n$ (in the subsequence)  the neighborhood as required by the
theorem exists. This contradicts our assumption that none of the
points $(x_n,t_n)$ have these neighborhoods. This contradiction
proves the result.
\end{proof}

\subsection{Completion of the proof of
Theorem~\protect{\ref{stdsoln}}}

The next proposition establishes the last of the conditions claimed
in Theorem~\ref{stdsoln}.

\begin{thm}
For the standard flow $T=1$.
\end{thm}

\begin{proof}
We have already seen in Corollary~\ref{1} that $T\le 1$. Suppose now
that $T<1$. Take $T_0<T$ sufficiently close to $T$. Then according
to Proposition~\ref{item5} there is a compact subset $X\subset
\Ar^3$ such that restriction of the flow to $(\Ar^3\setminus
X)\times [0,T_0]$ is $\epsilon$-close to the standard evolving flow
on $S^2\times (0,\infty),(1-t)h_0\times ds^2$, where $h_0$ is the
round metric of scalar curvature $1$ on $S^2$. In particular,
$R(x,T_0)\le (1+\epsilon)(1-T_0)^{-1}$ for all $x\in \Ar^3\setminus
X$. Because of Theorem~\ref{R_0cannbhd} and the definition of
$(C(\epsilon),\epsilon)$-canonical neighborhoods, it follows that at
any point $(x,t)$ with $R(x,t)\ge r^{-2}$, where $r>0$ is the
constant given in Theorem~\ref{R_0cannbhd}, we have $\partial
R/\partial t(x,t)\le C(\epsilon)R^2(x,t)$. Thus, provided that
$T-T_0$ is sufficiently small, there is a uniform bound to $R(x,t)$
for all $x\in \Ar^3\setminus X$ and all $t\in [T_0,T)$. Using
Theorem~\ref{shiw/deriv} and the fact that the standard flow is
$\kappa$-non-collapsed implies that the restrictions of the metrics
$g(t)$ to $\Ar^3\setminus X$ converge smoothly to a limiting
Riemannian metric $g(T)$ on $\Ar^3\setminus X$. Fix a non-empty open
subset $\Omega\subset \Ar^3\setminus X$ with compact closure. For
each $t\in [0,T)$ let $V(t)$ be the volume of
$(\Omega,g(t)|_{\Omega})$. Of course, ${\rm lim}_{t\rightarrow
T}V(t)={\rm Vol}_{g(T)}(\Omega)>0$.

Since the metric $g(T)$ exists in a neighborhood of infinity and has
bounded curvature there, if the limit metric $g(T)$ exists on all of
$\Ar^3$, then we can extend the flow keeping the curvature bounded.
This contradicts the maximality of our flow subject to the condition
that the curvature be locally bounded in time. Consequently, there
is a point $x\in \Ar^3$ for which the limit metric $g(T)$ does not
exist. This means that $\overline{\rm lim}_{t\rightarrow
T}R(x,t)=\infty$. That is to say,  there is a sequence of
$t_n\rightarrow T$ such that setting $Q_n=R(x,t_n)$, we have
$Q_n\rightarrow \infty$ as $n$ tends to infinity. By
Theorem~\ref{R_0cannbhd} the second hypothesis in the statement of
Theorem~\ref{smlmtflow} holds for the sequence
$(\Ar^3,g(t),(x,t_n))$. All the other hypotheses of this theorem as
well as the extra hypothesis in Theorem~\ref{kaplimit} obviously
hold for this sequence. Thus, according to Theorem~\ref{kaplimit}
the based flows $(\Ar^3,Q_ng(Q_n^{-1}t'+t_n),(x,0))$ converge
smoothly to a $\kappa$-solution. Since the asymptotic volume of any
$\kappa$-solution is zero (see Theorem~\ref{asympvol}), we see that
for all $n$ sufficiently large, the following holds:

\begin{claim}
For any $\epsilon>0$, there is $A<\infty$ such that for all $n$
sufficiently large we have
$${\rm Vol}(B_{Q_ng}(x,t_n,A))<\epsilon A^3.$$
\end{claim}

Rescaling, we see that for all $n$ sufficiently large we have
$${\rm Vol}\,B_g(x,t_n,A/\sqrt{Q_n})< \epsilon(A/\sqrt{Q_n})^3.$$
Since the curvature of $g(t_n)$ is non-negative and since the $Q_n$
tend to $\infty$, it follows from the Bishop-Gromov Inequality
(Proposition~\ref{BishopGromov}) that for any $0<A<\infty$ and any
$\epsilon>0$, for all $n$ sufficiently large we have
$${\rm Vol}\,B_g(x,t_n,A)<\epsilon A^3.$$

On the other hand, since $\Omega$ has compact closure, there is an
$A_1<\infty$ with $\Omega\subset B(x,0,A_1)$. Since the curvature of
$g(t)$ is non-negative for all $t\in [0,T)$, it follows from
Lemma~\ref{posdistdec} that the distance is a non-increasing
function of $t$, so that for all $t\in [0,T)$ we have $\Omega\subset
B(x,t,A_1)$. Applying the above, for any $\epsilon>0$ for all $n$
sufficiently large we have
$${\rm Vol}\,(\Omega,g(t_n))\le {\rm Vol}_g\,B(x,t_n,A_1)<\epsilon
(A_1)^3.$$ But this contradicts the fact that
$${\rm lim}_{n\rightarrow
\infty}{\rm Vol}\,\Omega,g(t_n)={\rm Vol}\,(\Omega,g(T))>0.$$ This
contradiction proves that $T=1$.
\end{proof}

This completes the proof of Theorem~\ref{stdsoln}.

\section{Some corollaries}

Now let us derive extra properties of the standard solution that
will be important in our applications.

\begin{prop}\label{Restim}
There is a constant $c>0$ such that for all $(p,t)$ in the standard
solution we have
$$R(p,t)\ge \frac{c}{1-t}.$$
\end{prop}

\begin{proof}
First, let us show that there is not a limiting metric $g(1)$
defined on all of $\Ar^3$. This does not immediately contradict the
maximality of the flow because we are assuming only that the flow is
maximal subject to having curvature locally bounded in time. Assume
that a limiting metric $(\Ar^3,g(1))$ exists.
 First,
notice that from the canonical neighborhood assumption and
Lemma~\ref{controlnbhd} we see that the curvature of $g(T)$ must be
unbounded at spatial infinity. On the other hand, by
Proposition~\ref{canonvary} every point of $(\Ar^3,g(1))$ of
curvature greater than $R_0$ has a $(2C,2\epsilon)$-canonical
neighborhood. Hence, since $(\Ar^3,g(1))$ has unbounded curvature,
it then has $2\epsilon$-necks of arbitrarily small scale. This
contradicts Proposition~\ref{narrows}. (One can also rule this
possibility out by direct computation using the spherical symmetry
of the metric.) This means that there is no limiting metric $g(1)$.

The next step is to see that for any $p\in \Ar^3$ we have ${\rm
lim}_{t\rightarrow 1}R(p,t)=\infty$. Let $\Omega\subset \Ar^3$ be
the subset of $x\in \Ar^3$ for which ${\rm liminf}_{t\rightarrow
1}R(x,t)<\infty$. We suppose that $\Omega\not=\emptyset$. According
to Theorem~\ref{omega} the subset $\Omega$ is open and the metrics
$g(t)|_\Omega$ converge smoothly to a limiting metric
$g(1)|_\Omega$. On the other hand, we have just seen that there is
not a limit metric $g(1)$ defined everywhere. This means that there
is $p\in \Ar^3$ with  ${\rm lim}_{t\rightarrow 1}R(p,t)=\infty$.
Take a sequence $t_n$ converging to $1$ and set $Q_n=R(p,t_n)$. By
Theorem~\ref{kaplimit} we see that, possibly after passing to a
subsequence, the based flows $(\Ar^3,Q_ng(t'-t_n),(p,0))$ converge
to a $\kappa$-solution. Then by Proposition~\ref{asympvol} for any
$\epsilon>0$ there is $A<\infty$ such that ${\rm
Vol}\,B_{Q_ng}(p,t_n,A)<\epsilon A^3$, and hence after rescaling we
have ${\rm Vol}\,B_g(p,t_n,A/\sqrt{Q_n})<\epsilon (A/\sqrt{Q_n})^3$.
By the Bishop-Gromov inequality (Proposition~\ref{BishopGromov}) it
follows that for any $0<A<\infty$, any $\epsilon>0$ and for all $n$
sufficiently large, we have ${\rm Vol}\,B_g(p,t_n,A)<\epsilon A^3$.
Take a non-empty subset $\Omega'\subset \Omega$ with compact
closure. Of course, ${\rm Vol}\,(\Omega',g(t))$ converges to ${\rm
Vol}\,(\Omega',g(T))>0$ as $t\rightarrow T$. Then there is
$A<\infty$ such that for each $n$, the subset $\Omega'$ is contained
in the ball $B(p_0,t_n,A)$. This is a contradiction since it implies
that for any $\epsilon>0$ for all $n$ sufficiently large we have
${\rm Vol}\,(\Omega',g(t))<\epsilon A^3$. This completes the proof
that for every $p\in \Ar^3$ we have ${\rm lim}_{t\rightarrow
1}R(p,t)=\infty$.

Fix $\epsilon>0$ sufficiently small and set $C=C(\epsilon)$. Then
for every $(p,t)$ with $R(p,t)\ge r^{-2}$ we have
$$\left|\frac{dR}{dt}(p,t)\right|\le CR^2(p,t).$$
 Fix $t_0=1-1/2r^2C$. Since the flow has curvature locally bounded in
 time, there is $2C\le C'<\infty$ such that $R(p,t_0)\le 1/(C'(1-t_0)$ for
 all $p\in \Ar ^3$.
 Since $R(p,t_0)=1/C'(1-t_0)$,
for all $t\in [t_0,1)$ we have
$$R(p,t)<{\rm
max}\left(\left[(C'-C)(1-t_0)\right]^{-1},\left[r^{-2}-C(1-t_0)\right]^{-1}\right).$$

This means that $R(p,t)$ is uniformly bounded as $t\rightarrow 1$,
contradicting what we just established. This shows that for $t\ge
1-1/2r^2C$ the result holds. For $t\le 1-1/2r^2C$ there is a
positive lower bound on the scalar curvature, and hence the result
is immediate for these $t$ as well.
\end{proof}

\begin{figure}[ht]
  \centerline{\epsfbox{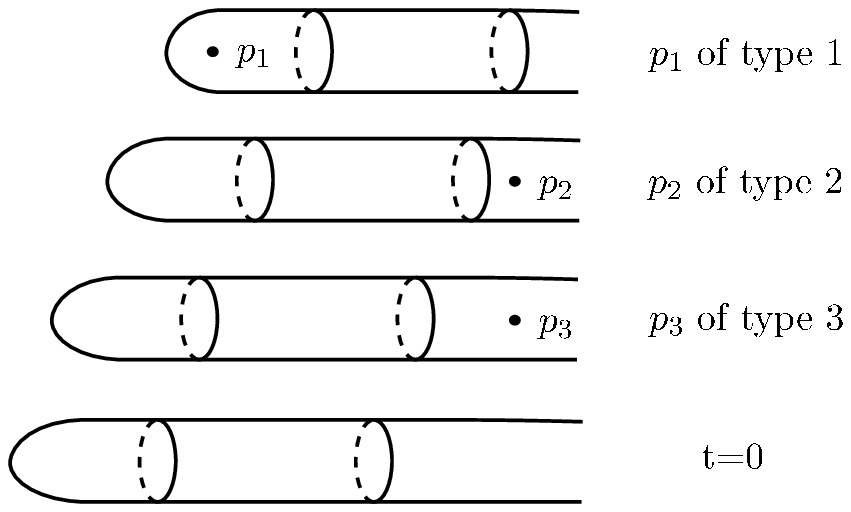}}
  \caption{Canonical neighborhoods in the standard solution.}\label{fig:types}
\end{figure}

\begin{thm}\label{stdsolncannbhd}
For any $\epsilon>0$ there is $C'(\epsilon)<\infty$ such that for
any point $x$ in the standard solution one of the following holds
(see {\sc Fig.}~\ref{fig:types}).
\begin{enumerate}
\item[(1)] $(x,t)$ is contained in the core of a  $(C'(\epsilon),\epsilon)$-cap.
\item[(2)] $(x,t)$ is the center of an evolving $\epsilon$-neck $N$ whose initial time-slice is
$t=0$ and whose initial time-slice is disjoint from the surgery cap.
\item[(3)] $(x,t)$ is the center of an evolving $\epsilon$-neck defined for rescaled time
$1+\epsilon$.
\end{enumerate}
\end{thm}

\begin{rem}
At first glance it may seem impossible for a point $(x,t)$ in the
standard solution to be the center of an evolving $\epsilon$-neck
defined for rescaled time $1+\epsilon$ since the standard solution
itself is only defined for time $1$. But this is indeed possible.
The reason is because the scale referred to for an evolving neck
centered at $(x,t)$ is $R(x,t)^{-1/2}$. As $t$ approaches one,
$R(x,t)$ goes to infinity, so that rescaled time $1$ at $(x,t)$ is
an arbitrarily small time interval measured in the scale of the
standard solution.
\end{rem}

\begin{proof}
By Theorem~\ref{R_0cannbhd}, there is $r_0$ such that if $R(x,t)\ge
r_0^{-2}$, then $(x,t)$ has a  $(C,\epsilon)$-canonical neighborhood
and if this canonical neighborhood is a strong $\epsilon$-neck
centered at $x$, then that neck extends to an evolving neck defined
for rescaled time $(1+\epsilon)$. By Proposition~\ref{Restim}, there
is $\theta<1$ such that if $R(x,t)\le r_0^{-2}$ then $t\le \theta$.
By Proposition~\ref{item5},  there is a compact subset $X\subset \Ar
^3$ such that if $t\le \theta$ and $x\notin X$, then there is an
evolving $\epsilon$-neck centered at $x$ whose initial time is zero
and whose initial time-slice is at distance at least one from the
surgery cap. Lastly, by compactness there is $C'<\infty$ such that
every $(x,t)$ for $x\in X$ and every $t\le \theta$ is contained in
the core of a $(C',\epsilon)$-cap.\end{proof}

\begin{cor}\label{stdsolncap}
Fix $\epsilon>0$. Suppose that $(q,t)$ is a point in the standard
solution with $t\le R(q,t)^{-1}(1+\epsilon))$ and with $(q,0)\in
B(p_0,0,(\epsilon^{-1}/2)+A_0+5)$. Then $(q,t)$ is contained in an
$(C'(\epsilon),\epsilon)$-cap.
\end{cor}

\begin{rem}
Recall that $p_0$ is the origin in $\Ar^3$ and hence is the tip of
the surgery cap. Also, $A_0$ is defined in Lemma~\ref{A_0}.
\end{rem}

\begin{cor}\label{stdsolnlimit}
For any $\epsilon>0$ let $C'=C'(\epsilon)$ be as in
Theorem~\ref{stdsolncannbhd}. Suppose that we have a sequence of
generalized Ricci flows $({\mathcal M}_n,G_n)$, points $x_n\in
{\mathcal M}_n$ with ${\bf t}(x_n)=0$, neighborhoods $U_n$ of $x_n$
in the zero time-slice of ${\mathcal M}_n$, and a constant
$0<\theta<1$. Suppose that there are embeddings $\rho_n\colon
U_n\times [0,\theta)\to {\mathcal M}_n$ compatible with time and the
vector field so that the Ricci flows $\rho_n^*G_n$ on $U_n$ based at
$x_n$ converge geometrically to the restriction of the standard
solution to $[0,\theta)$. Then for all $n$ sufficiently large, and
any point $y_n$ in the image of $\rho_n$ one of the following holds:
\begin{enumerate}
\item[(1)] $y_n$ is contained in the core of a
$(C'(\epsilon),\epsilon)$-cap
\item[(2)] $y_n$ is the center of a strong $\epsilon$-neck
\item[(3)] $y_n$ is the center of an evolving $\epsilon$-neck whose
initial time-slice is at time $0$.
\end{enumerate}
\end{cor}

\begin{proof}
This follows immediately from Theorem~\ref{stdsolncannbhd} and
Proposition~\ref{canonvary}.
\end{proof}

There is one property that we shall use later in proving the
finite-time extinction of Ricci flows with surgery for manifolds
with finite fundamental group (among others). This is a distance
decreasing property which we record here.

Notice that for the standard initial metric constructed in
Lemma~\ref{stdinit} we have the following:

\begin{lem}\label{decrease}
Let $S^2$ be the unit sphere in $T_0\Ar^3$. Equip it with the metric
$h_0$ that is twice the usual metric (so that the scalar curvature
of $h_0$ is $1$). We define a map $\rho\colon S^2\times
[0,\infty)\to \Ar^3$ by sending the point $(x,s)$ to the point at
distance $s$ from the origin in the radial direction from $0$ given
by $x$ (all this measured in the metric $g_0$). Then $\rho^*g_0\le
h_0\times ds^2$.
\end{lem}

\begin{proof}
Clearly, the metric $\rho^*g_0$ is rotationally symmetric and its
component in the $s$-direction is $ds^2$. On the other hand, since
each cross section $\{s\}\times S^2$ maps conformally onto a sphere
of radius $\le\sqrt{2}$ the result follows.
\end{proof}

\chapter{Surgery on a $\delta$-neck}

\section{Notation and the Statement of the Result}\label{necksurgery}

In this chapter we describe the surgery process. For this chapter we
fix:

\begin{enumerate}
\item[(1)] A $\delta$-neck $(N,g)$ centered at a point $x_0$.
We denote by $\rho\colon  S^2\times (-\delta^{-1},\delta^{-1})\to N$
the diffeomorphism that gives the $\delta$-neck structure.
\item[(2)] Standard initial conditions $(\Ar^3,g_0)$.
\end{enumerate}
We denote by $h_0\times ds^2$ the metric on $S^2\times \Ar$ which is
the product of the round metric $h_0$ on $S^2$ of scalar curvature
$1$ and the Euclidean metric $ds^2$ on $\Ar$. We denote by
$N^-\subset N$ the image $\rho((-\delta^{-1},0]\times S^2)$ and we
denote by $s\colon N^-\to (-\delta^{-1},0]$ the composition
$\rho^{-1}$ followed by the projection to the second factor.

Recall that the standard initial metric $(\Ar^3,g_0)$ is invariant
under the standard $SO(3)$-action on $\Ar^3$. We let $p_0$ denote
the origin in $\Ar^3$. It is  the fixed point of this action and is
called the  {\em tip}\index{tip, of the standard initial metric|ii}
of the standard initial metric. Recall from Lemma~\ref{A_0} that
there are $A_0>0$ and an isometry
$$\psi\colon (S^2\times
(-\infty,4],h_0\times ds^2)\to (\Ar^3\setminus B(p_0,A_0),g_0).$$ The
composition of $\psi^{-1}$ followed by projection onto the second
factor defines a map $s_1\colon \Ar^3\setminus B(p_0,A_0)\to
(-\infty,4]$. Lastly, there is $0<r_0<A_0$ such that on $B(p_0,r_0)$
the metric $g_0$ is of constant sectional curvature $1/4$. We extend
the map $s_1$ to a continuous map $s_1\colon \Ar^3\to
(-\infty,4+A_0]$ defined by $s_1(x)=A_0+4-d_{g_0}(p,x)$.  This map
is an isometry along each radial geodesic ray emanating from $p_0$.
It is smooth except at $p_0$ and sends $p_0$ to $4+A_0$. The
pre-images of $s_1$ on $(-\infty,4+A_0)$ are $2$-spheres with round
metrics of scalar curvature at least $1$.

The surgery\index{surgery, on a $\delta$-neck|ii} process is a local
one defined on the $\delta$-neck $(N,g)$. The surgery process
replaces  $(N,g)$ by a smooth Riemannian manifold $({\mathcal
S},\widetilde g)$. The underlying smooth manifold ${\mathcal S}$ is
obtained by gluing together $\rho(S^2\times(-\delta^{-1},4))$ and
$B(p_0,A_0+4)$ by identifying $\rho(x,s)$ with $\psi(x,s)$ for all
$x\in S^2$ and all $s\in (0,4)$. The functions $s$ on $ N^-$ and
$s_1$ agree on their overlap and hence together define a function
$s\colon {\mathcal S}\to (-\delta^{-1},4+A_0]$, a function smooth
except at $p_0$. In order to define the metric $\widetilde g$ we
must make some universal choices. We fix once and for all two bump
functions $\alpha\colon [1,2]\to [0,1]$, which is required to be
identically $1$ near $1$ and identically $0$ near $2$, and
$\beta\colon [4+A_0-r_0,4+A_0]\to [0,1]$, which is required to be
identically $1$ near $4+A_0-r_0$ and identically $0$ on
$[4+A_0-r_0/2,A_0]$. These functions are chosen once and for all and
are independent of $\delta$ and $(N,g)$. Next we set
$\eta=\sqrt{1-\delta}$.
The purpose of this choice is the following:

\begin{claim}\label{distdecrease}
Let $\xi\colon N\to \Ar^3 $ be the map that sends $\rho(S^2\times
[A_0+4,\delta^{-1}))$ to the origin $0\in \Ar^3$ (i.e., to the tip of
the surgery cap) and for every $s<A_0+4$ sends $(x,s)$ to the point
in $\Ar^3$ in the radial direction $x$ from the origin at
$g_0$-distance $A_0+4-s$. Then $\xi$ is a distance decreasing map
from $(N,R(x_0)g)$ to $(\Ar^3,\eta g_0)$.
\end{claim}

\begin{proof}
Since $R(x_0)g$ is within $\delta$ of $h_0\times ds^2$, it follows
that $R(x_0)g\ge \eta (h_0\times ds^2)$. But according to
Lemma~\ref{decrease} the map $\xi$ given in the statement of the
claim is a distance non-increasing map from $h_0\times ds^2$ to
$g_0$. The claim follows immediately.
\end{proof}

\begin{figure}[ht]
  \centerline{\epsfbox{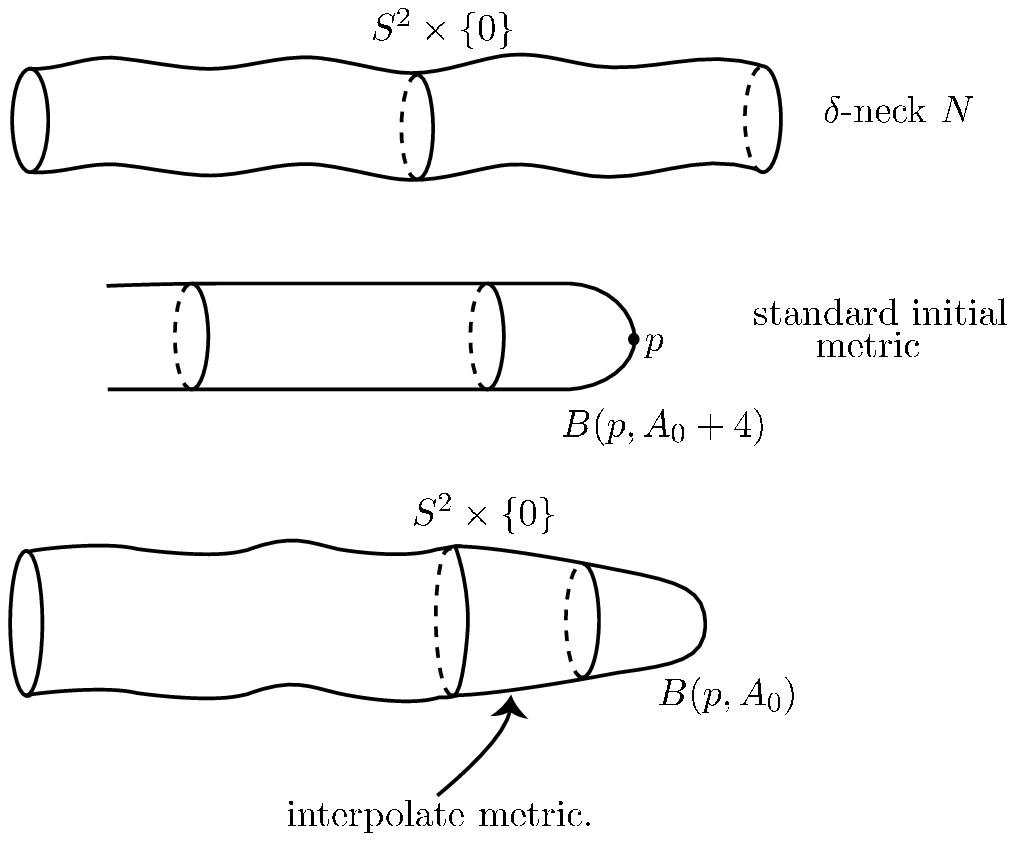}}
 \caption{Local Surgery}\label{fig:locsurg}
\end{figure}

 The last choices we need to make are of
constants $C_0<\infty$ and $q<\infty$, with $C_0 \gg q$, but both of these are
independent of $\delta$. These choices will be made later. Given all these
choices, we define a function
$$f(s)=\begin{cases} 0 & s\le 0 \\ C_0\delta e^{-q/s} & s>0,\end{cases}$$
and then we define the metric $\tilde g$ on ${\mathcal S}$ by first
defining a metric:
$$\hat g=\begin{cases} {\rm exp}(-2f(s))R(x_0)\rho^*g & \ \ {\rm on} \ \
s^{-1}(-\infty,1] \\
{\rm exp}(-2f(s))\left(\alpha(s)R(x_0)\rho^*g+(1-\alpha(s))\eta
g_0\right)  & \ \ {\rm on}\ \  s^{-1}([1,2]) \\
{\rm exp}(-2f(s))\eta g_0 & \ \ {\rm on}\ \ s^{-1}([2,A_{r_0}]\\
\left[\beta(s){\rm exp}(-2f(s))+(1-\beta(s)){\rm exp}(-2f(4+A_0))\right]\eta
g_0 & \ \ {\rm on} \ \ s^{-1}([A_{r_0},A']),\end{cases}$$ where
$A_{r_0}=4+A_0-r_0$ and $A'=A_0+4$.
 Then we define
$$\widetilde g=R(x_0)^{-1}\hat g.$$
See {\sc Fig.}~\ref{fig:locsurg}.

\begin{thm}\label{LOCALSURGERY}
There are constants $C_0,q,R_0<\infty$ and $\delta'_0>0$ such that the
following hold for the result $({\mathcal S},\widetilde g)$ of surgery on
$(N,g)$ provided that $R(x_0)\ge R_0$, $0<\delta\le \delta'_0$. Define $f(s)$
as above with the constants $C_0,\delta$ and then use $f$ to define surgery on
a $\delta$-neck $N$ to produce $({\mathcal S},\tilde g)$. Then the following
hold.
\begin{itemize}
\item Fix $t\ge 0$. For any $p\in N$,
let $X(p)={\rm max}(0,-\nu_{g}(p))$, where $\nu_ g(p)$ is the
smallest eigenvalue of $Rm_{ g}(p)$. Suppose that for all $p\in N$
we have:
\begin{enumerate}\item[(1)] $R(p)\ge \frac{-6}{1+4t}$, and
\item[(2)] $R(p)\geq 2X(p)\left({\rm log}X(p)+{\rm
log}(1+t)-3\right)$, whenever $0<X(p)$.
\end{enumerate}
Then the curvature of $({\mathcal S},\widetilde g)$ satisfies the
same equations at every point of ${\mathcal S}$ with the same value
of $t$.
\item The restriction of the metric $\widetilde g$ to $s^{-1}([1,4+A_0])$ has
positive sectional curvature.
\item Let $\xi\colon N\to {\mathcal S}$
be the map given in Claim~\ref{distdecrease}. Then it is a distance
decreasing map from $g$ to $\tilde g$.
\item For any $\delta''>0$ there is $\delta'_1=\delta'_1(\delta'')>0$  such that if
$\delta\le {\rm min}(\delta'_1,\delta'_0)$, then  the restriction of
$\hat g$ to $B_{\hat g}(p_0,(\delta'')^{-1})$ in $({\mathcal S},\hat
g)$ is $\delta''$-close in the $C^{[1/\delta'']}$-topology to the
restriction of the standard initial metric $g_0$ to
$B_{g_0}(p_0,(\delta'')^{-1})$.
\end{itemize}
\end{thm}

The rest of this chapter is devoted to the proof of this theorem.

Before starting the curvature computations let us make a remark about the surgery cap.

\begin{defn}\label{surgerycap}
The image in ${\mathcal S}$ of $B_{g_0}(p_0,0,A_0+4)$ is called the {\em
surgery cap}\index{surgery cap|ii}.
\end{defn}

The following is immediate from the definitions provided that
$\delta>0$ is sufficiently small.

\begin{lem}\label{capball}
The surgery cap in $({\mathcal S},\tilde g)$ has a metric that
differs from the one coming from a rescaled version of the standard
solution. Thus, the image of this cap is not necessarily a metric
ball. Nevertheless for $\epsilon<1/200$ the image of this cap will
be contained in the metric ball in ${\mathcal S}$ centered at $p_0$
of radius $R(x_0)^{-1/2}(A_0+5)$ and will contain the metric ball
centered at $p_0$ of radius $R(x_0)^{-1/2}(A_0+3)$. Notice also that
the complement of the closure of the surgery cap in ${\mathcal S}$
is isometrically identified with $N^-$.
\end{lem}

\section{Preliminary computations}

We shall compute  in a slightly more general setup.
 Let $I$ be an open interval contained in $(-\delta^{-1},4+A_0)$ and let
$h$ be a metric  on $S^2\times I$ within $\delta$ in the
$C^{[1/\delta]}$-topology of the restriction to this open
submanifold of the standard metric $h_0\times ds^2$.  We let $\hat
h=e^{-2f}h$. Fix local coordinates near a point $y\in S^2\times I$.
We denote by $\nabla$ the covariant derivative for $h$ and by
$\widehat\nabla$ the covariant derivative for $\hat h$.
 We also denote by
$(R_{ijkl})$ the matrix of the Riemann curvature operator of $h$
in the associated basis of $\wedge^2T(S^2\times I)$ and by $(\hat
R_{ijkl})$ the matrix of the Riemann curvature operator of $\hat
h$ with respect to the same basis. Recall the formula for the
curvature of a conformal change of metric (see, (3.34) on p.51 of
\cite{Sakai}):
\begin{eqnarray}\label{confchange}\hat R_{ijkl} & = &
e^{-2f}\left(R_{ijkl}-f_jf_kh_{il}+f_jf_lh_{ik}+f_if_kh_{jl}-f_if_lh_{jk}\right.
\\ & & \left.-
(\wedge^2h)_{ijkl}|\nabla
f|^2-f_{jk}h_{il}+f_{ik}h_{jl}+f_{jl}h_{ik}-f_{il}h_{jk}\right).\nonumber\end{eqnarray}
Here, $f_i$ means $\partial_if$,
$$f_{ij}={\rm Hess}_{ij}(f)=\partial_if_j-f_l\Gamma_{ij}^l,$$
and $\wedge^2h$ is the metric induced by $h$ on $\wedge^2TN$, so
that
$$\wedge^2h_{ijkl}=h_{ik}h_{jl}-h_{il}h_{jk}.$$

Now we introduce the notation $O(\delta)$. When we say that a
quantity is $O(\delta)$ we mean that there is some universal
constant $C$ such that, provided that $\delta>0$ is sufficiently
small, the absolute value of the quantity is $\le C\delta$. The
universal constant is allowed to change from inequality to
inequality.

In our case we take local coordinates adapted to the $\delta$-neck:
$(x^0,x^1,x^2)$ where $x^0$ agrees with the $s$-coordinate and
$(x^1,x^2)$ are Gaussian local coordinates on the $S^2$ such that
$dx^1$ and $dx^2$ are orthonormal at the point in question in the
round metric $h_0$. The function $f$ is a function only of $x^0$.
Hence $f_i=0$ for $i=1,2$. Also, $f_0=\frac{q}{s^2}f$. It follows
that
$$|\nabla f|_h=\frac{q}{s^2}f\cdot(1+O(\delta)),$$ so that
$$|\nabla
f|_h^2=\frac{q^2}{s^4}f^2\cdot (1+O(\delta)).$$

Because the metric $h$ is $\delta$-close to the product $h_0\times
ds^2$, we see that $h_{ij}(y)=(h_0)_{ij}(y)+O(\delta)$ and the
Christoffel symbols $\Gamma_{ij}^k(y)$ of $h$ are within $\delta$ in
the $C^{([1/\delta]-1)}$-topology of those of the product metric
$h_0\times ds^2$. In particular, $\Gamma_{ij}^0=O(\delta)$ for all
$ij$.
  The components $f_{ij}$ of the Hessian with
respect to $h$ are given by
$$f_{00}=\left(\frac{q^2}{s^4}-\frac{2q}{s^3}\right)f+\frac{q}{s^2}fO(\delta),$$
$$f_{i0}=\frac{q}{s^2}fO(\delta)\ \ \ {\rm for}\ \ 1\le i\le 2,$$
$$f_{ij}=\frac{q}{s^2}fO(\delta)\ \ \ {\rm for}\ \ 1\le i,j\le 2.$$
In the following $a,b,c,d$ are indices taking values $1$ and $2$.
Substituting in Equation~(\ref{confchange}) yields
\begin{eqnarray*}\hat
R_{0a0b} & = &
e^{-2f}\left(R_{0a0b}+\frac{q^2}{s^4}f^2h_{ab}-h_{ab}(\frac{q^2}{s^4})f^2(1+O(\delta))
+\left(\frac{q^2}{s^4}-\frac{2q}{s^3}\right)fh_{ab} \right. \\
& & \left. +\frac{q}{s^2}fO(\delta)\right)\\
& = &
e^{-2f}\left(R_{0a0b}+\left(\frac{q^2}{s^4}-\frac{2q}{s^3}\right)fh_{ab}+
h_{ab}(\frac{q^2}{s^4})f^2O(\delta)
+\frac{q}{s^2}fO(\delta)\right)\end{eqnarray*} Also, we have
\begin{eqnarray*}\hat
R_{ab0c} & = &
e^{-2f}\left(R_{ab0c}-(\wedge^2h)_{ab0c}(\frac{q^2}{s^4})f^2(1+O(\delta))
+\frac{q}{s^2}fO(\delta)\right) \\
& = & e^{-2f}\left(R_{ab0c}+(\frac{q^2}{s^4})f^2O(\delta)
+\frac{q}{s^2}fO(\delta)\right).\end{eqnarray*}
 Lastly,
\begin{eqnarray*} \hat R_{1212} & = &
e^{-2f}\left(R_{1212}-(\wedge^2h)_{1212}\frac{q^2}{s^4}f^2(1+O(\delta))+\frac{q}{s^2}fO(\delta)\right)
\\
& = &
e^{-2f}\left(R_{1212}-\frac{q^2}{s^4}f^2(1+O(\delta))+\frac{q}{s^2}fO(\delta)\right).
\end{eqnarray*}

Now we are ready to fix the constant $q$. We fix it so that for all $s\in
[0,4+A_0]$ we have \begin{equation}\label{qeqn}
 q \gg (4+A_0)^2\ \ \ {\rm and}\ \ \ \frac{q^2}{s^4}e^{-q/s} \ll 1.
 \end{equation}
It follows immediately that $q^2/s^4\gg q/s^3$ for all $s\in [0,4+A_0]$.
 We are not yet ready
to fix the constant $C_0$, but once we do we shall always require $\delta$ to
satisfy
$\delta\ll C_0^{-1}$ so that for all $s\in [0,4+A_0]$ we have
$$\frac{q}{s^2}f^2\ll \frac{q^2}{s^4}f^2\ll \frac{q}{s^2}f\ll 1.$$
(These requirements are not circular, since $C_0$ and $q$ are chosen
independent of $\delta$.)

  Using
these inequalities and putting our computations in matrix form show
the following.

\begin{cor}\label{delta0}
There is $\delta'_2>0$, depending on $C_0$ and $q$,
such that if $\delta\le \delta'_2$ then we have
\begin{equation}\label{maineqn}
\left(\hat
R_{ijkl}\right)=e^{-2f}\left[\left(R_{ijkl}\right)+\begin{pmatrix}
-\frac{q^2}{s^4}f^2 & 0  \\ 0 &
\left(\frac{q^2}{s^4}-\frac{2q}{s^3}\right)f\begin{pmatrix} 1 & 0
\\ 0 & 1\end{pmatrix}
\end{pmatrix}+\left(\frac{q}{s^2}fO(\delta)\right)\right].\end{equation}
\end{cor}

Similarly, we have the equation relating scalar curvatures
$$\hat R=e^{2f}\left(R+4\triangle f-2|\nabla f|^2\right),$$
and hence
$$\hat R=e^{2f}\left(R+4\left(\frac{q^2}{s^4}-\frac{2q}{s^3}\right)f-2\frac{q^2}{s^4}f^2+
\frac{q}{s^2}fO(\delta)\right).$$

\begin{cor}\label{Rcompar}
For any constant $C_0<\infty$ and any $\delta<{\rm min}(\delta_2',C_0^{-1})$ we
have $\hat R\ge R$.
\end{cor}

\begin{proof}
By our choice of $q$, since $C_0\delta<1$, then $f^2\ll f$ and $q^2/s^4\gg {\rm
max}(q/s^3,q/s^2)$ so that the result follows immediately from the above
formula.
\end{proof}

Now let us compute the eigenvalues of the curvature $ R_{ijkl}(y)$
for any $y\in S^2\times I$.

\begin{lem}\label{newbasis}
There is a $\delta_3'>0$ such that the following hold if $\delta\le
\delta_3'$. Let $\{e_0,e_1,e_2\}$ be an orthonormal basis for the
tangent space at a point $y\in S^2\times I$ for the metric
$h_0\times ds^2$ with the property that $e_0$ points in the
$I$-direction. Then there is a basis $\{f_0,f_1,f_2\}$ for this
tangent space so that the following hold:
\begin{enumerate}
\item[(1)] The basis is orthonormal in the metric $h$.
\item[(2)] The change of basis matrix expressing the $\{f_0,f_1,f_2\}$
in terms of $\{e_0,e_1,e_2\}$ is of the form ${\rm Id}+O(\delta)$.
\item[(3)] The Riemann curvature of $h$ in the basis $\{f_0\wedge f_1,f_1\wedge f_2,f_2\wedge f_0\}$
of $\wedge^2T_y(S^2\times I)$ is
$$\begin{pmatrix} 1/2 & 0 & 0 \\ 0 & 0 & 0 \\ 0 & 0 & 0
\end{pmatrix}+O(\delta).$$
\end{enumerate}
\end{lem}

\begin{proof}
Since $h$ is within $\delta$ of $h_0\times ds^2$ in the
$C^{[1/\delta]}$-topology, it follows that the matrix for $h(y)$ in
$\{e_0,e_1,e_2\}$ is within $O(\delta)$ of the identity matrix, and the matrix
for the curvature of $h$ in the associated basis of $\wedge^2T_y(S^2\times I)$
is within $O(\delta)$ of the curvature matrix for $h_0\times ds^2$, the latter
being the diagonal matrix with diagonal entries $\{1/2,0,0\}$. Thus, the usual
Gram-Schmidt orthonormalization process constructs the basis $\{f_0,f_1,f_2\}$
satisfying the first two items. Let $A=(A^{ab})$ be the change of basis matrix
expressing the $\{f_a\}$ in terms of the $\{e_b\}$, so that $A={\rm
Id}+O(\delta)$. The curvature of $h$ in this basis is then given by
$B^{\rm tr}(R_{ijkl})B$ where $B=\wedge^2A$ is the induced change of basis matrix
expressing the basis $\{f_0\wedge f_1,f_1\wedge f_2,f_2\wedge f_0\}$ in terms
of $\{e_0\wedge e_1,e_1\wedge e_2,e_2\wedge e_0\}$. Hence, in the basis
$\{f_0\wedge f_1,f_1\wedge f_2,,f_2\wedge f_0\}$ the curvature matrix for $h$
is within $O(\delta)$ of the same diagonal matrix. For $\delta$ sufficiently
small then the eigenvalues of the curvature matrix for $h$ are within
$O(\delta)$ of $(1/2,0,0)$.
\end{proof}

\begin{cor}\label{2ndnewbasis}
The following holds provided that $\delta\le \delta_3'$. It is
possible to choose the basis $\{f_0,f_1,f_2\}$ satisfying the
conclusions of Lemma~\ref{newbasis} so that in addition the
curvature matrix for $(R_{ijkl}(y))$ is of the form
$$\begin{pmatrix} \lambda & 0 & 0 \\ 0 & \alpha & \beta \\ 0 & \beta & \gamma\end{pmatrix}$$
with $|\lambda-\frac{1}{2}|\le O(\delta)$ and
$|\alpha|,|\beta|,|\gamma|\le O(\delta)$.
\end{cor}

\begin{proof}
We have an $h$-orthonormal basis $\{f_0\wedge f_1,f_1\wedge
f_2,f_2\wedge f_0\}$ for $\wedge^2T_y(S^2\times \Ar)$ in which the
quadratic form $(R_{ijkl}(y)$ is
$$\begin{pmatrix} 1/2 & 0 & 0 \\ 0 & 0 & 0 \\ 0 & 0 & 0\end{pmatrix}
+O(\delta).$$ It follows that the restriction to the $h$-unit sphere
in $\wedge^2T_y(S^2\times \Ar)$ of this quadratic form achieves its
maximum value at some vector $v$, which, when written out in this
basis, is given by $(x,y,z)$ with $|y|,|z|\le O(\delta)$ and
$|x-1|\le O(\delta)$. Of course, this maximum value is within
$O(\delta)$ of $1/2$. Clearly, on the $h$-orthogonal subspace to
$v$, the quadratic form is given by a matrix all of whose entries
are $O(\delta)$ in absolute value. This gives us a new basis of
$\wedge^2T_y(S^2\times I)$ within $O(\delta)$ of the given basis in
which $(R_{ijkl}(y))$ is diagonal. The corresponding basis for
$T_y(S^2\times\Ar)$ is as required.
\end{proof}

Now we consider the expression  $(\hat R_{ijkl}(y))$ in this basis.

\begin{lem}\label{Rhatmatrix} Set $\delta'_4={\rm min}(\delta_2',\delta_3')$.  Suppose
that $\delta\le {\rm min}(\delta'_4,C_0^{-1})$. Then in the basis
$\{f_0,f_1,f_2\}$ for $T_y(S^2\times I)$ as in Corollary~\ref{2ndnewbasis} we
have
$$(\hat R_{ijkl}(y)) =e^{-2f}\left[\begin{pmatrix} \lambda & 0 & 0
\\ 0 & \alpha & \beta \\ 0 & \beta  & \gamma\end{pmatrix}
+\begin{pmatrix} -\frac{q^2}{s^4}f^2 & 0  \\ 0 &
\left(\frac{q^2}{s^4}-\frac{2q}{s^3}\right)f\begin{pmatrix} 1 & 0
\\ 0 & 1\end{pmatrix}
\end{pmatrix}+\left(\frac{q^2}{s^4}fO(\delta)\right)\right]$$
where $\lambda,\alpha,\beta,\gamma$ are the constants in
Lemma~\ref{2ndnewbasis} and the first matrix is the expression for
$(R_{ijkl}(y))$ in this basis.
\end{lem}

\begin{proof}
We simply conjugate the expression in Equation~(\ref{maineqn}) by the change of
basis matrix and use the fact that by our choice of $q$ and the fact that
$C_0\delta<1$, we have $f\gg  f^2$ and $q/s^3\ll q^2/s^4$.
\end{proof}

\begin{cor}\label{13.10}
Assuming that $\delta\le {\rm min}(\delta'_4,C_0^{-1})$, there is an
$h$-orthonormal basis $\{f_0,f_1,f_2\}$ so that in the associated basis for
$\wedge^2T_y(S^2\times I)$ the matrix $(R_{ijkl}(y))$ is diagonal and given by
$$\begin{pmatrix} \lambda & 0 & 0 \\ 0 & \mu & 0 \\ 0 & 0 &
\nu\end{pmatrix}$$ with $|\lambda-1/2|\le O(\delta)$ and
$|\mu|,|\nu|\le O(\delta)$. Furthermore, in this same basis the
matrix  $(\hat R_{ijkl}(y))$ is
$$e^{-2f}\left[\begin{pmatrix}\lambda & 0 & 0 \\ 0 & \mu & 0 \\ 0 & 0 &
\nu\end{pmatrix}+\begin{pmatrix} -\frac{q^2}{s^4}f^2 & 0  \\ 0 &
\left(\frac{q^2}{s^4}-\frac{2q}{s^3}\right)f\begin{pmatrix} 1 & 0 \\
0 & 1\end{pmatrix}\end{pmatrix}+\frac{q^2}{s^4}fO(\delta)\right].$$
\end{cor}

\begin{proof}
To diagonalize $(R_{ijkl}(y))$ we need only rotate in the
$\{f_1\wedge f_2,f_2\wedge f_3\}$-plane. Applying this rotation to
the expression in Lemma~\ref{newbasis} gives the result.
\end{proof}

\begin{cor}\label{improve}
There is a constant $A<\infty$ such that the following holds for the given value of $q$
and any $C_0$ provided that $\delta$ is sufficiently small.
Suppose that the eigenvalues for the curvature matrix of $h$ at $y$
are $\lambda\ge \mu\ge \nu$.
Then the eigenvalues for the curvature of $\hat h$ at the point $y$
are given by $\lambda',\mu',\nu'$, where
$$\left|\lambda'-e^{2f}\left(\lambda-\frac{q^2}{s^4}f^2\right)\right|
\le \frac{q^2}{s^4}fA\delta$$
$$\left|\mu'-e^{2f}\left(\mu+\left(\frac{q^2}{s^4}-\frac{2q}{s^3}\right)f\right)\right|\le
\frac{q^2}{s^4}fA\delta$$
$$\left|\nu'-e^{2f}\left(\nu+\left(\frac{q^2}{s^4}-\frac{2q}{s^3}\right)f\right)\right|\le
\frac{q^2}{s^4}fA\delta.$$

In particular, we have
$$\nu'\ge e^{2f}\left(\nu+\frac{q^2}{2s^4}f\right)$$
$$\mu'\ge e^{2f}\left(\mu+\frac{q^2}{2s^4}f\right).$$
\end{cor}

\begin{proof}
Let $\{f_0,f_1,f_2\}$ be the $h$-orthonormal basis given in
Corollary~\ref{13.10}. Then $\{e^ff_0,e^ff_1,e^ff_2\}$ is orthonormal for $\hat
h=e^{-2f}h$. This change multiplies the curvature matrix by $e^{4f}$. Since
$f\ll 1$, $e^{4f}<2$ so that the expression for $(\hat R_{ijkl}(y))$ in this
basis is exactly the same as in Lemma~\ref{Rhatmatrix} except that the factor
in front is $e^{2f}$ instead of $e^{-2f}$. Now, it is easy to see that since
$((q^2/s^4)fA\delta)^2\ll (q^2/s^4)fA\delta$ the eigenvalues will differ from
the diagonal entries by at most a constant multiple of $(q^2/s^4)fA\delta$.

 The first three inequalities are  immediate
from the previous corollary. The last two follow since $q^2/s^4\gg q/s^3$ and
$\delta\ll 1$.
\end{proof}

One important consequence of this computation is the following:

\begin{cor}\label{smeigen}
For the given value of $q$ and for any $C_0$,
assuming that $\delta>0$ is sufficiently small, the smallest
eigenvalue of ${\rm Rm}_{\hat h}$ is greater than the smallest eigenvalue
of ${\rm Rm}_h$ at the same point. Consequently, at any point where $h$
has non-negative curvature  so does $\hat h$.
\end{cor}

\begin{proof}
Since $|\lambda-1/2|, |\mu|, |\nu|$ are all $O(\delta)$ and since
$\frac{q^2}{s^4}f\ll 1$, it follows that the smallest eigenvalue of $(\hat
R_{ijkl}(y))$ is either $\mu'$ or $\nu'$. But it is immediate from the above
expressions  that $\mu'>\mu$ and $\nu'>\nu$. This completes the proof.
\end{proof}

Now we are ready to fix $C_0$.
There is a universal constant $K$ such that for all $\delta>0$ sufficiently small
and for any $\delta$-neck $(N,h)$ of
scale one, every eigenvalue of ${\rm Rm}_h$ is at least $-K\delta$.
We set
$$C_0=2Ke^{q}.$$

\begin{lem}\label{nuprime>0}
With these choices of $q$ and $C_0$ for any $\delta>0$ sufficiently small we have
$\nu'>0$ and $\mu'>0$ for $s\in [1,4+A_0]$ and $\lambda'>1/4$.
\end{lem}

\begin{proof}
Then by the previous result we have
$$\nu'\ge e^{2f}\left(\nu+\frac{q^2}{2s^4}f\right).$$
It is easy to see that since $q\gg (4+A_0)$ the function $(q^2/2s^4)f$ is an
increasing function on $[1,4+A_0]$. Its value at $s=1$ is
$(q^2/2)e^{-q}C_0\delta>K\delta$. Hence $\nu+\frac{q^2}{2s^4}f>0$ for all $s\in
[1,4+A_0]$ and consequently $\nu'>0$ on this submanifold. The same argument
shows $\mu'>0$. Since $q^2/s^4 f^2\ll 1$ and $0<f$, the statement about
$\lambda'$ is immediate.
\end{proof}

\section{The proof of Theorem~\protect{\ref{LOCALSURGERY}}}

\subsection{Proof of the first two items for $s< 4$}

We consider the metric in the region $s^{-1}(-\delta^{-1},4)$ given by
$$h=\alpha(s)R_g(x_0)\rho^*g+(1-\alpha(s))\eta g_0.$$
There is a constant $K'<\infty$ (depending on the
$C^{[1/\delta]}$-norm of  $\alpha$) such that $h$ is within
$K'\delta$ of the product metric in the
$C^{[1/(K'\delta)]}$-topology. Thus, if $\delta$ is sufficiently
small, all of the preceding computations hold with the error term
$(q^2/s^4)fAK'\delta$. Thus, provided that $\delta$ is sufficiently
small,  the conclusions about the eigenvalues hold for $e^{-2f}h$ in
the region $s^{-1}(-\delta^{-1},4)$. But $e^{-2f}h$ is exactly equal
to $R(x_0)\tilde g$ in this region. Rescaling, we conclude that on
$s^{-1}(-\delta^{-1},4)$ the smallest eigenvalue of $\tilde g$ is
greater than the smallest eigenvalue of $g$ at the corresponding
point and that $R_{\tilde g}\ge R_{g}$ in this same region.

The first conclusion of Theorem~\ref{LOCALSURGERY}  follows by
applying the above considerations to the case of
$h=R_g(x_0)\rho^*g$. Namely, we have:

\begin{prop}
Fix $\delta>0$ sufficiently small. Suppose that for some $t\ge 0$
and every point $p\in N$ the curvature of $h$ satisfies:
\begin{enumerate}
\item[(1)] $R(p)\ge \frac{-6}{1+4t}$, and
\item[(2)] $R(p)\geq 2X(p)\left({\rm log}X(p)+{\rm
log}(1+t)-3\right)$ whenever $0<X(x,t)$.
\end{enumerate}

Then  the curvature $({\mathcal S},\widetilde g)$ satisfies the same
equation with the same value of $t$ in the region
$s^{-1}(-\delta^{-1},4)$. Also, the curvature of $\tilde g$ is
positive in the region $s^{-1}[1,4)$.
\end{prop}

\begin{proof}
According to Corollary~\ref{smeigen}, the smallest eigenvalue of $\hat
h$ at any point $p$ is greater than or equal to the smallest
eigenvalue of $h$ at the corresponding point. According to
Corollary~\ref{Rcompar}, $\hat R(p)\ge R(p)$ for every $p\in
{\mathcal S}$. Hence, $X_{\hat h}(p)\le X_h(p)$. If $X_h(p)\ge
e^3/(1+t)$, then we have
$$\hat R(p)\ge R(p)\ge 2X_h(p)({\rm log}(X_h(p)+{\rm log}(1+t)-3)\ge 2X_{\hat h}(p)
({\rm log}(X_{\hat h}(p)+{\rm log}(1+t)-3).$$
If $X_h(p)<e^3(1+t)$, then $X_{\hat h}(p)<e^3/(1+t)$. Thus, in this case
since we are in a $\delta$-neck, provided that $\delta$ is sufficiently small, we have
 $R(p)\ge 0$ and hence
$$\hat R(p)\ge R(p)\ge 0>2X_{\hat p}(p)({\rm log}(X_{\hat h}(p)+{\rm log}(1+t)-3).$$

This completes the proof in both cases.

This establishes the first item in the conclusion of
Theorem~\ref{LOCALSURGERY} for $\delta>0$ sufficiently small on
$s^{-1}(-\delta^{-1},4)$. As we have seen in Lemma~\ref{nuprime>0},
the curvature is positive on $s^{-1}[1,4)$.
\end{proof}

\subsection{Proof of the first two items for $s\ge 4$}

Now let us show that the curvature on $\tilde g$ is positive in the region
$s^{-1}([4,4+A_0])$. First of all in the preimage of the
interval $[4,4+A_0,-r_0]$ this follows
from Corollary~\ref{smeigen} and the fact that $\eta g_0$ has non-negative
curvature. As for the region $s^{-1}([4+A_0-r_0,4+A_0])$, as $\delta$ tends to
zero, the metric here tends smoothly to the restriction of the metric $g_0$ to
that subset. The metric $g_0$ has positive curvature on
$s^{-1}([4+A_0-r_0,4+A_0])$. Thus, for all $\delta>0$ sufficiently small the
metric $\tilde g$ has positive curvature on all of $s^{-1}([4+A_0-r_0,4+A_0])$.
This completes the proof of the first two items.

\subsection{Proof of the third item}

 By
construction the restriction of the metric $\widetilde g$ to
$s^{-1}((-\delta^{-1},0])$ is equal to the metric $\rho^*g$. Hence,
in this region the mapping  is an isometry. In the region
$s^{-1}([0,4])$ we have $R(x_0)\rho^*g\ge \eta g_0$ so that by
construction in this region $\rho^*g\ge \tilde g$. Lastly, in the
region $s^{-1}([4,A_0+4])$ we have $R(x_0)^{-1}\eta g_0\ge \tilde
g$. On the other hand, it follows from Lemma~\ref{distdecrease} that
the map from $([0,\delta^{-1}]\times S^2,R(x_0)\rho^* g)$ to
$(B(p_0,4+A_0),\eta g)$ is distance decreasing. This completes the
proof of the third item.

\subsection{Completion of the proof}

As $\delta$ goes to zero,  $f$ tends to zero in the
$C^\infty$-topology and $\eta$ limits to $1$. From this  the fourth
item is clear.

This completes the proof of Theorem~\ref{LOCALSURGERY}.

\section{Other properties of the result of surgery}

\begin{lem}\label{stdballs}
Provided that $\delta>0$ is sufficiently small the following holds. Let $(N,g)$
be a $\delta$-neck and let $({\mathcal S},\tilde g)$ be the result of surgery
along the cental $2$-sphere of this neck. Then for any $0<D<\infty$ the ball
$B_{\tilde g}(p,D+5+A_0)\subset {\mathcal S}$ has boundary contained in
$s_N^{-1}(-(2D+2),-D/2)$.
\end{lem}

\begin{proof}
The Riemannian manifold $({\mathcal S},\tilde g)$ is identified by a
diffeomorphism with the union of $s_{N}^{-1}(-\delta^{-1},0]$ to
$B_{g_0}(p_0,A_0+4)$ glued along their boundaries. Thus, we have a natural
identification of ${\mathcal S}$ with the ball $B_{g_0}(p,A_0+4+\delta^{-1})$
in the standard solution. This identification pulls back the metric $\tilde g$
to be within $2\delta$ of the standard initial metric. The result then follows
immediately for $\delta$ sufficiently small.
\end{proof}

\chapter{Ricci Flow with surgery: the definition}\label{sect:surgery}

In this chapter we introduce Ricci flows with surgery. These objects
are closely related to generalized Ricci flows but they differ
slightly. The space-time of a Ricci flow with surgery has an open
dense subset that is a manifold, and the restriction of the Ricci
flow with surgery to this open subset is a generalized Ricci flow.
Still there are other, more singular points allowed in a Ricci flow
with surgery.

\section{Surgery space-time}

\begin{defn}
By a {\em space-time}\index{space-time|ii} we mean a paracompact Hausdorff space
${\mathcal M}$ with a continuous function ${\bf t}\colon {\mathcal
M}\to \Ar$, called {\em time}. We require that the image of ${\bf
t}$ be an interval $I$, finite or infinite with or without
endpoints, in $\Ar$. The interval $I$ is called the {\em time-interval of
definition of space-time}\index{time-interval of definition}. The initial point of $I$, if there is
one,  is the {\em initial time}\index{initial time|ii} and the final point of $I$, if there
is one, is the {\em final time}\index{final time|ii}. The level sets of ${\bf t}$ are
called the {\em time-slices}\index{time-slices|ii} of space-time, and the preimage of the
initial (resp., final) point of $I$ is the {\em initial} (resp.,
{\em final}) time-slice.
\end{defn}

We are interested in a certain class of space-times, which we call
{\em surgery space-times}. These objects have a `smooth structure'
(even though they are not smooth manifolds). As in the case of a
smooth manifold, this smooth structure is given by local coordinate
charts with appropriate overlap functions.

\subsection{An exotic chart}

There is one exotic chart, and we begin with its description. To
define this chart we consider the open unit square $(-1,1)\times
(-1,1)$. We shall define a new topology, denoted by ${\mathcal P}$,
on this square. The open subsets of ${\mathcal P}$ are the open
subsets of the usual topology on the open square together with open
subsets of $(0,1)\times [0,1)$. Of course, with this topology the
`identity' map $\iota\colon {\mathcal P}\to (-1,1)\times (-1,1)$  is
a continuous map. Notice that the restriction of the topology of
${\mathcal P}$ to the complement of the closed subset $[0,1)\times
\{0\}$ is a homeomorphism onto the corresponding subset of the open
unit square. Notice that the complement of $(0,0)$ in ${\mathcal P}$
is a manifold with boundary, the boundary being $(0,1)\times\{0\}$.
(See {\sc Fig.}~\ref{fig:spacetime} in the Introduction.)

Next, we define a `smooth structure' on ${\mathcal P}$ by defining a
sheaf of germs of `smooth' functions. The restriction of this sheaf
of germs of `smooth functions'  to the complement of $(0,1)\times
\{0\}$ in ${\mathcal P}$ is the usual sheaf of germs of smooth
functions on the corresponding subset of the open unit square.
In particular, a function is smooth near $(0,0)$ if and only if its restriction
to some neighborhood of $(0,0)$ is the pullback under $\iota$
of a usual smooth function on a
neighborhood of the origin in the square.
 Now
let us consider the situation near a point of the form $x=(a,0)$ for
some $0<a<1$. This point has arbitrarily small neighborhoods $V_n$
that are identified under $\iota$ with open subsets of $(0,1)\times
[0,1)$. We say that a function  $f$ defined in a neighborhood of $x$
in ${\mathcal P}$ is {\em smooth at $x$} if its restriction to one
of these neighborhoods $V_n$ is the pullback via $\iota|_{V_n}$ of a
smooth function in the usual sense on the open subset $\iota(V_n)$
of the upper half space. One checks directly that this defines a
sheaf of germs of `smooth' functions on ${\mathcal P}$. Notice that
the restriction of this sheaf to the complement of $(0,0)$ is the
structure sheaf of smooth functions of a smooth manifold with
boundary. Notice that the map $\iota\colon {\mathcal P}\to
(-1,1)\times (-1,1)$ is a smooth map in the sense that it pulls back
smooth functions on open subsets of the open unit square to smooth
functions on the corresponding open subset of ${\mathcal P}$.

Once we have the notion of smooth functions on ${\mathcal P}$, there
is the categorical notion of a diffeomorphism between open subsets
of ${\mathcal P}$: namely a homeomorphism with the property that it
and its inverse pull back smooth functions to smooth functions. Away
from the origin, this simply means that the map is a diffeomorphism
in the usual sense between manifolds with boundary, and in a
neighborhood of $(0,0)$ it factors through a diffeomorphism of
neighborhoods of the origin  in the square. While $\iota\colon
{\mathcal P}\to (-1,1)\times (-1,1)$ is a smooth map, it is not a
diffeomorphism.

We define the {\sl tangent bundle} of ${\mathcal P}$ in the usual
manner. The tangent space at a point is the vector space of
derivations of the germs of smooth functions at that point. Clearly,
away from $(0,0)$ this is the usual ($2$-plane) tangent bundle of
the smooth manifold with boundary. The germs of smooth functions at
$(0,0)$ are, by definition,  the pullbacks under $\iota$ of germs
of smooth functions at the origin for the unit square, so that the
tangent space of ${\mathcal P}$ at $(0,0)$ is identified  with the
tangent space of the open unit square at the origin. In fact, the
map $\iota$ induces an isomorphism from the tangent bundle of
${\mathcal P}$ to the pullback under $\iota$ of the tangent bundle
of the square. In particular, the tangent bundle of ${\mathcal P}$
has a given trivialization from the partial derivatives $\partial_x$
and $\partial_y$ in the coordinate directions on the square. We use
this trivialization to induce a smooth structure on the tangent
bundle of ${\mathcal P}$: that is to say, a section of $T{\mathcal
P}$ is smooth if and only if it can be written as $\alpha
\partial_x+\beta
\partial_y$ with $\alpha$ and $\beta $ being smooth functions on
${\mathcal P}$. The smooth structure agrees off of $(0,0)\in
{\mathcal P}$ with the usual smooth structure on the tangent bundle
of the smooth manifold with boundary.
By a {\em smooth vector field on ${\mathcal P}$}  we mean a smooth
section of the tangent bundle of ${\mathcal P}$.
 Smooth vector fields act as derivations on the smooth
functions on ${\mathcal P}$.

 We let ${\bf t}_{\mathcal P}\colon {\mathcal P}\to \Ar$ be the pullback via
$\iota$ of the usual projection to the second factor on the unit
square. We denote by $\chi_{\mathcal P}$ the smooth vector field
$\iota^*\partial_2$. Clearly, $\chi_{\mathcal P}({\bf t}_{\mathcal
P}) =1$. Smooth vector fields on ${\mathcal P}$ can be uniquely
integrated locally to smooth integral curves in ${\mathcal P}$. (At
a manifold with boundary point, of course only vector fields
pointing into the manifold can be locally integrated.)

\subsection{Coordinate charts for a surgery space-time}

Now we are ready to introduce the types of coordinate charts\index{coordinate charts for space-time} that we
shall use in our definition of a surgery space-time. Each coordinate
patch comes equipped with a smooth structure (a sheaf of germs of
smooth functions) and a tangent bundle with a smooth structure, so
that smooth vector fields act as derivations on the algebra of
smooth functions. There is also a distinguished smooth function,
denoted ${\bf t}$, and a smooth vector field, denoted $\chi$,
required to satisfy $\chi({\bf t})=1$. There are three types of
coordinates:
\begin{enumerate}
\item[(1)] The coordinate patch is an open subset of the strip ${\Ar}^{n}\times I$,
where $I$ is an interval, with its usual smooth structure and
tangent bundle; the function ${\bf t}$ is the projection onto $I$;
and the vector field $\chi$ is the unit tangent vector in the
positive direction tangent to the foliation with leaves $\{x\}\times
I$. The initial point of $I$, if there is one, is the initial time
of the space-time and the final point of $I$, if there is one, is
the final time of the space-time.
\item[(2)] The coordinate patch an open subset of $\Ar^n\times [a,\infty)$, for some $a\in \Ar$,
 with its
usual smooth structure as a manifold with boundary and its usual
smooth tangent bundle; the function ${\bf t}$ is the projection onto
the second factor; and the vector field is the coordinate partial
derivative associated with the second factor. In this case we
require that $a$ not be the initial time of the Ricci flow.
\item[(3)] The coordinate patch is a product of ${\mathcal P}$ with an open subset of
$\Ar^{n-1}$ with the smooth structure (i.e., smooth functions and
the smooth tangent bundle) being the product of the smooth structure
defined above on ${\mathcal P}$ with the usual smooth structure of
an open subset of $\Ar^{n-1}$; the function ${\bf t}$ is, up to an
additive constant, the pullback of the function ${\bf t}_{\mathcal
P}$ given above on ${\mathcal P}$; and the vector field $\chi$ is
the image of the vector field $\chi_{\mathcal P}$ on ${\mathcal P}$,
given above, under the product decomposition.
\end{enumerate}

An ordinary Ricci flow is covered by coordinate charts of the first
type.  The second and third are two extra types of coordinate charts
for a Ricci flow with surgery that are not allowed in generalized
Ricci flows. Charts of the second kind are smooth
manifold-with-boundary charts, where the boundary is contained in a
single time-slice,  not the initial time-slice, and the flow exists
for some positive amount of forward time from this manifold.

All the structure described above for ${\mathcal P}$ -- the smooth
structure, the tangent bundle with its smooth structure, smooth
vector fields acting as derivations on smooth functions -- exist for
charts of the third type. In addition, the unique local integrability
of smooth vector fields  hold for coordinate charts of the third
type. Analogous results for coordinate charts of the first two types
are clear.

Now let us describe the allowable overlap functions between charts.
Between charts of the first and second type these are the smooth
overlap functions in the usual sense that preserve the functions
${\bf t}$ and the vector fields $\chi$ on the patches. Notice that
because the boundary points in charts of the second type are
required to be at times other than the initial and final times, the
overlap of a chart of type one and a chart of type two is disjoint
from the boundary points of each. Charts of the first two types are
allowed to meet a chart of the third type only in its manifold and
manifold-with-boundary points. For overlaps between charts of the
first two types with a chart of the third type, the overlap
functions are  diffeomorphisms between open subsets preserving the
local time functions ${\bf t}$ and the local vector fields $\chi$.
Thus, all overlap functions are diffeomorphisms in the sense given
above.

\subsection{Definition and basic properties of surgery space-time}

\begin{defn}\label{spacetime}
A {\em surgery space-time}\index{space-time, surgery|ii} is a
space-time  ${\mathcal M}$ equipped with a maximal atlas  of charts
covering ${\mathcal M}$, each chart being of one of the three types
listed above, with the overlap functions being diffeomorphisms
preserving the functions ${\bf t}$ and the vector fields $\chi$. The
points with neighborhoods of the first type are called {\em smooth
points}, those with neighborhoods of the second type but not the
first type are called {\em exposed points}\index{exposed points|ii},
and all the other points are called {\em singular
points}\index{singular points, space-time|ii}. Notice that the union of
the set of smooth points and the set of exposed points forms a
smooth manifold with boundary (possibly disconnected). Each
component of the boundary of this manifold is contained in a single
time-slice. The union of those components contained in a time
distinct from the initial time and the final time is called the {\em
exposed region}\index{exposed region|ii}. and the boundary points of
the closure of the exposed region form the set of the singular
points of ${\mathcal M}$. (Technically, the exposed points are
singular, but we reserve this word for the most singular points.) An
$(n+1)$-dimensional surgery space-time is by definition of
homogeneous dimension $n+1$.

By construction, the local smooth functions ${\bf t}$ are compatible
on the overlaps and hence fit together to define a global smooth
function ${\bf t}\colon {\mathcal M}\to \Ar$\index{${\bf t}$}, called the {\em time}\index{time}
function. The level sets of this function are called the {\em
time-slices}\index{time-slices|ii} of the space-time, and ${\bf t}^{-1}(t)$ is denoted
$M_t$. Similarly, the tangent bundles of the various charts are
compatible under the overlap diffeomorphisms and hence glue together
to give a global smooth tangent bundle on space-time.  The smooth
sections of this vector bundle, the smooth vector fields on space
time, act as derivations on the smooth functions on space-time. The
tangent bundle of an $(n+1)$-dimensional surgery space-time is a
vector bundle of dimension $(n+1)$. Also, by construction the local
vector fields $\chi$ are compatible and hence glue together to
define a global vector field, denoted $\chi$\index{$\chi$|ii}. The vector field and
time function satisfy
$$\chi({\bf t})=1.$$
At the manifold points (including the exposed points) it is a usual
vector field. Along the exposed region and the initial time-slice
the vector field points into the manifold; along the final
time-slice it points out of the manifold.
\end{defn}

\begin{defn}
Let ${\mathcal M}$ be a surgery space-time. Given a space $K$ and an
interval $J\subset \Ar$  we say that an embedding $K\times J\to
{\mathcal M}$ is {\em compatible with  time  and the vector field}\index{compatible with time and the vector field|ii}
if: (i) the restriction of ${\bf t}$ to the image agrees with the
projection onto the second factor and (ii) for each $x\in X$ the
image of $\{x\}\times J$ is the integral curve for the vector field
$\chi$.
 If in addition $K$ is a subset of $M_{t}$
we require that $t\in J$ and that the map $K\times\{t\}\to M_{t}$ be
the identity. Clearly, by the uniqueness of integral curves for vector fields, two such
embeddings agree on their common interval of definition, so that,
given $K\subset M_{t}$ there is a maximal interval $J_K$ containing
$t$ such that such an embedding is defined on $K\times J_K$. In the
special case when $K=\{x\}$ for a point $x\in M_{t}$ we say that
such an embedding is {\em the maximal flow line} through $x$. The
embedding of the maximal interval through $x$ compatible with time
and the vector field $\chi$ is called {\em the domain of definition}
of the flow line through $x$. For a more general subset $K\subset
M_{t}$ there is an embedding $K\times J$ compatible with time and
the vector field $\chi$ if and only if, for every $x\in K$, the
interval $J$ is contained in the domain of definition of the flow
line through $x$.
\end{defn}

\begin{defn}
Let ${\mathcal M}$ be a surgery space-time with $I$ as its time
interval of definition. We say that $t\in I$ is a {\em regular}\index{time!regular|ii} time
if there is an interval $J\subset I$ which is an open neighborhood
in $I$ of $t$, and a diffeomorphism $M_{t}\times J\to {\bf
t}^{-1}(J)\subset{\mathcal M}$ compatible with time and the vector
field. A time is {\em singular}\index{time!singular|ii} if it is not regular. Notice that if
all times are regular, then space-time is a product $M_t\times I$
with ${\bf t}$ and $\chi$ coming from the second factor.
\end{defn}

\begin{lem}
Let ${\mathcal M}$ be an $(n+1)$-dimensional surgery space-time, and fix $t$.
The restriction of the smooth structure on ${\mathcal M}$ to the time-slice
$M_{t}$ induces the structure of a smooth $n$-manifold on this time-slice. That
is to say, we have a smooth embedding of $M_{t}\to {\mathcal M}$. This smooth
embedding identifies the tangent bundle of $M_{t}$ with a codimension-one
subbundle of the restriction  of tangent bundle of ${\mathcal M}$ to $M_{t}$.
This subbundle is complementary to the line field spanned by $\chi$. These
codimension-one subbundles along the various time-slices fit together to form
a smooth, codimension-one subbundle of the tangent bundle of space-time.
 \end{lem}

\begin{proof}
These statements are immediate for any coordinate patch, and hence
are true globally.
\end{proof}

\begin{defn}
We call the codimension-one subbundle of the tangent bundle of
${\mathcal M}$ described in the previous lemma  the {\em horizontal
subbundle}\index{horizontal subbundle|ii}, and we denote it
${\mathcal HT}({\mathcal M})$\index{${\mathcal HT}({\mathcal M})$|ii}.
\end{defn}

\section{The generalized Ricci flow equation}

In this section we introduce the Ricci flow equation for surgery
space-times, resulting in an object that we call Ricci flow with
surgery.

\subsection{Horizontal metrics}

\begin{defn}
By a {\em horizontal metric $G$ on a surgery space-time ${\mathcal
M}$}\index{horizontal metric|ii} we mean a $C^\infty$ metric on ${\mathcal HT}{\mathcal M}$. For
each $t$, the horizontal metric $G$ induces a Riemannian metric,
denoted $G(t)$, on the time-slice $M_t$. Associated to a horizontal
metric $G$ we have the {\em horizontal covariant derivative},
denoted $\nabla$. This is a pairing between horizontal vector fields
$$X\otimes Y\mapsto \nabla_XY.$$
On each time slice $M_t$ it is the usual Levi-Civita connection associated to
the Riemannian metric $G(t)$. Given a function $F$ on space-time, by its
gradient $\nabla F$ we mean its horizontal gradient. The value of this gradient
at a point $q\in M_t$ is the usual $G(t)$-gradient of $F|_{M_t}$. In
particular, $\nabla F$ is a smooth horizontal vector field on space-time. The
horizontal metric $G$ on space-time has its (horizontal) curvatures ${\rm Rm}_G$.
These are smooth symmetric endomorphisms of the second exterior power of
${\mathcal HT}{\mathcal M}$. The value of ${\rm Rm}_G$ at a point $q\in M_t$ is
simply the usual Riemann curvature operator of $G(t)$ at the point $q$.
Similarly, we have the (horizontal) Ricci curvature ${\rm Ric}={\rm Ric}_G$, a
section of the symmetric square of the horizontal cotangent bundle, and the
(horizontal) scalar curvature denoted $R=R_G$.
\end{defn}

The only reason for working in ${\mathcal HT}{\mathcal M}$ rather
than individually in each slice is to emphasize the fact that all
these horizontal quantities vary smoothly over the surgery
space-time.

Suppose that $t\in I$ is not the final time and suppose that
$U\subset M_{t}$ is an open subset with compact closure. Then there
is $\epsilon>0$ and an embedding $i_U\colon
U\times[t,t+\epsilon)\subset {\mathcal M}$ compatible with time and
the vector field. Of course, two such embeddings agree on their
common domain of definition. Notice also that for each $t'\in
[t,t+\epsilon)$ the restriction of the map $i_U$ to $U\times \{t'\}$
induces an diffeomorphism from $U$ to an open subset $U_{t'}$ of
$M_{t'}$. It follows that the local flow generated by the vector
field $\chi$ preserves the horizontal subbundle. Hence, the vector
field $\chi$ acts by Lie derivative on the sections of ${\mathcal
HT}({\mathcal M})$ and on all associated bundles (for example the
symmetric square of the dual bundle).

\subsection{The equation}

\begin{defn}
A {\em  Ricci flow with surgery}\index{Ricci flow!with surgery|ii} is a pair $({\mathcal M},G)$
consisting of a surgery space-time ${\mathcal M}$ and a horizontal
metric $G$ on ${\mathcal M}$ such that for every $x\in {\mathcal M}$
 we have
\begin{equation}\label{GRFEq}
{\mathcal L}_\chi(G)(x)=-2{\rm Ric}_G(x)) \end{equation} as sections
of the symmetric square of the dual to ${\mathcal HT}({\mathcal
M})$. If space-time is $(n+1)$-dimensional, then we say that the
Ricci flow with surgery is $n$-dimensional (meaning of course that
each time-slice is an $n$-dimensional manifold).
\end{defn}

\begin{rem}
 Notice that at an exposed point and at points at the initial and the final time
   the Lie derivative is a one-sided derivative.
\end{rem}

\subsection{Examples of Ricci flows with surgery}

\begin{exam}
One example of a Ricci flow with surgery is ${\mathcal M}=M_0\times
[0,T)$ with time function ${\bf t}$ and the vector field $\chi$
coming from the second factor. In this case the Lie derivative
${\mathcal L}_\chi$ agrees with the usual partial derivative in the
time direction, and hence our generalized Ricci flow equation is the
usual Ricci flow equation. This shows that an ordinary Ricci flow is indeed
a Ricci flow with surgery.
\end{exam}

The next lemma gives an example of a  Ricci flow with surgery where
the topology of the time-slices changes.

\begin{lem}\label{gluingsoln}
Suppose that we have manifolds $M_1\times (a,b]$ and $M_2\times
[b,c)$ and compact, smooth codimension-$0$ submanifolds
$\Omega_1\subset M_1$ and $\Omega_2\subset M_2$ with open
neighborhoods $U_1\subset M_1$ and $U_2\subset M_2$ respectively.
Suppose we have a diffeomorphism $\psi\colon U_1\to U_2$ carrying
$\Omega_1$ onto $\Omega_2$. Let $(M_1\times (a,b])_0$ be the subset
obtained by removing $(M_1\setminus \Omega_1)\times \{b\}$ from
$M_1\times (a,b]$. Form the topological space
$${\mathcal M}=(M_1\times (a,b])_0\cup M_2\times [b,c)$$
where  $\Omega_1\times\{b\}$ in $(M_1\times (a,b])_0$ is identified
with $\Omega_2\times \{b\}$ using the restriction of $\psi$ to
$\Omega_1$. Then ${\mathcal M}$ naturally inherits the structure of
a surgery space-time where the time function restricts to
$(M_1\times (a,b])_0$ and to $M_2\times [b,c)$ to be the projection
onto the second factor and the vector field $\chi$ agrees with the
vector fields coming from the second factor on each of $(M_1\times
(a,b])_0$ and $M_2\times [b,c)$.

 Lastly, given Ricci
flows $(M_1,g_1(t)),\ a<t\le b$, and $(M_2,g_2(t)),\ b\le t<c$, if
$\psi\colon (U_1,g_1(b))\to (U_2,g_2(b))$ is an isometry, then these
families fit together to form a smooth horizontal metric $G$ on
${\mathcal M}$ satisfying the Ricci flow equation, so that
$({\mathcal M},G)$ is a  Ricci flow with surgery.
\end{lem}

\begin{proof}
As the union of Hausdorff spaces along closed subsets, ${\mathcal
M}$ is a Hausdorff topological space. The time function is the one
induced from the projections onto the second factor. For any point
outside the $b$ time-slice there is the usual smooth coordinate
coming from the smooth manifold $M_1\times (a,b)$ (if $t<b$) or
$M_2\times (b,c)$ (if $t>b$). At any point of
$(M_2\setminus\Omega_2)\times \{b\}$ there is the smooth manifold
with boundary coordinate patch coming from $M_2\times[b,c)$. For any
point in ${\rm int}(\Omega_1)\times\{b\}$ we have the smooth
manifold structure obtained from gluing $({\rm int}(\Omega_1))\times
(a,b]$ to ${\rm int}(\Omega_2)\times [b,c)$ along the $b$ time-slice
by $\psi$. Thus, at all these points we have neighborhoods on which
our data determine a smooth manifold structure. Lastly, let us
consider a point $x\in
\partial \Omega_1\times \{b\}$.
Choose local coordinates $(x^1,\ldots,x^n)$ for a neighborhood $V_1$
of $x$ such that $\Omega_1\cap V_1=\{x^n\le 0\}$. We can assume that
$\psi$ is defined on all of $V_1$. Let $V_2=\psi(V_1)$ and take the
 local coordinates on $V_2$ induced from the $x^i$ on $V_1$.
 Were we to identify $V_1\times (a,b]$
with $V_2\times [b,c)$ along the $b$ time-slice using this map, then
this union would be a smooth manifold.  There is a
neighborhood of the point $(x,b)\in {\mathcal M}$ which is obtained
from the smooth manifold $V_1\times (a,b]\cup_\psi V_2\times [b,c)$
by inducing a new topology where the open subsets are, in addition
to the usual ones, any open subset of the form $\{x^n>0\}\times
[b,b')$ where $b<b'\le c$. This then gives the coordinate charts of
the third type near the points of $\partial \Omega_2\times \{b\}$.
Clearly, since the function ${\bf t}$ and the vector field
$\partial/\partial t$ are smooth on $V_1\times (a,b]\cup_\psi
V_2\times [b,c)$, we see that these objects glue together to form
smooth objects on ${\mathcal M}$.

Given the Ricci flows $g_1(t)$ and $g_2(t)$ as in the statement,
they clearly determine a (possibly singular) horizontal metric on
${\mathcal M}$. This horizontal metric is clearly smooth except
possibly along the $b$ time-slice. At any point of
$(M_2\setminus\Omega_2)\times\{b\}$ we have a one-sided smooth
family, which means that on this set the horizontal metric is
smooth. At a point of ${\rm int}(\Omega_2)\times \{b\}$, the fact
that the metrics fit together smoothly is an immediate consequence
of Proposition~\ref{patch}. At a point $x\in
\partial \Omega_2\times \{b\}$ we have neighborhoods $V_2\subset M_2$
of $x$ and $V_1\subset M_1$ of $\psi^{-1}(x)$ that are isometrically
identified by $\psi$. Hence, again by Lemma~\ref{patch} we see that
the Ricci flows fit together to form a smooth family of metrics on
$V_1\times (a,b]\cup_\psi V_2\times [b,c)$. Hence, the induced
horizontal metric on ${\mathcal M}$ is smooth near this point.
\end{proof}

The following is obvious from the definitions.

\begin{prop}
Suppose that $({\mathcal M},G)$ is a Ricci flow with surgery. Let
${\rm int}{\mathcal M}$ be the open subset consisting of all
smooth $(n+1)$-manifold points, plus all manifold-with-boundary
points at the initial time and the final time. This space-time
inherits the structure of a smooth manifold with boundary. This
structure together with the restrictions to it of ${\bf t}$ and the
vector field $\chi$ and the restriction of the horizontal metric $G$
form a generalized Ricci flow whose underlying smooth manifold is ${\rm int}{\mathcal M}$.
\end{prop}

\subsection{Scaling and translating Ricci flows with surgery}

Suppose that $({\mathcal M},G)$ is a Ricci flow with surgery. Let
$Q$ be a positive constant. Then we can define a new  Ricci flow
with surgery by setting $G'=QG$, ${\bf t'}=Q{\bf t}$ and
$\chi'=Q^{-1}\chi$. It is easy to see that the resulting data still
satisfies the generalized Ricci flow equation,
Equation~(\ref{GRFEq}). We denote this new  Ricci flow with surgery by
$(Q{\mathcal M},QG)$ where the changes in ${\bf t}$ and $\chi$ are
indicated by the factor $Q$ in front of the space-time.

It is also possible to translate a Ricci flow with surgery $({\mathcal
M},G)$ by replacing the time function ${\bf t}$ by ${\bf t'}={\bf
t}+a$ for any constant $a$, and leaving $\chi$ and $G$ unchanged.

\subsection{More basic definitions}

\begin{defn}
Let $({\mathcal M},G)$ be a Ricci flow with surgery, and let $x$ be
a point of space-time. Set $t={\bf t}(x)$. For any  $r>0$ we define
$B(x,t,r)\subset M_t$ to be the metric ball of radius $r$ centered
at $x$ in the Riemannian manifold $(M_t,G(t))$.
\end{defn}

\begin{defn}
Let $({\mathcal M},G)$ be a Ricci flow with surgery, and let $x$ be
a point of space-time. Set $t={\bf t}(x)$. For any $r>0$ and $\Delta
t>0$ we say that the {\em backward parabolic neighborhood}\index{parabolic neighborhood|ii}
$P(x,t,r,-\Delta t)$ {\em exists} in ${\mathcal M}$ if there is an
embedding $B(x,t,r)\times (t-\Delta t,t]\to {\mathcal M}$ compatible
with time and the vector field. Similarly, we say
that the {\em forward parabolic neighborhood} $P(x,t,r,\Delta t)$
{\em exists} in ${\mathcal M}$ if there is an embedding
$B(x,t,r)\times [t,t+\Delta t)\to {\mathcal M}$ compatible with time
and the vector field. A {\em parabolic neighborhood} is either a
forward or backward parabolic neighborhood.
\end{defn}

\begin{defn}\label{defnkappa}
Fix $\kappa>0$ and $r_0>0$. We say that a Ricci flow with surgery
$({\mathcal M},G)$ is {\em $\kappa$-noncollapsed on scales $\le
r_0$}\index{Ricci flow!$\kappa$-non-collapsed|ii} if the following holds for every point $x\in {\mathcal M}$ and
for every $r\le r_0$. Denote ${\bf t}(x)$ by $t$. If the  parabolic
neighborhood $P(x,t,r,-r^2)$ exists in ${\mathcal M}$ and if
$|{\rm Rm}_G|\le r^{-2}$ on $P(x,t,r,-r^2)$, then  ${\rm
Vol}\,B(x,t,r)\ge \kappa r^3$.
\end{defn}

\begin{rem}
For $\epsilon>0$ sufficiently small,
an $\epsilon$-round component satisfies the first condition in the above definition
 for some $\kappa>0$ depending only on
the order of the fundamental group of the underlying manifold, but there is no universal
$\kappa>0$ that works for all $\epsilon$-round manifolds. Fixing an integer $N$ let ${\mathcal C}_N$
be the class of closed $3$-manifolds
 with the property that any finite free factor of $\pi_1(M)$
has order at most $N$. Then any $\epsilon$-round component of any
time-slice of any Ricci flow $({\mathcal M},G)$ whose initial
conditions consist of a manifold in ${\mathcal C}_N$ will have
fundamental group of order at most $N$ and hence will satisfy the first condition
in the above definition for some $\kappa>0$ depending only on $N$.
\end{rem}

We also have the notion of the curvature being  pinched toward
positive, analogous to the notions for Ricci flows and  generalized Ricci flows.

\begin{defn}
Let $({\mathcal M},G)$ be a $3$-dimensional Ricci flow with surgery,
whose time domain of definition is contained in $[0,\infty)$. For
any $x\in {\mathcal M}$ we denote the eigenvalues of ${\rm Rm}(x)$ by
$\lambda(x)\ge \mu(x)\ge \nu(x)$ and we set $X(x)={\rm
max}(0,-\nu(x))$.
 We say
that its {\em curvature is pinched toward positive}\index{curvature!pinched toward positive|ii} if the following
hold for every $x\in {\mathcal M}$:
\begin{enumerate}
\item[(1)] $R(x)\ge \frac{-6}{1+4{\bf t}(x)}$.
\item[(2)] $R(x)\geq 2X(x)\left({\rm log}X(x)+{\rm
log}(1+{\bf t}(x))-3\right)$, whenever $0<X(x)$.
\end{enumerate}

Let $(M,g)$ be a Riemannian manifold and let $T\ge 0$.
We say that $(M,g)$ has curvature {\em pinched toward positive up to time $T$} if
the above two inequalities hold for all $x\in M$ with ${\bf t}(x)$ replaced by $T$.
\end{defn}

Lastly, there is the definition of canonical neighborhoods for a
Ricci flow with surgery, there is the following extension of the
notion for a generalized Ricci flow.

\begin{defn}
Fix constants $(C,\epsilon)$ and a constant $r$. We say that a Ricci
flow with surgery $({\mathcal M},G)$ {\em satisfies the strong
$(C,\epsilon)$-canonical neighborhood assumption with parameter $r$}\index{parameter, canonical neighborhood|ii}
if every point $x\in {\mathcal M}$ with $R(x)\ge r^{-2}$ has a strong
$(C,\epsilon)$-canonical neighborhood in ${\mathcal M}$. In all
cases except that of the strong $\epsilon$-neck, the strong
canonical neighborhood of $x$ is a subset of the time-slice
containing $x$, and the notion of a $(C,\epsilon)$-canonical
neighborhood has exactly the same meaning as in the case of an
ordinary Ricci flow. In the case of a strong $\epsilon$-neck
centered at $x$ this means that there is an embedding
$\left(S^2\times (-\epsilon^{-1},\epsilon^{-1})\right)\times ({\bf
t}(x)-R(x)^{-1},{\bf t}(x)]\to {\mathcal M}$, mapping $(q_0,0)$ to $x$,
 where $q_0$ is the basepoint of $S^2$, an embedding
compatible with time and the vector field,
such that the pullback
of $G$ is a Ricci flow on $S^2\times (-\epsilon^{-1},\epsilon^{-1})$
which, when the time is shifted by $-{\bf t}(x)$ and then the flow
is rescaled by $R(x)$, is within $\epsilon$ in the
$C^{[1/\epsilon]}$-topology of the standard evolving round cylinder
$\left(S^2\times(-\epsilon^{-1},\epsilon^{-1}),h_0(t)\times
ds^2\right),\ -1<t\le 0$, where the scalar curvature of the $h_0(t)$
is $1-t$.
\end{defn}

Notice that $x$ is an exposed point or sufficiently close to an exposed point
then $x$ cannot be the center of a strong $\epsilon$-neck.

\chapter{Controlled Ricci flows with surgery}

 We do not
wish to consider all Ricci flows with surgery. Rather we shall
concentrate on $3$-dimensional flows (that is to say $4$-dimensional
space-times) whose singularities are closely controlled both
topologically and geometrically. We introduce the hypotheses that we
require these evolutions to satisfy. Then main result, which is
stated in this chapter and proved in the next two, is that these
controlled $3$-dimensional Ricci flows with surgery always exist for
all time with any compact $3$-manifold as initial metric.

\subsection{Normalized initial conditions}

Consider a compact connected Riemannian $3$-manifold $(M,g(0))$
satisfying
\begin{enumerate}
\item[(1)] $|{\rm Rm}(x,0)|\le 1$ for all $x\in M$ and
\item[(2)] for every $x\in M$ we have ${\rm Vol}\,B(x,0,1)\ge \omega/2$ where
$\omega$ is the volume of the unit ball in $\Ar^3$.
\end{enumerate}

Under these conditions we say that $(M,g(0))$ is {\em normalized}.
Also, if $(M,g(0))$ is the initial manifold of a Ricci flow with
surgery then we say that it is a {\em normalized initial metric}\index{initial metric, normalized|ii}. Of
course, given any compact Riemannian $3$-manifold $(M,g(0))$ there
is a positive constant $Q<\infty$ such that $(M,Qg(0))$ is
normalized.

Starting with a normalized initial metric implies that the flow exists
and has uniformly bounded curvature for a fixed amount of time. This is the
content of the following claim which is an immediate corollary of
Theorem~\ref{compuniq}, Proposition~\ref{patch}, Theorem~\ref{shi},
and  Proposition~\ref{kappa0r0t0}.

\begin{claim}\label{newkappa0r0t0}
 There is $\kappa_0$ such that the following holds. Let $(M,g(0))$ be a normalized initial
metric. Then the solution to the Ricci flow equation with these
initial conditions exists for $t\in [0,2^{-4}]$, and  $|R(x,t)|\le
2$ for all $x\in M$ and all $t\in [0,2^{-4}]$. Furthermore, for any
$t\in [0,2^{-4}]$ and any $x\in M$ and any $r\le \epsilon$ we have
${\rm Vol}\,B(x,t,r)\ge \kappa_0r^3$.
\end{claim}

\section{Gluing together evolving necks}

\begin{prop}\label{neckglue}
There is $0<\beta<1/2$ such that the following holds for any
$\epsilon<1$. Let $(N\times [-t_0,0],g_1(t))$ be an evolving
$\beta\epsilon$-neck centered at $x$ with $R(x,0)=1$. Let $
(N'\times (-t_1,-t_0],g_2(t))$ be a strong $\beta\epsilon/2$-neck.
Suppose we have an isometric embedding of $N\times \{-t_0\}$ with
$N'\times \{-t_0\}$ and the strong $\beta\epsilon/2$-neck structure on $N'\times
(-t_1,-t_0]$ is centered at the image of $(x,-t_0]$. Then the union
$$N\times [-t_0,0]\cup N'\times (-t_1,-t_0]$$ with the induced
one-parameter family of metrics contains a strong $\epsilon$-neck
centered at $(x,0)$.
\end{prop}

\begin{figure}[ht]
  \centerline{\epsfbox{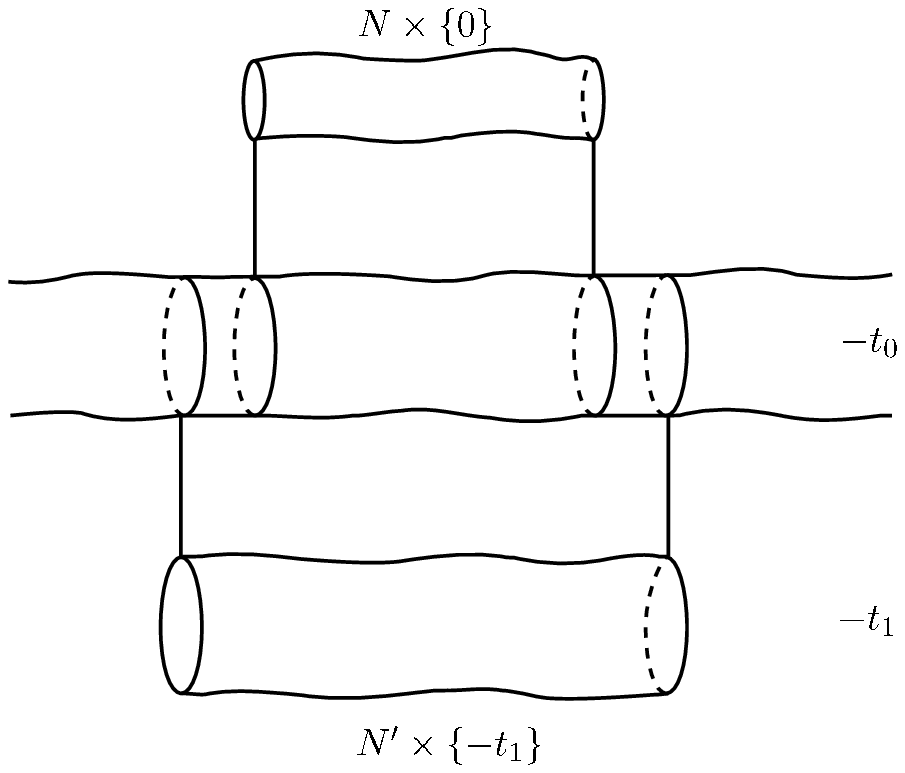}}
  \caption{Gluing together necks.}\label{fig:gluing}
\end{figure}

\begin{proof}
Suppose that the result does not hold. Take a sequence of $\beta_n$
tending to zero and counterexamples $(N_n\times
[-t_{0,n},0],g_{1,n}(t));\ (N'_n\times
(-t_{1,n},-t_{0,n}],g_{2,n}(t))$. Pass to a subsequence so that the
$t_{0,n}$ tend to a limit $t_{0,\infty}\ge 0$. Since $\beta_n$ tends
to zero, we can take a smooth limit of a subsequence and this limit
is an evolving cylinder $(S^2\times \Ar, h(t)\times ds^2)$, where
$h(t)$ is the round metric of scalar curvature $1/(1-t)$ defined for
some amount of backward time. Notice that, for all $\beta$
sufficiently small, on a $\beta\epsilon$-neck the derivative of the
scalar curvature is positive. Thus, $R_{g_{1,n}}(x,-t_{0,n})<1$.
Since we have a strong neck structure on $N'_n$ centered at
$(x,-t_{0,n})$, this implies that $t_{1,n}>1$ so that the limit is
defined for at least time $t\in [0, 1+t_{0,\infty})$. If
$t_{0,\infty}>0$, then, restricting to the appropriate subset of
this limit, a subset with compact closure in space-time, it follows
immediately that for all $n$ sufficiently large there is a strong
$\epsilon$-neck centered at $(x,0)$. This contradicts the assumption
that we began with a sequence of counterexamples to the
proposition.

Let us consider the case when $t_{0,\infty}=0$. In this case the
smooth limit is an evolving round cylinder defined for time
$(-1,0]$.
 Since $t_{1,n}>1$ we
see that for any $A<\infty$ for all $n$ sufficiently large the ball
$B(x_n,0,A)$ has compact closure in every time-slice and there are
uniform bounds to the curvature on $B(x_n,0,A)\times (-1,0]$. This
means that the limit is uniform for time $(-1,0]$ on all these
balls. Thus, once again for all $n$ sufficiently large we see that
$(x,0)$ is the center of a strong $\epsilon$-neck in the union. In
either case we have obtained a contradiction, and hence we have
proved the result. See {\sc Fig.}~\ref{fig:gluing}.
\end{proof}

\subsection{First assumptions}\label{firstassumpt}

\noindent{\bf Choice of $C$ and $\epsilon$:} The first thing we need
to do is fix for the rest of the argument $C<\infty$ and
$\epsilon>0$. We do this in the following way. We fix $0<\epsilon\le
{\rm
min}(1/200,\left(\sqrt{D}(A_0+5)\right)^{-1},\bar\epsilon_1/2,\bar
\epsilon'/2,\epsilon_0)$ where  $\bar\epsilon_1$ is the constant from
Proposition~\ref{narrows}, $\bar\epsilon'$ is the constant from
Theorem~\ref{kappasummary}, $\epsilon_0$ is the constant from Section~\ref{10.1},
and $A_0$ and $D$ are the constants from
Lemma~\ref{A_0}.
 We fix $\beta<1/2$, the constant from
Proposition~\ref{neckglue}.  Then we let $C$ be the maximum of the
constant $C(\epsilon)$ as in Corollary~\ref{kappacannbhd} and
$C'(\beta\epsilon/3)+1$ as in Theorem~\ref{stdsolncannbhd}.

For all such $\epsilon$, Theorem~\ref{bcbd} holds for $\epsilon$ and
Proposition~\ref{narrows}, Proposition~\ref{canonvary} and
Corollaries~\ref{kappacannbhd} and~\ref{limitcannbhd} and Theorems~\ref{smlmtflow} and~\ref{kaplimit}
 hold for $2\epsilon$. Also, all the topological results of the Appendix hold
for $2\epsilon$ and $\alpha=10^{-2}$.

Now let us turn to the assumptions we shall make on the Ricci flows
with surgery that we shall consider. Let ${\mathcal M}$ be a
space-time. Our first  set of assumptions are basically topological
in nature. They are:

\bigskip

\noindent {\bf Assumption (1). Compactness and dimension:} {\em Each
time-slice $M_t$ of space-time is a compact $3$-manifold
containing no embedded $\Ar P^2$ with trivial normal bundle.}

\bigskip

\noindent{\bf Assumption (2). Discrete singularities:} {\em The set of
singular times is a discrete subset of  $\Ar$.}

\bigskip

\noindent {\bf Assumption (3). Normalized initial conditions:} {\em
$0$ is the initial time of the Ricci flow with surgery and the
initial metric $(M_0,G(0))$ is normalized.}

\bigskip

 It follows
from Assumption (2) that for any time $t$ in the time-interval of
definition of a Ricci flow with surgery, with $t$ being distinct
from the initial and final times (if these exist), for all
$\delta>0$ sufficiently small, the only possible singular time in
$[t-\delta,t+\delta]$ is $t$. Suppose that $t$ is a singular time.
The singular locus at time $t$ is a closed, smooth subsurface
$\Sigma_{t}\subset M_{t}$\index{$\Sigma_t$|ii}. From the local
model, near every point of $x\in \Sigma_{t}$ we see that this
surface separates $M_{t}$ into two pieces:
$$M_{t}=C_{t}\cup_{\Sigma_{t}}E_{t},$$
where $E_{t}$\index{$E_t$|ii} is the exposed region at time $t$ and $C_{t}$\index{$C_t$|ii} is the
complement of the interior of $E_{t}$ in $M_{t}$. We call $C_t$ the
{\em continuing region}\index{continuing region|ii}.  $C_t\subset M_t$ is the maximal subset of
$M_t$ for which there is $\delta>0$ and an embedding $C_t\times
(t-\delta,t]\to {\mathcal M}$ compatible with time and the vector
field.

\noindent {\bf Assumption (4). Topology of the exposed regions:} {\em
At all singular times $t$  we require that $E_{t}$ be a finite
disjoint union of $3$-balls. In particular, $\Sigma_{t}$ is a
finite disjoint union of $2$-spheres.}

\bigskip

The next assumptions are geometric in nature.
Suppose that $t$ is a surgery time. Let ${\mathcal M}_{(-\infty,t)}$ be ${\bf t}^{-1}((-\infty,t))$
and let $(\widehat {\mathcal M}_{(-\infty,t)},\widehat G)$ be the maximal extension of
$({\mathcal M}_{(-\infty,t)},G)$ to time $t$, as given in Definition~\ref{hatM}.

\noindent{\bf Assumption (5). Boundary components of the exposed
regions:} {\em There is a surgery control parameter
function\index{surgery control parameter|ii},
$\delta(t)>0$\index{$\delta(t)$|ii}, a non-increasing function of
$t$, such that each component of $\Sigma_{t}\subset M_{t}$ is the
central $2$-sphere of a strong $\delta(t)$-neck in $(\widehat{\mathcal
M}_{(-\infty,t)},\widehat G)$.}

\bigskip

Suppose that $t$ is a singular time. Then for all $t^-<t$ with $t^-$
sufficiently close to $t$,  the manifolds $M_{t^-}$ are diffeomorphic and are
identified under the flow. Applying the flow (backward) to $C_{t}$ produces a
diffeomorphism from $C_{t}$  onto a compact submanifold with boundary
$C_{t^-}\subset M_{t^-}$. Our next assumption concerns the nature of the
metrics $G(t^-)$ on the {\em disappearing region}\index{disappearing region|ii}
$D_{t^-}=M_{t^-}\setminus C_{t^-}$\index{$D_{t^-}$|ii}. The following holds for
every $t^-<t$ sufficiently close to $t$.

\noindent{\bf Assumption (6). Control on the disappearing region:} {\em For any
singular time $t$, for all $t^-<t$ sufficiently close to $t$, each point of
$x\in D_{t^-}$ has a strong $(C,\epsilon)$-canonical neighborhood in
$M_{t^-}$.}

\noindent{\bf Assumption (7). Maximal flow intervals:} {\em Let $t$ be
the initial time or a singular time and let $t'$ be the first
singular time after $t$ if such exists, otherwise let $t'$ be the
least upper bound of the time-interval of definition of the Ricci
flow with surgery. Then the restriction of the Ricci flow with
surgery to $[t,t')$ is a maximal Ricci flow. That is to say, either
$t'=\infty$ or, as $t\rightarrow t'$ from below, the curvature of $G(t)$
is unbounded so that this restricted Ricci flow cannot be extended
{\bf as a Ricci flow} to any larger time.}

\bigskip

{\bf From now on $C$ and $\epsilon$ have fixed values as described
above and all Ricci flows with surgeries are implicitly assumed to
satisfy Assumptions (1) -- (7).}

\section{Topological consequences of  Assumptions (1) -- (7)}

Next we show that the topological control that we are imposing on
the $3$-dimensional Ricci flows with surgery are enough to allow
us to relate the topology of a time-slice $M_T$ in terms of a later
time-slice $M_{T'}$ and topologically standard pieces. This is the
result that will be used to establish the topological theorems
stated in the introduction.

\begin{prop}\label{surgerytoptype}
Suppose that $({\mathcal M},G)$ is a generalized Ricci flow satisfying Assumptions
(1) -- (7). Let $t$ be a singular time. Then the following holds
for any $t^-<t$ sufficiently close to $t$. The manifold $M_{t^-}$ is
diffeomorphic to a manifold obtained in the following way. Take the disjoint
union of $M_t$, finitely many $2$-sphere bundles over $S^1$, and finitely many
closed $3$-manifolds admitting  metrics of constant positive curvature. Then
perform connected sum operations between (some subsets of) these components.
\end{prop}

\begin{proof}
Fix $t'<t$ but sufficiently close to $t$.
By Assumption 4 every component of $E_t$ is a $3$-ball and hence every component
of $\partial E_t=\partial C_t$ is a $2$-sphere. Since $C_t$ is diffeomorphic to
 $C_{t'}\subset M_{t'}$ we see that every component of $\partial C_{t'}=\partial
 D_{t'}$ is a $2$-sphere. Since every component of $E_t$ is a $3$-ball,
 the passage from the smooth manifold $M_{t'}$ to the smooth manifold
 $M_t$ is effected by removing the
 interior of $D_{t'}$ from $M_{t'}$ and gluing a $3$-ball onto each component of $\partial C_{t'}$
 to form $M_t$.

By Assumption (5) every point of $D_{t'}$ has a strong
$(C,\epsilon)$-canonical neighborhood. Since $\epsilon$ is
sufficiently small it follows from Proposition~\ref{epstopology}
that every component of $D_{t'}$ that is also a component of
$M_{t'}$ is diffeomorphic either to a manifold admitting a metric of
constant positive curvature (a $3$-dimensional space-form), to
$\Ar P^3\#\Ar P^3$ or to a $2$-sphere bundle over $S^1$. In the passage from
$M_{t'}$ to $M_t$ these components are removed.

Now let us consider a component of $D_{t'}$ that is not a component
of $M_{t'}$. Such a component is a connected subset of $M_{t'}$ with
the property that every point is either contained in the core of a
$(C,\epsilon)$-cap or is the center of an $\epsilon$-neck and whose
frontier in $M_{t'}$ consists of $2$-spheres that are central
$2$-spheres of $\epsilon$-necks. If every point is the center of an
$\epsilon$-neck, then according to Proposition~\ref{epschains}
$D_{t'}$ is an $\epsilon$-tube and in particular is diffeomorphic to
$S^2\times I$. Otherwise $D_{t'}$ is contained in a capped or double
capped $\epsilon$-tube. Since the frontier of $D_{t'}$ is non-empty
and is the union of central $2$-spheres of an $\epsilon$-neck, it
follows that either $D_{t'}$ is diffeomorphic to a capped
$\epsilon$-tube or to an $\epsilon$-tube. Hence, these components of
$D_{t'}$ are diffeomorphic either to $S^2\times (0,1)$, to $D^3$, or
to $\Ar P^3\setminus B^3$. Replacing a $3$-ball component of
$D_{t'}$ by another $3$-ball leaves the topology unchanged.
Replacing a component of $D_{t'}$ that is diffeomorphic to
$S^2\times I$ by the disjoint union of two $3$-balls  has the
effect of doing a surgery along the core $2$-sphere of the cylinder
$S^2\times I$ in $M_{t'}$. If this $2$-sphere separates $M_{t'}$
into two pieces then doing this surgery effects a connected sum
decomposition. If this $2$-sphere does not separate, then the
surgery has the topological effect of doing a connected sum
decomposition into two pieces, one of which is diffeomorphic to
$S^2\times S^1$, and then removing that component entirely.
Replacing a component of $D_{t'}$ that is diffeomorphic to $\Ar
P^3\setminus B^3$ by a $3$-ball, has the effect of doing a
connected sum decomposition on $M_{t'}$ into pieces, one of which is
diffeomorphic to $\Ar P^3$, and then removing that component.

 From this description the proposition follows immediately.
\end{proof}

\begin{cor}\label{Tgood0good}
Let $({\mathcal M},G)$ be a generalized Ricci flow satisfying Assumptions
(1) -- (7)  with initial conditions $(M,g(0))$. Suppose
that for some $T$ the time-slice $M_T$ of this generalized flow
satisfies Thurston's Geometrization Conjecture\index{Thurston's
Geometrization Conjecture}. Then the same is true for the manifold
$M_t$ for any $t\le T$, and in particular $M$ satisfies Thurston's
Geometrization Conjecture. In addition,
\begin{enumerate}
\item[(1)] If for some $T>0$ the manifold $M_T$ is empty, then $M$ is a
connected sum of manifolds diffeomorphic to $2$-sphere bundles over
$S^1$ and $3$-dimensional space-forms, i.e., compact
$3$-manifolds that admit a metric of constant positive curvature.
\item[(2)]
If for some $T>0$ the manifold $M_T$ is empty and if $M$ is connected and simply connected,
then $M$ is diffeomorphic to $S^3$.
\item[(3)]
If for some $T>0$ the manifold $M_T$ is empty and if $M$ has finite fundamental group, then $M$
is a $3$-dimensional space-form.
\end{enumerate}
\end{cor}

\begin{proof}
Suppose that $M_T$ satisfies the Thurston Geometrization Conjecture and that $t_0$ is
the largest surgery time $\le T$. (If there is no such surgery time
then $M_T$ is diffeomorphic to $M$ and the result is established.)
Let $T'<t_0$ be sufficiently close to $t_0$ so that
$t_0$ is the only surgery time in the interval $[T',T]$. Then according to
the previous proposition $M_{T'}$ is obtained from $M_T$ by first taking the disjoint union
 of $M_T$ and copies of $2$-sphere bundles over $S^1$ and $3$-dimensional space forms.
 In the Thurston Geometrization Conjecture the first step is to decompose the manifold as
 a connected sum of prime $3$-manifolds and then to treat each prime piece independently.
 Clearly, the prime decomposition of $M_{T'}$ is obtained from the prime decomposition of $M_T$
 by adding a disjoint union with $2$-sphere bundles over $S^1$ and $3$-dimensional space forms.
 By definition
 any $3$-dimensional space-form satisfies Thurston's Geometrization Conjecture.
 Since any diffeomorphism of $S^2$ to itself is isotopic to either the identity
 or to the antipodal map, there are two diffeomorphism types of $2$-sphere bundles over
 $S^1$: $S^2\times S^1$ and the non-orientable $2$-sphere bundle over $S^1$.
 Each is obtained from $S^2\times I$ be gluing the ends together by an isometry of the round
 metric on $S^2$. Hence, each has a homogeneous geometry modeled on $S^2\times \Ar$, and hence
 satisfies Thurston's Geometrization Conjecture. This proves that if $M_T$ satisfies this
 conjecture, then so does $M_{T'}$. Continuing this way by induction, using the fact that there
 are only finitely many surgery times completes the proof of the first statement.

Statement (1) is proved analogously. Suppose that $M_T$ is a disjoint union of
connected sums of $2$-sphere bundles over $S^1$
and $3$-dimensional space-forms. Let $t_0$ be
the largest surgery time $\le T$ and let $T'<t_0$ be sufficiently close to $t_0$.
(As before, if there is no such $t_0$ then $M_T$ is diffeomorphic to $M$ and the result is
established.) Then it is clear from the previous proposition that $M_{T'}$ is also a disjoint
union of connected sums of $3$-dimensional space-forms and $2$-sphere bundles over $S^1$.
Induction as in the previous case completes the argument for this case.

The last two statements are immediate from this one.
\end{proof}

\section{Further conditions on surgery}

\subsection{The surgery parameters}

The process of doing surgery requires fixing the scale $h$ at which one does
the surgery. We shall have to allow this scale $h$ to be a function of time,
decreasing sufficiently rapidly with $t$. In fact, the scale is determined by
two other functions of time which also decay to zero as time goes to infinity--
a canonical neighborhood parameter $r(t)$ determining the curvature threshold
above which we have canonical neighborhoods and the surgery control parameter
$\overline\delta(t)$ determining how close to cylinders (products of the round $2$-sphere
with an interval) the regions where we do
surgery are. In addition to these functions, in order to prove inductively that
we can do surgery we need to have a non-collapsing result. The non-collapsing
parameter $\kappa>0$ also decays to zero rapidly as time goes to infinity.
Here then are the functions that will play the crucial role in defining
the surgery process.

\begin{defn}\label{hrestricted}
We have: (i) {\em a canonical neighborhood parameter}\index{canonical neighborhood parameter|ii},
 $r(t)>0$\index{$r(t)>0$|ii}, and
(ii) {\em a surgery control parameter}\index{surgery control
parameter} $\overline\delta(t)>0$\index{$\delta(t)$}. We use these
to define the surgery scale function $h(t)$. Set
$\rho(t)=\overline\delta(t)r(t)$. Let $h(t)=
h(\rho(t),\overline\delta(t))\le
\rho(t)\cdot\overline\delta(t)=\overline\delta^2(t)r(t)$ be the
function given by Theorem~\ref{hexist}. We require that $h(0)\le
R_0^{-1/2}$ where $R_0$ is the constant from
Theorem~\ref{LOCALSURGERY}.

 In addition, there is a
function $\kappa(t)>0$\index{$\kappa(t)$|ii} called the {\em non-collapsing parameter}\index{non-collapsing parameter|ii}.
All three functions $r(t)$, $\overline\delta(t)$ and $\kappa(t)$ are
required to  be  positive,  non-increasing functions of $t$.
\end{defn}

We shall consider Ricci flows with surgery $({\mathcal M},G)$ that
satisfy Assumptions (1) -- (7) and also satisfy:

\noindent {\bf For any singular time $t$ the surgery at time $t$ is
performed with control $\overline\delta(t)$ and at scale
$h(t)=h(\rho(t),\overline\delta(t))$, where
$\rho(t)=\overline\delta(t)r(t)$, in the sense that each boundary
component of $C_t$ is the central $2$-sphere of a strong
$\overline\delta(t)$-neck centered at a point $y$ with
$R(y)=h(t)^{-2}$.}

There is quite a bit of freedom in the choice of these parameters.
But it is not complete freedom. They must decay rapidly enough as
functions of $t$. We choose to make $r(t)$ and $\kappa(t)$  step
functions, and we require $\overline\delta(t)$ to be bounded above
by a step function of $t$. Let us fix the step sizes.

\begin{defn}
We set $t_0=2^{-5}$, and for any $i\ge 0$ we define $T_i=2^it_0$.
\end{defn}

The steps we consider are $[0,T_0]$ and then $[T_i,T_{i+1}]$ for
every $i\ge 0$. The first step  is somewhat special. Suppose that
$({\mathcal M},G)$ is a Ricci flow with surgery with normalized
initial conditions. Then
 according to
Claim~\ref{newkappa0r0t0} the flow exists on $[0,T_1]$ and the norm
of the Riemann curvature is bounded by $2$ on $[0,T_1]$, so that
by Assumption (7) there are no surgeries in this time interval. Also,
by Claim~\ref{newkappa0r0t0} there is a $\kappa_0>0$ so
that ${\rm Vol}\,B(x,t,r)\le \kappa_0r^3$ for every $t\le T_1$ and
$x\in M_t$ and every $r\le \epsilon$.

\begin{defn}\label{surgeryparams}
{\em Surgery parameter sequences}\index{surgery parameter sequence} are sequences
\begin{enumerate}
\item[(i)] ${\bf r}=r_0\ge r_1\ge r_2\ge \cdots >0$, with $r_0=\epsilon$\index{${\bf r}$\ii},
\item[(ii)] ${\bf K}=\kappa_0\ge \kappa_1\ge \kappa_2\ge \cdots >0$\index{${\bf K}$}
with $\kappa_0$ as in Claim~\ref{newkappa0r0t0}, and \item[(iii)] $
\Delta=\delta_0\ge \delta_1\ge \delta_2\ge \cdots
>0$ with $\delta_0={\rm min}(\beta\epsilon/3,\delta_0',K^{-1},D^{-1})$\index{$\Delta$|ii}
 where $\delta_0'$ is the constant
from Theorem~\ref{LOCALSURGERY} and $\beta<1/2$ is the constant from
Proposition~\ref{neckglue}, $\epsilon$ is the constant that we have
already fixed, and $K$ and $D$ are the constants from
Lemma~\ref{A_0}.
\end{enumerate}
We shall also refer to partial sequences defined for indices $0,\ldots, i$ for some $i>0$
as surgery parameter sequences if they are positive, non-increasing and if their initial
terms satisfy the conditions given above.
\end{defn}

We let $r(t)$ be the step function whose value on $[T_i,T_{i+1})$ is $r_{i+1}$
and whose value on $[0,T_0)$ is $r_0$. We say that a Ricci flow with surgery
satisfies the strong $(C,\epsilon)$-canonical neighborhood assumption with parameter ${\bf
r}$ if it satisfies this condition with respect to the step function $r(t)$
associated with ${\bf r}$. This means that any $x\in {\mathcal M}$ with
$R(x)\ge r^{-2}({\bf t}(x))$ has a strong $(C,\epsilon)$-canonical neighborhood in
${\mathcal M}$. Let $\kappa(t)$ be the step function whose value on
$[T_i,T_{i+1})$ is $\kappa_{i+1}$ and whose value on $[0,T_0)$ is $\kappa_0$.
Given $\kappa>0$, we say that a Ricci flow defined on $[0,t]$ is $\kappa$-non-collapsed
on scales $\le \epsilon$ provided that for every point $x$ not contained in a
component of its time-slice with positive sectional curvature,
if for some $r\le \epsilon$, the parabolic neighborhood $P(x,{\bf t}(x),r,-r^2)$ exists in ${\mathcal M}$
and the norm of the Riemann curvature is bounded on this backward
parabolic neighborhood by $r^{-2}$, then ${\rm Vol}\,B(x,{\bf t}(x),r)\ge
\kappa r^3$.
We say that a Ricci flow with surgery is ${\bf K}$-non-collapsed on scales
$\epsilon$ if for every $t\in [0,\infty)$ the restriction of the flow to
$[0,t]$ is $\kappa(t)$-non-collapsed on scales $\le \epsilon$.
  Lastly, we fix a non-increasing function $\overline\delta(t)>0$
with $\overline\delta(t)\le \delta_{i+1}$ if $t\in [T_i,T_{i+1})$ for all $i\ge
0$ and $\overline\delta(t)\le \delta_0$ for $t\in [0,T_0)$. We denote the fact
that such inequalities hold for all $t$ by saying $\overline\delta(t)\le
\Delta$.

Having fixed surgery parameter sequences ${\bf K}$, ${\bf r}$ and
$\Delta$, defined step functions $r(t)$ and $\kappa(t)$, and fixed
$\overline\delta(t)\le \Delta$ as above, we shall consider only
Ricci flows with surgery  where the surgery at time $t$ is defined
using the surgery parameter functions $r(t)$ and
$\overline\delta(t)$. In addition,  we require that these Ricci
flows with surgery satisfy Assumptions (1) -- (7).

What we shall  show is that there are surgery parameter
sequences ${\bf r}$, ${\bf K}$ and $\Delta$ with the property that
for any normalized initial metric and any positive,
non-increasing function $\overline\delta(t)\le \Delta$, it is
possible to construct a Ricci flow with surgery using the surgery
parameters $r(t)$ and $\overline\delta(t)$ with the given initial
conditions and furthermore that this Ricci flow with surgery
satisfies the Assumptions (1) -- (7), has curvature pinched toward
positive, satisfies the canonical neighborhood assumption, and satisfies the
non-collapsing assumption using these parameters.

In fact we shall prove this inductively, constructing the step
functions inductively one step at a time. Thus, given surgery
parameter sequences indexed by $0,\ldots,i$ we show that there are appropriate choices of
$r_{i+1}, \kappa_{i+1}$ and $\delta_{i+1}$ such that the following is true.
Given a Ricci flow with
surgery defined on time $[0,T_i)$ satisfying all the properties with
respect to the first set of data,  that Ricci
flow with surgery extends to one defined for time $[0,T_{i+1})$ and satisfies
Assumptions (1) -- (7), the canonical neighborhood assumption and the
non-collapsing assumption with respect to the extended surgery
parameter sequences, and has curvature pinched toward positive.
As stated this is not quite true; there is a slight twist: we must also assume
that $\overline\delta(t)\le \delta_{i+1}$ for all $t\in [T_{i-1},T_{i+1})$.
It is for this reason that we consider
pairs consisting of sequences $\Delta$ and a surgery control
parameter $\overline\delta(t)$ bounded above by $\Delta$.

\section{The process of surgery}\label{sect:process}

Suppose given surgery parameter sequences $\{r_0,\ldots,r_{i+1}\}$,
$\{\kappa_0,\ldots,\kappa_{i+1}\}$ and
$\Delta_i=\{\delta_0,\ldots,\delta_i\}$ and also given a positive,
decreasing function $\overline\delta(t)\le \Delta_i$, defined for
$t\le T_{i+1}$ with $\delta_0={\rm
min}(\alpha\epsilon/3,\delta_0',K^{-1},D^{-1})$ as above. Suppose
that $({\mathcal M},G)$ is a Ricci flow with surgery defined for
$t\in [0,T)$ that goes singular at time $T\in (T_i,T_{i+1}]$. We
suppose that it satisfies Assumptions (1) -- (7). Since the flow has
normalized initial conditions and goes singular at time $T$, it
follows that $i\ge 1$. We suppose that $({\mathcal M},G)$ satisfies
the $(C,\epsilon)$-canonical neighborhood assumption with parameter
$r_{i+1}$ and that its curvature is pinched toward positive. By
Theorem~\ref{omega} we know that there is a maximal extension
$(\widehat{\mathcal M},\widehat G)$ of this generalized flow to time
$T$ with the $T$ time-slice being $(\Omega(T),G(T))$. Set
$\rho=\overline\delta(T)r_{i+1}$, and set
$h(T)=h(\rho(T),\overline\delta(T))$ as in Theorem~\ref{hexist}.
Since $\overline\delta(T)\le \delta_0<1$, we see that
$\rho<r_{i+1}$. According to Lemma~\ref{2epshorns} there are
finitely many components of $\Omega(T)$ that meet $\Omega_\rho(T)$.
Let $\Omega^{\rm big}(T)$ be the disjoint union of all the
components of $\Omega(T)$ that meet $\Omega_\rho(T)$.
Lemma~\ref{2epshorns} also tells us that $\Omega^{\rm big}(T)$
contains a finite collection of disjoint $2\epsilon$-horns with
boundary contained in $\Omega_{\rho/2C}$, and the complement of the
union of the interiors of these horns is a compact submanifold with
boundary containing $\Omega_\rho$.  Let ${\mathcal
H}_1,\ldots,{\mathcal H}_j$ be a disjoint union of these
$2\epsilon$-horns. For each $i$ fix a point $y_i\in {\mathcal H}_i$
with $R(y_i)=h^{-2}(T)$. According to Theorem~\ref{hexist} for each
$i$ there is a strong $\delta(T)$-neck centered at $y_i$ and
contained in ${\mathcal H}_i$.
 We orient the $s$-direction of the
neck so that its positive end lies closer to the end of the horn
than its negative end. Let $S^2_i$ be the center of this strong
$\delta(T)$-neck. Let ${\mathcal H}_i^+$ be the unbounded
complementary component of $S^2_i$ in ${\mathcal H}_i$. Let $C_T$ be
the complement of $\coprod_{i=1}^j{\mathcal H}^+_i$ in $\Omega^{\rm
big}(T)$. Then we do surgery on these necks as described in
Section~\ref{necksurgery}, using the constant $q=q_0$ from
Theorem~\ref{LOCALSURGERY},
 removing the positive half of the neck,
and gluing on the cap from the standard solution. This creates a compact
$3$-manifold $M_T=C_T\cup_{\coprod_iS^2_i} B_i$, where each $B_i$ is a copy of the
metric ball of radius $A_0+4$ centered around the tip of the standard solution
 (with the metric scaled by
$h^2(T)$ and then perturbed near the boundary of $B_i$ to
match $g(T)$). Notice that in this process we have removed every
component of $\Omega(T)$ that does not contain a point of
$\Omega_\rho(T)$. The result of this operation is to produce a
compact Riemannian $3$-manifold $(M_T,G_T)$ which is the $T$
time-slice of our extension of $({\mathcal M},G)$. Let $(M_T,G(t)),\
T\le t<T'$, be the maximal Ricci flow with initial conditions
$(M_T,G_T)$ at $t=T$. Our new space-time is the union of
$M_T\times [T,T')$ and $({\mathcal M},G)\cup C_T\times \{T\}$  along $C_T\times \{T\}$.
Since the isometric embedding  $C_T\subset M_T$ extends to an
Here, we view $({\mathcal M},G)\cup C_T\times \{T\}$ as a subspace
of $(\widehat{\mathcal M},\widehat G)$ via the
isometric embedding of $C_T$ into $\Omega(T)$. The time functions
and vector fields glue to provide analogous data for this new
space-time. According to Lemma~\ref{gluingsoln} the horizontal
metrics glue together to make a smooth metric on space-time
satisfying the Ricci flow equation.

Notice that the continuing region at time $T$ is exactly $C_T$
whereas the exposed region is $\coprod_iB_i$, which is a disjoint
union of $3$-balls. The disappearing region is the complement of
the embedding of $C_T$ in $M_{t'}$ for $t'<T$ but sufficiently close to it obtained by
flowing $C_T\subset \Omega_T$ backward. The
disappearing region contains $M_{t'}\setminus \Omega(T)$ and also
contains all components of $\Omega(T)$ that do not contains points
of $\Omega_\rho(T)$, as well as the ends of those components of
$\Omega(T)$ that contain points of $\Omega_\rho(T)$.

\begin{defn}
The operation described in the previous paragraph is the {\em surgery
operation at time $T$} using the surgery parameters $\overline\delta(T)$ and $r_{i+1}$.
\end{defn}

\section{Statements about the existence of Ricci flow with
surgery}

What we shall establish is the existence of surgery satisfying Assumptions (1) -- (7)
above and also satisfying the curvature pinched toward positive assumption, the
strong canonical neighborhood assumption, and the non-collapsing assumption.
This requires first of all that we begin with a compact, Riemannian
$3$-manifold $(M,g(0))$ that is normalized, which we are assuming. It also
requires careful choice of upper bounds $\Delta=\{\delta_i\}$ for the surgery
control parameter $\overline\delta(t)$ and careful choice of the canonical neighborhood
parameter ${\bf r}=\{r_i\}$ and of the non-collapsing step function ${\bf
K}=\{\kappa_i\}$.

Here is the statement that we shall establish.

\begin{thm}\label{MAIN}
There are surgery parameter sequences
$${\bf
K}=\{\kappa_i\}_{i=1}^\infty,\Delta=\{\delta_i\}_{i=1}^\infty,{\bf
r}=\{r_i\}_{i=1}^\infty\index{${\bf K}$}\index{$\Delta$}\index{${\bf
r}$}$$
  such that the following holds. Let $r(t)$
be the step function whose value on $[T_{i-1},T_i)$ is $r_i$.
Suppose that $\overline\delta\colon [0,\infty)\to \Ar^+$ is any
non-increasing function with $\overline\delta(t)\le \delta_i$
whenever $t\in [T_{i-1},T_i)$. Then the following holds: Suppose
that $({\mathcal M},G)$ is a Ricci flow with surgery defined for
$0\le t<T$ satisfying Assumptions (1) -- (7). In addition, suppose that
the following conditions:
\begin{enumerate}
\item[(1)] the generalized flow has curvature
pinched toward positive,
\item[(2)] the flow satisfies the strong $(C,\epsilon)$-canonical neighborhood assumption
with parameter ${\bf r}$ on $[0,T)$, and \item[(3)] the flow is ${\bf K}$
non-collapsed on $[0,T)$ on scales $\le \epsilon$.
\end{enumerate}
 Then there is
an extension of $({\mathcal M},G)$ to a Ricci flow with surgery
defined for all $0\le t<\infty$ and satisfying Assumptions (1) -- (7) and
the above three conditions.
\end{thm}

This of course leads immediately to the existence result for Ricci
flows with surgery defined for all time with any normalized initial
conditions.

\begin{cor}\label{RFSexists}
Let ${\bf K}$, ${\bf r}$ and $\Delta$ be surgery parameter sequences
provided by the previous theorem. Let $\overline\delta(t)$ be a
non-increasing positive function with $\overline\delta(t)\le
\Delta$. Let $M$ be a compact $3$-manifold containing no $\Ar P^2$
with trivial normal bundle. Then there is a Riemannian metric $g(0)$
on $M$ and a Ricci flow with surgery defined for $0\le t<\infty$
with  initial metric $(M,g(0))$. This Ricci flow with surgery
satisfies the seven assumptions and is ${\bf K}$-non-collapsed on
scales $\le \epsilon$. It also satisfies the strong
$(C,\epsilon)$-canonical neighborhood assumption with parameter
${\bf r}$ and has curvature pinched toward positive. Furthermore,
any surgery at a time $t\in [T_i,T_{i+1})$ is done using
$\overline\delta(t)$ and $r_{i+1}$.
\end{cor}

\begin{proof} (Assuming Theorem~\ref{MAIN})
Choose a metric $g(0)$ so that $(M,g_0)$ is normalized. This is
possible by beginning with any Riemannian metric on $M$ and scaling
it by a sufficiently large positive constant to make it normalized.
According to Proposition~\ref{kappa0r0t0} and the definitions of
$T_i$ and $\kappa_0$ there is a Ricci flow $(M,g(t))$ with these
initial conditions defined for $0\le t\le T_2$ satisfying Assumptions (1) -- (7)
and the three conditions of the previous theorem. The assumption that
$M$ has no embedded $\Ar P^2$ with trivial normal bundle is needed so that Assumption (1) holds
for this Ricci flow.
 Hence, by the previous theorem we can extend this
Ricci flow to a Ricci flow with surgery defined for all $0\le
t<\infty$ satisfying the same conditions.
\end{proof}

Showing that after surgery Assumptions (1) -- (7) continue to hold and
that the curvature is pinched toward positive is direct and only requires that
$\overline\delta(t)$ be smaller than some universal positive constant.

\begin{lem}\label{14.26}
Suppose that $({\mathcal M},G)$ is a Ricci flow with surgery going
singular at time $T\in [T_{i-1},T_i)$. We suppose that $({\mathcal
M},G)$ satisfies Assumptions (1) - (7),  has curvature pinched toward
positive, satisfies the strong $(C,\epsilon)$-canonical neighborhood
assumption with parameter ${\bf r}$ and is ${\bf K}$ non-collapsed.
Then the result of the surgery operation  at time $T$ on $({\mathcal M},G)$ is a
 Ricci flow with surgery defined on $[0,T')$ for some $T'>T$.
 The resulting Ricci flow with
surgery satisfies Assumptions (1) -- (7). It also has curvature pinched
toward positive.
\end{lem}

\begin{proof}
It is immediate from the construction and  Lemma~\ref{gluingsoln} that
the result of performing the surgery operation at time $T$ on a Ricci flow with
surgery produces a new Ricci flow with surgery. Assumptions (1) -- (3)
clearly hold for the result. and Assumptions (4) and (5) hold because of
the way that we do surgery. Let us consider Assumption (6). Fix $t'<T$
so that there are no surgery times in $[t',T)$. By flowing backward
using the vector field $\chi$ we have an embedding $\psi\colon
C_t\times [t',T]\to \widehat {\mathcal M}$ compatible with time and
the vector field. For any $p\in M_{t'}\setminus\psi({\rm
int}\,C_T\times \{t'\})$ the limit as $t$ tends to $T$ from below of
the flow line $p(t)$ at time $t$ through $p$ either lies in
$\Omega(T)$ or it does not. In the latter case, by definition we
have
$${\rm lim}_{t\rightarrow T^-}R(p(t))=\infty.$$
In the former case, the limit point  either is contained in the end
of a strong $2\epsilon$-horn cut off by the central $2$-sphere of
the strong $\delta$-neck centered at one of the $y_i$ or is
contained in a component of $\Omega(T)$ that contains no point of
$\Omega_\rho(T)$. Hence, in this case we have
$${\rm lim}_{t\rightarrow T^-}R(p(t))>\rho^{-2}>r_i^{-2}.$$
Since $M_{t'}\setminus \psi({\rm int}\,C_T\times \{t'\})$ is compact
for every $t'$, there is $T_1<T$ such that $R(p(t))>r_i^{-2}$ for
all $p\in M_{t'}\setminus \psi({\rm int}C_T\times \{t'\})$ and all
$t\in [T_1,T)$. Hence, by our assumptions all these points have
strong $(C,\epsilon)$-canonical neighborhoods. This establishes that
Assumption (6) holds at the singular time $T$. By hypothesis
Assumption (6) holds at all earlier singular times. Clearly, from the
construction the Ricci flow on $[T,T')$ is maximal. Hence,
Assumption (7) holds for the new Ricci flow with surgery.

 From Theorem~\ref{LOCALSURGERY} the
fact that $\delta(T)\le \delta_i\le \delta_0\le\delta'_0$ and
$h(T)\le R_0^{-1/2}$ imply that the Riemannian manifold $(M_T,G(T))$
has curvature pinched toward positive for time $T$. It then follows
from Corollary~\ref{pincha} that  the Ricci flow defined on $[T,T')$
with $(M_T,G(T))$ as initial conditions has curvature pinched toward
positive. The inductive hypothesis is that on the time-interval
$[0,T)$ the Ricci flow with surgery has curvature pinched toward
positive. This completes the proof of the lemma.
\end{proof}

\begin{prop}\label{surgerydist}
Suppose that $({\mathcal M},G)$ is a Ricci flow with surgery
satisfying Assumptions (1) -- (7) in Section~\ref{firstassumpt}.
Suppose that $T$ is a surgery time, suppose  that the surgery
control parameter $\delta(T)$ is less than $\delta_0$ in
Definition~\ref{surgeryparams}, and suppose that the scale of the
surgery $h(T)$ is less than $R_0^{-1/2}$ where $R_0$ is the constant
from Theorem~\ref{LOCALSURGERY}. Fix $t'<T$  sufficiently close to
$T$. Then  there is an embedding $\rho\colon M_{t'}\times [t',T)\to
{\mathcal M}$ compatible with time and the vector field. Let $X(t')$
be a component of $M_{t'}$ and let $X(T)$ be a component obtained
from $X(t')$ by doing surgery at time $T$. We view $\rho^*G$ as a
one-parameter family of metrics $g(t)$ on $X(t')$. There is an open
subset $\Omega\subset X(t')$ with the property that ${\rm
lim}_{t'\rightarrow T^-} g(t')|_\Omega$ exists (we denote it by
$g(T)|_\Omega$) and with the property that $\rho|_{\Omega\times
[t',T)}$ extends to a map $\widehat \rho\colon \Omega\times
[t',T]\to {\mathcal M}$. This defines a map for $\Omega\subset
X(t')$ onto an open subset $\Omega(T)$ of $X(T)$ which is an
isometry from the limiting metric $g(T)$ on $\Omega$ to
$G(T)|_{\Omega}$. Suppose that all of the $2$-spheres along which we
do surgery are separating. Then this map extends to a map $X(t')\to
X(T)$. For all $t<T$ but sufficiently close to $T$ this extension is
a distance decreasing map from $(X(t')\setminus \Omega,g(t))$ to
$X(T)$.
\end{prop}

\begin{proof}
This is immediate from the third item in Theorem~\ref{LOCALSURGERY}.
\end{proof}

\begin{rem}
If we have a non-separating surgery $2$-sphere then there will a
component $X(T)$ with surgery caps on both sides of the surgery
$2$-sphere and hence we cannot extend the map even continuously
over all of $X(t')$.
\end{rem}

The other two inductive properties in Theorem~\ref{MAIN} -- that the result is
${\bf K}$-non-collapsed and also that it satisfies the strong
$(C,\epsilon)$-canonical neighborhood assumption with parameter ${\bf r}$ -
require appropriate inductive choices of the sequences. The arguments
establishing these are quite delicate and intricate. They are given in the next
two sections.

\section{Outline of the proof of Theorem~\protect{\ref{MAIN}}}

Before giving the proof proper of Theorem~\ref{MAIN} let us outline
how the argument goes. We shall construct the surgery parameter
sequences $\Delta$, ${\bf r}$, and ${\bf K}$ inductively. Because of
Lemma~\ref{kappa0r0t0} we have the beginning of the inductive
process. We suppose that we have defined sequences as required up to
index $i$ for some $i\ge 1$. Then we shall extend them one more step
to sequences defined up to $(i+1)$, though there is a twist: to do
this we must redefine $\delta_i$ in order to make sure that the
extension is possible. In Chapter~\ref{sectnoncoll} we establish the
non-collapsing result assuming the strong canonical neighborhood
result. More precisely, suppose that we have a Ricci flow with
surgery $({\mathcal M},G)$ defined for time $0\le t<T$ with $T\in
(T_i,T_{i+1}]$ so that the restriction of this flow to the
time-interval $[0,T_i)$ satisfies the inductive hypothesis with
respect to the given sequences. Suppose also that the entire Ricci
flow with surgery has
 strong $(C,\epsilon)$-canonical neighborhoods for some $r_{i+1}>0$.
Then there is $\delta(r_{i+1})>0$ and $\kappa_{i+1}>0$ such that,
provided that $\overline\delta(t)\le \delta(r_{i+1})$ for all $t\in
[T_{i-1},T)$, the Ricci flow with surgery $({\mathcal M},G)$ is
$\kappa_{i+1}$ non-collapsed on scales $\le \epsilon$.

In Section~\ref{sectcannbhd} we show that the strong $(C,\epsilon)$-canonical neighborhood
assumption extends for some parameter $r_{i+1}$, assuming again that
$\overline\delta(t)\le \delta(r_{i+1})$ for all $t\in [T_{i-1},T)$.

Lastly, in Section~\ref{sectcomplete} we complete the proof by
showing that the number of surgeries possible in $[0,T_{i+1})$ is
bounded in terms of the initial conditions and $\overline\delta(T)$.
The argument for this is a simple volume comparison argument --
under Ricci flow with normalized initial conditions, the volume
grows at most at a fixed exponential rate and under each surgery an
amount of volume, bounded below by a positive constant depending
only on $\overline \delta(T_{i+1})$, is removed.

\chapter{Proof of the non-collapsing}\label{sectnoncoll}

The precise statement of the non-collapsing result is given in the
next section. Essentially, the proof of non-collapsing in the
context of Ricci flow with surgery is the same as the proof in the
case of ordinary Ricci flows. Given a point $x\in {\mathcal M}$, one
finds a parabolic neighborhood whose size, $r'$, is determined by
the constants $r_i$, $C$ and $\epsilon$, contained in ${\bf
t}^{-1}([T_{i-1},T_i))$ and on which the curvature is bounded by
$(r')^{-2}$. Hence, by the inductive hypothesis, the final
time-slice of this neighborhood is $\kappa_i$-non-collapsed.
Furthermore, we can choose this neighborhood so that the reduced
${\mathcal L}$-length of its central point from $x$ is bounded by
$3/2$. This allows us to produce an open subset at an earlier time
whose reduced volume is bounded away from zero. Then using
Theorem~\ref{THM} we transfer this conclusion to a non-collapsing
conclusion for $x$. The main issue in this argument is to show that
there is a point in each earlier time-slice whose reduced length
from $x$ is at most $3/2$. We can argue as in the case of a Ricci
flow if we can show that any curve parameterized by backward time
starting at $x$ (a point where the hypothesis of
$\kappa$-non-collapsing holds) that comes close to  a surgery cap
either from above or below must have large ${\mathcal L}$-length. In
establishing the relevant estimates we are forced to require that
$\delta_i$ be sufficiently small.

\section{The statement of the non-collapsing result}

Here, we shall assume that after surgery the strong canonical
neighborhood assumption holds, and we shall establish the
non-collapsing result.

\begin{prop}\label{KAPPALIMIT}
Suppose that for some $i\ge 0$ we have surgery parameter sequences $\delta_0\ge
\delta_1\ge\cdots\ge \delta_i>0$, $\epsilon=r_0\ge r_1\ge \cdots\ge
r_i>0$ and $\kappa_0\ge \kappa_1\ge \cdots\ge \kappa_i>0$.
 Then there is $0<\kappa\le \kappa_i$ and for any $0<r_{i+1}\le r_i$ there is
$0<\delta(r_{i+1})\le \delta_i$ such that the following holds. Suppose that
$\overline\delta\colon [0,T_{i+1}]\to \Ar^+$ is a non-increasing function with
$\overline\delta(t)\le \delta_j$ for all $t\in [T_j,T_{j+1})$ and
$\overline\delta(t)\le \delta(r_{i+1})$ for all $t\in [T_{i-1},T_{i+1})$.
Suppose that $({\mathcal M},G)$ is a Ricci flow with surgery defined for $0\le
t<T$ for some $T\in (T_i,T_{i+1}]$ with surgery control parameter
$\overline\delta(t)$. Suppose that the restriction of this Ricci flow with
surgery  to the time-interval $[0,T_i)$ satisfies the hypothesis of
Theorem~\ref{MAIN} with respect to the given sequences. Suppose also that the
entire Ricci flow with surgery $({\mathcal M},G)$ satisfies Assumptions (1) --
(7) and the strong $(C,\epsilon)$-canonical neighborhood assumption with
parameter $r_{i+1}$. Then $({\mathcal M},G)$ is $\kappa$-non-collapsed on all
scales $\le \epsilon$.
\end{prop}

\begin{rem}
Implicitly,
$\kappa$ and $\delta(r_{i+1})$ are also allowed to depend on $t_0, \epsilon$,
and $C$, which are fixed, and also $i+1$.
Also recall that the non-collapsing condition allows for two outcomes: if $x$ is a point at which the
hypothesis of the non-collapsing hold, then
there is a lower bound on the volume of a ball centered at $x$, or  $x$
is contained in a component of its time-slice that has positive sectional curvature.
\end{rem}

\section{The proof of non-collapsing when $R(x)=r^{-2}$ with
$r\le r_{i+1}$}\label{rsmall}

Let us begin with an easy case of the non-collapsing result, where
non-collapsing follows easily from the strong canonical neighborhood
assumption, rather than from using ${\mathcal L}$-length and
monotonicity along ${\mathcal L}$-geodesics. We suppose that we have
a Ricci flow with surgery $({\mathcal M},G)$ defined for $0\le t<T$
with $T\in [T_i,T_{i+1})$, and a constant $r_{i+1}\le r_i$, all
satisfying the hypothesis of Proposition~\ref{KAPPALIMIT}. Here is
the result that establishes the non-collapsing in this case.

\begin{prop}
Let $x\in {\mathcal M}$ with ${\bf t}(x)=t$ and with $R(x)=r^{-2}\ge r_{i+1}^{-2}$. Then there is $\kappa>0$ depending only
on $C$  such that ${\mathcal M}$ is $\kappa$-non-collapsed at $x$; that is to say,
if $R(x)=r^{-2}$ with $r\le r_{i+1}$, then ${\rm
Vol}\,B(x,t,r)\ge \kappa r^3$, or $x$ is contained in a component of $M_t$ with positive sectional curvature.
\end{prop}

\begin{proof}
 Since
$R(x)\ge r_{i+1}^{-2}$, by assumption any such $x$ has a strong
$(C,\epsilon)$-canonical neighborhood. If this neighborhood is a
strong $\epsilon$-neck centered at $x$ then the result is clear for
a non-collapsing constant $\kappa$ which is universal. If the
neighborhood is an $\epsilon$-round component containing $x$, then
$x$ is contained in a component of positive sectional curvature. Likewise,
if $x$ is contained in a $C$-component then by definition it is contained in a component of its
time-slice with positive sectional curvature.

Lastly, we suppose that $x$ is contained in the core $Y$ of a
$(C,\epsilon)$-cap ${\mathcal C}$. Let $r'>0$ be such that the
supremum of $|{\rm Rm}|$ on $B(x,t,r')$ is $(r')^{-2}$.
Then, by the
definition
 of a $(C,\epsilon)$-cap,
${\rm vol}\,B(x,t,r')\ge C^{-1}(r')^3$.
Clearly, $r'\le r$ and there is a point $y\in \overline{B(x,t,r')}$
with $R(y)=(r')^{-2}$. On the other hand, by the
definition of a $(C,\epsilon)$-cap, we have $R(y)/R(x)\le C$, so that $r'/r\ge C^{-1/2}$.
Thus, the
volume of $B(x,t,r)$ is at least $C^{-5/2}r^3$.

This completes an examination of all cases and establishes the
proposition.
\end{proof}

\section{Minimizing ${\mathcal L}$-geodesics exist when
$R(x)\le r^{-2}_{i+1}$: the statement}

The proof of the non-collapsing result when $R(x)=r^{-2}$ with
$r_{i+1}<r\le \epsilon$ is much more delicate. As we indicated
above, it is analogous to the proof of non-collapsing for Ricci
flows given in Chapter~\ref{noncoll}. That is to say, in this case
the result is proved using the length function on the Ricci flow
with surgery and the monotonicity of the reduced volume. Of course,
unlike the case of Ricci flows treated in Chapter~\ref{noncoll},
here not all points of a Ricci flow with surgery ${\mathcal M}$ can
be reached by minimizing ${\mathcal L}$-geodesics, or rather more
precisely by minimizing ${\mathcal L}$-geodesics contained in the
open subset of smooth points of ${\mathcal M}$. (It is only for the
latter ${\mathcal L}$-geodesics that the analytic results of
Chapter~\ref{lengthfn} apply.) Thus, the main thing to establish in
order to prove non-collapsing is that for any Ricci flow with
surgery $({\mathcal M},G)$ satisfying the hypothesis of
Proposition~\ref{KAPPALIMIT} there are minimizing ${\mathcal
L}$-geodesics in the open subset of smooth points of ${\mathcal M}$
to `enough' of ${\mathcal M}$ so that we can run the same reduced
volume argument that worked in Chapter~\ref{noncoll}. Here is the
statement that tells us that there are minimizing ${\mathcal
L}$-geodesics to `enough' of ${\mathcal M}$.

\begin{prop}\label{delta0ri+1}
For each $r_{i+1}$ with $0<r_{i+1}\le r_i$, there is
$\delta=\delta(r_{i+1})>0$ (depending implicitly on $t_0$, $C$,
$\epsilon$, and $i$) such that if $\overline\delta(t)\le \delta$ for
all $t\in [T_{i-1},T_{i+1}]$ then the following holds. Let
$({\mathcal M},G)$ be a Ricci flow with surgery satisfying the
hypothesis of Proposition~\ref{KAPPALIMIT} with respect to the given
sequences and $r_{i+1}$, and let $x\in {\mathcal M}$ have ${\bf
t}(x)=T$ with $T\in [T_i,T_{i+1})$. Suppose that for some $r\ge
r_{i+1}$ the parabolic neighborhood $P(x,r,T,-r^2)$ exists in
${\mathcal M}$ and $|{\rm Rm}|\le r^{-2}$ on this neighborhood. Then
there is an open subset $U$ of ${\bf t}^{-1}[T_{i-1},T)$ contained
in the open subset of smooth manifold points  of ${\mathcal M}$ with
the following properties:
\begin{enumerate}
\item[(1)] For every $y$ in $U$ there is a minimizing ${\mathcal L}$-geodesic
connecting $x$ to $y$.
\item[(2)] $U_t=U\cap {\bf t}^{-1}(t)$ is non-empty for every $t\in [T_{i-1},T)$
\item[(3)] For each $t\in [T_{i-1},T)$ the restriction of ${\mathcal
L}$ to $U_t$ achieves its minimum and that minimum is at most
$3\sqrt{(T-t)}$.
\item[(4)] The subset of $U$ consisting of all $y$ with the property that
${\mathcal L}(y)\le {\mathcal L}(y')$ for all $y'\in {\bf
t}^{-1}({\bf t}(y))$ has the property that its intersection with
${\bf t}^{-1}(I)$ is compact for every compact interval $I\subset
[T_{i-1},T)$.
\end{enumerate}
\end{prop}

The basic idea in proving this result is to show that all paths
beginning at $x$ and parameterized by backward time that come close
to the exposed regions have large ${\mathcal L}$-length. If we can
establish this, then the existence of such paths will not be an
impediment to using the analytic estimates from
Chapter~\ref{lengthfn} to show that for each $t\in [T_{i-1},T)$
there is a point whose ${\mathcal L}$-length from $x$ is at most
$3\sqrt{T-t}$, and that the set of points that minimize the
${\mathcal L}$-length from $x$ in a given time-slice form a compact
set.

Given Proposition~\ref{delta0ri+1},  arguments from
Chapter~\ref{noncoll} will be applied to complete the proof of
Proposition~\ref{KAPPALIMIT}.

\section{Evolution of neighborhoods of surgery caps}

 We begin this analysis required to prove Proposition~\ref{delta0ri+1}
 by studying the
evolution of surgery caps\index{surgery cap!evolution of|ii}.
Proposition~\ref{altern} below is the main result along these lines.
Qualitatively, it says that if the surgery control parameter
$\delta$ is sufficiently small, then as a surgery cap evolves in a
Ricci flow with surgery it stays near the rescaled version of the
standard flow for any rescaled time less than one unless the entire
cap is removed (all at once) by some later surgery. In that case,
the evolution of the cap is close to the rescaled version of the
standard flow until it is removed. Using this result we will show
that if a path parameterized by backward time has final point  near
a surgery cap and has initial point with scalar curvature not too
large, then this path must enter this evolving neighborhood either
from the `top' or `size' and because of the estimates that we derive
in this chapter such a path  must have large ${\mathcal L}$-length.

\begin{figure}[ht]
  \centerline{\epsfbox{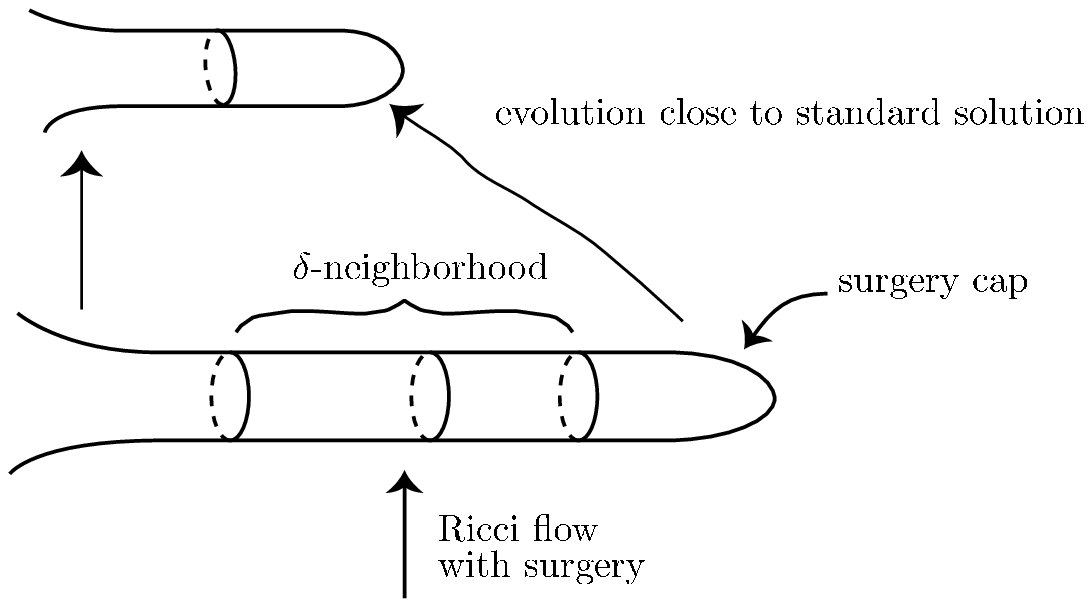}}
  \caption{Evolution of a surgery cap.}\label{fig:surgcap}
\end{figure}

\begin{prop}\label{altern}
Given $A<\infty$, $\delta''>0$ and $0<\theta<1$, there is
$\delta''_0=\delta''_0(A,\theta,\delta'')$ ($\delta_0''$ also
depends on $r_{i+1}$, $C$, and $\epsilon$, which are all now fixed)
such that the following holds. Suppose that $({\mathcal M},G)$ is a
Ricci flow with surgery defined for $0\le t<T$ with surgery control
parameter $\overline\delta(t)$. Suppose that it satisfies the strong
$(C,\epsilon)$-canonical neighborhood assumption at all points $x$
with $R(x)\ge r_{i+1}^{-2}$. Suppose also that $({\mathcal M},G)$ has
curvature that is pinched toward positive. Suppose that there is a
surgery at some time $\bar t$ with $T_{i-1}\le \bar t<T$ with $\bar
h$ as the surgery scale parameter. Set $T'={\rm min}(T,\bar t+\theta
\bar h^2)$. Let $p\in M_{\bar t}$ be the tip of the cap of a surgery
disk. Then, provided that $\bar\delta(\bar t)\le \delta''_0$ one of
the following holds:
\begin{enumerate}
\item[(a)] There is an embedding $\rho\colon B(p,\bar t,A\bar h)\times [\bar t,T')\to {\mathcal M}$
compatible with time and the vector field. Let $g'(t),\ \bar t\le
t<T'$, be the one-parameter family of metrics on $B(p,\bar t,A\bar
h)$ given by $\rho^*G$. Shifting this family by $-\bar t$ to make
the initial time $0$ and scaling it by $(\bar h)^{-2}$ produces a
family of metrics $g(t),\ 0\le t<{\rm min}((T-\bar t)\bar
h^{-2},\theta)$, on $B_{g}(p,0,A)$ that are within $\delta''$ in the
$C^{[1/\delta'']}$-topology of the standard flow on the ball of
radius $A$ at time $0$ centered at the tip of its cap.
\item[(b)] There is $\bar t_+\in (\bar t,T')$ and an embedding
$B(p,\bar t,A\bar h)\times [\bar t,\bar t_+)\to {\mathcal M}$
compatible with time and the vector field so that the previous item
holds with $\bar t_+$ replacing $T'$. Furthermore, for any $t<\bar
t_+$ but sufficiently close to $\bar t_+$ the image of $B(p,\bar
t,A\bar h)\times\{t\}$ is contained in the region $D_{t}\subset M_t$
that disappears at time $\bar t_+$.
\end{enumerate}
See {\sc Fig.}~\ref{fig:surgcap}.
\end{prop}

\index{surgery cap!evolution of}

\begin{proof} The method of proof is to assume that the result is false and take
a sequence of counterexamples with surgery control parameters
$\delta_n$ tending to zero. In order to derive a contradiction we
need to be able to take smooth limits of rescaled versions of these
Ricci flows with surgery, where the base points are the tips of the
surgery caps. This is somewhat delicate since the surgery cap is not
the result of moving forward for a fixed amount of time under Ricci
flow, and consequently Shi's theorem does not apply. Fortunately,
the metrics on the cap are bounded in the $C^\infty$-topology so
that Shi's theorem with derivatives does apply. Let us start
by examining limits of the sort we need to take.

\begin{claim}
Let $(N,g_N)$ be a strong $\delta'$-neck with $N_0$ its middle half.
Suppose that $({\mathcal S},g)$ is the result of doing surgery on
(the central $2$-sphere) of $N$, adding a surgery cap ${\mathcal C}$
to $N^-$. Let $h$ be the scale of $N$. Let $({\mathcal S}_0(N),g')$
be the union of $N^-_0\cup {\mathcal C}$ with its induced metric as
given in Section~\ref{necksurgery}, and let $({\mathcal
S}_0(N),\widehat g_0)$ be the result of rescaling $g_0$ by $h^{-2}$.
Then for every $\ell<\infty$ there is a uniform bound to
$|\nabla^\ell {\rm Rm}_{\widehat g_0}(x)|$ for all $x\in {\mathcal
S}_0(N)$.
\end{claim}

\begin{proof}
Since $(N,g_N)$ is a strong $\delta'$-neck of scale $h$, there is a
Ricci flow on $N$ defined for backward time $h^2$. After rescaling
by $h^{-2}$ we have a flow defined for backward time $1$.
Furthermore, the curvature of the rescaled flow is bounded on the
interval $(-1,0]$. Since the closure of $N_0$ in $N$ is compact, the
restriction of $h^{-2}g_N$ to $N_0\subset N$ at time $0$ is
uniformly bounded in the $C^\infty$-topology by Shi's theorem
(Theorem~\ref{shi}). The bound on the $k^{th}$-derivatives of the
curvature depends only on the curvature bound and hence can be taken
to be independent of $\delta'>0$ sufficiently small and also
independent of the strong $\delta'$-neck $N$. Gluing in the cap with
a  $C^\infty$-metric that converges smoothly to the standard initial
metric $g_0$ as $\delta'$ tends to zero using a fixed
$C^\infty$-partition of unity produces a family of manifolds
uniformly bounded in the $C^\infty$-topology.
\end{proof}

This leads immediately to:

\begin{cor}
Given a sequence of $\delta'_n\rightarrow 0$ and strong
$\delta'_n$-necks $(N(n),g_{N(n)})$ of scales $h_n$ and results of
surgery $(S_0(N(n)),g(n))$ with tips $p_n$ as in the previous claim,
then after passing to a subsequence there is a smooth limit
$(S_\infty,g_\infty,p_\infty)$ of a subsequence of the
$(S_0(N(n)),h_n^{-2}g_0(n)),p_n)$. This limit is the metric from Section~\ref{necksurgery}
that gives the standard initial conditions for a surgery cap.
\end{cor}

\begin{proof}
That there is a smooth limit of a subsequence is immediate from the
previous claim. Since the $\delta_n$ tend to zero, it is clear that the
limiting metric is the standard initial metric.
\end{proof}

\begin{lem}
Suppose that we have a sequence of $3$-dimensional Ricci flows
with surgeries $({\mathcal M}_n,G_n)$ that satisfy the strong
$(C,\epsilon)$-canonical neighborhood assumption with parameter
$r_{i+1}$, and have curvature pinched toward positive. Suppose that
there are surgeries in ${\mathcal M}_n$ at times $t_n$ with surgery
control parameters $\delta'_n$ and scales $h_n$. Let $p_n$ be the
tip of a surgery cap for the surgery at time $ t_n$. Also suppose
that there is $0\le\theta_n<1$ such that for every $A<\infty$,
for all $n$ sufficiently large
there are embeddings
$B(p_a, t_n,A h_n)\times [t_n, t_n+ h^2_n\theta_n)\to {\mathcal
M}_n$ compatible with time and the vector field. Suppose that
$\delta'_n\rightarrow 0$ and $\theta_n\rightarrow \theta<1$ as $n\rightarrow
\infty$.  Let $({\mathcal M}'_n,G'_n,p_n)$ be the Ricci flow with
surgery obtained by shifting time by $-t_n$ so that surgery occurs
at ${\bf t}=0$ and rescaling by $h_n^{-2}$ so that the scale of the
surgery becomes one. Then, after passing to a subsequence, the
sequence converges smoothly to a limiting flow
$(M_\infty,g_\infty(t),(p_\infty,0)),\ 0\le t<\theta$. This limiting flow is isomorphic to
the restriction of the standard flow to time $0\le t<\theta$.
\end{lem}

\begin{proof}
Let $Q<\infty$ be an upper bound for the scalar curvature of the
standard flow on the time interval $[0,\theta)$. Since
$\delta'_n\rightarrow 0$, according to the previous corollary, there is
a smooth limit at time $0$ for a subsequence, and this limit is the
standard initial metric. Suppose that, for some $0\le
\theta'<\theta$, we have established that there is a smooth limiting
flow on $[0,\theta']$. Since the initial conditions are the standard
solution, it follows from the uniqueness statement in
Theorem~\ref{stdsoln} that in fact the limiting flow is
isomorphic to the restriction of the standard flow to this time
interval. Then the scalar curvature of the limiting flow is bounded
by $Q$. Hence, for any $A<\infty$, for all $n$ sufficiently large,
the scalar curvature of the restriction of $G'_n$ to the image of
$B_{G'_n}(p_n,0,2A)\times[0,\theta']$ is bounded by $2Q$. According
to Lemma~\ref{controlnbhd} there is an $\eta>0$ and a constant
$Q'<\infty$, each depending only on $Q$, $r_{i+1}$, $C$ and
$\epsilon$, such that for all $n$ sufficiently large, the scalar
curvature of the restriction of $G_n'$ to $B_{G_n'}(p_n,0,A)\times
[0,{\rm min}(\theta'+\eta,\theta_n))$ is bounded by $Q'$. Because of
the fact that the curvature is pinched toward positive, this implies
that on the same set the sectional curvatures are uniformly bounded.
Hence, by  Shi's theorem with derivatives
(Theorem~\ref{shiw/deriv}), it follows that there are uniform bounds
for the curvature in the $C^\infty$-topology. Thus, passing to a
subsequence we can extend the smooth limit to the time interval
$[0,\theta'+\eta/2]$ unless $\theta'+\eta/2\ge \theta$. Since $\eta$
depends on $\theta$ (through $Q$), but is independent of $\theta'$,
we can repeat this process extending the time-interval of definition
of the limiting flow by $\eta/2$ until $\theta'+\eta/2\ge \theta$.
Now suppose that $\theta'+\eta/2\ge \theta$. Then the argument shows
that by passing to a subsequence we can extend the limit to any
compact subinterval of $[0,\theta)$. Taking a diagonal sequence
allows us to extend it to all of $[0,\theta)$. By the uniqueness of
the standard flow, this limit is the standard flow.
\end{proof}

\begin{cor}\label{goodtotheta}
With the notation and assumptions of the previous lemma, for all
$A<\infty$, and any $\delta''>0$, then  for all $n$ sufficiently
large, the restriction of $G_n'$ to the image
$B_{G_n'}(p_n,0,A)\times [0,\theta_n)$ is within $\delta''$ in the
$C^{[1/\delta'']}$-topology of the restriction of the standard
solution to the ball of radius $A$ about the tip for time $0\le
t<\theta_n$.
\end{cor}

\begin{proof}
Let $\eta>0$ depending on $\theta$ (though $Q$) as well as
$r_{i+1}$, $C$ and $\epsilon$ be as in the proof of the previous
lemma, and take $0<\eta'<\eta$. For all $n$ sufficiently large
$\theta_n>\theta-\eta'$, and consequently for all $n$ sufficiently
large there is an embedding $B_{G_n}(p_n,t_n,Ah_n)\times
[t_n,t_n+h_n^2(\theta-\eta')]$ into ${\mathcal M}_n$ compatible with
time and with the vector field. For all $n$ sufficiently large, we
consider the restriction of $G_n'$ to $B_{G_n'}(p_n,0,A)\times
[0,\theta-\eta']$. These converge smoothly to the restriction of the
standard flow to the ball of radius $A$ on the time interval
$[0,\theta-\eta']$. In particular, for all $n$ sufficiently large,
the restrictions to these time intervals are within $\delta''$ in
the $C^{[1/\delta'']}$-topology of the standard flow. Also, for all $n$ sufficiently
large, $\theta_n-(\theta-\eta')<\eta$. Thus, by
Lemma~\ref{controlnbhd}, we see that the scalar curvature of $G_n'$
is uniformly bounded (independent of $n$) on
$B_{G_n'}(p_n,0,A)\times [0,\theta_n)$. By the assumption that the
curvature is pinched toward positive, this means that the sectional
curvatures of the $G_n'$ are also uniformly bounded on these sets,
and hence so are the Ricci curvatures. (Notice that these bounds are
independent of $\eta'>0$.) By Shi's theorem with derivatives
(Theorem~\ref{shiw/deriv}), we see that there are uniform bounds on
the curvatures in the $C^\infty$-topology on these subsets, and
hence bounds in the $C^\infty$-topology on the Ricci curvature.
These bounds are independent of both $n$ and $\eta'$. Thus, choosing
$\eta'$ sufficiently close to zero, so that $\theta_n-\eta'$ is also
close to $\theta$ for all $n$ sufficiently large, we see that for all such large $n$
and all
$t\in [\theta-\eta',\theta)$,  the
restriction of $G_n'$ to $B_{G_n'}(p_n,0,A)\times \{t\}$ is
arbitrarily close in the $C^{[1/\delta'']}$-topology to
$G_n'(\theta-\eta')$. The same is of course true of the standard
flow. This completes the proof of the corollary.
\end{proof}

Now we turn to the proof proper of Proposition~\ref{altern}. We fix
$A<\infty$, $\delta''>0$ and $\theta<1$. We are free to make $A$
larger so we can assume by Proposition~\ref{item5} that for the
standard flow the restriction of the flow to $B(p_0,0,A)\setminus
B(p_0,0,A/2)$ remains close to a standard evolving $S^2\times
[A/2,A]$ for time $[0,\theta]$. Let $K<\infty$ be a constant with
the property that $R(x,t)\le K$ for all $x\in B(p_0,0,A)$ in the
standard flow and all $t\in [0,\theta]$. If there is no
$\delta''_0>0$ as required, then we can find a sequence
$\delta'_n\rightarrow 0$ as $n\rightarrow\infty$ and Ricci flows with
surgery $({\mathcal M}_n,G_n)$ with surgeries at time $t_n$ with
surgery control parameter $\delta_n(t_n)\le \delta_n'$  and surgery
scale parameter $ h_n=h(r_{i+1}\delta_n( t_n),\delta_n(t_n))$
satisfying the hypothesis of the lemma but not the conclusion. Let
$T'_n$ be the final time of $({\mathcal M}_n,G_n)$. Let $\theta_n\le
\theta$ be maximal subject to the condition that there is an
embedding $\rho_n\colon B_{G_n}(x, t_n,Ah_n)\times [ t_n,
t_n+h_n^2\theta_n)\to {\mathcal M}_n$ compatible with time and the
vector field. Let $G'_n$ be the result of shifting the time by
$-t_n$ and scaling the result by $h_n^{-2}$. According to
Lemma~\ref{goodtotheta}, for all $n$ sufficiently large, the
restriction of $G_n'$ to the image of $\rho_n$ is within $\delta''$
in the $C^{[1/\delta'']}$-topology of the standard flow restricted
to the ball of radius $A$ about the tip of the standard solution on
the time interval $[0,\theta_n)$. If $\theta_n={\rm
min}(\theta,(T_n'-t_n)/h_n^2)$, then the first conclusion of
Proposition~\ref{altern} holds for $({\mathcal M}_n,G_n)$ for all
$n$ sufficiently large which contradicts our assumption that the
conclusion of this proposition holds for none of the $({\mathcal
M}_n,G_n)$. If on the other hand $\theta_n<{\rm
min}(\theta,(T_n'-t_n)/h_n^2)$, we need only show that all of
$B(x_n,t_n,Ah_n)$ disappears at time $t_n+h_n^2\theta_n$ in order to
show that the second conclusion of Proposition~\ref{altern} holds
provided that $n$ is sufficiently large. Again this would contradict
the fact that the conclusion of this proposition holds for none of
the $({\mathcal M}_n,G_n)$.

So now let us suppose that $\theta_n<{\rm
min}(\theta,(T'_n-t_n)/h_n^2)$. Since there is no further extension
in forward time for $B(p_n,t_n,A h_n)$, it must be the case that $
t_n+ h_n^2\theta_n$ is a surgery time and there is some flow line
starting at a  point of $B(p_n,t_n,A h_n)$ that does not continue to
time $t_n+h_n^2\theta_n$. It remains to show that in this case that
for any $t< t_n+ h_n^2\theta_n$ sufficiently close to $t_n+
h_n^2\theta_n$ we have $\rho_n\left(B_{G_n}(x, t_n,Ah_n)\times
\{t\}\right)\subset D_t$, the region in $M_t$ that disappears at
time $t_n+h_n^2\theta_n$.

\begin{claim}
Suppose that $\theta_n<{\rm min}(\theta,(T'_n- t_n)/ h_n^2)$. Let
$\Sigma_{1},\ldots,\Sigma_{k}$ be the $2$-spheres along which
we do surgery at time $ t_n+h_n^2\theta_n$. Then for any $t< t_n+
h_n^2\theta_n$ sufficiently close to $ t_n+ h_n^2\theta_n$ the
following holds provided that $\delta_n'$ is sufficiently small. The image
$$\rho_n\left(B_{g_n}(x,t_n,A h_n)\times
\{t\}\right)$$ is disjoint from the images $\{\Sigma_{i}(t)\}$ of the $\{\Sigma_{i}\}$
under the
backward flow to time $t$ of the spheres $\Sigma_{i}$ along which
we do surgery at time $t_n+h_n^2\theta_n$.
\end{claim}

\begin{proof}
There is a constant $K'<\infty$ depending on $\theta$ such that for
the standard flow we have $R(x,t)\le K'$ for all $x\in B(p_0,0,A)$
and all $t\in [0,\theta)$ for the standard solution. Consider the
embedding $\rho_n\left(B(p_n,t_n,Ah_n)\times [t_n,t_n+h_n^2
\theta_n)\right)$. After time shifting by $-t_n$ and rescaling by
$h_n^{-2}$, the flow $G_n'$ on the image of $\rho_n$ is within
$\delta''$ of the standard flow. Thus, we see that for all $n$ sufficiently
large and for every point $x$ in the image of $\rho_n$ we have $
R_{G_n'}(x)\le 2K'$ and hence $R_{G_n}(x)\le 2K'h_n^{-2}$.

Let $h_n'$ be the scale of the surgery at time $t_n+ h_n^2\theta_n$.
(Recall that $h_n$ is the scale of the surgery at time $t_n$.)
Suppose that $\rho_n(B(p_n,t_n,Ah_n)\times \{t'\})$ meets one of the
surgery $2$-spheres $\Sigma_{i}(t')$ at time $t'$ at a point
$y(t')$. Then, for all $t\in [t',t_n+h_n^2\theta_n)$ we have
 the image $y(t)$ of $y(t')$ under the flow. All these points
$y(t)$ are points of intersection of $\rho_n(B(p,t_n,Ah_n)\times
\{t\})$ with $\Sigma_{i}(t)$. Since $y(t)\in \rho_n(B(p,t_n,A
h_n)\times\{t\})$, we have $ R(y(t))\le 2K' h_n^{-2}$. On the other
hand $R(y(t)) (h'_n)^2$ is within $O(\delta)$ of  $1$ as $t$ tends
to $t_n+h_n^2\theta_n$. This means that $h_n/h_n'\le \sqrt{3K'}$ for
all $n$ sufficiently large. Since the standard solution has
non-negative curvature, the metric is a decreasing function of $t$,
and hence the diameter of $B(p_0,t,A)$ is at most $ 2A$ in the standard
solution. Using Corollary~\ref{goodtotheta} we see that  for all $n$
sufficiently large, the diameter of $\rho_n\left(B(p,t_n,Ah_n)\times
\{t\}\right)$ is at most $Ah_n\le 4\sqrt{K'} Ah_n'$.
 This
means that for $\delta'_n$ sufficiently small the distance at time
$t$ from $\Sigma_{i}(t)$ to the complement of the $t$ time-slice
of the strong $\delta_n( t_n+ h_n^2\theta_n)$-neck $N_{i}(t)$
centered at $\Sigma_{i}(t)$ (which is at least
$(\delta'_n)^{-1}h_n'/2$) is much larger than the diameter of
$$\rho_n(B(p_n,t_n,Ah_n)\times \{t\}).$$  Consequently,  for all $n$
sufficiently large, the image $\rho_n(B(p_n,t_n,Ah_n)\times \{t\})$
is contained in  $N_{i}(t)$. But by our choice of $A$, and
Corollary~\ref{goodtotheta} there is an $\epsilon$-neck of rescaled
diameter approximately $Ah_n/2$ contained in $\rho_n(B(p_n,t_n,A
h_n)\times \{t\})$. By Corollary~\ref{s2isotopic} the spheres coming
from the neck structure in $$\rho_n(B(p_n,t_n,Ah_n)\times \{t\})$$ are
 isotopic in $N_{i}(t)$ to the central $2$-sphere of this neck.
This is a contradiction because in $N_{i}(t)$ the central
$2$-sphere is homotopically non-trivial whereas  the spheres in
$\rho_n(B(p_n,t_n,A\bar h_n)\times \{t\})$  clearly bound
$3$-disks.
\end{proof}

Since $\rho_n(B(p_n,t_n,Ah_n)\times\{t\})$ is disjoint from the
backward flow to time $t$ of all the surgery $2$-spheres
$\Sigma_{i}(t)$ and since $\rho_n(B(p_n,t_n,Ah_n)\times \{t\})$ is
connected,
 if there is a flow line starting at some point $z\in
B(p,t_n,Ah_n)$ that disappears at time $t_n+ h_n^2\theta_n$, then
the flow from every point of $B(p,t_n,Ah_n)$ disappears at time $
t_n+ h_n^2\theta_n$. This shows that if $\theta_n<{\rm
min}(\theta,T'_n- t_n/h_n^{2})$, and if there is no extension of
$\rho_n$ to an embedding defined at time $t_n+h^2_n\theta_n$, then
all forward flow lines beginning at points of $B(p,t_n,A h_n)$
disappear at time $ t_n+h_n^2\theta_n$, which of course means that
for all $t<t_n+h^2_n\theta_n$ sufficiently close to
$t_n+h_n^2\theta_n$ the entire image $\rho_n(B(p,t_n,Ah_n)\times
\{t\})$ is contained in the disappearing region $D_t$. This shows
that for all $n$ sufficiently large, the second conclusion of
Proposition~\ref{altern} holds,  giving a contradiction.

This completes the proof of Proposition~\ref{altern}.
\end{proof}

\begin{rem}
Notice that it is indeed possible that $B_{G}(x,t,Ah)$ is removed at
some later time, for example as part of a capped  $\epsilon$-horn
associated to some later surgery.
\end{rem}

\section{A length estimate}

We use the result in the previous section about the evolution of
surgery caps to establish the length estimate on paths parameterized
by backward time approaching a surgery cap from above.

\begin{defn}\label{defnoverlinedelta0}
Let $c>0$ be the constant from Proposition~\ref{Restim}. Fix
$0<\overline\delta_0<1/4$ such that if $g$ is within
$\overline\delta_0$ of $g_0$ in the $C^{[1/\overline
\delta]}$-topology then $|R_{g'}(x)-R_{g_0}(x)|<c/2$ and $|{\rm
Ric}_{g'}-{\rm Ric}_{g_0}|<1/4$.
\end{defn}

Here is the length estimate.

\begin{figure}[ht]
  \centerline{\epsfbox{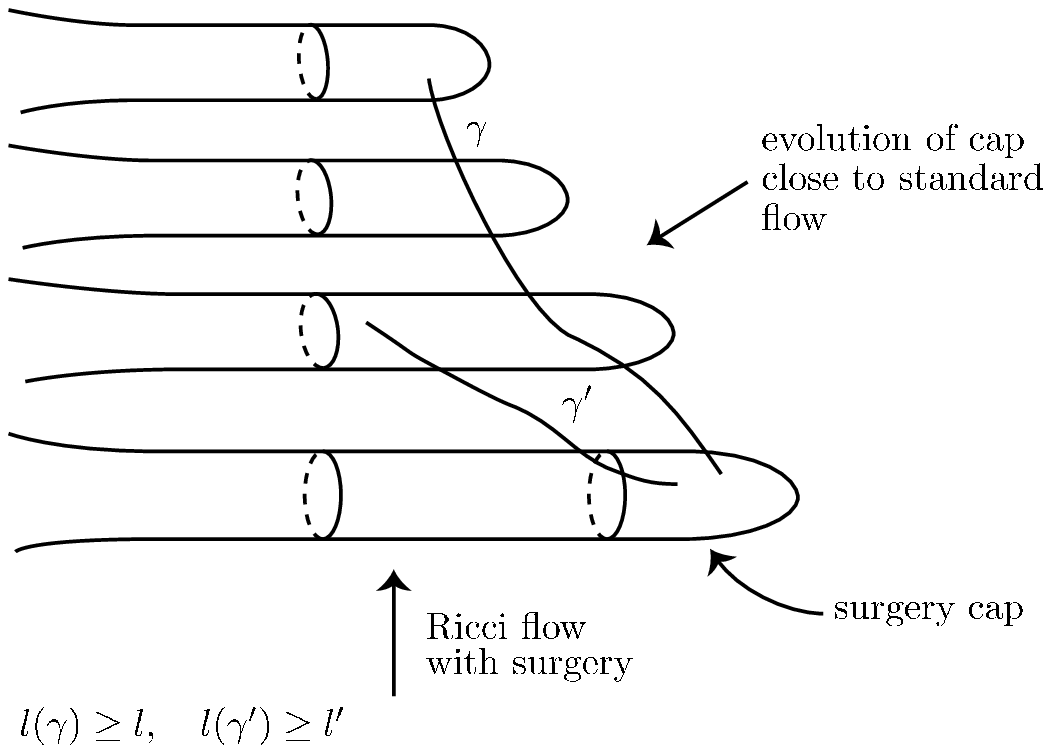}}
  \caption{Paths in evolving surgery caps are long.}\label{fig:long}
\end{figure}

\begin{prop}\label{bigell}
For any $\ell<\infty$ there is $A_0=A_0(\ell)<\infty$,
$0<\theta_0=\theta_0(\ell)<1$, and for any $A\ge A_0$ for the
constant
$\delta''=\delta''(A)=\delta''_0(A,\theta_0,\overline\delta_0)>0$
from Proposition~\ref{altern} the following holds. Suppose that
$({\mathcal M},G)$ is a Ricci flow with surgery defined for $0\le
t< T<\infty$. Suppose that it satisfies the strong
$(C,\epsilon)$-canonical neighborhood assumption at all points $x$
with $R(x)\ge r_{i+1}^{-2}$. Suppose also that the solution has
curvature pinched toward positive. Suppose that there is a surgery
at some time $\bar t$ with $T_{i-1}\le \bar t<T$ with
$\overline\delta(\bar t)$ as the surgery control parameter and with
$h$ as the surgery scale parameter. Then the following holds
provided that $\overline\delta(\bar t)\le \delta''$. Set $T'={\rm
min}(T,\bar t+h^2\theta_0 )$. Let $p\in M_{\bar t}$ be the tip of
the cap of a surgery disk at time $\bar t$. Suppose that $P(p,\bar
t,Ah,T'-\bar t)$ exists in ${\mathcal M}$. Suppose that we have
$t'\in [\bar t,\bar t+h^2/2]$ with $t'\le T'$, and suppose that we
have a curve $\gamma(\tau)$ parameterized by backward time $\tau\in
[0,T'- t']$ so that $\gamma(\tau)\in M_{T'-\tau}$ for all $\tau\in
[0,T'- t']$. Suppose that the image of $\gamma$ is contained in the
closure of $P(p,\bar t,Ah,T'-\bar t)\subset {\mathcal M}$. Suppose
further:
\begin{enumerate}
\item[(1)] either that $T'=\bar t +\theta_0 h^2\le T$ or
that $\gamma(0)\subset \partial B(p,\bar t,Ah)\times\{T'\}$; and
\item[(2)]  $\gamma(T'-t')\in B(p,\bar t,Ah/2)\times{t'}$.
\end{enumerate}
Then
$$\int_0^{T'- t'}\left(R(\gamma(t))+|X_\gamma(t)|^2\right)dt>\ell.$$
See {\sc Fig.}~\ref{fig:long}.
\end{prop}

\begin{proof}
The logic of the proof is as follows. We
fix $\ell<\infty$. We shall determine the relevant value of
$\theta_0$ and then of $A_0$ in the course of the argument. Then for any $A\ge A_0$
we define
 $\delta''(A)= \delta''_0(A,\theta_0,\overline\delta_0)$, as in
Proposition~\ref{altern}.

The integral expression is invariant under time translation and also
under rescaling. Thus, we can (and do) assume that $\bar t=0$ and
that the scale $h$ of the surgery is $1$. We use the embedding of
$P(p,0,A,T')\to {\mathcal M}$ and write the restriction of the flow
to this subset as a one-parameter family of metrics $g(t),\ 0\le
t\le T'$, on $B(p,0,A)$. With this renormalization, $0\le t'\le
1/2$, also $T'\le \theta_0$,  and $\tau=T'-t$.

Let us first consider the case when  $T'=\theta_0\le T$. Consider the
standard flow $(\Ar^3,g_0(t))$, and let $p_0$ be its tip. According
to Proposition~\ref{Restim},  for all $x\in \Ar^3$ and all $t\in
[0,1)$ we have $R_{g_0}(x,t)\ge c/(1-t)$. By Lemma~\ref{altern} and
since we are assuming that $\overline\delta(\bar t)\le \delta''=
\delta''_0(A,\theta_0,\overline\delta_0)$, we have that $R(a,t)\ge
c/2(1-t)$ for all $a\in B(p,0,A)$ and all $t\in [0,\theta]$. Thus,
 we have
\begin{eqnarray*}
\int_0^{\theta_0-t'}\left(R(\gamma(\tau))+|X_\gamma(\tau)|^2\right)d\tau
& \ge & \int_{t'}^{\theta_0}
\frac{c}{2(1-t)}dt \\
& = & \frac{-c}{2}\left({\rm log}(1-\theta_0)-{\rm
log}(1-t')\right)dt \\
& \ge &  \frac{-c}{2}\left({\rm log}(1-\theta_0)+{\rm
log}(2)\right).
\end{eqnarray*}
 Hence, if
$\theta_0<1$ sufficiently close to $1$, the  integral will be $> \ell$.
This fixes the value of $\theta_0$.

\begin{claim}
There is $A'_0<\infty$ with the property that for any $A\ge A'_0$
the restriction of the standard solution $g_0(t)$ to
$\left(B(p_0,0,A)\setminus B(p_0,0,A/2)\right)\times [0,\theta_0]$
is close to an evolving family $(S^2\times [A/2,A],h_0(t)\times
ds^2)$. In particular, for any $t\in [0,\theta_0]$, the
$g_0$-distance at time $t$ from $B(p_0,0,A/2)$ to the complement of
$B(p_0,0,A)$ in the standard solution is more than $A/4$.
\end{claim}

\begin{proof}
This is immediate from Proposition~\ref{item5} and the fact that
$\theta_0<1$.
\end{proof}

Now fix $A_0={\rm max}(A_0', 10\sqrt{\ell})$ and let $A\ge A_0$.

Since $\overline\delta_0<1/4$ and since $T'\le \theta_0$, for $\overline\delta(\bar t)\le
\delta''_0(A,\theta_0,\overline\delta_0)$
by Proposition~\ref{altern}
the $g(T')$-distance  between $B(p,0,A/2)$ and $\partial
B(p,0,A)$  is at least $A/5$.

 Since the
flow on $B(p,0,A)\times[0,T']$ is within $\overline\delta_0$ of
the standard solution, and since the curvature of the standard
solution is non-negative, for any horizontal tangent vector $X$ at
any point of $B(p,0,A)\times [0,T']$ we have that
$${\rm Ric}_g(X,X)\ge -\frac{1}{4}|X|_{g_0}^2\ge -\frac{1}{2}|X|^2_g,$$ and hence
$$\frac{d}{dt}|X|_g^2\le |X|_g^2.$$
Because  $T'\le 1$, we see that
$$|X|^2_{g(T')}\le e\cdot|X|^2_{g(t)}<3|X|^2_{g(t)}$$
for any $t\in [0,T']$.

Now suppose that $\gamma(0)\in \partial B(p,0,A)\times\{T'\}$. Since
the image of $\gamma$ is contained in the closure of $P(p,0,A,T')$
for every $\tau\in [0,T']$ we have
$\sqrt{3}|X_\gamma(\tau)|_{g(T'-\tau)}\ge |X_\gamma(\tau)|_{g(T')}$.
Since the flow $g(t)$ on $P(p,0,A, T')$ is within
$\overline\delta_0$ in the $C^{[1/\overline\delta_0]}$-topology of
the standard flow on the corresponding parabolic neighborhood,
$R(\gamma(t))\ge 0$ for all $t\in [0,T']$. Thus, because of these
two estimates we have
\begin{equation}\label{*}
\int_0^{T'-t'}\left(R(\gamma(\tau))+|X_\gamma(\tau)|^2\right)d\tau\ge
\int_0^{T'-t'}
\frac{1}{3}|X_\gamma(\tau)|_{g(T')}^2d\tau.\end{equation}
 Since $\gamma(0)\in
\partial B(p,0,A)\times \{T'\}$ and $\gamma(T')\in B(p,0,A/2)$, it follows from
Cauchy-Schwarz that
\begin{eqnarray*}(T'-t')^2\int_0^{T'}|X_\gamma(\tau)|_{g(T')}^2d\tau &
\ge & \left(\int_0^{T'-t'}|X_\gamma(\tau)|_{g(T')}d\tau\right)^2
\\
& \ge & \left(d_{g(T')}(B(p,0,A/2),\partial
B(p,0,A))\right)^2\ge \frac{A^2}{25}. \end{eqnarray*}
 Since $T'-t'<1$, it
immediately follows from this and Equation~(\ref{*}) that
$$\int_0^{T'-t'}\left(R(\gamma(\tau))+|X_\gamma(\tau)|^2\right)d\tau\ge
\frac{A^2}{75}.$$ Since $A\ge A_0\ge 10\sqrt{ \ell}$, this
expression is $>\ell$. This completes the proof of
Proposition~\ref{bigell}
\end{proof}

\subsection{Paths with short ${\mathcal L}_+$-length avoid the surgery caps}

Here we show that a path parameterized by backward time that ends in
a surgery cap (or comes close to it)  must have
long ${\mathcal L}$-length. Let $({\mathcal M},G)$ be a Ricci flow
with surgery, and let $x\in {\mathcal M}$ be a point with ${\bf
t}(x)=T\in (T_i,T_{i+1}]$. We suppose that these data satisfy the
hypothesis of Proposition~\ref{delta0ri+1} with respect to the given
sequences and $r\ge r_{i+1}>0$. In particular, the parabolic
neighborhood $P(x,T,r,-r^2)$ exists in ${\mathcal M}$ and $|{\rm Rm}|$ is
bounded on this parabolic neighborhood by $r^{-2}$.

Actually, here we do not work directly with the length function
${\mathcal L}$ defined from $x$, but rather with a closely related
function. We set $R_+(y)={\rm max}(R(y),0)$.

\begin{figure}[ht]
  \centerline{\epsfbox{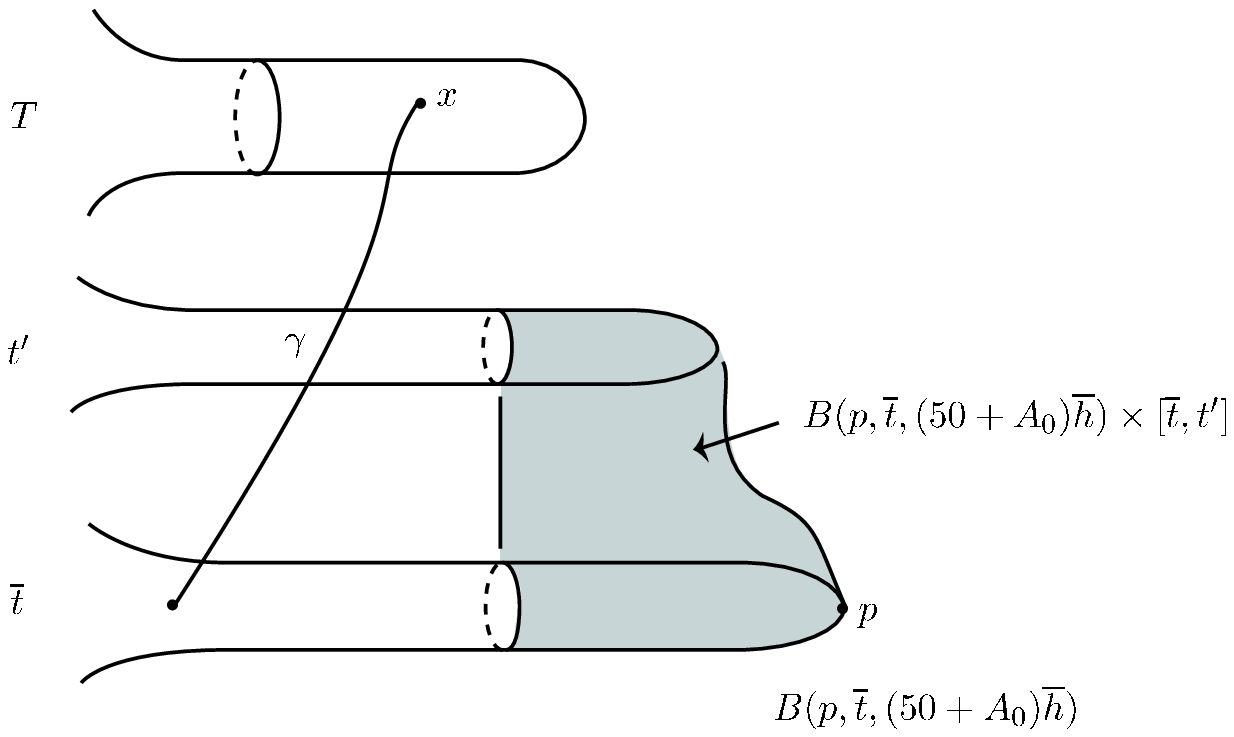}}
  \caption{Avoiding neighborhoods of surgery caps}\label{fig:avoid}
\end{figure}

\begin{lem}\label{L_0}
Given $L_0<\infty$, there is
$\overline\delta_1=\overline\delta_1(L_0,r_{i+1})>0$, independent of
$({\mathcal M},G)$ and $x$, such that if $\overline\delta(t)\le
\overline\delta_1$ for all $t\in [T_{i-1},T)$, then for any curve
$\gamma(\tau),\ 0\le \tau\le \tau_0$, with $\tau_0\le T-T_{i-1}$,
parameterized by backward time with $\gamma(0)=x$ and with
$${\mathcal L}_+(\gamma)=\int_0^{\tau_0}\sqrt{\tau}
\left(R_+(\gamma(\tau))+|X_\gamma|^2\right)d\tau< L_0$$\index{${\mathcal L}_+$|ii} the
following two statements hold: \begin{enumerate} \item[(2)] Set
$$\tau'={\rm min}\left(\frac{r_{i+1}^4}{(256)L_0^2},{\rm ln}(\root 3\of{2})r_{i+1}^2\right).$$
Then for all $\tau\le {\rm
min}(\tau',\tau_0)$ we have $\gamma(\tau)\in P(x,T,r/2,-r^2)$.
\item[(2)] Suppose that  $\bar t\in [T-\tau_0,T)$ is a surgery time with $p$
being the tip of the surgery cap at time $\bar t$ and with the scale
of the surgery being $\bar h$. Suppose   $t'\in [\bar t,\bar t+\bar
h^2/2]$ is such that there is an embedding
$$\rho\colon B(p,\bar t,(50+A_0)\bar h)\times [\bar t,t']\to {\mathcal M}$$
compatible with time and the vector field. Then the image of $\rho$
is disjoint from the image of $\gamma$. See {\sc Fig.}~\ref{fig:avoid}.
\end{enumerate}
\end{lem}

\begin{rem}
Recall that $(A_0+4)\bar h$ is the radius of the surgery cap (measured in the
rescaled version of the standard initial metric) that is glued in when
performing surgery with scale $\bar h$.
\end{rem}

\begin{proof}
We define $\ell=L_0/\sqrt{\tau'}$, then define $A={\rm
max}(A_0(\ell),2(50+A_0))$ and $\theta=\theta_0(\ell)$. Here,
$A_0(\ell)$ and $\theta_0(\ell)$ are the constants in
Proposition~\ref{bigell}.
  Lastly, we require $\overline\delta_1\le \delta''(A)$ from
Proposition~\ref{bigell}. Notice that, by construction,
 $\delta''(A)=\delta''_0(A,\theta,\overline\delta_0)$ from
Proposition~\ref{altern}. Thus, if $p$ is the tip of a surgery cap at time $\bar t$
with the scale of the surgery being $\bar h$, then it follows that
for any $\Delta t\le \theta$, if there is an embedding
$$\rho\colon B(p,\bar t,A\bar h)\times [\bar t,\bar t+\bar h^2\Delta t)\to {\mathcal M}$$
 compatible with time and the vector field, then the induced flow (after time shifting by $-\bar t$
 and scaling by $(\bar h)^{-2}$ is within $\overline\delta_0$ in the
 $C^{[1/\overline\delta_0]}$-topology of the standard solution.
 In particular, the scalar curvature
at any point of the image of $\rho$  is positive and is within a
multiplicative factor of two of the scalar curvature at the
corresponding point of the standard flow.

 Recall that we have $r\ge r_{i+1}$ and
that $P(x,T,r,-r^2)$ exists in ${\mathcal M}$ and that $|{\rm Rm}|\le
r^{-2}$ on this parabolic neighborhood.  We begin by proving by
contradiction that there is no $\tau\le  \tau'$ with the property
that $\gamma(\tau)\not\in P(x,T,r/2,-r^2)$. Suppose there is such a
$\tau\le \tau'$. Notice that by construction
$\tau'<r_{i+1}^2<r^2$. Hence,  for the first $\tau''$ with the
property that $\gamma(\tau'')\not\in P(x,T,r/2,-r^2)$ the point
$\gamma(\tau'')\in
\partial B(x,T,r/2)\times \{T-\tau''\}$.

\begin{claim}
$\int_0^{\tau''}|X_\gamma(\tau)|d\tau> r/2\sqrt{2}$.
\end{claim}

\begin{proof}
Since $|{\rm Rm}|\le r^{-2}$ on $P(x,T,r,-r^2)$, we have $|{\rm Ric}|\le
2r^{-2}$ on $P(x,T,r,-\tau'')$. Thus, for any tangent vector $v$ at
a point of $B(x,T,r)$ we have
$$\left|\frac{d(\langle v,v\rangle_{G(T-\tau)})}{d\tau}\right|\le 2r^{-2}\langle v,v\rangle_{G(T-\tau)}$$
for all
$\tau\in [0,\tau'']$. Integrating gives that for any $\tau\le
\tau''$ we have
$${\rm
exp}(-2r^{-2}\tau'')\langle v,v\rangle_{G(T)}\le \langle
v,v\rangle_{G(T-\tau)} \le {\rm exp}(2r^{-2}\tau'')\langle
v,v\rangle_{G(T)}.$$ Since  $\tau''\le \tau'$ and $r\ge r_{i+1}$ by
the assumption on $\tau'$ we have
$${\rm exp}(2r^{-2}\tau'')\le {\rm exp}(2\root 3\of{2})< 2.$$ This
implies that for all $\tau\le \tau''$ we have
$$\frac{1}{\sqrt{2}}| X_\gamma(\tau)|_{G(T)}<
|X_\gamma(\tau)|_{G(T-\tau)}<\sqrt{2}|X_\gamma(\tau)|_{G(T)},$$
and hence
$$\int_0^{\tau''}|X_\gamma(\tau)|d\tau>
\frac{1}{\sqrt{2}}\int_0^{\tau''}|X_\gamma(\tau)|_{G(T)}\ge
\frac{r}{2\sqrt{2}},$$ where we use the fact that
$d_T(\gamma(0),\gamma(\tau''))=r/2$.
\end{proof}

Applying Cauchy-Schwarz to $\tau^{1/4}|X_\gamma|$ and $\tau^{-1/4}$ on the
interval $[0,\tau'']$ yields
\begin{eqnarray*}\int_0^{\tau''}\sqrt{\tau}\left(R_+(\gamma(\tau))+|X_\gamma(\tau)|^2\right)d\tau &
\ge & \int_0^{\tau''}\sqrt{\tau}|X_\gamma(\tau)|^2d\tau
\\
& \ge &
\frac{\left(\int_0^{\tau''}|X_\gamma(\tau)|d\tau\right)^2}{\int_0^{\tau''}\tau^{-1/2}d\tau}
\\ & > &
\frac{r^2}{16\sqrt{\tau''}}\ge L_0.
\end{eqnarray*}
Of course, the integral from $0$ to $\tau''$ is less than or equal
the entire integral from $0$ to $\tau_0$ since the integrand is
non-negative, contradicting the assumption that ${\mathcal
L}_+(\gamma)\le L_0$. This completes the proof of the first
numbered statement.

We turn now to the second statement. We impose a further condition
on $\overline\delta_1$. Namely, require that
$\overline\delta_1^2<r_{i+1}/2$. Since $r_{i}\le r_0\le \epsilon<1$,
we have $\overline\delta_1^2r_i<r_{i+1}/2$. Thus, the scale of the
surgery, $\bar h$, which is $\le \overline\delta_1^2r_i$ by
definition, will also be less than $r_{i+1}/2$, and hence there is
no point of $P(x,T,r,-r^2)$ (where the curvature is bounded by
$r^{-2}\le r_{i+1}^{-2}$) in the image of $\rho$ (where the scalar
curvature is greater than $(\bar h)^{-2}/2>2r_{i+1}^{-2}$). Thus, if
$\tau'\ge \tau_0$ we have completed the proof. Suppose that
$\tau'<\tau_0$. It suffices to establish that for every $\tau_1\in
[\tau',\tau_0]$ the point $\gamma(\tau_1)$ is not contained in the
image of $\rho$ for any surgery cap and any $t'$ as in the
statement. Suppose that in fact there is $\tau_1\in [\tau', \tau_0]$
with $\gamma(\tau_1)$ contained in the image of $\rho(B(p,\bar
t,(A_0+50)\bar h)\times [\bar t,t'])$ where $\bar t\le t'\le \bar
t+(\bar h)^2/2$ and where $p$ is the tip
 of some surgery cap at time $\bar t$. We estimate
\begin{eqnarray}\nonumber
\lefteqn{\int_0^{\tau_0}\sqrt{\tau}\left(R_+(\gamma(\tau))+|X_\gamma(\tau)|^2\right)d\tau} & & \\
& \ge &
\int_{\tau'}^{\tau_0}\sqrt{\tau}\left(R_+(\gamma(\tau))+|X_\gamma(\tau)|^2\right)d\tau
\nonumber
\\ & \ge &
\sqrt{\tau'}\int_{\tau'}^{\tau_1}\left(R_+(\gamma(\tau))+|X_\gamma(\tau)|^2\right)d\tau.
\label{tau_1est}
\end{eqnarray}

 Let $\Delta t\le
T-\bar t$ be the supremum of the set of $s$ for which there is a
parabolic neighborhood $P(p,\bar t,A\bar h,s)$ embedded in ${\bf
t}^{-1}((-\infty,T])\subset {\mathcal M}$. Let $\Delta t_1={\rm
min}(\theta \bar h^2,\Delta t)$. We consider $P(p,\bar t,A\bar
h,\Delta t_1)$. First, notice that since $\bar h\le
\overline\delta_1^2r_{i}<r_{i+1}/2$,
the scalar curvature on $P(p,\bar
t,A\bar h,\Delta t_1)$ is larger than $(\bar
h)^{-2}/2>r_{i+1}^{-2}\ge r^{-2}$. In particular, the parabolic
neighborhood $P(x,T,r,-r^2)$ is disjoint from $P(p,\bar t,A\bar
h,\Delta t_1)$. This means that there is some $\tau''\ge \tau'$ such
that $\gamma(\tau'')\in
\partial P(p,\bar t,A\bar h,\Delta t_1)$ and $\gamma|_{[\tau'',\tau_1]}\subset
P(p,\bar t,A\bar h,\Delta t_1)$. There are two cases to consider.
The first is when $\Delta t_1=\theta \bar h^2$, $\tau''=T-(\bar
t+\Delta t_1)$ and $\gamma(\tau'')\in B(p,\bar t,A\bar
h)\times\{\bar t+\Delta t_1\}$. Then, according to
Proposition~\ref{bigell},
\begin{equation}\label{**}
\int_{\tau''}^{\tau_1}R_+(\gamma(\tau))d\tau> \ell.\end{equation}

Now let us consider the other case. If $\Delta t_1<\theta \bar h^2$,
this means that either $\bar t+\Delta t_1=T$ or, according to
Proposition~\ref{altern}, at the time $\bar t+\Delta t_1$ there is a
surgery that removes all of $B(p,\bar t,A\bar h)$. Hence, under
either possibility it must be the case that $\gamma(\tau'')\in
\partial B(p,\bar t,A\bar h)\times \{T-\tau''\}$. Thus, the remaining
case to consider is when, whatever $\Delta t_1$ is,
$\gamma(\tau'')\subset
\partial B(p,\bar t,A\bar h)\times \{T-\tau''\}$.
Lemma~\ref{bigell} and the fact that $R\ge 0$ on $P(p,\bar t,A\bar
h,\Delta t_1)$ imply that
$$\ell<\int_{\tau''}^{\tau_1}\left(R(\gamma(\tau))+|X_\gamma(\tau)|^2\right)d\tau =
\int_{\tau''}^{\tau_1}\left(R_+(\gamma(\tau))+|X_\gamma(\tau)|^2\right)d\tau.$$

Since $\ell=L_0/\sqrt{\tau'}$ and $\tau''\ge \tau'$, it follows from
Equation~(\ref{tau_1est}) that in both
cases
$${\mathcal L}_+(\gamma)\ge
\int_{\tau''}^{\tau_1}\sqrt{\tau}\left(R_+(\gamma(\tau))+|X_\gamma(\tau)|^2\right)d\tau>
\ell\sqrt{\tau'}=L_0,$$ which contradicts our hypothesis. This
completes the proof of Lemma~\ref{L_0}.
\end{proof}

\subsection{Paths with small energy avoid
 the disappearing regions}

At this point we have shown that paths of small energy do not approach
the surgery caps from above. We also need to rule out that they can be arbitrarily
close from below. That is to say,
 we need to see that paths whose ${\mathcal L}$-length is not
too large avoid neighborhoods of the disappearing regions at all
times just before the surgery time at which they disappear. Unlike
the previous estimates which were universal for all $({\mathcal
M},G)$ satisfying the hypothesis of Proposition~\ref{delta0ri+1}, in
this case the estimates will depend on the Ricci flow with surgery.
First, let us fix some notation.

\begin{figure}[ht]
  \centerline{\epsfbox{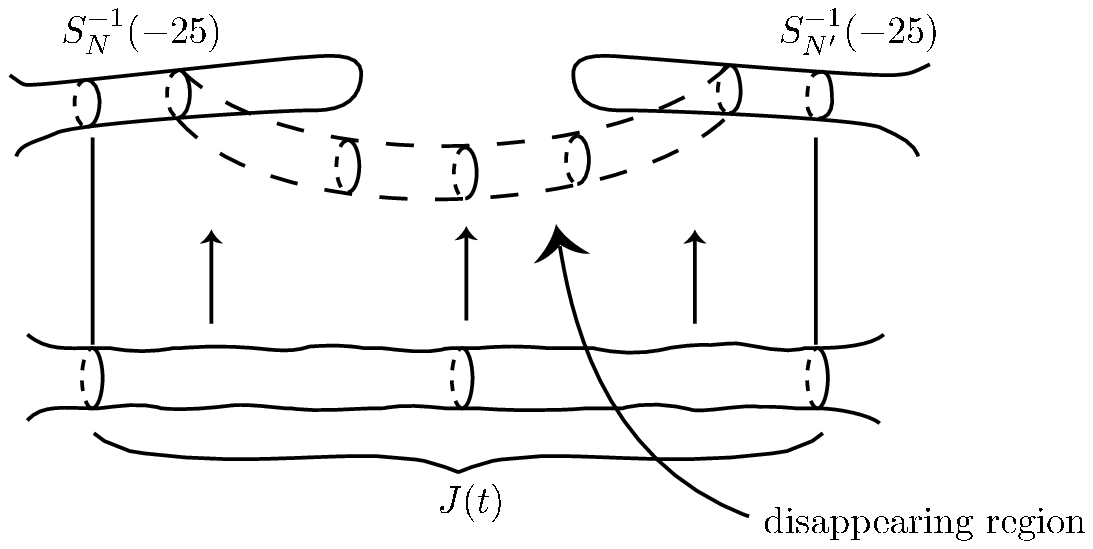}}
  \caption{Paths of short length avoid disappearing regions.}\label{fig:disappear}
\end{figure}

\begin{defn}\label{J}
 Suppose that $\bar t$ is a surgery
time, that $\tau_1>0$,  and that there are no other surgery times in
the interval $(\bar t-\tau_1,\bar t]$. Let $\{\Sigma_i(\bar t)\}_i$
be the $2$-spheres on which we do surgery at time $\bar t$. Each
$\Sigma_i$ is the central $2$-sphere of a strong $\delta$-neck
$N_i$. We can flow the cylinders $J_0(\bar
t)=\cup_is_{N_i}^{-1}(-25,0])$ backward to any time
 $t\in (\bar t-\tau_1,\bar t]$. Let $J_0(t)$ be
the result. There is an induced function, denoted
$\coprod_is_{N_i}(t)$, on $J_0(t)$. It takes values in $ (-25,0]$.
We denote the boundary of $J_0(t)$ by $\coprod_i\Sigma_i(t)$. Of
course, this boundary is the result of flowing
$\coprod_i\Sigma_i(\bar t)$ backward to time $t$. (These backward
flows are possible since there are no surgery times in $(\bar
t-\tau_1,\bar t)$.) For each $t\in [\bar t-\tau_1,\bar t)$ we also
have the disappearing region $D_t$: -- the region that disappears at
time $\bar t$. It is an open submanifold whose boundary is
$\coprod_i\Sigma_i(t)$. Thus, for every $t\in (\bar t-\tau_1,\bar
t)$ the subset $J(t)=J_0(t)\cup D_t$ is an open subset of $M_t$. We
define
$$J(\bar t-\tau_1,\bar t)=\cup_{t\in (\bar t-\tau_1,\bar t)}J(t).$$
Then $J(\bar t-\tau_1,\bar t)$ is an open subset of ${\mathcal M}$.
See {\sc Fig.}~\ref{fig:disappear}.
\end{defn}

\begin{lem}\label{15.17}
Fix a Ricci flow with surgery $({\mathcal M},G)$, a point $x\in
{\mathcal M}$ and constants $r\ge r_{i+1}>0$ as in the statement of
Proposition~\ref{delta0ri+1}. For any $1<\ell<\infty$ the following
holds. Suppose that $\bar t\in [T_{i-1},T)$ is a surgery time and
that  $\gamma(\tau)$ is a path with $\gamma(\tau)\in M_{\bar
t-\tau}$. Let $\{p_1,\ldots,p_k\}$ be the tips of all the surgery
caps at time $\bar t$ and let $\bar h$ be the scale of surgery at
time $\bar t$. Suppose that for some $0<\tau_1\le \ell^{-1} \bar
h^2$ there are no surgery times in the interval $(\bar t-\tau_1,\bar
t)$. We identify all $M_t$ for $t\in [\bar t-\tau_0,\bar t)$ with
$M_{\bar t-\tau_1}$ using the flow. Suppose that $\gamma(0)\in
M_{\bar t}\setminus \cup_{i=1}^kB(p_i,\bar t,(50+A_0)\bar h)$, and
lastly, suppose that
$$\int_{0}^{\tau_1}|X_\gamma(\tau)|^2d\tau\le \ell.$$
Then $\gamma$ is disjoint from the open subset $J(\bar t-\tau_1,\bar
t))$ of ${\mathcal M}$.
\end{lem}

\begin{proof}
Suppose that the lemma is false and let $\gamma\colon [0,\bar
\tau]\to {\mathcal M}$ be a path satisfying the hypothesis of the
lemma with $\gamma(\bar \tau)\in J(\bar t-\tau_1,\bar t)$.
 Since $$\gamma(0)\in
M_{\bar t}\setminus \cup_iB(p_i,\bar t,(50+A_0)\bar h),$$ if follows
that $\gamma(0)$ is separated from the boundary of
$s_{N_i}^{-1}(-25,0]$ by distance at least $20\bar h$. Since the $
J_0(t)$ are contained in the disjoint union of strong $\delta$-necks
$N_i$ centered at the $2$-spheres along which we do surgery, and
since $\tau_1\le \bar h^2/\ell<\bar h^2$, it follows that, provided
that $\delta$ is sufficiently small, for every $t\in [\bar
t-\tau_1,\bar t)$, the metric on $J_0(t)$ is at least
$1/2$ the metric on $ J_0(\bar t)$.
 It follows that, for $\delta$ sufficiently small, if there is a $\tau\in [0,\tau_1]$
with $\gamma(\tau)\in J(t)$ then $\int_0^{\tau_1}|X_\gamma|d\tau> 10\bar h$.
Applying Cauchy-Schwarz we see that
$$\int_0^{\tau_1}|X_\gamma(\tau)|^2d\tau\ge (10\bar h)^2/\tau_1.$$
Since $\tau_1\le \ell^{-1}(\bar h)^2$, we see that
$$\int_0^{\tau'}|X_\gamma(\tau)|^2d\tau>\ell,$$
contradicting our hypothesis.
\end{proof}

\subsection{Limits of a sequence of paths with short ${\mathcal L}$-length}

Now using Lemmas~\ref{L_0} and~\ref{15.17} we show that it is possible to take limits
of certain types of sequences of paths parameterized by backward time to create
minimizing ${\mathcal L}$-geodesics.

We shall work with a compact subset of ${\bf t}^{-1}([T_{i-1},T])$ that is
obtained by removing appropriate open neighborhoods of the exposed regions.

\begin{defn}
Fix $\ell<\infty$. Let $\theta_0=\theta_0(\ell)$ be as in
Proposition~\ref{bigell}. For each surgery time $\bar t\in
[T_{i-1},T]$, let $h(\bar t)$ be the scale of the surgery.
 Let $p_1,\ldots,p_k$ be
the tips of the surgery caps at time $\bar t$. For each $1\le j\le k$,
we consider $B_j(\bar
t)=B(p_j,\bar t,(A_0+10)h(\bar t))$, and we let $\Delta t_j\le {\rm
min}(\theta_0,(T-\bar t)/h^2(\bar t))$ be maximal subject to the condition that
there is an embedding $\rho_j\colon B_j(\bar t)\times [\bar t,\bar t+h^2(\bar
t)\Delta t_j)$ into ${\mathcal M}$ compatible with time and the vector field.
Clearly, $B'_j=B(p_j,\bar t,(10+A_0)h)\cap C_{\bar t}$ is contained in $J(\bar
t)$. Let $\bar t'$ be the previous surgery time if there is one, otherwise set
$\bar t'=0$. Also for each $\bar t$ we set $\tau_1(\ell,\bar t)={\rm
min}\left(h(\bar t)^2/\ell,\bar t-\bar t'\right)$. For each $t\in (\bar
t-\tau_1(\ell,\bar t),\bar t)$ let $\widetilde J(t)\subset J(t)$ be the union
of $D_t$, the disappearing region at time $t$, and $\coprod_iB_i'(t)$,  the
result of flowing $\coprod_iB'_i$ backward to time $t$. Then we set
$\widetilde J(\bar t-\tau_1(\ell,\bar t),\bar t)\subset J(\bar
t-\tau_1(\ell,\bar t),\bar t)$ equal to the union over $t\in (\bar
t-\tau_1(\ell,\bar t),\bar t)$ of $\widetilde J(t)$.

 By construction,
for each surgery time $\bar t$, the union  $$\nu_{\rm sing}(\ell,\bar
t)=\widetilde J(\bar t-\tau_1(\ell,\bar t),\bar t)\cup\bigcup
_iB_i\times [\bar t,\bar t+h^2(\bar t)\Delta t_i)$$ is an open
subset of ${\mathcal M}$ containing all the exposed regions and
singular points at time $\bar t$.

We define $Y(\ell )\subset {\bf t}^{-1}([T_{i-1},T])$ to be the
complement of the $\cup_{\bar t}\nu_{\rm sing}(\ell,\bar t)$ where
the union is over all surgery times $\bar t\in [T_{i-1},T]$.
Clearly, $Y(\ell )$ is a closed subset of ${\bf
t}^{-1}([T_{i-1},T])$ and hence $Y(\ell )$ is a compact subset
contained in the open subset of smooth points of ${\mathcal M}$.
(Notice that $Y(\ell )$ depends on $\ell $ because $\tau_1(\ell,\bar
t)$ and $\theta_0$ depend on $\ell$.)
\end{defn}

\begin{prop}\label{Lgeoexist}
Fix $0<L<\infty$. Set
$$L_0=L+4(T_{i+1})^{3/2}.$$ Suppose that for all
$t\in [T_{i-1},T_{i+1}]$, the surgery control parameter
$\overline\delta(t)\le \overline\delta_1(L_0,r_{i+1})$
where the right-hand side is the constant from
Lemma~\ref{L_0}. Suppose that $\gamma_n$ is a sequence of paths in
$({\mathcal M},G)$ parameterized by backward time $\tau\in
[0,\bar\tau]$ with $\bar \tau\le T-T_{i-1}$, with $\gamma_n(0)=x$ and
with
$${\mathcal L}(\gamma_n)\le L$$
for all $n$. Then:
\begin{enumerate}
\item[(1)]
After passing to a subsequence, there is a limit $\gamma$ defined on
$[0,\bar\tau]$. The limit $\gamma$ is a continuous path and is a
uniform limit of the $\gamma_n$. The limit is contained in the open subset
of smooth points of ${\mathcal M}$ and has finite
${\mathcal L}$-length satisfying
$${\mathcal L}(\gamma)\le {\rm liminf}_{n\rightarrow\infty}{\mathcal L}(\gamma_n).$$
\item[(2)]
If there is a point $y\in M_{T-\bar\tau}$ such that
$\gamma_n(\bar\tau)=y$ for all $n$, and if the $\gamma_n$ are a
sequence of paths parameterized by backward time from $x$ to $y$
with ${\rm lim}_{n\rightarrow\infty} {\mathcal L}(\gamma_n)$ being no greater than
the ${\mathcal L}$-length of any path from $x$ to $y$, then the limit
$\gamma$ of a subsequence is a minimizing ${\mathcal L}$-geodesic
connecting $x$ to $y$ contained in the open subset of smooth points of ${\mathcal M}$.
\item[(3)] There is $\ell<\infty$ depending only on $L$ such that
any path $\gamma$ parameterized by backward time from $x$ to a
point $y\in {\bf t}^{-1}([T_{i-1},T))$ whose ${\mathcal L}$-length
is at most $L$ is contained in the compact subset $Y(\ell)$ given in
the previous definition.
\end{enumerate}
\end{prop}

\begin{proof}
Given $L_0$, we set
$$\tau'={\rm min}\left(\frac{r_{i+1}^4}{(256)L_0^2},{\rm ln}(\root 3\of{2})r_{i+1}^2\right)$$
as in Lemma~\ref{L_0} and then define $\ell=L_0/\sqrt{\tau'}$. We
also let $A={\rm min}(2(50+A_0),A_0(\ell))$ and
$\theta_0=\theta_0(\ell)$ as in Proposition~\ref{bigell}. Lastly, we
let
$\overline\delta_1(L_0,r_{i+1})=\delta''(A)=\delta''(A,\theta_0,\overline\delta_0)$
from Propositions~\ref{bigell} and~\ref{altern}. We suppose that
$\delta(t)\le \overline\delta_1(L_0,r_{i+1})$ for all $t\in
[T_{i-1},T]$.

Let $\bar t\in [T_{i-1},T]$ be a surgery time, and let $\bar h$ be
the scale of the surgery at this time. For each surgery cap
${\mathcal C}$ with tip $p$ at a time $\bar t\in [T_{i-1},T]$ let
$\Delta t({\mathcal C})$ be the supremum of those $s$ with $0\le s\le \theta_0\bar
h^2$ for which there is an embedding
$$\rho_{\mathcal C}\colon
B(p,\bar t,2(A_0+50)\bar h)\times[\bar t,\bar t+s)\to{\mathcal M}$$
compatible with time and the vector field. We set
$$P_0({\mathcal C})=\rho_{\mathcal C}\left(B(p,\bar t,(A+50)\bar h)\times
[\bar t,\bar t+{\rm min}(\bar h^2/2,\Delta t({\mathcal C})))\right).$$

\begin{claim} Any path $\gamma$ beginning at $x$ and parameterized by backward
time misses $P_0({\mathcal C})$ if ${\mathcal L}(\gamma)<L$.
\end{claim}

\begin{proof}
Set $\tau_0=T-\bar t$. Of course, $\tau_0\le T-T_{i-1}\le T_{i+1}-T_{i-1}$.
Consider the restriction of $\gamma$ to $[0,\tau_0]$.
We have
\begin{eqnarray*}
\lefteqn{\int_0^{\tau_0}\sqrt{\tau}\left(R_+(\gamma_n(\tau))+|X_{\gamma_n}(\tau)|^2\right)d\tau
}  \\
& \le &
\int_0^{T-T_{i-1}}\sqrt{\tau}\left(R_+(\gamma_n(\tau))+|X_{\gamma_n}(\tau)|^2\right)d\tau \\
& \le & \int_0^{T-T_{i-1}}\sqrt{\tau}\left(R(\gamma_n(\tau))+|X_{\gamma_n}(\tau)|^2\right)d\tau+
\int_0^{T-T_{i-1}}
6\sqrt{\tau}d\tau \\
& = &
\int_0^{T-T_{i-1}}\sqrt{\tau}\left(R(\gamma_n(\tau))+|X_{\gamma_n}(\tau)|^2\right)d\tau
+4(T-T_{i-1})^{3/2}
\\
& \le &
\int_0^{\tau_0}\sqrt{\tau}\left(R(\gamma_n(\tau))+|X_{\gamma_n}(\tau)|^2\right)d\tau+4(T_{i+1})^{3/2}
\end{eqnarray*}
Thus, the hypothesis that ${\mathcal L}(\gamma_n)\le L$ implies that
\begin{equation}\label{L'ineq}
\int_0^{\tau_0}\sqrt{\tau}\left(R_+(\gamma_n(\tau))+|X_{\gamma_n}(\tau)|^2\right)d\tau\le
L_0.\end{equation}
The claim now follows immediately from Lemma~\ref{L_0}.
\end{proof}

Now set $t'$ equal to the last surgery time before $\bar t$ or set $t'=0$ if $\bar t$ is
the first surgery time. We set $\tau_1(\bar t)$ equal to the minimum of $\bar t-t'$ and
$\bar h^2/\ell$.

Assume that $\gamma(0)=x$ and that ${\mathcal L}(\gamma)\le L$. It follows from
Lemma~\ref{L_0} that the restriction of the path $\gamma$ to $[0,\tau']$ lies
in a region where the Riemann curvature is bounded above by $r^{-2}\le
r_{i+1}^{-2}$. Hence, since $\bar h<\bar\delta(t)^2r_{i+1}\ll r_{i+1}$, this
part of the path  is disjoint from all strong $\delta$-necks (evolving
backward for rescaled time $(-1,0]$). That is to say, $\gamma|_{[0,\tau']}$ is
disjoint from $J_0(t)$ for every $t\in (\bar t-\tau_1(\bar t),\bar t)$ for any
surgery time $\bar t\le T$. It follows immediately that $\gamma|_{[0,\tau']}$
is disjoint from $J(\bar t-\tau_1(\bar t),\bar t)$.

\begin{claim}
For every surgery time $\bar t\in [T_{i-1},T]$, the path $\gamma$
starting at $x$ with ${\mathcal L}(\gamma)\le L$ is disjoint from
$J(\bar t,\bar t-\tau_1(\bar t))$.
\end{claim}

\begin{proof}
By the remarks above, it suffices to consider surgery times $\bar
t\le T-\tau'$. It follows immediately from the previous claim that
for any surgery time $\bar t$, with the scale of the surgery being
$\bar h$ and with $p$ being the tip of a surgery cap at this time,
we have $\gamma$ is disjoint from $B(p,\bar t,(50+A_0)\bar h)$.
Also,
$$\int_{T-\bar t}^{T-\bar t+\tau_1(\bar t)}\sqrt{\tau}|X_\gamma(\tau)|^2d\tau
\le {\mathcal L}_+(\gamma)\le L_0.$$
Since we can assume $T-\bar t\ge \tau'$ this implies that
$$\int_{T-\bar t}^{T-\bar t-\tau_1(\bar t)}|X_\gamma(\tau)|^2d\tau\le L_0/\sqrt{\tau'}=\ell.$$
The claim is now immediate from Lemma~\ref{15.17}.
\end{proof}

From these two claims we see immediately that $\gamma$ is contained in the
compact subset $Y(\ell)$ which is contained in the open subset of smooth points
of ${\mathcal M}$. This proves the third item in the statement of the
proposition. Now let us turn to the limit statements.

Take a sequence of paths $\gamma_n$ as in the statement of
Proposition~\ref{Lgeoexist}.
 By Lemma~\ref{L_0} the
restriction of each $\gamma_n$ to the interval $[0,{\rm
min}(\bar\tau,\tau')]$ is contained in $P(x,T,r/2,-r^2)$. The
arguments in the proof of Lemma~\ref{geoexist} (which involve
changing variables to $s=\sqrt{\tau}$) show that, after passing to a
subsequence, the restrictions of the $\gamma_n$ to $[0,{\rm
min}(\bar\tau,\tau')]$ converge uniformly to a path $\gamma$ defined
on the same interval. Furthermore,
\begin{equation*}
\int_0^{{\rm
min}(\bar\tau,\tau')}\sqrt{\tau}|X_\gamma(\tau)|^2d\tau\le {\rm
liminf}_{n\rightarrow\infty}\int_0^{{\rm
min}(\bar\tau,\tau')}\sqrt{\tau}|X_{\gamma_n}(\tau)|^2d\tau,
\end{equation*}
so that
\begin{eqnarray}\label{gamliminf1}\lefteqn{\int_0^{{\rm
min}(\bar\tau,\tau')}\sqrt{\tau}\left(R(\gamma(\tau)+|X_\gamma(\tau)|^2\right)d\tau\le} \\
& & {\rm liminf}_{n\rightarrow\infty}\int_0^{{\rm
min}(\bar\tau,\tau')}\sqrt{\tau}\left(R(\gamma_n(\tau))+|X_{\gamma_n}(\tau)|^2\right)d\tau.\nonumber
\end{eqnarray}

If $\bar\tau\le \tau'$, then we have established the existence of a
limit as required. Suppose now that $\bar\tau>\tau'$. We turn our
attention to the paths $\gamma_n|_{[\tau',\bar\tau]}$. Let
$T_{i-1}<\bar t\le T-\tau'$ be either a surgery time or $T-\tau'$, and let $t'$ be the
maximum of the last surgery time before $\bar t$ and $T_{i-1}$. We consider the
restriction of the $\gamma_n$ to the interval $[T-\bar t,T-t']$. As
we have seen, these restrictions  are disjoint from $J(\bar
t-\tau_1(\bar t),\bar t)$ and also from the exposed region at time
$\bar t$, which is denoted $E(\bar t)$,  and from $J_0(\bar t)$. Let
$$Y={\bf t}^{-1}([T-t',T-\bar t])\setminus
\left(J(\bar t-\tau_1(\bar t),\bar t)\cup (E(\bar t)\cup J_0(\bar
t))\right).$$ This is a compact subset with the property that any
point $y\in Y$ is connected by a backward flow line lying entirely
in $Y$ to a point $y(t')$ contained in $M_{t'}$.

 Since $Y$ is compact there is a finite upper bound on the
Ricci curvature on $Y$, and hence to ${\mathcal L}_\chi(G)$ at any
point of $Y$. Since all backward flow lines from points of $Y$
extend all the way to $M_{t'}$, it follows that there is a constant
$C'$ such that
$$|X_{\gamma_n}(\tau)|_{G(t')}\le C'|X_{\gamma_n}(\tau)|_{G(t)}$$
for all $t\in [t',\bar t]$. Our hypothesis that the ${\mathcal
L}(\gamma_n)$ are uniformly bounded, the fact that the curvature is
pinched toward positive and the fact that there is a uniform bound
on the lengths of the $\tau$-intervals implies that the
$$\int_{T-\bar t}^{T-t'} \sqrt{\tau}|X_{\gamma_n}(\tau)|^2d\tau$$
are uniformly bounded. Because $T-\bar t$ is at least $\tau'>0$, it
follows that the $\int_{T-\bar t}^{T-t'}|X_{\gamma_n}|^2d\tau$ have
a uniform upper bound. This then implies that there is a constant
$C_1$ such that for all $n$ we have
$$\int_{T-\bar t}^{T-t'}|X_{\gamma_n}(\tau)|_{G(t')}^2d\tau\le
C_1.$$ Thus, after passing to a subsequence, the $\gamma_n$ converge
uniformly to a continuous $\gamma$ defined on $[T-\bar t,T-t']$.
Furthermore, we can arrange that the convergence is a weak
convergence in $W^{1,2}$. This means that $\gamma$ has a derivative
in $L^2$ and
$$\int_{T-\bar t}^{T-t'}|X_{\gamma}(\tau)|^2d\tau\le {\rm
liminf}_{n\rightarrow
\infty}\int_{T-\bar t}^{T-t'}|X_{\gamma_n}(\tau)|^2d\tau.$$

Now we do this simultaneously for all $\bar t=T-\tau'$ and for all the finite number of surgery times
in $[T_{i-1},T-\tau']$. This
gives a limiting path $\gamma\colon [\tau',\bar\tau]\to {\mathcal M}$.
Putting together the above inequalities we see that the limit
satisfies
\begin{equation}\label{gamliminf2}
\int_{\tau'}^{\bar\tau}\sqrt{\tau}\left(R(\gamma(\tau))+|X_\gamma(\tau)|^2\right)d\tau\le
{\rm
liminf}_{n\rightarrow\infty}\int_{\tau'}^{\bar\tau}\sqrt{\tau}\left(R(\gamma_n(\tau))+
|X_{\gamma_n}(\tau)|^2\right)d\tau.
\end{equation}
 Since we have already
arranged that there is a limit on $[0,\tau']$, this produces a
limiting path $\gamma\colon [0,\tau_0]\to {\mathcal M}$. By
Inequalities~\ref{gamliminf1} and~\ref{gamliminf2} we see that
$${\mathcal L}(\gamma)\le {\rm liminf}_{i\rightarrow \infty}{\mathcal
L}(\gamma_n).$$
The limit lies in the compact subset $Y(\ell)$ and hence is contained
in the open subset of smooth points of ${\mathcal M}$.
 This completes the proof of the first statement of the proposition.

Now suppose, in addition to the above, that all the $\gamma_n$
have the same endpoint $y\in M_{T-\tau_0}$ and that ${\rm
lim}_{n\rightarrow\infty}{\mathcal L}(\gamma_n)$ is less than or equal
to the ${\mathcal L}$-length of any path parameterized by backward
time connecting $x$ to $y$. Let $\gamma$ be the limit of a
subsequence as constructed in the proof of the first part of this
result. Clearly, by what we have just established, $\gamma$ is a path
parameterized by backward time from $x$ to $y$ and ${\mathcal
L}(\gamma)\le {\rm lim}_{n\rightarrow\infty}{\mathcal L}(\gamma_n)$.
This means that $\gamma$ is a minimizing ${\mathcal L}$-geodesic
connecting $x$ to $y$, an ${\mathcal L}$-geodesic contained in the open subset of smooth
points of ${\mathcal M}$.

This completes the proof of the proposition.
\end{proof}

\begin{cor}\label{bardelta1}
Given $L<\infty$, let
$\overline{\delta}_1=\overline\delta_1(L+4(T_{i+1}^{3/2},r_{i+1})$ be as given in Lemma~\ref{L_0}.
If
$\overline\delta(t)\le \overline\delta_1$ for all $t\in [T_{i-1},T_{i+1}]$,
then for any $x\in {\bf t}^{-1}([T_i,T_{i+1}))$ and for any $y\in
M_{T_{i-1}}$, if there is a  path $\gamma$ parameterized by
backward time connecting $x$ to $y$ with ${\mathcal L}(\gamma)\le L$,
then there is a minimizing ${\mathcal L}$-geodesic contained in the
open subset of smooth points of ${\mathcal M}$ connecting $x$ to
$y$.
\end{cor}

\begin{proof}
Choose an ${\mathcal L}$-minimizing sequence of paths from $x$  to
$y$ and apply the previous proposition.
\end{proof}

\subsection{Completion of the proof of Proposition~\protect{\ref{delta0ri+1}}}

Having found a compact subset of the open subset of smooth points of ${\mathcal M}$
that contains all paths parameterized by backward time whose ${\mathcal L}$-length is not
too large, we are in a position to prove Proposition~\ref{delta0ri+1}, which states
the existence of a minimizing ${\mathcal L}$-geodesics in ${\mathcal M}$ from $x$
and gives estimates on their ${\mathcal L}$-lengths.

\begin{proof}(of Proposition~\ref{delta0ri+1}).
Fix $r\ge r_{i+1}>0$. Let $({\mathcal M},G)$ and $x\in {\mathcal M}$ be
as in the statement of Proposition~\ref{delta0ri+1}. We set
$L=8\sqrt{T_{i+1}}(1+T_{i+1})$, and we set
$$\delta={\rm min}\left(\delta_i,\overline \delta_1(L+4(T_{i+1})^{3/2},r_{i+1})\right),$$
where $\overline\delta_1$ is  as given in Lemma~\ref{L_0}.
Suppose that $\overline\delta(t)\le
\delta$ for all $t\in [T_{i-1},T_{i+1})$. We set $U$ equal to the
subset of ${\bf t}^{-1}([T_{i-1},T))$ consisting of all points $y$
for which there is a path $\gamma$ from $x$ to $y$, parameterized by
backward time, with ${\mathcal L}(\gamma)<L$. For each $t\in
[T_{i-1},T)$ we set $U_t=U\cap M_t$. According to
Corollary~\ref{bardelta1} for any $y\in U$ there is a minimizing
${\mathcal L}$-geodesic connecting $x$ to $y$ and this geodesic lies
in the open subset of ${\mathcal M}$ consisting of all the smooth
points of ${\mathcal M}$; in particular, $y$ is a smooth
point of ${\mathcal M}$. Let ${\mathcal L}_x\colon U\to \Ar$ be the
function that assigns to each $y\in U$ the ${\mathcal L}$-length of
a minimizing ${\mathcal L}$-geodesic from $x$ to $y$. Of course,
${\mathcal L}_x(y)<L$ for all  $y\in U$. Now let us show that the
restriction of ${\mathcal L}_x$ to any time-slice $U_t\subset U$
achieves its minimum along a compact set. For this, let $y_n\in U_t$
be a minimizing sequence for ${\mathcal L}_x$ and for each $n$ let
$\gamma_n$ be a minimizing ${\mathcal L}$-geodesic connecting $x$ to
$y_n$. Since ${\mathcal L}(\gamma_n)<L$ for all $n$, according to
Proposition~\ref{Lgeoexist}, we can pass to a subsequence that
converges to a limit, $\gamma$, connecting $x$ to some point $y\in
M_t$ and ${\mathcal L}(\gamma)\le {\rm
inf}_{n}{\mathcal L}(\gamma_n)<L$. Hence, $y\in
U_t$, and clearly ${\mathcal L}_x|_{U_t}$ achieves its minimum at
$y$. Exactly the same argument with $y_n$ being a sequence of points
at which ${\mathcal L}_x|_{U_t}$ achieves its minimum shows  that
the subset of $U_t$ at which ${\mathcal L}_x$ achieves its minimum
is a compact set.

We set $Z\subset U$ equal to the set of $y\in U$ such that
${\mathcal L}_x(y)\le {\mathcal L}_x(y')$ for all $y'\in  U_{{\bf
t}(y)}$.

\begin{claim}
The subset  $Z'=\{z\in Z|{\mathcal L}_x(z)\le L/2\}$ has the property
that for any compact interval $I\subset [T_{i-1},T)$ the
intersection ${\bf t}^{-1}(I)\cap Z'$ is compact.
\end{claim}

\begin{proof}
Fix a compact interval $I\subset [T_{i-1},T)$. Let $\{z_n\}$ be a
sequence in $Z'\cap {\bf t}^{-1}(I)$. By passing to a subsequence we
can assume that the sequence ${\bf t}(z_n)=t_n$ converges to some
$t\in I$, and that ${\mathcal L}_x(z_n)$ converges to some $D\le
L/2$. Since the surgery times are discrete, there is a neighborhood
$J$ of $t$ in $I$ such that the only possible surgery time in $J$ is
$t$ itself. By passing to a further subsequence if necessary, we can
assume that $t_n\in J$ for all $n$. Fix $n$. First, let us consider
the case when $t_n\ge t$. Let $\gamma_n$ be a minimizing ${\mathcal
L}$-geodesic from $x$ to $z_n$. Then we form the path $\widehat
\gamma_n$ which is the union of $\gamma_n$ followed by the flow line
for the vector field $\chi$ from the endpoint of $\gamma_n$ to
$M_t$. (This flowline exists since there is no surgery time in the
open interval $(t,t_n]$.) If $t_n<t$, then we set $\widehat
\gamma_n$ equal to the restriction of $\gamma_n$ to the interval
$[0,T-t]$. In either case let $\hat y_n\in M_t$ be the endpoint of
$\widehat\gamma_n$. Since $M_t$ is compact, by passing to a
subsequence we can arrange that the $\hat y_n$ converge to a point $y\in
M_t$. Clearly,
 ${\rm lim}_{n\rightarrow\infty }z_n=y$.

It is also the case that ${\rm lim}_{n\rightarrow\infty}{\mathcal L}(\widehat\gamma_n)=
{\rm lim}_{n\rightarrow\infty}{\mathcal L}(\gamma_n)=D\le L/2$. This
means that $y\in U$ and that ${\mathcal L}_x(y)\le D\le L/2$. Hence,
the greatest lower bound of the values of ${\mathcal L}_x$ on $U_t$
is at most $D\le L/2$, and consequently $Z'\cap U_t\not=\emptyset$.
Suppose that the minimum value of ${\mathcal L}_x$ on $U_t$ is
$D'<D$. Let $z\in U_t$ be a point where this minimum value is
realized, and let $\gamma$ be a minimizing ${\mathcal L}$-geodesic
from $x$ to $z$. Then by restricting $\gamma$ to subintervals
$[0,t-\mu]$ shows that the minimum value of ${\mathcal L}_x$ on
$U_{t+\mu}\le (D'+D)/2$ for all $\mu>0$ sufficiently small. Also,
extending $\gamma$ by adding a backward vertical flow line from $z$
shows that the minimum value of ${\mathcal L}_x$ on $U_{t-\mu}$ is
at most $(D'+D)/2$ for all $\mu>0$ sufficiently small. (Such a
vertical flow line backward in time exists since $z\in U$ and hence
$z$ is contained in the smooth part of ${\mathcal M}$.) This
contradicts the fact that limit of the minimum values of ${\mathcal
L}_x$ on $U_{t_n}$ converge to $D$ as $t_n$ converges to $t$. This
contradiction proves  that the minimum value of ${\mathcal L}_x$ on
$U_t$ is $D$, and consequently the point $y\in Z'$. This proves that
$Z'\cap {\bf t}^{-1}(I)$ is compact, establishing the claim.
\end{proof}

At this point we  have established that Properties (1),(2), and (4);
So it remains only to prove Property (3) of Proposition~\ref{delta0ri+1}.
 To do this we define the reduced length function
$l_x\colon U\to \Ar$ by
$$l_x(q)=\frac{{\mathcal L}_x(q)}{2\sqrt{T-{\bf t}(q)}} \ \ {\rm and} \ \
l_x^{\rm min}(\tau)={\rm min}_{q\in M_t}l_x(q).$$ We consider the
subset ${\mathcal S}$ of $\tau'\in (0,T-T_{i-1}]$ with $l_x^{\rm
min}(\tau)\le L/2$ for all $\tau\le \tau'$. Recall that by the
choice of $L$, we have $3\sqrt{T-T_{i-1}}< L/2$. Clearly, the
minimum value of $l_x$ on $U_{T-\tau}$ converges to $0$
as $\tau\rightarrow 0$, implying that this set is non-empty. Also, from
its definition, ${\mathcal S}$ is an interval with $0$ being one
endpoint.

\begin{lem}\label{lvalue}
Let $l_x^{\rm min}(\tau')$ be the minimum value of $l_x$ on
$U_{T-\tau'}$. For any $\tau\in {\mathcal S}$ we have $l_x^{\rm
min}(\tau)\le 3/2$.
\end{lem}

\begin{proof}
Given that we have already established Properties 1,2 and 4 of
Proposition~\ref{delta0ri+1}, this is immediate from
Corollary~\ref{compactcoro}.
\end{proof}

Now let us establish that ${\mathcal S}=(0,T-T_{i-1}]$. As we remarked above,
${\mathcal S}$ is a non-empty interval with $0$ as one endpoint. Suppose that
it is of the form $(0,\tau]$ for some $\tau<T-T_{i-1}$. Then by the previous
claim, we have $l_x^{\rm min}(\tau)\le 3/2$ so that there is an ${\mathcal
L}$-geodesic $\gamma$ from $x$ to a point $y\in M_{T-\tau}$ with ${\mathcal
L}(\gamma)\le 3\sqrt{\tau}< L/2$.
 This implies that for all $\tau'>\tau$ but sufficiently close to
$\tau$, there is a point $y(\tau')\in U_{T-\tau'}$ with ${\mathcal
L}_x(y(\tau'))<L/2$. This shows that all $\tau'$ greater than and sufficiently close
to $\tau$ are contained in ${\mathcal S}$.
 This is a
contradiction of the assumption that ${\mathcal S}=(0,\tau]$.

Suppose now that ${\mathcal S}$ is of the form $(0,\tau)$, and set $t=T-\tau$.
Let $t_n\rightarrow t$ and $z_n\in Z'\cap U_{t'}$. The same argument as above shows
that for every $n$ we have ${\mathcal L}_x(z_n)\le 3\sqrt{T-t_n}$. For all $n$
sufficiently large, there are no surgery times in the interval $(t,t_n)$.
Hence, by passing to a subsequence, we can arrange that the $z_n$ converge to a
point $z\in M_t$. Clearly, $${\mathcal L}_x(z)\le {\rm
limsup}_{n\rightarrow\infty}{\mathcal L}_x(z_n)\le 3\sqrt{T-t},$$ so that $\tau\in
{\mathcal S}$. This again contradicts the assumption that ${\mathcal
S}=(0,\tau)$.

The only other possibility is that the set of $\tau$ is $(0,T-T_{i-1}]$ and the
minimum value of ${\mathcal L}$ on $U_t$ is at most $3\sqrt{T-t}$ for all $t\in
[T_{i-1},T)$. This is exactly the third property stated in
Proposition~\ref{delta0ri+1}. This completes the proof of that proposition.
\end{proof}

\section{Completion of the proof of
Proposition~\protect{\ref{KAPPALIMIT}}}

Now we are ready to establish Proposition~\ref{KAPPALIMIT}, the
non-collapsing result. We shall do this by finding a parabolic
neighborhood whose size, $r'$, depends only on $r_i$, $C$ and
$\epsilon$, on which the sectional curvature is bounded by
$(r')^{-2}$ and so that the ${\mathcal L}$-distance from $x$ to any
point of the final time-slice of this parabolic neighborhood is
bounded. Recall that in Section~\ref{rsmall} we established it when
$R(x)=r^{-2}$ with $r\le r_{i+1}<\epsilon$. Here we assume that
$r_{i+1}<r\le \epsilon$.  Fix $\delta=\delta(r_{i+1})$ from
Proposition~\ref{delta0ri+1} and set $L=8\sqrt{T_{i+1}}(1+T_{i+1})$.

First of all, in Claim~\ref{newkappa0r0t0} we have  seen that there
is $\kappa_0$ such that ${\bf t}^{-1}[0,T_1]$ is $\kappa_0$
non-collapsed on scales $\le \epsilon$. Thus, we may assume that
$i\ge 1$.

Recall that ${\bf t}(x)=T\in (T_{i},T_{i+1}]$. Let $\gamma$ be an
${\mathcal L}$-geodesic contained in the smooth part of ${\mathcal
M}$ from $x$ to a point in $ M_{T_{i-1}}$ with ${\mathcal
L}(\gamma)\le 3\sqrt{T-T_{i-1}}$. That such a $\gamma$ exists was
proved in Proposition~\ref{delta0ri+1}. We shall find a point $y$ on
this curve with  $R(y)\le 2r_i^{-2}$. Then  we find a backward
parabolic neighborhood centered at $y$ on which ${\mathcal L}$ is
bounded and so that the slices have volume bounded from below. Then
we can apply the results from Chapter~\ref{noncoll} to establish the
$\kappa$ non-collapsing.

\begin{claim}\label{tau0cl}
There is $\tau_0$  with ${\rm max}(\epsilon^2,T-T_i)\le \tau_0\le
T-T_{i-1}-\epsilon^2$ such that $R(\gamma(\tau_0))< r_i^{-2}$.
\end{claim}

\begin{proof}
Let $T'={\rm max}(\epsilon^{2},T-T_i)$ and let
$T''=T-T_{i-1}-\epsilon^2$, and suppose that $R(\gamma(\tau))\ge
r_i^{-2}$ for all $\tau\in [T',T'']$. Then we see that
$$\int_{T'}^{T''}\sqrt{\tau}\left(R(\gamma(\tau)+|X_\gamma(\tau)|^2\right)d\tau
\ge \frac{2}{3}r_i^{-2}\left((T'')^{3/2}-(T')^{3/2}\right).$$ Since
$R\ge -6$ because the curvature is pinched toward positive, we see that
\begin{eqnarray*}
{\mathcal L}(\gamma) & \ge & \frac{2}{3}r_i^{-2}\left((T'')^{3/2}-(T')^{3/2}\right)-\int_0^{T'}
6\sqrt{\tau}d\tau-\int_{T''}^{T-T_{i-1}}6\sqrt{\tau}d\tau \\
& = & \frac{2}{3}r_i^{-2}\left((T'')^{3/2}-(T')^{3/2}\right)-
4(T')^{3/2}-4\left((T-T_{i-1})^{3/2}-(T'')^{3/2}\right).
\end{eqnarray*}

\begin{claim} We have the following estimates:
\begin{eqnarray*}
(T'')^{3/2}-(T')^{3/2} & \ge &
\frac{1}{4}(T-T_{i-1})^{3/2} \\
4(T')^{3/2} & \le &  4(T-T_{i-1})^{3/2} \\
4\left((T-T_{i-1})^{3/2}-(T'')^{3/2}\right) & \le &
\frac{2t_0}{25}\sqrt{(T-T_{i-1})}.\end{eqnarray*}
\end{claim}

\begin{proof}
Since $T_i-T_{i-1}\ge t_0$ and $T\ge T_i$, we see that
$T''/(T-T_{i-1})\ge 0.9$. If $T'=T-T_i$, then since
$T<T_{i+1}=2T_i=4T_{i-1}$ we have $T'/(T-T_{i-1})\le 2/3$. If
$T'=\epsilon^2$, since $\epsilon^2\le t_0/50$, and $T-T_{i-1}\ge
t_0$, we see that $T'\le (T-T_{i-1})/50$. Thus, in both cases we
have $T'\le 2(T-T_{i-1})/3$. Since $(0.9)^{3/2}>0.85$ and
$(2/3)^{3/2}\le 0.6$, the first inequality follows.

The second inequality is clear since $T'<(T-T_{i-1})$.

The last inequality is clear from the fact that
$T''=T-T_{i-1}-\epsilon^2$ and $\epsilon\le \sqrt{t_0/50}$.
\end{proof}

 Putting these together yields
 $${\mathcal L}(\gamma)\ge
 \left[\left(\frac{1}{6}r_i^{-2}-4\right)\left(T-T_{i-1}\right)-\frac{2t_0}{25}\right]
 \sqrt{T-T_{i-1}}.$$

Since
$$r_i^{-2}\ge r_0^{-2}\ge \epsilon^{-2}\ge 50/t_0,$$
and $T-T_{i-1}\ge t_0$ we see that \begin{eqnarray*}{\mathcal
L}(\gamma) & \ge &
\left[\left(\frac{50}{6t_0}-4\right)t_0-\frac{2t_0}{25}\right]\sqrt{T-T_{i-1}}
\\
& \ge & (8-5t_0)\sqrt{T-T_{i-1}} \\
& \ge & 4\sqrt{T-T_{i-1}}.
\end{eqnarray*}

(The last inequality uses the fact that $t_0=2^{-5}$.) But this
contradicts the fact that ${\mathcal L}(\gamma)\le
3\sqrt{T-T_{i-1}}$.
\end{proof}

Now fix $\tau_0$ satisfying Claim~\ref{tau0cl}.
 Let $\gamma_1$ be the
restriction of $\gamma$ to the subinterval $[0,\tau_0]$, and let
$y=\gamma_1(\tau_0)$. Again using the fact that $R(\gamma(\tau))\ge
-6$ for all $\tau$, we see that
\begin{equation}\label{Lgamma1}
{\mathcal L}(\gamma_1)\le {\mathcal L}(\gamma)+4(T-T_{i-1})^{3/2}\le
3(T_{i+1})^{1/2}+4(T_{i+1})^{3/2}.\end{equation}

Set $t'=T-\tau_0$. Notice that from the definition we
have $t'\le T_i$. Consider $B=B(y,t',\frac{r_i}{2C})$, and
define $\Delta={\rm min}(r_i^2/16C,\epsilon^2)$. According to
Lemma~\ref{controlnbhd} every point $z$ on a backward flow line
starting in $B$ and defined for time at most $ \Delta$  has the
property that $R(z)\le 2r_i^{-2}$. For any surgery time $\bar t$ in
$[t'-\Delta,t')\subset [T_{i-1},T)$ the scale $\bar h$ of the
surgery at time $\bar t$ is $\le \delta(\bar t)^2 r_i$, and hence
every point of the surgery cap has scalar curvature at least
$D^{-1}\delta(\bar t)^{-4}r_i^{-2}$. Since $\overline\delta(\bar t)\le
\overline \delta\le \delta_0\le {\rm min}(D^{-1},1/10)$, it follows
that every point of the surgery cap has curvature at least
$\delta_0^{-3}r_i^{-2}\ge 1000r_i^{-2}$. Thus,
 no point
$z$ as above can lie in a surgery cap. This means that the entire
backward parabolic neighborhood $P(y,t',\frac{r_i}{2C},-\Delta)$
exists in ${\mathcal M}$, and the scalar curvature is bounded by $2
r_i^{-2}$ on this backward parabolic neighborhood. Because of the
curvature pinching toward positive assumption, there is $C'<\infty$
depending only on $r_i$ and such that the Riemann curvature is
bounded by $C'$ on $P(y,t',\frac{r_i}{2C},-\Delta)$.

Consider the one-parameter family of metrics $g(\tau),\ 0\le \tau\le
\Delta$, on $B(y,t',\frac{r_i}{2C})$ obtained by restricting the
horizontal metric $G$ to the backward parabolic neighborhood.
There is $0<\Delta_1\le \Delta /2$
depending only on $C'$ such that for every $\tau\in [0,\Delta_1]$ and
every non-zero tangent vector $v$ at a point of
$B(y,t',\frac{r_i}{2C})$ we have
$$\frac{1}{2}\le \frac{|v|^2_{g(\tau)}}{|v|^2_{g(0)}}\le 2.$$
Set $\hat r={\rm min}(\frac{r_i}{32C},\Delta_1/2)$, so that $\hat r$ depends
only on $r_i$, $C$, and $\epsilon$. Set $t''=t'-\Delta_1$. Clearly,
$B(y, t'',\hat r)\subset B(y,t',\frac{r_i}{2C})$ so that $B(y,t'',\hat r)\subset
P(y,t',\frac{r_i}{2C},-\Delta)$. Of course, it then follows that the parabolic neighborhood
$P(y,t'',\hat r,-\Delta_1)$ exists in ${\mathcal M}$ and
$$P(y,t'',\hat r,-\Delta_1)\subset P(y,t',\frac{r_i}{2 C},-\Delta),$$ so that the
Riemann curvature is bounded above by $C'$ on the parabolic
neighborhood $P(y,t'',\hat r,-\Delta_1)$. We set $r'={\rm
min}(\hat r,(C')^{-1/2},\sqrt{\Delta_1}/2)$, so that $r'$ depends only on
$r_i$, $C$, and $\epsilon$. Then  the parabolic neighborhood
$P(y,t'',r',-(r')^2)$ exists in ${\mathcal M}$ and $|{\rm Rm}|\le
(r')^{-2}$ on $P(y,t'',r',-(r')^2)$. Hence, by the inductive
non-collapsing assumption either $y$ is contained in a component of $M_{t''}$ of positive sectional curvature or
$${\rm Vol}\,B(y,t'',r')\ge
\kappa_i(r')^3.$$

If $y$ is contained in a component of $M_{t''}$ of positive sectional curvature, then by
  Hamilton's result, Theorem~\ref{flowtoround},
under Ricci flow the component of $M_{t''}$  containing $y$ flows forward as
a family of  components of positive sectional curvature until it disappears. Since there is path moving backwards in time
from $x$ to $y$, this means that the original point $x$ is contained
in a component of its time-slice with positive sectional curvature.

Let us consider the other possibility when ${\rm Vol}\,B(y,t'',r')\ge
\kappa_i(r')^3$.
 For each $z\in B(y,t'',r')$ let
 $$\mu_z\colon [T-t',T-t'']\to  B(y,t',\frac{r_i}{2C})$$ be
 the $G(t')$-geodesic connecting $y$ to $z$. Of course
 $$|X_{\mu_z}(\tau)|_{G(t')}\le \frac{r'}{\Delta_1}$$
 for every $\tau\in [0,\Delta_1]$. Thus,
 $$|X_{\mu_z}(\tau)|_{G(T-\tau)}\le
 \frac{\sqrt{2}r'}{\Delta_1}$$ for all $\tau\in [T-t',T-t'']$.
 Now we let $\tilde \mu_z$ be the resulting path parameterized by backward time
 on the time-interval $[T-t',T-t'']$.
 We estimate
 \begin{eqnarray*}{\mathcal
 L}(\tilde \mu_z) & = &
 \int_{T-t'}^{T-t''}\sqrt{\tau}\left(R(\tilde\mu_z(\tau))+|X_{\tilde\mu_z}(\tau)|^2\right)d\tau
 \\
& \le & \sqrt{T-t''}\int_{T-t'}^{T-t''}\left(2 r_i^{-2}+\frac{2(r')^2}{\Delta_1^2}\right)d\tau
\\
& \le & \sqrt{T-t''}(2r_i^{-2}\Delta_1+\frac{1}{2}) \le
\sqrt{T}\left(\frac{1}{16C}+\frac{1}{2}\right).
\end{eqnarray*}
In passing to the last inequality we use the fact, from the
definitions that $r'\le \sqrt{\Delta_1}/2$ and $\Delta\le
r_i^2/16C$, whereas $\Delta_1\le \Delta/2$.

Since $C>1$, we see that
$${\mathcal L}(\tilde \mu_z)\le \sqrt{T}.$$
Putting this together with the estimate, Equation~(\ref{Lgamma1}),
for ${\mathcal L}(\gamma_1)$ tells us that for each $z\in
B(y,t'',r')$ we have
$${\mathcal L}(\gamma_1*\tilde\mu_z)\le
4(T_{i+1})^{1/2}+4(T_{i+1})^{3/2}\le L/2.$$ Hence, by
Proposition~\ref{delta0ri+1} and the choice of $L$, there is a
minimizing ${\mathcal L}$-geodesic from $x$ to each point of
$B(y,t'',r')$ of length $\le L/2$, and these geodesics are contained
in the smooth part of ${\mathcal M}$. In fact, by
Proposition~\ref{Lgeoexist} there is a compact subset $Y$ of
the open subset of smooth points of ${\mathcal M}$ that contains all
the minimizing ${\mathcal L}$-geodesics from $x$ to points of
$B(y,t'',r')$.

Then, by Corollary~\ref{fullmeasure} (see also,
Proposition~\ref{lips}), the intersection, $B'$, of ${\mathcal U}_x$
with $B(y,t'',r')$
 is an open subset of full measure in $B(y,t'',r')$.
Of course, ${\rm Vol}\,B'={\rm Vol}\,B(y,t'',r')\ge \kappa_i(r')^3$
and the function $l_x$ is bounded by $L/2$ on $B'$. It now follows
from Theorem~\ref{THM} that there is $\kappa>0$ depending only on
$\kappa_i$, $r'$, $\epsilon$ and $L$ such that $x$ is $\kappa$
non-collapsed on scales $\le \epsilon$. Recall that $L$ depends only
on $T_{i+1}$, and $r'$ depends only on $r_i,C,C'$ and $\epsilon$,
whereas $C'$ depends only on $r_i$. Thus, in the final analysis,
$\kappa$ depends only on $\kappa_i$ and $r_i$ (and $C$ and
$\epsilon$ which are fixed). This entire analysis assumed that for
all $t\in [T_{i-1},T_{i+1})$ we have the inequality
$\overline\delta(t)\le
\overline\delta_1(L+4(t_{i+1})^{3/2},r_{i+1})$ as in
Lemma~\ref{L_0}. Since $L$ depends only on $i$ and $t_0$, this shows
that the upper bound for $\delta$ depends only on $r_{i+1}$ (and on
$i$, $t_0$, $C$, and $\epsilon$). This completes the proof of
Proposition~\ref{KAPPALIMIT}.

\chapter{Completion of the proof of Theorem~\protect{\ref{MAIN}}}

We have established the requisite non-collapsing result assuming the
existence of strong canonical neighborhoods. In order to complete
the proof of Theorem~\ref{MAIN} it remains for us to show the
existence of strong canonical neighborhoods. This is the result of
the next section.

\section{Proof of the strong canonical neighborhood assumption}\label{sectcannbhd}

 \begin{prop}\label{extend}
Suppose that for some $i\ge 0$ we have surgery parameter sequences
$\delta_0\ge \delta_1\ge\cdots\ge \delta_i>0$, $\epsilon=r_0\ge
r_1\ge \cdots\ge r_i>0$ and $\kappa_0\ge \kappa_1\ge \cdots\ge
\kappa_i>0$. For any $r_{i+1}\le r_i$, let $\delta(r_{i+1})>0$ be
the constant in Proposition~\ref{KAPPALIMIT} associated to these
three sequences and to $r_{i+1}$.
 Then there are positive
constants $r_{i+1}\le r_i$ and $\delta_{i+1}\le\delta(r_{i+1})$ such that the
following holds. Suppose that $({\mathcal M},G)$ is a Ricci flow with surgery
defined for $0\le t<T$ for some $T\in (T_i,T_{i+1}]$ with surgery control
parameter $\overline\delta(t)$. Suppose that the restriction of this Ricci flow
with surgery to ${\bf t}^{-1}([0,T_i))$ satisfies Assumptions (1) -- (7) and
also the five properties given in the hypothesis of Theorem~\ref{MAIN} with
respect to the given sequences. Suppose also that $\overline\delta(t)\le
\delta_{i+1}$ for all $t\in [T_{i-1},T]$. Then $({\mathcal M},G)$  satisfies
the strong $(C,\epsilon)$-canonical neighborhood assumption with parameter
$r_{i+1}$.
 \end{prop}

\begin{proof}
Suppose that the result does not hold.
  Then we can take a sequence of $r_a\rightarrow 0$ as
 $a\rightarrow \infty$, all less than $r_i$, and for each $a$ a sequence
 $\delta_{a,b}\rightarrow 0$ as $b\rightarrow \infty$ with each
 $\delta_{a,b}\le \delta(r_a)$, where $\delta(r_a)\le \delta_i$ is the constant in Proposition~\ref{KAPPALIMIT}
 associated to the three sequences given in the statement of this proposition and $r_a$, such that for
 each $a,b$ there is a Ricci flow with surgery $({\mathcal
 M}_{(a,b)},G_{(a,b)})$ defined for $0\le t<T_{(a,b)}$ with $T_i<
 T_{(a,b)}\le T_{i+1}$ with control parameter $\overline \delta_{(a,b)}(t)$
 such the flow satisfies the hypothesis of the proposition with
 respect to these constants but fails to satisfy the conclusion.

\begin{lem}\label{16.2}
For each $a$, and given $a$, for all $b$ sufficiently large there is $t_{(a,b)}\in
[T_i,T_{(a,b)})$ such that the restriction of $({\mathcal
M}_{(a,b)},G_{(a,b)})$ to ${\bf t}^{-1}\left([0,t_{(a,b)})\right)$
satisfies the strong $(C,\epsilon)$-canonical neighborhood
assumption with parameter $r_a$ and furthermore, there is
$x\in{\mathcal M}_{(a,b)}$ with ${\bf t}(x_{(a,b)})=t_{(a,b)}$ at
which the strong $(C,\epsilon)$-canonical neighborhood assumption
with parameter $r_a$ fails.
\end{lem}

\begin{proof}
Fix $a$.
 By supposition, for each $b$ there is a point $x\in {\mathcal M}_{(a,b)}$ at
which the strong $(C,\epsilon)$-canonical neighborhood assumption
fails for the parameter $r_a$. We call points at which this
condition fails {\em counterexample points}. Of course, since
$r_a\le r_i$ and since the restriction of $({\mathcal
M}_{(a,b)},g_{(a,b)})$ to ${\bf t}^{-1}([0,T_i))$ satisfies the
hypothesis of the proposition, we see that any counterexample point
$x$ has ${\bf t}(x)\ge T_i$. Take a sequence $x_n=x_{n,(a,b)}$   of
counterexample points with ${\bf t}(x_{n+1})\le {\bf t}(x_{n})$ for
all $n$ that minimizes ${\bf t}$ among all counterexample points in
the sense that for any $\xi>0$ and for any counterexample point
$x\in {\mathcal M}_{(a,b)}$ eventually ${\bf t}(x_{n})< {\bf
t}(x)+\xi$. Let $t'=t_{(a,b)}'={\rm lim}_{n\rightarrow\infty}{\bf
t}(x_{n})$. Clearly, $t'\in [T_i,T_{(a,b)})$, and by construction
the restriction of $({\mathcal M}_{(a,b)},G_{(a,b)})$ to ${\bf
t}^{-1}([0,t'))$ satisfies the $(C,\epsilon)$-canonical neighborhood assumption
with parameter $r_a$. Since the surgery times are discrete, there is
$t''=t_{(a,b)}''$ with $t'<t''\le T_{(a,b)}$ and a diffeomorphism
$\psi=\psi_{(a,b)}\colon M_{ t'}\times[t',t'')\to {\bf
t}^{-1}([t',t''))$ compatible with time and the vector field. We
view $\psi^*G_{(a,b)}$ as a one-parameter family of metrics
$g(t)=g_{(a,b)}(t)$ on $M_{t'}$ for $t\in [t',t'')$. By passing to a
subsequence we can arrange that ${\bf t}(x_{n})\in [t',t'')$ for all
$n$. Thus, for each $n$ there are $y_n=y_{n,(a,b)}\in M_{t'}$ and
$t_{n}\in [t',t'')$ with $\psi(y_{n},t_{n})=x_{n}$. Since
$M_{t'}$ is a compact $3$-manifold, by passing to a further
subsequence we can assume that $y_{n}\rightarrow x_{(a,b)}\in
M_{t'}$. Of course, $t_{n}\rightarrow t'$ as $n\rightarrow\infty$ and
 ${\rm lim}_{n\rightarrow
\infty}x_{n}=x_{(a,b)}$ in ${\mathcal M}_{(a,b,)}$.

 We claim that, for
all $b$ sufficiently large, the strong $(C,\epsilon)$-canonical neighborhood assumption
with parameter $r_a$ fails at $x_{(a,b)}$. Notice that since $x_{(a,b)}$
is the limit of a sequence where the strong
$(C,\epsilon)$-neighborhood assumption fails, the points in the
sequence converging to $x_{(a,b)}$ have scalar curvature at least
$r_a^{-2}$. It follows that $R(x_{(a,b)})\ge r_a^{-2}$.  Suppose
that $x_{(a,b)}$ satisfies the strong $(C,\epsilon)$-canonical
neighborhood assumption with parameter $r_a$. This means that there
is a neighborhood $U=U_{(a,b)}$ of $x_{(a,b)}\in M_{t'}$
which is a strong $(C,\epsilon)$-canonical neighborhood of
$x_{(a,b)}$. According to Definition~\ref{defncanonnbhd} there are
four possibilities. The first two we consider are that $(U,g(t'))$
is an $\epsilon$-round component or a $C$-component. In either of
these cases, since the defining inequalities given in
Definition~\ref{defnepsround} and~\ref{Ccomponent} are strong
inequalities,  all metrics on $U$ sufficiently close to $g(t')$
in the $C^\infty$-topology the satisfy these same inequalities. But as $n$ tends to
$\infty$, the metrics $g(t_n)|_U$ converge in the $C^\infty$-topology
to $g(t')|_U$. Thus, in these two cases, for all $n$ sufficiently
large, the metrics $g(t_n)$ on $U$ are $(C,\epsilon)$-canonical
neighborhood metrics of the same type as $g(t_{(a,b)}')|_U$. Hence, in
either of these cases, for all $n$ sufficiently large $x_{n,(a,b)}$
has a strong $(C,\epsilon)$-canonical neighborhood of the same type
as $x_{(a,b)}$, contrary to our assumption about the sequence
$x_{n,(a,b)}$.

Now suppose that there is
 a $(C,\epsilon)$-cap whose core contains $x_{(a,b)}$.
This is to say that $(U,g(t'))$ is a $(C,\epsilon)$-cap whose core
contains $x_{(a,b)}$. By Proposition~\ref{canonvary}, for all $n$
sufficiently large, $(U,g(t_{n}))$ is also a $(C,\epsilon)$-cap with
the same core. This core contains $y_{n}$ for all $n$ sufficiently
large, showing that $x_{n}$ is contained in the core of a
$(C,\epsilon)$-cap  for all $n$ sufficiently large.

Now let us consider the remaining case when $x_{(a,b)}$ is the
center of a strong $\epsilon$-neck. In this case we have an
embedding $\psi_{U_{(a,b)}}\colon
U_{(a,b)}\times({t_{(a,b)}'}-R^{-1}(x_{(a,b)}),t_{(a,b)}']\to
{\mathcal M}_{(a,b)}$ compatible with time and the vector field and
a diffeomorphism $f_{(a,b)}\colon S^2\times
(-\epsilon^{-1},\epsilon^{-1})\to U_{(a,b)}$ so that $
(f_{(a,b)}\times{\rm
Id})^*\psi_{U_{(a,b)}}^*(R({x_{(a,b)}})G_{(a,b)})$ is
$\epsilon$-close in the $C^{[1/\epsilon]}$-topology to the
evolving product metric $h_0(t)\times ds^2,\ -1<t\le 0$, where $h_0(t)$
is a round metric of scalar curvature $1/(1-t)$ on $S^2$ and $ds^2$
is the Euclidean metric on the interval. Here, there are two
subcases to consider.
\begin{enumerate}
\item[(i)] $\psi_{U_{(a,b)}}$ extends backward past $
t_{(a,b)}'-R^{-1}(x_{(a,b)})$.
\item[(ii)] There is a flow line through a point $y_{(a,b)}\in U_{(a,b)}$ that is defined
on the interval $[t_{(a,b)}'-R^{-1}(x_{(a,b)}),t_{(a,b)}']$ but with
the value of the flow line at $t_{(a,b)}'-R^{-1}(x_{(a,b)})$ an
exposed point.
\end{enumerate}

Let us consider the first subcase. The embedding $\psi_{U_{(a,b)}}$
extends forward in time because of the diffeomorphism
$\psi_{(a,b)}\colon M_{t_{(a,b)}'}\times [t_{(a,b)}',t_{(a,b)}'')\to
{\mathcal M}_{(a,b)}$ and, by assumption, $\psi_{U_{(a,b)}}$ extends
backward in time some amount. Thus, for all $n$ sufficiently large,
we can use these extensions of $\psi_{U_{(a,b)}}$ to define an
embedding $\psi_{n,(a,b)}\colon U_{(a,b)}\times ({\bf
t}(x_{n,(a,b)})-R^{-1}(x_{n,(a,b)}),{\bf t}(x_{n,(a,b)})]\to
{\mathcal M}_{(a,b)}$ compatible with time and the vector field.
Furthermore, since the $\psi_{n,(a,b)}$ converge in the
$C^\infty$-topology as $n$ tends to infinity to $\psi_{U_{(a,b)}}$,
the Riemannian metrics $(f_{(a,b)}\times{\rm
Id})^*\psi_{n,(a,b)}^*(R(x_{n,(a,b)})G_{a,b})$ converge in the
$C^\infty$-topology to the pullback  $(f_{(a,b)}\times{\rm
Id})^*\psi_{U_{(a,b)}}^*(R(x_{(a,b)})G_{a,b})$. Clearly then, for
fixed $(a,b)$ and for all $n$ sufficiently large the pullbacks of
the rescalings of these metrics by $R(x_{n,(a,b)})$ are within
$\epsilon$ in the $C^{[1/\epsilon]}$-topology of the standard
evolving flow $h_0(t)\times ds^2, -1<t\le 0$, on the product of
$S^2$ with the interval. Under these identifications the points
$x_{n,(a,b)}$ correspond to points $(p_{n,(a,b)},s_{n,(a,b)})\in
S^2\times (-\epsilon^{-1},\epsilon^{-1})$ where ${\rm
lim}_{n\rightarrow\infty}s_{n,(a,b)}=0$. The last thing we  do is
to choose diffeomorphisms $\varphi_{n,(a,b)}\colon
(-\epsilon^{-1},\epsilon^{-1})\to (-\epsilon^{-1},\epsilon^{-1})$
that are the identity near both ends, such that $\varphi_{n,(a,b)}$
carries $0$ to $s_{n,(a,b)}$ and such that the $\varphi_{n,(a,b)}$
converge to the identity in the $C^\infty$-topology for fixed
$(a,b)$ as $n$ tends to infinity. Then, for all $n$ sufficiently
large, the composition
$$S^2\times (-\epsilon^{-1},\epsilon^{-1})\stackrel{{\rm Id}\times
\varphi_{n,(a,b)}}{\longrightarrow}S^2\times
(-\epsilon^{-1},\epsilon^{-1})\stackrel{f_{(a,b)}}{\longrightarrow}
U\stackrel{\psi_{n,(a,b)}}{\longrightarrow} {\mathcal M}_{(a,b)}$$
is a strong $\epsilon$-neck centered at $x_{n,(a,b)}$. This shows that for any
$b$ for which the first subcase holds, for all $n$ sufficiently large, there is
a strong $\epsilon$-neck centered at $x_{n,(a,b)}$.

Now suppose that the second subcase holds for all $b$. Here, unlike
all previous cases, we shall have to let $b$ vary and we shall prove
the result only for $b$ sufficiently large.
We shall show that for all $b$ sufficiently large, $x_{(a,b)}$ is contained
in the core of a $(C,\epsilon)$-cap. This will establish the result by contradiction,
for as we showed in the previous case, if $x_{(a,b)}$ is contained in the core of
a $(C,\epsilon)$-cap, then the same is true for the $x_n$ for all $n$ sufficiently large,
contrary to our assumption.

For the moment fix $b$. Set $\bar
t_{(a,b)}=t_{(a,b)}'-R_{G_{(a,b)}}(x_{(a,b)})^{-1}$. Since, by supposition the
embedding $\psi_{U_{(a,b)}}$ does not extend  backwards past $\bar t_{(a,b)}$, it
must be the case that $\bar t_{(a,b)}$  is a surgery time and
furthermore that there is a surgery cap ${\mathcal C}_{(a,b)}$ at
this time with the property that there is a point $y_{(a,b)}\in
U_{(a,b)}$ such that $\psi_{U_{(a,b)}}(y_{(a,b)},t)$ converges to a
point $z_{(a,b)}\in {\mathcal C}_{(a,b)}$ as $t$ tends to $\bar
t_{(a,b)}$ from above. (See {\sc Fig.}~\ref{fig:neckcap}.) We denote by
$p_{(a,b)}$ the tip of ${\mathcal C}_{(a,b)}$, and we denote by
$\bar h_{(a,b)}$ the scale of the surgery at time $\bar t_{(a,b)}$.

Since the statement that $x_{(a,b)}$ is contained in the core of a $(C,\epsilon)$-cap is
a scale invariant statement, we are free to replace $({\mathcal M}_{(a,b)},G_{(a,b)})$ with
$(\widetilde{\mathcal M}_{(a,b)},\widetilde G_{(a,b)})$, which has been rescaled to make
$\bar h_{(a,b)}=1$ and shifted in time so that $\bar t_{(a,b)}=0$. We denote the new time function
by ${\bf \tilde t}$.
(Notice that this rescaling and time-shifting is different from what we usually do.
Normally, when we have a base point like $x_{(a,b)}$ we rescale to make its scalar curvature one and
we shift time to make it be at time $0$. Here we have rescaled based on the scale of the surgery cap
rather than $R(x_{(a,b)})$.) We set $\tilde Q_{(a,b)}=R_{\widetilde G_{(a,b)}}(x_{(a,b)})$ and we set
$\tilde t'_{(a,b)}={\bf \tilde t}(x_{(a,b)})$. Since the initial time of the strong $\epsilon$-neck
is zero,  $\tilde t'_{(a,b)}=\tilde Q_{(a,b)}^{-1}$. We denote the flow line backward in time from
 $y_{(a,b)}$ by $y_{(a,b)}(\tilde t),\ 0\le \tilde t\le \tilde t'_{(a,b)}$, so that
 $y_{(a,b)}(\tilde t'_{(a,b)})=y_{(a,b)}$.
 Since $U_{(a,b)}$ is a strong $\epsilon$-neck, by our choice of $\epsilon$,
 it follows from Lemma~\ref{neckcurv} and rescaling that
 $R(\psi(y_{(a,b)},\tilde t))$ is within $(0.01) \tilde Q_{(a,b)}$ of
 $\tilde Q_{(a,b)}/(1+\tilde Q_{(a,b)}( \tilde t_{(a,b)}'-\tilde t))$ for all
 $t\in (0,\tilde t_{(a,b)}']$. By taking limits as $t$
 approaches $0$, we see that
  $R_{\widetilde G_{(a,b)}}(z_{(a,b)})$ is within $(0.01)\tilde Q_{(a,b)}$ of $\tilde Q_{(a,b)}/2$.
Let $D$ be the universal constant given in  Lemma~\ref{A_0}, so that
the scalar curvature at any point of the standard initial metric is
at least $D^{-1}$ and at most $D$.  It follows from the third item
in Theorem~\ref{LOCALSURGERY} that, since we have rescaled to make
the surgery scale one, for all $b$ sufficiently large the scalar
curvature on the surgery ${\mathcal C}_{(a,b)}$ is at least
$(2D)^{-1}$ and at most $2D$. In particular, for all $b$
sufficiently large
$$(2D)^{-1}\le R_{\widetilde G_{(a,b)}}(z_{(a,b)})\le 2D.$$
Together with the above estimate relating $R_{\widetilde G_{(a,b)}}(z_{(a,b)})$ and
$\tilde Q_{(a,b)}$, this gives
\begin{equation}\label{DQ}
(5D)^{-1}\le
\tilde Q_{(a,b)}\le 5D.
\end{equation} Since the flow line from
$z_{(a,b)}$ to $y_{(a,b)}$ lies in the closure of a strong
$\epsilon$-neck of scale $\tilde Q_{(a,b)}^{-1/2}$, the scalar curvature is
less than
 $6D$ at every point of this flow line.
 According to Proposition~\ref{Restim} there is $\theta_1<1$
(depending only on $D$) such that $R(q,t)\ge 8D$ for all $(q,t)$
in the standard solution with $t\ge \theta_1$.

\begin{figure}[ht]
  \centerline{\epsfbox{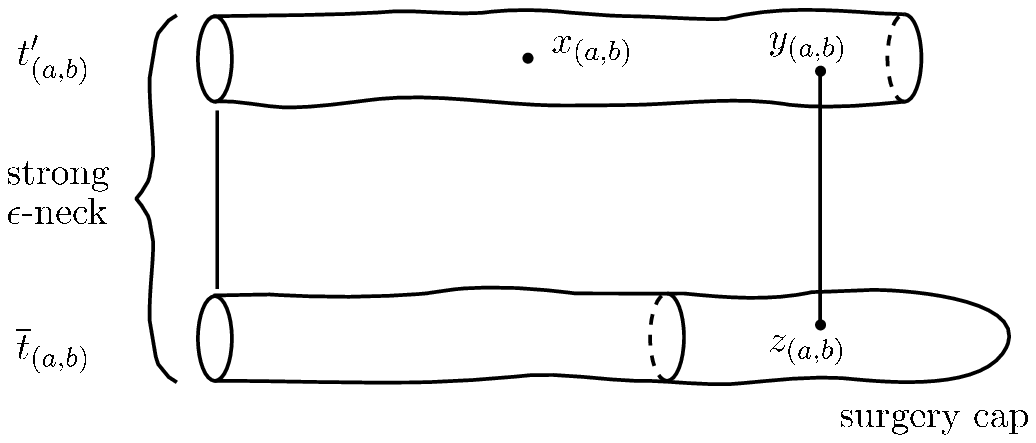}}
  \caption{A strong neck with initial time in a surgery cap}\label{fig:neckcap}
\end{figure}

By the fifth property of Theorem~\ref{stdsoln} there is
$A'(\theta_1)<\infty$ such that in the
 standard flow, $B(p_0,0,A)$ contains $B(p_0,\theta_1,A/2)$ for every $A\ge A'(\theta_1)$.
 We set $A$ equal to the maximum of $A'(\theta_1)$ and
 $$3\left((1.2)\sqrt{5D}\epsilon^{-1}+(1.1)(A_0+5)+C\sqrt{5D}\right).$$
Now for any $\overline\delta>0$ for all  $b$ sufficiently large, we
have $\delta_{(a,b)}\le \delta''(A,\theta_1,\overline\delta)$, where
$\delta''(A,\theta_1,\overline\delta)$
is the constant given in Proposition~\ref{bigell}.

\begin{claim}
 Suppose that $b$ is sufficiently large so that
$\delta_{(a,b)}\le \delta''(A,\theta_1,\overline\delta_0)$, where $\overline\delta_0$ is the constant
given in Definition~\ref{defnoverlinedelta0}. Then
$\tilde t_{(a,b)}'\le \theta_1$.
\end{claim}

\begin{proof}
In this proof we shall fix $(a,b)$, so we drop these indices from
the notation.
 Consider
$s\le \theta_1$ maximal so that there is an embedding
$$\psi=\psi_{(a,b)}\colon
B(p_0,0,A)\times [0,s)\to \widetilde{\mathcal
M}_{(a,b)}$$ compatible with time and the vector field. First
suppose that $s<\theta_1$. Then according to
Proposition~\ref{altern} either the entire ball $B(p,0,A)$
disappears at time $s$ or $s$ is the final time of the
time interval of definition for the flow $(\widetilde{\mathcal
M}_{(a,b)},\widetilde G_{(a,b)})$. Since we have the flow line from $z\in
B(p_0,0,A)$ extending to time $\tilde t'=\tilde t'_{(a,b)}$, in either case this
implies that $\tilde t'<s$, proving that $t'<\theta_1$ in this case.

Now suppose that $s=\theta_1$. By the choice of
$\theta_1$, for the standard solution  the scalar curvature at every
$(q,\theta_1)$ is at least $8D$. Since
$\delta_{(a,b)}\le \delta''(A,\theta_1,\overline\delta_0)$, by the
definition of $\overline\delta_0$ given in
Definition~\ref{defnoverlinedelta0} and by Proposition~\ref{altern}
the scalar curvature of the pullback of the metric under  $\psi$ is
within a factor of two of the scalar curvature of the rescaled
standard solution. Hence, the scalar curvature along the flow line
$(z,t)$ through $z$ limits to at least
$8D$ as $t$ tends to $\theta_1$.
Since the scalar curvature on $(z,t)$ for $t\in [0,\tilde t']$ is
bounded above by $6D$, it follows that $\tilde t'<\theta_1$ in this case as
well. This completes the proof of the
claim.
\end{proof}

Thus, we have maps
$$\psi_{(a,b)}\colon B(p_0,0,A)\times [0,\tilde t'_{(a,b)}]\to \widetilde {\mathcal M}_{(a,b)}$$
compatible with time and the vector field, with the property that for each $\delta>0$, for
all $b$ sufficiently large the pullback under this map of $\widetilde G_{(a,b)}$ is within $\delta$
in the $C^{[1/\delta]}$-topology of the restriction of the standard solution.
Let $w_{(a,b)}$ be the result of flowing $x_{(a,b)}$ backward to time $0$.

\begin{claim}
For all $b$ sufficiently large, $w_{(a,b)}\in \psi_{(a,b)}(B(p_0,0,A)\times\{0\})$.
\end{claim}

\begin{proof}
First notice that, by our choice of $\epsilon$,
 every point in the $0$ time-slice of the closure of the strong $\epsilon$-neck
centered at $x_{(a,b)}$ is within distance $(1.1)\tilde Q_{(a,b)}^{-1}\epsilon^{-1}$
of $w_{(a,b)}$. In particular,
$$d_{\widetilde G_{(a,b)}}(w_{(a,b)},y_{(a,b)})<(1.1)\tilde Q_{(a,b)}^{-1/2}\epsilon^{-1}.$$
Since $y_{(a,b)}$ is contained in the surgery cap and the scale of the surgery
at this time is $1$, $y_{(a,b)}$ is within distance $A_0+5$ of $p_{(a,b)}$.
Hence, by the triangle inequality and Inequality~(\ref{DQ}), we have
\begin{eqnarray*}
d_{\widetilde G_{(a,b)}}(w_{(a,b)},p_{(a,b)}) & < & (1.1)\tilde Q_{(a,b)}^{-1/2}\epsilon^{-1}+(A_0+5) \\
& < & (1.1)\sqrt{5D}\epsilon^{-1}+(A_0+5).
\end{eqnarray*}
For $b$ sufficiently large, the image  $\psi_{(a,b)}(B(p_0,0,A))$
contains the ball of radius $(0.95)A$
centered at $p_{(a,b)}$. Since by our choice of $A$ we have $(0.95)A>(1.1)\sqrt{5D}\epsilon^{-1}+(A_0+5)$,
the claim follows.
\end{proof}

We define $q_{(a,b)}\in B(p_0,0,A)$ so that
$\psi_{(a,b)}(q_{(a,b)},0)=w_{(a,b)}$. Of course,
$$\psi_{(a,b)}(q_{(a,b)},\tilde t'_{(a,b)})=x_{(a,b)}.$$ If follows from the
above computation that for all $b$ sufficiently large we have
$$d_0(q_{(a,b)},p_0)<(1.15)\tilde Q_{(a,b)}^{-1/2}\epsilon^{-1}+(1.05)(A_0+5).$$
Since the standard flow has non-negative curvature, it is a distance non-increasing flow. Therefore,
$$d_{\tilde t'_{(a,b)}}(q_{(a,b)},p_0)<(1.15)\tilde Q_{(a,b)}^{-1/2}\epsilon^{-1}+(1.05)(A_0+5).$$

Suppose that a point $(q,\tilde t'_{(a,b)})$ in the standard solution were the center of a
$\beta\epsilon/3$-neck, where $\beta$ is the constant from Proposition~\ref{neckglue}.
 Of course, for all $b$ sufficiently large, $R(q,\tilde
t'_{(a,b)})>(0.99)\tilde Q_{(a,b)}$. Since $\beta<1/2$ and
$\epsilon<\sqrt{5D}(A_0+5)/2$ and $\tilde Q_{(a,b)}\le 5D$, it follows from the
above distance estimate that this neck would contain $(p_0,\tilde t'_{(a,b)})$.
But this is impossible: since $(p_0,\tilde t'_{(a,b)})$ is an isolated fixed
point of an isometric $SO(3)$-action on the standard flow, all the sectional
curvatures at $(p_0,\tilde t'_{(a,b)})$ are equal, and this is in contradiction
with estimates on the sectional curvatures at any point of an $\epsilon$-neck
given in Lemma~\ref{neckcurv}. We can then conclude from Theorem~\ref{stdsolncannbhd}
that for all $b$ sufficiently large, the point $(p_0,\tilde t'_{(a,b)})$ is
contained in the core of a $(C(\beta\epsilon/3),\beta\epsilon/3)$-cap
$Y_{(a,b)}$ in the $\tilde t'_{(a,b)}$ time-slice of the standard solution.
Now note that for
all $b$ sufficiently large, the scalar curvature of $(q_{(a,b)},\tilde
t'_{(a,b)})$ is at least $(0.99)\tilde Q_{(a,b)}$, since the scalar curvature of $x_{(a,b)}$
is equal to $Q_{(a,b)}$. This implies that the
diameter of $Y_{(a,b)}$ is at most
$$(1.01)\tilde Q^{-1/2}_{(a,b)}C(\beta\epsilon/3)<(1.1)\sqrt{5D}C(\beta\epsilon/3).$$
Since $B(p_0,0,A)$ contains $B(p_0,\tilde t'_{(a,b)},A/2)$, and since
$C>C(\beta\epsilon/3)$, it follows from the definition of $A$, the above
distance estimate, and the triangle inequality that for all $b$ sufficiently
large $B(p_0,0,A)\times\{\tilde t'_{(a,b)}\}$ contains $Y_{(a,b)}$.

 Since $C>
C(\beta\epsilon/3)+1$ and since for $b$ sufficiently large
$\psi_{(a,b)}^*\widetilde G_{(a,b)}$ is arbitrarily close to the
restriction of the standard solution metric, it follows from
Lemma~\ref{canonvary} that for all $b$ sufficiently large, the image $\psi_{(a,b)}(Y_{(a,b)})$
 is a $(C,\epsilon)$-cap whose core contains $x_{(a,b)}$.
As we have already remarked, this contradicts the assumption that no $x_n$
has a strong $(C,\epsilon)$-canonical neighborhood.

 This completes the proof in the last case and establishes
Lemma~\ref{16.2}.
\end{proof}

\begin{rem}
Notice that even though $x_{(a,b)}$ is the center of a strong
$\epsilon$-neck, the canonical neighborhoods of the $x_n$ constructed
in the second case are not a strong $\epsilon$-necks but rather are
$(C,\epsilon)$-caps coming from applying the flow to a
neighborhood of the surgery cap ${\mathcal C}$.
\end{rem}

Now we return to the proof of Proposition~\ref{extend}. For each $a$, we pass
to a subsequence (in $b$) so that Lemma~\ref{16.2} holds for all $(a,b)$. For
each $(a,b)$, let $t_{(a,b)}$ be as in that lemma.  We fix a point
$x_{(a,b)}\in {\bf t}^{-1}(t_{(a,b)})\subset {\mathcal M}_{(a,b)}$ at which the
canonical neighborhood assumption with parameter $r_a$ fails.  For each $a$
choose $b(a)$ such that $\delta_{b(a)}\rightarrow 0$ as $a\rightarrow \infty$. For each
$a$ we set $({\mathcal M}_a,G_a)=({\mathcal M}_{(a,b(a))},G_{(a,b(a))})$, we
set $t_a=t_{(a,b(a))}$, and we let $x_a=x_{(a,b(a))}\in {\mathcal M}_a$. Let
$(\widetilde {\mathcal M}_a,\widetilde G_a)$ be the Ricci flow with surgery
obtained from $({\mathcal M}_a,G_a)$ by shifting $t_a$ to $0$ and rescaling the
metric and time by $R(x_a)$. We have the points $\tilde x_a$ in the $0$
time-slice of $\widetilde {\mathcal M}_a$ corresponding to $x_a\in {\mathcal
M}_a$. Of course, by construction $R_{\widetilde G_a}(\tilde x_a)=1$ for all
$a$.

We shall take limits of a subsequence of this sequence of based
Ricci flows with surgery. Since $r_a\rightarrow 0$ and $R(x_a)\ge
r_a^{-2}$, it follows that $R(x_a)\rightarrow \infty$. By
Proposition~\ref{KAPPALIMIT}, since $\delta_{b(a)}\le \delta(r_a)$
it follows that the restriction of $(\widetilde {\mathcal
M}_a,\widetilde G_a)$ to ${\bf t}^{-1}(-\infty,0)$ is
$\kappa$-non-collapsed on scales $\le \epsilon R_{G_a}^{1/2}(x_a)$.
By passing to a subsequence we arrange that one of the following two
possibilities holds:
\begin{enumerate}
\item[(i)] There is $A<\infty$ and $t'<\infty$ such that,
for each $a$ there is a flow line through a point $y_a$ of $B_{\widetilde
G_a}(\tilde x_a,0,A)$ that is not defined on all of $[-t',0]$. (See
{\sc Fig.}~\ref{fig:poss1}.)
\item[(ii)] For every $A<\infty$ and every $t'<\infty$, for all $a$
sufficiently large all flow lines through points of $B_{\widetilde
G_a}(\tilde x_a,0,A)$ are defined on the interval $[-t', 0]$.
\end{enumerate}

\begin{figure}[ht]
  \centerline{\epsfbox{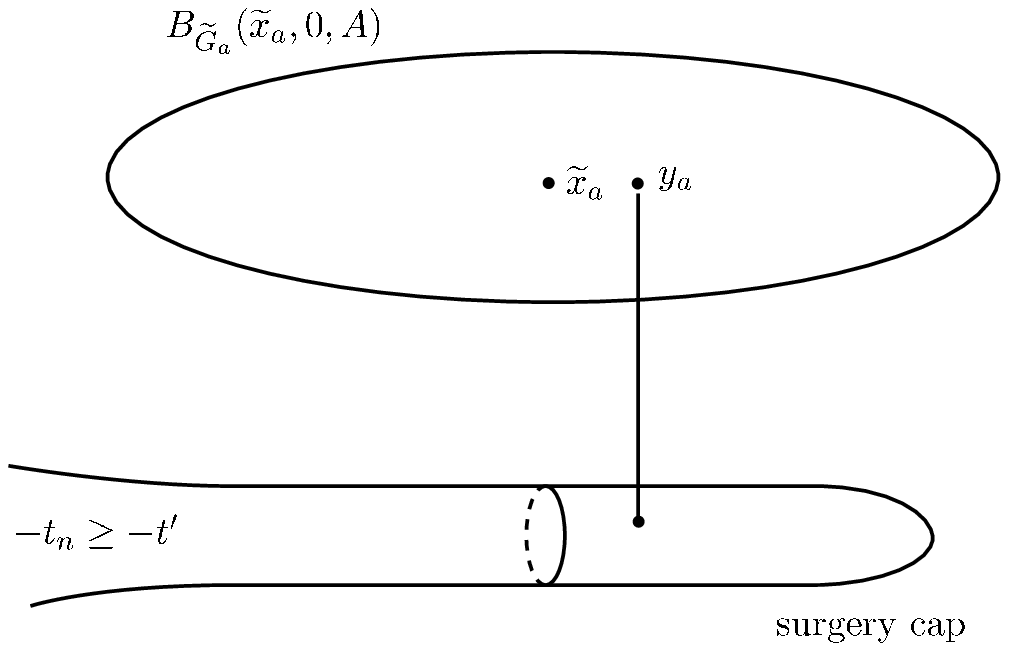}}
  \caption{Possibility (i).}\label{fig:poss1}
\end{figure}

Let us consider the second case. By Proposition~\ref{KAPPALIMIT}
these rescaled solutions are $\kappa$-non collapsed on scales $\le
\epsilon R_{G_a}(x_a)^{1/2}$ for all $t<0$. Since this condition is
a closed constraint, the same is true if $t=0$. Since $R(x_a)\ge
r_a^{-2}$, by construction every point $\tilde
x\in(\widetilde{\mathcal M}_a,\widetilde G_a)$ with $R(\tilde x)\ge
1$ and ${\bf t}(\tilde x)<0$ has a strong $(C,\epsilon)$-canonical
neighborhood.

\begin{claim} For all $a$ sufficiently large,
every point $\tilde x\in(\widetilde{\mathcal M}_a,\widetilde G_a)$ with $R(\tilde x)>1$ and
${\bf t}(\tilde x)=0$ has a $(2C,2\epsilon)$-canonical neighborhood.
\end{claim}

\begin{proof}
Assume that $\tilde x\in \widetilde{\mathcal M}$ has $R(\tilde
x)>1$. Suppose that $\tilde x$ is an exposed point. If $a$ is
sufficiently large, then $\delta_{b(a)}$ is arbitrarily close to
zero and hence by the last item in Theorem~\ref{LOCALSURGERY} and
the structure of the standard initial condition,  we see that
$\tilde x$ is contained in the core of a $(2C,2\epsilon)$-cap.

Suppose now that $\tilde x$ is not an exposed point.
 Then we can take a sequence of points
$\tilde y_n\in \widetilde{\mathcal M}_a$ all lying on the flow line for the vector field through $\tilde x$
converging to $\tilde x$ with ${\bf t}(\tilde y_n)<0$. Of course, for all $n$ sufficiently large
$R(\tilde y_n)>1$,
which implies that for all $n$ sufficiently large $\tilde y_n$ has a strong $(C,\epsilon)$-canonical
neighborhood. Passing to a subsequence, we can arrange that all of these canonical neighborhoods
are of the same type. If they are all $\epsilon$-round components, all $C$-components, or all
$(C,\epsilon)$-caps whose cores contain $y_n$, then by taking limits and arguing as in the proof of
Lemma~\ref{omegacanon} we see that $\tilde x$ has a strong
$(2C,2\epsilon)$-canonical neighborhood of the same type. On the other hand, if $\tilde y_n$ is the center
of a strong $\epsilon$-neck for all $n$, then according to Claim~\ref{epslimit},
the limit point $\tilde x$ is the center
of a strong $2\epsilon$-neck.
\end{proof}

Since we have chosen $\epsilon>0$ sufficiently  small so that Theorem~\ref{kaplimit}
applies with $\epsilon$ replaced by $2\epsilon$,  applying this theorem shows
that we can pass to a subsequence and take a smooth limiting
flow of a subsequence of the rescaled flows $(\widetilde {\mathcal
M}_a,\widetilde G_a)$ based at $\tilde x_a$ and defined for all
$t\in (-\infty,0]$. Because the $({\mathcal M}_a,G_a)$ all have
curvature pinched toward positive and since $R(x_a)\rightarrow \infty$
as $a$ tends to infinity, this result says that the limiting flow
has non-negative, bounded curvature and is $\kappa$-non-collapsed on
all scales. That is to say, the limiting flow is a
$\kappa$-solution. By Corollary~\ref{limitcannbhd} this contradicts
the fact that the strong $(C,\epsilon)$-canonical neighborhood
assumption fails at $x_a$ for every $a$. This contradiction shows
that in the second case there is a subsequence of the $a$ such that
$x_a$ has a strong canonical neighborhood and completes the proof of the second case.

Let us consider the first case. In this case we will arrive at a
contradiction by showing that for all $a$ sufficiently large, the
point $x_a$ lies in a strong $(C,\epsilon)$-canonical neighborhood
coming from a surgery cap. Here is the basic result we use to find
that canonical neighborhood.

\begin{lem}\label{lemstdcan}
Suppose that there are $A',D',t'<\infty$ such that the following
holds for all $a$ sufficiently large.  There is a point $y_a\in
B_{\widetilde G_a}(\tilde x_a,0,A')$ and a  flow line of $\chi$
 beginning at $y_a$, defined for backward time and ending at a
 point $z_a$ in a
surgery cap ${\mathcal C}_a$ at time $-t_a$ for some $t_a\le t'$. We
denote this flow line by $y_a(t), -t_a\le t\le 0$. Furthermore,
suppose that the scalar curvature on  the  flow line from $y_a$ to
$z_a$ is bounded by $D'$. Then for all $a$ sufficiently large, $x_a$
has a strong $(C,\epsilon)$-canonical neighborhood.
\end{lem}

\begin{proof}
  The proof is by contradiction. Suppose
the result does not hold. Then there are $A',D',t'<\infty $ and we
can pass to a subsequence (in $a$) such that  the hypotheses of the
lemma hold for every $a$ but no $x_a$ has a strong
$(C,\epsilon)$-canonical neighborhood. The essential point of the
argument is to show that {\bf in the units of the surgery scale} the
elapsed time between the surgery time and $0$ is less than $1$ and
the distance from the point $z_a$ to the tip of the surgery cap is
bounded independent of $a$.

 By Lemma~\ref{A_0}, the fact
that the scalar curvature at $z_a$ is bounded by $D'$ implies that
for all $a$ sufficiently large the scale $\bar h_a$ of the surgery
at time $-t_a$ satisfies
\begin{equation}\label{hequat}
\bar h_a^2\ge (2D'D)^{-1}.
\end{equation}
(Recall that we are working in the rescaled flow $(\widetilde {\mathcal M}_a,\widetilde G_a)$.)

Now we are ready to show that the elapsed time is bounded less than
one in the surgery scale.

\begin{claim}
There is $\theta_1<1$, depending on $D'$ and $t'$, such that for all
$a$ sufficiently large we have  $t_a< \theta_1\bar h_a^2$.
\end{claim}

\begin{proof}
We consider two cases: either $t_a\le \bar h_a^2/2$ or $\bar
h_a^2/2<t_a$. In the first case, the claim is obviously true with
$\theta_1$ anything greater than $1/2$ and less than one. In the
second case, the curvature everywhere along the flow line is at most
$ D'<(2t_aD')\bar h_a^{-2}\le (2t'D')\bar h_a^{-2}$. Using
Proposition~\ref{Restim} fix $1/2<\theta_1<1$ so that every point of
the standard solution $(x,t)$ with $t\ge (2\theta_1-1)$ satisfies
$R(x,t)\ge 6t'D'$. Notice that  $\theta_1$ depends only on $D'$ and
$t'$. If $t_a<\theta_1\bar h_a^2$, then the claim holds for this
value of $\theta_1<1$. Suppose $t_a\ge \theta_1\bar h_a^2$, so
that $-t_a+(2\theta_1-1)\bar h_a^2<0$. For all  $a$
sufficiently large we have $\delta_a\le \delta_0''(A_0+5,
\theta_1,\overline \delta_0)$ where $\overline\delta_0$ is the constant
from Definition~\ref{defnoverlinedelta0} and $\delta''_0$ is the constant
from Proposition~\ref{altern}.
This means that the scalar curvatures
at corresponding points of the rescaled standard solution and the
evolution of the surgery cap (up to time $0$) in $\widetilde{\mathcal M}_a$
differ by at most a factor of two. Thus, for these $a$, we have
$R(y_a,(-t_a+(2\theta_1-1)\bar h_a^{2}))\ge 3(t'D')\bar h_a^{-2}$
from the definition of $\overline\delta_0$ and
Proposition~\ref{altern}. But this is impossible since
$-t_a(2\theta_1-1)\bar h_a^3<0$ and $3t'D'/\bar h_a^2\ge 3t_aD'/\bar h_a^2>D'$ as
 $t_a\ge \bar h_a^2/2$. Hence,
$R(y_a,(-t_a+(2\theta_1-1)\bar h_a^2))\le 2t_aD'\bar
h_a^{-2}\le 2t'D'\bar h_a^{-2}$. This contradiction
shows that if $a$ is sufficiently large then $t_a<\theta_1\bar
h_a^2$.
\end{proof}

We pass to a subsequence so that $t_a\bar h_a^{-2}$ converges to
some $\theta\le \theta_1$. We define $\widetilde C$ to be the
maximum of $C$ and $3\epsilon^{-1}\beta^{-1}$. Now, using Part 5 of
Theorem~\ref{stdsoln} we set $A''\ge (9\widetilde
C+3A')\sqrt{2DD'}+6(A_0+5)$ sufficiently large so that in the
 standard flow $B(p_0,0,A'')$ contains
$B(p_0,t,A''/2)$ for any $t\le (\theta_1+1)/2$. This constant is chosen only to
depend on $\theta_1$, $A'$, and $C$. As $a$ tends to infinity, $\delta_a$ tends
to zero which means, by Proposition~\ref{altern}, that for all $a$ sufficiently
large there is an embedding $\rho_a\colon B(p_0,-t_a,A''\bar h_a)\times
[-t_a,0]\to \widetilde {\mathcal M}_a$ compatible with time and the vector
field such that (after translating by $t_a$ to make the flow start at time $0$
and scaling by $\bar h_a^{-2}$) the restriction of $\widetilde G_a$ to this
image is close in the $C^\infty$-topology to the restriction of the standard
flow to $B(p_0,0,A'')\times [0,\bar h_a^{-2}t_a]$. The image $\rho_a(p_0,-t_a)$
is the tip $p_a$ of the surgery cap ${\mathcal C}_a$ in $\widetilde {\mathcal
M}_a$. In particular, for all $a$ sufficiently large the image
$\rho_a\left(B(p_0,-t_a,A''\bar h_a)\times\{0\}\right)$ contains the $A''\bar
h_a/3$-neighborhood of the image $\rho_a(p_0,0)$ of the tip of the surgery cap
under the flow forward to  time $0$. By our choice of $A''$, and
Equation~(\ref{hequat}), this means that for all $a$ sufficiently large
$\rho_a\left(B(p_0,-t_a,A''\bar h_a)\times\{-t_a\}\right)$ contains the
$(3\widetilde C+A') +2(A_0+5)\bar h_a$ neighborhood of $p_a=\rho_a(p_0,-t_a)$.
Notice also that, since the standard solution has positive curvature and hence
the distance between points is non-increasing in time by
Lemma~\ref{posdistdec}, the distance at time $0$ between $\rho_a(p_0,0)$ and
$y_a$ is less than $2(A_0+5)\bar h_a$. By the triangle inequality, we conclude
that for all $a$ sufficiently large, $\rho_a\left(B(p_0,-t_a,A''\bar h_a)\times
\{0\}\right)$ contains the $3\widetilde C$-neighborhood of $x_a$. Since the
family of metrics on $\rho_a\left(B(p_0,-t_a,A''\bar h_a)\times
[-t_a,0]\}\right)$ (after time-shifting by $t_a$ and rescaling by $\bar
h_a^{-2}$) are converging smoothly to the ball $B(p_0,0,A'')\times [0,\theta]$
in the standard flow,
 for all $a$ sufficiently large then the flow from time $-t_a$ to $0$ on the
$3\widetilde C$-neighborhood of $x_a$ is, after rescaling by $\bar
h_a^{-2}$,
 very nearly isometric to the restriction of the standard flow from time $0$ to $\bar h^{-2}_at_a$ on
 the $3\widetilde C\bar h_a^{-1}$-neighborhood of some
point $q_a$ in the standard flow. Of course, since the scalar
curvature of $x_a$ is one,  $R(q_a,\bar h_a^{-2}t_a)$ in the standard flow is
close to $\bar h_a^{-2}$. Hence, by Theorem~\ref{stdsolncannbhd}
there is a neighborhood $X$ of $(q_a,\bar h_a^{-2}t_a)$ in the standard solution
that either is a $(C,\epsilon)$-cap, or is an evolving
$\beta\epsilon/3$-neck centered at $(q_a,\bar h_a^{-2}t_a)$. In the latter case
either the evolving neck is defined for backward time
$(1+\beta\epsilon/3)$ or its initial time-slice is the zero
time-slice and this initial time-slice lies at distance at least $1$
from the surgery cap. Of course, $X$ is contained in the
$CR(q_a,\bar h_a^{-2}t_a)^{-1/2}$ neighborhood of $(q_a,\bar h_a^{-2}t_a)$ in the standard
solution. Since $\widetilde C\ge C$ and $R(q_a,\bar h_a^{-2}t_a)$ is close to
$\bar h_a^{-2}$, the neighborhood $X$ is contained in the
$2\widetilde C\bar h_a^{-1}$-neighborhood of $(q_a,\bar h_a^{-2}t_a)$ in the
standard solution.  Hence, after rescaling, the corresponding
neighborhood of $x_a$ is contained in the $3\widetilde
C$-neighborhood of $x_a$. If either of the first two cases in Theorem~\ref{stdsolncannbhd} occurs
for a subsequence of $a$ tending to infinity, then  by
Lemma~\ref{canonvary} and the fact that $\widetilde C>{\rm
max}(C,\epsilon^{-1})$, we see that there is a subsequence of $a$
for which $x_a$ either is contained in the core of a
$(C,\epsilon)$-cap or is the center of a strong $\epsilon$-neck.

We must examine further the last case. We suppose that for every $a$
this last case holds. Then for all $a$ sufficiently large we have an
$\beta\epsilon/3$-neck $N_a$ in the zero time-slice of $\widetilde{\mathcal
M}_a$ centered at $x_a$. It is an evolving neck and there is an
embedding $\psi\colon N_a\times [-t_a,0]\to \widetilde{\mathcal
M}_a$ compatible with time and the vector field so that the initial
time-slice $\psi(N_a\times \{-t_a\})$ is in the surgery time-slice
$M_{-t_a}$ and is disjoint from the surgery cap, so in fact it is
contained in the continuing region at time $-t_a$. As we saw above,
the image of the central $2$-sphere $\psi(S^2_a\times\{-t_a\})$ lies
at distance at most  $ A''\bar h_a$ from the tip of the surgery cap
$p_a$ (where, recall, $A''$ is a constant independent of $a$). The
$2$-sphere, $\Sigma_a$, along which we do surgery,   creating the surgery
cap with $p_a$ as its tip, is the central $2$-sphere of a strong
$\delta_{b(a)}$-neck. As $a$ tends to infinity the surgery control
parameter $\delta_{b(a)}$ tends to zero. Thus, for $a$ sufficiently
large this strong $\delta_{b(a)}$-neck will contain a strong
$\beta\epsilon/2$- neck $N'$ centered at $\psi(x_a,-t_a)$. Since we
know that the continuing region at time $-t_a$ contains a
$\beta\epsilon/3$-neck centered at $(x_a,-t_a)$, it follows that
$N'$ is also contained in $C_{-t_a}$. That is to say, $N'$  is
contained in the negative half of the $\delta_{b(a)}$-neck centered
at $\Sigma_a$. Now we are in the situation of
Proposition~\ref{neckglue}. Applying this result tells us that $x_a$
is the center of a strong $\epsilon$-neck.

This completes the proof that for all $a$ sufficiently large, $x_a$
has a $(C,\epsilon)$-canonical neighborhood in contradition to our
assumption. This contradiction completes the proof of
Lemma~\ref{lemstdcan}.
\end{proof}

There are several steps required to complete the proof of
Proposition~\ref{extend}. The first step helps us apply the previous claim to
find strong $(C,\epsilon)$-canonical neighborhoods.

\begin{claim}\label{D(A)delta(A)}
Given any $A<\infty$ there is $D(A)<\infty$ and $\delta(A)>0$ such that for all $a$
sufficiently large, $|{\rm Rm}|$ is bounded by $D(A)$ along all backward flow lines
beginning at a point of $B_{\widetilde
G_a}(\tilde x_a,0,A)$ and defined for backward time at most $\delta(A)$.
\end{claim}

\begin{proof}
Since all points  $y\in ({\mathcal M}_a,G_a)$ with $R_{G_a}(y)\ge
r^{-2}_a$ and ${\bf t}(y)<{\bf t}(x_a)$ have strong
$(C,\epsilon)$-canonical neighborhoods, and since $R(x_a)=r^{-2}_a$,
we see that all points $y\in (\widetilde{\mathcal M}_a,\widetilde
G_a)$ with ${\bf t}(y_a)<0$ and with $R_{\tilde G_a}(y_a)\ge 1$ have
strong $(C,\epsilon)$-canonical neighborhoods. It follows that all
points in $(\widetilde {\mathcal M}_a,\widetilde G_a)$ with ${\bf
t}(y)\le 0$ and $R(y)>1$ have strong $(2C,2\epsilon)$-canonical
neighborhoods. Also, since $\delta_a\le \delta(r_a)$, where
$\delta(r_a)$ is the constant given in Proposition~\ref{KAPPALIMIT},
and since the condition of being $\kappa$-non-collapsed is a closed
constraint, it follows from Proposition~\ref{KAPPALIMIT} that these
Ricci flows with surgery are $\kappa$-non-collapsed for a fixed
$\kappa>0$. It is now immediate from Theorem~\ref{bcbd} that there
is a constant $D_0(A)$ such that $R$ is bounded above on
$B_{\widetilde G_a}(\tilde x_a,0,A)$ by $D_0(A)$. Since every point
$y\in ({\mathcal M}_a,G_a)$ with $R(y)>1$ of the sequence of with
scalar curvature at least $1$ has a $(C,\epsilon)$ canonical
neighborhood, it follows from the definition that for every such
point $y$ we have $\left|\partial R(y)/\partial t\right|<CR(y)^2$.
Arguing as in Lemma~\ref{controlnbhd} we see that there is a
constant $\delta(A)>0$ and a bound $D'(A)$, both depending only in
$D_0(A)$,
 for the scalar curvature at all points
of backward flow lines beginning in $B_{\widetilde
G_a}(\tilde x_a,0,A)$ and defined for backward time at most $\delta(A)$. Since the curvature is pinched toward positive,
it follows that there is a bound $D(A)$ depending only on $D'(A)$ to $|{\rm Rm}|$ on the same flow lines.
\end{proof}

\begin{claim} After passing to a subsequence (in $a$), either:
\begin{enumerate}
\item[(1)]
for each $A<\infty$ there are $D(A)<\infty$ and $t(A)>0$ such that for all $a$
sufficiently large $P_{\widetilde G_a}(\tilde x_a,0,A,-t(A))$ exists in
$\widetilde {\mathcal M}_a$ and $|{\rm Rm}|$ is bounded by $D(A)$ on this
backward parabolic neighborhood, or
\item[(2)] each $x_a$ has a strong
$(C,\epsilon)$-canonical neighborhood.
\end{enumerate}
\end{claim}

\begin{proof}
First notice that if there is $t(A)>0$ for which the backwards parabolic neighborhood
$P=P_{\widetilde G_a}(\tilde x_a,0,A,-t(A))$ exists,
then, by Claim~\ref{D(A)delta(A)},
there are constants $D(A)<\infty$ and $\delta(A)>0$ such that, replacing $t(A)$ by ${\rm min}(t(A),\delta(A))$,
 $|{\rm Rm}|$ is bounded by $D(A)$ on $P$.
Thus, either Item (1) holds or
passing to a subsequence, we can suppose that there is some $A<\infty$ for which no $t(A)>0$ as
required by Item (1) exists.
Then,  for each $a$ we find a point
$y_a\in B_{\widetilde G_a}(\tilde x_a,0,A)$ such that the backwards
flow line from $y_a$ meets a surgery cap at a time $-t_a$ where
${\rm lim}_{a\rightarrow \infty}(t_a)=0$. Then, by the previous claim, for all $a$ sufficiently large, the
sectional curvature along any backward flow line beginning in $B_{\widetilde G_a}(\widetilde x_a,0,A)$
and defined for backward time $t_a$ is bounded by
a constant $D(A)$ independent of $a$.  Under our assumption this
means that for all $a$ sufficiently large, there is a point $y_a\in
B_{\widetilde G_a}(\tilde x_a,0,A)$ and a backwards flow line
starting at $y_a$ ending at a point $z_a$ of a surgery cap, and the
sectional curvature along this entire flow line is bounded by $D(A)<\infty$. Thus,
applying Lemma~\ref{lemstdcan} produces the strong
$(C,\epsilon)$-canonical neighborhood around $x_a$, proving the claim.
\end{proof}

But we are assuming that no $x_a$ has a strong $(C,\epsilon)$-canonical
neighborhood. Thus, the consequence of the previous claim is that for each
$A<\infty$ there is a $t(A)>0$ such that for all $a$ sufficiently large
$P_{\widetilde G_a}(\tilde x_a,0,A,-t(A))$ exists in $\widetilde {\mathcal
M}_a$ and there is a bound, depending only on $A$ for $|{\rm Rm}|$ on
this backward parabolic neighborhood. Applying Theorem~\ref{sliceconv} we see that, after passing to a
subsequence, there is a smooth limit $(M_\infty,g_\infty,x_\infty)$ to the zero
time-slices $(\widetilde{\mathcal M}_a,\widetilde G_a,\tilde x_a)$.
Clearly, since the curvatures of the sequence are pinched toward positive, this
limit has non-negative curvature.

Lastly, we show that $(M_\infty,g_\infty)$ has bounded curvature. By Part 3 of Proposition~\ref{canonvary}
each point of $(M_\infty,g_\infty)$ with scalar curvature
greater than one has a $(2C,2\epsilon)$-canonical neighborhood. If a point lies in an $2\epsilon$-component
or in a $2C$-component, then $M_\infty$ is compact, and hence clearly has bounded curvature.
Thus, we can assume that each $y\in M_\infty$ with $R(y)>1$ is either the center of a $2\epsilon$-neck
or is contained in the core of a $(2C,2\epsilon)$-cap. According to Proposition~\ref{narrows}
 $(M_\infty,g_\infty)$ does not contain $2\epsilon$-necks of arbitrarily high curvature.
 It now follows then that $(M_\infty,g_\infty)$ there is a bound to the scalar curvature of any
 $2\epsilon$-neck and of any $(2C,2\epsilon)$-cap, and hence it follows that $(M_\infty,g_\infty)$
 has bounded curvature.

\begin{claim}
If the constant $t(A)>0$ cannot be chosen independent of $A$, then
after passing to a subsequence, the $x_a$ have strong
$(C,\epsilon)$-canonical neighborhoods.
\end{claim}

\begin{proof}
 Let $Q$ be the bound of the scalar curvature of
$(M_\infty,g_\infty,x_\infty)$. Then by Lemma~\ref{controlnbhd}
 there is a constant
$\Delta t>0$ such that if $R_{\widetilde G_a}(y,0)\le 2Q$, then the
scalar curvature is bounded by $16Q$ on the backward flow line from
$y$ defined for any time $\le \Delta t$. Suppose that there is
$A<\infty$ and a  subsequence of $a$ for which there is a flow line
beginning at a point $y_a \in B_{\widetilde G_a}(\tilde x_a,0,A)$
defined for backward time at most $\Delta t$ and ending at a point
$z_a$ of a surgery cap. Of course, the fact that the scalar
curvature of $(M_\infty,g_\infty)$ is at most $Q$ implies that for
all $a$ sufficiently large, the scalar curvature of $B_{\widetilde
G_a}(\tilde x_a,0,A)$ is less than $2Q$. This implies  that for all
$a$ sufficiently large the scalar curvature along the flow line
from $y_a$ to $z_a$ in a surgery cap is $\le 16Q$. Now invoking
Lemma~\ref{lemstdcan} we see that for all $a$ sufficiently large the
point $\widetilde x_a$ has a strong $(C,\epsilon)$-canonical
neighborhood. This is a contradiction, and this contradiction proves
that we can choose $t(A)>0$ independent of $A$.
\end{proof}

Since we are assuming that no $x_a$ has a strong
$(C,\epsilon)$-canonical neighborhood, this means that it is
possible to find a constant $t'>0$ such that $t(A)\ge t'$ for all
$A<\infty$. Now let $0<T'\le \infty$ be the maximum possible value for such
$t'$. Then for every $A$ and every $T<T'$ the
parabolic neighborhood $P_{\widetilde G_a}(\tilde x_a,0,A,T)$ exists
for all $a$ sufficiently large. According to
Theorem~\ref{kaplimit}, after passing to a subsequence, there is
a limiting flow $(M_\infty,g_\infty(t),x_\infty),\ -T'<t\le 0$, and
this limiting flow has bounded, non-negative curvature.
If $T=\infty$, this limit is a $\kappa$-solution, and hence the $x_a$ have
strong $(C,\epsilon)$-canonical neighborhoods for all $a$ sufficiently large,
which is a contradiction.

Thus, we can assume that $T'<\infty$.
 Let $Q$ be
the bound for the scalar curvature of this flow. Since $T'$ is
maximal, for every $t>T'$, after passing to a subsequence, for all
$a$ sufficiently large there is $A<\infty$ and a backwards flow
line, defined for a time less than $t$, starting at a point $y_a$ of
$B_{\widetilde G_a}(\tilde x_a,0,A)$ and ending at a point $z_a$ of
a surgery cap. Invoking Lemma~\ref{controlnbhd} again, we see that
for all $a$ sufficiently large, the scalar curvature is bounded on
the flow line from $y_a$ to $z_a$ by a constant independent of $a$.
Hence, as before, we see that for all $a$ sufficiently large $x_a$
has a strong $(C,\epsilon)$-canonical neighborhood; again this is  a contradiction.

Hence, we have now shown  that our assumption that the strong
$(C,\epsilon)$-canonical neighborhood assumption fails for all $r_a$
and all $\delta_{a,b}$ leads to a contradiction and hence is false.

This completes the proof of Proposition~\ref{extend}.
\end{proof}

\section{Surgery times don't accumulate}\label{sectcomplete}

Now we turn to the proof of Theorem~\ref{MAIN}.
 Given surgery parameter
sequences $$\Delta_i=\{\delta_0,\ldots,\delta_i\}$$
$${\bf r_i}=\{r_0,\ldots,r_i\}$$
$${\bf K_i}=\{\kappa_0,\ldots,\kappa_i\},$$
 we let
$r_{i+1}$ and $\delta_{i+1}$ be as in Proposition~\ref{extend} and
then set $\kappa_{i+1}=\kappa(r_{i+1})$ as in
Proposition~\ref{KAPPALIMIT}.
 Set $${\bf
r_{i+1}}=\{{\bf r_i},r_{i+1}\}$$
$${\bf K_{i+1}}=\{{\bf
K_i},\kappa_{i+1}\}$$
$$\Delta_{i+1}=\{\delta_0,\ldots,\delta_{i-1},\delta_{i+1},\delta_{i+1}\}.$$
Of course, these are also surgery parameter sequences.

Let $\overline \delta\colon [0,T]\to \Ar^+$ be any non-increasing positive
function and let $({\mathcal M},G)$ be a Ricci flow with surgery defined on
$[0,T)$ for some $T\in [T_i,T_{i+1})$ with surgery control parameter
$\overline\delta$. Suppose $\overline \delta\le \Delta_{i+1}$ and that this
Ricci flow with surgery satisfies the conclusion of Theorem~\ref{MAIN} with
respect to these sequences on its entire interval of definition. We wish to
extend this Ricci flow with surgery to one defined on $[0,T')$ for some $T'$
with $T<T'\le T_{i+1}$ in such a way that $\overline\delta$ is the surgery
control parameter and the extended Ricci flow with surgery continues to satisfy
the conclusions of Theorem~\ref{MAIN} on its entire interval of definition.

We may as well assume that the Ricci flow $({\mathcal M},G)$ becomes
singular at time $T$. Otherwise we would simply extend by Ricci flow
to a later time $T'$. By Proposition~\ref{KAPPALIMIT} and
Proposition~\ref{extend} this extension will continue to satisfy the
conclusions of Theorem~\ref{MAIN} on its entire interval of
definition. If $T\ge T_{i+1}$, then we have extended the Ricci flow
with surgery to time $T_{i+1}$ as required and hence completed the
inductive step. Thus, we may as well assume that $T<T_{i+1}$.

Consider the maximal extension of $({\mathcal M},G)$ to time $T$. Let $T^-$ be
the previous surgery time, if there is one, and otherwise be zero. If the $T$
time-slice, $\Omega(T)$, of this maximal extension is all of $M_{T^-}$, then
the curvature remains bounded as $t$ approaches $T$ from below. According to
Proposition~\ref{flowextend} this means that $T$ is not a surgery time and we
can extend the Ricci flow on $(M_{T^-},g(t)),\ T^-\le t<T$, to a Ricci flow on
$(M_{T^-},g(t)),\ T^-\le t<T'$ for the maximal time interval (i.e. so that the
flow becomes singular at time $T'$ or $T'=\infty$). But we are assuming that
the flow goes singular at  $T$. That is to say,  $\Omega(T)\not= M_{T^-}$. Then
we can do surgery at time $T$ using $\overline \delta(T)$ as the surgery
control parameter, setting $\rho(T)=r_{i+1}\delta(T)$. Let $(M_T,G(T))$ be the
result of surgery. If $\Omega_{\rho(T)}(T)=\emptyset$, then the surgery process
at time $T$ removes all of $M_{T'}$. In this case, the Ricci flow is understood
to exist for all time and to be empty for $t\ge T$. In this case we have
completed the extension to $T_{i+1}$, and in fact all the way to $T=\infty$,
and hence completed the inductive step in the proof of the proposition.

We may as well assume that $\Omega_{\rho(T)}(T)\not=\emptyset$ so
that the result of surgery is a non-empty manifold $M_T$.
 Then we use this compact Riemannian
$3$-manifold as the initial conditions of a Ricci flow beginning
at time $T$. According to Lemma~\ref{14.26} the union along
$\Omega(T)$ at time $T$ of this Ricci flow with $({\mathcal M},G)$
is a Ricci flow with surgery satisfying Assumptions (1) -- (7) and whose
curvature is pinched toward positive.

Since the surgery  control parameter $\overline\delta(t)$ is at most
$\delta(r_{i+1})$, the constant from Proposition~\ref{KAPPALIMIT},
 for all $t\in [T_{i-1},T]$, since $T\le T_{i+1}$,
and since the restriction of $({\mathcal M},G)$ to ${\bf
t}^{-1}([0,T_i))$ satisfies Proposition~\ref{KAPPALIMIT}, we see by
Proposition~\ref{extend} that the extended Ricci flow with surgery
satisfies the conclusion of Theorem~\ref{MAIN} on its entire time
interval of definition.

Either we can repeatedly apply this process, passing from one surgery
time to the next and eventually reach $T\ge T_{i+1}$, which would
prove the inductive step, or there is an unbounded number of
surgeries in the time interval $[T_i,T_{i+1})$. We must rule out the
latter case.

\begin{lem}\label{bddsurgery}
Given a Ricci flow with surgery $({\mathcal M},G)$ defined on $[0,T)$ with
$T\le T_{i+1}$ with surgery control parameter $\overline\delta$ a
non-increasing positive function defined on $[0,T_{i+1}]$ satisfying the
hypotheses of Theorem~\ref{MAIN} on its entire time-domain of definition, there
is a constant $N$ depending only on the volume of $(M_0,g(0))$, on $T_{i+1}$,
on $r_{i+1}$, and on $\overline\delta(T_{i+1})$ such that this Ricci flow with
surgery defined on the interval $[0,T)$ has at most $N$ surgery times.
\end{lem}

\begin{proof}
Let $(M_t,g(t))$  be the $t$ time-slice of $({\mathcal M},G)$. If $t_0$ is not
a surgery time, then ${\rm Vol}(t)={\rm Vol}(M_t,g(t))$ is a smooth function
of $t$ near $t_0$ and
$$\frac{d{\rm
Vol}}{dt}(t_0)=-\int_{M_{t_0}}Rd{\rm vol},$$ so that, because of the curvature
pinching toward positive hypothesis, we have $\frac{d{\rm Vol}}{dt}(t_0)\le
6{\rm Vol}(t_0)$. If $t_0$ is a surgery time, then either $M_{t_0}$ has fewer
connected components than  $M_{t_0^-}$ or we do a surgery in an $\epsilon$-horn
of $M_{t_0^-}$. In the latter case we remove the end of the $\epsilon$-horn,
which  contains the positive half of a $\overline\delta(t_0)$-neck of scale
$h(t_0)$. We then  sew in a ball with volume at most $(1+\epsilon)Kh^3(t_0)$,
where $K<\infty$ is the universal constant given in Lemma~\ref{A_0}. Since
$h(t_0)\le \overline\delta^2(t_0)r(t_0)\le \delta_0^2r(t_0)$ and since we have
chosen
 $\overline\delta(t_0)\le\delta_0<K^{-1}$, it follows that this
operation lowers volume by at least $\delta^{-1}h^2(t_0)/2$. Since
$\overline\delta(t_0)\ge \overline \delta(T_{i+1})>0$ and the canonical
neighborhood parameter $r$ at time $t_0$ is at least $r_{i+1}>0$, it follows
that $h(t_0)\ge h(T_{i+1})>0$. Thus, each surgery at time $ t_0\le T_{i+1}$
along a $2$-sphere removes at least a fixed amount of volume depending on
$\overline \delta(T_{i+1})$ and $r_{i+1}$. Since under Ricci flow the volume
grows at most exponentially, we see that there is a bound depending only on
$\overline\delta(T_{i+1})$, $T_{i+1}$, $r_{i+1}$ and ${\rm Vol}(M_0,g(0))$ to
the number of $2$-sphere surgeries that we can do in this time interval. On the
other hand, the number of components at any time $t$ is at most $N_0+S(t)-D(t)$
where $N_0$ is the number of connected components of $M_0$, $S(t)$ is the
number of $2$-sphere surgeries performed in the time interval $[0,t)$ and
$D(t)$ is the number of connected components removed by surgeries at times in
the interval $[0,t)$.  Hence, there is a bound on the number of components in
terms of $N_0$ and $S(T)$ that can be removed by surgery in the interval
$[0,T)$. Since the initial conditions are normalized, $N_0$ is bounded by the
volume of  $(M_0,g(0))$. This completes the proof of the result.
\end{proof}

This lemma completes the proof of the fact that for any $T\le T_{i+1}$, we
encounter only a fixed bounded number  surgeries in  the Ricci flow with
surgery from $0$ to $T$. The bound depends on the volume of the initial
manifold as well as the surgery constants up to time $T_{i+1}$. In particular,
for a given initial metric $(M_0,g(0))$ there is a uniform bound, depending
only on the surgery constants up to time $T_{i+1}$, on the number of surgeries
in any Ricci flow with surgery defined on a subinterval of $[0,T_{i+1})$. It
follows that the surgery times cannot accumulate in any finite interval. This
completes the proof of Theorem~\ref{MAIN}.

To sum up, we have sequences $\Delta$, ${\bf K}$ and ${\bf r}$ as given in
Theorem~\ref{MAIN}. Let $\overline \delta\colon [0,\infty)\to \Ar$ be a
positive, non-increasing function with $\overline \delta\le \Delta$. Let $M$ be
a compact $3$-manifold that contains no embedded $\Ar P^2$ with trivial normal
bundle. We have proved that for any normalized initial Riemannian metric
$(M_0,g_0)$ there is a Ricci flow with surgery with time-interval of definition
$[0,\infty)$ and with $(M_0,g_0)$ as initial conditions. This Ricci flow with
surgery is ${\bf K}$-non-collapsed and satisfies the strong
$(C,\epsilon)$-canonical neighborhood theorem with respect to the parameter
${\bf r}$. It also has curvature pinched toward positive. Lastly, for any $T\in
[0,\infty)$ if there is a surgery at time $T$ then this surgery is performed
using the surgery parameters $\overline\delta(T)$ and $r(T)$, where if $T\in
[T_i,T_{i+1})$ then $r(T)=r_{i+1}$. In this Ricci flow with surgery, there are
only finitely many surgeries on each finite time interval. As far as we know
there may be infinitely many surgeries in all.

\chapter{Finite-time extinction}\label{energy}

Our purpose in this chapter is to prove the following finite-time
extinction theorem for certain Ricci flows with surgery which, as we
shall show below, when combined with the theorem on the existence of
Ricci flows with surgery defined for all $t\in [0,\infty)$
(Theorem~\ref{MAIN}), immediately yields Theorem~\ref{Theorem1},
thus completing the proof of the Poincar\'e Conjecture and the
$3$-dimensional space-form conjecture.

\section{The result}

\begin{thm}\label{extinct}
Let $(M,g(0))$ be a compact, connected normalized Riemannian $3$-manifold.
Suppose that the fundamental group of $M$ is a free product of finite groups
and infinite cyclic groups. Then $M$ contains no $\Ar P^2$ with trivial normal
bundle.  Let $({\mathcal M},G)$ be the Ricci flow with surgery defined for all
$t\in [0,\infty)$ with $(M,g(0))$ as initial conditions given by
Theorem~\ref{MAIN}. This Ricci flow with surgery becomes extinct after a finite
time in the sense that the time-slices $M_T$ of ${\mathcal M}$ are empty for
all $T$ sufficiently large.
\end{thm}

Let us quickly show how this theorem implies our main result
Theorem~\ref{Theorem1}.

\begin{proof}(of Theorem~\ref{Theorem1} assuming
Theorem~\ref{extinct}). Fix a normalized metric $g(0)$ on $M$, and let
$({\mathcal M},G)$ be the Ricci flow with surgery defined for all $t\in
[0,\infty)$ produced by Theorem~\ref{MAIN} with initial conditions $(M,g(0))$.
According to Theorem~\ref{extinct} there is $T>0$ for which the time-slice
$M_T$ is empty. By Corollary~\ref{Tgood0good}, if there is $T$ for which $M_T$
is empty, then for any $T'<T$ the manifold $M_{T'}$ is a disjoint union of
connected sums of $3$-dimensional spherical space forms and $2$-sphere bundles
over $S^1$. Thus, the manifold $M=M_0$ is a connected sum of $3$-dimensional
space-forms and $2$-sphere bundles over $S^1$. This proves
Theorem~\ref{Theorem1}. In particular, if $M$ is simply connected, then $M$ is
diffeomorphic to $S^3$, which is the statement of the Poincar\'e Conjecture.
Similarly, if $\pi_1(M)$ is finite then $M$ is diffeomorphic to a connected sum
of a $3$-dimensional spherical space-form and $3$-spheres, and hence $M$ is
diffeomorphic to a $3$-dimensional spherical space-form.
\end{proof}

The rest of this chapter is devoted to the proof of
Theorem~\ref{extinct} which will then complete the proof of
Theorem~\ref{Theorem1}.

\subsection{History of this approach}

The basic idea for proving finite-time extinction is to use a min-max function
based on the area (or the closely related energy) of $2$-spheres or $2$-disks
in the manifold. The critical points of the energy functional are harmonic maps
and they play a central role in the proof. For a basic reference on harmonic
maps see \cite{SacksUhlenbeck}, \cite{SchoenYau}, and \cite{Jost}. Let us
sketch the argument. For a compact Riemannian manifold $(M,g)$ every non-zero
element $\bbeta\in \pi_2(M)$ has associated with it an area, denoted
$W_2(\bbeta,g)$, which is the infimum over all maps $S^2\to M$ in the free
homotopy class of $\bbeta$ of the energy of the map. We find it convenient to
set $W_2(g)$ equal to the minimum over all non-zero homotopy classes $\bbeta$
of $W_2(\bbeta,g)$. In the case of a Ricci flow $g(t)$ there is an estimate
(from above) for the forward difference quotient of $W_2(g(t))$ with respect to
$t$. This estimate shows that after a finite time $W_2(g(t))$ must go negative.
This  is absurd since $W_2(g(t))$ is  always non-negative. This means that the
Ricci flow cannot exist for all forward time. In fact, using the
distance-decreasing property for surgery in Proposition~\ref{surgerydist} we
see that, even in a Ricci flow with surgery, the same forward difference
quotient estimate holds for as long as  $\pi_2$ continues to be non-trivial,
i.e., is not killed by the surgery. The forward difference quotient estimate
means that eventually all of $\pi_2$ is killed in a Ricci flow with surgery and
we arrive at a time $T$ for which every component of the $T$ time-slice, $M_T$,
has trivial $\pi_2$. This result holds for all Ricci flows with surgery
satisfying the conclusion of Theorem~\ref{MAIN}.

Now we fix $T_0$ so that every component of $M_{T_0}$ has trivial
$\pi_2$. It follows easily from the description of surgery that the
same statement holds for all $T\ge T_0$. We wish to show that, under
the group-theoretic hypothesis of Theorem~\ref{extinct}, at some
later time $T'>T_0$ the time-slice $M_{T'}$ is empty. The argument
here is similar in spirit. There are two approaches. The first
approach is due to Perelman \cite{Perelman3}. Here, one represents a
non-trivial element in $\pi_3(M_{T_0},x_0)$ by  a non-trivial
element in $\pi_2(\Lambda M,*)$, where $\Lambda M$ is the free loop
space on $M$ and $*$ is the trivial loop at $x_0$. For any compact
family $\Gamma$ of homotopically trivial loops in $M$ we consider
the areas of minimal spanning disks for each of the loops in the
family and set $W(\Gamma)$ equal to the maximal area of these
minimal spanning disks. For a given element in $\gamma\in
\pi_2(\Lambda M)$ we set $W(\gamma)$ equal to
 the infimum over all representative $2$-sphere families $\Gamma$
for $\gamma$  of $W(\Gamma)$. Under Ricci flow, the forward
difference quotient of this invariant satisfies an inequality and
the distance-decreasing property of surgery
(Proposition~\ref{surgerydist}) says that the inequality remains
valid for Ricci flow with surgery. The inequality implies that the
value $W(\gamma)$ goes negative in finite time, which is impossible.

The other approach, by Colding-Minicozzi \cite{ColdingMinicozzi}, is
to represent a non-trivial element in $\pi_3(M_T)$ as a non-trivial
element in $\pi_1({\rm Maps}(S^2,M))$, and associate to such an
element the infimum over all representative families of the maximal
energy of the $2$-spheres in the family. Again, one shows that under
Ricci flow the forward difference quotient of this minimax satisfies
an inequality that implies that it goes negative in finite time. As
before, the distance-decreasing property of surgery
(Proposition~\ref{surgerydist}) implies that this inequality is
valid for Ricci flows with surgery. This tells us that the manifold
must completely disappear in finite time.

Our first reaction was that, of the two approaches, the one
considered by Colding-Minicozzi was preferable since it seemed more
natural and it had the advantage of avoiding the boundary issues
that occupy most of Perelman's analysis in \cite{Perelman3}.
 In the Colding-Minicozzi approach one must construct paths of $2$-spheres with the
 property that when the energy of the $2$-sphere is close to the
 maximum value along the path, then the $2$-sphere in question
 represents a point in the space ${\rm Maps}(S^2,M)$ that is
 close to a (usually) non-minimal critical point for the energy
 functional on this space. Such paths are needed in order to
 establish the forward difference quotient result alluded to above.
In Perelman's approach, one deals only with area-minimizing disks so
that one avoids having to deal with non-minimal critical points at
the expense of dealing with the technical issues related to the
boundary. Since the latter are one-dimensional in nature, they are
much easier to handle. In the end we decided to follow Perelman's
approach, and that is the one we present here.  In \cite{Perelman3}
there were two points that we felt required quite a bit of argument
beyond what Perelman presented. In \S 2.2 on page 4 of
\cite{Perelman3}, Perelman asserts that there is a local, pointwise
curvature estimate that can be obtained by adapting arguments in the
literature; see Lemmas~\ref{CSSHI} and~\ref{csshi1} for the precise
statement. To implement this adaption required further non-trivial
arguments. We present these arguments in Section~\ref{18.6}. In \S
2.5 on page 5 of \cite{Perelman3} Perelman asserts that an
elementary argument establishes a lower bound on the length of a
boundary curve of a minimal annulus; see
Proposition~\ref{annlengthest} for a precise statement. While the
statement seems intuitively clear, we found the argument, while
elementary, was quite intricate. We present this argument in
Section~\ref{18.5}.

The first use of these types of ideas to show that geometric objects
must disappear in finite time under Ricci flow is due to Hamilton
\cite{HamiltonNSRF3M}. He was considering a situation where a
time-slice $(M,g(t_0))$ of a $3$-dimensional Ricci flow had
submanifolds on which the metric was close to (a truncated version)
of a hyperbolic metric of finite volume. He wished to show that
eventually the boundary tori of the truncation were incompressible
in the $3$-manifold. If not, then there would be an immersed minimal
disk in $M$ whose boundary was a non-trivial loop on the torus. He
represented this relative homotopy class by a minimal energy disk in
$(M,g(t_0))$ and proved the same sort of forward difference quotient
estimate for the area of the minimal disk in the relative homotopy
class.  The same contradiction -- the forward difference quotient
implies that after a finite time the area would go negative if the
disk continued to exist --- implies that after a finite amount of
time this compressing disk must disappear. Using this he showed that
for sufficiently large time all the boundary tori of almost
hyperbolic submanifolds in $(M,g(t))$ were incompressible.

In the next section we deal with $\pi_2$ and, using $W_2$, we show
that given a Ricci flow with surgery as in Theorem~\ref{MAIN} there
is $T_1<\infty$ such that for all $T\ge T_1$ every connected
component of $M_T$ has trivial $\pi_2$. Then in the section after
that, by analyzing $W_3$, we show that, under the group-theoretic
hypothesis of Theorem~\ref{extinct}, there is a $T_2<\infty$ such
that $M_T=\emptyset$ for all $T\ge T_2$. In both these arguments we
need the same type of results -- a forward difference inequality for
the energy function;  the statement that away from surgery times
this function is continuous; and lastly, the statement that the
value of the energy function at a surgery time is at most the liminf
of its values at a sequence of times approaching the surgery time
from below.

\subsection{Existence of the Ricci flow with surgery}

Let $(M,g(0))$ be as in the statement of Theorem~\ref{extinct}, so that $M$ is
a compact, connected $3$-manifold whose fundamental group is a free product of
finite groups and infinite cyclic groups. By scaling $g(0)$ by a sufficiently
large constant, we can assume that $g(0)$ is normalized. Let us show that such
a manifold cannot contain an embedded $\RR P^2$ with trivial normal bundle.
First note that since $\RR P^2$ has Euler characteristic one, it is not the
boundary of a compact $3$-manifold. Hence, an $\RR P^2$ embedded with trivial
normal bundle does not separate the connected component of $M$ containing it.
Also, any non-trivial loop in $\Ar P^2$ has non-trivial normal bundle in $M$ so
that inclusion of $\RR P^2$ into $M$ induces an injection on fundamental
groups. Under the fundamental group hypotheses, $M$ decomposes as a connected
sum of $3$-manifolds with finite fundamental groups and $2$-sphere bundles over
$S^1$, see \cite{Hempel}.
 Given
an $\RR P^2$ with trivial normal bundle embedded in a connected sum,
it can be replaced by one contained in one of the connected factors.
[Proof: Let $\Sigma=\Sigma_1\cup\cdots\cup\Sigma_n$ be the spheres
giving the connected sum decomposition of $M$.  Deform the $\Ar P^2$
until it is transverse to $\Sigma$ and let $\gamma$ be a circle of
intersection of $\Ar P^2$ with one of the $\Sigma_i$ that is
innermost on $\Sigma_i$ in the sense that $\gamma$ bounds a disk $D$
in $\Sigma_i$ disjoint from all other components of intersection of
$\Sigma_i$ and $\Ar P^2$. The loop $\gamma$ also bounds a disk $D'$
in $\Ar P^2$. Replace $D'$ by $D$ and push $D$ slightly off to the
correct side of $\Sigma_i$. This will produce a new embedded $\Ar
P^2$ with trivial normal bundle in $M$ and at least one fewer
component of intersection with $\Sigma$. Continue inductively until
all components of intersection with $\Sigma$ are removed.]

Now suppose that we have an $\Ar P^2$ with trivial normal bundle
embedded disjointly from $\Sigma$, and hence embedded in one of the
prime factors of $M$.
 Since it
does not separate this factor, by the Mayer-Vietoris sequence (see p. 149 of
\cite{Hatcher}) the first homology of the factor in question maps onto $\Zee$
and hence the factor in question has infinite fundamental group. But this group
also contains the cyclic subgroup of order two, namely the image of $\pi_1(\RR
P^2)$ under the map induced by the inclusion. Thus, the fundamental group of
this prime factor is not finite and is not infinite cyclic. This is a
contradiction. (We have chosen to give a topological argument for this result.
There is also an argument using the theory of groups acting on trees which is
more elementary in the sense that it uses no $3$-manifold topology. Since it is
a more complicated, and to us, a less illuminating argument, we decided to
present the topological argument.)

Thus, by Theorem~\ref{MAIN},  for any compact $3$-manifold $M$ whose
fundamental group is a free product of finite groups and infinite cyclic groups
and for any normalized metric $g(0)$ on $M$ there is a Ricci flow with surgery
$({\mathcal M},G)$ defined for all time $t\in [0,\infty)$ satisfying the
conclusion of Theorem~\ref{MAIN} with $(M,g(0))$ as the initial conditions.

\begin{defn}
Let $I$ be an interval (which is allowed to be open or closed at
each end and finite or infinite at each end). By a {\em path of
components} of a Ricci flow with surgery $({\mathcal M},G)$ defined
for all $t\in I$ we mean a connected, open subset ${\mathcal
X}\subset {\bf t}^{-1}(I)$ with the property that for every $t\in I$
the intersection ${\mathcal X}(t)$ of ${\mathcal X}$ with each
time-slice $M_t$ is a connected component of $M_t$.
\end{defn}

Let ${\mathcal X}$ be a path of components in a Ricci flow with surgery
$({\mathcal M},G)$, a path defined for all $t\in I$. Let $I'$ be a subinterval
of $I$ with the property that no point of $I'$ except possibly its initial
point is a surgery time.  Then the intersection of ${\mathcal X}$ with ${\bf
t}^{-1}(I')$ is  the Ricci flow on the time interval $I'$ applied to ${\mathcal
X}(t)$ for any $t\in I'$.  Thus, for such intervals $I'$ the intersection,
${\mathcal X}(I')$, of ${\mathcal X}$ with ${\bf t}^{-1}(I')$ is determined by
the time-slice ${\mathcal X}(t)$ for any $t\in I'$. That is no longer
necessarily the case if some point of $I'$ besides its initial point is a
surgery time. Let $t\in I$ be a surgery time, distinct from the initial point
of $I$ (if there is one), and let $I'\subset I$ be an interval of the form
$[t',t)$ for some $t'<t$ sufficiently close to $t$ so that there are no surgery
times in $[t',t)$. Then, as we have just seen, ${\mathcal X}(I')$ is a Ricci
flow on the connected manifold ${\mathcal X}(t')$. There are several possible
outcomes of the result of surgery at time $t$ on this manifold. One possibility
is that the surgery leaves this connected component unchanged (affecting only
other connected components). In this case, there is no choice for ${\mathcal
X}(t)$: it is the continuation to time $t$ of the Ricci flow on ${\mathcal
X}(t')$. Another possibility is that ${\mathcal X}(t')$ is completely removed
by the surgery at time $t$. In this case the manifold ${\mathcal X}$ cannot be
continued to time $t$, contradicting the fact that the path of components
${\mathcal X}$ exists for all $t\in I$. The last possibility is that at time
$t$ surgery is done on ${\mathcal X}(t')$ using one or more $2$-spheres
contained in ${\mathcal X}(t')$. In this case the result of surgery on
${\mathcal X}(t')$ results in one or several connected components and
${\mathcal X}(t)$ can be any one of these.

\section{Disappearance of components with non-trivial $\pi_2$}

Let $({\mathcal M},G)$ be a Ricci flow with surgery satisfying the
conclusions of Theorem~\ref{MAIN}. We make no assumptions about the
fundamental group of the initial manifold $M_0$. In this section we
shall show that at some finite time $T_1$ every connected component
of $M_{T_1}$ has trivial $\pi_2$ and that this condition persists
for all times $T\ge T_1$. There are two steps in this argument.
First, we show that there is a finite time $T_0$ such that after
time $T_0$ every $2$-sphere surgery is performed along a
homotopically trivial $2$-sphere. (Using Kneser's theorem on
finiteness of topologically non-trivial families of $2$-spheres, one
can actually show by the same argument that after some finite time
all $2$-sphere surgeries are done along $2$-spheres that bound
$3$-balls. But in fact, Kneser's theorem will follow from what we do
here.)

After time $T_0$ the number of components with non-trivial $\pi_2$
is a weakly monotone decreasing function of time. The reason is the
following. Consider a path of components ${\mathcal X}$ defined for
$t\in [T_0,t']$ with the property that each time-slice ${\mathcal
X}(t)$ has non-trivial $\pi_2$. Using the fact that after time $T_0$
all the $2$-sphere surgeries are along homotopically trivial
$2$-spheres,  one shows easily that ${\mathcal X}$ is determined by
its initial time-slice ${\mathcal X}(T_0)$. Also, it is easy to see
that if there is a component of $M_{t'}$ with non-trivial $\pi_2$,
then it is the final time-slice of some path of components defined
for $t\in [T_0,t']$ with every time-slice of this path having
non-trivial $\pi_2$. This then produces an injection from the set of
connected components of $M_{t}$ with non-trivial $\pi_2$ into the
set of connected components of $M_{T_0}$ with non-trivial $\pi_2$.

The second step in the argument is to fix a path ${\mathcal X}(t)$,
$T_0\le t\le t'$, of connected components with non-trivial $\pi_2$
and to consider the function $W_2=W_2^{\mathcal X}$ that assigns to
each $t\in [T_0,t']$ the minimal area of a homotopically non-trivial
$2$-sphere mapping into ${\mathcal X}(t)$. We show that this
function is continuous except at the surgery times. Furthermore, we
show that if $t$ is a surgery time, then $W_2(t)\le {\rm
liminf}_{t'\rightarrow t^-}W_2(t)$. Lastly, we show that at any
point $t\ge T_0$ we have
$$\frac{dW_2}{dt}(t)\le -4\pi-\frac{1}{2}R_{\rm min}(t)W_2(t),$$
in the sense of forward difference quotients. It follows easily from
the bound $R_{\rm min}(t)\ge -6/(4t+1)$ that there is $T_1({\mathcal
X})$  such that $W_2$ with these three properties cannot be
non-negative for all $t\in [T_0,T_1({\mathcal X})]$ and hence
$t'<T_1$. Since there are only finitely many components with
non-trivial $\pi_2$ at time $T_0$ it follows that there is
$T_1<\infty$ such that every component of $M_T$ has trivial $\pi_2$
for every $T\ge T_1$.

\subsection{A group-theory lemma}

To bound the number of homotopically non-trivial $2$-spheres in a
compact $3$-manifold we need the following group theory lemma.

\begin{lem}\label{grouptheory} Suppose that $G$ is a finitely generated group,
say generated by $k$ elements. Let $G=G_1*\cdots*G_\ell$ be a free
product decomposition of $G$ with non-trivial free factors, i.e.,
with $G_i\not=\{1\}$ for each $i=1,\ldots,\ell$. Then $\ell\le k$.
\end{lem}

\begin{proof}
This is a consequence of Grushko's theorem \cite{Stallings}, which
says that given a map of a finitely generated free group $F$ onto
the  free product $G$, one can decompose the free group as a free
product of free groups $F=F_1*\cdots*F_\ell$ with $F_i$ mapping onto
$G_i$.
\end{proof}

\subsection{Homotopically non-trivial families of $2$-spheres}

\begin{defn}
Let $X$ be a compact $3$-manifold (possibly disconnected). An
embedded $2$-sphere in $X$ is said to be {\em homotopically
essential} if the inclusion of the $2$-sphere into $X$ is not
homotopic to a point map of the $2$-sphere to $X$. More generally,
let $F=\{\Sigma_1,\ldots,\Sigma_n\}$ be a family of disjointly
embedded $2$-spheres in $X$. We say that the family is {\em
homotopically essential} if
\begin{enumerate}
\item[(i)] each $2$-sphere in the family is homotopically essential, and
\item[(ii)] for any $1\le i<j\le n$, the inclusion of $\Sigma_i$ into $X$ is not
homotopic in $X$ to the inclusion of $\Sigma_j$ into $X$.
\end{enumerate}
\end{defn}

Notice that if $F=\{\Sigma_1,\ldots,\Sigma_n\}$ is a homotopically
essential family of disjointly  embedded $2$-spheres in $X$, then
any subset $F$ is also homotopically essential.

\begin{lem}
Let $X$ be a compact $3$-manifold (possibly disconnected). Then
there is a finite upper bound to the number of spheres in any
homotopically essential family of disjointly embedded $2$-spheres.
\end{lem}

\begin{proof}
Clearly, without loss of generality we can assume that $X$ is
connected. If $F$ is a homotopically essential family of $2$-spheres
in $X$, then by van Kampen's theorem, see p. 40 of \cite{Hatcher},
there is an induced graph of groups decomposition of $\pi_1(X)$ with
all the edge groups being trivial. Since the family is homotopically
essential, it follows that the group associated with each vertex of
order $1$ and each vertex of order $2$ is non-trivial group. The
rank of the first homology of the graph underlying the graph of
groups, denoted $k$, is bounded above by the rank of $H_1(X)$.
Furthermore, by the theory of graphs of groups there is a free
product decomposition of $\pi_1(X)$ with the free factors being the
vertex groups and then $k$ infinite cyclic factors. Denote by $V_i$
the number of vertices of order $i$ and by $E$ the number of edges
of the graph. The number $E$ is the number of $2$-spheres in the
family $F$.
 An elementary combinatorial argument shows that
$$2V_1+V_2\ge E+3(1-k).$$
Thus, we have a free product decomposition of $\pi_1(X)$ with at least
$E+3(1-k)$ non-trivial free factors. Since $k$ is bounded by the rank of
$H_1(X)$,
 applying
Lemma~\ref{grouptheory} and using the fact that the fundamental group of a
compact manifold is finitely presented establishes the result.
\end{proof}

\subsection{Two-sphere surgeries are trivial after finite
time}

\begin{defn}
Let $({\mathcal M},G)$ be a Ricci flow with surgery. We say that a
surgery along a $2$-sphere $S_0(t)$ at time $t$ in $({\mathcal
M},G)$ is a {\em homotopically essential surgery} if, for every
$t'<t$ sufficiently close to $t$,  flowing  $S_0(t)$ backwards from
time $t$ to time $t'$ results in a homotopically essential
$2$-sphere $S_0(t')$ in $M_{t'}$.
\end{defn}

\begin{prop}\label{homess}
Let $({\mathcal M},G)$ be a Ricci flow with surgery satisfying
Assumptions (1) -- (7) in Chapter~\ref{sect:surgery}. Then there can
be only finitely many homotopically essential surgeries along
$2$-spheres in $({\mathcal M},G)$.
\end{prop}

\begin{proof}
 Associate to each compact $3$-manifold $X$
 the invariant $s(X)$ which is the maximal number of spheres in
 any homotopically essential family of embedded $2$-spheres in $X$.
 The main step in establishing the corollary is the following:

\begin{claim}
Let $({\mathcal M},G)$ be a Ricci flow with surgery and for each $t$ set
$s(t)=s(M_t)$. If $t'<t$ then $s(t')\ge s(t)$. If  we do surgery at time $t$
along at least one homotopically essential $2$-sphere, then $s(t)<s(t')$ for
any $t'<t$.
\end{claim}

\begin{proof}
Clearly, for any $t_0$ we have $s(t)=s(t_0)$ for $t\ge t_0$
sufficiently close to $t_0$. Also, if $t$ is not a surgery time,
then $s(t)=s(t')$ for all $t'<t$ and sufficiently close to $t$.
According to Proposition~\ref{surgerytoptype}, if $t$ is a surgery
time then for $t'<t$ but sufficiently close to it, the manifold
$M_t$ is obtained from $M_{t'}$  by doing surgery on a finite number
of $2$-spheres and removing certain components of the result. We
divide the operations into three types: (i) surgery along
homotopically trivial $2$-spheres in $M_{t'}$, (ii) surgery along
homotopically non-trivial $2$-spheres in $M_{t'}$, (iii) removal of
components. Clearly, the first operation does not change the
invariant $s$ since it simply creates a manifold that is the
disjoint union of a manifold homotopy equivalent to the original
manifold with a collection of homotopy $3$-spheres. Removal of
components will not increase the invariant. The last operation to
consider is surgery along a homotopically non-trivial $2$-sphere.
Let $F_t$ be a homotopically essential family of disjointly embedded
$2$-spheres in $M_t$. This family of $2$-spheres in $M_t$ can be
deformed to miss the $3$-disks (the surgery caps) in $M_t$ that we
sewed in doing the surgery at time $t$ along a homotopically
non-trivial $2$-sphere. After deforming the spheres in the family
$F_t$ away from the surgery caps, they produce a disjoint family
$F'_{t'}$ of $2$-spheres in the manifold $M_{t'}$, for $t'<t$ but
$t'$ sufficiently close to $t$. Each $2$-sphere in $F'_{t'}$ is
disjoint from the homotopically essential $2$-sphere $S_0$ along
which we do surgery at time $t$. Let $F_{t'}$ be the family
$F'_{t'}\cup \{S_0\}$. We claim that $F_{t'}$ is a homotopically
essential family in $M_{t'}$.

First, suppose that one of the spheres $\Sigma$ in $F_{t'}$ is
homotopically trivial in $M_{t'}$. Of course, we are in the case
when the surgery $2$-sphere is homotopically essential, so $\Sigma$
is not $S_0$ and hence is the image of one of the $2$-spheres in
$F_t$. Since $\Sigma$ is homotopically trivial, it is the boundary
of a homotopy $3$-ball $B$ in $M_{t'}$. If $B$ is disjoint from the
surgery $2$-sphere $S_0$, then it exists in $M_t$ and hence $\Sigma$
is homotopically trivial in $M_t$, which is not possible from the
assumption about the family $F_t$. If $B$ meets the surgery
$2$-sphere $S_0$, then since the spheres in the family $F_{t'}$ are
disjoint, it follows that $B$ contains the surgery $2$-sphere $S_0$.
This is not possible since in this case $S_0$ would be homotopically
trivial in $M_{t'}$, contrary to assumption.

We also claim that no distinct members of $F_{t'}$ are homotopic.
For suppose that two of the members $\Sigma$ and $\Sigma'$ are
homotopic. It cannot be the case that one of $\Sigma$ or $\Sigma'$
is the surgery $2$-sphere $S_0$ since, in that case, the other one
would be homotopically trivial after surgery, i.e., in $M_t$. The
$2$-spheres $\Sigma$ and $\Sigma'$ are the boundary components of a
submanifold  $A$ in $M_{t'}$ homotopy equivalent to $S^2\times I$.
If $A$ is disjoint from the surgery $2$-sphere $S_0$, then $A$
exists in $M_t$ and $\Sigma$ and $\Sigma'$ are homotopic in $M_{t}$,
contrary to assumption. Otherwise, the surgery sphere $S_0$ must be
contained in $A$. Every $2$-sphere in $A$ is either homotopically
trivial in $A$ or is homotopic in $A$ to either boundary component.
If $S_0$ is homotopically trivial in $A$, then it would be
homotopically trivial in $M_{t'}$ and this contradicts our
assumption. If $S_0$ is homotopic in $A$ to each of $\Sigma$ and
$\Sigma'$, then each of $\Sigma$ and $\Sigma'$ is homotopically
trivial in $M_{t'}$, contrary to assumption. This shows that the
family $F_{t'}$ is homotopically essential. It follows immediately
that doing surgery on a homotopically non-trivial $2$-sphere
strictly decreases the invariant $s$.
\end{proof}

Proposition~\ref{homess} is immediate from this claim and the
previous lemma.
\end{proof}

\subsection{For all $T$ sufficiently large $\pi_2(M_T)=0$}

We have just established that given any Ricci flow with surgery $({\mathcal
M},G)$ satisfying the conclusion of Theorem~\ref{MAIN} there is $T_0<\infty$,
depending  on $({\mathcal M},G)$, such that all surgeries after time $T_0$
either are along homotopically trivial $2$-spheres or remove entire components
of the manifold. Suppose that $M_{T_0}$ has a component ${\mathcal X}(T_0)$
with non-trivial $\pi_2$, and suppose that we have a path of components
${\mathcal X}(t)$ defined for $t\in [T_0,T)$ with the property that each
time-slice has non-trivial $\pi_2$. If $T$ is not a surgery time, then there is
a unique extension of ${\mathcal X}$ to a path of components with non-trivial
$\pi_2$ defined until the first surgery time after $T$. Suppose that $T$ is a
surgery time and let us consider the effect of surgery at time $T$ on
${\mathcal X}(t)$ for $t<T$ but close to it. Since no surgery after time $T_0$
is done on a homotopically essential $2$-sphere there are three possibilities:
(i) ${\mathcal X}(t)$ is untouched by the surgery, (ii) surgery is performed on
one or more homotopically trivial $2$-spheres in ${\mathcal X}(t)$, or (iii)
the component ${\mathcal X}(t)$ is completely removed by the surgery. In the
second case, the result of the surgery on ${\mathcal X}(t)$ is a disjoint union
of components one of which is homotopy equivalent to ${\mathcal X}(t)$, and
hence has non-trivial $\pi_2$, and all others are homotopy $3$-spheres. This
implies that there is a unique extension of the path of components preserving
the condition that every time-slice has non-trivial $\pi_2$, unless the
component ${\mathcal X}(t)$ is removed by surgery at time $T$, in which case
there is no extension of the path of components to time $T$. Thus, there is a
unique maximal such path of components starting at ${\mathcal X}(T_0)$ with the
property that every time-slice has non-trivial $\pi_2$. There are two
possibilities for the interval of definition of this maximal path of components
with non-trivial $\pi_2$. It can be $[T_0,\infty)$ or it is of the form
$[T_0,T)$, where the surgery at time $T$ removes the component ${\mathcal
X}(t)$ for $t<T$ sufficiently close to it.

\begin{prop}\label{pi2}
Let $({\mathcal M},G)$ be a Ricci flow with surgery satisfying the conclusion
of Theorem~\ref{MAIN}. Then there is some time $T_1<\infty$ such that every
component of $M_T$ for any $T\ge T_1$ has trivial $\pi_2$. For every $T\ge
T_1$, each component of $M_T$ either has finite fundamental group, and hence
has a homotopy $3$-sphere as universal covering, or has  contractible universal
covering.
\end{prop}

If $M$ is a connected $3$-manifold with $\pi_2(M)=0$, then the
universal covering, $\widetilde M$, of $M$ is a $2$-connected
$3$-manifold. The covering $\widetilde M$ is compact if and only if
$\pi_1(M)$ is finite. In this case $\widetilde M$ is a homotopy
$3$-sphere. If $\widetilde M$ is non-compact then $H_3(\widetilde
M)=0$, so that all its homology groups and hence, by the Hurewicz
theorem, all its homotopy groups vanish. It follows from the
Whitehead theorem  that $\widetilde M$ is contractible in this case.
This proves the last assertion in the proposition modulo the first
assertion.

The proof of the first assertion of this proposition occupies the
rest of this subsection. By the above discussion we see that the
proposition holds unless there is a path of components ${\mathcal
X}$ defined for all $t\in [T_0,\infty)$ with the property that every
time-slice has non-trivial $\pi_2$. We must rule out this
possibility. To achieve this we introduce the area functional.

\begin{lem}
Let $X$ be a compact Riemannian manifold with $\pi_2(X)\not=0$. Then
there is a positive number $e_0=e_0(X)$ with the following two
properties:
\begin{enumerate} \item Any map $f\colon S^2\to X$ with area less
than $e_0$ is homotopic to a point map. \item There is a minimal
$2$-sphere $f\colon S^2\to X$, which is a branched immersion, with
the property that the area of $f(S^2)=e_0$ and with the property
that $f$ is not homotopic to a point map.
\end{enumerate}
\end{lem}

\begin{proof}
The first statement is Theorem 3.3 in \cite{SacksUhlenbeck}. As for
the second, following Sacks-Uhlenbeck, for any $\alpha>1$ we
consider the perturbed energy $E_\alpha$ given by
$$E_\alpha(s)=\int_{S^2}\left(1+|ds|^2\right)^\alpha da.$$
 According to
\cite{SacksUhlenbeck} this energy functional is Palais-Smale on the
space of $H^{1,2\alpha}$ maps and has an absolute minimum among
homotopically non-trivial maps, realized by a map $s_\alpha\colon
S^2\to X$. We consider a decreasing sequence of $\alpha$ tending to
$1$ and the minimizers $s_\alpha$ among homotopically non-trivial
maps. According to \cite{SacksUhlenbeck}, after passing to a
subsequence, there is a weak limit which is a strong limit on the
complement of a finite set of points in $S^2$. This limit extends to
a harmonic map of $S^2\to M$, and its energy is less than or equal
to the limit of the $\alpha$-energies of $s_\alpha$. If the result
is homotopically non-trivial then it realizes a minimum value of the
usual energy among all homotopically non-trivial maps, for were
there a homotopically non-trivial map of smaller energy, it would
have smaller $E_\alpha$ energy than $s_\alpha$ for all $\alpha$
sufficiently close to $1$. Of course if the limit is a strong limit,
then the map is homotopically non-trivial, and the proof is
complete.

We must examine the case when the limit is truly a weak limit. Let
$s_n$ be a sequence as above with a weak limit $s$. If the limit is
truly a weak limit, then there is bubbling. Let $x\in S^2$ be a
point where the limit $s$  is not a strong limit. Then according to
\cite{SacksUhlenbeck} pre-composing with a sequence of conformal
dilations $\rho_n$ centered at this point leads to a sequence of
maps $s_n'$ converging uniformly on compact subsets of $\Ar ^2$ to a
non-constant harmonic map $s'$ that extends over the one-point
compactification $S^2$. The energy of this limiting map $s'$ is at
most the limit of the $\alpha$-energies of the $s_\alpha$. If $s'$
is homotopically non-trivial, then, arguing as before, we see that
it realizes the minimum energy among all homotopically non-trivial
maps, and once again we have completed the proof. We rule out the
possibility that $s'$ is homotopically trivial. Let $\alpha$ be the
area, or equivalently the energy, of $s'$. Let $D\subset \Ar^2$ be a
disk centered at the origin which contains three-quarters of the
energy of $s'$ (or equivalently three-quarters of the area of $s'$),
and let $D'$ be the complementary disk to $D$ in $S^2$.  For all $n$
sufficiently large the area of $s_n'|D$ minus the area of $s_n'|D'$
is at least $\alpha/3$. The restrictions of $s_n'$ on $\partial D$
are converging smoothly to $s'|\partial D'$. Let $D_n\subset S^2$ be
$\rho_n^{-1}(D)$. Then the area of $s_n|D_n$ equals the area of
$s_n'|D$ and hence is at least the area of $s'|D'$ plus $\alpha/4$
for all $n$ sufficiently large. Also, as $n$ tends to infinity the
image $s_n(D_n)$ converges smoothly, after reparameterization, to
$s'(\partial D)$. Thus, for all $n$ large, we can connect
$s_n(\partial D_n)$ to $s'(\partial D')$ by an annulus $A_n$
contained in a small neighborhood of $s'(\partial D')$ and whose
area tends to $0$ as $n$ goes to infinity.  For all $n$ sufficiently
large, the resulting $2$-sphere $\Sigma_n$ made out of
$s_n|(S^2\setminus D_n)\cup A_n\cup S'(D')$ is homotopic to $s(S^2)$
since $s'$ is homotopically trivial. Also, for all $n$ sufficiently
large, the area of $\Sigma_n$ is less than the area of $s_n$ minus
$\alpha/5$. Reparameterizing this $2$-sphere by a conformal map
leads to a homotopically non-trivial map of energy less than the
area of $s_n$ minus $\alpha/5$. Since as $n$ tends to infinity, the
limsup of the areas of the $s_n$ converge to at most $e_0$, for all
$n$ sufficiently large we have constructed a homotopically
non-trivial map of energy less than $e_0$, which contradicts the
fact that the minimal $\alpha$ energy for a homotopically
non-trivial map tends to $e_0$ as $\alpha$ tends to $1$.

Of course, any minimal energy map of $S^2$ into $M$ is conformal because there
is no non-trivial holomorphic quadratic differential on $S^2$. It follows that
such a map is a branched immersion.
\end{proof}

Now suppose that ${\mathcal X}$ is a path of components defined for
all $t\in [T_0,\infty)$ with $\pi_2({\mathcal X}(t))\not= 0$ for all
$t\in [T_0,\infty)$.
 For
each $t\ge T_0$ we define $W_2(t)$ to be $e_0({\mathcal X}(t))$, where $e_0$ is
the invariant  given in the previous lemma. Our assumption on ${\mathcal X}$
means that $W_2(t)$ is defined and positive for all $t\in [T_0,\infty)$.

\begin{lem}\label{W_2}
$$\frac{d}{dt}W_2(t)\le -4\pi-\frac{1}{2}R_{\rm min}(t)W_2(t)$$
in the sense of forward difference quotients. If $t$ is not a
surgery time, then $W_2(t)$ is continuous at $t$, and if $t$ is a
surgery time, then $$W_2(t)\le {\rm liminf}_{t'\rightarrow
t^-}W_2(t').$$
\end{lem}

Let us show how this lemma implies Proposition~\ref{pi2}. Because
the curvature is pinched toward positive, we have
$$R_{\rm min}(t)\ge (-6)/(1+4t).$$
Let $w_2(t)$ be the function satisfying the differential equation
$$\frac{dw_2}{dt}=-4\pi+\frac{3w_2}{1+4t}$$
and $w_2(T_0)=W_2(T_0)$. Then by Lemma~\ref{fordiffquot} and
Lemma~\ref{W_2} we have $W_2(t)\le w_2(t)$ for all $t\ge T_0$. On
the other hand, we can integrate to find
$$w_2(t)=w_2(T_0)\frac{(4t+1)^{3/4}}{(4T_0+1)^{3/4}}+4\pi
(4T_0+1)^{1/4}(4t+1)^{3/4}-4\pi (4t+1).$$ Thus, for $t$ sufficiently
large, $w_2(t)<0$. This is a contradiction since $W_2(t)$ is always
positive, and $W_2(t)\le w_2(t)$.

This shows that to complete the proof of Proposition~\ref{pi2} we
need only establish Lemma~\ref{W_2}.

\begin{proof}(of Lemma~\ref{W_2})
Let $f\colon S^2\to (X(t_0),g(t_0))$ be a minimal $2$-sphere.

\begin{claim}\label{varareaform}
$$\frac{d{\rm Area}_{g(t)}(f(S^2))}{dt}(t_0)\le -4\pi -\frac{1}{2}
R_{\rm min}(g(t_0)){\rm Area}_{g(t_0)}f(S^2).$$
\end{claim}

\begin{proof}
 Recall that, for any immersed surface $f\colon S^2\to (M,g(t_0))$, we
have ([Ha])
 \begin{eqnarray}\label{ha1}
 \f1{d}{dt}{\rm Area}_{g(t)}(f(S^2))\bigl|_{t=t_0}\bigr. & = &
 \int_{S^2}\frac{1}{2}{\rm Tr}|_{S^2}\Bigl(\frac{\partial g}{\partial t}\Bigr)\Bigl|_{t=t_0}\Bigr.da \\
 & = &
 -\int_{S^2}(R-{\rm Ric}({\bf n},{\bf n}))da \nonumber
 \end{eqnarray}
 where $R$ denotes the scalar curvature of $M$, ${\rm Ric}$ is the
Ricci curvature of $M$, and ${\bf n}$ is the unit normal vector field of $\g1$
in $M$. Now suppose that $f(S^2)$ is minimal. We can rewrite this as
 \begin{eqnarray}
 \f1{d}{dt}{\rm
 Area}_{g(t)}(f(S^2))\big |_{t=0}
 &= &
 -\int_{S^2}K_{S^2}da-\f1{1}{2}\int_{S^2}(|A|^2+R)da,  \label{ha2}
 \end{eqnarray}
where $K_{S^2}$ is the Gaussian curvature of $S^2$ and $A$ is the
second fundamental form of $f(S^2)$ in $M$. (Of course, since
$f(S^2)$ is minimal, the determinant of its second fundamental form
is $-|A|^2/2$.) Even if $f$ is only a branched minimal surface,
(\ref{ha2}) still holds when  the integral on the right is replaced
by the integral over the immersed part of $f(S^2)$. Then by the
Gauss-Bonnet theorem we have
 \begin{eqnarray}
\f1{d}{dt}{\rm Area}_{g(t)}(f(S^2))\bigl|_{t=t_0}\bigr.  & \leq &
-4\pi-\f1{1}{2}{\rm
 Area}_{g(t_0)}(S^2)\min_{x\in M}\{R_{g}(x,t_0)\}.  \label{ha3}
 \end{eqnarray}
\end{proof}

Since $f(S^2)$ is a homotopically non-trivial sphere in ${\mathcal X}(t)$ for
all $t$ sufficiently close to $t_0$ we see that $W_2(t)\le {\rm
Area}_{g(t)}f(S^2)$. Since ${\rm Area}_{g(t)}f(S^2))$ is a smooth function of
$t$, the forward difference quotient statement in Lemma~\ref{W_2} follows
immediately from Claim~\ref{varareaform}.

We turn now to continuity at non-surgery times. Fix $t'\ge T_0$
distinct from all surgery times. We show that  the function
$e_0(t')$ is continuous at $t'$. If $f\colon S^2\to {\mathcal
X}(t')$ is the minimal area, homotopically non-trivial sphere, then
the area of $f(S^2)$ with respect to a nearby metric $g(t)$ is close
to the area of $f(S^2)$ in the metric $g(t')$. Of course, the area
of $f(S^2)$ in the metric $g(t)$ is greater than or equal to
$W_2(t)$. This proves that $W_2(t)$ is  upper semi-continuous at
$t'$. Let us show that it is lower semi-continuous at  $t'$.

\begin{claim}
Let $(M,g(t)),\ t_0\le t\le t_1$, be a Ricci flow on a compact
manifold. Suppose that $|{\rm Ric}_{g(t)}|\le D$ for all $t\in
[t_0,t_1]$ Let $f\colon S^2\to (M,g(t_0))$ be a $C^1$-map. Then
$${\rm Area}_{g(t_1)}f(S^2)\le {\rm
Area}_{g(t_0)}f(S^2)e^{4D(t_1-t_0)}.$$
\end{claim}

\begin{proof}
The rate of change of the area of $f(S^2)$ at time $t$ is
$$\int_{f(S^2)}\frac{\partial g}{\partial t}(t)da=-2\int_{f(S^2)}{\rm Tr}|_{TS^2}({\rm Ric}_{g(t)})da\le
4D{\rm Area}_{g(t)}f(S^2) .$$ Integrating from $t_0$ to $t_1$ gives
the result.
\end{proof}

Now suppose that we have a family of times $t_n$ converging to a
time $t'$ that is not a surgery time. Let $f_n\colon S^2\to
{\mathcal X}(t_n)$ be the minimal area non-homotopically trivial
$2$-sphere in ${\mathcal X}(t_n)$, so that the area of $f_n(S^2)$ in
${\mathcal X}(t_n)$ is $e_0(t_n)$. Since $t'$ is not a surgery time,
for all $n$ sufficiently large we can view the maps $f_n$ as
homotopically non-trivial maps of $S^2$ into ${\mathcal X}(t')$. By
the above claim, for any $\delta>0$ for all $n$ sufficiently large,
the area of $f_n(S^2)$ with respect to the metric $g(t')$ is at most
the area of $f_n(S^2)$ plus $\delta$. This shows that for any
$\delta>0$ we have $W_2(t')\le W_2(t_n)+\delta$ for all $n$
sufficiently large, and hence $W_2(t')\le {\rm
liminf}_{n\rightarrow\infty}W_2(t_n)$. This is the lower
semi-continuity.

The last thing to check is the behavior of $W_2$ near a surgery time
$t$. According to the description of the surgery process given in
Section~\ref{sect:process}, we write ${\mathcal X}(t)$ as the union
of a compact subset $C(t)$ and a finite number of surgery caps. For
every $t'<t$ sufficiently close to $t$ we have an embedding
$n_{t'}\colon C(t)\cong C(t')\subset {\mathcal X}(t')$ given by
flowing $C(t)$ backward under the flow to time $t'$. As
$t'\rightarrow t$ the maps $\eta_{t'}$ converge in the
$C^\infty$-topology to isometries, in the sense that the
$n_{t'}^*(g(t'))|_{C(t')}$ converge smoothly to $g(t)|_{C(t)}$.
Furthermore, since the $2$-spheres along which we do surgery are
homotopically trivial they separate $M_{t'}$. Thus, the maps
$n_{t'}^{-1}\colon C(t')\to C(t)$ extend to maps $\psi_{t'}\colon
{\mathcal X}(t')\to {\mathcal X}(t)$. The image under $\psi_{t'}$ of
${\mathcal X}(t')\setminus C(t')$ is contained in the union of the
surgery caps.  Clearly, since all the $2$-spheres on which we do
surgery at time $t$ are homotopically trivial, the maps $\psi_{t'}$
are homotopy equivalences. If follows from
Proposition~\ref{surgerydist} that for any $\eta>0$ for all $t'<t$
sufficiently close to $t$, the map $\psi_{t'}\colon {\mathcal
X}(t')\to {\mathcal X}(t)$ is a
 homotopy equivalence that is a $(1+\eta)$-Lipschitz map. Thus,
 given $\eta>0$ for all $t'<t$ sufficiently close to $t$,
 for any minimal $2$-sphere $f\colon S^2\to ({\mathcal X}(t'),g(t'))$ the area of
 $\psi_{t'}\circ f\colon S^2\to ({\mathcal X}(t),g(t))$ is at most $(1+\eta)^2$
 times the area of $f(S^2)$. Thus, given $\eta>0$ for all $t'<t$ sufficiently close to $t$
we see that $W_2(t)\le (1+\eta)^2W_2(t')$.
 Since this is true for every $\eta>0$, it follows that
 $$W_2(t)\le {\rm liminf}_{t'\rightarrow t^-}W_2(t').$$

This establishes all three statements in Proposition~\ref{pi2} and
completes the proof of the proposition.

As an immediate corollary of Proposition~\ref{pi2}, we obtain the
sphere theorem for closed $3$-manifolds.

\begin{cor}\label{spherethm}
Suppose that $M$ is a closed, connected $3$-manifold containing no
embedded $\Ar P^2$ with trivial normal bundle, and suppose that
$\pi_2(M)\not= 0$. Then either $M$ can be written as a connected sum
$M_1\# M_2$ where neither of the $M_i$ is homotopy equivalent to
$S^3$ or $M_1$ has a prime factor that is a $2$-sphere bundle over
$S^1$. In either case, $M$ contains an embedded $2$-sphere which is
homotopically non-trivial.
\end{cor}

\begin{proof}
Let $M$ be as in the statement of the corollary. Let $g$ be a
normalized metric on $M$, and let $({\mathcal M},G)$ be the Ricci
flow with surgery defined for all time with $(M,g)$ as initial
conditions. According to Proposition~\ref{pi2} there is $T<\infty$
such that every component of $M_T$ has trivial $\pi_2$. Thus, by the
analysis above, we see that there must be surgeries that kill
elements in $\pi_2$: either the removal of a component with
non-trivial $\pi_2$ or surgery along a homotopically non-trivial
$2$-sphere. We consider the first such surgery in $M$.
 The only
components with non-trivial $\pi_2$ that can be removed by surgery
are $S^2$-bundles over $S^1$ and $\Ar P^3\#\Ar P^3$. Since each of
these has homotopically non-trivially embedded $2$-spheres, if the
first surgery killing an element in $\pi_2$ is removal of such a
component, then, because all the earlier $2$-sphere surgeries are
along homotopically trivial $2$-spheres, the homotopically
non-trivial embedded $2$-sphere in this component deforms back to an
embedded, homotopically non-trivial $2$-sphere in $M$.  The other
possibility is that the first time an element in $\pi_2(M)$ is
killed it is by surgery along a homotopically non-trivial
$2$-sphere. Once again, using the fact that all previous surgeries
are along homotopically trivial $2$-spheres, deform this $2$-sphere
back to $M$ producing a homotopically non-trivial $2$-sphere in $M$.
\end{proof}

\begin{rem}
Notice that it follows from the list of disappearing components that
the only ones with non-trivial $\pi_2$ are those based on the
geometry $S^2\times \Ar$; that is to say, $2$-sphere bundles over
$S^1$ and $\Ar P^3\#\Ar P^3$. Thus, once we have reached the level
$T_0$ after which all $2$-sphere surgeries are performed on
homotopically trivial $2$-spheres the only components that can have
non-trivial $\pi_2$ are components of these types. Thus, for example
if the original manifold has no $\Ar P^3$ prime factors and no
non-separating $2$-spheres, then when we reach time $T_0$ we have
done a connected sum decomposition into components each of which has
trivial $\pi_2$. Each of these components is either covered by a
contractible $3$-manifold or by a homotopy $3$-sphere, depending on
whether its fundamental group has infinite or finite order.
\end{rem}

\section{Extinction}

Now we assume that the Ricci flow with surgery $({\mathcal M},G)$ satisfies the
conclusion of Theorem~\ref{MAIN} and also has initial condition $M$ that is a
connected $3$-manifold whose fundamental group satisfies the hypothesis of
Theorem~\ref{extinct}. The argument showing that  components with non-trivial
$\pi_3$ disappear after a finite time is, in spirit, very similar to the
arguments above, though the technical details are more intricate in this case.

\subsection{Forward difference quotient for $\pi_3$}

Let $M$ be a compact, connected $3$-manifold. Fix a base point
$x_0\in M$. Denote by $\Lambda M$ the free loop space of $M$. By
this we mean the space of $C^1$-maps of $S^1$ to $M$ with the
$C^1$-topology. The components of $\Lambda M$ are the conjugacy
classes of elements in $\pi_1(M,x_0)$. The connected component of
the identity of $\Lambda M$ consists of all homotopically trivial
loops in $M$. Let $*$ be the trivial loop at $x_0$.

\begin{claim}
Suppose that $\pi_2(M,x_0)=0$. Then $\pi_2(\Lambda M,*)\cong
\pi_3(M,x_0)$ and $\pi_2(\Lambda M,*)$ is identified with the free
homotopy classes of maps of $S^2$ to the component of $\Lambda M$
consisting of homotopically trivial loops.
\end{claim}

\begin{proof}
An element in $\pi_2(\Lambda M,*)$ is represented by a map
$S^2\times S^1\to M$ that sends $\{{\rm pt}\}\times S^1$ to $x_0$.
Hence, this map factors through the quotient of $S^2\times S^1$
obtained by collapsing $\{{\rm pt}\}\times S^1$ to a point. The
resulting quotient space is homotopy equivalent to $S^2\vee S^3$,
and a map of this space into $M$ sending the wedge point to $x_0$
is, up to homotopy, the same as an element of $\pi_2(M,x_0)\oplus
\pi_3(M,x_0)$. But we are assuming that $\pi_2(M,x_0)=0$. The first
statement follows. For the second, notice that since $\pi_2(M,x_0)$
is trivial, $\pi_3(M,x_0)$ is identified with $H_3$ of the universal
covering $\widetilde M$ of $M$.  Hence, for any map of $S^2$ into
the component of $\Lambda M$ containing the trivial loops, the
resulting map $S^2\times S^1\to M$ lifts to  $\widetilde M$. The
corresponding element in $\pi_3(M,x_0)$ is the image of the
fundamental class of $S^2\times S^1$ in $H_3(\widetilde
M)=\pi_3(M)$.
\end{proof}

\begin{defn}\label{Wdefn}
Fix a homotopically trivial loop $\gamma\in \Lambda M$. We set
$A(\gamma)$ equal to the infimum of the areas of any spanning disks
for $\gamma$, where by definition a spanning disk is a Lipschitz map
$D^2\to M$ whose boundary is, up to reparameterization, $\gamma$.
Notice that $A(\gamma)$ is a continuous function of $\gamma$ in
$\Lambda M$. Also, notice that $A(\gamma)$ is invariant under
reparameterization of the curve $\gamma$. Now suppose that
$\Gamma\colon S^2\to \Lambda M$ is given with the image consisting
of homotopically trivial loops.
 We define $W(\Gamma)$ to be equal to the maximum over all $c\in S^2$ of
$A(\Gamma(c))$. More generally, given a homotopy class $\xi\in \pi_2(\Lambda
M,*)$ we define $W(\xi)$ to be equal to the infimum over all (not necessarily
based) maps $\Gamma\colon S^2\to\Lambda M$ into the component of $\Lambda M$
consisting of homotopically trivial loops representing $\xi$ of $W(\Gamma)$.
\end{defn}

Now let us formulate the analogue of Proposition~\ref{pi2} for
$\pi_3$. Suppose that ${\mathcal X}$ is a path of components of the
Ricci flow with surgery $({\mathcal M},G)$ defined for $t\in
[t_0,t_1]$. Suppose that $\pi_2({\mathcal X}(t_0),x_0)=0$ and that
$\pi_3({\mathcal X}(t_0),x_0)\not=0$. Then,  the same two conditions
hold for ${\mathcal X}(t)$ for each $t\in [t_0,t_1]$. The reason is
that at a surgery time $t$, since all the $2$-spheres in ${\mathcal
X}(t')$ ($t'<t$ but sufficiently close to $t$) along which we are
doing surgery are homotopically trivial, the result of surgery is a
disjoint union of connected components: one connected component is
homotopy equivalent to ${\mathcal X}(t')$ and all other connected
components are homotopy $3$-spheres. This means that either
${\mathcal X}(t)$ is homotopy equivalent to ${\mathcal X}(t')$ for
$t'<t$ or ${\mathcal X}(t)$ is a homotopy $3$-sphere. In either case
both homotopy group statements hold for ${\mathcal X}(t)$. Even more
is true: The distance-decreasing map ${\mathcal X}(t')\to {\mathcal
X}(t)$ given by Proposition~\ref{surgerydist} is either a homotopy
equivalence or a degree one map of ${\mathcal X}(t')\to {\mathcal
X}(t)$. In either case, it induces an injection of $\pi_3({\mathcal
X}(t'))\to \pi_3({\mathcal X}(t))$. In this way a non-zero element
in $\xi(t_0)\in\pi_3({\mathcal X}(t_0))$ produces a family of
non-zero elements $\xi(t)\in \pi_3({\mathcal X}(t))$ with the
property that under Ricci flow these elements agree and at a surgery
time $t$ the degree one map constructed in
Proposition~\ref{surgerydist} sends $\xi(t')$ to $\xi(t)$ for all
$t'<t$ sufficiently close to it. Since $\pi_2({\mathcal X}(t))$ is
trivial for all $t$, we identify $\xi(t)$ with a homotopy class of
maps of $S^2$ to $\Lambda {\mathcal X}(t)$. We now define a function
$W_\xi(t)$ by associating to each $t$ the invariant $W(\xi(t))$.

Here is the result that is analogous to Lemma~\ref{W_2}.

\begin{prop}\label{W_3fordiff}
Suppose that $({\mathcal M},G)$ is a Ricci flow with surgery as in
Theorem~\ref{MAIN}. Let ${\mathcal X}$ be a path of components of ${\mathcal
M}$ defined for all $t\in[t_0,t_1]$ with $\pi_2({\mathcal X}(t_0))=0$. Suppose
that $\xi\in \pi_3( X(t_0),*)$ is a non-trivial element. Then the function
$W_\xi(t)$ satisfies the following inequality in the sense of forward
difference quotients:
$$\frac{d W_\xi(t)}{dt}\le -2\pi-\frac{1}{2}R_{\rm
min}(t)W_\xi(t).$$ Also, for every $t\in [t_0,t_1]$ that is not a
surgery time the function $W_\xi(t)$ is continuous at $t$. Lastly,
if $t$ is a surgery time then
$$W_\xi(t)\le {\rm liminf}_{t'\rightarrow t^-}W_\xi(t').$$
\end{prop}

In the next subsection we assume this result and use it to complete
the proof.

\subsection{Proof of Theorem~\protect{\ref{extinct}} assuming
Proposition~\protect{\ref{W_3fordiff}}}

 According to Proposition~\ref{pi2} there is $T_1$ such that
every component of $M_T$ has trivial $\pi_2$ for every $T\ge T_1$.
Suppose that Theorem~\ref{extinct} does not hold for this Ricci flow
with surgery. We consider a path of components ${\mathcal X}(t)$ of
${\mathcal M}$ defined for $[T_1,T_2]$.  We shall show that there is
a uniform upper bound to $T_2$.

\begin{claim}
 ${\mathcal X}(T_1)$ has non-trivial $\pi_3$.
\end{claim}

\begin{proof}
By hypothesis the fundamental group of $M_0$ is a free product of
infinite cyclic groups and finite groups. This means that the same
is true for the fundamental group of each component of $M_t$ for
every $t\ge 0$, and in particular it is true for ${\mathcal
X}(T_0)$. But we know that $\pi_2({\mathcal X}(T_0))=0$.

\begin{claim}
Let $X$ be a compact $3$-manifold. If $\pi_1(X)$ is a non-trivial
free product or if $\pi_1(X)$ is isomorphic to $\Zee$, then
$\pi_2(X)\not=0$.
\end{claim}

\begin{proof}
See \cite{Hempel}, Theorem 5.2 on page 56 (for the case of a copy of $\Zee$)
and \cite{Hempel} Theorem 7.1 on page 66 (for the case of a free product
decomposition).
\end{proof}

Thus, it follows that $\pi_1({\mathcal X}(T_1))$ is a finite group
(possibly trivial). But a $3$-manifold with finite fundamental group
has a universal covering that is a compact $3$-manifold with trivial
fundamental group. Of course, by Poincar\'e duality any simply
connected $3$-manifold  is a homotopy $3$-sphere. It follows
immediately that $\pi_3({\mathcal X}(T_1))\cong \Zee$. This
completes the proof of the claim.
\end{proof}

Now we can apply Proposition~\ref{W_3fordiff} to our path of components
${\mathcal X}$ defined for all $t\in [T_1,T_2]$. First recall by
Theorem~\ref{MAIN} that the curvature of $({\mathcal M},G)$ is pinched toward
positive which implies that $R_{\rm min}(t)\ge (-6)/(1+4t)$. Let $w(t)$ be the
function satisfying the differential equation
$$w'(t)=-2\pi+\frac{3}{1+4t}w(t)$$
with initial condition $w(T_1)=W_\xi(T_1)$. According to
Proposition~\ref{W_3fordiff} and Proposition~\ref{fordiffquot} we
see that $W_\xi(t)\le w(t)$ for all $t\in [T_1,T_2]$. But direct
integration shows that
$$w(t)=W_\xi(T_1)\frac{(4t+1)^{3/4}}{(4T_1+1)^{3/4}}+2\pi
(4T_0+1)^{1/4}(4t+1)^{3/4}-2\pi (4t+1).$$ This clearly shows that
$w(t)$ becomes negative for $t$ sufficiently large, how large
depending only on $W_\xi(T_1)$ and $T_1$. On the other hand, since
$W_\xi(t)$ is the infimum of areas of disks, $W_\xi(t)\ge 0$ for all
$t\in [T_1,T_2]$. This proves that $T_2$ is less than a constant
that depends only on $T_1$ and on the component ${\mathcal X}(T_1)$.
Since there are only finitely many connected components of
$M_{T_1}$, this shows that $T_2$ depends only on $T_1$ and the
Riemannian manifold $M_{T_1}$. This completes the proof of
Theorem~\ref{extinct} modulo Proposition~\ref{W_3fordiff}.
\end{proof}

Thus, to complete the argument for Theorem~\ref{extinct} it remains
only to prove Proposition~\ref{W_3fordiff}.

\subsection{Continuity of $W_\xi(t)$}

In this subsection we establish the two continuity conditions for
$W_\xi(t)$ stated in Proposition~\ref{W_3fordiff}.

\begin{claim}
If $t$ is not a surgery time, then $W_\xi(t)$ is continuous at $t$.
\end{claim}

\begin{proof}
Since $t$ is not a surgery time, a family $\Gamma(t)\colon S^2\to
\Lambda {\mathcal X}(t)$ is also a family $\Gamma(t')\colon
S^2\to\Lambda{\mathcal X}(t')$ for all nearby $t'$. The minimal
spanning disks for the elements of $\Gamma(t)(x)$ are  also spanning
disks in the nearby ${\mathcal X}(t')$ and their areas vary
continuously with $t$. But the maximum of the areas of these disks
is an upper bound for $W_\xi(t)$. This immediately implies that
$W_\xi(t)$ is upper semi-continuous at $t$.

The result for lower semi-continuity is the same as in the case of
$2$-spheres. Given a time $t$ distinct from a surgery time and a
family $\Gamma\colon S^2\to \Lambda {\mathcal X}(t')$ for a time
$t'$ near $t$ we can view the family $\Gamma$ as a map to $\Lambda
{\mathcal X}(t)$. The areas of all minimal spanning disks for the
loops represented by points $\Gamma$ measured in ${\mathcal X}(t)$
are at most $(1+\eta(|t-t'|))$ times their areas measured in
${\mathcal X}(t')$, where $\eta(|t-t'|)$ is a function going to zero
as $|t-t'|$ goes to zero. This immediately implies the lower
semi-continuity at the non-surgery time $t$.
\end{proof}

\begin{claim}
Suppose that $t$ is a surgery time. Then
$$W_\xi(t)\le {\rm
liminf}_{t'\rightarrow t^-}W_\xi(t').$$
\end{claim}

\begin{proof}
This is immediate from the fact from Proposition~\ref{surgerydist} that for any
$\eta>0$ for every $t'<t$ sufficiently close to $t$ there is a homotopy
equivalence ${\mathcal X}(t')\to {\mathcal X}(t)$ which is a
$(1+\eta)$-Lipschitz map.
\end{proof}

To prove Proposition~\ref{W_3fordiff} and hence
Theorem~\ref{extinct}, it remains to prove the forward difference
quotient statement for $W_\xi(t)$ given in
Proposition~\ref{W_3fordiff}.

\subsection{A further reduction of
Proposition~\protect{\ref{W_3fordiff}}}

Let $\Gamma\colon S^2\to \Lambda {\mathcal X}(t_0)$ be a family. We
must construct an appropriate deformation of the family of loops
$\Gamma$ in order to establish Proposition~\ref{W_3fordiff}. Now we
are ready to state the more technical estimate for the evolution of
$W(\Gamma)$ under Ricci flow that will imply the forward difference
quotient result for $W_\xi(t)$ stated in
Proposition~\ref{W_3fordiff}. Here is the result that shows a
deformation as required exists.

\begin{defn}
Let $(M,g(t)),\ t_0\le t\le t_1$, be a Ricci flow on a compact
$3$-manifold. For any $a$ and any $t'\in [t_0,t_1]$ let
$w_{a,t'}(t)$ be the solution to the differential equation
\begin{equation}\label{wdiffeqn}
\frac{dw_{a,t'}}{dt}=-2\pi-\frac{1}{2}R_{\rm min}(t)w_{a,t'}(t)
\end{equation}
with initial condition $w_{a,t'}(t')=a$. We also denote $w_{a,t_0}$
by $w_a$.
\end{defn}

\begin{prop}\label{XI}
Let $(M,g(t)),\ t_0\le t\le t_1$, be a Ricci flow on a compact
$3$-manifold.  Fix a map $\Gamma$ of $S^2$ to $\Lambda M$
 whose image consists of homotopically trivial loops and $\zeta>0$. Then there
is a continuous family $\widetilde\Gamma(t),\ t_0\le t\le t_1$, of
maps $S^2\to \Lambda M$ whose image consists of homotopically
trivial loops with $[\widetilde\Gamma(t_0)]=[\Gamma]$ in
$\pi_3(M,*)$ such that for each $c\in S^2$ we have
$|A(\widetilde\Gamma(t_0)(c))-A(\Gamma(c))|<\zeta$ and furthermore,
one of the following two alternatives holds:
\begin{enumerate}
\item[(i)] The length of $\widetilde\Gamma(t_1)(c)$ is less than $\zeta$.
\item[(ii)] $A(\widetilde\Gamma(t_1)(c))\le w_{A(\widetilde\Gamma(t_0)(c))}(t_1)+\zeta$.
\end{enumerate}
\end{prop}

Before proving this result we shall show it implies the forward difference
quotient result in Proposition~\ref{W_3fordiff}. Let ${\mathcal X}$ be a path
of components. Suppose that $\pi_2({\mathcal X}(t),x_0)=0$ for all $t$. Fix
$t_0$ and fix a non-trivial element $\xi\in \pi_3({\mathcal X}(t_0),x_0)$,
which we identify with a non-trivial element in $\xi\in \pi_2(\Lambda{\mathcal
X}(t_0),*)$. Fix an interval $[t_0,t_1]$ with the property that there are no
surgery times in the interval $(t_0,t_1]$. Restricting to this interval the
family ${\mathcal X}(t)$ is a Ricci flow on ${\mathcal X}(t_0)$. In particular,
all the ${\mathcal X}(t)$ are identified under the Ricci flow.
 Let
$w(t)$ be the solution to Equation~(\ref{wdiffeqn}) with value
$w(t_0)=W_\xi({\mathcal X},t_0)$. We shall show that $W_\xi(t_1)\le
w(t_1)$. Clearly, once we have this estimate, taking limits as $t_1$
approaches $t_0$ establishes the forward difference quotient result
at $t_0$.

\begin{defn}
Let $A(t)=\int_{t'}^t\frac{1}{2}R_{\rm min}(s)ds$.
\end{defn}

Direct integration shows the following:

\begin{claim}\label{vclaim} We have
$$w_{a,t'}(t'')={\rm exp}(-A(t''))\left(a-2\pi\int_{t'}^{t''}{\rm
exp}(A(t))dt\right).$$ If $a'>a$, then for $t_0\le t'<t''\le t_1$,
we have
$$w_{a',t'}(t'')=w_{a,t'}(t'')+(a'-a){\rm exp}(-A(t'')).$$
\end{claim}

 The next thing to
establish is the following.

\begin{lem}\label{shortpi2triv}
Given a compact Riemannian manifold $(X,g)$ with $\pi_2(X)=0$. Then
there is $\zeta>0$ such that if $\xi\in \pi_3({\mathcal X})$ is
represented by a family $\Gamma\colon S^2\to\Lambda X$ with the
property that for every $c\in S^2$ the length of the loop
$\Gamma(c)$ is less than $\zeta$, then $\xi$ is the trivial homotopy
element.
\end{lem}

\begin{proof}
We choose $\zeta$ smaller than the injectivity radius of $(X,g)$.
Then any pair of points at distance less than $\zeta$ apart are
joined by a unique geodesic of length less than $\zeta$.
Furthermore, the geodesic varies smoothly with the points. Given a
map $\Gamma\colon S^2\to \Lambda X$ such that every loop of the form
$\Gamma(c)$ has length at most $\zeta$, we consider the map $f\colon
S^2\to X$ defined by $f(c)=\Gamma(c)(x_0)$, where $x_0$ is the base
point of the circle. Then we can join each point $\Gamma(c)(x)$ to
$\Gamma(c)(x_0)$ by a geodesic of length at most $\zeta$ to fill out
a map of the disk $\widehat \Gamma(c)\colon D^2\to X$. This disk is
smooth except at the point $\Gamma(c)(x_0)$. The disks
$\widehat\Gamma(c)$ fit together as $c$ varies to make a continuous
family of disks parameterized by $S^2$ or equivalently a map
$S^2\times D^2$ into $X$ whose boundary is the family of loops
$\Gamma(c)$. Now shrinking the loops $\Gamma(c)$ across the disks
$\hat\Gamma(c)$ to $\Gamma(c)(x_0)$ shows that the family $\Gamma$
is homotopic to a $2$-sphere family of constant loops at different
points of $X$. Since we are assuming that $\pi_2(X)$ is trivial,
this means the family of loops is in fact trivial as an element of
$\pi_2(\Lambda X,*)$, which  means that the original element $\xi\in
\pi_3(X)$ is trivial.
\end{proof}

Notice that this argument also shows the following:

\begin{cor}\label{shlnsharea}
Let $(X,g)$ be a compact Riemannian manifold. Given $\eta>0$ there is a
$0<\zeta<\eta/2$ such that any $C^1$-loop $c\colon S^1\to X$ of length less
than $\eta$ bounds a disk in $X$ of area less than $\eta$.
\end{cor}

Now we return to the proof that Proposition~\ref{XI} implies
Proposition~\ref{W_3fordiff}. We consider the restriction of the
path ${\mathcal X}$ to the time interval $[t_0,t_1]$. As we have
already remarked, since there are no surgery times in $(t_0,t_1]$,
this restriction is a Ricci flow and all the ${\mathcal X}(t)$ are
identified with each other under the flow.
 Let $w(t)$ be the solution to Equation~(\ref{wdiffeqn})
with initial condition $w(t_0)=W_\xi(t_0)$. There are two cases to
consider: (i) $w(t_1)\ge 0$ and $w(t_1)< 0$.

Suppose that $w(t_1)\ge 0$. Let $\eta>0$ be given. Then by Claim~\ref{vclaim}
and  Corollary~\ref{shlnsharea}, there is $0<\zeta<\eta/2$ such that the
following two conditions hold:
\begin{enumerate}
\item[(a)] Any loop in ${\mathcal X}(t_1)$ of length less than $\zeta$
bounds a disk of area less than $\eta$.
\item[(b)] For every $a\in [0,W_\xi(t_0)+2\zeta]$
the solution $w_{a}$ satisfies $w_{a}(t_1)<w(t_1)+\eta/2$.
\end{enumerate}

Now fix a map $\Gamma\colon S^2\to \Lambda{\mathcal X}(t_0)$, whose
image consists of homotopically trivial loops, with $[\Gamma]=\xi$,
and with $W(\Gamma)<W_\xi(t_0)+\zeta$. According to
Proposition~\ref{XI} there is a one-parameter family
$\widetilde\Gamma(t),\ t_0\le t\le t_1$, of maps $S^2\to
\Lambda{\mathcal X}(t)$, whose images consist of homotopically
trivial loops, with $[\widetilde\Gamma(t_0)]=[\Gamma]=\xi$ such that
for every $c\in S^2$ we have
$A(\widetilde\Gamma(t_0)(c))<A(\Gamma(c))+\zeta$ and one of the
following holds
\begin{enumerate}
\item[(i)] the length of $\widetilde\Gamma(t_1)(c)$ is less than $\zeta$, or
\item[(ii)] $$A(\widetilde\Gamma(t_1)(c))<
w_{A(\widetilde\Gamma(t_0)(c))}(t_1)+\zeta.$$
\end{enumerate}
Since $A(\widetilde\Gamma(t_0)(c))<A(\Gamma(c))+\zeta<W_\xi(t_0)+2\zeta$, it
follows from our choice of $\zeta$ that for every $c\in S^2$ either
\begin{enumerate}
\item[(a)] $\widetilde\Gamma(t_1)(c)$ has length less than $\zeta$ and hence bounds a disk
of area less than $\eta$, or
\item[(b)]
$A(\widetilde\Gamma(t_1)(c))<w_{W_\xi(t_0)+2\zeta}(t_1)+\zeta<w(t_1)+\eta/2+\eta/2=
w(t_1)+\eta$.
\end{enumerate}
Since we are assuming that $w(t_1)\ge 0$, it now follows that for
every $c\in S^2$ we have $A(\widetilde\Gamma(t_1)(c))<w(t_1)+\eta$,
and hence $W(\widetilde\Gamma(t_1))< w(t_1)+\eta$. This shows that
for every $\eta>0$ we can find a family $\widetilde\Gamma(t)$ with
$\widetilde\Gamma(t_0)$ representing $\xi$ and with
$W(\widetilde\Gamma(t_1))< w(t_1)+\eta$. This completes the proof of
Proposition~\ref{XI} when $w(t_1)\ge 0$.

Now suppose that $w(t_1)<0$. In this case, we must derive a contradiction since
clearly it must be the case that for any one-parameter family
$\widetilde\Gamma(t)$ we have $W(\widetilde\Gamma(t_1))\ge 0$. We fix $\eta>0$
such that $w(t_1)+\eta<0$. Then using Lemma~\ref{shortpi2triv} and
Claim~\ref{vclaim}, we fix $\zeta$ with $0<\zeta<\eta/2$ such that:
\begin{enumerate}
\item[(i)] If $\Gamma\colon S^2\to \Lambda {\mathcal X}(t_1)$ is a family
of loops and each loop in the family is of length less than $\zeta$,
then the family is homotopically trivial.
\item[(ii)] For any $a\in [0,W_\xi(t_0)+2\zeta]$ we have $w_a(t_1)<w(t_1)+\eta/2$.
\end{enumerate}

We fix a map $\Gamma\colon S^2\to {\mathcal X}(t_0)$ with
$[\Gamma]=\xi$ and with $W(\Gamma)<W_\xi(t_0)+\zeta$. Now according
to Proposition~\ref{XI} there is a family of maps
$\widetilde\Gamma(t)\colon S^2\to \Lambda{\mathcal X}(t)$ with
$[\widetilde\Gamma(t_0)]=[\Gamma]=\xi$ and for every $c\in S^2$ we
have $A(\widetilde\Gamma(t_0)(c))<A(\Gamma(c))+\zeta$ and also
either $A(\widetilde\Gamma(t_1)(c))\le
w_{A(\widetilde\Gamma(t_0)(c))}(t_1)+\zeta$ or the length of
$\widetilde\Gamma(t_1)(c)$ is less than $\zeta$. It follows that for
every $c\in S^2$ we have $A(\widetilde\Gamma(t_0)(c))\le
W(\Gamma)+\zeta<W_\xi(t_0)+2\zeta$. From the choice of $\zeta$ this
means that
$$A(\widetilde\Gamma(t_1)(c))< w(t_1)+\eta/2+\zeta<w(t_1)+\eta<0$$ if
the length of $\widetilde\Gamma(t_1)(c)$ is at least $\zeta$. Of course, by
definition $A(\widetilde\Gamma(t_1)(c))\ge 0$ for every $c\in S^2$. This
implies  that for every $c\in S^2$ the loop $\widetilde\Gamma(t_1)(c)$ has
length less than $\zeta$. By Lemma~\ref{shortpi2triv} this implies that
$\widetilde\Gamma(t_1)$ represents the trivial element in
$\pi_2(\Lambda{\mathcal X}(t_1))$, which is a contradiction.

At this point, all that it remains to do in order to complete the
proof of Theorem~\ref{extinct} is to establish Proposition~\ref{XI}.
The rest of this chapter is devoted to doing that.

\section{Curve-shrinking flow}\label{sect:csf}

Given $\Gamma$, the idea for constructing the one-parameter family
$\widetilde \Gamma(t)$ required by Proposition~\ref{XI} is to evolve
an appropriate approximation $\widetilde \Gamma(t_0)$ of $\Gamma$ by
the curve-shrinking flow. Suppose that $(M,g(t)),\ t_0\le t\le t_1$,
is a Ricci flow of compact manifolds and that $c\colon S^1\times
[t_0,t_1]\to (M,g(t_0))$ is a family of parameterized, immersed
$C^2$-curves. We denote by $x$ the parameter on the circle. Let
$X(x,t)$ be the tangent vector $\partial c(x,t)/\partial x$ and let
$S(x,t)=X(x,t)/(|X(x,t)|_{g(t)})$ be the unit tangent vector to $c$.
We denote by $s$ the arc length parameter on $c$. We set
$H(x,t)=\nabla _{S(x,t)}S(x,t)$, the curvature vector of $c$ with
respect to the metric $g(t)$.
 We define the curve-shrinking
flow by
$$\frac{\partial c(x,t)}{\partial t}=H(x,t),$$
where $c(x,t)$ is a one-parameter family of curves and $H(x,t)$ is
the curvature vector of the curve $c(\cdot,t)$ at the point $x$ with
respect to the metric $g(t)$. We denote by $k(x,t)$ the curvature
function: $k(x,t)=|H(x,t)|_{g(t)}$. We shall often denote the
one-parameter family of curves by $c(\cdot,t)$. Notice that if
$c(x,t)$ is a curve-shrinking flow and if $x(y)$ is a
reparameterization of the domain circle then $c'(y,t)=c(x(y),t)$ is
also a curve-shrinking flow.

\begin{claim}\label{curvshrexist}
For any immersed $C^2$-curve $c\colon S^1\to (M,g(t_0))$ there is a
curve-shrinking flow $c(x,t)$ defined for $t\in [t_0,t_1')$ for some
$t_1'>t_0$ with the property that each $c(\cdot,t)$ is an immersion.
Either the curve-shrinking flow extends to a curve-shrinking flow
that is a family of immersions defined at $t_1'$ and beyond, or
${\rm max}_{x\in S^1}k(x,t)$ blows up as $t$ approaches $t_1'$ from
below.
\end{claim}

For a proof of this result, see Theorem 1.13 in \cite{AG}.

\subsection{The proof of Proposition~\protect{\ref{XI}} in a simple case}

The main technical hurdle to overcome is that in general the curve
shrinking flow may not exist if the original curve is not immersed
and even if the original curve is immersed the curve-shrinking flow
can develop singularities, where the curvature of the curve goes to
infinity. Thus, we may  not be able to define the curve-shrinking
flow as a flow defined on the entire interval $[t_0,t_1]$, even
though the Ricci flow is defined on this entire interval. But to
show the idea of the proof,  let us suppose for a moment that the
starting curve is embedded and that no singularities develop in the
curve-shrinking flow and show how to prove the result.

\begin{lem}\label{embcase}
Suppose that $c\in \Lambda M$ is a homotopically trivial, embedded
$C^2$-loop. and suppose that there is a curve-shrinking flow
$c(x,t)$ defined for all $t\in [t_0,t_1]$ with each $c(\cdot,t)$
being an embedded smooth curve. Consider the function $A(t)$ which
assigns to $t$ the minimal area of a spanning disk for $c(\cdot,t)$.
Then $A(t)$ is a continuous function of $t$ and
$$\frac{dA}{dt}(t)\le -2\pi-\frac{1}{2}R_{\rm min}(t)A(t)$$
in the sense of forward difference quotients.
\end{lem}

\begin{proof}
According to results of Hildebrandt and Morrey, \cite{Hildebrandt}
and \cite{Morrey}, for each $t\in [t_0,t_1]$, there is a smooth
minimal disk spanning $c(\cdot,t)$. Fix $t'\in [t_0,t_1)$ and
consider a smooth minimal disk $D\to (M,g(t'))$ spanning
$c(\cdot,t)$. It is immersed, see \cite{HardtSimon} or \cite{GL}.
The family $c(\cdot,t)$ for $t$ near $t'$ is an isotopy of
$c(\cdot,t')$. We can extend this to an ambient isotopy
$\varphi_t\colon M\to M$ with $\varphi_{t'}={\rm Id}$. We impose
coordinates $\{x_\alpha\}$ on $D$; we let $h_{\alpha\beta}(t')$ be
the metric induced  on $\varphi_{t'}(D)$ by $g(t')$, and we let $da$
be the area form induced by the Euclidean coordinates on $D$. We
compute
$$\frac{d}{dt}\bigl|_{t=t'}\bigr.{\rm Area}(\varphi_t(D))
=\frac{d}{dt}\bigl|_{t=t'}\bigr.\int_{\varphi_t(D)}\sqrt{{\rm
det}(h_{\alpha\beta})(t)}da.$$ Of course,
\begin{eqnarray*}
\frac{d}{dt}\bigl|_{t=t'}\bigr.\int_{\varphi_t(D)}\sqrt{{\rm
det}(h_{\alpha\beta}(t))}da & = & -\int_{\varphi_{t'}(D)}\Bigl({\rm
Tr}\,{\rm
Ric}^T\Bigr)\sqrt{{\rm det}(h_{\alpha\beta}(t))}da \\
& & +\int_{\varphi_{t'}(D)}{\rm div}\left(\frac{\partial
\varphi_{t'}}{dt}\right)^T\sqrt{{\rm det}(h_{\alpha\beta}(t))}da.
\end{eqnarray*}
 Here, ${\rm Ric}^T$
denotes the restriction of the Ricci curvature of $g(t')$ to the
tangent planes of $\varphi_{t'}(D)$ and
$\frac{\partial(\varphi_{t'})^T}{\partial t}$ is the component of
$\varphi_{t'}$ tangent to $\varphi_{t'}(D)$. Setting $\hat A$ equal
to the second fundamental form of $\varphi_{t'}(D)$, using the fact
that $\varphi_{t'}(D)$ is minimal and arguing as in the proof of
Claim~\ref{varareaform}, we have
\begin{eqnarray*}
\lefteqn{-\int_{\varphi_{t'}(D)}\Bigl({\rm Tr}\,{\rm
Ric}^T\Bigr)\sqrt{{\rm det}(h_{\alpha\beta}(t'))}da } \\
 & = &
-\int_{\varphi_{t'}(D)}K_{\varphi_{t'}(D)}da
-\f1{1}{2}\int_{\varphi_{t'}(D)}(|\hat A|^2+R)da \\
& \le & -\int_{\varphi_{t'}(D)}K_{\varphi_{t'}(D)}\sqrt{{\rm
det}(h_{\alpha\beta}(t'))}da-\f1{1}{2}{\rm
 Area}{\varphi_{t'}(D)}\min_{x\in M}\{R(x,t')\}.
 \end{eqnarray*}

Integration by parts shows that
$$\int_{\varphi_{t'}(D)}{\rm div}\left(\frac{\partial
\varphi_{t'}}{dt}\right)^T\sqrt{{\rm
det}(h_{\alpha\beta}(t'))}da=-\int_{ \varphi_{t'}(\partial
D)}\left(\frac{d\varphi_t}{dt}\bigl|_{t=t'}\bigr.\right)\cdot n
ds,$$ where $n$ is the inward pointing normal vector to
$\varphi_{t'}(D)$ along $\varphi_{t'}(\partial D)$. Of course, by
definition, if the variation along the boundary is given by the
curve-shrinking flow, then along $\varphi_{t'}(\partial D)$ we have
$$\left(\frac{d\varphi_t}{dt}\bigl|_{t=t'}\bigr.\right)\cdot n=k_{\rm geod}.$$

Thus, we have
\begin{eqnarray*}\lefteqn{\frac{d}{dt}\bigl|_{t=t'}\bigr.\int_{\varphi_t(D)}\sqrt{{\rm
det}(h_{\alpha\beta}(t))}da } \\\ & \le &
-\int_{\varphi_{t'}(D)}K_{\varphi_{t'}(D)}da-\int_{\varphi_{t'}(\partial
D)}k_{\rm geod}ds-\frac{1}{2}R_{\rm min}(t'){\rm
Area}(\varphi_{t'}(D)).\end{eqnarray*} Of course, the Gauss-Bonnet
theorem allows us to rewrite this as
$$\frac{d}{dt}\bigl|_{t=t'}\bigr.\int_{\varphi_t(D)}\sqrt{{\rm
det}(h_{\alpha\beta}(t))}da \le -2\pi-\frac{1}{2}R_{\rm min}(t'){\rm
Area}(\varphi_{t'}(D)).$$
\end{proof}

Let $\psi(t)$ be the solution to the ODE
$$\psi'(t)=-2\pi-\frac{1}{2}R_{\rm min}(t)\psi(t)$$
with $\psi(t^-)=A(t^-)$. The following is immediate from the
previous lemma and Lemma~\ref{fordiffquot}.

\begin{cor}
With notation and assumptions as above, if the curve-shrinking flow
is defined on the interval $[t^-,t^+]$ and if the curves
$c(\cdot,t)$ are embedded for all $t\in [t^-,t^+]$ then
$$A(t^+)\le \psi(t^+).$$
\end{cor}

Actually, the fact that the loops in the curve-shrinking flow are
embedded is not essential in dimensions $\ge 3$.

\begin{lem}\label{nonembeasycase}
Suppose that the dimension of $M$ is at least $3$ and that
$c(\cdot,t)$ is a $C^2$-family of homotopically trivial, immersed
curves satisfying the curve-shrinking equation defined for $t^-\le
t\le t^+$. For each $t$, let $A(t)$ be the infimum of the areas of
spanning disks for  $c(\cdot,t)$. Then $A(t)$ is a continuous
function and, with $\psi$ as above, we have
$$A(t^+)\le \psi(t^+).$$
\end{lem}

\begin{proof} We first remark that continuity has already been established.
To show the inequality, we begin  with a claim.

\begin{claim}
It suffices to prove the following for every $\delta>0$. There is a
$C^2$-family $\hat c(x,t)$ of immersions within $\delta$ in the
$C^2$-topology to $c(x,t)$ defined on the interval $[t^-,t^+]$ such
that
$$A(t^+)\le \psi_{\delta,\hat c}(t^+)$$
where $\psi_{\delta,\hat c}$ is the solution of the ODE
$$\psi_{\delta,\hat c}'(t)=-2\pi+2\delta L_{\hat c}(t)-\frac{1}{2}R_{\rm
min}(t)\psi_{\delta,\hat c}(t)$$ with value $A(\hat c(t^-))$ at
$t^-$, and where $L_{\hat c}(t)$ denotes the length of the loop
$\hat c(\cdot,t)$.
\end{claim}

\begin{proof} (of the claim) Suppose that for each $\delta$ there is such
a $C^2$-family as in the statement of the claim. Take a sequence
$\delta_n$ tending to zero, and let $\hat c_n(\cdot,t)$ be a family
as in the claim for $\delta_n$. Then by the continuity of the
infimum of areas of the spanning disk in the $C^1$-topology, we see
that
$${\rm lim}_{n\rightarrow \infty}A(\hat c_n(\cdot, t^\pm))=A(c(\cdot,t^\pm)).$$
Since the $\hat c_n(x,t)$ converge in the $C^2$-topology to
$c(x,t)$, the lengths $L(\hat c_n(t))$ are uniformly bounded and the
$A(\hat c_n(t^-))$ converge to $A(c(t))$. Thus, the
$\psi_{\delta_n,\hat c_n}$ converge uniformly to $\psi$ on
$[t^-,t^+]$, and taking limits shows the required inequality for
$A(c(t))$, thus proving the claim.
\end{proof}

Now we return to the proof of the lemma.  Let $\hat c(x,t)$ be a
generic $C^2$-immersion  sufficiently close to $c(x,t)$ in the
$C^2$-topology so that the following hold:
\begin{enumerate}
\item[(1)] the difference of the curvature of $\hat c$ and of $c$ at
every $(x,t)$ is a vector of length less than $\delta$,
\item[(2)] the difference of $\partial \hat c/\partial t$ and $\partial
c/\partial t$ is a vector of length less than $\delta$,
\item[(3)] the ratio of the arc lengths of $\hat c$ and $c$ at every
$(x,t)$ is between $(1-\delta)$ and $(1+\delta)$.
\end{enumerate}
The generic family $\hat c(x,t)$ consists of embedded curves for all
but a finite number of $t\in [t^-,t^+]$ and at the exceptional $t$
values the curve is immersed. Let $t_1<t_2<\cdots <t_k$ be the
values of $t$ for which $\hat c(\cdot,t)$ is not embedded. We set
$t_0=t^-$ and $t_{k+1}=t^+$. Notice that it suffices to show that
$$A(\hat c(t_{i+1}))-A(\hat c(t_{i}))\le
\psi_{\delta,\hat c}(t_{i+1})-\psi_{\delta,\hat c}(t_i)$$ for
$i=0,\ldots,k$. To establish this inequality for the interval
$[t_i,t_{i+1}]$,  by continuity it suffices to establish the
corresponding inequality for every compact subinterval contained in
the interior of this interval. This allows us to assume that the
approximating family is a family of embedded curves. Let the
endpoints of the parameterizing interval  be denoted $a$ and $b$.
Fix $t'\in [a,b]$ and let $D$ be a minimal disk spanning $\hat
c(\cdot, t')$, and let $\varphi_t$ be an isotopy as in the argument
given the proof of Lemma~\ref{embcase}. According to this argument
we have
$$\frac{d}{dt}A(\hat c(t))|_{t=t'}\le -2\pi -\frac{1}{2}R_{\rm
min}(t')A(c(t'))+\int_{c(x,t')}\left[k_{\rm geod}(\hat
c)-\left(\frac{d\varphi_t}{dt}|_{t=t'}\right)\cdot n\right] ds$$ in
the sense of forward difference quotients. The restriction of
$\frac{d\varphi_t}{dt}|_{t=t'}$ to the boundary of $D$ agrees with
$\partial \hat c(x,t)/\partial t$. Hence, by our conditions on the
approximating family, and since for $c(\cdot,t)$ the corresponding
quantities are equal,
$$\left|k_{\rm geod}(\hat c)-\left(\frac{d\varphi_t}{dt}|_{t=t'}\right)\cdot
n\right|<2\delta.$$ Integrating over the circle implies that
$$\frac{d}{dt}A(\hat c(t))|_{t=t'}\le -2\pi -\frac{1}{2}R_{\rm
min}(t')A(c(t'))+2\delta L_{\hat c}(t).$$
 The result is then
immediate from Lemma~\ref{fordiffquot}.
\end{proof}

\subsection{Basic estimates for curve-shrinking}

Let us establish some elementary formulas. To simplify the formulas we often
drop the variables $x,t$ from the notation, though they are understood to be
there.

\begin{lem}\label{1stcurvshr}
Assume that $(M,g(t)),\ t_0\le t\le t_1$, is a Ricci flow and that
$c=c(x,t)$ is a solution to the curve-shrinking flow. We have vector
fields $X=\partial /\partial x$ and $H=\partial /\partial t$ defined
on the domain surface. We denote by $|X|_{c^*g}^2$ the function on
the domain surface whose value at $(x,t)$ is $|(X(x,t))|_{g(t)}^2$.
We define $S=|X|_{c^*g}^{-1}X$, the unit vector in the $x$-direction
measured in the evolving metric. Then,
$$\frac{\partial}{\partial t}(|X|_{c^*g}^2)(x,t)=-2{\rm Ric}_{g(t)}(X(x,t),X(x,t))-2k^2|X(x,t)|_{g(t)}^2,$$
and
$$[H,S](x,t)=\left(k^2+{\rm Ric}_{g(t)}(S(x,t),S(x,t))\right)S(x,t).$$
\end{lem}

\begin{proof}
Notice that as $t$ varies $|X|_{c^*g}^2$ is not the norm of the
vector field $X$ with respect to the pullback of a fixed metric
$g(t)$. On the other hand, when we compute $\nabla_H X$ at a point
$(x,t)$ we are taking a covariant derivative with respect to the
pullback of a fixed metric $g(t)$ on the surface. Hence, in
computing $H(|X|_{c^*g}^2)$ the usual Leibniz rule does not apply.
In fact, there are two contributions to $H(|X|_{c^*g}^2)$: one, the
usual Leibniz rule differentiating in a frozen metric $g(t)$ and the
other coming from the effect on $|X|_{c^*g}^2$ of varying the metric
with $t$. Thus, we have
$$H(|X|_{c^*g}^2)(x,t)=-2{\rm Ric}_{c^*g(t)}(X(x,t),X(x,t))
+2\langle\nabla_HX(x,t),X(x,t)\rangle_{c^*g(t)}.$$ Since $t$ and $x$
are coordinates on the surface swept out by the family of curves,
$\nabla_HX=\nabla_XH$, and hence the second term on the right-hand
side of the previous equation can be rewritten as
$2\langle\nabla_XH(x,t),X(x,t)\rangle_{c^*g(t)}$. Since $X(x,t)$ and
$H(x,t)$ are orthogonal in $c^*g(t)$ and since $X=|X|_{c^*g}S$,
computing covariant derivatives in the metric $c^*g(t)$, we have
\begin{eqnarray*}
2\langle\nabla_XH,X\rangle_{c^*g(t)} & = & -2\langle H,\nabla_XX\rangle_{c^*g(t)}\\
 & = & -2\langle H,|X|_{g(t)}^2\nabla_SS\rangle_{c^*g(t)}-
 2\langle H,|X|_{g(t)}S(|X|_{c^*g})S\rangle_{c^*g(t)} \\
 & = & -2\langle H,H\rangle_{c^*g(t)}|X|_{g(t)}^2\\
 & = & -2k^2|X|_{c^*g}^2.
 \end{eqnarray*}
This proves the first inequality. As for the second, since $X$ and
$H$ commute we have
$$[H,S]=[H,|X|_{c^*g}^{-1}X]=H\left((|X|_{c^*g}^2)^{-1/2}\right)X=
\frac{-1}{2\left(|X|^2_g\right)^{3/2}}H(|X|_{c^*g})^2)X.$$ According to the
first equation, we can rewrite this as
$$[H,S](x,t)=\left(k^2+{\rm Ric}_{c^*g(t)}(S(x,t),S(x,t))\right)S(x,t).$$
\end{proof}

Now let us compute the time derivative of $k^2$. In what follows we drop the
dependence on the metric $c^*g(t)$ from all the curvature terms, but it is
implicitly there.

\begin{lem}\label{2ndcurvshr}
\begin{eqnarray*}
\frac{\partial }{\partial t}k^2 & = & \frac{\partial ^2}{\partial
s^2}(k^2)-2\langle(\nabla_XH)^\perp,(\nabla_SH)^\perp\rangle_{c^*g}+2k^4\\
& & -2{\rm Ric}(H,H) +4k^2{\rm Ric}(S,S)+2{\rm
Rm}(H,S,H,S),\end{eqnarray*} where the superscript $\perp$ means the
image under projection to the orthogonal complement of $X$.
\end{lem}

\begin{proof}
Using the same conventions as above for the function $|H|_{c^*g}$
and but leaving the metric implicit, we have
\begin{equation}\label{*1}
\frac{\partial }{\partial t}k^2=\frac{\partial }{\partial
t}(|H|^2_{c^*g})=-2{\rm
Ric}(H,H)+2\langle\nabla_HH,H\rangle_{c^*g}.\end{equation}
 Now we compute
(using the second equation from Lemma~\ref{1stcurvshr})
\begin{eqnarray*}
\nabla_HH & = & \nabla_H\nabla_SS \\
& = & \nabla_S\nabla_HS+\nabla_{[H,S]}S+
{\mathcal R}(H,S)S \\
& = & \nabla_S\nabla_SH+\nabla_S([H,S])+\nabla_{[H,S]}S+{\mathcal R}(H,S)S \\
& = & \nabla_S\nabla_SH+\nabla_S\left((k^2+{\rm
Ric}(S,S))S\right)+(k^2+{\rm Ric}(S,S))\nabla_SS+
{\mathcal R}(H,S)S \\
& = & \nabla_S\nabla_SH+2(k^2+{\rm Ric}(S,S))H+S(k^2+{\rm
Ric}(S,S))S+ {\mathcal R}(H,S)S.
\end{eqnarray*}
Using this, and the fact that $\langle H,S\rangle_{c^*g}=0$, we have
\begin{equation}\label{*2}
2\langle\nabla_HH,H\rangle_{c^*g}=2g(\nabla_S\nabla_SH,H)+4k^4+4k^2{\rm
Ric}(S,S))+2{\rm Rm}(H,S,H,S).\end{equation}
 On the other hand,
\begin{equation}\label{*3}
S(S(\langle
H,H\rangle_{c^*g}))=2\langle\nabla_S\nabla_SH,H\rangle_{c^*g}+2\langle\nabla_SH,\nabla_SH\rangle_{c^*g}.
\end{equation}

 We write
$$\nabla_SH=(\nabla_SH)^\perp+\langle\nabla_SH,S\rangle_{c^*g}S.$$
Since $H$ and $S$ are orthogonal, we have
$\langle\nabla_SH,S\rangle_{c^*g}=-\langle
H,\nabla_SS\rangle_{c^*g}=-\langle H,H\rangle$. Thus, we have
$$\nabla_SH=(\nabla_SH)^\perp-\langle H,H\rangle_{c^*g}S. $$
It follows that
$$-2\langle \nabla_SH,\nabla_SH\rangle_{c^*g}=
-2\langle(\nabla_SH)^\perp,(\nabla_SH)^\perp\rangle_{c^*g}-2k^4.$$ Substituting
this into Equation~(\ref{*3}) gives
\begin{equation}\label{*4}
2\langle\nabla_S\nabla_SH,H\rangle_{c^*g}=S(S(|H|_{c^*g}^2))
-2\langle(\nabla_SH)^\perp,(\nabla_SH)^\perp\rangle_{c^*g}-2k^4.
\end{equation}
Plugging this into Equation~(\ref{*2})  and using
Equation~(\ref{*1}) yields
\begin{eqnarray*}
\frac{\partial}{\partial t}k^2 & = & -2{\rm Ric}(H,H)+S(S\langle
H,H\rangle_{c^*g})-2\langle(\nabla_SH)^\perp,(\nabla_SH)^\perp\rangle_{c^*g}
\\
& & +2k^4+4k^2{\rm Ric}(S,S)+2{\rm Rm}(H,S,H,S).\end{eqnarray*} Of course,
$S(S(\langle H,H\rangle_{c^*g}))=(k^2)''$ so that this gives the result.
\end{proof}

Grouping together the last three terms in the statement of the
previous lemma, we can rewrite the result as
\begin{equation}\label{*4'}\frac{\partial }{\partial t}k^2\le
(k^2)''-2\langle(\nabla_SH)^\perp,(\nabla_SH)^\perp\rangle_{c^*g}+2k^4+\widehat
Ck^2,\end{equation} where the primes refer to the derivative with
respect to arc length along the curve and $\widehat C$ is a constant
depending only on an upper bound for the norm of the sectional
curvatures  of the ambient manifolds in the Ricci flow.

\begin{claim}\label{3rdcurvshr}
 There is a constant $C_1<\infty$ depending only on an
upper bound for the norm of the sectional curvatures of the ambient
manifolds in the Ricci flow $(M,g(t)),\ t_0\le t\le t_1$, such that
$$\frac{\partial}{\partial t} k\le k''+k^3+C_1k.$$
\end{claim}

\begin{proof}
We set $C_1=\widehat C/2$, where $\widehat C$ is as in
Inequality~\ref{*4'}. It follows from Inequality~(\ref{*4'}) that
\begin{equation}
\label{*5} 2k\frac{\partial k}{\partial t}\le
2kk''+2(k')^2+2k^4-2\langle(\nabla_SH)^\perp,(\nabla_SH)^\perp\rangle_{c^*g}+\hat
Ck^2.\end{equation} Since $k^2=\langle H,H\rangle_{c^*g}$, we see
that $(k^2)'=2\langle\nabla_SH,H\rangle_{c^*g}$. Since $H$ is
perpendicular to $S$, this can be rewritten as
$(k^2)'=2\langle(\nabla_SH)^\perp,H\rangle_{c^*g}$. It follows that
$$k'=\frac{\langle(\nabla_SH)^\perp,H\rangle_{c^*g}}{|H|_{c^*g}}.$$
Hence,
$$(k')^2\le \frac{\langle(\nabla_SH)^\perp,H\rangle_{c^*g}^2}{|H|_{c^*g}^2}\le
\langle(\nabla_SH)^\perp,(\nabla_SH)^\perp\rangle_{ c^*g}.$$ Plugging this into
Equation~(\ref{*5}) gives
$$\frac{\partial k}{\partial t}\le k''+k^3+C_1k.$$
\end{proof}

Now we define the {\em total length} of the curve $c(x,t)$,
$$L(t)=\int |X|_{c^*g}dx=\int ds.$$
 We also define the {\em total curvature} of the curve $c(x,t)$,
$$\Theta(t)=\int k|X|_{c^*g}dx=\int kds.$$

\begin{lem}\label{1835} There is a constant $C_2<\infty$ depending only on an
upper bound for the norm of the sectional curvatures of the ambient
manifolds in the Ricci flow such that
\begin{equation}\label{Lformula}
\frac{d}{dt}L\le \int (C_2-k^2)ds\end{equation} and
$$\frac{d}{dt}\Theta\le C_2\Theta.$$
\end{lem}

\begin{proof}
$$\frac{d}{dt}L=\int\frac{\partial}{\partial t}\sqrt{|X|_{c^*g}^2}dx.$$
By Lemma~\ref{1stcurvshr} we have $$\frac{d}{dt}L=\int
\frac{1}{2|X|_{c^*g}}\left(-2{\rm Ric}(X,X)-2k^2|X|_{c^*g}^2\right)dx.$$ Thus,
\begin{equation}\label{dLdteqn}\frac{d}{dt}L=\int (-{\rm Ric}(S,S)-k^2)|X|_{c^*g}dx=\int(-{\rm
Ric}(S,S)-k^2)ds.\end{equation} The first inequality in the lemma
then follows by taking $C_2$ to be an upper bound for the norm of
${\rm Ric}_{g(t)}$.

Now let us consider the second inequality in the statement.
$$\frac{d}{dt}\Theta=\int\frac{\partial}{\partial t}(k|X|_{c^*g})dx=\int\left(\frac{\partial k}{\partial
t}|X|_{c^*g} +k\frac{\partial |X|_{c^*g}}{\partial t}\right)dx.$$ Thus, using
Claim~\ref{3rdcurvshr} and the first equation in Lemma~\ref{1stcurvshr} we have
\begin{eqnarray*}
\frac{d}{dt}\Theta & \le &  \int (k''+k^3+C_1k)ds+\int
\frac{k}{2|X|_{c^*g}}(-2{\rm Ric}(X,X)-2k^2|X|_{c^*g}^2)dx \\
& = & \int(k''+k^3+C_1k)ds-\int k({\rm Ric}(S,S)+k^2)ds \\& = &
\int(k''+C_1k-k\,{\rm Ric}(S,S))ds.
\end{eqnarray*}
Since $\int k''ds=0$ by the fundamental theorem of calculus, we get
$$\frac{d}{dt}\Theta\le  C_2\Theta,
$$
 for an
appropriate constant $C_2$ depending only on an upper bound for the
norm of the sectional curvatures of the ambient family $(M,g(t))$.
\end{proof}

\begin{cor}\label{lncurvgrowth}  The following holds for the constant $C_2$ as in the previous
lemma. Let $c(x,t)$ be a curve-shrinking flow, let $L(t)$ be the
total length of $c(t)$ and let $\Theta(t)$ be the total curvature of
$c(t)$. Then for any $t_0\le t'<t''\le t_1$ we have $$L(t'')\le
L(t')e^{C_2(t''-t')}$$
$$\Theta(t'')\le \Theta(t')e^{C_2(t''-t')}.$$
\end{cor}

\subsection{Ramp solutions in $ M\times S^1$}

As we pointed out in the beginning of Section~\ref{sect:csf} the
main obstacle we must overcome is that the curve-shrinking flow does
not always exist for the entire time interval $[t_0,t_1]$. The
reason is the following: Even though, as we shall see, it is
possible to bound the total curvature of the curve-shrinking flow in
terms of the total curvature of the initial curve and the ambient
Ricci flow, there is no pointwise estimate on the curvature for the
curve-shrinking flow. The idea for dealing with this problem, which
goes back to \cite{AG}, is to replace the original situation of
curves in a manifold with graphs by taking the product of the
manifold with a circle and using ramps. We shall see that in this
context the curve-shrinking flow always exists. The problem then
becomes to transfer the information back from the flows of ramps to
the original manifold.

Now suppose that the Ricci flow is of the form $(M, g(t))\times
(S^1_\lambda,ds^2)$ where $(S^1_\lambda,ds^2)$ denotes the  circle
of length $\lambda$. Notice that the sectional curvatures of this
product flow depend only on the sectional curvatures of $(M, g(t))$
and, in particular, are independent of $\lambda$. Let $U$ denote
vector field made up of unit tangent vectors in the direction of the
circle factors. Let $u(x,t)=\langle S,U\rangle_{g(t)}$.

\begin{claim}\label{rampcomp}
$$\frac{\partial u}{\partial t}=u''+(k^2+{\rm Ric}(S,S))u\ge u''-C'u,$$
where $C'$ is an upper bound for the norm of the Ricci curvature of
$(M,g(t))$.
\end{claim}

\begin{proof}
Since $U$ is a constant vector field and hence parallel along all curves and
since ${\rm Ric}(V,U)=0$ for all tangent vectors $V$, by Lemma~\ref{1stcurvshr}
we have
\begin{eqnarray*}
\frac{\partial}{\partial t}\langle S,U\rangle_{g(t)} & = & -2{\rm
Ric}(S,U)
+\langle dc(\nabla_HS),U\rangle_{g(t)} \\
& = & \langle dc(\nabla_HS),U\rangle_{g(t)}
=  \langle dc([H,S]+\nabla_SH),U\rangle_{g(t)} \\
& = & (k^2+{\rm Ric}(S,S))u+\langle dc(\nabla_SH),U\rangle_{g(t)} \\
& = & (k^2+{\rm Ric}(S,S))u+S(dc(\langle H),U\rangle_{g}) \\
 & = & (k^2+{\rm Ric}(S,S))u+S(\langle dc(\nabla_SS),U\rangle_g) \\
&  = & (k^2+{\rm Ric}(S,S))u+S(S(u))  =  (k^2+{\rm Ric}(S,S))u+u''.
 \end{eqnarray*}
\end{proof}

\begin{defn}
A curve $c\colon S^1\to  M\times S^1_\lambda$ is said to be a {\em
ramp} if $u$ is strictly positive.
\end{defn}

The main results of this section show that the curve-shrinking flow
is much better behaved for ramps  than for the general smooth curve.
First of all, as the next corollary shows, the curve-shrinking flow
applied to a ramp produces a one-parameter families of ramps. The
main result of this section shows that for any ramp as initial
curve, the curve-shrinking flow does not develop singularities as
long as the ambient Ricci flow does not.

\begin{cor}\label{ubd}
If $c(x,t),\ t_0\le t< t'_1<\infty$, is a solution of the curve
shrinking flow in $( M, g(t))\times (S^1_\lambda,ds^2)$ and if
$c(t_0)$ a ramp, then $c(t)$ is a ramp for all $t\in [t_0,t'_1)$.
\end{cor}

\begin{proof}
From  the equation in Claim~\ref{rampcomp}, we see that for  $C'$ an
upper bound for the norm of the Ricci curvature, we have
$$\frac{\partial}{\partial t}\left(e^{C't}u\right)\ge \left(e^{C't}u\right)''.$$
It now follows from a standard maximum principle argument that the
minimum value of $e^{C't}u$ is a non-decreasing function of $t$.
Hence, if $c(t_0)$ is a ramp then each $c(t)$ is a ramp and in fact
$u(x,t)$ is uniformly bounded away from zero in terms of the minimum
of $u(x,t_0)$ and the total elapsed time $t_1-t_0$.
\end{proof}

\begin{lem}\label{csforramps}
Let $(M,g(t)),\ t_0\le t\le t_1$, be a Ricci flow. Suppose that
$c\colon S^1\to ( M\times S^1_\lambda, g(t)\times ds^2)$ is a ramp.
Then there is a curve-shrinking flow $c(x,t)$ defined for all $t\in
[t_0,t_1]$ with $c$ as the initial condition at time $t=t_0$. The
curves $c(\cdot,t)$ are all ramps.
\end{lem}

\begin{proof}
The real issue here is to show that the curve-shrinking flow exists
for all $t\in [t_0,t_1]$. Given this, the second part of the
statement follows from the previous corollary. If the curve
shrinking flow does not exist on all of $[t_0,t_1]$ then by
Claim~\ref{curvshrexist} there is a $t'_1\le t_1$ such that the
curve-shrinking flow exists on $[t_0,t_1')$ but $k$ is unbounded on
$S^1\times [t_0,t'_1)$. Thus, to complete the proof we need to see
that for any $t'_1$ for which the curve-shrinking flow is defined on
$[t_0,t'_1)$ we have a uniform bound on $k$ on this region.

Using Claim~\ref{3rdcurvshr} and Claim~\ref{rampcomp} we compute
\begin{eqnarray*} \frac{\partial}{\partial t}\left(\frac{k}{u}\right)
& = &\frac{1}{u}\frac{\partial
k}{\partial t}-\frac{k}{u^2}\frac{\partial u}{\partial t} \\
& \le & \frac{k''+k^3+C_1k}{u}-\frac{k}{u^2}\left(u''+(k^2+{\rm
Ric}(S,S))u\right)
\\
& = & \frac{k''}{u}-\frac{ku''}{u^2}+\frac{C_1k}{u}-\frac{k}{u}{\rm
Ric}(S,S).
\end{eqnarray*}
On the other hand,
$$\left(\frac{k}{u}\right)''=\frac{k''u-u''k}{u^2}-2\left(\frac{u'}{u}\right)
\left(\frac{k'u-u'k}{u^2}\right).$$ Plugging this in, and using the
curvature bound on the ambient manifolds we get
$$\frac{\partial}{\partial
t}\left(\frac{k}{u}\right)\le
\left(\frac{k}{u}\right)''+\left(\frac{2u'}{u}\right)\left(\frac{k}{u}\right)'+C'\frac{k}{u},$$
for a constant $C'$ depending only on a bound for the norm of the
sectional curvature of the ambient Ricci flow. A standard maximum
principle argument shows that the maximum of $k/u$ at time $t$ grows
at most exponentially rapidly in $t$. Since $u$ stays bounded away
from zero, this implies that for ramp solutions on a finite time
interval, the value of $k$ is bounded.
\end{proof}

Next let us turn to the growth rate of the area of a minimal annulus
connecting two ramp solutions.

\begin{lem}\label{mugrowth} Suppose that the dimension $n$ of $ M$
is at least three. Let $c_1(x,t)$ and $c_2(x,t)$ be ramp solutions
in $( M, g)\times (S^1_\lambda,ds^2)$ with the image under the
projection to $S^1_\lambda$ of each $c_i$ being of degree one. Let
$\mu(t)$ be the infimum of the areas of annuli in $( M\times
S^1_\lambda, g(t)\times ds^2)$ with boundary $c_1(x,t)\cup
c_2(x,t)$. Then $\mu(t)$ is a continuous function of $t$ and
$$\frac{d}{dt}\mu(t)\le (2n-1){\rm max}_{x\in M}|{\rm Rm}(x,t)|\mu(t),$$
in the sense of forward difference quotients.
\end{lem}

\begin{proof}
Fix a time $t'$. First assume that the loops $c_1(\cdot,t')$ and
$c_2(\cdot,t')$ are disjoint. Under Ricci flow the metrics on the
manifold immediately become real analytic (see \cite{Bando}) and
furthermore, under the curve-shrinking flow the curves $c_1$ and
$c_2$ immediately become analytic (see \cite{GageHamilton}).
[Neither of these results is essential for this argument because we
could approximate both the metric and the curves by real analytic
objects.] Establishing the results for these and taking limits would
give the result in general. Since $c_1(\cdot,t')$ and
$c_2(\cdot,t')$ are homotopic and are homotopically non-trivial
there is an annulus connecting them and there is a positive lower
bound to the length of any simple closed curve in any such annulus
homotopic to a boundary component. Hence, there is a minimal annulus
spanning $c_1(\cdot,t')\coprod c_2(\cdot, t')$ According to results
of Hildebrandt (\cite{Hildebrandt}) and Morrey (\cite{Morrey}) any
minimal annulus $A$ with boundary the union of these two curves is
real analytic up to and including the boundary and is immersed
except for finitely many branch points. By shifting the boundary
curves slightly within the annulus, we can assume that there are no
boundary branch points. Again, if we can prove the result for these
perturbed curves taking limits will give the result for the original
ones.  Given the deformation vector $H$ on the boundary of the
annulus, extend it to a deformation vector $\hat H$ on the entire
annulus. The first order variation of the area at time $t'$ of the
resulting deformed family of annuli is given by
$$\frac{d {\rm Area}\,A}{d t}(t')=\int_{A}(-{\rm Tr}({\rm Ric}^T(g(t'))))da
+\int_{\partial A}-k_{\rm geod}ds,$$ where ${\rm Ric}^T$ is the
Ricci curvature in the tangent directions to the annulus. (The first
term is the change in the area of the fixed annulus as the metric
deforms. The second term is the change in the area of the family of
annuli in the fixed metric. There is no contribution from moving the
annulus in the normal direction since the original annulus is
minimal.)
 If $A$ is embedded, then by the Gauss-Bonnet theorem, we have
$$\int_{\partial A}-k_{\rm geod}ds =\int_{A}Kda$$
where $K$ is the Gaussian curvature of $A$. More generally, if $A$
has interior branch points of orders $n_1,\ldots,n_k$ then there is
a correction term and the formula is
$$\int_{\partial A}-k_{\rm
geod}ds=\int_{A}Kda-\sum_{i=1}^k2\pi(n_i-1).$$ Thus, we see
$$\frac{d {\rm Area}\, A}{d t}(t')\le\int_{A}(-{\rm Tr}({\rm Ric}^T(g(t')))+K)da.$$
On the other hand, since $A$ is a minimal surface,  $K$ is at most
the sectional curvature of $( M, g(t'))\times (S^1_\lambda,ds^2)$
along the two-plane tangent to the annulus. Of course, the trace of
the Ricci curvature along $A$ is at most $2(n-1)|{\rm max}_{x\in
M}{\rm Rm}(x,t')|$. Hence,
$$\frac{d{\rm Area}\, A}{dt}(t')\le (2n-1)|{\rm max}_{x\in M}{\rm Rm}(x,t')|\mu(t').$$

This computation was done assuming that $c_2(\cdot,t')$ is disjoint
from $c_1(\cdot,t')$. In general, since the dimension of $ M$ is at
least three, given $c_1(\cdot,t')$ and $c_2(\cdot,t')$ we can find
$c_3(\cdot,t)$ arbitrarily close to $c_2(\cdot,t')$ in the
$C^2$-sense and disjoint from both $c_1(\cdot,t')$ and
$c_2(\cdot,t')$. Let $A_3$ be a minimal annulus connecting
$c_1(\cdot,t')$ to $c_3(\cdot,t')$ and $A_2$ be an minimal annulus
connecting $c_3(\cdot,t')$ to $c_2(\cdot,t')$. We apply the above
argument to these annuli to estimate the growth rate of minimal
annuli connecting the corresponding curve-shrinking flows. Of course
the sum of these areas (as a function of $t$) is an upper bound for
the area of a minimal annulus connecting the curve-shrinking flows
starting from $c_1(\cdot,t')$ and $c_2(\cdot,t')$. As we choose
$c_3(\cdot,t')$ closer and closer to $c_2(\cdot,t')$, the area of
$A_2$ tends to zero and the area of $A_3$ tends to the area of a
minimal annulus connecting $c_1(\cdot,t')$ and $c_2(\cdot,t')$. This
establishes the continuity of $\mu(t)$ at $t'$ and also establishes
the forward difference quotient estimate in the general case.
\end{proof}

\begin{cor}
Given  curve-shrinking flows $c_1(\cdot,t)$ and $c_2(\cdot,t)$ for
ramps of degree one in $(\bar M,\bar g(t))\times (S^1_\lambda,ds^2)$
the minimal area of an annulus connecting $c_1(\cdot,t)$ and
$c_2(\cdot,t)$ grows at most exponentially with time with an
exponent determined by an upper bound on the sectional curvature of
the ambient flow, which in particular is independent of $\lambda$.
\end{cor}

\section{Proof of Proposition~\protect{\ref{XI}}}

Now we are ready to use the curve-shrinking flow for ramps in
$M\times S_\lambda^1$ to establish Proposition~\ref{XI} for $M$. As
we indicated above, the reason for replacing the flow $(M,g(t))$
that we are studying with its product with $S^1_\lambda$ and
studying ramps in the product is that the curve-shrinking flow
exists for all time $t\in [t_0,t_1]$ for these. By this mechanism we
avoid the difficulty of finite time singularities in the curve
shrinking flow. On the other hand, we have to translate results for
the ramps back to results for the original Ricci flow $(M,g(t))$.
This requires careful analysis.

\subsection{Approximations to the original family}

The first step in the proof of Proposition~\ref{XI} is to identify
the approximation to the family $\Gamma$ that we shall use. Here is
the lemma that gives the needed approximation together with all the
properties we shall use.

Given a loop $c$ in $M$ and $\lambda>0$ we define a loop $c^\lambda$
in $M\times S^1_\lambda$. The loop $c^\lambda$ is obtained by
setting $c^\lambda(x)=(c(x),x)$ where we use  a standard
identification of the domain circle (the unit circle) for the free
loop space with $S^1_\lambda$, an identification that defines a loop
in $S^1_\lambda$ of constant speed $\lambda/2\pi$.

\begin{lem}\label{goodapprox}
Given a continuous map $\Gamma\colon S^2\to \Lambda M$ representing
an element of $\pi_3( M,*)$ and $0<\zeta<1$, there is a continuous
map $\widetilde \Gamma\colon S^2\to \Lambda M$ with the following
properties:
\begin{enumerate}
\item[(1)] $[\widetilde \Gamma]=[\Gamma]$ in $\pi_3(M,*)$.
\item[(2)] For each $c\in S^2$ the loop $\widetilde \Gamma(c)$ is a
$C^2$-loop.
\item[(3)] For each $c\in S^2$ the length of $\widetilde \Gamma(c)$ is
within $\zeta$ of the length of $\Gamma(c)$.
\item[(4)] For each $c\in S^2$, we have $|A(\widetilde
\Gamma(c))-A(\Gamma(c))|<\zeta$.
\item[(5)] There is a constant $C_0<\infty$ depending only on $\Gamma$,
on the bounds for the norm of the Riemann curvature operator of the
ambient Ricci flow, and on $\zeta$ such that for each $c\in S^2$ and
each $\lambda\in (0,1)$ the total length and the total curvature of
the ramp $\widetilde \Gamma(c)^\lambda$ are both bounded by $C_3$.
\end{enumerate}
\end{lem}

Before proving this lemma we need some preliminary definitions and
constructions.

\begin{defn}
Let $c\colon S^1\to M$ be a $C^1$-map.  Fix a positive integer $n$.
By a {\em regular $n$-polygonal approximation} to $c$ we mean the
following. Let $\xi_n={\rm exp}(2\pi i/n)$, and consider the points
$p_k=c(\xi^k_n)$ for $k=1,\ldots,n+1$. For each $1\le k\le n$, let
$A_k$ be a minimal geodesic in $M$ from $p_k$ to $p_{k+1}$. We
parameterize $A_k$ by the interval $[\xi_n^k,\xi_n^{k+1}]$ in the
circle at constant speed. This gives a  piecewise geodesic map
$c_n\colon S^1\to M$.
\end{defn}

The following is immediate from the definition.

\begin{claim}\label{onec}
Given $\zeta>0$ and a $C^1$-map $c\colon S^1\to M$ then for all $n$
sufficiently large the following hold for the $n$-polygonal
approximation $c_n$ of $c$.
\begin{enumerate}
\item[(a)] the length of $c_n$ is within $\zeta$ of the length of $c$.
\item[(b)] there is a map of the annulus $S^1\times I$ to $M$ connecting
$c_n$ to $c$ with the property that the image is piecewise smooth
and of area less than $\zeta$.
\end{enumerate}
\end{claim}

\begin{proof}
The length of $c$ is the limit of the lengths of the $n$-polygonal
approximations as $n$ goes to infinity. The first item is immediate
from this. As to the second, for $n$ sufficiently large, the
distance between the maps $c$ and $c_n$ will be arbitrarily small in
the $C^0$-topology, and in particular will be much smaller than the
injectivity radius of $M$. Thus, for each $k$ we can connect $A_k$
to the corresponding part of $c$ by a family of short geodesics.
Together, these form an annulus, and it is clear that for $n$
sufficiently large the area of this annulus is arbitrarily small.
\end{proof}

As the next result shows, for $\zeta>0$, the integer $n(c)$
associated by the previous claim to a $C^1$-map $c$ can be made
uniform as $c$ varies over a compact subset of $\Lambda M$.

\begin{claim}\label{uniformc}
Let $X\subset \Lambda M$ be a compact subset and let $\zeta>0$ be
fixed. Then there is $N$ depending only on $X$ and $\zeta$ such the
conclusion of the previous claim holds for every $c\in X$ and every
$n\ge N$.
\end{claim}

\begin{proof}
Suppose the result is false. Then for each $N$ there is $c_N\in X$
and $n\ge N$ so that the lemma does not hold for $c_N$ and $n$.
Passing to a subsequence, we can suppose that the $c_N$ converge to
$c_\infty\in X$. Applying Claim~\ref{onec} we see that there is $N$
such that the conclusion of Claim~\ref{onec} holds with $\zeta$
replaced by $\zeta/2$ for $c_\infty$ and all $n\ge N$. Clearly, then
by continuity for all $n\ge N$ the conclusion of Claim~\ref{onec}
holds for the $n$-polygonal approximation for every $c_l$ for all
$l$ sufficiently large. This is a contradiction.
\end{proof}

\begin{cor}\label{Gammaapprox}
Let $\Gamma\colon S^2\to \Lambda M$ be a continuous map  with the
property that $\Gamma(c)$ is homotopically trivial for all $c\in
S^2$. Fix $\zeta>0$. For any $n$ sufficiently large denote by
$\Gamma_n$ the family of loops defined by setting $\Gamma_n(c)$
equal to the  $n$-polygonal approximation to $\Gamma(c)$. There is
$N$ such that for all $n\ge N$ we have
\begin{enumerate}
\item[(1)] $\Gamma_n$ is a continuous family of $n$-polygonal loops in
$M$.
\item[(2)] For each $c\in S^2$, the loop
$\Gamma_n(c)$ is a homotopically trivial loop in $M$ and its length
is within $\zeta$ of the length of $\Gamma(c)$.
\item[(3)] For each $c\in S^2$, we have $|A(\Gamma_n(c))-A(\Gamma(c))|<\zeta$.
\end{enumerate}
\end{cor}

\begin{proof}
Given $\Gamma$ there is a uniform bound over all $c\in S^2$ on the
maximal speed of $\Gamma(c)$. Hence, for all $n$ sufficiently large,
the lengths of the sides in the $n$-polygonal approximation to
$\Gamma(c)$ will be uniformly small. Once this length is less than
the injectivity radius of $M$, the minimal geodesics between the
endpoints are unique and vary continuously with the endpoints. This
implies that for $n$ sufficiently large the family $\Gamma_n$ is
uniquely determined and itself forms a continuous family of loops in
$M$. This proves the first item. We have already seen that, for $n$
sufficiently large, for all $c\in S^2$ there is an annulus
connecting $\Gamma(c)$ and $\Gamma_n(c)$. Hence, these loops are
homotopic in $M$. The first statement in the second item follows
immediately. The last statement in the second item and third item
follow immediately from Claim~\ref{uniformc}.
\end{proof}

The next step is  to turn these $n$-polygonal approximations into
$C^2$-curves. We fix, once and for all, a $C^\infty$ function
$\psi_n$ from the unit circle to $[0,\infty]$ with the following
properties:
\begin{enumerate}
\item[(1)] $\psi_n$ is non-negative and vanishes to infinite order at the point $1$ on the unit circle.
\item[(2)] $\psi_n$ is periodic with period $2\pi i/n$.
 \item[(3)] $\psi_n$ is positive on the interior of the interval $[1,\xi_n]$ on the unit circle,
and the restriction of $\psi_n$ to this interval is symmetric about
${\rm exp}(\pi i/n)$, and is increasing from $1$ to ${\rm exp}(\pi
i/n)$.
\item[(4)] $\int_1^{\xi_n}\psi_n(s)ds =2\pi/n$.
\end{enumerate}

 Now we define a map $\tilde\psi_n\colon S^1\to S^1$ by
$$\tilde \psi_n(x)= \int_{1}^x\psi_n(y)dy.$$ It is easy to see that
the conditions on $\psi_n$ imply that this defines a $C^\infty$-map
from $S^1$ to $S^1$ which is a homeomorphism and is a diffeomorphism
on the complement of the $n^{th}$ roots of unity.

Now given an $n$-polygonal loop $c_n$ we define the smoothing
$\tilde c_n$ of $c_n$ by $\tilde c_n=c_n\circ \tilde \psi_n$. This
smoothing $\tilde c_n$ is a $C^\infty$-loop in $M$ with the same
length as the original polygonal loop $c_n$. Notice that the
curvature of $\tilde c_n$ is not itself a continuous function: just
like the polygonal map it replaces, it has a $\delta$-function at
the `corners' of $c_n$.

\begin{proof} (of Lemma~\ref{goodapprox})
Given a continuous map $\Gamma\colon S^2\to \Lambda M$ and $\zeta>0$
we fix $n$ sufficiently large so that Corollary~\ref{Gammaapprox}
holds for these choices of $\Gamma$ and $\zeta$. Let $\widetilde
\Gamma=\widetilde \Gamma_n$ be the family of smoothings of the
family $\Gamma_n$ of $n$-polygonal loops. Since this smoothing
operation changes neither the length nor the area of a minimal
spanning disk, it follows  immediately from the construction and
Corollary~\ref{Gammaapprox} that $\widetilde \Gamma$ satisfies the
conclusions of Lemma~\ref{goodapprox} except possibly the last one.

To establish the last conclusion we must examine the lengths and
total curvatures of the ramps $\widetilde \Gamma(c)^\lambda$
associated to this family of $C^2$-loops. Fix $\lambda$ with
$0<\lambda<1$, and consider the product Ricci flow $(M,g(t))\times
(S^1_\lambda,ds^2)$ where the metric on $S^1_\lambda$ has length
$\lambda$.

\begin{claim}
For any $0<\lambda<1$, the length of the ramp $\widetilde
\Gamma(c)^\lambda$ is at most $\lambda$ plus the length of
$\Gamma(c)$. The total curvature of $\widetilde\Gamma(c)^\lambda$ is
at most $n\pi$.
\end{claim}

\begin{proof}
The arc length element for $\widetilde \Gamma(c)^\lambda$ is
$\sqrt{a(x)^2+(\lambda/2\pi)^2}dx\le (a(x)+\lambda/2\pi)dx$ where
$a(x)dx$ is the arc length element for $\widetilde\Gamma(c)$.
Integrating gives the length estimate.

The total curvature of $\widetilde \Gamma(c)^\lambda$ is the sum
over the intervals $[\xi_n^k,\xi_n^{k+1}]$ of the total curvature on
these intervals. On any one of these intervals we have a curve in a
totally geodesic, flat surface: the curve lies in the product of a
geodesic arc in $M$ times $S^1_\lambda$. Let $u$ and $v$ be unit
tangent vectors to this surface, $u$ along the geodesic (in the
direction of increasing $x$) and $v$ along the $S^1_\lambda$ factor.
These are parallel vector fields on the flat surface.
 The tangent
vector $X(x)$ to the restriction of $\widetilde \Gamma(c)^\lambda$
to this interval is $L\psi_n(x) u+(\lambda/2\pi)v$, where $L$ is the
length of the geodesic segment we are considering. Consider the
first-half subinterval $[\xi_n^k,\xi_n^k\cdot{\rm exp}(\pi i/n)]$.
The tangent vector $X(x)$ is $(\lambda/2\pi) v$ at the initial point
of this subinterval and is $L\psi_n(\xi_n^k\cdot{\rm exp}(\pi
i/2))u+(\lambda/2\pi) v$ at the final point. Throughout this
interval the vector is of the form $a(x)u+(\lambda/2\pi) v$ where
$a(x)$ is an increasing function of $x$. Hence, the tangent vector
is always turning in the same direction and always lies in the first
quadrant (using $u$ and $v$ as the coordinates). Consequently, the
total turning (the integral of $k$ against arc-length) over this
interval is the absolute value of the difference of the angles at
the endpoints. This difference is less than $\pi/2$ and tends to
$\pi/2$ as $\lambda$ tends to zero, unless $L=0$ in which case there
is zero turning for any $\lambda>0$. By symmetry, the total turning
on the second-half subinterval $[\xi_n^k\cdot{\rm exp}(\pi
i/n),\xi_n^{k+1}]$ is also bounded above by $\pi/2$. Thus, for any
$\lambda>0$, the total turning on one of the segments is bounded
above by $\pi$. Since there are $n$ segments this gives the upper
bound of $n\pi$ on the total turning of $\widetilde
\Gamma(c)^\lambda$ as required.
\end{proof}

This claim completes the proof of the last property required of
$\widetilde\Gamma=\widetilde \Gamma_n$ and hence completes the proof
of Lemma~\ref{goodapprox}.
\end{proof}

Having fixed $\Gamma$ and $\zeta>0$, we fix $n$ and set
$\widetilde\Gamma=\widetilde \Gamma_n$. We choose $n$ sufficiently
large so that $\widetilde \Gamma$ satisfies Lemma~\ref{goodapprox}.
 Fix $\lambda\in (0,1)$ and
define $\widetilde \Gamma^\lambda\colon S^2\to (M\times S^1_\lambda)$, by
setting $\widetilde \Gamma^\lambda(c)=\widetilde \Gamma(c)^\lambda$.

Fix  $c\in S^2$, and let $\widetilde\Gamma_c^\lambda(t),\ t_0\le
t\le t_1$, be the curve-shrinking flow given in
Lemma~\ref{csforramps} with initial data the ramp
$\widetilde\Gamma^\lambda(c)$ . As $c$ varies over $S^2$ these fit
together to produce a one-parameter family
$\widetilde\Gamma^\lambda(t)$ of maps $S^2\to\Lambda(M\times
S^1_\lambda)$. Let $p_1$ denote the projection of $M\times
S^1_\lambda$ to $M$. Notice that for any $\lambda$ we have
$\widetilde\Gamma^\lambda_c(t_0)=\widetilde\Gamma^\lambda(c)$, so
that $p_1\widetilde \Gamma_c^\lambda(t_0)=\widetilde \Gamma(c)$. We
shall show that for $\lambda>0$ sufficiently small, the family
$p_1\widetilde\Gamma^\lambda(t)$ satisfies the conclusion of
Proposition~\ref{XI} for the fixed $\Gamma$ and $\zeta>0$. We do
this in steps. First, we show that fixing one $c\in S^2$, for
$\lambda$ sufficiently small (depending on $c$) an analogue  of
Proposition~\ref{XI} holds for the one-parameter family of loops
$p_1\widetilde\Gamma_c^\lambda(t)$.  By this we mean that either
$p_1\widetilde\Gamma^\lambda_c(t_1)$ has length less than $\zeta$ or
$A(p_1\widetilde\Gamma^\lambda_c(t_1))$ is at most the value
$v(t_1)+\zeta$, where $v$ is the solution to the
Equation~(\ref{wdiffeqn}) with initial condition
$v(t_0)=A(\widetilde\Gamma(c))$. (Actually, we establish a slightly
stronger result, see Lemma~\ref{xistronger}.) The next step in the
argument is to take a finite subset ${\mathcal S}\subset S^2$ so
that for every $c\in S^2$ there is $\hat c\in {\mathcal S}$ such
that $\widetilde\Gamma(c)$ and $\widetilde\Gamma(\hat c)$ are
sufficiently close. Then, using the result of a single $c$, we fix
$\lambda>0$ sufficiently small so that the analogue of
Proposition~\ref{XI} for individual curves (or rather the slightly
stronger version of it) holds for every $\hat c\in {\mathcal S}$.
Then we complete the proof of Proposition~\ref{XI} using the fact
that for every $c$ the curve $\widetilde\Gamma(c)$ is sufficiently
close to a curve $\widetilde\Gamma(\hat c)$ associated to an element
$\hat c\in {\mathcal S}$.

\subsection{The case of a single $c\in S^2$}

According to Lemma~\ref{goodapprox}, for all $\lambda\in (0,1)$  the
lengths and total curvatures of the $\widetilde\Gamma^\lambda(c)$
are uniformly bounded for all $c\in S^2$. Hence, by
Corollary~\ref{lncurvgrowth} the same is true for
$\widetilde\Gamma_c^\lambda(t)$ for all $c\in S^2$ and all $t\in
[t_0,t_1]$.

\begin{claim}\label{areachange}
There is a constant $C_4$ depending on $t_1-t_0$, on the curvature
bound of the sectional curvature of the Ricci flow $(M,g(t)),\
t_0\le t\le t_1$, on the original family $\Gamma$ and on $\zeta$
such that for any $c\in S^2$ and any $t_0\le t'<t''\le t_1$ we have
$$A(p_1\widetilde \Gamma^\lambda_c(t''))-A(p_1\widetilde
\Gamma^\lambda_c(t'))\le C_4(t''-t').$$
\end{claim}

\begin{proof}
All the constants in this argument are allowed to depend on
$t_1-t_0$, on the curvature bound of the sectional curvature of the
Ricci flow $(M,g(t)),\ t_0\le t\le t_1$, on the original family
$\Gamma$ and on $\zeta$ but are independent of $\lambda$, $c\in
S^2$, and $t'<t''$ with $t_0\le t'$ and $t''\le t_1$. First, let us
consider the surface $S^\lambda_c[t',t'']$ in  $M\times S^1_\lambda$
swept out by $c(x,t),\ t'\le t\le t''$. We denote by ${\rm
Area}\,(S^\lambda_c[t',t''])$ the area of this surface with respect
to the metric $g(t'')\times ds^2$. We compute the derivative of this
area for fixed $t'$ as $t''$ varies. There are two contributions to
this derivative: (i) the contribution due to the variation of the
metric $g(t'')$ with $t''$ and (ii) the contribution due to
enlarging the surface. The first is
$$\int_{S^\lambda_c[t',t'']}-{\rm
Tr}\, {\rm Ric}^Tda$$
 where ${\rm Ric}^T$ is the restriction of the
Ricci tensor of the ambient metric $g(t'')$ to the tangent planes to
the surface and $da$ is the area form of the surface in the metric
$g(t'')\times ds^2$. The second contribution is
$\int_{c(x,t'')}|H|ds$. According to Lemma~\ref{1835} there is a
constant $C'$ (depending only on the curvature bound for the
manifold flow, the initial family $\Gamma(t)$ and $\zeta$ and
$t_1-t_0$) such that the second term is bounded above by $C'$. The
first term is bounded above by $C''{\rm Area}\,S^\lambda_c[t',t'']$
where $C''$ depends only on the bound on the sectional curvatures of
the ambient Ricci flow. Integrating we see that there is a constant
$C_1'$ such that the derivative of the area function is at most
$C_1'$. Since its value at $t'$ is zero, we see that ${\rm
Area}\,S^\lambda_c[t',t'']\le C'_1(t''-t')$. It follows that the
area of $p_1S^\lambda_c[t',t'']$ with respect to the metric $g(t'')$
is at most $C_1'(t''-t')$.

Now we compute an upper bound for the forward difference quotient of
$A(p_1\widetilde \Gamma^\lambda_c(t))$ at $t=t'$. For any $t''>t'$
we have a spanning disk for $p_1\widetilde\Gamma^\lambda_c(t'')$
defined by taking the union of a minimal spanning disk for
$p_1\widetilde \Gamma^\lambda_c(t')$ and the annulus
$p_1S^\lambda_c[t',t'']$. As before, the derivative of the area of
this family of disks has two contributions, one coming from the
change in the metric over the minimal spanning disk at time $t'$ and
the other which we computed above to be at most $C_1'$. Thus, the
derivative is bounded above by
$C'_2A(p_1\widetilde\Gamma^\lambda_c(t'))+C_1'$. This implies that
the forward difference quotient of $A(p_1\widetilde
\Gamma^\lambda_c(t'))$ is bounded above by the same quantity. It
follows immediately that the areas of all the minimal spanning
surfaces are bounded by a constant depending only on the areas of
the minimal spanning surfaces at time $t_0$, the sectional curvature
of the ambient Ricci flow and $t_1-t_0$. Hence, there is a constant
$C_4$ such that the forward difference quotient of
$A(p_1\widetilde\Gamma^\lambda_c(t))$ is bounded above by $C_4$.
This proves the claim.
\end{proof}

Next, by the uniform bounds on total length of all the curves
$\widetilde\Gamma_c^\lambda(t)$, it follows from
Equation~(\ref{Lformula}) that there is a constant $C_5$ (we take
$C_5>1$) depending only on the curvature bound of the ambient
manifolds and the family $\Gamma$ such that for any $c\in S^2$ we
have
\begin{equation}\label{C5eqn}\int_{t_0}^{t_1}\int_{\widetilde\Gamma_c^\lambda(t)}
k^2dsdt\le C_5. \end{equation} Thus, for any constant $1<B<\infty$
there is a subset $I_B(c,\lambda)\subset [t_0,t_1]$ of measure at
least $(t_1-t_0)-C_5B^{-1}$ such that
$$\int_{\widetilde\Gamma_c^\lambda(t)}k^2ds\le B$$
for every $t\in I_B(c,\lambda)$. (Later, we shall fix $B$
sufficiently large depending on $\Gamma$ and $\zeta$.)

Now we need a result for curve-shrinking that in some ways is
reminiscent of Shi's theorem for Ricci flows.

\begin{lem}\label{CSSHI}
Let $(M,g(t)),\ t_0\le t\le t_1$, be a Ricci flow. Then there exist
constants $\delta>0$ and $\widetilde C_i<\infty$ for
$i=0,1,2,\ldots$, depending only on $t_1-t_0$  and a bound for the
norm of the curvature of the Ricci flow, such that the following
holds. Let $c(x,t)$ be a curve-shrinking flow that is an immersion
for each $t$. Suppose that at a time $t'$ for some $0<r<1$ such that
$t'+\delta r^2<t_1$, the length of $c(\cdot,t')$ is at least $r$ and
the total curvature of $c(\cdot,t')$ on any subarc of length $r$ is
at most $\delta$. Then for every $t\in [t',t'+\delta r^2)$ the
curvature $k$ and the higher derivatives satisfy
$$k^2\le\widetilde C_0(t-t')^{-1}$$
$$|\nabla_SH|^2\le \widetilde C_1(t-t')^{-2},$$
$$|\nabla_S^iH|^2\le \widetilde C_i(t-t')^{-(i+1)}.$$
\end{lem}

The first statement follows from arguments very similar to those in
Section 4 of \cite{AG}. Once $k^2$ is bounded by $\widetilde
C_0/(t-t')$ the higher derivative statements are standard, see
\cite{Altschuler}. For completeness we have included the proof of
the first inequality in the last section of this chapter.

We now fix $\delta>0$ (and also $\delta<1$) as described in the last
lemma for the Ricci flow  $(M,g(t)),\ t_0\le t\le t_1$. By
Cauchy-Schwarz it follows that for every $t\in I_B(c,\lambda)$,
  and for any arc $J$ in $\widetilde\Gamma_c^\lambda(\cdot,t)$ of
length at most $\delta^2B^{-1}$ we have
$$\int_{J\times\{t\}}k\le \delta.$$
Applying the previous lemma,  for each $a\in I_B(c,\lambda)$ with
$a\le t_1-B^{-1}-\delta^5B^{-2}$ we set
$J(a)=[a+\delta^5B^{-2}/2,a+\delta^5B^{-2}]\subset
[t_0,t_1-B^{-1}]$. Then for all $t\in \cup_{a\in
I_B(c,\lambda)}J(a)$ for which the length of $\widetilde
\Gamma_c^\lambda(\cdot,t)$ is at least $\delta^2B^{-1}$ we have that
$k$ and all the norms of spatial derivatives of $H$ are pointwise
uniformly bounded. Since $I_B(c,\lambda)$ covers all of $[t_0,t_1]$
except a subset of measure at most $C_5B^{-1}$, it follows that the
union $\widehat J_B(c,\lambda)$ of intervals $J(a)$ for $a\in
I_B(c,\lambda)\cap [t_0,t_1-B^{-1}-\delta^5B^{-2}]$ cover all of
$[t_0,t_1]$ except a subset of measure at most
$C_5B^{-1}+B^{-1}+\delta^5B^{-1}<3C_5B^{-1}$. Now it is
straightforward to pass to a finite subset of these intervals
$J(a_i)$ that cover all of $[t_0,t_1]$ except a subset of measure at
most $3C_5B^{-1}$.
 Once we
have a finite number of $J(a_i)$, we order them along the interval
$[t_0,t_1]$ so that  their initial points form an increasing
sequence. (Recall that they all have the same length.) Then if we
have $J_i\cap J_{i+2}\not=\emptyset$, then $J_{i+1}$ is contained in
the union of $J_i$ and $J_{i+2}$ and hence can be removed from the
collection without changing the union. In this way we reduce to a
finite collection of intervals $J_i$, with the same union, where
every point of $[t_0,t_1]$ is contained in at most $2$ of the
intervals in the collection. Once we have arranged this we have a
uniform bound, independent of $\lambda$ and $c\in S^2$, on the
number of these intervals. We let $J_B(c,\lambda)$ be the union of
these intervals. According to the construction and Lemma~\ref{CSSHI}
these sets $J_B(c,\lambda)$ satisfy the following:
\begin{enumerate}
\item[(1)] $J_B(c,\lambda)\subset [t_0,t_1-B^{-1}]$ is a union of a bounded number of intervals
(the bound being independent of $c\in S^2$ and of $\lambda$) of
length $\delta^5B^{-2}/2$.
\item[(2)] The measure of $J_B(c,\lambda)$ is at least
$t_1-t_0-3C_5B^{-1}$.
\item[(3)] For every $t\in J_B(c,\lambda)$ either the length of $\widetilde\Gamma_c^\lambda(t)$ is
less than $\delta^2 B^{-1}$ or there are uniform bounds, depending
only on the curvature bounds of the ambient Ricci flow and the
initial family $\Gamma$, on the curvature and its higher spatial
derivatives of $\widetilde\Gamma_c^\lambda(t)$.
\end{enumerate}

Now we fix $c\in S^2$ and $1<B<\infty$ and we fix a sequence of
$\lambda_n$ tending to zero. Since the number of intervals in
$J_B(c,\lambda)$ is bounded independent of $\lambda$, by passing to
a subsequence of $\lambda_n$ we can suppose that the number of
intervals in $J_B(c,\lambda_n)$ is independent of $n$, say this
number is $N$, and that their initial points (and hence the entire
intervals since all their lengths are the same) converge as $n$ goes
to infinity. Let $\hat J_1,\ldots, \hat J_N$ be the limit intervals,
and for each $i, 1\le i\le N$, let $J_i\subset \hat J_i$ be a
slightly smaller interval contained in the interior of $\hat J_i$.
We choose the $J_i$ so that they all have the same length. Let
$J_B(c)\subset [t_0,t_1-B^{-1}]$ be the union of the $J_i$. Then an
appropriate choice of the length of the $J_i$ allows us to arrange
the following:
\begin{enumerate}
\item[(1)] $J_B(c)\subset J_B(c,\lambda_n)$ for all $n$ sufficiently
large.
\item[(2)] $J_B(c)$ covers all of $[t_0,t_1]$ except a subset of length
$4C_5B^{-1}$.
\end{enumerate}

Now fix one of the intervals $J_i$ making up $J_B(c)$. After passing
to a subsequence (of the $\lambda_n$), one of the following holds:
\begin{enumerate}
\item[(3)] there are uniform bounds for the curvature and all its
derivatives for the curves $\widetilde\Gamma_c^{\lambda_n}(t)$, for
all $ t\in J_i$ and all $n$, or
\item [(4)] for each $n$ there is $t_n\in J_i$ such that the length of
$\widetilde\Gamma_c^{\lambda_n}(t_n)$ is less than $\delta^2B^{-1}$.
\end{enumerate}

By passing to a further subsequence, we arrange that the same one of
the Alternatives (3) and (4) holds for every  one of the intervals
$J_i$ making up $J_B(c)$.

The next claim is the statement that a slightly stronger version of
Proposition~\ref{XI} holds for $p_1\widetilde\Gamma_c^\lambda(t)$.

\begin{lem}\label{xistronger}
Given $\zeta>0$, there is   $1<B<\infty$, with $B>(t_1-t_0)^{-1}$,
depending only on $\Gamma$ and the curvature bounds on the ambient
Ricci flow $(M,g(t)),\ t_0\le t\le t_1$, such that the following
holds. Let $t_2=t_1-B^{-1}$. Fix $ c\in S^2$. Let $v_c$ be the
solution to Equation~(\ref{wdiffeqn}) with initial condition
$v_{c}(t_0)=A(p_1(\widetilde\Gamma(c)))$, so that in our previous
notation $v_c=w_{A(p_1(\widetilde\Gamma(c)))}$. Then for all
$\lambda>0$ sufficiently small, either
$A(p_1\widetilde\Gamma_c^\lambda(t_1))<v_c(t_1)+\zeta/2$ or the
length of $\widetilde\Gamma_c^\lambda(t)$ is less than $\zeta/2$ for
all $t\in [t_2,t_1]$.
\end{lem}

\begin{proof}
In order to establish this lemma  we need a couple of claims about
functions on $[t_0,t_1]$ that are approximately dominated by
solutions to Equation~(\ref{wdiffeqn}). In the first claim the
function in question is dominated on a finite collection of
subintervals by solutions to these equations and the subintervals
fill up most of the interval. In the second, we also allow the
function to only be approximately dominated by the solutions to
Equation~(\ref{wdiffeqn}) on these sub-intervals. In both claims the
result is that on the entire interval the function is almost
dominated by the solution to the equation with the same initial
value.

\begin{claim}\label{C1C4delta}
Fix $C_4$ as in Claim~\ref{areachange} and fix a constant
$\widetilde A>0$. Given $\zeta>0$ there is $\delta'>0$ depending on
$C_4$, $t_1-t_0$, and $\widetilde A$ as well as the curvature bound
of the ambient Ricci flow such that the following holds. Suppose
that $f\colon [t_0,t_1]\to \Ar$ is a function and suppose that
$J\subset [t_0,t_1]$ is a finite union of intervals. Suppose that on
each interval $[a,b]$ of $J$ the function $f$ satisfies
$$f(b)\le w_{f(a),a}(b).$$ Suppose further that for any $t'<t''$ we
have
$$f(t'')\le f(t')+C_4(t''-t').$$
Then, provided that the total length of $[t_0,t_1]\setminus J$ is at
most $\delta'$ and $0\le f(t_0)\le \widetilde A$, we have
$$f(t_1)\le w_{f(t_0),t_0}(t_1)+\zeta/4.$$
\end{claim}

\begin{proof}
We write $J$ as a union of disjoint intervals $J_1,\ldots,J_k$ so
that  $J_i<J_{i+1}$ for every $i$. Let $a_i$, resp. $b_i$, be the
initial, resp. final, point of $J_i$. For each $i$ let $\delta_i$ be
the length of the interval between $J_i$ and $J_{i+1}$. (Also, we
set $\delta_0=a_1-t_0$, and $\delta_{k}=t_1-b_k$.) Let $C_6\ge 0$ be
such that $R_{\rm min}(t)\ge -2C_6$ for all $t\in [t_0,t_1]$. Let
$V(a)$ be the maximum value of $|w_{a,t_0}|$ on the interval
$[t_0,t_1]$ and let $V={\rm max}_{a\in [0,\widetilde A]}V(a)$. let
$C_7=C_4+2\pi+C_6V$.
 We shall prove by
induction that
$$f(a_i)-w_{f(t_0),t_0}(a_i)\le \sum_{j=0}^{i-1}\left(C_7\delta_j
\prod_{\ell=j+1}^{i-1}e^{C_6|J_\ell|}\right)$$ and
$$f(b_i)-w_{f(t_0),t_0}(b_i)\le
\sum_{j=0}^{i-1}\left(C_7\delta_j\prod_{\ell=j+1}^ie^{
C_6|J_\ell|}\right).$$

We begin the induction by establishing the result at $a_1$. By
hypothesis we know that
$$f(a_1)\le f(t_0)+C_4\delta_0.$$
 On the other
hand, from the defining differential equation for $w_{f(t_0),t_0}$
and the definitions of $C_6$ and $V$ we have
$$w_{f(t_0),t_0}(a_1)\ge f(t_0) -(C_6V+2\pi)\delta_0.$$
Thus,
$$f(a_1)-w_{f(t_0),t_0}(a_1)\le (C_4+2\pi+C_6V)\delta_0=C_7\delta_0,$$
which is exactly the formula given in the case of $a_1$.

Now suppose that we know the result for $a_i$ and let us establish
it for $b_i$. Let $\alpha_i=f(a_i)-w_{f(t_0),t_0}(a_i)$, and let
$\beta_i=f(b_i)-w_{f(t_0),t_0}(b_i)$. Then by Claim~\ref{vclaim} we
have
$$\beta_i\le e^{C_6|J_i|}\alpha_i.$$
Given the inductive inequality for $\alpha_i$, we immediately get
the one for $\beta_i$.

Now suppose that we have the inductive inequality for $\beta_i$.
Then
$$f(a_{i+1})\le f(b_i)+C_4\delta_i.$$
 On the other hand, by
 the definition of $C_6$ and $V$ we have
$$w_{f(t_0),t_0}(a_{i+1})-w_{f(t_0),t_0}(b_i)\ge
-(C_6V+2\pi)\delta_{i}.$$ This yields
$$f(a_{i+1})-w_{f(t_0),t_0}(a_{i+1})\le \beta_i+C_7\delta_i.$$
Hence, the inductive result for $\beta_i$ implies the result for
$\alpha_{i+1}$. This completes the induction.

Applying this to $a_{k+1}=t_1$ gives
$$f(t_1)-w_{f(t_0),t_0}(t_1)\le
\sum_{j=0}^k\left(C_7\delta_j\prod_{\ell=j+1}^{k}e^{
C_6|J_\ell|}\right)\le C_7\sum_{j=0}^k\delta_je^{C_6(t_1-t_0)}.$$

Of course $\sum_{j=1}^k\delta_j=t_1-t_0-\ell(J)\le\delta'$, and
$C_7$ only depends on $ C_6,C_4$ and $V$, while $V$  only depends on
$\widetilde A$ and $C_6$ only depends on the sectional curvature
bound on the ambient Ricci flow. Thus, given $C_4, \widetilde A$ and
$t_1-t_0$ and the bound on the sectional curvature of the ambient
Ricci flow, making $\delta'$ sufficiently small makes
$f(t_1)-w_{f(t_0),t_0}(t_1)$ arbitrarily small. This completes the
proof of the claim.
\end{proof}

Here is the second of our claims:

\begin{claim}\label{deltaprime}
Fix $\zeta>0$, $A$ and $C_6, C_4$ as in the last claim, and let
$\delta'>0$ be as in the last claim. Suppose that we have $J\subset
[t_0,t_1]$ which is a finite disjoint union of intervals with
$t_1-t_0-|J|\le \delta'$.
 Then there
is $\delta''>0$ ($\delta''$ is allowed to depend on $J$) such that
the following holds. Suppose that we have a function $f\colon
[t_0,t_1]\to \Ar$ such that:
\begin{enumerate}
\item[(1)] For all $t'<t''$ in $[t_0,t_1]$ we have $f(t'')-f(t')\le
C_4(t''-t')$.
\item[(2)] For any interval $[a,b]\subset J$ we have
$f(b)\le w_{f(a),a}(b)+\delta''$.
\end{enumerate}
Then $f(t_1)\le w_{f(t_0),t_0}(t_1)+\zeta/2$.
\end{claim}

\begin{proof}
We define $C_7$ as in the previous proof. We use the notation
$J=J_1\coprod\cdots\coprod J_k$ with $J_1<J_2<\cdots<J_k$ and let
$\delta_i$ be the length of the interval separating $J_{i-1}$ and
$J_i$. The arguments in the proof of the previous claim work in this
context to show that
$$f(a_i)-w_{f(t_0),t_0}(a_i)\le
\sum_{j=0}^{i-1}\left(C_7\delta_j\prod_{\ell=j+1}^{i-1}(e^{C_6|J_\ell|}+\delta'')\right).$$
Applying this to $a_{k+1}$ and taking the limit as $\delta''$ tends
to zero, the right-hand side tends to a limit smaller than
$\zeta/4$. Hence, for $\delta''$ sufficiently small the right-hand
side is less than $\zeta/2$.
\end{proof}

Now let us return to the proof of Lemma~\ref{xistronger}. Recall
that $c\in S^2$ is fixed. We shall apply the above claims to the
curve-shrinking flow $\widetilde\Gamma_c^\lambda(t)$ and thus prove
Lemma~\ref{xistronger}. Now it is time to fix $B$. First, we fix
$\widetilde A=W(\Gamma)+\zeta$,  we let $C_2$ be as in
Corollary~\ref{lncurvgrowth}, $C_4$ be as in Claim~\ref{areachange},
$C_5$ be as in Equation~(\ref{C5eqn}), and $C_6$ be as in the proof
of Claim~\ref{C1C4delta}. Then we have $\delta'$ depending on
$C_6,C_4,\widetilde A$ as in Claim~\ref{C1C4delta}. We fix $B$ so
that:
\begin{enumerate}
\item[(1)]  $B\ge 3C_5(\delta')^{-1}$,
\item[(2)] $B\ge 3e^{C_2(t_1-t_0)}\zeta^{-1}$, and
\item[(3)] $B>C_2/({\rm log}4-{\rm log}3)$.
\end{enumerate}
The first step in the proof of Lemma~\ref{xistronger} is the
following:

\begin{claim}\label{limitcs}
After passing to a subsequence of $\{\lambda_n\}$, either:
\begin{enumerate}
\item[(1)] for each $n$ sufficiently large there is $t_n\in J_B(c)$ with
the length of $\widetilde\Gamma_c^{\lambda_n}(t_n)<\delta^2B^{-1}$,
or
\item[(2)] for each component $J_i=[t_i^-,t_i^+]$ of $ J_B(c)$, after
composing $\widetilde\Gamma_c^{\lambda_n}(x,t)$ by a
reparameterization of the domain circle (fixed in $t$ but a
different reparameterization for each $n$) so that the $\widetilde
\Gamma_c^{\lambda_n}(t_i^-)$ have constant speed, there is a smooth
limiting curve-shrinking flow denoted $\widetilde\Gamma_c(t)$, for
$t\in J_i$ for the sequence $p_1\widetilde\Gamma_c^{\lambda_n}(t),\
t_i^-\le t\le t_i^+$. The limiting flow consists of immersions.
\end{enumerate}
\end{claim}

\begin{proof}
Suppose that the first case does not hold for any subsequence. Fix a
component $J_i$ of $J_B(c)$. Then, by passing to a subsequence, by
the fact that $J_B(c)\subset J_B(c,\lambda_n)$ for all $n$, the
curvatures and all the derivatives of the curvatures of
$\widetilde\Gamma_c^{\lambda_n}(t)$ are uniformly bounded
independent of $n$ for all $t\in J_B(c)$. We reparameterize the
domain circle so that the $\widetilde\Gamma_c^{\lambda_n}(t_i^-)$
have constant speed. By passing to a subsequence we can suppose that
the lengths of the $\widetilde\Gamma_c^{\lambda_n}(t_i^-)$ converge.
The limit is automatically positive since we are assuming that the
first case does not hold for any subsequence. Denote by
$S_n=S_c^{\lambda_n}(t_i^-)$ the unit tangent vector to $\widetilde
\Gamma_c^{\lambda_n}(t_i^-)$ and by $u_n$ the inner product $\langle
S_n,U\rangle$. Now we have a family of loops with tangent vectors
and all higher derivatives bounded. Since $u_n$ is everywhere
positive, since $\int u_nds=\lambda_n$, since the length of the loop
$\widetilde \Gamma_c^{\lambda_n}(t_i^-)$ is bounded away from $0$
independent of $n$, and since
$|(u_n)'|=|\langle\nabla_{S_n}S_n,U\rangle|$ is bounded above
independent of $n$, we see that $u_n$ tends uniformly to zero as $n$
tends to infinity. This means that the $|p_1(S_n)|$ converge
uniformly to one as $n$ goes to infinity.
 Since
the ambient manifold is compact, passing to a further subsequence we
have a smooth limit of the
$p_1\widetilde\Gamma_c^{\lambda_n}(t_i^-)$. The result is an
immersed curve in $(M,g(t_i^-))$ parameterized at unit speed. Since
all the spatial and time derivatives of the $p_1\widetilde
\Gamma_c^{\lambda_n}(t)$ are uniformly bounded, by passing to a
further subsequence, there is a  smooth map $f\colon S^1\times
[t_i^-,t_i^+]\to M$ which is a smooth limit of the sequence
$\widetilde\Gamma_c^{\lambda_n}(t),\ t_i^-\le t\le t_i^+$. If for
some $t\in [t_i^-,t_i^+]$ the curve $f|_{S^1\times\{t\}}$ is
immersed, then this limiting map along this curve agrees to first
order with the curve-shrinking flow. Thus, for some $t>t_i^-$ the
restriction of $f$ to the interval $[t_i^-,t]$ is a curve-shrinking
flow. We claim that $f$ is a curve-shrinking flow on the entire
interval $[t_i^-,t_i^+]$. Suppose not. Then there is a first $t'\le
t_i^+$ for which $f|_{S^1\times \{t'\}}$ is not an immersion.
According to Lemma~\ref{curvshrexist} the maximum of the norms of
the curvature of the curves $f(t)$ must tend to infinity as $t$
approaches $t'$ from below. But the curvatures of $f(t)$ are the
limits of the curvatures of the family $p_1\widetilde
\Gamma_c^{\lambda_n}(t)$ and hence are uniformly bounded on the
entire interval $[t_i^-,t_i^+]$.
 This contradiction shows that the
entire limiting surface
$$f\colon S^1\times [t_i^-,t_i^+]\to (M,g(t))$$
is a curve-shrinking flow of immersions.
\end{proof}

\begin{rem}\label{xicomp}
Notice that if the first case holds then by the choice of $B$ we
have a point $t_n\in J_B(c)$ for which the length of
$\widetilde\Gamma_c^{\lambda_n}(t_n)$ is less than
$e^{-C_2(t_1-t_0)}\zeta/3$.
\end{rem}

For each $n$,  the family of curves
$p_1\widetilde\Gamma_c^{\lambda_n}(t)$  in $M$ all have
$p_1\widetilde\Gamma_c^{\lambda_n}(t_0)=\widetilde\Gamma(c)$ as
their initial member. Thus, these curves are all homotopically
trivial. Hence, for each $t\in J_B(c)$ the limiting curve
$\widetilde\Gamma(c)(t)$ of the
$p_1\widetilde\Gamma_c^{\lambda_n}(t)$ is then also homotopically
trivial. It now follows from Lemma~\ref{nonembeasycase},
Claim~\ref{limitcs} and Remark~\ref{xicomp} that one of the
following two conditions holds:
\begin{enumerate}
\item[(1)] for some $t\in J_B(c)$ the length of $\widetilde\Gamma(c)(t)$ is less than or equal to
 $e^{-C_2(t_1-t_0)}\zeta/3$ or
\item[(2)] the function $A(t)$ that assigns to each $t\in
J_B(c)$ the area of the minimal spanning disk for
$p_1\widetilde\Gamma(c)(t)$ satisfies
$$\frac{dA(t)}{dt}\le -2\pi-\frac{1}{2}R_{\rm min}(t)A(t)$$
in the sense of forward difference quotients.
\end{enumerate}
By continuity,  for any $\delta''>0$ then for all $n$
 sufficiently large one of the following two conditions holds:
\begin{enumerate}
\item[(1)] there is $t_n\in J_B(c)$ such that the length of $\widetilde\Gamma_c^{\lambda_n}(t_n)$
is less than $e^{-C_2(t_1-t_0)}\zeta/2$, or
\item[(2)]  for every $t\in J_B(c)$, the areas of the minimal
spanning disks for $p_1(\widetilde\Gamma_c^{\lambda_n}(t))$ satisfy
$$\frac{dA(p_1\widetilde\Gamma_c^{\lambda_n}(t))}{dt}\le -2\pi-\frac{1}{2}R_{\rm
min}(t)A(p_1\widetilde\Gamma_c^{\lambda_n}(t))+\delta''$$ in the
sense of forward difference quotients.
\end{enumerate}

Suppose  that for every $n$ sufficiently large, for every  $t\in
J_B(c)$ the length of $\widetilde\Gamma_c^{\lambda_n}(t)$ is at
least $e^{-C_2(t_1-t_0)}\zeta/2$. We have already seen in
Claim~\ref{areachange} that for every $t'<t''$ in $[t_0,t_1]$ the
areas satisfy
$$A(p_1(\widetilde\Gamma_c^{\lambda_n}(t'')))-A(p_1(\widetilde\Gamma_c^{\lambda_n}(t')))\le
C_4(t''-t').$$ Since the total length of the complement $J_B(c)$ in
$[t_0,t_1]$ is at most $3C_5B^{-1}$, it follows from our choice of
$B$ that this total length is at most the constant $\delta'$ of
Claim~\ref{C1C4delta}. Invoking Claim~\ref{deltaprime} and the fact
that $A(\widetilde\Gamma(c))\le W(\Gamma(c))+\zeta\le
W(\Gamma)+\zeta=\widetilde A$, we see that for all $n$ sufficiently
large we have
$$A(p_1(\widetilde\Gamma^c_{\lambda_n}(t_1)))-v_c(t_1)<\zeta/2.$$

The other possibility to consider is that for each $n$ there is
$t_n\in J_B(c)$ such that the length of
$\widetilde\Gamma_c^{\lambda_n}(t_n)<e^{-C_2(t_1-t_0)}\zeta/2$.
Since $J_B(c)\subset [t_0,t_1-B^{-1}]$, in this case we invoke the
first inequality in Corollary~\ref{lncurvgrowth} to see that the
length of $\widetilde\Gamma_c^{\lambda_n}(t)<\zeta/2$ for every
$t\in [t_1-B^{-1},t_1]$. This completes the proof of
Lemma~\ref{xistronger}.
\end{proof}

\subsection{The completion of the proof of
Proposition~\protect{\ref{XI}}}

Now we wish to pass from  Lemma~\ref{xistronger} which deals with an
individual $c\in S^2$ to a proof of Proposition~\ref{XI} which deals
with the entire family $\widetilde\Gamma$. Let us introduce the
following notation. Suppose that $\omega\subset S^2$ is an arc. Then
$\widetilde\Gamma(\omega)=\cup_{c\in \omega}\widetilde \Gamma(c)$ is
an annulus in $M$ and for each $t\in [t_0,t_1]$ we have
 the annulus $\widetilde\Gamma_\omega^\lambda(t)$ in $M\times
S^1_\lambda$.

A finite set ${\mathcal S}\subset S^2$ with the property that for
$c\in S^2$ there is $\hat c\in {\mathcal S}$ and an arc $\omega$ in
$S^2$ joining $c$ to $\hat c$ so that the area of the annulus
$\widetilde\Gamma(\omega)$ is less than $\nu$ is called {\em a
$\nu$-net} for $\widetilde\Gamma$. Similarly, if for every $c\in
S^2$ there is $\hat c\in {\mathcal S}$ and an arc $\omega$
connecting them for which the area of the annulus
$\widetilde\Gamma_\omega^\lambda(t_0)$ is less than $\nu$, we say
that ${\mathcal S}$ is {\em a $\nu$-net} for
$\widetilde\Gamma^\lambda$. Clearly, for any $\nu$ there is a subset
${\mathcal S}\subset S^2$ that is a $\nu$-net for $\widetilde\Gamma$
and for $\widetilde\Gamma^\lambda$ for all $\lambda$ sufficiently
small.

\begin{lem}\label{muarea}
There is a $\mu>0$ such that the following holds. Let $c,\hat c\in
S^2$. Suppose that there is an arc $\omega$ in $S^2$ connecting $c$
to $\hat c$ with the area of the annulus
$\widetilde\Gamma_\omega^\lambda(t_0)$ in $M\times S^1_\lambda$ less
than $\mu$. Let $v_{\hat c}$, resp., $v_c$, be the solution to
Equation~(\ref{wdiffeqn}) with initial condition $v_{\hat
c}(t_0)=A(\widetilde\Gamma(\hat c))$, resp.,
$v_c(t_0)=A(\widetilde\Gamma(c))$. If
 $$A(p_1\widetilde\Gamma^\lambda_{\hat c}(t_1))\le v_{\hat c}+\zeta/2,$$
then
$$A(p_1(\widetilde\Gamma^\lambda_c(t_1))\le v_c+\zeta.$$
\end{lem}

\begin{proof}
 First of all we require that $\mu<e^{-(2n-1)C'(t_1-t_0)}\zeta/4$
where $C'$ is an upper bound for the norm of the Riemann curvature
tensor at any point of the ambient Ricci flow.  By
Lemma~\ref{mugrowth} the fact that the area of the minimal annulus
between the ramps $\widetilde\Gamma_c^\lambda(t_0)$ and
$\widetilde\Gamma_{\hat c}^\lambda(t_0)$ is less than $\mu$ implies
that the area of the minimal annulus between the ramps
$\widetilde\Gamma_c^\lambda(t_1)$ and $\widetilde\Gamma_{\hat
c}^\lambda(t_1)$ is less than $\mu e^{(2n-1)C'(t_1-t_0)}=\zeta/4$.
The same estimate also holds for the image under the projection
$p_1$ of this minimal annulus. Thus, with this condition on $\mu$,
and for $\lambda$ sufficiently small, we have
$$\left|A(p_1\widetilde\Gamma_c^\lambda(t_1))-A(p_1\widetilde\Gamma_{\hat c}^\lambda(t_1))\right|<\zeta/4.$$
The other condition we impose upon $\mu$ is that if $a,\hat a$ are
positive numbers at most $W(\Gamma)+\zeta$ and if $a<\hat a+\mu$
then
$$w_{a,t_0}(t_1)<w_{\hat a,t_0}(t_1)+\zeta/4.$$
Applying this with $a=A(\widetilde\Gamma(c))$ and $\hat
a=A(\widetilde\Gamma(\hat c))$ (both of which are at most
$W(\widetilde \Gamma)<W(\Gamma)+\zeta$), we see that these two
conditions on $\mu$ together imply the result.
\end{proof}

 We must also examine what
happens if the second alternative holds for $\widetilde\Gamma_{\hat
c}^\lambda$. We need the following lemma to treat this case.

\begin{lem}\label{BARMU}
There is  $\delta>0$ such that for any $r>0$  there is $\bar\mu>0$,
depending on $r$ and on the curvature bound for the ambient Ricci
flow such that the following holds. Suppose that $\gamma$ and $\hat
\gamma$ are ramps in $(M,g(t))\times S^1_\lambda$. Suppose that the
length of $\gamma$ is at least $r$ and suppose that on any
sub-interval $I$ of $\gamma$ of length $r$ we have
$$\int_Ikds<\delta.$$
 Suppose also that there is an annulus connecting $\gamma$ and
$\hat\gamma$ of area less than $\bar\mu$. Then the length of
$\hat\gamma$ is at least $3/4$ the length of $\gamma$.
\end{lem}

We give a proof of this lemma in the next section. Here we finish
the proof of Proposition~\ref{XI} assuming it.

\begin{claim}\label{mulength}
There is $\mu>0$ such that the following holds. Suppose that $c,\hat
c\in S^2$ are such that there is an arc $\omega$ in $S^2$ connecting
$c$ and $\hat c$ such that the area of the annulus
$\widetilde\Gamma^\lambda(\omega)$ is at most $\mu$.  Set
$t_2=t_1-B^{-1}$. If the length of $\widetilde\Gamma_{\hat
c}^\lambda(t)$ is less than $\zeta/2$ for all $t\in [t_2,t_1]$, then
the length of $p_1\widetilde\Gamma_c^\lambda(t_1)$ is less than
$\zeta$.
\end{claim}

\begin{proof}
The proof is by contradiction: Suppose that the length of
$p_1\widetilde\Gamma_c^\lambda(t_1)$ is at least $\zeta$ and the
length of $\widetilde\Gamma_{\hat c}^\lambda(t)$ is less than
$\zeta/2$ for all $t\in [t_2,t_1]$. Of course, it follows that the
length of $\widetilde\Gamma_c^\lambda(t_1)$ is also at least
$\zeta$. The third condition on $B$ is equivalent to
$$e^{C_2B^{-1}}<4/3.$$
It then follows from
 Corollary~\ref{lncurvgrowth} that for every $t\in
 [t_2,t_1]$ the length of $\widetilde\Gamma^c_\lambda(t)$ is at least $3\zeta/4$.
On the other hand, by hypothesis for every such $t$, the length of
$\widetilde\Gamma_{\hat c}^\lambda(t)$ is less than $\zeta/2$. It
follows from Equation~(\ref{Lformula}) that
$$\int_{t_2}^{t_1}\left(\int k^2ds\right)dt\le
C_2\left(\int_{t_2}^{t_1}L(\widetilde\Gamma_c^\lambda(t))dt\right)-
L(\widetilde\Gamma_c^\lambda(t_1))+L(\widetilde\Gamma_c^\lambda(t_2)).$$
(Here $L$ is the length of the curve.) From this and
Corollary~\ref{lncurvgrowth} we see that there is a constant $C_8$
depending on the original family $\Gamma$ and on the curvature of
the ambient Ricci flow such that
$$\int_{t_2}^{t_1}\left(\int_{\widetilde\Gamma_c^\lambda(t)} k^2ds\right)dt\le C_8.$$
Since $t_1-t_2=B^{-1}$, this implies that there is $t'\in [t_2,t_1]$
with
$$\int_{\widetilde\Gamma_c^\lambda(t')}k^2ds\le C_8B.$$
By Cauchy-Schwarz, for any subinterval $I$ of length $\le r$ in
$\widetilde\Gamma_c^\lambda(t')$ we have
$$\int_Ikds\le \sqrt{C_8Br}.$$
We choose $0<r\le \zeta$ sufficiently small so that $\sqrt{C_8Br}$
is less than or equal to the constant $\delta$ given in
Lemma~\ref{BARMU}. Then we set $\bar \mu$ equal to the constant
given by that lemma for this value of $r$.

Now suppose that $\mu$ is sufficiently small so that the solution to
the equation
$$\frac{d\mu(t)}{dt}=(2n-1)|{\rm Rm}_{g(t)}|\mu(t)$$
with initial condition $\mu(t_0)\le \mu$ is less than $\bar\mu$ on
the entire interval $[t_0,t_1]$. With this condition on $\mu$,
Lemma~\ref{mugrowth} implies that for every $t\in [t_0,t_1]$ the
ramps $\widetilde\Gamma_c^\lambda(t)$ and $\widetilde\Gamma_{\hat
c}^\lambda(t)$ are connected by an annulus of area at most
$\bar\mu$. In particular, this is true for
$\widetilde\Gamma_c^\lambda(t')$ and $\widetilde\Gamma_{\hat
c}^\lambda(t')$. Now we have all the hypotheses of Lemma~\ref{BARMU}
at time $t'$. Applying this lemma we conclude that
$$L(\widetilde\Gamma_{\hat c}^\lambda(t'))\ge \frac{3}{4}L(\widetilde\Gamma_c^\lambda(t')).$$
But this is a contradiction since by assumption
$L(\widetilde\Gamma_{\hat c}^\lambda(t'))<\zeta/2$ and the
supposition that $L(p_1\widetilde\Gamma_c^\lambda(t_1))\ge \zeta$
led to the conclusion that $L(\widetilde\Gamma_c^\lambda(t'))\ge
3\zeta/4$. This contradiction shows that our supposition that
$L(p_1\widetilde\Gamma_c^\lambda(t_1))\ge \zeta$ is false.
\end{proof}

Now we complete the proof of Proposition~\ref{XI}.

\begin{proof} (of Proposition~\ref{XI}.)
Fix $\mu>0$ sufficiently small so that Lemma~\ref{muarea} and
Claim~\ref{mulength} hold. Then we choose a $\mu/2$-net $X$ for
$\widetilde\Gamma$. We take $\lambda$ sufficiently small so that
Lemma~\ref{xistronger} holds for every $\hat c\in {\mathcal S}$. We
also choose $\lambda$ sufficiently small so that $X$ is a $\mu$-net
for $\widetilde\Gamma^\lambda$. Let $c\in S^2$. Then there is $\hat
c\in {\mathcal S}$ and an arc $\omega$ connecting $c$ and $\hat c$
such that the area of $\widetilde\Gamma^\lambda(\omega)<\mu$. Let
$v_{\hat c}$, resp., $v_c$ be the solution to
Equation~\ref{wdiffeqn} with initial condition $v_{\hat
c}(t_0)=A(\widetilde\Gamma({\hat c}))$, resp.,
$v_c(t_0)=A(\widetilde\Gamma(c))$. According to
Lemma~\ref{xistronger} either $A(p_1\widetilde\Gamma_{\hat
c}^\lambda(t_1))<v_{\hat c}(t_1)+\zeta/2$ or the length of
$\widetilde\Gamma_{\hat c}^\lambda(t)$ is less than $\zeta/2$ for
every $t\in [t_2,t_1]$ where $t_2=t_1-B^{-1}$. In the second case,
Claim~\ref{mulength} implies that the length of
$p_1\widetilde\Gamma_c^\lambda(t_1)$ is less than $\zeta$. In the
first case, Lemma~\ref{muarea} tells us that
$A(p_1\widetilde\Gamma_c^\lambda(t_1))<v_c(t_1)+\zeta$. This
completes the proof of Proposition~\ref{XI}.
\end{proof}

\section{Proof of Lemma~\protect{\ref{BARMU}}: annuli of small area}\label{18.5}

Except for the brief comments that follow, our proof involves
geometric analysis that takes place on an abstract annulus with
bounds on its area, upper bounds on its Gaussian curvature, and on
integrals of the geodesic curvature on the boundary.
Proposition~\ref{annlengthest} below gives the precise result along
these lines. Before stating that proposition, we show that its
hypotheses hold in the situation that arises in Lemma~\ref{BARMU}.
Let us recall  the situation of Lemma~\ref{BARMU}. We have ramps
$\gamma$ and $\hat \gamma$ in which are real analytic embedded
curves in the real analytic Riemannian manifold $(M\times
S^1_\lambda,g\times ds^2)$. By a slight perturbation we can assume
they are disjointly embedded. These curves that are connected by an
annulus $A_0\to M\times S^1_\lambda$ of small area, an area bounded
above by, say,  $\mu$. We take an energy minimizing map of an
annulus $\psi\colon A\to M\times S^1_\lambda$ spanning
$\gamma\coprod \hat \gamma$. According to \cite{Hildebrandt}, $\psi$
is a real analytic map and the only possible singularities
(non-immersed points) of the image come from the {\em branch points}
of $\psi$, i.e., points where $d\psi$ vanishes. There are finitely
many branch points. If there are branch points on the boundary, then
the restriction of $\psi$ to $\partial A$ will be a homeomorphism
rather than a diffeomorphism onto $\gamma\coprod \hat \gamma$.
Outside the branch points, $\psi$ is a conformal map onto its image.
The image is an area minimizing annulus spanning $\gamma\coprod \hat
\gamma$. Thus, the area of the image is at most $\mu$. According to
\cite{White} the only branch points on the boundary are false branch
points, meaning that a local smooth reparameterization of the map on
the interior of $A$ near the boundary branch point removes the
branch point. These reparameterizations produce a new smooth
structure on $A$, identified with the original smooth structure on
the complement of the boundary branch points. Using this new smooth
structure on $A$ the map $\psi$ is an immersion except at finitely
many interior branch points. From now on the domain surface $A$ is
endowed with this new smooth structure. Notice that, after this
change, the domain is no longer real analytic; it is  only smooth.
Also, the original annular coordinate is not smooth at the finitely
many boundary branch points.

The pullback of the metric $g\times ds^2$ is a smooth symmetric
two-tensor on $A$. Off the finite set of interior branch points it
is positive definite and hence a Riemannian metric, and in
particular, it is a Riemannian metric near the boundary. It vanishes
at each interior branch point.
  Since the geodesic
curvature $k_{\rm geod}$ of the boundary of the annulus is given by $k\cdot n$
where $n$ is the unit normal vector along the boundary pointing into $A$, we
see that the restriction of the geodesic curvature to $\gamma$, $k_{\rm
geod}\colon \gamma\to \Ar$ has the property that for any sub-arc $I$ of
$\gamma$ of length $r$ we have
$$\int_I|k_{\rm geod}|ds<\delta.$$ Lastly, because the map of $A$ into $M\times
S^1_\lambda$ is minimal, off  the set of interior branch points, the
Gaussian curvature of the pulled back metric is bounded above by the
upper bound for the sectional curvature of $M\times S^1_\lambda$,
which itself is bounded independent of $\lambda$ and $t$, by say
$C'>0$.

Next, let us deal with the singularities of the pulled back metric on $A$
caused by the interior branch points. As the next claim shows, it is an easy
matter to deform the metric slightly near each branch point without increasing
the area much and without changing the upper bound on the Gaussian curvature
too much. Here is the result:

\begin{claim}
Let $\psi\colon A\subset M\times S^1_\lambda$ be an area-minimizing
annulus of area at most $\mu$ with smoothly embedded boundary as
constructed above. Let $h$ be the induced (possibly singular) metric
on $A$ induced by pulling back $g\times ds^2$ by $\psi$, and let
$C''>0$ be an upper bound on the Gaussian curvature of $h$ (away
from the branch points). Then  there is a deformation $\tilde h$ of
$h$, supported near the interior branch points, to a smooth metric
with the property that the area of the deformed smooth metric is at
most $2\mu$ and where the upper bound for the curvature of $\tilde
h$ is $2C''$.
\end{claim}

\begin{proof} Fix an interior branch point $p$.
 Since $\psi$  is smooth and
conformal onto its image, there is a disk in $A$ centered at $p$ in
which $h=f(z,\bar z)|dz|^2$ for a smooth function $f$ on the disk.
The function $f$ vanishes at the origin and is positive on the
complement of the origin. Direct computation shows that the Gaussian
curvature $K(h)$ of $h$ in this disk is given by
$$K(h)=\frac{-\triangle f}{2f^2}+\frac{|\nabla f|^2}{2f^3}\le C,$$
where $\triangle$ is the usual Euclidean Laplacian on the disk and $|\nabla
f|^2=(\partial f/\partial x)^2+(\partial f/\partial y)^2$.  Now consider the
metric $(f+\epsilon)|dz|^2$ on the disk. Its Gaussian curvature is
$$\frac{-\triangle f}{2(f+\epsilon)^2}+\frac{|\nabla f|^2}{2(f+\epsilon)^3}.$$
\begin{claim} For all $\epsilon>0$ the Gaussian curvature of $(f+\epsilon)|dz|^2$
is at most $2C''$. \end{claim}

\begin{proof}
We see that $-\triangle f\le C''f^2$, so that
\begin{eqnarray*}\frac{-\triangle
f}{(f+\epsilon)^2}+\frac{|\nabla f|^2}{(f+\epsilon)^3} & = &
\frac{(f+\epsilon)(-\triangle
f)+|\nabla f|^2}{(f+\epsilon)^3} \\
& \le & \frac{C''f^3-\epsilon\triangle f}{(f+\epsilon)^3}\le
C''+\frac{\epsilon f^2C''}{(f+\epsilon)^3}\le 2C''.
\end{eqnarray*}
\end{proof}

Now we fix a smooth function $\rho(r)$ which is identically one on a subdisk
$D'$ of $D$ and vanishes near $\partial D$ and we replace the metric $h$ on the
disk by
$$h_\epsilon=(f+\epsilon \rho(r))|dz|^2.$$ The above computation shows that the
Gaussian curvature of $h_\epsilon$ on $D'$ is bounded above by
$2C''$. As $\epsilon$ tends to zero the restriction of the metric
$h_\epsilon$ to $D\setminus D'$ converges uniformly in the
$C^\infty$-topology to $h$. Thus, for all $\epsilon>0$ sufficiently
small the Gaussian curvature of $h_\epsilon$ on $D\setminus D'$ will
also be bounded by $2C''$. Clearly, as $\epsilon$ tends to zero the
area of the metric $h_\epsilon$ on $D$ tends to the area of $h$ on
$D$.

Performing this construction near each of the finite number of interior branch
points and taking $\epsilon$ sufficiently small gives the perturbation $\tilde
h$ as required.
\end{proof}

Thus, if $\gamma$ and $\hat \gamma$ are ramps as in
Lemma~\ref{BARMU}, then replacing $\hat \gamma$ by a close $C^2$
approximation we have an abstract smooth annulus with a Riemannian
metric connecting $\gamma$ and $\hat \gamma$. Taking limits shows
that  establishing the conclusion of Lemma~\ref{BARMU} for a
sequence of better and better approximations to $\hat \gamma$ will
also establish it for $\hat \gamma$. This allows us to assume that
$\gamma$ and $\hat \gamma$ are disjoint. The area of this annulus is
bounded above by a constant arbitrarily close to $\mu$. The Gaussian
curvature of the Riemannian metric is bounded above by a constant
depending only on the curvature bounds of the ambient Ricci flow.
Finally, the integral of the absolute value of the geodesic
curvature over any interval of length $r$ of $\gamma$ is at most
$\delta$.

With all these preliminary remarks, we see that Lemma~\ref{BARMU}
follows from:

\begin{prop}\label{annlengthest}
Fix $0<\delta<1/100$.  For each $0<r$ and $C''<\infty$ there is a
$\mu>0$ such that the following holds. Suppose that $A$ is an
annulus with boundary components $c_0$ and $c_1$. Denote by $l(c_0)$
and $l(c_1)$ the lengths of $c_0$ and $c_1$, respectively. Suppose
that the Gaussian curvature of $A$ is bounded above by $C''$.
Suppose that  $l(c_0)>r$ and that for each sub-interval $I$ of $c_0$
of length $r$, the integral of the absolute value of the geodesic
curvature along $I$ is less than $\delta$. Suppose that the area of
$A$ is less than $\mu$. Then
$$l(c_1)\ge \frac{3}{4}l(c_0).$$
\end{prop}

To us, this statement was intuitively extremely reasonable but we
could not find a result along these lines stated in the literature.
Also, in the end, the argument we constructed is quite involved,
though elementary.

The intuition is that we exponentiate in from the boundary component
$c_0$ using the family of geodesics perpendicular to the boundary.
The bounds on the Gaussian curvature and local bounds on the
geodesic curvature of $c_0$ imply that the exponential mapping will
be an immersion out to some fixed distance $\delta$ or until the
geodesics  meet the other boundary, whichever comes first.
Furthermore, the metric induced by  this immersion will be close to
the product metric. Thus, if there is not much area, it must be the
case that, in the measure sense, most of the geodesics in this
family must meet the other boundary before distance $\delta$. One
then deduces the length inequality. There are two main difficulties
with this argument that must be dealt with. The first is due to the
fact that we do not have a pointwise bound on the geodesic curvature
of $c_0$, only an integral bound of the absolute value over all
curves of short length. There may be points of arbitrarily high
geodesic curvature. Of course, the length of the boundary where the
geodesic curvature is large is very small. On these small intervals
the exponential mapping will not be an immersion out to any fixed
distance. We could of course, simply omit these regions from
consideration and work on the complement. But these small regions of
high geodesic curvature on the boundary can cause focusing (i.e.,
crossing of the nearby geodesics). We must estimate out to what
length along the boundary this happens. Our first impression was
that the length along the boundary where focusing occurred would be
bounded in terms of the total turning along the arc in $c_0$. We
were not able to establish this. Rather we found a weaker estimate
where this focusing length is bounded in terms of the total turning
 and the area bounded by the triangle cut out by the two geodesics
that meet. This is a strong enough result for our application. Since
the area is small and the turning on any interval of length $r$ is
small, a maximal collection of focusing regions will meet each
interval of length $r$ in $c_0$ in a subset of small total length.
Thus, on the complement (which is most of the length of $c_0$) the
exponential mapping will be an immersion out to length $\delta$ and
will be an embedding when restricted to each interval of length one.
The second issue to face is to show that the exponential mapping on
this set is in fact an embedding, not just an immersion. Here one
uses standard arguments invoking the Gauss-Bonnet theorem to rule
out various types of pathologies, e.g., that the individual
geodesics are not embedded or geodesics that end on $c_0$ rather
than $c_1$, etc. Once these are ruled out, one has established that
the exponential map on this subset is an embedding and the argument
finishes as indicated above.

\subsection{First reductions}

Of course, if the hypothesis of the proposition holds for $r>0$ then
it holds for any $0<r'<r$. This allows us to assume that $r<{\rm
min}((C'')^{-1/2},1)$. Now let us scale the metric by $4r^{-2}$. The
area of $A$ with the rescaled metric is  $4r^{-2}$ times the area of
$A$ with the original metric. The Gaussian curvature of $A$ with the
rescaled metric is less than $(r^2C''/4)\le 1$. Furthermore, in the
rescaled metric $c_0$ has length greater than $2$ and the total
curvature along any interval of length $1$ in $c_0$ is at most
$\delta$. This allows us to assume (as we shall) that  $r=1$, that
$C''\le 1$, and that $l(c_0)\ge 2$. We must find a $\mu>0$ such that
the proposition holds provided that the area of the annulus is less
than $\mu$.

The function $k_{\rm geod}\colon c_0\to \Ar$ is smooth. We choose a
regular value $\alpha$ for $k_{\rm geod}$ with $1<\alpha<1.1$. In
this way we divide $c_0$ into two disjoint subsets, $Y$ where
$k_{\rm geod}>\alpha$, and $X$ where $k_{\rm geod}\le \alpha$. The
subset $Y$ is a union of finitely many disjoint open intervals and
$X$ is a disjoint union of finitely many closed intervals.

\begin{rem}\label{Ylengthrem}
The condition on $k_{\rm geod}$ implies that for any arc $J$ in
$c_0$ of length $1$ the total
 length of $J\cap Y$ is less than $\delta$.
\end{rem}

Fix $\delta'>0$. For each $x\in X$ there is a geodesic $D_x$ in $A$
whose initial point is $x$ and whose initial direction is orthogonal
to $c_0$. Let $f(x)$ be the minimum of $\delta$ and the distance
along $D_x$ to the first point (excluding $x$) of its intersection
with $\partial A$. We set
$$S_X(\delta') =\{(x,t)\in X\times [0,\delta']\bigl|\bigr.\  t\le f(x)\}.$$
The subset $S_X(\delta')$ inherits a Riemannian metric from the
product of the  metric on $X$ induced by the embedding $X\subset
c_0$ and the standard metric on the interval $[0,\delta']$.

\begin{claim}\label{expcomp}
There is $\delta'>0$ such that the following holds.
 The
exponential mapping defines a map ${\rm exp}\colon S_X(\delta')\to
A$ which is a local diffeomorphism and the pullback of the metric on
$A$ defines a metric on $S_X(\delta')$ which is at least
$(1-\delta)^2$ times the given product metric.
\end{claim}

\begin{proof}
This is a standard computation using the Gaussian curvature upper
bound and the geodesic curvature bound.
\end{proof}

Now we fix $0<\delta'<1/10$ so that Claim~\ref{expcomp} holds, and
we set $S_X=S_X(\delta')$. We define
$$\partial_+S_X=\{(x,t)\in S_X \bigl|\bigr.\, t=f(x)\}.$$
Then the boundary of $S_X$ is made up of $X$, the arcs $\{x\}\times
[0,f(x)]$ for $x\in \partial X$ and $\partial_+(S_X)$. For any
subset $Z\subset X$ we denote by $S_Z$ the intersection $(Z\times
[0,\delta])\cap S_X$, and we denote by $\partial_+(S_Z)$ the
intersection of $S_Z\cap \partial_+S_X$.

Lastly, we fix $\mu>0$ with $\mu<(1-\delta)^2(\delta')/10$. Notice
that this implies that $\mu<1/100$. We now assume that the area of
$A$ is less than this value of $\mu$ (and recall that $r=1$, $C''=1$
and $l(c_0)\ge 2$). We must show that $l(c_1)>3l(c_0)/4$.

\subsection{Focusing triangles}

\begin{figure}[ht]
  \relabelbox{
  \centerline{\epsfbox{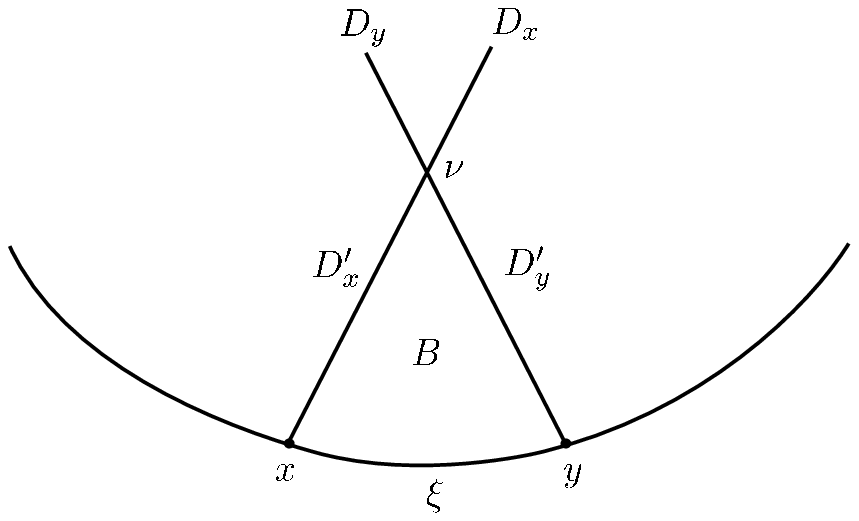}}}
  \relabel{1}{$x$}
  \relabel{2}{$y$}
  \relabel{3}{$D_y$}
  \relabel{4}{$D_x$}
  \relabel{5}{$D_x'$}
  \relabel{6}{$D_y'$}
   \relabel{e}{$\xi$}
  \endrelabelbox
  \caption{Focusing triangle.}\label{fig:focus}
\end{figure}

By a {\em focusing triangle} we mean the following. We have distinct
points $x,y\in X$ and sub-geodesics $D'_x\subset D_x$ and
$D'_y\subset D_y$ that are embedded arcs with $x$, respectively $y$,
as an endpoint. The intersection $D'_x\cap D'_y$ is a single point
which is the other endpoint of each of $D'_x$ and $D'_y$. Notice
that since $D'_x\subset D_x$ and $D'_y\subset D_y$, by construction
both $D'_x$ and $D'_y$ have lengths at most $\delta'$. We have an
arc $\xi$ in $c_0$ with endpoints $x$ and $y$ and the loop
$\xi*D'_y*(D'_x)^{-1}$ bounds a disk $B$ in $A$. The arc $\xi$ is
called the {\em base} of the focusing triangle and with, respect to
an orientation of $c_0$, if $x$ is the initial point of $\xi$ then
$D'_x$ is called the {\em left-hand side} of the focusing triangle
and $D'_y$ is called its {\em right-hand side}. See {\sc
Fig.}~\ref{fig:focus}.

Our main goal here is the following lemma which gives an upper bound
for the length of the base, $\xi$, of a focusing triangle in terms
of the turning along the base and the area of the region $B$
enclosed by the triangle.

\begin{lem}\label{laest}
 Suppose that we have a focusing triangle ${\mathcal T}$ with base $\xi$ bounding a disk $B$
 in $A$.
 Suppose that the length of $\xi$ is at most one.
 Then
$$l(\xi)\le \left(\int_\xi k_{\rm geod}ds+{\rm Area}(B)\right).$$
\end{lem}

\begin{proof}
We begin with a preliminary computation. We denote by $a(B)$ the
area of $B$. We define
$$t_{\xi}=\int_{\xi}k_{\rm geod}ds \ \ \ {\rm and} \ \ \ T_\xi=\int_{\xi}|k_{\rm geod}|ds.$$

Recall that given a piecewise smooth curve, its {\em total turning}
is the integral of the geodesic curvature over the smooth part of
the boundary plus the sum over the break points of $\pi$ minus the
interior angle at the break point.  The Gauss-Bonnet theorem tells
us that for any compact surface with piecewise smooth boundary the
integral of the Gaussian curvature over the interior of the surface
plus the total turning around the boundary equals $2\pi$ times the
Euler characteristic of the surface.

\begin{claim}\label{theta/Kbd}
The angle $\theta_B$ between $D'_x$ and $D'_y$ at the vertex $v$
satisfies $$\theta_B\le t_\xi+a(B)$$ and for any measurable subset
$B'\subset B$ we have
$$\theta_B-t_\xi-a(B)\le \int_{B'}Kda<a(B).$$
\end{claim}

\begin{proof}
 Since $D'_x$ and
$D'_y$ meet $\partial A$ in right angles, the total turning around
the boundary of $B$  is
$$t_\xi+2\pi-\theta_B.$$
Thus, by Gauss-Bonnet, we have
$$\theta_B=\int_BKda+ t_\xi.$$
But $K\le 1$, giving the first stated inequality. On the other hand
$\int_BKda=\int_BK^+da+\int_BK^-da$, where $K^+={\rm max}(K,0)$ and
$K^-=K-K^+$. Since $0\le \int_BK^+da\le a(B)$ and $\int_BK^-\le 0$,
the second string of inequalities follows.
\end{proof}

In order to make the computation we need to know that this triangle
is the image under the exponential mapping of a spray of geodesics
out of the vertex $v$. Establishing that requires some work.

\begin{claim}\label{avgeo}
Let $a\in {\rm int}\,\xi$. There is a shortest path in $B$ from $a$
to $v$. This shortest path is a geodesic meeting $\partial B$ only
in its end points. It has length $\le(1/2)+\delta'$.
\end{claim}

\begin{proof}
The length estimate is obvious: Since $\xi$ has length at most $1$,
a path along $\partial A\cap B $ from $a$ to the closest of $x$ and
$y$ has length at most $1/2$. The corresponding side has length at
most $\delta'$. Thus, there is a path from $a$ to $v$ in $B$ of
length at most $(1/2)+\delta'$.

Standard convergence arguments show that there is a shortest path in
$B$ from $a$ to $v$. Fix $a\in {\rm int}(\partial A\cap B)$. It is
clear that the shortest path cannot meet either of the  `sides'
$D'_x$ and $D'_y$ at any point other than $v$. If it did, then there
would be an angle at this point and a local shortcut, cutting off a
small piece of the angle, would provide a shorter path. We must rule
out that the shortest path from $a$ to $v$ meets $\partial A\cap B$
in another point. If it does, let $a'$ be the last such point
(parameterizing the geodesic starting at $a$). The shortest path
from $a$ then leaves $\partial A$ at $a'$ in the direction tangent
at $a'$ to $\partial A$. (Otherwise, we would have an angle which
would allow us to shorten the path just as before.) This means that
we have a geodesic $\gamma$ from $v$ to $a'$ whose interior is
contained in the interior of $B$ and which is tangent to $\partial
A$ at $a'$. We label the endpoints of $\partial A\cap B$ so that the
union of $\gamma$ and the interval on $\partial A\cap B$ from $a'$
to $y$ gives a $C^1$-curve. Consider the disc $B'$ bounded
 by $\gamma$, the arc of $\partial A$ from $a'$
to $y$, and $D'_y$. The total turning around the boundary is at most
$3\pi/2+\delta$, and the integral of the Gaussian curvature over
$B'$ is at most the area of $B$, which is less than
$\mu<1/20<(\pi/4)-\delta$. This contradicts the Gauss-Bonnet
theorem.
\end{proof}

\begin{claim}
For any $a\in (\partial A\cap B)$ there is a unique minimal geodesic
in $B$ from $a$ to $v$.
\end{claim}

\begin{proof}
Suppose not; suppose there are two $\gamma$ and $\gamma'$ from $v$
to $a$. Since they are both minimal in $B$, each is embedded, and
they must be disjoint except for their endpoints. The upper bound on
the curvature and the Gauss-Bonnet theorem implies that the angles
that they make at each endpoint are less than $\mu<\pi/2$. Thus,
there is a spray of geodesics (i.e. geodesics determined by an
interval $\beta$ in the circle of directions at $v$) coming out of
$v$ and moving into $B$ with extremal members of the spray being
$\gamma$ and $\gamma'$. The geodesics $\gamma$ and $\gamma'$ have
length at most $(1/2)+\delta'$, and hence the exponential mapping
from $v$ is a local diffeomorphism on all geodesics of length at
most the length of $\gamma$. Since the angle they make at $a$ is
less than $\pi/2$ and since the exponential mapping is a local
diffeomorphism near $\gamma$, as we move in from the $\gamma$ end of
the spray we find geodesics from $v$ of length less than the length
of $\gamma$ ending on points
 of $\gamma'$. The same
Gauss-Bonnet argument shows that the angles that each of these
shorter geodesics makes with $\gamma'$ is at most $\mu$. Consider
the subset $\beta'$ of $\beta$ which are directions of geodesics in
$B$ of length $<(1/2)+\delta'$ that end on points of $\gamma'$ and
make an angle less than $\mu$ with $\gamma'$. We have just seen that
$\beta'$ contains an open neighborhood of the end of $\beta$
corresponding to $\gamma$. Since the Gaussian curvature is bounded
above by $1$, and these geodesics all have length at most
$1/2+\delta$, it follows that the exponential map is a local
diffeomorphism near all such geodesics. Thus, $\beta'$ is an open
subset of $\beta$. On the other, hand if the direction of
$\gamma''\not=\gamma'$ is a point $b''\in \beta$ which is an
endpoint of an open interval  $ \beta'$, and if this interval
separates $b''$ from the direction of $\gamma$ then the length of
$\gamma''$ is less than the length of each point in the interval.
Hence, the length of $\gamma''$ is less than $(1/2)+\delta'$.
Invoking Gauss-Bonnet again we see that the angle between $\gamma''$
and $\gamma'$ is $<\mu$.

This proves that if $U$ is an open  interval in $\beta'$ then the
endpoint of $U$ closest to the direction of $\gamma'$ is also
contained in $\beta'$ (unless that endpoint is the direction of
$\gamma'$). It is now elementary to see that $\beta'$ is all of
$\beta$ except the endpoint corresponding to $\gamma'$. But this is
impossible. Since the exponential mapping is a local diffeomorphism
out to distance $(1/2)+\delta'$, and since $\gamma'$ is embedded,
any geodesic from $v$ whose initial direction is sufficiently close
to that of $\gamma'$ and whose length is at most $(1/2)+\delta'$
will not cross $\gamma'$.
\end{proof}

See {\sc Fig.}~\ref{fig:spray}

\begin{figure}[ht]
  \centerline{\epsfbox{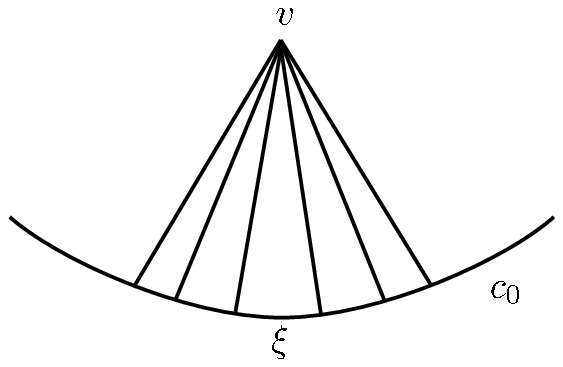}}
  \caption{Spray of geodesics from $v.$}\label{fig:spray}
\end{figure}

\begin{rem}
The same argument shows that from any $a\in (\partial A\cap B)$
there is a unique embedded geodesic in $B$ from $v$ to $a$ with
length at most $(1/2)+\delta'$. (Such geodesics may cross more than
once, but the argument given in the lemma applies to sub-geodesics
from $v$ to the first point of intersection along $\gamma$.)
\end{rem}

Let $E$ be the sub-interval of the circle of tangent directions at
$v$ consisting of all tangent directions of geodesics pointing into
$B$ at $v$. The endpoints of $E$ are the tangent directions for
$D'_x$ and $D'_y$. We define a function from $\xi$ to the interval
$E$ by assigning to each $a\in \xi$ the direction at $v$ of the
unique minimal geodesic in $B$ from $v$ to $a$. Since the minimal
geodesic is unique, this function is continuous and, by the above
remark, associates to $x$ and $y$ the endpoints of $E$. Since
geodesics are determined by their initial directions, this function
is one-to-one. Hence it is a homeomorphism from $\xi$ to $E$. That
is to say the spray of geodesics coming out of $v$ determined by the
interval $E$ produces a diffeomorphism between a wedge-shaped subset
of the tangent space at $v$ and $B$. Each of the geodesics in
question ends when it meets $\xi$.

Now that we have shown that the region enclosed by the triangle is
the image under the exponential map from the vertex $v$ of a
wedge-shaped region in the tangent space at $v$, we can make the
usual computation relating length and geodesic curvature. To do this
we pull back to the tangent space at $v$, and, using polar
coordinates, we write $\xi$ as $\{s=h(\psi); \psi\in E\}$ where $s$
is the radial coordinate and $\psi$ is the angular coordinate.
Notice that $h(\psi)\le (1/2)+\delta'$ for all $\psi\in E$. (In
fact, because the angles of intersection at the boundary are all
close to $\pi/2$ we can give a much better estimate on $h$ but we do
not need it.)
 We consider the one-parameter family of arcs $\lambda(t)$ defined
 to be the graph of the function $t\mapsto s(t)=th(\psi)$, for $0\le t\le
 1$. We set $l(t)$ equal to the length of
 $\lambda(t)$.

 \begin{claim}
$$\frac{dl}{dt}(t)\le {\rm max}_{\psi\in E}h(\psi)\int_{\lambda(t)}k_{\rm
geod}ds.$$
 \end{claim}

\begin{proof}
First of all notice that, by construction,  the curve $\xi$, which
is defined by $\{s=h(\psi)\}$, is orthogonal to the radial geodesics
to the endpoints. As a consequence, $h'(\psi)=0$ at the endpoints.
Thus, each of the curves $\lambda(t)$ is orthogonal to the radial
geodesics through its end points. Therefore, as we vary the family
$\lambda(t)$ the formula for the derivative of the length is
$$l'(t)=\int_{\lambda(t)}k_{\rm geod}(\psi)h(\psi)|{\rm
cos}(\theta(\psi,t))|ds$$ where $\theta(\psi,t)$ is the angle at
$(th(\psi),\psi)$ between the curve $s=th(\psi)$ and the radial
geodesic. The result follows immediately.
\end{proof}

Next, we must bound the turning of $\lambda(t)$. For this we invoke the
Gauss-Bonnet theorem once again. Applying this to the wedge-shaped disk $W(t)$
cut out by $\lambda(t)$ gives
$$\int_{W(t)}Kda+\int_{\lambda(t)}k_{\rm geod}ds=\theta_B.$$
From Claim~\ref{theta/Kbd} we conclude that
$$\int_{\lambda(t)}k_{\rm geod}ds\le t_\xi+a(B).$$

Of course, by Claim~\ref{avgeo} we have ${\rm max}_{\psi\in
E}h(\psi)\le (1/2)+\delta'$. Since $l(0)=0$, this implies that
$$l(\xi)=l(1)\le (a(B)+t_B)((1/2)+\delta')<a(B)+t_B.$$
This completes the proof of Lemma~\ref{laest}.
\end{proof}

\begin{cor}\label{muepsbd}
Suppose that ${\mathcal T}$ is a focusing triangle with base $\xi$
of length at most one. Then the length of $\xi$ is at most
$\delta+\mu$. More generally, suppose we have a collection of
focusing triangles ${\mathcal T}_1,\ldots,{\mathcal T}_n$ whose
bases all lie in a fixed interval of length one in $c_0$. Suppose
also that the interiors of disks bounded by these focusing triangles
are disjoint. Then the sum of the lengths of the bases is at  most
$\delta+\mu$.
\end{cor}

\begin{proof}
The first statement is immediate from the previous lemma. The second
comes from the fact that the sum of the areas of the disks bounded
by the ${\mathcal T}_i$ is at most $\mu$ and the sum of the total
turnings of the $\xi_i$ is at most $\delta$.
\end{proof}

This completes our work on the local focusing issue. It remains to
deal with global pathologies that would prevent the exponential
mapping from being an embedding out to distance $\delta'$.

\subsection{No $D_x$ is an embedded arc with both ends in $c_0$}

One thing that we must show  is that the geodesics $D_x$ are
embedded. Here is a special case that will serve some of our
purposes.

\begin{lem}\label{noloops}
For each $x\in X$, there is no non-trivial sub-geodesic of $D_x$
which is a homotopically trivial embedded loop in $A$.
\end{lem}

\begin{proof}
Were there such a loop, its  total turning would be $\pi$ minus the
angle it makes when the endpoints of the arc meet. Since $K\le 1$
and the area of the disk bounded by this loop is less than the area
of $A$ which in turn is less than $\mu<\pi$,  one obtains a
contradiction to the Gauss-Bonnet theorem.
\end{proof}

Next, we rule out the possibility that  one of the geodesics $D_x$
has both endpoints contained in $c_0$. This is the main result of
this section. In a sense, what the argument we give here shows that
if there is a $D_x$ with both ends on $c_0$, then under the
assumption of small area, $D_x$ cuts off a thin tentacle of the
annulus. But out near the end of this thin tentacle there must be a
short arc with large total turning, violating our hypothesis on the
integrals of the geodesic curvature over arcs of length at most one.

\begin{lem}\label{Dxacross}
There is no $x$ for which $D_x$ is an embedded arc  with both
endpoints on $c_0$ and otherwise disjoint from $\partial A$.
\end{lem}

\begin{proof}
Suppose that there were such a $D_x$. Then $D_x$ separates $A$ into
two components, one of which, $B$, is a topological disk. Let $c_0'$
be the intersection of $c_0$ with $B$. We consider two cases: Case
(i): $l(c'_0)\le 1$ and Case (ii): $l(c_0')>1$.

Let us show that the first case is not possible. Since $D_x$ is a
geodesic and $D_x$ is perpendicular to $c_0$ at one end, the total
turning around the boundary of $B$ is at most
$$3\pi/2+\int_{c_0'}k_{\rm geod}ds<3\pi/2+\delta,$$
where the last inequality uses the fact that the length of $c_0'$ is
at most one. On the other hand, $\int_{B}Kda<\mu$, and
$\mu<1/20<(\pi/2)-\delta$. This contradicts the Gauss-Bonnet
theorem.

Now let us consider the second case. Let $J$ be the subinterval of
$c_0'$ with one end point being $x$ and with the length of $J$ being
one. We orient $J$ so that $x$ is its initial point. We set
$X_J=J\cap X$. We define $S_{X_J}(B)\subset S_{X_J}$ as follows. For
each $y\in X_J$ we let $f_B(y)$ be the minimum of $\delta'$ and the
distance along $D_y$ to the first point (excluding $y$) of $D_y$
contained in $\partial B$ and let $D_y(B)$ be the sub-geodesic of
$D_y$ of this length starting at $y$. Then $S_{X_J}(B)\subset
S_{X_J}$ is the union over $y\in X_J$ of $[0,f_B(y)]$.  Clearly, the
exponential mapping defines an immersion of $S_{X_J}(B)$ into $B$.
We need to replace $X_ J$ by a slightly smaller subset in order to
make the exponential mapping be an embedding. To do this we shall
remove bases of a maximal focusing triangles in $B$.

First notice that for each $y\in  X_J$ the exponential mapping is an
embedding on $D_y(B)$. The reason is that the image of $D_y(B)$ is a
geodesic contained in the ball $B$. Lemma~\ref{noloops} then shows
that this geodesic is embedded. This leads to:

\begin{claim}\label{locemb}
For any component $c$ of $X_J$, the restriction of the exponential
mapping to $S_c(B)=\left(c\times [0,\delta')\right)\cap S_{X_J}(B)$
is an embedding.
\end{claim}

\begin{proof}
Since the geodesics that make up $S_c(B)$ have length at most
$\delta'<1/10$ and since the curvature of the annulus is bounded
above by $1$, the restriction of the exponential mapping to $S_c(B)$
is a local diffeomorphism. The restriction  to each $\{y\}\times
[0,f_B(y)]$ is an embedding onto $D_y(B)$. If the restriction of the
exponential mapping to $S_c(B)$ is not an embedding, then there are
$y\not= y'$ in $c$ such that the geodesics $D_y(B)$ and $D_{y'}(B)$
meet. When they meet, they meet at a positive angle and by the
Gauss-Bonnet theorem this angle is less than $\mu+\delta$.  Thus,
all the geodesics starting at points sufficiently close to $y'$ and
between $y$ and $y'$ along $c$ must also meet $D_y(B)$. Of course,
if a sequence of $D_{y_i}(B)$ meet $D_y(B)$, then the same is true
for the limit. It now follows that $D_{y''}(B)$ meets $D_y(B)$ for
all $y''$ between $y$ and $y'$. This contradicts the fact that
$D_y(B)$ is embedded.
\end{proof}

\begin{claim}\label{limittri}
Any focusing triangle for $J$ must contain a component of $J\setminus X_J $. If
$\{{\mathcal T}_n\}$ is an infinite sequence of focus triangles for $J$, then,
after passing to a subsequence, there is a limiting focusing triangle for $J$.
\end{claim}

\begin{proof}
The first statement is immediate from Claim~\ref{locemb}.
 Since $X\cap J$ is compact, it is clear that after passing to a
 subsequence each of the sequence of left-hand sides and the
 sequence of right-hand sides  converge to a geodesic arc orthogonal to
 $J$ at points of $X$. Furthermore, these limiting geodesics meet in
 a point at distance at most $\delta'$ from the end of each. The only
 thing remaining to show is that the limiting left- and right-hand
 sides do not begin at the same point of $X$.
This is clear since each focusing triangle contains one of the finitely many
components of $J\setminus X_J$.
\end{proof}

Using Claim~\ref{limittri} we see that if there is  a  focusing
triangle for $J$ there is a first point $x_1$ in $ X_J$ whose
associated geodesic contains the left-hand side of a focusing
triangle for $J$. Then since the base length of any focusing
triangle is bounded by a fixed constant, invoking again
Claim~\ref{limittri}, that there is a focusing triangle ${\mathcal
T}_1$ for $J$ that has left-hand side contained in the geodesic
$D_{x_1}$ and has a maximal base among all such focusing triangles,
Maximal in the sense that the base of this focusing triangle
contains the base of any other focusing triangle with left-hand side
contained in $D_{x_1}$. Denote its base by $\xi_1$ and denote the
right-hand endpoint of $\xi_1$ by $y_1$. For the triangle we take
the geodesic arcs to the first point of intersection measured along
$D_{y_1}$. Set $J_1=\overline{J\setminus \xi}$, and repeat the
process for $J_1$. If there is a focusing triangle for $J_1$ we find
the first left-hand side of such focusing triangles and then find
the maximal focusing triangle ${\mathcal T}_2$ with this left-hand
side.

\begin{claim}
The interior of ${\mathcal T}_2$ is disjoint from the interior of
${\mathcal T}_1$.
\end{claim}

\begin{proof}
Since by construction the interiors of the bases of ${\mathcal T}_1$
and ${\mathcal T}_2$ are disjoint, if the interior of ${\mathcal
T}_2$ meets ${\mathcal T}_1$, then one of the sides of ${\mathcal
T}_2$ crosses the interior of one of the sides of ${\mathcal T}_1$.
But since ${\mathcal T}_1$ is a maximal focusing triangle with its
left-hand side, neither of the sides of ${\mathcal T}_2$ can cross
the interior of the left-hand side of ${\mathcal T}_1$. If one of
the sides of ${\mathcal T}_2$ crosses the interior of the right-hand
side of ${\mathcal T}_1$, then the right-hand side of ${\mathcal
T}_1$ is the left-hand side of a focusing triangle for $J_1$. Since
by construction the left-hand side of ${\mathcal T}_2$ is the first
such, this means that the left-hand side of ${\mathcal T}_2$ is the
right-hand side of ${\mathcal T}_1$. This means that the right-hand
side of ${\mathcal T}_2$ terminates when it meets the right-hand
side of ${\mathcal T}_1$ and hence the right-hand side of ${\mathcal
T}_2$ ends the first time that it meets the right-hand side of
${\mathcal T}_1$ and hence does not cross it.
\end{proof}

We continue in this way constructing focusing triangles for $J$ with
disjoint interiors. Since each focusing triangle for $J$ contains a
component of $J\setminus X_J$, and as there are only finitely many
such components, this process must terminate after a finite number
of steps. Let ${\mathcal T}_1,\ldots,{\mathcal T}_k$ be the focusing
triangles so constructed, and denote by $\xi_i$ the base of
${\mathcal T}_i$. Let $X_J'= X_J\setminus \cup_{i=1}^k\xi_i$.

\begin{defn}
We call the triangles ${\mathcal T}_1,\ldots,{\mathcal T}_k$
constructed above, the {\em maximal set of focusing triangles} for
$J$ relative to $B$.
\end{defn}

\begin{claim}
The length of $X'_J$ is at least $1-2\delta-\mu$.
\end{claim}

\begin{proof}
Since the interiors of the ${\mathcal T}_i$ are disjoint, according
to Corollary~\ref{muepsbd}, we have $\sum_il(\xi_i)<\delta+\mu$. We
also know by Remark~\ref{Ylengthrem} that the length of $X_J$ is at
least $(1-\delta)$. Putting these together gives the result.
\end{proof}

We define $S_{X'_J}(B)$ to be the intersection of $S_{X_J}(B)$ with
$S_{X'_J}$.

\begin{claim}
The restriction of the exponential mapping to $S_{X'_J}(B)$ is an
embedding.
\end{claim}

\begin{proof}
Suppose that we have distinct points $x',y'$ in $X'_J$ such that
$D_{x'}(B)\cap D_{y'}(B)\not=\emptyset$. We assume that $x'<y'$ in
the orientation on $J$. Then there is a focusing triangle for $J$
whose base is the sub-arc of $J$ with endpoints $x'$ and $y'$, and
hence the left-hand side of the focusing triangle is contained in
$D_{x'}(B)$. Since $x'$ is not a point of $\cup_{i}\xi_i$ either it
lies between two of them, say $\xi_j$ and $\xi_{j+1}$ or it lies
between the initial point $x$ of $J$ and the initial point of
$\xi_1$ or it lies between the last $\xi_n$ and the final point of
$J$.

But $x'$ cannot lie before $\xi_1$, for this would contradict the
construction which took as the left-hand endpoint of $\xi_1$ the
first point of $J$ whose geodesic contained the left-hand side of a
focusing triangle for $J$. Similarly, $x'$ cannot lie between
$\xi_j$ and $\xi_{j+1}$ for any $j$ since the left-hand endpoint of
$\xi_{j+1}$ is the first point at or after the right-hand endpoint
of $\xi_j$ whose geodesic contains the left-hand side of a focusing
triangle for $J$. Lastly, $x'$ cannot lie to the right of the last
$\xi_k$, for then we would not have finished the inductive
construction.
\end{proof}

We pull back the metric of $A$ to the $S_{X'_J}(B)$ by the
exponential mapping. Since this pullback metric is at least
$(1-\delta)^2$ times the product of the metric on $X'_J$ induced
from $c_0$ and the usual metric on the interval, and since the map
on this subset is an embedding, we see that the area of the region
of the annulus which is the image under the exponential mapping of
this subset is at least
$$(1-\delta)^2\int_{X'_J}f_B(x)ds,$$
where $s$ is arc length along $X'_J$. Of course, the area of this
subset is at most $\mu$. This means that, setting $Z$ equal to the
subset of $ X'_J(\delta')$ given by
$$Z=\{z\in X'_J|f_B(z)<\delta'\},$$ the total length of $X'_J\setminus Z$
satisfies
$$l(X'_J\setminus Z)\le (1-\delta)^{-2}(\delta')^{-1}\mu<\frac{1}{10},$$
where the last inequality is an immediate consequence of our choice
of $\mu$. Thus, the length of $Z$ is at least $(0.9-2\delta-\mu)\ge
0.87$. Let $\partial_+S_Z(B)$ be the union of the final endpoints
(as opposed to the initial points) of the $D_x(B)$ as $x$ ranges
over $Z$.
 Of course, since $f_B(z)<\delta'$
for all $z\in Z$, it must be the case that  the exponential mapping
embeds $\partial_+S_Z(B)$ into $\partial B$. Furthermore, the total
length of the image of $\partial_+S_Z(B)$ is at least
$(1-\delta)l(Z)\ge 0.86$. The boundary of $B$ is made up of two
pieces: $D_x$ and an arc on $c_0$. But the length of $D_x$ is at
most $\delta'<1/20$ so that not all of $\partial_+S_Z(B)$ can be
contained in $D_x$. Thus, there is $z\in Z$, distinct from $x$ such
that $D_z$ has both endpoints in $c_0$. It then follows  that all
points of $Z$ that are separated (along $J$) from $x$ by $z$ have
the same property. Since the length of $Z$ is at least $0.86$, it
follows that there is a point $z\in X$ at least distance $0.85$
along $J$ from $x$ with the property that $D_z$ has both endpoints
in $c_0$. The complementary component of $D_z$ in $A$, denoted $B'$,
 is a disk that is contained in $B$ and the length of $B'\cap c_0$ is at
least $0.85$ less than the length of $B\cap c_0$.

The length of $B'\cap c_0$ cannot be less than $1$, for that gives a
contradiction. But if the length of $B'\cap c_0$ is greater than one, we now
repeat this construction replacing $B$ by $B'$. Continuing in this way we
eventually we cut down the length of $B\cap c_0$ to be less than one and hence
reach a contradiction.
\end{proof}

\subsection{For every $x\in X$, the geodesic $D_x$ is embedded}

The steps in the above argument inductively constructing disjoint
maximal focusing triangles and showing that their bases have a small
 length and that off of them the map is an embedding will be
repeated in two other contexts. The next context is to rule out the
case when a sub-arc of $D_x$  forms a homotopically non-trivial loop
in $A$.

\begin{lem}
For any $x\in X$ there is no sub-geodesic of $D_x$ that is an
embedded loop in $A$.
\end{lem}

\begin{proof}
We have already treated the case when the loop bounds a disk. Now we
need to treat the case when the loop is homotopically non-trivial in
$A$. Let $D'_x\subset D_x$ be the minimal compact sub-geodesic
containing $x$ that is not an embedded arc. Let ${\rm int}\,B$ be
the complementary component of $D'_x$ in $A$ that contains
$c_0\setminus\{x\}$. There is a natural compactification of ${\rm
int}\, B$ as a disk and an immersion of this disk into $A$, an
immersion that is two-to-one along the shortest sub-geodesic of
$D'_x$ from $x$ to the point of intersection of $D_x$ with itself.
We do exactly the same construction as before. Take a sub-arc $J$ of
length one with $x$ as an endpoint and construct $S_{X_J}(B)$
consisting of the union of the sub-geodesics of $D_z$, for $z\in
J\cap X$ that do not cross the boundary of $B$. We then construct a
sequence of maximal focusing triangles along $J$ relative to $B$
just as in the previous case. In this way we construct a subset $Z$
of $X\cap J$ of total length at least $0.87$ with the property that
for every $z\in Z$ the final end of $D'_z(B)$ lies in $\partial B$.
Furthermore, the length of the arcs that these final ends sweep out
is at least $0.86$. Hence, since the total length of the part of the
boundary of $B$ coming from $D_x'$ is at most $2\delta'<0.2$, there
must be a $z\in Z$ for which $D_z(B)$ has both ends in $c_0$. This
puts us back in the case ruled out in Lemma~\ref{Dxacross}.
\end{proof}

\subsection{Far apart $D_x$'s don't meet}

Now the last thing that can prevent the exponential mapping in the
complement of the focusing triangles from being an embedding is that
geodesics $D_x$ and $D_y$ might meet even though $x$ and $y$ are far
apart along $c_0$. Our next goal is to rule this out.

\begin{lem}
Let $x,y$ be distinct points of $X$. Suppose that there are
sub-geodesics $D'_x\subset D_x$ and $D'_y\subset D_y$ with a common
endpoint. Then the arc $D'_x*(D'_y)^{-1}$ cuts $A$ into two
complementary components, exactly one of which is a disk, denoted
$B$. Then it is not possible for  $B\cap c_0$ to contain an arc of
length $1$.
\end{lem}

\begin{proof}
The proof is exactly the same as in Lemma~\ref{Dxacross} except that
the part of the boundary of $B$ that one wants to avoid has length
at most $2\delta'<0.2$ instead of $\delta'$. Still, since (in the
notation of the proof of Lemma~\ref{Dxacross}) the total length of
$Z$ is at least $0.87$ so that the lengths of the other ends of the
$D_z$ as $z$ ranges over $Z$ is at least $0.86$, there is $z\in Z$
for which both ends of $D_z$ lie in $c_0$. Again this puts us back
in the case ruled out by Lemma~\ref{Dxacross}.
\end{proof}

As a special case of this result we have the following.

\begin{cor}\label{noholes}
Suppose that we have an arc $\xi$ of length at most $1$ on $c_0$.
Denote the endpoints of $\xi$ by $x$ and $y$ and  suppose that
$D_x\cap D_y\not=\emptyset$. Let $D'_x$ and $D'_y$ be sub-geodesics
containing $x$ and $y$ respectively ending at the same point, $v$,
and otherwise disjoint. Then the loop $\xi*D'_y*(D'_x)^{-1}$ bounds
a disk in $A$.
\end{cor}

\begin{proof}
If not, then it is homotopically non-trivial in $A$ and replacing
$\xi$ by its complement, $c_0\setminus{\rm int}\,\xi$, gives us
exactly the situation of the previous lemma. (The length of
$c_0\setminus{\rm int}\,\xi$ is at least one since the length of
$c_0$ is at least $2$.)
\end{proof}

Let us now summarize what we have established so far about the intersections of
the geodesics $\{D_x\}_{x\in X}$.

\begin{cor}\label{focussummary}
For each $x\in X$, the geodesic $D_x$ is an embedded arc in $A$.
Either it has length $\delta'$ or its final point lies on $c_1$.
Suppose there are   $x\not=x'$ in $X$ with $D_x\cap
D_{x'}\not=\emptyset$. Then there is an arc $\xi$ on $c_0$
connecting $x$ to $x'$ with the length of $\xi$ at most
$\delta+\mu$. Furthermore, for sub-geodesics $D'_x\subset D_x$,
containing $x$, and $D'_{x'}\subset D_{x'}$, containing $x'$, that
intersect exactly in an endpoint of each, the loop
$\xi*D'_{x'}*(D'_x)^{-1}$ bounds a disk $B$ in $A$, and the length
of $\xi$ is at most the turning of $\xi$ plus the area of $B$.
\end{cor}

\subsection{Completion of the proof}

We have now completed all the technical work on focusing and we have
also shown that  the restriction of the exponential mapping to the
complement of the bases of the focusing regions is an embedding. We
are now ready to complete the proof of
Proposition~\ref{annlengthest}.

Let $J$ be an interval of length one in $c_0$. Because of
Corollary~\ref{focussummary} we can  construct the maximal focusing
triangles for $J$ as follows. Orient $J$, and begin at the initial
point of $J$. At each step we consider the first $x$ (in the
subinterval of $J$ under consideration) which intersects a $D_y$ for
some later $y\in J$. If we have such $y$, then we can construct the
sides of the putative triangle for sub-geodesics of $D_x$ and $D_y$.
But we need to know that we have a focusing triangle. This is the
content of Corollary~\ref{noholes}. The same reasoning works when we
construct the maximal such focusing triangle with a given left-hand
side, and then when we show that in the complement of the focusing
triangles the map is an embedding. Thus, as before, for an interval
$J$  of length $1$, we construct a subset $X_J'\subset X\cap J$ of
length at least $0.97$ such that the restriction of the exponential
mapping to $S_{X'_J}$ is an embedding. Again the area estimate shows
that there is a subset $Z\subset X'_J$ whose length is at least
$0.87$ with the property that for every $z\in Z$ the geodesic $D_z$
has both endpoints in $\partial A$. By Lemma~\ref{Dxacross}, the
only possibility for the final endpoints of all these $D_z$'s is
that they lie in $c_1$.

In particular, there are $x\in X$ for which $D_x$ spans from $c_0$
to $c_1$. We pick one such, $x_0$, contained in the interior of $X$,
and use it as the starting point for a construction of maximal
focusing triangles all the way around $c_0$. What we are doing at
this point actually is cutting the annulus open along $D_{x_0}$ to
obtain a disk and we construct a maximal family of focusing
triangles of the interval $[x'_0,x_0'']$ obtained by cutting $c_0$
open at $x_0$ relative to this disk. Here $x_0'$ and $x_0''$ are the
points of the disk that map to $x_0$ when the disk is identified to
form $A$. Briefly, having constructed a maximal collection of
focusing triangles for a subinterval $[x'_0,x]$, we consider the
first point $y$ in the complementary interval $[x,x''_0]$ with the
property that there is $y'$ in this same interval, further along
with $D_y\cap D_{y'}\not=0$. Then, using
Corollary~\ref{focussummary} we construct the maximal focusing
triangle on $[x,x''_0]$ with left-hand side being a sub-geodesic of
$D_y$. We then continue the construction inductively until we reach
$x''_0$. Denote by $\xi_1,\ldots,\xi_k$ the bases of these focusing
triangles and let $X'$ be $X\setminus\cup_i\xi_i$.

The arguments above show that the exponential mapping is an
embedding of $S_{X'}$ to the annulus.

\begin{claim}
For every subinterval $J$ of length one in $c_0$ the total length of
the bases $\xi_i$ that meet $J$ is at most $2\delta+\mu<0.03$.
\end{claim}

\begin{proof}
Since, by Corollary~\ref{muepsbd}, every base of a focusing triangle
has length at most $\delta+\mu$, we see that the union of the bases
of focusing triangles meeting $J$ is contained in an interval of
length $1+2(\delta+\mu)<2$. Hence, the total turning of the bases of
these focusing triangles is at most $2\delta$ whereas the sum of
their areas is at most $\mu$. The result now follows from
Corollary~\ref{muepsbd}.
\end{proof}

By hypothesis there is an integer $n\ge 1$ such that the length
$l(c_0)$ of $c_0$ is greater than $n$ but less than or equal to
$n+1$. Then it follows from the above that the total length of the
bases of all the focusing triangles in our family is at most
$$(n+1)(2\delta+\mu)< 0.03(n+1)\le 0.06n\le  0.06 l(c_0).$$

Since the restriction of the exponential mapping to $S_{X'}$ is an
embedding, it follows from Claim~\ref{expcomp} and the choice of
$\delta'$ that, for any open subset $Z$ of $X'$, the area of the
image under the exponential mapping of $S_Z$ is at least
$(1-\delta)^2\int_Zf(x)ds$, where $ds$ is the arc length along $Z$.
Also, the image under the exponential mapping of $\partial_+(S_Z)$
is an embedded arc in $A$ of length at least $(1-\delta)l(Z)$. Since
the length of $X'$ is at least $(0.94)l(c_0)$ and since the area of
$A$ is less than $\mu<(1-\delta)^2\delta'/10$, it follows that the
subset of $X'$ on which $f$ takes the value $\delta'$ has length at
most $0.10<(0.10)l(c_0)$. Hence, there is a subset $X''\subset X'$
of total length at least $(0.84)l(c_0)$ with the property that
$f(x)<\delta'$ for all $x\in X''$. This means that for every $x\in
X''$ the geodesic $D_x$ spans from $c_0$ to $c_1$, and hence the
exponential mapping embeds $\partial_+S_{X''}$ into $c_1$. But we
have just seen that the length of the image under the exponential
mapping of $\partial_+S_{X''}$  is at least
$$(1-\delta)l(X'')>(0.99)l(X'')>(0.83)l(c_0).$$
It follows  that the length of $c_1$ is at least
$(0.83)l(c_0)>3(l(c_0))/4$.

This completes the proof.

\section{Proof of the first inequality in Lemma~\protect{\ref{CSSHI}}}\label{18.6}

Here is the statement that we wish to establish when the manifold
$(W,h(t))$ is the product of $(M,g(t))\times (S^1_\lambda,ds^2)$.

\begin{lem}\label{csshi1}
Let $(W,h(t)),\ t_0\le t\le t_1$, be a Ricci flow and fix
$\Theta<\infty$. Then there exist constants $\delta>0$ and $0<r_0\le
1$   depending only on the curvature bound for the ambient Ricci
flow and $\Theta$ such that the following holds. Let $c(x,t),\
t_0\le t\le t_1$, be a curve-shrinking flow with $c(\cdot,t)$
immersed for each $t\in [t_0,t_1]$ and with the total curvature of
$c(\cdot,t_0)$ being at most $\Theta$. Suppose that there is $0<r\le
r_0$ and at a time $t'\in [t_0,t_1-\delta r^2]$ such that the length
of $c(\cdot,t')$ is at least $r$ and the total curvature of
$c(\cdot,t')$ on any sub-arc of length $r$ is at most $\delta$. Then
for every $t\in [t',t'+\delta r^2]$ the curvature $k$ satisfies
$$k^2\le  \frac{2}{(t-t')}.$$
\end{lem}

The rest of this section is devoted to the proof of this lemma. In
\cite{AG} such a local estimate was established when the ambient
manifold was Euclidean space and the curve in question is a graph. A
related result for hypersurfaces that are graphs appears in
\cite{EckerHuisken}. The passage from Euclidean space to a general
Ricci flow is straightforward, but it is more delicate to use the
bound on total curvature on initial sub-arcs of length $r$ to show
that in appropriate coordinates the evolving curve can be written as
an evolving graph, so that the analysis in [2] can be applied.

We fix $\delta>0$ sufficiently small. We fix $t'\in [t_0,t_1-\delta
r^2]$ for which the hypotheses of the lemma hold.  The strategy of
the proof is to first restrict to the maximum subinterval of
$[t',t_2]$ of $[t',t'+\delta r^2]$ on which $k$ is bounded by
$\sqrt{2/(t-t')}$. If $t_2<t'+\delta r^2$, then $k$ achieves the
bound $\sqrt{2/(t-t')}$ at time $t_2$. We show that in fact on this
subinterval $k$ never achieves the bound. The result then follows.
To show that $k$ never achieves the bound, we show that on a
possibly smaller interval of time $[t',t_3]$ with $t_3\le t_2$ we
can write the restriction of the curve-shrinking flow to any
interval whose length at time $t'$ is $(0.9)r$ as a family of graphs
in a local coordinate system so that the function $f$ (of arc and
time) defining the graph has derivative along the arc bounded in
norm by $1/2$. We take $t_3\le t_2$ maximal with respect to these
conditions. Then with both the bound on $k$ and the bound on the
derivative of $f$ one shows that the spatial derivative of $f$ never
reaches $1/2$ and also that the curves do not move too much so that
they always remain in the coordinate patch. The only way that this
can happen is that if $t_2=t_3$, that is to say, on the entire time
interval where we have the curvature bound, we also can write the
curve-shrinking flow as a flow of graphs with small spatial
derivatives. Then it is convenient to replace the curve-shrinking
flow by an equivalent flow, introduced in \cite{AG}, called the
graph flow. Applying a simple maximum principle argument to this
flow we see that $k$ never achieves the value $\sqrt{2/(t-t')}$ on
the time interval $[t',t_2]$ and hence the curvature estimate
$k<\sqrt{2/(t-t')}$ holds throughout the interval $(t',t'+\delta
r^2]$.

\subsection{A bound for $\int kds$}

Recall that $k$ is the norm of the curvature vector $\nabla_SS$, and
in particular, $k\ge 0$. For any sub-arc $\gamma_{t'}$ of
$c(\cdot,t')$  at time $t'$ we let $\gamma_t$ be the result at time
$t$ of applying the curve-shrinking flow to $\gamma_{t'}$. The
purpose of this subsection is to show that $\int_{\gamma_t} kds$ is
small for all $t\in [t',t'+\delta r^2]$ and all initial arcs
$\gamma_{t'}$ of length at most $r$.

\begin{claim}\label{C0claim}
There is a constant $D_0<\infty$, depending only on $\Theta$ and the
curvature bound of the ambient Ricci flow such that for every $t\in
[t',t'+\delta r^2]$ and every sub-arc $\gamma_{t'}$ whose length is
at most $r$, we have $\int_{\gamma_t} kds<D_0$ and $l(\gamma_t)\le
D_0r$, where $l(\gamma_t)$ is the length of $\gamma_t$.
\end{claim}

\begin{proof}
This is immediate from Corollary~\ref{lncurvgrowth} applied to all
of $c(\cdot,t)$.
\end{proof}

Now we fix $t_2\le t'+\delta r^2$ maximal subject to the condition
that $k(x,t)\le \sqrt{\frac{2}{t-t'}}$ for all $x$ and all $t\in
[t',t_2]$. If $t_2<t'+\delta r^2$ then there is $x$ with
$k(x,t_2)=\sqrt{\frac{2}{(t_2-t')}}$.

Now consider a curve $\gamma_{t'}$ of length $r$. From the integral estimate in
the previous claim and the assumed pointwise estimate on $k$, we see that
$$\int_{\gamma_t}k^2 ds\le {\rm max}_{x\in\gamma_t}k(x,t)\int_{\gamma_{t}} kds
<\sqrt{\frac{2}{t-t'}}\cdot D_0.$$ Using Equation~(\ref{dLdteqn}),
it follows easily that, provided that $\delta>0$ is sufficiently
small, the length of $\gamma_t$ is at least $(0.9)r$ for all $t\in
[t',t_2]$, and more generally for any subinterval $\gamma'_{t'}$ of
$\gamma_{t'}$ and for any $t\in [t',t_2]$ the length of the
corresponding interval $\gamma'_t$ is at least $(0.9)$ times the
length of $\gamma'_{t'}$. We introduce a cut-off function on
$\gamma_{t'}\times [t',t_2]$ as follows. First, fix a smooth
function $\psi\colon [-1/2,1/2]\to [0,1]$ which is identically zero
on $[-0.50,-0.45]$ and on $[0.45,0.50]$, and is identically $1$ on
$[-3/8,3/8]$. There is a constant $D'$ such that $|\psi'|\le D'$ and
$|\psi''|\le D'$. Now we fix the midpoint $x_0\in \gamma_{t'}$ and
define the signed distance from $(x_0,t)$, denoted
$$s\colon \gamma_{t'}\times [t',t_2]\to \Ar,$$
as follows:
$$s(x,t)=\int _{x_0}^x|X(y,t)|dy.$$

We define the cut-off function
$$\varphi(x,t)=\psi\left(\frac{s(x,t)}{r}\right).$$

\begin{claim}\label{1887}
There is a constant $D_1$ depending only on the curvature bound for
 the ambient Ricci flow such that for any sub-arc $\gamma_{t'}$ of
 length $r$, defining $\varphi(x,t)$ as above, for all $x\in \gamma_t$
 and all $t\in [t',t_2]$ we have
$$\left|\frac{\partial \varphi(x,t)}{\partial t}\right|\le
\frac{D_1}{r\sqrt{t-t'}}+D_1.$$
\end{claim}

\begin{proof}
Clearly,
$$\frac{\partial \varphi(x,t)}{\partial
t}=\psi'\left(\frac{s(x,t)}{r}\right)\cdot \frac{1}{r}\frac{\partial
s(x,t)}{\partial t}.$$ We know that $|\psi'|\le D'$ so that
$$\left|\frac{\partial \varphi(x,t)}{\partial
t}\right|\le \frac{D'}{r}\left|\frac{\partial s(x,t)}{\partial
t}\right|.$$ On the other hand,
$$s(x,t)=\int_{x_0}^x|X(y,t)|dy,$$
so that  $$ \left|\frac{\partial s(x,t)}{\partial t}\right| =
\left|\int_{x_0}^x\frac{\partial |X(y,t)|}{\partial t}dy\right|,$$
By Lemma~\ref{1stcurvshr} we have $$\frac{\partial
|X(y,t)|}{\partial y}dy  = \left(-{\rm
Ric}(S(y,t),S(y,t))-k^2(y,t)\right)ds,$$
 so that there is a constant $D$ depending only on the bound of the
sectional curvatures of the ambient Ricci flow with
\begin{eqnarray*}\left|\frac{\partial s(x,t)}{\partial t}\right| &
\le &  \int_{x_0}^x(D+k^2)ds\le
Dl(\gamma_t)+\int_{x_0}^xk^2(y,t)ds(y,t), \end{eqnarray*} and hence
by Claim~\ref{C0claim}
$$\left|\frac{\partial s(x,t)}{\partial t}\right|\le
  DD_0r+
\int_{x_0}^xk^2(y,t)ds(y,t).$$
  Using
the fact that $k^2\le 2/(t-t')$,  we have
$$\int_{x_0}^xk^2(y,t)ds(y,t)\le\sqrt{\frac{2}{t-t'}}\int_{x_0}^xkds\le
\frac{\sqrt{2}D_0}{\sqrt{t-t'}}.$$ Putting all this together, we see
that there is a constant $D_1$ such that
$$\left|\frac{\partial \varphi(x,t)}{\partial t}\right|\le
D_1(\frac{1}{r\sqrt{t-t'}}+1).$$
\end{proof}

\begin{claim}\label{claim1885}
There is a constant $D_2$ depending only on the curvature bound of
the ambient Ricci flow and $\Theta$ and a constant $D_3$ depending
only on the curvature bound of the ambient Ricci flow, such that for
any $t\in [t',t_2]$ and any sub-arc $\gamma_{t'}$ of length $r$, we
have
$$\left|\frac{d}{dt}\int_{\gamma_t}\varphi
kds\right|\le
D_2\left(1+\frac{1}{r\sqrt{t-t'}}\right)+\frac{D_2}{r^2}+D_3\int_{\gamma_t}\varphi
kds.$$
\end{claim}

\begin{proof}
We have
$$\left|\frac{d}{dt}\int_{\gamma_t}\varphi
kds\right|\le\left|\int_{\gamma_t}\frac{\partial
\varphi(x,t)}{\partial
t}kds\right|+\left|\int_{\gamma_t}\varphi\frac{\partial
(kds)}{\partial t}\right|.$$ Using Claim~\ref{1887} for the first
term and Claim~\ref{3rdcurvshr} and arguing as in the proof of
Lemma~\ref{1835} for the second term, we have
$$\left|\frac{d}{dt}\int_{\gamma_t}\varphi
kds\right|\le
D_1\left(1+\frac{1}{r\sqrt{t-t'}}\right)\int_{\gamma_t}kds+\left|\int_{\gamma_t}
\varphi k''ds\right|+\int_{\gamma_t}C'_1\varphi kds,$$ where $C_1'$
depends only on the ambient curvature bound. We bound the first term
by
$$D_1D_0\left(\frac{1}{r\sqrt{t-t'}}+1\right),$$
where $D_0$ is the constant depending on $\Theta$ and the ambient
curvature bound from Claim~\ref{C0claim}. Since the ends of
$\gamma_t$ are at distance at least $(0.45)r$ from $x_0$ all $t\in
[t',t_2]$, we see that for all $t\in [t',t_2]$
$$\int_{\gamma_t}\varphi k''=\int_{c(\cdot,t)}\varphi k''.$$
Integrating by parts we have
$$\int_{c(\cdot,t)} \varphi k'' ds=\int_{c(\cdot,t)}\varphi'' k ds,$$
where the prime here refers to the derivative along $c(\cdot,t)$
with respect to arc length. Of course $|\varphi''|\le
\frac{D'}{r^2}$. Thus, we see that
$$\left|\int_{\gamma_t}\varphi k''ds\right|\le
\frac{D'}{r^2}\int_{c(\cdot,t)}kds\le \frac{D'D_0}{r^2}.$$ Putting
all this together, we have
$$\left|\frac{d}{dt}\int_{\gamma_t}\varphi
kds\right|\le
D_2\left(1+\frac{1}{r\sqrt{t-t'}}\right)+\frac{D_2}{r^2}+D_3\int_{\gamma_t}\varphi
kds$$ for $D_2=D_0\,{\rm max}(D',D_1)$ and $D_3=C_1'$. This gives
the required estimate.
\end{proof}

\begin{cor}
For any $t\in [t',t_2]$  and any sub-arc $\gamma_{t'}$ of length $r$
we have
$$\int_{\gamma_t}\varphi kds
\le  D_4\sqrt{\delta}$$ for a constant $D_4$ that depends only on
the sectional curvature bound of the ambient Ricci flow and
$\Theta$.
\end{cor}

\begin{proof}
This is immediate from the previous result by integrating from $t'$
to $t_2\le t'+\delta r^2$, and using the fact that $\delta<1$ and
$r<1$, and using the fact that $$\int_{\gamma_{t'}}\varphi kds\le
\int_{\gamma_{t'}}kds<\delta$$ since $\gamma_{t'}$ has length at
most $r$.
\end{proof}

This gives:

\begin{cor}\label{D_4cor}
For $\gamma_{t'}\subset c(\cdot,t')$ a sub-arc of length at most $r$
 and for any $t\in [t',t_2]$, we have
$$\int_{\gamma_t}kds\le
2D_4\sqrt{\delta}.$$ For any $t\in [t',t_2]$ and any sub-arc
$J\subset c(\cdot,t)$ of length at most $r/2$ with respect to the
metric $h(t)$, we have
$$\int_Jk(x,t)ds(x,t)\le 2D_4\sqrt{\delta}.$$
\end{cor}

\begin{proof}
We divide an  interval $\gamma_{t'}\subset c(\cdot,t')$ of length at
most $r$ into two subintervals $\gamma'_{t'}$ and $\gamma''_{t'}$ of
lengths at most $r/2$. Let $\hat \gamma'_{t'}$ and $\hat
\gamma''_{t'}$ be intervals of length $r$ containing $\gamma'_{t'}$
and $\gamma''_{t'}$ respectively as  middle subintervals. We then
apply the previous corollary to  $\hat \gamma'_{t'}$ and $\hat
\gamma''_{t'}$ using the fact that $\varphi k\ge 0$ everywhere and
$\varphi k=k$ on the middle subintervals of $\hat \gamma'_{t'}$ and
$\hat \gamma ''_{t'}$. For an interval $J\subset \gamma_t$ of length
$r/2$, according to Lemma~\ref{lncurvgrowth} the length of
$\gamma_{t'}|_J$ with respect to the metric $h(t')$ is at most $r$,
and hence this case follows from the previous case.
\end{proof}

\subsection{Writing the curve flow as a graph}

Now we restrict attention to $[t',t_2]$, the maximal interval in
$[t',t'+\delta r^2]$ where $k^2\le 2/(t-t')$. Let $\gamma_{t'}$ be
an arc of length $r$ in $c(\cdot,t')$ and let $x_0$ be the central
point of $\gamma_{t'}$. Denote $\gamma_{t'}(x_0)=p\in W$. We take
the $h(t')$-exponential mapping from $T_pW\to W$. This map will be a
local diffeomorphism out to a distance determined by the curvature
of $h(t')$.   For an appropriate choice of the ball (depending on
the ambient curvature bound) the metric on the ball induced by
pulling back $h(t)$ for all $t\in [t',t_2]$ will be within $\delta$
in the $C^1$-topology to the Euclidean metric $h'=h(t')_p$. By this
we mean that
\begin{enumerate}
\item[(1)] $\left|\langle X,Y\rangle_{h(t)}-\langle
X,Y\rangle_{h'}\right|<\delta |X|_{h'}|Y|_{h'}$ for all tangent vectors in the
coordinate system, and \item[(2)] viewing the connection $\Gamma$ as a bilinear
map on the coordinate space with values in the coordinate space we have $
|\Gamma(X,Y)|_{h'}<\delta |X|_{h'}|Y|_{h'}$.
\end{enumerate}

We choose $0<r_0\le 1$ so that it is much smaller than this
distance, and hence  $r$ is also much smaller than this distance. We
lift to the ball in $T_pW$.

 We fix orthonormal coordinates with respect to the metric $h'$ so
that the tangent vector of $\gamma_{t'}(x_0)$ points in the positive
$x^1$-direction. Using these coordinates we decompose the coordinate
patch as a product of an interval in the $x^1$-direction and an open
ball, $B$, spanned by the remaining Euclidean coordinates. From now
on we shall work in this coordinate system using this product
structure.  To simplify the notation in the coming computations, we
rename the $x^1$-coordinate the $z$-coordinate. Ordinary derivatives
of a function $\alpha$ with respect to $z$ are written $\alpha_z$.
When we write norms and inner products without indicating the metric
we implicitly mean that the metric is $h(t)$. When we use the
Euclidean metric on these coordinates we denote it explicitly. Next,
we wish to understand how $\gamma_t$ moves in the Euclidean
coordinates under the curve-shrinking flow. Since we have
$|\nabla_SS|_{h}=k$, it follows that $|\nabla_SS|_{h'}\le
\sqrt{1+\delta}k\le 2/\sqrt{t-t'}$, and hence, integrating tells us
that for any $x\in \gamma_{t'}$ we have
$$\left|\gamma_t(x)-\gamma_{t'}(x)\right|_{h'}\le 4\sqrt{t-t'}\le 4\sqrt{\delta}r.$$
This shows that for every $t\in [t',t_2]$, the curve $\gamma_t$ is
contained in the coordinate patch that we are considering. This
computation also implies that the $z$-coordinate of $\gamma_t$
changes by at most $4\sqrt{\delta} r$ over this time interval.

Because the total curvature of $\gamma_{t'}$ is small and the metric
is close to the Euclidean metric, it follows that the tangent vector
at every point of $\gamma_{t'}$ is close to the positive
$z$-direction. This means that we can write $\gamma_{t'}$ as a graph
of a function $f$ from a subinterval in the $z$-line to $Y$ with
$|f_z|_{h'}<2\delta$. By continuity, there is $t_3\in (t',t_2]$ such
that all the curves $\gamma_{t}$ are written as graphs of functions
(over subintervals of the $z$-axis that depend on $t$) with
$|f_z|_{h'}\le 1/10$. That is to say, we have an open subset $U$ of
the product of the $z$-axis with $[t',t_3]$, and the evolving curves
define a map $\tilde\gamma$ from $U$ into the coordinate system,
where the slices at constant time are graphs $z\mapsto (z,f(z,t))$
and are the curves $\gamma_t$. Using the coordinates $(z,t)$ gives a
new flow of curves by moving in the $t$-direction. This new flow is
called the {\sl graph-flow}. It is
 a reparameterization of the curve shrinking flow in such a way that
the $z$-coordinate is preserved. We denote by $Z=Z(z,t)$ the image
under the differential of the map $\widetilde\gamma$ of the tangent
vector in the $z$-direction and by $Y(z,t)$ the image under the
differential of $\widetilde \gamma$ of the tangent vector in the
$t$-direction. Notice that $Z$ is the tangent vector along the
parameterized curves in the graph flow. Since we are now using a
different parameterization of the curves from the one determined by
the curve-shrinking flow, the tangent vector $Z$ has the same
direction but not necessarily the same length as the tangent vector
$X$ from the curve-shrinking parameterization. Also, notice that in
the Euclidean norm we have $|Z|_{h'}^2=1+|f_z|^2_{h'}$. It follows
that on $U$ we have
$$(1-\delta)(1+|f_z|_{h'}^2)\le |Z(z,t)|^2_{h(t)}\le (1+\delta)(1+|f_z|_{h'}^2).$$
In particular, because of our restriction to the subset where $|f_z|_{h'}\le
1/10$ we have  $(1-\delta)\le |Z(z,t)|^2_{h(t)}\le (1.01)(1+\delta)$.

Now we know that $\gamma_{t'}$ is a graph of a function $f(z,t')$
defined on some  interval $I$ along the $z$-axis.  Let $I'$ be the
subinterval of $I$ centered in $I$ with $h'$-length $(0.9)$ times
the $h'$-length of $I$. By the above estimate on $|Z|$ it follows
that  the restriction of $\gamma_{t'}$ to $I'$ has length between
$(0.8)r$ and $r$, and also that the $h'$-length of $I'$ is between
$(0.8)r$ and $r$. The above estimate means that, provided that
$\delta>0$ is sufficiently small, for every $t\in [t',t_3]$ there is
a subinterval of $\gamma_t$ that is the graph of a function defined
on all of $I'$. We now restrict attention to the family of curves
parameterized by $I'\times [t',t_3]$. For every $t\in[t',t_3]$ the
curve $\gamma_t|_{I'}$ has length between $(0.8)r$ and $r$. The
curve-shrinking flow is not defined on this product because under
the curve-shrinking flow the $z$-coordinate of any given point is
not constant. But the graph flow defined above, and studied in
\cite{AG} (in the case of Euclidean background metric), is defined
on $I'\times [t',t_3]$ since this flow preserves the $z$-coordinate.
The time partial derivative in the curve-shrinking flow is given by
\begin{equation}\label{1814}
\nabla_SS=\frac{\nabla_ZZ}{|Z|^2}-\frac{1}{|Z|^4}\langle \nabla_ZZ,Z\rangle Z.
\end{equation} The time partial derivative in the graph-flow is given by
$Y=\partial\widetilde \gamma/\partial t$. The tangent vector $Y$ is
characterized by being $h'$-orthogonal to the $z$-axis and differing from
$\nabla_SS$ by a functional multiple of $Z$.

\begin{claim}\label{curvflowcompare}
$$Y=\frac{\nabla_ZZ-\langle\Gamma(Z,Z),\partial_z\rangle_{h'} Z}{|Z|^2}
=\nabla_SS+\left(\frac{\langle \nabla_ZZ,Z\rangle}{|Z|^4}
-\frac{\langle\Gamma(Z,Z),\partial_z\rangle_{h'} }{|Z|^2}\right)Z.$$
\end{claim}

\begin{proof}
In our Euclidean coordinates, $Z=(1,f_z)$ so that
$\nabla_ZZ=(0,f_{zz})+\Gamma(Z,Z)$. Thus,
$$\langle \nabla_ZZ,\partial_z\rangle_{h'}=\langle
\Gamma(Z,Z),\partial_z\rangle_{h'}.$$ Since $\langle
Z,\partial_z\rangle_{h'}=1$, it follows that
$$\frac{ \nabla_ZZ-\langle\Gamma(Z,Z),\partial_z\rangle_{h'} Z}{|Z|^2}$$ is
$h'$-orthogonal to the $z$-axis and hence is a multiple of $Y$. Since it
differs by a multiple of $Z$ from $\nabla_SS$, it follows that it is $Y$. This
gives the first equation; the second follows from this and
Equation~(\ref{1814}).
\end{proof}

To simplify the notation we set
$$\psi(Z)=\frac{\langle\Gamma(Z,Z),\partial_z\rangle_{h'}}{|Z|^2}.$$
Notice that from the conditions on $\Gamma$ and $h'$ it follows
immediately that $|\psi(Z)|<(1.5)\delta$.

\subsection{Proof that $t_3=t_2$}

At this point we have a product coordinate system on which the
metric is almost the Euclidean metric in the $C^1$-sense, and we
have the graph flow given by
$$Y=\frac{\partial\widetilde \gamma}{\partial
t}=\frac{\nabla_ZZ}{|Z|^2}-\psi(Z)Z$$ defined on $[t',t_3]$ with
image always contained in the given coordinate patch and written as
a graph over a fixed interval $I'$ in the $z$-axis. For every $t\in
[t',t_3]$ the length of $\gamma_{t'}|_{I'}$ in the metric $h(t')$ is
between $(0.8)r$ and $r$. The function $f(z,t)$ whose graphs give
the flow satisfies $|f_z|_{h'}\le 1/10$. Our next goal is to
estimate $|f_z|_{h'}$ and show that it is always less than  $1/10$
as long as $k^2\le 2/(t-t')$ and $t-t'\le \delta r^2$ for a
sufficiently small $\delta$, i.e., for all $t\in [t',t_2]$; that is
to say, our next goal is to prove that $t_3=t_2$. In all the
arguments that follow $C'$ is a constant that depends only on the
curvature bound for the ambient Ricci flow, but the value of $C'$ is
allowed to change from line to line.

The first step in doing this is to consider the angle between
$\nabla_ZZ$ and $Z$.

\begin{claim}\label{1891} Provided that $\delta>0$ is sufficiently small,
the angle (measured in $h(t)$) between $Y$ and $Z=(1,f_z)$ is
greater than $\pi/4$. Also,
\begin{enumerate}
\item[(1)] $$k\le |Y|<\sqrt{2}k.$$
\item[(2)] $$\left|\langle \nabla_ZZ,Z\rangle\right|<(k+ 2\delta) |Z|^3.$$
\item[(3)] $$\left|\nabla_ZZ\right|<2(|Y|+\delta).$$
\item[(4)] $$|\langle Y,Z\rangle |\le |Y||f_z|(1+3\delta).$$
\end{enumerate}
\end{claim}

\begin{proof}
Under the hypothesis that $|f_z|_{h'}\le 1/10$, it is easy to see
that the Euclidean angle between $(0,f_{zz})$ and $(1,f_z)$ is at
most $\pi/2-\pi/5$. From this, the first statement follows
immediately provided that $\delta$ is sufficiently small. Since $Y$
is the sum of $\nabla_SS$ and a multiple of $Z$ and since
$\nabla_SS$ is $h(t)$-orthogonal to $Z$, it follows that
$|Y|=|\nabla_SS|\left({\rm cos}(\theta)\right)^{-1}$, where $\theta$
is the angle between $\nabla_SS$ and $Y$. Since $Y$ is a multiple of
$(0,f_{zz})$, it follows from the first part of the claim that the
$h(t)$-angle between $Y$ and $\nabla_SS$ is less than $\pi/4$. Item
(1) of the claim then follows from the fact that by definition
$|\nabla_SS|=k$.

Since
$$\frac{\nabla_ZZ}{|Z|^2}=Y+\psi(Z)Z,$$ and $|Z|^2\le (1.01)(1+\delta)$, the third item is
immediate. For Item (4), since $Y$ is $h'$-orthogonal to the $z$-axis, we have
$$\left|\langle Y,Z\rangle_{h'}\right|=\left|\langle Y,(0,f_z)\rangle_{h'}\right|\le |Y|_{h'}|f_z|_{h'}.$$
From this and the comparison of $h(t)$ and $h'$, the Item (4) is
immediate. Lastly, let us consider Item (2). We have
$$\langle Y,Z\rangle =\frac{\langle \nabla_ZZ,Z\rangle}{|Z|^2}-\langle
\Gamma(Z,Z),\partial_z\rangle_{h'}.$$ Thus, from Item (4) we have
$$\frac{\langle \nabla_ZZ,Z\rangle}{|Z|^2}\le |Y||f_z|(1+3\delta)+(1.5)\delta.$$
Since $Y<\sqrt{2}k$ and $|f_z|<1/10$, the Item (2) follows.
\end{proof}

\begin{claim}\label{psiclaim}
The following hold provided that $\delta>0$ is sufficiently small:
\begin{enumerate}
\item[(1)] $|Z(\psi(Z))|<C'(1+\delta |Y|)$, and
\item[(2)] $|Y(\psi(Z))|<C'(|Y|+\delta|\nabla_ZY|)$.
\end{enumerate}
(Recall that $C'$ is a constant depending only on the curvature bound of the
ambient Ricci flow.)
\end{claim}

\begin{proof}
For the first item, we write $Z(\psi(Z))$ as a sum of terms where the
differentiation by $Z$ acts on the various. When the $Z$-derivative acts on
$\Gamma$ the resulting term has norm bounded by a constant depending only on
the curvature of the ambient Ricci flow. When the $Z$-derivative acts on one of
the $Z$-terms in $\Gamma(Z,Z)$ the norm of the result is bounded by
$2\delta|\nabla_ZZ||Z|$. Action on each of the other $Z$-terms gives a term
bounded in norm by the same expression. Lastly, when the $Z$-derivative acts on
the constant metric $h'$ the norm of the result is bounded by $2\delta^2$.
Since $|\nabla_ZZ|\le 2(|Y|+\delta)$, the first item follows.

We compute $Y(\psi(Z))$ in a similar fashion. When the $Y$-derivative acts on
the $\Gamma$, the norm of the result is bounded by $C'|Y|$. When the
$Y$-derivative acts on one of the $Z$-terms the norm of the result is bounded
by $2\delta|\nabla_YZ|$. Lastly, when the $Y$-derivative acts on the constant
metric $h'$, the norm of the result is bounded by $\delta^2|Y|$. Putting all
these terms together establishes the second inequality above.
\end{proof}

Now we wish to compute $\int_{I'\times\{t\}}|Z|^2dz.$ To do this we
first note that using the definition of $Y$, and arguing as in the
proof of the first equation in of Lemma~\ref{1stcurvshr} we have we
have
\begin{eqnarray*}\frac{\partial}{\partial t}|Z|^2
& = & -2{\rm Ric}(Z,Z)+2\langle \nabla_YZ,Z\rangle \\
& = & -2{\rm Ric}(Z,Z)+2\langle \nabla_ZY,Z\rangle \end{eqnarray*}
Direct computation shows that
$$2\langle
\nabla_Z\left(\frac{\nabla_ZZ}{|Z|^2}\right),Z\rangle
=Z\left(\frac{Z(|Z|^2)}{|Z|^2}\right)-2\frac{|\nabla_ZZ|^2}{|Z|^4}|Z|^2.$$ Thus
from the Claim~\ref{curvflowcompare}, we have
\begin{eqnarray}\label{Z_t}
\frac{\partial}{\partial t}|Z|^2 & = & 2\langle\nabla_ZY,Z\rangle -2{\rm
Ric}(Z,Z)\nonumber \\
& = &
Z\left(\frac{Z(|Z|^2)}{|Z|^2}\right)-2\frac{|\nabla_ZZ|^2}{|Z|^4}|Z|^2-2\langle
\nabla_Z(\psi(Z)Z),Z\rangle-2{\rm Ric}(Z,Z) \nonumber \\
& = & Z\left(\frac{Z(|Z|^2)}{|Z|^2}\right)-2|Y|^2|Z|^2 +V,
\end{eqnarray}
where $$V=-4|Z|^2\langle Y,\psi(Z)Z\rangle-2\psi^2(Z)|Z|^4-2\langle
\nabla_Z(\psi(Z)Z),Z\rangle-2{\rm Ric}(Z,Z).$$ By Item (1) in
Claim~\ref{psiclaim} and Item (4) in Claim~\ref{1891} we have
\begin{equation}\label{Vest}
|V|<C'(1+\delta|Y|).
\end{equation}

Using this and the fact that $|Y|\le \sqrt{2}k$ we compute:
\begin{eqnarray*}
\frac{d}{dt}\int_{I'\times\{t\}}|Z|^2dz & \le & \int_{I'\times\{t\}}
Z\left(\frac{Z(|Z|^2)}{|Z|^2}\right)dz +\int_{I'\times \{t\}}\left(C'(1+\delta
k)\right) dz
 \\
& = & \frac{Z(|Z|^2)}{|Z|^2}\Bigl|_0^a\Bigr.+\int_{I'\times\{t\}}\left(C'(1+\delta k)\right) dz \\
& = & 2\frac{\langle \nabla_ZZ,Z\rangle}{|Z|^2}\Bigl|_0^a\Bigr.
+\int_{I'\times\{t\}}\left(C'(1+\delta k)\right) dz,
\end{eqnarray*}
where we denote the endpoints of $I'$ by $\{0\}$ and $\{a\}$. By Item (2) in
Claim~\ref{1891}, the first term is at most
$2(k+2\delta)\sqrt{(1.01)(1+\delta)}$, which is at most $\frac{8}{\sqrt{t-t'}}$
and the second term is at most $C'(1+\delta k) r$. Now integrating from $t'$ to
$t$ we see that for any $t\in [t',t_3]$ we have
$$\int_{I'\times\{t\}}|Z|^2dz\le\int_{I'\times\{t'\}}|Z|^2dz+16\sqrt{\delta}r+C'\delta
r^3+C'\delta^{3/2}r^2.$$  Since $|f_z(z,t')|_{h'}\le 2\delta$ and $|Z|^2$ is
between $(1-\delta)(1+|f_z|_{h'}^2)$ and $(1+\delta)(1+|f_z|_{h'}^2)$,  we see
that $\int_{I'\times\{t'\}}|Z|^2dz\le (1+3\delta)\ell_{h'}(I')$. It follows
that for any $t\in [t',t_3]$ we have
$$\int_{I'\times\{t\}}|Z|^2dz\le (1+3\delta)l_{h'}(I')+C'(\sqrt{\delta}r+\delta r^3+\delta^{3/2}r^2).$$
Since $|Z|^2$ is between $(1-\delta)(1+|f_z|_{h'}^2)$ and
$(1+\delta)(1+|f_z|_{h'}^2)$, we see that there is a constant
$C_1''$ depending only on the ambient curvature bound such that for
any $t\in [t',t_3]$, denoting by $\ell_{h'}(I')$ the length of $I'$
with respect to $h'$, we have
$$\int_{I'\times\{t\}}|f_z|_{h'}^2dz\le 4\delta\ell_{h'}(I')
+C_1''(\sqrt{\delta}r+\delta r^3+\delta^{3/2}r^2).$$ Since
$(0.8)r\le \ell_{h'}(I')\le r<1$, we see  that provided that
$\delta$ is sufficiently small, for each $t\in [t',t_3]$ there is
$z(t)\in I'$ such $|f_z(z(t),t)|^2_{h'}\le 2C_1''\sqrt{\delta}$. If
we have chosen $\delta$ sufficiently small, this means that for each
$t\in [t',t_3]$ there is $z(t)$ such that $|f_z(z(t),t)|_{h'}\le
1/20$. Since by Corollary~\ref{D_4cor}
$\int_{I'\times\{t\}}kds<2D_4\sqrt{\delta}$,
 provided that $\delta$ is sufficiently small, it follows that for all
 $t\in[t',t_3]$ the
curve $\gamma_t|_{I'}$ is a graph of $(z,t)$ and $|f_z|_{h'}<1/10$. But by
construction either $t_3=t_2$ or there is a point in $(z,t_3)\in I'\times
\{t_3\}$ with $|f_z(z,t_3)|_{h'}=1/10$. Hence, it must be the case that
$t_3=t_2$, and thus our graph curve flow is defined for all $t\in [t',t_2]$ and
satisfies the derivative bound $|f_z|_{h'}< 1/10$ throughout the interval
$[t',t_2]$.

\subsection{Proof that $t_2=t'+\delta r^2$}

The last step is to show that the inequality $k^2< 2/(t-t')$ holds
for all $t\in [t',t'+\delta r^2]$.

 We fix a point $x_0$. We continue all the notation, assumptions and results of the previous section.
 That is to say, we lift the evolving family of curves to the tangent space $T_{x_0}M$
 using the exponential mapping, which is a
 local diffeomorphism. This tangent space is split as the
product of the $z$-axis and $B$. On this coordinate system we have
the evolving family of Riemannian metrics $h(t)$ pulled back from
the Ricci flow and also we have the Euclidean metric $h'$ from the
metric $h(t')$ on $T_{x_0}M$. We fix an interval $I'$ on the
$z$-axis
  of $h'$-length between $(0.8)r$ and $r$.    We choose $I'$ to be centered at
 $x_0$ with respect to the $z$-coordinate. On $I'\times [t',t_2]$
  we have the graph-flow  which is reparameterization
  of the pull back of the curve-shrinking flow. The graph-flow is given as the graph
  of a function $f$ with $|f_z|_{h'}<1/10$. The vector fields $Z$
  and $Y$ are as in the last section.

 We follow closely the discussion in Section 4 of \cite{AG} (pages 293 -294).
 Since we are not working in a flat background, there are two
 differences:
 (i) we take covariant derivatives instead of ordinary derivatives
 and (ii) there are various correction terms from curvature, from covariant derivatives,
 and from the fact that $Y$ is not equal to $\nabla_ZZ/|Z|^2$.

Notice that
\begin{eqnarray*}
Z\left(\frac{Z(|Z|^2)}{|Z|^2}\right) & = & \frac{|Z|^2_{zz}}{|Z|^2}-
\frac{\left(|Z|^2_z\right)^2}{|Z|^4} \\
& = & \frac{|Z|^2_{zz}}{|Z|^2}-4\langle Y,Z\rangle^2 -8\langle
Y,Z\rangle \psi(Z)|Z|^2-4\psi^2(Z)|Z|^4 \end{eqnarray*}
 Thus, it
follows from Equation~(\ref{Z_t}) that we have
$$\frac{\partial}{\partial
 t}|Z|^2=\frac{(|Z|^2)_{zz}}{|Z|^2}-2|Z|^2|Y|^2-4\langle Z,Y\rangle^2+V,$$
where $|V|\le C'(1+\delta |Y|)$ for a constant $C'$ depending only
on the curvature bound of the ambient flow.

Similar computations show that
 \begin{eqnarray*}\frac{\partial}{\partial
 t}|Y|^2 & = &
\frac{(|Y|^2)_{zz}}{|Z|^2}-\frac{2|\nabla_ZY|^2}{|Z|^2}
 -\frac{4}{|Z|^2}\langle\frac{\nabla_ZZ}{|Z|^2},Y\rangle\langle
 \nabla_ZY,Z\rangle \\
 &  & -2{\rm Ric}(Y,Y)+2\frac{{\rm
 Rm}(Y,Z,Y,Z)}{|Z|^2}-2\langle \nabla_Y(\psi(Z)Z),Y\rangle.
 \end{eqnarray*}
Of course, $$\langle \frac{\nabla_ZZ}{|Z|^2},Y\rangle=|Y|^2+\psi(Z)\langle
Z,Y\rangle.$$ Hence, putting all this together and using Claim~\ref{psiclaim}
we have
$$\frac{\partial}{\partial
 t}|Y|^2  =
\frac{(|Y|^2)_{zz}}{|Z|^2}-\frac{2|\nabla_ZY|^2}{|Z|^2}
 -\frac{4|Y|^2}{|Z|^2}\langle
 \nabla_ZY,Z\rangle+W,$$
 where
 $$|W|\le C'|Y|(|Y|+\delta|\nabla_ZY|).$$

 Now let us consider
$$Q=\frac{|Y|^2}{2-|Z|^2}.$$
Notice that since $|f_z|_{h'}<1/10$, it follows that $1-\delta\le
|Z|^2<(1.01)(1+\delta)$ on all of $[t',t_2]$. We now make computation following
the computations on p. 294 of \cite{AG} and adding in the error terms.
\begin{eqnarray*}
Q_t & = & \frac{|Y|^2_t}{(2-|Z|^2)}+\frac{|Y|^2|Z|^2_t}{(2-|Z|^2)^2} \\
& = & \frac{|Y|^2_{zz}}{|Z|^2(2-|Z|^2)}-\frac{2|\nabla_ZY|^2}{|Z|^2(2-|Z|^2)}
 -\frac{4|Y|^2}{|Z|^2(2-|Z|^2)}\langle
 \nabla_ZY,Z\rangle+\frac{W}{(2-|Z|^2)} \\
 & &
 +\frac{|Y|^2|Z|^2_{zz}}{|Z|^2(2-|Z|^2)^2}-\frac{2|Z|^2||Y|^4}{(2-|Z|^2)^2}
 -\frac{4|Y|^2}{(2-|Z|^2)^2}\langle Z,Y\rangle^2+\frac{|Y|^2}{(2-|Z|^2)^2}V.
\end{eqnarray*}

On the other hand,
\begin{eqnarray*}
\frac{Q_{zz}}{|Z|^2}=\frac{|Y|^2_{zz}}{|Z|^2(2-|Z|^2)}+
\frac{|Y|^2|Z|^2_{zz}}{|Z|^2(2-|Z|^2)^2}+\frac{2|Y|^2_z|Z|^2_z}{|Z|^2|(2-|Z|^2)^2}
+\frac{2|Y|^2\left(|Z|^2_z\right)^2}{|Z|^2(2-|Z|^2)^3}.
\end{eqnarray*}

From Claim~\ref{1891} we have
$$|Z|^2_z=2\langle\nabla_ZZ,Z\rangle=2|Z|^2\langle Y,Z\rangle+2\psi(Z)|Z|^4.$$
Plugging in this expansion gives
\begin{eqnarray*}
\frac{Q_{zz}}{|Z|^2} & = & \frac{|Y|^2_{zz}}{|Z|^2(2-|Z|^2)}+
\frac{|Y|^2|Z|^2_{zz}}{|Z|^2(2-|Z|^2)^2} \\
& & +\frac{8\langle \nabla_ZY,Y\rangle\langle
Y,Z\rangle}{(2-|Z|^2)^2}+\frac{8\psi(z)|Z|^2\langle
\nabla_ZY,Y\rangle}{(2-|Z|^2)^2} \\
& & +\frac{8|Z|^2|Y|^2\langle
Y,Z\rangle^2}{(2-|Z|^2)^3}+\frac{16\psi(Z)|Z|^4|Y|^2\langle
Y,Z\rangle}{(2-|Z|^2)^3}+\frac{8\psi^2(Z)|Z|^6|Y|^2}{(2-|Z|^2)^3}.
\end{eqnarray*}
Expanding, we have
\begin{eqnarray*}
\frac{Q_{zz}}{|Z|^2} & = & \frac{|Y|^2_{zz}}{|Z|^2(2-|Z|^2)}+
\frac{|Y|^2|Z|^2_{zz}}{|Z|^2(2-|Z|^2)^2}+\frac{8\langle
\nabla_ZY,Y\rangle\langle Y,Z\rangle}{(2-|Z|^2)^2}  \\
& & +\frac{8|Y|^2|Z|^2\langle Y,Z\rangle^2}{(2-|Z|^2)^3}+U,
\end{eqnarray*}
where
$$|U|\le C'(|Y|^2+\delta|\nabla_ZY||Y|+\delta|Y|^3).$$

Comparing the formulas yields
\begin{eqnarray*}
Q_t & = & \frac{Q_{zz}}{|Z|^2}-\frac{8\langle \nabla_ZY,Y\rangle\langle
Y,Z\rangle}{(2-|Z|^2)^2}-\frac{8|Y|^2|Z|^2\langle Y,Z\rangle^2}{(2-|Z|^2)^3} \\
& & -\frac{2|\nabla_ZY|^2}{|Z|^2(2-|Z|^2)}
 -\frac{4|Y|^2}{|Z|^2(2-|Z|^2)}\langle
 \nabla_ZY,Z\rangle \\
 & & -\frac{2|Z|^2||Y|^4}{(2-|Z|^2)^2}
 -\frac{4|Y|^2}{(2-|Z|^2)^2}\langle Z,Y\rangle^2+A,
\end{eqnarray*}
where
$$|A|\le C'(|Y|^2+\delta|\nabla_ZY||Y|+\delta|Y|^3).$$

Using Item (4) of Claim~\ref{1891} this leads to
\begin{eqnarray*}
Q_t & \le & \frac{Q_{zz}}{|Z|^2}+\frac{8(1+3\delta)|Y|^2|f_z|
|\nabla_ZY|}{(2-|Z|^2)^2} -\frac{4|Y|^2\langle
 \nabla_ZY,Y\rangle }{|Z|^2(2-|Z|^2)} \\
 & &
-\frac{2|\nabla_ZY|^2}{|Z|^2(2-|Z|^2)}
 -\frac{2|Z|^2||Y|^4}{(2-|Z|^2)^2}
 +|A|
 \end{eqnarray*}
Next, we have
\begin{claim}
$$|\langle \nabla_ZY,Z\rangle|\le (|f_z|(|\nabla_ZY|+2\delta|Y|)+\delta|Z|^2|Y|)(1+\delta).$$
\end{claim}

\begin{proof}
Since $Y=(0,\phi)$ for some function $\phi$, we have
$\nabla_ZY=(0,\phi_z)+\Gamma(Z,Y)$ and hence
$$|\langle\nabla_ZY,Z\rangle_{h'}| = |\langle \nabla_ZY,(1,f_z)\rangle_{h'}|\le |\langle
f_z ,\phi_z\rangle_{h'}|+|\langle \Gamma(Z,Y),Z\rangle_{h'}|\le
|f_z|_{h'}|\phi_z|_{h'}+\delta|Z|^2|Y|.$$ On the other hand
$\nabla_ZY=(0,\phi_z)+\Gamma(Z,Y)$ so that $|\phi_z|_{h'}\le|\nabla_ZY|+\delta
|Z||Y|$. From this the claim follows.
\end{proof}

Now for $\delta>0$ sufficiently small, using the fact that
$1-\delta<|Z|^2<(1+\delta)(1.01)$ we can rewrite this as
\begin{eqnarray*}
Q_t & \le & \frac{Q_{zz}}{|Z|^2}+\frac{8(1+3\delta)|Y|^2|f_z|
|\nabla_ZY|}{(2-|Z|^2)^2} +\frac{4|Y|^2|(1+\delta)
 |\nabla_ZY||f_z|}{|Z|^2(2-|Z|^2)} \\
 & & -\frac{2|\nabla_ZY|^2}{|Z|^2(2-|Z|^2)}-\frac{(1.95)|Y|^4}{(2-|Z|^2)^2} +\tilde
 A,
 \end{eqnarray*}
where $\tilde A\le C'(|Y|^2+\delta|Y||\nabla_ZY|+\delta|Y|^3)$. Of
course, $|Y||\nabla_ZY|+|Y|^3\le 2|Y|^2+|\nabla_ZY|^2+|Y|^4$. Using
this, provided that $\delta$ is sufficiently small, we can rewrite
this as

\begin{eqnarray*}
Q_t & \le & \frac{Q_{zz}}{|Z|^2} +\frac{1}{(2-|Z|^2)}\cdot
\left[\frac{8(1+3\delta)|\nabla_ZY||f_z||Y|^2-(0.9)|Y|^4}{(2-|Z|^2)} \right.\\&
&
 +\left. \frac{4|Y|^2|(1+\delta)
 |\nabla_ZY||f_z|-(1.9)|\nabla_ZY|^2}{|Z|^2}\right]-Q^2+\tilde A''
\end{eqnarray*}
where $\tilde A''\le C'(|Y|^2)$. We denote the quantity within the
brackets by $B$ and we estimate
\begin{eqnarray*}
B & \le &
8(1+3\delta)\frac{|Y|^2|\nabla_ZY|(1/10)(1+\delta)}{(2-(1.01)(1+\delta))}
+\frac{4(1/10)(1+\delta)|Y|^2|\nabla_ZY|}{(1-\delta)} \\
& & -\frac{(1.9)}{(1.01)(1+\delta)}|\nabla_ZY|^2
-\frac{(0.9)|Y|^4}{2-(1.01)(1+\delta)} \\
& \le  & (1.6)|Y|^2|\nabla_ZY|-(0.8)|\nabla_ZY|^2-(0.8)|Y|^4 \\
& \le  & 0.
\end{eqnarray*}
Therefore,
$$Q_t\le \frac{Q_{zz}}{|Z|^2}-Q^2+|\tilde A|\le \frac{Q_{zz}}{|Z|^2}-(Q-C_1')^2+(C_1')^2,$$
for some constant $C'_1>1$ depending only on the curvature bound for the
ambient Ricci flow.

Denote by $l$ the length of $I'$ under $h'$. As we have already
seen, $(0.8)r\le l\le r$. We translate the $z$-coordinate so that
$z=0$ is one endpoint of $I'$ and $z=l$ is the other endpoint; the
point $x_0$ then corresponds to $z=l/2$. Consider the function
$g=l^2/(z^2(l-z)^2)$ on $I'\times [t',t_2]$. Direct computation
shows that $g_{zz}\le 12g^2$.
 Now set
$$\widetilde Q=Q-C_1'$$
 and
$$h=\frac{1}{t-t'}+\frac{4(1-\delta)^{-1}l^2}{z^2(l-z)^2}+C_1'.$$
 Then
$$ -h_t+(1-\delta)^{-1}h_{zz}+(C_1')^2\le h^2,$$
so that
$$(\widetilde Q-h)_t\le \frac{\widetilde Q_{zz}}{|Z|^2}-\frac{h_{zz}}{1-\delta}-\widetilde
Q^2+h^2.$$ Since both $h$ and $h_{zz}$ are positive, at any point
where $\widetilde Q -h\ge 0$ and $\widetilde Q_{zz}<0$, we have
$(\widetilde Q-h)_t<0$. At any point where $Q_{zz}\ge 0$, using the
fact that $|Z|^2\ge (1-\delta)$ we have
$$(\widetilde Q-h)_t\le (1-\delta)^{-1}(\widetilde Q-h)_{zz}-\widetilde
Q^2+h^2.$$ Thus, for any fixed $t$, at any local maximum for
$(\widetilde Q-h)(\cdot,t)$ at which $(\widetilde Q-h)$ is $\ge 0$
we have $(\widetilde Q-h)_t\le 0$. Since $\widetilde Q-h$ equals
$-\infty$ at the end points of $I'$ for all times, there is a
continuous function $f(t)={\rm max}_{z\in I'}(\widetilde Q-h)(z,t)$,
defined for all $t\in (t',t_2]$ approaching $-\infty$ uniformly as
$t$ approaches $t'$ from above. By the previous discussion, at any
point where $f(t)\ge 0$ we have $f'(t)\le 0$ in the sense of forward
difference quotients. It now follows that $f(t)\le 0$ for all $t\in
(t',t_2]$.
 This means that for all $t\in (t',t_2]$ at the $h'$-midpoint $x_0$ of
$I'$ (the point where $z=l/2$) we have
$$Q(x_0,t)\le \frac{1}{t-t'}+\frac{16\cdot 4(1-\delta)^{-1}}{l^2}+C_1'.$$
Since $l\ge (0.8)r$ and since $t-t'\le \delta r^2$, we see that
provided $\delta$ is sufficiently small (depending on the bound of
the curvature of the ambient flow) we have
$$Q(x_0,t)< \frac{3}{2(t-t')}$$ for all $t\in [t',t_2].$
Of course, since $|Z|^2\ge 1-\delta$ everywhere, this shows that
$$k^2(x_0,t)\le |Y(x_0,t)|^2=
(2-|Z(x_0,t)|^2)Q(x_0,t)< \frac{2}{(t-t')}$$ for all $t\in
[t',t_2]$. Since $x_0$ was an arbitrary point of $c(\cdot,t')$,
 this shows that
 $k(x,t)<\sqrt{\frac{2}{t-t'}}$ for all $x\in c(\cdot,t)$ and all $t\in
 [t',t_2]$. By the definition of $t_2$ this implies that $t_2=t'+\delta r^2$ and completes the
 proof of  Lemma~\ref{csshi1}.

\chapter{Appendix: Canonical neighborhoods}\label{sect:canonnbhd}

Recall that an $\epsilon$-neck structure on a Riemannian manifold
$(N,g)$ centered at a point $x\in N$ is a diffeomorphism $\psi\colon
S^2\times (-\epsilon^{-1},\epsilon^{-1})\to N$ with the property
that $x\in \psi(S^2\times\{0\})$ and the property that $R(x)\psi^*g$
is within $\epsilon$ in the $C^{[1/\epsilon]}$-topology of the
product metric $h_0\times ds^2$, where $h_0$ is the round metric on
$S^2$ of scalar curvature $1$ and $ds^2$ is the Euclidean metric on
the interval. Recall that the {\em scale} of the $\epsilon$-neck is
$R(x)^{-1/2}$. We define $s=s_N\colon N\to
(-\epsilon^{-1},\epsilon^{-1})$ as the composition of $\psi^{-1}$
followed by the projection to the second factor.

\section{Shortening curves}

\begin{lem}\label{curveshort}
The following holds for all $\epsilon>0$ sufficiently small. Suppose
that $(M,g)$ is a Riemannian manifold and that $N\subset M$ is an
$\epsilon$-neck centered at $x$. Let $S(x)$ be the central
two-sphere of this neck and suppose that $S(x)$ separates $M$. Let
$y\in M$. Orient $s$ so that $y$ lies in the closure of the positive
side of $S(x)$. Let $\gamma\colon [0,a]\to M$ be a rectifiable curve
from $x$ to $y$. If $\gamma$ contains a point of
$s^{-1}(-\epsilon^{-1},-\epsilon^{-1}/2)$ then there is a
rectifiable curve from $x$ to $y$ contained in the closure of the
positive side of $S(x)$ whose length is at most the length of
$\gamma$ minus $\frac{1}{2}\epsilon^{-1}R(x)^{-1/2}$.
\end{lem}

\begin{proof}
Since $\gamma$ contains a point on the negative side of $S(x)$ and
it ends on the positive side of $S(x)$, there is a $c\in (0,a)$ such
that $\gamma(c)\in S(x)$ and $\gamma|_{(c,a]}$ is disjoint from
$S(x)$. Since $\gamma|_{[0,c]}$ has both endpoints in $S(x)$ and
also contains a point of $s^{-1}(-\epsilon^{-1},-\epsilon^{-1}/2)$,
it follows that for $\epsilon$ sufficiently small, the length of
$\gamma|_{[0,c]}$ is at least $3\epsilon^{-1}R(x)^{-1/2}/4$. On the
other hand, there is a path $\mu$ in $S(x)$ connecting $x$ to
$\gamma(c)$ of length at most $2\sqrt{2}\pi(1+\epsilon)$. Thus, if
$\epsilon$ is sufficiently small, the concatenation of $\mu$
followed by $\gamma|_{[c,a]}$ is the required shorter path.
\end{proof}

\section{The geometry of an $\epsilon$-neck}

\begin{lem}\label{neckcurv}
For any $0<\alpha<1/8$  there is $\epsilon_1=\epsilon_1(\alpha)>0$ such that
the following two conditions hold for all $0<\epsilon\le \epsilon_1$.
\begin{enumerate} \item[(1)] If $(N,g)$ is an $\epsilon$-neck centered at $x$ of
scale one (i.e., with $R(x)=1$) then the principal sectional
curvatures at any point of $N$ are within $\alpha/6$ of
$\{1/2,0,0\}$. In particular, for any $y\in N$ we have
$$(1-\alpha)\le R(y)\le (1+\alpha).$$
\item[(2)] There is unique two-plane of maximal sectional curvature at
every point of an $\epsilon$-neck, and the angle between the
distribution of two-planes of maximal sectional curvature and the
two-plane field tangent to the family of two-spheres of the
$\epsilon$-neck structure is everywhere less than $\alpha$.
\end{enumerate}
\end{lem}

\begin{proof}
The principal curvatures and their directions are continuous
functions of the metric $g$ in the space of metrics with the
$C^2$-topology. The statements follow immediately.
\end{proof}

\begin{cor}\label{s2isotopic}
The following holds for any $\epsilon>0$ sufficiently small. Suppose
that $(N,g)$ is an $\epsilon$-neck and we have  and an embedding
$f\colon S^2 \to N$ with the property that the restriction of $g$ to
the image of this embedding is within $\epsilon$ in the
$C^{[1/\epsilon]}$-topology to the round metric $h_0$ of scalar
curvature one on $S^2$ and with the norm of the second fundamental
form less than $\epsilon$. Then the two-sphere $f(S^2)$ is isotopic
in $N$ to any member of the family of two-spheres coming from the
$\epsilon$-neck structure on $N$.
\end{cor}

\begin{proof}
By the previous lemma, if $\epsilon$ is sufficiently small for every
$n\in N$  there is a unique two-plane, $P_n$, at each point on which
the sectional curvature is maximal. The sectional curvature on this
two-plane is close to $1/2$ and the other two eigenvalues of the
curvature operator at $n$ are close to zero. Furthermore, $P_n$
makes small $g$-angle with the tangent planes to the $S^2$-factors
in the neck structure. Under the condition that the restriction of
the metric to $f(S^2)$ is close to the round metric $h_0$ and the
norm of the second fundamental form is small, we see that for every
$p\in S^2$ the two-plane $df(T_pS^2)$ makes a small $g$-angle with
$P_n$ and hence with the tangent planes to the family of two-spheres
coming from the neck structure. Since $g$ is close to the product
metric, this means that the angle between $df (T_nS^2)$ and the
tangents to the family of two-spheres coming from the neck
structure, measured in the product metric, is also small. Hence, the
composition of $f$ followed by the projection mapping $N\to S^2$
induced by the neck structure determines a submersion of $S^2$ onto
itself. Since $S^2$ is compact and simply connected, any submersion
of $S^2$ onto itself is a diffeomorphism. This means that $f(S^2)$
crosses each line $\{x\}\times (-\epsilon^{-1},\epsilon^{-1})$
transversely and in exactly one point. Clearly then, it is isotopic
in $N$ to any two-sphere of the form $S^2\times \{s\}$.
\end{proof}

\begin{lem}\label{directions}
 For any $\alpha>0$  there is $\epsilon_2=\epsilon_2(\alpha)>0$
such that the following hold for all $0<\epsilon\le \epsilon_2$.
Suppose that $(N,g)$ is an $\epsilon$-neck centered at $x$ and
$R(x)=1$. Suppose that $\gamma$ is a minimal geodesic in  $N$ from
$p$ to $q$. We suppose that $\gamma$ is parameterized by arc length,
is of length $\ell>\epsilon^{-1}/100$, and that $s(p)<s(q)$. Then
for all $s$ in the domain of definition of $\gamma$ we have
$$|\gamma'(s)-(\partial/\partial s)|_g<\alpha.$$
 In particular, the angle between $\gamma'$ and $\partial/\partial
s$ is less than $2\alpha$.
 Also,  any
member $S^2$ of the family of two-spheres in the $N$ has intrinsic
diameter at most $(1+\alpha)\sqrt{2}\pi$.
\end{lem}

\begin{proof}
Let us consider a geodesic $\mu$ in the product Riemannian manifold
$S^2\times \Ar$ with the metric on $S^2$ being of constant Gaussian
curvature $1/2$, i.e., radius $\sqrt{2}$. Its projections, $\mu_1$
and $\mu_2$, to $S^2$ and to $\Ar$, respectively,
 are also geodesics, and $|\mu|=\sqrt{|\mu_1|^2+|\mu_2|^2}$.
 For $\mu$ to be a minimal geodesic, the same is true of each of its projections.
In particular, when $\mu$ is minimal, the length of $\mu_1$ is at
most $\sqrt{2}\pi$. Hence, for any $\alpha'>0$, if $\mu$ is
sufficiently long and if the final endpoint has a larger $s$-value
than the initial point, then the angle between the tangent vectors
$\mu'(s)$ and $\partial/\partial s$ is less than $\alpha'$. This
establishes the result for the standard metric on the model for
$\epsilon$-necks.

The first statement now follows for all $\epsilon$ sufficiently
small and all $\epsilon$-necks because minimal geodesics between a
pair of points in a manifold vary continuously in the $C^1$-topology
as a function of the space of metrics with the $C^k$-topology, since
$k\ge 2$. The second statement is obvious since the diameter of any
member of the family of two-spheres in the standard metric is
$\sqrt{2}\pi$.
\end{proof}

\begin{cor}\label{neckgeo}
For any $\alpha>0$ there is $\epsilon_3=\epsilon_3(\alpha)>0$ such
that the  following hold for any $0<\epsilon\le \epsilon_3$ and any
$\epsilon$-neck $N$ of scale $1$ centered at $x$.
\begin{enumerate}
\item[(1)] Suppose that $p$ and $q$ are points of  $N$
with either $|s(q)-s(p)|\ge \epsilon^{-1}/100$ or $d(p,q)\ge
\epsilon^{-1}/100$. Then  we have
$$(1-\alpha)|s(q)-s(p)|\le d(p,q)\le (1+\alpha)|s(q)-s(p)|.$$
\item[(2)]
$$B(x,(1-\alpha)\epsilon^{-1})\subset
N\subset B(x,(1+\alpha)\epsilon^{-1}).$$
\item[(3)] Any geodesic that exits from both ends of $N$ has length at least $2(1-\alpha)
\epsilon^{-1}$.\end{enumerate}
\end{cor}

\begin{cor}\label{compint}
The following holds for all $\epsilon>0$ sufficiently small. Let $N$
be an $\epsilon$-neck centered at $x$. If $\gamma$ is a shortest
geodesic in $N$ between its endpoints and if
$|\gamma|>R(x)^{-1/2}\epsilon^{-1}/100$, then $\gamma$ crosses each
two-sphere in the neck structure on $N$ at most once.
\end{cor}

There is a closely related lemma.

\begin{figure}[ht]
  \relabelbox {
  \centerline{\epsfbox{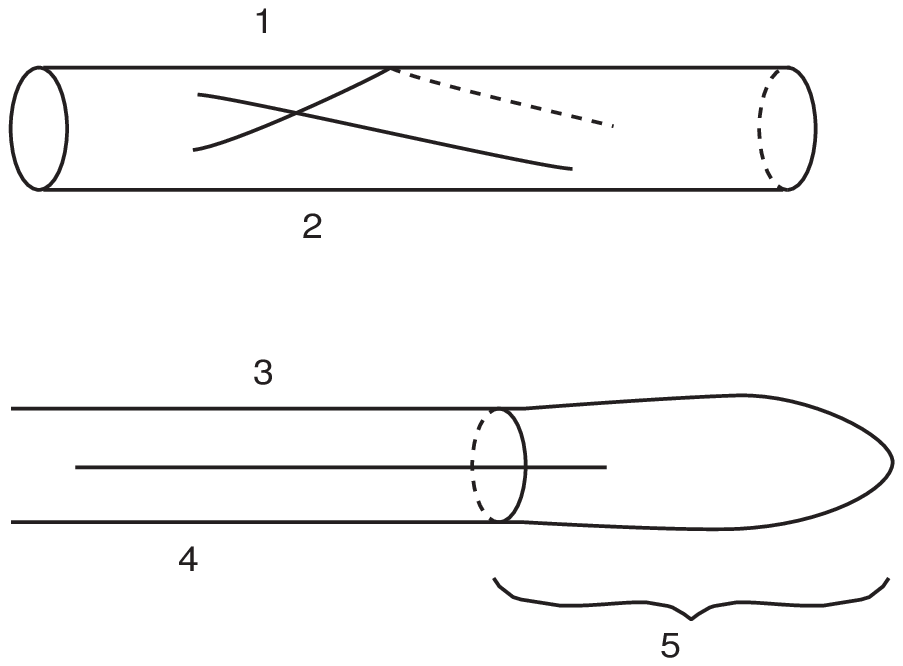}}}
  \relabel{1}{shortest geodesics}
  \relabel{2}{$\epsilon$-neck}
  \relabel{3}{shortest geodesic}
  \relabel{4}{$\epsilon$-cap}
  \relabel{5}{core}
  \endrelabelbox
  \caption{Shortest geodesics in necks and caps}\label{fig:short}
\end{figure}

\begin{lem}\label{genint}
The following holds for every $\epsilon>0$ sufficiently small.
Suppose that $(M,g)$ is a Riemannian manifold and that $N\subset M$
is an $\epsilon$-neck centered at $x$ and suppose that $\gamma$ is a
shortest geodesic in $M$ between its endpoints and that the length
of every component of $N\cap|\gamma|$ has length at least $
R(x)^{-1/2}\epsilon^{-1}/8$. Then $\gamma$ crosses each two-sphere
in the neck structure on $N$ at most once; see {\sc
Fig.}~\ref{fig:short}.
\end{lem}

\begin{proof}
We parameterize $\gamma$ as a map from $[a,b]\to M$. By
Corollary~\ref{compint}, provided that $\epsilon>0$ is sufficiently
small, each component of $\gamma\cap N$ crosses each two-sphere of
the neck structure at most once. Suppose that there is some
two-sphere $S^2\times\{x\}$ that is crossed by two different
components of $\gamma$. Let $c<d$ be two points of intersection of
$\gamma$ with $S^2\times\{s\}$.

There are two cases to consider. Suppose that the two components of
$\gamma\cap N$ cross $S^2\times \{x\}$ in opposite directions. In
this case, since each component of $\gamma\cap N$ has length at
least $\epsilon^{-1}/8$, then applying Corollary~\ref{neckgeo} we
can take the two-sphere that they both cross to be $S^2\times \{s\}$
for some $s\in (-(0.9)\epsilon^{-1},(0.9)\epsilon^{-1})$. Applying
Corollary~\ref{neckgeo} again we see that the distance from this
sphere to the complement of $N$ is  at least
$R(x)^{-1/2}\epsilon^{-1}/20$. Let $c<d$ be the points of
intersection. Remove $\gamma([c,d])$ from $\gamma$ and replace it by
a path in $S^2\times\{s\}$ between $\gamma(c)$ and $\gamma(d)$. If
$\epsilon$ is sufficiently small, by Lemma~\ref{directions} we can
choose this path to have length at most $2\pi$, and hence the result
will be a shorter path.

The other possibility is that $\gamma$ crosses $S^2\times\{s\}$
twice in the same direction. In this case the central two-sphere of
$N$ does not separate $M$ and $\gamma$ makes a circuit transverse to
the two-sphere. In particular, by Corollary~\ref{neckgeo} the length
of $\gamma([c,d])$ is bounded below by
$2(1-\alpha)R(x)^{-1/2}\epsilon^{-1}$ where we can take $\alpha>0$
as close to zero as we want by making $\epsilon$ smaller. Clearly,
then in this case as well, replacing $\gamma([c,d])$ with a path of
length less than $2\pi R(x)^{-1/2}$ on $S^2\times\{s\}$ will shorten
the length of $\gamma$.
\end{proof}

\begin{cor}\label{coreint}
The following holds for all $\epsilon>0$ sufficiently small and any
$C<\infty$. Let $X$ be an $(C,\epsilon)$-cap in a complete
Riemannian manifold $(M,g)$, and let $Y$ be its core and let $S$ be
the central two-sphere of the $\epsilon$-neck $N=X-\overline Y$. We
orient the $s$-direction in $N$ so that $Y$ lies off the negative
end of $N$. Let $\widehat Y$ be the union of $\overline Y$ and the
closed negative half of $N$ and let $S$ be the boundary of $\widehat
Y$. Suppose that $\gamma$ is a minimal geodesic in $(M,g)$ that
contains a point of the core $Y$. Then the intersection of $\gamma$
with $\widehat Y$ is an interval containing an endpoint of $\gamma$;
see {\sc Fig.}~\ref{fig:short}.
\end{cor}

\begin{proof}
If $\gamma$ is completely contained in $\widehat Y$ then the result
is clear. Suppose that the path is $\gamma\colon [a,b]\to M$ and
$\gamma(d)\in Y$ for some $d\in [a,b]$. Suppose that there are
$a'<d<b'$ with $\gamma(a')$ and $\gamma(b')$ contained in $S$. Then,
by Corollary~\ref{neckgeo}, replacing $\gamma|_{[a',b']}$ with a
path on $S$ joining $\gamma(a')$ to $\gamma(b')$ creates a shorter
path with the same endpoints. This shows that at least one of the
paths $\gamma|_{[a,d]}$ or $\gamma|_{[d,b]}$, let us say
$\gamma|_{[a,d]}$, is contained in $\widehat Y$. The other path
$\gamma|_{[d,b]}$ has an endpoint in $Y$ and exits from $\hat Y$,
hence by Corollary~\ref{compint} there is a subinterval $[d,b']$
such that either $\gamma(b')$ is contained in the frontier of $X$ or
$b=b'$ and furthermore $\gamma([d,b'])$ crosses each two-sphere of
the $\epsilon$-neck structure on $N$ at most once. Since $\gamma$ is
not contained in $\hat Y$, there is $b''\in [d,b']$ such that
$\gamma(b'')\in S$. We have constructed a subinterval of the form
$[a,b'']$ such that $\gamma([a,b''])$ is contained in $\widehat Y$.
If $b'=b$, then it follows from the fact that $\gamma|_{[d,b]}$
crosses each two-sphere of $N$ at most once that $\gamma|_{[b'',b]}$
is disjoint from $Y$.  This establishes  the result  in this case.
Suppose that $b'<b$.
 If there is $c\in [b',b]$ with $\gamma(c)\in \widehat
Y$ then the length of $\gamma([b'',c])$ is at least twice the
distance from $S$ to the frontier of the positive end of $N$. Thus,
we could create a shorter path with the same endpoints by joining
$\gamma(b'')$ to $\gamma(c)$ by a path of $S$. This means that
$\gamma|_{[b',b]}$ is disjoint from $S$ and hence from $\widehat Y$,
proving the result in this case as well.
\end{proof}

We also wish to compare distances from points outside the neck with
distances in the neck.

\begin{lem}\label{neckdistcomp}
Given $0<\alpha<1$ there is $\epsilon_4=\epsilon_4(\alpha)>0$ such
that the following holds for any $0<\epsilon\le \epsilon_4$. Suppose
that $N$ is an $\epsilon$-neck centered at $x$ in a connected
manifold $M$ (here we are not assuming that $R(x)=1$). We suppose
that the central $2$-sphere of $N$ separates $M$. Let $z$ be a point
outside of  the middle two-thirds of $N$ and lying on the negative
side of the central $2$-sphere of $N$. (We allow both the case when
$z\in N$ and when $z\not\in N$.) Let $p$ be a point in the middle
half of $N$. Let $\mu\colon [0,a]\to N$ be a straight line segment
(with respect to the standard product metric) in the positive
$s$-direction in $N$ beginning at $p$ and ending at a point $q$ of
$N$. Then
$$(1-\alpha)(s(q)-s(p))\le d(z,q)-d(z,p)\le (1+\alpha)(s(q)-s(p)).$$
\end{lem}

\begin{proof}
This statement is clearly true for the product metric on an infinite
cylinder, and hence by continuity, for any given $\alpha$, the
result holds for all $\epsilon>0$ sufficiently small.
\end{proof}

\noindent{\bf N.B.} It is important that the central two-sphere of
$N$ separates the ambient manifold $M$. Otherwise, there may be
shorter geodesics from $z$ to $q$ entering the other end of $N$.

\begin{figure}[ht]
  \relabelbox {
  \centerline{\epsfbox{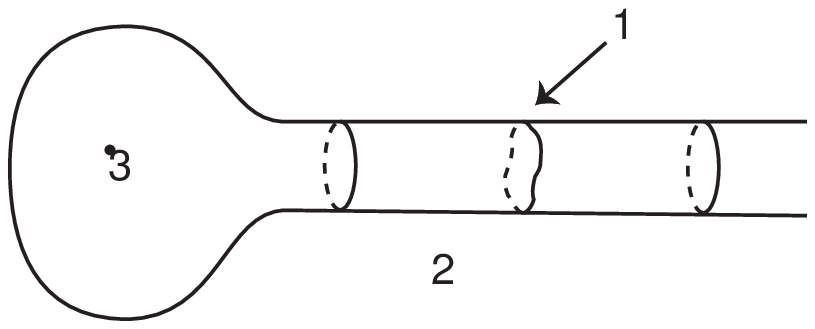}}}
  \relabel{1}{$\partial B(z,d)$}
  \relabel{2}{$\epsilon$-neck}
  \relabel{3}{$z$}
  \endrelabelbox
  \caption{Intersection of metric balls and
  necks}\label{fig:intball}
\end{figure}

\begin{lem}\label{topiso}
Given any $\alpha>0$  there is $\epsilon(\alpha)>0$ such that  the
following holds for any $0<\epsilon\le \epsilon(\alpha)$. Suppose
that $N$ is an $\epsilon$-neck centered at $x$ in a connected
manifold $M$ (here we are not assuming that $R(x)=1$) and that $z$
is a point outside the middle two-thirds of $N$. We suppose that the
central two-sphere of $N$ separates $M$. Let $p$ be a point in the
middle sixth of $N$ at distance $d$ from $z$. Then the intersection
of the boundary of the metric ball $B(z,d)$ with $N$ is a
topological $2$-sphere contained in the middle quarter of $N$ that
maps homeomorphically onto $S^2$ under the projection mapping $N\to
S^2$ determined by the $\epsilon$-neck structure. Furthermore, if
$p'\in \partial B(z,d)$ then $|s(p)-s(p')|<\alpha
R(x)^{-1/2}\epsilon^{-1}$; see {\sc Fig.}~\ref{fig:intball}.
\end{lem}

\begin{proof}
The statement is scale-invariant, so we can assume that $R(x)=1$.
 Denote by $S(z,d)$ the boundary of the metric ball $B(z,d)$. We
orient $s$ so that $z$ lies to the negative side of the central
two-sphere of $N$. It follows immediately from the previous result
that, provided that $\epsilon>0$ is sufficiently small,  $S(z,d)$
intersects any line ${y}\times (-\epsilon^{-1}/3,\epsilon^{-1}/3)$
in at most one point. To complete the proof we need only show that
$S(z,d)$ is contained
$s^{-1}((s(p)-\alpha\epsilon^{-1},s(p)+\alpha\epsilon^{-1}))$. The
distance from $d$ to any point in the two-sphere factor of $N$
containing $p$ is contained in the interval $[d-2\pi,d+2\pi]$.
Provided that $\epsilon$ is sufficiently small depending on
$\alpha$, the result follows immediately from
Lemma~\ref{neckdistcomp}.
\end{proof}

\section{Overlapping $\epsilon$-necks}

The subject of this section is the internal geometric properties of
$\epsilon$-necks and of intersections of $\epsilon$-necks. We assume that
$\epsilon\le 1/200$.

\begin{figure}[ht]
  \centerline{\epsfbox{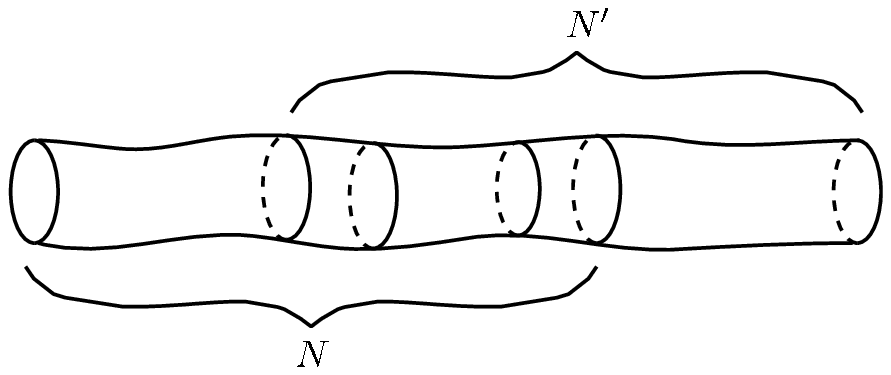}}
  \caption{Overlapping $\epsilon$-necks.}\label{fig:overlap}
\end{figure}

\begin{prop}\label{S2intersection}
Given $0<\alpha\le 10^{-2}$, there is
$\epsilon_5=\epsilon_5(\alpha)>0$ such that the following hold for
all $0<\epsilon\le \epsilon_5$. Let $N$ and $N'$ be $\epsilon$-necks
centered at $x$ and $x'$, respectively, in a Riemannian manifold
$X$:
\begin{enumerate}
\item[(1)] If $N\cap N'\not=\emptyset$ then
$1-\alpha<R(x)/R(x')<1+\alpha$. In particular, denoting the scales
of $N$ and $N'$ by $h$ and $h'$ we have
$$1-\alpha<\frac{h}{h'}<1+\alpha.$$
\item[(2)] Suppose $y\in N\cap N'$ and $S$ and $S'$ are the two-spheres in the $\epsilon$-neck
structures on $N$ and $N'$, respectively, passing through $y$. Then
the angle between $TS_y$ and $TS'_y$ is less than $\alpha$.
\item[(3)] Suppose that $y\in N\cap N'$. Denote by $\partial/\partial s_N$ and $\partial/\partial s_{N'}$
the tangent vectors in the $\epsilon$-neck structures of $N$ and
$N'$, respectively. Then at the point $y$, either
$$|R(x)^{1/2}(\partial /\partial s_N)-R(x')^{1/2}\partial/\partial s_{N'}|<\alpha$$
or
$$|R(x)^{1/2}(\partial /\partial s_N)+R(x')^{1/2}\partial/\partial s_{N'}|<\alpha.$$
 \item[(4)] Suppose that
one of the two-spheres $S'$ of the $\epsilon$-neck structure on $N'$
is completely contained in $N$. Then $S'$ is a section of the
projection mapping on the first factor
$$p_1\colon S^2\times(-\epsilon^{-1},\epsilon^{-1})\to S^2.$$
In particular, $S'$ is isotopic in $N$ to any one of the two-spheres
of the $\epsilon$-neck structure on $N$ by an isotopy that moves all
points in the interval directions.
\item[(5)] If $N\cap N'$ contains a point $y$ with $(-0.9)\epsilon^{-1}\le s_N(y)\le
(0.9)\epsilon^{-1}$, then there is a point $y'\in N\cap N'$ such
that
$$-(0.96)\epsilon^{-1}\le s_N(y')\le (0.96)\epsilon^{-1}$$
$$-(0.96)\epsilon^{-1}\le s_{N'}(y')\le (0.96)\epsilon^{-1}.$$
The two-sphere $S(y')$ in the neck structure on $N$ through $y'$ is contained
in $N'$ and the two-sphere $S'(y')$ in the neck structure on $N'$ through $y'$
is contained in $N$. Furthermore, $S(y')$ and $S'(y')$ are isotopic in $N\cap
N'$. Lastly, $N\cap N'$ is diffeomorphic to $S^2\times (0,1)$ under a
diffeomorphism mapping $S(y)$ to $S^2\times \{1/2\}$, see {\sc
Fig.}~\ref{fig:overlap}.
\end{enumerate}
\end{prop}

\begin{proof}
Fix $0<\epsilon_5(\alpha)\le{\rm
min}(\epsilon_1(\alpha_1),\epsilon_2(\alpha/3),\epsilon_3(\alpha),\alpha/3)$
sufficiently small so that Corollary~\ref{s2isotopic} holds. The
first two items are then immediate from Lemma~\ref{neckcurv}. The
third statement is immediate from Lemma~\ref{directions}, and the
fourth statement from Corollary~\ref{s2isotopic}. Let us consider
the last statement. Let $y\in N\cap N'$ have $-(0.9)\epsilon^{-1}\le
s_N(y)\le (0.9)\epsilon^{-1}$. By reversing the $s$-directions of
$N$ and/or $N'$ if necessary, we can assume that $0\le s_N(y)\le
(0.9)\epsilon^{-1}$ and that $\partial_{s_N}$ and
$\partial_{s_{N'}}$ almost agree at $y$. If $-(0.96)\epsilon^{-1}\le
s_{N'}(y)\le (0.96)\epsilon^{-1}$, we set $y'=y$. Suppose that
$s_{N'}(y)>(0.96)\epsilon^{-1}$. We move along the straight line
though $y$ in the neck structure on $N$ in the negative direction to
a point $y'$ with $(0.96)\epsilon^{-1}=s_{N'}(y')$ According to Item
3 of this result we have $-(0.96)\epsilon^{-1}\le s_{N}(x')\le
(0.96)\epsilon^{-1}$. There is a similar construction when
$s_{N'}(y)<-(0.96)\epsilon^{-1}$. In all cases this allows us to
find $y'$ such that both the following hold:
$$-(0.96)\epsilon^{-1}\le s_N(y')\le (0.96)\epsilon^{-1}$$
$$-(0.96)\epsilon^{-1}\le s_{N'}(y')\le (0.96)\epsilon^{-1}.$$

Let $y'$ be any point satisfying both these inequalities. According
to Lemma~\ref{directions} and Part (1) of this result, the diameter
of $S(y')$ is at most $2\pi h$, where $h$  is the scale of $N$ and
$N'$. Since $\epsilon^{-1}\ge 200$, it follows from
Corollary~\ref{neckgeo} that $S(y')$ is contained in $N'$.
Symmetrically $S'(y')$ is contained in $N$.

Now consider the intersection of any straight line in the neck
structure on $N$ with $N'$. According to Part (3), this intersection
is connected. Thus, $N\cap N'$ is a union of open arcs in the
$s_N$-directions thought the points of $S(y')$. These arcs can be
used to define a diffeomorphism from  $N\cap N'$  to $S^2\times
(0,1)$ sending $S(y')$ to $S^2\times \{1/2\}$. Also, we have the
straight line isotopy from $S'(y')$ to $S(y')$ contained in $N\cap
N'$.
\end{proof}

\section{Regions covered by $\epsilon$-necks and
$(C,\epsilon)$-caps}

Here we fix $0<\epsilon\le 1/200$ sufficiently small so that all the
results in the previous two sections hold with $\alpha=10^{-2}$.

\subsection{Chains of $\epsilon$-necks}

\begin{defn}
Let $(X,g)$ be a Riemannian manifold. By a {\em finite chain} of
$\epsilon$-necks in $(X,g)$, we mean a sequence $N_a,\ldots,N_{b}$,
of $\epsilon$-necks in $(X,g)$ such that: \begin{enumerate} \item
for all $i,\ a\le i<b$, the intersection $N_i\cap N_{i+1}$ contains
the positive-most quarter of $N_i$ and the negative-most quarter of
$N_{i+1}$ and is contained in the positive-most three-quarters of
$N_i$ and the negative-most three-quarters of $N_{i+1}$, and
\item  for all $i,\ a<i\le b$, $N_i$ is disjoint from the negative end of
$N_a$.
\end{enumerate}
By an {\em infinite chain} of $\epsilon$-necks in $X$ we mean a
collection $\{N_i\}_{i\in I}$ for some interval $I\subset \Zee$,
infinite in at least one direction, so that for each finite
subinterval $J$ of $I$ the subset of $\{N_i\}_{i\in J}$ is a chain
of $\epsilon$-necks.
\end{defn}

Notice that in an $\epsilon$-chain $N_i\cap N_j=\emptyset$ if
$|i-j|\ge 5$.

\begin{lem}\label{S2timesI}
The union $U$ of the $N_i$ in a finite or infinite chain of
$\epsilon$-necks is diffeomorphic to $S^2\times (0,1)$. In
particular, it is an $\epsilon$-tube.
\end{lem}

\begin{proof}
Let us first prove the result for finite chains. The proof that $U$
is diffeomorphic to $S^2\times (0,1)$ is by induction on $b-a+1$. If
$b=a$, then the result is clear. Suppose that we know the result for
chains of smaller cardinality. Then $N_a\cup\cdots\cup N_{b-1}$ is
diffeomorphic to $S^2\times (0,1)$. Hence by Part (5) of
Proposition~\ref{S2intersection}, $U$ is the union of two manifolds
each diffeomorphic to $S^2\times (0,1)$ meeting in an open subset
diffeomorphic to $S^2\times (0,1)$. Furthermore, by the same result
 in the intersection there is a two-sphere isotopic to each of the
two-sphere factors from the two pieces.  It now follows easily that
the union is diffeomorphic to $S^2\times (0,1)$. Now consider an
infinite chain. It is an increasing union of finite chains each
diffeomorphic to $S^2\times (0,1)$ and with the two-spheres of one
isotopic to the two-spheres of any larger one. It is then immediate
that the union is diffeomorphic to $S^2\times (0,1)$.
\end{proof}

Notice that the frontier of the union of the necks in a finite
chain, $U=\cup_{a\le i\le b}N_i$, in $M$ is equal to the frontier of
the positive end of $N_b$ union the frontier of the negative end of
$N_a$. Thus, we have:

\begin{cor}\label{chainfrontier}
Let $\{N_a,\ldots,N_b\}$ be a chain of $\epsilon$-necks. If a
connected set $Y$ meets both $U=\cup_{a\le i\le b}N_i$ and its
complement, then $Y$ either contains  points of the frontier of the
negative end $N_a$ or of the positive  end of $N_b$.
\end{cor}

The next result shows there is no frontier at an infinite end.

\begin{lem}\label{inffrontier}
Suppose that $\{N_0,\cdots\}$ is an infinite chain of
$\epsilon$-necks in $M$. Then the frontier of $U=\cup_{i=0}^\infty
N_i$ is the frontier of the negative end of $N_0$.
\end{lem}

\begin{proof}
Suppose that $x$ is a point of the frontier of $U$. Let $x_i\in U$
be a sequence converging to $x$. If the $x_i$ were contained in a
finite union of the $N_k$, say $N_0\cup\cdots\cup N_\ell$, then $x$
would be in the closure of this union and hence by the previous
comment would be either be in the frontier of the negative end of
$N_0$ or the frontier of the positive end of $N_\ell$. But the
frontier of the positive end of $N_\ell$ is contained in
$N_{\ell+1}$ and hence contains no points of the frontier of $U$.
Thus, in this case $x$ is a point of the frontier of the negative
end of $N_0$. If $\{x_i\}$ is not contained in any finite union,
then after passing to a subsequence, we can suppose that $x_i\in
N_{k(i)}$ where $k(i)$ is an increasing sequence tending to
infinity. Clearly $R(x_i)$ converges to $R(x)<\infty$. Hence, there
is a uniform lower bound to the scales of the $N_{k(i)}$. For all
$i$ sufficiently large $x_i\not\in N_0$. Thus, for such $i$ any path
from $x_i$ to $x$ must traverse either $N_0$ or $N_{k(j)}$ for all
$j\ge i+5$. The length of such a path is at least the minimum of the
width of $N_0$ and the width of $N_{k(j)}$ for some $j$ sufficiently
large. But we have just seen that there is a positive lower bound to
the scales of the $N_{k(j)}$ independent of $j$, and hence by
Corollary~\ref{neckgeo} there is a positive lower bound, independent
of $j$, to the widths of the $N_{k(j)}$. This shows that there is a
positive lower bound, independent of $i$, to the distance from $x_i$
to $x$ .This is impossible since $x_i$ converges to $x$.
\end{proof}

 In fact, there is a geometric version of Lemma~\ref{S2timesI}.

\begin{lem}
There is $\epsilon_0>0$ such that  the following holds for all
$0<\epsilon\le \epsilon_0$. Suppose that $\{N_j\}_{j\in J}$ is a
chain of $\epsilon$-necks in a Riemannian manifold $M$. Let
$U=\cup_{j\in J}N_j$. Then there exist an interval $I$ and a smooth
map $p\colon U\to I$ such that every fiber of $p$ is a two-sphere,
and if $y$ is in the middle $7/8$'s of $N_j$ then the fiber
$p^{-1}(p(y))$ makes a small angle at every point with the family of
two-spheres in the $\epsilon$-neck $N_j$.
\end{lem}

\begin{proof}
Since according to Lemma~\ref{neckcurv} the two-spheres for $N_j$
and $N_{j+1}$ almost line up, it is an easy matter to interpolate
between the projection maps to the interval to construct a fibration
of $U$ by two-spheres with the given property. The interval $I$ is
simply the base space of this fibration.
\end{proof}

A finite or infinite chain $\{N_j\}_{j\in J}$ of $\epsilon$-necks is
{\em balanced} provided that for every $j\in J$, not the largest
element of $J$, we have
\begin{equation}\label{distineq}
(0.99)R(x_j)^{-1/2}\epsilon^{-1}\le d(x_j,x_{j+1})\le
(1.01)R(x_j)^{-1/2}\epsilon^{-1},
\end{equation}
where, for each $j$, $x_j$ is the central point of $N_j$.

Notice that in a balanced chain $N_j\cap N_{j'}=\emptyset$ if
$|j-j'|\ge 3$.

\begin{lem}\label{firstchain}
There exists $\epsilon_0>0$ such that for all $0<\epsilon\le
\epsilon_0$ the following is true.  Suppose that $N$ and $N'$ are
$\epsilon$-necks centered at $x$ and $x'$, respectively, in a
Riemannian manifold $M$. Suppose that $x'$ is not contained in $N$
but is contained in the closure of $N$ in $M$. Suppose also that the
two-spheres of the neck structure on $N$ and $N'$ separate $M$.
Then, possibly after reversing the $\epsilon$-neck structures on $N$
and/or $N'$, the pair $\{N,N'\}$ forms a balanced chain.
\end{lem}

\begin{proof} By Corollary~\ref{neckgeo},
Inequality~(\ref{distineq}) holds for $d(x,x')$. Once we have this
inequality, it follows immediately from the same corollary that,
possible after reversing, the $s$-directions $\{N,N'\}$ makes a
balanced chain of $\epsilon$-necks. (It is not possible for the
positive end of $N_b$ to meet $N_a$ for this would allow us to
create a loop meeting the central two-sphere of $N_b$ transversely
in a single point, so that this two-sphere would not separate $M$.)
\end{proof}

\begin{lem}\label{addaneck}
There exists $\epsilon_0>0$ such that for all $0<\epsilon\le
\epsilon_0$ the following is true.  Suppose that
$\{N_a,\ldots,N_b\}$ is a balanced chain in a Riemannian manifold
$M$ with $U=\cup_{i=a}^bN_i$. Suppose that the two-spheres of the
neck structure of $N_a$ separate $M$.
 Suppose that $x$
is a point of the frontier of $U$ contained in the closure of the
plus end of $N_b$ that is also the center of an $\epsilon$-neck $N$.
Then possibly after reversing the direction of $N$, we have that
$\{N_a,\ldots,N_b,N\}$ is a balanced chain. Similarly, if $x$ is in
the closure of the minus end of $N_a$, then (again after possibly
reversing the direction of $N$) we have that $\{N,N_a,\ldots,N_b\}$
is a balanced $\epsilon$-chain.
\end{lem}

\begin{proof}
 The two cases are
symmetric; we consider only the first. Since $x$ is contained in the
closure of $N_b$, clearly $N_b\cap N\not=\emptyset$. Also, clearly,
provided that $\epsilon>0$ is sufficiently small, $d(x_b,x)$
satisfies Inequality~(\ref{distineq}) so that Lemma~\ref{firstchain}
the pair $\{N_b,N\}$ forms an $\epsilon$-chain, and hence a balanced
$\epsilon$-chain. It is not possible for $N$ to meet the negative
end of $N_a$ since the central two-sphere of $N_a$ separates $M$.
Hence $\{N_a,\ldots,N_b,N\}$ is a balanced chain of
$\epsilon$-necks.
\end{proof}

\begin{prop}\label{epschains}
There exists $\epsilon_0>0$ such that for all $0<\epsilon\le
\epsilon_0$ the following is true. Let $X$ be a connected subset of
a Riemannian manifold $M$ with the property that every point $x\in
X$ is the center of an $\epsilon$-neck $N(x)$ in $M$. Suppose that
the central two-spheres of these necks do not separate $M$. Then
there is a subset $\{x_i\}$ of $X$ such that the necks $N(x_i)$
(possibly after reversing their $s$-directions) form a balanced
chain of $\epsilon$-necks $\{N(x_i)\}$ whose union $U$ contains $X$.
The union $U$ is diffeomorphic to $S^2\times (0,1)$. It is an
$\epsilon$-tube.
\end{prop}

\begin{proof}
According to Lemma~\ref{addaneck} for $\epsilon>0$ sufficiently
small the following holds. Suppose that we have a balanced chain of
$\epsilon$-necks $N_a\ldots,N_{b}$, with $N_i$ centered at $x_i\in
X$, whose union $U$  does not contain $X$.
 Then one of the following holds:
\begin{enumerate}
\item[(1)] It is possible to find an $\epsilon$-neck $N_{b+1}$ centered
at a point of the intersection of $X$ with the closure of the
positive end of $N_{b}$ so that $N_a,\ldots,N_{b+1}$ is a balanced
$\epsilon$-chain.
\item[(2)] It is possible to find an $\epsilon$-neck $N_{a-1}$ centered at a
point of the intersection of $X$ with the closure of the negative
end of $N_a$ so that $N_{a-1},N_a,\ldots,N_b$ is a balanced
$\epsilon$-chain.
\end{enumerate}

Now assume that there is no finite balanced chain of
$\epsilon$-necks $N(x_i)$ containing $X$. Then we can repeatedly
lengthen a balanced chain of $\epsilon$-necks centered at points of
$X$ by adding necks at one end or the other. Suppose that we have a
half-infinite balanced chain $\{N_0,N_1,\ldots,\}$. By
Lemma~\ref{inffrontier} the frontier of this union is the frontier
of the negative end of $N_0$. Thus, if we can construct a balanced
chain which is infinite in both directions, then the union of the
necks in this chain is a component of $M$ and hence contains the
connected set $X$. If we can construct a balanced chain that is
infinite at one end but not the other that cannot be further
extended, then the connected set is disjoint from the frontier of
the negative end of the first neck in the chain and, as we have see
above, the `infinite' end of the chain has no frontier. Thus, $X$ is
disjoint from the frontier of $U$ in $M$ and hence is contained in
$U$. Thus, in all cases we construct a balanced chain of
$\epsilon$-necks containing $X$. By Lemma~\ref{S2timesI} the union
of the necks in this chain is diffeomorphic to $S^2\times (0,1)$ and
hence is an $\epsilon$-tube.
\end{proof}

\begin{lem}\label{S2times01}
The following holds for every $\epsilon>0$ sufficiently small. Let
$(M,g)$ be a connected Riemannian manifold. Suppose that every point
of $M$ is the center of an $\epsilon$-neck. Then either $M$ is
diffeomorphic to $S^2\times (0,1)$ and is an $\epsilon$-tube, or $M$
is diffeomorphic to an $S^2$-fibration over $S^1$.
\end{lem}

\begin{proof}
If the two-spheres of the $\epsilon$-necks do not separate $M$, then
it follows from the previous result that $M$ is an $\epsilon$-tube.
If one of the two-spheres does separate, then take the universal
covering $\widetilde M$ of $M$.  Every point of $\widetilde M$ is
the center of an $\epsilon$-neck (lifting an $\epsilon$-neck in $M$)
and the two-spheres of these necks separate $\widetilde M$. Thus the
first case applies, showing that $\widetilde M$ is diffeomorphic to
$S^2\times (0,1)$. Every point is the center of an $\epsilon$-neck
that is disjoint from all its non-trivial translates under the
fundamental group. This means that the quotient is fibered by
$S^2$'s over $S^1$, and the fibers of this fibration are isotopic to
the central two-spheres of the $\epsilon$-necks.
\end{proof}

\section{Subsets of the union of cores of $(C,\epsilon)$-caps
and $\epsilon$-necks.} In this section we fix $0<\epsilon\le 1/200$
so that all the results of this section hold with $\alpha=0.01$.

\begin{figure}[ht]
  \centerline{\epsfbox{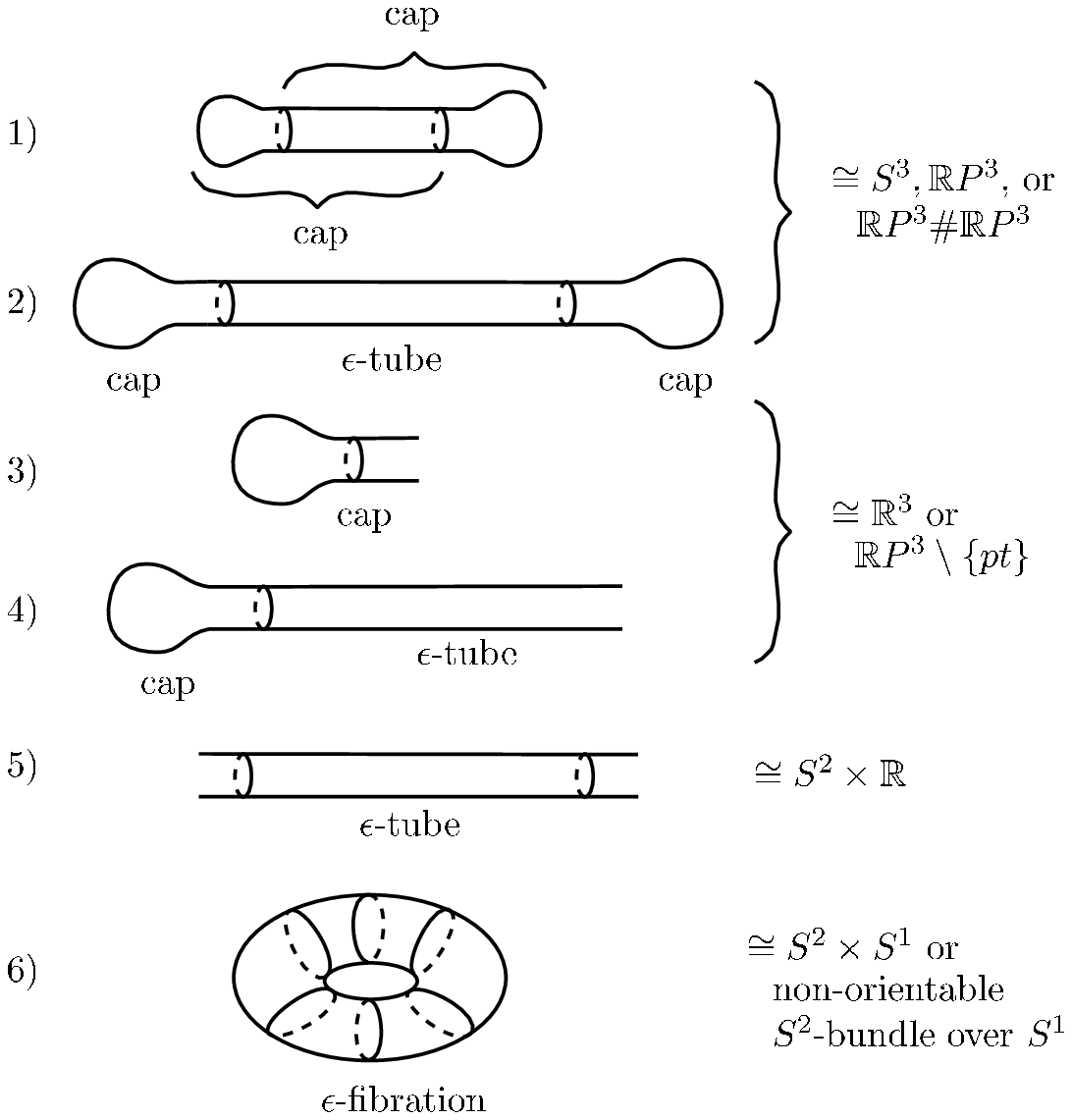}}
  \caption{Components covered by $\epsilon$-necks and
  $\epsilon$-caps.}\label{fig:comp}
\end{figure}

\begin{prop}\label{Xcontainedin}
For any $C<\infty$ the following holds. Suppose that $X$ is a
connected subset of a Riemannian three-manifold $(M,g)$. Suppose
that every point of $X$ is either the center of an $\epsilon$-neck
or is contained in the core of a $(C,\epsilon)$-cap. Then one of the
following holds:
\begin{enumerate}
\item[(1)] $X$ is contained in a component of $M$ that is the union of
two $(C,\epsilon)$-caps. This component is diffeomorphic to $S^3$,
$\Ar P^3$ or $\Ar P^3\#\Ar P^3$.
\item[(2)] $X$ is contained in a component of $M$ that is a double $C$-capped $\epsilon$-tube.
 This component is diffeomorphic to $S^3$, $\Ar P^3$ or $\Ar
P^3\#\Ar P^3$.
\item[(3)] $X$ is contained in a single $(C,\epsilon)$-cap.
\item[(4)] $X$ is contained in a $C$-capped $\epsilon$-tube.
\item[(5)] $X$ is contained in an $\epsilon$-tube.
\item[(6)] $X$ is contained in a component of $M$ that is an
$\epsilon$-fibration, which itself is a union of $\epsilon$-necks.
\end{enumerate}
(See {\sc Fig.}~\ref{fig:comp}.)
\end{prop}

\begin{proof}
We divide the proof into two cases: Case I: There is a point of $X$
contained in the core of a $(C,\epsilon)$-cap. Case II: Every point
of $X$ is the center of an $\epsilon$-neck.

\noindent{\bf Case I:} We begin the study of this case with a claim.

\begin{claim}
It is not possible to have an infinite chain of $(C,\epsilon)$-caps
$C_0\subset C_1\subset \cdots$ in $M$ with the property that for
each $i\ge 1$, the closure of the core of $C_i$ contains a point of
the frontier of $C_{i-1}$
\end{claim}

\begin{proof}
We argue by contradiction. Suppose there is such an infinite chain.
Fix a point $x_0\in C_0$ and let $Q_0=R(x_0)$. For each $i\ge 1$ let
$x_i$ be a point in the frontier of $C_{i-1}$ that is contained in
the closure of the core of $C_i$. For each $i$ let $N_i$ be the
$\epsilon$-neck in $C_{i}$ that is the complement of the closure of
its core. We orient the $s_{N_{i}}$-direction so that the core of
$C_{i}$ lies off the negative end of $N_{i}$. Let $S'_i$ be the
boundary of the core of $C_i$. It is the central two-sphere of an
$\epsilon$-neck $N'_i$ in $C_i$. We orient the $s$-direction of
$N'_i$ so that the non-compact end of $C_i$ lies off the positive
end of $N'_i$. We denote by $h_{i-1}$ the scale of $N_{i-1}$ and by
$h_i'$ the scale of $N'_i$. By Lemma~\ref{neckcurv} the ratio
$h_{i-1}/h'_i$ is between $0.99$ and $1.01$. Suppose that $S'_i$ is
disjoint from $C_{i-1}$. Then one of the complementary components of
$S'_i$ in $M$ contains $C_{i-i}$, and of course, one of the
complementary components of $S'_i$ is the core of $C_i$. These
complementary components must be the same, for otherwise $C_{i-1}$
would be disjoint from the core of $C_i$ and hence the intersection
of $C_{i-1}$ and $C_i$ would be contained in $N_i$. This cannot
happen since  $C_{i-1}$ is contained in $C_i$. Thus, if $S'_i$ is
disjoint from $C_{i-1}$, then the core of $C_i$ contains $C_{i-1}$.
 This means that the distance from $x_0$ to the
complement of $C_i$ is greater than the distance of $x_0$ to the
complement of $C_{i-1}$ by an amount equal to the width of $N_i$.
Since the scale of $N_i$ is at least $C^{-1/2}R(x_0)^{-1/2}$ (see
(5) of Definition~\ref{epscap}), it follows from
Corollary~\ref{neckgeo} that this width is at least
$2(0.99)C^{-1/2}R(x_0)^{-1/2}\epsilon^{-1}$.

Next suppose that $S'_i$ is contained in $C_{i-1}$. Then one of the
complementary components $A$ of $S'_i$ in $M$ has closure contained
in $C_{i-1}$. This component cannot be the core of $C_i$ since the
closure of the core of $C_i$ contains a point of the frontier of
$C_{i-1}$ in $M$. Thus, $A$ contains $N_i$. Of course, $A\not= N_i$
since the frontier of $A$ in $M$ is $S'_i$ whereas $N_i$ has two
components to its frontier in $M$. This means that $C_i$ does not
contain $A$, which is a contradiction since $C_i$ contains $C_{i-1}$
and $A\subset C_{i-1}$.

Lastly, we suppose that $S'_i$ is neither contained in $C_{i-1}$ nor
in its complement. Then $S'_i$ must meet $N_{i-1}$.  According to
Proposition~\ref{S2intersection} the $s$-directions in $N_{i-1}$ and
$N'_i$ either almost agree or are almost opposite. Let $x\in
S'_i\cap \partial N_{i-1}$ so that $s_{N'_i}(x)=0$. Move from $x$
along the $s_{N'_i}$-direction that moves into $N_{i-1}$ to a point
$x'$ with $|s_{N_i}(x')|=(0.05)\epsilon^{-1}$. According to
Proposition~\ref{S2intersection}
$(0.94)\epsilon^{-1}<s_{N_{i-1}}(x')<(0.96)\epsilon^{-1}$.
 Let $S'(x')$ be the two-sphere in the neck structure for $N'_i$
 through this point. According to Proposition~\ref{S2intersection},
 $S'(x')\subset N_{i-1}$, and $S'(x')$ is isotopic in $N_{i-1}$ to its
 central two-sphere. One of the complementary components of $S'(x')$
 in $C_i$, let us call it $A'$, is diffeomorphic to $S^2\times (0,1)$. Also, one of the
 complementary components $A$ of $S'(x')$ in $M$ contains the core of
 $C_{i-1}$. As before, since $C_{i-1}\subset C_i$, the complementary
 component $A$ cannot meet $C_i$ in $A'$. This means that
 the $s_{N_{i-1}}$- and $s_{N'_i}$-directions almost line up along
 $S'(x')$. This means that
 $S'(x')=s_{N'_i}^{-1}(-(0.05)\epsilon^{-1})$. Since the diameter of
 $S'(x')$ is less than $2\pi h_{i-1}$, and since
  $s_{N_{i-1}}(x')\ge (0.94)\epsilon^{-1}$, it follows
 that $S'(x')\subset
 S_{N_{i-1}}^{-1}((0.9\epsilon^{-1},\epsilon^{-1}))$.
Since the distance from $S'_i$ to the central two-sphere is at least
$(0.99)\epsilon^{-1}h'_i$,  It follows from Corollary~\ref{neckgeo}
that the central two-sphere of $N_i$ is disjoint from $C_{i-1}$ and
lies off the positive end of $N_{i-1}$. This implies that the
distance from $x_0$ to the complement of $C_i$ is greater than the
distance from $x_0$ to the complement of $C_{i-1}$ by an amount
bounded below by the distance from the central two-sphere of $N_i$
to its positive end. According to Corollary~\ref{neckgeo} this
distance is at least $(0.99)\epsilon^{-1}h_i$, where $h_i$ is the
scale of $N_i$. But we know that  $h_i\ge C^{-1/2}R(x_0)^{-1/2}$.

Thus, all cases either lead to a contradiction or to the conclusion
that the distance from $x_0$ to the complement of $C_i$ is at least
a fixed positive amount (independent of $i$) larger than the
distance from $x_0$ to the complement of $C_{i-1}$.
 Since the diameter of any $(C,\epsilon)$-cap is
uniformly bounded, this contradicts the existence of an infinite
chain $C_0\subset C_1\subset \cdots$ contrary to the claim. This
completes the proof of the claim.
\end{proof}

Now let us turn to the proof of the proposition. We suppose first
that there is a point $x_0\in X$ that is contained in the core of a
$(C,\epsilon)$-cap. Applying the previous claim, we can find a
$(C,\epsilon)$-cap $C_0$ containing $x_0$ with the property that no
point of $X$ contained in the  frontier of $C_0$ is contained in the
closure of the core of a $(C,\epsilon)$-cap $C_1$ that contains
$C_0$.

There are three possibilities to examine:
\begin{enumerate}
\item[(i)] $X$ is disjoint from the frontier of $C_0$.
\item[(ii)] $X$ meets the frontier of $C_0$ but every point of this
intersection is the center of an $\epsilon$-neck.
\item[(iii)] There is a point of the intersection of $X$
with the frontier of $C_0$ that is contained in the core of
$(C,\epsilon)$-cap.
\end{enumerate}

In the first case, since $X$ is connected, it is contained in $C_0$.
In the second case we let $N_1$ be an $\epsilon$-neck centered at a
point of the intersection of $X$ with the frontier of $C_0$, and we
replace $C_0$ by $C_0\cup N_1$ and repeat the argument at the
frontier of $C_0\cup N_1$. We continue in this way creating $C_0$
union a balanced chain of $\epsilon$-necks $C_0\cup N_1\cup N_2\cup
\cdots\cup N_k$. At each step it is possible that either there is no
point of the frontier containing a point of $X$, in which case the
union, which is a $C$-capped $\epsilon$-tube, contains $X$. Another
possibility is that we can repeat the process forever creating a
$C$-capped infinite $\epsilon$-tube. By Lemma~\ref{inffrontier} this
union is a component of $M$ and hence contains $X$.

We have shown that one of following holds:
\begin{enumerate}
\item[(a)] There is a $(C,\epsilon)$-cap that contains $X$.
\item[(b)] There is a finite or infinite $C$-capped $\epsilon$-tube that
contains $X$.
\item[(c)] There is a $(C,\epsilon)$-cap or a finite $C$-capped
$\epsilon$-tube $\widetilde C$ containing a point of $X$ and there
is a point of the intersection of $X$ with the frontier of
$\widetilde C$ that is contained in the core of a
$(C,\epsilon)$-cap.
\end{enumerate}

In the first two cases we have established the proposition. Let us
examine the third case in more detail. Let $N_0\subset C_0$ be the
$\epsilon$-neck that is the complement of the closure of the core of
$C_0$.
 First notice that by
Lemma~\ref{S2times01} the union $N_0\cup N_1\cup\cdots\cup N_k$ is
diffeomorphic to $S^2\times (0,1)$, with the two-spheres coming from
the $\epsilon$-neck structure of each $N_i$ being isotopic to the
two-sphere factor in this product structure. It follows immediately
that $\widetilde C$ is diffeomorphic to $C_0$. Let $C'$ be a
$(C,\epsilon)$-cap whose core contains a point of the intersection
of $X$ with the frontier of $\widetilde C$. We use the terminology
`the core of $\widetilde C$' to mean $\widetilde C\setminus N_k$.
Notice that if $k=0$, this is exactly the core of $C_0$. To complete
the proof of the result we must show that the following hold:

\begin{claim}\label{cccomponent}
If $C'$ is a $(C,\epsilon)$-cap whose core contains a point of the
frontier of $\widetilde C$, then $\widetilde C\cup C'$ is a
component of $M$ containing $X$.
\end{claim}

\begin{proof}
We suppose that $\widetilde C$ is the union of $C_0$ and a balanced
chain $N_0,\ldots,N_k$ of $\epsilon$-necks. We orient this chain so
that $C_0$ lies off the negative end of each of the $N_i$.
 Let $S'$ be the boundary of the core of $C'$ and let $N'$ be an
$\epsilon$-neck contained in $C'$ whose central two-sphere is $S'$.
We orient the direction $s_{N'}$ so that the positive direction
points away from the core of  $C'$. The first step in proving this
claim is to establish the following.

\begin{claim}\label{sinside}
Suppose that there is a two-sphere $\Sigma\subset N'$ contained in
the closure of the positive half of $N'$ and also contained in
$\widetilde C$. Suppose that $\Sigma$ is isotopic in $N'$ to the
central two-sphere $S'$ of $N'$. Then $\widetilde C\cup C'$ is a
component of $M$, a component containing $X$.
\end{claim}

\begin{proof}
$\Sigma$ separates $\widetilde C$ into two components: $A$, which
has compact closure in $\widetilde C$, and $B$, containing the end
of $\widetilde C$. The two-sphere $\Sigma$ also divides $C'$ into
two components. Since $\Sigma$ is isotopic in $N'$ to $S'$, the
complementary component $A'$ of $\Sigma$ in $C'$ with compact
closure contains the closure of the core of $C'$. Of course, the
frontier of $A$ in $M$ and the frontier of $A'$ in $M$ are both
equal to $\Sigma$. If $A=A'$, then the closure of the core of $C'$
is contained in the closure of $A$ and hence is contained in
$\widetilde C$, contradicting our assumption that $C'$ contains a
point of the frontier of $\widetilde C$. Thus, $A$ and $A'$ lie on
opposite sides of their common frontier. This means that $\overline
A\cup \overline A'$ is a component of $M$. Clearly, this component
is also equal to $\widetilde C\cup C'$. Since $X$ is connected and
this component contains a point $x_0$ of $X$, it contains $X$. This
completes the proof of Claim~\ref{sinside}.
 \end{proof}

Now we return to the proof of Claim~\ref{cccomponent}. We consider
three cases.

\noindent{\bf First Subcase: $S'\subset \widetilde C$.} In this case
we apply Claim~\ref{sinside} to see that $\widetilde C\cup C'$ is a
component of $M$ containing $X$.

 \noindent {\bf Second Subcase: $S'$ is disjoint from
$\widetilde C$.} Let $A$ be the complementary component of $S'$ in
$M$ containing $\widetilde C$. The intersection of $A$ with $C'$ is
either the core of $C'$ or is a submanifold of $C'$ diffeomorphic to
$S^2\times (0,1)$. The first case is not possible since it would
imply that the core of $C'$ contains $\widetilde C$ and hence
contains $C_0$, contrary to the way we chose $C_0$. Thus, the core
of $C'$ and the the complementary component $A$ containing
$\widetilde C$ both have $S'$ as their frontier and they lie on
opposite sides of $S'$. Since the closure of the core of $C'$
contains a point of the frontier of $\widetilde C$, it must be the
case that $S'$ also contains a point of this frontier. By
Proposition~\ref{S2intersection},  the neck $N'\subset C'$ meets
$N_{k}$ and there is a two-sphere $\Sigma\subset N'\cap N_{k}$
isotopic in $N'$ to $S'$ and isotopic in $N_{k}$ to the central
two-sphere of $N_{k}$. Because $N_{k}\subset \widetilde C$ and
$\widetilde C$ is disjoint from the core of $C'$, we see that
$\Sigma$ is contained in the positive half of $N'$. Applying
Claim~\ref{sinside} we see that $\widetilde C\cup C'$ is a component
of $M$ containing $X$.

\noindent{\bf Third Subcase: $S'\cap \widetilde C\not=\emptyset$ and
$S'\not\subset \widetilde C$.} Clearly, in this case $S'$ contains a
point of the frontier of $\widetilde C$ in $M$, i.e., a point of the
frontier of the positive end of $N_k$ in $M$. Since $N_k\cap
N'\not=\emptyset$, by Lemma~\ref{neckcurv} the scales of $N_k$ and
$N'$ are within $1\pm 0.01$ of each other, and hence the diameter of
$S'$ is at most $2\pi$ times the scale of $N_k$. Since the central
two-sphere $S'$ of $N'$ contains a point in the frontier of the
positive end of $N_k$, it follows from Lemma~\ref{neckgeo} that $S'$
is contained on the positive side  of the central two-sphere of
$N_k$ and that the frontier of the positive end of $N_k$ is
contained in $N'$. By Proposition~\ref{S2intersection} there is a
two-sphere $\Sigma$ in the neck structure for $N'$ that is contained
in $N_k$ and is isotopic in $N_k$ to the central two-sphere from
that neck structure. Let $A$ be the complementary component of
$\Sigma$ in $M$ that contains $\widetilde C\setminus N_k$.
 If the complementary component of $\Sigma$ that contains
$C'\setminus N'$ is not $A$, then $\widetilde C\cup C'$ is a
component of $M$ containing $X$. Suppose that $A$ is also the
complementary component of $\Sigma$ in $M$ that contains
$C'\setminus N'$. Of course, $A$ is contained in the core of $C'$.
If $k\ge 1$, we see that $A$ and hence the core of $C'$ contains
$\widetilde C\setminus N_k$, which in turn contains the core of
$C_0$. This contradicts our choice of $C_0$. If $k=0$, then
$C_0=A\cup (N_0\cap (M\setminus A))$. Of course, $A\subset C'$.
Also, the frontier of $N_0\cap(M\setminus A)$ in $M$ is the union of
$A$ and the frontier of the positive end of $N_0$ in $M$. But we
have already established that the frontier of the positive end of
$N_0$ in $M$ is contained in $N'$. Since $A\subset C'$, it follows
that all of $C_0$ is contained in $C'$. On the other hand, there is
a point of the frontier of $C_0$ contained in the closure of the
core of $C'$. This then contradicts our choice of $C_0$.

 This
completes the analysis of all the cases and hence completes the
proof of Claim~\ref{cccomponent}.
\end{proof}

The last thing to do in this case in order to prove the proposition
in Case I is to show that $\widetilde C\cup C'$ is diffeomorphic to
$S^3$, $\Ar P^3$, or $\Ar P^3\#\Ar P^3$. The reason for this is that
$\widetilde C$ is diffeomorphic to $C_0$;  hence $\widetilde C$
either is diffeomorphic to an open three-ball or to a punctured $\Ar
P^3$. Thus, the frontier of $C'$ in $\widetilde C$ is a two-sphere
that bounds either a compact three-ball or the complement of an open
three-ball in $\Ar P^3$. Since $C'$ itself is diffeomorphic either
to a three-ball or to a punctured $\Ar P^3$, the result follows.

\noindent{\bf Case II:} Suppose that every point of $X$ is the
center of an $\epsilon$-neck. Then if the two-spheres of these necks
separate $M$, it follows from Proposition~\ref{epschains} that $X$
is contained in an $\epsilon$-tube in $M$.

It remains to consider the case when the two-spheres of these necks
do not separate $M$. As in the case when the two-spheres separate,
we begin building a balanced chain $\epsilon$-necks with each neck
in the chain centered at a point of $X$. Either this construction
terminates after a finite number of steps in a finite
$\epsilon$-chain whose union contains $X$, or it can be continued
infinitely often creating an infinite $\epsilon$ chain containing
$X$ or at some finite stage (possibly after reversing the indexing
and the $s$-directions of the necks)  we have a balanced
$\epsilon$-chain $N_a\cup\cdots\cup N_{b-1}$ and a point of the
intersection of $X$ with the frontier of the positive end of
$N_{b-1}$ that is the center of an $\epsilon$-neck $N_b$ with the
property that $N_b$ meets the negative end of $N_a$. Intuitively,
the chain wraps around on itself like a snake eating its tail. If
the intersection of $N_a\cap N_b$ contains a point $x$ with
$s_{N_a}(x)\ge -(0.9)\epsilon^{-1}$, then according to
Proposition~\ref{S2intersection} the intersection of $N_a\cap N_b$
is diffeomorphic to $S^2\times (0,1)$ and the two-sphere in this
product structure is isotopic in $N_a$ to the central two-sphere of
$N_a$ and is isotopic in $N_b$ to the central two-sphere of $N_b$.
In this case it is clear that $N_a\cup \cdots \cup N_b$ is a
component of $M$ that is an $\epsilon$-fibration.

We examine the possibility that the intersection $N_a\cap N_b$
contains some points in the negative end of $N_a$ but is  contained
in $s_{N_a}^{-1}((-\epsilon^{-1},-(0.9)\epsilon^{-1}))$. Set $A=
s_{N_a}^{-1}((-\epsilon^{-1},-(0.8)\epsilon^{-1}))$. Notice that
since $X$ is connected and $X$ contains a point in the frontier of
the positive end of $N_a$ (since we have added at least one neck at
this end), it follows that $X$ contains points in $s_{N_a}^{-1}(s)$
for all $s\in [0,\epsilon^{-1})$.
 If there are
no points of $X$ in  $A$, then we replace $N_a$ by an
$\epsilon$-neck $N_a'$ centered at a point of
$s_{N_a}^{-1}\left((0.15)\epsilon^{-1}\right)\cap X$. Clearly, by
Lemma~\ref{neckgeo} $N'_a$ contains
$s_{N_a}^{-1}(-(0.8)\epsilon{-1},\epsilon^{-1})$ and is disjoint
from $s_{N_a}^{-1}((-\epsilon^{-1},-(0.9)\epsilon^{-1})$, so that
$N_a',N_{a+1}, \ldots, N_b$ is a chain of $\epsilon$-necks
containing $X$. If there is a point of $X\cap A$, then we let
$N_{b+1}$ be a neck centered at this point. Clearly,
$N_a\cup\cdots\cup N_{b+1}$ is a component, $M_0$, of $M$ containing
$X$. The preimage in the universal covering of $M_0$ is a chain of
$\epsilon$-necks infinite in both directions. That is to say, the
universal covering of $M_0$ is an $\epsilon$-tube. Furthermore, each
point in the universal cover of $M_0$ is the center of an
$\epsilon$-neck that is disjoint from all its non-trivial covering
translates. Hence, the quotient $M_0$ is an $\epsilon$-fibration.

We have now completed the proof of Proposition~\ref{Xcontainedin}.
\end{proof}

As an immediate corollary we have:

\begin{prop}\label{epstopology}
For all $\epsilon>0$ sufficiently small the following holds. Suppose
that $(M,g)$ is a connected Riemannian manifold such that every
point is either contained in the core of a $(C,\epsilon)$-cap in $M$
or is the center of an $\epsilon$-neck in $M$. Then one of the
following holds:
\begin{enumerate}
\item[(1)] $M$ is diffeomorphic to $S^3$, $\Ar P^3$ or $\Ar P^3\#\Ar
P^3$,
and $M$ is either a double $C$-capped $\epsilon$-tube or is the
union of two $(C,\epsilon)$-caps.
\item[(2)] $M$ is diffeomorphic to  $\Ar^3$ or $\Ar P^3\setminus \{{\rm point}\}$,
and $M$ is either a $(C,\epsilon$-cap or a $C$-capped
$\epsilon$-tube.
\item[(3)] $M$ is diffeomorphic to $S^2\times \Ar$ and is an
$\epsilon$-tube.
\item[(4)] $M$ is diffeomorphic to an $S^2$-bundle over $S^1$ and is an $\epsilon$-fibration.
\end{enumerate}
\end{prop}

\bibliographystyle{plain}
\bibliography{bibliography}

\printindex

\end{document}